\documentclass[dvips,oneside,10pt]{memoir}
 
\DeclareSymbolFont{symbolsC}{U}{txsyc}{m}{n}
\SetSymbolFont{symbolsC}{bold}{U}{txsyc}{bx}{n}
\DeclareFontSubstitution{U}{txsyc}{m}{n}
\makeatletter
\def\re@DeclareMathSymbol#1#2#3#4{%
    \let#1=\undefined
    \DeclareMathSymbol{#1}{#2}{#3}{#4}}
\re@DeclareMathSymbol{\leftsquigarrow}{\mathrel}{symbolsC}{102}
\makeatother

\cftpagenumbersoff{part} 
\usepackage{xcolor}
\ifdraftdoc
 \usepackage[color]{showkeys}
 \colorlet{refkey}{rgb:blue,1}
 \colorlet{labelkey}{rgb:red,1}
 \changemarkstrue           
\else
 \changemarksfalse
\fi

\makeatletter
\newcommand{\corrected}[1]{%
  \@bsphack
  \ifchangemarks
    \v@rid{\small$\Uparrow$ #1}{\small$\Uparrow$ #1}%
  \fi
  \@esphack}
\makeatother

\vfuzz \hfuzz
\newcommand{\version}{2}

\usepackage{amsmath,amsthm,amssymb}

\usepackage{palatino,eulervm}

\usepackage[backrefs]{amsrefs}

\input{chapter-style-blue.tex}

\newcommand{\checked}{\text{\ \checkmark\ }}
\renewcommand{\checked}{}

\usepackage[final]{graphicx}
\usepackage{array}
\usepackage{pstricks}
\usepackage{pstcol}
\usepackage{fancybox}
\usepackage[all]{xy}
\usepackage{ulem}
\usepackage{here}
\usepackage{epsfig}
\usepackage{verbatim}
\usepackage{enumerate}
\usepackage{longtable}
\usepackage{supertabular}
\usepackage{booktabs}
\usepackage[overload]{textcase}  

  \newcommand{\TimeStampStart}{}
  \newcommand{\mytoday}{\footnote{\TimeStampStart}}

 \renewcommand{\mytoday}{}

\numberwithin{equation}{section}

\newtheorem{nada}{}[chapter]

\newtheorem{theorem}[nada]{Theorem}
\newtheorem{lemma}[nada]{Lemma}

\newtheorem{cor}[nada]{Corollary}
\newtheorem{prop}[nada]{Proposition}
\newtheorem{conj}[nada]{Conjecture}

\theoremstyle{definition}
\newtheorem{definition}[nada]{Definition}
\newtheorem{example}[nada]{Example}
\newtheorem{construction}[nada]{Construction}

\theoremstyle{remark}
\newtheorem{remark}[nada]{Remark}
\newtheorem{question}{Question}

\newcommand{\sltwo}{\text{SL}_2(\reals[x])}
\newcommand{\sltwop}{\text{SL}_2(\reals[x])_P}

\newcommand{\matint}{\looparrowleft}

\newcommand{\uhp}{\textsf{UHP}}
\newcommand{\uhpc}{\overline{\textsf{UHP}}}
\newcommand{\rhp}{\textsf{RHP}}

\newcommand{\fdb}{Fa\`{a} di Bruno}
\newcommand{\dersub}{derivative-substitution}

\newcommand{\multiaff}[1]{{\textsf{MA}}_{#1}}

  \newcommand{\grace}[1]{\textsf{Grace}_{#1}}
  \newcommand{\homog}[1]{\textsf{homog}_{#1}}
  \newcommand{\cone}{\mathcal{C}}
  
  \newcommand{\hbc}[1]{\nv{#1}{H}}

  \newcommand{\nv}[2]{\textsf{NV}_{#1}(#2)}
  \newcommand{\nvr}[2]{\textsf{rNV}_{#1}(#2)}
  \newcommand{\rnv}[2]{\nvr{#1}{#2}}
  \newcommand{\nvf}[2]{\widehat{\textsf{NV}_{#1}}(#2)}
  \newcommand{\nvlace}{\leftsquigarrow}  

\newcommand{\iposlace}{\ensuremath{\stackrel{\imag}{\thicksim}}}
\newcommand{\pposlace}{\ensuremath{\stackrel{P}{\thicksim}}}
\newcommand{\hposlace}{\ensuremath{\stackrel{H}{\thicksim}}}
\newcommand{\hlace}{\ensuremath{\stackrel{H}{\longleftarrow}}}
\newcommand{\place}{\ensuremath{\stackrel{P}{\longleftarrow}}}
\newcommand{\ulace}{\ensuremath{\stackrel{U}{\longleftarrow}}}

\newcommand{\calC}{\mathcal{C}}
\newcommand{\calD}{\mathcal{D}}
\newcommand{\calE}{\mathcal{E}}
\newcommand{\calS}{\mathcal{S}}
\newcommand{\calM}{\mathcal{M}}

\newcommand{\init}[2]{#1\hspace*{-3pt}{\mathbf{\bigl/}}#2}

\newcommand{\up}[1]{\text{\texttt{U}}_{#1}(\complexes)}
\newcommand{\rup}[1]{{\text{\texttt{U}}_{#1}}}
\newcommand{\rupint}[1]{{\text{\texttt{P}}_{#1}}}  

\newcommand{\hb}[1]{\textsf{HP}_{#1}}

\newcommand{\rhb}[1]{\textsf{rHP}_{#1}}

\newcommand{\isknown}[1]{\text{\color{red}\uwave{$#1$}}}
\newcommand {\eqques}{\ensuremath {\stackrel {\mathrm{?}}{=}}}

\newcommand{\bfsf}[1]{\textbf{\textsf{#1}}}

  \newcommand{\quada}{\textsf{Q}_I}
  \newcommand{\quadb}{\textsf{Q}_{II}}
  \newcommand{\quadc}{\textsf{Q}_{III}}
  \newcommand{\quadd}{\textsf{Q}_{IV}}

  \newcommand{\sector}[1]{\mathcal{S}_{#1}}

\newcommand{\diag}{\operatorname{diag}}

\newcommand{\intmod}[1]{{\mathbb{Z}}/{#1}{\mathbb{Z}}}
\newcommand{\septimes}{{\,\mathbf{\otimes}\,}}
\newcommand{\hadprod}{\circledast}
\newcommand{\ppoly}{\textsf{ poly}}

\newcommand{\updn}{{\mathcal{U}\mathcal{L}}}
\newcommand{\roots}{{\textsf{ roots}}}
\newcommand{\norm}[1]{\left\Vert#1\right\Vert} 

\newcommand{\polycpx}{\mathfrak{I}}
\newcommand{\polycpxpos}{\mathfrak{I}^{pos}}

\newcommand{\polycpxf}{\widehat{\mathfrak{I}}}
\newcommand{\polycpxclose}{\overline{\mathfrak{I}}}
\newcommand{\plane}{{\mathsf H}}

\newcommand{\stabled}[1]{{\mathcal{H}_{#1}}} 
\newcommand{\stabledopen}[1]{\mathcal{H}_{#1}^o}
\newcommand{\stabledf}[1]{{\widehat{\mathcal{H}_{#1}}}}
\newcommand{\stabledcf}[1]{{\widehat{\mathcal{H}_{#1}}(\complexes)}}
\newcommand{\polypos}[1]{P^{+}_{#1}}

\newcommand{\stabledc}[1]{{\mathcal{H}_{#1}(\complexes)}}

\newcommand{\reals}{{\mathbb R}}
\newcommand{\complexes}{{\mathbb C}}
\newcommand{\imag}{\boldsymbol{\imath}}
\newcommand{\cals}{\mathcal{S}}

\newcommand{\ttt}{\mathcal{T}}
\newcommand{\cpxlace}{\curlyeqprec}
\newcommand{\maj}{{\,\prec_m\,}}

\newcommand{\poslace}{\ensuremath{\stackrel{+}{\thicksim}}}

\newcommand{\algless}{\operatorname{{\raisebox{-.5mm}{\rotatebox[origin=c]{90}{\textsf{A}}}}}}

\newcommand{\intoper}[1]{{\textstyle\int}^{#1}}
\newcommand{\expoper}{\textsc{exp}}

\newcommand{\intexp}[1]{\intoper{#1}\,\expoper\,}
\newcommand{\rev}[1]{{#1}^{\textsc{ rev}}}

\newcommand{\diffd}{{\text{{\textsf{D}}}}} 
\newcommand{\diffi}{{\text{{\textbf{I}}}}} 
\newcommand{\diffj}{{\text{{\textbf{J}}}}} 
\newcommand{\diffk}{{\text{{\textbf{K}}}}} 
\newcommand{\diffl}{{\text{{\textbf{L}}}}} 
\newcommand{\diffr}{{\text{{\textbf{R}}}}} 
\newcommand{\sdiffi}{{\text{{\textbf{\tiny{I}}}}}} 
\newcommand{\sdiffj}{{\text{{\textbf{\tiny{J}}}}}} 

\newcommand{\image}[1]{\mathfrak{P}_{#1}}
\newcommand{\imagef}[1]{\widehat{\mathfrak{P}}_{#1}}
\newcommand{\imagepos}[1]{\mathfrak{P}^{pos}_{#1}}
\newcommand{\imageposf}[1]{\widehat{\mathfrak{P}}^{pos}_{#1}}

\newcommand{\xdxydy}{x\partialx+y\partialy}
\newcommand{\xdxydypd}{x_1\frac{\partial}{\partial x_1} + \cdots + %
x_d\frac{\partial}{\partial x_d}}
\newcommand{\partialy}{\frac{\partial}{\partial y}} 
\newcommand{\partialx}{\frac{\partial}{\partial x}} 
\newcommand{\polar}[1]{\text{{\textsf{ D}}}_{#1}}
\newcommand{\polard}[1]{\text{{\textsf{ D}}}_{#1}}
\newcommand{\qderiv}{\polard{q}}
\newcommand{\daffine}{\polard{\affa}}

\newcommand{\charpoly}[1]{\text{{\textsc{ cp}}}(#1)}

\newcommand{\lesseq}{{\,\gg\,}}
\newcommand{\greateq}{{\,\ll\,}}
\newcommand{\lesseqeq}{{\,\underline{\gg}\,}}
\newcommand{\greateqeq}{{\,\underline{\ll}\,}}

\newcommand{\lessless}{\operatorname{\lessdot}}
\newcommand{\lessgreat}{\operatorname{\gtrdot}}
\newcommand{\lesslesseq}{{\,\operatorname{\underline{\lessdot}\,}}}
\newcommand{\lessgreateq}{{\operatorname{\underline{\gtrdot}}}}

\newcommand{\weaki}{\blacktriangleleft}

\newcommand{\posplus}{^\mathbf{+}}

\newcommand{\pgreateqeq}{{\,\underline{\ll}\posplus\,}}

\newcommand{\plessless}{\poslace}
\newcommand{\plesslesseq}{\poslace}
\newcommand{\lace}{{\,\thicksim\,}}

\newcommand{\rising}[2]{{\left<\underline{#1}\right>_{#2}}} 
\newcommand{\arising}[3]{{\left<\underline{#1;#2}\right>_{#3}^{\affa}}}
\newcommand{\falling}[2]{{(\underline{#1})}_{#2}}     

\newcommand{\aaa}{\mathbf{a}}
\newcommand{\bbb}{\mathbf{b}}
\newcommand{\ccc}{\mathbf{c}}
\newcommand{\ddd}{\mathbf{d}}

\newcommand{\eee}{\mathbf{e}}

\newcommand{\xx}{\mathbf{x}}
\newcommand{\yy}{\mathbf{y}}
\newcommand{\zz}{\mathbf{z}}
\newcommand{\mydiamond}[1]{\underset{#1}{\,\Diamond}\,}
\newcommand{\mydiamondA}[2]{\underset{#1}{\overset{#2}{\Diamond}}\,}

\newcommand{\bx}{{\mathbf{x}}}
\newcommand{\by}{{\mathbf{y}}}

\newcommand{\mycone}[1]{\textbf{ C}\,\textsc{one}(#1)}

\newcommand{\affa}{{\mathbb{A}}}
\newcommand{\affb}{{\mathbb{B}}}
\newcommand{\alessless}{\,\lessdot^{\affa}\,}

\newcommand{\alesslesseq}{\,\underline{\lessdot}^{\affa}\,}

\newcommand{\epolyint}[2]{\text{\textbf{{P}}}^{#2}/\text{\textbf{ Z}}_{#1}}
\newcommand{\epoly}[1]{\epolyint{#1}{}}
\newcommand{\epolypos}[1]{\epolyint{#1}{+}}
\newcommand{\epolyneg}[1]{\epolyint{#1}{-}}

\newcommand{\allpolyint}[1]{\text{\textbf{P}}^{#1}}

\newcommand{\allpoly}{\allpolyint{}}
\newcommand{\allpolyalt}{\allpolyint{alt}}
\newcommand{\allpolypos}{\allpolyint{pos}}
\newcommand{\allpolyneg}{\allpolypos}
\newcommand{\allpolypm}{\allpolyint{\pm}}
\newcommand{\allpolypmclose}{\overline{\allpolypm}}
\newcommand{\allpolyaltclose}{\overline{\allpolyalt}}
\newcommand{\allpolyposclose}{\overline{\allpolyneg}}
\newcommand{\allpolysep}{\allpolyint{sep}}

\newcommand{\allpolyaffine}{\allpolyint{\affa}}
\newcommand{\allpolyaffinef}{\allpolyintf{\affa}}

\newcommand{\allpolypsd}{$\allpoly$-positive semi-definite}  

\newcommand{\hypergeo}[2]{{\ _{#1}F_{#2}}}
\newcommand{\foneone}{\hypergeo{1}{1}}  
\newcommand{\fonetwo}{\hypergeo{1}{2}} 

\newcommand{\hyper}[1]{{\textsf{ Hyp}}_{#1}}
\newcommand{\hyperclose}[1]{{\overline{\textsf{ Hyp}}}_{#1}}
\newcommand{\hyperpos}[1]{{\textsf{ Hyp}}_{#1}^{pos}}
\newcommand{\smp}{stable matrix polynomial}

\newcommand{\gsubpos}{\gsub} 

\newcommand{\gsubalt}{\gsub^{\text{alt}}}
\newcommand{\gsubsep}{\gsub^{\text{sep}}}

\newcommand{\gsubgen}{\gsub^{\text{gen}}}
\newcommand{\gsubcloseneg}{\overline{\gsubalt}}

\newcommand{\gsubclosepos}{\overline{\gsubpos}}
\newcommand{\gsubposclose}{\overline{\gsubpos}}
\newcommand{\gsub}{\mathbf{P}}
\newcommand{\gsubclose}{\overline{\gsub}}

\newcommand{\gsubplus}{\gsub^{\text{pos}}} 
\newcommand{\gsubcloseplus}{\gsubclose^{\text{pos}}} 
\newcommand{\gsubplusclose}{\gsubcloseplus}

\newcommand{\gsubint}[1]{\gsub^{#1}}

\newcommand{\gsubf}{\allpolyintf{}}
\newcommand{\gsubaltf}{\allpolyintf{alt}}
\newcommand{\gsubposf}{\allpolyintf{pos}}

\newcommand{\xsub}{\mathbf{Sub}}

\newcommand{\mywp}{\raisebox{1mm}{\ensuremath{\boldsymbol{\wp}}}}
\newcommand{\partialpoly}[1]{\mywp_{#1}}
\newcommand{\partialpolypos}[1]{\mywp^{pos}_{#1}}
\newcommand{\partialpolyclose}[1]{\overline{\mywp}_{#1}}
\newcommand{\partialpolyposclose}[1]{\overline{\mywp}^{pos}_{#1}}
\newcommand{\partialpolyf}[1]{\widehat{\mywp}_{#1}}

\newcommand{\partX}[1]{\frac{\partial^{#1}}{\partial x^{#1}}} 
\newcommand{\partY}[1]{\frac{\partial^{#1}}{\partial y^{#1}}}

\newcommand{\allpolyintf}[1]{\widehat{{\allpoly}^{#1}}}
\newcommand{\allpolyf}{\allpolyintf{}}
\newcommand{\allpolyposf}{\allpolyintf{pos}}
\newcommand{\allpolynegf}{\allpolyposf}

\newcommand{\allpolyaltf}{\allpolyintf{alt}}
\newcommand{\allpolysepf}{\allpolyintf{sep}}

\newcommand{\pospoly}{\mathcal Q}

\newcommand{\hpoly}{{\mathfrak C}}

\newcounter{numprob}

\newcommand{\problem}[1]{%
\addtocounter{numprob}{1}%
\colorbox{yellow}{#1}%
\index{Z-problem-chap-\thechapter-\thenumprob}%
\setlength{\fboxrule}{.4pt}}

\newcommand{\noproblem}[1]{}

\newcounter{mynoteC}
\newcommand{\mynote}[1]{%
\addtocounter{mynoteC}{1}%
\index{Z-note-chap-\thechapter-\themynoteC}%
\marginpar{%
\vskip-\baselineskip 
\raggedright\footnotesize
\itshape{\hrule\smallskip\colorbox{yellow}{\parbox{1in}{\themynoteC: #1}}\par\smallskip\hrule}}}

\newcommand{\tpp}{\mathcal{S}}
\newcommand{\nsd}{negative subdefinite}

\newcommand{\smalltwo}[4]{%
\left(\begin{smallmatrix} #1&#2 \\ #3&#4\end{smallmatrix}\right)}
\newcommand{\smalltwodet}[4]{%
\left|\begin{smallmatrix} #1&#2 \\ #3&#4\end{smallmatrix}\right|}
\newcommand{\smalltwobyone}[2]{%
\left(\begin{smallmatrix} #1 \\ #2\end{smallmatrix}\right)}
\newcommand{\partmatrix}[4]{%
\begin{pmatrix}
  #1 & \vdots       & #2 \\
       \hdotsfor{3}      \\
  #3 & \vdots       & #4
\end{pmatrix}
}
\newcommand{\twodet}[4]{\begin{vmatrix}#1 & #2 \\ #3 & #4 \end{vmatrix}}
\newcommand{\smallsquare}[4]{
\raisebox{0.5\depth}{\xymatrix@-1pc{
#1 \ar@{<-}[r]&  #2 \ar@{<-}[d] \\
#3 \ar@{->}[u]&  #4 \ar@{->}[l] 
}}
}
\renewcommand{\smallsquare}[4]{
\raisebox{4ex}{\xymatrix@=.2pc{
&& #2 \ar@{->}[dll] \ar@{<-}[drr] && \\
#1 &&&& #4\\
&& #3 \ar@{->}[ull] \ar@{<-}[urr] &&
}}}

\renewcommand{\smallsquare}[4]{
\raisebox{4ex}{\xymatrix@=.2pc{
&& #2 \ar@{->}[dll] \ar@{<-}[drr] && \\
#1 &&&& #4\\
&& #3 \ar@{->}[ull] \ar@{<-}[urr] &&
}}}

\newcommand{\smallsquareLL}[4]{
\raisebox{4ex}{\xymatrix@=.2pc{
&& #2  \ar@{}[dll]|{\myrot{\lessless}{45}}
\ar@{}[drr]|{\myrot{\lessless}{315}}  && \\
#1 &&&& #4 \ar@{}[dll]|{\myrot{\lessless}{45}}\\
&& #3 \ar@{}[ull]|{\myrot{\lessless}{315}}  &&
}}}

  \newcommand{\myrot}[2]{\rotatebox[origin=c]{#2}{$#1$}}
  \newcommand{\myup}[1]{\myrot{#1}{270}}
  
  \newcommand{\arbsquare}[6]{
    \begin{matrix}
#1 & #5 & #2 \\
#6 && #6 \\
#3 & #5 & #4 
    \end{matrix}
}
\newcommand{\arbupsquare}[5]{\arbsquare{#1}{#2}{#3}{#4}{#5}{\myup{#5}}}

\newcommand{\chap}[1]{Chapter~\ref{cha:#1}}
\newcommand{\chapsec}[2]{$\mathsection$\,\ref{cha:#1}.\ref{sec:#2}}

\newcommand{\slem}[1]{\pageref{lem:#1}}
\newcommand{\sthm}[1]{\pageref{thm:#1}}
\newcommand{\scorr}[1]{\pageref{cor:#1}}
\newcommand{\spropp}[1]{\pageref{prop:#1}}

\newcommand{\genfct}[4]{%
\text{#1} &%
\makebox[.1in]{}  &\displaystyle\sum_{i=0}^\infty & #2  &\dfrac{(-y)^i}{i!} &%
\quad=\quad #3 & #4\\
}
\newcommand{\genfctpp}[4]{%
\text{#1} &%
\makebox[.1in]{}  f(x,y) &  \mapsto & #2 &
\quad\quad #3 & #4\\
}

\newcommand{\genfctpd}[4]{%
\text{#1} &%
\makebox[.1in]{}  &\displaystyle\sum_{\diffi} & #2  &\dfrac{(-\yy)^\sdiffi}{\diffi!} &%
\quad=\quad #3 & #4\\
}

\newcommand{\mypage}[1]{(p.~\pageref{#1})}
\newcommand{\seepage}[1]{(See p.~\pageref{#1}.)}

\newcommand{\Szego}{Szeg\"{o}}
\newcommand{\Polya}{P\'{o}lya}
\newcommand{\Mobius}{M\"{o}bius}
\newcommand{\bbs}{Borceau, Branden, Shapiro}

\newcommand{\onepoly}[2]{
\noindent\hrulefill\\[.2cm]
\noindent{\textbf{\large #1}}\\[.2cm]
\centerline{#2}\\[.2cm]
\noindent%
$
\begin{array}{r>{$}p{4in}<{$}}
}

\newcommand{\workhome}{somewhere}

\newcounter{aquery}

\title{Polynomials, roots, and interlacing \\[.2cm] \textsf{ Version\ \version}}
\author{Steve Fisk \\ Bowdoin College \\ Brunswick, Me 04011}
\date{\today}

\makeindex
\CompileMatrices

\usepackage{memhfixc}

\begin{document}

\chapterstyle{section}

\frontmatter

\begin{titlingpage}
\aliaspagestyle{titlingpage}{plain}
\setlength{\droptitle}{30pt} 
\maketitle
\end{titlingpage}

\setcounter{tocdepth}{1}
\tableofcontents
\listoffigures
\listoftables
\pagebreak
\include{preface}


\mainmatter
\chapterstyle{BlueBox}

\part{Polynomials with all real roots}

\chapter{Polynomials in One Variable}

\label{cha:polynomials}

\renewcommand{\TimeStampStart}{Friday, January 18, 2008: 09:57:49}
\mytoday  


\section{Notation and basic definitions}
\label{sec:intro-notation}
We begin with the definitions of the various kinds of interlacing, the 
classes of polynomials of interest, and the properties of linear
transformations that we shall be using.

\index{interlacing}

If $f$ is a polynomial of degree $n$ with all real roots, then define 
$$
\roots(f) = (a_1,\dots,a_n)$$
where $a_1 \le a_2 \le \dots \le
a_n$ are the roots of $f(x)=0$.  If we write $\roots(f)$ then we are
assuming that $f$ actually has all real roots.  Given two polynomials
$f,g$ with $\roots(f)=(a_1,\dots,a_n)$ and $\roots(g) =
(b_1,\dots,b_m)$ then we say that $f$ and $g$ \emph{interlace} if
these roots alternate. If all the roots are  distinct then there are
only four ways that this can happen:

\begin{alignat}{2}
 g &\greateq &f \quad\quad& \label{eqn:defn-2} 
 a_1 < b_1 < a_2 < b_2 < \cdots < a_n < b_n \\
 g& \lessless& f \quad\quad& 
 b_1 < a_1 < b_2 < a_2 < \cdots < b_n  \nonumber \\
 f& \greateq &g \quad\quad& 
 b_1 < a_1 < b_2 < a_2 < \cdots < b_n < a_n \nonumber \\
 f &\lessless &g \quad\quad& 
 a_1 < b_1 < a_2 < b_2 < \cdots < a_n.  \nonumber 
\end{alignat}
We define $g \lesseq f$ to mean that $f \greateq g$, and $g\lessgreat
f$ to mean $f \lessless g$. Notice that when we write $f\greateq g$ we
assume that both $f$ and $g$ have all real roots, which are
necessarily distinct.  If two polynomials interlace then their degrees
differ by at most one.


If the inequalities in \eqref{eqn:defn-2} are relaxed, so that all
inequalities ``$<$'' are replaced by ``$\le$'', then we write $g
\greateqeq f$, $g\lesslesseq f$, $f\greateqeq g$ and $f\lesslesseq g$
respectively.  We call interlacings with $\lessless$ or $\greateq$
\emph{strict interlacings}.  It's easy to see that non-strict
interlacings are limits of strict interlacings - see
Remark~\ref{rem:strict-limit} for a non ad hoc argument.

\index{strict interlacings}
\index{interlacing!strict}

\index{interlacing}
\index{\ zzlesseq@$\lesseq$}
\index{\ zzlessless@$\lessless$}
\index{\ zzgreateq@$\greateq$}
\index{\ zzlesslesseq@$\lesslesseq$}
\index{\ zzlesseqeq@$\lesseqeq$}

If $f$ and $g$ are two polynomials such that there is a polynomial $h$
with $f \greateqeq h$, $g \greateqeq h$, (or $h \lesslesseq g$, $h \lesslesseq f$)
then we say that $f,g$ have a \emph{common interlacing}.  Unless
stated otherwise, we assume that $f$ and $g$ have the same degree.
\index{common interlacing}
\index{interlacing!common}

The two interlacings $f\greateqeq g$ and $f \lesslesseq
g$ have very similar properties. To eliminate many special cases, we
introduce the notation
\begin{alignat*}{2}
  f \longleftarrow g \quad\quad\text{which means }\quad\quad&
  f\greateqeq g& \quad\text{ or }\quad & f\lesslesseq g
\end{alignat*}
\noindent%
One way to remember this notation is to notice that $f\longleftarrow
g$ means that
\begin{itemize}
\item $f$ and $g$ interlace
\item The largest root belongs to $f$.
\item The arrow points to the polynomial with the largest root.
\end{itemize}

There are several important subclasses of polynomials with all real roots.
\index{\ Ppoly@$\allpoly$}
\index{\ Ppolypos@$\allpolyalt$}
\index{\ Ppolyneg@$\allpolypos$}
\index{\ Ppolyint@$\allpolyint{I}$}
\index{\ Ppolyposneg@$\allpolypm$}
\begin{align*}
  \allpoly &=  \{\text{ all polynomials with all real roots}\} \\
  \allpolypos &= \{\text{ all polynomials with all real roots, all negative,} \\
& \qquad \text{and all coefficients are the same sign}\}\\
  \allpolyalt &= \{\text{ all polynomials with all real roots, all
    positive,} \\
  &\qquad \text{and the coefficients alternate in sign}\}\\
  \allpolypm &= \allpolyalt \cup \allpolypos \\
  \allpolyint{I} &= \{\text{ all polynomials with all real roots, all
    roots lying in a set $I$}\} \\
  \allpolypmclose &= \allpolyint{(-\infty,0]} \cup
  \allpolyint{[0,\infty)}
\end{align*}
It is helpful to remember that that $\allpolypos$ has a superscript
\emph{pos} because whenever the leading coefficient is positive all
the coefficients are positive.  If we want to only consider
polynomials of degree $n$ then we write $\allpoly(n)$,
$\allpolypos(n)$, and so on. 

 We are particularly interested in properties of
linear transformations of polynomials. Suppose that $T$ is a
linear transformation that maps polynomials to polynomials.

\index{interlacing!preserves}
\index{linear transformation!preserves roots}
\index{linear transformation!preserves interlacing}
\index{linear transformation!preserves degree}
\index{linear transformation!preserves sign of leading coefficient}

\begin{itemize}
\item  $T$ \emph{preserves roots} if whenever
  $f$ has all real roots, so does $Tf$. Using the notation above, we
  can write this as $T\colon{}\allpoly\longrightarrow\allpoly$.
\item  $T$ \emph{preserves interlacing} if
  whenever $f$ and $g$ interlace, so do $Tf$ and $Tg$.
\item  $T$ \emph{preserves degree} if the
  degrees of $f$ and $Tf$ are the same.
\item $T$ \emph{preserves the sign of the leading coefficient} if $f$
  and $Tf$ have leading coefficients of the same sign.

\end{itemize}

The quintessential interlacing theorem is a well known consequence of
Rolle's Theorem. The proof given below analyzes the graph of $f$ and
$f^\prime$ - see Remark~\ref{rem:rolle-2} and Remark~\ref{rem:rolle-3}
for proofs that do not involve any analysis of the graph.

\index{Rolle's Theorem}

\begin{theorem} \label{thm:rolle}
  If $f$ is a polynomial with all real roots, and all these roots are
  distinct, then $f \lessless f^\prime$.  
\end{theorem}

\begin{proof}
  By Rolle's Theorem, there is a root of $f^\prime$ between
  consecutive roots of $f$.  If the degree of $f$ is $n$, then there
  are $n-1$ intervals between roots of $f$, so $f^\prime$ has at least
  $n-1$ roots.  Since $f^\prime$ has degree $n-1$, it has exactly
  $n-1$ roots.
\end{proof}

If $f$ has all real roots, but all the roots of $f$ are not distinct,
then $f$ and $f^\prime$ have roots in common.  In this case we have
that $f \lesslesseq f^\prime$.

\index{\ Deriv@$\diffd{}$} The derivative operator $\diffd f =
f^\prime$ is of fundamental importance. It is a linear transformation
$\allpoly\longrightarrow\allpoly$,
$\allpolyalt\longrightarrow\allpolyalt$, and
$\allpolypos\longrightarrow\allpolyneg$. We just observed that
$\diffd$ preserves roots.  We will see in Theorem~\ref{thm:rolle-2} that it also
preserves interlacing.

We will occasionally make use of the stronger version of Rolle's theorem that
states that there are an \emph{odd} number of roots (counting
multiplicity) of $f^\prime$ between any two consecutive distinct roots of $f$.

\begin{figure}[htbp]
  \begin{center}
    \leavevmode
    \epsfig{file=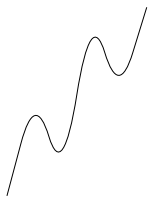,height=1in}
    \caption{A quintic with 4 critical points, and at most 3 roots}
    \label{fig:quintic}
  \end{center}
\end{figure}

\begin{example} \label{ex:no-integral}
  If $f$ is a polynomial with all real roots, then it is not
  necessarily true that some anti-derivative of $f$ has all real
  zeros.  Equivalently, the derivative does not map $\allpoly(n)$ \emph{onto}
  $\allpoly(n-1)$.  A quintic $p(x)$ such as in
  Figure~\ref{fig:quintic} has four critical points, so $p^\prime$ has
  four zeros.  However, no horizontal line intersects $p$ in more than
  three points, and so no integral of $p^\prime$ will have more than
  three roots.  For instance, $(x+1)(x+2)(x+4)(x+5)$ is not the
  derivative of a quintic with all real roots.
\end{example}
                                

\index{polynomials!integral not in $\allpoly$}

We introduce a relation that is weaker than interlacing. See
Section~\ref{sec:relation-weaker-than} for more properties.
  \begin{definition}
    Suppose that 

\begin{align*}
  f(x) &= (x-r_1)\cdots(x-r_n) & \quad & r_1\le r_2\cdots\le r_n \\
  g(x) &= (x-s_1)\cdots(x-s_m) & \quad & s_1\le s_2\cdots\le s_m 
\end{align*}

We say that $g\weaki f$ if $m\le n$, and $s_i\le r_i$ for $i=1,\dots,m$.

In other words, each root of $g$ is at most the corresponding root of
$f$. Note that the degrees of $f$ and $g$ do not need to be equal.
  \end{definition}

Common interlacing and $\weaki$ imply interlacing.
The proof of the Lemma is immediate.

\begin{lemma} \label{lem:ci}
\index{interlacing!common}
\index{common interlacing}
If
\begin{enumerate}
\item $f\weaki g$
\item $f$ and $g$ have the same degree
\item $f$ and $g$ have a common interlacing
\end{enumerate}
 then $f \greateqeq g$.
\end{lemma}

Common interlacings can be split into two interlacings.
\index{common interlacing}\index{interlacing!common}

\begin{lemma}
  If $f,g$ have the same degree and have a common interlacing then
  there are $f_1,f_2,g_1,g_2\in\allpoly$ such that
  \begin{align*}
    f_1 & \longleftarrow g_1 & f & = f_1 f_2 \\
    g_2 & \longleftarrow f_2 &  g & = g_1 g_2
      \end{align*}
\end{lemma}
\begin{proof}
  Choose $h\lessless g$ and $h \lessless f$. There is one root $r_i$
  of $f$ and one root $s_i$ of $g$ between the $i$th and $(i+1)$st
  roots of $h$. We may assume that the roots of $f$ and $g$ are
  distinct. The desired polynomials are defined by 
  \begin{align*}
    \roots(f_1) &= \{r_i\,\mid\, r_i< s_i \}&     
    \roots(g_1) &= \{s_i\,\mid\, r_i< s_i \}\\
    \roots(f_2) &= \{r_i\,\mid\, r_i> s_i \}&
    \roots(g_2) &= \{s_i\,\mid\, r_i> s_i \}
  \end{align*}
\end{proof}

\begin{remark}
  Here are a few important properties about polynomials. See
  \cite{rahman}.  First of all, the roots of a polynomial are
  continuous functions of the coefficients. Consequently, if a
  sequence of polynomials $\{f_n\}$ converges to the polynomial $f$
  and all $f_n$ have all real roots, then so does $f$.  In addition if
  there is a polynomial $g$ such that $g\longleftarrow f_n$ for all
  $n$ then $g\longleftarrow f$. We should also mention that roots of
  polynomials with real coefficients come in conjugate pairs. It
  follows that if a polynomial of degree $n$ has $n-1$ real roots,
  then it has all real roots. Lastly, if $fg$ has all real roots,
  then so do $f$ and $g$.
\end{remark}

\section{Sign interlacing}
\label{sec:intro-sign-interlacing}

\index{sign interlacing}

The easiest way to show that two polynomials with all real roots
interlace is to show that one \emph{sign interlaces} the other.
The \emph{sign} of a number $c$ is defined as

\begin{eqnarray*}
sgn\ c &=& 
\begin{cases}
  1 & \text{ if } c>0 \\
  -1 & \text{ if } c<0 \\
  0 & \text{ if } c=0 \\
\end{cases}
\end{eqnarray*}

We say that \emph{$f$ sign interlaces $g$} if
the sign of $f$ alternates on the roots of $g$.  More precisely, there
is an $\epsilon=\pm1$ so that if $r$ is the $i$-th largest root of $g$
then the sign of $f(r)$ is $\epsilon(-1)^i$.  

We usually assume that polynomials have positive leading coefficient.
In this case, if $f\greateq g$ then there is exactly one root of $f$
that is larger than the largest root $r$ of $g$. Since $f$ has
positive leading coefficient, it must be the case that $f(r)$ is
negative. Consequently, if $g$ has degree $n$ then the sign of $f$ at
the $i$-th largest root of $g$ is $(-1)^{i+n+1}$.

\begin{figure}[htbp]
  \begin{center}
    \leavevmode
    \epsfig{file=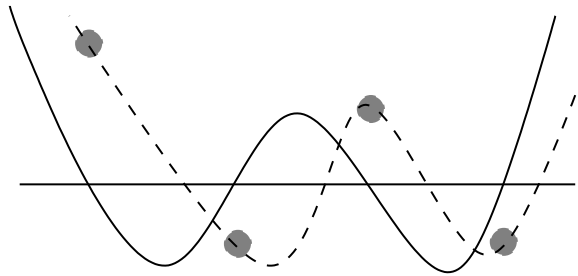,width=3in}
    \caption{Two interlacing polynomials}
    \label{fig:interlace-1}
  \end{center}
\end{figure}

It should be clear that if $f$ and $g$ interlace then $f$ sign
interlaces $g$, and $g$ sign interlaces $f$.  The converse is not
always true, and is discussed in the first three lemmas.
Figure~\ref{fig:interlace-1} shows two interlacing polynomials
$f\greateq g$, where $f$ is the dashed line, and $g$ the solid line.
The solid circles indicate the value of $f$ at the zeros of $g$ -
notice that these circles alternate above and below the $x$-axis.

\index{sign interlacing}
\index{interlacing!sign}

\begin{lemma}
  \label{lem:int-1}
  Suppose $f,g$ are polynomials satisfying
  \begin{itemize}
  \item $f$ has all real roots, all distinct.
  \item $\deg(f) = \deg(g)+1$
  \item $g$ sign interlaces $f$
  \end{itemize}
  then $f \lessless g$.
\end{lemma}

\begin{proof}
  Assume that $f$ has $n$ roots.  Since the sign of $g$ alternates
  between consecutive roots of $f$, we find $n-1$ roots of $g$, one
  between each pair of consecutive roots of $f$.
\end{proof}

If we knew that $f$ sign interlaced $g$ in Lemma~\ref{lem:int-1}, then
it does not follow that $f\lessless g$.  For instance, take
$f=x^3$ and $g = (x-1)(x+1)$, so that the roots of $f$ and $g$ are $\{0,0,0\}$ and
$\{-1,1\}$.  The reason that $g$ and  $f$ do not interlace is
that although we have roots of $f$ between roots of $g$, there aren't
any roots of $f$ outside of the roots of $(g)$.  In order to find these
roots, we  must control the behaviors of $f$ and $g$  by specifying
appropriate signs at the largest and smallest roots of $g$.

\begin{figure}[htbp]
  \begin{center}
    \leavevmode
    \epsfig{file=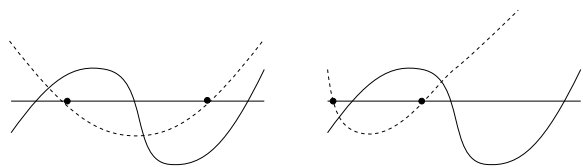,width=3in}
    \caption{Sign interlacing doesn't imply interlacing.}
    \label{fig:possibility-1}
  \end{center}
\end{figure}

Figure~\ref{fig:possibility-1} shows graphically that if $f$ (solid) sign
interlaces $g$ (dashed)  then it is possible that $f$ is not interlaced by $g$.

\begin{lemma}
  \label{lem:int-2}
  Suppose $f,g$ are polynomials satisfying
\begin{itemize}
\item $f$ has all real roots, all distinct.
\item $\deg(f) = \deg(g)-1$
\item $g$ has positive leading coefficient.
\item $g$ sign interlaces $f$.
\item Let $\alpha$ be the smallest root of $f$, and  $\beta$ the largest
  root of $f$.   If either
  \begin{enumerate}
  \item $sgn\ f(\beta) = -1$ or
  \item $sgn\ f(\alpha) = 1$
  \end{enumerate}
  then $g \lessless f$.
\end{itemize}
  
\end{lemma}
\begin{proof}
Since $g$ sign interlaces $f$ there is a root of $g$ between adjacent
roots of $f$.  This determines $\deg(f)-2$ roots, so it remains to find
where the last  two roots are. If $sgn(g(\beta)) =  -1$ then
since $\lim_{x\rightarrow\infty} sgn(g(x)) = 1$ there are an odd
number of roots (hence exactly one) to the right of $\beta$.
Similarly there is a root to the left of $\alpha$. See Figure~
\ref{fig:possibility-1}. 
\end{proof}

\begin{figure}[htbp]
  \begin{center}
    \leavevmode
    \epsfig{file=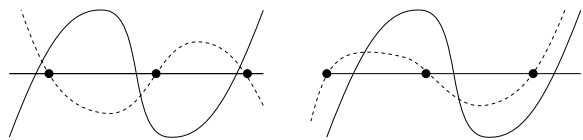,width=3in}
    \caption{Interlacing when the degrees are equal.}
    \label{fig:possibilities-2}
  \end{center}
\end{figure}

If $f$ and $g$ both have the same degree and they sign interlace, then
either $f \greateq g$ or $g \greateq f$. Again, we need extra
information about the behavior on the extremes to determine which case
holds.  Figure~\ref{fig:possibilities-2} shows the two possibilities.

\begin{lemma}
  \label{lem:int-0}
  Suppose $f,g$ are polynomials satisfying
\begin{itemize}
\item $f$ has all real roots, all distinct.
\item $g$ has positive leading coefficient.
\item $\deg(f) = \deg(g) = n$
\item $g$ sign interlaces $f$.
\item Let $\alpha$ be the smallest root of $f$ and $\beta$ the largest
  root.  Then
  \begin{enumerate}
  \item If $sgn\ g(\beta) = 1$ then $f \greateq g$.
  \item If $sgn\ g(\beta) = -1$ then $g \greateq f$.
  \item If $sgn\ g(\alpha) = - (-1)^n$ then $g \greateq f$.
  \item If $sgn\ g(\alpha) = (-1)^n $ then $f \greateq g$.
  \end{enumerate}
\end{itemize}
\end{lemma}

\begin{proof}
  Assume that $f$ has $n$ roots.  Since the sign of $g$ alternates
  between consecutive roots of $f$, we find $n-1$ roots of $g$, one
  between each pair of consecutive roots of $f$.  Now $g$ has an odd
  number of roots in each of these intervals, and there cannot be
  three roots of $g$ in any of them, for then $g$ would have at least
  $n+2$ roots, so there is one more root of $g$, to either the left of
  the smallest root of $f$, or to the right of the largest root.
  Thus, we have either $f \greateq g$ or $g \greateq f$.

  All four parts are similar, so we will only establish the first one.
  Assume that $sgn\ g(\beta) = 1$. Counting the number of
  alternations of sign of $g$ gives
  \begin{eqnarray*}
    \lim_{x\rightarrow-\infty} sgn\,g(x) &=& (-1)^n  \\ sgn\ 
    g(\alpha) &=& (-1)^{n-1} sgn\ \beta = (-1)^{n-1} 
  \end{eqnarray*}
  and hence there is a root of $g$ in $(-\infty,\alpha)$. This shows
  that $f \greateq g$.
\end{proof}

As a first consequence of these elementary but fundamental results, we 
choose

\begin{lemma} \label{lem:sum-1}
  If $f,g$ have positive leading coefficients and there is an $h$ so
  that $h\longleftarrow f,g$ then for all positive $\alpha$ and $\beta$ the
  linear combination $\alpha f + \beta g$ has all real roots and $h
  \longleftarrow \alpha f+\beta g$.
\end{lemma}

\begin{proof}
  $f,g$ both sign interlace $h$, and since they have positive leading
  coefficients, they are the same sign on all the roots of $h$.
  Consequently, $\alpha f + \beta g$ has the same sign on a root of $h$
  as does $f$ or $g$. The result follows from Lemmas~\ref{lem:int-0}
  and \ref{lem:int-2}. 
\end{proof}

  Although the sum of two arbitrary polynomials with all real roots
  does not necessarily have all real roots, they might if one of them
  is sufficiently small.
  \begin{lemma}\label{lem:add-small}
    Suppose $f\in\allpoly$ has degree $n$ and all distinct roots. If
    $g\in\allpoly$ has degree $n+1$ then $f+\epsilon g\in\allpoly$ for
    sufficiently small $\epsilon$.
  \end{lemma}
  \begin{proof}
    The roots of $f+\epsilon g$ converge to the roots of $f$ as
    $\epsilon$ goes to zero. Since the $n$ roots of $f$ are all
    distinct, $f+\epsilon g$ has $n$ real roots for $\epsilon$
    sufficiently small. Since there is only one root unaccounted for,
    it must be real.
  \end{proof}

The next well-known lemma is easy to see geometrically.\added{12/2/7}

  \begin{lemma}\label{lem:f+c}\ 
    \begin{enumerate}
    \item If there is an $\alpha$ such that $f(x)+c\in\allpoly$ for
      all $|c|>\alpha$ then $f$ is linear.
    \item If there is an $\alpha$ such that $f(x)+c\in\allpoly$ for
      all $c>\alpha$ then $f$ is quadratic with negative leading coefficient.
    \item If there is an $\alpha$ such that $f(x)+c\in\allpoly$ for
      all $c<\alpha$ then $f$ is quadratic with positive leading coefficient.
    \end{enumerate}
  \end{lemma}

\begin{figure}[htbp]
  \begin{center}
    \leavevmode
    \epsfig{file=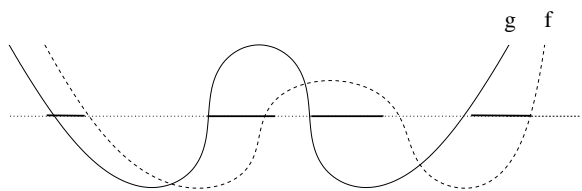,width=3in}
    \caption{Intervals containing roots of $f+\alpha g$, $\alpha>0$.}
    \label{fig:intervals}
  \end{center}
\end{figure}

Sign interlacing determines the intervals where sums of polynomials
have roots. In general the roots lie in alternating intervals, but
which set of intervals depends on the behavior at infinity. This is in
turn determined by the sign of the leading coefficient of the sum.
For instance Figure~\ref{fig:intervals} shows two polynomials $f$
(dashed) and $g$ (solid).  The intervals on the $x$-axis where
$f+\alpha g$ has roots for positive $\alpha$ are solid; the intervals where
they do not have a root (since both functions have the same sign) are
dotted.

\index{determinants!of coefficients}
\begin{lemma} \label{lem:interlace-and-log-concave}
  Assume that polynomials $f,g,h$  have positive
  leading coefficients. If $f$ and $g$ interlace, and
  $\smalltwodet{f}{g}{g}{h}<0$ for all $x$ then $g$ and $h$ interlace. 
  In more detail we can say that
  \begin{itemize}
  \item If $f\lessless g$ then $g\lessless h. $
  \item If $f\greateq g$ then $g\longleftarrow h$.
  \item If $f\lessgreat g$ then $g\longrightarrow h$.
  \item If $f\lesseq g$ then $g\lesseq h$.
  \end{itemize}
\end{lemma}
\begin{proof}
  Assume that the degree of $g$ is $n$.  If $r$ is a root of $g$ then
  upon evaluating the determinant at $r$ we find that $f(r)h(r)<0$.
  Since $f$ sign interlaces $g$, so also does $h$ sign interlace $g$.
  This implies that the degree of $h$ is at least $n-1$. It is not
  possible that the degree of $fh$ is greater than $2n$, for then
  $fh-g^2$ will be positive for large values of $x$, since the leading
  coefficient of $fh$ is positive.  We now consider the four cases in
  turn.
  
  If $f\lessless g$ then the degree of $f$ is $n+1$. Since the degree
  of $h$ is at least $n-1$, and the degree of $fh$ is at most $2n$,
  the degree of $h$ is exactly $n-1$.  Since $h$ sign interlaces $g$,
  we can apply Lemma~\ref{lem:int-1} to conclude that $g\lessless h$.
  
  If $f\lessgreat g$ then the degree of $f$ is $n-1$, and the degree
  of $h$ is either $n+1$, $n-1$ or $n$.  The sign of $f$ at the
  largest root $r$ of $g$ is positive, so the sign of $h$ at $r$ is
  negative. Since the leading coefficient of $h$ is positive, there is
  a root of $h$ in $(r,\infty)$, and so the degree of $h$ can not be
  $n-1$. Lemma~\ref{lem:int-0} now implies that $h\lessless g$ if the degree of
  $h$ is $n+1$. If the degree of $h$ is $n$ then $h\greateq g$.
  
  When $f\lesseq g$ then the degree of $h$ is $n$ or $n-1$. Again the
  sign of $f$ on the largest root $r$ of $g$ is positive, so there is
  a root of $h$ in $(r,\infty)$. Thus $h$ has exactly $n$ roots and
  $g\lesseq h$.

  If $f\greateq g$ then a similar argument applies.
\end{proof}

\begin{cor} \label{cor:fbxcg}
  If $f\lessless g$ have positive leading coefficients and
  $h = af +(bx+c)g$ then $h\lessless g$ if and only if the sign of the 
  leading coefficient of $h$ is $sign(a)$.
\end{cor}
\begin{proof}
  Since $f$ sign interlaces $g$ it follows that $h$ sign interlaces
  $g$. If $\beta$ is the largest root of $g$ then $f(\beta)<0$ so the
  sign of $h(\beta)$ is $-sign(a)$.  Use Lemma~\ref{lem:int-2}.
\end{proof}

As a consequence we get a special case that will prove useful later.
\begin{cor} \label{cor:agxgp}
  If $g\in\allpoly(n)$ has positive leading coefficient and 
  $a\not\in(0,n)$ then $xg^\prime-ag\in\allpoly$. 
  
  If $a\in(0,n)$ and $g$ has $p$ positive roots and $m$ negative
  roots then $xg^\prime-ag$ has at least $p-1$ positive roots and at
  least 
  $m-1$ negative roots. It is possible that $xg^\prime-ag\not\in\allpoly$.
\end{cor}
\begin{proof}
  We may assume that the leading coefficient of $g$ is $1$.  The
  leading coefficient of $xg^\prime-ag$ is $(n-a)$. If $a>n$
  then we can apply Corollary~\ref{cor:fbxcg}. If $a<0$ then note that $xg
  \lesslesseq g$ and $xg\lesslesseq xg^\prime$ and therefore $xg
  \lesslesseq  xg^\prime + (-a)g$.

  If $a\in(0,n)$ then $xg^\prime-ag$ sign interlaces the positive
  roots of $g$ and the negative roots of $g$. This implies that it has 
  at least $p-1$ positive roots, and $m-1$ negative roots, but it does 
  not imply that the remaining two roots are real.  For example, if
  $g=x^2-1$ and $a=1$ then $xg^\prime-g = x^2+1$ which has no real
  roots. 
\end{proof}

The locations of the roots of an arbitrary linear combination can be 
nearly any value. For instance, 

\begin{lemma} \label{lem:pos-neg-1}
  If $f,g$ interlace, are not multiples of one another, and $(a,b)$ is
  any non-empty interval then there is an $\alpha$ such that
\begin{itemize}
\item  $ f + \alpha g \in\allpoly$.
\item  $f+\alpha g$ has a root in $(a,b)$.
\end{itemize}
\end{lemma}
\begin{proof}
  Choose an $r\in(a,b)$ that is not a root of $g$.  If we define
  $\alpha = -\dfrac{f(r)}{g(r)}$ then $r$ is a root of $f+\alpha g$.
  Now $f+\alpha g$ is not the zero polynomial since $f$ and $g$ are
  not multiples of one another, so $f+\alpha g$ is in $\allpoly$.
\end{proof}

However, in the case of the derivative, we can give precise
information about the location of the roots of a linear combination.

\begin{lemma}\label{lem:pos-neg-3}
  If $f\in\allpolyint{(a,b)}(n)$ and $\alpha>0$ then $
  f+\alpha f'\in\allpolyint{(a-\alpha{n},b)}$. 
\end{lemma}
\begin{proof}
  The only question is the location of the smallest root. Assume $f$
  has positive leading coefficient. Since $sgn\,f(a-\alpha{n}) =
  (-1)^n$ it suffices to show that
\begin{equation}
  \label{eqn:pos-neg-3}
\frac{1}{\alpha} > \frac{-f'}{f}(a-\alpha{n})
\end{equation}
since this implies $sgn\,(f+\alpha f')(a-\alpha{n}) = (-1)^n$.
Suppose that the roots of $f$ are $r_1,\dots,r_n$. Then
\begin{align*}
    \frac{-f'}{f}(a-\alpha{n}) &= 
\sum \frac{1}{\alpha{n}-a +r_i} \\
\intertext{Now since $r_i>a$ it follows that 
$ \alpha{n}-a+r_i   > \alpha{n}$ and so}
\frac{1}{\alpha} = n \,\frac{1}{\alpha n} &> \sum \frac{1}{\alpha{n}-a+r_i}
  \end{align*}
\end{proof}

The following lemma combines addition and subtraction, since
$f^2-g^2=(f+g)(f-g)$. It is not a coincidence that the terms come from
the binomial expansion $(f+g)^2$ - the result will be generalized in
Lemma~\ref{cor:fg-imag}. 

\begin{lemma}\label{lem:fg-2}
  Suppose that $f\greateq g$, the leading coefficient of $f$ is
  $c_f$, and the leading coefficient of $g$ is $c_g$. 

\begin{tabular}{crcllc}
 If& $|c_f| > |c_g|$ & then & $f^2-g^2$ & $\greateq $&$ 2\,fg$\\
 If& $|c_f| < |c_g|$ & then & $f^2-g^2$ & $\lesseq $&$ 2\,fg$ \\
 If& $|c_f| = |c_g|$ & then & $f^2-g^2$ & $\lessgreat $&$ 2\,fg$
\end{tabular}
\end{lemma}
\begin{proof}
  The proof is a standard sign interlacing argument. For example, if
  the degree is $3$, then the diagram below shows where the roots of
  $f+g$ (marked ``+'') and $f-g$ (marked ``--'') are:

\centerline{\xymatrix{
g \ar@{-}[r]^{+} & f\ar@{-}[r]_{-} & g\ar@{-}[r]^{+} & f \ar@{-}[r]_{-} &
 g\ar@{-}[r]^{+} & f
}}

In the third case this is all of the roots. In the first case there is
a root on the far right, and in the second case the remaining root is
on the far left. 
\end{proof}

The following convolution lemma is a non-obvious consequence of sign interlacing.
It will be generalized to a wider class of polynomials in
Lemma~\ref{lem:convolution}. 
\index{convolution}

\begin{lemma} \label{lem:fg-ab}
  Choose $a_1\le \dots \le a_n$, $b_1\le \dots\le b_n$ and define
  $f=\prod(x-a_i)$, $g=\prod(x-b_i)$. The following polynomial has
  all real roots
$$h = \sum_{i=1}^n
\left(\frac{f}{x-a_i}\right)\,\cdot\,\left(\frac{g}{x-b_{n+1-i}}\right)
$$
\end{lemma}
\begin{proof}
  Without loss of generality we may assume that the roots of $fg$ are
  $r_1< \cdots < r_{2n}$. We will show that the sign of $h(r_m)$ only
  depends on $m$. More precisely, we show that
  $sgn(h(r_m))=(-1)^{m+1}$ if $m\le n$, and $(-1)^{m}$ if $m>n$.  This
  shows that there are $2n-2$ intervals between roots where the sign
  at the endpoints is  positive and negative, so $h$ has $2n-2$
  roots, and hence is in $\allpoly$ since its degree is $2n-2$.

  We may assume that $r_m = a_i$. Almost all terms of $h$ are zero
  when evaluated at $r_m$:
  $$ h(r_m) = \left(
    \frac{f}{x-a_i}\,\cdot\,g\,\cdot\,\frac{1}{x-b_{n+1-i}}\right)(a_i)$$
  We investigate each of the three factors in turn. The first factor
  is simply $f^\prime(a_i)$, and the sign of this  is $(-1)^{i+n}$.
  
  Next, since $g(a_i) = \prod_j(a_i-b_j)$, the sign of $g(a_i)$
  depends on the number of $b_j$'s that are greater than $a_i$. Since
  $a_i$ is the $m$-th largest root, there are exactly $m-i$ roots of
  $g$ that are smaller than $a_i$, and hence there are $n-(m-i)$ roots
  of $g$ that are greater than $a_i$. Therefore, $sgn(g(a_i) =
  (-1)^{n-m+i}$.  
  
  Finally, consider the sign $\epsilon$ of $a_i - b_{n+1-i}$. The
  largest $j$ such that $b_j$  is less than $a_i$ is $j=m-i$.
  Consequently, the sign $\epsilon$ is negative whenever
  $(n+1-i)>m-i$, and so $\epsilon = 1$ if $m\le n$ and $-1$ if $m>n$.
  Combining these calculations shows that the sign of $h(r_m)$ is
  $(-1)^{i+n}(-1)^{n-m+i} \epsilon $ which is as claimed.
\end{proof}

\section{A quantitative version of sign interlacing}
\label{sec:sign-quant}

If $f,g$ are two polynomials of the same degree $n$ that sign
interlace, then we know the signs of $g$ evaluated at the $n$ roots of
$f$.  This is similar to interpolation, where we are given the
\emph{values} of $g$, not just the \emph{signs}, at these roots and use this
information to determine $g$. We can combine these two ideas  to
give a quantitative version of sign interlacing.

\index{interpolation}
\index{quantitative sign interlacing}

Let $f$ be a polynomial with distinct roots $\roots(f) =
(a_1,\dots,a_n)$.  The $n$ polynomials
$$
\frac{f(x)}{x-a_1}, \ \cdots, \ \frac{f(x)}{x-a_n}$$
are of degree
$n-1$ and form a basis for the polynomials of degree $n-1$.   If we add
any polynomial of degree $n$ (e.g. $f$ itself) then we get a basis for
the polynomials of degree $n$.

\begin{lemma} \label{lem:sign-quant}
  Assume that $f$ is a polynomial of degree $n$, with positive leading
  coefficient, and with roots $\{a_1,\dots,a_n\}$. Suppose that $g$ is
  a polynomial with positive leading coefficient.

\index{partial fractions}
\begin{enumerate} 
\item If $g$ has degree $n-1$ and we write
    \begin{equation}
      \label{eqn:quant-1}
      g(x) = c_1 \, \frac{f(x)}{x-a_1} + \dots + c_n
      \,\frac{f(x)}{x-a_n}      
  \end{equation}
  then $f \lessless g$ if and only if all $c_i$ are positive.
\item If $g$ has degree $n$ and we write
    \begin{equation}
      \label{eqn:quant-2}
      g(x) = c f(x) + c_1 \, \frac{f(x)}{x-a_1} + \dots + c_n
    \,\frac{f(x)}{x-a_n}
    \end{equation}
    
    then $c$ is positive, and
    \begin{itemize}
    \item $f \greateq g$ if and only if all $c_i$ are positive.
    \item $f \lesseq g$ if and only if all $c_i$ are negative.
    \end{itemize}
  \item If $g$ has degree $n+1$ and we write
    \begin{equation}
      \label{eqn:quant-3}
      g(x) = (d x + c) f(x) + c_1 \, \frac{f(x)}{x-a_1} + \dots + c_n
    \,\frac{f(x)}{x-a_n}
    \end{equation}
    
    then $d$ is positive, and $g \lessless f$ if and only if all $c_i$
    are negative.

  \end{enumerate}
\end{lemma}

\begin{proof}
  Evaluate $g$ (in all three cases) at the roots of $f$:
  $$
  g(a_i) = c_i (a_i-a_1)\cdots(a_i-a_{i-1})
  (a_i-a_{i+1})\cdots(a_i-a_n)$$
  so $sgn\,g(a_i) = sgn(c_i)
  (-1)^{i+n}$.  In the first case we know that $\deg(g)=\deg(f)-1$.
  If all the $c_i$ have the same sign then $g$ sign interlaces
  $f$, and the result follows from Lemma~\ref{lem:int-1}.

  In the second case, if all $c_i$ are
  positive, and $\deg(g)=\deg(f)$, then $g$ sign interlaces $f$. Since
  the sign of $g$ at the largest root $a_n$ of $f$ is the sign of the
  leading coefficient of $f$, and the leading coefficient of $g$ is
  positive, the result follows from Lemma~\ref{lem:int-0}.
  
  In the third case, $sgn\ g(a_n) = sgn\ c_n$ and since the leading
  coefficient of $g$ is positive, we find $g \lessless f$ by Lemma~\ref{lem:int-2}.
  
  Conversely, if $f \lessless g$, then $g$ sign interlaces $f$, and so
  the signs of $g(a_i)$ alternate.  Thus, all $c_i$'s have the same
  sign, which must be positive since the leading coefficient of $g$ is
  positive.  The converses of the second and third cases are similar,
  and omitted. 
\end{proof} 

\begin{remark}
  We can extend Lemma~\ref{lem:sign-quant} by noticing that we do not require the
  coefficients $c_i$ to be constant.   
   If $c_i(x)$ is positive (or negative) on an interval
    containing the roots of $f$, then $g$ sign interlaces $f$.

  One special case of \eqref{eqn:quant-1} is when all the $c_1$ are 
  $1$.  In this case, $g$ is exactly the derivative of $f$.

  If we relax the positivity and negativity conditions of the
  parameters $c_i$'s in Lemma~\ref{lem:sign-quant}, then we still have
  interlacing.  For example, when the degree of $g$ is $n-1$, then $f
  \lesslesseq g$ if and only if all $c_i$ are non-negative.  This case
  where some of the coefficients $c_i$ are zero follows by continuity.
\end{remark}

\begin{cor} \label{cor:xbf}
  If $f\lesslesseq g$ have positive leading coefficients then\\
  $xf-bg\lesslesseq f$ for any positive $b$.  In addition, $xg
  + b f \lesslesseq g$.  If $b$ is negative then it is possible
  that $xf-bg$ and $xg+bf$ are not in $\allpoly$.
\end{cor}
\begin{proof}
  The first is immediate from Lemma~\ref{lem:sign-quant} if we write\\ $g = \sum
  c_i\dfrac{f}{x-a_i}$. In the second case we can write $f = \alpha x
  g + \beta g - r$ where $g\lesslesseq r$. Since $(xg+b f) + \gamma g
  = (\alpha+b) g + (b\beta+\gamma) g -br$ the conclusion follows.
\end{proof}

We can also get interlacings such as $xg+bf  \lesslesseq g$ if we
restrict $g$ to $\allpolypos$.  See Lemma~\ref{lem:xfg}. 


The equation
\begin{equation}
  \label{eqn:sign-quant}
  (-1)^{i+n} sgn\,g(a_i) = sgn\,c_i
\end{equation}
shows that the signs of the $c_i$'s encode the sign of $g$ evaluated
at the roots of $f$. We can easily determine the $c_i$'s in terms of
$g$ and $f$.  For instance, if $f \lessless g$ then $c_i =
g(a_i)/f^\prime(a_i)$.

There is  another condition on the signs of the coefficients in
\eqref{eqn:quant-1} that implies that $g$ has all real roots.

\begin{lemma} \label{lem:sign-1}
  Suppose that the roots of the monic polynomial $f$ are
  $a_1<a_2<\cdots<a_n$. If $c_1,\dots c_n$ in \eqref{eqn:quant-1} has
  exactly one sign change then $g(x)$ in \eqref{eqn:quant-1} has all
  real roots.  If $c_1,\dots,c_n$ has exactly two sign changes,
  $c_1>0$,  and
  $\sum c_i<0$ then $g(x)$ has all real roots.\footnote{The condition
    $\sum c_i<0$ is necessary. If $f(x)=(x+1)(x+2)(x+3)$ then a quick
    computation shows that
    $\frac{f}{x+1}-\frac{f}{x+2}+\frac{f}{x+2}\not\in\allpoly$, and
    $\sum c_i = 1$.}
\end{lemma}
\begin{proof}
  If the $c_i$'s have exactly one sign change, then assume that the
  terms of the sequence $c_1,\dots,c_r$ are positive, and
  $c_{r+1},\dots,c_n$ are negative.  \eqref{eqn:quant-1} shows that
  $g$ sign interlaces $f$ between $a_1,\dots,a_r$ giving $r-1$ roots,
  and also between $a_{r+1}\dots,a_n$ giving $n-r-1$ roots for a total
  of $n-2$ roots.  Since $g$ has degree $n-1$ it follows that $g$ has
  all real roots.
  
  If there are exactly two sign changes then the above argument shows
  that there are at least $n-3$ real roots.  The $c_i$'s are initially
  positive, then negative, and then positive again. The leading
  coefficient of $g$ is $\sum c_i$ which is negative, and since $sgn\
  g(a_n)= sgn\ c_n >0$ there is a root of $g$ in the interval
  $(a_n,\infty)$, which again implies that $g$ has all real roots.
\end{proof}

  Common interlacing puts no restrictions on the signs of the
  coefficients.

  \begin{lemma}
    Any sign pattern can be realized by a pair of polynomials with a
    common interlacing.
  \end{lemma}
  \begin{proof}
    This means that if we are given any sequence $\epsilon_i$ of
    pluses and minuses then we can find two polynomials $g_1,g_2$ with
    a common interlacing such that if $g_2 = \prod(x-r_i)$ and
\[
g_1 = b g_2 + \sum a_i \frac{g_2(x)}{x-r_i}
\]
then $sign(a_i) = \epsilon_i$.

We prove this constructively; the initial cases are trivial. Assume
we are given $g_1,g_2$ as above, and choose $\alpha>\beta$ to both be
smaller than the smallest root of $g_1g_2$. We will show that the sign
sequence of 
\begin{gather*}
  f_1 = (x-\alpha)g_1,\qquad f_2=(x-\beta)g_2\\
\intertext{is}
+, \epsilon_1,\dots,\epsilon_n
\end{gather*}
By construction $f_1$ and $f_2$ have a common interlacing.
Now $a_i = g_1(r_i)/g_2'(r_i)$. Thus the sequence of signs
determined by $f_1,f_2$ is
\[
\frac{f_1(\beta)}{f_2'(\beta)}, \frac{f_1(r_1)}{f_2'(r_1)}, \dots,
\frac{f_1(r_n)}{f_2'(r_n)}
\]
Now
\begin{align*}
  \frac{f_1(\beta)}{f_2'(\beta)} &=
(\beta-\alpha)\frac{g_1(\beta)}{g_2'(\beta)}\\
  \frac{f_1(r_i)}{f_2'(r_i)} &=
  \frac{(r_i-\alpha)g_1(r_i)}{(r_i-\beta)g_2(r_i)+g_2'(r_i)} =
(r_i-\alpha)\frac{f_1(r_i)}{f_2'(r_i)}\\
\end{align*}
Since $\beta$ is smaller than all the roots of $g_1$ and $g_2$, the
signs of $g_1(\beta)$ and $g_2'(\beta)$ are opposite. Also
$r_i-\alpha$ is positive, so the signs are as claimed.

If we choose $\alpha>\beta$ instead, we get a minus sign for the sign,
so all possibilities can happen. 
  \end{proof}

The next result is a simple consequence of quantitative sign
interlacing, and is useful for establishing many interlacing results.
\index{quantitative sign interlacing}

\begin{lemma} \label{lem:fghjk}
  Assume $f,g,j,k\in\allpoly$ have positive leading coefficients, and
$f\lessless g$, $k\greateq f \greateq j$.  Define
$$ h = \ell_1f + \ell_2 j + \ell_3k + q\,g$$
where $\ell_1,\ell_2,\ell_3$ are linear, and $q$ is at most quadratic. The
interlacing $h\lessless f$ holds if these three conditions are true:
\begin{enumerate}
\item The leading coefficients of $\ell_1,\ell_2,\ell_3$ are positive.
\item $q$  and $\ell_2$  are negative on the roots of $f$.
\item $\ell_3$ is positive on the roots of $f$.
\end{enumerate}
\end{lemma}
\begin{proof}
  Assume that $\roots(f)=(a_1,\dots,a_n)$ and write
  \begin{alignat*}{2}
    g &= & \sum \alpha_i \frac{f}{x-a_i} & \\
    j &= \beta f +&  \sum \beta_i \frac{f}{x-a_i} & \\
    k &= \gamma f + & \sum \gamma_i \frac{f}{x-a_i} & 
  \end{alignat*}
  {where $\alpha_i,\beta_i,\beta,\gamma$ are positive, and
  $\gamma_i$ is negative. From the definition }
 $$ h = (\ell_1 + \beta\ell_2 +\gamma\ell_3)f + \sum(\beta_i\ell_2
+\gamma_i\ell_3+\alpha_i q)\frac{f}{x-a_i}. 
$$
The hypotheses guarantee that the leading coefficient of
$\ell_1+\beta\ell_2+\gamma\ell_3$ is positive, and that 
$\beta_i\ell_2
+\gamma_i\ell_3+\alpha_i q$ is negative on the roots of $f$. The
result now follows from Lemma~\ref{lem:sign-quant}.
\end{proof}

We will see consequences of the theorem in
\chapsec{positive}{elem-interlacing}.  As a quick example,
\begin{equation}
  \label{eqn:xffp}
  xf -f -f^\prime\in\allpoly \quad\text{ if } f\in\allpoly
\end{equation}
because $\ell_1=x-1,\ell_2=\ell_2=0,q=-1$ satisfy conditions $1,2,3$.
The hypotheses on $\ell_1,\ell_2,\ell_3$ can be weakened by only
requiring that $\ell_1 +\beta\ell_2+\gamma\ell_3$ has positive leading
coefficient.

  \label{eqn:xffp-1}

If we assume that all roots of $f$ are positive (or negative) we can get more
specific interlacing results.  These interlacing facts will be useful
in showing that various linear transformations preserve interlacing.
See \chapsec{positive}{elem-interlacing}.

\index{Leibnitz}
If polynomials $f,g$ are in the relation $f \longleftarrow g$, then $g$ behaves
like a  derivative of $f$.

\begin{lemma}[Leibnitz Rule] \label{lem:leibnitz}
  Suppose that $f,f_1,g,g_1$ are polynomials with positive leading
  coefficients, and with all real roots. Assume that
  $f$ and $g$ have no common roots.
  \begin{itemize}
  \item If $f \longleftarrow f_1$ and $g \longleftarrow g_1$ then $f_1g,fg_1,fg
    \longleftarrow fg_1 + f_1g\longleftarrow f_1g_1$
  \end{itemize}
  \noindent In particular, $fg_1 + f_1g$ has all real
  roots. Conversely, if $f,g,f_1,g_1$ are polynomials that satisfy $fg
  \longleftarrow f_1g + fg_1$ then
  \begin{enumerate}
  \item $f,g,f_1,g_1\in\allpoly$
  \item $f\longleftarrow f_1$
  \item $g\longleftarrow g_1$
  \end{enumerate}
\end{lemma}

  \begin{proof}
  We only consider the case where the left hand side is $fg$, the
  other cases are similar.
  Let $f = a \prod_1^n(x-a_i)$ , $g = b\prod_1^n(x-b_i)$
  If we write

\begin{align} 
  f_1(x) & = c f(x) + c_1\frac{f(x)}{x-a_1} + \dots +
  c_n\frac{f(x)}{x-a_n}\label{eqn:pf-leib-1}\\
  g_1(x) & = d g(x) + d_1\frac{g(x)}{x-b_1} + \dots +
  d_m\frac{g(x)}{x-b_m}\label{eqn:pf-leib-2}
\end{align}
where all the $c_i$'s and $d_i$'s are positive and $c,d$ are
non-negative then

\begin{multline} \label{eqn:pf-leib-3}
  fg_1+f_1g = (c+d) f(x)g(x) + \\ c_1\frac{f(x)g(x)}{x-a_1} + \dots +
  c_n\frac{f(x)g(x)}{x-a_n} + d_1\frac{f(x)g(x)}{x-b_1} + \dots +
  d_m\frac{f(x)g(x)}{x-b_m}
\end{multline}
and the result follows from Lemma~\ref{lem:sign-quant}.

Conversely, if $fg \longleftarrow fg_1+f_1g$ then we first observe
that $f,g\in\allpoly$.  We can write $f_1$ in terms of $f$ as in
\eqref{eqn:pf-leib-1}, and $g_1$ in terms of $g$ as in
\eqref{eqn:pf-leib-2}, but we do not know the signs of the $c_i$'s and
the $d_i$'s.  However, from \eqref{eqn:pf-leib-3} we know all the
coefficients are negative.  Thus the $c_i$'s and $d_i$'s are negative.
The coefficients $c,d$ are non-negative since all the polynomials have
positive leading coefficients.  Thus $f\longleftarrow f_1$ and
$g\longleftarrow g_1$.
\end{proof}

\begin{cor}\label{cor:leibnitz-linear}
  Suppose that $f\lesslesseq g$ have positive leading coefficients,
  and $a,c$ are positive. 

If $ad\le bc$ then $(ax+b)f+(cx+d)g\in\allpoly$.
  \end{cor}
  \begin{proof}
    Note that $ad\le bc$ is equivalent to $cx+d\greateqeq ax+b$. 
    Apply Lemma~\ref{lem:leibnitz}.
  \end{proof}

\index{Van Vleck polynomials}
\index{Heine-Stieltjes polynomials}
\begin{remark}\label{rem:lame}
  The Lam\'{e} differential equation is
  \begin{equation}
    \label{eqn:lame}
    f(x)y^{\prime\prime} + g(x)y^\prime + h(x)y = 0
  \end{equation}
  Heine and Stieltjes proved that if $g\in\allpoly(p)$, and
  $f\lesslesseq g$ then there are exactly $\binom{n+p-1}{n}$
  polynomials $h(x)$ such that the equation has a solution $y$ which
  is a polynomial of degree $n$. Moreover, both $h$ and $y$ are in
  $\allpoly$. The polynomials $h$ are called Van Vleck polynomials,
  and the corresponding solution is a Heine-Stieltjes polynomial.
  
Although $h$ and  $f$ do not interlace, we do have \cite{shah} that
  $fy^\prime \lesslesseq h y$. To see this, note that $y^\prime\lesslesseq
  y^{\prime\prime}$ and $f\lesslesseq g$, so by Leibnitz we have
$$ f y^\prime \lesslesseq fy^{\prime\prime} + g y^{\prime} = - hy$$
This also shows that if $y\in\allpoly$ then $h\in\allpoly$. 

\end{remark}

\index{cone of interlacing} Recall that a \emph{cone}  is a set
in a vector space that is closed under non-negative linear
combinations. Lemma~\ref{lem:sign-quant} describes all polynomials that
interlace a given polynomial. If we restrict ourselves to polynomials
with positive leading coefficients then this set of polynomials is a
cone. This is also true for polynomials in several variables - see
\chap{topology}.\ref{sec:cone}.

\section{Elementary interlacings} 
\label{sec:elem-interlace}

In this section we use  an algebraic definition of $f\greateqeq g$
to derive many fundamental interlacings. The following is a simple
consequence of Lemma~\ref{lem:sign-quant}.

\begin{cor} \label{cor:rewrite}
  If $f,g$ have positive leading coefficients then $f\longleftarrow g$
  if and only there is a polynomial $h$ with positive leading
  coefficient and non-negative $\alpha$ such that $f\lesslesseq h$ and
  $g = \alpha f+ h$. In addition, if $g\longleftarrow f$ then there is
  a positive $\alpha$ and an $h$ with positive leading coefficient
  satisfying $f\lesslesseq h$ such that $g = \alpha f -h$. If the
  interlacing $f\longleftarrow g$ is strict then $f\lessless h$.
\end{cor}

If $f\longleftarrow g$ then for any $\alpha,\beta$ the linear
combination $\alpha f + \beta g$ interlaces $f$ and $g$.  The exact
direction depends on $\alpha$ and $\beta$. The converse is also true
and of fundamental importance -- see Proposition~\ref{prop:pattern}. This result is
useful, for it gives us a simple way to show that the sum of
certain polynomials with all real roots has all real roots.

\begin{cor} \label{cor:where-roots}
Suppose $f\longleftarrow g$, both $f$ and $g$ are monic, and
$\alpha,\beta$ are non-zero.
\begin{alignat*}{3}
\intertext{If $f\greateqeq g$}
&&\alpha f + \beta g & \longrightarrow & f &
\text{\quad if $\beta$ and $\alpha+\beta$ have the same sign}\\
&&\alpha f + \beta g & \longleftarrow & f &
\text{\quad if $\beta$ and $\alpha+\beta$ have opposite signs or
  $\alpha+\beta=0$}\\
&&\alpha f + \beta g & \longleftarrow & g &
\text{\quad if $\alpha$ and $\alpha+\beta$ have the same sign}\\
&&\alpha f + \beta g & \longrightarrow & g &
\text{\quad if $\alpha$ and $\alpha+\beta$ have opposite signs or $\alpha+\beta=0$}\\
\intertext{If $f\lesslesseq g$ and $\alpha,\beta$ have the same sign}
f & \longleftarrow & \alpha f + \beta g & \longleftarrow & g &
\\
f &\longrightarrow & \alpha f - \beta g & \longleftarrow & g &
\end{alignat*}
\end{cor}

\begin{proof}
  Since $f,g$ are monic we can write $g = f+h$ as in
  Corollary~\ref{cor:rewrite}.  We prove the first two assertions; the
  remainder are similar. Substitution yields $ \alpha f + \beta g =
  (\alpha+\beta)f + \beta h$. Since $(\alpha+\beta)$ and $\beta$ are
  the same sign, the desired interlacing follows since $f\rightarrow
  f$ and $h\rightarrow f$.

  For the second one, write in terms of $g$: $\alpha f+ \beta g =
  (\alpha+\beta)g-\alpha h$. If $\alpha+\beta$ and $\alpha$ have
  opposite signs then $g\rightarrow f$ and $h\rightarrow f$.
  \end{proof}

   
  

We have seen that we can add interlacings. This is a fundamental way
of constructing new polynomials with all real roots.
Here are some variations.

 \begin{lemma}
   \label{lem:add-interlace}
   Suppose $f,g,h$ have positive leading coefficients, and
   $\alpha,\beta$ are positive.
\begin{enumerate}
\item If $f\longleftarrow g$ and $h\longleftarrow g$ then $\alpha
  f+\beta h\longleftarrow g$.
\item If $f\longleftarrow g$ and $f\longleftarrow h$ then $f
  \longleftarrow\alpha g+\beta h$.
\item If $f\longleftarrow g \lessless h$ then $\alpha f -
  \beta h \longleftarrow g$.
\end{enumerate}
 \end{lemma}
 \begin{proof}
   Immediate from the representations of Lemma~\ref{lem:sign-quant} and Corollary~\ref{cor:rewrite}.
 \end{proof}

\index{interlacing square} \index{square!interlacing}
An \emph{interlacing square} is a collection of four polynomials
$f,g,h,k$ that interlace as follows: $\smallsquare{f}{g}{h}{k}$.  If
we have such a square then there are interlacings between certain
products. It is interesting to note that although we do not have the
interlacing $fg \longleftarrow gk$ nor $fg \longleftarrow fh$ it is
true that the sum interlaces: $fg \longleftarrow gk+fh$.

\begin{lemma} \label{lem:int-4}
  If $f,g,h,k$ have positive leading coefficients and
\smallsquare{f}{g}{h}{k} is an interlacing square then 

\centerline{
\xymatrix@-1pc{
fh \ar@{<-}[rd] & & fg & & fk \\
& fk+gh \ar@{->}[ru] \ar@{->}[rd] & & fh+gk \ar@{->}[lu] \ar@{->}[rd] \ar@{->}[ld] \ar@{->}[ru] \\
kg \ar@{<-}[ru] & & hk & & hg
}}

\end{lemma}

\begin{proof} 
  The interlacings involving $fk+gh$ follow from
  Corollary~\ref{cor:where-roots}, as do the right arrows from $fh+gk$. 
  We now show that $fg$ and $gk+fh$ interlace.  Let $f$ and $g$ have roots
  $\{a_1,\dots,a_n\}$ and $\{b_1\le \cdots\le b_n\}$. Interlacing
  follows from these evaluations at the roots of $f$ and $g$:
  \begin{align*}
    sgn(gk+fh)(a_i) & = sgn\,g(a_i)k(a_i) = 1 \\
    sgn(gk+fh)(b_i) & = sgn\,f(b_i)h(b_i) = -1 \\
  \end{align*}
  Since $f\longleftarrow g$  the sign of $gk+fh$ at the largest root $b_n$
  of $fg$ is negative, and so $fg \longleftarrow  gk+fh$.  The remaining
  case is similar.
\end{proof}

\begin{remark} 
  We can use the representation of Corollary~\ref{cor:rewrite} to
  derive a simple relationship between the coefficients of interlacing
  polynomials that will be generalized later. Suppose $f\greateq g$
  have positive leading coefficients and $f = a_n x^n +
  a_{n-1}x^{n-1}+\cdots$, $ g=b_nx^n + b_{n-1}x^{n-1}+\cdots$. We can
  write $g = \alpha f + r$ where $\alpha = b_n/a_n$, and the leading
  coefficient of $r$ is $b_{n-1} - \frac{b_n}{a_n}a_{n-1}$. As long as
  $g$ is not a multiple of $f$ we know that the leading coefficient of
  $r$ is positive. Using the fact that $a_n$ is positive, we conclude
  that $\smalltwodet{a_n}{a_{n-1}}{b_{n}}{b_{n-1}}>0$. In
  Corollary~\ref{cor:log-con-coef} we will show that all such
  determinants of consecutive coefficients are positive.
\end{remark}

Squares don't interlace.

\begin{lemma}\label{lem:int-squares}
\index{interlacing!of squares}
  If $f,g\in\allpoly$, and $f^2$ and $g^2$ interlace then $f$ is a
  constant multiple of $g$.
\end{lemma}
\begin{proof}
  If $g$ has a root that is not a root of $f$ then $g^2$ has two
  consecutive roots that are not separated by a root of $f$. Thus $g$
  divides $f$, and $f$ divides $g$, so they are constant multiples.
\end{proof}

\section{Linear combinations of polynomials}
\label{sec:linear-combinations}

When do all linear combinations of two polynomials with all real roots
have all real roots?  The answer is quite simple: they interlace. This
important result is central to the generalizations of interlacing to
polynomials in several variables that is developed in
Chapter~\ref{cha:p2}, and of polynomials with complex roots.  The
following proposition has been (re)discovered many times, see
\cite{rahman}*{page 210}.

\begin{prop} \label{prop:pattern}
  Suppose relatively prime polynomials $f,g$ have all real roots, all
  distinct, with $\deg(f) \ge \deg(g)$.  Suppose that for all positive
  $\alpha,\beta$ the linear combination $\alpha f + \beta g$ has all
  real roots.  Then either
  \begin{enumerate}
  \item  $\deg(f) = \deg(g)$ and there is a polynomial $k$ such that
    $k\lessless f$, $k\lessless g$.
  \item  $\deg(f) = \deg(g)+1$ and then there is a polynomial $k$ such that
    $k\lessless g$, $k\lesseq f$ or $k\greateq f$. 
  \item  $\deg(f) = \deg(g)+2$ and there is a $k$ so that 
    $f \lessless k \lessless g$.\\[.1cm]

\noindent\hspace*{-1cm}{If  $\alpha f + \beta g$ has all real roots
    for all \emph{positive and negative} $\alpha,\beta$, then
    either} \\[.1cm]
\item  $\deg(f) = \deg(g)$ and either $f \greateq g$ or $g \greateq f$.
\item  $\deg(f) = \deg(g)+1$ and  $f \lessless g$.
  \end{enumerate}
\noindent%
In addition, the converses of all the cases are valid.
\end{prop}
\begin{proof}
  Assume that the roots of $f$a nd $g$ are $a_1\le\dots\le a_n$ and $
  b_1\le\dots\le b_m$. These $n+m$ points determine $n+m+1$ open
  intervals, and $f$ and $g$ have constant sign on the interior of each
  of these intervals.  Enumerate these intervals from the right, so
  that $I_1$ is either $(a_n,\infty)$ or $(b_m,\infty)$.  On adjacent
  intervals, one of $f,g$ has the same 
  sign on both intervals, and one has different sign on the two
  intervals.  
  
  If $f,g$ have leading coefficients with the same sign then they have
  the same sign on $I_1$, and hence on all odd numbered intervals.
  Thus, the roots of $\alpha f + \beta g$ lie in the even numbered
  intervals $I_2,I_4,\dots$. By the continuity of the roots as a
  function of $\alpha,\beta$, there are the same number of roots of
  $\alpha f + \beta g$ in each interval.  By choosing $\alpha=1$ and
  $\beta$ small, we see that each interval has at most one root.
  Thus, since $\alpha f + \beta g$ has $n$ roots, the number of even
  numbered intervals, $\lfloor (n+m+1)/2\rfloor$, is at least $n$.
  This implies that $n=m$ or $n=m+1$.
  
  If we assume that $f,g$ have leading coefficients with opposite
  signs, then the roots lie in the odd intervals of which there are
  $\lfloor (n+m+2)/2\rfloor$.  Consequently, $n=m, n=m+1$ or $n=m+2$.
  
  If both endpoints of an interval are roots of $f$ (say) then $g$ is
  the same sign on the closed interval.  If we choose $\beta$
  sufficiently small, then we can force $f+\beta g$ to have the same
  sign on the interval, and hence there are no roots there.  If we can
  choose $\beta$ to be any sign, then choosing $\beta$ large gives at
  least two roots of $\alpha f + \beta g$ on the interval.
  Consequently, when $\alpha,\beta$ are unrestricted in sign,
  intervals have one endpoint a root of $f$ and the other a root of
  $g$. This implies (4) and (5).

  If  $\deg(f) = \deg(g)$, and  their leading coefficients have
  the same sign then the $n$ roots of $\alpha f + \beta g$ lie in the
  $n$ even intervals.  By the above paragraph each interval has one
  endpoint a root of  $f$ and one a root of $g$. Choose $k$ to have
  one root in each of the odd intervals.
  
  If the leading coefficients differ in sign, then the the roots lie in the $n+1$
  odd intervals.  Since $I_1$ and $I_{2n+1}$ each have only one root
  of $f,g$, and every other interval has one root of $f$ and one root
  of $g$, choose $k$ to have one root in each even interval. 

  If $\deg(f)=\deg(g)+2$ then there are $n$ odd intervals and $n-1$
  even intervals, so all the roots are in the odd intervals. It
  follows that from the right we have a root of f, then a root of $f$
  and  of $g$ in some order, a root of $f$ and of $g$ in some order
  and so on until the leftmost root which is a root of $f$. We choose 
  $k$ to have roots in all the even intervals.

  The case $\deg(f) = \deg(g)+1$ is similar and omitted.

\end{proof}

\begin{figure}[htbp]
  \begin{center}
    \leavevmode
    \epsfig{file=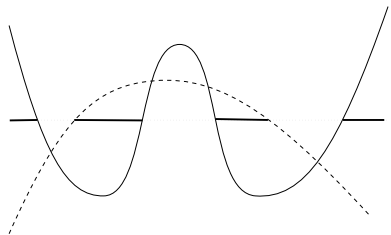,width=2in}
    \caption{Roots of $f+g$ where the degrees of $f,g$ differ by $2$}
    \label{fig:fplus2}
  \end{center}
\end{figure}

The only unfamiliar case is when the degree differs by two.  The thick
lines in Figure~\ref{fig:fplus2} shows where the roots of $f+\beta g$
are located in case $f$ has degree 4 and $g$ degree 2, and $\beta$ is
positive. A similar discussion shows that if $f-\alpha g\in\allpoly$
for all positive $\alpha$, then $f$ and $g$ have a common interlacing.

\index{common interlacing}
\index{interlacing!common}
\begin{cor} \label{cor:af-common}
  Assume that the monic polynomials $f_1,\dots f_n$ all have degree
  $n$, and have all real roots.  The following are equivalent:
  \begin{enumerate}
  \item For all positive $\alpha_1,\dots,\alpha_n$, $\alpha_1 f_1 +
    \dots + \alpha_n f_n$ has all real roots.
  \item $f_1,\dots,f_n$ have a common interlacing.
  \end{enumerate}
\end{cor}

\begin{proof}
  By choosing all but two of the $\alpha$'s sufficiently small, we
  find from Proposition~\ref{prop:pattern} that for all $i\ne j$,
  $f_i,f_j$ has a common interlacing.  It is easy to see that this
  implies that they all have a common interlacing. \seepage{cor:root-zones}
\end{proof}

If we have three polynomials whose linear combinations always have roots,
then the third is a linear combination of the  other two.

\begin{cor} \label{cor:combinations-1}
  Suppose that $f,g,h$ are three distinct polynomials such that for
  all $\alpha,\beta,\gamma$ the polynomial $\alpha f+\beta g+\gamma h$
  has all real roots. Then two of them interlace, and the third is a
  linear combination of the other two.

  Otherwise stated, the largest vector space contained in $\allpoly$ has
  degree $2$.
\end{cor}
\begin{proof}
  If we set one of $\alpha,\beta,\gamma$ to zero then by
  Proposition~\ref{prop:pattern} we know that any two of $f,g,h$ interlace.  After
  renaming, we may assume that $f\greateqeq g \greateqeq h$.  Fix an
  $\alpha$ and consider the two polynomials $k_\alpha = f+\alpha h$
  and $g$.  Since every linear combination of $k_\alpha$ and $g$ has
  all real roots, we know that $k_\alpha$ and $g$ interlace.
  Moreover, $k_0\greateqeq g$ and $ k_\alpha \lesseqeq g$ when
  $\alpha$ is large.  There is an $\alpha$ for which the interlacing
  switches from $\greateqeq$ to $\lesseqeq$; at this value $k_\alpha$
  and $g$ have the same roots.  Thus, for this choice of $\alpha$ the
  linear combination $f+\alpha h$ is a multiple of $g$.
\end{proof}


\index{interlacing!common}
\index{common interlacing}

Next we consider the case of linear combinations where each coefficient is
a polynomial. Notice that the order of interlacing is reversed.

\begin{lemma} \label{lem:pattern2}
  Suppose $h,k$ are polynomials of the same degree with positive
  leading coefficients.  If for every pair of polynomials $f \greateqeq
  g$ with positive leading coefficients the polynomial $fh+gk$ has all
  real roots, then $k \greateqeq  h$. If $k$ and $h$ have different
  degrees then $k\longleftarrow h$.
\end{lemma}

\begin{proof}
  We may assume that $k$ and $h$ have no factors in common.  We will
  only need to assume that the conclusion only holds for \emph{linear}
  $f,g$.  By continuity, $fh+gk$ holds when $f$ and $g$ are multiples
  of one another.  Choosing $f = ax$, $g=bx$ shows that $x(ah+bk)$ has
  all real roots for all positive $a,b$.  From
  Proposition~\ref{prop:pattern} we conclude that $h,k$ have a common
  interlacing.
  
  Next, choose $f=a(x-r)$ and $g=b(x-s)$ where $r$ is larger than all
  roots of $hk$, and $s$ is smaller than all the roots of $hk$. Since
  $f\greateq g$ we know from the hypotheses that all positive linear
  combinations of $fh$ and $gk$ lie in $\allpoly$, and hence $fh$ and
  $gk$ have a common interlacing. Combined with the fact that $h$ and
  $k$ have a common interlacing, it follows that $k\greateqeq h$.

  The proof of the second statement is similar and omitted.
\end{proof}

\section{Basic topology of polynomials}
\label{sec:topology}

In this section we note a few facts about the topology of polynomials
and of interlacing.  The set $V$ of all polynomials of degree $n$ is a
vector space of dimension $n+1$, and is homeomorphic to
$\reals^{n+1}$.  The set $\allpoly(n)$ of all polynomials of degree
$n$ in $\allpoly$ is a subset of $V$, and so has a topological
structure.  We will determine the boundary and interior of
$\allpoly(n)$.

Since roots are continuous functions of the coefficients, if a
polynomial $f$ has all real roots and they are all distinct, then
there is an open neighborhood of $f$ in which all the polynomials have
all real roots, all distinct.  This is half of the next theorem:

\begin{theorem}\label{thm:topology-1} \ 
\begin{itemize}
\item The interior of $\allpoly(n)$ consists of  polynomials with
  all distinct roots.
\item The boundary of $\allpoly(n)$ consists of  polynomials with
  repeated roots.
\end{itemize}
\end{theorem}

\begin{proof}
Suppose that $f$ has a repeated root, so that we can write $f =
(x-a)^2 g$ for some polynomial $g$ in $\allpoly$. For any positive
$\epsilon$ the polynomial $((x-a)^2+\epsilon)g$ has two complex
roots. Since $\epsilon$ can be arbitrarily small, every  neighborhood  of
$f$ has a polynomial that is not in $\allpoly$. Thus, $f$ is in the
boundary of $\allpoly(n)$. Since polynomials with distinct roots are in the
interior of $\allpoly(n)$, this establishes the theorem.
\end{proof}




$\allpoly$  satisfies a closure property. 

\begin{lemma} \label{lem:poly-closure-1}
  If $f$ is a polynomial, and there are polynomials $f_i\in\allpoly$
  such that $\lim f_i = f$ then $f\in\allpoly$.
\end{lemma}
\begin{proof}
  Since the roots of a polynomial are continuous functions of the
  coefficients, if $f$ had a complex root then some $f_i$ would have
  a complex root as well. Thus, all roots of $f$ are real, and $f\in\allpoly$.
\end{proof}

If we allow sequences of polynomials whose degree can increase without
bound then we can get analytic functions such as $e^x$ as the limits
of polynomials in $\allpoly$. The proof of the lemma shows that such
limits will never have complex roots.  Such functions are considered in 
\chap{analytic}.

We can use Proposition~\ref{prop:pattern} to show 
\begin{cor}\emph{Limits preserve interlacing.}

  Suppose that $f_1,f_2,\dots$ and $g_1,g_2,\dots$ are sequences of
  polynomials with all real roots that converge to polynomials $f$ and $g$
  respectively. If $f_n$ and $g_n$ interlace for $n=1,2,\dots$ then
  $f$ and $g$ interlace.
\end{cor}
\begin{proof}
  Since $f_n$ and $g_n$ interlace, $f_n+\alpha g_n\in\allpoly$ for all 
  $\alpha$.  Taking the limit, we see that $f+\alpha g$ has all real
  roots for all $\alpha$. We apply Proposition~\ref{prop:pattern} to conclude that
  $f$ and $g$ interlace.
\end{proof}

The nature of interlacing can change under limits. For instance,  as
$n\rightarrow\infty$ the interlacing
\begin{align*}
  (1-\frac{x}{n})x(1+\frac{x}{n}) & \lessless (x-1)(x+1) \\
\intertext{converges to}
   x & \lessgreat (x-1)(x+1)
\end{align*}

However, if the degrees don't change then the direction is preserved.
\begin{lemma}\label{lem:lace-limit}
  Suppose that $\lim f_i=f$, $\lim g_i = g$ are all polynomials of
  degree $n$. If $f_i\greateqeq g_i$ for all $i$ then $f\greateqeq g$.
\end{lemma}
\begin{proof}
  Write $g_i = \alpha_i f_i + h_i$ where $f_i\lesslesseq h_i$, $h_i$
  has positive leading coefficient and $\alpha_i>0$. Clearly the limits $\lim
  h_i=h$ and $\lim \alpha_i = \alpha$ exist. Since $g = \alpha f + h$
  where $f\lesslesseq h$ and $h$ has non-negative leading coefficient
  we have $f\greateqeq g$. 

  An alternative argument is to consider the location of the roots.
\end{proof}

\section{Linear transformations that preserve interlacing}
\label{sec:poly-linear-trans}

\index{interlacing!preserves}
\index{linear transformation!preserves interlacing}
\index{degree!preserves}
\index{linear transformation!preserves degree}
\index{degree!has constant}
\index{linear transformation!has constant degree}

In this section we establish some properties of a linear
transformation that preserves interlacing.  We begin by observing that if a
linear transformation maps $\allpoly$ to itself then it preserves
interlacing. This is an important result, for it frees us from having
to try to prove that a linear transformation preserves interlacing -
all we need to do is to show that the  linear transformation maps
$\allpoly$ to itself.

\begin{theorem} \label{thm:only-roots}
If  a linear transformation $T$  maps
$\allpoly$ to itself then $T$ preserves interlacing. 

If we also assume  that $T$  preserves degree, then  $f\lesslesseq g$ implies
$Tf\lesslesseq Tg$. In addition, let  $h\greateqeq k$. 
\begin{enumerate}
\item If the leading coefficients of $T(x),T(x^2),T(x^3),\dots$ all
  have the same sign, then  $Th \greateqeq Tk$.
\item If the leading coefficients of $T(x),T(x^2),T(x^3),\dots$ have
  alternating signs, then  $Th \lesseqeq Tk$.
\end{enumerate}
If $T$ maps $\allpoly$ bijectively to $\allpoly$ then all interlacings
$\greateqeq$ above can be replaced with $\greateq$.  More generally,
if $T$ is a $1-1$ open map then all $\greateqeq$ and $\lesslesseq$ can be
\index{open map}
replaced with $\greateq$ and $\lessless$.
\end{theorem}

\begin{proof}
  If $f$ and $g$ are in $\allpoly$ then for all $\alpha$ the polynomial
  $f+\alpha g$ lies in $\allpoly$. Since $T$ maps $\allpoly$ to itself,
  $T(f+\alpha g)=T(f)+\alpha T(g)$ is in $\allpoly$. We now apply Proposition~\ref{prop:pattern}
  to conclude that $Tf$ and $Tg$ interlace.
 
  Now assume $T$ preserves degree. We may assume that $h$ and $k$ have
  positive leading coefficients. We write $k = bh + m$ where $b$ is
  non-negative, $h\lesslesseq m$, and $m$ has positive leading
  coefficient.  Since $Tk = bTh + Tm$ and $Th\lesslesseq Tm$, the
  results (1) and (2) follow from Corollary~\ref{cor:rewrite}.
  
  It remains to show that if $T$ maps $\allpoly$ onto itself and
  preserves degree then $T$ preserves strict interlacing.  $T$ is an
  onto map from the set of polynomials of degree $n$ to itself.
  Choose disjoint opens sets $f\in\mathcal{O}_1, g\in\mathcal{O}_2$ of
  polynomials so that $h\in\mathcal{O}_1$ and $k\in\mathcal{O}_2$
  implies $h\greateqeq k$. Then $T(\mathcal{O}_1)$ and
  $T(\mathcal{O}_2)$ are disjoint open sets whose elements interlace,
  so $Tf\greateq Tg$. Note that this only used $1-1$ and open, so the
  last assertion follows.  
\end{proof}
\index{open map}

\begin{remark}
  The two conditions for the signs of the leading coefficients of the
  sequence $T(x), T(x^2),T(x^3),\dots$ are the only
  possibilities. See Corollary~\ref{cor:leading-coef} below.
\end{remark}

\begin{remark}\index{linear transformation!doesn't preserve interlacing}
    If a transformation merely maps $\allpolypos$ to $\allpoly$, then
    $T$ does not necessarily preserve interlacing. We will see
    (Lemma~\ref{lem:affine-x}) that the linear transformation $T(x^n) =
    (-1)^{n(n+1)/2}\,x^n$ maps
    $\allpolypos\longrightarrow\allpoly$. However, 
    \begin{alignat*}{2}
      T(x+1) &= 1-x & & \quad \text{ has root  } (1)\\
      T(x+1)^2 &= 1-2x-x^2 & &\quad  \text{ has roots  } (-2.4,.4)\\
      T(x+1)^2 &  \not\hspace*{-4pt}\lesslesseq \  T(x+1)
    \end{alignat*}
  \end{remark}

\begin{cor} \label{cor:quant-transform}
  If\, $T$ is a linear transformation such that $T(xf)\lesslesseq T(f)$
  for all $f\in\allpolyint{\sdiffi}$ (where $\diffi$ is an interval
  containing $0$) and maps polynomials with positive leading
  coefficients to polynomials with positive leading coefficients then
  $T$ defines a linear transformation
  $\allpolyint{\sdiffi}\longrightarrow\allpoly$ that preserves
  interlacing.
\end{cor}

\begin{proof}
We first show that if
$h,(x-\alpha)h,(x-\beta)h\in\allpolyint{\sdiffi}$ where $\beta<\alpha$
then $T(x-\alpha)h \greateqeq T(x-\beta)h$. By hypothesis we see that
$T(x-\alpha)h \lesslesseq Th$, and so the conclusion follows from
$$ T(x-\beta)h = T(x-\alpha)h + (\alpha-\beta) T(h)$$

Next, assume $f\greateqeq g$ are in $\allpolyint{\sdiffi}$,   let
$r$ and $s$ be the largest roots of $f$ and $g$, and define $h =
f/(x-r)$. Since we have the interlacings $ (x-r)h \greateqeq g$ and $g
\greateqeq (x-s)h$ it follows from the remark above that we can find a
chain of interlacing polynomials 

\begin{gather*}
 (x-r)h \greateqeq \cdots \greateqeq g \greateqeq \cdots \greateqeq
(x-s)h \\
\intertext{where consecutive polynomials satisfy the conditions of the
  first paragraph. Thus}
 T(x-r)h \greateqeq \cdots \greateqeq T(g) \greateqeq \cdots \greateqeq
T(x-s)h \\
\end{gather*}
Since $(x-r)h \greateq (x-s)h$ it follows that 
$T(x-r)h \greateq T(x-s)h$ and hence the sequence of polynomials is
mutually interlacing. In particular, $T(f) = T(x-\alpha)h \greateqeq T(g)$.
\end{proof}

The next result is fundamental. 

\begin{theorem} \label{thm:rolle-2}
  Differentiation preserves interlacing. 
\end{theorem}
\begin{proof}
  The linear transformation $T(f)=f^\prime$ maps $\allpoly$ to itself,
  and therefore preserves interlacing.  The derivative does not map
  $\allpoly(n)$ onto $\allpoly(n-1)$. (see
  \chap{operators}.\ref{sec:integration-exponential}). However, it's
  easy to see that the derivative is an open map, so strict
  interlacing is preserved.
\end{proof}

  The converse to Theorem~\ref{thm:rolle-2} does not hold\footnote{
If $f =x(x-1)(x-2)$ and $g=(x-.5)^2$ then $f^\prime \lesslesseq
g^\prime$ but it is not true that $f\lesslesseq g$.
}.

The simplest linear transformations mapping $\allpoly$ to itself are
of the form $x^i\mapsto a_i x^i$. Such linear transformations are well
understood, see \chapsec{analytic}{anal-gen-fct}.  Here we determine
their image on $\allpolypos$. 

\begin{lemma} \label{lem:multiplier-image}
  Suppose the linear transformation $T(x^i)=a_ix^i$ maps $\allpoly$ to
  itself. Then
  \begin{enumerate}[(1)]
  \item All non-zero terms of the sequence $a_0,a_1,\dots$ are
    consecutive.
  \item Either all non-zero $a_i$ have the same sign, or the signs
    alternate.
  \item $T$ maps $\allpolypos$ to either $\allpolyneg$ or to
    $\allpolyalt$.
  \end{enumerate}
\end{lemma}
\begin{proof}
  If the non-zero $a_i$'s aren't consecutive and there are at least
  two consecutive zeros then there are $r,s$ such that $s-r>2$, 
  $a_r\ne0$, $a_{r+1}=\cdots=a_{s-1}=0$, $a_s\ne0$. In this case
  $T(x^r(1+x)^{s-r}) = x^r(a_r + a_s x^{s-r})$, and the latter
  polynomial is not in $\allpoly$. If there is an isolated zero, then
  $s-r=2$, and so $a_r + a_{r+2}x^2$ has all real roots. This implies
  that $a_r$ and $a_{r+2}$ have opposite signs. 
  
  Next, $T(x^r(x^2-1)) = (a_{r+2} - a_r x^2)x^r$, so if $a_r$ and
  $a_{r+2}$ are non-zero then they are the same sign. Consequently,
  there are no isolated zeros. Thus, the signs of the first two
  non-zero terms determine the signs of the rest of the terms, so we
  can conclude (1) and (2).  If all the signs are the same then
  $T(\allpolypos)\subseteq \allpolypos$, and if they alternate then
  $T(\allpolyneg)\subseteq \allpolyalt$
\end{proof}

If a linear transformation maps $\allpoly$ to itself, then the leading
coefficients of the images of $x^i$ also determine a linear
transformation that maps $\allpoly$ to itself.

\index{coefficients!leading}

\begin{lemma} \label{lem:leading-coef}
  If $T\colon{}\allpoly\longrightarrow\allpoly$ is a linear transformation
  that preserves degree, and the leading coefficient of $T(x^i)$ is
  $c_i$, then the linear transformation $S(x^i)= c_ix^i$ also
  maps $\allpoly$ to itself.
\end{lemma}
\begin{proof}
  Let $T(x^i) = g_i(x)$, and choose $f=\sum_0^n a_i x^i\in\allpoly$.
  If we substitute $\alpha x$ for $x$, and apply $T$, we find that
  $\sum_0^n a_i \alpha^i g_i(x)\in\allpoly$. Next, substitute
  $x/\alpha$ for $x$ to find that $\sum_0^n a_i \alpha^i
  g_i(x/\alpha)\in\allpoly$. Since $g_i(x)$ is a polynomial of degree
  $i$ we see that $\lim_{\alpha\rightarrow0} \alpha^i g_i(x/\alpha) =
  c_i x^i$, and so 
$$ \lim_{\alpha\rightarrow0} \sum_{i=0}^n a_i \alpha^i g_i(x/\alpha) =
\sum a_i c_i x^i = S(f) \in\allpoly.$$
\end{proof}

\index{coefficients!leading}

\begin{cor} \label{cor:leading-coef}
  If $T\colon{}\allpoly\longrightarrow\allpoly$ is a linear transformation
  that preserves degree then  the leading coefficients of $T(x^i)$
  either all have the same sign, or their signs alternate.
\end{cor}

\begin{proof}
  Use Lemma~\ref{lem:multiplier-image} and Lemma~\ref{lem:leading-coef}.
\end{proof}

  \section{Linear combinations and determinants }
\label{sec:lin-comb-and-det}
\index{determinants!and interlacing squares} \index{interlacing
  square} \index{square!interlacing} Suppose that we have an
interlacing square $\smallsquare{f}{g}{h}{k}$. There are close
connections between the determinant $\smalltwodet{f}{g}{h}{k}$ and the
interlacings \\ 
$\alpha f+\beta g \lessless \alpha h+\beta k$ for all
$\alpha,\beta\in\reals$.

\begin{lemma} \label{lem:inequality-4}
  Assume that $f,g,h,k$ are polynomials with positive leading
  coefficients whose degrees satisfy $deg(f) =
  deg(g)+1=deg(h)+1=deg(k)+2$.  If for all $\alpha,\beta$, not both
  zero, we have that $\alpha f+\beta g \lessless \alpha h+\beta k$
  then
  \begin{enumerate}
\item $\smallsquareLL{f}{g}{h}{k}$  
  \item     $\smalltwodet{f}{g}{h}{k} < 0$ for all $x$.
\end{enumerate}
\index{determinants!and interlacing squares}
Similar results hold if we  replace ``$\lessless$'' by ``$\greateq$''.
\end{lemma}
\begin{proof}
  We know that $h\lesslesseq k$ and $f \lesslesseq g$ by
  Proposition~\ref{prop:pattern}, and the assumptions on degree.
  Setting $\alpha=0$ or $\beta=0$ shows that $f \lessless h$ and $g
  \lessless k$.  If the determinant is zero then there is an $r$ such
  that $f(r)k(r)=g(r)h(r)$.  

  There are two cases. First, assume that $f(r)=0$. Since $h$ and
  $f$ have no roots in common, we know $h(r)\ne0$. If we set $\alpha =
  k(r)$, $\beta=-h(r)$ then for this choice we get the contradiction
  that
$$ (\alpha f + \beta g)(r) = 0 = (\alpha h + \beta k)(r).$$

  The second case is $f(r)\ne0$. The choices $\alpha=-g(r)$,
  $\beta=f(r)$ yield another contradiction.

If we evaluate the determinant at a root $r$ of $k$ we find
$h(r)g(r)\ne0$ Thus $k$ and $h$ have no roots in common, and so 
$h\lessless k$. Similarly, $f \lessless g$. 

Since the determinant is never zero, it always has the same sign. We
determine the sign of the determinant by evaluating it at the largest
root $r$ of $f$.  The result is $-g(r)h(r)$ which is negative since
$f\lessless g$ and $f\lessless h$. 
\end{proof}

There is a converse.

\begin{lemma} \label{lem:inequality-4b}
Suppose $f,g,h,k$ are polynomials in $\allpoly$ with  positive leading 
coefficients. If they satisfy
\begin{enumerate}
\item $f\lessless g$
\item $h\lessless k$
\item $f\lessless h$
\item $\smalltwodet{f}{g}{h}{k}$ is never zero.
\end{enumerate}
then $f+\alpha g \lessless h+\alpha k$ for all $\alpha$. A similar
result holds if $\lessless$ is replaced by $\greateq$.
\end{lemma}
\index{determinants!and interlacing squares}
\begin{proof}
  Assumptions (1) and (2) imply that $f+\alpha g$ and $h+\alpha k$
  have all real roots for all $\alpha$ and that the degree of
  $f+\alpha g$ is always one more than the degree of $h+\alpha k$.
  Assumption (3) implies that the conclusion is true for $\alpha=0$.
  If $(f+\alpha g)(r) = (h+\alpha k)(r)$ then eliminating $\alpha$
  shows that $(fk-gh)(r)=0$. By assumption (4) the determinant is
  never zero, so as $\alpha$ varies $f+\alpha g$ and $h+\alpha k$
  never have a root in common. Consequently, at all times they must
  satisfy $f+\alpha g\lessless h+\alpha k$. The case with $\greateq$
  is similar. 
\end{proof}

\begin{example}
  It is easy to construct examples of interlacing squares whose
  determinants are always negative. We will see more of these determinants
  in Chapter~\ref{cha:affine}. If we apply
  Theorem~\ref{thm:product-4} to the product
$$ f(x,y,z)=(1 + x + y + z) (2 + x + 3 y + z) (5 + x + 2 y + 2 z) $$
then we get the interlacing square

\[
\begin{matrix}
(1 + x) (2 + x) (5 + x)  & \lessless 
 & 19 + 19\, x + 4\, x^2  \\[-.4cm]
\rotatebox{270}{$\lessless$} & & \rotatebox{270}{$\lessless$} \\
29 + 31\, x + 6\, x^2  & \lessless & 
36 + 16 x
\end{matrix}
\]

\noindent%
and the determinant is always negative.  If we apply the theorem to
the determinant
\[
\begin{vmatrix}
   1 + x + y + z & 0 & 0 \\
 0 & 2 + y + 2\,z & 1 \\
 0 & 1 & 1 + 3\,y + z
\end{vmatrix} 
\]

\noindent%
we get an interlacing square where all interlacings are $\greateq$:

\[
\begin{matrix}
1+x & \greateq & 8 + 7x \\[-.4cm]
\rotatebox{270}{$\greateq$} & & \rotatebox{270}{$\greateq$} \\
5+4x & \greateq & 18+7x
\end{matrix}
\]

\noindent%
and the determinant arising from this square is always negative.

  \end{example}

  \begin{remark}
    If we differentiate the entries of the first square
    \begin{align*}
      M & =
      \begin{pmatrix}
(1 + x) (2 + x) (5 + x)    & 19 + 19\, x + 4\, x^2  \\
29 + 31\, x + 6\, x^2  & 36 + 16 x
              \end{pmatrix}\\
M' &=
\begin{pmatrix}
   3 x^2+16 x+17 & 8 x+19 \\
 12 x+31 & 16
\end{pmatrix}
    \end{align*}
then we have that
\begin{enumerate}
\item $|M|>0$ for all $x$.
\item $|M'|>0$ for all $x$.
\item $|M+\alpha M'|>0$ for all $x$ and $\alpha\ge0$.
\end{enumerate}
The last condition suggests that $|M|$ and $|M'|$ ``positively
interlace'' in some sense. Indeed, $M$ and $M'$ arise from polynomials
$f$ and $\frac{\partial f}{\partial x}$ that do interlace, as
polynomials of three variables. The reason why we require $\alpha$ to
be non-negative is that $M-\alpha M'$ can have negative
coefficients, and the result we are using requires that all the
coefficients are positive.
  \end{remark}

The following lemma arises naturally when considering the linear
transformation $f\mapsto f+\alpha f'$ for polynomials with complex
coefficients (Corollary~\ref{cor:cpxpoly-fp}).

\begin{lemma} \label{lem:inequality-1}
  If $f \lessless  g$ (or $f\greateq g$) have positive leading coefficients then \\
  $\smalltwodet{f}{g}{f^\prime}{g^\prime} < 0$.  If $f\lesslesseq g$
  then the determinant is also negative except at common roots of $f$
  and $g$ where it is zero. Thus, if $f\lessless g$ then $f+\alpha f'
  \lessless g+\alpha g'$ for all $\alpha\in\reals$.
\end{lemma}
\index{determinants!involving $f$ and $f^\prime$}

\begin{proof}
We give two proofs. 
Since $f\lessless g$ it follows from Corollary~\ref{cor:rewrite} that
$f+\alpha g$ has all distinct roots for all $\alpha$. Since
differentiation preserves $\lessless$ we can apply
Lemma~\ref{lem:inequality-4} to the interlacing $ f + \alpha g
\lessless f' + \alpha g'$. 

\index{partial fractions}
If we don't use Lemma~\ref{lem:inequality-4} then we can use partial
fractions.  If $f\lessless g$ we  write $g=\sum
b_i\displaystyle\frac{f}{x-a_i}$ where $f(x)=a\prod(x-a_i)$, and all
$b_i$ are positive. Now
$$ \frac{fg^\prime - gf^\prime}{f^2} = \left(\frac{g}{f}\right)^\prime = \left( \sum
  \frac{b_i}{x-a_i}\right)^\prime = - \sum \frac{b_i}{(x-a_i)^2}< 0$$
If $x$ is not a root of $f$ then
$\smalltwodet{f}{g}{f^\prime}{g^\prime}<0$. If $x$ is a root of $f$
then the determinant is equal to $-f^\prime(x)g(x)$ which is negative
since $f^\prime$ and $g$ interlace $f$.

If the interlacing is not strict, then write $f=hf_0$ and $g=hg_0$
where $f_0\lessless g_0$. Since $\smalltwodet{f}{g}{f^\prime}{g^\prime}
= h^2\smalltwodet{f_0}{g_0}{f_0^\prime}{g_0^\prime}$ the result
follows from the first part.

If $f\greateq g$ then there is a positive $\alpha$ and polynomial $h$
with positive leading coefficient such that $g = \alpha f +h$ and
$f\lessless h$. The result now follows from the first case since
$$ \smalltwodet{f}{g}{f^\prime}{g^\prime} = 
\smalltwodet{f}{\alpha f+h}{f^\prime}{\alpha f^\prime+h^\prime} =
\smalltwodet{f}{h}{f^\prime}{h^\prime} 
$$
\end{proof}

An immediate corollary (using $f \lessless f^\prime)$ is a well known
property of derivatives called the Laguerre inequality.
\index{Laguerre inequality}

\begin{cor} \label{cor:inequality-3}
  If $f$ has all real roots, all distinct, then 
$\smalltwodet{f}{f^\prime}{f^\prime}{f^{\prime\prime}} < 0$.
\end{cor}

There is a converse to Lemma~\ref{lem:inequality-1}.

\index{determinants!involving $f$ and $f^\prime$}
\begin{lemma}
  If $f\in\allpoly$, the degree of $g$ is at least one less than the degree the
  of $f$,  and $\smalltwodet{f}{g}{f^\prime}{g^\prime}\ne 0$
for all $x$ then $f$ and $g$ interlace.
\end{lemma}
\begin{proof}
  It suffices to show that for all $\alpha$ the polynomial $f+\alpha
  g$ is in $\allpoly$. When $\alpha=0$ we know that $f\in\allpoly$.
  If there is an $\alpha$ such that $f+\alpha g\notin\allpoly$ then
  there is an $\alpha$ for which $f+\alpha g$ has a double root, $r$.
  In this case we have that
  \begin{align*}
    0 &= f(r) + \alpha g(r) \\
    0 &= f^\prime(r) + \alpha g^\prime(r)
  \end{align*}
  so the determinant $\smalltwodet{f}{g}{f^\prime}{g^\prime}$ is zero
  for this value of $r$.
\end{proof}

If all linear combinations of $f$ and $g$ have all real roots, then
the corollary tells us what direction they interlace.

\begin{cor}\label{cor:inequality-3a}
  If the following three conditions hold then $f\greateqeq g$:
  \begin{enumerate}
  \item $f,g$ have the same degree.
  \item $f+\alpha g\in\allpoly$ for all $\alpha\in\reals$.
  \item $\smalltwodet{f}{g}{f'}{g'}\le0$ for all $x$.
  \end{enumerate}

\end{cor}
The next lemma allows us to conclude that two polynomials interlace
from interlacing information, and knowledge of a determinant.

\begin{lemma} \label{lem:square-1} 
  Suppose that $f\longleftarrow g$, $h\longleftarrow k$ are
  polynomials with positive leading coefficients whose degrees satisfy
  $\deg(f)\ge \deg(h)$ and $\deg(g)\ge \deg(k)$ where all
  ``$\longleftarrow$'' are ``$\greateq$'' or all are ``$\lessless$''.
  The following are equivalent
  \begin{enumerate}
  \item $\smallsquare{f}{g}{h}{k}$
  \item The determinant $\smalltwodet{f}{g}{h}{k}$ is negative at
    roots of $fghk$.
  \end{enumerate}
\end{lemma}
\index{determinants!and interlacing squares}

\begin{proof}
  If $r$ is a root of $g$ then the determinant is $f(r)k(r)$ which is
  negative by the interlacing assumptions. The argument is similar for 
  the roots of $g,h,k$.
  
  Conversely, if the determinant is negative at the $i$-th root $r$ of
  $g$ then we have that $f(r)k(r)<0$.  Since $f$ sign interlaces $g$, so does $k$
  sign interlace $g$.  Similarly we see that $h$ sign interlaces $f$.
  The degree assumptions imply that these sign interlacings are
  actually interlacings.
\end{proof}

A determinant condition gives information about the signs of the
leading coefficients.

\index{determinants!and signs of coefficients}
\begin{lemma} \label{lem:log-concave-coef}
  Assume that $f\longleftarrow h \longleftarrow k$, and 
  that there is a polynomial $g$ such that  $\smalltwodet{f}{g}{h}{k}<0$.
  Then the signs of the leading coefficients of $f$ and $k$ are the
  same.
\end{lemma}
\begin{proof}
  Without loss of generality we may assume that the leading
  coefficient of $h$ is positive.  Suppose that $r$ is the largest
  root of $h$, the leading coefficient of $f$ is $c_f$, and the
  leading coefficient of $k$ is $c_k$. Since the determinant is
  negative we know that  $f(r)k(r)$ is
  negative. Since $f\longleftarrow  h$ the sign of $f(r)$ is $-sgn(c_f)$, and
  similarly the sign of $k(r)$ is $sgn(c_k)$. This implies that $c_f$ and
  $c_k$ have the same sign.
\end{proof} 

\section{Polynomials from determinants}
\label{sec:poly-from-det}

We can construct polynomials from determinants. Although the
constructions yield nothing new in one variable, they are an essential
ingredient for the construction of polynomials in two variables. 
We emphasize the properties of principle submatrices that correspond
to properties of $f(x)/(x-r_i)$ where $f(x) = \prod(x-r_k)$.

We consider determinants of matrices of the form $xI+C$ where $C$ is a
symmetric matrix and $I$ is the identity matrix. More generally we
consider $xD+C$ where $D$ is a positive definite matrix, but these
are really no different. Since $D$ has a positive definite square
root, 

$$ |xD+C| = |D^{-1}|\, |xI + D^{-1/2}CD^{-1/2}| $$
and $D^{-1/2}CD^{-1/2}$ is symmetric.

The determinants $|xI+C|$ generalize the product representation of a
polynomial. If $f(x) = (x-r_1)\cdots(x-r_n)$ then $f(x) = |xI+C|$
where $C$ is the diagonal matrix whose diagonal is
$\{-r_1,\dots,-r_n\}$. 

One important well known fact is that $|xI+C|$ has all real
roots. Indeed, since $C$ is symmetric, its eigenvalues are all real,
and the roots of $|xI+C|$ are the negatives of these eigenvalues. 

If $r$ is a root of $f$ then the analog of $f(x)/(x-r)$ is a
principle submatrix.  If $M$ is a $d$ by $d$ matrix, and
$\alpha\subset\{1,2,\dots,d\}$ then $M\{\alpha\}$ is the (sub)matrix of
$M$ formed by the entries in the rows and columns indexed by $\alpha$.
We let $M[i]$ denote the submatrix resulting from the deletion of the
$i$th row and column of $M$. Such a submatrix is called a
\emph{principle submatrix}.  \index{principle submatrix} 

The maximal factor $f(x)/(x-r_i)$ interlaces $f(x)$. The analogous
statement is also true for principle submatrices.  The following is
well known, but the proof below is elementary.
\begin{theorem} \label{thm:principle-1}
  If $A$ is an $n$ by $n$ Hermitian matrix, and $B$ is any
  principle submatrix of order $n-1$ then their characteristic
  polynomials interlace.
\end{theorem}

\begin{proof}
Choose $\alpha\in\mathbb R$,  partition $
A =\left(\begin{array}{c|c}
B & c \\\hline c^\ast & d
\end{array}\right)
$
and consider the following equation that follows from linearity of the
determinant:

$$
\left|\begin{array}{c|c}
  B - xI  & c \\\hline c^\ast & d-x+\alpha
\end{array}\right|
=
\left|\begin{array}{c|c}
  B - xI  & c \\\hline c^\ast & d-x
\end{array}\right|
+
\left|\begin{array}{c|c}
  B - xI  & c \\\hline 0 & \alpha
\end{array}\right|
$$

Since the matrix on the left hand side is the characteristic
polynomial of a Hermitian matrix, $|A-xI| + \alpha|B-xI|$ has all real
roots for any $\alpha$, and hence the eigenvalues interlace.
\end{proof}

The derivative of $f(x) = \prod_1^n(x-r_n)$ can be written as

\begin{align}
  \frac{d}{dx}\,f(x) &= \frac{f(x)}{x-r_1} + \cdots +
  \frac{f(x)}{x-r_n}\notag \\
  \intertext{There is a similar formula for determinants.  Suppose
    that $C$ is a $n$ by $n$ matrix. }
  \label{eqn:log-der-mat}
  \frac{d}{dx}\, |xI+C| &= |M[1]| + \cdots  + |M[n]|
\end{align}

The factors of a polynomial with all distinct roots have all distinct
roots. The same holds for principle submatrices. This is immediate
since the roots of the submatrix interlace the roots of the original
matrix.

Finally, we have some ways in which principle submatrices don't
generalize properties of maximal factors. 

\begin{example}
  If the roots of a polynomial $f(x)=\prod_1^n(x-r_i)$ are distinct, then
  the polynomials $f(x)/(x-r_i)$ are linearly independent, and span
  the space of all polynomials of degree $n-1$. The matrix
  $\smalltwo{x}{1}{1}{x}$ has two distinct roots $\{1,-1\}$, yet
  the two principle submatrices are equal.  
\end{example}

\begin{example}
  If the roots of $f(x)$ are ordered $r_1\le \cdots\le r_n$ then the
  factors are mutually interlacing:
$$ \frac{f(x)}{x-r_1} \greateqeq \frac{f(x)}{x-r_2} \greateqeq \cdots
\greateqeq \frac{f(x)}{x-r_n}
$$
Although the $n$ principle submatrices have characteristic
polynomials that interlace the characteristic polynomial of $M$, it is
not the case that these submatrices even have interlacing
characteristic polynomials, much less mutually interlacing ones.  They
only have a common interlacing.
\end{example}

We can interpret \eqref{eqn:quant-3} in terms of determinants. This
representation of $g(x)$ will prove useful when we consider
polynomials in two variables.  Since the $c_i$ in \eqref{eqn:quant-3}
are negative the vector
$$v=(\sqrt{-c_1},\dots,\sqrt{-c_n}\,)$$
is real. We assume that $f(x)$
is monic so that $f(x)=\prod(x-a_i)$, and define the diagonal matrix
$A$ whose diagonal entries are $a_1,\dots,a_n$.  Let $I$ be the $n$ by
$n$ identity matrix, and $|W|$ the determinant of a matrix $W$. Expanding by
the first row shows that

\index{determinants!and sign interlacing}
\begin{align}
  \label{eqn:quant-det}
\begin{vmatrix} x+c & v \\ v^t & x I - A\end{vmatrix} &=
\begin{vmatrix}
  x+c & \sqrt{-c_1} & \sqrt{-c_2} & \dots & \sqrt{-c_n} \\
\sqrt{-c_1} & x-a_1 & 0 & \dots & 0 \\
\sqrt{-c_2} & 0 & x-a_2  & \dots & 0 \\
\vdots & \vdots & \vdots & \ddots & \vdots\\
\sqrt{-c_n} & 0 & 0 & \dots & x-a_n
\end{vmatrix} \\[.3cm]
&= ( x + c) f(x) + c_1 \, \frac{f(x)}{x-a_1} + \dots +
c_n \frac{f(x)}{x-a_n}\notag\\
&= g(x)\notag
\end{align} 
Equation \eqref{eqn:quant-det} shows that $g\in\allpoly$ since $g$ is
the characteristic polynomial of the symmetric matrix
$\smalltwo{c}{v}{v^t}{-A}$.
\index{characteristic polynomial}

  \section{A relation weaker than interlacing}
  \label{sec:relation-weaker-than}

  The relation $\weaki$ is weaker than interlacing, yet still has some
  of the properties of interlacing. In particular, it is sometimes
  preserved by differentiation, and there are 
  inequalities for $2$ by $2$ determinants of coefficients.

Consider some simple examples and properties.
\begin{enumerate}
\item $(x+1)^n \weaki x^n$
\item If $h\weaki g$ and $g\weaki f$ then $h\weaki f$.
\item If $g\weaki f$ and $s\le r$ then $(x-s)\,g\weaki (x-r)\,f$.
\item If $g\weaki f$ and $h\in\allpoly$ then $hg\weaki hf$.
\item If $g\weaki f$ and $f\in\allpolypos$ then $g\weaki x^k\,f$ for
  any positive integer $k$. 
\item If $f \greateqeq g$ then $g\weaki f$.
\item If $f\lesslesseq g$ then we generally do  not have $g\weaki
f$ since the smallest root belongs to $f$.
\end{enumerate}

$g\weaki f$ determines a chain of interlacings.

\begin{lemma}\label{lem:weaki-chain}
If $g\weaki f$ have the same degree then there are $h_1,\dots,h_k$ such that

\[
\begin{matrix}
  g &= h_1 & \lesseqeq & h_2 & \lesseqeq & \cdots &
  \lesseqeq &h_k = f \\
  g &= h_1 & \weaki & h_2 & \weaki & \cdots &
  \weaki &h_k = f 
\end{matrix}
\]
If $f$ and $g$ have different degrees then we can write $f = hk$ where
$g\weaki h$, and $h$ has the same degree as $g$.
\end{lemma}
\begin{proof}
  Let $g = \prod_1^m(x-s_i)$ and $f = \prod_1^n(x-r_i)$. The following
  interlacings give the desired sequence of polynomials.
  
\begin{multline*}
g \lesseqeq \left(\prod_2^m (x-s_i)\right)\left(\prod_1^1
  (x-r_i)\right)
\lesseqeq \left(\prod_3^m (x-s_i)\right)\left(\prod_1^2
  (x-r_i)\right) \lesseqeq \cdots \\
\lesseqeq \left(\prod_m^m (x-s_i)\right)\left(\prod_1^{m-1}
  (x-r_i)\right) \lesseqeq \prod_1^{m}(x-r_i) = f
\end{multline*}

For the second part, we take $k$ to be the polynomial whose roots are
the $n-m$ largest roots of $f$. 
\end{proof}

Differentiation does not preserve $\weaki$ if the degrees are
different. For example, if $f=(x+1)(x+2)$ and $g=x(x+1)(x+2)$, then
$f\weaki g$ but $f' \weaki g'$ fails.  We do have 

\begin{lemma}\label{lem:weak-diff}
  If $g\weaki f$ have the same degree then $g'\weaki f'$.
\end{lemma}
\begin{proof}
  The result follows by differentiating the chain of interlacings in
  the previous lemma. 
\end{proof}

\begin{lemma}\label{lem:weak-ineq}
  Suppose that $g\weaki f$, $f$ and $g$ have all positive coefficients, and
$f = \sum_1^n a_i x^i$, $g=\sum _1^m b_i x^i$. Then
\[
\begin{vmatrix}
  b_i & b_j \\ a_i & a_j 
\end{vmatrix} \ge 0\quad \text{for}\ 1\le i < j \le n
\]
\end{lemma}
\begin{proof} 
  We first assume that $f$ and $g$ have the same degree.
  In order to show that $b_ia_j\ge b_ja_i$, we note that this is
  equivalent to $\frac{b_i}{a_i} \ge \frac{b_j}{a_j}$. Thus it
  suffices to take $j=i+1$ since then
\[
\frac{b_i}{a_i} \ge\frac{b_{i+1}}{a_{i+1}} \ge\cdots\ge \frac{b_j}{a_j}. 
\]
  Since $\weaki$ is preserved by differentiation, we can differentiate
  $i$ times, and so assume that $i=0$ and $j=1$. Write
  $f=\prod_1^m(x+r_i)$, $g=\prod_1^m(x+s_i)$ where $0\le r_i\le s_i$
  for $1\le i\le m$. Then
\begin{multline*}
\begin{vmatrix}
  b_0 & b_1 \\ a_0 & a_1 
\end{vmatrix}
=
\begin{vmatrix}
  \prod_1^m s_i & \prod_1^m \left(\sum_1^m \frac{1}{s_i}\right) \\
  \prod_1^m r_i & \prod_1^m \left(\sum_1^m \frac{1}{r_i}\right) \\
\end{vmatrix} = 
\left(\prod s_i\right)\left(\prod r_i\right) \sum_1^m \left(\frac{1}{r_i}
-\frac{1}{s_i}\right)
\end{multline*}

Now assume that the degrees are different. Note that if $p_1 = \sum
\alpha_ix^i$, $p_2 = \sum \beta_ix^i$ and $p_3 = \sum \gamma_ix^i$
then
\[
\begin{vmatrix} \beta_i & \beta_j \\ \alpha_i &
  \alpha_j\end{vmatrix}>0\qquad\text{and}\qquad
\begin{vmatrix} \gamma_i & \gamma_j \\\beta_i &
  \beta_j\end{vmatrix}>0\qquad\implies\qquad
\begin{vmatrix} \alpha_i & \alpha_j \\\beta_i & \beta_j\end{vmatrix}>0
\]
By Lemma~\ref{lem:weaki-chain} it suffices to prove the lemma for $g$
and $(x+\alpha)g$ where $\alpha>0$. If $g=\sum b_ix^i$ then the
determinant is
\[
\begin{vmatrix}
  b_i & b_j\\ \alpha b_i + b_{i-1} & \alpha b_j + b_{j-1}
\end{vmatrix}
=
\begin{vmatrix}
  b_i & b_j \\ b_{i-1} & b_{j-1}
\end{vmatrix}
\]
As above, it suffices to show that
$\smalltwodet{b_i}{b_{i+1}}{b_{i-1}}{b_i}\ge0$, but this is Newton's
  inequality (Theorem~\ref{thm:newton}). 
\end{proof}

A matrix is TP$_2$ is all $2$ by $2$ submatrices have non-negative
determinant. \index{\ TP2@TP$_2$}

\begin{lemma}\label{lem:weak-tp2}
  If $f_1\weaki f_2\weaki \cdots \weaki f_k$ are in $\allpolypos$ then
  the matrix of coefficients of 
  \[
\begin{pmatrix}
f_1\\f_2\\\vdots\\f_k    
\end{pmatrix}
\]
is TP$_2$. 
\end{lemma}
\begin{proof}
  Since transitivity of $\weaki$ implies that $f_i\weaki f_j$ for $i<j$
  the result follows from Lemma~\ref{lem:weak-ineq}.
\end{proof}

Here's a typical corollary:
\begin{cor}\label{cor:weak-ggff}
If $g\weaki f$ are in $\allpolypos$ then  
\[
\begin{vmatrix}
 g \\ .& g \\.&.& g \\
.&.& f \\.&.&.& f \\.&.&.&.& f \\
\end{vmatrix}=
\begin{vmatrix}
  b_0 & b_1 & b_2 &  \cdots \\
0 &   b_0 & b_1 & b_2  & \cdots \\
0 & 0 &   b_0 & b_1 & b_2  & \cdots \\
0 & 0&   a_0 & a_1 & a_2  & \cdots \\
0 & 0 & 0 &   a_0 & a_1 & a_2  & \cdots \\
0 & 0& 0 & 0 &   a_0 & a_1 & a_2  & \cdots \\
\end{vmatrix}
\]
is  TP$_2$.
\end{cor}
\begin{proof}
  The rows are the coefficients of the sequence of polynomials that
  satisfy the relations
\[
g \weaki x g \weaki x^2 g \weaki x^2\,f \weaki x^3\, f \weaki x^4\,f
\]
The result now follows from Lemma~\ref{lem:weak-tp2}.

\end{proof}

\begin{example}\label{ex:roots-ineq}
  Here is a  example that uses $\weaki$ to get an inequality for the
  location of the roots of the derivative. Suppose that $f =
  \prod_1^n(x-r_k)$. Since
    $f \weaki (x-r_k)^k\,(x-r_n)^{n-k}$ we differentiate to find that
\[
   f' \weaki  (x-r_k)^{k-1}\,(x-r_n)^{n-k-1}\,(nx - kr_n - (n-k)r_k)  
\]
It follows that the $k$th root of $f'$ lies in the interval
\[   \bigl[r_k,r_k + \frac{k}{n}(r_n-r_k)\bigr].\]
\end{example}

  \section{Polytopes and polynomials}
  \label{sec:polyt-polyn}

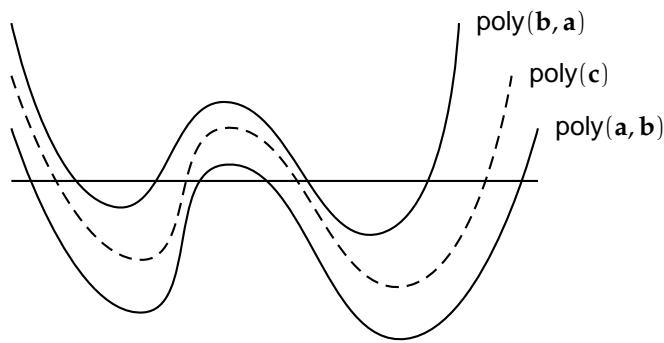
\begin{figure}[htbp]
  \centering
  \begin{pspicture}(0,-3)(11,3)
\psset{xunit=.7cm,yunit=.7cm}
    \psline(0,0)(10,0)
    \pscurve[linestyle=dashed](0,2)(2.5,-1.5)(4,1)(7.5,-2)(9.5,2)
    \pscurve(0,1)(2.5,-2.5)(4,.3)(7.5,-3)(10,1)
    \pscurve(0,3)(2,-.5)(4,1.5)(7,-1)(8.5,3)
    \rput[l](8.7,3){$\ppoly(\bbb,\aaa)$}
    \rput[l](10.2,1){$\ppoly(\aaa,\bbb)$}
    \rput[l](9.7,2){$\ppoly(\ccc)$}
  \end{pspicture}

  \caption{Graphs of interval polynomials}
  \label{fig:poly-interval}
\end{figure}

  A polynomial $a_0+a_1x+ \cdots+a_nx^n$ corresponds to the point
  $(a_0,\dots,a_n)$ in $\reals^{n+1}$. Thus, we can relate geometric
  properties to algebraic properties. In particular, we study two questions:
  \begin{enumerate}
  \item When does a polytope $\mathcal{P}$ in $\reals^{n+1}$ have the property that
    all points in $\mathcal{P}$ determine polynomials with all real roots?
  \item When do polytopes $\mathcal{P}$ and $\mathcal{P'}$ have the
    property that every polynomial determined by points in
    $\mathcal{P}$ interlaces every polynomial determined by points in
    $\mathcal{P'}$?
  \end{enumerate}
  There are simple answers when the polytopes are products of
  intervals, and the proofs are based on sign interlacing. If we have
  vectors $\aaa=(a_i)$ and $\bbb=(b_i)$ then we define
\begin{align*}
 \ppoly(\aaa) &= a_0 + a_1 x+ a_2 x^2 + a_3x^3 + a_4x^4 + \cdots \\
  \ppoly(\aaa,\bbb) &= a_0 + b_1 x+ a_2 x^2 + b_3x^3 + a_4x^4 + \cdots 
\end{align*}
\begin{lemma}\label{lem:polytope-box}
  Suppose that $0\le\aaa\le\bbb$, and that there exists an $f$
  satisfying the interlacing $f\greateqeq \ppoly(\aaa,\bbb)$ and $f\greateqeq
  \ppoly(\bbb,\aaa)$. Then,  for all $\ccc$ satisfying $\aaa\le\ccc\le\bbb$
  we have that  $f\greateqeq \ppoly(\ccc)$.

\end{lemma}
\begin{proof}
If the roots of $f$ are $r_1\le\cdots\le r_n$, then
$\ppoly(\aaa,\bbb)$ and $\ppoly(\bbb,\aaa)$ sign interlace $f$, and so
   the sign of either of them at $r_i$ is $(-1)^{n+i}$. The roots are
   negative since all the coefficients are positive, so the conclusion
  will follow from 
\[ \ppoly(\aaa,\bbb)(-\alpha) \le \ppoly(\ccc)(-\alpha) \le
\ppoly(\bbb,\aaa)(-\alpha)\qquad\text{for}\ \alpha\ge0 
\]
Expanded, this is
\begin{align*}
  a_0 - b_1\,\alpha + a_2\,\alpha^2 - b_3\,\alpha^3 + a_4\,\alpha^4
  \cdots & \le
  c_0 - c_1\,\alpha + c_2\,\alpha^2 - c_3\,\alpha^3 + c_4\,\alpha^4
  \cdots \\ &  \le
  b_0 - a_1\,\alpha + b_2\,\alpha^2 - a_3\,\alpha^3 + b_4\,\alpha^4
  \cdots 
\end{align*}
{which follows from}
\begin{gather*}
  a_{2i}\,\alpha^{2i} \le c_{2i}\,\alpha^{2i} \le b_{2i}\,\alpha^{2i} \\
  -b_{2i+1}\,\alpha^{2i+1} \le -c_{2i+1}\,\alpha^{2i+1} \le -a_{2i+1}\,\alpha^{2i+1}. \\
\end{gather*}
Figure~\ref{fig:poly-interval} shows the configuration of
$\ppoly(\aaa,\bbb)$, $\ppoly(\bbb,\aaa)$ and $\ppoly(\ccc)$. 

\end{proof}

Next, we show that we don't need to assume common interlacing.

\begin{lemma}\label{lem:polytope-box-2}
  Suppose that $0\le\aaa\le\bbb$, and that both $\ppoly(\aaa,\bbb)$
  and $\ppoly(\bbb,\aaa)$ have all real roots. Then
  $\ppoly(\aaa,\bbb)$ and $\ppoly(\bbb,\aaa)$ have a common
  interlacing.
\end{lemma}
\begin{proof}
  We prove this by induction on the degree $n$. We assume that
  $\aaa\ne\bbb$. From the proof of Lemma~\ref{lem:polytope-box} we
  know that the graph of $\ppoly(\bbb,\aaa)$ lies strictly above the
  graph of $\ppoly(\aaa,\bbb)$ for negative $x$. When $n=2$ there are
  two possibilities for the graphs of $\ppoly(\aaa,\bbb)$
  and $\ppoly(\bbb,\aaa)$ (Figure~\ref{fig:ppoly-1})

  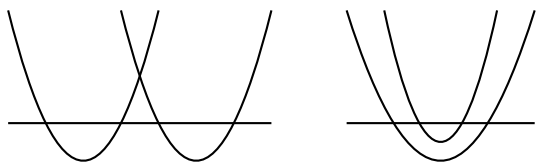
\begin{figure}[htbp]
    \centering
    \begin{pspicture}(0,-1)(14,3)
      \psset{xunit=.5cm,yunit=.5cm}
      \parabola(0,3)(2,-1)
      \parabola(3,3)(5,-1)
      \psline(0,0)(7,0)
      \parabola(9,3)(11.5,-1)
      \parabola(10,3)(11.5,-.5)
      \psline(9,0)(14,0)
    \end{pspicture}
    \caption{The possibilities for a quadratic}
    \label{fig:ppoly-1}
  \end{figure}
Both of the roots of one quadratic can not be to either the left or
the right of the roots of the other quadratic, for then there is an
intersection. Thus, they must be arranged as in the right hand picture,
and there is clearly a common interlacing.

We now consider the general case. The key point is that we can
differentiate. If we define
\[
\aaa' = (a_1,2a_2,3a_3,\cdots,na_{n-1})\qquad 
\bbb' = (b_1,2b_2,3a_3,\cdots,na_{n-1})
\]
then we have the relations
\[
\ppoly(\aaa,\bbb)' = \ppoly(\bbb',\aaa')\qquad 
\ppoly(\bbb,\aaa)' = \ppoly(\aaa',\bbb').
\]
The inductive hypothesis implies that 
$\ppoly(\aaa,\bbb)'$ and $\ppoly(\bbb,\aaa)'$ have a common
interlacing. 

Suppose that $\ppoly(\aaa,\bbb)$ has roots $r_1\cdots<r_n$. Both
$\ppoly(\aaa,\bbb)$ and $\ppoly(\bbb,\aaa)$ are positive on the
intervals
\[
(r_{n-2},r_{n-1}),(r_{n-4},r_{n-3}),\dots
\]
and so the only place that $\ppoly(\bbb,\aaa)$ can have roots is in
the intervals
\[
(r_{n-1},r_{n}),(r_{n-3},r_{n-2})\dots.
\]
We now show that none of these intervals can have more than two roots
of $\ppoly(\bbb,\aaa)$. Figure ~\ref{fig:ppoly-3} is an example of
what can not happen.
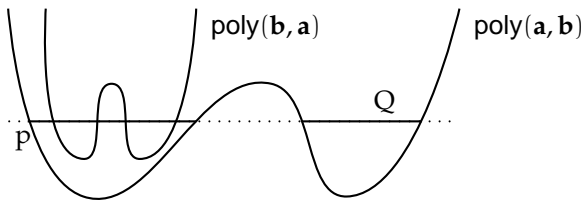
\begin{figure}[htbp]
  \centering
\begin{pspicture}(0,-2)(12,3)
\psset{xunit=.5cm,yunit=.5cm}
  \psline[linestyle=dotted](0,0)(12,0)
  \rput(.4,-.5){p}
  \rput(10,.5){Q}
  \rput[l](12.2,2.5){$\ppoly(\aaa,\bbb)$}
  \rput[l](5.2,2.5){$\ppoly(\bbb,\aaa)$}
  \psline(.52,0)(5,0)
  \psline(7.8,0)(11,0)
  \pscurve(0,3)(2,-2)(7,1)(9,-2)(12,3)
  \pscurve(1,3)(2,-1)(2.75,1)(3.5,-1)(5,3)
\end{pspicture}
    \caption{An impossible configuration}
  \label{fig:ppoly-3}
\end{figure}
We observe the simple property of polynomials $g,h$ with a common
interlacing:
\begin{quote}
  For any $x\in\reals$ the number of roots of $g$ greater than $x$ and
  the number of roots of $h$ greater than $x$  differ by at most one.
\end{quote}

We now see that Figure~\ref{fig:ppoly-3} is impossible, since the two
largest roots of $\ppoly(\aaa,\bbb)'$ are greater than all the roots
of $\ppoly(\bbb,\aaa)'$. Moreover, if there were two roots of
$\ppoly(\bbb,\aaa)$ in the rightmost solid interval labeled $Q$,
then there are $5$ roots of $\ppoly(\bbb,\aaa)'$ to the right of the
point $p$, but only $3$ roots of $\ppoly(\aaa,\bbb)'$ to the right.

Continuing to the left, suppose up to some
point that we find all intervals with $\ppoly(\aaa,\bbb)$ negative
have two roots of $\ppoly(\bbb,\aaa)$. The next interval to the left
can't have $0$ roots, for the same reason as above. Also, it can't
have more than two roots, for then there would be too many roots
of $\ppoly(\bbb,\aaa)'$.
\end{proof}

Combining the last two lemmas we get
\begin{prop}\label{prop:polytope-box}
  Suppose that $0\le\aaa\le\bbb$, and that both $\ppoly(\aaa,\bbb)$
  and $\ppoly(\bbb,\aaa)$ have all real roots. Then, for all $\ccc$
  satisfying $\aaa\le\ccc\le\bbb$ we have that $\ppoly(\ccc)\in\allpoly$.

\end{prop}

\begin{remark}
  If we drop the hypotheses that that the coefficients are positive,
  then there might be no common interlacing, and the conclusion fails.
  If we take $\aaa=(1,-2,1)$, $\bbb=(1,2,1)$ and $\ccc=(1,0,1)$ then
  $\aaa\le \ccc\le\bbb$.  However, there is no common interlacing, and
  \begin{align*}
    \ppoly(\aaa,\bbb) &= 1+2x+x^2 &\in\allpoly\\
    \ppoly(\bbb,\aaa) &= 1-2x+x^2 &\in\allpoly\\
    \ppoly(\ccc) &= 1+x^2 &\not\in\allpoly\\
  \end{align*}
\end{remark}
\begin{cor}
  If $0\le\aaa\le\bbb$, $\ppoly(\aaa,\bbb)$ and $\ppoly(\bbb,\aaa)$
  have all real roots then all points in the product of
  intervals
\[ \aaa\times\bbb=(a_0,b_0)\times(a_1,b_1) \times\cdots \times  (a_{n},b_{n})\]
determine polynomials in $\allpoly$.
\end{cor}

\begin{lemma}\label{lem:polytope-2}
  Suppose that $0\le \aaa \le \ccc \le \bbb$ and $0\le \aaa' \le \ccc'
  \le \bbb'$. If
\[
\begin{matrix} \ppoly(\aaa,\bbb) \\ \ppoly(\bbb,\aaa) \end{matrix}\biggr\} 
\quad\greateqeq\quad
\biggl\{\begin{matrix} \ppoly(\aaa',\bbb') \\ \ppoly(\bbb',\aaa') \end{matrix}
\]
then $\ppoly(\ccc) \greateqeq \ppoly(\ccc')$.
\end{lemma}
\begin{proof}
  Two applications of Lemma~\ref{lem:polytope-box} show that
\[ \ppoly(\aaa,\bbb) \greateqeq \ppoly(\ccc') \qquad \ppoly(\bbb,\aaa)
\greateqeq \ppoly(\ccc')\]
The proof of the lemma holds if we replace $\greateqeq$ by
$\lesseqeq$, so we conclude that 
$ \ppoly(\ccc) \greateqeq \ppoly(\ccc')$.
\end{proof}
The Lemma shows that all points in the box $\aaa\times\bbb$ determine
polynomials that interlace all the polynomials corresponding to points in
$\aaa'\times\bbb'$.

\section{The \fdb\ problem}
\label{sec:fdb-problem}

The \fdb\ formula is an expression for the $m$'th derivative of a
composition. We conjecture that certain polynomials associated with
this formula have all real roots. We are able to establish this in the
first few cases,  for $e^x$, and for $x^d$. If $f$ and $g$ are polynomials then
the \fdb\ formula\cite{wjohnson} is
\begin{multline*}
 \frac{d^m}{dx^m} (f(g(x)) = \\
\sum \frac{m!}{b_1!\cdots b_m!} f^{(k)}(g(x))
\left(\frac{g'(x)}{1!}\right)^{b_1}
\left(\frac{g^{(2)}(x)}{2!}\right)^{b_2}
\cdots
\left(\frac{g^{(m)}(x)}{m!}\right)^{b_m}
\end{multline*}
where the sum is over all ${b_1+2b_2+\cdots+mb_m=n}$ and $b_1+\cdots+b_m=k$.

For example, the first few are

\[
\begin{array}{rl}
\toprule
m &\\
\midrule
 0 & f(g(x)) \\
1 & f'(g(x))\, g'(x)\\
2 & f''(g(x))\, g'(x)^2+f'(g(x))\, g''(x)\\
3 & f^{(3)}(g(x))\, g'(x)^3+3 f''(g(x))\,
   g''(x) g'(x)+f'(g(x))\, g^{(3)}(x)\\
4 & f^{(4)}(g(x))\, g'(x)^4+6 g''(x)\, f^{(3)}(g(x))\, g'(x)^2+ \\ 
 & \  4 f''(g(x))\, g^{(3)}(x)\, g'(x)+3 f''(g(x))\,
 g''(x)^2+f'(g(x))\, g^{(4)}(x)\\
\bottomrule
\end{array}
\]

We can simplify these expressions by writing
\begin{align*}
  \frac{d^m}{dx^m} (f(g(x)) &=
\sum f^{(k)}(g(x)) A_{m,k}(x) \\
\intertext{and defining}
F_{m}(x,y) & = \sum  y^k A_{m,k}(x) 
\end{align*}

Here are the first few values of $F_{m}(x,y)$

\[
\begin{array}{rl}
\toprule
m & F_m(x,y) \\
\midrule
 0 & 1 \\
1 &  g'(x) \,y \\
2 &  g'(x)^2 \,y^2 + g''(x) \,y\\
3 &  g'(x)^3 \,y^3 +3    g''(x) g'(x) \,y^2 + g^{(3)}(x) \,y\\
4 &  g'(x)^4 \,y^4 +6 g''(x)  g'(x)^2 \,y^3 +   
     \bigl[ 4 g^{(3)}(x) g'(x)\,   + 3 g''(x)^2 \bigr] y^2 +   g^{(4)}(x)\,y \\
\bottomrule
\end{array}
\]


Note that $F_m$  only depends on $g$, and not on $f$. If $\alpha$
and $m$ are fixed the transformation $g\mapsto F_m(\alpha,y)$ is \emph{not} a
linear transformation. It follows from \cite{johnson} that $F_m$
satisfies the recurrence

\begin{equation}\label{eqn:fdb-rec}
F_{m+1}(x,y) = \frac{d}{dx}F_m(x,y) + g'(x) \,y\, F_m(x,y)
\end{equation}
If we iterate this we get the simple formula
\[
F_m = \bigl(\diffd_x + g'(x)y\bigr)^m\,(1)
\]

We can recover $ \frac{d^m}{dx^m} (f(g(x))$ from $F_m$ as follows:

\[  \frac{d^m}{dx^m} (f(g(x)) = 
\left(\sum \diffd_y^k A_{m,k}(x)\right) \ 
\left(\sum f^{(k)}(g(x)) \frac{y^k}{k!}\right)
\]
Note that the term on the left is $F_{m}(x,\diffd_y)$ and the term on
the right is the Taylor series of $f(g(x)+y)$.

\begin{conj}\label{conj:fdb-1}
Suppose that $g\in\allpolyf$. For any $\alpha\in\reals$ and
$m=0,1,\dots,$
\begin{enumerate}
\item  $F_m(\alpha,y)$ has all real roots.
\item $F_{m+1}(\alpha,y)\lesslesseq F_m(\alpha,y)$.
\end{enumerate}
  
\end{conj}

We can verify this in the first few cases.
\begin{description}
\item[m=1] We need to show that $g'(\alpha)y\in\allpoly$, which is trivial.
\item[m=2] This is again trivial, since 
\[ F_2(\alpha,y) = y(g''(\alpha) + g(\alpha)^2 y \]
\item[m=3] Factoring out $y$ we have a quadratic whose  discriminant
  is
\[
(3\,g''(\alpha)\,g'(\alpha))^2 - 4 (g'(\alpha))^3 g'''(\alpha) =
g'(\alpha)^2 \bigl[ 9 g''(\alpha) ^2 - 4g'(\alpha)g'''(\alpha)\bigr]
\]
and this is non-negative by Newton's inequalities.

It takes a bit more work to verify interlacing $F_3\lesslesseq
F_2$. Factoring out $y$ we need to show that for all $\alpha$
\[
(g')^3 y^2 + 3 g'' g' y + g''' \lesslesseq (g')^2 y + g'' \]
If we substitute the root of the right hand side into the left hand
side and multiply by $(g')^2$ we get

\[
g'(g'')^2 - 3 (g'')^2 g' + g''' (g')^2 = g' \bigl( g' g''' -
2(g'')^2\bigr)
\]
By Newton's inequalities this has sign opposite from the leading
coefficient of $F_3$, so $F_3\lesslesseq F_2$.
\end{description}

The conjecture holds for $g(x)=e^x$.
\index{Bell polynomial}
\begin{lemma}\label{lem:fdb-ex}
  If $g=e^x$ then Conjecture \ref{conj:fdb-1} holds.
\end{lemma}
\begin{proof}
  The key fact, which follows from \cite{wjohnson}, is that $F_m$ has a
  simple form if $g=e^x$. That is,
\[
F_m(x,y) = B_m(e^xy)
\]
where $B_m$ is the Bell polynomial. We know that $B_m(x)\lesslesseq
B_{m-1}(x)$, and therefore for all $\alpha$ we have that $B_m(\alpha
x)\lesslesseq B_{m-1}(\alpha x)$.
\end{proof}

The conjecture also holds for $x^d$.

\begin{lemma}\label{lem:fdb-xd}
  If $g=x^d$ then Conjecture \ref{conj:fdb-1} holds.
\end{lemma}
\begin{proof}
  To prove the lemma it suffices to show that there are polynomials
  $H_n\in\allpoly$ satisfying
  \[ x^n F_n(x,y) = H_n(x^dy) \] and $H_n\lesslesseq H_{n-1}$.  For
  example, if $d=3$ we compute using \eqref{eqn:fdb-rec} that

\[
\begin{array}{rrcl}
  \toprule
  n &x^n F_n(x,y) && H_n(z) \\
  \midrule
0 & 1 && 1 \\
1 &x( 3x^2y) &&   3z \\
2 & x^2(9 y^2 x^4+6 y x) && 9z^2 + 6z \\
3 & x^3(27 y^3 x^6+54 y^2 x^3+6 y) && 27z^3 + 54 z^2 + 6z \\
4 & x^4(81 y^4 x^8+324 y^3 x^5+180 y^2 x^2) && 81 z^4 + 324z^3+180z^2\\
\bottomrule
\end{array}
\]

We define $H_m$ by the recursion
\begin{align}
  H_0 &= 1 \notag \\
  H_{m+1} &= -m H_m + d\,z\,( H_m + H_m') \label{eqn:fdb-rec-h}
\end{align}

We will prove that $F_m = H_m(x^dy)$ by induction, and it holds for
$m=0$. 

\begin{align*}
  x^{m+1} F_m &= x^{m+1}(F_m' + d x^{d-1} \, y\, F_m) \\
  &= x^{m+1}F_m' + d x^{d+m} \,y\,F_m \\
  H_{m+1}(x^dy) &= -mH_m(x^dy) + dx^dy( H_m(x^dy) + H'(z)[x^dy]) \\
  &= -mH_m(x^dy) + dx^dy H_m(x^dy) + x H_m'(x^dy) \\
  &= -mx^mF_m + d x^dy x^mF_m + x(x^mF_m)'\\
  &= d x^dy x^mF_m + x^{m+1} F_m'
\end{align*}

Since $F_m$ has all positive coefficients $H_m$ has all positive
coefficients, and by induction $H_m\in\allpolypos$. The recursion
\eqref{eqn:fdb-rec-h} shows that $H_m \lesslesseq H_{m-1}$.
\end{proof}

\section{Root Trajectories}
\label{sec:poly-trajectories}

If we are given polynomials $p(x,\alpha)$ where $\alpha$ is a
parameter, then we can investigate the movement of the roots of
$p(x,\alpha)$ as $\alpha$ varies. There is a large body of literature
in the case that $p(x,\alpha)$ is an orthogonal polynomial
\cites{driver01,muldoon}. We limit ourselves to two simple cases.  We
first look at the effect of small perturbations of the roots of a
polynomial on the roots of the derivative, and then consider the root
trajectories of $f+\alpha g$ where $f \longleftarrow g$. See
Figure~\ref{fig:i-trajectory} for trajectories of $f + \imag\, t\, g$
in $\complexes$.

\begin{lemma}
  If $f\in\allpoly$, and we increase the roots of $f$ then the roots
  of $f^\prime$ also increase. \cite{anderson}
\end{lemma}
\begin{proof}
  Without loss of generality we may assume that the roots of $f$ are
  distinct.  If $g$ is the polynomial resulting from increasing the
  roots of $f$ by a sufficiently small amount, then $g\greateq f$.
  Consequently $g^\prime \greateq f^\prime$, so each root of
  $g^\prime$ is larger than the corresponding root of $f^\prime$.
\end{proof}

\begin{lemma} \label{lem:trajectories}
  If $f\longleftarrow g$ have positive leading coefficients then the
  roots of $f+\alpha g$ decrease as $\alpha$ increases.
\end{lemma}
\begin{proof}
  We'll give three proofs!  The simplest proof requires results from
  \chap{poly-matrices} (Corollary~\ref{cor:lin-comb-new}): if $\alpha>\beta$ then
  the matrix $\smalltwo{1}{\beta}{\alpha}{1}$ preserves interlacing,
  and hence $f+\beta g \greateqeq f+\alpha g$.

  Next,  a proof that uses the product representation of $f$. We
  assume that $f,g$ have no roots in common. 
  Since $f\longleftarrow g$ we can write
\begin{align*}
  g &= \beta f + \sum_i b_i \frac{f}{x-a_i}\\
  \intertext{where $\beta\ge0$, the $b_i$ are non-negative, and the
    $a_i$ are the roots of $f$. A root $r$ of $f+\alpha g=0$
    satisfies}
  0 &= f(r)\left(1+\alpha\beta + \alpha \sum_i b_i \frac{1}{r-a_i}\right).\\
  \intertext{ Since $f(r)\ne0$, we remove the factor of $f(r)$, and
    differentiate (remember $r$ is a function of $\alpha$)} 0 &=
  \beta+\sum_i b_i \frac{1}{r-a_i} - \alpha \sum_i
  \frac{b_i}{(r-a_i)^2}\frac{dr}{d\alpha}\\
  \intertext{and solving for $r^\prime$ shows that $r^\prime$ is
    negative} r^\prime(\alpha) &=
  -\left(\alpha^2\sum\frac{b_i}{(r-a_i)^2}\right)^{-1}
\end{align*}

Finally, a conceptual proof. We assume that $\alpha$ is
non-negative; the case where it is non-positive is similar. When
$\alpha=0$ the roots of $f+\alpha g$ are the roots of $g$. As
$\alpha\rightarrow\infty$ the roots go the roots of $g$.  The
roots are continuous functions of $\alpha$. Since every
real number $r$ is the root of at most one equation $f+\alpha g=0$
(choose $\alpha = -f(r)/g(r)$) it follows that the roots start at
roots of $f$ for $\alpha=0$, decrease to the roots of $g$, and never
back up.
\end{proof}

\section{Recursive sequences of polynomials}
\label{sec:log-concave}

We give two different ways of constructing sequences of interlacing
polynomials. The first is to use a recurrence; the second is the use
of log concavity.  
 The interlacing constructions of previous sections can be iterated to
create infinite sequences of interlacing polynomials.  

\begin{lemma} \label{lem:recursive}
  Suppose $\{a_i\}$, $\{b_i\}$, $\{c_i\}$ are sequences where all
  $a_i$ and $c_i$ are positive, and the $b_i$ are
  unrestricted. Define a sequence of polynomials recursively by 
  \begin{align*}
    p_{-1} &= 0 \\
    p_0 & = 1 \\
    p_n & = (a_n x + b_n)p_{n-1} - c_n p_{n-2} \quad\quad\text{for n
    $\ge 1$}
  \end{align*} 
  These polynomials form an infinite interlacing sequence:
$$ p_1 \lessgreat p_2 \lessgreat p_3 \lessgreat \cdots $$ 
\end{lemma}
\begin{proof}
  We first note that $p_2$ evaluated at the root of $p_1$ is $-c_2$
  which is negative. This shows that $p_2\lessless p_1$; the general
  case now follows by induction from Lemma~\ref{lem:fghjk}.
\end{proof}

It is important to note that the sequences of Lemma~\ref{lem:recursive} are
exactly the sequences of orthogonal polynomials; see
\cite{szego}*{page 106}.  These polynomials have the property that the
degree of $p_n$ is $n$; we can also construct sequences of interlacing
polynomials that all have the same degree.

We can sometimes replace the constant $c_n$ by a quadratic:

\begin{lemma}\label{lem:recursive-q}
  Suppose that $f_n\in\allpolypos$ is a sequence of polynomials
  satisfying the recurrence
\[
f_{n+1}(x) = \ell_n(x) \,f_n(x) - q_n(x) \,f_{n-1}(x)
\]
where
\begin{enumerate}
\item $\ell_n$ is linear with positive leading coefficient.
\item $q_n$ is quadratic with positive leading coefficient.
\item $q_n$ is positive at the roots of $f_n$.
\end{enumerate}
then $f_{n+1} \lessless f_n$.
\end{lemma}
\begin{proof}
  This follows from Lemma~\ref{lem:fghjk}.
\end{proof}
Here's an example where it is easy to verify the third condition.
\begin{cor}\label{cor:orthog-quad}
  Suppose that $f_n\in\allpolypos$ is a sequence of polynomials
  satisfying the recurrence
\[
f_{n+1}(x) = (a_nx+b_n) \,f_n(x) + x(c_nx + d_n) \,f_{n-1}(x)
\]
where $c_n\le0$ and $d_n\ge0$. Then $f_{n+1} \lessless f_n$.
  
\end{cor}
\begin{proof}
  Since $f_n\in\allpolypos$ its roots are negative, and when $r<0$ we
  have that 
  $r(rc_n+d_n)<0$. Now apply the lemma.
\end{proof}

\begin{example}

\index{polynomials!Narayana}\index{Narayana polynomials}

The Narayana polynomials $N_n$ are defined by the formula
\[
N_n = \sum_{i=0}^n \frac{1}{n} \binom{n}{k}\binom{n}{k-1}x^k.
\]
 and satisfy the recurrence
\[
N_n = \frac{1}{n+1} \biggl( (2n-1)(x+1)\, N_{n-1} - (n-2) (x-1)^2 N_{n-2}\biggr).
\]
It follows from Lemma~\ref{lem:recursive-q} that $N_n\lessless N_{n-1}$. 

  \end{example}

\begin{example}\label{ex:comm-narayan}
\index{polynomials!Narayana}\index{Narayana polynomials}
\index{Commutator}
    We show that the sequence of commutators determined by the
    Narayana polynomials is an interlacing sequence. Let $T(x^n) =
    N_n$, and define
\[
f_n = [\diffd,T] (x^n) = N_n' - nN_{n-1}
\]
The $f_n's$ satisfy the recurrence relation
\begin{align*}
f_n & = \frac{2n-1}{(n+1)(n-1)^2}\, \bigl( 
n^2-n-1 +\left(n^2-n+1\right) x\bigr)
 \,f_{n-1}  \\
& -\frac{(n-2) n^2 }{(n-1)^2 (n+1)}\,(x-1)^2\, f_{n-2}
\end{align*}
Now $f_3 = 1+3x$, so using induction  and the recurrence we conclude that
$f_n(1)\ne0$. 
It follows from Lemma~\ref{lem:recursive-q} that 
\[ \cdots \lessless f_n \lessless f_{n-1} \lessless f_{n-2} \cdots \]
\end{example}

\begin{example}
  In \cite{little} they considered the recurrence relation
\[
p_{n+1} = 2 x \,p_n - \bigl(x^2 + (2n-1)^2\bigr)\,p_{n-1}.
\]
This satisfies all the hypotheses of Lemma~\ref{lem:recursive-q}, so
we conclude that all $p_n$ have all real roots, and consecutive ones
interlace. 
\end{example}

A small modification to the recurrence in Lemma~\ref{lem:recursive}
determines a series of polynomials that are known as 
orthogonal Laurent polynomials\index{orthogonal polynomials!Laurent}
\index{Laurent orthogonal polynomials}.\cite{dimitrov-ranga}

\begin{lemma}\label{lem:recursive-laurent}
If all coefficients $a_n,b_n,c_n$ are positive, $a_n>c_n$, and
  \begin{align*}
    p_{-1} &= 0 \\
    p_0 & = 1 \\
    p_n & = (a_n x + b_n)p_{n-1} - c_n x\, p_{n-2} \quad\quad\text{for n
    $\ge 1$}
  \end{align*} 
then $p_n\in\allpolyalt$, and 
$ p_1 \lessgreat p_2 \lessgreat p_3 \lessgreat \cdots $  
\end{lemma}
\begin{proof}
  We proceed by induction. The case $p_2\lessless p_1$ is easy, so
  assume that $p_{n-1} \lessless p_{n-2}$. The hypotheses on the
  coefficients implies that the coefficients of all $p_{k}$ alternate
  in sign, so $p_{n-1}$ and $p_{n-2}$ have all positive
  roots. Consequently 
  \begin{gather*}
    a_n x\, p_{n-1} \lesslesseq p_{n-1} \greateq c_nx\, p_{n-2}\\
\intertext{and so}
   p_n = a_nx \,p_{n-1} - c_nx\,p_{n-2} \lessless p_{n-1}
  \end{gather*}
Note that from $f\lessless g\greateq h$ we can not always conclude
that $f-h\lessless g$. We can do so here because the leading
coefficient of $a_nxp_{n-1}$ is greater than the leading coefficient of
$c_nxp_{n-2}$ by the hypothesis that $a_n>c_n$. 

\end{proof}

\begin{remark}
  The reverse of a polynomial $f(x)$ of degree $n$ is $x^nf(1/x)$. 
  The reverse of a (Laurent) orthogonal polynomials also satisfies a
  simple recurrence. For orthogonal polynomials as in
  Lemma~\ref{lem:recursive} the recurrence is
  \begin{align*}
    p_n^{rev} &= (a_n + b_n x)\,p_{n-1}^{rev} - c_n\,x^2\,p_{n-2}^{rev} \\
\intertext{and for Laurent polynomials as in
  Lemma~\ref{lem:recursive-laurent}} 
p_n^{rev} &= (a_n + b_n x)\,p_{n-1}^{rev} - c_n\, x\, p_{n-2}^{rev}
  \end{align*}
\end{remark}

If we modify the signs of coefficients in Lemma~\ref{lem:recursive}
then we still know the location of the roots. See page \pageref{eqn:pq} for
more results.

\index{purely imaginary roots}
\begin{lemma} \label{lem:recursive-i}
  Suppose $\{a_i\}$,  $\{c_i\}$ are sequences of positive terms. 
  Define a sequence of polynomials recursively by 
  \begin{align*}
    p_{-1} &= 0 \\
    p_0 & = 1 \\
    p_n & = a_n x \,p_{n-1} + c_n p_{n-2} \quad\quad\text{for n
    $\ge 1$}
  \end{align*} 
All the roots of $p_n$ are purely imaginary.
\end{lemma}
\begin{proof}
It is easy to see that $p_n$ only has  terms of degree with the same
parity as $n$. Consequently, if we define $q_n(x) = \imag^n p_n(\imag
x)$ then $q_n$ has all real coefficients, and satisfies the recurrence
$$ q_n = a_n\, x\, q_{n-1} - c_n\,q_{n-2}$$
Thus $q_n$ has all real roots, and so $p_n$ has purely imaginary roots.
\end{proof}

\begin{example}
    The $d$-dimensional \index{cross polytope}cross polytope is 
 $$\mathcal{O}_d = \{(x_1,\dots,x_d)\in\reals^d\mid\ |x_1|+\cdots+|x_d|\le n\}$$
The number of integer points $p_d$ in ${\mathcal O}_d$ having exactly $k$
non-zero  entries is $2^k\binom{d}{k}\binom{n}{k}$ since there are 
\begin{itemize}
\item $2^k$ choices of sign for the non-zero entries.
\item $\binom{d}{k}$ choices for the non-zero coordinates
\item $\binom{n}{k}$ positive integer solutions to $x_1+\cdots x_k\le n$
\end{itemize}
and therefore $p_d = \sum_{k=0}^n 2^k\binom{d}{k}\binom{n}{k}$.  $p_d$
is known as the \index{Ehrhart polynomial}Ehrhart polynomial of
${\mathcal O}_d$. If we define $q_d = \sum_{k=0}^n
2^k\binom{d}{k}\binom{x-1/2}{k}$ then we can verify that $q_d$
satisfies the recurrence relation $$
q_d = \frac{2}{d}x q_{d-1} +
\frac{d-1}{d} q_{d-2}$$
and thus $q_d$ has all imaginary roots. It
follows from Lemma~\ref{lem:recursive-i} that the real part of all the
roots of $p_d$ is $-1/2$. (See \cite{bump}.)
  \end{example}

\begin{example}
  Assume that $f_0$ has all positive roots and define a sequence by
  the recurrence $f_{n+1} = f_n + x f_n^\prime$.  All these
  polynomials have the same degree, all real roots, and
  \eqref{eqn:int-3} is satisfied. (See
  Lemma~\ref{lem:inequality-6a}.)
\end{example}

We now briefly discuss log-concavity.

\index{log concave}
\index{log concave!strict}

A function $f$ is \textit{strictly concave} if for every choice of
distinct $a,b$ in the domain of $f$ we have
$$\frac{f(a)+f(b)}{2} < f\left(\frac{a+b}{2}\right).$$
We say that $f$ is
\textit{strictly log concave} if its logarithm is strictly concave.
This is equivalent to
$$ f(a)f(b) < f\left(\frac{a+b}{2}\right)^2.$$ If we restrict $f$
to integer values (i.e., a sequence), then at the points \\ $i-1,i,i+1$ a
strictly log concave function $g$ satisfies
$$ g(i-1)\,g(i+1) < g(i)^2.$$ A sequence of functions $f_1,f_2,\dots$ is
strictly log concave if the sequence determined by evaluating at any
point is strictly log concave.  This is equivalent to
$$ \text{For all } x, \quad\quad f_{n-1}(x)f_{n+1}(x) < f_n(x)^2.$$



We can extend Lemma~\ref{lem:interlace-and-log-concave}  to chains of
interlacings.

\begin{cor} \label{cor:int-3}
  Assume that $f_1,f_2,\dots$ is a sequence of polynomials with
  positive leading coefficients satisfying
  \begin{enumerate}
  \item The sequence is strictly log concave.
  \item The leading coefficients  have the same sign, or
    alternating signs.
  \item $f_1 \longleftarrow f_2$
  \end{enumerate}
  then 
  \begin{equation}
    \label{eqn:int-3}
  f_1 \longleftarrow f_2 \longleftarrow f_3 \longleftarrow f_4
  \longleftarrow \cdots.    
  \end{equation}
\end{cor}

The property of log-concavity is not much stronger than interlacing.
\begin{lemma} \label{lem:make-log-concave}
  Suppose $f,g,h$ are polynomials of the same degree with positive
  leading coefficients and $f\greateq g$.  Then,  $g\greateq h$ if and
  only if there is a positive constant $c$ such that $f,g,ch$ is
   strictly log-concave.
\end{lemma}

\begin{proof}
  If $f,g,ch$ is strictly log-concave then $g\greateq h$ by
  Corollary~\ref{cor:int-3}.  Conversely, if $g\greateq h$ then there is a $c$
  such that $f,g,ch$ is log-concave if and only if $fh/g^2$ is bounded
  above.  Since $f,g,h$ have the same degree, the limit of $fh/g^2$ is
  finite as $x$ goes to $\pm\infty$.  Consequently, it suffices to
  show that $fh/g^2$ is bounded above at each root of $g$.  Now $g^2$
  is positive near a root $a$ of $g$, so in order for $fh/g^2$ to be
  bounded above near $a$ it must be the case that $f(a)h(a) < 0$.  If
  the roots of $g$ are $\{a_1<\cdots<a_n\}$, then $sgn\, f(a_i) =
  (-1)^{n+i}$, and $sgn\, h(a_i) = (-1)^{n+i+1}$, so $sgn\,
  f(a_i)h(a_i) = -1$. We conclude that $fh/g^2$ is bounded above.
\end{proof}

\section{The positive and negative parts}

We consider properties of the positive and negative parts of a
polynomial, and apply the properties to a recurrence.

If $f\in\allpoly$, then $f^{neg}$ is the monic
polynomial whose roots are the negative roots of $f$, and $f^{pos}$
the monic polynomial whose roots are the positive roots of $f$.  For
example, if we know that $f\greateq g$ then exactly one of these
possibilities holds:
\begin{alignat*}{3}
  f^{neg} & \greateq g^{neg} &\  &\And\ & f^{pos} & \greateq g^{pos} \\
  f^{neg} & \lessgreat g^{neg} &\  &\And\ & f^{pos} & \lessless g^{pos} \\
\end{alignat*}
Similarly, if we know that $f\lessless g$ then there are exactly two
possibilities:
\begin{alignat*}{3}
  f^{neg} & \lesseq g^{neg} &\  &\And\ & f^{pos} & \lessless g^{pos} \\
  f^{neg} & \lessless g^{neg} &\  &\And\ & f^{pos} & \greateq g^{pos} \\
\end{alignat*}

\index{negative part of a polynomial} \index{polynomial!negative part}
\index{positive part of a polynomial} \index{polynomial!positive part}
\index{f$^{neg}$}
\index{f$^{pos}$}

\begin{lemma}\label{lem:rec-pm}
  The notation {\framebox{$\lesseq\ \lessless$}} describes the
  interlacing of the positive and negative parts of a pair of
  polynomials. If the polynomials are $f,g$, then this example means
  that $f^{neg}\lesseq g^{neg}$ and $f^{pos}\lessless g^{pos}$. An
  arrow between boxes means the interlacing of $f,g$ is described in
  the first box, and the interlacing of $f+xg,f$ is described in the
  target box. Then, we have

\centerline{
\xymatrix{
{\framebox{$\lesseq\ \lessless$}} \ar@{->}[r] & 
{\framebox{$\greateq \  \lesseq$}} \ar@{<->}[r] & 
{\framebox{$\lessless \  \greateq$}}\\
{\framebox{$\lessgreat\ \lessless$}} \ar@{->}[r] & 
{\framebox{$\lessless \  \lesseq$}} \ar@{<->}[r] & 
{\framebox{$\greateq \  \greateq$}}\\
}}

\end{lemma}
\begin{proof}
  There are six assertions to be verified. We describe in detail how
  to prove one of them, namely

\centerline{
\xymatrix{
{\framebox{$\lessless \  \lesseq$}} \ar@{<-}[r] & 
{\framebox{$\greateq \  \greateq$}}.
}}

The diagram to be proved is the following statement:
$$
f^{neg}\greateq g^{neg} \And  f^{pos}\greateq g^{pos}
\implies
h^{neg}\lessless f^{neg} \And h^{pos} \lesseq f^{pos} 
$$
where $h = f + xg$. We use a sign interlacing argument, which we
illustrate with an example where $f$ has two positive and three
negative roots. These assumptions imply the arrangement of the roots
in the diagram below. We can determine the sign of $g$ and $f+xg$ at
the roots of $f$; they are given in the second and third lines of the
diagram. There is a root in the leftmost segment since $f$ and $f+xg$
have asymptotes of opposite sign as $x\rightarrow-\infty$. It follows that
there is a root of $f+xg$ in every segment of the first line of the
diagram. This tells us where all the roots of $f+xg$ are, and implies
the conclusion.

\centerline{
\xymatrix@1@R=1pt{
 \ar@{-}[r]^g & f \ar@{-}[r]^g & f \ar@{-}[r]^g & f \ar@{-}[r] & 0 \ar@{-}[r]^g & 
f \ar@{-}[r]^g & f  & \\
&+ &- &+ &+ &- &+ & {\text{sign of $g$}}\\
&- &+ &- &+ &- &+ & {\text{sign of $f+xg$}}\\
}}

\end{proof}

We can apply this result to a recurrence.

\index{recurrence}
\begin{cor}
  Consider a sequence satisfying the recursion $p_k = p_{k-1} +x
  p_{k-2}$. If $\{p_1,p_2\}$ satisfy any of the six interlacing
  conditions of Lemma~\ref{lem:rec-pm} then all $p_k$ are in
  $\allpoly$. 
\end{cor}

The recurrence $q_n = x q_{n-1} + q_{n-2}$ looks like the recurrence
for orthogonal polynomials, but the sign is wrong. We show that if
$q_0=1,q_1=x$ then the roots are all on the imaginary axis, and the
roots of consecutive terms interlace. (See Lemma~\ref{lem:recursive-i}.)

By induction we see that $q_{2n}$ only has even powers of $x$, and
$q_{2n+1}$ only has odd powers. Write
$$ q_{2n}(x) = p_{2n}(x^2) \quad\quad q_{2n+1}(x) =x\,p_{2n+1}(x^2).$$
Here are the first few values:

$$
\begin{array}{ccc}
k & q_k & p_k \\[.1cm]
0 & 1 & 1 \\
1 & x & 1\\
2 & 1+x^2 & 1+x \\
3 & 2x+x^3 & 2+x \\
4 & 1+3x^2+x^4 & 1+3x+x^2
\end{array}
$$
The recurrence for $q$ translates into  two recurrences for $p$:
\begin{align*}
  p_{2k+1}(x) & = p_{2k}(x) +  p_{2k-1}(x)\\
  p_{2k}(x) & = x p_{2k-1}(x) +  p_{2k-2}(x)\\
\intertext{and as matrix equations these are}
\begin{pmatrix}
  1&1\\1&0
\end{pmatrix}
\begin{pmatrix}
  p_{2k}\\p_{2k-1}
\end{pmatrix}
&=
\begin{pmatrix}
  p_{2k+1}\\p_{2k}
\end{pmatrix}\\
\begin{pmatrix}
  x&1\\1&0
\end{pmatrix}
\begin{pmatrix}
  p_{2k+1}\\p_{2k}
\end{pmatrix}
&=
\begin{pmatrix}
  p_{2k+2}\\p_{2k+1}
\end{pmatrix}
\end{align*}

We claim that we have the interlacings
\begin{equation}
  \label{eqn:pq}
  p_{2k} \greateq p_{2k+1} \lessgreat p_{2k+2}
\end{equation}

The first is trivial since $p_{2k}\lessless p_{2k-1}$. The second
follows from Lemma~\ref{lem:sign-quant}. In particular, all $p_k$ are
in $\allpolypos$.

Now $q_{2n}(x)=p_{2n}(x^2)$, so all roots of $q_{2n}$ lie on the
imaginary axis. The roots of $q_{2n+1}$ also lie on the imaginary
axis, and include $0$. The interlacings \eqref{eqn:pq} imply that the
roots of the $q_k$ interlace.

If we had that $q_n = x q_{n-1} + \alpha_n\,q_{n-2}$ where $\alpha_n$
is positive then the same argument applies, and the roots interlace on
the imaginary axis.

\section{Polynomials without all real roots}
\label{sec:poly-no-roots}

We have seen many ways of showing that polynomials have all real
roots; here we note that certain polynomials do \emph{not} have all real
roots.

\begin{enumerate}
\item $p_n = 1+x+\frac{x^2}{2}+\dots+\frac{x^n}{n!}$ does not have all
  real roots if $n>1$. Since $p_n^\prime = p_{n-1}$, if $p_n$ had all
  real roots, so would $p_{n-1}$.  The quadratic formula shows that
  $p_2$ does not have all real roots, so $p_n$ does not either.

\item $q_{2n} = 1+\frac{x^2}{2}+\dots+\frac{x^{2n}}{(2n)!}$ and 
$q_{2n+1} = x+\frac{x^3}{3!}+\dots+\frac{x^{2n+1}}{(2n+1)!}$ do not have all
real roots for $n>0$.  Again, note that $q_2$ does not have all real
roots, and $q_n^\prime = q_{n-1}$.

\item $p_n = 1!x + 2!x^2 + \cdots + n!x^n$ does not have all real roots. If
  it did then its $(n-2)$th derivative would have two real roots, but  
$$ \left(\frac{d}{dx}\right)^{n-2}\, p_n = (n-1)!(n-2)!\left(
  \frac{1}{n-1} + x + \frac{n^2}{2}x^2\right)
$$
and this polynomial has no real roots. A more conceptual proof is
that the exponential map (see
\chap{operators}.\ref{sec:integration-exponential}) takes $p_n$ to
$1+x+x^2+\cdots+x^n$. This polynomial isn't in $\allpoly$, yet the
exponential map sends $\allpoly$ to itself.
\index{exponential operator}

\item We can use Lemma~\ref{lem:inequality-1} to show that there are no real
  roots. For example, 
\begin{quote}
  If $f$ and $g$ interlace and have no roots in common then the
  polynomial $ f\,g^\prime - f^\prime\,g$ has no real roots.
\end{quote}
\end{enumerate}

If a monic polynomial does not have \emph{any} real roots then it is
positive for all real values. Such polynomials are often stable - see
Chapter~\ref{cha:stable}. 


\chapter{Polynomials with all positive coefficients}
\label{cha:positive}

\renewcommand{\TimeStampStart}{Tuesday, August 21, 2007: 21:12:23}
\mytoday  

In this chapter we consider $\allpolypos$, the collection of all
polynomials that have all negative roots. If the leading coefficient
is positive then all coefficients are positive.  We introduce an
interlacing relation $\poslace$ for such polynomials that is weaker
than interlacing for $\allpoly$, but is the natural definition for
$\allpolypos$. We call this relation \emph{positive interlacing}.
\index{positive interlacing}
\index{interlacing!positive}

\section{Basic Properties}
\label{sec:pos-signs}

We establish some simple properties of $\allpolypos$.

\index{coefficients!sign conditions}
\begin{lemma} \label{lem:all-neg-converse}
Suppose that $f\in\allpolypos$.

\begin{enumerate}
\item If  $g$ in $\allpoly$ has all positive coefficients 
  then $g\in\allpolypos$.
\item If $g\in\allpoly$ and  $f\lesslesseq g$ then $g\in\allpolypos$.
\item If $g(x)\in\allpoly$ then there is a positive $\alpha$ such that
  $g(x+\alpha)\in\allpolypos$. 
\item $f(x)\in\allpolypos$ iff $f(-x)\in\allpolyalt$.
\item $f(-x^2)$ has all real roots if and only if
  $f(x)\in\allpolypos$.
\item  $f(x^2)$ has all real roots if and only if $f\in\allpolyalt$.

\end{enumerate}

\end{lemma}

\begin{proof}
  If $f\in\allpoly$ has all positive coefficients, and $r$ is
  positive, then $f(r)$ is a sum of positive terms, and hence $r$ is
  not a root of $f$. In addition, $0$ is not a root since $f(0)$ is
  the constant term which is  positive by assumption. 
  
  If $f\lesslesseq g$, the roots of $f$ are negative, and are
  interlaced by the roots of $g$, so the roots of $g$ are also all
  negative. Thus, $g\in\allpolypos$.  For the next part, take $\alpha$
  to be larger than the absolute value of the smallest root of $f$.
  
  The roots of $f(-x^2)$ are the square roots of the negatives of the
  roots of $f(x)$, so the roots of $f(x)$ are all negative. The case
  that $f\in\allpolyalt$ is similar.

\end{proof}

The next result constructs a polynomial in $\allpoly$ from two
interlacing polynomials in $\allpolypos$. See Theorem~\ref{thm:hurwitz} for a
converse. 

\begin{lemma} \label{lem:hurwitz}
  If $f\lesslesseq g$, $f,g\in\allpolypos$, then 

  \begin{enumerate}
  \item $  f(-x^2) \lesslesseq xg(-x^2)$.
  \item $f(-x^2)+xg(-x^2)\in\allpoly$.
  \end{enumerate}
\end{lemma}
\begin{proof}
  If the roots of $f(x)$ are $(-r_1^2,\dots,-r_n^2)$ then the roots of
  the new polynomial $f(-x^2)$ are $(\pm r_1,\dots,\pm r_n)$. Consequently,
  $f(-x^2)\lesslesseq xg(-x^2)$, which implies the lemma.
\end{proof}

Here's a different way to tell if a polynomial is in $\allpolypm$.

\begin{lemma} \label{lem:ispm}
  If $f$ and $g$ have no factors in common, $g(0)\ne0$, and $f$
  interlaces both $g$ and $xg$ then $g$ is in $\allpolyalt$ or
  $\allpolypos$.
\end{lemma}
\begin{proof}
  If we remove a zero from the roots of $xg$ then there will be two
  consecutive roots of $f$ with no root of $g$ in between unless $0$ is
  the largest or smallest of the roots of $xg$.
\end{proof}


Recall (Corollary~\ref{cor:leading-coef}) that if
$T\colon{}\allpoly\longrightarrow\allpoly$ then the leading coefficients are
all the same sign, or alternate in sign. If
$T\colon{}\allpolypos\longrightarrow\allpoly$ then there are more complicated
sign patterns. See Question~\ref{ques:signs}.

\begin{lemma} \label{lem:signs-T}
  If either $T\colon{}\allpolypos\longrightarrow\allpolyneg$ or
$T\colon{}\allpolyalt\longrightarrow\allpolyalt$, $T$ preserves degree, and
$T(1)$ is positive then the leading coefficient of $T(x^n)$ is positive.

  If $T\colon{}\allpolypos\longrightarrow\allpolyalt$ or 
$T\colon{}\allpolyalt\longrightarrow\allpolypos$, and $T$ preserves degree
then the signs of the leading coefficients of $T(x^n)$ alternate.
\end{lemma}
\begin{proof}
  If $T(x^n) = \sum c_{i,n}x^i$ then for any $\alpha$ 
  \begin{align*}
    T(\,(x+\alpha)^n\,) &= T(x^n) + n\alpha T(x^{n-1}) + \cdots \\
    &= c_{n,n}x^n + (c_{n-1,n}+n\alpha c_{n-1,n-1})x^{n-1} + \cdots
  \end{align*}
  If $T\colon{}\allpolypos\longrightarrow\allpolyneg$ then for any $\alpha>0$
  the polynomial $(x+\alpha)^n$ is in $\allpolypos$, and so
  $T((x+\alpha)^n)$ is also in $\allpolypos$. It follows that for any
  positive $\alpha$ both $c_{n,n}$ and $c_{n-1,n}+n\alpha c_{n-1,n-1}$
  have the same sign. Choosing $\alpha$ sufficiently large shows that
  the leading coefficients $c_{n,n}$ and $c_{n-1,n-1}$ have the same
  sign.

  If  $T\colon{}\allpolyalt\longrightarrow\allpolyalt$ then we choose
  $\alpha<0$ and the argument is similar. For the last part replace
  $x$ by $-x$.
\end{proof}

\added{9/2/7}
  If all the roots are positive then the coefficients
  alternate. Here's another alternation result.
\index{coefficients!alternating signs}

  \begin{lemma}\label{lem:even-alternate}
    If  a polynomial has all real roots and all  the
    coefficients of odd degree are zero then the signs of the
    non-zero coefficients alternate.
  \end{lemma}
  \begin{proof}
    If $f(x)$ is such a polynomial then the non-zero terms all have
    even degree, so $f(x)=f(-x)$. Thus, the roots of $f$ come in pairs
    $r_i,-r_i$, which implies 
\[ f(x) = \prod_i (x^2 - r_i^2) \]
and therefore the non-zero coefficients alternate in sign.
  \end{proof}

There is an interesting condition that guarantees that the signs of
the coefficients alternate in pairs.

\index{signs!patterns of coefficients}
\begin{lemma}\label{lem:sign-patterns}
  If $f(x)$ and $f(-x)$ interlace then the sign pattern of the
  coefficients of $f$ is $\cdots++--++--\cdots$.
\end{lemma}
\begin{proof}
We note that if
  \begin{align*}
    f(x) &= a_0 + a_1x + a_2 x^2 + \cdots \\
\intertext{then}
    f(x)+f(-x) &= 2( a_0 + a_2x^2 + \cdots) = 2g(x^2) \\
\intertext{Since the left hand side is in $\allpoly$ it follows that
  $g(x)$ has all positive roots. Similarly}
    f(x)-f(-x) &= 2x( a_1 + a_3x^2 + \cdots) = 2xh(x^2) 
  \end{align*}
  implies $h$ has all positive roots. Since $g$ and $h$ have all positive
  roots, their coefficients alternate, and so the coefficients of $f$
  have pattern $\cdots++--++\cdots$ since it is the intertwining of two
  alternating sequences.
\end{proof}

\section{Interlacing in $\allpolypos$}
\label{sec:pos-interlacing}

Our definition of interlacing for $\allpoly$ requires that all linear
combinations of a pair of interlacing polynomials lie in $\allpoly$.
Such a definition can not work if we replace $\allpoly$ by
$\allpolypos$, because if we choose a negative coefficient then the
resulting polynomial could have both positive and negative
coefficients. This motivates our definition:

\index{interlacing!in $\allpolypos$}
\index{\ zzplessless@$\plessless$}

\begin{definition}
  If $f,g\in\allpoly$ and $\deg(f)>\deg(g)$ then we say that
  $f\poslace g$ if and only if $f + \alpha g$ is in $\allpoly$ for
  all \emph{positive} $\alpha$.  We say that $f\plessless g$ if there
  are open neighborhoods $\mathcal{O}_f,\mathcal{O}_g$ of $f$ and $g$
  such that for $f_0\in\mathcal{O}_f$ and $g_o\in\mathcal{O}_g$ we
  have $f_0\poslace g_0$.

  We write $f \pgreateqeq g$ if there is a positive $\beta$ and
  $h\in\allpoly$ with positive leading coefficient such that
  $f\plesslesseq h$ and $g = \beta f + h$.
\end{definition}

It is clear from the definition that if $f,g\in\allpolypos$ satisfy
$f\lesslesseq g$ then $f \plesslesseq g$.  We can express
$\plesslesseq$ in terms of $\lesslesseq$.

\begin{lemma}
  Suppose that $f,g\in\allpoly$ have positive leading coefficients.
  The following are equivalent
  \begin{enumerate}
  \item $f \plesslesseq g$
  \item There is $k\in\allpoly$ such that $k\greateqeq f$ and $k
    \lesslesseq g$. That is, $f$ and $g$ have a common interlacing.
\index{common interlacing}
\index{interlacing!common}
  \end{enumerate}
\end{lemma}

\begin{proof}
 If there is such a polynomial  $k$ then Lemma~\ref{lem:add-interlace} implies that all positive
 linear combinations are in $\allpolypos$.  The converse follows from
 Proposition~\ref{prop:pattern}. 
\end{proof}

Linear transformations that map $\allpolypos$ to itself preserve
$\plesslesseq$. The proof is the same as in Theorem~\ref{thm:only-roots}.

\begin{lemma} \label{lem:pos-1}
  Suppose that $T$ is a linear transformation that maps $\allpolypos$ to itself, 
  and maps polynomials with all positive coefficients to polynomials
  with all positive coefficients. Then  $f \plesslesseq g$ implies
  that $Tf \plesslesseq Tg$.
\end{lemma}

Unlike $\lesslesseq$, it is not true that $f\plesslesseq g,h$ implies
that $g+h\in\allpoly$. For example, if
\begin{align*}
 f & = (x + 2)(x + 6)(x+7) & g & =  (x + 1)(x + 3) & h & =  (x + 4)(x + 5)
\end{align*}
then $f \plesslesseq g,h$ but $g+h$ has imaginary
roots%
\footnote{But see p. \pageref{lem:pos-1-stab}}%
.
In addition, even though
$f\plesslesseq h$, $f$ and $f+h$ do not interlace.

Next, we have a simple property of the derivative that  
will be generalized to homogeneous polynomials in Corollary~\ref{cor:fagp2}, and to
polynomials in two variables in Corollary~\ref{cor:fpp}.

\begin{cor} \label{cor:fagp}
  If $f\lesslesseq g$ both have positive leading coefficients then \\ {$f
  \plesslesseq -g^\prime$}.  Equivalently, $f-ag^\prime\in\allpoly$ for
  $a>0$.
\end{cor}
\begin{proof}
 Apply Lemma~\ref{lem:add-interlace} to the interlacings $f\lesslesseq g
 \lesslesseq g^\prime$.
\end{proof}

\begin{lemma} \label{lem:inequality-6a}
  If $f\in\allpoly$ then $f \pgreateqeq xf^\prime$. If $f\in\allpolyalt$ 
  then $f \greateqeq xf^\prime$, and if $f\in\allpolypos$ then
  $xf^\prime \greateqeq f$. 
\end{lemma}
\begin{proof}
  Since the leading coefficient of $xf^\prime$ is positive it suffices
  to show that $f+\alpha xf^\prime\in\allpoly$ for all
  positive $\alpha$.  So write $f =
  \prod(x-r_i)$ and let $g = xf$. Then
  \begin{align*}
    f+\alpha xf^\prime & = \frac{g}{x-0} + \sum \alpha 
    \frac{g}{x-r_i}\\
\intertext{where we recall that}
f^\prime & = \sum  \frac{f}{x-r_i}
  \end{align*}
This shows that $f+\alpha xf^\prime$ interlaces $xf$ for positive
$\alpha$, and hence is in $\allpoly$.

The second part follows since $f \lesslesseq f^\prime$, and $0$ is
either greater or less than all the roots of $f$. 
\end{proof}

\section{Linear combinations in $\allpolypos$}
\label{sec:pos-lin-comb}

In Proposition~\ref{prop:pattern} we found that if $\alpha f+ \beta g$ have all real
roots for all non-negative $\alpha,\beta$ then $f$ and $g$ have a
common interlacing.  The assumptions were not strong enough to conclude
that $f$ and $g$ interlace.  The next result gives assumptions
about non-constant combinations that lead to interlacing. 

\begin{lemma} \label{lem:axbfcg}
  Suppose that polynomials $f,g$ are in $\allpolypos$ and have
  positive leading coefficients. In addition, assume that for all
  non-negative $\alpha,\beta,\gamma$ the polynomial $(\alpha
  x+\beta)f+\gamma g$ has all real roots.  Then $g \longleftarrow f$.
\end{lemma}
\begin{proof}
  If we apply Proposition~\ref{prop:pattern} we see that for all positive $\beta$ the
  polynomials $(x+\beta)f$ and $g$ have a common interlacing.  If we
  choose $-\beta$ smaller than the smallest root of $f$ and $g$ then since
  $(x+\beta)f$ and $g$ have a common interlacing, it follows that the
  smallest root of $g$ is less than the smallest root of $f$. If the
  combined roots of $f$ and $g$ have two consecutive roots of $f$ (or
  $g$) then upon choosing $-\beta$ in between these roots we find that
  there is no common interlacing.  It is possible to make these
  choices with positive $\beta$'s since all roots of $f$ and $g$ are
  negative.
\end{proof}

\begin{lemma} \label{lem:ispm-2}
  If $f+\alpha g$ and $xf+\alpha g$ are in $\allpolypos$ for all positive
  $\alpha$ then $g\longleftarrow f$.  Equivalently, if $f \plesslesseq
  g$ and $xf \plesslesseq g$ where $f,g\in\allpolypos$ then
  $g\longleftarrow f$.
\end{lemma}
\begin{proof}
  The only was that common interlacing of $f$ and $g$, and of $xf$ and
  $g$, can occur is if $f$ and $g$ interlace.
\end{proof}

It is not easy for linear combinations of polynomials to remain in
$\allpolyaltclose$ or $\allpolyposclose$.

\begin{lemma} \label{lem:pos-neg-2}
Suppose that $f,g$ have the same degree, and that $f+\alpha
g\in\allpolyaltclose\cup\allpolyposclose$ for all $\alpha\in\reals$. 
Then at least one of the following is true
\begin{enumerate}
\item $g$ is a constant multiple of $f$.
\item There are constants $a,b,c,d$ such that $f=(ax+b)x^r$,\\ $g=(cx+d)x^r$.
\end{enumerate}
\end{lemma}
\begin{proof}
  We may assume that the degree is greater than $1$. 
  If it is not the case that $g$ is a constant multiple of $f$ then
  by Lemma~\ref{lem:pos-neg-1} we can find $\alpha,\beta$ such that $f+\alpha g$ 
  has a positive root, and $f+\beta g$ has a negative root. Since by
  hypothesis we have that $f+\alpha g\in\allpolyaltclose$, and
  $f+\beta g\in\allpolyposclose$, it follows that for some
  $\alpha\le\gamma\le\beta$ we have that $f+\gamma
  g\in\allpolyaltclose\cap\allpolyposclose$. Now the only
  polynomials in $\allpolyaltclose\cap\allpolyposclose$ are multiples
  of powers of $x$, so we find that $f+\gamma g = e x^r$ for some
  $e,r$. Substituting shows that $e x^r +
  (\alpha-\gamma)g\in\allpoly$ for all $\alpha$. Consequently, $g$ and 
  $x^r$ interlace. This implies that $g = (cx+d)x^r$ or
  $g=cx^{r-1}$. The latter is not possible since $f$ has degree at least 
  $r$, so $g = (cx+d)x^r$. The result now follows by substituting for $g$.
\end{proof}

If $f\greateqeq g$ then the largest root belongs to $f$, and the
smallest to $g$. Thus, if $f,g\in\allpolypos$ we have $xg\lesslesseq
f$, and if $f,g\in\allpolyalt$ we know $xf\lesslesseq g$. 

\begin{lemma} \label{lem:xfg}
  Suppose $f\greateq g$ have positive leading coefficients.  If
  $f,g\in\allpolypos$ then $xg+f\lessless f$. If $f,g\in\allpolyalt$
  then $xf - g \lessless g$.
\end{lemma}
\begin{proof}
  Since $g\in\allpolyalt$ we know that $xf\lessless g$, and therefore
  $xf-g\lessless g$ by Corollary~\ref{cor:where-roots}. The second case follows by
  replacing $x$ by $-x$. 
\end{proof}

\begin{cor} \label{cor:tf-3}
  If $a\not\in(0,n)$ then 
  \begin{align*}
    \text{If $f\in\allpolypos(n)$ then } & -af + xf^\prime 
\begin{cases} 
  \greateqeq f & a > n \\
  \lessgreateq f & a=n \\
  \lesseqeq f & a < 0
\end{cases} \\
    \text{If $f\in\allpolyalt(n)$ then } & -af + xf^\prime 
\begin{cases} 
  \lesseqeq f & a > n \\
  \lessgreateq f & a=n \\
  \greateqeq f & a < 0
\end{cases}
\end{align*}
\end{cor}

\begin{proof}
  If $f\in\allpoly$ and $a\not\in(0,n)$ then
  $-af+xf^\prime\in\allpoly$ by Corollary~\ref{cor:agxgp}. The remaining parts
  follow the usual arguments.
\end{proof}

\section{Various linear combinations}
\label{sec:elem-interlacing}

In this section we first look at linear combinations of the form
$$h = (ax+b) f + (cx^2+dx+e) g$$
where $f\lesslesseq g$. We are
interested in what restrictions we may place on $a,b,c,d,e$ to
conclude that $h$ interlaces $f$, or that $h$ is in $\allpolypos$ or
$\allpolyalt$. The  results are all easy consequences of
Lemma~\ref{lem:fghjk}, or or Lemma~\ref{lem:sign-quant}.

\begin{theorem} \label{thm:tff} Let $Tf = (ax+b)f + (cx^2+dx+e)g$
  where $c\ne0$, $f\in\allpolypos(n)$ has positive leading
  coefficient, and $f\lesslesseq g$. Assume that $cx^2+dx+e$ has
  constant sign $\epsilon$ on all roots of $f$ and $(a+nc)\epsilon <0$
  for all positive integers $n$.
  \begin{enumerate}
  \item   $Tf\lessless f$.
  \item  If  $f\in\allpolyint{(-\infty,\alpha)}$ and
    $(a+nc)\cdot(Tf)(\alpha)>0$ then $Tf\in\allpolyint{(-\infty,\alpha)}$.
  \item If  $f\in\allpolyint{(\alpha,\infty)}$ and
    $(-1)^n(a+nc)\cdot(Tf)(\alpha)>0$ then $Tf\in\allpolyint{(\alpha,\infty)}$.
  \item If  $f\in\allpolypos$ and $(a+nc), b, e$ have
    the same sign then $Tf\in\allpolypos$.
  \item If  $f\in\allpolyalt$ and $(a+nc),b,e$
    have the sign then $Tf\in\allpolyalt$.
  \end{enumerate}
\end{theorem}
\begin{proof}
  In order to show that $Tf \lessless f$ we need to show that $Tf$
  sign interlaces $f$, and that $Tf$ evaluated at the largest root of
  $f$ has sign opposite to the leading coefficient of $Tf$
  (Lemma~\ref{lem:int-2}).  If $z$ is a root of $f$, then $(Tf)(z) =
  g(z)\cdot(cz^2+dz+e)$. Since $g$ sign interlaces $f$,
  we must have that $cz^2+dz+e$ has constant sign on any interval
  containing $\roots(f)$.
  
  If $z$ denotes the largest root of $f$, then the sign of
  $(Tf)(z)$ equals \\ $sgn\, g(z)\,sgn(az^2+bz+c)$, and since $f$ has positive
  leading coefficient, $sgn\,g(z)$ is positive.  The leading
  coefficient of $Tf$ is $a+nc$, so a sufficient condition for
  $Tf\lessless f$ is that $sgn(a+nc)\,sgn(az^2+bz+c)$ is negative.
  
  Now assume that $f\in\allpolyint{(-\infty,\alpha)}$.  Since
  $Tf\lessless f$, the only root of $Tf$ that could possibly be
  greater than $\alpha$ is the largest root. If $(Tf)(\alpha)$ has the
  same sign as the leading coefficient $a+nc$ of $Tf$ then there are
  no roots greater than $\alpha$.  

  In case $f\in\allpolyint{(\alpha,\infty)}$, we follow the preceding
  argument; the condition is that $(Tf)(\alpha)$ has the same behavior
  as $(Tf)(x)$ as $x\longrightarrow-\infty$.

  For last two parts we use  the fact that $(Tf)(0) = bf(0) + e
  f^\prime(0)$.   If $f\in\allpolypos$ then $f(0)$ and $f^\prime(0)$
  are positive.   If $f\in\allpolyalt$ then $f(0)$ and $f^\prime(0)$
  have sign $(-1)^n$.

\end{proof}

The corollary below follows immediately from Lemma~\ref{lem:sign-quant}.

\begin{cor} \label{cor:combinations}
  Suppose that $f,g,h,k$ have positive leading coefficients and
  satisfy $f\greateq g \greateq h$, $g\lessless k$.  Then for all $a$
  and positive $\alpha,\beta,\gamma$ we have
  \begin{enumerate}
  \item If $g\in\allpolypos$ then $xf + ag + \alpha h - \beta f -
    \gamma k   \lessless g$

  \item If $g\in\allpolyalt$ then $xh + ag - \alpha h + \beta f +
    \gamma  k  \lessless g$
  \item $xg + \beta h \greateq g$
  \end{enumerate}
\end{cor}

  \section{The positivity hierarchy}
  \label{sec:positivity-hierarchy}

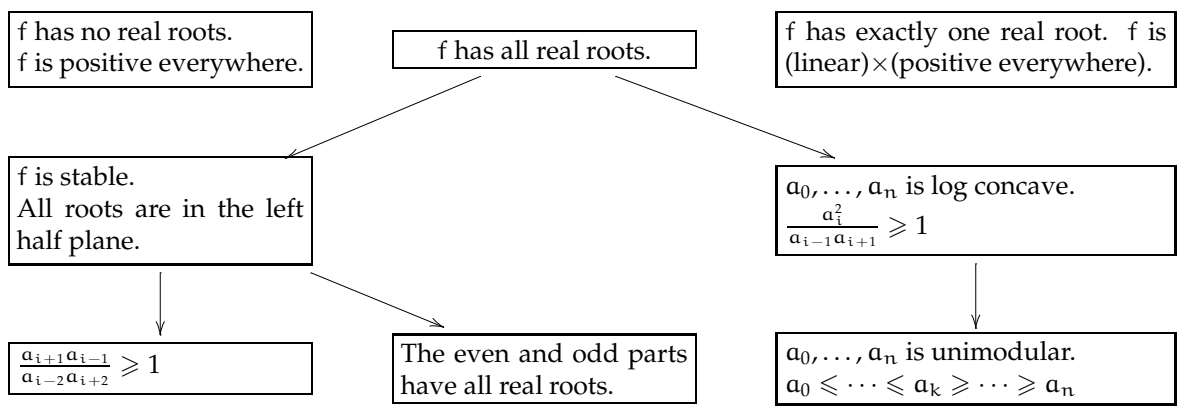
\begin{figure}[h]

\centerline{\xymatrix{
{\fbox{\parbox{1.5in}{$f$ has no real roots.\\ $f$ is positive everywhere.}}} &
{\fbox{\parbox{1.5in}{\centerline{$f$ has all real roots.}}}}
  \ar@{->}[dl] \ar@{->}[dr]   &
{\fbox{\parbox{2in}{$f$ has exactly one real root. 
$f$ is (linear)$\times$(positive everywhere).}}}\\
{\fbox{\parbox{1.5in}{$f$ is stable. \\ All roots are in the left half plane.}}}
\ar@{->}[d] \ar@{->}[dr]&
&
{\fbox{\parbox{2in}{$a_0,\dots,a_n$ is log concave.\\
    $\frac{a_i^2}{a_{i-1}a_{i+1}}\ge1$}}} \ar@{->}[d]\\
{\fbox{\parbox{1.5in}{$\frac{a_{i+1}a_{i-1}}{a_{i-2}a_{i+2}}\ge1$}}} &
{\fbox{\parbox{1.5in}{The even and odd parts have all real roots.}}} &
{\fbox{\parbox{2in}{$a_0,\dots,a_n$ is unimodular. \\ $a_0\le\cdots\le a_k
    \ge \cdots \ge a_n$}}}
   }}

    \caption{The positivity hierarchy}
    \label{tab:pos-hier}
  \end{figure}

\added{6/3/7}
  If we are given a polynomial with all positive coefficients then
  there are many properties that it might have, and  we should
  always check to see which properties it does have. These properties
  are given in Figure~\ref{tab:pos-hier}. We assume that $f =
  a_0+\cdots + a_nx^n$ where all $a_i$ are positive. The conditions
  get weaker as we move downward.

\index{polynomial!stable}
\index{stable!polynomial}
\index{log concave}
\index{unimodular}
\index{even and odd parts}

\chapter{Matrices that preserve interlacing}
\label{cha:poly-matrices}

\renewcommand{\TimeStampStart}{Tuesday, March 11, 2008: 09:48:00}
\mytoday  

This chapter is concerned with the question:
\begin{quote}
  Suppose that $M$ is a matrix, either a constant matrix, or with
  polynomial entries. Suppose that $v$ is a vector of polynomials
  lying in $\allpoly$, or perhaps in $\allpolypos$. When is $Mv$ also
  a vector of polynomials in $\allpoly$?
\end{quote}

The answer when $M$ is constant is roughly that $M$ satisfy a property
called total positivity. Such matrices are discussed in the first
section, and the characterization is proved in the following
section. The next section discusses two by two matrices that preserve
interlacing, and Section 4 covers some results about general $M$. The
final section has some examples of matrices that preserve interlacing
for polynomials in $\allpolypos$.

\section{Totally positive matrices}
\label{sec:poly-totally-positive}

In this section we recall properties of totally positive matrices.  We
will see that these matrices occur whenever we consider linear
combinations of interlacing polynomials. See \cite{ando} and
\cite{fomin} for the proofs of the theorems in this section.

\index{strictly totally positive} 
\index{totally positive} 
\index{matrix!totally positive} 
\index{matrix!strictly totally positive} 

A matrix is \emph{strictly totally positive} if all its minors are positive. If
all the minors are  non-negative we say it is \emph{totally
  positive}.  A simple consequence of the Cauchy-Binet formula for
the determinant of the product of two matrices is that products of
totally positive matrices are totally positive, and products of
strictly totally  matrices are strictly totally positive.

A weakening of this definition is often sufficient.  We say
that a matrix is is \emph{strictly totally positive$_2$} if all
elements and all two by two determinants are positive. It is
\emph{totally positive$_2$} if all elements and all two by two
determinants are non-negative.

The first theorem is due to Whitney:

\index{Whitney theorem}
\begin{theorem}[Whitney]
  Every  totally positive matrix is the limit of a
  sequence of strictly totally positive matrices.
\end{theorem}

The next theorem describes a useful decomposition of totally
non-negative matrices:

\index{Loewner-Whitney theorem}
\begin{theorem}[Loewner-Whitney]
  Any invertible totally positive matrix is a product of elementary
  Jacobi matrices with non-negative matrix entries.
\end{theorem}

There are three kinds of elementary Jacobi matrices. Let $E_{i,j}$ be
the matrix whose $i,j$ entry is $1$ while all other entries are $0$.
If $I$ is the identity matrix then the elementary Jacobi matrices are
$I+t E_{i,i+1}$, $I+t E_{i+1,i}$, and $I+(t-1)E_{i,i}$ where $t$ is
positive. For example, in the case $i=2$ for $4$ by $4$ matrices they
are

\begin{equation}\label{eqn:jacobi-matrix}
\begin{pmatrix}
  1&0&0&0\\
  0&1&t&0\\
  0&0&1&0\\
  0&0&0&1
\end{pmatrix}\ 
\hspace*{.5cm}
\begin{pmatrix}
  1&0&0&0\\
  0&1&0&0\\
  0&t&1&0\\
  0&0&0&1
\end{pmatrix}\ 
\hspace*{.5cm}
\begin{pmatrix}
  1&0&0&0\\
  0&t&0&0\\
  0&0&1&0\\
  0&0&0&1
\end{pmatrix}\ 
\end{equation}

If $\smalltwo{a}{c}{b}{d}$ is strictly totally positive, then $a,b,c,d>0$, and
$ad-bc>0$. An explicit factorization into Jacobi matrices is

\begin{equation}\label{eqn:jacobi-factor}
 \begin{pmatrix}{a}& {c}\\{b}&{d} \end{pmatrix}=
\begin{pmatrix}{1}& {0}\\{0}&{b} \end{pmatrix}
\begin{pmatrix}{a}& {0}\\{0}&{1}\end{pmatrix}
\begin{pmatrix}{1}& {0}\\{0}&{\frac{ad-bc}{ab}}\end{pmatrix}
\begin{pmatrix}{1}& {0}\\{\frac{ab}{ad-bc}}&{1}\end{pmatrix}
\begin{pmatrix}{1}& {\frac{c}{a}}\\{0}&{1}\end{pmatrix}
\end{equation}

\section{Constant matrices preserving mutual interlacing}

Suppose that we have a vector $v$ of polynomials in $\allpoly$. If $A$
is a constant matrix when is $Av$ also a vector of polynomials in
$\allpoly$?
\index{matrix!preserving interlacing}

Even if we are only considering three polynomials, a condition such as
\\ $f_1 
\greateq f_2 \greateq f_3$ is not strong enough to draw any conclusions.
The reason is that when we form linear combinations, we will need to 
know the relationship between $f_1$ and $f_3$ as well.  Thus, we
assume that $f_1,f_2,f_3$ is mutually interlacing:

\index{interlacing!mutually}
\index{mutually interlacing}

\begin{definition}
  A sequence of polynomials $f_1,\dots,f_n$ is \emph{mutually
    interlacing} if and only if for all $1\le i < j\le n$ we have $f_i
  \greateqeq f_j$, and each $f_i$ has positive leading coefficient.

The roots of a mutually interlacing have a simple ordering. If we
denote the ordered roots of $f_i$ by $r_i^j$ then the roots of $f_1,\dots,f_n$
 are ordered
$$ r_1^1\le r_2^1 \le\dots\le r_n^1\ \le\  r_1^2\le r_2^2 \le\dots\le r_n^2\
\le\  r_1^3\le r_2^3 \le\dots\le r_n^3\ \le\  \cdots $$
\end{definition}

Here are a few examples of mutually interlacing sequences of polynomials:
\begin{lemma} \label{lem:ti} \ 

  \begin{enumerate}
  \item If $g(x)$ is a polynomial with   roots $ a_1\le 
    \cdots \le a_n$, then the sequence below is  mutually
    interlacing:
    $$
    \frac{g}{x-a_1},\ \frac{g}{x-a_{2}},\ \dots,\ 
    \frac{g}{x-a_{n-1}},\ \frac{g}{x-a_n}$$
  \item Any subsequence of a mutually interlacing sequence is mutually
    interlacing.
    
  \item  If $f_1,f_2,\dots,f_n$ is a mutually interlacing sequence then
     the sequences below are mutually interlacing
\begin{align*} 
&
f_1,\,f_1+f_2,\,f_2,\,f_2+f_3,\,\dots,\,f_{n-1},\,f_{n-1}+f_n,\,f_n\\
&f_1,f_1+f_2,f_1+f_2+f_3,\dots,f_1+f_2+\cdots+f_n\\
    & f_1,f_2,f_2,f_2,\dots,f_2
  \end{align*}
\end{enumerate}
\end{lemma}

\begin{proof}
The roots of the first sequence of
  polynomials are
    $$
  (a_2,a_3,\dots,a_{n}),\ (a_1,a_3,\dots,a_{n}), \dots,\ 
  (a_1,\dots,a_{n-2},a_n),\ (a_1,\dots,a_{n-2},a_{n-1}) 
  $$
  which are interlacing.  The remaining assertions follow easily from
  Corollary~\ref{cor:where-roots}. 
\end{proof}

We first determine when a linear combination of mutually interlacing
polynomials has all real roots. The answer is simple:  the coefficients have
at most one change of sign.

\begin{lemma} \label{lem:faf}
  Suppose that $a_1,\dots,a_n$ are non-zero constants, \\ and $f_1
  \greateqeq \cdots \greateqeq f_n$ is a mutually interlacing sequence
  of polynomials. The following are equivalent:
  \begin{enumerate}
  \item For all possible $f_i$ as above $\sum a_i f_i$ has all real roots.
  \item There is an integer $k$ so that all the $a_i$ where $i\le k$
    have the same sign, and all $a_i$ where  $i>k$ have the same sign.
  \end{enumerate}
\end{lemma}

\begin{proof}
   We can write $f_i = f_k +\epsilon_i r_i$ where $f_i\lesslesseq
   r_i$, $r_i$ has positive leading coefficient, and 
  $\epsilon_i=-1$ if $i<k$ and $\epsilon_i=1$ if $i>k$. Now
  $$
f=  \sum a_i f_i  = (\sum_i a_i)f_k + \sum_i \left(a_i \epsilon_i r_i\right)$$
  and by hypothesis
  the coefficients $a_i \epsilon_i$ all have the same sign. It follows
  from Lemma~\ref{lem:sum-1} that $f\in\allpoly$.

  Conversely, consider the three mutually interlacing polynomials 
  \begin{align*}
    g_1 &= (x-8)(x-24) \\
    g_2 &=  (x-5)(x-21) \\
    g_3 &= (x-2)(x-13) \\
\intertext{where $g_1\greateq g_2 \greateq g_3$. Observe that }
 g_1 - g_2 +  g_3 &= x^2 -21x +13 \not\in\allpoly.
  \end{align*}

  Now embed $g_1,g_2,g_3$ in a mutually interlacing sequence so that
  $f_1=(1/|a_1|)g_1$, $f_i=(1/|a_i|)g_2$, and
  $f_n=(1/|a_n|)g_3$. Since 
  $$
  a_1 f_1 + a_if_i + a_nf_n = sgn(a_1) g_1 +sgn(a_i)g_2 +
  sgn(a_n)g_3$$
  we know that if the left hand side has all real roots then  we can not
  have that $sgn(a_1)=sgn(a_n)\ne sgn(a_i)$. This implies that there
  can be at most one sign change.
\index{sign change}
\end{proof}

The following obvious corollary will be generalized in the next
section where sums are replaced by integrals.

\begin{cor}\label{cor:faf}
  Suppose that $a_1,\dots,a_n$ are positive constants, \\ and $f_1
  \greateqeq \cdots \greateqeq f_n$ is a mutually interlacing sequence
  of polynomials. Then
  \begin{equation}
    \label{eqn:faf-1}
    f_1 \greateqeq a_1 f_1 + \cdots + a_n f_n \greateqeq f_n
  \end{equation}
\end{cor}

Our first theorem is a simple consequence of the general theorems of
the previous section.

\begin{theorem} \label{thm:tnn}
  Suppose $f_1,\dots,f_n$ are mutually interlacing, and $A = (a_{i,j})$
  is an $n$ by $n$ matrix of constants such that 
  $$A \cdot(f_1,\dots,f_n) = (g_1,\dots,g_n)$$
  If $A$ is totally positive, then $g_1,\dots,g_n$ is mutually interlacing.
\end{theorem}

\begin{proof}
  Since the limit of a mutually interlacing sequence is mutually
  interlacing, and any non-invertible totally positive matrix is
  the limit of strictly totally positive matrices, we may assume that $A$ is
  invertible.  By the Loewner-Whitney theorem we can assume that A is
  an elementary Jacobi matrix. The three classes of matrices applied
  to $f_1\greateqeq \cdots \greateqeq f_n$ give three easily verified mutually
  interlacing sequences:
$$
\begin{matrix}
  f_1 & \cdots & f_{i-1} & f_i+tf_{i+1} & f_{i+1} & \cdots & f_n \\[.2cm]
  f_1 & \cdots & f_{i-1} & f_i & tf_i + f_{i+1} & \cdots & f_n \\[.2cm]
  f_1 & \cdots & f_{i-1} & tf_i & f_{i+1} & \cdots & f_n \\
\end{matrix}
$$
\end{proof}

The two by two case for all positive entries follows from the
theorem. However, we can determine the action of an arbitrary two by
two matrix.

\begin{cor} \label{cor:lin-comb-new}
  Suppose that $f_1 \greateq f_2$ (or $f_1 \lessless f_2$) have
  positive leading coefficients, and
  $\smalltwo{a}{b}{c}{d}\smalltwobyone{f_1}{f_2} =
  \smalltwobyone{g_1}{g_2}$, where without loss of generality we may
  assume that $a$ is positive.  If the determinant $ad-bc$ is positive
  then $g_1\greateq g_2$, and if it is negative $g_2 \greateq g_1$.
\end{cor}
\begin{proof}
  Assume that  the determinant $ad-bc$ is
    non-zero. If $r$ is a common root of $g_1$ and $g_2$ then since the
    matrix is invertible we conclude that $f_1(r)=f_2(r)=0$ which
    contradicts the assumption that $f_1$ and $f_2$ have no common
    factors. Now any linear combination of $g_1$ and $g_2$ is a linear
    combination of $f_1$ and $f_2$, and thus has all real
    roots. Consequently, $g_1$ and $g_2$ interlace. The direction of
    interlacing will not change if we continuously change the matrix
  $\smalltwo{a}{b}{c}{d}$ as long as the determinant never becomes
  zero. Consequently, if the determinant is positive we can deform
  the matrix to the identity, in which case $g_1 \greateq g_2$. If the
  determinant is negative we can deform the matrix to
  $\smalltwo{1}{0}{1}{-1}$, and since $f-g \greateq f$ we are done.
\end{proof}

\begin{example}
  Consider $\smalltwo{1}{b}{1}{d}\smalltwobyone{f}{g} =
  \smalltwobyone{f+bg}{f+dg}$.  If $f\greateqeq  g$ and $0<b<d$ then
  $f+bg \greateqeq f+dg$. 

  The examples of mutually interlacing sequences in Lemma~\ref{lem:ti}(3) are
  given by the following totally positive matrices, where we take
  $n=3$ for clarity, and the ``.'' are zeros.

$$
\begin{pmatrix}
  1 & . & . \\ 1 & 1 & . \\  . & 1 & . \\ . & 1 & 1 \\ .&.&1
\end{pmatrix}
\quad\quad
\begin{pmatrix}
  1&.&.\\1&1&. \\1&1&1
\end{pmatrix}
\quad\quad
\begin{pmatrix}
  1&.&.\\ .&1&. \\.&1&.\\ .&1&.
\end{pmatrix}
$$
\end{example}

  An immediate corollary of Corollary~\ref{cor:lin-comb-new} is that the
  inverse of a totally positive matrix sometimes preserves
  interlacing.

  \begin{cor} \label{cor:inv-2by2}
If
\begin{enumerate}
\item $M$ is a $2$ by $2$ totally positive invertible matrix with
  non-negative entries
\item $f,g,u,v \in\allpoly$ have positive leading coefficients
\item $M\smalltwobyone{f}{g} = \smalltwobyone{u}{v}$
\end{enumerate}
then
$ f\greateqeq g \text{ if and only if } u \greateqeq v $.
\end{cor}


\begin{remark}
  The inverses of the $3$ by $3$ Jacobi matrices (except the diagonal
  ones) do not satisfy the conclusions of the lemma. However, we do have that
  the matrices of the form
\[ \begin{pmatrix} 1 & 0 & -a \\ 0 & 1 & 0 \\ -b & 0 & 1
\end{pmatrix}
\]
satisfy the conclusions if $0 \le ab < 1$. That is, if $0\le a,b$ and
$0\le ab<1$,
\[ \begin{pmatrix} 1 & 0 & -a \\ 0 & 1 & 0 \\ -b & 0 & 1
\end{pmatrix}
\begin{pmatrix}
  f \\ g \\ h
\end{pmatrix} =
\begin{pmatrix}
  u \\ v \\ w
\end{pmatrix}
\]
$f,g,h$ is mutually interlacing, and $f,g,h,u,v,w$ have positive
leading coefficients then $u,v,w$ is mutually interlacing. This is
equivalent to showing that
\begin{xalignat*}{3}
    f - a h & \greateqeq g & g & \greateqeq -bf + h & 
f - a h & \greateqeq  h - bf
\end{xalignat*}
The first follows from $f \greateqeq g \greateqeq h$. For the second,
write $f = \alpha g-r, h=\beta g + s$ where $g\lesslesseq r,s$. Then
$g\greateqeq (\beta - b \alpha)g + s + b r$ since the leading
coefficient of the right hand side is positive by hypothesis. Finally,
the last one is a two by two matrix $
\smalltwo{1}{-a}{-b}{1}\smalltwobyone{f}{h} =
\smalltwobyone{f-ah}{h-bf}$. The lemma shows that the interlacing is
satisfied if the determinant is positive, which is true since $0\le
ab<1$.

\end{remark}

Surprisingly, we do not need the full force of total positivity to
conclude that a matrix preserves interlacing, but just the
non-negativity of all two by two submatrices.

\begin{theorem} \label{thm:totally-preserves}
  Suppose that $A=(a_1,\dots,a_n)$ and $B=(b_1,\dots,b_n)$ are vectors
  of positive real numbers, and $f=(f_1,\dots,f_n)$ is a vector of
  mutually interlacing polynomials.
  Then
$$
Af \greateqeq Bf \text{ for all $f$ if and only if }
\smalltwobyone{A}{B} \text{ is } \text{ totally positive}_2.$$
\end{theorem}

\begin{proof}
  Consider the mutually interlacing set of polynomials
  $$c_1f_1,c_2f_2,\dots,c_nf_n$$ where the $c_i$ are positive. By continuity of
  roots, if we set $c_i=c_j=1$, and all other $c$'s to zero, then
$$ a_if_i + a_j f_j \greateqeq b_i f_i + b_j f_j.$$
We can rewrite this as
$\smalltwo{a_i}{a_j}{b_i}{b_j}\smalltwobyone{f_i}{f_j} =
\smalltwobyone{h}{k}$ where $h\greateqeq k$. By Corollary~\ref{cor:lin-comb-new} we must 
have $\smalltwodet{a_i}{a_j}{b_i}{b_j}\ge0$.

Conversely, given $A$ and $B$ construct a totally positive $n$ by
$n$ matrix $C$ whose first row is $A$, second row is $B$, and whose
remaining rows are all $0$. If we apply $C$ to the mutually interlacing
sequence $f,g,g,\cdots,g$ the conclusion follows from Theorem~\ref{thm:tnn}.
\end{proof}

A small modification of this argument shows

\index{totally positive$_2$}
\begin{cor} \label{cor:tnn}
  Suppose $f_1,\dots,f_n$ are mutually interlacing, and $A = (a_{i,j})$
  is an $n$ by $m$ matrix of constants such that 
  $$A \cdot(f_1,\dots,f_n) = (g_1,\dots,g_m)$$
  If $A$ is totally
 positive$_2$, then $g_1,\dots,g_m$ is mutually interlacing.
\end{cor}

If we fix the constants rather than the polynomials, then we need a
determinant condition on the polynomials. 

  \begin{lemma}\label{lem:int-fam-int-1}
    Suppose that we have the following interlacing diagram 

  $$ \begin{matrix}
f_1 & \greateq & f_2 & \greateq & \dots & \greateq & f_n \\[-.1cm]
\rotatebox{270}{$\greateq$} &&  \rotatebox{270}{$\greateq$}
&&&&\rotatebox{270}{$\greateq$} \\ 
g_1 & \greateq & g_2 & \greateq & \dots & \greateq & g_n 
    \end{matrix}
$$
where 
\begin{enumerate}
\item The $f_i$ are mutually interlacing.
\item The $g_i$ are mutually interlacing.
\item $\smalltwodet{f_i}{f_{i+1}}{g_i}{g_{i+1}}<0$ for $1 \le i < n$
\end{enumerate}
then for any non-negative $\alpha_1,\dots,\alpha_n$ we have
$
\sum \alpha_i\,f_i \greateqeq \sum \alpha_i\, g_i
$
  \end{lemma}
  \begin{proof}
    We may assume that all $\alpha_i$ are positive.  We proceed by
    induction. The case $n=2$ is Lemma~\ref{lem:inequality-4}. 
    let $h_1 = \alpha f_1 + \alpha_2 f_2$ and $k_1 = \alpha g_1 +
    \alpha_2 g_2$ then by the lemma we know that $h_1 \greateq k_1$ and
    $\{h_1,f_2,\dots,f_n\}$ and $\{k_1,g_2,\dots,g_n\}$ are mutually
    interlacing. The determinant condition is satisfied since
$$
\begin{vmatrix}
  h_1 & f_2 \\ k_1 & g_2 
\end{vmatrix} =
\begin{vmatrix}
  \alpha_1 f_1 + \alpha_2 f_2 & f_2 \\ \alpha_1 g_1 + \alpha_2 g_2 & g_2
\end{vmatrix}
=
\alpha_1\begin{vmatrix}
  f_1 & f_2 \\ g_1 & g_2 
\end{vmatrix} < 0
$$
We continue combining terms  until done.
  \end{proof}

Mutually interlacing sequences have a nice interpretation when
expressed in the interpolation basis.  Suppose $f_1,\dots,f_n$ is a
sequence of mutually interlacing polynomials of degree $d$ with
positive leading coefficients.  All the smallest roots of
$f_1,\dots,f_n$ are smaller than all the second smallest roots of
$f_1,\dots,f_n$, and so on. Thus, we can always find a $g$ such that
$g\lesslesseq f_i$ for all $i$. It follows that there are non-negative
constants $a_{i,j}$ such that
$$f_i =
  \sum a_{i,j}\dfrac{g(x)}{x-b_j}.$$

If we set $A=(a_{i,j})$, $F=(f_1,\dots,f_n)$ and
$G=(\dfrac{g}{x-b_j})$ then $AG=F$. The next lemma characterizes
those $A$'s that always yield mutually interlacing sequences.

\begin{cor} \label{cor:tot-pos-interp}
  Suppose $G$ and $A$ are as above.  For all choices of $G$ the
  sequence $AG$ is mutually interlacing if and only if $A$ is totally
  positive$_2$.
\end{cor}

\begin{proof}
We simplify what needs to be proved. First of all, it suffices to
assume that $A$ has only two rows since the conclusion only involves
two by two matrices. Thus we will show that if for all $g$ we have
$$ \sum a_i \dfrac{g(x)}{x-b_i} \greateqeq \sum c_i
\dfrac{g(x)}{x-b_i}$$
then $\smalltwodet{a_i}{a_{i+1}}{c_i}{c_{i+1}}\ge0$. Next, notice that we
can multiply any column of $A$ by any positive number, and the
determinants will still be non-negative. By continuity we can set all
other columns to zero, and so we  assume
that 
$$ a_i  \dfrac{g(x)}{x-b_i} + a_{i+1}\dfrac{g(x)}{x-b_{i+1}}
\greateqeq 
c_i  \dfrac{g(x)}{x-b_i} + c_{i+1}\dfrac{g(x)}{x-b_{i+1}}
$$
Next, we choose $g = x^{n-2}(x-\alpha)(x-\beta)$ where $b_i=\alpha$
and $b_{i+1}=\beta$, and $\beta>\alpha$. We can cancel the powers of
$x$, and so we have
$$ a_i  (x-\beta) + a_{i+1}(x-\alpha)
\greateqeq c_i  (x-\beta) + c_{i+1}(x-\alpha)
$$
The root of the left hand side is not less than the right hand side,
so
$$ \frac{c_{i+1}\alpha + c_i\beta}{c_i+c_{i+1}} \le
\frac{a_{i+1}\alpha + a_i\beta}{a_i+a_{i+1}} $$
and from this we find that
$$ (\alpha-\beta)(a_{i+1}c_i-a_ic_{i+1})\ge0$$
which is the desired conclusion.

The converse follows from Corollary~\ref{cor:tnn} since the set
$\left\{\frac{g(x)}{x-b_i}\right\}$ is mutually interlacing.
\end{proof}

The following property of mutually interlacing polynomials generalizes
Lemma~\ref{lem:fg-ab}, and will be extended to integrals in the next
section (Lemma~\ref{lem:convolution-int}). 

\index{convolution!mutually interlacing}
\index{mutually interlacing!convolution}
\index{mutually interlacing!sum in reverse order}

  \begin{lemma}\label{lem:convolution}
    If $f_1,\dots,f_n$ and $g_1,\dots,g_n$ are two sequences of
    mutually interlacing polynomials with positive leading
    coefficients, then
$$ f_1\,g_n + f_2\,g_{n-1} + \cdots + f_n\,g_1 \in\allpoly.$$
Equivalently, $\sum f_i(x) g_i(-x)\in\allpoly$.
  \end{lemma}
  \begin{proof}
    Assume the $f$'s have degree $r$ and the $g$'s have degree $s$. We
    may assume that the roots are distinct.
It is helpful to visualize the location
    of these roots in two dimensions, so we plot the roots of $f_i$
    and $g_{n+1-i}$ with $y$ coordinate $i$. For instance, if
    $r=3,s=2,n=4$ then a possible diagram of the roots is

    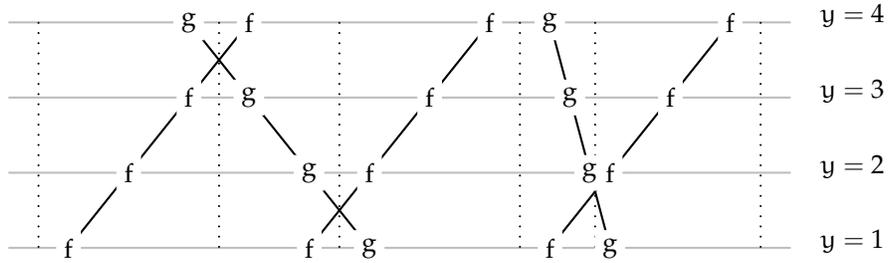
\begin{figure}[htbp]
      \centering
            \label{fig:convolution}
    \begin{pspicture}(-1,1)(13,4)
      \psset{xunit=.8cm,yunit=1cm}
      \psline[linecolor=lightgray](-1,1)(12,1)
      \psline[linecolor=lightgray](-1,2)(12,2)
      \psline[linecolor=lightgray](-1,3)(12,3)
      \psline[linecolor=lightgray](-1,4)(12,4)
      \psline(0,1)(3,4)
      \psline(4,1)(7,4)
      \psline(8,1)(11,4)
      \psline(2,4)(5,1)
      \psline(8,4)(9,1)
      \psline[linestyle=dotted](-.5,1)(-.5,4)
      \psline[linestyle=dotted](11.5,1)(11.5,4)
      \psline[linestyle=dotted](7.5,1)(7.5,4)
      \psline[linestyle=dotted](2.5,1)(2.5,4)
      \psline[linestyle=dotted](4.5,1)(4.5,4)
      \psline[linestyle=dotted](8.75,1)(8.75,4)
      \rput*(0,1){f}
      \rput*(1,2){f}
      \rput*(2,3){f}
      \rput*(3,4){f}
      \rput*(4,1){f}
      \rput*(5,2){f}
      \rput*(6,3){f}
      \rput*(7,4){f}
      \rput*(8,1){f}
      \rput*(9,2){f}
      \rput*(10,3){f}
      \rput*(11,4){f}
      \rput*(2,4){g}
      \rput*(3,3){g}
      \rput*(4,2){g}
      \rput*(5,1){g}
      \rput*(8,4){g}
      \rput*(8.33,3){g}
      \rput*(8.66,2){g}
      \rput*(9,1){g}
      \rput*[Bl](12.5,1){$y=1$}
      \rput*[Bl](12.5,2){$y=2$}
      \rput*[Bl](12.5,3){$y=3$}
      \rput*[Bl](12.5,4){$y=4$}
    \end{pspicture}
    \caption{Roots of mutually interlacing polynomials}
  \end{figure}

The roots of the $f$'s lie in three increasing groups, because
each $f_i$ has three roots, and the $f_i$ are mutually
interlacing. Since the $g$'s occur in the opposite order, their groups
are decreasing. In addition, since  the $f$'s are mutually
interlacing, each group of four roots on a line does not vertically
overlap any other group of four roots.

In this diagram the dotted lines represent values of $x$ for which
each of the terms $f_1\,g_n,\cdots,f_n\,g_1$ have the same sign
slightly to the left of $x$. To see this, consider a term
$f_i\,g_{n+1-i}$. At any value $x$ and sufficiently small $\epsilon$ the
sign of $f_i(x-\epsilon)\,g_{n+1-i}(x-\epsilon)$ is $(-1)$ to the
number of roots to the right of $x$. The doted lines arise at the
intersection of upward and downward lines, or at gaps where there are
no up or down lines. 

The general case is no different. There are $n$ levels, $r$ upward
lines, and $s$ downward lines. It is easy, but notationally
complicated, to see that there are $r+s+1$ vertical lines where all
signs are the same. Between each pair of signs there is a root of the
sum, accounting for all the roots.

The final statement follows from the first since the interlacing
direction of $\{g_i(-x)\}$ is the reverse of the interlacing direction
of $ \{g_i(x)\}$.
  \end{proof}

  \section{Determinants and mutual interlacing }
  \label{sec:det-mutual-int}

\index{matrix of coefficients}
\index{mutually interlacing!determinants of coefficients}
\index{mutually interlacing!determined by a polynomial}
\index{totally positive}

If $f_1,\dots,f_d$ is a sequence of mutually interlacing polynomials
of degree $n$, and $f_i = \sum a_{i,j}x^j$ then the matrix
$(a_{d-i,j})$ is called the \emph{matrix of coefficients} of
$f_1,\dots,f_d$. If the matrix of coefficients of a sequence of
mutually interlacing polynomials is TP$_k$ (or totally positive), we
say that the sequence is TP$_k$ (or totally positive). We show that
any mutually interlacing sequence is TP$_2$, but not necessarily
TP$_3$. In addition, we show that certain mutually interlacing set of
polynomials are totally positive.

  \begin{example}\label{ex:3-mi-not-tp}
Here are  three mutually interlacing polynomials
\begin{align*}
f_1 &=    ( 2 + x ) \,( 9 + x ) \,( 14 + x ) \\
f_2&=  ( 5 + x ) \,( 10 + x ) \,( 15 + x ) \\
f_3 &=  ( 6 + x ) \,( 12 + x ) \,( 18 + x ) 
\end{align*}
The matrix of coefficients is 
\[
\begin{pmatrix}
1296&396&36&1\\ 750&275&30&1\\252&172&25&1
\end{pmatrix}
\]
It's easy to see that it is TP$_2$, but the four three by three
determinants are all negative, so it is not TP$_3$.

\end{example}

We first have two useful facts.

\begin{lemma}\
  
\begin{enumerate}
\item Suppose that $C$ is
the matrix of coefficients of $V=(f_1,f_2,\dots,f_n)$, and $M$ is a
compatible matrix, then the matrix of coefficients of $MV$ is $MC$.

\item If $g\in\allpolypos$, and the mutually interlacing sequence
  $V=(f_1,\dots,f_d)$ is totally positive, then the sequence
  $(g\,f_1,g\,f_2,\cdots,g\,f_n)$ is also totally positive.
\end{enumerate}

\end{lemma} 

\begin{proof}
  The first is immediate from the definition. For the second, since
  we can factor $g$ into positive factors, it suffices to assume that
  $g = x+\alpha$, where $\alpha>0$. If $f_1=\sum a_ix^i$ and $f_2=\sum
  b_ix^i$ then the matrices of coefficients satisfy

\begin{multline*}
\begin{pmatrix}
  a_0 & a_1 & a_2  & \cdots \\
  b_0 & b_1 & b_2  & \cdots \\
  \vdots
\end{pmatrix}
\begin{pmatrix}
  \alpha & 1 & & & \\  & \alpha & 1 & & \\  &  & \alpha &1  &  \\ & &  &  & \ddots
\end{pmatrix}
= \\
\begin{pmatrix}
  \alpha\,a_0 & \alpha\,a_1+a_0 & \alpha\,a_2+a_1 & \cdots\\
  \alpha\,b_0 & \alpha\,b_1+b_0 & \alpha\,b_2+b_1 & \cdots\\
\vdots
\end{pmatrix}
\end{multline*}

The first matrix is the matrix of coefficients of $V$, the second
is totally positive, and the third is the matrix of
coefficients of $(x+\alpha)V$. The conclusion follows.
\end{proof}

If $f=\prod(x+r_i)$ where $r_1<r_2<\cdots<r_n$ then we say that the
sequence 
\[
\frac{f(x)}{x+r_n},\frac{f(x)}{x+r_{n-1}},\cdots,\frac{f(x)}{x+r_1}
\]
of mutually interlacing polynomials is \emph{determined by $f$}.

  \begin{prop}
    If $\,V=(f_1,\dots,f_n)$ is a sequence of mutually interlacing
    polynomials determined by $f\in\allpolypos$, and $M$ is a totally
    positive matrix, then the matrix of coefficients of $MV$ is
    totally positive.
  \end{prop}
  \begin{proof}
    Since the matrix of coefficients of $MV$ is $M$ times the matrix of
    coefficients of $V$, it suffices to show that $V$ is totally positive.
    We thus may assume that $f_i = f/(x+r_{n+1-i})$, where $f(x) =
    (x+r_1)\cdots(x+r_n)$ and $0<r_1<\cdots<r_n$. The matrix of
    coefficients is an $n$ by $n$ matrix. Consider the submatrix
    determined by a set $R$ of $d$  rows and a set  of $d$
    columns. Define
\[
g(x) = \prod_{k\not\in R} (x+r_k) \qquad 
h(x) = \prod_{k\in R} (x+r_k) 
\]
Note $f_i = g(x)\cdot h(x)/(x+r_i)$. We want to prove that all the
determinants of size $d$ by $d$ are positive for the
following set of polynomials
\[
\left\{ g(x) \cdot h(x)/(x+r_k)\mid k\in R\right\}
\]
By the lemma, it suffices to show that the determinant of the matrix
of coefficients of
\[
\left\{ h(x)/(x+r_k)\mid k\in R\right\}
\]
is positive. For instance, if $d=4$ and $R=\{1,2,3,4\}$ then the matrix
of coefficients is
\[
\begin{pmatrix}
  r_2\,r_3\,r_4&r_2\,r_3 + r_2\,r_4 + r_3\,r_4&r_2 + r_3 +
  r_4&1\\
  r_1\,r_3\,r_4&r_1\,r_3 + r_1\,r_4 +
  r_3\,r_4&r_1 + r_3 + r_4&1 \\
  r_1\,r_2\,r_4&r_1\,r_2 + r_1\,r_4 +
  r_2\,r_4&r_1 + r_2 + r_4&1 \\
  r_1\,r_2\,r_3& r_1\,r_2 + r_1\,r_3 + r_2\,r_3&
  r_1 + r_2 + r_3&1 \\
\end{pmatrix}
\]

First of all, the determinant is a polynomial of degree $\binom{k}{2}$
in the $r_i$. If $r_i=r_j$ then two rows are equal, and the
determinant is zero. Thus, $r_i-r_j$ divides the determinant.
Consequently,
\[
\Delta =\prod_{i<j} (r_j - r_i)
\]
divides the determinant. But this is a polynomial of degree
$\binom{k}{2}$, so the determinant is a constant multiple of
$\Delta$. We can check that that the constant of proportionality is
$1$. Since $r_i<r_j$, all terms of the product in $\Delta$ are
positive, so the determinant is positive.

  \end{proof}

If the matrix of coefficients of a sequence $V$ of mutually interlacing
polynomials is not totally positive, then $V$ is not obtained by
multiplying the sequence determined by the factors of some polynomial
by a totally positive matrix. For instance, the polynomials in
Example~\ref{ex:3-mi-not-tp} form such a sequence.

  \begin{lemma}
    If $f_1,\dots,f_d$ is a sequence of mutually interlacing
    polynomials in $\allpolypos$ then the matrix of coefficients is
    TP$_2$.
  \end{lemma}
  \begin{proof}
    It suffices to prove that all two by two determinants are
    non-negative for $f\greateq g$. Write $g = a f+ \sum
    \alpha_i f/(x+r_i)$ where $f = \prod(x+r_i)$, and all $a_i$ are
    non-negative. If we write  this equation in terms of matrices
\[
\begin{pmatrix}
  1 & 0 & 0 & \cdots \\ a & \alpha_1 & \alpha_2 & \cdots
\end{pmatrix}
\begin{pmatrix}
  f \\ f/(x+r_n) \\ \vdots \\ f/(x+r_1) 
\end{pmatrix}
=
\begin{pmatrix}
  f\\g
\end{pmatrix}
\]
we see the matrix on the left is totally positive, and the middle
matrix is totally positive by the Proposition and taking limits, so
the conclusion follows.

  \end{proof}

  We will prove this same result in a different way in
  Corollary~\ref{cor:log-con-coef}. If we apply the Lemma to $xf$ and
  $f$ where $f=\sum a_ix^i\in\allpolypos$, then the matrix of
  coefficients is
\[
\begin{pmatrix}
  a_0 & a_1 & a_2 & \cdots \\ 0 & a_0 & a_1 & \cdots
\end{pmatrix}
\] 
The two by two determinants show that $a_k^2\ge a_{k-1}a_{k+1}$. This
is a special case of Newton's inequality \mypage{thm:newton}.

  \section{Interlacing polynomials and the Hurwitz matrix}
  \label{sec:total-posit-inter}

\index{totally positive!coefficients of $\allpolypos$}
  A pair of interlacing polynomials determines a totally positive
  matrix. We prove it inductively; the key step is given in the
  following lemma.

  \begin{lemma}\label{lem:hurwitz-factor}
    Define $T_c(f,g) = (xg+cf,f)$, and 
\[
\mathcal{D} = \left\{ (f,g)\mid f\lessless g \text{ or } f\greateq g
  \text{ in } \allpolypos \text{ with positive leading coefficients }\right\}.
\]
\begin{description}
\item[1-1:] $T_c:\mathcal{D}\longrightarrow\mathcal{D} \text{ is one to one}.$
\item[onto:] Given $\alpha\in\mathcal{D}$ there is a positive $c$ and
  $\beta\in\mathcal{D}$ such that $T_c(\beta)=\alpha$.
\end{description}
\end{lemma}
  \begin{proof}
    The first part is clear. For the second choose
    $(h,k)\in\mathcal{D}$, and set $c = h(0)/k(0)$ and $g = (h-ck)/x$.
    Since $h$ and $ck$ have the same constant term, $g$ is a
    polynomial. If $h\lessless k$ then $gx =h-ck\lessless k$ so
    $T_c(k,g) = (h,k)$ and $k\greateq g$. 

    If $h\greateq k$ then write $h = a k- \ell$ where $a$ is positive
    and $h\lessless \ell$ and follow the previous argument to conclude
    that $k\lessless g$.

  \end{proof}

  If we iterate the lemma, we see that we can find positive $c_i$ and
  positive $b$ such that 
\begin{equation}\label{eqn:hurwitz-factor}
(f,g) = T_{c_1}\, \cdots\, T_{c_n}\, H(b,0).
\end{equation}
This factorization can be turned into a factorization of matrices. 
If $f = \sum a_ix^i$, $g = \sum b_ix^i$ and we define

    \begin{equation}
      \label{eqn:tp-h-0}
J(c) =     \begin{pmatrix}
        c & 1 & 0 & 0 &0 &0& \hdots \\
        0 & 0 & 1 & 0 &0 &0&  \hdots\\
        0 & 0 & c & 1 &0 &0&  \hdots \\
        0 & 0 & 0 & 0 & 1 &0& \hdots\\
0&0&        0 & 0 & c & 1 &  \hdots \\
        \vdots & \vdots & \vdots & \vdots & \vdots & \ddots
      \end{pmatrix}
\quad
H(f,g) =     \begin{pmatrix}
        a_0 & a_1 & a_2 & \hdots \\
       0& b_0 & b_1 & b_2 & \hdots \\
       0&  a_0 & a_1 & a_2 & \hdots \\
       0&0& b_0 & b_1 & b_2 & \hdots \\
       0&0&  a_0 & a_1 & a_2 & \hdots \\
\vdots & \vdots & \vdots & \vdots & \vdots & \ddots
      \end{pmatrix}
    \end{equation}
then $J(c)H(f,g) = H(xg+cf,f)$. The matrix $H(f,g)$ is sometimes
called the \emph{Hurwitz matrix}.

\index{Hurwitz matrix}

    \begin{prop}\label{prop:hurwitz-totally-pos}
      If $f\lessless g$ in $\allpolypos$ have positive leading
      coefficients then
      \begin{enumerate}
      \item $H(f,g)$ is totally positive.
      \item There are positive $c_i$ and $b$ such that 
\[H(f,g)= J(c_1)\cdots J(c_n)H(b,0).\]
      \end{enumerate}
    \end{prop}
    \begin{proof}
      Since all $J(c)$ and $H(b,0)$ are totally positive, the first
      statement is a consequence of the second. The second is a
      consequence of the factorization \eqref{eqn:hurwitz-factor} and
      the lemma.
    \end{proof}

\begin{remark}
    If we write out the matrix determined by $f\lesslesseq f'$
  \[  \begin{pmatrix}
      a_1 & 2 a_2 & 3a_3 & \dots & \\
      a_0 & a_1 & a_2 & a_3 & \dots \\
0 &      a_1 & 2 a_2 & 3a_3 & \dots & \\
0&       a_0 & a_1 & a_2 & a_3 & \dots \\
\vdots & \vdots & \vdots & \vdots & \ddots
    \end{pmatrix}
    \]
then we see that
$\smalltwodet{ka_k}{(k+1)a_{k+1}}{a_{k-1}}{a_k}\ge0$. This is
equivalent to
\[
a_k^2 \ge \frac{k+1}{k}\ a_{k-1}a_{k+1}
\] and this is a weak form of Newton's inequality \seepage{thm:newton}.
  \end{remark}

We will prove the following corollary in a very different way later -
see Theorem~\ref{thm:tp}.

\begin{cor}\label{cor:pos-is-tp}
  If $f=\sum a_i x^i$ is in $\allpolypos$ then the matrix below is
  totally positive.
\[
\begin{pmatrix}
  a_0 & a_1 & a_2 & . & . \\
. & a_0 & a_1 & a_2 & . \\
. & . & a_0 & a_1 & . \\
\vdots & & & & \ddots
\end{pmatrix}
\]
\end{cor}
\begin{proof}
  Apply the proposition to $f\lesslesseq f'$ and select the submatrix
  whose rows come from $f$.
\end{proof}

\section{The division algorithm}
\label{sec:division-algorithm}
\index{division algorithm}
\index{Euclidean algorithm}

We study the division algorithm applied to interlacing
polynomials. The basic result is the following lemma:

\begin{lemma}\label{lem:div-alg}
  Suppose that $f,g$ have positive leading coefficients, and
  $f\lessless g$ in $\allpolypos$. If we divide $f$ by $g$
\[
f = (ax+b)g -r \quad\quad\text{where $deg(r)<deg(g)$}
\]
then
\begin{enumerate}
\item $g\lessless r$.
\item $r\in\allpolypos$
\item $r$ has positive leading coefficient.
\item ${deg}(r)=n-2$
\end{enumerate}
\end{lemma}
\begin{proof}
  Since $f$ sign interlaces $g$, it follows that $r$ sign interlaces
  $f$.  The degree of $r$ is less than the degree of $g$ so we have
  that $g\lesslesseq r$.  The interlacing is strict, since any common
  factor of $g$ and $r$ is a common factor of $f$ and $g$. Since $f$
  strictly interlaces $g$, they have no common factors. Since
  $g\in\allpolypos$, so is $r$. 

  Now we determine the leading coefficient of $r$. We may assume that
  $f,g$ are monic, so
  \begin{align*}
    f &= x^n + a_1 x^{n-1} + a_2 x^{n-2} + \cdots \\
    g &= x^{n-1} + b_1 x^{n-2} + b_2 x^{n-3} + \cdots \\
    r &= (x + a_1 - b_1) g -f \\
    &= (a_1b_1 + b_2 - b_1^2-a_2)x^{n-2} + \cdots
  \end{align*}
  \index{totally positive} \index{Hurwitz matrix}
  Using Proposition~\ref{prop:hurwitz-totally-pos} and reversing $f,g$
  we see that
\[
0 \le
\begin{vmatrix}
  1 & b_1 & b_2 \\ 1 & a_1 & a_2 \\ 0 & 1 & b_1
\end{vmatrix}
=
a_1b_1 + b_2 -b_1^2-a_2
\]
Because $r$ interlaces $g$  the degree of $r$ is $n-2$, and so the
leading coefficient is positive.
\end{proof}

\begin{remark}\added{3/11/8}
  If we apply the division algorithm to $f\lesslesseq f'$ where
  $f\in\allpolypos$ then all coefficients of $r$ must be positive. The
  resulting inequalities are either Newton's inequalities, or simple
  consequences of them. \index{Newton's inequalities}
\end{remark}

It's curious that the lemma  does not depend on the positivity of $f$ and
$g$.

\begin{lemma}\label{lem:div-alg-2}
  Suppose that $f,g$ have positive leading coefficients, and
  $f\lessless g$ in $\allpoly$. If we divide $f$ by $g$
\[
f = (ax+b)g -r \quad\quad\text{where $deg(r)<deg(g)$}
\]
then
\begin{enumerate}
\item $g\lessless r$.
\item $r$ has positive leading coefficient.
\item ${deg}(r)=n-2$
\end{enumerate}
\end{lemma}
\begin{proof}
  If we write
\[
f(x) = \prod_1^n (x-r_i) \qquad g = \sum _1^n \alpha_i\, \frac{f}{x-r_i}
\]
where $\alpha_i\ge0$ then a computation shows that the coefficient of
$r$ is
\[
- \,\frac{1}{\bigl(\alpha_1+\cdots+\alpha_n)^2}
\sum_{i<j} \alpha_i\alpha_j (r_i-r_j)^2
\]
\end{proof}

\begin{example}
  We can apply this in the usual way to solve $Ff-Gg=1$. For example,
  \begin{align*}
    f &= (x+1)(x+3)(x+6)\\
    g &= (x+2)(x+4)\\
    f &= (x+4)\bigl[g\bigr] - \bigl[
    5x+14\bigr]\\
    g &= \frac{5x+16}{25}\bigl[(5x+4)\bigr] - \frac{24}{25}\\
\intertext{and recursively solving we get}
1 &= \underbrace{\frac{39+36x+5x^2}{24}}_{G}\bigl[(x+2)(x+4)\bigr] -
\underbrace{\frac{5x+16}{24}}_{F} \bigl[(x+1)(x+3)(x+6)\bigr]
  \end{align*}
Notice that in this example $F$ and $G$ are in $\allpolypos$. In
general, if $f$ has degree $n$ then the algorithm yields $F,G$ with
$deg(F)= deg(G)-1=deg(f)-2$.
\end{example}

\begin{lemma}
  If $f\lessless g$ in $\allpolypos$, $deg(F)= deg(G)-1=deg(f)-2$, and
  $fF-gG=1$ then
  \begin{enumerate}
  \item $F,G\in\allpolypos$.
  \item $f \lessless g,G$ and $g,G \lessless F$
  \end{enumerate}
\end{lemma}
\begin{proof}
  If $\alpha$ is a root of $f$ then $g(\alpha)(G(\alpha)=1$ so $G$ and
  $g$ have the same sign at roots of $f$. Since $deg(G)=deg(f)-1$ it
  follows that $f\lesslesseq G$ - in particular, $G\in\allpoly$. The
  same arguments show $G\lesslesseq F$ and $g\lesslesseq F$. All
  interlacings are strict since any common divisor divides $1$.
\end{proof}

\index{Bezout identity}
The identity $fF-gG=1$ is known as the \emph{Bezout identity}.

\begin{cor}\label{cor:division}
  If $f\lessless g$ in $\allpolypos$ then there exist
  $F,G\in\allpolypos$ such that 
  \begin{enumerate}
  \item $fF-gG=1$
  \item  $f \lessless g,G$ and $g,G \lessless F$
  \end{enumerate}
\end{cor}

Restated in terms of the matrix $\smalltwo{f}{g}{F}{G}$

\begin{cor}
  A strict interlacing pair $(f,g)$ in $\allpolypos$ can be extended to an
  interlacing matrix of determinant $1$. 
\end{cor}

\begin{example}
  If we assume that the anti-diagonal elements are $(x-a)(x-b)$ and
  $1$ then there are exactly two solutions to the Bezout identity:
\[
\begin{pmatrix}
  x-\frac{a+b\pm \sqrt{(a-b)^2-4}}{2} & (x-a)(x-b) \\
    1 &   x-\frac{a+b\pm \sqrt{(a-b)^2-4}}{2}
\end{pmatrix}
\]
Here is a family of matrices satisfying the Bezout identity
\[
\left(
\begin{array}{ll}
 x+2 & k x^2+2 k x+2 x+3 \\
 1 & k x+2
\end{array}
\right)
\]
\end{example}

\begin{remark}
  If $M$ is a matrix as in the lemma, then $M$ and $M^{-1}$ both map
  interlacing polynomials to interlacing polynomials.
  \seepage{prop:2by2} This is different from linear transformations on
  \emph{one} polynomial, where the only linear transformations $T$ for
  which $T$ and $T^{-1}$ map $\allpoly$ to itself are affine
  transformations.
\end{remark}

\index{\ SL2R@$SL_2(\reals)$}
Next we show that all elements of $SL_2(\reals)$ can be realized by
evaluating polynomials.

\begin{lemma}
  If $\smalltwodet{a}{b}{c}{d}=1$ then there are
  $f,g,h,k\in\allpoly(n)$ such that
\[
\begin{vmatrix}
  f&g\\h&k
\end{vmatrix}=1 \quad\quad
\begin{pmatrix}
  f(0) & g(0) \\ h(0) & k(0)
\end{pmatrix}=
\begin{pmatrix}
  a & b \\ c & d
\end{pmatrix}
\]
\end{lemma}
\begin{proof}
Define
\begin{quote}
  $S$ is the set of all $\smalltwo{f}{g}{h}{k}$ such that
  $\smalltwodet{f}{g}{h}{k}=1$ and all polynomials in the same row or
  column interlace.
\end{quote}
Note that if $\smalltwodet{f}{g}{h}{k}=1$ and polynomials in one row
or column interlace then $\smalltwo{f}{g}{h}{k}\in S$. We now claim
that multiplication by $\smalltwo{\alpha}{0}{\beta}{1/\alpha}$ and
$\smalltwo{\alpha}{\beta}{0}{1/\alpha}$ map $S$ to itself. For
instance,
\[
\begin{pmatrix}
  \alpha & 0 \\ \beta & 1/\alpha
\end{pmatrix}
\begin{pmatrix}
  f&g\\h&k
\end{pmatrix}
=
\begin{pmatrix}
  \alpha f & \alpha g \\ \beta f + h/\alpha & \beta g + k/\alpha
\end{pmatrix}
\]
The resulting matrix has determinant one, and one interlacing row, so
it is in $S$.

From Corollary~\ref{cor:division} we know that $S$ is non-zero, so we choose
$\smalltwo{f}{g}{h}{k}$ in $S$. The element
$\smalltwo{f(0)}{g(0)}{h(0)}{k(0)}$ is in $SL_2(\reals)$. Since
$SL_2(\reals)$ is generated by upper and lower triangular matrices as
above, all elements of $SL_2(\reals)$ have the desired representation.
  \end{proof}

Following the usual matrix representation of the Euclidean algorithm
we get a unique representation for interlacing polynomials.
\index{Euclidean algorithm}

\begin{lemma}
  If $f\lesslesseq g$ have positive leading coefficients then there
  are unique positive $a_1,\dots,a_n$ and positive $\alpha$ such that
\[
\begin{pmatrix}
  f\\g
\end{pmatrix}
=
\alpha
\begin{pmatrix}
  a_nx+b_n & -1 \\ 1 & 0
\end{pmatrix}
\cdots
\begin{pmatrix}
  a_1x+b_1 & -1 \\ 1 & 0
\end{pmatrix}
\begin{pmatrix}
  1\\0
\end{pmatrix}
\]
\end{lemma}
\begin{proof}
  Simply note that if $f = (ax+b)g-r$ then
\[
\begin{pmatrix}
  f\\g 
\end{pmatrix}
=
\begin{pmatrix}
  ax+b& -1 \\1 & 0
\end{pmatrix}
\begin{pmatrix}
  g\\r
\end{pmatrix}
\]
\end{proof}

  \section{Integrating families of polynomials}
\label{sec:int-fam}

We can extend the concept of a mutually interlacing sequence
$f_1,\dots,f_n$ of polynomials to an interlacing family $\{f_t\}$,
where $t$ belongs to some interval. Most results of  section~\ref{sec:det-mutual-int}
have straightforward generalizations, where summation is replaced by
integration, and totally positive matrices are replaced by totally
positive functions. 

We first introduce a quantitative measure of root separation, and use it
to integrate families of polynomials. 
If
$f\in\allpoly$, then $\delta(f)$ is the minimum distance between roots
of $f$. If $f$ has repeated roots then $\delta(f)=0$, and so $f$ has all
distinct roots if and only if $\delta(f)>0$.  It is easy to see that an
equivalent definition is
\begin{align}
 \delta(f) & = \text{ largest $t$ such that }  f(x)\greateqeq f(x+t)\notag\\
\intertext{It easily follows from the definition of $\delta(f)$ that}
  0\le s \le t \le \delta(f) &\implies f(x+s)\greateqeq f(x+t)
\label{eqn:sep-prop}
\end{align}

This last property is the motivation for the following definitions.

\index{family!of interlacing polynomials}
\index{interlacing family}

\begin{definition}
  Suppose that $I$ is an interval, and for each $t\in I$ we have a
  polynomial $f_t(x)$ with all real roots such that the assignment $t
  \mapsto f_t$ is continuous.  We say $\{f_t\}$ is a \emph{continuous
    family of polynomials}.  If $f_t \greateqeq f_s$ for every $t \le
  s$ in $I$, then $\{f_t\}$ is an \emph{interlacing family}. The
  \emph{interval of interlacing at a point $v$}, written
  $\rho(f_t,v)$, is the largest $\epsilon$ such that $f_v\greateqeq f_s$
  for all $s\in(v,v+\epsilon)$. If there is no such $\epsilon$ we set
  $\rho(f_t,v)=0$.  A continuous family of polynomials $\{f_t\}$ is
  \emph{locally interlacing} if $\rho(f_t,v)$ is positive for all $v$
  in the interior of $I$. The interval of interlacing of a locally
  interlacing family, $\rho(f_t)$, is the largest $\epsilon$ such that
  if $s,v\in I$ and $v < s < v+\epsilon$, then $f_v \greateqeq f_s$. If
  there is no such $\epsilon$ then the interval of interlacing is $0$.
  
\end{definition}

The following is a generalization of Corollary~\ref{cor:faf}.

\begin{prop} \label{prop:family-int}
  Suppose $\{f_t\}$ is a continuous family of interlacing polynomials
  on $[a,b]$, and that $w(t)$ is a positive function on $[a,b]$. The
  polynomial
  \begin{equation}
    \label{eqn:integral-1}
    p(x) = \int_a^b f_t(x)\ w(t)dt    
  \end{equation}
has all real roots, and $f_a \greateqeq p \greateqeq f_b$.
\end{prop}

\begin{proof}
  Choose $a = t_0 < t_1 < \dots < t_n = b$, and consider the 
  function 
$$  s(x,t) = f_{t_0}(x)w(t)\chi_{[t_0,t_1]} + \dots +
f_{t_{n-1}}(x)w(t)\chi_{[t_{n-1},t_n]} 
$$
where $\chi_{[c,d]}$ is the characteristic function of the interval
$[c,d]$.  It suffices to show that the polynomial $\int_a^b
s(x,t)\,dt$ has all real roots, since \eqref{eqn:integral-1} is the limit of
such functions.  Note that
$$
\int_a^b s(x,t)\,dt = f_{t_0}\cdot \int_{t_0}^{t_1} w(t)\,dt + \dots +
f_{t_{n-1}}\cdot \int_{t_{n-1}}^{t_n} w(t)\,dt
$$
Each of the coefficients of the polynomials $f_{t_i}$ is positive,
and the set of polynomials $f_{t_0},\dots,f_{t_n}$ is mutually
interlacing since $\{f_t\}$ is an interlacing family, so $\int_a^b
s(x,t)\,dt$ has all real roots by Lemma~\ref{lem:ti}.
\end{proof}

\begin{example}
  Since $\delta\falling{x}{n} = 1$, if we let $w(t)=1$ then $ \int_0^1
  \falling{x+t}{(n)}\,dt$ is in $\allpoly.$ More generally, for any
  $w(t)$ that is positive on $(0,1)$ we have

$$ \falling{x}{n} \greateqeq \int_0^1\falling{x+t}{(n)}\,w(t)\,dt
\greateqeq \falling{x+1}{(n)}.$$ 
\end{example}

\begin{example}\label{ex:interpolating-mutual}
  Given any mutually interlacing sequence of polynomials we can find
  an  interlacing family that contains it. Suppose that $f_1\greateqeq
  \cdots\greateqeq f_n$ is a mutually interlacing family. Define
an interlacing family on $(1,n)$ such that  $F_i = f_i$ by 
$$ F_t = (i+1-t)f_i + (t-i)f_{i+1} \text{ for } i \le t \le i+1.$$
We can easily compute the integral 
\begin{align*}
  \int_1^n F_t\ dt &= \sum_{i=1}^{n-1} \int_i^{i+1}
  (i+1-t)f_i+(t-i)f_{i+1}\,dt\\
&= \sum_{i=1}^{n-1} (\frac{1}{2} f_i + \frac{1}{2}f_{i+1})\\
&= \frac{1}{2} f_1 + f_2 + \cdots + f_{n-1} + \frac{1}{2}f_n
\end{align*}

Here is a similar example, where the sequence of mutually interlacing
polynomials is determined by the roots of a fixed polynomial. Suppose
that $f(x) = (x-r_1)\cdots(x-r_n)$ where $r_1\le\cdots\le r_n$, and
let $f_i = f/(x-r_i)$. We know that $ f_n \greateqeq \cdots \greateqeq
f_1$ is a mutually interlacing sequence. Define an interlacing family
on $(r_1,r_n)$ by
$$ F_t = \frac{f(x)}{(x-r_i)(x-r_{i+1})}\,(x-r_i-r_{i+1}+t) \text{ on
} r_i\le t \le r_{i+1}$$
This family satisfies $F_{r_i} = f_i$. The integral is equal to
\begin{align*}
  \int_{r_1}^{r_n} F_t\,dt &= \sum_{i=1}^{n-1} \int_{r_i}^{r_{i+1}}
   \frac{f(x)}{(x-r_i)(x-r_{i+1})}\,(x-r_i-r_{i+1}+t) \\
&= \sum_{i=1}^{n-1}
(r_{i+1}-r_i)\,\frac{f(x)}{(x-r_i)(x-r_{i+1})}(x-(r_i+r_{i+1})/{2})
\\
\end{align*}
The integral interlaces $f$, since the integral is a positive linear
combination of polynomials interlacing $f$.

In neither of these examples is $F_t$ differentiable; the
constructions involve piecewise linear functions. In
Lemma~\ref{lem:interp-mutual} we will find a differentiable family for
the sequence in the second example.
\end{example}

The next result is the continuous analog of
Lemma~\ref{lem:int-fam-int-1}. The proof is omitted.

\begin{lemma}\label{lem:int-fam-int-2}
  Suppose that $f_t$ and $g_t$ are  interlacing families on
  $[a,b]$. If
  \begin{enumerate}
  \item $f_t(x) \greateq g_t(x)$ for all $a\le t \le b$
  \item $\smalltwodet{f(s)}{f(t)}{g(s)}{g(t)} <0 $ for $a\le s < t \le
    b$
  \end{enumerate}
then
$$ \int_a^b f_t(x)\,dt \greateqeq \int_a^b g_t(x)\,dt$$
\end{lemma}

The conclusion of Proposition~\ref{prop:family-int} was only that $p$
weakly interlaced $f_a$ and $f_b$, since the limit of interlacing
polynomials is only weakly interlacing polynomials.  The next lemma is
used to deduce strict interlacing for integrals.

\begin{lemma}
  Suppose $\{f_t\}$ is a continuous family of interlacing polynomials
  on $[a,b]$, and that $w(t)$ is a positive function on $[a,b]$.
  Choose $c_1,c_2$ such that $a < c_1 < c_2 < b$ and consider the the
  polynomials
  \begin{eqnarray*}
    p_1  & = \displaystyle \int_a^{c_1} f_t(x)\ w(t)dt \\
    p_2  & = \displaystyle \int_{c_2}^{b} f_t(x)\ w(t)dt 
  \end{eqnarray*}
Then, $f_a \greateqeq p_1 \greateq p_2 \greateqeq f_b$.
\end{lemma}

\begin{proof}
  From the previous theorem and the hypothesis we know that $$ f_a
  \greateqeq p_1 \greateqeq f_{c_1} \greateq f_{c_2} \greateqeq p_2
  \greateqeq f_b$$
Since $\{f_t\}$ is interlacing on $[a,b]$, we see that $p_1 \greateq p_2$.
\end{proof}

We can integrate locally interlacing families to get new ones. We
first hold the interval constant. 

\begin{prop}\label{prop:family-int-2}
  If $\{f_t\}$ is a locally interlacing family on $[a,b]$, $w(x)$
  positive on $[a,b]$, and $r \le \rho(f)/2$, then 
$$ g_s(x) = \int_s^{s+r} f_t(x)w(t)\,dt$$
is a locally interlacing family on $[a,b-r]$ with $\rho(g_s)\ge r$.
\end{prop}

\begin{proof}
Certainly all members $g_s$ have roots, so it remains to show that
there is local interlacing. Consider $g_s$ and $g_{s+\epsilon}$, where
$0< \epsilon\le r$.  By the previous lemma, the three integrals 

\begin{align*}
  p_1 &= \int_s^{s+\epsilon} f_t(x)w(t)\,dt \\
  p_2 &= \int_{s+\epsilon}^{s+r} f_t(x)w(t)\,dt \\
  p_3 &= \int_{s+r}^{s+r+\epsilon} f_t(x)w(t)\,dt \\
\end{align*}
\noindent%
satisfy $p_1 \greateqeq p_2 \greateqeq p_3$ and $p_1 \greateq p_3$.
Since $g_s = p_1 + p_2$, and $g_{s+\epsilon} =p_2+p_3$, it follows
that $g_s \greateq g_{s+\epsilon}$.

\end{proof}

Next, we vary the domain of integration.

\begin{lemma}\label{lem:family-int-2a}
  If $\{f_t(x)\}$ is a locally interlacing family on $\reals$, $w(x)$
  a non-negative function, and
  $$
  g_s(x) = \begin{cases} \displaystyle\int_0^s f_t(x)w(t)\,dt & s>
    0 \\ f_0(x) & s=0\end{cases}$$
  then $\{g_s\}$ is a locally
  interlacing family with $\rho(g)\ge \rho(f)$.
\end{lemma}
\begin{proof}
We need to show that $g_u \greateqeq g_v$ for $0\le u \le v \le \rho(f)$.  
  If we consider an approximating sum for $g_v$ then all the terms are
  mutually interlacing, which easily implies the result. 
\end{proof}

If the weight is a function of two variables, then integrating a
locally interlacing family against the weight gives a new family whose
members all have roots.  However, we need an extra condition to
guarantee that it is locally interlacing, and total positivity is what
we need.

\begin{lemma}
  Suppose that $\{f_t\}$ is a family of interlacing polynomials on
  $[a,b]$, and $u(x),w(x)$ are positive on $[a,b]$. Set
  \begin{align*}
    p_u &= \int_a^b f_t(x)\,u(t)\,dt \\
    p_w &= \int_a^b f_t(x)\,w(t)\,dt. 
  \end{align*}
If $\smalltwodet{u(s)}{u(t)}{w(s)}{w(t)}$ is positive for all $a\le
s<t\le b$, then $p_u \greateqeq p_w$ 
\end{lemma}

\begin{proof}
  For any $a < t_1 < \dots < t_n < b$, the approximations to $p_u$ and
  $p_w$ are
  \begin{align*}
    \sum_{i=1}^{n-1} \left(\int_{t_i}^{t_{i+1}} u(t)\,dt\right)\
    f_{t_i} \\
    \sum_{i=1}^{n-1} \left(\int_{t_i}^{t_{i+1}} w(t)\,dt\right)\
    f_{t_i} 
  \end{align*}
Now since 

\begin{multline*}
\left(\int_{t_i}^{t_{i+1}} w(t)\,dt\right)
\left(\int_{t_j}^{t_{j+1}} u(t)\,dt\right) - 
\left(\int_{t_j}^{t_{j+1}} w(t)\,dt\right)
\left(\int_{t_i}^{t_{i+1}} u(t)\,dt\right) = \\
\int_{t_i}^{t_{i+1}}\int_{t_j}^{t_{j+1}}(w(t)u(s) - u(t)w(s))\,dt\,ds
> 0
\end{multline*}
and the sequence of polynomials $f_{t_1},\dots,f_{t_n}$ is mutually
interlacing, the result follows from Theorem~\ref{thm:totally-preserves}.  Note
that we can allow $w$ and $u$ to have infinitely many zeros on $[a,b]$
since all the integrals will still be positive.
\end{proof}

\index{mutually interlacing}

\index{totally positive}
\index{totally positive!strongly, of order 2}
\begin{definition}
  A function $w(s,t)$ is \emph{strongly totally positive of order 2 (or $STP_2$)
  } on an interval $[a,b]$ if for all $t<t^\prime,s<s^\prime$ in
  $[a,b]$ the determinant
$$ 
\begin{vmatrix}
  w(s,t) & w(s,t^\prime) \\
  w(s^\prime,t) & w(s^\prime,t^\prime) 
\end{vmatrix}
$$
and all its entries are positive. For example, $e^{xy}$ is $STP_2$.
\end{definition}

\begin{cor}
  Suppose that $\{f_t\}$ is an interlacing  family on $[a,b]$, and
  that $w(s,t)$ is $STP_2$.  The family of polynomials
$$ g_s = \int_a^b f_t(x)w(t,s)\,dt$$
is a  interlacing family.
\end{cor}  
\begin{proof}
  The lemma shows that it is locally interlacing, and it follows that
  it is an interlacing family since $f_a \greateqeq g_s \greateqeq f_b$.
\end{proof}

As an example, consider the family
$$
g_t(x) = \int_0^1 e^{ty} \falling{x+y}{n}\,dy$$
Since
$\delta\falling{x+y}{n}=1$  and $e^{ty}$ is $STP_2$ we know
that $g_t$ is an interlacing family.  Since $g_t$ is defined for all
$t$, we have an interlacing family on $\reals$.

\section{Families, root zones, and common interlacings}
  \label{sec:families-root-zones}

  In this section we consider families of polynomials that are indexed
  by a set $S$. We say that a family $\{f_t\}_{t\in S}$ has a common
  interlacing if there is a $g\in\allpoly$ such that $f_t\lesslesseq
  g$ for all $t\in\ S$. The basic question is
  \begin{quote}
    What properties do we need in order to have a common interlacing?
  \end{quote} \index{common interlacing}

  Here's an  example.  Consider the family $\{F_t\}_{t\ge0}$ where
  $F_t(x) = f(x) + t g(x)$ and $deg(f)=deg(g)+1$. If
  $F_t\in\allpoly$ for all positive $t$ then there is a common
  interlacing. This is Proposition~\ref{prop:pattern}.

  \index{root zone} It is useful to consider all the zeros of a
  family. Suppose that $\{f_t\}_{t\in S}$ is a family of polynomials
  of degree $n$. The $k$'th \emph{root zone} is
\[ \Delta_k = \bigl\{ w\in\reals\mid \exists t\in S \And \text{$w$
is the k'th largest root of $f_t$}\bigr\}
\]

Since we assume that families are continuous, it follows that if $S$
is connected then all the root zones are intervals.  The following is
elementary, but points out how root zones can be useful.

\begin{lemma}\label{lem:root-zones}
  If a family $F$ consists of polynomials in $\allpoly$ of
  degree $n$,  and all root zones
  $\Delta_1,\dots,\Delta_n$ are disjoint then $F$ has a common
  interlacing. The family  has a common interlacing
  iff the intersection of any two root zones has at most one element.
\end{lemma}
\begin{proof}
  Pick the roots of a common interlacing to lie between consecutive
  root zones.
\end{proof}

\begin{cor}\label{cor:root-zones}
  If any two polynomials of a family have a common interlacing, then
  the family has a common interlacing.
\end{cor}
\begin{proof}
  The hypotheses imply that all elements in $\Delta_k$ are less than
  or equal to all elements of $\Delta_{k+1}$. 
\end{proof}

\begin{example}
  The root zones of the family $\{f+tg\}_{t\ge0}$ mentioned above have a
  simple description. Suppose that $f$ and $g$ have degree $n$, and
  let $\roots(f) = (r_i)$ and $\roots(g) = (s_i)$. Since $f,g$ have a
  common interlacing, it follows that the root zones are
\[
\Delta_k = \bigl[ \min(r_k,s_k),\,\max(r_k,s_k)\bigr].
\]

\end{example}

If all members of the family $\{f+\alpha g\}_{\alpha\ge0}$ have all
real roots then the family has a common interlacing. We can rephrase
this in a way that leads to an important generalization:
\begin{quote}
  If $\beta^2 f + g\in\allpoly$ for all real $\beta$ then
  $\{\beta^2f+g\}_{\beta\in\reals}$ has a common interlacing.
\end{quote}
The next lemma adds a linear term.

\begin{lemma}\label{lem:mat-poly-2}
  Suppose that $f,g,h$ are polynomials with positive leading
  coefficients such that $deg(f)=deg(g)=deg(h)+1$
  \begin{quote}
    
    If $\beta^2 f + \beta h + g\in\allpoly$ for all real $\beta$ then
    $\{\beta^2f+\beta h + g\}_{\beta\in\reals}$ has a common
    interlacing.
  \end{quote}
\end{lemma}
\begin{proof}
  Define $F(x,\beta) = \beta^2f+\beta h + g$. If $deg(f)=n$ then for
  each $\beta$ the polynomial $F(x,\beta)$ has exactly $n$ roots,
  since $F(x,\beta)$ has degree $n$ for all $\beta$. If $r_i(\beta)$
  is the $i$'th largest root of $F(x,\beta)$ then the continuous
  function $\beta\mapsto r_i(\beta)$ is called the $i$'th solution
  curve. 

  The $i$'th solution curve meets the $x$-axis in the $i$'th largest
  root of $g$ since $F(x,0)=g$. In addition,
  $\lim_{|\beta|\rightarrow\infty} r_i(\beta)$ is the $i$'th largest
  root of $f$. Thus, the $i$'th solution curve is asymptotic to the
  vertical line through the $i$'th largest root of $f$, and meets the
  $x$-axis in the $i$'th largest root of $g$.

  Consideration of the geometry will finish the proof. For any $x_0$
  there are at most two points on the graph with $x$-coordinate
  $x_0$. There can be no vertical line that meets two solution curves,
  for then there would be at least four intersections on some
  vertical line. See Figure~\ref{fig:ci-2}, where the dashed lines are
  the vertical asymptotes. Thus, the root zones are separated, and we can apply
  Lemma~\ref{lem:root-zones}.

  \begin{figure}[htbp]
    \centering
  \begin{pspicture}(0,-2)(3,2)
    \psline[linestyle=dashed](0,-2)(0,2)
    \psline[linestyle=dashed](3,-2)(3,2)
    \pscurve(.1,-2)(.5,-1)(2,0)(.5,1)(.1,2)
    \pscurve(2.9,-2)(2.5,-1)(1,0)(2.5,1)(2.9,2)
  \end{pspicture}
    
    \caption{Impossible solution curves}
    \label{fig:ci-2}
  \end{figure}
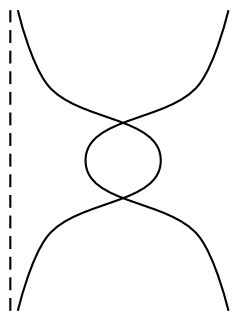
\end{proof}

  \begin{cor}
    Suppose that $f,h$ are two polynomials satisfying
    \begin{enumerate}
    \item $f\in\allpoly$ has positive leading coefficient.
    \item $h$ has degree less than $n$, and does not necessarily have
      all real roots.
    \item $f+\alpha h\in\allpoly$ for $0 \le\alpha\le1$
    \end{enumerate}
then $\{f+\alpha h\}_{0\le\alpha\le1}$ has a common interlacing.
  \end{cor}
  \begin{proof}
    Consider the following
    \begin{align*}
      f + \alpha h\in\allpoly &\qquad \text{if }\ |\alpha|\le 1 \\
f + \frac{2\beta}{\beta^2+1}h\in\allpoly &\qquad \text{if }\ \beta\in\reals\
\text{since}\ |\frac{2\beta}{\beta^2+1}|\le 1\\
\beta^2 f + 2\beta h + f \in\allpoly &\qquad \text{if }\ \beta\in\reals
    \end{align*}
We can now apply the lemma.
  \end{proof}

  \section{Two by two matrices preserving interlacing}
\label{sec:small-preserve}

\index{matrix!preserving interlacing} When does a two by two matrix
with polynomial entries preserve interlacing?  First, a precise
definition.

\begin{definition}
  Suppose that $M$ is a two by two matrix with polynomial entries, and
  $M\smalltwobyone{f}{g} = \smalltwobyone{r}{s}$. Then
  \begin{itemize}
  \item $M$ preserves interlacing if 
$f\lessless g \text{ or } f \greateq g \implies r\lessless s \text{ or } r \greateq s.$
  \item $M$ preserves interlacing for $\lessless$ if 
$f\lessless g  \implies r\lessless s \text{ or } r \greateq s.$
  \item $M$ preserves interlacing for $\greateq$ if 
$ f \greateq g \implies r\lessless s \text{ or } r \greateq s.$
  \item $M$ weakly preserves interlacing if all interlacings are
    $\lesslesseq$ or $\greateqeq$. 
  \end{itemize}
\end{definition}

A number of earlier results can be cast as matrix statements. Assume
that $\alpha,\beta$ are positive.

\begin{xalignat*}{2}
  f\longleftarrow g &\implies \alpha f \longleftarrow \beta g &
  \smalltwo{\alpha}{0}{0}{\beta} & \text{ preserves interlacing.}\\
  f\longleftarrow g &\implies x f \longleftarrow x g &
  \smalltwo{x}{0}{0}{x} & \text{ weakly preserves interlacing.}\\
  f\lessless g &\implies xf-\alpha g \lessless \beta f&
  \smalltwo{x}{-\alpha}{\beta}{0} & \text{ preserves interlacing for $\lessless$.}\\
  f\lessless g &\implies xg+\alpha f \lessless \beta g&
  \smalltwo{\alpha}{x}{0}{\beta} & \text{ preserves interlacing for $\lessless$.}
\end{xalignat*}

We can combine these matrices using the following lemma:

\begin{lemma}
  Suppose that $u,v,w$ are vectors of length two. If
  $\smalltwobyone{u}{v}$ and $\smalltwobyone{u}{w}$ preserve
  interlacing, then for any positive $\alpha,\beta$ the matrix
  $\smalltwobyone{u}{\alpha v+\beta w}$ preserves interlacing. If one
  of them weakly preserves interlacing, then the conclusion still
  holds.
\end{lemma}
\begin{proof}
  Let $f\longleftarrow g$, $\phi=\smalltwobyone{f}{g}$,
  $\smalltwobyone{u}{v}\phi= \smalltwobyone{r}{s}$ and
  $\smalltwobyone{u}{w}\phi=\smalltwobyone{r}{t}$. Then
  $\smalltwobyone{u}{\alpha v+\beta w}\phi=\smalltwobyone{r}{\alpha
    s+\beta t}$. Since $r\longleftarrow s$ and $r\longleftarrow t$ and
  at least one of these is strict, the result follows by adding
  interlacings. 
\end{proof}

\begin{cor} \label{cor:2by2}
  Suppose that $\alpha,\beta$ are positive.
  \begin{enumerate}
  \item $\smalltwo{\alpha}{x}{0}{\beta}$ preserves interlacing.
  \item $\smalltwo{x}{0}{\alpha}{x}$ preserves interlacing.
  \item $\smalltwo{x}{-\alpha}{0}{x}$  preserves interlacing. 
  \item $\smalltwo{x}{-\alpha}{\beta}{x}$  preserves interlacing.
  \end{enumerate}
\end{cor}
\begin{proof}
  The two matrices on the left satisfy the conditions of the lemma, so
  the matrix on the right preserves interlacing:

  \begin{alignat*}{2}
    \smalltwo{\alpha}{0}{0}{\beta} & \And \smalltwo{0}{x}{0}{\beta} &
    \implies & \smalltwo{\alpha}{x}{0}{\beta}\\
    \smalltwo{x}{0}{0}{x} & \And \smalltwo{x}{0}{\alpha}{0} &
    \implies & \smalltwo{x}{0}{\alpha}{x}
\intertext{In the next one, the second matrix represents the
  interlacing $\epsilon xg-\alpha g\greateqeq g$, which follows from
  $\epsilon x g \greateqeq xg \lesslesseq g$. We then take the limit as $\epsilon\rightarrow0^+$.}
    \smalltwo{x}{0}{0}{x} & \And \smalltwo{0}{\epsilon x - \alpha}{0}{x} &
    \implies & \smalltwo{x}{-\alpha}{0}{x}\\
\intertext{The next one follows from previous ones:}
    \smalltwo{x}{-\alpha}{\beta}{0} & \And
    \smalltwo{x}{-\alpha}{0}{x} &
    \implies & \smalltwo{x}{-\alpha}{\beta}{x}
  \end{alignat*}

\end{proof}

The condition that a matrix preserves interlacing is not much more
than invertibility. Note that the signs of the leading coefficients
are not relevant.

\begin{lemma}
  If the two by two matrix $M$ is invertible for all $x$, and there
  are polynomials $r\lessless s$ and $u\lessless v$ such that 
$M\smalltwobyone{r}{s} = \smalltwobyone{u}{v}$ then $M$ preserves
interlacing for $\lessless$.
\end{lemma}
\begin{proof}
  Assume that there is an $\alpha$ such that
  $u(\alpha)=v(\alpha)=0$. Since $M$ evaluated at $\alpha$ is
  invertible, it follows that $r(\alpha)=s(\alpha)=0$. Since $r$ and
  $s$ have no roots in common, it follows that $u$ and $v$ also have
  no roots in common. Thus, whatever the order the zeros of $u$ and
  $v$, it is the same order for all $u$ and $v$. Since there is one
  example where they interlace, they always interlace.
\end{proof}

\begin{prop} \label{prop:2by2}
  Suppose that $M=\smalltwo{f}{g}{h}{k}$ is a matrix of polynomials
  with positive leading coefficients.  $M$ preserves interlacing if
  and only if 
  \begin{enumerate}
  \item $f\longrightarrow g$, $h\longrightarrow k$, and both
    interlacings are strict.
  \item The determinant $\smalltwodet{f}{g}{h}{k}$ is never zero.
  \end{enumerate}
\end{prop}
\begin{proof}
  Assume that $M$ preserves interlacing, and let
  $M\smalltwobyone{r}{s} = \smalltwobyone{u}{v}$. If $r\greateq s$
  then $fr + gs\in\allpoly$. It follows from Lemma~\ref{lem:pattern2}
  that $g\greateqeq f$. If $g$ and $f$ had a common factor, then $u$
  and $v$ would have a common factor. Thus $g\greateq f$. Similarly,
  $k\greateq h$. Suppose that $\rho$ is a root of $|M|$. It follows
  that $\smalltwobyone{u(\rho)}{v(\rho)}=\smalltwobyone{0}{0}$ which
  contradicts the hypothesis that $u$ and $v$ have no roots in
  common. 
  
  Conversely, by Leibnitz (Lemma~\ref{lem:leibnitz}) both $u=fr+gs$ and
  $v=hr+ks$ have all real roots.  By the preceding lemma, it suffices
  to show that there is one example where the roots of $u$ and $v$
  interlace.  For this purpose, choose $r=g$ and $s=f$ and apply
  Lemma~\ref{lem:int-4}.

  The cases with $\lessless$ are similar and omitted.
\end{proof}

\begin{cor}\label{cor:2by2-a}
  If $f\lessless g$ and the absolute value of the leading coefficient
  of $f$ is greater than the absolute value of the leading coefficient
  of $g$ then $\smalltwo{f}{-g}{g}{f}$ preserves interlacing for
  $\lessless$. 
\end{cor}
\begin{proof}
  Since $f\lessless g$ the determinant is always positive.  The
  following interlaces by Lemma~\ref{lem:fg-2}, and so proves the
  corollary:
$$ 
\begin{pmatrix}
  f&-g\\g&f
\end{pmatrix}
\begin{pmatrix}
  f\\g
\end{pmatrix}=
\begin{pmatrix}
  f^2-g^2 \\2\,fg
\end{pmatrix}
$$
\end{proof}

\begin{example} \label{ex:preserve-3}
  If we pick $f\greateq g$ then $\smalltwodet{g}{f}{g'}{f'}\ne0$.
  Consequently, $\smalltwo{g}{f}{g'}{f'}$ preserves interlacing. For
  instance, $\smalltwo{2x}{x^2}{1}{2x}$ weakly preserves interlacing.
\end{example}

\begin{example}\label{ex:powerofi}
    If we expand $(x+\imag)^n$ we can write it as
    $f_n(x)+g_n(x)\imag$ where $f_n$ and $g_n$ have all real
    coefficients. Multiplication by $x+\imag$ leads to the recurrence
$$
\begin{pmatrix}
  x & -1 \\ 1 & x 
\end{pmatrix}
\begin{pmatrix}
  f_n \\g_n
\end{pmatrix}=
\begin{pmatrix}
  f_{n+1} \\ g_{n+1}
\end{pmatrix}
$$
By Corollary~\ref{cor:2by2} it follows that $f_n \lessless g_n$. This
is an elementary fact about polynomials whose roots all have negative
imaginary part. See Chapter~\ref{cha:complex-coef}.
  \end{example}

\begin{example}\label{ex:even-odd}
  Here is another example where matrix multiplication naturally arises
  (see Theorem~\ref{thm:hurwitz}).  Suppose that $f(x)$ is a
  polynomial.  Write $f(x) = f_e(x^2) + xf_o(x^2)$. The even part of
  $f$ is $f_e(x)$, and the odd part is $f_o(x)$. For example, consider
  powers of $x+1$: \\[.2cm]

\index{even part}
\index{odd part}

\begin{center}
\begin{tabular}[here]{cccc}
\toprule
& $(x+1)^3$ & $(x+1)^4$ & $(x+1)^5$ \\[.1cm]
\midrule
even part & $3x+1$ & $x^2+6x+1$ & $5x^2+10x+1$ \\[.1cm]
odd part & $x+3$ & $4x+4$ & $x^2+10x+5$\\
\bottomrule
\end{tabular}
\end{center}

It's clear from the definition that if $g(x) = (x+\alpha)\,f(x)$ then

\begin{align*}
  g(x) &= \alpha\,f_e(x^2) + x^2\,f_o(x^2) + x\,f_e(x^2) + \alpha\,x\,f_o(x^2)\\
  g_e(x) & =  \alpha\, f_e(x) + x f_o(x) \\
  g_o(x) & =   f_e(x) + \alpha\, f_o(x) 
\end{align*}
{which we can express as the matrix product}
\[
\begin{pmatrix}
  \alpha & x \\ 1 & \alpha
\end{pmatrix}
\begin{pmatrix}
  f_e \\ f_o
\end{pmatrix}=
\begin{pmatrix}
  g_e \\g_o
\end{pmatrix}\\
\]

We will see in Example~\ref{ex:hri} that if $h \longleftarrow  k$
have positive leading coefficients, and $h,k$ lie in $\allpolypos$ then  

$$
\begin{pmatrix}
  1 & x \\ 1 & 1
\end{pmatrix}
\begin{pmatrix}
  h \\ k
\end{pmatrix}=
\begin{pmatrix}
  r \\s
\end{pmatrix}
$$
implies that $r \longleftarrow s$. However, the matrix does not
preserve arbitrary pairs of interlacing polynomials:
$$
\begin{pmatrix}
  1 & x \\ 1 & 1
\end{pmatrix}
\begin{pmatrix}
(x + 2)(x - 2)\\ (x - 1)(x + 3)
\end{pmatrix}=
\begin{pmatrix}
  -4 - 3\ x + 3\ x^2 + x^3\\ -7 + 2\ x + 2\ x^2
\end{pmatrix}
$$
The latter two polynomials are in $\allpoly$, but they don't
interlace. However, it is worth noticing that 
$$
\begin{pmatrix}
  1 & x \\ 1 & 1
\end{pmatrix}^2
\begin{pmatrix}
(x + 2)(x - 2)\\ (x - 1)(x + 3)
\end{pmatrix}=
\begin{pmatrix}
-4 - 10\ x + 5\ x^2 + 3\ x^3 \\ -11 - x + 5\ x^2 + x^3
\end{pmatrix}
$$
and the latter two polynomials do interlace.  This is a general fact:

\begin{lemma}
  If $\alpha>0$ then the matrix $\smalltwo{\alpha}{x}{1}{\alpha}^2$
  weakly preserves interlacing.
\end{lemma}
\begin{proof}
  First of all, $\smalltwo{\alpha}{x}{1}{\alpha}^2 =
  \smalltwo{x+\alpha^2}{2\alpha x}{2\alpha}{x+\alpha^2}$. We notice
  that 
$$\begin{pmatrix}{x+\alpha^2}&{2\alpha
    x}\\{2\alpha}&{x+\alpha^2}
\end{pmatrix}
\begin{pmatrix}
  r\\s
\end{pmatrix}=
  \begin{pmatrix}{(x+\alpha^2)r+2\alpha xs}\\{2\alpha r+(x+\alpha^2)s}
  \end{pmatrix}=
  \begin{pmatrix}
    u\\v
  \end{pmatrix}
$$  
and both $u$ and $v$
  are in $\allpoly$.
  
  If we follow the proof of Proposition~\ref{prop:2by2} then
  $\smalltwodet{x+\alpha^2}{2\alpha x}{2\alpha}{x+\alpha^2} =
  \smalltwodet{\alpha}{x}{1}{\alpha}^2 = (x-\alpha)^2$ has a zero at
  $x=\alpha$, and is otherwise positive. As long as $u$ and $v$ do not
  have $\alpha$ for a zero then the argument of the proposition
  applies.  This is the case if and only if $r(\alpha)+s(\alpha)=0$.
  If we choose $\epsilon$ so that $s_\epsilon(x)=s(x+\epsilon)$
  satisfies $r\greateq s_\epsilon$,
  $r(\alpha)+s_\epsilon(\alpha)\ne0$, then we can apply the argument
  of the proposition. Taking limits finishes the proof.

\end{proof}

\end{example}

\subsection{Polynomials of equal degree}

For the rest of this section  we will always assume 
that
$$
deg(f) = deg(g) = deg(h) = deg(k) 
$$

$$\begin{pmatrix}{f}&{g}\\{h}&{k}
\end{pmatrix}
\begin{pmatrix}{r}\\{s}
\end{pmatrix}
= \begin{pmatrix}{u}\\{v}
\end{pmatrix}
$$

Lemma~\ref{lem:pattern-2} in the next section shows that in order to get
non-trivial results we must assume that $f,g,h,k,r,s$ all have
positive leading coefficients.  We consider two different ways that
interlacing can be preserved.

\begin{description}
\item[($\greateqeq$)] $\quad r\greateqeq s$ implies $u \greateqeq v$
\item[($\lesslesseq$)]   $\quad r\lesslesseq s$ implies $u \lesslesseq v$
\end{description}

These conditions put constraints on $f,g,h,k$.

\begin{lemma}
  Suppose $\smalltwo{f}{g}{h}{k}\smalltwobyone{r}{s} =
  \smalltwobyone{u}{v}$, and consider the two interlacings
\begin{align}
  \label{eqn:fag-1}
  f+\alpha g & \lesslesseq h+\alpha k  \text{ for $\alpha>0$} \\
  xf+\alpha g & \lesslesseq xh+\alpha k \text{ for $\alpha>0$} \label{eqn:fag-2}
\end{align}

\begin{tabular}{ccl}
 $r\greateqeq s$ implies $u \greateqeq v$ & $\implies$ &
 \eqref{eqn:fag-1} \\
$r\lesslesseq s$ implies $u \lesslesseq v$ & $\implies$ &
 \eqref{eqn:fag-2} 
\end{tabular}

\end{lemma}
\begin{proof}
  These cases are all similar. When  \textbf{($\greateqeq$)} holds then we take
  $r=x$ and $s=\alpha x$. When \textbf{($\lesslesseq$)} holds then we set $r=x^2$
  and $s=\alpha x$. Since $s$ must have a positive leading coefficient
  we require that $\alpha$ is positive. The results follow by
  computing $u,v$ and dividing the interlacing $u\longleftarrow v$ by
  $x$.
\end{proof}

If \eqref{eqn:fag-1} holds then we see from considering the left side
that $f$ and $g$ have a common interlacing.  If both \eqref{eqn:fag-1}
and \eqref{eqn:fag-2} then Lemma~\ref{lem:ispm-2} shows that we have true
interlacing.

\begin{remark}
  If we only require that $f\lesslesseq s$ implies that $u \lesslesseq
  v$ and do not require that the degrees of $f,g,h,k$ are equal then
  we do not get \eqref{eqn:fag-2}, and so have more possibilities. In
  particular note that if $r\lesslesseq s$ and
  $\smalltwo{x}{-b}{1}{0}\smalltwobyone{r}{s}=\smalltwobyone{u}{v}$
  then $u\lesslesseq v$. However, then entries of
  $\smalltwo{x}{-b}{1}{0}$ do not satisfy \eqref{eqn:fag-2}.
\end{remark}

When \eqref{eqn:fag-1} and \eqref{eqn:fag-2} hold we can combine the
restrictions into one polynomial.

\begin{lemma} \label{lem:fgzhk}
  If \eqref{eqn:fag-1} and \eqref{eqn:fag-2} hold then for any
$y,z\in\reals$ we have that
$$ g + zf + y k + yz h \in \allpoly$$
\end{lemma}
\begin{proof}
  For any positive $\alpha$ and any $y\in\reals$ we know that
  \begin{align*}
    g + \alpha f + y(k+\alpha h) & \in\allpoly \\
    g + \alpha xf + y(k+\alpha xh) & \in\allpoly \\
\intertext{Rewriting these yields}
(g+yk) + \alpha(f+yh) &\in\allpoly\\
(g+yk) + \alpha x(f+yh) &\in\allpoly\\
\intertext{Since $g+yk\in\allpoly$ and $f+yh\in\allpoly$ we apply
  Lemma~\ref{lem:ispm-2} to conclude that for any $z\in\reals$}
(g+yk) + z(f+yh) & \in\allpoly
   \end{align*}
\end{proof}
 
We will study such polynomials in \chap{extending}.

\section{Matrices of polynomials that preserve interlacing}
\label{sec:mat-poly-preserve}

\index{matrix!preserving interlacing}
We investigate the set $\tpp$ of all matrices $M$ with
polynomial entries such that if $v$ is a vector of mutually interlacing
polynomials then so is $Mv$.    $\tpp$ is a semigroup
since it's closed under multiplication, and $\tpp$ contains all
totally positive matrices.  In this section we give a few general
constructions.

We do not have a complete analog to Lemma~\ref{lem:faf}, but just a necessary condition.
Notice that the order of interlacing is reversed.

\begin{cor}  \label{cor:pattern}
  Suppose $f_1,\dots,f_n$ is a sequence of polynomials of the same degree
  and all positive leading coefficients.  If for all mutually
  interlacing polynomials $g_1\greateqeq \cdots\greateqeq g_n$ with positive leading
  coefficients the polynomial $f_1g_1 + \dots + f_ng_n$ has all real
  roots, then the $f_i$ are mutually interlacing, and 
  $$ f_n \greateqeq \cdots \greateqeq f_2 \greateqeq f_1$$
\end{cor}

\begin{proof}
  By taking limits we may assume that all but  $g_i$ and $g_j$, $i<j$, are zero.
  Now apply Lemma~\ref{lem:pattern2} to find that $g_j \greateqeq g_i$.
\end{proof}

We have seen an example of this reverse order in Lemma~\ref{lem:fg-ab}. Note
that we assumed that all the polynomials had positive leading
coefficients. This is not a convenience:

\begin{lemma} \label{lem:pattern-2}
  Suppose that $f\greateq g$ have positive leading coefficients. If
  for every $r\greateq s$ we have that $fs+gr\in\allpoly$, then $r$
  and $s$ have leading coefficients of the same sign.
\end{lemma}
\begin{proof}
    Choose $r=f^\prime$ and $s=-g^\prime$.  From
  Lemma~\ref{lem:inequality-1} we know that $fs+gr =
  \smalltwodet{f}{g}{f^\prime}{g^\prime}$ is never zero, and so has no
  real roots.  
\end{proof}

We do not have a satisfactory answer to the problem of when mutual
interlacing is preserved.  We will reconsider this question in
\chapsec{ideal}{line-transf} where we find relationships
with polynomials in three variables. Here are some necessary
conditions; they are easily seen not to be sufficient.

\begin{lemma}
  Suppose that $A = (p_{i,j})$ is a $m$ by $n$ matrix of polynomials
  of the same degree with the property that for all vectors $f =
  (f_1,\dots,f_n)$ of mutually interlacing polynomials, the polynomials
  $Af = (g_1,\dots,g_m)$ are mutually interlacing. Then
  \begin{itemize}
  \item The polynomials in any column of $A$ are mutually interlacing. 
    $$ p_{1,j} \greateqeq p_{2,j} \greateqeq \dots \greateqeq p_{n,j}$$
  \item The polynomials in any row of $A$ are mutually interlacing in
    the reverse order:
    $$ p_{j,1} \lesseqeq p_{j,2} \lesseqeq \dots \lesseqeq p_{j,n}$$
  \end{itemize}
\end{lemma}

\begin{proof}
  By continuity, we can apply $A$ to the vector $(f,0,\dots,0)$.  This
  implies that all columns are mutually interlacing.  The interlacing
  in the rows follows from Corollary~\ref{cor:pattern}.
\end{proof}

The next lemma determines some necessary conditions that depend on the
coefficients.

\index{totally positive$_2$}  
\begin{lemma}\label{lem:pattern-3}
  Suppose that $M$ is a matrix of polynomials of degree $n$ and
  positive leading coefficients that preserves mutually interlacing
  polynomials. The matrices of constant terms and of leading
  coefficients of $M$ are both totally positive$_2$.
  \end{lemma}
  \begin{proof}
    It suffices to prove this for two by two matrices. So, assume that
    $M=\smalltwo{f_1}{f_2}{f_3}{f_4}$ where $f_i$ has constant term
    $c_i$ and leading coefficient $d_i$. Let $M_r$ be the result of
    substituting $rx$ for $x$ in $M$. If $r$ is positive then $M_r$
    preserves interlacing. Consequently, $\lim_{r\rightarrow0} M_r =
    \smalltwo{c_1}{c_2}{c_3}{c_4}$ preserves mutually interlacing
    polynomials. Similarly, $\lim_{r\rightarrow\infty} r^{-n}M_r =
    \smalltwo{d_1}{d_2}{d_3}{d_4}$ preserves mutually interlacing
    polynomials. The conclusion now follows from
    Theorem~\ref{thm:totally-preserves}.
  \end{proof}

We can use the Leibnitz lemma (Lemma~\ref{lem:leibnitz}) to get a some classes of
elementary matrices that are analogous to the Jacobi matrices
\eqref{eqn:jacobi-matrix}.  Choose $g \greateqeq f$ with positive
leading coefficients. Assume that $r\greateqeq s$ and compute 
\[\smalltwo{f}{g}{0}{f}\smalltwobyone{r}{s}=\smalltwobyone{fr+gs}{fs}.\]
Lemma~\ref{lem:leibnitz} implies that $fr+gs\greateqeq fs$, so
$\smalltwo{f}{g}{0}{f}$ preserves interlacing. Similarly, the matrix
$\smalltwo{g}{0}{f}{g}$ also preserves interlacing.  We can generalize
this construction.

The first elementary matrix $E_i$ has $f$ on its diagonal, $g$ in
position $i,i+1$ for some $i$, and the remaining entries are $0$. The
second elementary matrix $E^0_i$ is not the transpose, but has $g$ on
its diagonal, $f$ in position $i+1,i$, and the remaining entries are
$0$. The last elementary matrix is the diagonal matrix $fI$. The
previous paragraph shows that these matrices preserve vectors of
mutually interlacing polynomials. For example with $n=4$ and $i=2$
these elementary matrices are

$$
E_2=\begin{pmatrix}
  f & 0 & 0 & 0 \\
  0 & f & g & 0 \\
  0 & 0 & f & 0 \\
  0 & 0 & 0 & f
\end{pmatrix}
\ \text{\ } \
E_2^0=\begin{pmatrix}
  g & 0 & 0 & 0 \\
  0 & g & 0 & 0 \\
  0 & f & g & 0 \\
  0 & 0 & 0 & g
\end{pmatrix}
\ \text{ } \
f\,I=\begin{pmatrix}
  f & 0 & 0 & 0 \\
  0 & f & 0 & 0 \\
  0 & 0 & f & 0 \\
  0 & 0 & 0 & f
\end{pmatrix}
$$

With these matrices we can construct examples of transformations that
preserve mutually interlacing sequences. For instance, if $n$ is $4$ then
$$ 
(E_1E_2E_3)\,(E_1E_2)\,E_1 =
\begin{pmatrix}
 f^6 & 3\,f^5\,g & 3\,f^4\,g^2 & f^3\,g^3 \\
 0 & f^6 & 2\,f^5\,g & f^4\,g^2 \\
0 & 0 & f^6 & f^5\,g \\
 0 & 0 & 0 & f^6 
\end{pmatrix}
$$
and consequently if $h_1\greateq h_2 \greateq h_3 \greateq h_4$ is
mutually interlacing then (after factoring out $f^3$) the following
are mutually interlacing.
\begin{align*}
   f^3 h_1 + 3\,f^2\,g h_2 + 3\,f\,g^2 h_3 +  g^3 h_4 &
   \greateqeq
  f^3 h_2 +  2\,f^2\,g h_3 +  f\,g^2 h_4 \\
&\greateqeq f^3 h_3+ f^2\,g h_4 \\
&\greateqeq f^3 h_4
\end{align*}

In general,  
\[ (E_1\cdots E_r)(E_1\cdots E_{r-1})\cdots(E_2E_1)E_1 =
\left(g^{j-i}f^{d+i-j}\binom{d+1-i}{j-i}\right)\] 
and consequently the
following sequence of polynomials is mutually interlacing
\[ \sum_{j=i}^d g^{j-i}f^{d+i-j}\binom{d+1-i}{j-i} h_i \quad\quad i=1,\dots,d\]

\begin{example}
  By taking coefficients of polynomials in $\gsubclose_3$ we can
  construct matrices of polynomials such that all rows and columns are
  interlacing, and all two by two determinants are positive. Here's an
  example of a $3\times3$ matrix of linear terms

\[
  \begin{pmatrix}
     86 x+1 & 625 x+76 & 1210 x+260 \\
 876 x+61 & 3953 x+666 & 4707 x+1351 \\
 1830 x+204 & 5208 x+1128 & 3216 x+1278
  \end{pmatrix}
\]

\noindent%
and one with quadratic terms

\[
\begin{pmatrix}
 644 x^2+114 x+1 & 9344 x^2+3019 x+81 & 30483 x^2+13579 x+466 \\
 7816 x^2+2102 x+78 & 52036 x^2+24710 x+1206 & 96288 x^2+73799 x+4818 \\
 19820 x^2+7918 x+417 & 72132 x^2+58752 x+4114 & 59610 x^2+110860 x+11777
\end{pmatrix}
\]

\end{example}

  \section{Matrices preserving interlacing in $\allpolypos$}
  \label{sec:matr-pres-interl}

  The class of matrices preserving interlacing in $\allpolypos$
  contains more than the totally positive matrices. For example, if
  $f\leftarrow g$ in $\allpolypos$, then $xg\leftarrow f$.  The
  corresponding matrix $\smalltwo{0}{x}{1}{0}$ preserves interlacing
  in $\allpolypos$, but not in $\allpoly$. In this section we will
  restrict ourselves to the following matrices:
  
\begin{quote}
  An NX matrix is a matrix whose entries are either non-negative
  constants, or positive multiples of $x$.
\index{NX matrix}
\end{quote}

We begin by determining which NX vectors map to $\allpolypos$.

\begin{lemma}\label{lem:nx-row}
   $v$ is an NX $1$ by $n$ matrix that maps all mutually interlacing
  sequences of polynomials in $\allpolypos$ to $\allpolypos$ if and
  only if all
  multiples of $x$ in $v$ occur to the right of all the positive constants.
\end{lemma}

\begin{proof}
If a multiple of $x$ is to the left of a positive constant, then we
may ignore all other coordinates, and just take $v=(x,1)$.  In this case
\[
f=x \qquad g=x+1 \qquad f\lessless g\qquad xf+g =
x^2+x+1\not\in\allpoly
\]

Conversely, assume that $v=(a_1,\dots,a_r,x\,a_{r+1},\dots,x\,a_n)$,
and let $F=(f_1,\dots,f_n)$ be mutually interlacing. Define
\[
g_1 = \sum_{k=1}^r a_k\,f_k\qquad g_2 = \sum_{k=r+1}^n  a_k\,f_k
\]
Note that $g_1\greateqeq g_2$ since the $f_i$ are mutually
interlacing. Consequently, $xg_2\lesslesseq g_1$ since all $f_i$ are
in $ \allpolypos$, and so $vF = g_1 + xg_2\in\allpolypos$.
\end{proof}

Next we look at two by two matrices in NX.

\begin{lemma}\label{lem:2x2-pos}
  The two by two matrices listed below  preserve 
  mutually interlacing polynomials in $\allpolypos$ if they satisfy 

  \begin{align*}
    1:\,\left(\begin{smallmatrix} \bullet & \bullet \\ \bullet & \bullet \end{smallmatrix}\right)
&&
  2:\,\left(\begin{smallmatrix} {\mathsf{ X}} & {\mathsf{ X}} \\ {\mathsf{ X}} & {\mathsf{ X}} \end{smallmatrix}\right)
&&
\text{determinant is $\ge0$}  \\
  3:\,\left(\begin{smallmatrix} {\mathsf{ X}} & {\mathsf{ X}} \\ \bullet & \bullet \end{smallmatrix}\right)
&&
    4:\,\left(\begin{smallmatrix} \bullet & {\mathsf{ X}} \\ \bullet & {\mathsf{ X}} \end{smallmatrix}\right)
&&
\text{determinant is $\le 0$} \\
    5:\,\left(\begin{smallmatrix} \bullet & {\mathsf{ X}} \\ \bullet & \bullet \end{smallmatrix}\right)
&&
  6:\;\left(\begin{smallmatrix} {\mathsf{ X}} & {\mathsf{ X}} \\ \bullet & {\mathsf{ X}} \end{smallmatrix}\right)
&&
\text{no restriction} 
  \end{align*}
\end{lemma}
\begin{proof}
  The $\bullet$ stands for a non-negative constant, and {\textsf{ X}} stands
  for a positive multiple of $x$. Since two by two matrices with
  positive determinant preserve interlacing, the first two matrices
  preserve interlacing.  The third and fourth cases follow from the
  factorizations 
  \begin{xalignat*}{2}
    \begin{pmatrix}       ax & bx \\ c & d     \end{pmatrix} &= 
    \begin{pmatrix}   0 & x \\ 1 & 0         \end{pmatrix}
    \begin{pmatrix}   c & d \\ a & b         \end{pmatrix} &
    \begin{pmatrix}   a & cx \\ b & dx          \end{pmatrix} &=
    \begin{pmatrix}   c & a \\ d & b           \end{pmatrix}
    \begin{pmatrix}  0 & x \\ 1 & 0          \end{pmatrix}
  \end{xalignat*}
\noindent
The last two cases  are easy and omitted.
\end{proof}

\begin{prop}
 \label{prop:NX} 
An NX matrix preserves mutually interlacing sequences of polynomials
in $\allpolypos$ if and only if
\begin{enumerate}
\item All two by two submatrices satisfy the
conditions of Lemma~\ref{lem:2x2-pos}.
\item All elements  that lie above or to the right of a
  multiple of $x$ are also multiples of $x$.
\end{enumerate}
\end{prop}
\begin{proof}
  It suffices to take two rows, and show that they determine
  interlacing polynomials. Choose mutually interlacing polynomials
  $(f_1,\dots,f_n)$. Let $a_r$ be the rightmost nonnegative constant in
  the first row, and $a_s$ the rightmost nonnegative constant in the
  second row. By the second hypothesis we know that $r\le s$. We
  denote the two rows by
  \begin{gather*}
    (a_1,a_2,\dots,a_r,x\,a_{r+1},\dots,x\,a_n) \qquad
    (b_1,b_2,\dots,b_s,x\,b_{s+1},\dots,x\,b_n) 
  \end{gather*}
  Define 
  \begin{align*}
    g_1 &= \sum_{i=1}^r a_i\,f_i & 
    g_2 &= \sum_{i=r+1}^s a_i\,f_i & 
    g_3 &= \sum_{i=s+1}^n a_i\,f_i \\
    h_1 &= \sum_{i=1}^r b_i\,f_i & 
    h_2 &= \sum_{i=r+1}^s b_i\,f_i & 
    h_3 &= \sum_{i=s+1}^n b_i\,f_i \\
  \end{align*}
We need to show that
\[
 g_1 + x\,g_2 + x\,g_3 \greateqeq h_1 + h_2 + x\,h_3
\]

The interlacings 
$
h_1\greateqeq g_2,\quad
g_1\greateqeq h_2,\quad
g_2\greateqeq g_3,\quad
g_2\greateqeq h_3\quad
$
follow from the mutual interlacing of the $f_i$'s. The interlacing
$h_2\greateqeq g_2$ follows from the hypothesis that 
$\left(\begin{smallmatrix} {\mathsf{ X}} & {\mathsf{ X}} \\ \bullet & \bullet
  \end{smallmatrix}\right)$ has negative determinant.
We now show that $g_1 + x\,g_3\greateqeq h_1 + x\,h_3$. Note that 
\[
F = (x\,f_{s+1},\dots,x\,f_n,f_1,\dots,f_r)
\]
is a mutually interlacing sequence. The
assumptions on the two by two matrices imply that
\[
M =
\begin{pmatrix}
  a_{s+1} & \dots & a_n & a_1 & \dots & a_r \\
  b_{s+1} & \dots & b_n & b_1 & \dots & b_r \\
\end{pmatrix}
\]
is TP$_2$. Since $MF = (g_1 + x\,g_3, h_1 + x\,h_3)$, we
find
\begin{align*}
  g_1 + x\,g_3& \greateqeq h_1 + x\,h_3 \\
x\,g_2 & \greateqeq h_1 + x\, h_3 \\
g_1 + x\,g_2 + x\,g_3 & \greateqeq h_1 + x\,h_3 \\
g_1 + x\,g_2 + x\,g_3 & \greateqeq h_2
\end{align*}
Adding these last two yields the conclusion.

Conversely, we know from Lemma~\ref{lem:nx-row} that all the $x$'s 
occur to the right of the constants. If (1) is satisfied, and (2) is
not, then $r< s$. It follows from the first part that the interlacing
is in the wrong direction. 

Suppose that a two by two matrix listed in Lemma~\ref{lem:2x2-pos} does
not satisfy the condition. In the first three cases we know the
interlacing is reversed. In the fourth case we take the matrix to be
$\smalltwo{1}{x}{1}{dx}$ and assume that $0<c<d$. Applying this to
$x\greateqeq x$ yields $x(x+1)\lesseq x(c+dx)$ since $0<c<d$.
\end{proof}

\begin{example}\label{ex:nx-1s}
  The matrix below preserves mutually interlacing sequences in
  $\allpolypos$ since all the determinants are zero, and the $x$'s
  satisfy the second condition.
 \[
 \begin{pmatrix}
    1 & x & x & \dots & x &x \\
    1 & 1 & x & \dots & x &x\\
    1 & 1 & 1 & \dots & x &x\\
    \vdots & & & \ddots & & \vdots\\
    1 & 1 & 1 & \dots & 1&x \\
    1 & 1 & 1 & \dots & 1 &1
  \end{pmatrix}
\]
\end{example}

\begin{example}\label{ex:hri}
    Suppose that $M$ is a totally positive matrix, and we form a new
    matrix $N$ by replacing the upper right corner element by $x$. All
    the hypothesis are met, so $N$ preserves mutual interlacing in
    $\allpolypos$. 

\end{example}

\begin{example}
  We can construct matrices with more complicated polynomials by
  multiplying these examples together. For instance, if
\[
m_1 =\begin{pmatrix}
  1 & x & x \\ 1& 1 & x \\ 1& 1 & 1
\end{pmatrix}
\qquad\qquad
m_2 =
\begin{pmatrix}
  1 & 0 & x \\ 1 & 1 & 0 \\ 1&2&1
\end{pmatrix}
\]
then the product $m_1\,m_2\,m_1$ below preserves mutually interlacing
polynomials in $\allpolypos$:

\[
\begin{pmatrix}
  1 + 7\,x & 6\,x + 2 \,x^2 & 3\,x + 5\,x^2 \\
 3 + 5\,x & 1 + 6\,x + x^2 & 5\,x + 3\,x^2 \\
 7 + x & 4 + 4\,x &1 + 7\,x  
\end{pmatrix}
\]

\end{example}

\begin{example}
  Here's an example where the matrix is not
  immediately obvious \cite{goodman-sun}*{Theorem 2.8}. Suppose $r$
    is a positive integer and that $f_0(x),\dots,f_{r-1}(x)$ are
    mutually interlacing. If 
\begin{multline*}
(1+x+\cdots+ x^{r-1})\,\bigl(f_0(x^r)+ x f_1(x^r) + \cdots+
x^{r-1}f_{r-1}(x^r)\bigr)
\\ = 
g_0(x^r)+ x g_1(x^r) + \cdots + x^{r-1}g_{r-1}(x^r)
\end{multline*}
then $g_0(x),\dots,g_{r-1}(x)$ are also mutually
interlacing. This follows from the observation that

\[
\begin{pmatrix}
  1 & x & x &  \hdotsfor{2} & x \\
  1 & 1 & x &  \hdotsfor{2} &x \\
  1 & 1 & 1 &  \hdotsfor{2} &x \\
  \vdots & & & \ddots&&\vdots\\
  1 & 1 & 1 & \hdots & 1 & x \\
  1 & 1 & 1 & \hdots & 1 & 1 \\
\end{pmatrix}
\begin{pmatrix}
  f_0(x) \\ f_1(x) \\ f_2(x) \\ \vdots \\ f_{r-2}(x) \\ f_{r-1}(x)
\end{pmatrix}
=
\begin{pmatrix}
  g_0(x) \\ g_1(x) \\ g_2(x) \\ \vdots \\ g_{r-2}(x) \\ g_{r-1}(x)
\end{pmatrix}
\]

\end{example}

\begin{example}\label{ex:hurwitz-x}
  Here's a similar example. Choose a positive integer $r$, and suppose
  that $g(x) = (x+\alpha)f(x)$. Write
  \begin{align*}
    f(x) &= \sum_{k=0}^{r-1} x^i f_i(x^r) &
    g(x) &= \sum_{k=0}^{r-1} x^i g_i(x^r) &
  \end{align*}
If $\{f_0,\dots,f_{r-1}\}$ are mutually interlacing then so are
$\{g_0,\dots,g_{r-1}\}$. The proof is the observation that
\[
\begin{pmatrix}
  \alpha & 0 & 0 & \hdots & x \\
  1 & \alpha & 0 & \hdots & 0 \\
0 & 1 & \alpha & \hdots & 0 \\
\vdots & & & \ddots & \vdots \\
0 & 0 & 0  & \hdots & \alpha
\end{pmatrix}
\begin{pmatrix}
  f_0(x) \\ f_1(x)  \\ \vdots  \\ f_{r-2}(x) \\ f_{r-1}(x)
\end{pmatrix}
=
\begin{pmatrix}
  g_0(x) \\ g_1(x)  \\ \vdots  \\ g_{r-2}(x)\\ g_{r-1}(x)
\end{pmatrix}
\]
This example can be used to give another proof of
Theorem~\ref{thm:hurwitz-gen}. 
\end{example}

  \section{Linear $2$ by $2$ matrices preserving interlacing }
  \label{sec:2-2-matrices}

  We are interested in determining when the matrix
\begin{equation}\label{eqn:2by2-1}
\begin{pmatrix}
  a_1x+b_1 &   a_2x+b_2 \\   a_3x+b_3 &   a_4x+b_4 
\end{pmatrix}
\end{equation}
preserves interlacing in $\allpolypos$.  We give a nice construction
due to Liu and Wang \cite{liu-wang}. There are constructions for
matrices of higher degree that preserve interlacing 
\seepage{lem:pres-int-det}.

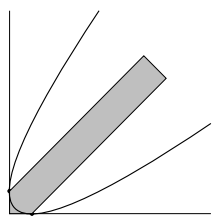
\begin{figure}[htbp] 
  \resizebox{3cm}{3cm}{%
    \begin{pspicture}(-1,-1)(9,9)
    \psline(0,9)(0,0)(9,0)
    \qdisk(1,0){.1cm}
    \qdisk(0,1){.1cm}
    \pspolygon[fillstyle=solid,fillcolor=lightgray](0,0)(1,0)(7,6)(6,7)(0,1)
    \psbezier[showpoints=false,linewidth=.5mm](4,9)(-2.7,-1)(-1,-2.7)(9,4)
  \end{pspicture}
}
  \caption{The graph of $(a-b)^2 - 2(a+b)+1 = 0$}
\label{fig:parabola}
\end{figure}

We begin by specializing to the matrix $M=\smalltwo{x+a}{x}{1}{x+b}$. If
we apply $M$ to the vector $(x,\alpha x)$ where $\alpha>0$
  then the result is a pair of interlacing polynomials. After factoring out $x$
\begin{gather*}
   x + a  + x \alpha  \greateqeq 1 + (b + x) \alpha  \\
 {-a}/(1 + \alpha )   \ge 
(-{1} -{b}\,\alpha )/{\alpha }\\
P=\alpha^2 b + \alpha( b+1-a) +1 
  \ge 0 \quad\text{for all } \alpha\ge0\\
\intertext{Next, apply $M$ to $(x^2,\alpha x)$. Factoring out $x$
  yields}
x(x+a+\alpha)  \lesslesseq x(1+\alpha)+b\alpha \\
\intertext{Substituting the root of the second term into the first
  gives}
\frac{-b\alpha}{1+\alpha}( \frac{-b\alpha}{1+\alpha} + a + \alpha ) 
\le 0 \\
R = \alpha^2 + (a+1-b)\alpha+a  \ge 0 \quad \text{for all $\alpha\ge0$}
\end{gather*}
Note that $R$ and $P$ have the same discriminant $Q = (a-b)^2 -
2(a+b)+1$. The graph of $Q=0$ is a parabola, and is given in Figure~\ref{fig:parabola}.

\begin{lemma}\label{lem:2x2-ab}
  $\smalltwo{x+a}{x}{1}{x+b}$ preserves interlacing in $\allpolypos$
  if and only $(a,b)$ is either in the parabola, or in the finite
  region bounded by the axes and parabola. See Figure~\ref{fig:parabola}.
\end{lemma}
\begin{proof}
  $Q$ intersects the axes at $(1,0)$ and $(0,1)$ and is negative in
  the region inside the parabola, so $P$ and $R$ are non-negative for
  those $(a,b)$. If $(a,b)$ is in the unbounded region below the
  parabola, then $Q$ is positive, and $b+1-a$ is negative. Thus, $P$
  has two positive roots, and so is not positive for positive
  $\alpha$. $R$ has two positive roots in the region above the
  parabola. Thus, if $M$ preserves interlacing then $(a,b)$ can't be
  in either unbounded region outside the parabola.
  
  Lemma~\ref{lem:2by2-2} shows that all points in the gray region of
  Figure~\ref{fig:parabola} determine matrices preserving interlacing.
  We can use Example~\ref{ex:ma-xyz} to show that all $(a,b)$ in the parabola
  determine interlacing preserving matrices.
\end{proof}

\begin{lemma}\label{lem:2by2-2}
If $a$ and $b$ are positive, and $|a-b|\le 1$ then
$\smalltwo{x+a}{x}{1}{x+b}$ preserves interlacing in $\allpolypos$.   
\end{lemma}
\begin{proof}
  We first claim that $\smalltwo{x+t}{x}{1}{x+t}$ preserves
  interlacing for positive $t$.   If we write
\[  \begin{pmatrix}
    x+t & x \\ 1 & x+t 
  \end{pmatrix}
  \begin{pmatrix}
    f \\ g
  \end{pmatrix} =
  \begin{pmatrix}
    (x+t)f+xg \\ f + (x+t)g
  \end{pmatrix}
\]
then each of the two terms on the bottom interlaces each of the two
terms of the top. Similarly, $\alpha(f+x g) \greateqeq f + (x+t)g$, so
$\smalltwo{\alpha}{\alpha x}{1}{x+t}$ preserves interlacing for any
positive $\alpha$.

Since $\smalltwo{\alpha}{\alpha x}{1}{x+t}$ and
$\smalltwo{x+t}{x}{1}{x+t}$ have the same second row, the matrix
$\smalltwo{x+t + \alpha}{(\alpha+1)x}{1}{x+t}$
preserves interlacing. Multiplying the first row and the second column
by $1/(\alpha+1)$, and replacing $x$ by $(\alpha+1)x$ shows that 
\begin{equation}\label{eqn:parabola-1}
\begin{pmatrix}
  x+ \frac{t+\alpha}{\alpha+1} & x \\ 1 &  x + \frac{t}{\alpha+1}
\end{pmatrix}
\end{equation}
preserves interlacing.
Suppose that $b+1\ge a \ge b$. If we choose $t = b/(b+1-a)$ and
$\alpha = (a-b)/(b+1-a)$ then $t$ and $\alpha$ are both non-negative,
and \eqref{eqn:parabola-1} reduces to
$\smalltwo{x+a}{x}{1}{x+b}$. The case $a+1\ge b \ge a$ is similar. 
\end{proof}

We have seen that the criterion for a $2\times2$ matrix of polynomials
to preserve interlacing is that the determinant is always positive. If
we restrict ourselves to $\allpolypos$ then we only need that the
determinant is positive for negative $x$.
The following result is due to Liu and Wang\cite{liu-wang}.

  \begin{lemma}\label{lem:preserve-2}
    Let $M = \smalltwo{\ell_1}{q}{c}{\ell_2}$ where $c>0$,
    $\ell_1,\ell_2\in\allpolypos(1)$ and $q\in\allpolypos(2)$, where
    $\ell_1,\ell_2,q$ have positive leading coefficients. If $|M|\ge0$
    for all $x\le0$ then $M$ preserves interlacing in $\allpolypos$.
  \end{lemma}
  \begin{proof}
    Let $M(f,g)^t = (F,G)^t$ where $f\lessless g$ in $\allpolypos$.
    The root $\alpha$ of $\ell_1$ is negative, and since the
    determinant is non-negative at $\alpha$, it follows that
    $q(\alpha)\le0$ and so $q\lesslesseq \ell_1$. From Leibnitz
    \mypage{lem:leibnitz} 
     we know that $F$ and $G$ are in $\allpolypos$.

    Since $G$ has positive coefficients and negative roots, $|M|$ is
    non-negative at the roots of $G$. In addition, we have
    $G=cf+\ell_2g\lesslesseq g$.  The identity
\[ cF = \ell_1 \,G - |M|\,g
\]
shows that $F$ sign interlaces $G$. Since the degree of $F$ is greater
than the degree of $G$ we conclude that $F\lesslesseq G$. 
  \end{proof}

If $q$ is linear then we get
\begin{cor}
  Let $M = \smalltwo{\ell_1}{\ell_3}{c}{\ell_2}$ where $c>0$ and
  $\ell_1,\ell_2,\ell_3\in\allpolypos(1)$ have positive leading
  coefficients. If $|M|\ge0$ for all $x\le0$ then $M$ preserves
  interlacing in $\allpolypos$.
  \end{cor}

Note that the corollary provides another proof of
Lemma~\ref{lem:2x2-ab}, since $\smalltwodet{x+a}{x}{1}{x+b}$ is
$(x+a)(x+b)-x$, and this determines the same region that is described
in the proposition.

Here is a simple consequence.

\begin{cor}
  $\smalltwo{x}{x}{1}{x+\gamma}$ and $\smalltwo{x+\gamma}{x}{1}{x}$
  preserve interlacing in $\allpolypos$ if and only if $0\le \gamma \le 1$.
\end{cor}
\begin{proof}
  Existence follows from the lemma. If $\gamma>1$ then $(0,\gamma)$
  and $(\gamma,0)$ lie outside the parabola.
\end{proof}

\begin{lemma}
  Suppose that $a,b,c,d$ are positive, and
  $\smalltwodet{a}{b}{c}{d}\ge0$.
  \begin{enumerate}
  \item $\smalltwo{a}{x+b}{c}{d}$ preserves interlacing in $\allpolypos$.\\[.1cm]
  \item $\smalltwo{a}{x+b}{c}{x+d}$ preserves interlacing
    $\allpolypos$ if and only if $c\ge a$, $d\ge b$.\\[.1cm]
\item $\smalltwo{x+a}{x+b}{c}{d}$ preserves interlacing $\allpolypos$ if and only
  if $a\ge b$, $c\ge d$.\\[.1cm]
\item $\smalltwo{x+a}{x+b}{x+c}{x+d}$ preserves interlacing $\allpolypos$ if and
  only if \\ $c\ge a\ge b$, and $(c-a)(a-b)\ge ad-bc\ge0$.
  \end{enumerate}
\end{lemma}
\begin{proof}
  The first one follows from the decomposition $\smalltwo{0}{x}{c}{d}
  + \smalltwo{a}{b}{c}{d}$.  If  $c\ge a$ and $d\ge b$ then 
\[
\begin{pmatrix}
  a&x+b\\c&x+d
\end{pmatrix}
=
\begin{pmatrix}
  1&0\\1&1
\end{pmatrix}
\begin{pmatrix}
  a & x+b \\ c-a & d-b
\end{pmatrix}
\]
shows that $\smalltwo{a}{x+b}{c}{x+d}$ preserves
interlacing. Conversely, if interlacing is preserved, then applying
$\smalltwo{a}{x+b}{c}{x+d}$ to $(1,\alpha)$ yields 
\[
a + \alpha( x +  b) \greateqeq c+ \alpha(x+d)
\]
Taking roots shows that $\alpha b + a \le \alpha d+c$ for all
positive $\alpha$  which implies $c\ge a$ and $d\ge b$.

The third one is similar: sufficiency follows from
\[
\begin{pmatrix}
  x+a & x+b \\ c&d
\end{pmatrix}
= 
\begin{pmatrix}
  a-b & x+b \\ c-d & d
\end{pmatrix}
\begin{pmatrix}
  1&1\\0&1
\end{pmatrix}
\]
and necessity is a consequence of applying $\smalltwo{x+a}{x+b}{c}{d}$
to $(x,\alpha)$. 
For the last one, we first compute
\[
\begin{pmatrix}  1 & 0 \\ 1 & 1 \end{pmatrix}
\begin{pmatrix}  x+a & x \\ c-a & \frac{ad-bc}{a-b} \end{pmatrix}
\begin{pmatrix}  1 & 0 \\ 0 & a-b \end{pmatrix}
\begin{pmatrix}  a & b \\ 0 & 1 \end{pmatrix}
=
\begin{pmatrix}  x+a & x+b \\ x+c & x+d \end{pmatrix}
\]
By Corollary~\ref{cor:inv-2by2} we see that
$\smalltwo{x+a}{x+b}{x+c}{x+d}$ preserves interlacing if and only if  
$\smalltwo{x+a} {x }{ c-a }{ \frac{ad-bc}{a-b}}$ does. From (3)
this is the case if and only if ${ c-a \ge \frac{ad-bc}{a-b}}$.
\end{proof}

  \section{Polynomial sequences from matrices}
  \label{sec:high-order-recur}

\index{totally positive!recurrence relations}
\index{NX matrix!recurrence relations}
\index{recurrence relations!from matrices}

A totally positive or NX $n\times n$ matrix determines sequences of
polynomials obeying an $n$'th order recurrence relation. We first
recall how to get sequences satisfying a recurrence relation from a
matrix. Suppose that $M$ is a matrix, $v$ is a vector, and the
characteristic polynomial of $M$ is
\[
a_0 + a_1 y + a_2 y^2 + \cdots + a_n y^n.
\]
Since $M$ satisfies its characteristic polynomial, we have
\[
a_0 (v) + a_1 (Mv) + a_2 (M^2v) + \cdots + a_n (M^nv) = 0.
\]
If we fix an integer $k$ and define $w_i$ to be the $k$th coordinate
of $M^iv$ then the $w_i$ also satisfy the same recurrence relation
\begin{equation}\label{eqn:rec-tp-1}
a_0 \,w_0 + a_1\,w_1 + a_2\,w_2 + \cdots + a_n\,w_n = 0.
\end{equation}

We now let $M$ be a totally positive matrix, and let $v$ be a matrix
of mutually interlacing polynomials. Since $M^iv=v_i$ is a vector of
mutually interlacing polynomials, all the coordinates of $v_i$ are are
in $\allpoly$.  These polynomials satisfy the recurrence
relation~\eqref{eqn:rec-tp-1}.  All these polynomials have the same
degree but  consecutive $v_i$ do not necessarily interlace.

\begin{example}
  Consider an example where $M= \left(\begin{smallmatrix}1&2&1\\ 1&3&3
      \\ 1&4 & 6\end{smallmatrix}\right)$ and we choose $v_0$ to be
  three mutually interlacing polynomials of degree $2$.

\[
\begin{array}{rccccc}
v_0 &= \bigl( &
( 1 + x )( 4 + x ) & ( 2 + x )( 5 + x ) & ( 3 + x )( 6 + x ) &
\bigr)\\[.1cm]
v_1 &= \bigl( &
2( 21 + 14x + 2x^2 )  &
88 + 53x + 7x^2 &
152 + 87x + 11x^2 &
\bigr)\\[.1cm]
v_2 &= \bigl( &
 370 + 221x + 29x^2 &
2( 381 + 224x + 29x^2 ) &
   2( 653 + 381x + 49x^2 )  &
\bigr)
 \end{array}
\]

If we just consider the second column, and observe that the
characteristic polynomial is $1 - 12 y + 10 y^2 - y^3$, then we
conclude that all polynomials defined by the sequence below are in
$\allpolypos$.

\begin{align*}
    p_{n+3} &= 10 p_{n+2} - 12 p_{n+1} + p_n \\
    p_0 &= ( 2 + x )( 5 + x )\\
    p_1 &= 88 + 53x + 7x^2 \\
    p_2 & = 2( 381 + 224x + 29x^2 ) 
\end{align*}
\end{example}

We now consider NX matrices. These are more interesting, since the
degrees of the polynomials can increase. 

\begin{example}
  We  take the matrix $M = 
  \left(\begin{smallmatrix}
    1 & x  \\ 1  & 1 
  \end{smallmatrix}\right)$. If we begin with $v_0=(x+1,x+2)$ then
the second column yields that  $p_n\in\allpolypos$, where

\begin{align*}
    p_{n+2} &=  2 p_{n+1} + (x-1)p_n \\
    p_0 &= x+2\\
    p_1 &= 2x+3
\end{align*}
  
\end{example}

\begin{example}
  We now take the matrix $M = 
  \left(\begin{smallmatrix}
    1 & x & x \\ 1 & 1 & x \\ 1 & 1 & 1
  \end{smallmatrix}\right)$. If we begin with $v_0=(x+1,x+2,x+3)$ then
the third column yields that all $p_n\in\allpolypos$, where

\begin{align*}
    p_{n+3} &= 3p_{n+2} + 3(x-1)p_{n+1} +(x-1)^2p_n\\
    p_0 &= x+3\\
    p_1 &= 3x+6\\
    p_2 & = 3x^2+14x+10
\end{align*}

\end{example}

\begin{example}
  In this example we let $M$ be the cube of the previous $M$. 
\[
M =
\begin{pmatrix}
   1 + 7\,x + x^2 & 6\,x + 3\,x^2 & 3\,x + 6\,x^2 \\
 3 + 6\,x & 1 + 7\,x + x^2 & 6\,x + 3\,x^2 \\
 6 + 3\,x & 3 +    6\,x & 1 + 7\,x + x^2 
\end{pmatrix}
\]
Since every row has a quadratic term, the degree of the $p_n$'s goes
up by two at each step, so consecutive $p_n$'s don't interlace. If we
again start with $v_0=(x+1,x+2,x+3)$ then all $p_n$ are in
$\allpolypos$:

\begin{align*}
    p_{n+3} &= 3(1+7x+x^2)p_{n+2} - 3 (x-1)^4p_{n+1} + (x-1)^6p_n\\
    p_0 &= x+3\\
    p_1 &= 5 + 42 x + 30 x^2 + 4 x^3\\
    p_2 & =8 + 260 x + 947 x^2 + 784 x^3 + 181 x^4 + 7 x^5
\end{align*}

\end{example}

\section{Interlacing via the complexes}
\label{sec:interl-via-compl}

In this section we show how we can use properties of complex numbers
to establish interlacing. Let $\uhp$ be the open upper half plane
$\{z\mid \Im(z)>0\}$. Note that
\[
f\in\allpoly \Leftrightarrow f(\alpha)\ne0 \ \text{for all
}\alpha\in\uhp.
\]
Our first result connects interlacing and quotients.

\begin{lemma}
  Suppose that $f,g$ are relatively prime and have positive leading
  coefficients. The following are equivalent
  \begin{enumerate}
  \item $f\longleftarrow g$
  \item $\frac{f}{g}\colon\uhp\longrightarrow\uhp$
  \end{enumerate}
\end{lemma}
\begin{proof}
  Suppose 1) holds. If we write
\[ g = a_0 f + \sum a_i \frac{f}{x-r_i} \]
where $\roots{(f)} = (r_i)$ then
\[
\frac{g}{f}(\alpha) = a_0 + \sum a_i \frac{1}{\alpha - r_i}
\]
If $\alpha\in\uhp$ then $\alpha-r_i\in\uhp$, so
$1/(\alpha-r_i)\in-\uhp$, and therefore
$(g/f)(\alpha)\in-\uhp$. Taking inverses shows that
$(f/g)(\alpha)\in\uhp$.

Conversely, if 2) holds and $f+ t g\not\in\allpoly$ for some
$t\in\reals$ then there is $\alpha\in\uhp$ such that $f(\alpha) +
tg(\alpha)=0$. Thus, $(f/g)(\alpha)\not\in\uhp$. Thus $f$ and $g$
interlace, and the previous paragraph shows that the direction must be
$f\longleftarrow g$. 
\end{proof}

The upper half plane is a cone \index{cone} - it is closed under
positive sums.  This additivity explains two fundamental results about
interlacing.

\begin{lemma}\label{lem:interlace-geom} Suppose that $f,g,h,k$ in $\allpoly$ have positive
  leading coefficients.
  \begin{enumerate}
  \item If $f\lesslesseq g,h$ then $f\lesslesseq g+h$.
  \item If $f\lesslesseq g$, $h\lesslesseq k$ then $ fh \lesslesseq fk
+ gh \lesslesseq gk$.
  \end{enumerate}
\end{lemma}
\begin{proof} If $f\lesslesseq g,h$ and $\sigma$ is in $\uhp$, then
  $\frac{g+h}{f}(\sigma)\in-\uhp$ since $\frac{g}{f}(\sigma)$
  and $\frac{h}{f}(\sigma)$ are, and $-\uhp$ is a cone. It follows
  that $f/(g+h)(\sigma)\in\uhp$.

  The second one is similar; it's just the addition of fractions:
\[ \frac{g}{f}(\sigma) + \frac{k}{h}(\sigma) =
\frac{gh+fk}{fh}(\sigma). \]
\end{proof}

Note that the region $\complexes\setminus(-\infty,0)$ is \emph{not} a
cone, so we should not expect $f\plessless g,h$ to imply
$g+h\in\allpoly$. Indeed, it's false.

  All points in the upper half plane can be realized by ratios of
interlacing polynomials.

  \begin{lemma}\label{lem:find-ratio} If $\alpha$ and $\sigma$ are in
$\uhp$ then there are $f\lessless g$ in $\allpoly$ with positive
leading coefficients such that $\frac{f(\sigma)}{g(\sigma)}=\alpha$.
  \end{lemma}
  \begin{proof}

We will show that the set
\[ S = \left\{
  \begin{pmatrix} f(\sigma) \\g(\sigma)
  \end{pmatrix}\,\mid\, f\lessless g\right\}
\] equals $\uhp$, considered as a subset of complex projective
space. We know all elements of $S$ lie in $\uhp$. Choose any
$f\lessless g$ and let $f(\sigma)/g(\sigma) = \beta$. Since
$SL_2(\reals)$ is transitive on $\uhp$ we can find
$M=\smalltwo{a}{b}{c}{d}$ in $SL_2(\reals)$ such that
$M(\beta)=\alpha$. Then
\[
\begin{pmatrix} af+bg\\cf+dg
\end{pmatrix}(\sigma) = M\,
\begin{pmatrix} f\\g
\end{pmatrix}(\sigma) = M(\beta)=\alpha
\] Since $M$ has determinant one we know that $af+bg\lessless cf +
dg$.
  \end{proof}

Matrices preserving interlacing are the same as matrices mapping
$\uhp$ to $\uhp$. The next lemma provides an alternative proof of
Corollary~\ref{cor:lin-comb-new}, and clarifies why there is a simple
determinant condition. Recall that if $M$ has all real coefficients
then $M$ maps $\uhp$ to itself if and only if $|M|>0$.

  \begin{lemma} Suppose that $M$ is a $2$ by $2$ matrix of polynomials
with positive leading coefficients. The following are equivalent
    \begin{enumerate}
    \item $M$ preserves interlacing.
    \item $M(\sigma)$ maps $\uhp$ to $\uhp$ for all
$\sigma\in\uhp$.
    \end{enumerate}
  \end{lemma}
  \begin{proof} Assume that $M$ preserves interlacing, and choose any
$\sigma$ and $\alpha$ in $\uhp$. By Lemma~\ref{lem:find-ratio} we
can find $f\lessless g$ such that $f(\sigma)/g(\sigma)=\alpha$. Then
\[ \bigl[M(\sigma)\bigr](\alpha) =
\bigl[M(\sigma)\bigr]\smalltwobyone{f(\sigma)}{g(\sigma)} =
\bigl[M\smalltwobyone{f}{g}\bigr](\sigma)\in\uhp
\] since $M\smalltwobyone{f}{g}$ is a pair of interlacing polynomials.

Conversely, if $f\lessless g$ and
$M\smalltwobyone{f}{g}=\smalltwobyone{F}{G}$ then for all
$\sigma\in\uhp$ we know
$\smalltwobyone{f(\sigma)}{g(\sigma)}\in\uhp$. Thus
$\smalltwobyone{F(\sigma)}{G(\sigma)}\in\uhp$, which implies that
$F\lesslesseq G$.
  \end{proof}

\section{Mutually interlacing polynomials}
\label{sec:mutu-interl-polyn}

We will generalize the following three equivalent statements that hold
for a real $2$ by $2$ matrix $M$ to mutually interlacing polynomials and
strictly totally positive matrices$_2$.

\begin{enumerate}
\item $M$ has positive determinant.
\item $M$ maps the upper half plane to itself.
\item $M$ preserves interlacing polynomials.
\end{enumerate}

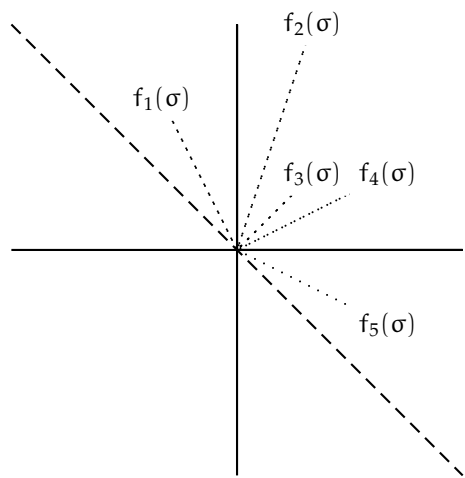
\begin{figure} \centering

\newpsobject{showgrid}{psgrid}{subgriddiv=1,griddots=10,gridlabels=6pt}
 \begin{pspicture}(-3,-3)(3,3)
\psline(-3,0)(3,0)(0,0)(0,3)(0,-3)
\psline[linestyle=dashed](-3,3)(3,-3)
\psline[linestyle=dotted](0,0)(-1,2)(0,0)(1,3)(0,0)(1,1)
(0,0)(2,1)(0,0)(2,-1) \rput*(-1,2){$f_1(\sigma)$}
\rput*(1,3){$f_2(\sigma)$} \rput*(1,1){$f_3(\sigma)$}
\rput*(2,1){$f_4(\sigma)$} \rput*(2,-1){$f_5(\sigma)$}
  \end{pspicture}
  \caption{Mutually interlacing polynomials evaluated at $\sigma$}
  \label{fig:mutual}
\end{figure}

First of all, mutually interlacing polynomials have a pretty geometric
interpretation. Suppose that $f_1,\dots,f_n$ are mutually interlacing.
For any $\sigma\in\uhp$ we know that
$f_i(\sigma)/f_j(\sigma)\in\uhp$ for $i<j$. In terms of angles this
means that the angle from $f_j(\sigma)$ to $f_i(\sigma)$ is less than
$\pi$. It follows that
\begin{quote} All the points $\{f_i(\sigma)\}$ lie in a half plane,
and are , in counterclockwise order, $f_1(\sigma),\dots,f_n(\sigma)$.
(Figure~\ref{fig:mutual})
\end{quote}

A $2$ by $2$ matrix with real entries preserves the upper half plane
iff it preserves interlacing. We know that an $n$ by $n$ matrix
preserves interlacing if and only if it is totally positive$_2$. We
now introduce the structure preserved by strictly totally positive$_2$
matrices.

\[
H^+_d = \{(u_1,\dots,u_d) \,\mid\, \frac{u_i}{u_j}\in\uhp\ \text{for }
1\le i < j \le d \}.  \] 

For instance,
\[ H_3^+ = \{(u,v,w)\,\mid\,
\frac{u}{v}\in\uhp, \frac{u}{w}\in\uhp,
\frac{v}{w}\in\uhp\} 
\]

The space $H_2^+$ is a subset of $P^2(\complexes)$, and is usually identified
with $\uhp$.  A particular point in $H_d^+$ is
$u=(\zeta^{d-1},\dots,\zeta,1)$ where $\zeta = e^{\pi\,\imag/d}$.  The
characterization of mutual interlacing is the same as for interlacing,
and follows easily from the case $d=2$.

\begin{lemma}
Assume $(f_1,\dots,f_d)$ is a sequence of polynomials. The following
are equivalent:
\begin{enumerate}
\item $(f_1,\dots,f_d)$ is strictly mutually interlacing.
\item $(f_1(\sigma),\dots,f_d(\sigma))\in H^+_d$  for all
  $\sigma\in\uhp$. 
\end{enumerate}
\end{lemma}

Let $\ttt_d$ be the semigroup of all invertible matrices with non-negative
entries for which the determinants of  all two by two submatrices are
positive. For example, $\ttt_2$ consists of all two by two
matrices with non-negative entries and positive determinant. 

\begin{lemma}
  $\ttt_d$ acts transitively on $H_d^+$. 
\end{lemma}
\begin{proof}
  We know that $\ttt_2$ acts transitively on $\uhp$. 
Choose two elements $(u_1,\dots,u_d)$ and $(v_1,\dots,v_d)$ of
$H_d^+$. The block matrix 
\[
\begin{pmatrix}
  I & 0 & 0 \\ 0 & M & 0 \\ 0 & 0 & I
\end{pmatrix}
\]
where $M$ is a $2$ by $2$ matrix in $\ttt_2$ preserves $H^+_d$. 
Since $\ttt_2$ is transitive we first carry $u_1,u_2$ to $v_1,v_2$,
then leaving the first coordinate fixed we carry $v_2,u_3$ to
$v_2,v_3$ and so on.
\end{proof}

\begin{lemma}
  For any $u\in H^+_d$ and $\sigma\in\uhp$ there is a vector $f$ of
  strictly mutually interlacing polynomials such that $f(\sigma)=u$.
\end{lemma}
\begin{proof}
  Choose any vector $g$ of mutually interlacing polynomials. Since
  $g(\sigma)\in H^+_d$ there is an $M\in\ttt_d$ such that
  $M(g(\sigma))=u$. Thus
$u = (Mg)(\sigma)$, so $Mg$ is the desired vector. 
\end{proof}

  \begin{prop}\label{prop:miHplus}
 Suppose that $M$ is a $d$ by $d$ matrix of positive real
numbers. The following are equivalent
\begin{enumerate}
\item $M$ preserves  mutually interlacing polynomials.
\item $M$ is strictly totally positive$_2$.
\item $M$ maps $H^+_d$ to itself.
\end{enumerate}
\end{prop}

\begin{proof}
  We have already seen the first two equivalences.
  Assume $M$ maps $H_d^+$ to itself. For any $\sigma\in\uhp$ and
  mutually interlacing $f_1\lessless f_2\cdots\lessless f_d$ we know
  that
\[
\bigl[M\cdot(f_1,\dots,f_d)\bigr](\sigma)=
M\cdot (f_1(\sigma),\dots,f_d(\sigma))\in H_d^+
\]
since $(f_1(\sigma),\dots,f_d(\sigma))\in H_d^+$.
Thus, $M$ preserves interlacing, and so is totally positive$_2$.

Conversely, suppose that $M$ is totally positive$_2$, and $u =
(u_1,\dots,u_d)$ is in $H_d^+$. Choose mutually interlacing
$f=(f_1,\dots,f_d)$ such that $f(\sigma)=u$ for some
$\sigma\in\uhp$. Then 
\[
Mu = M(f(\sigma)) = (Mf)(\sigma) \in H^+_d
\]
since $Mf$ consists of strictly totally interlacing polynomials.
\end{proof}

\begin{remark}\label{rem:tp2}
  We give a direct proof that (2) implies (3).  It suffices to take
  $M$ to be the $2$ by $d$ matrix $\smalltwobyone{a}{b}$.  Assume
  $a=(a_i)$, $b = (b_i)$ and $v = (v_i)$. We show that $Mv\in\uhp$
  by computing the imaginary part of the quotient of $Mv$.

  \begin{gather*}
    \Im\left(\frac{\sum a_i v_i}{\sum b_iv_i}\right)
    = \frac{\Im\,\bigl(\sum a_i v_i\bigr)\bigl({\sum
        b_i \overline{v_i}}\bigr)}
    {\bigl|{\sum b_i v_i}\bigr|^2} \\
   \intertext{Ignoring the positive denominator, the imaginary part is}
 \Im\biggl[ \sum a_ib_i|v_i|^2 + \sum_{i<j} a_ib_j v_i\overline{v_j} +
a_j b_i v_j \overline{v_i}\biggr]
=
\sum_{i<j} \begin{vmatrix}a_i&a_j\\b_j & b_j\end{vmatrix}
\Im(v_i\overline{v_j})
  \end{gather*}
All the two by two determinants are positive by hypothesis, and $\Im
(v_i\overline{v_j})$ is positive since $v_i/v_j\in\uhp$.
\end{remark}


\chapter{Homogeneous Polynomials}
\label{cha:homog}

\renewcommand{\TimeStampStart}{Tuesday, March 11, 2008: 09:47:15}
\mytoday  
We begin the study of polynomials in two variables with a very special
case: homogeneous polynomials. These results are generally well known,
and will be generalized in later chapters.

\section{Introduction} 
A homogeneous polynomial of degree $n$ is a polynomial in $x,y$ such
that the degree of every term is $n$.   A polynomial $f$ of
degree $n$ in one variable determines  a homogeneous polynomial
of degree $n$:
\begin{align}
  F(x,y) &= y^n f(x/y).\notag \\ 
\intertext{The coefficients of $F$ are those of $f$. If} f(x)
  &= a_0 + a_1x + \dots + a_n x^n \label{eqn:homog:1} \\ 
\index{homogenized  polynomial}
\intertext{then the homogenized  polynomial is} 
F(x,y) &= a_0 y^n + a_1xy^{n-1} + \dots + a_n x^n \label{eqn:homog:2}
\end{align}
We can recover $f$ from $F$ since $f(x) = F(x,1)$.

The key relation for us between homogeneous polynomials and
polynomials in one variable is that factorization for homogeneous
polynomials corresponds to the property of having all real roots.

\begin{lemma}
  Let $F(x,y)$ be homogeneous, and set $f=F(x,1)$.  We can  factor $F$ into
  linear factors if and only if $f\in\allpoly$.
\end{lemma}
\begin{proof}
  If $f=\prod(a_ix+b_i)$ then $F=\prod(a_ix+b_iy)$ and conversely. 
\end{proof}

\index{reverse!of a polynomial}
The reversal of $f$ , written $\rev{f}$,
is
\begin{equation}
  \label{eqn:homog-3}
  a_n + a_{n-1}x + \cdots + a_0x^n.
\end{equation}
We can express $\rev{f}$ in terms of $F$ since $\rev{f}(x) =
F(1,x)$. This shows that if $f$ is in $\allpoly$ then $\rev{f}$ is
also in $\allpoly$.

\section{Polar derivatives} \label{sec:polar}
Since a polynomial in $x$ and $y$ has derivatives with respect to both
$x$ and $y$ we can construct a new derivative operator for polynomials
in one variable. From \eqref{eqn:homog:2} we find

\begin{align}
 \partialy{} F &= na_0y^{n-1} + (n-1)a_1y^{n-2}x +
    \dots + a_nx^{n-1}.\notag \\ 
\intertext{If we substitute $y=1$, then the right hand side becomes} 
& na_0 + (n-1)x + \dots + a_{n-1}x^{n-1}\label{eqn:polar-eqn} \\ 
\intertext{which can  be written as} 
&n(a_0 + a_1x + \dots + a_nx^n) - x(a_1 +2a_2x +
    \dots + nx_{n-1})\notag \\ &= nf - xf^\prime \notag
\end{align}
  
\index{polar derivative} This is called the \emph{polar derivative} of
$f$ and is written $\partialy f$.  \index{polar derivative} The degree
of $\partialy f$ is $n-1$, and the formula for $\partialy f$ has an
$n$ in it.  The polar derivative is \emph{not} not defined on
$\allpoly$, but only on $\allpoly(n)$.
It is useful to see that the polar derivative can be defined as
follows:
$$ \text{reverse $f$; }\quad\quad\text{differentiate; }
\quad\quad\text{reverse}
$$
Since these operations preserve $\allpoly$, the polar derivative is in
$\allpoly$. 

\index{interpolation!basis} If we assume $f = \prod(x-a_i)$ then 
\begin{align}\label{eqn:polar-series}
 \partialy f &=  \sum (-a_i)\, \frac{f(x)}{x-a_i}
 \\
\intertext{which is similar to the representation for the derivative}
 \partialx f &=  \sum  \frac{f(x)}{x-a_i}\notag
\end{align}
If $f\in\allpolypm$ then the polar derivative has all real
roots\footnote{It fails to have all real roots if there are positive
  and negative roots - e.g. ${(x+1)^2(x-2)^2}$} by Lemma~\ref{lem:inequality-6a}.
Since polar derivatives act like derivatives, it is not surprising
that we have this result (due to Laguerre). In fact, in
Chapter~\ref{cha:p2} we will see that they are exactly derivatives.

\begin{lemma} \label{lem:polar-1}
  If $f\in\allpolypm$ and $f$ has all distinct roots,
  then $f\lessless \partialy f$. The polar derivative maps
  $\allpolypos$ to itself, and $\allpolyalt$ to itself.
\end{lemma}
\begin{proof}
The result follows from Lemma~\ref{lem:sign-quant} and
\eqref{eqn:polar-series}. The mapping properties follow from \eqref{eqn:polar-eqn}.
\end{proof}

\begin{cor} \label{cor:fagp2}
  If $f\in\allpolypm$ then $f-a\, \partialx \partialy\, f\in\allpoly$
  for $a>0$.
\end{cor}
\begin{proof}
  Use Lemma~\ref{lem:polar-1} and Corollary~\ref{cor:fagp}.
\end{proof}

\begin{lemma} \label{lem:polar-3}
  If $f = a_0 + \dots + a_nx^n$ is in $\allpoly$ then $\partialx f$
  and $\partialy f$ interlace. If $f$ has positive leading coefficient 
  then 
  \begin{itemize}
  \item If $a_{n-1}$ is negative then $\frac{\partial}{\partial
      x}f \greateqeq \frac{\partial}{\partial y}f$.
  \item If $a_{n-1}$ is positive then $\frac{\partial}{\partial
      x}f \lesseqeq \frac{\partial}{\partial y}f$.
  \end{itemize}
\end{lemma}
\begin{proof}
  After dividing  by $a_n$   we may assume that $f$ is monic.
  Since $\frac{\partial}{\partial y}f = -xf^\prime+nf$ we see that
  $\frac{\partial}{\partial y}f$ 
  sign interlaces $f^\prime$.  The sign of $\frac{\partial}{\partial
  y}f$ at the largest root of $f^\prime$   is negative.  Since
  the leading coefficient of $\frac{\partial}{\partial y}f$ is
  $a_{n-1}$, the result follows.
\end{proof}

Polar derivatives preserve interlacing, but the direction depends on
certain of the coefficients.

\begin{lemma} \label{lem:polar-2}
Suppose that $f = \sum_0^n a_ix^i$ and $g = \sum_0^n b_i x^i$ satisfy
$f \greateq g$.
\begin{itemize}
\item If $a_{n-2}/a_{n-1} > b_{n-2}/b_{n-1} $ then $\partialy f \greateq
  \partialy g$
\item If $a_{n-2}/a_{n-1} < b_{n-2}/b_{n-1} $ then $\partialy g \greateq
  \partialy f$
\item If $f,g\in\allpolyalt$ or $f,g\in\allpolypos$ then $\partialy f
  \greateq \partialy g$
\end{itemize}
\end{lemma}
\begin{proof}
  Since $\partialy$ maps $\allpoly$ to $\allpoly$, it follows from
  Theorem~\ref{thm:only-roots} that $\partialy f$ and $\partialy g$ interlace.
If we compute $\partialy f$ by reversing $f$, differentiating, and
then reversing, we see that if $f$ has roots all the same sign then
$\partialy $ preserves interlacing. It remains to find the direction
of interlacing.

In general, if $\partialy f \greateq \partialy g$, then we can
differentiate $n-2$ times to find that
$$ n!a_{n-1}x + (n-1)! a_{n-2} \greateq n!b_{n-1}x + (n-1)! b_{n-2} $$
If $\partialy g \greateq \partialy f$ then this interlacing would be
reversed.  This implies the result. 
\end{proof}

We can consider the poset of all partial derivatives of a polynomial
$f$. In case $f$ is in  $\allpolypm$ then we
have predictable interlacings. For instance, 
\begin{cor} \label{cor:polar-5}
  If $f\in\allpolypos$ is a polynomial of degree $n$, then whenever  $i+j\le n$
  we have $$ 
\frac{\partial^i}{\partial x^i} \frac{\partial^j}{\partial y^j} f
\greateq 
\frac{\partial^{i-1}}{\partial x^{i-1}} \frac{\partial^{j+1}}{\partial
  y^{j+1}} f
$$
\end{cor}

The poset of partial derivatives is log concave.  
\begin{lemma} \label{lem:polar-6}
  If $f\in\allpolypos$ has all distinct roots then 
$ \smalltwodet{f}
{\frac{\partial}{\partial x}f}
{\frac{\partial}{\partial y}f}
{\frac{\partial}{\partial x}\frac{\partial}{\partial y}f} 
<0
$
\end{lemma}
\begin{proof}
  From the definition we have
$$\smalltwodet{f}
{\frac{\partial}{\partial x}f}
{\frac{\partial}{\partial y}f}
{\frac{\partial}{\partial x}\frac{\partial}{\partial y}f} 
 = 
\smalltwodet{f}{f^\prime}{nf-xf}{nf^\prime-(xf^\prime)^\prime}
=-
\smalltwodet{f}{f^\prime}{xf^\prime}{(xf^\prime)^\prime}
$$
Since $f\in\allpolypos$ we know $f\lesseq xf^\prime$, and the result
follows from Lemma~\ref{lem:inequality-1}.
\end{proof}

\section{Newton's inequalities}
\label{sec:newton}

We are going to derive relations among the coefficients of one or two
polynomials whose roots have appropriate properties.  The idea is to
differentiate  the polynomials until we get a polynomial of
degree 1 or 2, and then use properties of linear and quadratic
polynomials.

It is immediate from the quadratic formula that if that $f(x) =
ax^2+bx+c$ has two distinct real roots then $b^2 > 4ac$.  We can use
the quadratic formula to derive inequalities for the coefficients of a
polynomial that has all real roots.  Suppose that we are given a
fourth degree polynomial $f = a_0 + a_1x + a_2x^2+a_3x^3 +a_4x^4$
with all distinct roots.  Then, $\partialx f$ has all distinct roots.
Differentiating once more, the second
derivative
$ 2a_2 + 6 a_3x + 12a_4x^2$ has two distinct roots,  so we
conclude that $ 36a_3^2  > 96a_2a_4$, or $a_3^2 > (4/3)a_2a_4$.
Newton generalized this argument.

\index{Newton's inequalities}
\index{log concave}

\begin{theorem}[Newton's Inequalities] \label{thm:newton}
  Suppose that  the polynomial \\ $a_nx^n + \dots + a_0$ has all
  real roots.
  \begin{enumerate}
  \item If the roots are all distinct then 
    \begin{equation}
      \label{eqn:newton-1}
  a_k^2 > (1 + \frac{1}{k})(1+\frac{1}{n-k}) a_{k-1}a_{k+1} ,\text{ for } k=1,\dots,n-1.
    \end{equation}

    \item If the roots are not necessarily distinct, then
    \begin{align}\label{eqn:newton}
      a_k^2 > a_{k-1}a_{k+1} \text{ for } k=1,\dots,n-1. 
      \end{align} 
{unless $a_k = a_{k-1} = 0$ or $a_k=a_{k+1}=0$}
  \end{enumerate} 
\end{theorem} 
 
\begin{proof} 
  Let $F$ be the homogeneous polynomial corresponding to $f$.  Since
  partial derivatives preserve distinct real roots, if $f$ has all
  distinct roots then the polynomial
  $$
  \partX{k-1}\partY{n-k-1}F =
  n!\left(\frac{a_{k-1}}{\binom{n}{k-1}}x^2 +
    2\frac{a_k}{\binom{n}{k}}xy +
    \frac{a_{k+1}}{\binom{n}{k+1}}y^2\right).$$
  has all real roots, and they are all distinct.  Consequently,
  \begin{equation}\label{eqn:newton-2}
 \left(\frac{a_k}{\binom{n}{k}}\right)^2 >
  \frac{a_{k-1}}{\binom{n}{k-1}} \cdot \frac{a_{k+1}}{\binom{n}{k+1}}
\end{equation}
which reduces to (\ref{eqn:newton-1}).  

  If the roots of $f$ aren't distinct, then the inequality in
  \eqref{eqn:newton-1} is replaced by $\ge$.  If both $a_{k-1}$ and
  $a_{k+1}$ are non-zero, then \eqref{eqn:newton} follows.  If one of
  them is zero then the only way for \eqref{eqn:newton} to fail is if $a_k=0$.
\end{proof}

There is a kind of converse to Newton's inequalities, see
\chapsec{analytic}{rapid}.  If all the coefficients are non-zero then
we can restate the conclusion in terms of the \index{Newton quotient}\emph{Newton quotients}:
\begin{equation}
  \label{eqn:newton-1-quotient}
  \frac{a_{k}^2}{a_{k-1}a_{k+1}} \ge \frac{k+1}{k}\,\frac{n-k+1}{n-k}
\end{equation}

\begin{example}\label{ex:newton-even-part} 
  We can use the Newton quotient to get inequalities for more widely
  separated coefficients. Suppose $f(x)=\sum a_ix^i\in\allpolypos(n)$.
  From \eqref{eqn:newton-1-quotient} we get
    \begin{multline*}
      a_k^4 \ge \biggl(\frac{k+1}{k}\frac{n-k+1}{n-k}\biggr)^2
      a_{k-1}^2a_{k+1}^2 \\
 \ge \biggl(\frac{k+1}{k}\frac{n-k+1}{n-k}\biggr)^2
\frac{k}{k-1}\frac{n-k+2}{n-k+1}a_ka_{k-2}\,
\frac{k+2}{k+1}\frac{n-k}{n-k-1}a_{k+2}a_{k}
\end{multline*}
\[
\frac{a_k^2}{a_{k-2}a_{k+2}} \ge \sqrt{\frac{(k+2)(k+1)}{k(k-1)}
\frac{(n-k+2)(n-k+1)}{(n-k)(n-k-1)}}
  \]
If we iterate this inequality, then for any $r$ satisfying $k\ge 2^r$,
$n-k\ge 2^r$
\begin{equation}\label{eqn:newton-3}
\frac{a_k^2}{a_{k-2^r}a_{k+2^r}} \ge \sqrt[2^r]{\frac{\rising{k+1}{2^r}}{\falling{k}{2^r}}
\frac{\rising{n-k+1}{2^r}}{\falling{n-k}{2^r}}}
\end{equation}

  \end{example}

A similar argument to the above puts an interesting restriction on
consecutive coefficients of a polynomial with all real roots.

\begin{lemma} \label{lem:newton-zeros}
  If a polynomial $f$ has all real roots, and  two consecutive
  coefficients of $f$ are zero, then we can write $f=x^2g$ where $g$ is
  a polynomial with all real roots.
\end{lemma}

\begin{proof}
  Write $f = \sum a_ix^i$, and let $F$ be the corresponding
  homogeneous polynomial.  Since partial derivatives preserve roots,
  the polynomial
  \begin{align*}
    G(x,y) &= \frac{\partial^{k-3}}{\partial x^{k-3}}
    \frac{\partial^{n-k}}{\partial y^{n-k}}\,F    \\
    & = c_1a_kx^3 + c_2a_{k-1}x^2y + c_3a_{k-2}xy^2 + c_4a_{k-3}y^3 
  \end{align*}
  has all real roots, where $c_1,c_2,c_3,c_4$ are non-zero multiples of the
  original coefficients only depending on $n$ and $k$.  If consecutive
  coefficients $a_{k-1}=a_{k-2}=0$, then $G(x,1) = c_1a_kx^3 +
  c_3a_{k-3}$.  This polynomial does not have all real roots unless
  $a_{k-3}=0$ or $a_k=0$.  A similar argument shows that if we have
  $r$ consecutive coefficients that are $0$, then if they do not
  include the constant term we can find a block of $r+1$ consecutive
  zero coefficients. Continuing, we see that at the very least we must
  have that the constant term and the coefficient of $x$ are $0$.
\end{proof}

\index{coefficients!Newton's inequalities}

\begin{cor} \label{cor:newton}
  If $f(x) = a_0 + \cdots + a_nx^n$ has all real roots and not both
  $a_0$ and $a_1$ are zero then $a_{i}^2-a_{i+1}a_{i-1}>0$ for
$i=1,\dots,n-1$.
  \end{cor}

  We can read Lemma~\ref{lem:newton-zeros} to say that if
  two consecutive coefficients are zero for $f\in\allpoly$, then
  all earlier ones are also zero. We can slightly extend this to two
  variables. 

  \begin{lemma}\label{lem:same-degree}
    Suppose $f(x,y) =  \cdots + g(x)y^r + by^{r+1} + cy^{r+2}+\cdots$
    satisfies 
\[
f(\alpha,y)\in\allpoly \ \text{for all}\ \alpha\in\reals
\]
then $g$ has even degree, and $cg(x)\le0$ for all $x$.
  \end{lemma}
  \begin{proof}
As is usual in these arguments, we differentiate $r$ times with
respect to $y$, reverse, and differentiate and reverse again until we
have a quadratic. It has the form $\alpha g(x) + \beta by + \gamma
cy^2$ for positive $\alpha,\beta,\gamma$. Since this has roots for all
values of $x$ we must have that $g(x)c\le0$, which implies $g(x)$ must
have even degree.
  \end{proof}
 This phenomenon is possible - for instance, take $f=x^4-y^2$. The
corollary below is a two variable analog of Lemma~\ref{lem:newton-zeros}.

\begin{cor}\label{cor:same-degree}
  If $f(x,y) = \sum f_i(x)y^i$ satisfies
  \begin{enumerate}
  \item $f(\alpha,y)\in\allpoly$ for all $\alpha\in\reals$.
  \item the degree of any $f_i$ is at most $1$.
  \item Some two consecutive $f_i$ are constant.
  \end{enumerate}
  then all earlier $f_i$ are constant.
\end{cor}

\begin{proof}
  Using the lemma inductively shows that all earlier coefficients are
  constant. 
\end{proof}

  \section{More inequalities and bounds}
  \label{sec:more-inequalities}

  There are many inequalities beyond the Newton inequalities. It's
  known that the set of polynomials of fixed degree with real
  coefficients is a semi-algebraic set, and so is determined by a
  finite set of inequalities. Such a set is given by
  Theorem~\ref{thm:tp}. In this section we show how to get
  inequalities using \index{majorization}majorization and the
  Arithmetic-Geometric inequality. We apply these to determining a
  bound for polynomials on an interval.

  The following well-known inequality uses the Arithmetic-Geometric
  inequality. I don't know if it follows from Newton's inequalities.
  Lemma~\ref{lem:abs-bound} uses a weaker inequality that holds for
  all the coefficients.

\begin{lemma}\label{lem:agi}
  If $a_0+\cdots + a_nx^n\in\allpolypos$ then
\begin{enumerate}
\item $a_{n-1}^{n} \ge n^n a_0 a_n^{n-1}$
\item $a_1 \displaystyle  \biggl(\frac{a_1}{n\,a_0}\biggr)^{n-1}\ge
  a_n$
\end{enumerate}

The inequality is strict unless $f = a(x+b)^n$ for constants $a,b$.
\end{lemma}
\begin{proof}
  If the roots are $-r_1,\dots,-r_n$  then
  \begin{align*}
    \frac{a_0}{a_n} &= r_1  \cdots  r_n \\
    \frac{a_{n-1}}{a_n} &= r_1+\cdots+ r_n
  \end{align*}
Since all $r_i$ are positive the Arithmetic-Geometric inequality
implies that 
\[
\frac{1}{n}\frac{a_{n-1}}{a_n} = \frac{r_1+\cdots+r_n}{n} \ge
(r_1\cdots r_n)^{1/n} = \biggl(\frac{a_0}{a_n}\biggr)^{1/n}
\]
The second inequality follows from the first by reversing $f$.
Equality holds in the AGI if and only if all roots are equal. 
\end{proof}

\index{arithmetic-geometric inequality}
\index{inequalities!arithmetic-geometric}

We can use bounds on the coefficients to bound the size of a
polynomial with all real roots.  We begin with a simple consequence of
MacLaurin's inequality. The important point is that the upper bound
does not depend on the degree of the polynomial.

\begin{lemma}\label{lem:abs-bound}
  If $f(x) = \sum_0^n a_ix^i\in\allpolypos$ then
$\displaystyle
\sup_{|x|\le r} |f(x)| \le |a_0| \, e^{ra_1/a_0}
$
\end{lemma}
\begin{proof}
  Since $f\in\allpolypos$ we can apply  Newton's inequalities, and so
\[
\frac{a_1}{a_0}> 2\frac{a_2}{a_1} >
2\biggl(\frac{3}{2}\frac{a_3}{a_2}\biggr) >
3\biggl(\frac{4}{3}\frac{a_4}{a_3}\biggr) >
\cdots k \frac{a_k}{a_{k-1}}
\]
\index{MacLaurin's inequality}

This gives a version of MacLaurin's inequality:
\[
a_k < a_1\,\frac{1}{k!}\biggl(\frac{a_1}{a_0}\biggr)^{k-1}
\]
and therefore for $r>1$
\[
\sup_{|x|\le r} |f(x)| \le r^n \sum_{k=0}^n a_k < a_0 \sum_0^n
\biggl(\frac{a_1}{a_0}\biggr)^k r^n/k! < a_0 \exp(ra_1/a_0)
\]
\end{proof}

To find bounds for arbitrary polynomials in $\allpoly$ we need a
result  of \cite{szasz}.

\begin{lemma}\label{lem:abs-bound-2}
  If $ f(x) = a_mx^m + \cdots + a_nx^n\in\allpoly$ and $a_ma_n\ne0$ then
\[
\sup_{|x|\le r} |f(x)| \le |a_m|r^m \exp\biggl(
r\biggl|\frac{a_{m+1}}{a_m}\biggr|
+
r^2\,\biggl|\frac{a_{m+1}^2}{a_m^2}\biggr|
+
3\,r^2\,\biggl|\frac{a_{m+2}}{a_m}\biggr|\biggr)
\]
\end{lemma}

\begin{definition}
  Given two sequences $\aaa=(a_1\le a_2 \le \cdots \le a_n)$ and
  $\bbb= \\ (b_1\le b_2 \le \cdots \le b_n)$ then we say $\aaa$
  majorizes $\bbb$, written $\bbb\maj\aaa$, if
  \begin{align*}
    b_n & \le a_n \\
b_{n-1} +     b_n & \le a_{n-1} + a_n \\
\vdots \qquad& \le \qquad\vdots \\
b_2 + \cdots + b_{n-1} +     b_n & \le a_2 + \cdots + a_{n-1} + a_n \\
b_1 + b_2 + \cdots + b_{n-1} +     b_n & = a_1 + a_2 + \cdots + a_{n-1} + a_n \\
  \end{align*}
\end{definition}

There is an important inequality between monomial symmetric functions
due to Muirhead \cites{steele,hlpolya}. 

\index{Muirhead's inequality}
\index{inequality!Muirhead}

\begin{theorem}[Muirhead]
  Given two sequences $\aaa\maj\bbb$ where \\ $\aaa=(a_1\le a_2 \le
  \cdots \le a_n)$ and $\bbb=  (b_1\le b_2 \le \cdots \le b_n)$ then for
  all positive $x_1,\dots,x_n$
  \begin{equation}
    \label{eqn:muirhead}
    \sum_{\sigma\in Sym_n}
    x^{a_1}_{\sigma(1)}\,x^{a_2}_{\sigma(2)}\,\cdots
    x^{a_n}_{\sigma(n)} \le 
    \sum_{\sigma\in Sym_n}
    x^{b_1}_{\sigma(1)}\,x^{b_2}_{\sigma(2)}\,\cdots
    x^{b_n}_{\sigma(n)} 
  \end{equation}
\end{theorem}



\begin{example}
  Consider a cubic $f=(x+r_1)(x+r_2)(x+r_3)$ where $r_1,r_2,r_3$ are
  positive. We also write $f = c_0+c_1x+c_2x^2+c_3x^3$. We list the four
  index sequences that sum to $4$ and the expansion (after
  homogenizing) of the sum in Muirhead's
  theorem in terms of the coefficients:\\[.2cm]

  \begin{tabular}{rcc}
    \toprule
    Name &  & expansion \\
    \midrule
    $M_{004}$ & 0,0,4 & $2\bigl(c_2^4 - 4 c_1 c_2^2\ c_3 + 2\ c_1^2\ c_3^2 + 
      4\ c_0\ c_2\ c_3^2\bigr)$\\
    $M_{013}$ & 0,1,3 & $-c_3 \left(2 c_3 c_1^2-c_2^2 c_1+c_0 c_2
      c_3\right)$ \\
    $M_{022}$ & 0,2,2 & $2 \left(c_1^2-2 c_0 c_2\right) c_3^2$\\
    $M_{112}$ & 1,1,2 & $2 c_0 c_2 c_3^2$\\
    \bottomrule
  \end{tabular}\\[.2cm]

For example, 
\begin{align*}
  M_{112} &= r_1r_2r_3^2+  r_1r_2^2 r_3  + r_1^2r_2r_3 \\
  &= r_1r_2r_3(r_1+r_2+r_3) \\
  &= c_0c_2
\end{align*}

To homogenize we form $c_3^4\frac{c_0}{c_3}\frac{c_2}{c_3} =
c_0c_2c_3^2$.  There are four majorizations that yield (after
simplifying) inequalities

\begin{align*}
  M_{013}\maj M_{004} &\quad
  2 c_2^4-9 c_1 c_3 c_2^2+9 c_0 c_3^2
  c_2+6 c_1^2 c_3^2 &\ge 0 \\
  M_{022} \maj M_{004} &\quad
  c_2^4-4 c_1 c_3 c_2^2+6 c_0 c_3^2 c_2+c_1^2
    c_3^2&\ge0 \\
  M_{112} \maj M_{013} & \quad
  -3 c_0 c_2+c_1^2 &\ge0\\
  M_{112} \maj M_{022} & \quad
  -2 c_3 c_1^2+c_2^2 c_1-3 c_0 c_2 c_3&\ge0
\end{align*}
Although the third is Newton's inequality, calculations show that the
first one is not a consequence of the two Newton inequalities.
\end{example}

\section{The matrix of a transformation}
\label{sec:poly-det}

If $T$ is a transformation on polynomials, then we define $\varphi(T)$
to be the matrix of $T$ in the basis $\{1,x,x^2,x^3,\dots\}$. We
investigate the relation between the properties of $T$ and the
properties of $\varphi(T)$.

It is often the case that $\varphi(T)$ is totally positive. However,
$T\colon{}\allpoly\longrightarrow\allpoly$ does not imply $\varphi(T)$ is
totally positive, nor is the converse true. Consider the two examples:
  \begin{example}
    Suppose that $T(x^i) = a_ix^i$, and $D$ is the diagonal matrix
    whose $i,i$ entry is $a_i$. Clearly $D$ is totally positive, yet
    most transformations of this form do not  preserve roots - see
    Theorem~\ref{thm:polya-schur}. 
  \end{example}
  \begin{example}
    If $T(f) = f(\diffd)x^2$ then we know that
    $T\colon{}\allpoly\longrightarrow\allpoly$. However,
\[
\varphi(T) =
\begin{pmatrix}
  0 & 0 & 1 & \hdots \\ 0 & 2 & \hdots \\ 2 & 0 & \hdots \\
\vdots 
\end{pmatrix}
\]
where all the remaining entries  are zero. This matrix is not totally positive.
  \end{example}

Although a transformation that preserves roots does not necessarily
determine a totally positive matrix, we do have a result for two by
two submatrices if the direction of interlacing is preserved.

\begin{lemma}
  If the linear transformation $T$ satisfies
  \begin{enumerate}
  \item   $T\colon{}\allpolypos\longrightarrow\allpolypos$
  \item $T$ preserves the sign of the leading coefficient.
  \item $f\lesslesseq g$ implies $T(f)\lesslesseq T(g)$ for $f\in\allpolypos$.
  \end{enumerate}
then all two
  by two submatrices of $\varphi(T)$ are non-negative. (Thus
  $\varphi(T)$ is TP$_2$.)
\end{lemma}
\begin{proof}
  It suffices to show that all two by two submatrices composed of
  adjacent entries have non-negative determinants. \cite{fomin}
  Since $T(x^n)\lesslesseq T(x^{n-1})$, it follows from
  Corollary~\ref{cor:log-con-coef}  that
  the determinants are non-negative. 
\end{proof}

  We now consider several examples where $\varphi(T)$ is totally
  positive. If $M_f(g) = f\cdot g$ is the multiplication
  transformation, then we write $\varphi(f)$ for $\varphi(M_f)$. One
  important case is when $f=a+bx$.

\begin{equation}\label{eqn:tp-axb}
\phi(a+bx) = 
\begin{pmatrix}
  a & b & 0 &  0 & \dots \\
  0 & a & b &  0 & \dots \\
  0 & 0 & a &  b & \dots \\
  0 & 0 & 0 &  a & \dots \\
  \vdots & \vdots & \vdots &  \vdots & \ddots
\end{pmatrix}
\end{equation}


\noindent
For another example,  take $ f(x)=a_0+a_1x+a_2x^2+a_3x^3$.

\begin{equation}\label{eqn:tp}
\varphi(f) = 
\begin{pmatrix}
  a_0 & a_1 & a_2 & a_3 & 0 & 0 & 0 & \hdots \\
0&  a_0 & a_1 & a_2 & a_3 & 0 & 0 &  \hdots \\
0&0&  a_0 & a_1 & a_2 & a_3 & 0 &   \hdots \\
0&0&0&  a_0 & a_1 & a_2 & a_3 &    \hdots \\
\vdots&\vdots&&&   & \ddots & \ddots &  \ddots \\
\end{pmatrix}
\end{equation}

\noindent
The matrix of a derivative transformation is similar:

\begin{equation}\label{eqn:tp-d}
\phi(f\mapsto \alpha f+f') = 
\begin{pmatrix}
  \alpha & 0 & 0 &  0 & \dots \\
  1 & \alpha & 0 &  0 & \dots \\
  0 & 2 & \alpha &  0 & \dots \\
  0 & 0 & 3 &  \alpha & \dots \\
  \vdots & \vdots & \vdots &  \vdots & \ddots
\end{pmatrix}
\end{equation}

\noindent
If $ f(x)=a_0+a_1x+a_2x^2+a_3x^3$ then 

\[
\varphi(g\mapsto f(\diffd)g) =
\begin{pmatrix}
  a_0 & 0& \hdots& & & & \\
  a_1 & a_0 & 0& \hdots & & & \\
  2a_2 & 2a_1 & a_0 &0  &\hdots & & \\
  6a_3 & 6a_2 & 3 a_1 & a_0 &0 &\hdots & \\
  24a_4 & 24 a_3 & 12 a_2 & 4 a_1 & a_0 &0 & \\
\vdots & \vdots & \vdots&\vdots &\vdots &\vdots & \ddots
\end{pmatrix}
\]

\noindent
Finally, the map of an affine transformation leads to binomial
coefficients:

\[
\phi(f(x)\mapsto f(x+1)) = 
\begin{pmatrix}
  1 & 0 & 0 &  0 & \dots \\
  1 & 1 & 0 &  0 & \dots \\
  1 & 2 & 1 &  0 & \dots \\
  1 & 3 & 3 &  1 & \dots \\
  \vdots & \vdots & \vdots &  \vdots & \ddots
\end{pmatrix}
\]

The first result is well known \cite{brenti}.

\index{totally positive}
\begin{theorem}\label{thm:tp}
  If a polynomial $f$ is in $\allpolypos$ then
  $\varphi(f)$ is totally positive.  \cite{brenti}
\end{theorem}
\begin{proof}
  We show by induction on the degree of $f$ that $\varphi(f)$ is
  totally positive.  If $M_f$ is multiplication by $f$ then the
  composition of $M_f$ and $M_g$ is $M_{fg}$, and consequently
  $\varphi(fg) = \varphi(f)\varphi(g)$. From \eqref{eqn:tp-axb} it is
  immediate that $\phi(a+bx)$ is a totally positive matrix if $a$ and
  $b$ are non-negative.  Since every polynomial in $\allpolypos$ can
  be written as a product of factors of the form $a+bx$ where $a,b$
  are positive, the conclusion follows from the fact that the product
  of two totally positive matrices is totally positive.
\end{proof}

\begin{theorem}
  If $f\in\allpolypos$ and $T(g) = f(\diffd)g$ then $\varphi(T)$ is
  totally positive.
\end{theorem}
\begin{proof}
  Same as above. Just note that $f(\diffd)$ is a composition of
  operators of the form $f\mapsto \alpha f+ \beta f'$. From
  \eqref{eqn:tp-d} it is clear that $\varphi(f\mapsto \alpha f+\beta
  f')$ is totally positive. 
\end{proof}

\begin{remark}
  These last two results are equivalent since the two matrices
  $\varphi(f)$ and $\varphi(g\mapsto f(\diffd)g)$ are similar. This
  isn't surprising since differentiation and multiplication are
  adjoint. If $F$ is the infinite diagonal matrix with $i$th
  diagonal entry  $i!$ then it is easy to see that
\[
\varphi(g\mapsto f(\diffd)g) = F \,\varphi(f)^t\,F^{-1}
\]

\end{remark}

\begin{remark}
  It is actually the case that the converse of Theorem~\ref{thm:tp} is
  also true: if $\varphi(f)$ is totally positive, then
  $f\in\allpolypos$. This is known as the Aissen-Schoenberg-Whitney
  Theorem  \cite{asw}.  Also, note that the two by
  two determinants give Newton's inequalities.
\end{remark}

\index{Cauchy-Binet identity}
\index{Hurwitz's theorem!for $\intmod{d}$}
\index{totally non-negative}

The general form of Hurwitz's theorem  (see Theorem~\ref{thm:hurwitz}) follows
immediately from this remark.

\begin{cor} \label{cor:tp-hurwitz}
  If $\,\sum a_ix_i$ is in $\allpolypos$ then for any positive $d$ the
  polynomial
$\sum a_{di}x^i$ is also in $\allpolypos$.
\end{cor}
\begin{proof}
  The totally positive matrix corresponding to $\sum a_{di}x^i$ is
  formed by taking the submatrix of \eqref{eqn:tp} whose rows and
  columns are divisible by $d$.
\end{proof}

We can apply Theorem~\ref{thm:tp} to determinants involving
evaluations at integers.

\begin{cor} \label{cor:tp-2}
  Suppose $f\in\allpolypos$, $d$ is a positive integer, and
  $\alpha>0$. Then
$$
\begin{vmatrix}
  \frac{1}{d!} f(\alpha+d) & \frac{1}{(d-1)!}f(\alpha+d-1) & \hdots &
  f(\alpha)\\
\vdots & & & \vdots \\
\vdots & & & \vdots \\
  \frac{1}{(2d)!} f(\alpha+2d) & \frac{1}{(2d-1)!}f(\alpha+2d-1) & \hdots &
  \frac{1}{d!}f(\alpha+d)\\
\end{vmatrix}\ge0
$$
\end{cor}
\begin{proof}
  Note that $f(x+\alpha)\in\allpolypos$. From Lemma~\ref{lem:xnfi} we know
  that
$$ \sum_{k=0}^\infty \frac{f(i+\alpha)}{i!}x^i \in\allpolyposf$$
We can now apply Theorem~\ref{thm:tp}.
\end{proof}

If $f = \sum \alpha_i x^i$ then we can define an infinite matrix which
is just $\varphi(f)$ (see \eqref{eqn:tp}) bordered by powers of $x$.
For instance, if $f$ has degree three then

\begin{equation}
  \label{eqn:tp-newt}
  N(f) = 
\begin{pmatrix}
1&  \alpha_0 & \alpha_1 & \alpha_2 & \alpha_3 & 0 & 0 & 0 & \hdots \\
x& 0&  \alpha_0 & \alpha_1 & \alpha_2 & \alpha_3 & 0 & 0 &  \hdots \\
x^2 &0&0&  \alpha_0 & \alpha_1 & \alpha_2 & \alpha_3 & 0 &   \hdots \\
x^3 &0 &0&0&  \alpha_0 & \alpha_1 & \alpha_2 & \alpha_3 &    \hdots \\
x^4 &0 &0&0& 0& \alpha_0 & \alpha_1 & \alpha_2 &     \hdots \\
\vdots&\vdots&&&   & \ddots & \ddots &  \ddots \\
\end{pmatrix}
\end{equation}
  
\index{Newton's inequalities}
\begin{prop}\label{prop:newton-det}
  Suppose $f(x) = \sum \alpha_i x^i\in\allpolypos$. If $d$ is odd then all
  $d$ by $d$ submatrices of $N(f)$ are non-negative for positive
  $x$. 
\end{prop}
\begin{proof}
We proceed by induction on the degree  $n$ of $f$. If $n$ is $0$ then 
\[
  N(f) = 
\begin{pmatrix}
1&  \alpha_0 & 0 & 0 & 0 &     \hdots \\
x& 0&  \alpha_0 & 0  & 0 &     \hdots \\
x^2 &0&0&  \alpha_0 & 0  &       \hdots \\
x^3 &0 &0&0&  \alpha_0 &        \hdots \\
\vdots&\vdots&&&   &   \ddots \\
\end{pmatrix}
\]
We need to check that every $d$ by $d$ submatrix has non-negative
determinant. If the first column isn't part of the submatrix then this is
trivial. If the first column is part of the submatrix then it's easy
to see that the only way the determinant is not zero is if the
submatrix looks like (for $d=5$):

\[
\begin{pmatrix}
  x^m & \alpha_0 & 0 & 0 &0\\
x^{m+1} & 0 & \alpha_0 & 0 &0\\
x^{m+2} & 0& 0 & \alpha_0  &0\\
x^{m+3} &0 &  0 & 0  &\alpha_0\\
x^{m+4} &0&  0 & 0  &0\\
\end{pmatrix}
\]
The determinant is $(-1)^{d-1}\,x^{m+4}\alpha_0^{d-1}$. Since $d$ is odd, this is
positive for positive $x$. 

  If we multiply $N(f)$ on the left by $\varphi(a+bx)$ (see
  \eqref{eqn:tp}) then the result is
\begin{equation}\label{eqn:tp-3}
\begin{pmatrix}
  a+bx & a \alpha_0 & a\alpha_1+b\alpha_0 & a\alpha_2+b\alpha_1 &
  \hdots \\
  (a+bx)x & 0 & a \alpha_0 & a\alpha_1+b\alpha_0  &   \hdots \\
  (a+bx)x^2 & 0 & 0 & a \alpha_0 &   \hdots \\
\vdots & \vdots & & & 
\end{pmatrix}
\end{equation}
This is nearly $N((a+bx)f)$, except that the first column of
\eqref{eqn:tp-3} is the first column of $N((a+bx)f)$ multiplied by
$a+bx$. 

Assume that the proposition holds true for $N(f)$. Any submatrix of
$N((a+bx)f)$ not containing the first column has non-negative
determinant. If a submatrix contains the first column, then its
determinant is $(a+bx)$ times the corresponding determinant in
$\varphi(a+bx)N(f)$, which is positive by
\index{Cauchy-Binet identity}Cauchy-Binet and the induction hypothesis.

\end{proof}

\begin{remark}
  The case of $d=3$ yields the cubic Newton inequalities
  (Corollary~\ref{cor:newton-cubic}) for positive $x$.
  \index{Newton's inequalities!cubic}
\end{remark}

\section{More coefficient inequalities }
\label{sec:homog-coef}

An important consequence of Lemma~\ref{lem:inequality-1} concerns the
coefficients of interlacing polynomials.

\index{coefficients!of interlacing polynomials}
\index{interlacing!coefficients of}
\begin{cor} \label{cor:log-con-coef} If $f = a_0 + \cdots + a_nx^n$,
  $g = b_0 + \cdots + b_nx^n$, $f \longleftarrow g$, there is no index
  $k$ such that $a_k=b_k=0$, and $f,g$ have positive leading
  coefficients then for $i=0,1,\dots,n-1$ we have
  \begin{equation}
    \label{eqn:log-con-coef}
      \smalltwodet{b_{i}}{b_{i+1}}{a_{i}}{a_{i+1}} > 0
  \end{equation}
\end{cor}

\begin{proof}
 From Lemma~\ref{lem:inequality-1} we have that 
\[
0 >\smalltwodet{f(0)}{g(0)}{f^\prime(0)}{g^\prime(0)} =
 \smalltwodet{a_0}{a_1}{b_0}{b_1}.
\]
since zero is not a common root.  In the general case we differentiate
$f\greateq g$ exactly $i-1$ times. The result is the interlacing $$
i!a_i + \falling{i+1}{i}a_{i+1} + \cdots +\falling{n}{i} x^{n-i+1}
\greateq i!b_i + \falling{i+1}{i}b_{i+1} + \cdots+ \falling{n}{i}
x^{n-i+1} $$
We can now apply the case $i=0$ to complete the argument.
\end{proof}

\begin{remark}\label{rem:interlace-picture}
There is a useful picture to help remember this result that will reappear in
higher dimensions. If we write the coefficients of $f$ and $g$ in the
diagram

\centerline{\xymatrix{
b_0 \ar@{-}[r]\ar@{-}[d]\ar@{-}[dr] &
b_1 \ar@{-}[r]\ar@{-}[d]\ar@{-}[dr] &
b_2 \ar@{-}[r]\ar@{-}[d]\ar@{-}[dr] &
b_3 \ar@{..}[r]\ar@{-}[d]\ar@{..}[dr] & \\
a_0 \ar@{-}[r] &
a_1 \ar@{-}[r] &
a_2 \ar@{-}[r] &
a_3 \ar@{..}[r] &
}}

\noindent%
then in each pair of adjacent triangles of a square the product of the
coefficients at the diagonally adjacent vertices is greater than the product of
the coefficients at the remaining two vertices.
\end{remark}

 \begin{remark}
    We can  derive these inequalities from Newton's inequalities. In
    fact, if $f = \prod_i(x-a_i)$, and $f_j = f/(x-a_j)$, then
    Newton's inequalities imply that the determinant
    \eqref{eqn:log-con-coef} is positive for the functions $f$ and
    $f_j$. Now the determinant is linear in the coefficients of $g$,
    so the coefficients of any convex combination of $f,f_1,\dots$
    also satisfy \eqref{eqn:log-con-coef}. These convex combinations
    describe all the possible $g$'s by Lemma~\ref{lem:sign-quant}.
  \end{remark}

If we restrict $f,g$ to lie in $\allpolypos$ then

\begin{lemma} \label{lem:log-con-coef-2}
  If $f=a_0+\cdots+a_nx^n$, $g=b_0+\cdots+b_nx^n$, $f\greateq g$,
  $f,g\in\allpolypos$ and $f,g$ have positive leading coefficients
  then for all $0\le r < s \le n$:
$$ 
\begin{vmatrix}   b_{r} & b_s\\ a_{r} & a_s   \end{vmatrix} > 0
$$
\end{lemma}
\index{determinants!of coefficients}
\begin{proof}
If $f,g\in\allpolypos$ then all coefficients are positive, so from
Corollary~\ref{cor:log-con-coef} we have the inequalities
\[
\frac{a_{r+1}}{a_r}>\frac{b_{r+1}}{b_r}\qquad
\frac{a_{r+2}}{a_{r+1}}>\frac{b_{r+2}}{b_{r+1}}\qquad\cdots\qquad
\frac{a_s}{a_{s-1}}>\frac{b_s}{b_{s-1}}
\]
Multiplying these inequalities yields the conclusion.

  
\end{proof}

In Lemma~\ref{lem:log-con-coef-2} we extended the Corollary~\ref{cor:log-con-coef} to
determinants of the form $\smalltwodet{a_{i+d}}{a_i}{b_{i+d}}{b_i}$.
Next, we observe that the coefficients of interlacing polynomials can
not agree very often:

\index{interlacing!coefficients of}
\begin{cor}
  If $f\greateq g$ are both in $\allpolypos$, then $f$ and $g$ can not
  agree at any two coefficients.
\end{cor}
\begin{proof}
  If $f,g$ have coefficients as above, and $a_k=b_k$,
  $a_{k+d}=b_{k+d}$, then  Lemma~\ref{lem:log-con-coef-2} is contradicted,
  for the determinant is 0.
\end{proof}


Here's another variation.

\begin{lemma}\label{lem:det-not-0}
  If $f=\sum a_ix^i$ and $g = \sum b_ix^i$ interlace, are linearly
  independent, and $\smalltwodet{a_0}{b_0}{a_1}{b_1}=0$ then
  $a_0=b_0=0$.
\end{lemma}
\begin{proof}
  Since the determinant is zero there is a $\lambda$ such that
  $\lambda(a_0,a_1) = (b_0,b_1)$. If $\lambda\ne0$ then by linear
  independence $\lambda f-g\ne0$, and $f$ and $\lambda f - g$
  interlace. Since $\lambda f - g$ has a factor of $x^2$, $f$ has a
  factor of $x$. Considering $g$ and $\lambda f-g$ shows that $g$ has
  a factor of $x$.

  If $\lambda=0$ then $g$ has a factor of $x^2$, so $f$ has a factor
  of $x$.
\end{proof}

\begin{remark}
  If we form the matrix of coefficients of three mutually interlacing
  quadratics $f\greateq g\greateq h$, then we know that all two by two
  determinants are positive (Corollary~\ref{cor:log-con-coef}), but
  the three by three determinant can be either positive or negative.
  In particular, the matrix of coefficients is not totally positive.
  \index{totally positive}
    \begin{align*}
      h &= (4 + x) (12 + x) & h &=  (8 + x) (12 + x) \\
      g &= (3 + x) (9 + x) & g&= (7 + x) (10 + x)\\
      f & = (2 + x) (8 + x) & f &=(6 + x) (9 + x)   \\
      -2 & = \begin{vmatrix} 
48& 16& 1\\
27& 12& 1 \\ 
16& 10& 1  
\end{vmatrix} & 4 &=
\begin{vmatrix}
  96& 20& 1\\
  70& 17& 1\\
  54& 15& 1
\end{vmatrix}
    \end{align*}

  \end{remark}

  Next we have a kind of reverse Newton's inequalities for two by two
  determinants of interlacing polynomials.

\index{Newton's inequalities!for interlacing polynomials}
  \begin{lemma}\label{lem:newton-interlacing}
    Suppose that $f=\sum a_ix^i$ and $g = \sum b_ix^i$ are in
    $\allpoly(n)$. If $f\greateqeq g$ then
\[
{ \begin{vmatrix}  b_{k} & b_{k+2} \\ a_{k} & a_{k+2}  \end{vmatrix}^2 }
\le {4} \,\cdot\,\frac{(k+2)(n-k)}{(k+1)(n-k-1)}\,\cdot
{
  \begin{vmatrix}  b_k & b_{k+1} \\ a_k & a_{k+1} \end{vmatrix}
  \begin{vmatrix}  b_{k+1} & b_{k+2} \\ a_{k+1} & a_{k+2} \end{vmatrix}}
\]
  \end{lemma}
  \begin{proof}
    Since $f\greateqeq g$ we know that $f+\alpha g\in\allpoly$ for all
    $\alpha$.  If we let $\beta = (k+2)(n-k)/((k+1)(n-k-1))$ then Newton's
    inequalities state that for all $\alpha$ we have
\[
(a_{k+1}+\alpha b_{k+1})^2 \ge \beta\,(a_k+\alpha b_k)(a_{k+2}+\alpha
b_{k+2})
\]
Since the difference is non-negative for all $\alpha$ the discriminant is
non-positive. A computation of the discriminant shows that
\[ 
\beta^2 
 \begin{vmatrix}  b_{k} & b_{k+2} \\ a_{k} & a_{k+2}  \end{vmatrix}^2
 - 4\,\beta 
  \begin{vmatrix}  b_k & b_{k+1} \\ a_k & a_{k+1} \end{vmatrix}
  \begin{vmatrix}  b_{k+1} & b_{k+2} \\ a_{k+1} & a_{k+2}
  \end{vmatrix} 
\le 0
\]
and  the conclusion follows. 
  \end{proof}

An equivalent formulation is that the determinant below has all real
roots.

\[
\begin{vmatrix}
  1 & \sqrt{\beta}\,x & x^2 \\
b_k & b_{k+1} & b_{k+2} \\
a_k & a_{k+1} & a_{k+2} 
\end{vmatrix}
\]

Rosset \cite{rosset} called the following corollary   the cubic Newton inequality 
\begin{cor}\label{cor:newton-cubic}
  If $f(x) = \sum a_ix^i\in\allpolypos(n)$ and $0\le k\le n-3$ then 
  \begin{equation*}
 4\,\left( {a_{k+1}}^2 - a_{k}\,a_{k+2} \right) \,\left( {a_{k+2}}^2 -
   a_{k+1}\,a_{k+3} \right) \ge
{\left( a_{k+1}\,a_{k+2} - a_{k}\,a_{k+3} \right) }^2
  \end{equation*}
\end{cor}
\begin{proof}
  Apply the lemma to $xf$ and $f$. 
\end{proof}

  \section{Log concavity in $\allpolypos$}
  \label{sec:log-concave-p}

  The logarithms of the coefficients of a polynomial in $\allpolypos$
  have useful and well-known concavity properties. We will revisit
  these ideas when we consider polynomials in two variables.

\index{Newton graph}
\index{log concave}

To begin, suppose that $f = a_0 + \cdots +a_nx^n$ has all positive
coefficients and is in $\allpolypos$. We define the Newton graph of
$f$ to be the collection of segments joining the points
\[
(0,\log a_0),\,(1,\log a_1),\cdots,\,(n,\log a_n).
\]

\noindent%
For instance, Figure~\ref{fig:newton-graph} is the Newton graph for
$n=4$. 

\begin{figure}[htbp]\label{fig:newton-graph}
  \centering
\begin{pspicture}(-1,0)(8,3)
 \psline{o-o}(0,0)(2,2)(4,2.5)(6,2.2)(8,1)  
\rput[r](0,0){$(0,\log a_0)\quad$}
\rput[r](2,2){$(1,\log a_1)\quad$}
\rput[b](4,2.5){$(2,\log a_2)$}
\rput[l](6,2.2){$\quad(3,\log a_3)$}
\rput[l](8,1){$\quad(4,\log a_4)$}
\end{pspicture}
  
  \caption{The Newton graph of a polynomial in $\allpolypos$}
\end{figure}

The Newton graph is concave. Indeed, the slope of the segment between
$(k,\log a_k)$ and $(k+1,\log a_{k+1})$ is $\log(a_{k+1}/a_k)$. Taking
logs of Newton's inequality yields
\[ \log\frac{a_{k+1}}{a_k} \ge \log\frac{a_{k+2}}{a_{k+1}},
\]
which means that the slopes are decreasing. This implies that the
Newton graph is concave.

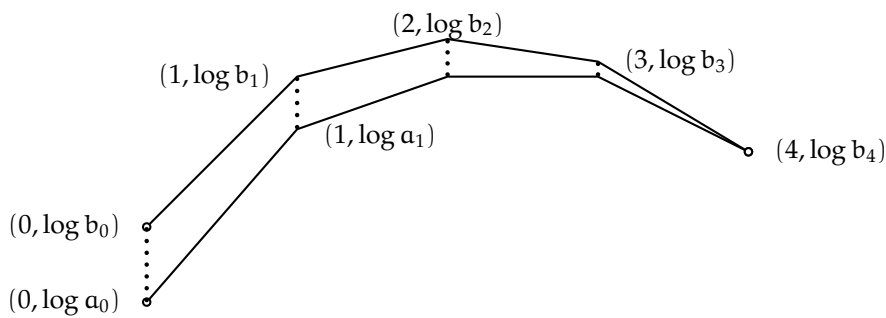
\begin{figure}[htbp]\label{fig:newton-graph-2}
  \centering
\begin{pspicture}(-1,-1)(8,3)
 \psline{o-o}(0,0)(2,2)(4,2.5)(6,2.2)(8,1)  
 \psline{o-o}(0,-1)(2,1.3)(4,2)(6,2)(8,1)  
\psline[linewidth=.6mm,linestyle=dotted](0,0)(0,-1)
\psline[linewidth=.6mm,linestyle=dotted](2,2)(2,1.3)
\psline[linewidth=.6mm,linestyle=dotted](4,2.5)(4,2)
\psline[linewidth=.6mm,linestyle=dotted](6,2.2)(6,2)
\rput[r](0,0){$(0,\log b_0)\quad$}
\rput[r](0,-1){$(0,\log a_0)\quad$}
\rput[r](2,2){$(1,\log b_1)\quad$}
\rput[l](2,1.2){$\quad(1,\log a_1)$}
\rput[b](4,2.5){$(2,\log b_2)$}
\rput[l](6,2.2){$\quad(3,\log b_3)$}
\rput[l](8,1){$\quad(4,\log b_4)$}
\end{pspicture}
\caption{The Newton graphs of two interlacing polynomials}
\end{figure}

Next, we consider the Newton graphs of two interlacing polynomials.
Assume that $f\greateq g$ in $\allpolypos$. In order to compare
coefficients we will assume that both are monic, so that  their
Newton graphs agree at the right hand endpoint. The Newton graph of
$f$ lies below the Newton graph of $g$, and the gap decreases as the
index increases. For example, Figure~\ref{fig:newton-graph-2} shows
$f$ and $g = \sum b_i x^i$. The dotted lines decrease in length as we
move to the right.

We use  Lemma~\ref{lem:log-con-coef-2} which says that
$\smalltwodet{b_k}{b_{k+1}}{a_k}{a_{k+1}}\ge0$. Taking logs shows that
\begin{align*}
  \log\frac{b_{k+1}}{b_k} & \le \log\frac{a_{k+1}}{a_k} & 
  \log\frac{b_{k+1}}{a_{k+1}} & \le \log\frac{b_{k}}{a_k} 
\end{align*}
and these inequalities imply the assertions of the last paragraph.

\section{A pair of recursions}
\label{sec:homog-pair}

In \cite{zhedanov} there is an interesting pair of recursions

\begin{align*}
p_{n+1}(z)&=z\,p_n(z)-\alpha_n\,z^n\,q_n(1/z) \\
q_{n+1}(z)&=z\,q_n(z)-\beta_n\,z^n\,p_n(1/z).
\end{align*}
If we homogenize these equations then we can find recurrences for
$p_n$ and $q_n$. Using these these recurrences we can easily show that
$p_n$ and $q_n$ are in $\allpoly$ for various choices of $\alpha_n$
and $\beta_n$. So consider $P(x,y)$ and $Q(x,y)$ satisfying the
recurrences
 
\begin{align*}
   P_{n+1}(x,y) &= x\,P_n(x,y) -\alpha_n\,y\,Q_n(y,x) \\
   Q_{n+1}(x,y) &= x\,Q_n(x,y) -\beta_n\,y\,P_n(y,x) \\
\intertext{with the initial conditions $P_0(x,y)=Q_0(x,y)=1$. We
  compute $Q_{n+1}(y,x)$ using each equation}
\frac{1}{x\,\alpha_{n+1}}\left(y\,P_{n+1}(y,x)-P_{n+2}(y,x)\right) &=
\frac{1}{x\,\alpha_{n}}\left(y\,P_{n}(y,x)-P_{n+1}(y,x)\right) - \\
& \quad \beta_n\,y\,P_n(y,x)
\end{align*}
\noindent{ Interchange $x$ and $y$ and simplify}
\begin{align*}
  P_{n+2}(x,y) &=
  \left(x+\frac{\alpha_{n+1}}{\alpha_n}\,y\right)P_{n+1}(x,y) +
  \left(\alpha_{n+1}\beta_n -
    \frac{\alpha_{n+1}}{\alpha_n}\right)\,xy\,P_n(x,y)\\
  \intertext{Since $P_n$ is homogeneous $P_n(x,1)=p_n(x)$ and }
  p_{n+2}(x) &=
  \left(x+\frac{\alpha_{n+1}}{\alpha_n}\,\right)p_{n+1}(x) +
  \left(\alpha_{n+1}\beta_n -
    \frac{\alpha_{n+1}}{\alpha_n}\right)\,x\,p_n(x)\\
  \intertext{If the coefficient of $p_n$ is negative then this defines
    a sequence of orthogonal Laurent polynomials
    \mypage{lem:recursive-laurent}. In the special case that
    $\alpha_n=\beta_n=\alpha$ for all $n$ then $p_n = q_n = f_n$ and
    we get the simple recurrence}
 \end{align*}

\vspace*{-1cm}

 \begin{equation}
   \label{eqn:bi-1}
f_{n+2} = (x+1)f_{n+1} + x(\alpha^2-1)f_n   
 \end{equation}

\begin{lemma}
  If $|\alpha|\le1$, $f_0=1,f_1=x-\alpha$ and $f$ satisfies
  \eqref{eqn:bi-1} then $f_n\in\allpolyalt$ and $f_{n+1}\lesslesseq f_n$.
\end{lemma}

 Another interesting case is when $\alpha_n=\beta_n = -q^{(n+1)/2}$.
 The resulting polynomials are known as the Rogers-Szego polynomials.
\index{Rogers-Szego polynomials}
\index{polynomials!Rogers-Szego}
The recurrence is 
 \begin{equation}
   \label{eqn:bi-2}
f_{n+2} = (x+1)f_{n+1} + x(q^{n+1}-1)f_n   
 \end{equation}

It follows that
\begin{lemma}
 If $q\ge1$, $f_0=1,f_1=x+q$ and $f$ satisfies
  \eqref{eqn:bi-2} then $f_n\in\allpolyalt$ and $f_{n+1}\lesslesseq f_n$.
  \end{lemma}


\chapter{Analytic Functions}
\label{cha:analytic}

\renewcommand{\TimeStampStart}{Monday, December 17, 2007: 17:22:10}
\mytoday
 
We introduce some classes of analytic functions that are the uniform
closure of certain sets of polynomials with all real roots.  We are
interested in the relationship between these classes, and the set of
generating functions of linear transformations of the form
$x^i\mapsto a_ix^i$ that preserve roots.

In later chapters we will generalize these results to more general
linear transformations, and to polynomials in several variables. 

\section{The uniform closure of $\allpolypos$}
\label{sec:anal-neg}

\index{\ F@$\allpolyposf$}

We begin with the simplest class of analytic functions, the closure of
polynomials with all positive coefficients and all real roots.

\begin{definition}
  $\allpolyposf$ is the set of all functions that are in the closure of
  $\allpolypos$ under uniform convergence on compact domains.
\end{definition}

From basic complex analysis we know
\begin{enumerate}
\item All functions in $\allpolyposf$ are analytic.
\item All functions in $\allpolyposf$ are \emph{entire} - this 
means that they are defined on the entire complex plane.
\item All derivatives of functions in $\allpolyposf$ are in
  $\allpolyposf$.
\item The product of two functions in $\allpolyposf$ is in $\allpolynegf$.
\item Every function in $\allpolyposf$ has all real roots.
\end{enumerate}

It is obvious that $\allpolypos\subset\allpolynegf$. The limit of
$\allpolyalt$ is $\allpolyaltf$, and satisfies the above properties,
except that $\allpolyalt\subset \allpolyaltf$. 
We now introduce  a few important non-polynomial functions in
$\allpolyposf$.

\begin{example} \emph{The exponential function.}
 $e^x$ is the most
important non-polynomial function in $\allpolyposf$. To see
that it is in $\allpolyposf$ note that
$$ \lim_{n\rightarrow\infty} \,\left(1+\frac{x}{n}\right)^n \ = e^x$$ 
so $e^x$ is the limit of polynomials  that have $n$-fold zeros at 
$-n$. These zeros go off to infinity, so $e^x$ has no zeros. 
More generally, $e^{ax}$ is in $\allpolyposf$ for any positive $a$. If
$a$ is negative then $e^{ax}\in\allpolyaltf$. 

Also, since 
$$ \lim_{n\rightarrow\infty} \left(1+x+\frac{x^2}{2!} + \dots +
\frac{x^n}{n!}\right) \ = e^x$$ 
we see that $e^x$ is also a limit of polynomials that are \emph{not}
in $\allpoly$, so $e^x$ is on the boundary of $\allpolyposf$.

\end{example}

\begin{example} \label{ex:q-exponential} \emph{The $q$-exponential function.}
  There is a $q$-analog of the exponential function. We first need to
  introduce the Pochhammer symbols \index{Pochhammer symbol}
  \begin{align*}
    (a;q)_n &= (1-a)(1-aq)\cdots(1-aq^{n-1})\\
    (a;q)_\infty &= (1-a)(1-aq)\cdots(1-aq^{n-1})\cdots\\
  \end{align*}
The following  basic identity \cite{gasper}*{(1.3.16)} is valid for
$|q|<1$:
\begin{equation} \label{eqn:q-exponential}
 \sum_{n=0}^\infty \frac{q^{\binom{n}{2}}}{(q;q)_n} z^n =
 (-z;q)_\infty
\end{equation}
We denote this function by $Exp_q(z)$. If we write out the right hand
side of \eqref{eqn:q-exponential}
$$ (-z;q)_\infty = (1+z)(1+qz)(1+q^2z)\cdots$$
then we see that for $0<q<1$ all the roots of $Exp_q(z)$ are negative, and so
$Exp_q(z)$ is in $\allpolyposf$. 
We will later use this variation (see \cite{askey}) on the
$q$-exponential function, where $[n] = (q^n-1)(/(q-1)$:
\begin{equation}
  \label{eqn:q-exponential-2}
  E_q(z)  = \sum_{n=0}^\infty
    \frac{q^{\binom{n}{2}}}{[n]!} z^n = ((1-q)z;q)_\infty \quad\quad |q|<1\\
\end{equation}

\end{example}

\begin{remark}
  $\allpolyposf$ is closed under products, so we can construct new
  elements in $\allpolyposf$ by multiplying appropriate factors. First
  of all, $x^r$ and $e^{ax}$ are always in $\allpolyposf$ for $r$ a
  non-negative integer, and $a$ positive. If the infinite product
  $\prod_{1}^\infty(1+ c_ix)$ converges, where all $c_i$ are positive,
  then the following is in $\allpolyposf$:

  \begin{equation}\label{eqn:type-1}
  c x^r
  e^{ax} \prod_{i=1}^\infty (1+c_i x)
\end{equation}

A famous result due to \Polya-Schur \cites{bilodeau,polya-szego2}
states that this representation describes all the elements of
$\allpolyposf$.
\end{remark}

\section{The uniform closure of $\allpoly$}
\label{sec:analytic-allpoly}

$\allpolyf$ is the uniform closure of $\allpoly$ on compact subsets. It
satisfies the same general analytic properties that $\allpolyposf$
does. Here are some examples of analytic functions in $\allpolyf$.

\begin{example} \emph{The Gaussian density.}
Since 
$$ e^{-x^2} = \lim_{n\rightarrow\infty}\left(1-\frac{x^2}{n}\right)^n
= \lim_{n\rightarrow\infty} 
\left(1-\frac{x}{\sqrt{n}}\right)^n
\left(1+\frac{x}{\sqrt{n}}\right)^n
$$ we see that 
 $$e^{-{x^2}} =
 1 - \frac{x^2}{1!} + \frac{x^4}{2!} - \frac{x^6}{3!} +  \dots 
$$
is in $\allpolyf\setminus\allpolyposf$ since it is the limit of
polynomials with both positive and negative roots.
\end{example}

\begin{example} \emph{Sine and Cosine.}
  The sine and cosine  are in $\allpolyf$ since they both have
  infinite products with positive and negative roots:
\index{sin and cos}
\begin{align*}
  \sin(x) &= x\,\prod_{k=1}^\infty 
\left(1 - \frac{x^2}{k^2\pi^2}\right) \\
 & =  \frac{x}{1} -  \frac{x^3}{3!} + \frac{x^5}{5!} - \cdots \\
\cos(x) & = \prod_{k=0}^\infty
\left(1 - \frac{4x^2}{(2k+1)^2\pi^2}\right) \\
& = 1 - \frac{x^2}{2!} + \frac{x^4}{4!} - \cdots 
\end{align*}
\end{example}

\index{Bessel function}
\index{\ aaaJi@$J_i$|see{Bessel function}}
\begin{example} \emph{Bessel Functions.}
The Bessel function of the first kind, $J_i(z)$, for $i$ a positive
integer, is defined as 

$$ J_i(z) = 
\frac{z^i}{2^i} \sum_{k=0}^\infty (-1)^k \frac{z^{2k}}{2^{2k}k!(k+i)!}
$$

There is  an infinite product representation \cite{grads}*{page 980}, where
we denote the zeros of $J_i(z)$ by $z_{i,m}$:

$$ J_i(z) = 
\left(\frac{z}{2}\right)^i \frac{1}{i!} \prod_{m=1}^\infty \left(1 -
\frac{z^2}{z_{i,m}^2}\right) 
$$

The product representation shows that $J_i(z)$ is in $\allpolyf$. 

Inequalities for elements in $\allpolyf$ follow by taking
  limits of inequalities for $\allpoly$. For instance, Newton's
  inequalities \eqref{eqn:newton} become $a_k^2 \ge a_{k-1}a_{k+1}$.
\index{Turan's inequality}\index{Newton's inequalities}
\index{Legendre polynomials}\index{polynomials!Legendre}
Using this observation,  we derive Turan's inequality:
  If $P_n(x)$ is the $n$-th Legendre polynomial, then 
\begin{equation}
  \label{eqn:turan}
 P_n(x)^2 - P_{n+1}(x)P_{n-1}(x) \ge 0\quad \text{ for } |x|<1 
\end{equation}

To see this, we need the identity
$$
  \sum_{n=0}^\infty P_n(x)\frac{y^n}{n!} = e^{xy} J_0(y\sqrt{1-x^2})
$$
Thus, if $|\alpha|<1$ the series $\sum P_n(\alpha)\frac{y^n}{n!}$ is
in $\allpolyf$. An application of Newton's inequalities
\eqref{eqn:newton-2} yields \eqref{eqn:turan}. See \cite{szego-turan}.

\end{example}

\index{Gamma function}\index{Gamma function}
\begin{example} \label{ex:gamma} \emph{The Gamma function}
  Although the Gamma function is meromorphic, its reciprocal is an
  entire function with the product \cite{polya-szego2}*{page 242}:
$$
\frac{1}{\Gamma(z+1)} \ =\ 
e^{\gamma z}\,\left(1+\frac{z}{1}\right)e^{-\frac{z}{1}}
\left(1+\frac{z}{2}\right)e^{-\frac{z}{2}} \cdots
\left(1+\frac{z}{n}\right)e^{-\frac{z}{n}} \cdots
$$
where $\gamma$ is Euler's constant.  The roots are at the negative
integers $-1,-2,\dots$. Moreover, if we choose a positive integer $k$
then the quotient
$$G(z) = \frac{\Gamma(z+1)}{\Gamma(kz+1)}$$
is an entire function since the poles of the numerator are contained
among the poles of the denominator. In addition, $G(z)$ is in $\allpolyf$ since
it has a product formula.
\end{example}

\begin{remark}
\index{P\'{o}lya-Schur}
There is an explicit description \cites{bilodeau,polya-szego2} of the
elements of $\allpolyf$ that is slightly more complicated than the
description of the elements of $\allpolyposf$. 
The elements of $\allpolyf$ have the product representation
\begin{equation} 
 \label{eqn:type-2}
  c x^r e^{ a x-b x^2} \prod_{i=1}^\infty
  (1+c_ix)e^{-c_ix}
\end{equation}
where $a,c_1,\dots$ are real, $r$ is a non-negative integer,
and $b$ is non-negative. Notice that  $e^{x^2}$ is not in
$\allpolyf$. 
\end{remark}

    We can use trigonometric identities to find some non-obvious members
    of $\allpolyf$. Recall that
$$\frac{1}{2} \sin \frac{x+y}{2} \,\cos \frac{x-y}{2} = \sin x+\sin
y$$
If we choose $a$ so that $|a|\le1$, then we can find $y$ so that $\sin
y = a$, and hence 
\[
\sin(x) + a = \frac{1}{2}  \sin \frac{x+a}{2} \,\cos
\frac{x-a}{2}
\]
 Since the right hand side is a product of members of
$\allpolyf$, it follows that $a+\sin x\in\allpolyf$. Multiplying
several of these together  yields half of
\index{sin and cos}
\begin{lemma}\label{lem:f-of-sin}
  Suppose that $f(x)$ is a polynomial. Then,
  $f(x)\in\allpolyint{[-1,1]}$ if and only if $f(\sin x)\in\allpolyf$.
\end{lemma}
\begin{proof}
  Conversely, let $r_1,\dots,r_n$ be the roots of $f(x)$. Then the
  roots of $f(\sin x)$ are solutions to $\sin x = r_i$. If any $r_i$
  has absolute value greater than $1$ or is not real, then $\arcsin
  r_i$ is also not real, so $f(\sin x)$ has non-real roots. 
\end{proof}

\section{The Hadamard Factorization Theorem}
\label{sec:hadamard-fac}

How can we show that a function $f$ is in $\allpolyf$ or
$\allpolyposf$? It is not enough to know that all the roots are real -
for example $e^{x^2}$ has no complex roots and yet is not in
$\allpolyf$.  If all roots of $f$ are real we can determine if $f$ is
in $\allpolyf$ using the coefficients of the Taylor series.

The Hadamard factorization theorem represents an arbitrary
entire function as an infinite product times an exponential factor. We
first summarize the necessary background (see \cite{rubel} or
\cite{boas}), and then use the machinery to exhibit a collection of
functions whose exponential part satisfies the necessary conditions to
belong to $\allpolyf$.

We begin with definitions. If $p$ is a non-negative integer,
we define the \emph{Weierstrass primary factor of order $p$} by
\begin{align*}
  E(z,0) &= (1-z) \\
  E(z,p) &= (1-z)\,exp(z+\frac{z^2}{2} + \cdots + \frac{z^p}{p})
\end{align*}
For a set $Z$ of complex numbers, \emph{the genus of $Z$, $p=p(Z)$} is
defined by
$$
p(Z) = \inf\left\{ q: q \text{ is an integer, } \sum_{z_i\in Z}
  \frac{1}{|z_i|^{q+1}} < \infty \right\}$$
The \emph{canonical  product  over $Z$} is defined by
$$ P_Z(z) = \prod_{z_i\in Z}  E\left( \frac{z}{z_i},p\right)$$
where $p = p(Z)$.  The maximum modulus of $f(z)$ on the circle $|z|=r$
is denoted $M(r)$. The \emph{order $\rho=\rho(f)$ of $f(z)$} is
$$\rho = \limsup_{r\rightarrow\infty} \frac{\log\log M(r)}{\log r}$$
We can now state 
\begin{theorem}[Hadamard Factorization Theorem]
If $f(z)$ is an entire function of order $\rho$ with an $m$-fold zero
at the origin, we have 
$$ f(z) = z^m e^{Q(z)}\, P_Z(z)$$
where $Q(z)$ is a polynomial of degree at most $\rho$, and $Z$ is the
set of zeros of $f$.  
\end{theorem}

The particular case we are interested in is $\rho<1$ and $p=0$. In
this case
\begin{equation} \label{eqn:desired-rep}
f(z) = z^m e^{\alpha z} \prod_{i}\left(1-\frac{z}{z_i}\right)
\end{equation}
If all the roots $z_i$ are negative then $f\in\allpolyposf$, and
$f\in\allpolyf$ if there are positive and negative roots. The main
difficulty is showing that all the roots are real.

How can we determine the two  parameters $\rho$ and $p$? First of all,
it is known that $p\le \rho$. 
Second, it is possible to find $\rho$ from the Taylor series of
$f$. Suppose that
$$ f(z) = \sum_{i=0}^\infty a_i z^i$$
If the limit 
\begin{equation} \label{eqn:find-rho}
\limsup_{n\rightarrow\infty} \frac{ n \log n}{\log (1/|a_n|)}
\end{equation}
is
finite then its value is $\rho$. 


We are in a position to determine the exponential part. We begin with
the Bessel function $J_0(-i\sqrt{x})$.

\begin{example}\textbf{The Bessel function.} \label{ex:entire-bessel}  %
  \index{Bessel function}
  $$J_0(-i\sqrt{z})\,=\,\sum_{i=0}^\infty \frac{z^i}{i!\,i!} \in \allpolyposf$$
  We first recall Stirling's formulas: \index{Stirling's formula} %
  \begin{align*}
    n! & \sim \ \sqrt{2\pi n}\ n^n e^{-n} \\
\log n! & \sim \   \frac{1}{2} \log(2\pi n) +n \log n -n
  \end{align*}
The limit \eqref{eqn:find-rho} is
$$ \limsup_{n\rightarrow\infty} \frac{n \log n}{\log (n!\, n!)} =
\limsup_{n\rightarrow\infty} \frac{n \log n}{2(\frac{1}{2} \log(2\pi n) +n \log n -n)} =
\frac{1}{2}
$$
Thus $\rho=\frac{1}{2}$, and hence $p=0$. This shows that we have
the desired representation \eqref{eqn:desired-rep}, and since all
terms are positive there are no positive zeros. Since the Bessel
function is known to have all real zeros the series is in
$\allpolyposf$. 

The general Bessel function has a series expansion
\begin{equation}
  \label{eqn:bessel-v}
  J_\nu(z) = \sum_{k=0}^\infty \frac{(-1)^k}{\Gamma(k+\nu+1)
    \,k!}\,\left(\frac{z}{2}\right)^{\nu+2k} 
\end{equation}
$J_\nu(z)$ has infinitely many real roots for all $\nu$, and all of
its roots are real for $\nu\ge-1$. Since a similar computation shows
that $p=0$ for $J_\nu(z)$, it follows that $J_\nu\in\allpolyf$ for
$\nu\ge-1$.

\end{example}

We now give some examples of computations of $\rho$; since we don't know
if the functions have all real roots, we can't conclude that they are
in $\allpolyf$. 




\begin{example}\label{lem:entire-fract}
  Suppose that $a_1,\dots,a_r$ and $b_1,\dots,b_s$ are positive. Then
$$ \sum_{n=0}^\infty
\frac{(a_1n)!\cdots(a_rn)!}{(b_1n)!\cdots(b_sn)!}z^n \text{ has } p<1 $$ 
if $\left(\sum_1^s b_k\right) - \left(\sum_1^r a_k\right)>1$.

  Compute the limit using Stirling's formula:
\begin{align*}
\limsup \frac{n \log n}{\log \prod (b_kn)! - \log \prod (a_k n)!} &=
\limsup \frac{n \log n}{\sum b_kn \log b_kn - \sum a_k n \log a_kn}
\\
&= \frac{1}{\sum b_k - \sum a_k}
\end{align*}
If the hypothesis are satisfied, then the last expression is less than
$1$, and $p<1$.

\end{example}

\begin{example}
  Here is a similar example.  Suppose that $a_1,\dots,a_r$ and
  $b_1,\dots,b_s$ are positive.  We define the hypergeometric
  series \index{hypergeometric series}
\[
\hypergeo{r}{s}(a_1,\dots,a_r;b_1,\dots,b_s;z) =
 \sum\frac{\rising{a_1}{n}\cdots\rising{a_r}{n}}{
\rising{b_1}{n}\cdots\rising{b_s}{n}}\frac{z^n}{n!} 
\]
  If we compute the limit as above, we find that it is
  $\frac{1}{1+s-r}$. 
 If $s-r>0$ then   $p<1$. Hurwitz showed that the series

 \begin{equation}
   \label{eqn:hypergeo}
\hypergeo{0}{s}(;b_1,\dots,b_s;z) =
 \sum\frac{1}{
\rising{b_1}{n}\cdots\rising{b_s}{n}}\frac{z^n}{n!} 
 \end{equation}
has all real roots if $b_1\ge -1$ and all the other $b_i$ are positive.

\end{example}

\begin{example} \label{ex:entire-q2} \textbf{A $q$-series.} Consider
  $$ \sum_{i=0}^\infty \frac{z^i}{q^{\binom{i}{2}}\,i!} 
 $$

The series converges if $|q|>1$, or $q=-1$. The
limit~\eqref{eqn:find-rho} in either case is $0$:
$$ \limsup_{n\rightarrow\infty} \frac{n \log n}{\log (n!\, q^{\binom{n}{2}})} = 
\limsup_{n\rightarrow\infty} \frac{n \log n}{n \log n + \frac{n(n-1)}{2}\log q} = 0
$$
and so $p=\rho=0$ and again the representation \eqref{eqn:desired-rep}
holds. If $q>1$ all terms are positive, and if $q=-1$ there are
positive and negative signs.
\end{example}

\begin{example} \label{ex:qfact} \textbf{A $q$-factorial.} 
 Consider the two series for  $q>1$.
$$
 \sum_{i=0}^\infty \frac{z^i}{(q;q)_i\, i!}    \quad\quad
 \text{ and }\quad\quad
 \sum_{i=0}^\infty \frac{z^i}{[i]!\, i!}  
$$
The limit for the first series is 
$$ \limsup_{n\rightarrow\infty} \frac{n \log n}{n \log n+ \sum_{i=0}^n
  \log(1-q^i)} =
\limsup_{n\rightarrow\infty} \frac{n \log n}{n \log n+
  \frac{n(n-1)}{2}\log q} = 0.
$$
The terms alternate in sign the first series.  If $q<-1$ then all
terms of the second series are positive.

\end{example}

  \begin{example}
Suppose $r>0$. If we compute \eqref{eqn:find-rho} for 
$$
(1-x)^{-r} = \sum_{n=0}^\infty \binom{r+n-1}{r} x^n$$
we find that
$\rho$ is infinite, which isn't surprising since this is not an entire
function. However,
$$
\expoper{}(1-x)^{-r} = \sum_{n=0}^\infty \binom{r+n-1}{r}
\frac{x^n}{n!}$$
has $\rho=1$.  A result of \Polya-Schur
(Lemma~\ref{lem:f-pos}) says that if $f\in\allpolyf$, and all
coefficients are non-negative then $f\in\allpolyposf$. If we knew that
$\expoper{} (1-x)^{-r}\in\allpolyf$ then we could conclude that
$\expoper{}(1-x)^{-r} \in\allpolyposf$.
  \end{example}

\begin{remark}
Consider $f = \displaystyle \sum \falling{\alpha}{i} {x^i}$. Since
$\expoper{}f = (1+x)^\alpha$ we see that
$\expoper{}f\not\in\allpolyf$.  However, when $\alpha=1/2$ we  have
the identity that $\expoper{}^2 f = e^{x/2} I_0(x/2)$ where $I_0$ is
the modified Bessel function of the first kind. Since
$$ I_0(x) = \sum _{i=0}^\infty \frac{x^i}{i!\,i!} = \expoper{}\,e^x$$
  we have that $\expoper{}^2\, f\in\allpolyf$.
  \index{Bessel function}
\end{remark}

\section{Polynomials with rapidly decreasing coefficients}
  \label{sec:rapid}

\begin{figure}[htbp]
  \centering
  \includegraphics*[width=3in]{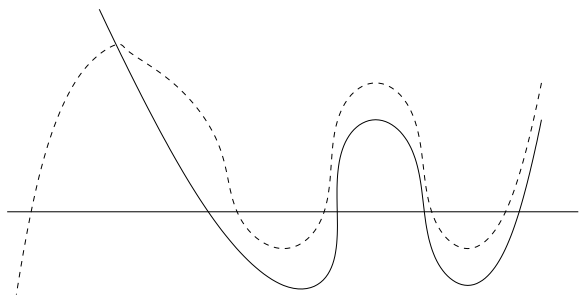}
  \caption{A polynomial with rapidly decreasing coefficients}
  \label{fig:rapid}
\end{figure}

If the coefficients of a polynomial are decreasing sufficiently rapidly
then the polynomial has all real roots. To motivate this, suppose that
$f\in\allpoly(n)$, and define $g(x) = f + \alpha x^{n+1}$. If $\alpha$
is very small then the graph of $g$ looks like Figure~\ref{fig:rapid}.
There are $n$ roots of $g$ that are close to the roots of $f$, so $g$
has all real roots.

\index{Newton's inequalities} 
\index{Newton's inequalities!converse}

The following theorem of Kurtz \cite{kurtz} gives a sufficient condition for
the existence of all real roots.

\begin{theorem} \label{thm:kurtz}
  Let $f=\sum a_ix^i$ be a polynomial of degree $n>1$ with all
  positive coefficients. If 
\begin{equation} \label{eqn:kurtz}
  a_i^2 - 4a_{i-1}a_{i+1}>0 \quad\quad  i=1,2,\dots,n-1
\end{equation}
then all roots of $f$ are real and distinct.
\end{theorem}
\index{Newton's inequalities!converse}

If we have equality in \eqref{eqn:kurtz} and $a_0=1$ then the
polynomial is
\[ f(x) = \sum_{k=0}^n x^k 4^{-\binom{n}{2}} \]
and the coefficients are rapidly decreasing. 
\begin{cor}\label{cor:kurtz-4}
  If $\alpha\ge4$ then the series below is in $\allpolyf$ and all of
  its partial sums are in $\allpoly$.
  \[\sum_{k=0}^\infty x^k
\alpha^{-\binom{n}{2}}\]
\end{cor}
  
\begin{cor}
  If $\alpha\ge2$ then then the series below is in $\allpolyf$ and all of
  its partial sums are in $\allpoly$.
  \[\sum_{k=0}^\infty x^k \alpha^{-k^2}\]
\end{cor}
\begin{proof}
Use the fact that $\alpha\ge2$ to check that \[
\left({a^{-k^2}}\right)^2 \ge 4
{a^{-(k-1)^2}}{a^{-(k+1)^2}} \]
\end{proof}

\begin{cor}
  If $q\ge4$ then then the series below is in $\allpolyf$ and all of
  its partial sums are in $\allpoly$.
  \[\sum_{k=0}^\infty \frac{x^k}{[k]!} \]
\end{cor}
\begin{proof}
  The \index{Newton quotient}Newton quotient is 
\[
 \frac{[n+2]!\,[n]!}{([n+1]!)^2} = \frac{[n+2]}{[n+1]} =
\frac{q^{n+2}-1}{q^{n+1}-1} > q \ge 4.
\]
\end{proof}

  Corollary~\ref{cor:kurtz-4} is in \cite{polya-szego2}*{V 176}.
  In \cite{polya-szego1}*{III 200} they  show that there is exactly one
  root in each annulus  $\alpha^{2n-2} < |z| < \alpha^{2n}$ for
  $n=1,2,\dots$. 
  There is a similar result due to Kuperberg \cite{kuperberg}. 
  \begin{lemma}
    If $q\ge4$ then the series below is in $\allpolyf$, and all
   partial sums are in $\allpoly$.
\[
\sum_{k=0}^\infty \frac{x^k}{(q;q)} = 
\sum_{k=0}^\infty \frac{x^k}{(1-q)(1-q^2)\cdots(1-q^k)}
\]
  \end{lemma}
\begin{proof}
  Kuperberg proves this without reference to Kurtz's result. He shows
  that if $F_{q,n}$ is the sum of the first $n$ terms of the series
  and $q\ge4$ then $F_{q,n}(q^{k+\frac{1}{2}})$ alternates in sign as $k$
  ranges from $0$ to $n$.
  \end{proof}

If all the partial sums of an infinite series have all real roots,
then the coefficients satisfy a quadratic exponential bound, but the
constant is $\sqrt{2}$, not $4$ as we would like.

  \begin{lemma}\label{lem:rapid-bound}
    If $\sum a_i x^i$ has all positive terms, and all partial sums
    have all real roots then
\begin{equation}\label{eqn:rapid-bound-1}
 a_n \le a_0\,\left(\frac{a_1}{a_0}\right)^n \,
 \left({2^{\binom{n}{2}}n!}\right)^{-1} 
\end{equation}
\end{lemma}
  \begin{proof}
    Since $a_0+a_1x+\cdots+ a_nx^n$ has all real roots we can apply
    Newton's inequalities \eqref{eqn:newton} 
    \begin{align*}
      \frac{a_{n-1}^2}{a_na_{n-2}} & \ge 2\frac{n}{n-1} \\
\frac{a_{n-1}}{a_n} &\ge 2 \,\frac{n}{n-1} \frac{a_{n-2}}{a_{n-1}} 
\ge 2 \,\frac{n}{n-1}\,2 \,\frac{n-1}{n-2}\, \frac{a_{n-3}}{a_{n-2}}
\ge \cdots \\
\frac{a_{n-1}}{a_n} &\ge2^{n-1} n \, \frac{a_0}{a_1} \\
a_n & \le \left(\frac{a_1}{a_0}\right) 2^{-n-1}\frac{1}{n} a_{n-1} \le 
\left(\frac{a_1}{a_0}\right)^2 2^{-n-1-(n-2)}\frac{1}{n(n-1)} a_{n-2}
\le \cdots \\
\intertext{from which \eqref{eqn:rapid-bound-1} follows.}
    \end{align*}
  \end{proof}

We end the section with a result where the terms are increasing.

\index{sign interlacing}
\index{rising factorial} \index{falling factorial}

\begin{lemma} \label{lem:aixi}
  Suppose $0\le a_0\le a_1 \le a_2\cdots\le\cdots.$ Then the polynomials 
\[
  f_n=\sum_{i=0}^n a_i\,(-1)^i \,\falling{x}{i}\quad\text{ and }\quad
g_n=\sum_{i=0}^n a_i \rising{x}{i}
\]
have all real roots. In addition, $f_{n+1}\lesslesseq f_n$ and
$g_{n+1}\lesslesseq g_n$. 
\end{lemma}
\begin{proof}
  The second case follows from the first upon substituting $-x$ for
  $x$ in the first case.  Note that if $k$ is an integer satisfying
  $1\le k\le n$ then $f_n(k) = \sum_0^{k-1} a_i (-1)^i \falling{k}{i} .$ As a
  function of $i$ the sequence $\{a_i\}$ is non-decreasing, and
  $\{\falling{k}{i} \}$ is increasing, so the sign of $f_n(k)$ is the sign of
  the largest term, and hence is $(-1)^{k-1}.$ Consequently, $f_n$
  sign interlaces $\falling{x}{n}$, and so in in $\allpoly.$ Since $f_{n+1} =
  f_n + a_{n+1}(-1)^{n+1} \falling{x}{n+1}$ it follows that
  $f_{n+1}\lesslesseq f_n$.
\end{proof}

\section{Interlacing of analytic functions}
\label{sec:analytic-interlace}

We use Proposition~\ref{prop:pattern} to motivate our definition of interlacing for
analytic functions in $\allpolyf$. 

\begin{definition}
  If $f,g\in\allpolyf$ then we say that $f$ and $g$ \emph{interlace} if
  $f+\alpha g$ is in $\allpolyf$ for all real $\alpha$.
\end{definition}

Later in this section we will relate interlacing of analytic functions
to the interlacing properties of their zeros.  Although it might be
hard to verify directly that a linear combination is in $\allpolyf$ we
can establish interlacing using limits.

\begin{lemma} \label{lem:anal-int}
  Suppose that $f,g\in\allpolyf$ and we have sequences of polynomials
  $(f_n),(g_n)$ where the degree of $f_n$ and $g_n$ is $n$ such that 
  $f_n\longrightarrow f$ and $g_n\longrightarrow g$. If $f_n$ and
  $g_n$ interlace for every $n$ then $f$ and $g$ interlace.
\end{lemma}
\begin{proof}
  Since $f_n$ and $g_n$ interlace we know that $f_n+\alpha g_n$ is in
  $\allpoly$ for all $\alpha$. Taking limits, we see that $f+\alpha g$
  is in $\allpolyf$ and hence $f,g$ interlace.
\end{proof}

Just as for polynomials, the derivative is interlacing.

\begin{lemma}
  If $f\in\allpolyf$ then $f^\prime\in\allpolyf,$ and $f$ and
  $f^\prime$ interlace. 
\end{lemma}
\begin{proof}
  If $f_n$ is a sequence of polynomials that converge to $f$, then the
  derivatives $f_n^\prime$ converge to $f^\prime$. Since $f_n+\alpha
  f_n^\prime$ converges to $f+\alpha f^\prime$ it follows that $f$ and
  $f^\prime$ interlace.
\end{proof}

\begin{remark}
  Exponential functions $e^{\alpha x}$ and $e^{\beta x}$ only
  interlace if $\alpha=\beta$. Indeed, if they did interlace, then
  $e^{\alpha x} + e^{\beta x}\in\allpolyf$. However,
\[ e^{\alpha x} + e^{\beta x} = e^{\alpha x}(1 + e^{(\beta-\alpha)
  x})\]
 has complex roots if $\alpha\ne\beta$.  The next three results
 generalize this simple observation.
\end{remark}

\index{interlacing!of exponentials}

We first see that no linear combination of polynomials and
exponentials is in $\rupint{2}$. 

\begin{lemma}
  If $f(x),g(x)$ are polynomials and $g(x)e^x+f(x)\in\allpolyf$ then
  $f=0$ and $g\in\allpoly$. 
\end{lemma}
\begin{proof}
  We first show that if $c\in\reals$ then $e^x+c\not\in\allpoly$. If
  it were then $(e^\diffd+c)(x^n) = (x+1)^n + c x^n\in\allpoly$, but
  the roots are complex for $n\ge3$.

  Next, if $e^x+f(x)\in\allpolyf$ where $f$ has degree $d$ then the
  $d$'th derivative is in $\allpolyf$ and is of the form $e^x+c$. This
  contradicts the first paragraph.

  If $g(x)e^x+c\in\allpolyf$ then replacing $x$ by $-x$ and
  multiplying by $e^x$ shows that $g(-x) + c\,e^x\in\allpolyf$, which
  we know is not possible. 

  Finally, if $g(x)e^x+f(x)\in\allpolyf$ then the $d$'th derivative is of
  the form $G(x)e^x + c$, so $f$ must be zero. Multiplying by $e^{-x}$
  shows that $g\in\allpoly$.
\end{proof}

Next, it is easy to describe the polynomials interlacing $e^x$.

  \begin{lemma}
    If $f\lace e^x$ then $f = (\alpha\, x+\beta)e^x$.
  \end{lemma}
  \begin{proof}
    The hypothesis means that $f(x) + \alpha e^x\in\allpolyf$ for all
    real $\alpha$. Thus \\ $ (f+\alpha e^x)(\diffd)x^n\in\allpoly$
    for all positive integers $n$. It follows that
    \[
 f(\diffd)x^n \lace e^\diffd x^n = (x+1)^n .
\] 
We can determine $f$ since we know all polynomials interlacing
$(x+1)^n$. Since $f(\diffd)x^n$ has degree at most $n$ it follows that
there are $a_n,b_n$ such that 
\[  f(\diffd)x^n = a_n(x+1)^n + b_n(x+1)^{n-1}.
\]
If we write  $f(x) = \sum d_ix^i$ then equating terms yields
\[
d_i\diffd^ix^n = d_i x^{n-i}\frac{n!}{(n-i)!} = x^{n-i}\biggl[
a_n\binom{n}{n-i} + b_n \binom{n-1}{n-i}\biggr]
\]
and so 
\begin{align*}d_i &= \frac{a_n}{i!} + \frac{b_n}{n}\frac{1}{(i-1)!}.
\\
\intertext{Setting $i=0$ shows that $a_n=d_0$, and $i=1$ yields $b_n/n =
d_1-d_0$. Thus}
  d_i &= d_0\frac{1}{i!} + (d_1-d_0)\frac{1}{(i-1)!}\\
\intertext{and consequently}
f(x) &= d_0 e^x + (d_1-d_0) \,x\,e^x.
\end{align*}

  \end{proof}

More generally we have

\begin{lemma}
  If $f(x) \lace g(x)e^x$ where $g(x)$ is a polynomial then there is
  a polynomial $h\lace g$ such that $f(x) = h(x)e^x$. 
\end{lemma}
\begin{proof}
  We first show that it suffices to show that $f(x)=H(x)e^x$ for some
  polynomial. If so, then since $f+\alpha ge^x\in\allpolyf$ it follows
  that $(H+\alpha g)e^x\in\allpolyf$, and therefore $g+\alpha
  H\in\allpoly$. Thus, $g$ and $H$ interlace.

  Assume that $g$ has degree $r$, and write $f = \sum d_ix^i$. If
  $n>r+1$ then $g(\diffd)e^\diffd x^n = g(\diffd)(x+1)^n$ has a factor
  of $(x+1)^2$. Since $f(\diffd)x^n\lace g(\diffd)(x+1)^n$ it follows
  that $-1$ is a root of $f(\diffd)x^n$. This gives us the equations,
  valid for $n>r+1$:
\[
\biggl(\sum d_i\diffd^i x^n\biggr)(-1) = n!\sum d_i
\frac{(-1)^{n-i}}{(n-i)!} =0
\]
that we view as an infinite matrix equation. For example, if $r=1$ we
have
\begin{equation}\label{eqn:flacex}
\begin{pmatrix}
  -1/6 & 1/2 & -1 & 1 & 0 & 0 & \dots \\
  1/24 & -1/6 & 1/2 & -1 & 1 & 0 & \dots \\
  -1/120& 1/24 & -1/6 & 1/2 & -1 & 1 &  \dots \\
  \vdots &&&&&&\ddots
\end{pmatrix}
\begin{pmatrix}
  d_0 \\d_1 \\d_2 \\ \vdots
\end{pmatrix}=
\begin{pmatrix}
  0 \\ 0 \\ 0 \\ \vdots
\end{pmatrix}
\end{equation}
In this example it is clear that we can solve for $d_3,d_4,\dots$ in
terms of $d_0,d_1,d_2$. If $M$ is the lower triangular matrix whose
$ij$ entry is  $(-1)^{i-j}/(i-j)!$ then the matrix above if formed by
removing the first three rows of $M$. Now the inverse of $M$ has entries
$1/(i-j)!$, and if we multiply \eqref{eqn:flacex} by $M^{-1}$ we
recove the identity matrix after $r+1$ columns. In the example, the
result of multiplying by $M^{-1}$ is
\[
\begin{pmatrix}
-1/6 & 1/2 & -1 & 1 & 0 & 0 & 0 & \dots \\
-1/8 & 1/3 & -1/2 & 0 & 1 & 0 & 0 & \dots\\
-1/20 & 1/8 & -1/6 & 0 & 0 & 1 & 0 & \dots\\
-1/72 & 1/30 & -1/24 & 0 & 0 & 0 & 1 & \dots\\
\vdots && \vdots &&\vdots &\ddots
\end{pmatrix}
\begin{pmatrix}
  d_0 \\d_1 \\d_2 \\d_3 \\ \vdots
\end{pmatrix}=
\begin{pmatrix}
  0 \\ 0 \\ 0 \\ 0 \\ \vdots
\end{pmatrix}
\]
Consequently we find
\[
\begin{array}{rrccrccrc}
  d_3 =& 1& d_2 &-&\frac{1}{2}& d_1 &+& \frac{1}{6}& d_0 \\
  d_4 =& \frac{1}{2}& d_2 &-&\frac{1}{3}& d_1 &+&  \frac{1}{8}&d_0 \\
\dots \\
  d_n =& \frac{1}{(n-2)!}&d_2 &-&  \frac{n-2}{(n-1)!}&d_1 &+&  \frac{(n-1)(n-2)/2}{n!}&d_0
\end{array}
\]
Adding up the series determines $f$:
\begin{multline*}
  d_0 + d_1 x + d_2 x^2 + d_2 \sum_{n=3}^\infty \frac{x^n}{(n-2)!}
- d_1 \sum_{n=3}^\infty  x^n\frac{n-2}{(n-1)!}
+ d_0 \sum_{n=3}^\infty x^n\frac{\binom{n-1}{2}}{n!} \\
= e^x\biggl( d_2x^2 + d_1(x-x^2) + d_0(1-x+\frac{1}{2}x^2)\,\biggr)
\end{multline*}
which shows that $f$ is $e^x$ times a quadratic.

In the general case we follow the same approach, and find that $f$ equals
\[
e^x\,\sum_{i=0}^{r+1} d_i \, x^i \biggl(\, \sum_{k=0}^{r+1-i} \,(-1)^k \frac{x^k}{k!}\biggr)
\]
which is a polynomial in $x$ of degree $r+1$ times an exponential, as desired.

\end{proof}

We now define interlacing of zeros:
\begin{definition}
  If $f$ and $g$ are in $\allpolyf$ then we say that \emph{the zeros
    of $f$ and $g$ interlace} if there is a zero of $f$ between every
  pair of zeros of $g$, and a zero of $g$ between every two zeros of
  $f$. If the set of zeros of $f$ and $g$ have least zeros $z_f$ and
  $z_g$ then we say that $f\lesslesseq g$ if $f$ and $g$ interlace,
  and $z_f \le z_g$. 
\end{definition}


If the zeros of functions in $\allpolyposf$ interlace and the
exponential factors are the same, then the
functions interlace. We would like to prove the converse. 

\begin{lemma}
  Suppose that 
  $$
  f = z^re^{\gamma  x}\prod_{i=1}^\infty (1+c_i x) \quad\quad g =
  x^re^{\gamma x}\prod_{i=1}^\infty (1+d_i x) $$
  are in $\allpolyposf$ where $c_1 > c_2
  > \dots > 0$ and $d_1 > d_2 > \dots > 0$. Then $f$ and $g$ 
  interlace if 
$$ c_1 > d_1 > c_2 > d_2 \dots \quad\text{ or } \quad
d_1> c_1 > d_2 > c_2 \dots$$
\end{lemma}
\begin{proof}
  Since  the $c_i$'s and the $d_i$'s interlace then the partial products
  $$
  f_n=\prod_{i=1}^n (1+ c_i x) \quad\quad g_n=\prod_{i=1}^n (1+ d_i
  x)$$
  interlace, and converge to interlacing functions. If we then
  multiply these interlacing functions by $x^re^{\gamma x}$ we find
  that $f$ and $g$ interlace.
\end{proof}

We can generalize part of Lemma~\ref{lem:sign-quant}. 
\begin{lemma}
  Suppose that $f(x)=\prod(1-x/a_i)$ is in $\allpolyposf$. Choose
  positive $\alpha_i$ such that $\sum \alpha_i<\infty.$ Then 
  \begin{enumerate}
  \item $g=\displaystyle\sum \alpha_i \frac{f}{1-x/a_i} \in\allpolyposf$
  \item $f$ and $\displaystyle\sum \alpha_i \frac{f}{1-x/a_i}$ interlace.
  \end{enumerate}
\end{lemma}
\begin{proof}
  Consider $$f_n = \prod_{i=1}^n (1-x/a_i)\quad\text{and}\quad g_n =
  \sum_{i=1}^n \alpha_i \frac{f}{1-x/a_i}$$
  We know that $f_n  \lesslesseq g_n$. The assumption on $\alpha_i$ implies that $g_n$
  converges uniformly to $g$. The result now follows by taking limits.
\end{proof}

\section{Characterizing $\allpolyf$ }

The following theorem shows that the
functions in $\allpolyf$ are exactly the right functions to generalize
polynomials.

\begin{theorem} \label{thm:polya-schur-2}
  An analytic function $f$ satisfies
  $f(\polar{})\allpoly\subset\allpoly$ if and only if $f$ is in
  $\allpolyf$.  It satisfies
  $f(\diffd{})\allpolypos\subset\allpolyneg$ if and only if it is in
  $\allpolyposf$.
\end{theorem}
\begin{proof}
  If $f$ is in $\allpolyf$, then it is the uniform limit of
  polynomials. The containment is true for each of the polynomials by
  Corollary~\ref{prop:fofd}, so it is true for $f$.  The converse is more
  interesting, and follows an argument due to \Polya\ and
  Schur\cite{polya-schur}.  \label{bib:polya-schur} First write
  $f(x) = \sum a_ix^i$, and set $q_n = f(\diffd{})x^n$.  Expanding
  the series
  \begin{align*}
    q_n & = \sum_{i=0}^n a_i \falling{n}{i}  x^{n-i} \\
    \intertext{and reversing $q_n$ gives a polynomial in $\allpoly$}
    p_n & = \sum_{i=0}^n a_i \falling{n}{i}  x^{i} \\
    \intertext{If we now replace $x$ by $x/n$ we find } p_n(x/n) & =
    \sum_{i=0}^n a_i x^{n-i} \prod_{k=1}^{i}(1-\frac{k}{n})
  \end{align*}
  The polynomials $p_n(x/n)$ converge uniformly to $f$, and so $f$ is
  in $\allpolyf$. The second case is similar.
\end{proof}

  \index{Laguerre inequality}

\section{\allpolypsd\ matrices} \changed{1/18/7}
  We derive some simple inequalities for quadratic forms involving the
  coefficients of a polynomial with all real roots. It is interesting
  that these inequalities are valid for a broader class of functions
  than just $\allpolyf$.

  Define an $n$ by $n$ matrix $Q$ to be \emph{\allpolypsd}\ if for every
  $f=\sum_0^n a_ix^i\in\allpoly$ and $A = (a_0,\dots,a_{n-1})$, we have $A
  Q A^t \ge 0$.  Of course, any positive semi-definite matrix is also
  \allpolypsd, but the converse is not true. In particular,
  \allpolypsd\ matrices can have negative entries. Indeed, Newton's
  inequalities imply that the matrix $ \left(\begin{smallmatrix}
      0&0&-2\\0&1&0\\-2&0&0
  \end{smallmatrix}\right)
$ is \allpolypsd.

We first note that it is easy to construct new
  \allpolypsd\ matrices from old ones. It's obvious that \allpolypsd\
  matrices form a cone. In addition, 

\begin{lemma}
  If $Q$ is a \allpolypsd\ matrix, and $B$ is the
  diagonal matrix of coefficients of a polynomial in $\allpolypos$,
  then $BQB$ is also \allpolypsd.
\end{lemma}
\begin{proof}
  If $A$ is the vector of coefficients of $f\in\allpoly$, then $AB$ is
  the vector of coefficients of the Hadamard product of $f$ and $g$.
  By Theorem~\ref{thm:hadamard-1} the Hadamard product is in
  $\allpoly$, it follows that $ A(BQB)A^t = (AB)Q(AB)^t \ge 0$.
\end{proof}

Here is a  construction of a \allpolypsd\ matrix. Note that it applies
to polynomials in $\allpoly$. 

  \begin{prop}\label{prop:gen-laguerre}
    Suppose that $$f(y) = \displaystyle\sum_{i=0}^\infty a_i {y^i} =
    e^{-g(y)+h(y)}\prod_j\left(1-\dfrac{y}{r_j}\right)$$
    where $h(y)$ is
    a polynomial with only terms of odd degree, $g(y)$ has all
    non-negative coefficients, and all exponents of $g$ are $\equiv 2
    \pmod{4}$.  If $m$ is a positive integer then
$$\sum_{j=-m}^{m}
(-1)^{m+j}{a_{m-j}}{a_{m+j}} \ge 0
$$
Since all real rooted polynomials satisfy the hypothesis, the $2m+1$
by $2m+1$ matrix whose anti-diagonal is alternating $\pm1$ beginning
and ending with $1$ is
\allpolypsd.
  \end{prop}
  \begin{proof}
    The coefficient of $y^{2m}$ in
$$ f(y)f(-y) = \left(\sum a_{s}{y^s}\right)\left(\sum
  a_{r}(-1)^r{y^r}\right)$$ 
is 
$$
\sum_{j=-m}^{m} (-1)^{j}{a_{m-j}}{a_{m+j}} 
$$
On the other hand, since $h(y)+h(-y)=0$ and $g(y)=g(-y)$ by hypothesis, we have 
\begin{align}
  f(y)f(-y) &= \left(e^{-g(y)+h(y)}\prod
    \left(1-\dfrac{y}{r_i}\right) \right)\, 
\left(e^{-g(-y)+h(-y)}\prod \left(1+\dfrac{y}{r_i}\right)\right) \notag\\
&= e^{-2g(y)}\,\prod \left( 1 - \frac{y^2}{r_i^2}\right) \label{eqn:lag-ineq}
\end{align}

Since $r_i^2$ is positive, the sign of the coefficient of $y^{2m}$ in
the infinite product is $(-1)^m$. The coefficient of $y^{2m}$ in
$e^{-2g(y)}$ is either $(-1)^m$ or $0$.  Consequently, the sign of the coefficient of
$y^{2m}$ in \eqref{eqn:lag-ineq} is also $(-1)^m$ or $0$, which finishes the
proof.
  \end{proof}

The case $m=1$ in the following corollary is the Laguerre inequality. 

  \begin{cor}
    Suppose that $f\in\allpoly$, and $m$ is a positive integer. Then
$$\sum_{j=-m}^{m}
(-1)^{m+j}\frac{f^{(m-j)}(x)}{(m-j)!}\frac{f^{(m+j)}(x)}{(m+j)!} \ge 0
$$
More generally, if $Q=(q_{ij})$ is a \allpolypsd\ matrix, then
$$ \sum_{i,j} q_{ij} \frac{f^{(i)}}{i!} \frac{f^{(j)}}{j!} \ge 0 $$
  \end{cor}
  \begin{proof}
    Use the Taylor series
$$ f(x+y) = \sum_{s=0}^\infty f^{(s)}(x)\frac{y^s}{s!}$$
  \end{proof}

  \begin{cor}
    If $P_n(x)$ is the Legendre polynomial,  $m$ is a positive
    integer, and $|x|<1$ then
$$\sum_{j=-m}^{m}
(-1)^{m+j}\frac{P_{m-j}(x)}{(m-j)!}\frac{P_{m+j}(x)}{(m+j)!} \ge 0
$$
  \end{cor}
  \begin{proof}
    The following  identity is valid for $|x|<1$. 
$$
  \sum_{n=0}^\infty \frac{P_n(x)}{n!}\,y^n  = e^{xy} J_0(y\sqrt{1-x^2})
$$
Now apply Proposition~\ref{prop:gen-laguerre}.
  \end{proof}

  \begin{remark}
    Krasikov\cite{krasikov} essentially shows that the matrix below (derived from
    his $V_4(f)$) is \allpolypsd.
$$    \begin{pmatrix}
      0 & 0 & 0 & 0 & -2 \\ 0 & 0 & 0 & -1 & 0 \\ 0 & 0 & 2 & 0 & 0 \\
      0 & -1 & 0 & 0 & 0 \\ -2 & 0 & 0 & 0 & 0
    \end{pmatrix}
$$
To see this directly,  notice we can write it as a sum
of \allpolypsd\ matrices.
\[
\begin{pmatrix}
      0 & 0 & 0 & 0 & 1 \\ 0 & 0 & 0 & -1 & 0 \\ 0 & 0 & 1 & 0 & 0 \\
      0 & -1 & 0 & 0 & 0 \\ 1 & 0 & 0 & 0 & 0
\end{pmatrix}
+
\begin{pmatrix}
      0 & 0 & 0 & 0 & -3 \\ 0 & 0 & 0 & 0 & 0 \\ 0 & 0 & 1 & 0 & 0 \\
      0 & 0 & 0 & 0 & 0 \\ -3 & 0 & 0 & 0 & 0
\end{pmatrix}
\]

The first one is \allpolypsd\ by the proposition. From equation
\eqref{eqn:newton-3} we have $a_2^2 \ge 6 a_0a_4$ which shows the
second matrix is \allpolypsd. More generally, if $0<k<n$ then Newton's
inequalities \index{Newton's inequalities} imply that the matrix with
$1$ in the $(k,k)$ entry, and
$-\frac{1}{2}\,\frac{k+1}{k}\,\frac{n-k+1}{n-k}$ in the $(k+1,k-1)$
and $(k-1,k+1)$ entries is a \allpolypsd\ matrix.
  \end{remark}


\chapter{Linear Transformations of Polynomials}
\label{cha:linear}

\renewcommand{\TimeStampStart}{Monday, October 15, 2007: 12:12:43}
\mytoday  

The aim of this chapter is to establish properties of linear
transformations that preserve roots.  In
Section~\ref{sec:characterize} we show that if we put assumptions on
$T$ such as $Tf \lesslesseq f$ for all $f\in\allpoly$ then $Tf$ is the
derivative. If $T$ is bijection on $\allpoly$, then $T$ is affine. We
see how to determine the possible domain and range of a linear
transformation.  In the next section we study transformations of the
form $f\mapsto f(T)(1)$ where $T$ is a linear transformation. The
following two sections consider properties of linear transformations
that are defined by recurrence relations. We then consider the effect
of \Mobius\ transformations on the the roots of polynomials and on
transformations. In the final section we begin the study of
transformations of the form $x^n \mapsto \prod_1^n (x+a_i)$.

\section{Characterizing transformations}
\label{sec:characterize}

Many linear transformations can be characterized by their mapping and
interlacing properties. 
The derivative can be characterized by its fundamental interlacing
property $f \lesslesseq f^\prime$. This result relies on the fact that
the only monic polynomial $f$ such that $(x+a)^n\lesslesseq f$ is
$f=(x+a)^{n-1}$.

\begin{theorem} \label{thm:fd}
  If $T$ is a linear transformation on polynomials such that $f
  \lesslesseq T(f)$ for all polynomials $f$ in $\allpolyalt$ (or for all
  polynomials in $\allpolypos$), then $T$ is a
  multiple of the derivative.
\end{theorem}
\begin{proof}
  Choosing $f = x^m$, we find $x^m \lesslesseq T(x^m)$, and hence
  there is a $c_m$ so that $T(x^m) = c_m x^{m-1}$. From
  $(x+1)^m\lesslesseq T(x+1)^m$ we see that there an $\alpha_m$ such
  that
\[  
\sum_{i=1}^{m}  \binom{m}{i}\, c_i\,x^{i-1} = 
T(x+1)^m =  
\alpha_m\,(x+1)^{m-1} =
\alpha_m\, \sum_{i=0}^{m-1}   \binom{m-1}{i}x^{i} 
\]
Equating coefficients shows that $c_{i} = (\alpha_m/m)\,i$,
 and so $T(f)$ is a multiple of $f'$.
\end{proof}

  
  


\begin{theorem} \label{thm:fTf}
  If $T$ is a linear transformation with the property that $Tf$ and
  $f$ have the same degree, and  interlace for all
  $f$ in $\allpolyalt$ (or for all $f$ in $\allpolypos$) then there are
  constants $a,b,c$ such that
$$ Tf = af + bxf^\prime + cf^\prime$$
If $Tf$ and $f$ interlace for all $f\in\allpoly$  then $Tf = a f + c f^\prime$.

\end{theorem}

\begin{proof}
By continuity $x^n$ and $ T(x^n)$ interlace, and so there are constants
$a_n,b_n$ such that $T(x^n) = a_nx^n + b_nx^{n-1}$.
Apply $T$ to 
the polynomial $f = x^r(x-\alpha)^2$.  We compute
\begin{align} \label{eqn:double-root}
 Tf &= T(x^{r+2} - 2\alpha x^{r+1} + \alpha^2x^r) \\
& = \notag 
a_{r+2}x^{r+2} + b_{r+2}x^{r+1} - 
2\alpha a_{r+1}x^{r+1} - \\
& \text{\hspace*{.2cm} } \notag
2\alpha b_{r+1}x^{r} +
\alpha^2 a_{r}x^{r} + \alpha^2 b_{r}x^{r-1} 
\end{align}

Now $f$ has a double root at $\alpha$, and  $Tf$ interlaces $f$, so $Tf$
must have $\alpha$ as a root.  Substituting $x=\alpha$ in
\eqref{eqn:double-root} gives the  recurrence relation
\begin{align*}
  0 & = \alpha^{r+2}(a_{r+2}-2a_{r+1} + a_r) + 
\alpha^{r+1}(b_{r+2}-2b_{r+1} + b_r) \\
\intertext{This equation holds for infinitely many  $\alpha$, so we get the two
  recurrences}
0 & = a_{r+2}-2a_{r+1} + a_r \\
0 & = b_{r+2}-2b_{r+1} + b_r 
\end{align*}
Since these recurrences hold for $1\le r \le n-2$, we solve these
equations and find that constants $a,b,c,d$ so that $a_r = a + br$,
$b_r = d + cr$. Substituting these values gives
$$  T(x^r)  = (a+br)x^r + (d+cr)x^{r-1} $$

Since $rx^r = x\diffd{}x^r$, it remains to show that $d=0$.  We know
that for any negative $\alpha$ the polynomials $(x+\alpha)^3$ and
$T(x+\alpha)^3$ interlace, and hence $(x+\alpha)^2$ divides
$T(x+\alpha)^3$. A computation shows that the remainder of
$T(x+\alpha)^3$ upon division by $(x+\alpha)^2$ is
$d\alpha(x+2\alpha)$. As this remainder is $0$ we conclude that $d=0$.
The case where $f\in\allpolypos$ is similar and omitted.

Finally, assume that $Tf$ and $f$ interlace for all $f\in\allpoly$. If
$b$ is not zero then we may divide $Tf$ by $b$, and thus assume that
$Tf = a f + x f^\prime + cf^\prime$. Choose $\alpha>|c|$ and consider
the polynomial
\begin{align*}
f & = (x-\alpha)(x+\alpha)(x-\alpha-1)(x+\alpha+1) \\
\intertext{A computation shows that}
(Tf)(-\alpha) &= 2(c-\alpha)\alpha(1+2\alpha)\\
(Tf)(\alpha) &= -2(c+\alpha)\alpha(1+2\alpha)
\end{align*}
The signs of the left hand side must alternate since $f$ and $Tf$ sign 
interlace, and so $(c-\alpha)(c+\alpha)>0$ which implies $\vert
c\vert>\alpha$. Since we chose $\alpha >|c|$ this is a
contradiction, and hence $b=0$.
\end{proof}

The remaining case is $Tf \lesslesseq f$.  Using an argument entirely
similar to that of Theorem~\ref{thm:fTf} we find

\begin{theorem} \label{thm:Tff}
  If $T$ is a linear transformation such that $Tf \lesslesseq f$ for all
  $f\in\allpoly$ then there are constants $a,b,c$ where $a$ and $c$
  have the same sign  such that
$$ Tf = axf + bf + cf^\prime$$
\end{theorem}

\begin{remark}
  If we restrict the domain then we can find  transformations that are
  not the derivative that satisfy $Tf \lesslesseq f$. Here are three
  examples:
  \begin{align*}
    T(f) & = (x^2-1)f' &\text{on}\ \allpolyint{(-1,1)}\\
    T(f) & = \diffd((x^2-1)f) &\text{on}\ \allpolyint{(-1,1)}\\
    T(x^n) &= \dfrac{x^{n+1}}{n+1}&  \text{on}\ \expoper{}(\allpoly)
  \end{align*}

\end{remark}

We can generalize Theorem~\ref{thm:fd} by restricting the domain of the
transformation to polynomials of degree $n$. Under these restrictions
the polar derivative makes an appearance.
\index{polar derivative}

\begin{theorem} \label{fd-2}
  Fix a positive integer $n$ and suppose that $T$ is a linear
  transformation that maps homogeneous polynomials of degree $n$ to
  homogeneous polynomials of degree $n-1$ with the property that
  $F\lesslesseq T(F)$ for all $f\in\allpoly(n)$ for which $F(x,1)=f(x)$.
  There are $b,c$ where $bc\ge0$ such that
$$ T(F) = b \frac{\partial F}{\partial x} + c \frac{\partial F}{\partial 
  y}$$
\end{theorem}
\begin{proof}
  Since $T$ is continuous we observe that if $F$ is the limit of
  polynomials of degree $n$ then $F$ and $TF$ interlace.  Considering
  that
$$ \lim_{\epsilon\rightarrow0} x^r(1+\epsilon x)^{n-r} = x^r$$
we conclude that $x^r$ and $T(x^r)$ interlace for all $r\le n$. Thus
there are $a_r,b_r,c_r$ such that
\begin{align*}
  T(x^r) &= x^{r-1}(a_r+b_rx+c_rx^2) \\
\intertext{Since $T(x^r(x-\alpha)^2)$ is interlaced by
  $x^r(x-\alpha)^2$, it has $\alpha$ for a root, and  we find}
0 =& \alpha^{r+2}(a_{r+2}+b_{r+2}\alpha+c_{r+2}\alpha^2) \\
&-2\alpha^{r+1}(a_{r+1}+b_{r+1}\alpha+c_{r+1}\alpha^2)\\
 &+\alpha^{r}(a_{r}+b_{r}\alpha+c_{r}\alpha^2) \\
\intertext{Since this is true for all $\alpha$ we find that}
0 &= a_{r+2} - 2a_{r+1} + a_r \\
0 &= b_{r+2} - 2b_{r+1} + b_r \\
0 &= c_{r+2} - 2c_{r+1} + c_r \\
\intertext{Solving these recurrences yields $a_r=A_0+A_1r$,
  $b_r=B_0+B_1r$, $c_r=C_0+C_1r$ which  implies}
T(x^r) &= (A_0+A_1r)x^{r-1} + (B_0+B_1r)x^r + (C_0+C_1r)x^{r+1}\\
\intertext{Since $T(x^n)$ is a polynomial of degree $n-1$ we find
  $B_0+B_1n=0$ and $C_0+C_1n=0$ so that }
T(x^r) &= (A_0+A_1r)x^{r-1} + B_1(r-n)x^r + C_1(r-n)x^{r+1}\\
\intertext{Since}
 T(x^{n-1}(x+1)) &= (A_0+A_1n)(x^{n-1}) +(A_0+A_1(n-1))x^{n-2}  - \\
&\quad B_1x^{n-1} - C_1x^n\\
\intertext{has degree $n-1$ we see that $C_1=0$. It remains to see that
  $A_0=0$. If not, we may assume $A_0=1$ and consider}
T((x-\alpha)^3) &= (x^2-3x\alpha+3) + 3A_1(x-\alpha)^2 -
B_1(-3x^2\alpha+6x\alpha^2-3\alpha^3) 
\end{align*}
{whose discriminant is $3\alpha^2(-1-4A_1-4B_1\alpha)$.}
Since $\alpha$ is arbitrary this can be negative unless $B_1=0$. Thus
$T(x^r) = x^{r-1} +A_1rx^{r-1}$. We see that $TF$ has degree $r-1$ so
$F \lessless TF$. The proof of Theorem~\ref{thm:fd} shows that this can not
happen, so $A_0$ is $0$, and
$ T(x^r) = A_1 (rx^{r-1}) + B_1 (r-n)x^r $ which implies the
conclusion by linearity. 
\end{proof}

  If $f$ is a polynomial then $f(\diffd)(x+a)^n$ is a sum of constants
  times powers of $x+a$. The next lemma is a  converse.
  \begin{lemma}
    Suppose $T\colon\allpoly\longrightarrow\allpoly$ satisfies 
    \[ T(x+a)^k = \sum_{i=0}^k \alpha_{k,i}\,(x+a)^i \]
    for all $a\in\reals$, where  $\alpha_{k,i}$ are constants not
    depending on  $a$. Then
    there is an $f\in\allpolyf$ such that $T = f(\diffd)$.
  \end{lemma}
  \begin{proof}
 Continuity of $T$ implies that $\diffd$ and $T$
    commute:
    \begin{align*}
      T\diffd(x+a)^k &
=
      T\lim_{\epsilon\rightarrow0} \frac{(x+a+\epsilon)^k -
        (x+a)^k}{\epsilon} \\
      &= \lim_{\epsilon\rightarrow0} 
      \frac{\sum_i \alpha_{k,i}(x+a+\epsilon)^i -
        \sum_i\alpha_{k,i}(x+a)^i}{\epsilon}\\
      &= \sum \alpha_{k,i}(x+a)^{k-1} = DT(x+a)^k
    \end{align*}

Now, any linear operator commuting with differentiation is a power
series $f(x)$ in $\diffd$. Since $f(\diffd)$ maps $\allpoly$ to
itself, it follows from  Theorem~\ref{thm:polya-schur-2} that
$f\in\allpolyf$. 
  \end{proof}

  The proof of the following lemma (due to Carncier-Pinkus-Pe{\~n}a\cite{pinkus-pena})
  clearly shows the significance of the monic hypothesis. The proof
  reduces to the fact that we know all the polynomials interlacing
  $(x+1)^n$. 
  \index{diagonals!interlacing}

  \begin{lemma}\label{lem:cpp-2}
    If $T:\allpoly(n)\longrightarrow\allpoly(n)$ satisfies
    \begin{enumerate}
    \item $T(x^k) = x^k + $ terms of lower degree for $0\le k \le n$.
    \item $T(x^k)\in\allpolypos$ for $0\le k \le n$.
    \end{enumerate}
then there is a polynomial $f$ such that
\begin{enumerate}
\item T = f(\diffd)
\item $\expoper{}(f)\in\allpoly$
\end{enumerate}
  \end{lemma}
  \begin{proof}
    Let $T(x^k) = \sum_0^ka_{k,i}x^o$ where $a_{k,k}=1$. We first
    describe the argument for $n=3$. Consider $T_\ast(x+y)^3$:\\[.2cm]

\centerline{\xymatrix@=.2cm{
    y^3 \\
    3a_{02}y^2 & 3xy^2 \\
    4a_{01}y & 6a_{11} xy & 3x^2y \\
    a_{00} & a_{10} x & a_{20}x^2 & x^3
}}

The monic hypothesis implies that the rightmost diagonal is
$(x+y)^3$. Since $T_\ast(x+y)^3\in\gsubplus_2$ the diagonals
interlace, so the adjacent diagonal is a multiple of $(x+y)^2$, and so
on. Thus $T_\ast(x+y)^3=$\\[.2cm]

\centerline{\xymatrix@=.2cm{
    y^3 \\
    a_{20}y^2 & 3xy^2 \\
    a_{10}y & 2a_{20} xy & 3x^2y \\
    a_{00} & a_{10} x & a_{20}x^2 & x^3
}}

If $f = \sum b_i x^i$ then $f(\diffd)(x+y)^3=$\\[.2cm]

\centerline{\xymatrix@=.2cm{
y^3\,b_0 \\
3y^2\,b_1 & xy^2\,b_0 \\
6y\,b_2 & 6xy\,b_1 & x^2y \, b_0 \\
6\,b_3 & 6x\, b_2 & 3x^2\,b_1 & x^3\,b_0
}}

 Equating coefficients shows that if
 $g=a_{00}+a_{10}x+a_{20}x^2+a_{30}x^3$ then
\[
 T = \frac{1}{6} \bigl( \exp^{-1}\,g^{rev}\bigr)(\diffd)
\]

In general, the rightmost diagonal of $T_\ast(x+y)^n$ is $(x+y)^n$, so
all the parallel diagonals are powers of $x+y$. Thus
\[
T_\ast(x+y)^n = \sum_{i=0}^n a_{i0}(x+y)^i
\]
If $g = \sum_{i=0}^n a_{i0}x^i$ then
\begin{align*}
  \biggl(\frac{1}{n!} \exp^{-1} \,g^{rev}\biggr)(\diffd)(x+y)^n &=
\biggl(\sum_{i=0}^n \frac{(n-i)!}{n!} a_{n-i,0}\diffd^i\biggr)(x+y)^n \\
&= \sum_{i=0}^n a_{i0}(x+y)^i = T_\ast(x+y)^n
\end{align*}
  \end{proof}

Affine transformations are the only transformations
whose inverses also preserve roots.

\index{affine transformation}
\begin{theorem} If a degree preserving linear transformation and its
  inverse both preserve roots then the linear transformation is an
  affine transformation. More precisely, if $T$ is a linear
  transformation such that
  \begin{itemize}
  \item $T(\allpoly)\subset\allpoly$
  \item $T^{-1}(\allpoly)\subset\allpoly$
  \item $deg(T(x^n)) = n$
  \end{itemize}
  then there are constants $a,b,c$ such that $Tf(x) = cf(ax+b)$.
\end{theorem}
\begin{proof}
  Using the three constants, we may assume that $T(1)=1$ and $T(x) =
  x$.  With these normalizations, it suffices to show that $Tf=f$.
  Given a polynomial $f\in\allpoly$, we can use Theorem~\ref{thm:only-roots} to
  apply $T^{-1}$ to the interlacing $ Tf \lesslesseq (Tf)^\prime$ to
  find that $f\lesslesseq T^{-1}(Tf)^\prime$. Since
  $T^{-1}(Tf)^\prime$ is a linear transformation defined on all of
  $\allpoly$, there is a constant $w$ such that $T^{-1}(Tf)^\prime =
  wf^\prime$, or equivalently $wT(f^\prime) = (Tf)^\prime$.

  We now can show that $T(x^2) = x^2$.  Assume that $T(x^2) = \alpha
  x^2 + \beta x + \gamma$.  Since
  \begin{align*}
    T(x^2)^\prime &= 2\alpha x + \beta \\
    wT(2x) &= 2wx
  \end{align*}
  we find that $T(x^2) = wx^2+\gamma$.  Solving for $T^{-1}$ gives
  $T^{-1}(x^2) = (1/w)(x^2-\gamma)$.  If $\gamma$ is non-zero, then it
  is not possible that both $T(x^2)$ and $T^{-1}(x^2)$ have real
  roots, so $T(x^2) = wx^2$. Since $T((x-1)^2) = wx^2 - 2x + 1$ has
  all real roots, the discriminant is non-negative which implies that
  $w\ge 1$. Similarly, since $T^{-1}(x^2)$ has all real roots then
  $(1/w)\ge 1$, and so $w=1$.

  If we assume by induction that $T(x^{n-1}) = x^{n-1}$ then applying 
  $T(f^\prime) = (Tf)^\prime$ shows that $T(x^n) = x^n + \alpha$.
  This polynomial does not have all real roots unless $\alpha=0$, so
  $T(x^n) = x^n$ and thus $Tf=f$.
\end{proof}

\begin{remark}
  There is an alternate argument for the end of the proof. If $T$
  commutes with differentiation $\diffd$ then $T = f(\diffd)$ for some
  $f\in\allpolyf$. Similarly $T^{-1} = f^{-1}(\diffd)$, so $f$ is
  invertible and has no zeros.  The only functions of this form in
  $\allpolyf$ are $f(x) = ae^{bx}$, and $ae^{b\diffd} g(x) = a
  g(x+b)$.
\end{remark}

\begin{cor}
  If $T$ is an invertible linear transformation such that
  $T\allpoly=\allpoly$ then $Tf = cf(ax+b)$.
\end{cor}

\begin{remark}
  Although the only bijections on $\allpoly$ are trivial, there are
  non-trivial bijections on $\allpolyalt$. If
  $\tilde{L}_n^{(\alpha)}(x)$ is the monic Laguerre polynomial, %
\index{Laguerre polynomials!monic}\index{polynomials!Laguerre}%
(\chapsec{operators}{laguerre}) then \cite{roman} the linear
transformation $x^n \mapsto \tilde{L}_n^\alpha(x)$ satisfies
$T=T^{-1}$. Since $T$ maps $\allpolyalt$ to itself (Corollary~\ref{cor:laguerre-x})
it follows that $T$ is a bijection on $\allpolyalt$.

\end{remark}

  \begin{lemma}
    Suppose that $T$ is a degree preserving linear transformations
    from polynomials to polynomials. If $f$ and $Tf$ have a common
    root for all $f$ then $Tf$ is a multiple of $f$. 
  \end{lemma} 
  \begin{proof}
    We can assume that $T(1)=1$, and let $T(x^n) = \sum
    a_{n,i}x^i$. Since $T(x^n)$ and $x^n$ have zero for a common root the
    constant terms $a_{n,0}$ are zero for $n>0$. Next, $x+t$ and
    $T(x+t) = a_{1,1}x+t$ have $-t$ for a common root, so $a_{1,1}=1$
    and $T(x)=x$. 

    Now assume that $T(x^k)=x^k$ for $k<n$. 
\[
      T(x+t)^n = T(x^n) + \sum _{k=1}^{n-1} x^k\binom{n}{k}t^k 
      = T(x^n) + (x+t)^n - x^n
\]
Since $-t$ is the only root of $(x+t)^n$ we see that
$T(x^n)(-t)=(-t)^n$ for all $t$ and therefore $T(x^n) = x^n$.
  \end{proof}

If $f$ and $Tf$ have a common interlacing then $Tf$ has a simple form.

\begin{lemma}
 If $T$ is a degree preserving linear transformation such that $Tf$ and $f$
  have a common interlacing for all $f\in\allpoly$
   then  there are $\alpha,\beta,\gamma$ so that
\[
    T(f) = (\alpha x+\beta)f' + \gamma f
\]
\end{lemma}
\begin{proof}
  Since $x^{n-1}$ is the only polynomial that interlaces $x^n$, we see
  that if $n\ge2$ then $T(x^n)$ has a $n-1$ fold root at
  zero. Write
\begin{align*}
T(1) &= b_0 \\
T(x^n)  &= a_nx^n + b_n x^{n-1} \qquad n\ge1.
\end{align*}

Now we   use the fact that $T(x-t)^n$ must have a root at $t$ if
$t\ge0$. For $n=2$ this means that
\begin{align*}
  0 = T(x-t)^2[t] &= T(x^2)[t] - 2t T(x)[t] + t^2 T[1] \\
\intertext{which implies}
  T(x^2) &= 2xT(x) - x^2 T(1)
\end{align*}
Thus $a_2 = 2a_1-b_0$ and $b_2 = 2b_1$. We now follow the same
argument, and show by induction that
\begin{align*}
  a_n &= na_1 - (n-1)b_0 & b_n &= n b_{n-1}
\end{align*}
which implies that
\begin{align*}
  T(x^n) &= x^n(na_1 - (n-1)b_0) + nx^{n-1}b_1 \\
  &= a_1 x( x^n)' - b_0(x(x^n)'- x^n)  + b_1 (x^n)'\\
\intertext{and by linearity}
T(f) &= a_1\,xf' -b_0\, xf' + b_0\, f + b_1\, f' \\
&= ((a_1-b_0)x + b_1)\,f' + b_0\, f.
\end{align*}
\end{proof}

Compositions mapping $\allpoly$ to itself are very simple.

\begin{lemma}\label{lem:char-affine}
  Suppose $g(x)$ is a polynomial.
  \begin{enumerate}[(1)]
  \item If $f(g(x))\in\allpoly$ for all $f\in\allpoly$ then $g$ is linear.
  \item If $f(g(x))\in\allpoly$ for all $f\in\allpolypos$ then $g$ is
    either linear, or is quadratic with negative leading coefficient.
  \item If $g(x)+c\in\allpoly$ for all $c\in\reals$ then $g$ is linear.
  \item If $g(x)+c\in\allpoly$ for all $c>0$ then $g$ is
    quadratic with negative leading coefficient.
  \end{enumerate}
\end{lemma}
\begin{proof}
  If we take $f$ to be linear, then (1) and (2) follow from  (3)
  and (4). If the degree of $g$ is odd then the horizontal line $y=-c$
  will intersect the graph of $g(x)$ in only one point if $c$ is
  sufficiently large. Thus, in this case the degree of $g$ is $1$.

  If the degree of $g$ is even then the horizontal line $y=-c$ with
  $c$ large and positive will intersect the graph of $g$ in at most
  two points, and hence the degree of $g$ is two. Since these
  horizontal lines are below the $x$-axis, the parabola must open
  downward, and so the leading coefficient is negative.
\end{proof}

  Here is another characterization of affine transformations. 

  \begin{lemma}\label{lem:affine-char}
    Suppose that $T\colon{}\allpoly\longrightarrow\allpoly$ is a linear
    transformation. $T$ distributes over multiplication (\,$T(f\cdot g)=Tf\cdot Tg$\,)
    if and only if $T$ is an affine transformation.
  \end{lemma}
  \begin{proof}
    If $p=T(x)$ then by induction  $T(x^n)=p^n$. Since
    $T(x+c)=p+c\in\allpoly$ for all $c\in\reals$, we see that $p$ is
    linear, so write $p = qx+b$. If $f(x) = \sum a_ix^i$ then the
    conclusion follows from
$$ T(f) = \sum a_i T(x^i) = \sum a_i p^i = f(p) = f(qx+b).$$
  \end{proof}

The next lemma describes those transformations whose translations also
preserve roots.

\begin{lemma} \label{lem:axb}
  If $T$ is a linear transformation on polynomials such that 
  \begin{enumerate}
  \item $(T+a)\allpoly\subset\allpoly$ for all $a\in\reals$.
  \item $deg(T(x^n)) = n$
  \end{enumerate}
then there are constants $c,d$ such that  $Tf = cf + df^\prime$.
\end{lemma}
\begin{proof}
  For any polynomial $g\in\allpoly$, $(T+a)g = Tg + ag\in\allpoly$.
  By Proposition~\ref{prop:pattern} $Tg$ and $g$ interlace.  The conclusion now
  follows from  Theorem~\ref{thm:fTf}.
\end{proof}

A similar argument shows
\begin{lemma}
  If $T$ is a linear transformation on polynomials such that 
  \begin{enumerate}
  \item $(T+a)\allpolyalt \subset \allpoly$ for all $a\in\reals$
  \item $deg(T(x^n)) = n$
  \end{enumerate}
then  $Tf(x) = cf + (d+ex)f^\prime$ for
certain $c,d,e\in\reals$.
\end{lemma}

It is difficult to tell when two linear transformations $T,S$ have
the property that $Tf$ and $Sf$ interlace for all $f\in\allpoly$, but
it is easy if $T,S$ are  multiplier transformations.
\index{multiplier~transformations} Such transformations map $x^n$ to a
constant multiple of $x^n$.

\begin{lemma} \label{lin-multiplier}
  Suppose $T(x^n) = t_nx^n$ and $S(x^n) = s_nx^n$ both map $\allpoly$
  to itself.  Assume that $t_n$ is always positive and
  $\smalltwodet{t_n}{t_{n+1}}{t_{n+1}}{t_{n+2}}\le0$.  If $Tf$ and $Sf$
  interlace for all $f\in\allpoly$ then there are $a,b\in\reals$ so
  that $Sf = a Tf + bx(Tf)^\prime$.
\end{lemma}
\begin{proof}
  We use the same approach as we have been following, and choose a
  test polynomial $f$ such that $Tf$ has a multiple root - this root
  must also be a root of $Sf$, and this will give us a relation
  between $T$ and $S$. We let $f = x^n(x-a)(x-1)$ so that 
  \begin{align*}
    Tf(x) &= \left(a\,{t_n} - \left( 1 + a \right) \,x\,{t_{1 + n}} +
    x^2\,{t_{2 + n}}\right)x^n\\
\intertext{Ignoring the factor of $x^n$, the discriminant of $f$ is }
\Delta &= {\left( -1 - a \right) }^2\,t_{1 + n}^2 - 4\,a\,{t_n}\,{t_{2 + n}}\\
\intertext{so if we choose $a$ to make $\Delta$ zero} 
a &= \frac{-t_{1 + n}^2 + 2\,{t_n}\,{t_{2 + n}} - 2\,{\sqrt{-\left( {t_n}\,t_{1 + n}^2\,{t_{2 + n}} \right)  + 
         t_n^2\,t_{2 + n}^2}}}{t_{1 + n}^2}\\
\intertext{then we can solve $Tf=0$ with this value of a. The
  hypotheses on $t_n$ ensure that $a\in\reals$. The root is}
r &= \frac{{t_n}\,{t_{2 + n}} - {\sqrt{-\left( {t_n}\,t_{1 + n}^2\,{t_{2 + n}} \right)  + t_n^2\,t_{2 + n}^2}}}  t_{1 + n\,{t_{2 + n}}}\\
\intertext{and if we compute $(Sf)(r)$ we get}
0 &=-\left( \frac{{s_{2 + n}}\,{t_n}\,{t_{1 + n}} + \left( -2\,{s_{1 + n}}\,{t_n} + {s_n}\,{t_{1 + n}} \right) \,{t_{2 + n}}}
    {t_{1 + n}^3\,{t_{2 + n}}} \right) \ \times\\
& \left(t_{1 + n}^2 - 2\,{t_n}\,{t_{2 + n}} +
  2\,{\sqrt{{t_n}\,{t_{2 + n}}\,\left( -t_{1 + n}^2 + {t_n}\,{t_{2
            + n}} \right) }}\right)\\
\intertext{If the second factor is zero then simple algebra shows that
  $t_n=0$, contradicting our hypothesis. Consequently}
0 &= {s_{2 + n}}\,{t_n}\,{t_{1 + n}} + \left( -2\,{s_{1 + n}}\,{t_n} +
  {s_n}\,{t_{1 + n}} \right) \,{t_{2 + n}}\\
\intertext{If we define $u_n = s_n/t_n$ then this relation simplifies
  to}
0 &= u_{n+2} - 2u_{n+1} + u_n
\end{align*}
This implies that $u_n = a+bn$ and hence $s_n = (a+bn)t_n$. It follows
easily that $Sf = aTf + b(Tf)^\prime$.

\end{proof}



  The transformation $x^n\mapsto \frac{\rising{x}{n}}{n!}$ was shown
  in \cite{redmond} to map polynomials whose imaginary part is $1/2$
  to polynomials whose imaginary part is $1/2$. Such a transformation can
  not map $\allpoly$ to itself, as the next lemma shows.

  \begin{lemma}\label{lem:cpx-lines}
    Suppose that $T$ maps polynomials with complex coefficients to
    polynomials with complex coefficients. Suppose that $L_1$ and $L_2$
    are intersecting lines in the complex plane with the property that
    $T\colon{}\allpolyint{L_1}\longrightarrow\allpolyint{L_1}$ and 
    $T\colon{}\allpolyint{L_2}\longrightarrow\allpolyint{L_2}$. Then, $T$ is
     a shifted multiplier transformation. That is, if $\alpha=L_1\cap
     L_2$ then there are constants $a_n$ such that $$T(x-\alpha)^n  =
     a_n \,(x-\alpha)^n \text{ for all integers $n$.}$$
  \end{lemma}
  \begin{proof}
    All polynomials in $\allpolyint{L_1}\cap\allpolyint{L_2}$ are
    constant multiples of powers of $(x-\alpha)$. Since $T$ maps 
$\allpolyint{L_1}\cap\allpolyint{L_2}$ to itself, the conclusion follows.
  \end{proof}

If the transformation considered in \cite{redmond} also mapped
$\allpoly$ to itself, then $T(x-1/2)^n$ would be a multiple of $(x-1/2)^n$,
for all $n$.  This is false for $n=2$.

\section{The domain and range of a linear transformation}
\label{sec:linear-domain}
\index{domain!of a linear transformation}
\index{range!of a linear transformation}

If $T$ is a linear transformation that maps $\allpolyint{\diffi}$ to
$\allpolyint{\diffj}$ where $\diffi,\diffj$ are intervals, then there
are restrictions on the kinds of intervals that $\diffi,\diffj$ may
be. Up to affine transformations there are three types of intervals:
all of $\reals$, half infinite intervals (e.g $(-\infty,1)$), and
finite intervals (e.g. $(0,1)$).  The only restriction in the next
lemma is that $T$ is non-trivial - that is, the image of $T$ is not
all constant multiples of a fixed polynomial.

\begin{lemma} \label{lem:pipj}
  If $\diffj$ is an interval, and
  $T\colon{}\allpolyint{\reals}\longrightarrow \allpolyint{\diffj}$ is
  non-trivial then  $\diffj=\reals$.
\end{lemma}
\begin{proof}
  Suppose that $T\colon{}\allpoly\longrightarrow\allpolyint{\diffj}$. Choose
  $f\greateqeq g$ such that $Tf$ and $Tg$ are not constant multiples
  of one another. Since $f$ and $g$ interlace, so do $Tf$ and $Tg$,
  and hence $Tf + \alpha Tg\in\allpolyint{\diffj}$ for all
  $\alpha\in\reals$. 

  For any $r$ that is not a root of $Tg$ the polynomial $h = Tf +
  \frac{-(Tf)(r)}{(Tg)(r)} \, Tg$ has $r$ for a root. $h$ is not the
  zero  polynomial since $Tf$ and $Tg$ are not constant multiples of
  one another. Consequently there are polynomials in
  $\allpolyint{\diffj}$ with arbitrarily large and small roots, so
  $\diffj$ must be $\reals$.
  
\end{proof}

Table~\ref{tab:domain-and-range} shows that all cases not eliminated
by the lemma are possible.  For instance, the reversal map $f\mapsto
f^{rev}$ satisfies $\allpolyint{(1,\infty)}(n)\longrightarrow
\allpolyint{(0,1)}$.

\begin{table}[htbp] \label{tab:domain-and-range}
$$
    \begin{array}{c|c|c|c|}
\multicolumn{1}{c}{} & \multicolumn{1}{c}{\text{all of }\reals} &
\multicolumn{1}{c}{\text{half infinite}} &
\multicolumn{1}{c}{\text{finite}} \\\cline{2-4}
\text{all of } \reals & x^i\mapsto H_i & - & - \\\cline{2-4}
\text{half infinite} & x^i \mapsto (-1)^{\binom{i}{2}}x^i & x^i\mapsto
\falling{x}{i} & f\mapsto f^{rev}\\\cline{2-4}
\text{finite} & x^i \mapsto {}_qH_i & x^i \mapsto A_i & x^i \mapsto
T_i\\\cline{2-4} 
    \end{array}
$$    
\caption{Realizing possible domains and ranges}
  \end{table}

If $T\colon{}\allpolyint{[a,b]}(n)\longrightarrow\allpoly$ then we can easily
find  upper and lower bounds for the range of $T$.

\begin{lemma} \label{lem:find-range} If
  $T\colon{}\allpolyint{[a,b]}(n)\longrightarrow\allpoly(n)$ and preserves
  interlacing then
  $T\colon{}\allpolyint{[a,b]}(n)\longrightarrow\allpolyint{[r,s}](n)$ where
  $s$ is the largest root of $T(x-b)^n$ and $r$ is the smallest root
  of $T(x-a)^n$.  In addition, for any $f\in\allpolyint{(a,b)}(n)$ we
  have
\[ T(x-a)^n \weaki T(f) \weaki T(x-b)^n \]
\end{lemma}
\begin{proof}
  Choose $f\in\allpolyint{(a,b)}(n)$.  If we join $(x-b)^n$ and
  $(x-a)^n$ by a sequence of interlacing polynomials containing $f$
  (Lemma~\ref{lem:weaki-chain}) then applying $T$ yields
  
  $$
  T((x-b)^n) \greateqeq \cdots \greateqeq Tf \greateqeq \cdots
  \greateqeq T(x-a)^n
$$

This shows that the largest root of $T(x-b)^n$ is at least as large as
the largest root of $T(f)$, and the smallest root of $T(x-a)^n$ is at
most the smallest root of $Tf$.
\end{proof}

  If $f\longleftarrow g$ then the $k$'th root of $f$ is bounded above
  by the $k$'th root of $g$, for all $k$ up to the degree of
  $g$. Here's a simple consequence of this idea.

  \begin{lemma}
    Suppose that $T\colon\allpoly\longrightarrow\allpoly$ preserves
    degree and the sign of the leading coefficient. If
    $f\in\allpoly(n)$ and $a_k$ is the $k$'th largest root of $f$
    ($1\le k\le n$) then
\[
\text{$k$'th largest root of $T(f)$}\ \le\ 
\text{$k$'th largest root of $T((x-a_k)^k)$}
\]
  \end{lemma}
  \begin{proof}
    Suppose $\roots(f)=(a_i)$ and define the polynomials
    \begin{align*}
      g_r &= (x-a_1)(x-a_2)\cdots(x-a_r) \\
      h_r &= (x-a_r)(x-a_{r+1})\cdots(x-a_k)^r \qquad(1\le r\le k)
    \end{align*}
We have interlacings
\[
f=g_n\lesslesseq g_{n-1}\lesslesseq \cdots \lesslesseq
g_k=h_1\lesslesseq h_2\cdots \lesslesseq h_k=(x-a_k)^k
\]
Since $T$ preserves the direction of interlacing, and the degree,
\[
T(f)\lesslesseq  \cdots \lesslesseq
T(g_k)=T(h_1)\lesslesseq \cdots \lesslesseq T((x-a_k)^k)
\]
The conclusion now follows from the observation above.
  \end{proof}

If we apply this kind of argument to the derivative, then we can use
the extra information that $f\lesslesseq f'$. When we replace $h_r$ with
\[
h_r = (x-a_r)(x-a_{r+1})\cdots(x-a_k)^k(x-a_{k+1})
\]
then a similar argument shows
\begin{lemma}[Peyson\cite{peyser}]
  If $f\in\allpoly(n)$ and $1\le k <n$ then the $k$'th root of $f'$
  lies in the interval
\[
\left[ a_k + \frac{a_{k+1}-a_k}{n-k+1}\ ,\
  a_{k+1}-\frac{a_{k+1}-a_k}{k+1}\right]
\]
\end{lemma}

  We can integrate the image of a linear transformation over an
  interval if the images of the endpoints interlace.

  \begin{lemma}\label{lem:int-ab}
    If $T\colon{}\allpolyint{(a,b)}\longrightarrow\allpoly$ and
    \begin{align*}
      1) & \quad T\ \text{preserves interlacing} \\
      2) & \quad T(x-b)^n \greateqeq T(x-a)^n\quad\quad n=1,2,\cdots \\
      \text{then}&
      \quad\quad\displaystyle \int_a^b T(x-t)^n\,dt \in\allpoly
    \end{align*}
  \end{lemma}
  \begin{proof}
From Lemma~\ref{lem:weaki-chain} we can join $(x-b)^n$ and $(x-a)^n$
by a sequence of interlacing polynomials that contains $f$ and $g$. 
Since $T$ preserves interlacing we know that
$$ T(x-b)^n  \greateqeq\cdots \greateqeq T(g) \greateqeq
\cdots \greateqeq T(f) \greateqeq \cdots \greateqeq T(x-a)^n
$$
Now the endpoints of this interlacing sequence interlace, so it
follows that the sequence is a mutually interlacing sequence. In
particular, we know that $T(x-t)^n$ is a family of mutually
interlacing polynomials on $(a,b)$, so the conclusion follows from
Proposition~\ref{prop:family-int}.
  \end{proof}

  It is not easy to have a finite image. Some multiplier
  transformations have finite images. For instance, the identity is a
  multiplier transformation. Less trivially, the map $f\mapsto xf'$
  maps $\allpolyint{[0,1]}$ to itself, and is the multiplier
  transformation $x^i\mapsto ix^i$. There's a simple restriction on
  such transformations.

  \begin{lemma}
    Suppose $T\colon{}x^i\mapsto a_i x^i$ maps $\allpolyint{[0,1]}$ to
    itself. Then $$a_0 \le a_1 \le a_2 \le \cdots$$
  \end{lemma}
  \begin{proof}
    We may assume that all  $a_i$ are positive. Consequently, the roots
    of $T(\,x^n-x^{n-1}\,)$ lie in $[0,1]$, and so 
$$ 0 \le T(\,x^n-x^{n-1}\,)(1) = (a_n x^n - a_{n-1}x^{n-1})(1) = a_n-a_{n-1}.$$
  \end{proof}

If the constants decrease sufficiently rapidly, then no finite
interval can be preserved.

\begin{lemma}
  Suppose that the multiplier transformation $T\colon{}x^i\mapsto a_i x^i$
  maps $\allpolypos$ to itself, and $\limsup |a_n|^{1/n}=0$.  Then
  there are no finite intervals
  $\diffi,\diffj$ such that
  $T(\allpolyint{\diffi})\subset\allpolyint{\diffj}$.
\end{lemma}
\begin{proof}
  It suffices to show that the absolute value of the largest root of
  $T(x+1)^n$ goes to infinity as $n\rightarrow\infty$. If $T(x+1)^n =
  a_nx^n+\cdots+a_0$, then the product of the roots is
  $a_0/a_n$. Since there is a root of absolute value at least
  $\left|{a_0}/{a_n}\right|^{1/n}$, the conclusion follows since 
  $|a_n|^{1/n}$ goes to $0$. 
\end{proof}

For instance, since $\displaystyle\limsup
\left({1}/{n!}\right)^{1/n}=0$, it follows that the exponential
  transformation $x^n\mapsto x^n/n!$ does not preserve any finite
  interval. 

\begin{remark}
  By Corollary~\ref{cor:leading-coef}, the leading coefficients of a
  linear transformation that preserves degree are either all the same
  sign, or they alternate.  If $T$ alternates the signs of the leading
  coefficients, and we define $S(f)(x) = T(f)(-x)$ then $S$ preserves
  the sign of the leading coefficient. Thus, we usually only consider
  linear transformations $T\colon{}\allpoly\longrightarrow\allpoly$ that
  preserve the sign of the leading coefficient.

\index{degree!possible sequences}%

Table~\ref{tab:leading-coef-deg} shows that there are many
possibilities for the degrees of $T(1),T(x),T(x^2),T(x^3),\dots$.
Since $x^n+\alpha x^{n-1}\in\allpoly$
we know $T(x^n)$ and $T(x^{n-1})$ interlace, and so the degrees of
$T(x^n)$ and $T(x^{n-1})$ differ by at most one, provided neither is zero.

  \begin{table}[h]
    \centering
  \begin{tabular}{cll}
\toprule
\multicolumn{2}{c}{Transformation} & degree of
$T(1),T(x),T(x^2),\dots$\\
\midrule
$g\mapsto$ & $ f g$ & $n,n+1,n+2,n+3,\dots$\\
$g\mapsto$ & $ f\ast g$ & $0,1,\dots,n-1,n,0,0,0,\dots$ \\
$g\mapsto$ & $ \diffd^k f(\diffd)g$ & $0,\dots,0,1,2,3,\dots$ \\
$g\mapsto$ & $ g(\diffd)f$ & $n,n-1,\dots,3,2,1,0,0,\dots$\\
$x^k\mapsto$ & $ c_k f+ d_k f'$ & $0, n \text{ or } n-1$\\
\bottomrule
  \end{tabular}\\[.2cm]

    \caption{Degrees of a linear transformation where $f\in\allpoly(n)$.}
    \label{tab:leading-coef-deg}
  \end{table}
  
  It is important to note that the transformation $T\colon{}g\mapsto
  g(\diffd)f$ maps $\allpoly\longrightarrow\allpoly$, yet it does not
  preserve the sign of the leading coefficient, nor do the signs of
  the leading coefficient alternate. Instead, the sign of the leading
  coefficient of $T(g)$ depends on $g$, for it is the sign of
  $g(0)$. This doesn't contradict Lemma~\ref{cor:leading-coef} since
  $T$ doesn't preserve degree.

  The last transformation maps $x^k$ to a linear combination of $f$
  and $f'$, which
  might be constant, so the degree is $0,n$ or $n-1$. 
\end{remark}

\section{Transformations determined by products}
\label{sec:tprod}

When does a linear transformation of the form

\begin{equation}
  \label{eqn:tprod}
  T(x^n) = \prod_{i=1}^n (x+a_i)
\end{equation}
\noindent%
preserve real roots? We can not answer this question but we can show
that such a non-trivial $T$ can not map $\allpoly$ to itself. Next we
will observe that if $T$ in \eqref{eqn:tprod} maps $\allpolypos$ to
$\allpoly$ then there are constraints on the parameters $a_i$.  There
are examples of such root preserving transformations, for in
Corollary~\ref{cor:rising-factorial} we will see that the choice $a_i=i-1$ maps
$\allpolypos$ to itself.

\begin{lemma} \label{lem:tprod}
  If $T$ is given in \eqref{eqn:tprod} and $T$ maps $\allpoly$ to
  itself then all $a_i$ are equal.
\end{lemma}
\begin{proof}
  Choose $b\in\reals$, assume $T$ preserves roots, and consider
  \begin{align}
    T((x-b)^2) &= b^2 -2b(x+a_1)+(x+a_1)(x+a_2) \label{eqn:tprod-2} \\
\intertext{The discriminant is}
& -(a_2-a_1)(4b+a_1-a_2) \notag
  \end{align}
  If $a_1\ne a_2$ then we can choose $b$ to make the discriminant
  negative, so such a linear transformation $T$ does not preserve
  roots. More generally
\begin{equation}
  \label{eqn:tprod-3}
  T(x^r(x-b)^2) = 
\prod_{i=1}^r(x+a_i) \,\cdot\,\left(b^2-2b(x+a_{r+1}) + (x+a_{r+1})(x+a_{r+2})\right)
\end{equation}
We can apply the same argument as above to conclude that
$a_{r+1}=a_{r+2}$, and hence if $T$ preserves roots then all $a_i$ are
equal.
\end{proof}

Note that the conclusion of the lemma implies that there is an $a$
such that $T(f) = f(x+a)$. We  restate the conclusion:

\begin{cor}
  If $T\colon{}\allpoly\longrightarrow\allpoly$ maps monic polynomials
  to monic polynomials, and $T(x^i)$ divides $T(x^{i+1})$ for all
  $i\ge0$ then $T(x^i) = (x+a)^i$ for some $a$.
\end{cor}

\begin{lemma} \label{lem:tprod-4}
  If $T$ is given in \eqref{eqn:tprod} and
  $T\colon{}\allpolypos\longrightarrow\allpoly$ then
$$ a_1 \le a_2 \le a_2 \le \cdots $$
Moreover, if $a_{n+1}=a_n$ for some $n$ them $a_m=a_n$ for all $m\ge n$.
\end{lemma}
\begin{proof}
  Choose negative $b$ so that $T((x-b)^2)$ is in $\allpoly$. If we
  choose $|b|$ sufficiently large in  \eqref{eqn:tprod-2} then the sign 
    of the discriminant is the sign of $a_2-a_1$ and hence $a_1\le
    a_2$. The same argument applied to \eqref{eqn:tprod-3} shows that
    $a_{r+1}\le a_{r+2}$.

    Next, suppose that $a_{n+1}=a_n$. If $f =  (x+a)^3 x^{n-1}$ then 
    \begin{align*}
      T(f) &= \prod_1^{n+2}(x+a_i) + 3a \prod_1^{n+1}(x+a_i) + 3a^2
      \prod_1^{n}(x+a_i) + a^3 \prod_1^{n-1}(x+a_i)\\
      &= \prod_1^{n-1}(x+a_i)\quad \times \\ &
      (x+a_{n+2})(x+a_{n+1})(x+a_{n}) +
      3a (x+a_{n+1})(x+a_{n}) + 3a^2 (x+a_{n}) + a^3 \\
      \intertext{and without loss of generality we may set
        $a_n=a_{n+1}=0$. The resulting cubic polynomial}
      & (x+a_{n+2})x^2 + 3a x^2 + 3a^2x + a^3\\
    \end{align*}
\noindent%
is in $\allpoly$ for all positive $a$. If we substitute $a=a_{n+2}$
then the roots are $(-.42\pm.35i)a_{n+2}$ and $-3.1 a_{n+2}$.
Consequently, $a_{n+2}=0$. We can now continue by induction.

\end{proof}

For other examples of such transformations: see
Lemma~\ref{lem:factorial-variant}.

\section{Composition}
\label{sec:oper-composition}
\index{composition}

If $T$ is a linear transformation, and $f,g$ are polynomials then we
can construct a new polynomial $f(T)g$.  In this section observe that
such composition is ubiquitous, \index{composition!is ubiquitous}
and give conditions for $f(T)g$ to define a map
$\allpolyint{\diffi}\times\allpoly\longrightarrow\allpoly$. 

\begin{lemma}
  If $T$ is a linear transformation that preserves degree then there
  is a linear transformation $W$ and a constant $c$ so that $T(f) =
  cf(W)(1)$ for every polynomial $f$.
\end{lemma}
\begin{proof}
  Since $T$ preserves degree the polynomials
  $T(x^0),T(x^1),T(x^2),\dots$ form a basis for the space of all
  polynomials. Define the linear transformation $W$ by $W(T(x^i)) =
  T(x^{i+1}),$ and set $c=T(1)$. We show by induction that for this
  choice of $W$ and $c$ we have $Tf = cf(W)(1)$. It suffices to
  establish this for the basis above.  We note that $T(x^0) = c =
  cW^0(1)$. The inductive step is
  $$ T(x^{n+1}) = W(T(x^n)) = W(cW^n(1)) = cW^{n+1}(1).$$
\end{proof}

We can also define $W$ by $W(Tf) = T(xf)$. Since $T$ is invertible we
also have that $W(g) = T(x(T^{-1}g))$ This representation shows
that if $T$ is the identity then $W$ is multiplication by $x$.
The inverse of a transformation expressed as a  composition has an
inverse that is easily expressed as a composition.

\begin{lemma}
  If $W$ is a linear transformation that increases degree by one,
  $T(1)=1$, and $T(f) = f(W)(1)$ then
$$ T^{-1}(f) = f(T^{-1}WT)(1)$$
\end{lemma}
\begin{proof}
  It is enough to verify the conclusion on a basis, and a good choice
  for basis is $\{W^n(1)\}$. We note
$$ T^{-1}W^n(1) = (T^{-1}WT)^n T^{-1}(1) = (T^{-1}WT)^n(1)$$
\end{proof}

We now look at some general properties of compositions.

\begin{theorem} \label{thm:compose-2}
  Suppose that the linear transformation $S$ has the property that
  $(S-\alpha)h\in\allpoly$ for all $h\in\allpoly$ and all $\alpha$ in an
  interval $I$. The bilinear map $f\times g\mapsto f(S)g$ defines a
  bilinear transformation
  $\allpolyint{I}\times\allpoly\longrightarrow\allpoly$.
\end{theorem}
\begin{proof}
  We may write $f=(x-a_1)\cdots(x-a_n)$ where all $a_i$ are in $I$.
  Since
  $$
  T(g) = f(S)g = (S-a_1)\cdots(S-a_n)g$$
  we see that the hypothesis on $S$ guarantees that $Tg\in\allpoly$.
\end{proof}

A similar argument yields
\begin{cor} \label{cor:compose-2}
  Suppose that $T\colon{}\allpolypos\longrightarrow\allpolyneg$ and for any
  $f\in\allpolypos$ we have that $f \longleftarrow Tf$. The map
  $f\times g\mapsto f(T)g$ maps $\allpolypos\times\allpolyneg$ to
  $\allpolypos$. If $g\in\allpolyneg$ is fixed the map $f\mapsto
  f(T)g$ preserves interlacing.
\end{cor}

\section{Recurrence relations for orthogonal polynomials}
\label{sec:linear-recur}

All orthogonal polynomials satisfy a \emph{three term recurrence}
\eqref{eqn:recur-1}.  We use this recurrence  to establish
recurrence relations for linear transformations based on orthogonal
polynomials. We also get recurrence relations for the inverses of
these linear transformations. Finally, we specialize to polynomials of
Meixner type.
\index{recurrence!orthogonal polynomials}

Assume that $p_1,p_2,\dots$ is a sequence of orthogonal polynomials
satisfying the three term recurrence

\begin{equation}
  \label{eqn:recur-1}
  p_{n+1} = (a_n x+ b_n) p_n + c_n p_{n-1}
\end{equation}

and define the linear transformations
\begin{align*}
  T\ : &\ x^n \mapsto p_n \\
  T^{-1}\ : &\ p_n \mapsto x^n\\
  A\ : &\ x^n \mapsto a_nx^n \\
  B\ : &\ x^n \mapsto b_nx^n \\
  C\ : &\ x^n \mapsto c_nx^{n-1} \\
\end{align*}

From the recurrence we find that

\begin{align}
  T(x\,x^n) & = a_n x p_n + b_n p_n + c_n p_{n-1}  \nonumber \\
  & = x\ T(A(x^n)) + T(B(x^n)) + T(C(x^n))  \nonumber \\
\intertext{Since this holds for all $n$, and $T$ is linear, we have for
  all polynomials $f$}
T(xf) & = xTA(f) + TB(f) + TC(f) \label{eqn:recur-2} \\
\intertext{Similarly we have a recurrence for the inverse of $T$}
T^{-1}(xp_n) &= T^{-1}((1/a_n) p_{n+1} - (b_n/a_n) p_n -
  (c_n/a_n)p_{n-1}) \nonumber \\
 &= (1/a_n) x^{n+1} - (b_n/a_n) x^n - (c_n/a_n) x^{n-1} \nonumber\\
& = xA^{-1}T^{-1}(p_n) - BA^{-1}T^{-1}(p_n) - CA^{-1}T^{-1}(p_n)
  \nonumber \\
\intertext{and consequently we have for all polynomials $f$}
T^{-1}(xf) &=
 xA^{-1}T^{-1}(f) - BA^{-1}T^{-1}(f) - CA^{-1}T^{-1}(f)
 \label{eqn:recur-3} \\
 &= (xA^{-1} - BA^{-1} -CA^{-1})\, T^{-1}f \nonumber \\
\intertext{If $T(1) = 1$ then we have an explicit composition
  representation of $T^{-1}$}
T^{-1}(f) &= f(xA^{-1} - BA^{-1} -CA^{-1})(1) \label{eqn:recur-1a}
\end{align}

\index{recurrence!orthogonal polynomials inverse}
Notice that this last representation is not well defined. However, we
take \eqref{eqn:recur-1a} to mean that we do not multiply $f$ out, and then
substitute, but rather we use the product representation that follows
from \eqref{eqn:recur-3}.

These recurrences are too general to be useful. If our goal is to
establish that a certain $T$ preserves interlacing, then we must have that
$T(xf)\lesslesseq T(f)$. This necessitates that there is some sort of
interlacing relationship between $A,B,C$. We have seen that if $A$ is
the identity then it's probable that $B,C$ are functions of
derivatives. There is a well known class of polynomials for which this
is true.

\index{Meixner class of polynomials} 
\index{polynomials!of Meixner type}

The Meixner class of orthogonal polynomials \cite{godsil}*{page~165}
consists of those orthogonal polynomials $p_n$ that satisfy a
recurrence

\index{recurrence!Mexier}
\begin{equation}
  \label{eqn:meixner}
  p_{n+1} = (x - a -\alpha n)p_n - (bn + \beta n(n-1))p_{n-1}
\end{equation}

Upon specializing \eqref{eqn:recur-3} and \eqref{eqn:recur-1a}, where
$A$ is the identity, $B = -(a+\alpha x\diffd)$, $C = -(b\diffd + \beta
x \diffd^2)$, we find the recurrences

\begin{align}
  \label{eqn:meixner-1}
  T(xf) & = xTf - T(\,(a+\alpha x \diffd+b\diffd+\beta x \diffd^2)\,f) \\
  T^{-1}(f) & =  f(x+a+\alpha x \diffd+b\diffd+\beta x \diffd^2)(1)
\label{eqn:meixner-2}
\end{align}

We have seen several examples of orthogonal polynomial
families. Although \emph{families} of orthogonal polynomials have many
special properties, there is nothing special about any particular
orthogonal polynomial, as the next theorem shows.

\begin{theorem}
  Every $f\in\allpoly$ is an orthogonal polynomial. In other words,
  there is a measure  that determines an orthogonal family of
  polynomials, and $f$ is a member of that family.
\end{theorem}
\begin{proof}
  Assume $f\in\allpoly(n)$. We will show that there is a sequence of
  polynomials $p_0,p_1,\dots,p_n=f$ where the degree of $p_i$ is $i$,
  and there are non-negative constants $a_i$ such that
  $p_{i+1} = (a_ix+b_i)p_i - p_{i-1}$ for $1\le i < n$. It follows
  from Favard's Theorem \cite{szego} that the sequence $p_0,\dots,p_n$
  is a set of orthogonal polynomials determined by a measure.
  
  We construct the $p_i$ inductively. Assume that $f$
  has positive leading coefficient. Let $p_n=f$, and
  $p_{n-1}=f^\prime$. Since $p_n   \lesslesseq p_{n-1}$ we can find a
  $p_{n-2}$ such that $p_n =  (a_nx+b_n)p_{n-1} - p_{n-2}$. Since
  $p_{n-1}\lesslesseq p_{n-2}$ we continue inductively.
\end{proof}

\section{Linear transformations satisfying a recursion}
\label{sec:linear-interlace}

If a linear transformation satisfies a recursion then we
can sometimes conclude that the transformation maps $\allpolyalt$ or
$\allpolypos$ to itself.

\begin{theorem} \label{thm:interlace}
   Suppose that linear transformations $A,B,C,D,E,F,H$ map polynomials
   with positive leading coefficients to polynomials with positive
   leading coefficients. Suppose that for all $f\in\allpolypos$ we
   have
   \begin{equation}
     \label{eqn:thm-interlace}
     Af \greateqeq f \greateqeq Bf\hspace*{1cm}
     Cf \greateqeq f \longleftarrow Df\hspace*{1cm}
     Ef \greateqeq f \longleftarrow Ff\hspace*{1cm}
     Hf \lesslesseq f     
   \end{equation}
   then the linear transformation defined by
   $$
   T(xf) = (\,HT + xTB + T(C-D) + (E-F)T\,)\,f$$
   maps $\allpolypos$ to
   itself and preserves interlacing. If the interlacing assumptions
   \eqref{eqn:thm-interlace} hold for all $f\in\allpolyalt$ then
   $$
   S(xf) = (\,HS + xSA + S(C-D) + (E-F)S\,)\,f$$
   maps $\allpolyalt$ to    itself and preserves interlacing.
\end{theorem}
\begin{proof}
  We prove the theorem by induction on the degree of $f$, so assume
  that $T$ maps $\allpolypos(n)$ to itself and preserves interlacing.
  By Corollary~\ref{cor:quant-transform} it suffices to show that $T(xf)\greateqeq Tf$ for all
  $f\in\allpolypos(n)$. If $f\in\allpolyneg(n)$ then $f\greateqeq Bf$
  and so $Tf\greateqeq TBf$. Now all roots of $TBf$ are negative, so
  $xTBf\lesslesseq Tf$. Next, we apply our  assumptions
  \eqref{eqn:thm-interlace} to the polynomial $Tf\in\allpolypos(n)$
  and find 
  \begin{xalignat*}{2}
    HTf & \lesslesseq Tf & \\
    TCf & \greateqeq Tf & Tf \greateqeq TDf \\
    ETf & \greateqeq Tf & Tf \greateqeq FTf 
  \end{xalignat*}
  Finally, use Lemma~\ref{lem:add-interlace} to conclude
   $$T(xf) = (\,HT + xTB + T(C-D) + (E-F)T\,)\,f \lesslesseq Tf$$ 
   and hence $T(xf)\lesslesseq Tf$ which implies that $T$ maps
   $\allpolypos(n+1)$ to itself and preserves interlacing. The case
   for $f\in\allpolyalt$ is similar.
\end{proof}

\section{\Mobius\ transformations}
\label{sec:mobius}

\index{M\"{o}bius transformation}

A \Mobius\ transformation determines a  linear transformation on
polynomials of fixed degree, and also an action of the space of all linear
transformations.  Recall that a \Mobius\ transformation (or linear
fractional transformation) is a transformation of the form
\begin{equation} \label{eqn:mobius-1}
  M:z \mapsto \frac{az+b}{cz+d}.
\end{equation}
\noindent%
To construct a map between polynomials we  restrict ourselves
to polynomials of a fixed degree $n$, and homogenize:
\begin{align} \label{eqn:mobius}
 \widetilde{M} : x^i & \mapsto (ax+b)^i (cx+d)^{n-i}.
\intertext{We can express $\widetilde{M}$ in terms of $M$:}
\widetilde{M}(f) &= (cx+d)^n f(Mx). \nonumber
\end{align}

\index{reverse!of a polynomial} We have already seen an example of a
\index{\Mobius\ transformation}\Mobius\ transformation in
\eqref{eqn:homog-3}. If we take $M(z) = 1/z$ then $\widetilde{M}(f) =
z^nf(1/z)$ which is the reverse of $f$.  If $I$ is an interval, then
we define ${M}(I)$ to be the image of $I$. With this notation, we have
the elementary lemma

\begin{lemma} \label{lem:mobius}
  If $M$ is a \Mobius\ transformation \eqref{eqn:mobius} then
  $\widetilde{M}$ is a linear transformation that maps
$\allpolyint{I}(n)$ to $\allpolyint{M^{-1}(I)}(n)$ bijectively.
\end{lemma}

\begin{proof}
  If $f\in\allpolyint{I}$ then $\widetilde{M}(f) = (cx+d)^n f(Mx)$. 
  A root $r$ of $\widetilde{M}(f)$ satisfies $f(Mr)=0$ so $r=M^{-1}s$ where 
  $s$ is a root of $f,$ and  thus $s\in I$. The inverse of
  $\widetilde{M}$ is $\widetilde{M^{-1}}$ and so $\widetilde{M}$ is a
  bijection. 
\end{proof}



\begin{lemma}
  If $\diffi,\diffj$ are intervals, and $T$ is a linear transformation
  $T\colon{}\allpolyint{\diffi}\longrightarrow\allpolyint{\diffj}$ and
  $\alpha\in\reals$ then the linear transformation $S(x^i) = \alpha^i
  T(x^i)$ maps $\allpolyint{\diffi/\alpha}\longrightarrow\allpolyint{\diffj}.$
\end{lemma} 
\begin{proof}
 The result follows from the commuting diagram 
 
 \centerline{ \xymatrix{
     \allpolyint{\diffi/\alpha} \ar@{->}[d]_{x\mapsto \alpha x}
     \ar@{..>}[rr]^S &&  \allpolyint{\diffj}\\ 
     \allpolyint{\diffi} \ar@{->}[rru]_T
}}
\end{proof}

The next lemma shows that only the type of interval matters.
\begin{lemma}\label{lem:interval-bijection}
  Assume $\diffi$ and $\diffj$ are intervals, and that $\affa_1$ and
  $\affa_2$ are affine transformations. There is a 1-1 correspondence
  between linear transformations mapping $\allpolyint{\diffi}$ to
  $\allpolyint{\diffj}$ and linear transformations mapping
  $\allpolyint{\affa_1\diffi}$ to $\allpolyint{\affa_2\diffj}$
\end{lemma}
\begin{proof}
  The correspondence is given by the diagram

\centerline{ \xymatrix@-1pc{
\allpolyint{\diffi} \ar@{->}[rrr]  \ar@{->}[dd]_{\affa_1 x\mapsto x}& & &
    \allpolyint{\diffj} \ar@{<-}[dd]^{x\mapsto\affa_2 x}\\
&&& \\
\allpolyint{\affa_1\diffi} \ar@{->}[rrr] & & & 
    \allpolyint{\affa_2\diffj} \\
}} 

\end{proof}

As another application, the next lemma shows  two linear transformations 
that are related by a \index{\Mobius\ transformation}\Mobius\ transformation.

\begin{lemma} \label{lem:xdi}
  Choose a positive integer $n$, and define linear transformations
  $T(x^i)=\binom{x+n-i}{n}$, $S(x^i)=\binom{x}{i}$, and the \Mobius\
  transformation $M(z) = \frac{z}{z+1}$. Then $T=S\widetilde{M}$ on
  $\allpoly(n)$. 
\end{lemma}
\begin{proof}
  It suffices to verify that $T=S\widetilde{M}$ for a basis of
  $\allpoly(n)$, where $\widetilde{M}(x^i) = x^i(x+1)^{n-i}$, so
  consider
  \begin{align*}
    S\widetilde{M}(x^i) &= S\left(x^i(1+x)^{n-i}\right) \\
    & = S\left(\sum_{j=0}^{n-i} \binom{n-i}{j} x^{i+j}\right) \\
    &= \sum_{j=0}^{n-i} \binom{n-i}{j}\binom{x}{i+j}\\
    &= \sum_{j=0}^{n-i} \binom{n-i}{n-(i+j)}\binom{x}{i+j}\\
    \intertext{and the last expression equals $T(x^i)$ using the
      Vandermonde identity} \binom{a+b}{c} &= \sum_{r+s=c}
    \binom{a}{r}\binom{b}{s} 
  \end{align*}
\end{proof}

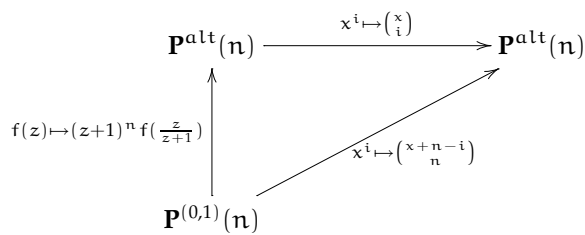
\begin{figure}\label{fig:binomial-diagram}
\centerline{ \xymatrix@-1pc{
    \allpolyalt(n) \ar@{->}[rrrrr]^{x^i\mapsto \binom{x}{i}} & & & & &
    \allpolyalt(n) \\
& & & & & \\
& & & & & \\
\allpolyint{(0,1)}(n) \ar@{->}[uuurrrrr]_{x^i\mapsto \binom{x+n-i}{n} }
\ar@{->}[uuu]^{f(z)\mapsto(z+1)^nf(\frac{z}{z+1})}
}} 
\caption{A \Mobius\ transformation acting on a linear transformation}
\end{figure}

We will later show that $S\colon{}\allpolyalt\longrightarrow\allpoly$. $M$
maps $(0,1)$ bijectively to $(0,\infty)$, so we have a commuting
diagram (Figure~\ref{fig:binomial-diagram}) of maps that preserve roots.

The second action determined by a \index{\Mobius\ transformation}\Mobius\ transformation is in the
space of linear transformations.  Let $T$ be  a linear transformation
that  preserves degree.  If $M$ is given in \eqref{eqn:mobius-1}
then we define a new linear transformation by
\begin{align}
  \label{eqn:mobius-2}
  T_M(x^r) &= (cx+d)^r\  (T(x^r))(Mx) \\
\intertext{If $T(x^r) = p_r$ then we can write this as}
  T_M(x^r) &= (cx+d)^r \,p_r(\frac{ax+b}{cx+d})\nonumber
\intertext{Note that $T_M\ne \widetilde{M}T$, since on $\allpoly(n)$}
 \widetilde{M}T\,(x^r) &= (cx+d)^n \,p_r(\frac{ax+b}{cx+d})\nonumber
\end{align}
If $p_r(0)\ne0$ then  both $T(x^r)$ and $T_M(x^r)$ have degree
$r$.  


We first note that this action respects composition.

\begin{lemma} \label{lem:mobius-comp}
  If $M,N$ are \Mobius\ transformations then $ (T_M)_{N} = T_{NM}$.
\end{lemma}
\begin{proof}
The proof follows from the commuting diagram

\centerline{
\xymatrix{
x^n \ar@{->}[rrrr]^T \ar@{->}[drr]^{T_M} \ar@{->}[dd]^{T_{NM}} &&&&
p_n \ar@{->}[dll]^{{M}}  \ar@/^3pc/[ddllll]^{NM}\\
&& M\,p_n \ar@{->}[dll]^{{N}}\\
{NM}\,p_n
}}
\end{proof}

Consider some examples.

\begin{enumerate}
\item If $T$ is the identity then $T_M\,(x^n) = (ax+b)^n$.
  Equivalently, $T_M\,(f) = f(ax+b)$.
\item If $T$ is any degree preserving linear transformation,
  $Mz=az+b,$ and we let $T(x^n)=p_n$ then $T_M\,(x^n) = p_n(ax+b)$.
  Consequently $T_M\,(f) = (Tf)(ax+b)$. Moreover, if $\diffi,\diffj$
  are intervals and
  $T\colon{}\allpolyint{\diffi}\longrightarrow\allpolyint{\diffj}$ then
  $T_M:\allpolyint{M^{-1}\diffi}\longrightarrow\allpolyint{\diffj}$.
\item The most interesting (and difficult) case is when the
  \index{reverse!of a polynomial} denominator is not constant.
  Because of the composition property (Lemma~\ref{lem:mobius-comp}), and the
  examples above, we only consider the case $M(z) = 1/z$. Since
  $x^nf(1/x)$ is the reverse of $f$ when $f$ has degree $n$ we see
  that
\[
  T_{1/z}(x^n) = \rev {\ T(x^n)}.
\]

\index{T@$T_{1/z}$}

  For example, suppose that \index{falling factorial}
  \begin{align*}
    T(x^n) &= (x-1)(x-2)\cdots(x-n).\\
    \intertext{Applying $M(z)=1/z$ yields}
    T_{1/z}\,(x^n) &= (1-x)(1-2x)\cdots(1-nx)\\
  \end{align*}
\end{enumerate}

\begin{example}
  We can combine \Mobius\ transformations with other root preserving
  linear transformations. For instance, we will see
  (Corollary~\ref{cor:rising-factorial}) that $\rising{x}{n}\mapsto
  x^n$ maps $\allpolyalt$ to itself and
  (Lemma~\ref{lem:falling-non-int}) $x^k\mapsto \falling{\alpha+n}{k}
  x^k$ maps $\allpolyalt$ to itself. If $M(z)=z/(1-z)$ the composition

\centerline{
\xymatrix{
\allpolyalt(n) \ar@{->}[rr]^{\rising{x}{k}\mapsto x^k} &&
\allpolyalt(n) \ar@{->}[rr]^{x^k\mapsto \falling{\alpha+n}{k}x^k} &&
\allpolyalt(n) \ar@{->}[rr]^{\widetilde{M}} &&
\allpolyint{(0,1)}
}}
\noindent%
is the linear transformation $\rising{x}{k}\mapsto
\falling{\alpha+n}{k}x^k(1-x)^{n-k}$ and maps $\allpolyalt(n)$ to
$\allpolyint{(0,1)}$. This is Example 1 in \cite{iserles}.
\end{example}

\section{Transformations and commutativity properties}
\label{sec:commute-trans}
  
  We first characterize transformations that commute with $x\mapsto
  \alpha x$, and those that commute with $x\mapsto x+1$. We also
  characterize some transformations defined by commuting diagrams.
  Finally, we list some diagrams that are satisfied by various linear
  transformations.

  \begin{lemma} \label{lem:commute-xpi}
    Suppose that the linear transformation $T$ preserves degree. Then
    $T$ commutes with $x\mapsto
    \alpha x$ where $\alpha\ne1$
%
%
%
if and only if there is a $g(x)=\sum a_ix^i$ such that
$Tf = f\ast g$. Equivalently, $T$ commutes with $x\mapsto \alpha x$
if and only if $T$ is a multiplier transformation. 

  \end{lemma}

  \begin{proof}
    Simply consider the action of $T$ on $x^n$. Assume that $T(x^n) =
    \sum_{k=0}^n b_i x^i$. Commutativity implies that 
$$ \alpha^n T(x^n) =    \alpha^n  \sum_{k=0}^n b_i x^i  = T(
x^n)(\alpha x)
= \sum_{k=0}^n b_i \alpha^i x^i $$
This implies that all $b_i$ are zero except for $b_n$, so $T$ is a
multiplier transformation. 
  \end{proof}

  \begin{lemma}
    Suppose that the linear transformation $T$ preserves degree. Then
    $T$ commutes with $\affa:x\mapsto x+1$
%
%
%
if and only if there is a $g(x)=\sum a_ix^i$ such that
$\displaystyle T(x^n) = \sum_{i=0}^n a_{n-i}\binom{n}{i}x^i$. Equivalently, 
$$ T(f) = (g \ast \affa f)^{rev}.$$

  \end{lemma}
  \begin{proof}
    That $T(x^n)$ has that form can be proved by induction. Conversely,
    note that the conclusion is linear in $g$, so we may take
    $g(x)=x^s$. It is easy to see that the diagram commutes.

\centerline{\xymatrix{ x^n \ar@{->}[d]_{{x\mapsto
            x+1} } 
        \ar@{->}[rrr]^{T } 
        &&& {\sum_{i=0}^n \binom{n}{i}a_{n-i}x^i}
        \ar@{->}[d]^{{x\mapsto x+1} } \\ 
        {\sum_{j=0}^n \binom{n}{j}x^j} \ar@{->}[rrr]^{T } 
        &&& {\sum_{i,j=0}^n \binom{n}{i}a_{n-i}\binom{i}{k}x^k} }}

  \end{proof}

  \begin{remark}
    If we ask that $T$ commutes with the seemingly more general
    transformation $x\mapsto x+y$, we find that $T$ is the same as is
    described in Lemma~\ref{lem:commute-xpi}. For example, the
    Bernoulli polynomials satisfy this commutativity property. However,
    this is not useful for us as most Bernoulli polynomials do not
    have all real roots. 

If we choose
    $g(x)=e^{-x}$ then 
$$ T(x^n) = \sum \binom{n}{k}\frac{-1)^k}{k!} x^k = L_n^{rev}(x)$$ so
we have a commuting square in two variables

\index{Laguerre polynomials!identities}
\begin{equation}\label{eqn:commute-laguerre-1}
  \checked \xymatrix{
.
      \ar@{->}[d]_{{x\mapsto x+y} }           
      \ar@{->}[rrr]^{{ x^n\mapsto L_n^{rev}(x)} }         
      &&&
 .
      \ar@{->}[d]^{{ x\mapsto x+y} } \\        
  .
      \ar@{->}[rrr]^{{ x^n\mapsto L_n^{rev}(x)} }         
      &&&
   .
}
\end{equation}

  \end{remark}

\begin{example} \label{ex:branden}
If $T(x^k) = \frac{\rising{x}{k}}{k!}$ and $\affa\, f(x) = f(1-x)$ then 
$\affa T = T \affa$.
\checked  

  \end{example}

\begin{example}
If $T(x^k) = \frac{\falling{x}{k}}{k!} = \binom{x}{k}$ and $\affa\,f(x) = f(-x-1)$
then $\affa T = T \affa $. 

  \end{example}

\begin{example}
\index{Hermite polynomials!identities}
$T\colon{}x^n\mapsto H_n$ satisfies
  \begin{equation}
    \checked\label{eqn:iden-diag-7}
    T(f)(x+y) = T_\ast f(x+2y) 
  \end{equation}
which leads to the diagram

\centerline{\xymatrix{
. \ar@{->}[rr]^T \ar@{->}[d]_{x\mapsto  x+2y} && . \ar@{->}[d]^{x\mapsto x+y} \\
. \ar@{->}[rr]^T  &&  
.
}} 

\end{example}

\begin{example}\index{Hermite polynomials!identities}
  $T\colon{}x^k\mapsto H_k(x)x^{n-k}$ satisfies
  \begin{equation}
    \label{eqn:iden-diag-5}
 T(x-2+\alpha)^n = H_n(\frac{\alpha x}{2}) =
      (-1)^n T(x-2-\alpha)^n
  \end{equation}

If we define $S(f) = T(f(x-2))$ then 

  \centerline{\xymatrix{
.
      \ar@{->}[d]_{x\mapsto \alpha{+x} }           
      \ar@{->}[rrr]^{S }         
      &&&
.
      \ar@{->}[d]^{x\mapsto \alpha{-x} } \\        
.
      \ar@{->}[rrr]^{S }         
      &&&
.
}}

\end{example}

\begin{example}
    $ T\colon{}x^k \mapsto \rising{x}{k}\falling{x-\alpha}{n-k}$ satisfies
    the identity on polynomials of degree $n$
\begin{equation}\label{eqn:iden-diag-1}
 (Tf)(\alpha-x) = (-1)^{n}\,T(f^{rev})
\end{equation}

  \centerline{\xymatrix{
.
      \ar@{->}[d]_{reverse }           
      \ar@{->}[rrr]^{T }         
      &&&
.
      \ar@{->}[d]^{x\mapsto\alpha-x } \\        
.
      \ar@{->}[rrr]^{(-1)^nT }         
      &&&
.
}}

  \end{example}

\begin{example}
 $T\colon{}x^i\mapsto \falling{x+n-i}{n}$ satisfies
 \begin{equation}
   \label{eqn:iden-diag-2}
   Tf = (T f^{rev})(-1-x)
 \end{equation}

since $\falling{-x-1+i}{n} = (-1)^n\falling{x+n-i}{n}$. Equivalently, 

  \centerline{\xymatrix{
.
      \ar@{->}[d]_{reverse }           
      \ar@{->}[rrr]^{T }         
      &&&
.
      \ar@{->}[d]^{x\mapsto-1-x } \\        
.
      \ar@{->}[rrr]^{T }         
      &&&
.
}}

  \end{example}

\begin{example} \index{Laguerre polynomials!identities} The
  Laguerre transformation $T\colon{}x^n\mapsto L_n(x)$ satisfies several
  identities. This square is the Laguerre identity
  \eqref{eqn:laguerre-addition}.

\begin{equation}\label{iden-diag-laguerre}
    \checked\centerline{\xymatrix{
x^n
      \ar@{->}[d]_{{x\mapsto x+y-1} }           
      \ar@{->}[rrr]^{T }         
      &&&
      L_n(x)
      \ar@{->}[d]^{homogenize } \\        
      (x+y-1)^n
      \ar@{->}[rrr]^{T }         
      &&&
      y^n\,L_n(x/y)
}}
\end{equation}

$T$ also satisfies

  \centerline{\xymatrix{
  .
      \ar@{->}[d]_{{x\mapsto \frac{x+y-1}{y}    }}           
      \ar@{->}[rrr]^{{ T  }}         
      &&&
  .
      \ar@{->}[d]^{{ x\mapsto \frac{x}{y}   }} \\        
   .
      \ar@{->}[rrr]^{{ T   }}         
      &&&
    .
}}

\index{Laguerre polynomials!identities}
\index{Hermite polynomials!identities}
We can combine the Laguerre and Hermite polynomials.
  Define the linear transformation $T$ acting on $\allpoly(n)$, and
  two induced transformations:
   \begin{align*}
     T(z^k) & = H_k(z)L_{n-k}(z)\\
     T_x(x^ry^s) & = H_r(x)L_{n-r}(x)\,y^s\\
     T_y(x^ry^s) & = x^r\,H_s(y)L_{n-s}(y)
   \end{align*}
We have the following commutative diagram

  \centerline{\xymatrix{
      & (x+y)^n 
      \ar@{->}[dl]_{T_x } 
      \ar@{->}[dr]^{T_y } 
      \\
      T_x(x+y)^n &&
      T_y(x+y)^n
      \ar@{->}[ll]^{x\mapsto 2x,y\mapsto y/2 }         
}}

\noindent
which if written out is
$$ \sum_{k=0}^n \binom{n}{k}\,H_k(x)L_{n-k}(x)\,y^{n-k} =
 \sum_{k=0}^n \binom{n}{k}\, H_{n-k}(y/2)L_{k}(y/2)\,(2x)^k
$$

\end{example}

\begin{example}\label{ex:bin-type}
\index{polynomials! of binomial type}
\index{binomial type}

  The falling factorials $\falling{x}{k}$ and $L_n/n!$ are polynomials
  of binomial type, and satisfy

  \begin{gather*}
    \falling{x+y}{n} = \sum_{k=0}^n \binom{n}{k}
    \falling{x}{k}\falling{y}{n-k} \\
    \frac{L_n(x+y)}{n!} = \sum_{k=0}^n \binom{n}{k}
    \frac{L_k(x)}{k!}    \frac{L_{n-k}(y)}{(n-k)!}
\intertext{A polynomial family $\{p_n\}$ of binomial type satisfies}
    p_n(x+y) = \sum_{k=0}^n \binom{n}{k}
    p_k(x)p_{n-k}(y)
  \end{gather*}

For a general family of polynomials of binomial type, the
transformation $T\colon{}x^n\mapsto p_n(x)$ satisfies

  \centerline{\xymatrix{
x^n
      \ar@{->}[d]_{{ x\mapsto x+y    }}           
      \ar@{->}[rrr]^{{ T  }}         
      &&&
p_n(x)
      \ar@{->}[d]^{{ x\mapsto x+y    }} \\        
      (x+y)^n
      \ar@{->}[rrr]^{{ T_x\,T_y   }}         
      &&&
      p_n(x+y)
}}
 where the transformations $T_x,T_y$ are the induced
 transformations. Considering the inverse transformations yields

  \centerline{\xymatrix{
x^n
      \ar@{->}[d]_{{ x\mapsto x+y    }}           
      \ar@{<-}[rrr]^{{ T^{-1}  }}         
      &&&
p_n(x)
      \ar@{->}[d]^{{ x\mapsto x+y    }} \\        
      (x+y)^n
      \ar@{<-}[rrr]^{{ T_x^{-1}\,T_y^{-1}   }}         
      &&&
      p_n(x+y)
}}

\end{example}

\section{Singular points of transformations}

\index{linear transformation!singular value} If $T$ is a
transformation, and $f\in\allpoly$ is a polynomial for which the
degree of $Tf$ is different from the degree of $Tg$ where $g$ is
arbitrarily close to $f$, then we say that $f$ is a singular value for
$T$. We are interested in transformations for which $(x+a)^n$ is a
singular value for infinitely many $n$. It is often the case that the
behavior of such a transformation changes at $a$, and that there is
some sort of symmetry around $x=a$. These symmetries are captured in
the commuting diagrams of the previous section. If $T$ is a linear
transformation that preserves degree, then there can be no
singularities, but there can be symmetries.  We know of the following
examples:

\newcommand{\symmetry}[2]{\item
  \makebox[2in]{#1\hfill}\hspace*{.2cm}\makebox[2in]{#2\hfill}\\[.3cm]}

\begin{enumerate}
\symmetry{$T\colon{}x^k \mapsto \rising{x}{k}\falling{x-\alpha}{n-k}$}%
{$T(x-1)^n = \rising{\alpha}{n}$}
\symmetry{$T\colon{}x^i\mapsto \falling{x+n-i}{n}$}%
{$T(x+1)^n = (-1)^n$}
\symmetry{$x^k \mapsto H_k(x)H_{n-k}(x)$}%
{$T(x-1)^n = 2^{n/2}H_n(0)$}
\symmetry{$x^k \mapsto H_k(x)H_{n-k}(x)$}%
{$T(x-\imath)^n = 2^n (1+\imath)^n$}
\symmetry{$x^k \mapsto L_k(x)L_{n-k}(x)$}%
{$T(x-1)^n =  0  \text{ if $n$ is  odd}$}
\symmetry{$T\colon{}x^k\mapsto H_k(x)x^{n-k}$}%
{$T\,(x-2)^n = H_n(0)$}
\symmetry{ $T\colon{}x^k \mapsto (x)_k(x)_{n-k}$}%
{$T(x-1)^n =  0  \text{ if $n$ is  odd}$}
\end{enumerate}

\section{Appell sequences}
\label{sec:appell}

Some of the properties of special functions and commuting diagrams are
easy for the general class of  Appell sequences.
A sequence of polynomials $\{p_n\}$ is an%
\index{Appell sequence}\emph{Appell sequence} \cite{rota} if the degree of $p_n$
is $n$, and 
\begin{equation}
  \label{eqn:appell-seq-1}
  p_n'(x) = n p_{n-1}.
\end{equation}
Such sequences are in $1-1$ correspondence with formal power series $g(x)$
\cite{roman} and satisfy
\begin{equation}
  \label{eqn:appell-seq-2}
   g(\diffd)x^n = p_n(x) \text{ for } n=0,1,\dots
\end{equation}

\begin{prop}
  Suppose $\{p_n\}$ is an Appell sequence determined by $g(x)$. 
  Then, $g(x)\in\allpolyf$  iff the linear transformation $x^n\mapsto p_n$
  maps $\allpoly$ to itself.
\end{prop}
\begin{proof}
The linear transformation is simply $f\mapsto g(\diffd)f$.
\end{proof}

If a sequence of polynomials satisfies $p_n = \alpha_n n p_{n-1}$,
then the polynomials $q_n = p_n/(\alpha_0\alpha_1\cdots\alpha_n)$ form
an Appell sequence since 

$$ q_n' = p_n'/(\alpha_0\cdots\alpha_n) =
\alpha_n n p_{n-1}/(\alpha_0\cdots\alpha_n) = 
 n p_{n-1}/(\alpha_0\cdots\alpha_{n-1}) = n q_{n-1}.
$$ 

The Hermite polynomials are nearly Appell sequences.  Since they
satisfy $H_n' = 2n H_{n-1}$, $H_n(x/2)$ is an Appell sequence.  For a
more interesting example, consider the transformation $T\colon{}f\mapsto
f(x+\imag) + f(x-\imag)$.  Now $Tf
= 2\cos(\diffd)f$, and so the sequences
$$ p_n = (x+\imag)^n + (x-\imag)^n$$ are Appell sequences with
$g(x)=\cos x$.

The Euler polynomials  form an Appell sequence with $g(x) =
\frac{1}{2}(e^x+1)$. This function is not in $\allpolyf$, so we can
not conclude that the linear transformation determined by the Euler
polynomials maps $\allpoly$ to itself. Indeed, it doesn't.
\index{Euler polynomials}
Similarly, the \index{Bernoulli polynomials}Bernoulli polynomials have
$g(x) = (e^x-1)/x$, and the corresponding linear transformation
doesn't map $\allpoly$ to itself.

Any Appell sequence satisfies the identity

\begin{equation}
  \label{eqn:appell-seq-3}
   p_n(x+y) = \sum_{k=0}^n \binom{n}{k} p_k(x) y^{n-k}
\end{equation}
which we have already seen for the Hermite polynomials in \eqref{eqn:iden-diag-7}.
Equation \eqref{eqn:appell-seq-3} implies that we have the commutative
diagram

\centerline{\xymatrix{
x^n \ar@{->}[d]_{{x\mapsto x+y}} \ar@{->}[rrrr]^{{x^n\mapsto p_n}} &&&&
  p_n(x) 
\ar@{->}[d]^{{x\mapsto x+y}}\\
(x+y)^n &&  && p_n(x+y) \ar@{<-}[llll]_{{x^n\mapsto p_n}}
}}

If $T\colon{}x^k\mapsto p_k(x)x^{n-k}$ then the lemma below shows that
$T_\ast(x-y)^r$ has degree $n$ for all $y$ other than $y=1$, so $1$ is a
singular point.

\begin{lemma}
  If $\{p_n\}$ is Appell and $T\colon{}x^k\mapsto p_k(x)x^{n-k}$, then 
$$ T_\ast(x-1+y)^r = x^{n-r} p_r(xy)$$
\end{lemma}
\begin{proof}
  Since the conclusion is a polynomial identity, it suffices to prove
  it for $y$ an integer. It's trivially true for $y=1$, and the case
  $y=0$ follows from

$$ T(x-1)^r = \sum_{k=0}^r \binom{r}{k} (-1)^{r-k} p_k(x) x^{r-k} =
x^{n-r} p_n(x-x)=x^{n-r} p_n(0).$$

Assume that it is true for  $y$,
  and consider $y+1$
\begin{gather*} T_\ast(x-1+y+1)^r = \sum_{k=0}^r \binom{r}{k} T_\ast(x-1+y)^k = 
  \sum_{k=0}^r \binom{r}{k} x^{n-k} p_k (xy) = \\
  x^{n-r} \sum_{k=0}^r \binom{r}{k} x^{r-k} p_k (xy) = x^{n-r}
  p_r(xy+x)= x^{n-r} p_r(x(y+1)).
\end{gather*}
\end{proof}

\section{The diamond product}
\label{sec:bilinear-diamond-1}

\index{diamond product!general}
\index{diamond product!original}
\index{product!diamond}

If we are given an invertible linear transformation $T$ on polynomials
then we can form a new bilinear transformation.

\centerline{\
\xymatrix{
\allpoly \times \allpoly \ar@{->}[rr]^{T\times T} \ar@{..>}[d]_{f\mydiamondA{}{}g}&&
  \allpoly\times\allpoly \ar@{->}[d]^{multiplication}\\
\allpoly \ar@{<-}[rr]_{T^{-1}} && \allpoly
}}

We call $f\mydiamond{}{} g = T^{-1}\left(Tf\cdot Tg\right)$ a diamond
product.  The product can be defined for any number of terms in a
simple way. For instance, observe that

$$ f\mydiamond{}{}(g\mydiamond{}{}h) =
f\mydiamond{}{} [ T^{-1}( Tg\cdot Th) ] =
T^{-1}( Tf \cdot T[T^{-1}( Tg\cdot Th) ]) = 
T^{-1}( Tf \cdot Tg \cdot Th)
$$

Consequently, we  see that the diamond product is
well defined for any number of factors, is associative and commutative, and equals

$$ f_1 \mydiamond{}{} f_2 \mydiamond{}{} \cdots \mydiamond{}{} f_n =
 T^{-1}\left(\, Tf_1 \cdot Tf_2  \cdots Tf_n \,\right)
$$

We can now prove a general result:

\begin{prop}\label{prop:diamond-gen}
  Suppose that $T$ is a linear transformation such that 
  \begin{enumerate}
  \item $T\colon{}\allpolyint{\diffi}\longrightarrow\allpolyint{\diffj}$
  \item
    $T^{-1}:\allpolyint{\diffj}(1)\longrightarrow\allpolyint{\diffk}(1)$
  \end{enumerate}
  Then every $f\in\allpolyint{\diffi}$ can be written the form
$$c(x-a_1)\mydiamond{}{}(x-a_2)\mydiamond{}{}\cdots\mydiamond{}{}(x-a_n)$$
where $c$ is a constant and each $x-a_i$ is in $\allpolyint{\diffk}$.
\end{prop} 
\begin{proof} 
From (1) we see that $T(f)\in\allpolyint{\diffj}$, and consequently
can be factored $T(f) = c_0(x-b_1)\cdots(x-b_n)$. Write
$T^{-1}(x-b_i)=c_i(x-a_i)$ where each $x-a_i$ is in
$\allpolyint{\diffk}$ by (2). If $c=c_0\cdots c_n$ then

\begin{align*}
c(x-a_1)\septimes(x-a_2)\septimes\cdots\septimes(x-a_n) & = 
c_0T^{-1}(c_1 T(x-a_1)\cdots c_n T(x-a_n)) \\
& = T^{-1}( c_0 (x-b_1)\cdots (x-b_n)) \\ & = T^{-1}(Tf) = f
\end{align*}

\end{proof} 

Here are some examples of diamond products.

\begin{description}
\item[binomial] Let $T^{-1}(x^n) = \binom{x}{n}$. Then 
\begin{align*}
\binom{x}{r}\mydiamond{} \binom{x}{s} & =
T^{-1}\,(\,T(\binom{x}{r})\,T(\binom{x}{s})\,) \\ &  = T^{-1}(x^{r+s})
\\ & =
\binom{x}{r+s}
\end{align*}

This was considered by Wagner \cite{wagner} who introduced the name diamond
product. See Proposition~\ref{prop:diamond-2}.

\item[falling factorial] Let $T^{-1}\falling{x}{n} = x^n$. Then 
$$ \falling{x}{n}\mydiamond{}{}\falling{x}{m} = \falling{x}{n+m}$$

This product  is discussed in \chapsec{affine}{affine-diamond}.

\item[exponential] If $T(x^n) = x^n/n!$ then
$$ x^n \mydiamond{}{} x^m = \binom{n+m}{n} x^{n+m}$$

See Question~\ref{ques:diamond-exp}.

\end{description}

We can generalize the diamond product construction by replacing
multiplication by some other bilinear map. Given
$m:\allpoly\times\allpoly\longrightarrow\allpoly$ then we form a new
bilinear transformation

$$f\mydiamondA{}{}g=  T^{-1}\left(m(T(f), T(g))\right)$$

See Proposition~\ref{prop:diamond-3} for an application. 

  Although affine transformations distribute over multiplication, they
  do not generally distribute over diamond multiplication. We need a
  connection between the affine transformation and the defining
  transformation.

  \begin{lemma}
    Suppose $f\mydiamond{}g = T^{-1}(Tf\cdot Tg)$ and $\affa$ is an
    affine transformation. If $\affa$ commutes with $T$ then $\affa$
    distributes over $\mydiamond{}$. Conversely, if $\affa$
    distributes over $\mydiamond{}$ and $\affa $ commutes with $T$ for
    polynomials of degree $1$, then $\affa$ commutes with $T$.
  \end{lemma}
  \begin{proof}
    The first part is just definition:
    \begin{gather*}
      \affa(f\mydiamond{}g) = \affa T^{-1}(Tf\cdot Tg) 
= T^{-1}\affa(Tf\cdot Tg) 
= T^{-1}(\affa Tf \cdot \affa Tg) \\
= T^{-1}(T\affa f\cdot T\affa g) 
= \affa f \mydiamond{} \affa g.
\end{gather*}

Conversely, using induction 
\begin{gather*}
  \affa T(fg) = \affa(Tf\mydiamond{}Tg) = \affa Tf\mydiamond{} \affa Tg =
  T\affa f\mydiamond{}T\affa g\quad\text{[induction]}\\
= T(\affa f\cdot \affa g) = T\affa(fg)
\end{gather*}
\end{proof}

\section{Other root preserving transformations}
\label{sec:other-transformations}

There are linear transformations that increase the degree of a
polynomial by more than $1$ and still preserve roots.  Such
transformations $W$ will not preserve interlacing if the degrees of
$Wf$ and $Wg$ differ by more than $1$.

A trivial example is $Tf = p\cdot f$ where $p\in\allpoly$. This
transformation does preserve interlacing. A more
interesting example that does not preserve interlacing  is
$$
T\colon{} x^n \mapsto (x^2-1)^n$$
so that $deg(Tf) = 2 deg(f)$.  This
transformation can also be expressed as $Tf = f(x^2-1)$.  If
$f\in\allpolyint{(-1,\infty)}$ and $\roots(f) = (a_1,\dots,a_n)$ then
$\roots(Tf) = (\dots,\pm(a_i+1)^{1/2},\dots)$.  Thus, $T$ preserves
roots for polynomials in $\allpolyint{(-1,\infty)}$.

We can generalize this idea.

\begin{lemma} \label{lem:sub-roots}
  Suppose that $g\in\allpoly(r)$ has the property that no relative
  maximum or minimum of $g$ has value lying in $(-1,1)$.  If
  $f\in\allpolyint{(-1,1)}(n)$ then
  $f(g(x))\in\allpolyint{(-1,1)}(nr)$.
\end{lemma}
\begin{proof}
  Since $g$ has no extrema with value in $(-1,1)$ it follows that
  every line $y=s$ where $-1\le s\le1$ meets the graph of $g$ in $r$
  distinct points.
  
  Without loss of generality we may assume that $f$ has all distinct
  roots.  If $f(t)=0$ then $-1\le t\le1$ and we get $r$ roots of
  $f(g(x))$ from the $r$ solutions to the equation $g(x)=t$. This
  accounts for all $nr$ of the roots of $f(g(x))$.

\end{proof}

\index{Chebyshev polynomials}

\begin{cor} \label{cor:cheby-2}
  If $f\in\allpolyint{(-1,1)}(n)$ and $T_r$ is the Chebyshev polynomial,
  then $f(T_r(x))$ is in $\allpolyint{(-1,1)}(nr)$.
\end{cor}

\begin{figure}[htbp]
  \begin{center}
    \leavevmode
    \epsfig{file=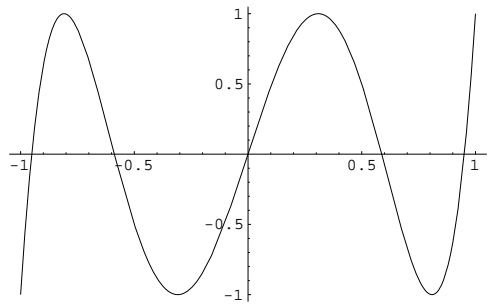,width=2.5in}
    \caption{The graph of the Chebyshev polynomial $ T_5 (x)$}
    \label{fig:cheby}
  \end{center}
\end{figure}

\begin{proof}
The graph of the Chebyshev polynomial $T_r$ (see Figure~\ref{fig:cheby})
on the interval $(-1,1)$ oscillates between its $r-1$ relative extrema 
of $1$ and $-1$. Now apply the lemma.
\end{proof}

For example, since $T_2=2x^2-1$, the transformation $f(x)\mapsto
f(2x^2-1)$ maps $\allpolyint{(-1,1)}$ to itself.

\chapter[Linear transformations that preserve roots ][Root preserving linear transformations]{Linear transformations that preserve roots }
\label{cha:operators}
 
\renewcommand{\TimeStampStart}{Thursday, January 17, 2008: 19:40:22}
\mytoday  

 We apply the results of the previous chapters to show that
particular linear transformations preserve roots.

\section{Multiplier transformations}
\label{sec:multiplier}
\index{multiplier~transformations}

The simplest linear transformations are \emph{multiplier
  transformations}, which are of the form $x^i\mapsto a_i x^i$. We
will see in Theorem~\ref{thm:polya-schur} that if the series $\sum a_i
\frac{x^i}{i!}$ is in $\allpolyposf$ then the multiplier
transformation maps $\allpoly$ to itself.  In
\chapsec{analytic}{hadamard-fac} we showed that various series were in
$\allpolyaltf$ or $\allpolyposf$. Consequently we have

  \begin{theorem} \label{thm:multiplier} The following linear
   transformations map $\allpoly$ to $\allpoly$ and $\allpolypos$ to
    $\allpolypos$. 

    \begin{enumerate}
    \item $\displaystyle x^i \mapsto \frac{x^i}{i!}$
    \item $\displaystyle x^i \mapsto \frac{x^i}{\Gamma(i+\alpha)}$ for
      $\alpha\ge0$
    \item $\displaystyle x^i \mapsto \frac{x^i}{\rising{\alpha}{i}}$ for
      positive $\alpha$. 
    \item $\displaystyle x^i \mapsto \frac{x^i}{(\alpha i)!}$ for
      positive integer $\alpha$.
    \item $\displaystyle x^i \mapsto q^{\binom{i}{2}}{x^i}$ for
      $|q|\le1/4$.
    \item $\displaystyle x^i \mapsto \frac{x^i}{[i]!}$ for
      $|q|>4$. \\[.1cm]
    \item $x^i \mapsto \falling{d}{i}  x^i$ for $d$ a positive integer.\\[.1cm]
    \item $x^i \mapsto \binom{d+i-1}{i} x^i$ for $d$ a positive
      integer.\\[.2cm]
    \item $x^i\mapsto \dfrac{i!q^{\binom{i}{2}}}{(q;q)_i}x^i$ for
      $|q|<1$.\\[.2cm]
    \item $x^i\mapsto \dfrac{q^{\binom{i}{2}}}{(q;q)_i}x^i$ for
      $|q|<1$.\\[.2cm]
    \end{enumerate}
    
  \end{theorem}
  \begin{proof}
    The first two follow from \eqref{eqn:bessel-v}, the third one
    from \eqref{eqn:hypergeo}, and the fourth from
    \cite{polya-szego2}*{62,\#162}.  The generating function of
    $\falling{d}{i}$ is in $\allpolypos$: $$\sum_{i=0}^\infty
    \falling{d}{i} \dfrac{x^i}{i!} = \sum_{i=0}^\infty \binom{d}{i}
    x^i = (1+x)^d.$$
    The generating function for $x^n/[n]!$ is in
    $\allpolyf$ - see \ref{sec:rapid}.  The
    generating function for $\binom{d+i-1}{i}$ is $$
    \sum_{i=0}^\infty
    \binom{d+i-1}{i} \dfrac{x^i}{i!} = e^x L_{d-1}(-x)$$
    where
    $L_{d-1}$ is the Laguerre polynomial (see
    Section~\ref{sec:laguerre}).%
    \index{Laguerre polynomials}\index{polynomials!Laguerre} Since
    $e^x\in\allpolyposf$ and $L_{d-1}(-x)\in\allpolypos$ their product
    is in $\allpolyposf$.

The last two use \eqref{eqn:q-exponential} and part (1).
  \end{proof}
  
  We can extend $(7)$ of the last result to non-integer values of $n$,
  but surprisingly there are restrictions.

\begin{lemma} \label{lem:falling-non-int}
  If $n$ is a positive integer and $ \alpha>n-2 \ge 0$ then 
  \begin{enumerate}
  \item the map $x^i\mapsto \falling{\alpha}{i}x^i$ maps $\allpolypos(n)$ to
    itself, and $\allpolyalt(n)$ to itself.
  \item the map $x^i\mapsto \binom{\alpha}{i}x^i$ maps $\allpolypos(n)$ to
    itself, and $\allpolyalt(n)$ to itself.
  \end{enumerate}
\end{lemma}
\begin{proof}

\index{recurrence!for $x^i\mapsto\falling{\alpha}{i}x^i$}

We find a recurrence and then proceed by induction. First, let $T(x^n) =
\falling{\alpha}{n}x^n$. Then
  \begin{align*}
    T(x\cdot x^i) &= \falling{\alpha}{i}(\alpha-i)x^{i+1} \\
    &= \alpha x\,\falling{\alpha}{i}\,x^i - x\,i\,\falling{\alpha}{i}\,x^i \\
    &= \alpha x \, T(x^i) - x\,T(\,x(x^i)^\prime\,) \\
    T(xf) &= \alpha x\, T(f) - x\,T(\,x f^\prime\,). \\
  \end{align*}
  
  We now show by induction on $m$ that $T$ maps $\allpolypos(m)$ to
  $\allpolypos(m)$ for $m<n$ and to $\allpoly(m)$ for $m=n$. For $m=1$ we
  take $f=x+r$ where $r>0$ and then $Tf = r+\alpha x$ is also in
  $\allpolypos(1)$.
  
  Assume that $T$ maps $\allpolypos(m)$ to itself and $m<n$. If we
  choose $g=(x+r)f$ where $f$ has positive leading coefficient and
  $r>0$ then
  $$
  T(g) = \alpha x T(f) - xT(xf^\prime) + rT(f).$$
  Since
  $f\in\allpolypos$ we know that $xf^\prime\greateqeq f$ and hence by
  induction $T(xf^\prime)\greateqeq T(f)$. It follows that $T(g)=(\alpha
  x+r)T(f) -xT(xf^\prime) \greateqeq T(f)$. If $m=n-1$ this shows that 
  $T$ maps $\allpolypos(n)$ to $\allpoly(n)$.
  
  Now assume that $m<n-1$. In order to show that $T(g)$ is in
  $\allpolypos$ we need to show that $(Tg)(0)$ and the leading
  coefficient of $T(g)$ have the same sign. First, $(Tg)(0) =
  r\,(Tf)(0)$ is positive since $r>0$ and $T(f)$ has positive leading
  coefficient and is in $\allpolypos$. Next, the leading coefficient
  of $T(g)$ is $s\falling{\alpha}{m+1}$ where $s>0$ is the leading
  coefficient of $f$. Consequently, the leading coefficient of $T(g)$
  is positive since $\alpha > n-2 \ge m$ and the conclusion now follows.
  
  The second part follows from the first by applying the exponential
  transformation (Theorem~\ref{thm:multiplier}(1)). The results for $\allpolyalt$
  follow from applying the above to the composition $f(x)\mapsto
  T(f(-x))(-x)$.
\end{proof}

If the generating function lies in $\allpolyf$ but not in
$\allpolyaltf$ or $\allpolyposf$ then the corresponding linear
transformation maps $\allpolypm$ to $\allpoly$ by
Theorem~\ref{thm:polya-schur}.

\begin{lemma} \label{gf-in-f}
The following linear transformations map $\allpolypm$ to
$\allpoly$.

\begin{enumerate}
\item $x^i \mapsto \begin{cases} 0 & i \text{ odd}\\
(-1)^{n} \dfrac{(2n)!}{n!}\,x^{2n} & i = 2n
\end{cases} $\\[.2cm]
\item   $    x^i  \mapsto \begin{cases} 0 & i \text{ even} \\ 
(-1)^nx^i & i = 2n+1 
\end{cases}$ \\[.2cm]
\item $    x^i  \mapsto \begin{cases} 0 & i \text{ odd} \\ 
(-1)^nx^i & i = 2n 
\end{cases} $
\end{enumerate}
\end{lemma}
\begin{proof}
The correspond the the functions $e^{-x^2}$, $\sin(x)$ and $\cos(x)$.  
\index{sin and cos}
\end{proof}

\index{q-exponential}

The $q$-exponential function was defined in Example~\ref{ex:q-exponential}.
We define the $q$-exponential transformation ($0<q<1$)
\begin{equation}
  \label{eqn:q-expo-operator}
  Exp_q(x^i): x^i \mapsto \frac{q^{\binom{i}{2}}}{(q;q)_i}x^i
\end{equation}

Theorem~\ref{thm:multiplier} shows that  $Exp_q$ maps $\allpoly$ to itself.

\section{Homogeneous transformations}
\label{sec:replacement}

Some linear transformations are only defined on polynomials of bounded
degree. For instance, if
$$T(x^i) = x^i(1-x)^{n-i}$$
then $T$ is a linear transformation on all
polynomials of degree at most $n$. 

\begin{lemma}
  If $a$ and $c$ are not both zero then the linear transformation $T$
  given by $$T(x^i) = (ax+b)^i(cx+d)^{n-i}$$
  maps $\allpoly(n)$ to
  itself.
\end{lemma}

\begin{proof}
  Let $F$ be the homogeneous polynomial corresponding to $f$.  If
  $f(x)= (x-a_1)\cdots(x-a_n)$ then $F(x,y) =
  (x-a_1y)\dots(x-a_ny)$.  The image  of $f$ under the
  transformation $T$ is simply $F(ax+b,cx+d)$.  In order to see this, we
  only need to check it for a basis, namely $x^i,$ in which case the
  result is clear.

  We therefore have the factorization $$T(f)=F(ax+b,cx+d) = 
  (ax+b-a_1(cx+d))\cdots(ax+b-a_n(cx+d))$$
  which shows that $Tf$ has all real roots. 
\end{proof}

This lemma also follows from the properties of \Mobius\ 
transformations in \chapsec{linear}{mobius} - the argument here shows
the factorization of the transformation. 
\begin{cor} \label{cor:homog-3}
  The linear transformation $T\colon{}x^i\mapsto x^i(1-x)^{n-i}$ maps
  $\allpolyalt(n)$ to $\allpolyint{(0,1)}(n)$.
\end{cor}
\begin{proof}
  This is  a special case of Lemma~\ref{lem:mobius} where $M=z/(1-z)$
  since $M$ maps $(0,1)$ to $(0,\infty$.
  Alternatively, if $\alpha$ is a root of $f\in\allpolyalt(n)$ then
  the corresponding root of $Tf$ is ${\alpha}/(\alpha+1)$.
\end{proof}

\index{operator!reversal}
\index{reverse!of a polynomial}

The reversal operator is the special case of this construction where
we set $a=0,b=1,c=1,d=0$.  See \eqref{eqn:homog-3}.

\index{reverse!of a polynomial}
\begin{lemma} \label{lem:reverse}
  If the polynomial $f(x)$ has all real
  roots  then $\rev{\,f}$ has all real
  roots.  If $f\in\allpolypos$ (resp. $\allpolyalt$)  then
  $\rev{f}\in\allpolypos$ (resp. $\allpolyalt$).
  
  If $f \greateq g$ and $0$ is a not root of $fg$ then $\rev{g}
  \greateq \rev{f}$.  If $f\lessless g$ and $f\in\allpolypm$ then
  $\rev{f} \lessless \rev{g}$.
\end{lemma}
\begin{proof}
  Without loss of generality we assume that $0$ is not a root.  The
  reverse of $f$ is $x^nf(\frac{1}{x})$, so its roots are the
  reciprocals of the non-zero roots of $f$.  The rest of the
  properties are immediate from this observation.
\end{proof}

If $f$ has $0$ as a root, then the degree of $\rev{\,f}$ is less than
$n$.  If $0$ is not a root of $f$, then  $\rev{(\,\rev{f})}=f$. However,
if $0$ is a $d$-fold root, then $x^d\rev{(\,\rev{f})}=f$

If we have a linear transformation $T$ on $\allpoly(n)$ then we can
conjugate $T$ by the reversal map to create a new linear
transformation $T^{rev}$ on $\allpoly(n)$. \index{reverse!conjugated by}

\centerline{
\xymatrix{
\allpoly(n) \ar@{.>}[rr]^{\rev{T}} \ar@{->}[d]_{\rev{}}&&  \allpoly(n) \ar@{<-}[d]^{\rev{}} \\
\allpoly(n) \ar@{->}[rr]^T&&  \allpoly(n)
}}
Here are some examples:

\begin{enumerate}
\item If $\diffd$ is the derivative, then 
  $\diffd^{\rev{}}$ is the polar derivative.
\item If $T$ is the transformation 
  $x^i\mapsto c_ix^i$ then $\rev{T}(x_i) = c_{n-i}x^i$.
\item If $T(x^i)=p_i$ where $p_i$ is symmetric ($p_i = \rev{(p_i)}$)
  then $\rev{T} (x^i) = p_{n-i}$. This is the case for the Hermite
  polynomials (see \chapsec{operators}{hermite}). 
\end{enumerate}

  \section{The Hadamard product}

\label{sec:hadamard-operators}
\index{Hadamard product}
\index{Hadamard product!generalized}

The generalized Hadamard product of two polynomials is the coordinate wise
product with an extra coefficient per coordinate:
\begin{multline*}
 (a_0 + a_1x + a_2x^2 + \dots + a_nx^n) \hadprod
 (b_0 + b_1x + b_2x^2 + \dots + b_nx^n) \\
=  a_0b_0c_0 + a_1b_1c_1x + a_2b_2c_2x^2 + \dots + a_nb_nc_nx^n
\end{multline*}

Equivalently, we can define it in terms of monomials:

\[
   x^i \hadprod x^j  \mapsto 
  \begin{cases}
    c_i\,x^i & i=j \qquad\text{and}\ c_i>0 \\ 0 & i\ne j
  \end{cases}
\]

With the help of the method of generating 
functions \mypage{prop:hadamard-prod-gen}, we will prove that
\begin{equation}\label{eqn:had-prod-prop}
\text{\parbox{4in}{  The generalized Hadamard product $\hadprod$ determines a map
  $\allpolypos\times\allpoly\mapsto\allpoly$ if and only if
  $\sum_1^\infty c_i\, \frac{x^i}{i!i!}\in\allpolyposf$.}}
\end{equation}

If all $c_i$ are $1$ we call it the Hadamard product and write $f\ast
g$ in place of $f\hadprod g$.  If $c_i=i!$ then we write $f\ast' g$.
It follows from the above result that $\ast$ and $\ast'$ map
$\allpolypos\times\allpoly$ to $\allpoly$.  There is a an alternative
proof for the Hadamard product (see \chapsec{p2}{hadamard}) that uses
an identification of it as a coefficient of a polynomial in two
variables.

The generalized Hadamard preserves interlacing.
\begin{lemma}\label{lem:had-prod-pres}
  If $\hadprod$ satisfies \eqref{eqn:had-prod-prop},
  $f\greateqeq g\in\allpolypos$, $h\in\allpoly$ then
$f\hadprod h \greateqeq g\hadprod h$.
\end{lemma}
\begin{proof}
  If $ g = \alpha f+ r$ where $f\lesslesseq r$ then $f\hadprod h
  \greateqeq \alpha f\hadprod h$ and $f\hadprod h \lesslesseq r
  \hadprod h$ since all coefficients are positive.  Adding these
  interlacings gives the result.

\end{proof}

Our goal is to prove a multiplicative property of $\hadprod$.
We begin with a special case.
  \begin{lemma}
    If $f\in\allpolypos$, $h\in\allpoly$, $\alpha,\beta>0$ and $\hadprod$ satisfies
    \eqref{eqn:had-prod-prop} then
\[
(x+\alpha)\,f\hadprod(x+\beta)\,h \lesslesseq f\hadprod h
\]
  \end{lemma}
  \begin{proof}
    We have the interlacings
    \begin{align*}
      xf  \hadprod xh & \lesslesseq f\hadprod h & \text{since $xf\hadprod xh =
        x(f\hadprod h)$}\\
xf \hadprod \beta h & \lesslesseq f \hadprod h & Lemma~\ref{lem:had-prod-pres} \\
\alpha f \hadprod x h & \lesslesseq f \hadprod h & Lemma~\ref{lem:had-prod-pres} \\
\alpha f \hadprod \beta h & \greateqeq f \hadprod h & 
    \end{align*}
and the conclusion follows upon adding the interlacings.
  \end{proof}

  It's surprising that the generalized Hadamard product allows us to
  multiply interlacings. In the case of $\ast$ there is a natural
  interpretation using polynomials with complex coefficients 
   \mypage{cor:had-prod-i}.

  \begin{lemma}\label{lem:hadamard-interlace}
    If $f\lesslesseq g$ in $\allpoly$,  $h \lesslesseq k$  in $\allpolypos$ and
    $\hadprod$ satisfies \eqref{eqn:had-prod-prop} then
\[
 f\hadprod h \lesslesseq g \hadprod k
\]
  \end{lemma}
  \begin{proof}
    Using Lemma~\ref{lem:sign-quant} we write
\[
g = \sum_i \frac{a_i\,f}{x-r_i}\ \text{and}\   k = \sum_j \frac{b_j\,
  h}{x-s_j}  
\]
where the $r$'s and $s$'s are
    negative, and the $a$'s and $b$'s are positive. From the above
    lemma we know that for all relevant $i,j$
\[
 f\hadprod h \lesslesseq  \frac{a_i\,f}{x-r_i} \,\hadprod \,  \frac{b_j\,h}{x-s_j} 
\]
Adding these interlacings gives the conclusion.
  \end{proof}

\begin{lemma} \label{lem:binom-2}
  If $m$ is a positive integer then 
  \begin{enumerate}
  \item the map $x^i\mapsto \binom{m}{i}x^i$ maps $\allpoly$ to
    itself.
  \item the map $x^i\mapsto \falling{m}{i} x^i$ maps $\allpoly$ to
    itself.
  \end{enumerate}
\end{lemma}
\begin{proof}
  The first statement is the map $f\mapsto (1+x)^m \ast f$. The second statement
  follows from the first statement and Theorem~\ref{thm:hadamard-2}. Another proof
  uses generating functions -see Theorem~\ref{thm:multiplier}. 
\end{proof}

It is obvious that the Hadamard product extends to a bilinear map
$\allpolyf\times\allpolypm \longrightarrow \allpolyf$. Since 
$$ e^{x^2} \ast (x+1)^2 = x^2 + 1$$
it follows that $e^{x^2}\not\in\allpolyf$.  Of course, we know that
$e^{-x^2}\in\allpolyf$. 

The effect of scaling on the Hadamard product is easy to describe.
  \begin{lemma}\label{lem:hp-scaling} Suppose that $\hadprod$ is a
generalized Hadamard product, and choose constants $\alpha,\beta$.  If
$S,T,U$ are regions satisfying
    \begin{align*} \hadprod\colon\allpolyint{S}\times\allpolyint{T}
&\longrightarrow \allpolyint{U}\\ \intertext{then}
\hadprod\colon\allpolyint{\alpha S}\times\allpolyint{\beta T}
&\longrightarrow \allpolyint{\alpha\beta U}\\
    \end{align*}
  \end{lemma}
  \begin{proof} It follows from the definition of generalized Hadamard
product that
\[ f(\alpha x)\hadprod g(\beta x) = (f\hadprod g)(\alpha\beta x)
\] Consideration of the diagram below yields the conclusion.

\centerline{\xymatrix{
\allpolyint{\alpha S}\times \allpolyint{\beta T} 
\ar@{->}[d]_{x\times y\mapsto \alpha x\times \beta y}
\ar@{->}[r]^{\hadprod} &
\ar@{<-}[d]^{x\mapsto x/(\alpha\beta)}
\allpolyint{\alpha\beta U} \\
\allpolyint{S}\times \allpolyint{T} 
\ar@{->}[r]_{\hadprod} & 
\allpolyint{U} 
}}

  \end{proof}

\added{20/1/7}
\begin{remark}
  Differentiating a Hadamard product has no nice properties in
  general. However, we do have the following properties that are
  easily verified, where $m$ is non-negative.
  \begin{align}
\notag 
    \frac{d}{dx}\bigl[ (1+x)^m \hadprod{}x\, f \bigr] &= m
    (1+x)^{m-1}\hadprod{} f\\
\notag 
    \frac{d}{dx}\bigl[ e^x \hadprod{}x\, f \bigr] &= 
    e^x \hadprod{} f\\
\notag 
    \left(\frac{d}{dx}\right)^m \bigl[e^x \hadprod{} x^m\, f\bigr] & =
    e^x\hadprod{} f
  \end{align}

Thus, if $f\in\allpoly$ then we can find polynomials $p_0,p_1,\dots$
so that
\begin{enumerate}
\item $p_0=e^x\hadprod{}f$.
\item all $p_i$ are in $\allpoly$.
\item $(d/dx) p_i = p_{i-1}$ for $i=1,2,\dots.$
\end{enumerate}
\end{remark}

Such polynomials have been called \index{very hyperbolic} \emph{very
  hyperbolic}. Here is a different construction of very hyperbolic
polynomials due to Chebeterov.

\begin{lemma}\label{lem:very-hyperbolic}
  If $f=\sum_0^\infty a_ix^i\in\allpolyf$ and $f_n = \sum_0^n a_i
  \frac{x^{n-i}}{(n-i)!}$ then all $f_n$ have all real roots and $f_n'
  = f_{n-1}$. 
\end{lemma}
\begin{proof}
  Clearly $f_n' = f_{n-1}$; we need to see that $f_n\in\allpoly$. Now
  \[
    (1+x)^n\ast' f =
\sum_0^n \binom{n}{i}\,i!\, a_i\, x^i 
= n! \sum_0^n a_i\, \frac{x^i}{(n-i)!}
\]
is in $\allpoly$, and taking the reverse shows $f_n\in\allpoly$. 
\end{proof}

For more properties of the Hadamard product, see
\chapsec{p2}{hadamard}.

\section{Differential operators}
\label{sec:easy}

The most general differential operators on polynomials that we consider
have the form
\begin{equation}
  \label{eqn:diff-op-1}
  g \mapsto \sum_{i=0}^n f_i(x) g^{(i)}(x)
\end{equation}
where the $f_i(x)$ are polynomials. We first determine some
restrictions on the coefficient functions for operators that map
$\allpoly$ to $\allpoly$. We then look at some general composition
properties of polynomials, and apply these to some particular
operators. In Proposition~\ref{prop:diff-op-poly} we will characterize, modulo some
positivity and degree conditions, those differential operators
\eqref{eqn:diff-op-1} that map $\allpoly$ to $\allpoly$. 

It is useful to view the differential operator 
as being determined by the two variable polynomial $f(x,y) = \sum
f_i(x)y^i$. These polynomials satisfy a substitution condition:

\begin{lemma} \label{lem:diff-op-a}
  If $f(x,y) = \sum f_i(x)y^i$ determines a differential operator
  \eqref{eqn:diff-op-1} that maps $\allpoly$ to itself then
  $f(x,\alpha)\in\allpoly$ for all $\alpha\in\reals$. 
\end{lemma}
\begin{proof}
  If the corresponding differential operator $T$ maps $\allpoly$ to
  $\allpoly$ then $T$ maps $\allpolyf$ to $\allpolyf$. If we apply $T$
to $e^{\alpha x} \in\allpolyf$ we find
$$ \left(\sum_{i=0}^n f_i(x)\diffd^i\right)\, e^{\alpha x} =
\left(\sum_{i=0}^n f_i(x)\alpha^i\right)\, e^{\alpha x} 
$$
Since the right hand side is in $\allpolyf$ we can multiply by
$e^{-\alpha x}$ and remain in $\allpolyf$. The result follows.
\end{proof}

\index{differential operators} \index{operator!differential}

We next consider two 
particular differential operators. Choose a polynomial $h(x) = \sum
a_i x^i$ and define $h(\diffd)g = \sum a_i g^{(i)}$. In this case the
function $f_i(x)$ of \eqref{eqn:diff-op-1} is just the constant $a_i$.
We also define
$$ h(x\diffd)g = \sum_{i=0}^n a_i \left(x\,\frac{d}{dx}\right)^i\,
g(x)$$
If we expand $(x\diffd)^ig$ we see we get a differential operator of
the form \eqref{eqn:diff-op-1}.

\begin{prop} \label{prop:fofd}
  The map $f\times g\mapsto f(\diffd{})g$ defines a bilinear map \\
   $\allpoly\times\allpoly\longrightarrow\allpoly$.
 
\end{prop}
\begin{proof}
It follows from  Corollary~\ref{cor:tf-1} that $(\diffd{}-\alpha)f\in\allpoly$ for all $f$. 
Apply Corollary~\ref{cor:compose-2}.
\end{proof}

\begin{example}
  If we choose $f=x^2-1$ then for any $g\in\allpoly$ it follows  that
  $g-g^{\prime\prime}\in\allpoly$. Moreover, if $g\lesslesseq h$ then
$ g-g^{\prime\prime} \lesslesseq h - h^{\prime\prime}$. 

If  $n$ is a positive integer then the linear transformation $f\mapsto
  \sum_{i=0}^n \binom{n}{r}f^{(r)}$ maps $\allpoly$ to itself - 
  simply choose $f=(1+x)^n$.
\end{example}

Suppose $f=x+1$ and $g=x+2$. Although $f\greateq g$, substituting the
derivative reverses the direction, e.g. $x+1 = f(\diffd)x \lesseq
g(\diffd)x=2x+1$. This is not true in general, but depends on the sign
of the constant terms.

\begin{prop} \label{prop:fofd-direction}
  If $f\greateqeq g$  and   $h\in\allpoly$ and if
\[
 \begin{cases}
f(0)g(0)>0 &\\
f(0)g(0)<0 &\\
\end{cases}
then
\begin{cases}
  f(\diffd)h  \lesseqeq g(\diffd)h &\\
  f(\diffd)h  \greateqeq g(\diffd)h &\\
\end{cases}
\]
\end{prop} 
\begin{proof} 
  Linearity implies that $f(\diffd)h$ and $g(\diffd)h$ interlace, but
  we don't know which direction.  If we write $f(\diffd)h = \beta
  g(\diffd)h + k$ then we need to determine the sign of the leading
  coefficient of $k$, and the signs of the leading coefficients of
  $f(\diffd)h$ and $g(\diffd)h$.  Compute

  \begin{align*}
h(x) &= x^n + a_{n-1}x^{n-1} + \cdots \\
f(x) &= x^m+\cdots + b_1 x + b_0 \\
g(x) &= x^m+\cdots + c_1x + c_0\\
f(\diffd)h &= b_0 x^n + (b_0a_{n-1}+nb_1)x^{n-1} + \cdots \\
g(\diffd)h &= c_0 x^n + (c_0a_{n-1}+nc_1)x^{n-1} + \cdots \\
f(\diffd)h  &= \frac{b_0}{c_0} g(\diffd)h + n \frac{b_1c_0 -
  b_0c_1}{c_0} x^{n-1}+\cdots
  \end{align*}
  Without loss of generality we assume that $f\greateq g$.
  Consequently, we know the numerator ${b_1c_0 -b_0c_1}$ is positive.

  If $b_0c_0>0$ then assume that $b_0>0$ and $c_0>0$. The leading
  coefficients of $f(\diffd)h$ and $g(\diffd)h$ are positive, and the
  leading coefficient of $k$ is also positive. This implies
  $g(\diffd)h\greateqeq f(\diffd)h$. 

  If $b_0c_0<0$ then assume that $b_0>0$ and $c_0<0$. The leading
  coefficient of $f(\diffd)h$ is positive,  the leading coefficient
  of $g(\diffd)h$  is negative, and the
  leading coefficient of $k$ is negative. This implies
  $g(\diffd)h\lesseqeq f(\diffd)h$. 

  The remaining cases are similar.
\end{proof}

\begin{example}
  We know that $f\lesslesseq T(g)$ for all $g\in\allpoly$ implies that
  $T$ is the derivative. However, we can find many different linear
  transformations $S,T$ so that $S(g)\lesslesseq T(g)$ for all
  $g\in\allpoly$. Choose $f_1\lessless f_2$ where $f_2(0)=0$ and
  define
  $$
  S(g) = f_1(\diffd)\,g\quad\quad T(g) = f_2(\diffd)\,g. $$
  The
  proposition above implies that $S(g)$ and $T(g)$ interlace, and
  since $f_2(0)=0$ the degree of $T(g)$ is one less than the degree of
  $S(g)$. Thus, $S(g)\lesslesseq T(g)$.
\index{interlacing!transformations}
\end{example}

\begin{cor}%
  The transformation\footnote{
It is not true that the transformation maps
$\allpolyalt\times\allpolypm$ to $\allpoly$. 
For example, if
$f=(x-.7)(x-.1)$ and $g=(x+.7)(x+.1)$ then $f(x\diffd)g = .0049 +
.216x + 2.47 x^2$ has no real roots.
} \label{cor:fxdg} 
 $f\times g\mapsto f(x\diffd)g$ defines a bilinear 
  map
\begin{enumerate}
\item  $\allpolypos \times \allpolypm \longrightarrow
  \allpolypm$. 
\item   $\allpolyint{I}(n) \times \allpoly\longrightarrow
    \allpoly$ if  $I = (-\infty,0)\cup(n,\infty).$
\item $\allpolypos \times \allpolyint{(0,1)} \longrightarrow
  \allpolyint{(0,1)}.$
\item $\allpolypos \times \allpolyint{(-1,0)} \longrightarrow
  \allpolyint{(-1,0)}.$
\end{enumerate}
\end{cor}
\begin{proof}
  To prove the first part for $f\in\allpolypos$ it suffices to show
  that $f(x\diffd)g\in\allpoly$ when $f= x+a$ and $a$ is positive.
  Since $f(x\diffd)g = ag + xg^\prime$ this follows from
  Corollary~\ref{cor:tf-3}. The second part also follows from the
  corollary.  The third part holds since
  $ag+xg^\prime\in\allpolyint{(0,1)}$ when $a$ is positive and
  $g\in\allpolyint{(0,1)}.$ The fourth part is similar.
\end{proof}

We can extend this to analytic functions:

\begin{cor} \label{cor:fxdg-analytic}
  If $f\in\allpolyf$ has no  roots in $[0,n]$ then $T(g)=f(x\diffd)g$
  defines a linear transformation $T\colon{}\allpoly(n)\longrightarrow\allpoly(n)$.
\end{cor}
\begin{proof}
  The reason this is not immediate is that $f$ is not necessarily
  presented as a limit of polynomials with no roots in $[0,n]$.
  However, if a sequence $f_i$ of polynomials converges uniformly to
  $f$, then the zeros of the $f_i$ converge uniformly to the zeros of
  $f$. If $f$ has no roots in $[0,n]$, then the roots of $f_i$ must go
  to to roots that are not in $[0,n]$ as $i$ goes to infinity. If
  $g\in\allpoly(n)$ then for $i$ sufficiently large no roots of $f_i$
  are in $[0,n]$, and we can apply Corollary~\ref{cor:tf-3a} as above.
\end{proof}

It's easy to express $f(x\diffd)g$ in terms of coefficients.
  If $g(x) = b_0 + \cdots + b_n x^n$ then 
  \begin{equation}
    \label{eqn:fxdg}  
f(x\diffd)g = f(0)b_0  + f(1)b_1x + \cdots f(n)b_nx^n.
  \end{equation}
  
  To see this, note that the conclusion is linear in both $f$ and $g$
  so it suffices to establish it for the case $f = x^r$ and $g=x^n$.
  In this case we see that
$$ f(x\diffd)g = (x\diffd)^r x^n = n^r x^n = f(n) g.$$

\index{Hadamard product}
  We can interpret \eqref{eqn:fxdg} as a Hadamard product, but the
  factor corresponding to $f$ is not in $\allpolyf$.  Given a
  polynomial $f$, we form the infinite series
  $$
  F(x) = f(0) + f(1)x + f(2)x^2 + \cdots $$
  We see from
  \eqref{eqn:fxdg} that $f(x\diffd)g = F\ast g$.  Moreover, we will
  see in \chapsec{operators}{euler} that this series is a
  rational function of the form $h(x)/(1-x)^r$, and so $F$ is not in
  $\allpolyf$. It is interesting that there are series that map
  $\allpoly$ to $\allpoly$ under the Hadamard product, but these
  series are not uniform limits of polynomials. This phenomenon is
  explained in Theorem~\ref{thm:hadamard-2}.

Using Corollary~\ref{cor:fxdg} and \eqref{eqn:fxdg} we get these corollaries.

\begin{cor}
  If $g = b_0 + b_1x + \cdots b_nx^n$ is in $\allpoly$ and $r$ is a
  positive integer then the following polynomials are in $\allpoly$.

$$ b_0 + 1^rb_1x + \cdots + n^r b_n x^n$$
$$ \falling{r}{r}  b_0 + \falling{r+1}{r}b_1x + \cdots + \falling{n+r}{r} b_n x^n$$
\end{cor}
\begin{proof}
  Choose $f = x^r$ in the first case, and  $f =
  (x+1)\cdots(x+r)$ in the second.
\end{proof}


\index{Bernstein polynomials}

We end this section with a discussion of the 
Bernstein polynomials. Given a function, the $n$-th Bernstein polynomial
is defined by
\begin{gather*} B_n(f) = \sum_{k=0}^n f\left(\frac{k}{n}\right) \binom{n}{k}
x^k(1-x)^{n-k}\\
\intertext{and in particular the polynomial corresponding to $x^r$ is}
B_n(x^r) = 
\sum_{k=0}^n \frac{k^r}{n^r} \binom{n}{k} x^k(1-x)^{n-k} \\
\intertext{Note that}
\sum_{k=0}^n \frac{k^r}{n^r} \binom{n}{k} x^k 
= \left(\frac{x\diffd}{n}\right)^r (1+x)^n
\end{gather*}
We therefore have the diagram

  \centerline{\xymatrix{
f(x)
      \ar@{.>}[drrr]_{{ f\mapsto B_n(f)    }}           
      \ar@{->}[rrr]^{{   }}         
      &&&
      f(\frac{x\diffd}{n})(1+x)^n
      \ar@{->}[d]^{{ z\mapsto \frac{z}{1-z}   }} \\        
      &&&
B_n(x)
}}

Now from Corollary~\ref{cor:fxdg} we know that $T\colon{}f\mapsto
f(\frac{x\diffd}{n})(1+x)^n$ determines a map
$\allpolypos\longrightarrow\allpolyint{(-1,0)}$. In addition, since
$f\in\allpolyint{(1,\infty)}$ implies that the roots of $f(x/n)$ are
greater than $n$, we find that $T$ also maps
$\allpolyint{(1,\infty)}$ to $\allpolyint{(-\infty,-1)}$. Moreover,
it follows from the Corollary that we actually have
$\allpolyint{\reals\setminus(0,1)}\longrightarrow
\allpoly$. 
Combined with the \Mobius\ transformation $z\mapsto z/(1-z)$, we get
the commuting diagram of spaces

  \centerline{\xymatrix{
\allpolyint{\reals\setminus(0,1)}
      \ar@{.>}[drrr]_{{ f\mapsto B_n(f)    }}           
      \ar@{->}[rrr]^{{ T  }}         
      &&&
\allpolyint{(-1,0)}
      \ar@{->}[d]^{{ z\mapsto \frac{z}{1-z}   }} \\
&&& 
\allpoly
}}

In summary,

\begin{lemma}
  The transformation $f\mapsto B_n(f)$ satisfies
  \begin{enumerate}
  \item $\allpolyint{\reals\setminus(0,1)}\longrightarrow \allpoly$
  \item $\allpolyint{(-\infty{,0)}} \longrightarrow
    \allpolyint{(-\infty,0)}$
  \item $\allpolyint{(1,-\infty)} \longrightarrow
    \allpolyint{(1,\infty)}$
  \end{enumerate}
\end{lemma}

\section{Transformations based on differentiation}
\label{sec:linear-diff}

 We now look at specific transformations that are based on
 differentiation. The proofs are all easy consequences of the results
 in \chap{polynomials} and \chap{linear}. 

\begin{cor} \label{cor:tf-1}
   If $Tf = (\alpha+\diffd{})f$  then
   $T\colon{}\allpoly\longrightarrow\allpoly$. If  $\alpha>0$ then
   $T\colon{}\allpolypos\longrightarrow\allpolyneg$. If $\alpha<0$ then 
   $T\colon{}\allpolyalt\longrightarrow\allpolyalt$.
\end{cor}

\begin{cor} \label{cor:x2-x}
  If $\alpha,\gamma>0$, $Tf = (\alpha x - \beta - \gamma
  x\diffd{})f, $ then \\ $T\colon{}\allpoly\longrightarrow\allpoly$. If 
  $f\in\allpolyalt$ then $h\lessless f$.  If we know in addition that
  $\beta>0$ then $T\colon{}\allpolyalt\longrightarrow\allpolyalt$.
\end{cor}

\begin{cor} \label{cor:fxd}
  If $a$ and $e$ are positive then the transformation \\  \hbox{$f\times g\mapsto
  f(ax-e\diffd)g$} defines a bilinear map
  $\allpoly\times\allpoly\longrightarrow\allpoly$. 
\end{cor}
\begin{proof}
  Apply Corollary~\ref{cor:x2-x}. 
\end{proof}

\begin{cor} \label{cor:fxxd2}
  The map $f\times g\mapsto f(x+x\diffd)g$ defines a bilinear map 
$\allpolypos\times\allpolyalt\longrightarrow\allpolyalt$.
\end{cor}
\begin{proof}
  Apply Corollary~\ref{cor:x2-x}.
\end{proof}

See Corollary~\ref{cor:tf-3} for the details about the next transformation.

\begin{cor} \label{cor:tf-3b} \label{cor:tf-3a}
  If $a\not\in(0,n)$ and $Tf= (-a+x\diffd)f$ then \\
  $T\colon{}\allpolypm(n)\longrightarrow\allpolypm(n)$ and
  $T\colon{}\allpoly(n)\longrightarrow\allpoly(n)$.
\end{cor}

\begin{cor} \label{cor:x2-b}
    If $Tf = (\beta  +(1- x^2)\diffd{})f $ where $\beta>0$
  and   $f\in\allpolyint{(-1,1)}$ then $Tf\lesslesseq f$ and
  $T\colon{}\allpolyint{(-1,1)}\longrightarrow\allpolyint{(-1,\infty)}$.
\end{cor}
\begin{proof}
  Note that $1-x^2$ is positive on the roots of $f$. Let \\ 
  \hbox{$h=(\beta +(1- x^2)\diffd{})f$.} Since the leading
  coefficient of $h$ is negative it follows from Lemma~\ref{lem:sign-quant} that
  $h\lesslesseq f$. At $x=1$ we see that $h(1) = \beta f(1)$ , so $h$
  has at exactly one root greater than $1$.  There are no roots less
  that $-1$, so $h\in\allpolyint{(-1,\infty)}$.
\end{proof}

\begin{cor} \label{cor:tf-euler}
  If $Tf = axf+x(1-x)f^\prime$ then $T$ maps
  $\allpolyint{(-\infty,0]}$ to itself.
\end{cor}
\begin{proof}
  Apply Lemma~\ref{lem:sign-quant} after factoring out $x$.
\end{proof}

A similar argument shows that

\begin{cor} \label{cor:euler-frob}
If $Tf = ((2x+\alpha)- (1-x^2)\diffd)f$ where $|\alpha|<1$ then $T$
maps $\allpolyint{(-1,1)}$ to itself.   
\end{cor}

\section{The falling factorial}
\label{sec:falling-factorial}
\index{falling factorial}
\index{\ Cbb@$\falling{x}{n}$}

Recall \index{falling factorial} $\falling{x}{n} = x(x-1)\cdots(x-n+1)$.  Since
$\falling{x}{n+1} = (x-n)\falling{x}{n}$ the falling factorials are of Meixner type
with parameters $(a,\alpha,b,\beta)=(0,1,0,0)$. If $T\colon{}x^n \mapsto
\falling{x}{n}$ then it follows from \eqref{eqn:meixner-1} and
\eqref{eqn:meixner-2} that 
  \begin{align}
    T(xf) &= xT(f) - T(xf^\prime) \label{eqn:recur-5} \\
T^{-1}(xf) & = xT^{-1}(f) + xT^{-1}(f)^\prime
\nonumber
  \end{align}

\index{falling factorial}
\index{Meixner class of polynomials} 
\index{polynomials!of Meixner type}

\noindent%
The polynomials $T^{-1}(x^n)$ are also known as the \emph{exponential
polynomials} \cite{rota}*{page 69} or the \emph{Bell polynomials}
\index{exponential polynomials}
\index{Bell polynomials}

\begin{prop} \label{prop:recur-5}
  If $T(x^n) = \falling{x}{n} = x(x-1)\dots(x-n+1)$, then $T$ maps $\allpolyalt$ to
  itself. $T^{-1}$ maps $\allpolypos$ to itself. 
\end{prop}

\begin{proof}
  From \eqref{eqn:recur-5}, we can write $T(xf) = xT(f) - TB(f)$ where
  $B(f) = xf^\prime$.  Since $Bf \lesslesseq f$ for any $f$ with all
  positive roots, Theorem~\ref{thm:interlace} shows that $T$ preserves interlacing
  for polynomials with all positive roots.

  We can also write $T^{-1}(xf) = xT^{-1}(f) + BT^{-1}f$, so that the
  second statement now follows from Theorem~\ref{thm:interlace}.
\end{proof}


\begin{remark}
  We can use the falling factorial to explicate a subtlety of 
  Theorem~\ref{thm:Tff}.  If $T$ is a transformation that preserves
  interlacing, then we can define $S$ by $T(xf) = S(Tf)$, and $Sg
  \lesslesseq g$ for $g=Tf$.  From Theorem~\ref{thm:Tff} we see
  that if $Sg\lesslesseq g$ then $S$ has a certain form, namely $Sg =
  (ax+b)g + (cx^2+dx+e)g^\prime$. This would seem to contradict the
  fact that $T$ is arbitrary.
  
  In the case of the falling factorial, where $T(x^n) = \falling{x}{n}$, we can
  show that $S(x^n) = x(x-1)^n$.  In particular, $S(x^n)$ is
  \emph{not} interlaced by $x^n$.  Moreover, $S$ does not have the
  simple form given above.  The reason is that we have ignored the
  domains of definition of these results.  The explicit form of $Sg$
  was determined by substituting $x^n$ for $g$.  However, we only know
  that $Sg\lesslesseq g$ when $g=Tf$, and $f\in\allpolyalt$, yet
  $T^{-1}x^n$ is a polynomial in $\allpolypos$.
  
  To establish that $S(x^n) = x(x-1)^n$, define $L(f(x)) = f(x-1)$.
  The map $S$ is defined by \index{falling factorial} $S(\falling{x}{n}) =
  \falling{x}{n+1}$, and this can be rewritten as $S(\falling{x}{n}) =
  x L( \falling{x}{n}))$.
  By linearity, $Sf = xLf$, and substituting $x^n=f$ shows that
  $S(x^n) = x(x-1)^n$.
\end{remark}

\index{falling factorial}

The falling factorial transformation is connected to a series identity.

\begin{lemma} \label{lem:xnfi}
  If $T(x^n) = \falling{x}{n}$, then for any polynomial $f$ we have
$$ \sum_{i=0}^\infty \frac{f(i)}{i!} x^i = e^x\,T^{-1}(f).$$
\end{lemma}

\begin{proof}
  By linearity we may assume that $f(x) = x^n$.  If we let $h_n(x) =
  \sum \frac{i^n}{i!} x^i$, then it is easy to verify that
  $xh^\prime_n(x) = h_{n+1}$.  Since induction and integration by
  parts show that there is a polynomial $g_n$ such that $h_n = g_n
  e^x$, we find that $g_n$ satisfies $g_{n+1} = x(g_n + g^\prime_n)$.

  We can now show by induction that $g_n = T^{-1}(x^n)$.  It is clearly
  true for $n=0$, and the recursion \eqref{eqn:recur-5} shows that $T^{-1}(x^{n+1}) =
  x(g_n + g_n^\prime) = g_{n+1}$.

  This is an ad hoc argument. See \cite{rota}*{page 70} for a proof
  based on the general principles of the operator calculus.
\end{proof}

From \eqref{eqn:meixner-2} we can express $T^{-1}$ as a composition:
\index{composition!falling factorial}
\begin{equation}\label{eqn:falling-1}
T^{-1}(f) = f(x+x\diffd{})(1).
\end{equation}
We could have applied Corollary~\ref{cor:fxxd2} to deduce that $T^{-1}$ preserves
interlacing.

\begin{cor} \label{cor:fi}
  If a polynomial $f$ has all negative roots, then  for all
  positive  $x$ 
 $$f(x)^2\ge f(x-1)f(x+1).$$ 
\end{cor}

\begin{proof}
  First, assume that $x$ is a positive integer.  Using
  Lemma~\ref{lem:xnfi} we see if \index{falling factorial}
  $T(x^n)=\falling{x}{n}$ and $f\in\allpolypos$ then
  $T^{-1}f\in\allpolyneg$, and so $e^x\,T^{-1}(f)$ is in
  $\allpolyposf$. The inequality is Newton's inequality
  (Theorem~\ref{thm:newton}) for $\allpolyf$.

  In general, let $\alpha = r+\epsilon$ where $r$ is an integer and
  $0\le \epsilon<1$. If $g(x) = f(x+\epsilon)$ then the desired
  inequality follows by applying the previous paragraph to $g(x)$
  since $g(r) = f(\alpha)$.
\end{proof}

\section{The rising factorial}
\label{sec:rising-factorial}

\index{rising factorial}
\index{\ Caa@$\rising{x}{n}$}

Recall $\rising{x}{n} = x(x+1)\cdots(x+n-1)$. If
$T\colon{}x^n \mapsto \rising{x}{n}$ then
  \begin{align}
    T(xf) &= xT(f) + T(xf^\prime) \label{eqn:recur-5a} \\
    T^{-1}(xf) & = xT^{-1}(f) - xT^{-1}(f)^\prime \nonumber
  \end{align}

\begin{cor} \label{cor:rising-factorial}
  If $T(x^n) = x(x+1)\dots(x+n-1)$, then $T^{-1}$ maps $\allpolyalt$
  to itself. $T$ maps $\allpolypos$ to itself. 
\end{cor}

\begin{proof}
  The proofs of the first two assertions are like the proofs for the
  falling factorial, and uses the recurrences of \eqref{eqn:recur-5a}.


\end{proof}

\index{recurrence!rising factorial}

The rising factorial satisfies a different recurrence than the falling
factorial, and it leads to a composition formula. Let\footnote{see
  \chapsec{affine}{affine-difference} for the motivation.} $\affa f(x)
= f(x+1)$. Then
\begin{align*}
  \rising{x}{n+1} &= x(x+1) \cdots (x+n) \\
&= x\cdot \affa (x\cdots (x+n-1)) \\
&= x \affa\rising{x}{n}
\end{align*}
 Consequently, if we define $T(x^n)=\rising{x}{n}$  then
$$ T(x^n) = \rising{x}{n} = (x\affa)\rising{x}{n-1} = \cdots = (x\affa)^n(1)$$
and so we have the useful composition formula
\begin{equation}
  \label{eqn:rising-affine}
T(f) = f(x\affa)(1)
\end{equation}

For more information about the rising factorial, see
\chapsec{affine}{affine-difference}.

Sometimes the rising factorial transformation is hidden. The following
lemma is based on the definition of $L$ found in \cite{al-salam}.

\begin{lemma}
  If $g=\sum a_ix^i$ then define $L(f,g)=\sum x^i a_i \Delta^if(0)$
where $\Delta g(x)=g(x+1)-g(x)$. 

If $f\in\allpolypos$ and $g\in\allpoly$ then $L(f,g)\in\allpoly$.
\end{lemma}

\begin{proof}
  Write $f=\sum b_i \rising{x}{i}$. Since $\Delta\rising{x}{i} =
  i\rising{x}{i-1}$ and $\Delta^i\rising{x}{i} = i!$ we see that $b_i
  = \frac{(\Delta^if)(0)}{i!}$ and hence $L(f,g) = \sum x^i a_ib_ii!$.
  Define $T\colon{}\rising{x}{i}\mapsto x^i$ and set $h=T(f)=\sum b_i x^i$. If
  $h\in\allpolypos$ and $g\in\allpoly$ then by Theorem~\ref{thm:hadamard-2} $\sum
  x^ia_ib_ii!\in\allpoly$. From Corollary~\ref{cor:rising-factorial} we know
  $T\colon{}\allpolypos\longrightarrow\allpolyneg$, so $f\in\allpolyneg$
  implies $h\in\allpolypos$. Thus we are done.
\end{proof}

The maps $x^n\mapsto \rising{x}{n}$ and $x^n\mapsto \falling{x}{n}$ are
conjugates by the map $x\mapsto -x$, as expressed in the diagram

\centerline{
\xymatrix{
\allpolyalt \ar@{|->}[rr]^{x^n\mapsto\rising{x}{n}} \ar@{|->}[d]_{x\mapsto-x}
&& \allpolyalt  \\
\allpolypos \ar@{|->}[rr]^{x^n\mapsto\falling{x}{n}} && \allpolyneg \ar@{|->}[u]_{x\mapsto-x}
}}
\noindent%
The diagram commutes since 
$$ x^n \mapsto (-1)^nx^n \mapsto (-1)^n\falling{x}{n}\mapsto \falling{-1}{n}(-x)^n = \rising{x}{n}.$$

If we compose the falling and rising factorial transformations we get
\begin{lemma} \label{lem:fall-rise} \ 
  \begin{enumerate}
  \item The transformation $\falling{x}{n}\mapsto \rising{x}{n}$ maps
    $\allpolypos\longrightarrow\allpolyneg$. 
  \item The transformation $\rising{x}{n}\mapsto \falling{x}{n}$ maps
    $\allpolyalt\longrightarrow\allpolyalt$. 
  \item The transformation $\falling{x}{k}\mapsto \rising{x}{n-k}$ maps
    $\allpolypos(n)\longrightarrow\allpolyneg(n)$.
  \item The transformation $x^k\mapsto \rising{x}{n-k}$ maps
    $\allpolyalt(n)\longrightarrow\allpolyalt(n)$.
  \item The transformation $\falling{x}{k}\mapsto \falling{\alpha}{k}x^k$ maps
    $\allpolypos(n)\longrightarrow\allpolypos(n)$ for $\alpha>n-2$.
  \end{enumerate}
\end{lemma}
\begin{proof}
  The first one follows from the diagram

\centerline{
\xymatrix{
\allpolypos \ar@{.>}[rr]^{(x^n)\mapsto\rising{x}{n}} 
\ar@{->}[d]_{\falling{x}{n}\mapsto x^n} && 
\allpolypos \\
 \allpolypos \ar@{->}[rru]_{x^n\mapsto \rising{x}{n}}
}}
\noindent%
and the second is similar. The next follows from the diagram

\centerline{
\xymatrix{
\allpolypos(n) \ar@{.>}[rr]^{\falling{x}{k}\mapsto\rising{x}{n-k}} &&
\allpolypos(n)\\
\allpolypos(n) \ar@{->}[rr]^{reverse} \ar@{<-}[u]^{\falling{x}{k}\mapsto x^k} &&
\allpolypos(n) \ar@{->}[u]_{x^k\mapsto\rising{x}{k}}
}}
\noindent%
The next two use the diagrams

\centerline{ \xymatrix{ \allpolyalt(n) \ar@{->}[d]_{reverse}
    \ar@{.>}[rr]^{x^i\mapsto\rising{x}{n-i}} && \allpolyalt(n) &&
    \allpolypos(n) \ar@{->}[d]_{\falling{x}{k}\mapsto x^k} 
    \ar@{.>}[rr]^{\falling{x}{k}\mapsto\falling{\alpha}{k}x^k}  && \allpolypos(n) \\
    \allpolyalt(n) \ar@{->}[urr]_{x^k\mapsto\rising{x}{k}} &&& &
    \allpolypos(n) \ar@{->}[urr]_{x^k\mapsto \falling{\alpha}{k}x^k}
}}

\end{proof}

Given any linear transformation $T$, the map $T(x^i)\mapsto
T(x^{i+1})$ usually does not preserve roots. The only non-trivial
cases known are when $T(x^n)=\falling{x}{n}$, $\rising{x}{n}$, or $H_n$
(Corollary~\ref{cor:hermite}). 

\begin{lemma} \label{lem:increment}
  The linear transformations $\falling{x}{i}\mapsto\falling{x}{i+1}$ and
  $\rising{x}{i}\mapsto\rising{x}{i+1}$ map
  $\allpoly\longrightarrow\allpoly$. 
\end{lemma}
\begin{proof}
  If $T(f)=xf$ and $S(f)=f(x+1)$ then the fact that the diagram 

\centerline{
\xymatrix{
{\falling{x-1}{i}} \ar@{|->}[r]^S \ar@{|->}[dr]_T & {\falling{x}{i}} \ar@{|->}[d] \\
& {\falling{x}{i+1}}
}}
\noindent%
commutes shows that the desired transformation is $TS^{-1}$. Now both
$S^{-1}$ and $T$ map $\allpoly$ to itself, and so the same is true for
$TS^{-1}$. Incidentally, this shows $T(x^n) = x(x-1)^n$. The proof for
$\rising{x}{i}$ is similar.

\end{proof}

\begin{lemma}\label{lem:rising-fact} 
Let $T(x^n) = \rising{x}{n}/n!.$
  \begin{align*}
T&: \allpolypos  \longrightarrow\allpolypos\\
T&: \allpolyint{(1,\infty)}  \longrightarrow\allpolyint{(1,\infty)}.     
  \end{align*}

\end{lemma}
\begin{proof}
  We have the identity $$(T(f))(1-x) = T(f(1-x))$$
which  is proved by taking  $f(x)=x^k$, and using

$$ \dfrac{\rising{1-x}{k}}{k!} = 
\sum_{i=0}^k \dfrac{\rising{x}{i}}{i!} \binom{k}{i}
$$

Now we know that $x^n\mapsto \rising{x}{n}$ maps
$\allpolypos\longrightarrow\allpolypos$, and upon application of the
exponential map, so does $T$.  
The identity above implies that we have the following commuting
diagram

\centerline{\xymatrix{
\allpolyint{(-\infty,0)} \ar@{<-}[d]_{x\mapsto 1-x} \ar@{->}[rr]^{T}
&& \allpolyint{(-\infty,0)} 
\ar@{->}[d]^{x\mapsto 1-x}  \\
\allpolyint{(1,\infty)} &&  \allpolyint{(1,\infty)} \ar@{<.}[ll]^{T} 
}}
\noindent%

\end{proof}

\begin{remark}
  $T$ above does not $\allpolyint{(0,1)}$ to $\allpoly$. Indeed,
  $T(x-1/2)^2$ has complex roots. It is known \cite{bump} that $T$
  maps $\allpolyint{\imath+\reals}$ to itself. See also 
  Question~\ref{ques:risingxk}. 
\end{remark}

If we need to prove a result about polynomials without all real roots,
then sign interlacing is  useful. Consider

\begin{lemma}\label{lem:num-pos-rts}
  If $g$ has $p$ positive roots, and $\alpha$ is positive, then\\
  $h=(x-\alpha)g-xg'$ has at least $p+1$ positive roots.
\end{lemma}
\begin{proof}
  Suppose that the positive roots of $g$ are $r_1<\cdots<r_p$. The
  sign of $h(r_i)$ is the sign of $-g'(r_i)$ and so  is $(-1)^{p+i+1}$.
  This shows that $h$ has at least $p$ roots, one between each
  positive root of $g$, and one to the right of the largest root of
  $g$. The sign of $h(0)$ is the sign of $-\alpha h(r_1)$. Since
  $\alpha$ is positive, there is a root between $0$ and $r_1$, for a
  total of at $p+1$ roots. 
\end{proof}

  \begin{lemma}\label{lem:pos-rts-rising-inv}
    If $T^{-1}:\rising{x}{n}\mapsto x^n$, and $f\in\allpoly$ has $p$
    positive roots, then $T^{-1}(f)$ has at least $p$ positive roots.
  \end{lemma}
  \begin{proof}
    We prove this by induction on the number $p$ of positive roots of $f$.
    There is nothing to prove if $p=0$.  Let $g(x) = T^{-1}(f)$.  By
    induction we assume that $g$ has $p$ positive roots.  If $h(x) =
    T^{-1}(x-\alpha)f$ where $\alpha$ is positive, then the recurrence
    \eqref{eqn:recur-5a} shows that $h = (x-\alpha)g -
    xg'$. Lemma~\ref{lem:num-pos-rts} shows that $h$ has at least
    $p+1$ positive roots.
  \end{proof}

\section{Hermite polynomials} 
\label{sec:hermite}

The Hermite \index{Hermite polynomials} polynomials are the orthogonal
polynomials for the weight function $e^{-x^2}$ on $(-\infty,\infty)$.
\index{Rodrigues' formula!Hermite}
There is an explicit formula (the Rodrigues' formula) for the Hermite
polynomials \cite{szego}*{page 25}:

\begin{gather} \label{eqn:hermite-4}
 H_n = (-1)^n  e^{x^2} \left(\frac{d}{dx}\right)^n\,
e^{-x^2} \\
\intertext{an expansion in a series}
H_n = \sum_{k=0}^{n/2} \frac{n!}{(n-2k)!}\,(-1)^k (2x)^{n-2k} 
\label{eqn:hermite-series}
\intertext{and a recurrence relation}
\label{eqn:hermite-recur}
H_n = 2x H_{n-1} - 2n H_{n-2}
\end{gather}
  An addition formula for Hermite polynomials is
  \begin{equation}\label{eqn:hermite-addition}
H_n(x+y) = \sum_{k=0}^n \binom{n}{k} H_k(x)(2y)^{n-k}
  \end{equation}

The Hermite polynomials are not of
\index{Meixner class of polynomials} 
\index{polynomials!of Meixner type}
Meixner type since they aren't monic, but if we rescale them to make
them monic then we get a recursion of the form $p_{n+1} = xp_n -
np_{n-1}$ which is Meixner with parameters
$(a,\alpha,b,\beta)=(0,0,1,0)$. From the recurrence
\eqref{eqn:hermite-recur} it is easy to verify that
$H_n^\prime=2nH_{n-1}$.  There is also a differential recurrence
\begin{equation}\label{eqn:hermite-5}
 H_{n+1} = 2xH_n - (H_n)^\prime = (2x-\diffd{})H_n
\end{equation}
and consequently
\begin{equation}
H_n = (2x-\diffd{})^n(1) \label{eqn:hermite-3}
\end{equation}

If we set $T(x^n) = H_n$ then we have the composition
\index{composition!Hermite}
\begin{align}
 T(f) = f(2x-\diffd{})(1) \label{eqn:hermite-1} \\
\intertext{It is easy to verify that}
  T^{-1}(f) &= f(x/2+\diffd{})(1) \label{eqn:hermite-2}
\end{align}

We need to be careful here because multiplication by $x$ does not
commute with $\diffd$. We define $(2x-\diffd)^n(1)$ inductively:
\begin{align*}
  (2x-\diffd)^n(1) &=
  (2x-\diffd)\,\left((2x-\diffd)^{n-1}(1)\right)\\
\intertext{and for $f=\sum a_ix^i$ we use linearity:}
   f(2x-\diffd)(1) &= \sum a_i \,(2x-\diffd)^i(1)
\end{align*}

In this light the following lemma establishes an interesting relation
between these two definitions. The proof is by induction and omitted.
\begin{lemma} \label{lem:xandD}
  If $f(x) = (ax+b\,\diffd)^n(1)$ then
$$ f(cx+d\,\diffd)(1) = \left(acx + (ad+{b}/{c})\diffd\right)^n(1)$$
\end{lemma}

The linear recurrences satisfied by the transformation $T(x^n)=H_n$
follow from \eqref{eqn:hermite-recur}, and are
\begin{align}
  T(xf) & = 2xT(f) - T(f^\prime) \label{eqn:recur-7} \\
T^{-1}(xf) & = (1/2) xT^{-1}(f) + (T^{-1}f)^\prime
\label{eqn:recur-7inv} 
\end{align}

\index{Hermite polynomials}
\index{orthogonal polynomials!Hermite}
\begin{cor} \label{cor:hermite}
  Let $H_n$ be the Hermite polynomial. 

  \begin{enumerate}
  \item If $T(x^n) = H_n$ then $T\colon{}\allpoly\longrightarrow\allpoly$. 
  \item The linear  transformation $x^i\mapsto H_{n-i}$ maps
    $\allpoly(n)$ to itself.
  \item The linear transformation $H_i\mapsto H_{i+1}$ maps
  $\allpoly(n)$ to itself.
\end{enumerate}
\end{cor}

\begin{proof}
  For the first part apply Corollary~\ref{cor:fxd} to
  \eqref{eqn:hermite-1}.  The next map is the conjugate of $T$ by the
  reversal map -- see \chapsec{operators}{replacement}.
  \index{reverse!conjugated by} The map $H_i\mapsto H_{i+1}$ is
  equivalent to $f\mapsto (2x-\diffd)f$ by \eqref{eqn:hermite-5}, and
  this maps $\allpoly$ to itself by Lemma~\ref{lem:sign-quant}. As an
  aside, the map $H_n\rightarrow H_{n+1}$ evaluated at $x^n$ is
  $x^{n-1}(2x^2-n)$.
\end{proof}

There are also generalized Hermite \index{Hermite polynomials}
polynomials \cite{roman}*{page 87} that are defined for any positive
$\nu$.  They satisfy $H_0^\nu=1$ and
$$
H_{n+1}^\nu = xH_n^\nu - \nu nH_{n-1}^\nu$$
The corresponding linear transformation preserves roots just as the
ordinary Hermite transformation does.

Different representations of orthogonal polynomials can lead to
different proofs of Corollary~\ref{cor:hermite}. For instance, using
\eqref{eqn:hermite-series} the Hermite 
\index{Hermite polynomials} polynomials can be defined in terms of the
derivative operator. 
\cite{fernandez98}
\begin{equation}\label{eqn:hermite-quad-diff}
H_n = 2^n exp(-\diffd{}^2/4) (x^n)
\end{equation}
The transformation $T\colon{}x^n\longrightarrow H_n$ can be factored as $T_1
T_2$ where $T_1(x^n) = (2x)^n$ and $T_2(x_n) =
exp(-\diffd{}^2/4)(x^n)$. Since $exp(-x^2/4)$ is in $\allpolyf$, we
can apply Theorem~\ref{thm:polya-schur-2}, and conclude that $T_2(f)\in\allpoly$.
Thus, $T$ maps $\allpoly$ to itself.

The action of $T\colon{}x^n\longrightarrow H_n(x)$ on $\sin$ and $\cos$ is
surprisingly simple.  It's easiest to compute $T$ on $e^{\imag x}$. 
\begin{align*}
  T(e^{\imag x}) &= \sum_{k=0}^\infty H_k(x) \frac{\imag^k}{k!} \\
&= e^{-2\imag x+1} \quad\quad\text{(from Table~\ref{tab:gen-fct-1})}\\
&= e(\cos(2x) - \imag \sin(2x)) \\
&= T(\cos(x)+\imag\sin(x))
\intertext{and so we conclude that }
T(\cos(x)) &= e\cos(2x) \\
T(\sin(x)) &= -e\sin(2x)
\end{align*}

\section{Charlier polynomials}
\label{sec:charlier}

\index{Charlier polynomials}
\index{orthogonal polynomials!Charlier}

\index{recurrence!monic Charlier}

The monic Charlier polynomials with parameter $\alpha$ are orthogonal
polynomials that satisfy the recurrence formula

\begin{equation}
  \label{eqn:charlier}
  C_{n+1}^{\alpha}(x) = (x-n-\alpha)C_{n}^{\alpha}(x) - \alpha n
  C_{n-1}^{\alpha}(x)
\end{equation}
 
\index{Meixner class of polynomials} 
\index{polynomials!of Meixner type}

The first five Charlier polynomials are in Appendix~\ref{cha:tables}.
Their recurrence shows them to be of Meixner type, with parameters
$(a,\alpha,b,\beta) = (\alpha,1,\alpha,0)$. They satisfy the
difference equation (analogous to the differential equation for the
Hermite \index{Hermite polynomials} \index{difference equation}
polynomials)
\begin{equation}\label{eqn:charlier-delta}
 \Delta C_{n}^{\alpha}(x) = nC_{n-1}^{\alpha}(x) 
\end{equation}
\noindent%
where $\Delta(f) = f(x+1)-f(x).$ Consequently, the Charlier
polynomials give an infinite sequence of polynomials with all real
roots that is closed under differences, just as the infinite sequence
of Hermite polynomials is closed under differentiation. 

The corresponding linear transformation is $T(x^n) = C^\alpha_n(x)$.
$T$ satisfies the recurrences 
\index{composition!Charlier}
\begin{align} 
  T(xf) & = (x-\alpha)Tf - T(xf^\prime) - \alpha T(f^\prime)
  \label{eqn:recur-8}\\
 T^{-1}(xf) & = (x+\alpha)T^{-1}(f) + x(T^{-1}f)^\prime + \alpha
 (T^{-1}f)^\prime \nonumber \\
T^{-1}(f) &= f(\,(x+\alpha)(\diffd{}+1)\,)(1)\nonumber
\end{align}

\begin{cor} \label{cor:charlier-x}
  If $T(x^n) = C^\alpha_n(x)$, the $n$-th Charlier polynomial, then $T$
  maps $\allpolyalt$ to itself.   $T^{-1}$ maps $\allpolypos$ to
  itself.
\end{cor}

\begin{proof}
  If we set $A(f) = xf^\prime$ and $\diffd{}(f) = f^\prime$, then the
  Charlier polynomials satisfy
\begin{align*}
  T(xf) & = (x-\alpha)T(f) - TA(f) - \alpha T\diffd{}(f) \\
  T^{-1}(xf) & = (x+\alpha)T^{-1}(f) + AT^{-1}(f) + \alpha \diffd{}T^{-1}(f) 
\end{align*}
Since $A(f) \greateq f$ and $f \lessless \diffd{}f$ the results follow
from Theorem~\ref{thm:interlace}.  
\end{proof}

Since $T^{-1}$ maps $\allpolypos$ to itself, we know that
$B_{n+1}\plesslesseq B_n$. Using the identity
\[ B_{n+1} = x(B_n + B_n')\]
we see that actually $B_{n+1}\lesslesseq B_n$. \seepage{ques:bell-polys}

\begin{remark} \label{rem:charlier-comp}
  The Charlier polynomials satisfy a formula that is analogous to the
  composition formula for Hermite polynomials. This will be used in
  \chapsec{affine}{affine-linear}.  If we define
  $\affa^{-1}\,f(x)=f(x-1)$, then the transformation $T$ above
  satisfies
  \begin{align}
    T(x^n) &= C_n^\alpha(x) \notag \\
    &= (x\affa^{-1}-\alpha)C_{n-1}^\alpha(x) \quad\text{by
      \eqref{eqn:charlier} and \eqref{eqn:charlier-delta}} \notag \\
    &= (x\affa^{-1}-\alpha)^n (1) \quad \text{by induction} \notag \\
\intertext{and so}
T(f) &= f(x\affa^{-1}-\alpha)(1) \label{eqn:charlier-sa}
  \end{align}
\end{remark}
\index{recurrence!Charlier polynomials}
\index{Bell polynomials}
\index{polynomials!Bell}
\index{Charlier polynomials}
\index{orthogonal polynomials!Charlier}

\section{Laguerre polynomials}
\label{sec:laguerre}

As with so many special functions there are variations in the
definition of the Laguerre polynomials.  The usual notation is
$L_n^\alpha(x)$; we will see that the $\alpha$ is quite important, and
so we also write $L_n(x;y) = L_n^y(-x)$. The definition of that we use
is in \cite{grads},  but we set the leading coefficient of
$L_n(x;y)$ to $1/n!$. The Rodrigures formula is

\index{Rodrigues' formula!Laguerre}
\begin{gather}
  \label{eqn:laguerre-rod}
  L_n^a(x) = \frac{1}{x^a e^{-x}} \frac{1}{n!} \diffd^n ( e^{-x} x^{n+a})\\
  L_n(x;y) = \frac{1}{x^y e^{x}} \frac{1}{n!} \diffd^n (
  e^{x} x^{n+y})\notag \\
\intertext{and a series expansion is  }
 \label{eqn:laguerre-gen}
  L_n(x;y) = \sum_{k=0}^n \binom{n+y}{n-k} \frac{x^k}{k!}
\end{gather}

The Hermite polynomials \eqref{eqn:hermite-quad-diff} and the Laguerre
polynomials have  nice representations in terms of the exponential of
a quadratic differential.

  \begin{equation}
    \label{eqn:laguerre-quad-diff}
    e^{\alpha \partial_x\partial_y}\,(x^ny^n) = n!\, \alpha^n\,
  L_n(\frac{-xy}{\alpha})
  \end{equation}

  The  homogenized Laguerre polynomials have an expansion
  \begin{gather}
    \label{eqn:laguerre-addition}
 (y+1)^nL_n(\frac{x}{y+1}) =    \sum_{k=0}^n L_k(x) \binom{n}{k}y^{n-k} 
\end{gather}

The expansion \eqref{eqn:laguerre-gen}
shows that $L_n(x;y)$ is a polynomial in $x$ and $y$ with 
positive coefficients.  

\index{Rodrigues' formula!Laguerre}
\index{Laguerre polynomials!Rodrigues' formula}

\index{Laguerre polynomials!monic}
An alternative definition is found in
\cite{roman} with leading coefficient $(-1)^n$. We denote these
alternative polynomials by $\tilde{L}^\alpha_n(x)$, and note that we have the
simple relationship $\tilde{L}^\alpha_n(x) = n! L_n^\alpha(x)$.
\index{Laguerre polynomials}
\index{polynomials!Laguerre}
The three term recurrence is not useful for us, but we will use the
composition formula \cite{roman}*{page 110}
\index{composition!Laguerre} \index{rising factorial}
\index{falling factorial}

\begin{align}
  \tilde{L}_n^{(\alpha)}(x) &= (x\diffd-x+\alpha+1)\tilde{L}_n^{\alpha+1}\notag\\
&= (x\diffd-x+\alpha+1)(x\diffd-x+\alpha+2)\cdots(x\diffd-x+\alpha+n)(1)\notag\\
& = \rising{x\diffd{} - x +
    \alpha+1}{n}(1).  \label{eqn:laguerre-carlitz} 
\end{align}

\begin{cor} \label{cor:laguerre-x}
  If $T$ is the Laguerre transformation $T(x^n) =
  \tilde{L}_n^{(\alpha)}(x)$ where $\alpha\ge-1$ then
  $T\colon{}\allpolypos\longrightarrow\allpolyalt$.
\end{cor}

\index{rising factorial}
\begin{proof}
  If we define $S(f) = (x\diffd{}-x+\alpha+1)(f)$ then from
  \eqref{eqn:laguerre-carlitz} the composition of maps
  $$
  x^n \mapsto \rising{x}{n} \mapsto \rising{S}{n} (1)$$
  is the
  Laguerre transformation.  If $U(x^n) = \rising{x}{n}$ and $V(f) =
  f(S)(1)$, then $T = VU$.  From Corollary~\ref{cor:rising-factorial} $U$ maps
  $\allpolypos$ to $\allpolypos$. By Corollary~\ref{cor:x2-x} $V$ maps
  $\allpolypos$ to $\allpoly$.  Consequently, $VU$ maps $\allpolypos$
  to $\allpoly$. Since the coefficients of $T(f)$ are alternating in sign, we see
  that $T$ maps to $\allpolyalt$. 
\end{proof}

If we define $T(x^n) = \tilde{L}_n^\alpha(-x)$, and $S(x^n) =
L_n(x;\alpha)$ then $S= T \circ \expoper{}.$ It follows that $S$ maps
$\allpolypos$ to itself. We will see in Lemma~\ref{lem:laguerre-2} that $S$
actually maps $\allpoly$ to itself.

The Laguerre polynomials are also related to the linear transformation 
$$R: f \mapsto \diffd^n(x^nf)$$
where $n$ is a fixed positive integer,
whose action on monomials is $R(x^i) = \falling{n+i}{i}\, x^i$. The exponential
generating function corresponding to $R$ is

\index{exponential generating function}

\begin{equation}
  \label{eqn:laguerre}
  \sum_{i=0}^\infty \,\frac{R(x^i)}{i!} \ = \ 
  \sum_{i=0}^\infty \frac{(n+i)!}{i!\,i!}\,x^i \ =\  e^x\, \tilde{L}_n(-x).
\end{equation}
Since both $e^x$ and $\tilde{L}_n(-x)$ are in $\allpolyposf$, so is
their product. The fact that the generating function of $R$ is in
$\allpolyposf$ also follows from Theorem~\ref{thm:polya-schur} since $R$ is a
multiplier transformation that maps $\allpolypos$ to itself.
\index{multiplier~transformations}

\section{The integration and exponential transformations}
\label{sec:integration-exponential}

\index{exponential transformation}
\index{integral transformation}
\index{transformation!integral}
\index{transformation!exponential}
\index{\ aaaint@$\intoper{}$}
\index{\ aaaexp@$\expoper{}$}

The exponential transformation $x^n\mapsto x^n/n!$ is a fundamental
linear transformation that maps $\allpoly$ to $\allpoly$.  Although
not all polynomials in $\allpoly$ have integrals that are in
$\allpoly$ (see Example~\ref{ex:no-integral}) we can find some polynomials for
which \emph{all} integrals have a root.  These polynomials are exactly
the image of the exponential transformation. We first define the
integral ($\intoper{}$) and exponential ($\expoper$) transformations.
Define
\begin{alignat*}{2}
  \exp:& \ x^n & \mapsto & x^n/n! \\
\intoper{}: &\ x^n & \mapsto & x^{n+1}/{n+1} \\
\intoper{k}: &\ x^n & \mapsto & x^{n+k}/\rising{n+1}{k} \\
\intoper{k}\exp: &\ x^n & \mapsto & x^{n+k}/(n+k)!
\end{alignat*}

Since $e^x\in\allpolyf$, we can apply Theorem~\ref{thm:polya-schur} to conclude
that $\expoper{}$ maps $\allpoly$ to itself. The integral doesn't
preserve all real roots as Example~\ref{ex:no-integral} shows, but we
do have

\begin{lemma} \label{lem:integral}
  If $f$ is a polynomial with all real roots, and  for every $k\ge0$ 
  we define $T(f) = \intexp{k}\,f$ then 
  \begin{enumerate}
  \item $T$ maps $\allpoly$ to $\allpoly$.
  \item $T$ maps $\allpolypos$ to $\allpolyneg$.
  \item $T$ maps $\allpolyalt$ to $\allpolyalt$.
  \item $T$ preserves  interlacing.
  \end{enumerate}
\end{lemma}

\begin{proof}
  It follows from the definition of $\intexp{k}$ that $$\intexp{k}f =
  e^x\ast (x^kf)$$ \index{Hadamard product}
  where ``$\ast$'' is the Hadamard product. The
  results now follow from the fact that $e^x\in\allpolyposf$ and
  Hadamard products preserve interlacing.
\end{proof}
\index{Hadamard product!and exponential}

When we apply various exponential transformations, we get some interesting
transformations that we  express in terms of coefficients:

\begin{cor} \label{cor:exp}
  If $f=\sum_{i=0}^n a_ix^i$ is in $\allpoly$ then the following
  polynomials also have all real roots:
  \begin{enumerate}
  \item $\displaystyle\sum_{i=0}^n \frac{a_i}{i!}\, x^i$
  \item $\displaystyle\sum_{i=0}^n \frac{a_i}{(n-i)!}\, x^i$
  \item $\displaystyle\sum_{i=0}^n \frac{a_i}{i!(n-i)!}\, x^i$
  \item $\displaystyle\sum_{i=0}^n \binom{n}{i}\, a_i\, x^i$
  \item $\displaystyle\sum_{i=0}^n \falling{n}{i} \, a_i\, x^i$
  \end{enumerate}
\end{cor}
\begin{proof}
  The first is $\expoper{}$ applied to $f$,  the second  is
  $\expoper{}$ applied to the reverse, and the third is $\expoper{}$
  applied to the second. Multiplying the third by $n!$ gives the fourth,
  and the fifth is $n!$ times the second.
\end{proof}

We can  prove a more general exponential result using the gamma function:

\begin{lemma} \label{lem:exp-ki}
  For any positive integer $k$ the
  transformations $x^i \mapsto \frac{i!}{(ki)!}x^i$ and $x^i \mapsto
  \frac{x^i}{(ki)!}$ map $\allpoly$ to itself.
\end{lemma}
\begin{proof}
\index{Gamma function}
  In Example~\ref{ex:gamma} we saw that $\Gamma(z+1)/\Gamma(kz+1)$ is in
  $\allpolyf$, and has all negative roots. The result now follows from
  Corollary~\ref{cor:fxdg-analytic} and \eqref{eqn:fxdg} since
  $\Gamma(k+1)/\Gamma(ki+1) = i!/(ki)!$.
\end{proof}

It is not  an accident that $\expoper{}$ occurs along with
$\intoper{}$.

\begin{theorem} \label{thm:exp-uniq}
  If $f$ is a polynomial such that $\intoper{k}f$ has all real roots
  for infinitely many positive integers $k$, then there is a
  polynomial $g$ with all real roots such that $f = \expoper(g)$.
\end{theorem}

\begin{proof}
{If we write} $    f(x) = \displaystyle\sum_{i=0}^n a_i \frac{x^i}{i!} $ then
  \begin{align}
    \intoper{k}f(x) & = 
\sum_{i=0}^n a_i \frac{x^{i+k}}{(i+k)!} \label{eqn:int-exp-1}\\
\intertext{The polynomial on the right hand side of  (\ref{eqn:int-exp-1}) has all real roots.  Replace
  $x$ by $kx$:}
\sum_{i=0}^n a_i \frac{(kx)^{i+k}}{(i+k)!} 
&= \frac{k^kx^k}{k!}\sum_{i=0}^n a_i
\frac{k}{k+1}\cdots\frac{k}{k+i}x^k\nonumber\\
\intertext{so the n-th degree polynomial}
&\sum_{i=0}^n a_i
\frac{k}{k+1}\cdots\frac{k}{k+i}x^k\nonumber
\end{align}
{has all real roots.  Taking the limit as
  $k\rightarrow\infty$ we see that}\quad
$\displaystyle \sum_{i=0}^n a_i x^i$
has all real roots.  This last polynomial is the desired $g$.
\end{proof}

The Laplace transform, \index{Laplace transform} when restricted to
polynomials, is closely related to the exponential transformation. If
$f(t)$ is any function then the Laplace transform of $f(t)$ is
$$ \mathcal{L}(f)(x) = \int_0^\infty e^{-xt} f(t) \,dt.$$

If $x$ is positive then $\mathcal{L}(t^n)(x) = \frac{n!}{x^{n+1}}$.
We say that the Laplace transform $\mathcal{L}(f)$ of a polynomial
$f=\sum a_i x^i$ is $\sum a_i \frac{i!}{x^{i+1}}$. Consequently, we
can express $\mathcal{L}(f)$ in terms of $\expoper{}$:
$$\mathcal{L}(f)(x) = \frac{1}{x} \expoper{}^{-1}(f)(\frac{1}{x})$$
since $\expoper{}^{-1}(x^n)(\frac{1}{x}) = \frac{n!}{x^{n}}.$ If
$f$ is a polynomial of degree  $n$ and $\mathcal{L}(f)(x)$ has $n$
real roots, then $0$ is not a root, so
$\expoper{}^{-1}(f)(\frac{1}{x})$ has all real roots. Taking the
reverse shows that $\expoper{}^{-1}(f)(x)$ has all real roots, and
hence applying $\expoper{}$ yields that $f\in\allpoly$. Summarizing,

\begin{lemma}
  If $f$ is a polynomial of degree $n$ and $\mathcal{L}(f)$ has $n$
  real roots then $f\in\allpoly(n)$.
\end{lemma}

\section{Binomial transformations}
\label{sec:binomial-trans}

In this section we study the transformation $T\colon{}x^n\mapsto
\binom{x}{n}$ and its inverse. The many identities satisfied by 
binomial coefficients yield recurrences for these transformations.

\begin{cor} \label{cor:binomial}
  The linear transformation $x^i\mapsto \binom{x}{i}$ maps

\begin{alignat*}{2}
 \allpolyint{(0,\infty)}& \longrightarrow& \allpolyint{(0,\infty)} & \\
 \allpolyint{(-\infty,-1)}&\longrightarrow&\ \hfill \allpolyint{(-\infty,-1)}&
  \end{alignat*}

\end{cor}
\begin{proof}
  The map $T$ is the composition $x^n \longrightarrow x^n/{n!}
  \longrightarrow \binom{x}{n}$ where the second map is the falling
  factorial (see Proposition~\ref{prop:recur-5}) and the first is the exponential
  map. Since the falling factorial maps $\allpolyalt$ to itself we
  have proved the first part. For the second one we use the identity $ \affa T =
  T \affa$ where $\affa f(x) = f(-x-1)$. This yields the communicative
  diagram

  \centerline{\xymatrix{
\allpolyint{(0,\infty)}
      \ar@{<-}[d]_{{\affa    }}           
      \ar@{->}[rrr]^{{T   }}         
      &&&
\allpolyint{(0,\infty)}
      \ar@{->}[d]^{{\affa   }} \\        
            \allpolyint{(-\infty,-1)}
      \ar@{->}[rrr]^{{ T   }}         
      &&&
            \allpolyint{(-\infty,-1)}
}}

since we know the top row, and $\affa^{-1} = \affa$.

\end{proof}

\begin{cor} \label{cor:Txdi}
  The map $T(x^i)=\binom{x+n-i}{n}$ maps $\allpolyint{(0,1)}$ to
  $\allpolyalt$.
\end{cor}
\begin{proof}
This follows from Lemma~\ref{lem:xdi} and the diagram following. 
See Lemma~\ref{lem:binomial-xdk} for a more general result.
\end{proof}

It is not the case that $T$ maps $\allpolyint{-1,0}$ to itself - even
$T(x+1/2)^2$ has imaginary roots.
We are now going to show that $T^{-1}$ maps $\allpolyint{-1,0}$ to
itself. 

We  use a simple binomial  identity to get a recurrence for $T^{-1}$.

\begin{align*}
  x\binom{x}{n} &= (n+1)\binom{x}{n+1} + n \binom{x}{n} \\
T^{-1}(  x\binom{x}{n}) &= (n+1)x^{n+1} + n x^n \\
&= x\cdot x^n +x(x+1)\cdot nx^{n-1} \\
&= x\cdot T^{-1}\binom{x}{n} + x(x+1)
\left(T^{-1}\binom{x}{n}\right)'\\
\intertext{By linearity it follows that}
T^{-1}(xf) &= xT^{-1}f + x(x+1)(T^{-1}f)'
\end{align*}

\begin{lemma}\label{lem:binomial-trans}
  The map $x^n\mapsto \binom{x}{n}$ satisfies
  $T^{-1}:\allpolyint{(-1,0)}\longrightarrow \allpolyint{(-1,0)}$ and
  preserves interlacing.
\end{lemma}
\begin{proof}
  If $T^{-1}f\in\allpolyint{(-1,0)}$ then $T^{-1}f \lesslesseq
  (T^{-1}f)'$ and so it follows that $x(x+1)T^{-1}f \lesslesseq T^{-1}f$ since all
  roots lie in $(-1,0)$. Using the recurrence yields
  $T^{-1}(xf)\lesslesseq Tf$ and we can follow the argument of
  Corollary~\ref{cor:quant-transform} to finish the proof.
\end{proof}

The behavior of $T^{-1}$ at the endpoints of $(-1,0)$  is
unusual:

\begin{lemma}
  $T^{-1}(x^n) \greateqeq T^{-1}(x+1)^n$. In fact, 
$$ (x+1)T^{-1}(x^n) = x T^{-1}(x+1)^n$$
\end{lemma}
\begin{proof}
  If we define the affine transformation $\affb f(x) = f(x+1)$ then we
  will prove the more general result $ (x+1) T^{-1}f = x T^{-1}\affb f$. It
  suffice to prove this for a basis, so choose $f = \binom{x}{n}$. The
  computation below establishes the lemma.
\begin{gather*} xT^{-1}\affb f= x T^{-1}\affb \binom{x}{n} = xT^{-1} \binom{x+1}{n} =
xT^{-1}(\binom{x}{n}+\binom{x}{n-1}) \\
= x(x^n + x^{n-1}) = (x+1) x^n = (x+1)T^{-1}f
\end{gather*}
\end{proof}

\begin{cor}
  The map $T\colon{}x^n\mapsto \binom{x}{n}$ satisfies 
\[ \int_0^1 T^{-1}(x+t)^n\,dt\in\allpoly\quad\text{for }n=1,2,\cdots.\]
\end{cor}
\begin{proof}
  We can apply Proposition~\ref{prop:family-int} once we show that
  $T^{-1}$ satisfies $T^{-1}(x+a)^n
  \greateqeq T^{-1}(x+b)^n$ if $0\le b \le a \le 1$. We show more
  generally that if the roots of $g$ are all at least as large as the
  corresponding roots of $f$ then $T^{-1}g\greateqeq T^{-1}f$. Since $T^{-1}$
  preserves interlacing we know that there is a sequence of
  interlacing polynomials from $T^{-1}x^n$ to $T^{-1}(x+1)^n$ that
  includes $f$ and $g$. Since the endpoints of this sequence interlace
  by the preceding lemma, we have a mutually interlacing sequence, and
  therefore $T^{-1}g\greateqeq T^{-1}f$.
\end{proof}

We now consider the transformation $T\colon{}x^n \mapsto \rising{x}{n}/n!$.
Since $\rising{-x}{n} = (-1)^n\binom{x}{n}$ we can apply the results
of this section to conclude

\begin{lemma}\label{lem:binomial-rising}
  The map $T\colon{}x^n \mapsto \rising{x}{n}/n!$ satisfies
  \begin{enumerate}
  \item
    $T\colon{}\allpolyint{(-\infty,0)}\cup\allpolyint{(1,\infty)}\longrightarrow \allpoly$. 
  \item $T^{-1}:\allpolyint{(0,1)}\mapsto \allpoly$.
  \end{enumerate}
\end{lemma}

If $\affa f(x) = f(1-x)$ then
$(x-1/2)^k$ is invariant up to sign under $T$ and   its roots are
invariant under $\affa$. The following shows that there are
complex roots, and they all have real part $1/2$.

\begin{lemma}
  If $T\colon{}x^n \mapsto \rising{x}{n}/n!$ then all roots of $T(x-1/2)^n$
  have real part $1/2$. If $n>1$ then there is a complex root.
\end{lemma}
\begin{proof}
If we define 
$$ S(f(x)) = (T(x-1/2))(x+1/2)$$
then  $S$  satisfies the recurrence
\begin{align*}
  S(x^{n+1}) &= \frac{1}{n+1} S(x^n) + \frac{n}{4(n+1)}S(x^{n-1})\\
\end{align*}
Since $S(1)=1$ and $S(x)=x$, it follows that $S(x^n)$ has purely
imaginary roots. The first part now follows from $S(x^n)(x-1/2) =
T(x-1/2)^n$.

If $T(x-1/2)^n$ had all real roots then they would all be $1/2$, but
it is easy to check that 
$T(x-1/2)^n$ is not a multiple of $(x-1/2)^n$. 
\end{proof}

\section{Eulerian polynomials}
\label{sec:euler}
\index{Eulerian polynomials}
\index{polynomials!Eulerian}

The Eulerian polynomials can be defined as the numerator of the sum of
an infinite series.  If $r$ is a non-negative integer, then the
Eulerian polynomial $A_r(x)$ is the unique 
polynomial  of degree $r$ such that
\begin{align}
        \sum_{i=0}^\infty i^r x^i &= \frac{A_r(x)}{(1-x)^{r+1}}.
\label{eqn:sum-of-powers} \\
\intertext{If we differentiate (\ref{eqn:sum-of-powers}) and multiply by $x$, we
get}
    \sum_{i=0}^\infty i^{r+1} x^i &= x\, \frac{(r+1)A_r +
      (1-x)A^\prime_r}{(1-x)^{r+2}}.
\nonumber \\
\intertext{and so}
A_0 &=1 \qquad A_1 = x \nonumber\\
A_{r} &= x\left(\,rA_{r-1} + (1-x)A^\prime_{r-1}\, \right)
\label{eqn:sum-of-powers-2}
\end{align}
\index{recurrence!Euler}

$A_n$ has all non-negative coefficients and Corollary~\ref{cor:tf-euler} shows
that the roots of $A_r$ are all non-positive, and $A_{r} \lesslesseq
A_{r-1}$. Using \eqref{eqn:sum-of-powers-2}, we can write the
recurrence for the transformation $T\colon{}x^n\mapsto A_n$:

\begin{equation}
  \label{eqn:euler-recur}
  T(xf) = xTf +  xT(xf^\prime) + x(1-x)(Tf)^\prime 
\end{equation}

There is a  \index{\Mobius\ transformation}\Mobius\ transformation
that simplifies this recurrence.

\begin{lemma} \label{lem:euler-mobius}
  If $Mz = \frac{z+1}{z}$ then $T_M: \allpolyint{(0,1)}\mapsto
  \allpolyint{(-1,0)}$ where $T_M$ is defined in \eqref{eqn:mobius-2}.
\end{lemma}
\begin{proof}
  Let $B_n = x^n A_n(\frac{x+1}{x})$. A bit of calculation using
  \eqref{eqn:sum-of-powers-2} shows that $B_0=1$ and
\begin{align*}
 B_{n+1} &= (x+1)(B_n+xB_n^\prime). \\
\intertext{This leads to the recurrence for $T_M:$}
 T_M(xf) &= (x+1)(T_M(f) + xT_M(f)^\prime)
\end{align*}
The conclusion now follows easily by induction using this recurrence.

\end{proof}

\begin{cor} \label{cor:euler}
  If $T\colon{}x^n\mapsto A_n$ then $T$ maps $\allpolyint{(-1,0)}$ to $\allpolypos$.
\end{cor}
\begin{proof}
  If $V$ is a linear transformation and $M$ is a \index{\Mobius\ transformation}\Mobius\ transformation then write $(M)\circ V=V_M$.  Setting $S=T_{(z+1)/z}$
  then with this notation we have that
  $$
  T = (z-1) \circ (\frac{1}{z})\circ T_{(z+1)/z}$$
  since the
  composition $(z-1)\circ(1/z)\circ((z+1)/z)$ is the identity. We now
  apply Lemma~\ref{lem:euler-mobius} and the general results of
  \chapsec{linear}{mobius} to deduce the commutative diagram

\centerline{
\xymatrix{
\allpolyint{(0,1)}  \ar@{->}[d]^{S}  \ar@{-->}[rr]^T &  & \allpolypos
\\
 \allpolyint{(-1,0)} \ar@{->}[rr]^{(1/z)}&&  \allpolyint{(-\infty,-1)}
\ar@{->}[u]^{(z-1)}&&  
}}

\end{proof}

It appears that the roots of $A_n$ go to minus infinity, and hence the
roots of polynomials in the image of $\allpolyint{(-1,0)}$ are not 
contained in any finite interval.

If $f$ is a polynomial of degree $n$, then  we
define the transformation $W$ by 
 \begin{align} \sum_{i=0}^\infty f(i)x^i &=
   \frac{(Wf)(x)}{(1-x)^{n+1}} \label{eqn:dfx-1}\\
\intertext{If $f = \sum_1^n b_i x^i$ then }
Wf & = \sum _1^n b_i (1-x)^{n-i}A_i \nonumber
\end{align}
This shows that $Wf$ is   a
polynomial of degree $n$.  $W$ can be realized as a composition with
$T$ and a homogeneous transformation
$$ W: x^k \mapsto A_k \mapsto (1-x)^{n-k}A_k.$$

\begin{cor}
  The map $W$ defines a linear transformation from $\allpolyint{(-1,0)}$
  to $\allpoly$.
\end{cor}

\begin{proof}
  If $f\in\allpolyint{(-1,0)}$ then $Tf\in\allpolyalt$ where $T$ is
  given in Corollary~\ref{cor:euler}. Since $W(x^n)=A_n\in \allpolypos$ has
  leading coefficient $1$ we can apply Lemma~\ref{lem:homog-xy-1} to conclude
  that $Wf\in\allpoly$.
\end{proof}

\section{Euler-Frobenius polynomials}

\index{Euler-Frobenius polynomials}
\index{polynomials!Euler-Frobenius}
\index{recurrence!Euler-Frobenius}

The Euler-Frobenius polynomials \cite{dubeau98} arise in interpolation 
problems.  These polynomials satisfy the recurrence
\begin{align}
  P_{n+1} &= 2xP_n - (1-x^2)P_n^\prime \label{eqn:euler-frob-1}\\
  &= (2x+(x^2-1)\diffd)P_n \nonumber\\
  &= \diffd ((x^2-1)P_n)\nonumber\\
\intertext{If we define the linear transformations $T(x^n)=P_n$ and
  $S(f)=\diffd\,((x^2-1)f)$ then}
T(x^n) &= S^n(1) \nonumber\\
\intertext{and for any polynomial $f$ we have}
T(f) &= f(S)\,(1) \nonumber
\end{align}

By Corollary~\ref{cor:euler-frob} the transformation $2x+(x^2-1)\diffd+\alpha$
maps $\allpolyint{(-1,1)}$ to itself for $|\alpha|<1$ so we conclude

\begin{lemma} \label{lem:euler-frob}
  If $T(x^n) = P_n$ then $T$ maps $\allpolyint{(-1,1)}$ to itself.
\end{lemma}

There is a modification of the Euler-Frobenius polynomials that are
also called Euler-Frobenius polynomials \cite{sampling-zeros}. They
are defined by applying a \index{\Mobius\ transformation}\Mobius\
transformation to $P_n$: 
\index{M\"{o}bius transformation}
\begin{align}
  E_n(x) &= (x-1)^n P_n\left(\frac{x+1}{x-1}\right)\nonumber\\
\intertext{These polynomials satisfy the recurrence}
  \label{eqn:euler-mobius-2}
E_n &= (1+nx)E_{n-1}+x(1-x)E_{n-1}^\prime
\end{align}

It is not obvious that all the coefficients of $E_n$ are positive, but
this can be shown \cite{sampling-zeros} by explicitly determining the
coefficients. The recurrence \eqref{eqn:euler-mobius-2} then shows that
all the roots of $E_n$ are negative, and that $E_n \lesslesseq
E_{n-1}$. We can use the fact that $P_n=S^n(1)$ to derive a
\index{Rodrigues' formula!Euler}
``Rodrigues' formula'' for $E_n$:
$$ E_n = \frac{(1-x)^{n+2}}{x}\,(x\diffd)^n \, \frac{x}{(1-x)^2}$$

\section{Even and odd parts}
\label{sec:hurwitz}

When does a polynomial with all real coefficients have roots whose
real parts are all negative?  The answer is given by  Hurwitz:
write a polynomial $f$ in terms of its even and odd parts. The
polynomial $f$ has all roots with negative real part iff the
polynomials corresponding to the even and odd parts interlace. See
\cite{polya-szego2}*{V171.4} or Proposition~\ref{prop:hb}.

In this section we generalize the transformations that assign either the
even or odd part to a polynomial.  Under appropriate assumptions,
these transformations preserve roots and form mutually interlacing sequences.

Write $f(x) = f_e(x^2) + xf_o(x^2)$. The even part of $f$ is $f_e(x)$, and
the odd part is $f_o(x)$ so we define two linear transformations by
$T_e(f) = f_e(x)$ and $T_o(f) = f_o(x)$ (see Page~\pageref{ex:even-odd}).
It is clear that $T_e$ and $T_o$ are linear, and satisfy the
recurrences

\index{even part}
\index{odd part}

\begin{align}
  \label{eqn:hurwitz}
  T_e(xf) & = xT_o(f) \\
  T_o(xf) &= T_e(f). \notag
\end{align}  

\index{Hurwitz's theorem!for $\intmod{2}$}
\begin{theorem}[Hurwitz's theorem] \label{thm:hurwitz} 
  If $f\in\allpolypos$ then $T_e(f)$ and $T_o(f)$ are in
  $\allpolypos$. In addition, $xT_o(f)$ is interlaced by $T_e(f)$.
  Thus,
$$ \text{if deg($f$) is } 
\begin{cases}
  \text{even} & \text{then } T_e(f) \lesslesseq T_o(f) \\ 
  \text{odd} & \text{then } T_e(f) \greateqeq T_o(f). 
\end{cases}
$$
\end{theorem}



We will derive this theorem from the following more general result. If
$f(x)$ is a polynomial and $d$ is a positive integer, then we can
uniquely write

\begin{equation}
  \label{eqn:hurwitz-1}
  f(x) = f_0(x^d) + x f_1(x^d) + \cdots + x^{d-1}f_{d-1}(x^d)
\end{equation}

\begin{theorem} \label{thm:hurwitz-gen}
\index{Hurwitz's theorem!for $\intmod{d}$} \index{mutually interlacing}
  If $f\in\allpolypos$,  $d$ is a positive integer, and
  $f_0,\dots,f_{d-1}$ are given in \eqref{eqn:hurwitz-1} then
  $f_0,f_1,\dots,f_{d-1}$ is a mutually interlacing sequence. In
  particular, all $f_i$ are in $\allpolypos$.
\end{theorem}

\begin{proof}
  The idea is to use a certain matrix $H$ such that $f(H)$ contains
  all the $f_i$'s as entries. We then show that $f(H)$ preserves
  mutually interlacing sequences. This yields $f_i\in\allpolypos$, and
  the mutual interlacing of the sequence.

  Assume that $f$ is monic, and write 
  $$
  f(x) = \sum_{i=0}^n a_i x^i = \prod_{i=1}^n(x+r_i)$$
  where all
  $r_i$ and $a_i$ are positive. Define the $d$ by $d$ matrix $H$ as
  follows. The upper right corner is $x$, the lower diagonal is all
  $1$'s, and the remaining entries are zero. For instance, if $d$ is
  three then $H,H^2,\dots,H^6$ are

\begin{multline*}
\begin{pmatrix}0& 0& x\\  1& 0& 0\\  0& 1& 0\end{pmatrix}\ ,
\begin{pmatrix}0& x& 0\\  0& 0& x\\  1& 0& 0 \end{pmatrix}\ ,
\begin{pmatrix}x& 0& 0\\  0& x& 0\\  0& 0& x \end{pmatrix}\ ,\\
\begin{pmatrix}0& 0& x^2\\  x& 0& 0\\  0& x& 0 \end{pmatrix}\ ,
\begin{pmatrix}0& x^2& 0\\  0& 0& x^2\\  x& 0& 0 \end{pmatrix}\ ,
\begin{pmatrix}x^2& 0& 0\\  0& x^2& 0\\  0& 0& x^2\end{pmatrix}
\end{multline*}

If $I$ is the identity matrix then 
$$\prod_{i=1}^n(H+r_iI) = \sum_{i=0}^n a_i H^i $$

An important property of $H$ is that $H^{d} = xId$, as we can see
from the above matrices when $d=3$. Upon consideration of the patterns
in the powers of $H$ we see that 
$$ \sum_{i=0}^n a_i H^i =
\begin{pmatrix}
  f_0 & xf_{d-1} & xf_{d-2} &\dots & xf_1 \\
f_1 & f_0 & xf_{d-1} & \dots & xf_2 \\
f_2 & f_1 & \ddots & & \vdots \\
\vdots & \vdots & & \ddots & xf_{d-1} \\
f_{d-1} & f_{d-2} & \dots & f_1 & f_0
\end{pmatrix}
$$
where all terms are either of the form $f_i(x)$ or $xf_i(x)$.

Since $(\sum a_i H^i)(1,0,\dots,0)^t = (f_0,f_1,\dots,f_{d-1})$ it
suffices to show that $H+rI$ preserves mutually interlacing
sequences. This follows from Example~\ref{ex:hri}.
\end{proof}

\deleted{5/21/07}

\begin{remark}\index{Newton's inequalities}
\index{even part!not arbitrary}
    The even part of a polynomial in $\allpolypos$ has all real
    roots. Is every polynomial in $\allpolypos$ the even part of a
    polynomial in $\allpolypos$? We use Newton's inequalities to show
    that the answer is no.

    Suppose that $f = f_e(x^2) + x f_o(x^2)\in\allpolypos(2n)$, and
    let $f_e = \sum a_ix^i$. From \eqref{eqn:newton-3}
    with $k=2$ and $r=1$ we see that
    \begin{equation}
      \label{eqn:newton-even-part}
      \frac{a_1^2}{a_0a_2} \ge
      \frac{4\cdot3}{2\cdot1}\frac{(2n)(2n-1)}{(2n-2)(2n-3)} >
      3\biggl(2\cdot\frac{n}{n-1}\biggr) 
    \end{equation}
Now Newton's inequality for $f_e$ is
\begin{gather*}
      \frac{a_1^2}{a_0a_2} \ge 2\frac{n}{n-1}
\end{gather*}
Thus, if $f_e$ is the even part of a polynomial in $\allpolypos$, then
it must satisfy \eqref{eqn:newton-even-part} which is stronger than
the usual Newton inequality. In particular, $(x+1)^n$ is not the even
part of a polynomial in $\allpolypos$ for $n\ge2$.
  \end{remark}

\begin{example}
    \index{polynomial!stable} If we begin with a polynomial in
    $\allpolypos$ (or a stable polynomial, see
    Chapter~\ref{cha:stable}) then we can form a tree by recursively
    splitting into even and odd parts. For instance, if we begin with
    $(x+1)^8$ we get Figure~\ref{fig:even-odd-tree}. The $2^k$
    polynomials in the $k$th row are mutually interlacing (but not
    left to right).

    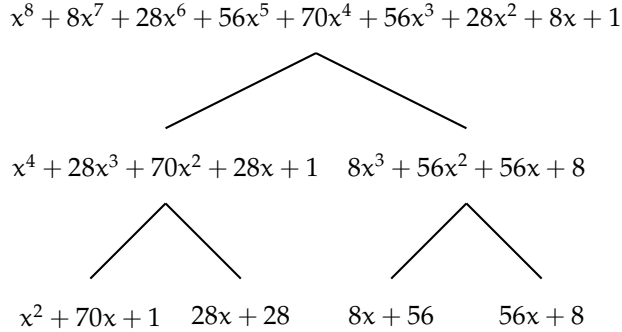
\begin{figure}
      \centering
 
\begin{pspicture}(0,3)(8,9)
  \rput(4,8){$x^8+8x^7+28x^6+56x^5+70x^4 + 56 x^3 + 28 x^2 +8x + 1$}
  \rput(2,6){$  x^4+28x^3+ 70x^2  + 28 x+ 1$}
  \rput(6,6){$8x^3+56x^2+56x+8$}
  \rput(1,4){$x^2+70x+1$}
  \rput(3,4){$28x+28$}
  \rput(5,4){$8x+56$}
  \rput(7,4){$56x+8$}
  \psline(4,7.5)(2,6.5)
  \psline(4,7.5)(6,6.5)
  \psline(2,5.5)(1,4.5)
  \psline(2,5.5)(3,4.5)
  \psline(6,5.5)(5,4.5)
  \psline(6,5.5)(7,4.5)
\end{pspicture}
     
      \caption{The even-odd tree of $(x+1)^8$}
      \label{fig:even-odd-tree}
    \end{figure}

  \end{example}

  There is an elementary reason why the even and odd parts never have
  a common factor.

  \begin{lemma}
    If $f\in\allpolypos$ then $f_e$ and $f_o$ have no common factor.
  \end{lemma}
  \begin{proof}
    Assume that they do. If $f_e = (x+a)g(x)$ and $f_o=(x+a)h(x)$ then
\[f(x) = f_e(x^2)+x f_o(x^2) = (x^2+a)(g(x^2) + x h(x^2))
\]
{This is not possible, since $f\in\allpolypos$ implies that $a$ is
positive, and  $x^2+a\not\in\allpoly$. }
  \end{proof}

\section{Chebyshev and Jacobi Polynomials}
\label{sec:chebyshev}
\index{Chebyshev polynomials}
\index{polynomials!Chebyshev}
\index{recurrence!Chebyshev}

The Chebyshev polynomials $T_n$ are orthogonal polynomials given by
the recurrence formula $T_0=1,T_1=x,$ and $T_{n+1}=2xT_n-T_{n-1}$.  We
study some \Mobius\ transformations associated to the Chebyshev and
Jacobi polynomials.

The transformations associated to the reversal of the Chebyshev and
Legendre polynomials satisfy some simply stated identities. We can use
these to establish mapping properties for the reversal of Chebyshev
polynomials.

For our purposes we need an explicit formula:
\begin{equation}
  \label{eqn:cheby}
  T_n(x) = 2^{-n} \sum_{i=0}^n \binom{2n}{2i}(x+1)^i (x-1)^{n-i}.
\end{equation}

\begin{lemma} \label{lem:chebyshev}
  If $T$ is the linear transformation  $x^n\mapsto T_n$ and
  $Mz=\frac{z+1}{z-1}$ then
  $T_M:\allpolyint{(0,1)}\rightarrow \allpolypos$.
\end{lemma}
\begin{proof}
Using \eqref{eqn:cheby} we find a simple expression for $T_M\,(x^n):$
\begin{align*}
  T_{M}\,(x^n) &= (x-1)^n T_n\left(\frac{x+1}{x-1}\right) \\
  &= (x-1)^n2^{-n} \sum_{i=0}^n \binom{2n}{2i} \left(
    \frac{x+1}{x-1}+1\right)^i
  \left(\frac{x+1}{x-1}-1\right)^{n-i} \\
  &= \sum_{i=0}^n\binom{2n}{2i} x^i \\
  &= \text{ even part of } (1+x)^{2n}\\
  \intertext{where the even part is defined in \chapsec{operators}{hurwitz}.
    Consequently} T_{M}\,(f) &= \text{ even part of } f(\,(1+x)^2\,).
\end{align*}

If $f\in\allpolyint{(0,1)}$ then $f(\,(1+x)^2\,)$ is in
$\allpolyint{(-2,0)}$ and by Theorem~\ref{thm:hurwitz} the even part is in
$\allpolypos$ as well. Thus $T_{M}$ maps $\allpolyint{(0,1)}$ to
$\allpolypos$.
\end{proof}

Formulas analogous to \eqref{eqn:cheby}  hold for all  Jacobi polynomials.
\index{Jacobi polynomials}
\index{polynomials!Jacobi}
From \cite{szego}*{(4.3.2)} 
\begin{align}
  P_n^{\alpha,\beta} (x) &= \sum_{k=0}^n
  \binom{n+\alpha}{n-k}\binom{n+\beta}{k} \left(\frac{x-1}{2}\right)^k
  \left(\frac{x+1}{2}\right)^{n-k} \notag \\
  \intertext{If we apply the \Mobius\ transformation ${M}$ of
    Lemma~\ref{lem:chebyshev} to the linear transformation $J:x^n\mapsto
    P^{\alpha,\beta}_n$ we get} J_M\,(x^n) &= \sum_{k=0}^n
  \binom{n+\alpha}{n-k}\binom{n+\beta}{k}x^k \label{eqn:jacobi} \\
  \intertext{In particular, the Legendre polynomials are $P^{0,0}_n$
    so defining $L(x^n)=P^{0,0}_n$ we have the elegant formula}
  L_M\,(x^n) &= \sum_{k=0}^n \binom{n}{k}^2x^k. \label{eqn:legendre}
\end{align}
\index{Legendre polynomials}
\index{polynomials!Legendre}

It follows from Corollary~\ref{cor:fx1y1} that
\begin{lemma} \label{lem:legendre}
  The map $x^n\mapsto L_M(x^n)$ given in \eqref{eqn:legendre} maps
  $\allpolyalt$ to itself.
\end{lemma}

We now state two identities for reverse polynomials, where $T_n$ and
$P_n$ are the Chebyshev and Legendre polynomials respectively.

\begin{xalignat}{2}
  \label{eqn:rev-iden-1}
  T(x^n) &= (T_n)^{rev} & T(1-x)^n &= 
  \begin{cases}
    0 & \text{$n$ odd}\\
    (1-x^2)^{n/2} & \text{$n$ even}
  \end{cases}\\
  S(x^n) &= (P_n)^{rev} & S(1-x)^n &=
  \begin{cases}
    0 & \text{$n$ odd}\\
    \binom{n-1}{n/2}\,2^{1-n}\,(1-x^2)^{n/2} & \text{$n$ even}
  \end{cases}
\end{xalignat}

\begin{lemma}\label{lem:cheby-rev}
  If $T$ is given above, then
  $T\colon{}\allpolyalt\cup\allpolypos\longrightarrow \allpoly$.
\end{lemma}
\begin{proof} \index{even part}
  The proof relies on the commutative diagram that holds for $(x-1)^n$
  which is a consequence of the identity \eqref{eqn:rev-iden-1}, and by
  linearity it holds for all polynomials. The map \emph{even} returns
  the even part of $f$.

  \centerline{\xymatrix{
      .
      \ar@{->}[d]_{x\mapsto x+1 }           
      \ar@{->}[rrr]^{T }         
      &&&
      .
      \ar@{<-}[d]^{x\mapsto 1-x^2 } \\        
      .
      \ar@{->}[rrr]^{even }         
      &&&
      .
}}

Since we know that the even part determines a map from $\allpolypos$
to itself, the claim in the top line of the following diagram follows
by commutativity.

  \centerline{\xymatrix{
      \allpolyint{(-\infty,0)}
      \ar@{->}[d]_{x\mapsto x+1 }           
      \ar@{->}[rrr]^{T }         
      &&&
      \allpoly
      \ar@{<-}[d]^{x\mapsto 1-x^2 } \\        
      \allpolyint{(-\infty,-1)}
      \ar@{->}[rrr]^{even }         
      &&&
      \allpolyint{(-\infty,0)}
}}

A similar argument applies for $\allpolyalt$.

\end{proof}

\section{Commuting diagrams of spaces}

The earlier commuting diagrams of transformations in
\chapsec{linear}{commute-trans} give rise to
transformations between spaces of polynomials. These are useful to us
because we generally understand three of the four transformations, and
so can deduce information about the fourth.

\index{commuting diagrams}

\begin{example}
$T\colon{}x^n\mapsto H_n$ satisfies (Corollary~\ref{cor:hermite})

\centerline{\xymatrix{
\allpolyint{\reals+\alpha{\imath}} \ar@{.>}[d]_T \ar@{->}[rr]^{x\mapsto
  x+\alpha{\imath}} && \allpoly \ar@{->}[d]^T  \\
\allpolyint{\reals+\alpha{\imath}/2} &&  \ar@{->}[ll]^{x\mapsto x+\alpha{\imath/2}}
\allpoly
}} 

\end{example}

\begin{example}
    $T\colon{}x^k \mapsto H_k(x)H_{n-k}(x)$ satisfies (Corollary~\ref{cor:hermite-xy})

  \centerline{\xymatrix{
\allpolyint{(0,1)}
      \ar@{->}[d]_{reverse }           
      \ar@{->}[rrr]^{T }         
      &&&
\allpoly
      \ar@{=}[d]^{ } \\        
\allpolyint{(1,\infty)}
      \ar@{->}[rrr]^{T }         
      &&&
\allpoly
}}

\end{example}
\begin{example}
  $T\colon{}x^k \mapsto L_k(x)L_{n-k}(x)$ satisfies
  (Corollary~\ref{cor:hermite-xy}) a similar diagram as the Hermite
  transformation above, but the image is smaller.

  \centerline{\xymatrix{
\allpolyint{(0,1)}
      \ar@{->}[d]_{reverse }           
      \ar@{->}[rrr]^{T }         
      &&&
\allpolyint{(0,\infty{})}
      \ar@{=}[d]^{ } \\        
\allpolyint{(1,\infty)}
      \ar@{->}[rrr]^{T }         
      &&&
\allpolyint{(0,\infty{})}
}}

\end{example}

\begin{example}
  For the transformation $ T\colon{}x^k \mapsto
  \rising{x}{k}\falling{x-\alpha}{n-k}$
the following diagram commutes at the function level, but we are
    not able to prove either the top or the bottom. If one were true,
    then the other follows from the commutativity.

  \centerline{\xymatrix{
\allpolyint{(0,1)}
      \ar@{->}[d]_{reverse }           
      \ar@{->}[rrr]^{T }         
      &&&
\allpolyint{(-\infty,0)}
      \ar@{->}[d]^{x\mapsto\alpha-x } \\        
\allpolyint{(1,\infty)}
      \ar@{->}[rrr]^{T }         
      &&&
\allpolyint{(\alpha,\infty)}
}}

  \end{example}

\begin{example}
 The transformation $T\colon{}x^i\mapsto \falling{x+d-i}{d}$ actually maps
 $\allpolyalt\longrightarrow\allpolysep$ (Lemma~\ref{lem:binomial-xdk}).

  \centerline{\xymatrix{
\allpolyint{(0,1)}
      \ar@{->}[d]_{reverse }           
      \ar@{->}[rrr]_{T }         
      &&&
\allpolyint{(0,\infty)}
      \ar@{->}[d]^{x\mapsto-1-x } \\        
\allpolyint{(1,\infty)}
      \ar@{->}[rrr]^{T }         
      &&&
\allpolyint{(-\infty,-1)}
}}

  \end{example}

  \begin{example}
 $T(x^k) = \frac{\rising{x}{k}}{k!}$
satisfies (Lemma~\ref{lem:rising-fact})

  \centerline{\xymatrix{
\allpolyint{(1,\infty)}
      \ar@{->}[d]_{{x\mapsto 1-x    }}           
      \ar@{->}[rrr]^{{T   }}         
      &&&
\allpoly
      \ar@{->}[d]^{{x \mapsto 1-x    }} \\        
            \allpolyint{(-\infty,0)}
      \ar@{->}[rrr]^{{ T   }}         
      &&&
      \allpoly
}}

  \end{example}

  \begin{example}
    $T\colon{}x^k\mapsto H_k(x)x^{n-k}$ acts (Lemma~\ref{lem:hermite-hxn}) on
    polynomials of degree $n$:

  \centerline{\xymatrix{
\allpolyint{(2,\infty{})}
      \ar@{->}[d]_{x\mapsto 4-t }           
      \ar@{->}[rrr]^{T }         
      &&&
      \allpoly
      \ar@{=}[d]^{ } \\        
      \allpolyint{(-\infty{},2)}
      \ar@{->}[rrr]^{T }         
      &&&
      \allpoly
}}

  \end{example}

\section{The Pincherlet\ derivative}
\index{Pincherlet derivative}

If we are given a linear transformation $T$ on polynomials we can form
the \emph{Pincherlet} derivative
$$ T^\prime(f) = T(xf) - xT(f) = [M,T](f)$$
where $M(f)=xf$, and $[M,T]$ is the commutator $TM-MT$. In some
circumstances $T^\prime$ will map $\allpoly$ to itself.

\begin{lemma}
  If $f\lesslesseq g$, $f$ has leading coefficient $a$, $g$ has
  leading coefficient $b$, then $f-xg\in\allpoly$ if $a\ge b$.
\end{lemma}
  \begin{proof}
    Immediate from Lemma~\ref{lem:sign-quant}.
  \end{proof}

  \begin{cor}
    Suppose $T\colon{}\allpoly\longrightarrow \allpoly$ maps polynomials with
    positive leading coefficients to polynomials with positive leading
    coefficients. If the sequence of leading coefficients of $T(x^i)$
    is strictly increasing then $T^\prime:\allpoly\longrightarrow\allpoly$.
  \end{cor}

  \begin{proof}
    Choose monic $f$. Then since $xf\lesslesseq f$ and since $T$
    preserves interlacing, $T(xf)\lessless T(f)$. By hypothesis
    $T(xf)$ and $T(f)$ meet the conditions of the lemma, so $T(xf) -
    xT(f) = T^\prime(f)\in\allpoly$.
  \end{proof}

  \begin{cor}
    The linear transformation $x^n \mapsto H_{n+1} - xH_n$ maps
    $\allpoly\longrightarrow\allpoly$. 
  \end{cor}
  \begin{proof}
    The leading coefficient of $H_n$ is $2^n$.
  \end{proof}

If $T$ maps $\allpolypos$ to itself, then we have weaker assumptions.

\begin{lemma}
  Suppose $T\colon{}\allpolypos\longrightarrow \allpolypos$ maps polynomials
  with positive leading coefficients to polynomials with positive
  leading coefficients, and preserves degree.  Then,
  $T^\prime:\allpolypos\longrightarrow\allpolypos$.
\end{lemma}
\begin{proof}
  Choose $f\in\allpolypos$. Then $T(xf)\lesslesseq Tf$ and since
  $T(xf)\in\allpolyposclose$ it follows that $xT(f) \greateqeq T(xf)$
  and hence $T(xf) - xT(f)\in\allpolypos$.
\end{proof}

\begin{cor}
  The linear transformation $x^n\mapsto n \rising{x}{n}$ maps
  $\allpolypos$ to itself.
\end{cor}
\begin{proof}
  If $T(x^n) = \rising{x}{n}$ then $T^\prime(x^n) =
  \rising{x}{n+1}-x\rising{x}{n} = n\rising{x}{n}$. Alternatively, this
  is the composition

\centering{
\xymatrix{ x^n \ar@{->}[rrr]^{differentiation} &&& n x^{n-1}
  \ar@{->}[rr]^{\text{mult. by } x} && n x^n \ar@{->}[rr]^{x^n\mapsto \rising{x}{n}}
  && n\rising{x}{n} 
}}

\end{proof}

\section{Hypergeometric Polynomials}
\label{sec:hypergeo}

\index{Hypergeometric series}
\index{Hypergeometric polynomials}
Hypergeometric series are a natural class of functions that are
sometimes polynomials. We show that for certain values of the
parameters these polynomials are in $\allpoly$. We begin with the
confluent hypergeometric function $\foneone$ which is defined by

$$
  \foneone(a,b;z) \ =\  1 + \frac{a}{b}\frac{z}{1!} +
  \frac{a(a+1)}{b(b+1)}\frac{z^2}{2!} + \cdots 
\ =\  \sum_{i=0}^\infty \frac{\rising{a}{i}}{\rising{b}{i}}\frac{z^i}{i!} 
$$

We are interested in the case that $a$ is a negative integer, and $b$
is a positive integer. Since  $\rising{a}{i}$ is zero if $a$ is a negative integer
and $i>|a|$, it follows that $\foneone(a,b;z)$ is a polynomial in
these cases, of degree $|a|$. For example, $ \foneone(-4,b;z)$ is
equal to

$$1 - \frac{4\,z}{b} + \frac{6\,z^2}{b\,\left( 1 + b \right) } - \frac{4\,z^3}{b\,\left( 1 + b \right) \,\left( 2 + b \right) } + 
  \frac{z^4}{b\,\left( 1 + b \right) \,\left( 2 + b \right) \,\left( 3
      + b \right) }$$
{and we can manipulate this into a more familiar form}
 $$\frac{1}{\falling{b+4}{4}}
\ \left(\falling{b+4}{4} - 4z\falling{b+4}{3} + 6z^2\falling{b+4}{2} - 4z^3\falling{b+4}{1}  + z^4\falling{b+4}{0}\right)$$
{In general when  $d$ a positive integer we can write}
$$
\foneone(-d,b;-z) = \frac{1}{(b+d)_d}
\sum_{i=0}^a \falling{b+d}{i}\binom{d}{i}z^i
$$

Applying Lemma~\ref{lem:falling-non-int} twice to this representation shows that 

\begin{lemma}
If $a$ is a positive integer, and $b$ is positive  then
$$\foneone(-a,b;-z)\in\allpolypos$$
\end{lemma}

There are many relationships involving hypergeometric functions - the
one below shows that the above result also follows from the fact that
the generalized Laguerre polynomial $L_n^\lambda$ has all real roots
for $\lambda>-1$.

$$
L_n^\lambda(z) = \frac{\rising{1+\lambda}{n}} {\Gamma(1 + \nu )}\,
\foneone( -n  , 1 + \lambda  ,z)
$$
\index{recurrence!Hypergeometric}

\noindent%
$\foneone$ satisfies recurrences in each parameter:

\begin{align*}
\foneone(a, b, z) &= \frac{ 2 + 2\,a - b + z}{1+a-b}\foneone(1+a,b,z)-
  \frac{a+1}{1+a-b}\foneone(2 + a, b, z)\\
\foneone(a,b;z) &= \frac{b+z}{b}\foneone(a,b+1;z) + \frac{(a-b-1)z}{b(b+1)}\foneone(a,b+2;z)
\end{align*}

These recurrences can be used to establish this interlacing square
for positive integral $a$, and positive $b$:

\centerline{
\xymatrix@-1pc{
\foneone(-(a+1),b;z) \ar@{->}[dd]^\greateq &
\ar@{}[r]^{\lessless}_{\ }
 & &
\ar@{->}[dd]^\greateq \foneone(-a,b;y)\\
& & & \\
\foneone(-(a+1),b+1;z) & \ar@{}[r]^{\lessless}_{\ } & & \foneone(-a,b+1;y)
}}

If the hypergeometric series has only one factor in the numerator then
we can apply the same argument. For instance

\begin{align*}
\fonetwo(a;b,c;z) &= \sum_{i=0}^\infty
\frac{\rising{a}{i}}{\rising{b}{i}\rising{c}{i}} \frac{z^i}{i!}\\
\intertext{and if we take $d$ to be a positive integer then}
\fonetwo(-d;b,c;-z) &=
\frac{1}{\rising{b+d}{d}\rising{c+d}{d}}
\sum_{i=0}^d \binom{d}{i}\falling{b+d}{i} \falling{c+d}{i} z^i
\end{align*}

and we can again apply Lemma~\ref{lem:falling-non-int} to yield that
$\fonetwo(-d;b,c;z)\in\allpolyalt$.

\section{Eigenpolynomials of linear transformations}
\label{sec:eigenpolynomials}

\index{eigenpolynomial} If $T\colon{}\allpoly\longrightarrow\allpoly$ is a
linear transformation then we can look for eigenpolynomials of $T$.
These are polynomials $e$ for which there is a $\lambda\in\reals$ so
that $Te =\lambda e$. We consider several questions: does $T$ have
eigenpolynomials? If so, are they in $\allpoly$? Do successive
eigenpolynomials interlace? What are explicit examples? 

\begin{lemma}\label{lem:eigen-1}
  Suppose $T\colon{}\allpoly\longrightarrow\allpoly$ maps polynomials of
  degree $n$ to polynomials of degree $n$. Let $c_n$ be the leading
  coefficient of $T(x^n)$. 
  \begin{enumerate}
  \item If all $c_n$ are distinct then $T$ has a
  unique eigenpolynomial of degree $n$ with eigenvalue $c_n$. 
\item If $|c_n|> |c_r|$ for $0\le r < n$ then $T$ has an
  eigenpolynomial of degree $n$ that is in $\allpoly$. 
  \end{enumerate}
\end{lemma}
\begin{proof}
  
  If we let $V$ be the vector space of all polynomials of degree at
  most $n$ then the assumption on $T$ 
  implies that $T$ maps $V$ to itself. In addition, the
  matrix representation $M$ of $T$ is upper triangular.  The
  $r$-th diagonal element of $M$ is the coefficient $c_r$ of $x^r$ in
  $T(x^r)$.  If they are all distinct then $T$ has an eigenvector
  corresponding to each eigenvalue.

\index{power method}
\index{eigenpolynomial}

We now use the power method to find the eigenvalue.  Choose any
initial value $v_0$ and define
$$ v_{k+1} = \frac{1}{|v_k|} \, Mv_k$$
   where $|v_k|$ is any vector norm. Since $M$ has a dominant
   eigenvalue, the $v_k$ converge to an eigenvector $v$ of $M$,
   provided that $v_0$ is not orthogonal to $v$. 

   Recasting this in terms of polynomials, we choose an initial
   polynomial $p_0$ in $\allpoly(n)$. Since $T$ maps
   $\allpoly$ to itself, all the $p_k$ are in $\allpoly$. If this
   happens to converge to zero, then $p_0$ was orthogonal to
   $p$. Simply perturb $p_0$, and we get an eigenpolynomial in
   $\allpoly$. 

\end{proof}

The sequence of possible eigenvalues is quite restricted.

\begin{cor}
  Suppose that $T\colon{}\allpoly\longrightarrow\allpoly$ preserves degree
  and has positive eigenvalues $\lambda_n$ for $n=0,1,\dots$. Then
\[ \sum_{n=0}^\infty \frac{\lambda_n}{n!}x^n\in\allpolyf.\]
\end{cor}
\begin{proof}
  The $\lambda_n$ are the leading coefficients of $T(x^n)$. Apply
  Lemma~\ref{lem:pd-lead-coef}.
\end{proof}

\begin{lemma}\label{lem:eigen-2}
  Suppose that $T$ is a linear transformation
  $\allpoly\longrightarrow\allpoly$, $p_n$ and $p_{n+1}$ are
  eigenpolynomials of degree $n$ and $n+1$ with eigenvalues
  $\lambda_n$ and $\lambda_{n+1}$. If $|\lambda_{n+1}|>|\lambda_n|$
  and the roots of $p_n$ are distinct 
  then $p_{n+1}\lesslesseq p_n$. 
\end{lemma}
\begin{proof}
  We show that $g(x)=p_n+\alpha p_{n+1}\in\allpoly$ for any $\alpha$. If we
  define
  \begin{align*}
g_k(x) &= \lambda_n^{-k} + \alpha \lambda_{n+1}^{-k} p_{n+1} \\
\intertext{then $T^k g_k = g$ so it suffices to show that
  $g_k\in\allpoly$. Now}
g_k(x) &= \lambda_n^{-k} \left( p_n +
  \left(\lambda_n/\lambda_{n+1}\right)^kp_{n+1}\right)      
  \end{align*}
Since $|\lambda_n/\lambda_{n+1}|<1$ we can apply
Lemma~\ref{lem:add-small} for $k$ sufficiently large since the roots
of $p_n$ are distinct. It follows that
$g_k\in\allpoly$, and hence $g\in\allpoly$. 

\end{proof}

\begin{remark}
  $T(f) = f - f''$ is a simple example of a linear transformation
  mapping $\allpoly\longrightarrow\allpoly$ that has only two
  eigenpolynomials.  If $T(f)=\lambda f$ then clearly $\lambda=1$ and
  $f''=0$, so $f$ has degree $0$ or $1$. The matrix representing $T$
  has all $1$'s on the diagonal, so the lemma doesn't apply. 
\end{remark}

We next show that multiplier transformations
\index{multiplier~transformations} generally have no interesting
eigenpolynomials, and then we determine the eigenpolynomials for the
Hermite transformation.

\begin{lemma}
  If $T(x^n) = t_nx^n$ where all the $t_i$ are distinct and non-zero,
  then the only eigenpolynomials are multiples of $x^i$.
\end{lemma}

\begin{proof}
 If $Te=\lambda e$ where $e = \sum a_ix^i$ then $a_i = \lambda t_i
 a_i$. Every non-zero $a_i$ uniquely determines $\lambda= 1/t_i$, so
 there is exactly one non-zero $a_i$.
\end{proof} 

Some simple linear transformations that map $\allpoly$ to itself have
no eigenpolynomials. For instance, consider $Tf = f+\alpha f^\prime$.
If $f+\alpha f^\prime = \lambda f$ then since the degree of $f^\prime$
is less than the degree of $f$ we see $\lambda=1$, and hence
$\alpha=0$.

The eigenpolynomials of the Hermite transformation $T(x^n) = H_n$
\index{Hermite polynomials} are given by a composition, and are also
orthogonal polynomials.

\begin{lemma}
  If $T(x^n)=H_n$ then the eigenpolynomials of $T$ are given by the
  composition $$ g_n = (3x-2\diffd)^n(1)$$
  and the corresponding eigenvalue is $2^n$.
\end{lemma}
\begin{proof}
  Since $T(f)=f(2x-\diffd)(1)$ we can apply Lemma~\ref{lem:xandD} to conclude
  that
  \begin{align*}
    T(g_n) & = g_n(2x-\diffd)(1) \\
&= (6x-4\diffd)^n(1) \\
&= 2^n g_n
  \end{align*}
\end{proof}

It isn't known how to determine if a linear transformation $T$ has all its
eigenpolynomials in $\allpoly$. If $T$ eventually maps into $\allpoly$
then all the  eigenpolynomials are in $\allpoly$.

\begin{lemma}\label{lem:find-eigen}
  If $T$ is a linear transformation such that for any polynomial $f$
  there is an integer $n$ such that $T^n\,f\in\allpoly$, then all
  eigenpolynomials corresponding to non-zero eigenvalues are in $\allpoly$. 
\end{lemma}
\begin{proof}
  If $Tf = \lambda f$ then choose $n$ as in the hypothesis and note
  that $T^n f = \lambda^n f$. If $\lambda$ is non-zero it follows that
  $f\in\allpoly$. 
\end{proof}

\section{Classical multiple orthogonal polynomials}

It is possible to construct  polynomials that are orthogonal
to several weight functions \cite{vanassche}.  The resulting
polynomials are called \index{multiple orthogonal polynomials}
\emph{multiple orthogonal polynomials}  and all have all
real roots. In this section we look at the corresponding differential
operators, and also observe some interlacing properties.

\begin{lemma} \
  \begin{enumerate}
  \item Suppose $\alpha>-1$. The operator $T_{n,\alpha}:f\mapsto x^{-\alpha}
    \diffd^n x^{n+\alpha} \,f$ maps $\allpoly\longrightarrow\allpoly$
    and $\allpolypos$ to $\allpolypos$.
  \item If $f\in\allpolypm$ then $T_{n+1,\alpha}f \greateqeq T_{n,\alpha}f$.
  \item Suppose $b\ge0$. The operator $S_{n,b,c}:f\mapsto e^{bx^2-cx} \diffd^n
    \, e^{-bx^2+cx}f$ maps $\allpoly\longrightarrow\allpoly$.
  \end{enumerate}
\end{lemma}
\begin{proof}
The first one follows from a more general result (Proposition \ref{prop:rod-laguerre}) in
several variables. To check interlacing we note that
\begin{align*}
  \beta T_{n,\alpha}f + T_{n+1,\alpha}f & =
x^{-\alpha}\diffd^n\left( \beta x^{\alpha+n} f + \diffd x^{n+1+\alpha}f\right)
\\
&= x^{-\alpha} \diffd^n x^{n+\alpha} \left( (\beta+n+\alpha+1)f + x
  f^\prime\right)\\
\intertext{Since $f\in\allpolypm$ we know that 
$h(x) = (\beta+n+\alpha+1)f + x  f^\prime \in\allpoly$
which implies that}
  \beta T_{n,\alpha}f + T_{n+1,\alpha}f & = x^{-\alpha} \diffd^n
  x^{n+\alpha} h
\end{align*}
where $h\in\allpoly$.  By the first part this is in $\allpoly$ for all
$\beta$ and so the interlacing follows.

For the Hermite case it suffices to show that 
$e^{bx^2-cx} D( e^{-bx^2+cx} f)$ is in $\allpoly$ if $f$ is in
$\allpoly$. Note that we know that $ D( e^{-bx^2+cx} f)$ is in $\allpolyf$
and has the form $e^{-bx^2+cx} g$, but we can't multiply by $e^{bx^2-cx}$
since $e^{bx^2}$ is not in $\allpolyf$. We just directly compute
\begin{equation*}
e^{bx^2-cx} D( e^{-bx^2+cx} f) =  ( f^\prime - (2bx-c) f)
\end{equation*}
and $f^\prime-(2bx-c)f$ is in $\allpoly$ by Lemma~\ref{lem:sign-quant}. 
\end{proof}

If $f=1$ in (1) then we have \index{Laguerre polynomials}Laguerre
polynomials, and if $f=1,c=0,b=1$ in (2) we have 
\index{Hermite polynomials} Hermite polynomials.

If we iterate these operators we get, up to a constant factor, the
classical multiple orthogonal polynomials of \cite{vanassche}. Let $n
= (n_1,\dots,n_r)$ be positive integers, and
$\alpha=(\alpha_0,\dots,\alpha_r)$ where $\alpha_1,\dots,\alpha_r$ are
all at least $-1$. The following polynomials are all in $\allpoly$. 
The interlacing property above shows that the Jacob-Pi\~neiro
polynomials $P_n^\alpha$ and $P_m^\alpha$ interlace if $n$ and $m$
are equal in all coordinates, and differ by one in the remaining
coordinate. 

\index{Jacobi-Pi\~neiro polynomials}
\index{polynomials!Jacobi-Pi\~neiro}
$$
\begin{array}{ll}
  \text{Jacobi-Pi\~neiro}  \quad &
\displaystyle  (1-x)^{-\alpha}\prod_{i=1}^r\left[x^{-\alpha_i} 
\frac{d^{n_i}}{dx^{n_i}}\,x^{n_i+\alpha_i}\right]\,
(1-x)^{\alpha_0+n_1+\cdots+n_r}\\
\text{Laguerre-I}  \quad &
\displaystyle  e^x \prod_{i=1}^r\left[x^{-\alpha_i} 
\frac{d^{n_i}}{dx^{n_i}}\,x^{n_i+\alpha_i}\right]\,e^{-x}\\
\text{Laguerre-II}   \quad & 
\displaystyle  x^{-\alpha_0} \prod_{i=1}^r\left[e^{\alpha_ix}
\frac{d^{n_i}}{dx^{n_i}} e^{-\alpha_ix}\right] \,x^{\alpha_0+n_1+\cdots+n_r}\\
\text{Hermite}  \quad  &\displaystyle \prod_{i=1}^r\left(e^{bx^2-cx} \diffd^n_i \,e^{-bx^2+cx}\right) \,1
\end{array}
$$

  \section{Laurent Orthogonal polynomials}
\index{orthogonal polynomials!Laurent}
\index{Laurent orthogonal polynomials}
  
Laurent orthogonal polynomials satisfy the recurrence
  \begin{align*}
    p_{-1} &= 0 \\    p_0  & = 1 \\
    p_n  &= (a_n x + b_n)p_{n-1} - c_n x\, p_{n-2} \qquad\text{for n
    $\ge 1$}
  \end{align*} 
  where all $a_n,b_n,c_n$ are positive. There are three simple
  recurrences that appear to have many mapping properties (see
  Appendix~\ref{chap:mma}). We are able to prove some mapping properties
  for two of them. Notice that the recurrence implies that the
  coefficients are alternating, and hence if  the roots are all real
  they are all positive.

\begin{lemma}
  If $p_n$ is the sequence of Laurent orthogonal polynomials
  arising from the recurrence $p_{n+1} = x \,p_n - n\,x\,p_{n-1}$ then 
  the transformation $T\colon{}x^n\mapsto p_n$ satisfies
  $T\colon{}\allpolyalt\longrightarrow\allpolyalt$ and preserves interlacing. 

  The transformation $T\colon{}x^n\mapsto p_n^{rev}$ maps $\allpolyalt$ to $\allpolyalt$.
\end{lemma}
\begin{proof}
  Using linearity 
  \begin{align*}
    T(x\cdot x^n) &= x p_n - nxp_{n-1} \\
    &= x\,T(x^n) - x T(x^n\,') \\
    T(xf) &= x\,T(f-f').
  \end{align*}
  We now prove the result by induction - the base case is easy and
  omitted. 
  Since $f-f'\greateqeq f$ it follows by induction that
  $T(f-f')\greateqeq T(f)$.  Since the roots are positive by
  induction, we know that 
\[
T(xf) = xT(f-f') \lesslesseq T(f)
\]
and the conclusion follows from Corollary~\ref{cor:quant-transform}. 

For the second part, we have the recurrence $T(xf) = T(f) -
x\,T(f')$. We can't apply Corollary~\ref{cor:quant-transform} since $T(f)$
does not have positive leading coefficients. However, if we let $S(f) =
T(f)(-x)$ then we can apply the lemma to $S$. Since $S$ satisfies the
recurrence 
\begin{align*}
  S(xf) &= S(f) + x\,S(f')\\
\intertext{we see that by induction}
S(xf) & = S(f)+x\,S(f) \lesslesseq S(f)
\end{align*}
and so $S\colon{}\allpolyalt\longrightarrow\allpolypos$.
\end{proof}

\begin{lemma}
  If $p_n$ is the sequence of Laurent orthogonal polynomials
  satisfying the recurrence $p_{n+1} = (x+n) \,p_n - n\,x\,p_{n-1}$ then 
  the transformation $T\colon{}x^n\mapsto p_n$ satisfies
  $T\colon{}\allpolyalt\longrightarrow\allpolyalt$ and preserves interlacing. 
  \end{lemma}
  \begin{proof}
    The proof is similar to the previous one. We check that
\begin{gather*}
T(xf) = x\,T(f + xf'-f')\\
\intertext{Now we know that}
 f \greateqeq f + x f' \lesslesseq f' \\
\text{so }\quad  f+ xf' - f' \greateqeq f \\
\intertext{By induction}
T(xf) = xT(f+ xf' - f') \lesslesseq T(f)
\end{gather*}
and we apply Corollary~\ref{cor:quant-transform} again to complete the proof.
\end{proof}

\section{Biorthogonal polynomials}

In this section I would like to sketch a completely different approach
to the problem of showing that a linear transformation preserves real
roots. The main results are by Iserles and Norsett, and can be found
in \cite{iserles}. The idea is to begin with a parametrized measure,
and to use it to associate to a polynomial a \emph{biorthogonal}
polynomial that is determined by the roots of the polynomial. There
are a number of non-degeneracy conditions that must be met. In general
the resulting transformation is highly non-linear, but there are quite
a number of measures for which we get linearity.

Choose a function $\omega(x,\alpha)$ and a measure $d\mu$, and define
a linear functional $$\mathcal{L}_\alpha(f) = \int
f(x)\omega(x,\alpha)\,d\mu$$
We start with a monic polynomial $g(x)$
with all real roots $r_1,\dots,r_n$.  Each root determines a linear
functional $ \mathcal{L}_{r_i}$. We next try to construct a monic
polynomial $p(x)$ of degree $n$ that satisfies
$$
\mathcal{L}_{r_1}(p)=0,\quad\mathcal{L}_{r_2}(p)=0,\quad\cdots\quad\mathcal{L}_{r_n}(p)=0$$
This is a set of linear equations. There are $n$ unknown coefficients,
and $n$ conditions, so as long as the determinant of the system is
non-zero we can uniquely find $p$, which is called a \emph{biorthogonal} polynomial.

Next, we want to show that $p$ has all real roots. The proof of this
is similar to the proof that a single orthogonal polynomial has all real
roots, and requires that the $w(x,\alpha)$ satisfy an interpolation
condition.

Finally, we define the transformation $T(g)=p$.  It is quite
surprising that there are \emph{any} choices that make $T$ linear. The
coefficients of $p$ are symmetric functions of the roots, and so can
be expressed in terms of the coefficients of $g$, but in general there
is no reason to suppose that there is linearity.



\chapter{Affine Transformations of Polynomials}
 
\renewcommand{\TimeStampStart}{Friday, January 18, 2008: 09:48:31}
\mytoday  

\label{cha:affine}    

In this chapter we investigate properties of polynomials associated to
affine transformations. Recall that an affine transformation $\affa$
satisfies $\affa(x) = qx+b$, and the action on polynomials is
$\affa(f(x)) = f(qx+b)$.  In order not to burden the reader (and the
writer) with lots of similar but different cases, we will generally
confine ourselves to the affine transformations of the form
$\affa(x)=qx+b$ where $q\ge1$.  We are interested in those polynomials
that interlace their affine transformations. In particular we will
study the set of polynomials
$$
\allpolyaffine = \{f\in\allpoly \mid f \greateqeq  \affa f\}.$$

\section{Introduction}
\label{sec:affine-intro}



We begin with general properties of these affine transformations for
$q>1$.  $\affa$ acts on the real line, and preserves order.  Since
$q>1$ we see that $\affa$ has a unique fixed point $\frac{b}{1-q}$.
Consequently $\affa$ is a bijection restricted to the two intervals
\begin{xalignat*}{2}
  \affa^- &= (-\infty,\frac{b}{1-q})&
  \affa^+ &= (\frac{b}{1-q},\infty)  
\end{xalignat*}
Since $\affa$ and $\affa^{-1}$ preserve order, it is easy to see that
$\affa$ is increasing on $\affa^+$ and it is decreasing on $\affa^-$.
More precisely,

$$
\begin{array}{ccc}
\alpha < \affa\alpha & \text{iff} & \alpha \in  (\frac{b}{1-q},\infty)\\
\affa\alpha < \alpha & \text{iff} & \alpha \in
(-\infty,\frac{b}{1-q})
\end{array}
$$

The action of $\affa$ on the roots of a polynomial is easy to
determine. If $f(\alpha)=0$ then $(\affa f)(\affa^{-1}\alpha) =
f(\alpha)=0$ so that the roots of $\affa f$ are $\affa^{-1}$ applied to
the roots of $f$. Since $\affa^{-1}$ preserves order it follows that
$f\longleftarrow g$ implies $\affa f\longleftarrow \affa g$.
If $f\greateqeq \affa f$ and $\alpha$ is a root of $f$ then we must
have that $\affa^{-1}\alpha \le \alpha$, and hence $\alpha \le \affa
\alpha$. We conclude that if $f\in\allpolyaffine$ then the roots of
$f$ must lie in $(\frac{b}{1-q},\infty)$. For the three examples
considered above we have

$$
\begin{array}{lccl}
\affa(x)=2x \quad&\quad f\in\allpolyaffine & \quad \implies\quad &
\text{roots of } f\in (0,\infty)\\
\affa(x)=2x+1 \quad&\quad f\in\allpolyaffine & \quad \implies\quad &
\text{roots of } f\in (-1,\infty)\\
\affa(x)=x+1 \quad&\quad f\in\allpolyaffine & \quad \implies\quad &
\text{roots of } f\in(-\infty,\infty)
\end{array}
$$

Before we proceed, we should show  that $\allpolyaffine$ is non-empty. Define
\begin{align}
  \label{eqn:affine-1}
  \arising{x}{a}{n} &= =\,\prod_{i=0}^{n-1} (\affa^i x - a) \\
\intertext{In case that $\affa (x) =x+1$  we have}
  \label{eqn:affine-1a}
 \arising{x}{0}{n} &= (x)(x+1)(x+2)\cdots(x+n-1)\ =\ \rising{x}{n}, \\
\intertext{if $\affa(x) = qx$ then }
 \label{eqn:affine-1b}
\arising{x}{-1}{n} &= (1+x)(1+qx)\cdots(1+q^{n-1}x), \\
\intertext{and again if $\affa(x) = qx$ then }
\label{eqn:affine-1c}
\arising{x}{0}{n} &= q^{\binom{n}{2}}\,x^n.
\end{align}

\noindent%
Since  $\affa^{-1}\alpha < \alpha$  for $\alpha\in\affa^+$, if
$a\in\affa^+$ then 
the roots of $\arising{x}{a}{n}$  in increasing  order are
\begin{gather*}
\affa^{-n+1}(a),\ \dots,\ \affa^{-2}(a),\
\affa^{-1}(a), a\\
\intertext{and the roots of $\affa \arising{x}{a}{n}$ are }
\affa^{-n}(a),\ \dots,\ \affa^{-3}(a),\ \affa^{-2}(a),\ \affa^{-1}(a)\\
\end{gather*}
Consequently if $a\in\affa^+$ then $\arising{x}{a}{n} \greateqeq \affa
\arising{x}{a}{n} $.

\begin{lemma}
  If $f\in\allpolyaffine$ then $\affa f\in\allpolyaffine$.
\end{lemma}
\begin{proof}
 If $f\in\allpolyaffine$ then $f\greateqeq \affa f$. Since $\affa$
 preserves interlacing $\affa f \greateqeq \affa^2 f$ and hence $\affa
 f\in\allpolyaffine$.
\end{proof}

Polynomials in $\allpolyaffine$ satisfy a hereditary property.

\begin{lemma} \label{lem:order-ideal} \index{order ideal}
  $\allpolyaffine$ is an order ideal.  That is, if
  $f\in\allpolyaffine$ and $g$ divides $f$ then
  $g\in\allpolyaffine$.
\end{lemma}
\begin{proof}
  It suffices to show that if we remove one linear factor from $f$ the
  resulting polynomial $g$ is still in $\allpolyaffine$. Since affine
  transformations preserve order, the effect of removing a root from
  $f$ is to also remove the corresponding root in $\affa f$. The
  remaining roots interlace, and so $g\in\allpolyaffine$.
\end{proof}

The next lemma gives a useful criterion for determining if a polynomial
is in $\allpolyaffine$.

\begin{lemma}\label{lem:affine-why}
  The following are equivalent
  \begin{enumerate}
  \item $f\in\allpolyaffine$
  \item There is a $g$ such that $f\longleftarrow g$ and
    $\affa f\longleftarrow g$.
  \item There is an $h$ such that $f\longleftarrow h$ and
    $f\longleftarrow \affa h$.
  \end{enumerate}
\end{lemma}
\begin{proof}
  If (1) holds then $f\greateqeq \affa f$.  Choosing $g=\affa f$ shows
  that (1) implies (2), and taking $h=f$ shows that (1) implies (3).
  If (2) holds and $f(x)=\prod(x-r_i)$ then $\affa^{-1} r_i \le r_i$, so
  (1) follows from Lemma~\ref{lem:ci}. If (3) holds then we know that
  $f\longleftarrow h$ and $\affa^{-1}f \longleftarrow h$, so (2)
  implies that $\affa^{-1}f \greateqeq f$. Applying $\affa$ to each
  side shows that (1) holds.
\end{proof}

\section{The affine derivative}
\label{sec:affine-derivative}

We  define the affine derivative  $\daffine f$ of $f$ as follows:

\begin{definition}
  If $f$ is a polynomial then
$  \daffine f =
  \displaystyle\frac{\affa f - f}{\affa x - x}.  
$
\end{definition}

\index{difference operator}
In case $\affa f = f(x+1)$ the affine derivative is the difference
operator $\Delta(f) = f(x+1)-f(x)$. If $\affa f = f(qx)$ then
$\daffine f$ is the q-derivative:
\begin{equation} \label{eqn:affine-product-3}
\qderiv(f) =
\begin{cases}
  \displaystyle   \frac{ f(qx) - f(x)}{ qx - x} & \text{ if $ x \ne 0$} \\
  f^\prime(x) & \text{ if $ x = 0$}
\end{cases}
\end{equation}

The affine derivative is well defined since $q\ne1$ implies $\affa(x)\ne
x$. It is immediate from the Taylor series that $\daffine $ is
continuous. If we let $\affa$ approach the identity, then $\daffine $
converges to the derivative. The degree of $\daffine f$ is one less
than the degree of $f$. As is to be expected for a derivative,
$\daffine x=1$, and there is a product rule
\begin{multline}
  \label{eqn:affine-product-1}
  \daffine (fg) = \frac{\affa f \cdot \affa g - f g}{\affa x - x} =
  \\
 \frac{\affa f \cdot \affa g - \affa f \cdot g + \affa f \cdot g -f
  g}{\affa x - x} = \affa f\cdot \daffine  g + g \daffine f
\end{multline}
In particular we have
\begin{equation}
  \label{eqn:affine-product-2}
  \daffine (xf) = \affa f + x \daffine f 
= f + (\affa x)\daffine f 
\end{equation}

  The formula for $\daffine$ on the product of three terms is
  \begin{align}
    \daffine(f_1f_2f_3) &= \daffine f_1\cdot(f_2f_3) + 
(\affa f_1)\cdot \daffine f_2 \cdot (f_3) +
\affa (f_1 f_2) \cdot \daffine f_3\notag\\
\intertext{and the general formula is}
\daffine(f_1\cdots f_n) &= \sum_{k=1}^n
 \affa(f_1\cdots f_{k-1}) \cdot \daffine f_k \cdot(f_{k+1}\cdots f_n)
\label{eqn:daffine-n}
  \end{align}

The next lemma is the heart of the results about the affine
derivative.

\begin{lemma} \label{lem:affine-faf}
  If $f\in\allpolyaffine$ and $\alpha$ is positive then $f-\alpha\affa
  f\in\allpolyaffine$.
\end{lemma}
\begin{proof}
  We may suppose that $f$ is monic of degree $n$, in which case the
  leading coefficient of $\affa f$ is $q^n$. Since $f \greateqeq \affa
  f$ we may write $\affa f = q^nf +h$ where $f\lesslesseq h$ so that
\begin{align*}
  f-\alpha\affa f &= (1-\alpha q^n)f - \alpha h 
\intertext{Writing $f = q^{-n}\affa f - k$ where $\affa f\lesslesseq
  k$ yields}
  f- \alpha\affa f &= (q^{-n}-\alpha)\affa f - k
\end{align*}
There are two cases, depending in the size of $\alpha$.

\textbf{Case 1.} If $0<\alpha<q^{-n}$ then $(1-\alpha q^{-n})$ and
$(q^{-n}-\alpha)$ are positive, so by Corollary~\ref{cor:where-roots} we have
\begin{align*}
f - \alpha \affa f =  (1-\alpha q^{-n})f-\alpha h & \longleftarrow f\\
f - \alpha \affa f =  (q^{-n}-\alpha)\affa f-k & \longleftarrow \affa f\\
\end{align*}

\textbf{Case 2} If $\alpha > q^{-n}$ then $(1-\alpha q^{-n})$ and
$(q^{-n}-\alpha)$ are negative, so by Corollary~\ref{cor:where-roots} we have
\begin{align*}
f - \alpha \affa f =  (1-\alpha q^{-n})f- \alpha h & \longrightarrow f\\
f - \alpha \affa f =  (q^{-n}-\alpha)\affa f-k & \longrightarrow \affa f\\
\end{align*}

In either case we can apply the affine criteria of Lemma~\ref{lem:affine-why} to
conclude that $f - \alpha \affa f \in\allpolyaffine$.

\end{proof}

\begin{remark}\label{rem:f+aAf}
     The proof of the lemma did not use an analysis of the
    roots. If we look at the roots then we get another proof of the
    lemma, as well as information about the behavior when $\alpha$ is
    positive. 

    If the roots of $f$ are $r_1<\cdots<r_n$ then we know that
$$ \affa^{-1}r_1 < r_1 < \affa^{-1}r_2 < r_2 < \cdots < \affa^{-1}r_n
< r_n.$$ 
If $\alpha$ is negative then the roots $\{s_i\}$ of $f + \alpha\affa f$ lie in
the intervals
$$ (-\infty,\affa^{-1}r_1), ( r_1 , \affa^{-1}r_2), \cdots ,(r_{n-1}, \affa^{-1}r_n)
$$ 
and consequently $s_i \le \affa^{-1}r_i < r_i\le s_{i+1}$. This shows
that $\affa s_i < s_{i+1}$, and hence $f+\alpha \affa
f\in\allpolyaffine$.

If $\alpha$ is positive the roots lie in the intervals
$$ (\affa^{-1}r_1 , r_1), ( \affa^{-1}r_2 , r_2),  \cdots ,( \affa^{-1}r_n,
 r_n)$$ 
The roots $\{s_i\}$ satisfy $s_i \le r_i \le \affa^{-1} r_{i+1} \le
s_{i+1}$. We can not conclude (indeed it's not true) that $\affa s_i \le
s_{i+1}$ from these inequalities. 

If we make stronger assumptions about the roots of $f$ then we can
draw conclusions about $f+\alpha \affa f$. Suppose that $f\greateq
\affa^2 f$. This means that the roots of $f$ satisfy $\affa^2 r_i <
r_{i+1}$. The roots $s_i$ now satisfy $s_i \le r_i, \affa^{-2}
r_{i+1}\le s_{i+1}$. Since the interlacing hypothesis implies that
$\affa r_i \le \affa^{-1}r_{i+1}$ we conclude that $\affa s_i \le
s_{i+1}$. A similar argument establishes the following result:
  \end{remark}

\begin{lemma} \label{lem:affine-faf-a}
  If $k$ is a positive integer greater than $1$, $f \greateq
  \affa^k f$ 
  and $\alpha$ is positive then $g = f+\alpha\affa f$ satisfies $g
  \greateq \affa^{k-1}g$.
\end{lemma}

The affine derivative has all the expected properties.

\begin{lemma}  \label{lem:affine-d}
  Suppose  $f\in\allpolyaffine$.
  \begin{enumerate}
  \item $\daffine f\in\allpolyaffine$.
  \item $f \lesslesseq \daffine f$.
  \item $\affa f \lesslesseq \daffine f$
  \item If $f,g\in\allpolyaffine$ satisfy $f\greateqeq
    g$ then $\daffine f \greateqeq \daffine g$.
  \end{enumerate}
\end{lemma}

\begin{proof}
  Lemma~\ref{lem:affine-faf} implies that $f-\affa f\in\allpolyaffine$, so the
  first result follows from Lemma~\ref{lem:order-ideal}.
  Since $f\greateqeq f - \affa f$ and the smallest root of $f-\affa f$
  is the fixed point, we can divide out by $x-\affa x$ and conclude
  from Lemma~\ref{lem:order-ideal} that $f\lesslesseq\daffine f$. The third part
  is similar.
  
  In order to show that $\daffine $ preserves interlacing it suffices
  to show that \\ $\daffine (x+a)f\lesslesseq \daffine f$ where
  $(x+a)f\in\allpolyaffine$. From the product rule
  \eqref{eqn:affine-product-2} we find
  \begin{equation*} 
    \daffine (x+a)f = (x+a)\daffine f + \affa f
  \end{equation*}
  Since $(x+a)\daffine f \lesslesseq \daffine f$ and $\affa f
  \lesslesseq \daffine f$ we can add these interlacings to conclude
  $\daffine (x+a)f\lesslesseq \daffine f$.
  \end{proof}

$\allpolyaffine$ is closed under  ordinary  differentiation.

\begin{lemma}
  If $f\in\allpolyaffine$ then $f^\prime\in\allpolyaffine$
\end{lemma}
\begin{proof}
  Since $f\in\allpolyaffine$ we know that $f + \alpha \affa
  f\in\allpoly$ for all $\alpha$. Since 
  $\frac{d}{dx}\affa f = q\,\affa(f^\prime)$ we have
  that $f^\prime + \alpha q \affa f^\prime\in\allpoly$ for all
  $\alpha$ and so $f^\prime\in\allpolyaffine$.
\end{proof}

Only in certain special cases do $\affa$ and $\polard{\affa}$ commute.
\begin{lemma} \label{lem:affine-commute2}
  For any positive integer $n$ we have that $q^n \affa^n \qderiv(f) =
  \qderiv(\affa^n f).$ Consequently, if $\affa \polard{\affa} =
  \polard{\affa} \affa$ then $\affa x = x+b$.
\end{lemma}
\begin{proof}
First note that 
$$ \affa(\affa x-x) = \affa(qx+b-x) = (q-1)(qx+b)+b = q(qx+b-x) =
q(\affa x-x).$$
This implies that for all positive $n$ we have $\affa^n(\affa x-x) = q^n(\affa
x-x)$, and hence $\affa^{-n}(\affa x -x) = q^{-n}(\affa x-x)$. We can
now compute:
\begin{align*}
  \qderiv(\affa^n f) &= \frac{ \affa^{n+1} f - \affa^n f}{\affa x - x}
  \\
&= \affa^n \, \frac{\affa f - f}{\affa^{-n}(\affa x-x)} \\
&= q^n \affa^n \qderiv f
\end{align*}
\end{proof}

We next consider properties of $f(\affa)g$ and $f(\daffine)g$. We first
observe this immediate consequence of Lemma~\ref{lem:affine-faf}.

\begin{lemma} \label{lem:affine-B1}
    If $g\in\allpolyalt$ and $f\in\allpolyaffine$ then
  $g(\affa)f\in\allpolyaffine$. 
\end{lemma}

\begin{cor}\label{cor:affine-faf-a}
  If $g\in\allpolypos(r)$ and $f \greateq \affa^{r+1}f$ then
  $g(\affa)f\in\allpolyaffine$. 

 If $f\greateq \affa^{r}$ then $\sum_0^r
    \binom{r}{k}\affa^k\,f\in\allpoly$. 

\end{cor}
\begin{proof}
  Factor $g$, and apply  Lemma inductively. For the second
  part, apply the first part to $g=(1+x)^r$. 
\end{proof}


Next we have an additive closure property of the affine derivative.
\begin{lemma} \label{lem:affine-B}
  If $g\in\allpolyaffine$ and $\alpha$ is positive then
  $g-\alpha\daffine g\in\allpolyaffine$.
\end{lemma}
\begin{proof}
  We know that $g\lesslesseq \daffine g$ so $g - \alpha \daffine
  g\greateqeq g$. Also, $g \greateqeq \affa g \lesslesseq \daffine g$,
   so $g-\alpha\daffine g\greateq \affa g$. Now apply
  Lemma~\ref{lem:affine-why}.
\end{proof}

\begin{lemma} \label{lem:affine-fg}
  If $g\in\allpolyalt$ and $f\in\allpolyaffine$ then
  $g(\daffine)f\in\allpolyaffine$. 
\end{lemma} 
\begin{proof} 
    Write $g(x) = (x+a_1)\dots(x+a_n)$ where all $a_i$ are negative.
    Since $$g(\daffine)f = 
  (\daffine+a_1)\dots(\daffine+a_n)f$$ we apply Lemma~\ref{lem:affine-B} inductively.
\end{proof}

  Here's a simple property about determinants.

  \begin{lemma}\label{lem:affine-det}
     Suppose $f,\affa f\in\allpolyaffine\cap\allpolypos$.

$$ \text{If } 
\begin{cases}
  \affa x > x \text{ for all positive $x$}& \\
  \affa x < x \text{ for all positive $x$}& 
\end{cases}
\text{ then } 
\twodet{ f}{\affa f}{\affa f }{\affa^2 f} \text{ is }
\begin{cases}
  \text{ negative for positive $x$.} & \\
  \text{ positive for positive $x$.} & 
\end{cases}
$$     

The same inequalities hold for 
$\twodet{f}{\daffine f}{\daffine f}{\daffine^2 f}$.
  \end{lemma}
  \begin{proof}
    We begin with the first statement.  We need to show that $ f \cdot
    \affa^2 f < \affa f \cdot \affa f$ for positive $x$. Since $\affa
    f$ is in $\allpolypos$, $\affa f$ is positive for positive $x$, so
    it suffices to show that $(\affa f/ \affa^2 f) > (f/\affa f)$.

    Since $f\in\allpolysep$ we know that $f \greateqeq \affa f$, and
    so there is a positive constant $c$ and non-negative $a_i$ such that
\begin{align*}
    \affa f &= c f + \sum a_i\frac{f}{x+r_i}
\intertext{where the roots of  $f$ are $\{r_i\}$. Thus }
 \frac{\affa f}{f} &= c + \sum \frac{a_i}{x+r_i} \\
\left(\frac{\affa f}{f}\right)' &= - \sum a_i \frac{1}{(x+r_i)^2}
\end{align*}
and consequently $\affa f/f$ is decreasing. Since $\affa x > x$, it
follows that
$$ \frac{f}{\affa f} (x) < \frac{f}{\affa f} (\affa x) = \frac{\affa
  f}{\affa^2 f} (x)   
$$

If $\affa x < x$ then the inequality is reversed. The last statement
follows from the above since 
$$ \twodet{ f}{\affa f}{\affa f }{\affa^2 f} =
\frac{1}{(\affa x - x)^2}
\twodet{f}{\daffine f}{\daffine f}{\daffine^2 f}.$$
\index{partial fractions}
\end{proof}

\begin{remark}
  $\arising{x}{a}{n}$ in \eqref{eqn:affine-1} is uniquely determined
by the three conditions below. 

$
\begin{array}{crl}
(1)\ &\ \arising{x}{a}{0}& =\ 1\\[.2cm]
(2)\ &\ \arising{a}{a}{n}& =\ 0\\[.2cm]
(3)\ &\ \daffine\,\arising{x}{a}{n}& =\  [n]\, \affa\, \arising{x}{a}{n-1}
\end{array}
$

This characterization of $\arising{x}{a}{n}$ 
is motivated a similar result by \cite{kac} in the case $\affa
x = qx$.

\end{remark}

\section{Sign interlacing for $\allpolyaffine$}
\label{sec:sec-sign-a}

\begin{figure}[htbp] \label{fig:aff-rts}
  \centering
      \includegraphics*[width=1.5in]{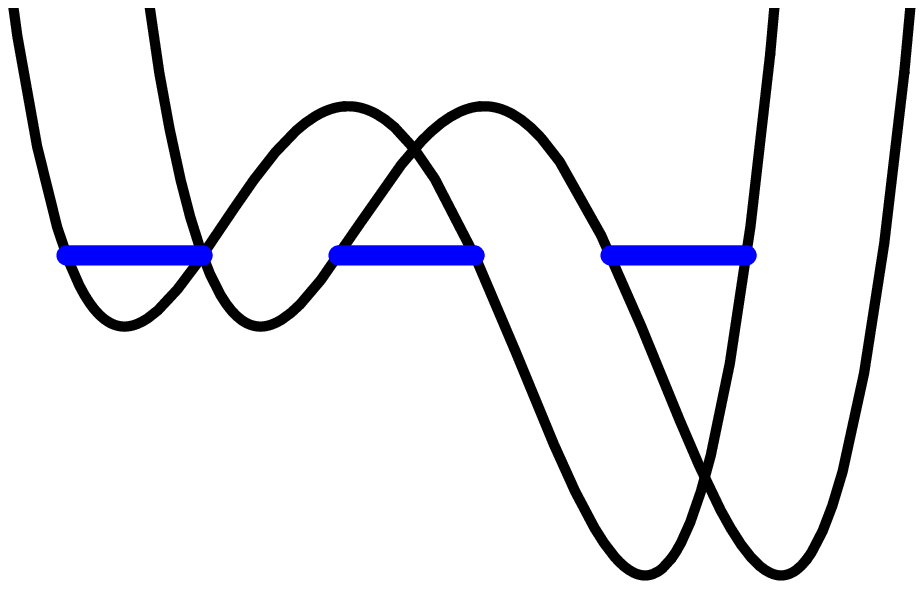}
  \caption{Understanding $\affa$}
\end{figure}

Let's look at some graphs   in order to better understand
$\allpolyaffine$. Suppose that $f\greateqeq \affa f$. The graph of
$\affa f$ is the graph of $f$, with some shifting to the left and
scaling. For instance, if $\affa x = x+1$ and $f=x(x-1)(x-3)(x-5)$
then the graph of $f$ and $\affa f$ is in Figure~\ref{fig:aff-rts}.

If $f\lesslesseq g$ and $\affa f \lessless g$ then the roots of $g$
must lie in the thickened segments of Figure~\ref{fig:aff-rts}. It's
clear that the roots of $g$ in the figure differ by at least one, as
is guaranteed by Lemma~\ref{lem:affine-why}. 

If we translate the qualitative aspects of this graph, we get a
version of sign interlacing. Assume that the roots of $f$ are
\begin{gather}
  \label{eqn:aff-rts-1}
  r_1 < r_2 < \cdots < r_n \\
\intertext{The roots of $\affa f$ are}
\label{eqn:aff-rts-2}
\affa^{-1} r_1 < \affa^{-1} r_2 < \cdots < \affa^{-1}r_n\\
\intertext{Since $f\greateqeq \affa f$ we can combine these orderings}
\label{eqn:aff-rts-3}
\affa^{-1} r_1 < r_1 \le  \affa^{-1} r_2 < r_2 \le 
\cdots \le \affa^{-1} r_n < r_n
\end{gather}

If $f\lesslesseq g$ and $\affa f \lesslesseq g$ all have positive
leading coefficients then it is clear from the graph that $g(r_n)>0$,
$g(\affa^{-1} r_n)>0$, $g(r_{n-1}<0$, $g(\affa^{-1}r_{n-1})<0$, and
generally
\begin{equation}
  \label{eqn:aff-si-1}
  (-1)^{n+i} g(r_i)>0 \quad\quad (-1)^{n+i} g(\affa^{-1}r_i) >0
\end{equation}
This is sign interlacing for $\affa$. We clearly have
\index{sign interlacing!for $\affa$}

\begin{lemma}\label{lem:sign-interlace-a}
If $g$ is a polynomial of degree $n-1$ that satisfies
\eqref{eqn:aff-si-1} then $g\in\allpolyaffine$.
  
\end{lemma}

The following polynomials play the role of $f(x)/(x-r_i)$ in
Lemma~\ref{lem:sign-quant}. For $1\le k \le n$ define

\begin{equation}
  \label{eqn:aff-rts-4}
  f_{[k]} = \affa\bigl[ (x-r_1)\cdots(x-r_{k-1})\bigr] \cdot \bigl[
  (x-r_{k+1}) \cdots (x-r_n)\bigr]
\end{equation}

These are polynomials of degree $n-1$ whose roots are 
\begin{equation}
  \label{eqn:aff-rts-5}
  \affa^{-1}r_1 < \affa^{-1}r_2 \cdots 
\affa^{-1}r_{k-1} <  r_{k+1} < \cdots < r_n
\end{equation}

Since the $r_i$ are all distinct, it is easy to see that the $f_{[k]}$ are
independent. Thus, they form a basis for the
polynomials of degree $n-1$.

\begin{lemma} \label{lem:affine-qsi}
  Suppose that $f\in\allpolyaffine$ and $ g = \sum a_k \, f_{[k]}$
  where all $a_i$ are non-negative.  Then $g\in\allpolyaffine$ and
  $f\lesslesseq g$ and $\affa f\lesslesseq g$.
\end{lemma}

\begin{proof}
  The key point is that $f \lesslesseq f_{[k]}$ and $\affa f \lesslesseq
  f_{[k]}$ - this follows from \eqref{eqn:aff-rts-1},
  \eqref{eqn:aff-rts-2}, \eqref{eqn:aff-rts-3},
  \eqref{eqn:aff-rts-4}. If all the $a_i$ are non-negative then we can
  add these interlacings to finish the proof.
\end{proof}

The fact that $ f \lesslesseq \daffine f$ follows from
Lemma~\ref{lem:affine-qsi}.  From \eqref{eqn:daffine-n} we see that
if $f(x) = (x-r_1)\cdots(x-r_n)$ then 

\begin{align*}
  \daffine(f) &= \daffine(x-r_1)\cdots(x-r_n) \\ &=
\sum_{k=1}^n \affa( \,(x-r_1)\cdots (x-r_{k-1})\,) \cdot 1 \cdot  (\,
(x-r_{k+1})\cdots(x-r_n)\,) \\
&= f_{[1]} + \cdots + f_{[n]}
\end{align*}
Consequently, this is another proof that $f\lesslesseq \daffine f$.

\section{The affine exponential}

  We can formally construct an affine exponential $E_\affa$ function,
  but the series might not converge. We want two simple 
  conditions

\begin{enumerate}[(1)]
\item $E_\affa(0) = 1$
\item $\daffine{}(E_\affa) = \affa\,E_\affa$
\end{enumerate}

The relation $\daffine{} E_\affa = \affa E_\affa$
determines $E_\affa$ up to a constant factor that is settled by
assuming $E_\affa(0)=1$.  If we write
\begin{align*}
  E_\affa(x) &= 1 + a_1 x + a_2 x^2 + \cdots \\
  \intertext{ and apply $\daffine{} E_\affa = \affa E_\affa$ we find}
  a_{n-1} \affa(x^{n-1}) &= a_n \daffine{}(x^n) \\
  \intertext{and recursively solving for the coefficients yields}
  E_\affa(x) &= \sum_{n=0}^\infty \prod_{k=1}^n
  \frac{\affa^{k-1}x}{[k]_\affa} \\
\intertext{where} [k]_\affa &=
  (\affa^k x-x)/(\affa x -x) = [k] = q^{k-1} + q^{k-2} +\cdots+ q+1.
\intertext{and}
  \affa^k x &= q^k x + b [k+1]
\end{align*}

\noindent%

Now this derivation is not valid if $b\ne0$, but we still can use this
definition since we have that
$$
\daffine \prod_{k=1}^n   \frac{\affa^{k-1}x}{[k]_\affa} = 
\affa \prod_{k=1}^{n-1}   \frac{\affa^{k-1}x}{[k]_\affa} 
$$
and thus term by term differentiation establishes condition (2). 
In particular

\begin{xalignat}{2}
  \affa x &= x \hspace*{1cm} & E_\affa(x)  &= \sum_{n=0}^\infty
  \frac{x^n}{n!}\label{eqn:exp-affa-1}\\
  \affa x &= qx  & E_\affa(x) &=
  \sum_{n=0}^\infty\frac{q^{\binom{n}{2}}}{[n]!} x^n  \label{eqn:exp-affa-3}\\
  \affa x &= q(x+1)  & E_\affa(x) &=
  \sum_{n=0}^\infty\frac{\prod_1^n (q^{k-1}x+q[k])}{[n]!}  \label{eqn:exp-affa-4}\\
  \affa x & = x+1  & E_\affa(x)  &= \sum_{n=0}^\infty
  \frac{\rising{x}{n}}{n!}\label{eqn:exp-affa-2}
\end{xalignat}

In the first case $E_\affa$ is the usual exponential.  The second one is the $q$-exponential of
\eqref{eqn:q-exponential} and has a representation as an infinite
product. The last two 
have  product representations for their finite sums:

\begin{xalignat}{2}
  \affa x & = qx & 
\sum_{n=0}^\infty \frac{q^{\binom{n}{2}}}{[n]!}x^n & =
\prod_{k=0}^\infty\left(1 - q^k(1-q)x\right)\label{eqn:af-prod-1}\\
  \affa x & = q(x+1) & 
\sum_{n=0}^N \left( \prod_{k=1}^n \frac{\affa^{k-1}x}{[k]}\right) & =
\prod_{k=1}^N\left( 1 + \frac{q^{k-1} x}{[k]}\right) \checked\label{eqn:af-prod-2}\\
\affa x &= x+1 & 
\sum_{n=0}^N  \frac{\rising{x}{n}}{n!} &= 
 \prod_{k=1}^N\left( 1 + \frac{x}{k}\right) = \frac{\rising{x+1}{N}}{N!}
\label{eqn:af-prod-3}
\end{xalignat}

Unlike the partial sums of the exponential function, the partial sums
of \eqref{eqn:exp-affa-2} and 
\eqref{eqn:exp-affa-4} have all real roots, as the product
representation shows. Since the sum $\sum q^{k-1}/[k]$ converges for
$|q|<1$ and diverges for $|q|>1$, it follows that the product
\eqref{eqn:af-prod-2} converges for $|q|< 1$ and diverges for $|q|\ge
1$. Consequently, equation \eqref{eqn:exp-affa-4} is only valid for $|q|<1$.
The series in \eqref{eqn:exp-affa-2} does not converge. 

The partial sums in the last two cases satisfy the exponential
condition (1) and nearly satisfy (2). Let $E_{\affa,N}$ denote the
partial sums in \eqref{eqn:af-prod-2} and \eqref{eqn:af-prod-3}. Then

\begin{xalignat*}{2}
  \affa x & = q(x+1) & \daffine(E_{\affa,N}) &= \affa\, E_{\affa,N-1}\\
  \affa x &= x+1 & \Delta(E_{\affa,N}) &= \affa\, E_{\affa,N-1}
\end{xalignat*}

\section{$\affa$-interlacing and the symmetric derivative}
  \label{sec:affine-int}
\index{interlacing!$\affa$-}

We say that $f$ and $g$ $\affa$-interlace if $f+\alpha
g\in\allpolyaffine$ for all $\alpha\in\reals$.  If the degree of $f$
is greater than the degree of $g$ and $f$ and $g$ $\affa$-interlace
then we write  $f\alesslesseq g$.  We give a criterion for affine
interlacing, and use it to determine when the symmetric affine derivative
is in $\allpolyaffine$. 
\index{symmetric derivative}

\begin{lemma}\label{lem:fandAf}
  If $f\in\allpolyaffine$ and $f\greateqeq \affa^2 f$ then $f$ and
  $\affa f$ $\affa$-interlace.
\end{lemma}
\begin{proof}
  From the interlacings $f\greateqeq \affa f$, $f\greateqeq \affa^2 f$,
  $\affa f\greateqeq \affa^2 f$ we conclude that if $\alpha>0$
  \begin{align*}
    f & \greateqeq f + \alpha \affa f\\
    f & \greateqeq \affa f + \alpha \affa^2 f\\
\intertext{and if $\alpha<0$}
    f & \lesseqeq f + \alpha \affa f\\
    f & \lesseqeq \affa f + \alpha \affa^2 f
  \end{align*}
Thus $f+\alpha \affa f$ and $\affa(f+\alpha \affa f)$ always have a
common interlacing, so we find that $f+\alpha \affa f\in\allpolyaffine$.
\end{proof}

Here is  an example of  polynomial $Q$ that $\affa$-interlaces
$\affa Q$. If $\affa$ is increasing on $\reals$ then define
\begin{align*}
  Q & = (x+1)(x+\affa^2 1)(x+\affa^4 1) \cdots (x+\affa^{2n}1)\\
\intertext{and if $\affa$ is decreasing define}
  Q & = (x-1)(x-\affa^2 1)(x-\affa^4 1) \cdots (x-\affa^{2n}1)\\
\end{align*}
If we apply $\affa^2$ to $Q$ we have even powers from $\affa^2$ to
$\affa^{2n+2}$, so clearly $Q\greateqeq \affa^2 Q$.

Here is a variation on the lemma for two polynomials.

\begin{lemma}
  If $f\lesslesseq g$ and $fg\in\allpolyaffine$ then $f\alesslesseq g$.
\end{lemma}
\begin{proof}
Suppose that the roots of $f(x)$ are
$f_1<\cdots$ and the roots of $g(x)$ are $g_1<\cdots$. The hypotheses
imply that the roots appear in the order
\[
\cdots g_i \le \affa g_i \le f_i \le \affa f_i \le g_{i+1} \le \affa g_{i+1}
\le \affa f_1 \cdots 
\]
from which it follows that $f\greateqeq \affa g$. For positive
$\alpha$ (first column) and negative $\alpha$ (second column) we have
\begin{align*}
  f \greateqeq & f + \alpha g & f \lesseqeq & f + \alpha g \\
  f \greateqeq & \affa f + \alpha \affa g &  f \lesseqeq & \affa
  f + \alpha \affa g \\ 
\end{align*}
Consequently, $f+\affa g\in\allpolyaffine$.
\end{proof}

Clearly $f\alesslesseq g$ implies that $f\lesslesseq g$ but
$\alesslesseq$ is a more restrictive condition than $\lesslesseq$, and
the converse does not usually hold.  We say that $f \alessless g$ if
$f \alesslesseq g$, and $f$ and $g$ have no common factors. Using sign
interlacing we can replace $\lesslesseq$ and $\alesslesseq$ in the
lemma by $\lessless$ and $\alessless$.  There is a simple criterion
for such interlacing.

\begin{lemma}
  Suppose that $f\lessless g$. Then $f \alessless g$ if and only if
  $\smalltwodet{f}{g}{\affa f}{\affa g}$ is never zero.
\end{lemma}
\begin{proof}
   From Lemma~\ref{lem:inequality-4b} we have that $f+\alpha g \greateq \affa f
   + \alpha \affa g$ for all $\alpha\in\reals$.  This means that
   $f+\alpha g\in\allpolyaffine$ so the conclusion follows.
\end{proof}

  The \index{symmetric derivative}symmetric derivative is 
  \begin{align*}
\daffine^{sym}(f) &= \frac{\affa(f) - \affa^{-1}(f)}{\affa(x) - \affa^{-1}x}
\intertext{For $\affa x = x+1$ this is the symmetric difference}
\daffine^{sym} &= \frac{f(x+1) - f(x-1)}{2} \\
\intertext{and for $\affa x = qx$ it is}
\daffine^{sym}(f) &= \frac{f(qx) - f(x/q)}{(q-1/q)x}.
  \end{align*}

  \begin{lemma}
    If $f\in\allpolyaffine$ then $\daffine^{sym}f\in\allpoly$. If in
    addition we have that $f\greateqeq \affa^2 f$ then 
$\daffine^{sym}f\in\allpolyaffine$.
  \end{lemma}
  \begin{proof}
Since $f\in\allpolyaffine$ we know that $f\greateqeq \affa f$ and
consequently $f-\affa f\in\allpoly$, so $\daffine^{sym}f\in\allpoly$.
Next, assume that $f\greateqeq \affa^2 f$. We have the interlacings
\begin{gather*}
  \affa^{-1} f  \greateqeq f  \greateqeq \affa f \\
  \affa^{-1} f  \greateqeq \affa^{-1} f  \greateqeq \affa f \\
\intertext{which imply that}
 f \lesslesseq \affa f - \affa^{-1} f \\
 \affa^{-1} f \lesslesseq \affa f - \affa^{-1} f 
\end{gather*}
The last two interlacings imply $\affa f - \affa^{-1}
f\in\allpolyaffine$, and so the conclusion follows.
  \end{proof}

Here's an example that shows that the symmetric derivative is not in
$\allpolyaffine$ if the second condition is not satisfied. If $\affa x
=x+1$ then
\[
\daffine^{sym} \falling{x}{n} = n(x-\frac{n-1}{2})\falling{x-1}{n-2}.
\]
If $n$ is odd then the symmetric derivative has a double root, and so
is not in $\allpolysep$.

\section{Linear transformations}
\label{sec:affine-linear}

In this section we investigate transformations that  map some
subset of $\allpoly$ to $\allpolyaffine$. The key idea is the
iteration of the linear transformation $f\mapsto
(x-\alpha)\affa^{-1}f-\beta f$.

\begin{lemma} \label{lem:affine-trans}
  The linear transformation $S\colon{}f\mapsto (x-\alpha)\affa^{-1}f - \beta
  \,f$ where $\beta\ge0$ maps
$$ S\colon{} \allpolyaffine\cap\allpolyint{(\alpha,\infty)}\longrightarrow
\allpolyaffine\cap\allpolyint{(\alpha,\infty)}
$$

\end{lemma}
\begin{proof}
{Since $f\in\allpolyaffine$ we know}
\begin{align}
  \affa^{-1}f & \greateqeq f\notag\\
  \intertext{and using $f\in\allpolyint{(\alpha,\infty)}$ we get}
  (x-\alpha)\affa^{-1}f & \lesslesseq \affa^{-1}f.\label{eqn:sep-3}\\
  (x-\alpha)\affa^{-1}f & \lesslesseq f.\label{eqn:sep-1}\\
  \intertext{By hypothesis $\beta\ge0$, so we can apply
    Corollary~\ref{cor:where-roots} to \eqref{eqn:sep-1} and
    Lemma~\ref{lem:add-interlace} to 
    \eqref{eqn:sep-3} yielding}
  (x-\alpha)\affa^{-1}f - \beta f & \lesslesseq f \notag\\
  (x-\alpha)\affa^{-1}f - \beta f & \lesslesseq \affa^{-1}f \notag\\
  \intertext{which shows that $S(f)\in\allpolyaffine$. Since
    $\affa^{-1}f\in\allpolyint{(\alpha,\infty)}$ we can apply
    Corollary~\ref{cor:where-roots} to \eqref{eqn:sep-1}} (x-\alpha)\affa^{-1}f
  -\beta f & \greateqeq (x-\alpha)\affa^{-1}f \notag
\end{align}
{and conclude that $S(f)\in\allpolyint{(\alpha,\infty)}$}.
\end{proof}

Now, define polynomials by $Q_n^{\alpha,\beta}(x) =
\left(\,(x-\alpha)\affa^{-1}-\beta\,\right)^n(1)$.  Some particular
examples are given in Table~\ref{tab:qn}. They are immediate from
their definitions, except for the Charlier polynomials, where we apply
the recurrence \eqref{eqn:charlier-sa}.


\begin{table}[htbp]
  \centering
  \begin{tabular}{|cccc|} \hline
  $\affa x$ &$ \alpha $&$ \beta $&$ Q_n^{\alpha,\beta} $ \\[.2cm]\hline
$  x+1 $&$ 0 $&$ 0 $&$ \falling{x}{n} $\\[.2cm]
$  x+1 $&$ 0 $&$ \beta $&$ C_n^\beta(x)$\\[.2cm]
$  qx $&$ 0 $&$ 0 $&$ q^{-\binom{n}{2}}x^n$ \\[.2cm]
$  qx $&$ -1 $&$ 0 $&$ q^{-\binom{n}{2}}(-x;q)_n$ \\\hline
\end{tabular}

\caption{Polynomials of the form $Q_n^{\alpha,\beta}$}
  \label{tab:qn}
\end{table}

\begin{prop}\label{prop:affine-trans}
  The map $T\colon{}x^n\mapsto Q^{\alpha,\beta}_n(x)$ is a linear
  transformation
  $$\allpolyint{(-\beta,\infty)}\longrightarrow
  \allpolyint{(\alpha,\infty)}\cap\allpolyaffine$$.
\end{prop}
\begin{proof}
  If we define $S(f) = (x-\alpha)\affa^{-1}f - \beta \,f$, then from
  the definition of $T$ we have that $T(x^n) = S^n(1) =
  Q^{\alpha,\beta}_n$, and so by linearity  $T(f) = f(S)(1)$. If we
  choose $f = \prod(x-r_i)$ in $\allpolyint{(-\beta,\infty)}$ then
  each root $r_i$ satisfies $r_i\ge -\beta$.  Notice that
$$
    T(f) = \prod(S-r_i) 
= \prod((x-\alpha)\affa^{-1} - \beta -r_i)(1)
$$
  Since $\beta+r_i\ge0$, if follows from the lemma that
  $T(f)\in\allpolyint{(\alpha,\infty)}\cap\allpolyaffine$.
\end{proof}

Here is a general criterion to tell if a transformation preserves
$\allpolyaffine$. See Lemma~\ref{lem:affine-cos} for an application. 

  \begin{lemma}\label{lem:affine-commute}
    If $T\colon{}\allpoly\longrightarrow\allpoly$ preserves interlacing, and
    $T\affa = \affa T$, then
    $T\colon{}\allpolyaffine\longrightarrow\allpolyaffine$. 
  \end{lemma}
  \begin{proof}
    Assume that $f\greateqeq \affa f$. Since $T$ preserves
    interlacing we have $Tf \greateqeq T \affa f$. Commutativity
    implies that $T \affa f =\affa T f$, and so $T(f)\greateqeq
    \affa T(f)$, and hence $T(f)\in\allpolyaffine$.
  \end{proof}

\section{The case $\affa x = x+1$}
\label{sec:affine-difference}

In this section we study the affine operator\footnote{This operator
  is often denoted by \textsf{E}.} $\affa x=x+1$.  The space of all
polynomials $f$ satisfying $f(x) \greateqeq f(x+1)$ is denoted
$\allpolysep$, where \emph{sep} stands for \emph{separated}. Here are
some necessary and sufficient conditions for a polynomial $f$ to
belong to $\allpolysep$:
\begin{enumerate}
\item $f(x) \greateqeq f(x+1)$
\item The distance between any two roots of $f$ is at least $1$.
\item $f(x)$ and $f(x+1)$ have a common interlacing.
\item For all positive $\beta$ we have $f(x) + \beta f(x+1)\in\allpoly$.
\end{enumerate}

For example, both \index{falling factorial} $\falling{x}{n}$ and
\index{rising factorial} $\rising{x}{n}$ are in
$\allpolysep$. The affine derivative is the difference operator
$\Delta(f) = f(x+1)-f(x)$. The difference operator commutes with
$\affa$.  Since the derivative maps $\allpolysep$ to itself, the roots
of the derivative can't get too close:

\begin{lemma} \label{lem:deriv-diff}
  If the distance between any two roots of a polynomial $f$ in
  $\allpoly$ is at least $1,$ then the distance between any two roots
  of the derivative $f^\prime$ is also at least $1$. Similarly, the
  distance between any two roots of $f(x+1)-f(x)$ is at least one.
\end{lemma}
\begin{proof}
We know that $f\in\allpolysep$ implies that both $f^\prime$ and
$\Delta f$ are also in $\allpolysep$.
\end{proof}

If we express a polynomial in the basis \index{falling factorial}
$\falling{x}{n}$, and use the fact that $\Delta \falling{x}{n} = n\falling{x}{n-1}$ we conclude
\begin{quote}
  If $\sum a_i \falling{x}{i}$ has roots that are separated by at least one,
  then the same is true for $\sum i a_i \falling{x}{i-1}$.
\end{quote}

From Lemma~\ref{lem:affine-d} we find that differences preserve interlacing:

\begin{lemma}
  If $f,g\in\allpolysep$ satisfy $f\greateqeq g$ then $\Delta f
  \greateqeq \Delta g$. 
\end{lemma}

\begin{cor} 
  If $g\in\allpolyalt$ and $f\in\allpolysep$ then
  $g(\affa)f\in\allpolysep$ where $\affa(f(x))=f(x+1)$. Similarly, if
  $g\in\allpolyalt$ then
  $g(\Delta)f\in\allpolysep$. 
\end{cor}
\begin{proof}
  Use Lemma~\ref{lem:affine-fg} for the second part.
\end{proof}

The second part of this corollary does not hold for all $f$. For
instance, if $f=x(x+1)$ then $f-\alpha \Delta f = (x+1)(x-2\alpha)$,
and the distance between roots is $|1+2\alpha|$. Thus $f-\alpha \Delta
f\in\allpolysep$ if and only if $\alpha\not\in (-1,0)$. This is true
generally.

\begin{lemma}
  If $f\in\allpolyint{\reals\setminus(-1,0)}$ then
    $f(\Delta):\allpolysep\longrightarrow\allpolysep.$ 
\end{lemma}
\begin{proof}
  It suffices to show that $f-\alpha \Delta f\in\allpolysep$ if
  $f\in\allpolysep$ and $\alpha\not\in(-1,0)$. Note that
\[f-\alpha\Delta f = -(\alpha+1)\bigl( \frac{\alpha}{\alpha+1}
f(x+1) - f(x)\bigr)
\]
Since $\alpha\not\in(-1,0)$ implies $\frac{\alpha}{\alpha+1}>0$ we can
  apply Corollary~\ref{cor:diff-sep} to finish the
  proof.
\end{proof}

Here is a nice consequence of composition.  

\begin{lemma}\label{lem:binomial-xdk}
  Suppose that $d$ is a positive integer.  The linear transformation
  $T\colon{}x^k \mapsto \falling{x+d-k}{d}$ maps $\allpolyalt$ to $\allpolysep$.
\end{lemma}
\begin{proof}
  Let $\affb x = x-1$. The important observation is  that $\affb^k
  \falling{x+d}{d} = T(x^k)$. Thus, $T(f) = f(\affb)\falling{x+d}{d}$.
  Thus we only have to verify the conclusion for one factor. Since
  $f\in\allpolyalt$, we have to show that $(\affb - \alpha)$ maps
  $\allpolysep$ to itself if $\alpha$ is positive. This follows from
  Lemma~\ref{lem:affine-faf} 
\end{proof}

Simple geometry forces a relationship between $\Delta f$ and $f^\prime$.

\begin{lemma}
  If $f\in\allpolysep$ then $f^\prime \greateq \Delta f$.
\end{lemma}
\begin{proof}
  Consider Figure~\ref{fig:affine-delta}. The points $a,b$ are
  consecutive zeros of $f$, and satisfy $b-a\ge1$. Points $d,e$ are on
  $f$, and the horizontal line between them has length exactly 
  $1$. Finally, $c$ is a root of $f^\prime$. 

  The $x$-coordinate of $d$ is a root of $f(x+1)-f(x)=0$, and so is a
  root of $\Delta f$. Since the $x$ coordinate of $d$ is obviously
  less than $c$, it follows that $f^\prime \greateq \Delta f$.

\begin{figure}[htbp]
  \begin{center}
    \leavevmode
    \epsfig{file=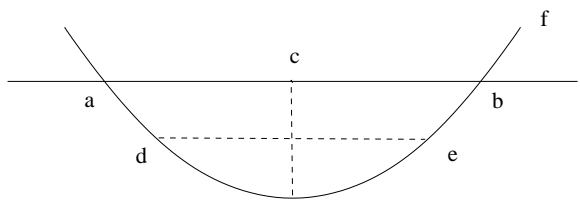,width=3in}
    \caption{The difference and the derivative}
    \label{fig:affine-delta}
  \end{center}
\end{figure}

\end{proof}

\begin{remark}
  If $f\in\allpolysep$ then we know that $f\lesslesseq f^\prime$,
  $f\lesslesseq \Delta f$, $f^\prime\in\allpolysep$, and $\Delta
  f\in\allpolysep$. Even though $f\lesslesseq f^\prime + \Delta f$, it
  is not necessarily the case that $f^\prime + \Delta
  f\in\allpolysep$. A simple counterexample is given by
  $f=x(x-1)(x-2)$.
\end{remark}

An important class of polynomials in $\allpolysep$ are the
Charlier polynomials $C_n^\alpha(x)$. %
\index{polynomials!Charlier}\index{Charlier polynomials}%
A simple consequence of the following result is that the Charlier polynomials
are in $\allpolysep$. 

\begin{cor} Assume $\alpha\le0$. \label{cor:affine-trans}
\begin{enumerate}
\item  The map $x^i\mapsto C_n^\alpha$ maps $\allpolyalt \longrightarrow
  \allpolyalt\cap\allpolysep$. 
\item  The map $x^i\mapsto \falling{x}{n}$ maps $\allpolyalt \longrightarrow
  \allpolyalt\cap\allpolysep$. 
\item  The map $x^i\mapsto \rising{x}{n}$ maps $\allpolypos \longrightarrow
  \allpolypos\cap\allpolysep$. 
\end{enumerate}
\end{cor}
\begin{proof}
  Apply Proposition~\ref{prop:affine-trans} to Table~\ref{tab:qn}. Switch signs for
  the rising factorial in the result for the falling factorial.
\end{proof}

Expressing this result in terms of coefficients yields
\begin{quote} \index{rising factorial}
  If $\sum a_i x^i$ has all negative roots, then 
     $\sum a_i \rising{x}{i}$ has roots whose mutual
  distances are all at least one.
\end{quote}

The following are consequences of  Corollary~\ref{cor:affine-faf-a}.

\begin{cor}
 If   $g\in\allpolypos(n-1)$ and all roots of
  $f$ are at least $n$ apart then $g(\affa)f\in\allpolysep$.
\end{cor}

\begin{cor}
  If the roots of $f$ are separated by at least $n+1$ and
  $f\in\allpolysep(n)$ then
$$\sum_{i=0}^n \binom{n}{i}f(n+i)\in\allpoly$$

If $f\in\allpolysep$ then
$$\sum_{i=0}^n (-1)^i\binom{n}{i}f(n+i)\in\allpoly$$
\end{cor}

\begin{proof}
  Apply the previous corollary to $(1+\affa)^nf$ and $(\affa-1)^nf$.
\end{proof}

The next result is an interesting special case of Lemma~\ref{lem:affine-fg}.

\begin{cor}
  If $g\in\allpolysep$ and $m$ is a positive integer then
$$ \sum_{i=0}^m (-1)^i\,\binom{m}{i}\,g(x+i) \in\allpolysep.$$
\end{cor}
\begin{proof}
  Apply Lemma~\ref{lem:affine-fg} with $f(x)=x^m$. 
\end{proof}

If we apply this result to \index{falling factorial} $\falling{x}{n}$ we get
$$ \sum_{i=0}^m (-1)^i\binom{m}{i}\falling{x+i}{n} \in\allpolysep.$$

The following corollary of Lemma~\ref{lem:affine-det} is a slight
strengthening of Corollary~\ref{cor:fi}.

  \begin{lemma}\label{lem:affine-det-2}
     Suppose $f\in\allpolysep\cap\allpolypos$. For positive $x$

     $$\smalltwodet{f(x+1)}{f(x)}{ f(x+2)}{f(x+1)} > 0$$
     
      \end{lemma}

      \begin{remark}
        This does not necessarily hold for all $x\in\reals$. If
        $f=x(x+1)$ then the determinant is  $2(x+1)(x+2)$ which is
        negative for $x\in(-2,-1)$. However, if $f=x(x+2)$ then the
        determinant is positive for all $x$. More generally, it
        appears that if $f(x)\greateqeq f(x+2)$ then the determinant
        is always positive. However, there are examples,
        e.g. $x(x+2)(x+3.5)$, where $f(x)$ and $f(x+2)$ don't
        interlace, yet the determinant is always positive.
      \end{remark}

\begin{example} \index{Rodrigues' formula!affine}
  Here is a family of polynomials whose definition is similar to the
  Rodrigues definition of the Legendre polynomials (See
  Question~\ref{ques:legendre}). Define  polynomials by the Rodrigues type formula
  $$
  p_n = \Delta^n\, \prod_{i=1}^n (x^2-i^2) $$
  Since $\prod_{i=1}^n
  (x^2-i^2) $ is in $\allpolysep$, the polynomial $p_n$ is in
  $\allpolysep$, and has degree $n$. Moreover, we claim that $p_{n+1}
  \lesslesseq p_n$. Since $\Delta$ preserves interlacing in
  $\allpolysep$ it suffices to show that
$$ \Delta\, \prod_{i=1}^{n+1} (x^2-i^2) \lesslesseq \prod_{i=1}^{n}
(x^2-i^2) $$ 
Now the roots of the left hand side are seen to be 
$$-n,-n+1,\cdots,-2,-1,-.5,1,2,\cdots,n$$ and the roots of the right
side are
$$ -n,-n+1,\cdots,-2,-1,1,2,\cdots,n$$
and so  they interlace. The $p_n$'s are \emph{not} orthogonal
polynomials since they do not satisfy a three term
recurrence. However, the $p_n$ appear to satisfy a recurrence
involving only three terms:

\begin{multline*} p_{n+2} = (n+2)^2 p_{n+1} +  x(2n+3) p_{n+1}
  \, + \\
(n+1)(n+2)\, \affa p_{n+1} + x(2n+3)\, \affa p_{n+1} -2(2n+3)(n+1)^3\,
\affa p_n
\end{multline*}
\end{example}

\index{interpolation!basis}

If we write a polynomial in the interpolation basis, then 
there is a simple estimate on the separation of roots. 
One consequence is that if $\delta(f)>0$ then $\delta(f')>\delta(f)$. 

\begin{lemma} \label{lem:interp-sep}
Suppose that $f = (x-a_1)\dots(x-a_n)$ where $a_1<\cdots<a_n$ and that 
$$
g(x) = bf(x) + \sum_i b_i\frac{f(x)}{x-a_i}$$
where $b\ge0$. If
$b_1\ge b_2 \ge \cdots \ge b_n>0$ then $\delta(g) > \delta(f)$.
\end{lemma}

\begin{proof}
  Choose consecutive roots $a_{j-1},a_j,a_{j+1}$, and points $x_1,x_2$
  satisfying $ a_{j-1} < x_1 < a_j < x_2 < a_{j+1}$. We will show that
  $\dfrac{g}{f}(x_2) - \dfrac{g}{f}(x_1) > 0$ if $|x_1-x_2| \le
  \delta(f)$.  It follows that for $|x_1-x_2|\le\delta(f)$, $x_1$ and
  $x_2$ can not be adjacent roots of $g$, so $\delta(g) > \delta(f)$.

We calculate
\begin{align*}
  \frac{g}{f}(x_2) - \frac{g}{f}(x_1) & = 
\sum_{i=1}^n \frac{b_i}{x_2-a_i} - \sum_{i=1}^n \frac{b_i}{x_1-a_i} \\
&=
\frac{b_1}{x_2-a_1} + \frac{b_n}{a_n-x_1} +
\sum_{i=1}^{n-1} \left(\frac{b_{i+1}}{x_2-a_{i+1}} -
  \frac{b_i}{x_1-a_i}\right)\\
\intertext{Now since $b_i \ge b_{i+1}$, we have}
&\ge 
\frac{b_1}{x_2-a_1} + \frac{b_n}{a_n-x_1} +
\sum_{i=1}^{n-1} b_{i+1}\left(\frac{1}{x_2-a_{i+1}} -
  \frac{1}{x_1-a_i}\right)\\
& =
\frac{b_1}{x_2-a_1} + \frac{b_n}{a_n-x_1} +
\sum_{i=1}^{n-1} b_{i+1}\frac{(a_{i+1}-a_i) - (x_2-x_1)}{(x_2-a_{i+1})(x_1-a_i)}
\end{align*}
The first two terms are positive; the numerator is positive since
$x_2-x_1 \le \delta \le a_{i+1}-a_i$, and the denominator is positive
since $x_1$ and $x_2$ lie in consecutive intervals determined by the
roots of $f$.
\end{proof}

\section{A multiplication in $\allpolysep$}
\label{sec:affine-diamond}

\index{product!diamond}

The product of two polynomials in $\allpolysep$ is not in
$\allpolysep$. We can remedy this situation by using a different
multiplication. We define the bilinear transformation $f \septimes g$
on a basis:

$$ \falling{x}{n} \septimes \falling{x}{m} =
 \falling{x}{n+m}$$

$\septimes$ is a diamond product - see \chapsec{linear}{bilinear-diamond-1}.
$\Delta$ and $\septimes$ interact as expected.

\begin{lemma}  
 $\Delta(f\septimes g) = \Delta f \septimes g + f \septimes
    \Delta g$
\end{lemma}  
  \begin{proof}
 We only need to check it on a basis. 
\begin{multline*}\Delta(\falling{x}{n}\septimes \falling{x}{m}) = \Delta
\falling{x}{n+m} = (n+m) \falling{x}{n+m-1} = \\ n
\falling{x}{n-1}\septimes\falling{x}{m} + m \falling{x}{n}\septimes\falling{x}{m-1}
\end{multline*}
\end{proof}

It follows by induction that

$$ \Delta(f_1\septimes\cdots\septimes f_n) = \sum_0^n
f_1\septimes\cdots \septimes \Delta{f_i}\septimes\cdots\septimes f_n$$

and in particular the difference of a product of linear terms is 

\[ \Delta((x-a_1)\septimes\cdots\septimes (x-a_n)) = \sum_0^n
(x-a_1)\septimes\cdots \septimes
\widehat{(x-a_i)}\septimes\cdots\septimes (x-a_n)
\]

\noindent
where the hat means ``omit''.

Every polynomial in $\allpolypos$ can be factored in terms of $\septimes$.

\begin{lemma}\label{lem:diam-sep-1}
  If $f\in\allpolypos(n)$ is monic then there are positive $a_1,\dots,a_n$ such that
  \[f(x) = (x+a_1)\septimes(x+a_2)\septimes\cdots\septimes(x+a_n).\]
\end{lemma}
\begin{proof}
  This follows from Proposition~\ref{prop:diamond-gen}
since the linear transformation \\ $T\colon{}\falling{x}{n}\mapsto x^n$ maps
$\allpolypos$ to itself.

\end{proof}

\begin{remark}
Here's an example of the Lemma:
$$ (x+3)(x+4) = (x+2)\septimes(x+6)$$
The  Lemma shows that the product of certain terms
  with positive coefficients is in $\allpoly$; there are products of linear
  terms with positive coefficients that are not in $\allpoly$:
$$ (x+1)\septimes(x+2) = (x+1-\imag)(x+1+\imag)$$
  If $f\in\allpolyalt$ then there might be no such representation. For
  example, the following is a unique representation:
  $$
  (x-1-\imag)\septimes(x-1+\imag) = (x-1)(x-2)$$
\end{remark}

\begin{lemma}\label{lem:sep-alt-times}
  If $f\in\allpolysep\cap\allpolyalt$ and $\beta\ge0$ then
$$ (x-\beta)\septimes f \in\allpolysep\cap\allpolyalt.$$  
\end{lemma}
\begin{proof}
  Let $T(f) = (x-\beta)\septimes f$. Since 
$$T(\falling{x}{n}) = \falling{x}{n+1} - \beta \falling{x}{n} = (x
\affa^{-1}-\beta) \falling{x}{n}$$
we see that $T(f) = (x\affa^{-1}-\beta)f$. 
Lemma~\ref{lem:affine-trans}  implies that
$T(f)\in\allpolysep\cap\allpolyalt$. 
\end{proof}

By iterating the Lemma, we get:

\begin{cor}
  If $a_1,\dots,a_n$ are positive then
$$ (x-a_1)\septimes(x-a_2)\septimes\cdots\septimes(x-a_n)
\in\allpolysep\cap\allpolyalt.$$ 
\end{cor}

We have seen that the following transformation preserved
$\allpolyalt$.

\begin{cor}
  The transformation $T\colon{}\falling{x}{n}\mapsto\falling{x}{n+1}$ maps
$$\allpolysep\cap\allpolyalt\longrightarrow\allpolysep\cap\allpolyalt.$$
\end{cor}
\begin{proof}
  Since $T(f) = x\septimes f$ we can apply Lemma~\ref{lem:sep-alt-times}.
\end{proof}

We would like to show that $\septimes$ preserves
$\allpolysep\cap\allpolyalt$. See Question~\ref{ques:sep-alt}.  This
follows easily from Question~\ref{ques:sep-bijection}, that the map
$x^n\longrightarrow\rising{x}{n}$ is a bijection between $\allpolypos$
and $\allpolysep\cap\allpolypos$.

\section{The analytic closure of $\allpolysep$}
\label{sec:allpolysepf}

The analytic closure $\allpolysepf$ is the uniform closure on compact
subsets of $\allpolysep$. Obviously $\allpolysepf\subsetneq \allpolyf$.
Since $\allpolysep$ is closed under differentiation and difference, so
is $\allpolysepf$. 

\begin{lemma} \label{lem:analytic-sep}
  The following functions are in $\allpolysepf$:
  \begin{enumerate}
  \item $\sin\pi x$, $\cos \pi x$.
  \item $e^{\alpha x}$
  \item $e^{\alpha x}f(x)$ for any $f\in\allpolysep$.
  \end{enumerate}
  $e^{-\alpha x^2}$ is not in $\allpolysepf$ for any  $\alpha$.
\end{lemma}
\begin{proof}
  From the infinite product
  \begin{align*}
    \sin \pi x &= \pi x \prod_{k=1}^\infty
    \left(1-\frac{x^2}{k^2}\right)\\
\intertext{we see that the polynomials}
 & \pi x \prod_{k=1}^n
    \left(1-\frac{x^2}{k^2}\right)\\
  \end{align*}
converge to $\sin \pi x$ and are in $\allpolysep$ since the roots are
$-n,-n+1,\dots,n-1,n$. A similar argument applies to $\cos \pi x$. 

Next, we show that for $\alpha>0$ we can find polynomials $p_n$ such
that
\begin{enumerate}
\item $p_n\in\allpolysep$
\item For every $r$ there is an $N$ such that $n>N$ implies that the
  roots of $p_n$ have absolute value at least $r$.
\item $\lim p_n = e^{\alpha x}$
\end{enumerate}
We start with the identity, valid for positive $\beta$:
\begin{equation}
  \label{eqn:exp-identity}
  (\beta+1)^x = \lim_{n\rightarrow\infty} \prod_{k=1}^n \left( 1 +
    \frac{\beta x}{n+\beta k}\right)
\end{equation}
Let $p_n$ be the product. The roots of $p_n$ are
$-(n/\beta)-n,\dots,-(n/\beta)-1$ so $p_n\in\allpolysep$. 
This shows that (1) and (2) are satisfied. If we choose  positive
$\beta$ such that $\beta+1=e^\alpha$ then $p_n$ also satisfies
(3). Thus, $e^{\alpha x}\in\allpolysepf$ for positive $\alpha$. Since
$g(x)\in\allpolysep$ if and only if $g(-x)\in\allpolysep$, it follows
that $e^{\alpha x}\in\allpolysepf$ for all $\alpha$. 

If $f\in\allpolysep$ then we can find an $N$ such that $n>N$ implies
that $p_n(x)f(x)\in\allpolysep$. Taking limits show that $e^{\alpha
  x}f(x)\in\allpolysepf$. 

\index{Rodrigues' formula!Hermite}

Finally, if $e^{-\alpha x^2}$ belongs to $\allpolysepf$, then so do
all its derivatives. The Rodrigues' formula for Hermite polynomials
\eqref{eqn:hermite-4} implies that the Hermite polynomials
$H_n(x/\sqrt{\alpha})$ would be in
$\allpolysep$. This means that all roots of $H_n(x)$ would be
separated by at least $\sqrt{\alpha}$. However, from
\cite{szego}{(6.31.21)} we see that the minimum distance between roots
of $H_n$ goes to zero as $n$ goes to infinity. 

Alternatively, if $e^{-\alpha x^2}\in\allpolysepf$ and $f\in
\allpolysep$ then differentiating $e^{-\alpha x^2}f$ and dividing out by
the exponential factor would imply that $f' - 2\alpha x
f\in\allpolysep$.  But if $f=x(x+1)(x+2)$ then $\delta(f' - 2 x f) =
.95$.
\end{proof}

The following consequence is due to Riesz \cite{stoyanoff}.
\begin{cor} 
  If $f\in\allpolysep$ and $\alpha\in\reals$ then $\alpha f + f'
  \in\allpolysep$. 
\end{cor}
\begin{proof}
  Since $e^{\alpha x}f(x)\in\allpolysep$ we differentiate and find
  that $(\alpha f+f')e^{\alpha x}\in\allpolysep$. The conclusion
    follows. 
\end{proof}

\begin{cor}\label{cor:fdf}
  If $g(x)\in\allpoly$ and $f\in\allpolysep$ then
  $g(\diffd)f\in\allpolysep$. 
\end{cor}
The next result is not an independent proof of
Lemma~\ref{lem:affine-faf} since that Lemma was used to show that
$\allpolysep$ is closed under differences. 

\begin{cor} \label{cor:diff-sep}
  If $f\in\allpolysep$ and $\beta>0$ then $\beta\,f(x+1) -
  f(x)\in\allpolysep$. 
\end{cor}
\begin{proof}
  Note that
$ \Delta ( e^{\alpha x} f) = e^{\alpha x}\left(e^\alpha f(x+1) -
  f(x)\right)
$.
\end{proof}

As long as function in $\allpolyf$ has no factor of $e^{-\alpha x^2}$,
and all roots are at least one apart then $f$ is in $\allpolysepf$. 

\begin{lemma}
  Assume that $a_{n+1}-a_n\ge 1$ for $n=1,2,\dots$. 
  \begin{enumerate}
  \item If $ f(x) = e^{\alpha x} \prod_{i=1}^\infty
        \left(1-\frac{x}{a_i}\right)\in\allpolyf$ then
        $f\in\allpolysepf$. 
  \item If $ f(x) = e^{\alpha x} \prod_{i=1}^\infty
        \left[\left(1-\frac{x}{a_i}\right)e^{-x/a_i}\right]\in\allpolyf$ then
        $f\in\allpolysepf$. 
      \item $1/\Gamma(x)\in\allpolysepf$.
  \end{enumerate}
\end{lemma}
\begin{proof}
  The assumptions in (1) and (2)  guarantee that the product
  converges. Part (1) follows from
  Lemma~\ref{lem:analytic-sep}(3). From the proof of this lemma we let
  $p(a,r)\in\allpolysep$ be a polynomial whose roots have absolute
  value at least $|r|$, and $\lim_{r\rightarrow\infty}p(a,r)=e^{ax}$.
  Consider the polynomial
$$ q_n(x) = p(\alpha,r_0)\, \prod_{i=1}^n
\left(\left(1-\frac{x}{a_i}\right) p(-1/a_i,r_i)\right)$$
We can choose the $r_i$ sequentially so that
$q_n\in\allpolysep$. Taking limits yields (2).
The Weierstrass infinite product representation
$$
\frac{1}{\Gamma(x)} = x e^{\gamma x}\prod_{n=1}^\infty
\left[\left(1+\frac{x}{n}\right)e^{-x/n}\right]$$
shows that $\frac{1}{\Gamma(x)}\in\allpolyf$. Since the roots are at
negative integers, part (2) implies that
$1/\Gamma(x)\in\allpolysepf$. 
\end{proof}

Since $\frac{1}{\Gamma(x)}$ is in $\allpolysepf$, so are all the
differences $\Delta^n\frac{1}{\Gamma(x)}$. We can find a recursive
formula for these differences. If we write
\begin{align*}
  \Delta^n\frac{1}{\Gamma(x)} &= \frac{f_n(x)}{\rising{x}{n} \Gamma(x)}\\
\intertext{then }
  \Delta^{n+1}\frac{1}{\Gamma(x)} &= \frac{f_n(x+1)}{\rising{x+1}{n}
    \Gamma(x+1)} - \frac{f_n(x)}{\rising{x}{n} \Gamma(x)}\\
&= \frac{f_n(x+1)-(x+n)f_n(x)}{\rising{x}{n+1} \Gamma(x)} 
\intertext{and consequently the numerators satisfy the recurrence}
f_{n+1}(x) &= f_n(x+1)-(x+n)f_n(x)
\end{align*}

Since $f_1(x)=1-x$, the next lemma shows that all of the numerators are in
$\allpolysep$.

  \begin{lemma}\label{lem:affine-recur-1}
    Suppose that $g_1(x)\in\allpolyint{(0,\infty)}\cap\allpolysep$ and
    the sequence $\{g_n\}$ is defined by the recurrence
$$ g_{n+1}(x) = g_n(x+1) - (x+n)g_n(x)$$

then $g_n(x)\in\allpolyint{(-n+1,\infty)}\cap\allpolysep.$
  \end{lemma}
  \begin{proof}
    
    The function $\dfrac{1}{\rising{x}{n+1}\Gamma(x)}$ is in
    $\allpolysepf$, as can be seen by removing initial factors from
    the infinite product;  its  zeros are $-n-1,-n-2,\cdots$.
    Since $\Delta^n\dfrac{g_1}{\Gamma(x)}=
    \dfrac{g_{n+1}}{\rising{x}{n+1}\Gamma(x)}$ is in $\allpolysepf$,
    the roots of $g_{n+1}$ must lie in $(-n,\infty)$.

  \end{proof}

If we consider $f/\Gamma$ we get
\begin{lemma}
  If $f\in\allpolyint{(1,\infty)}\cap \allpolysep$ then
$f(x+1)-x f(x)\in\allpolyint{(0,\infty)}\cap \allpolysep$.
\end{lemma}
\begin{proof}
  The largest root of $1/\Gamma(x)$ is $0$, so
  $f(x)/\Gamma(x)\in\allpolysepf$. Taking differences yields
\[
\frac{f(x+1)}{\Gamma(x+1)} - \frac{f(x)}{\Gamma(x)}, 
=
\frac{f(x+1) - x f(x)}{x \Gamma(x)} \in\allpolysepf
\]
and the numerator is in $\allpolysep$. The largest zero of $x
\Gamma(x)$ is $-1$, so the conclusion follows.
\end{proof}

Here's a result that we get from using $\falling{x}{n}$ in place of $\frac{1}{\Gamma(x)}$.
\begin{lemma}
If $f\in\allpolysep \cap \allpolyint{(n,\infty)}$ then
  $nf+(x+1)\Delta f\in\allpolysep \cap \allpolyint{(n-1,\infty)}$.
\end{lemma}
\begin{proof}
  Since the largest root of $\falling{x}{n}$ is $n-1$ it follows that
  $\falling{x}{n}f\in\allpolysep$. Now taking the difference
  preserves $\allpolysep$, so
  \begin{align*}
    \Delta \falling{x}{n}f &= f\, \Delta\falling{x}{n} + \left(\Delta f\right)
    \affa \falling{x}{n} \\
&= \falling{x}{n-1}\left( nf + (x+1) \Delta f\right)
  \end{align*}
is in $\allpolysep$. The result follows since factors of polynomials
in $\allpolysep$ are in $\allpolysep$.
\end{proof}

Motivated by the analogy of $\rising{x}{n}/n!$ and the exponential, we
define a Laguerre like polynomial, and show that consecutive ones
interlace.

\begin{lemma}
  If $f_n (x) = \Delta^n \rising{x+1}{n}\falling{x-1}{n}$ then
  $f_n\in\allpolysep$ and $f_n\lesslesseq f_{n-1}$. 
\end{lemma}
\begin{proof}
 If we let $p_n(x) = \rising{x+1}{n}\falling{x-1}{n}$ then
  $$
  \Delta p_n(x) = n \rising{x+2}{n-1}\falling{x-1}{n-1}(1+2x).$$
Observe that that $\Delta p_n \lesslesseq p_{n-1}$. Since
  $p_n\in\allpolysep$ the conclusion follows by applying $\Delta^{n-1}$.
\end{proof}

  \begin{lemma}\label{lem:affine-cos}
    If $T(f) = f(x+\imath) + f(x-\imath)$ then
    $T\colon{}\allpolysep\longrightarrow\allpolysep$.  \cite{walker97}
  \end{lemma}
  \begin{proof}
    Observe that $T(f) = 2 \cos(\diffd)f$, so that $T$ maps
    $\allpoly\longrightarrow\allpoly$ and preserves interlacing.
    Since $T$ commutes with $\affa$, the result follows from
    Lemma~\ref{lem:affine-commute}.
  \end{proof}

  \section{Integrating families of polynomials in $\allpolysep$ and $\allpolyaffine$}
\label{sec:affine-sep}

We now apply the integration results of
\chapsec{poly-matrices}{int-fam} to $\allpolysep$.  Since polynomials
in $\allpolysep$ are exactly those polynomials whose roots are
separated by at least one,  $f\in\allpolysep$ iff $\delta(f)\ge1$.
We will see more applications of integration when we consider the
generalization of $\allpolysep$ to two variables.

One useful way of constructing locally interlacing families is to
start with $f(x)\in\allpolysep$, and form the family $\{f(x+t)\}$. The
interval of interlacing is at least one. 

Here is a simple application. Note that the hypothesis
$f'\in\allpolysep$ doesn't imply $f\in\allpoly$, yet we can conclude
that $\Delta f\in\allpoly$.

\begin{cor} \label{cor:family-trans}
  Suppose that $w(x)$ is positive on $(0,1)$. The mapping
\begin{equation} \label{eqn:family-trans}
 f \mapsto \int_0^1\,f(x+t)w(t)\,dt
\end{equation}
defines a linear transformation $\allpolysep\longrightarrow \allpoly$.
\end{cor}

\begin{cor}
  If $f^\prime\in\allpolysep$ then $\Delta f\in\allpoly$ where $\Delta
  f = f(x+1)-f(x)$.
\end{cor}
\begin{proof}
  Choose $w(t)=1$, $a=0$, $b=1$, and $f_t(x) = f'(x+t)$ in
  Proposition~\ref{prop:family-int}. 
\end{proof}

\begin{remark}
  \index{interlacing!doesn't preserve}

  Integration is an example of a linear transformation \\
  $T\colon{}\allpolysep\longrightarrow\allpoly$ that does not preserve
  interlacing.  For example
  \begin{xalignat*}{2}
    f(x) & = x(x+1) & g(x) &= x\\
    T(f) &= \int_0^1 f(x+t)\,dt & T(g) &= \int_0^1 g(x+t)\,dt\\
    & = (x+.59)(x+1.4) & & = x+.5
  \end{xalignat*}
  Thus, $f\lesslesseq g$ yet $Tf$ and $Tg$ don't interlace.
\end{remark}

The following result is  convolution of polynomials in $\allpolysep$. 

\index{convolution}

\begin{lemma}\label{lem:convolution-int}
If  $f,g\in\allpolysep$   then
$$ \int_0^1 f(x+t)g(x-t)\, dt\, \in\allpoly$$
\end{lemma}
\begin{proof}
  If we choose $0<t_1<\cdots<t_n$ then let $f_i(x) =  f(x+t_i)$
  and $g_{n+1-i}(x) = g(x-t_i)$. By Lemma~\ref{lem:convolution} we
  know that $\sum f_i\, g_{n+1-i}\in\allpoly$. Taking the limit shows
  the integral is in $\allpoly$. 
  
  Note: the diagram in the proof of Lemma~\ref{lem:convolution} is
  essentially the graph of $f(x+t)g(x-t)$. Consideration of this
  diagram shows that we can find $f,g$ so that if $h$ is their
  convolution then $\delta(h)$ is as small as desired.
\end{proof}

Since we only need separation in one variable, we phrase the next
result in terms of mutually interlacing families.

\index{polynomials!mutually interlacing}
\index{interlacing!mutually}

\begin{cor}\label{cor:int-fam-3}
  If $\{f_t\}$ and $\{g_t\}$ are locally interlacing families on
  $\reals$ with $\rho(f)\ge1$ and $\rho(g)\ge 1$, and we define
$$ \Phi = \int_0^1 \,\int_0^{1-t} \,f_t(x)g_s(x)\,ds\,dt$$
then $\Phi\in\allpoly$.
\end{cor}
\begin{proof}
  If we define $h_t = \int_0^t g_s(x)\,dt$ then by Lemma~\ref{lem:family-int-2a}
  $\{h_t\}$ is a locally interlacing family with $\rho(h)\ge1$. Since
  $\Phi = \int_0^1 f_t(x)h_{1-t}(x)\,dt$ we see that $\Phi$ is a
  convolution of mutually interlacing polynomials. ($\Phi$ is
  generally not in $\allpolysep$.)
\end{proof}

\begin{cor}
  If $f\in\allpolysep$ then $\displaystyle\int_0^1 f(x+t)e^{-t}\,dt\in\allpoly$.
\end{cor}
\begin{proof}
  If we choose $g(x) = e^{-x}$ which is in $\allpolysepf$, then we can
  apply the lemma to find that $$\int_0^1 f(x+t)e^{x-t}\,dt = e^x
  \int_0^1 f(x+t)e^{-t}\,dt \in\allpoly.$$ If we multiply by $e^{-x}$
  we get the conclusion. 

\end{proof}

  We now to assume that $\affa x = qx$ where $q\ge1$.  We can
  integrate polynomials in $\allpolyaffine$ since they determine
  sequences of interlacing polynomials. Assume that
  $f(x)\in\allpolyaffine$. It's clear that
\begin{equation}\label{eqn:aff-int-q}
f(x) \greateqeq f(tx) \greateqeq f(qx) \text{ if } 1 \le t \le q
\end{equation}
Consequently, we have an interlacing family $\{f(tx)\}$.  Applying the
results of \chapsec{poly-matrices}{int-fam} we conclude that

\begin{cor} \label{cor:family-trans-q}
  Suppose that $w(x)$ is positive on $(1,q)$. The mapping
\begin{equation} \label{eqn:family-trans-q}
 f \mapsto \int_1^q\,f(tx)w(t)\,dt
\end{equation}
defines a linear transformation $\allpolyaffine\longrightarrow \allpoly$.
\end{cor}

\begin{cor}
If $f' \in\allpolyaffine$ then $\qderiv(f)\in\allpoly$.
\end{cor}
\begin{proof}
  Immediate from the computation 
$$\displaystyle\int_1^q f'(tx)\,dt = 
(f(qx)-f(x))/x = (q-1) \qderiv(f)$$
\end{proof}

\index{convolution}
We also have convolution. If $f(x)\in\allpolyaffine$ then replacing
$x$ by $x/q$ in \eqref{eqn:aff-int-q} shows 

\begin{equation}\label{eqn:aff-int-q-2}
f(x/q) \greateqeq f(x/t) \greateqeq f(x) \text{ if } 1 \le t \le q
\end{equation}

Equations \eqref{eqn:aff-int-q} and \eqref{eqn:aff-int-q-2} are
reversed if $0<q<1$. We can now apply
Lemma~\ref{lem:convolution-int} to conclude

\begin{lemma}\label{lem:convolution-int-q}
  If $f,g\in\allpolyaffine$ where $\affa x = qx$, and $w(t)$ is
  positive on $(1,q)$ if $q>1$ and on $(q,1)$ if $0<q<1$ then
$$ \int_1^q w(t) f(xt)g(x/t)\,dt \in\allpoly$$.
\end{lemma}

\section{Determinants and $\allpolyaffine$}
\label{sec:det-and-pa}

In this section we look at connections between $\allpolyaffine$ and
polynomials of the form $|xD+C|$ with symmetric $C$ and diagonal $D$.
We show that there is an analog of the differentiation formula
\eqref{eqn:daffine-n} where the factors are replaced by determinants
of principle submatrices. The rest of the section is concerned only
with $\affa x = x+1$.

We assume that $D$ is an $n$ by $n$ diagonal matrix with diagonal
entries $a_i$, and $C=(c_{ij})$ is a symmetric matrix. We are
interested in the matrix $M = x D + C$. If $M[i]$ is the $i$th
principle submatrix, then
  \begin{equation}
    \label{eqn:det-dif-1}
    \frac{d}{dx}\,det(M) = a_1\, det(M[1]) + \cdots + a_n\, det(M[n])
  \end{equation}
  The formula for the affine derivative is similar:

\begin{lemma} \label{lem:det-def}
Assume all $a_i$ are positive.
If $M_{[r]}$ is constructed from $M[r]$ by applying $\affa$ to the
first $r-1$ diagonal elements of $M[r]$ then 

\begin{equation}
    \label{eqn:det-dif-2}
\daffine\,det(M) = a_1\, det(M_{[1]}) + \cdots + a_n\, det(M_{[n]})
  \end{equation}
\end{lemma}

For example,
\begin{align*}
 M&= \begin{pmatrix}
  a_1x + c_{11} & c_{12} & c_{13} \\
 c_{12} & a_2x + c_{22}& c_{23} \\
 c_{13} & c_{23} &  a_3x + c_{33} 
\end{pmatrix} \\
M_{[1]}&= \begin{pmatrix}
a_2x + c_{22} & c_{23} \\
 c_{23} &  a_3x + c_{33}
\end{pmatrix} \\
M_{[2]}&= \begin{pmatrix}
\affa(a_1x + c_{11}) & c_{13} \\
 c_{13} &  a_3x + c_{33}
\end{pmatrix} \\
M_{[3]}&= \begin{pmatrix}
\affa(a_1x + c_{11}) & c_{12} \\
 c_{12} & \affa(  a_2x + c_{22})
\end{pmatrix} 
\end{align*}

$$\daffine\, det(M) = a_1\, det(M_{[1]}) + a_2\, det(M_{[2]}) + a_3\, det(M_{[3]})
$$

\begin{proof}[Proof (of \ref{lem:det-def})]
  Note that the term $a_i det(m_{[i]})$ is the determinant of the $n$
  by $n$ matrix formed from $M$ by setting all elements of the $i$-th
  row and column to zero, except for the main diagonal element, which
  is $a_i$. Since all terms of the equation are now $n$ by $n$
  matrices, we can multiply by $D^{-1/2}$ on each side, and so we may
  assume that all $a_i$ are equal to $1$. We now apply the following
  lemma, with $y=\affa x$, to finish the proof.
\end{proof}

\begin{lemma}\label{lem:det-def-new}
  If $C$ is an $n$ by $n$ symmetric matrix, and $D_k$ is the $n$ by
  $n$ diagonal matrix whose first $k$ diagonal elements are $y$ and
  the remaining ones are $x$ then
  \begin{multline}
    \label{eqn:det-def-new}
    \frac{det(C+x I) - det(C+y I)}{x-y} = \\
det\,(C+D_1)[1] + det\,(C+D_2)[2] + \cdots + \,det\,(C+D_n)[n]
\end{multline}
\end{lemma}

Notice that in the special case that $M$ is a diagonal matrix the
formula \eqref{eqn:det-dif-2} reduces to \eqref{eqn:daffine-n}. 
The limit of \eqref{eqn:det-def-new} as $y\rightarrow x$ is
\eqref{eqn:det-dif-1}. 

For the remainder of this section we restrict ourselves to $\affa
x=x+1$ and $\allpolysep$. If $f(x) = \prod(x-r_i)$ then the polynomials
$f/(x-r_i)$ are mutually interlacing. It is easy to find examples
where the determinants of the principle submatrices of $|xI+C|$ are not
mutually interlacing. However, if $|xI+C|$ is in $\allpolysep$ then we
can find a family of mutually interlacing polynomials that generalize
the $f_{[k]}$ of \eqref{eqn:aff-rts-4}.

Let $M = xI+C$ and write
$M=\smalltwodet{x+c_{11}}{v}{v^t}{M[1]}$. Then
\begin{align*}
  \begin{vmatrix}
    \affa(x+c_{11}) & v \\ v^t & M[1]
  \end{vmatrix}
&= 
\begin{vmatrix}
  x+c_{11} & v \\ v^t & M[1]
\end{vmatrix}+
\begin{vmatrix}
  1 & 0 \\ v^t & M[1]
\end{vmatrix}\\
&= |M| + |M[1]|\\
\intertext{Since $|M|\lesslesseq |M[1]|$ we see that }
|M| & \greateqeq 
\begin{vmatrix}
    \affa(x+c_{11}) & v \\ v^t & M[1]
\end{vmatrix}
\end{align*}
If we let $M_{[r]}$ be the result of applying $\affa$ to the first $r$
diagonal elements of $M$ then inductively we see that
$$
|M| \greateqeq |M_{[1]}| \greateqeq |M_{[2]}| \greateqeq \cdots
\greateqeq |M_{[n]}| = \affa |M|
$$
Thus, if $|M|\in\allpolysep$ then $|M|\greateqeq \affa|M|$ and
consequently the following polynomials are mutually interlacing
$$
\left\{ |M|,|M_{[1]}|,|M_{[2]}|, \cdots, |M_{[n-1]}| , \affa
  |M|\right\}
$$

Every factor of a polynomial in $\allpolysep$ is in $\allpolysep$, but
the analogous property does not hold for principle submatrices. 
For instance, the matrix $M$ below has the property that
$\delta(|M[k]|) < \delta(|M|)$ for $k=1,2,3,4,5$. 

$$
\begin{pmatrix}
x+2 & 4 & 6 & 6 & 5 \\
4 & x+2 & 5 & 8 & 6 \\
6 & 5 & x & 7 & 6 \\
6 & 8 & 7 & x+2 & 4 \\
5 & 6 & 6 & 4 & x  
\end{pmatrix}
$$

However, the roots are always distinct.

\begin{lemma}
  If $M=xI+C$ is a matrix whose determinant has all distinct roots,
  then the determinants of all principle submatrices have all distinct
  roots. 
\end{lemma}
\begin{proof}
  The determinant of a principle submatrix interlaces the determinant
  of $M$.
\end{proof}

\section{The case $\affa x = qx$ with $q>1$}
\label{sec:affine-q}

In this section we consider the affine transformation $\affa f(x) =
f(qx)$ where we first assume that that $q>1$. We know that $f
\greateqeq \affa f$ implies that $f\in\allpolyalt$. There are simple
conditions on the roots that characterize such polynomials.  If the
roots of $f$ are $\{a_i\}$ then the roots of $\affa f$ are
$\{a_i/q\}$. If $f\in\allpolyalt$ then since $q>1$ the ordering of the
roots is
$$
a_1/q < a_1 < a_2/q < a_2 < \cdots $$

The ratio of consecutive roots is at least $q$. Conversely, if
$f\in\allpolyalt$ and the ratio of consecutive roots is at least $q$
then $f\greateqeq \affa f$.

\noindent%
The affine derivative is called the $q$-derivative, and is given in
\eqref{eqn:affine-product-3}. The next lemma specializes the results
in Section \ref{sec:affine-derivative}.

\begin{cor}\label{cor:q>1}
  Suppose $\affa x = qx$ where $q>1$. The following are true:
  \begin{enumerate}
  \item If $f\in\allpolyaffine$ then $f\in\allpolyalt$. 
  \item If $f\in\allpolyaffine$ then $f^\prime\in\allpolyaffine$.
  \item If $f\in\allpolyaffine$ then $\qderiv f\in\allpolyaffine$.
  \item If $f\greateqeq g\in\allpolyaffine$ then
    $\qderiv{}f\greateqeq\qderiv{}g$. 
  \item If $g\in\allpolyalt$ and $f\in\allpolyaffine$ then
    $g(\affa)f\in\allpolyaffine$. 
  \item If $g\in\allpolyalt$ and $f\in\allpolyaffine$ then
    $g(\qderiv{})f\in\allpolyaffine$. 
  \end{enumerate}
\end{cor}

For example, if we take $g=(x-1)^n$ then the following polynomials are
in $\allpolyaffine$:
\begin{align*}
  (\affa-1)^n f &= \sum_{i=0}^n \binom{n}{i} f(q^i x)\\
  (\qderiv{}-1)^n f &= \sum_{i=0}^n \binom{n}{i} \qderiv{}^i f( x)
\end{align*}

In order to apply Lemma~\ref{lem:affine-fg} we
need the following formula for higher derivatives that is easily
proved by induction:

$$ \qderiv{}^m\,g(x) = \frac{1}{(q-1)^m q^{\binom{n}{2}}x^m}
\quad \sum_{i=0}^m (-1)^i \binom{m}{i}_q\,q^{\binom{i}{2}} \,g(q^ix).
$$
 $\binom{n}{i}_q$ is the \emph{Gaussian binomial coefficient}
\index{Gaussian binomial coefficient}
defined by
$$ \binom{n}{i}_q = \frac{[n]!}{[k]!\,[n-k]!}$$ where as usual $[n]! =
[n][n-1]\cdots[1]$.

\begin{cor}
  If $g\in\allpolyaffine$ and $m$ is a positive integer then 
$$
\frac{1}{x^m}
\quad \sum_{i=0}^m (-1)^i \binom{m}{i}_q\,q^{\binom{i}{2}}
\,g(q^ix)\in\allpolyaffine. 
$$
\end{cor}

Lemma~\ref{lem:affine-det} implies that
\begin{lemma}\label{lem:affine-det-3}
     Suppose $f\in\allpolyaffine\cap\allpolypos$ and $q>1$.
     For all positive $x$ 

     $$\smalltwodet{f(qx)}{f(x)}{f(q^2x)}{f(qx)}>0$$
  \end{lemma}

\section{The case $\affa x = qx$ with $0<q<1$}
\label{sec:affine-q-2}

The cases $q>1$ and $0<q<1$ are closely related. Let $\affa_q(x)=qx$,
and $\affa_{1/q}(x)= x/q$.  Consequently, if $f\in\allpoly$ then
$$
f(x) \greateq \affa_q \, f(x) \Longleftrightarrow f(-x) \greateq
\affa_{1/q}f(-x)$$

The transformation $x\mapsto -x$ is a bijection between
$\allpolyaffine$ for $q>1$ and $\allpolyaffine$ for $0<1/q<1$. 

For the rest of this section we assume that $0<q<1$.

The fundamental polynomials in $\allpolyaffine$ are 
$$ (-x;q)_n = (1+x)(1+qx)(1+q^2x)+\cdots(1+q^{n-1}x)$$
which have roots
$$ -q^{-(n-1)} < - q^{-(n-2)} < \cdots < -q^{-1} < -1$$
Note that $\qderiv(x;q)_n = [n](x;q)_{n-1}$.

We can translate the results of Lemma~\ref{cor:q>1}:
\begin{cor}\label{cor:0<q<1}
  Suppose $\affa x = qx$ where $0<q<1$. The following are true:
  \begin{enumerate}
  \item If $f\in\allpolyaffine$ then $f\in\allpolypos$. 
  \item If $f\in\allpolyaffine$ then $f^\prime\in\allpolyaffine$.
  \item If $f\in\allpolyaffine$ then $\qderiv f\in\allpolyaffine$.
  \item If $f\greateqeq g\in\allpolyaffine$ then
    $\qderiv{}f\greateqeq\qderiv{}g$. 
  \item If $g\in\allpolypos$ and $f\in\allpolyaffine$ then
    $g(\affa)f\in\allpolyaffine$. 
  \item If $g\in\allpolypos$ and $f\in\allpolyaffine$ then
    $g(\qderiv{})f\in\allpolyaffine$. 
  \end{enumerate}
\end{cor}

Lemma~\ref{lem:affine-det} implies that
\begin{lemma}\label{lem:affine-det-4}
     Suppose $f\in\allpolyaffine\cap\allpolypos$ and $0<q<1$.
     For all positive $x$ 

     $$\smalltwodet{f(qx)}{f(x)}{f(q^2x)}{f(qx)}<0$$
  \end{lemma}

\section{The $q$-exponential}
\label{sec:affine-exp}

In this section we continue to assume that $0<q<1$. 
We will use the $q$-exponential $E_\affa$ of \eqref{eqn:exp-affa-3}
to derive properties of the $q$-Laguerre polynomial. We have already
seen $E_q(x)$  in \eqref{eqn:q-exponential}. $E_q(x)$ is an
entire function if $|q|<1$.  Using the product representation it is
easy to verify that $\qderiv{} E_q(\alpha x) = \alpha\affa E_q(\alpha
x)$.  \index{q-exponential}

From \eqref{eqn:q-exponential} we see that all the terms have positive
coefficients, and 
 the roots of the $q$-exponential $Exp_q(x)$ are 
$\left\{-q^{-n}\right\}_{n=1}^\infty$, so clearly we have
$$
E_q(x) \greateqeq \affa E_q(x)$$
Thus, $E_q$ is an entire
function that is in the closure of $\allpolyaffine$.




\index{Rodrigues' formula!Laguerre}
\index{Laguerre polynomials!Rodrigues' formula}
\index{Laguerre polynomials!$q$}

\index{Rodrigues' formula!$q$-Laguerre}
The Rodrigues' formula for Laguerre polynomials is
$$ L_n(x) = \frac{1}{e^{-x}} \frac{1}{n!} \diffd^n ( e^{-x} x^n)$$
We define the $q$-Laguerre polynomials by the formula
$$
L_n^\affa(x) = \frac{1}{\affa^n E_q(-x)} \qderiv{}^n ( E_q(-x)
(-x;q)_n)$$
However, since $E_q(-x)\not\in\allpolyaffinef$, it is by
no means obvious that this is even a polynomial, let alone that it is
in $\allpolyaffine$. To this end, define

$$ S(n) = \left\{f\cdot \affa^n(E_q(-x))\mid
  f\in\allpolyaffine\right\}$$

We can control the action of $\qderiv{}$ in $S(n)$:

\begin{lemma} For any non-negative integer $n$, 
  $\qderiv: S(n) \longrightarrow S(n+1)$
\end{lemma}
\begin{proof}
  Apply the $q$-Leibnitz formula to the product:
  \begin{align*}
    \qderiv ( f \cdot \affa^n E_q(-x)) &= 
\qderiv f \cdot \affa(\affa^n E_q(-x)) + f \cdot \qderiv(\affa^n E_q(-x)) \\
&= \qderiv f \cdot \affa^{n+1} E_q(-x) - q^n f \cdot \affa^n \qderiv
E_q(-x) \\
&= \left( \qderiv f - q^n f\right) \cdot \affa^{n+1}E_q(-x)
  \end{align*}
Since $\qderiv f - q^n f\in\allpolyaffine$ by Lemma~\ref{lem:affine-faf} it
follows that $\qderiv ( f \cdot \affa^n E_q(-x))$ is in $S(n+1)$.
\end{proof}

\begin{cor}
  The $q$-Laguerre polynomials are in $\allpolyaffine$.
\end{cor}
\begin{proof}
  Since $(-x;q)_n\in\allpolyaffine$, it follows that $ (-x;q)_n\cdot
  E_q(-x)$ is in $S(0)$. The Lemma above shows that
  $\qderiv^n(\left[-x;q)_n\cdot E_q(-x)\right]$ is in $S(n)$. This
  means that it has the form $g\cdot \affa^n E_q(-x)$ where $g$ is in
  $\allpolyaffine$, which proves the corollary.

\end{proof}

There is a simple formula for the $q$-Laguerre polynomials.

\[ 
L_n^\affa(x) = \sum_{i=0}^n (-1)^i \, x^i \, \binom{n}{i}\,  \frac{ [n]!}{[i]!}
\]

\index{q-Hermite polynomials}
\index{polynomials!q-Hermite}

There are several $q$-analogs of the Hermite polynomials \cite{koekoek}.
The variant that we use is motivated by the desire for a recursive
formula similar to  \eqref{eqn:hermite-3}. We set $H^q_0=1, H^q_1=x$
and
$$
H^q_{n+1} = x \,H^q_n - q^n\,[n]\,H^q_{n-1}.
$$
Since the coefficient of $H^q_{n-1}$ is negative whenever $q$ is
positive, we see that $H^q_n$ has all real roots for all positive $q$.
From the definition it is easy to show by induction that
$$ 
\qderiv{}\,H^q_n = [n]\,H^q_{n-1}
$$
which is the desired analog of \eqref{eqn:hermite-3}.  
If we combine these last two equations we find
$$ 
H^q_n = (x-q^{n-1}\, \qderiv{})(x-q^{n-2}\, \qderiv{q}) \dots (x-
\qderiv{})(1).
$$
which is the $q$-analog of the identity $H_n = (2x-\diffd)^n(1)$.

The $q$-hermite polynomials can not be in $\allpolyaffine$ since they
have both positive and negative roots. However, they appear to be
nearly so. See Question~\ref{ques:q-hermite}.

\section{q-transformations}
\label{sec:qseries}

In this section we show that certain linear transformations associated
to q-series map appropriate subsets of $\allpoly$ to $\allpoly$.
From earlier results specialized to the case at hand we have

\begin{cor} \label{affine-q-trans}
  Let $\affa x = qx$ where $q>1$. 

  \begin{enumerate}
  \item The linear transformation $x^n\mapsto q^{\binom{n}{2}}x^n$
    maps $\allpolyalt\longrightarrow\allpolyalt\cap \allpolyaffine$.
  \item The linear transformation $x^n\mapsto (-x;q)_n$
    maps $\allpolyalt\longrightarrow\allpolyint{(1,\infty)}\cap \allpolyaffine$.
  \end{enumerate}
\end{cor}

We next study the latter transformation and its inverse for all values of $q$.
 Define

\begin{equation} \label{eqn:qseries}
 T_q(x^i) = (x;q)_i = (1-x)(1-qx)\cdots(1-q^{i-1}x)
\end{equation}

We will show the following:

\begin{theorem} \label{thm:qseries}
  If $T_q$ is defined in \eqref{eqn:qseries}, then
  \begin{enumerate}
  \item If $0<q<1$ then 
    $T_q:\allpolyint{\reals\setminus(0,1)}\longrightarrow\allpolyint{\reals\setminus(0,1)}$.
  \item If $0<q<1$ then $T_q^{-1}:\allpolyint{(0,1)}\longrightarrow\allpolyint{(0,1)}$.
  \item If $1<q$ then $T_q:\allpolyint{(0,1)}\longrightarrow\allpolyint{(0,1)}$.
  \item If $1<q$ then 
    $T_q^{-1}:\allpolyint{\reals\setminus(0,1)}\longrightarrow\allpolyint{\reals\setminus(0,1)}$.
  \end{enumerate}
\end{theorem}

The idea is to discover the recursion relations satisfied by $T_q$ and
$T_q^{-1}$. We will show that $T_{(1/q)} = T_q^{-1}$, and this allows
us to prove the results about $T_q^{-1}$.  From the definition

\begin{align}
  T_q(x\cdot x^{n-1})=  T_q(x^n) &= (1-x)(1-qx)(1-q^2x)\cdots(1-q^{n-1}x) \notag\\
  &= (1-x)\affa\left((1-x)(1-qx)\cdots(1-q^{n-2})\right) \notag\\
  &= (1-x)\affa T_q(x^{n-1})\notag \\
  \intertext{which by linearity implies}
  T_q((x+b)f) &= (1-x)\affa T_q(f)+ b T_q(f)\label{eqn:qseries-1}\\
  \intertext{Similarly,} x\,T_q(x^n) &=
  \frac{1}{q^n}\left(T_q(x^n)-T_q(x^{n+1})\right)\ =\ 
  T_q(\,(1-x)\affa^{-1}\,x^n)\notag\\
  T_q^{-1}\left(xT_q(x^n)\right) &= (1-x)\affa^{-1} x^n\notag\\
  T^{-1}_q\left(xT_q(g)\right) &= (1-x)\affa^{-1}\,g \notag\\
  \intertext{Setting $g=T_q^{-1}(f)$ yields}
  T_q^{-1}(xf) &= (1-x)\affa^{-1}T_q^{-1}(f)\notag \\
  \intertext{Applying this last recursion to $x^n$ yields}
  T_q^{-1}(x^n) &=
  \left(1-\frac{x}{1}\right)\left(1-\frac{x}{q}\right)
  \cdots\left(1-\frac{x}{q^{n-1}}\right)\ =\ T_{1/q}(x^n) \notag
\end{align}

Consequently, we have the equality  $T_q^{-1} = T_{1/q}$. This relation shows that part 1
of  Theorem~\ref{thm:qseries} implies part 4, and part 3 implies part 2.

The next Lemma establishes properties  of the recurrence
\eqref{eqn:qseries-1}.

\begin{lemma} \label{lem:qseries}
  Suppose that $0<q<1$, $b$ is positive, $\affa(x)=qx$, and define
  $$S(f) = (1-x)\affa f + b f.$$ Then  if
  $f\in\allpolyint{(1,\infty)}$ and $\affa(f)\greateqeq f$ then
  $S(f)\in\allpolyint{(1,\infty)}$, and $\affa S(f) \greateqeq S(f)$.
\end{lemma}
\begin{proof}
  Since $\affa (f)\greateqeq f$ and all roots of $f$ are greater than
  $1$ it follows that $(1-x)\affa(f)\greateqeq f$. By linearity, we
  know that $(1-x)\affa(f) + b f + \gamma f\in \allpoly$ for any
  $\gamma$, and consequently $S(f)\lesslesseq f$. This also shows that
  $S(f)\in\allpoly$. Since $f\in\allpolyint{(1,\infty)}$ we know that
  $S(f)$ has at most one root less than $1$. It's easy to see that the
  coefficients of $S(f)$ alternate, so $S(f)\in\allpolyalt$. We need
  to show that $S(f)$ has no roots in $(0,1)$. The leading
  coefficients of $f$ and $S(f)$ alternate in sign, and it's easy to
  check that $f$ is positive as $x\rightarrow-\infty$. Since all roots
  of $f$ are greater than $1$, it follows that $f$ is positive on
  $(0,1)$. Similarly, $\affa(f)$ is positive on $(0,1)$. Finally $1-x$
  is positive on $(0,1)$, and so $S(f)$ is positive there. This
  implies that $S(f)\in\allpolyint{(1,\infty)}.$

  It remains to show that $\affa(S(f))\greateqeq S(f)$, and we do this
  by analyzing the location of the roots of $S(f)$.
  Suppose that the roots of $f$ are $r_1<\cdots<r_n$. The roots of
  $\affa f$ are $\frac{r_1}{q}<\cdots<\frac{r_n}{q}$, and since $\affa
  f \greateqeq f$ we have the ordering
  $$
  1 < r_1 < \frac{r_1}{q} < r_2 < \frac{r_2}{q}< \cdots < r_n <
  \frac{r_n}{q}$$
  The signs of $f$ and $\affa f$ are constant on each
  of the intervals defined by consecutive terms of the above sequence.
  Now $1-x$ is negative for $x>1$, so the roots of $S(f)$ lie in the
  intervals where $f$ and $\affa f$ have the same sign. Since both $f$
  and $\affa f$ have the same sign on $(\frac{r_n}{q},\infty)$ the
  roots of $S(f)$ are located in the intervals

  $$(\frac{r_1}{q},r_2),\quad(\frac{r_2}{q},r_3),\quad\cdots\quad(\frac{r_{n-1}}{q},r_n),
  \quad(\frac{r_n}{q},\infty)$$
  Since $S(f)\lesslesseq f$ there is
  exactly one root in each interval. The smallest possible ratio
  between the first two roots of $S(f)$ is found by taking the largest
  possible first root, $r_2$, and dividing it into the smallest
  possible second root, $\frac{r_2}{q}$. This ratio is $\frac{1}{q}$,
  and the same argument applies to all the possible ratios. Since all
  consecutive ratios of the roots of $S(f)$ are at least
  $\frac{1}{q}$, we find that $\affa S(f)\greateqeq S(f)$.
\end{proof}

\begin{proof}[Proof of Theorem \ref{thm:qseries}]
We begin with part 1, so assume that $0<q<1$ and $f\in\allpolyint{(1,\infty)}.$
Now $T_q((x+b)f) = S(T_q(f))$ by \eqref{eqn:qseries-1}, so by induction
it suffices to prove our result for $f$ linear, which is trivial.
The case for $f\in\allpolypos$ is similar and omitted.

The proof of part 3 is entirely similar to part 2 - we need a modified
version of Lemma~\ref{lem:qseries} - and is omitted. As observed above, the
remaining parts follow from 1 and 3.
\end{proof}

  \section{Newton's inequalities for $\affa x = qx$ and $\affa x = x+1$}
  \label{sec:newton-qx}
\index{Newton's inequalities!for $\allpolyaffine$}
  Suppose $\affa x = q x$ where $q>0$.  The coefficients of a
  polynomial in $\allpolyaffine$ satisfy constraints that are the
  $q$-generalizations of Newton's inequalities. The key observation is
  that $\allpolyaffine$ is closed under reversal.

  \begin{lemma}
    If $\affa x = qx,q>0$ and $f\in\allpolyaffine$ then
    $f^{rev}\in\allpolyaffine$. 
  \end{lemma}
  \begin{proof}
    If the roots of $f$ are 
\[r_1< r_2< \cdots < r_n \]
then the roots of $f^{rev}$ are
\[ \frac{1}{r_n} < \frac{1}{r_{n-1}} < \cdots < \frac{1}{r_1}. \]
Since $f\in\allpolyaffine$ the ratio of consecutive roots of $f$
satisfies $r_{i+1}/r_i\ge q$ for $1\le i <n$.
The ratio of consecutive roots of $f^{rev}$ satisfies
\[ \frac{1/r_{k}}{1/r_{k+1}} \ge q \]
and so $f^{rev}\in\allpolyaffine$.
  \end{proof}

\begin{prop}\label{prop:newton-q}
        If $\affa x = q x$ and we choose $f(x) = a_0+\cdots+
    a_nx^n$ in $\allpolyaffine(n)$ then 
\begin{equation}\label{eqn:q-newton} 
\frac{a_{k+1}^2}{a_{k}a_{k+2}}\ge
\frac{4}{(q+1)^2}\,\frac{[k+2][n-k]}{[k+1][n-k-1]}
\end{equation}
In addition,
\begin{equation}\label{eqn:q-newton-2} 
\frac{a_{k+1}^2}{a_{k}a_{k+2}}\ge
\frac{1}{q}\frac{4}{(q+1)^2}\,\frac{[k+2][n-k]}{[k+1][n-k-1]}
\end{equation}
%
\end{prop}
\begin{proof}
  We follow the usual proof of Newton's inequalities.  Write
\[ f(x) = \sum_{k=0}^n b_k \binom{n}{k}_q\,x^k \]
where $\binom{n}{k}_q = [n]!/([k]![n-k]!)$. The $q$-derivative of $f$
satisfies
\[ \qderiv{f} = [n]\sum_{k=1}^{n} b_k \binom{n-1}{k-1}_q \,x^{k-1}\]
Consequently, if we then apply the $q$-derivative $k$ times, reverse,
apply the $q$-derivative $n-k-2$ times, and reverse again, the
resulting polynomial is in $\allpolyaffine$:
\begin{equation}\label{eqn:q-newton-1}
 [n]\cdots[n-k+1]\biggl( b_k \binom{2}{0}_q + b_{k+1}\binom{2}{1}_qx
+ b_{k+2}\binom{2}{2}_qx^2\biggr)
\end{equation}
After removing the constant factor, the discriminant satisfies
\[ (q+1)^2 \,b_{k+1}^2 - 4\, b_kb_{k+2} \ge 0\]
Substituting $b_k = a_k /\binom{n}{k}_q$ and simplifying yields the
first conclusion.

Now suppose that $0<q<1$. The roots of \eqref{eqn:q-newton-1} are
negative, and the coefficients are positive. Let
\begin{align*}
  r_1 &= \frac{-(q+1)b_{k+1} - \sqrt{(q+1)^2b_{k+1}^2 - 4 b_k
      b_{k+2}}}{2b_{k+2}}\\
  r_2 &= \frac{-(q+1)b_{k+1} + \sqrt{(q+1)^2b_{k+1}^2 - 4 b_k
      b_{k+2}}}{2b_{k+2}}\\
  \intertext{Since \eqref{eqn:q-newton-1} is in $\allpolyaffine$ we
    know $  r_1 \le q\,r_2$ so} 0 \le r_2 - q\,r_1 &=
  (q+1)\biggl((q-1)b_{k+1} + \sqrt{(q+1)^2b_{k+1}^2 - 4 b_k
    b_{k+2}}\biggr)/2b_{k+2} \\
  \intertext{Simplifying yields $qb_{k+1}^2 \ge b_kb_{k+2}$ which establishes
    the second part}
\end{align*}
 
\end{proof}

\begin{remark}
The product 
\begin{equation}\label{eqn:q-bin-thm}
 \prod_{k=0}^{n-1}(x+q^k) = \sum_{k=0}^n x^k 
\binom{n}{k}_q q^{\binom{n-k}{2}} 
\end{equation}
has Newton quotients 
\[
\frac{1}{q}\,\frac{[k+2][n-k]}{[k+1][n-k-1]}.
\]
Numerical computation suggests that this is probably the best possible bound. 
\end{remark}

\begin{example}
  If $q=1$ then $4/(q+1)^2=1$ and $[n]=n$, so  \eqref{eqn:q-newton}
  becomes the usual Newton inequalities \eqref{eqn:newton-1}.
  If we take  $f(x) =
    a_0+a_1x+a_2x^2+a_3x^3\in\allpolyaffine$ then
\begin{align*}
k&=1 & \frac{a_1^2}{a_{0}a_{2}}&\ge  \frac{4}{(q+1)^2}(1+q+q^2) \,\ge\, 3 \\
k&=2& \frac{a_2^2}{a_{1}a_{3}}&\ge   \frac{4}{(q+1)^2}(1+q+q^2) \,\ge\, 3 \\
\end{align*}
The bound $3$ is the same as the usual Newton's inequality bound, and is
realized if $q=1$.  If $q\ne1$ then the bound is better. For example,
if $q=2$ then the ratio is at least $3\frac{1}{9}$.

\end{example}

  Although the reverse of a polynomial in $\allpolysep$ generally
  isn't in $\allpolysep$, there are also inequalities for the
  coefficients of polynomials in $\allpolysep\cap\allpolypos$. The idea is to apply
  Proposition~\ref{prop:hurwitz-totally-pos} to the interlacing $f\lesslesseq
  \Delta f$. For example, if
  $f=a_0+\cdots+a_3x^3\in\allpolysep\cap\allpolypos$ then from the
  Proposition the following is totally positive:
\[
\begin{pmatrix}
  a_1+a_2+a_3 & 2a_2 + 3a_3 & 3a_3 & 0 & \dots & \dots \\
  a_0 & a_1 & a_2 & a_3 & 0 & \dots  \\
0&  a_1+a_2+a_3 & 2a_2 + 3a_3 & 3a_3 &0 &   \dots  \\
0&   a_0 & a_1 & a_2 & a_3 & \dots\\
\vdots & \vdots & \vdots & \vdots & \vdots & \ddots
\end{pmatrix}
\]

For example, from the submatrix $
\begin{pmatrix}
  a_2 & a_2 \\ 2a_2+3a_3 & 3a_3
\end{pmatrix}
$ we get the inequality $a_2 \ge 3a_3$. The corresponding matrix for a
polynomial of degree $n$ is
\[
\begin{pmatrix}
  a_{n-1} & a_n \\ (n-1)a_{n-1} + \binom{n}{2}a_n & na_n
\end{pmatrix}
\]
from which we conclude that $a_{n-1}\ge \binom{n}{2}a_n$.

\section{The case $\affa x = -x$}
\label{sec:affine-minus-one}
In this section we look at the transformation $\affa(x)=-x$. This is a
degenerate case, but the associated linear transformation is
interesting. We can characterize the elements of $\allpolyaffine$ in
terms of pairs of interlacing polynomials.

\begin{lemma}
  If $\affa (x)=-x$ and $f(x)\greateqeq f(-x)$ then there are
  $f_0,f_1$ in $\allpolypos$ such that $f_0
  \longleftarrow f_1$ and $f = f_0(x) \,f_1(-x)$.
  Conversely, given $f_0,f_1$ satisfying these conditions then
  $f_0(x)\,f_1(-x)$ is in $\allpolyaffine$.
\end{lemma}
\begin{proof}
   Write $f(x) = f_0(x)
  f_1(-x)$ where $f_0(x)$ contains all the negative roots, and
  $f_1(-x)$ all the positive roots. Since $f(-x) = f_0(-x)f_1(x)$ the
  result follows easily from consideration of the interlacing
  of $f(x)$ and $f(-x)$. 
\end{proof}

The affine derivative of $f(x)$ is $(f(x) - f(-x))/2x$. If we write
\begin{align*}
  f(x) & = a_0 + a_1x + a_2x^2 + a_3x^3 + \cdots  \\
\intertext{then the affine derivative is}
 \daffine(f) & = \frac{1}{2}\, \,(a_1 + a_3x^2 + a_5 x^4 + \cdots ). \\
\intertext{If we recall the \emph{odd part} of $f$}
   f_o(z) & = a_1 + a_3z + a_5 z^2 + \cdots
\intertext{then we can write}
  \daffine(f) &= \frac{1}{2}\,f_o(x^2).\\
\intertext{An unexpected  consequence is that }
  \daffine(\,\daffine(f)\,) &=0.
\end{align*}
\index{odd part} \index{even part} We saw the odd part (and the even
part) of $f$ in \chapsec{operators}{hurwitz}. Since the roots of
$\daffine(f)$ are the square roots of the roots of $f_o$ we see that
in order to have $\daffine(f)\in\allpoly$ we must have that the signs
of $\daffine(f)$ alternate. See Lemma~\ref{lem:sign-patterns}.

One result of these computations is that if the degree of $f$ is even
then the degree of $\daffine(f)$ is two less than the degree of $f$,
and so does \emph{not} interlace $f$. It is also easy to verify that
we do not have interlacing if the degree is odd.

There is a linear transformation associated with $\affa$ that maps
$\allpolypos$ to $\allpolyaffine$. Define
$$T(x^n) = x\,\cdot\,\affa x\,\cdot\,\affa^2 x\,\cdots\,\affa^{n-1}x =
(-1)^{\binom{n}{2}}x^n$$
Unlike linear transformations $\allpoly\longrightarrow\allpoly$
(Corollary~\ref{cor:leading-coef}), the leading coefficients of this transformation
neither alternate, nor have the same sign.

\begin{lemma} \label{lem:affine-x}
If $T(x^n)=(-1)^{\binom{n}{2}}x^n$ then
\begin{enumerate}
\item $T\colon{}\allpolypos\cup\allpolyneg\longrightarrow \allpoly$
\item If $f\in\allpolypos(2n)$ then $(Tf)(x)\greateqeq (Tf)(-x)$.
\item If $f\in\allpolypos(2n+1)$ then $(Tf)(-x)\greateqeq (Tf)(x)$.
\end{enumerate}
\end{lemma}

\begin{proof}
  We first get the recurrence relation for $T$. 
\begin{align*}
    T(x\,x^n) &= (-1)^{\binom{n+1}{2}}x^{n+1} =
    x\,(-1)^{\binom{n}{2}}(-x)^n = x\, (Tx^n)(-x)\\
\intertext{and by linearity}
T(xf) &= x (Tf)(-x) \\
T((x+\alpha)f) &= x\,(Tf)(-x) + \alpha\,(Tf)(x)
\end{align*}

We use this recursion to prove the lemma by induction on the degree of
$f$. For $n=1$ we chose $a>0$. Now $T(x+a) = x+a$ so we have that
$-x+a\greateq x+a$ since $a$ is positive. 

We will do one of the cases; the others are similar. Suppose that $g =
(x+a)f\in\allpolypos(2n)$. Then $f\in\allpolyneg(2n-1)$ and by
induction hypothesis $(Tf)(-x) \greateqeq (Tf)(x)$. We want to show
that $(Tg)(x) \greateqeq (Tg)(-x)$. From the recurrence relation

\begin{align*}
  \begin{pmatrix}{Tg(x)}\\{A(Tg(x))}\end{pmatrix} = 
\begin{pmatrix}{(Tg)(x)}\\{(Tg)(-x)} \end{pmatrix}
&= \begin{pmatrix}{x(Tf)(-x) + a (Tf)(x)}\\{-(x(Tf)(x) + a (Tf)(-x))}\end{pmatrix}\\
&= \begin{pmatrix}{x}&{a}\\{a}&{-x}\end{pmatrix}
 \begin{pmatrix}{(Tf)(-x)}\\{(Tf)(x)}\end{pmatrix}
\end{align*}

$Tg \greateqeq \affa Tg$ follows from Corollary~\ref{cor:2by2} since
$\affa Tf \greateqeq Tf$.

\end{proof}

\section{The case $\affa x=0$}
\label{sec:affa0}

It is surprising that the degenerate case $\affa x=0$ leads to
interesting results. The class $\allpolyaffine$ is empty, since $\affa
f = f(0)$ is a constant, and can not interlace $f$. However, there is
a class of polynomials $\affa_0$ that are based on the affine
derivative and its conjugate.  \index{conjugate}

In this case the affine derivative is also called the 
\emph{lower truncation operator} since the effect of $\diffd_\affa$ is to remove the lowest
term: \index{truncation!operator} \index{lower truncation operator}
\begin{align*}
  \diffd_L(f) &= \frac{f(x)-f(0)}{x-0} \\
\diffd_L(a_0+a_1x+\dots+a_nx^n) &= a_1 + a_2x + \dots+a_nx^{n-1}
\intertext{The conjugate of the lower truncation
operator by the reversal operator is the upper truncation operator}
\diffd_U(a_0+a_1x+\dots+a_nx^n) &= a_0 + a_1x + \dots+a_{n-1}x^{n-1}
\end{align*}
Neither of these operators preserve roots \seepage{lem:trunc}. This
leads us to consider a set of polynomials for which both operators
preserve roots.  \index{conjugate}

\begin{definition}
  $\affa_0$ is the largest set of polynomials in $\allpolypos$ with
  all distinct roots such that $\diffd_U(\affa_0)\subset\affa_0$ and
$\diffd_L(\affa_0)\subset\affa_0$.
\end{definition}

It is unexpected that there is a characterization of the members of
$\affa_0$ in terms of coefficients, and that this characterization is
a kind of converse of Newton's inequalities (Theorem~\ref{thm:newton}).

\begin{cor} \label{cor:kurtz}
  A polynomial $f=\sum a_ix^i$ in $\allpolypos$ with all distinct
  roots is in $\affa_0$ if and only if \eqref{eqn:kurtz} holds.
\end{cor}
\begin{proof}
  If $f\in\affa_0$ then by applying the right combinations of
  $\diffd_U$ and $\diffd_L$ we find that $a_{i-1}+a_ix+a_{i+1}x^2$ is
  in $\affa_0$. Since all polynomials in $\affa_0$ have all distinct
  roots the inequalities follow.  
  
  Conversely, suppose that the coefficients of $f$ satisfy
  \eqref{eqn:kurtz}. If we apply either $\diffd_L$ or $\diffd_U$ to
  $f$ then the resulting polynomial $g$ has the same coefficients as
  $f$, except for either $a_0$ or $a_n$. Consequently, $g$ also
  satisfies the conditions of Theorem~\ref{thm:kurtz} and so has all real roots,
  all distinct.
\end{proof}

The class $\affa_0$ satisfies the usual closure properties

\index{exponential operator}
\index{operator!exponential}

\begin{cor}
  $\affa_0$ is closed under differentiation and the exponential operator.
\end{cor}
\begin{proof}
  If $f=\sum a_ix^i$ is in $\affa_0$ then we only need check that the
  coefficients of $f^\prime$ and $\expoper{f}$ satisfy
  \eqref{eqn:kurtz}. This is immediate from \eqref{eqn:kurtz} and
  \begin{align*}
    i^2a_i^2 & \ge 4\,(i-1)(i+1)\,a_{i-1}a_{i+1} \\
    \frac{a_i^2}{(i!)^2}  & \ge 4\,\frac{a_{i-1}}{(i-1)!}\,\frac{a_{i+1}}{(i+1)!} 
  \end{align*}
\end{proof}

It is not the case that $\diffd_L(f)$ and $f$ interlace. In fact, they
are as far from sign interlacing as possible. If $r$ is a non-zero root of
$f\in\allpolypm$, then $\diffd_L(f)(r) = -f(0)/r$, and so
$\diffd_L(f)$ has the same sign at all roots of $f$.


\part{Polynomials in several variables}

\chapter{Polynomials in two variables }
\label{cha:p2}

\renewcommand{\TimeStampStart}{Tuesday, March 11, 2008: 09:48:06}
\mytoday    

In this chapter we generalize our results from polynomials in one
variable to polynomials in two variables. Our  goal is to
generalize $\allpoly$ and $\allpolypos$, to define interlacing even
though there are no roots, and to use these facts to deduce properties
of polynomials in one variable.

\section{The substitution property and  determinants}
\label{sec:substitution-defn}

All the polynomials that we will consider in this chapter satisfy a
property called \emph{substitution}. This is the analog of ``all real
roots'' for polynomials in two variables. However, this alone is not
sufficient; we will need another condition to generalize $\allpoly$.

\index{substitution}
\index{substitution!x-substitution}
\begin{definition}
  If $f$ is a polynomial in variables $x,y$ then $f$ satisfies
  \emph{$x$-substitution} if for every choice of $a$ in $\reals$ the
  polynomial $f(x,a)$ has all real roots, and the number of roots is
  the same for all $a$.  We say that $f$ satisfies \emph{substitution}
  if it satisfies $x$-substitution and $y$-substitution.  We let
  $\xsub_2$ be the set of all polynomials in two variables that
  satisfy substitution.
\end{definition}

\index{\ sub2@$\xsub_2$}

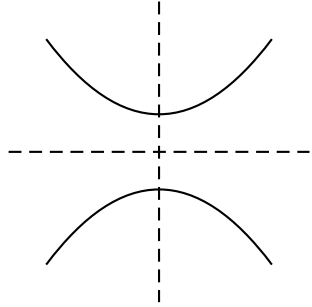
\begin{figure}[htbp]
  \begin{center}
    
\psset{unit=1cm}
\begin{pspicture}(-2,-2)(2,2)
\psline[linestyle=dashed](-2,0)(2,0)
\psline[linestyle=dashed](0,2)(0,-2)
\parabola(1.5,1.5)(0,.5)
\parabola(-1.5,-1.5)(0,-.5)
\end{pspicture}

    \caption{The graph of a polynomial  satisfying $y$ but not $x$-substitution. }
    \label{fig:not-sub}
  \end{center}
\end{figure}

A polynomial can satisfy
$x$-substitution and not $y$-substitution.  Figure~\ref{fig:not-sub}
shows the graph $f(x,y)=0$ of such a polynomial of total degree $2$.
Every vertical line meets the graph in two points, so $f$ satisfies
$y$-substitution. Since there are some horizontal lines that do
not meet the graph, $f$ does not satisfy $x$-substitution.

In the remainder of this section we discuss the substitution
properties of polynomials that are defined by determinants.
\index{substitution!arising from determinants}

\begin{example}\label{ex:lax-example}
  We start with an  example of a polynomial that satisfies a different kind
  of substitution.  Recall that the eigenvalues of a symmetric matrix
  are all real.  Choose symmetric matrices $A,B$ and consider
  the polynomial 
\begin{equation}\label{eqn:sub-from-det}
f(x,y) = |I+xA+yB|
\end{equation}
where $I$ is the identity
  matrix.  We claim that $f(x,y)$ satisfies the property:
  \begin{itemize}
  \item For any $\alpha,\beta$ the polynomial $f(\alpha z,\beta z)$
    has all real roots.
  \end{itemize}
  Indeed, we see that
  $$
  g(z)=f(\alpha z,\beta z) = | I + z(\alpha A + \beta B)|.$$  Since
  $\alpha A +\beta B$ is symmetric, the roots of $g(z)$ are given by
  $-r^{-1}$ where $r$ is an eigenvalue of $\alpha A+\beta B$.
  
  It is not the case that $|I+x A + yB|$ satisfies substitution. Here's
  a small example:
  \begin{gather*}
    A = 
    \begin{pmatrix}
      2 & 4 \\ 4 & 4
    \end{pmatrix} \quad\quad
B = 
\begin{pmatrix}
  2 & 4 \\ 4 & 2
\end{pmatrix}\\
f(x,y) = |I+xA+yB| = 1 + 6\,x - 8\,x^2 + 4\,y - 20\,x\,y - 12\,y^2\\
f(x,1) = -7 - 14\,x - 8\,x^2
  \end{gather*}
and the latter polynomial has two complex
roots. Figure~\ref{fig:lax-1} shows why every line through the origin
will meet the graph, yet there are horizontal and vertical lines that
miss it. 
\begin{figure}[htbp]
  \centering
  \includegraphics*[width=1.5in]{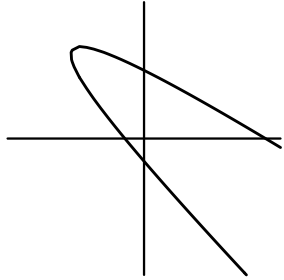}
  \caption{Meets all lines through the origin, fails substitution}
  \label{fig:lax-1}
\end{figure}
\end{example}

\begin{example}
  We again start with equation~\eqref{eqn:sub-from-det}.  Now assume that
    \begin{gather*}
      A,B,A^{-1}B \text{ are invertible and symmetric}\\
    \end{gather*}
    We claim that 
\begin{itemize}
\item For any $\alpha,\beta$ the polynomial $f(\alpha z,\beta z)$
    has all real roots.
  \item $f(x,y)$ satisfies $x$ and $y$ substitution.
\end{itemize}
    
        Choose $\alpha\in\reals$ and observe
\begin{gather*}
f(x,\alpha) = \bigl| A\bigr| \,\cdot\, \bigl| xI + (\alpha A^{-1}B+
A^{-1})\bigr|
    \end{gather*}
    Since $A^{-1}B$ and $A$ are symmetric by hypothesis, $\alpha
    A^{-1}B+A^{-1}$ is symmetric. Thus the roots of $f(x,\alpha)$ are
    all real since they are the negative of the eigenvalues of a
    symmetric matrix. For $f(\alpha,x)$, notice that $B^{-1}A$ is
    symmetric, since it's the inverse of a symmetric matrix.
    
    It is easy to find matrices that satisfy the assumptions of this
    example: take $A$ and $B$ to commute. For example, we take
$$ A = 
\
\begin{pmatrix}
 1 & -1 & -1 & 1 \\ -1 & -1 & 0 & -3 \\ -1 & 0 & -2 & 1 \\ 1 & -3 & 1 & 1   
\end{pmatrix}
$$
and let $B=A^2/10$. Figure~\ref{fig:matrix-not-p2} shows the graph of
the determinant, where the dot is the origin. The segments are not
linear, even though they appear to be. It is clear that it satisfies
substitution, and every line through the origin meets the graph in 4
points. 

\begin{figure}[htbp]
  \begin{center}
    \leavevmode
    \includegraphics*[height=2in]{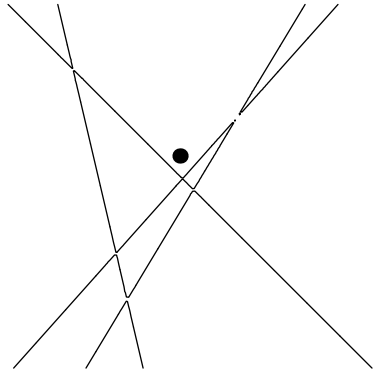}
    \caption{The graph of a determinant}
    \label{fig:matrix-not-p2}
  \end{center}
\end{figure}

  \end{example}

\begin{example}
    Next, we assume that $A,B$ are symmetric and positive definite $n$
    by $n$ matrices. We claim that $f(x,y)$ in
    \eqref{eqn:sub-from-det} satisfies these properties:
    \begin{itemize}
\item For any $\alpha,\beta$ the polynomial $f(\alpha z,\beta z)$
    has all real roots.
    \item $f$ satisfies $x$ and $y$ substitution.
    \item All coefficients of $f$ are positive.
    \end{itemize}

The first part follows as before since $A,B$ are symmetric. Next, 
$A$ is positive definite, and so it has a positive definite
 square root. Let $E^2=A$ where $E$ is positive definite. Factoring
 out $|A|$ yields
$$ f(x,y) = |A|\cdot|xI + y\,E^{-1}BE^{-1}+A^{-1}|$$
This representation shows that $f(x,\alpha)$ has all real
roots. Factoring out $B$ shows that $f(\alpha,y)$ has all real roots.
The fact that all the coefficients are positive follows from a more
general fact \cite{speyer}. 

\begin{lemma}\label{lem:pos-def-sum}
  Suppose that $A_1,\dots,A_d$ are positive definite. All the
  coefficients of $|x_1A_1+\cdots+x_dA_d|$ are positive.
\end{lemma}
\begin{proof}
  We prove it by induction, and it is immediate for $d=1$.  Assume all
  matrices are $n$ by $n$. If $\diffi$ is a subset of
  $\{1,2,\dots,n\}$, and $M$ is a matrix, then $M[\diffi]$ denotes the
  matrix composed from the rows and columns of $M$ listed in $\diffi$.
  Since $A_1$ is positive definite, we may conjugate by $A_1^{-1/2}$,
  and so assume that $A_1=I$, the identity. Upon expanding the
  determinant we see that

$$|x_1I+(x_2A_2+\cdots+x_dA_d)| = \sum_\diffi x_1^{n-|\diffi|}(x_2A_2+\cdots
x_dA_d)[\diffi]$$
Now any principle submatrix of a positive definite matrix is positive
definite, so all the terms in the sum are positive definite. Thus, by
the inductive hypothesis, all the summands have positive coefficients,
and so the lemma is proved.
\end{proof}

The following is a sort of converse to the above, in that it allows us
to determine that the eigenvalues are positive from properties of the
polynomials. 

\begin{lemma}\label{lem:pos-eigen}
 If $A$ is a matrix and  $|I+xA|\in\allpolypos$ then $A$ has all
 positive eigenvalues.
\end{lemma}
\begin{proof}
  If the eigenvalues of $A$ are $r_1,\dots,r_n$ then $|I+xA| =
  \prod(1-xr_i)$. Since $|I+xA|$ is in $\allpolypos$ it has all
  negative roots, which implies that all $r_i$ are positive.
\end{proof}

A consequence of this lemma is a kind of converse to
Lemma~\ref{lem:pos-def-sum}. 

\begin{cor}\label{cor:pos-eigen}
 If $A,B$ are matrices, and $f(x,y) = |I+xA+yB|$ is a polynomial
  that satisfies
  \begin{enumerate}
  \item All coefficients of $f$ are positive.
  \item $f(x,0)\in\allpoly$.
  \item $f(0,x)\in\allpoly$.
  \end{enumerate}
then $A$ and $B$ have positive eigenvalues.
\end{cor}
\end{example}

  \begin{example}\label{ex:det-aic}
    Here is another way of constructing matrices that satisfy
    substitution that we will study in detail later. We assume that
    $A$ and $C$ are symmetric $n$ by $n$ matrices, $A$ is positive
    definite, and define
 \begin{equation}
   \label{eqn:p2-from-det}
   f(x,y) = \det(xA+yI+C)
 \end{equation}
We claim that 
 \begin{itemize}
 \item $f(x,y)$ satisfies $x$ and $y$ substitution.
 \item The coefficients of terms of degree $n$ are positive.
 \end{itemize}

 By assumption $\alpha A+C$ is symmetric, so $f(\alpha,y)$ has all
 real roots, since its roots are the negative eigenvalues of $\alpha
 A+C$.  Let $A=E^2$ where $E$ is positive definite. Since
$$ f(x,y) = |A|\cdot|xI+yA^{-1}+E^{-1}CE^{-1}|$$
and $E^{-1}CE^{-1}$ is symmetric, it follows that $f(x,\alpha)$ has
all real roots.

The polynomial $f^H(x,y)$ formed by the terms of degree $n$ is
 $\det(xA+yI)$, and the roots of $f^H(1,y)$ are the negative
 eigenvalues of $A$. Thus, all roots of $f(1,y)$ are negative since $A$ is
 positive definite, and hence $f^H$ has all positive coefficients.

  \end{example}

We summarize the different matrices in
Table~\ref{tab:prop-of-mat}. 

\added{16/12/6}

\begin{table}[htbp]
  \centering
  \begin{tabular}{cccc|ccc}
\toprule
$A$ & $B$ & $A^{-1}B$ &  $C$ &  Sub. & $f^H$ & Coefficients \\
&&&&& positive & positive \\
\midrule
sym. & sym. & & I &  & & \\
sym. & sym. & sym. & I &  \bfsf X & & \\
posdef. & posdef & & I &  \bfsf X & \bfsf X & \bfsf X \\
posdef. & I & & sym. &  \bfsf X &  \bfsf X & \\
\bottomrule
  \end{tabular}
  \caption{Properties of $|Ax+By+C|$}
  \label{tab:prop-of-mat}
\end{table}

  We can compute an explicit example of Lemma~\ref{lem:pos-def-sum}
  that shows all terms are non-negative; this properly belongs in a
  later chapter since the matrices are only positive semi-definite.

  \begin{lemma}
    Suppose that $X = diag(x_i)$, $Y = diag(y_i)$ are $d$ by $d$
    diagonal matrices, and $v_1,\dots,v_d$ are $d$-vectors. Let $M =
    (v_1^t,\dots,v_d^t)^t$. Then
\[
|X + y_1 v_1^tv + \cdots + y_d v_d^tv_d| = \sum_{I\subset\{1,\dots,d\}}
X_{I'}\,Y_I\,|M[I])|^2
\]
  \end{lemma}
  \begin{proof}
    \index{Neil White} By induction 
    it\footnote{thanks to Neil White for this argument} 
    suffices to evaluate $|y_1v_1^tv + \cdots +
    y_d v_d^tv_d|$. But this can be written as
\[
|M^tYM| = |M^t||Y||M| = y_1\cdots y_d\,|M|^2
\]
  \end{proof}

  \begin{cor}
    If $A$ is symmetric, $|I+xA|=\prod_1^d(1+r_ix)$, $v_1\dots,v_d$
    are $d$-vectors then
    \begin{equation*}
      \label{eqn:psd-det}
      |I+xA+y_1v_1^tv_1+\cdots +y_dv_d^tv_d| =
\sum_{I\subset\{1,\dots,d\}} \frac{f(x)}{\prod_{i\in I}(1+r_ix)}\,Y_I \, |M[I]|^2
    \end{equation*}
  \end{cor}
  \begin{proof}
    Diagonalize $A$ and apply the lemma.
  \end{proof}
If $d=2$ we have
\begin{multline}\label{eqn:det-psd-2}
  |I + xD + y_1v_1v_1^t + y_2v_2v_2^t| =\\
= f + y_1 \sum v_{1,i}^2\frac{f}{1+xr_i} +y_2 \sum
v_{2,i}^2\frac{f}{1+xr_i} \\ + 
y_1y_2 \sum_{i<j} \begin{vmatrix} v_{1,i}& v_{1,j}
  \\ v_{2,i} & v_{2,j}
\end{vmatrix}^2 \frac{f}{(1+xr_i)(1+xr_j)} 
\end{multline}
The statement for $n=3$ is
\begin{multline}\label{eqn:det-psd-3}
  |I + xD + y_1v_1v_1^t + y_2v_2v_2^t + y_3v_3v_3^t|  \\
=f+y_1 \sum v_{1,i}^2\frac{f}{1+xr_i} +
y_2 \sum v_{2,i}^2\frac{f}{1+xr_i} + 
y_3 \sum v_{3,i}^2\frac{f}{1+xr_i} 
\\
\ + y_1y_2 \sum_{i<j} 
\begin{vmatrix} 
v_{1,i}& v_{1,j}  \\ v_{2,i} & v_{2,j}
\end{vmatrix}^2 
\frac{f}{(1+xr_i)(1+xr_j)}  \\
\ + y_1y_3 \sum_{i<j} 
\begin{vmatrix} v_{1,i}& v_{1,j}
  \\ v_{3,i} & v_{3,j}
\end{vmatrix}^2 
\frac{f}{(1+xr_i)(1+xr_j)} \\  
\ + y_2y_3 \sum_{i<j} 
\begin{vmatrix} 
v_{2,i}& v_{2,j}
  \\ v_{3,i} & v_{3,j}
\end{vmatrix}^2 
\frac{f}{(1+xr_i)(1+xr_j)} \\ 
\ +y_1y_2y_3 \sum_{i<j<k} 
\begin{vmatrix} 
v_{1,i}& v_{1,j} & v_{1,k} \\
v_{2,i}& v_{2,j} & v_{2,k} \\
v_{3,i}& v_{3,j} & v_{3,k} 
\end{vmatrix}^2 
\frac{f}{(1+xr_i)(1+xr_j)(1+xr_k)} 
\end{multline}

\section{Interlacing in   $\xsub_2$}
\label{sec:sub-interlacing}

Interlacing can be easily defined for polynomials satisfying
substitution. 

\index{\ zzlesslesseq@$\lesslesseq$!in $\xsub_d$}
\index{interlacing!in $\xsub_d$}
\begin{definition}
  If $f$ and $g$ are polynomials that satisfy substitution, then we
  say that $f$ and $g$ \emph{interlace} if for every real $\alpha$ the polynomial
  $f+\alpha g$ satisfies substitution. Using the definition of
  substitution, we can restate this in terms of one variable
  polynomials:
\begin{quote}
  If $f,g\in\xsub_2$ then $f$ and $g$ interlace iff for all
  $a\in\reals$ we have that $f(x,a)$ and $g(x,a)$ interlace, as do
  $f(a,y)$ and $g(a,y)$.
\end{quote}
  We  define $\lesslesseq$ in terms of these one variable
  polynomials.
  \begin{quote}
  If $f,g\in\xsub_2$ then 
  $f \lesslesseq g$  iff    for all $a\in\reals$  we have
  $f(x,a)\lesslesseq g(x,a)$ and $f(a,y)\lesslesseq g(a,y)$.
  \end{quote}
\end{definition}


For homogeneous polynomials in two variables we can verify interlacing
by reducing to a one variable problem.
\begin{lemma}
  Suppose $f,g\in\xsub_2$ are homogeneous and that the coefficients of $f$ 
  are all positive.  The following are equivalent
  \begin{enumerate}
  \item $f\lesslesseq g$
  \item $f(1,y) \lesslesseq g(1,y)$
  \item $f(x,1) \lesslesseq g(x,1)$
  \end{enumerate}
\end{lemma}
\begin{proof}
  Note that the first interlacing is in $\xsub_2$, and the last two in
  $\allpoly$. We can substitute $x=1$ in the first interlacing to
  deduce the second, so assume $f(1,y) \lesslesseq g(1,y)$. By
  assumption $f(1,y)\in\allpolypm$. If the roots of $f(1,y)$ are
  $r_1,\dots,r_n$ then the roots of $f(a,y)$ are $\{r_i/a\},$ and the
  roots of $f(x,b)$ are $\{b/r_i\}$. There are similar expressions for
  the roots of $g$. Since the roots are either all positive, we have
  $f(a,y)\lesslesseq g(a,y)$ and $f(x,a) \lesslesseq g(x,a)$. Thus $f$
  and $g$ interlace in $\xsub_2$.  Similarly we can show that 1 and 3
  are equivalent.
\end{proof}






The relation of interlacing is preserved under limits.

\begin{lemma} \label{lem:sub2-closure-2}
  Suppose $f_i,g_i$ are in $\xsub_2$, and $f_i\rightarrow f$,
  $g_i\longrightarrow g$ where $f,g$ are polynomials. If $f_i$ and
  $g_i$ interlace for all $i$ then $f$ and $g$ interlace.
\end{lemma}

\begin{proof}
  Since $f_i+\alpha g_i$ is in $\xsub_2$ for all $i$, and $f_i+\alpha
  g_i$ converges to $f+\alpha g$ it follows that $f+\alpha
  g\in\xsub_2$.  (See Lemma~\ref{lem:sub-closure-1}.) Consequently $f$
  and $g$ interlace.
\end{proof}




Substitution is \emph{not} preserved under the operations of
differentiation\footnote{ The partial derivative $\frac{\partial
    f}{\partial x_1}$ of $ f = (x_1+x_2+1)(x_1+x_2)(x_1-x_2)$ has
  imaginary roots for $x=-1/4$.}.  Consequently, we must restrict
ourselves to a subset of $\xsub_2$.

\section{Polynomials in two variables}
\label{sec:leading-coefficient}

We continue our journey toward the generalization of $\allpoly$ to
polynomials in two variables.  For polynomials of one variable it is
often important to know the sign of the leading coefficient.  For two
variables the homogeneous part is the analog of the leading
coefficient. The homogeneous part will determine the asymptotic
behavior of the graph, which is central to our generalization.

If $f$ is a polynomial in two variables, and we let $x$ and $y$ get
large simultaneously then the behavior of $f(x,y)$ is determined by
the terms of highest total degree. If the maximum total degree is $n$
then

\index{homogeneous part}
\index{polynomials!homogeneous part}
\index{\ aaaH@$f^H$}

 \begin{equation}
   \label{eqn:O-2}
   f(x,y) = \sum_{i+j=n} c_{ij}x^iy^j \quad+\quad
   \sum_{i+j<n} c_{ij}x^iy^j 
 \end{equation}
 
 The {homogeneous part} $f^H$ is the first sum in \eqref{eqn:O-2}.


\index{positivity condition}
\index{condition!positivity}

\begin{definition} \label{defn:positivity-condition-2}
A homogeneous polynomial $f(x,y)=c_0x^n + \cdots + c_ny^n$  satisfies
the \emph{positivity   condition} iff all $c_i$ are positive.


\end{definition}

Using the homogeneous part, we define the
class $\gsubpos_2$ of polynomials that forms the $2$-dimensional
analog of ``polynomials with all real roots''.

\index{\ Pgsub2@$\gsubpos_2$}

\begin{definition} \label{def:p2}
  
  $\gsubpos_2(n)$ consists of all polynomials $f$ of degree $n$ such that $f$
  satisfies $x$-substitution and $y$-substitution, and $f^H$ satisfies
  the positivity condition.

$$\gsubpos_2  = \gsubpos_2(1) \cup \gsubpos_2(2) \cup \gsubpos_2(3)
\cup \dots$$


Sometimes the homogeneous part of a polynomial might be negative. For
instance, if $f(x,y)\in\rupint{2}(n)$, then the homogeneous part of
$f(-x,-y)$ has all negative signs if $n$ is odd. In this case,
$-f(-x,-y)\in\rupint{2}$. If $n$ is even then $f(-x,-y)\in\gsub_2$. To
simplify exposition, we just write $f(-x,-y)\in\pm\rupint{2}$, which is
true for every $n$. We  also express this as $\pm f(-x,-y)\in\rupint{2}$.
\end{definition}

It's useful to note that $f^H$ is not arbitrary.

\begin{lemma} \label{lem:sub-fH2}
  If $f\in\gsubpos_2$ then $f^H\in\allpolypos$.
\end{lemma}
\begin{proof}
  See the proof of the more general result (Lemma~\ref{lem:sub-fH}).
\end{proof}

Here are a few elementary facts about $\rupint{2}$ and $f^H$.

\begin{lemma}\ 
  \begin{enumerate}
  \item   If $f(x,y)\in\gsubpos_2$ then $f(y,x)\in\gsubpos_2$.
  \item   If also $g(x,y)\in\rupint{2}$, then $f(x,y)g(x,y)\in\gsub_2$.
  \item If $a_1,\dots,a_n$ are positive, and $b_1,\dots,b_n$ are
    arbitrary, then 
$$ (x+a_1y+b_1)(x+a_2y+b_2)\cdots(x+a_ny+b_n)\in\rupint{2}$$
\item If $g\in\gsubpos_2$ then $(fg)^H = f^H\,g^H$.
  \item If $g\in\gsubpos_2$ has degree $n-1$ then $(f+g)^H = f^H$.
  \item $(f^H)^H = f^H$.
  \item $\left(\frac{\partial f}{\partial x}\right)^H = 
\frac{\partial}{\partial x}\,f^H$
  \end{enumerate}
\end{lemma}

\begin{proof}
  $f(y,x)^H$ is the reverse of $f(x,y)^H$ and $f(y,x)$ satisfies
  substitution since $f(x,y)$ does. The second part follows from the
  $(fg)^H= f^H\cdot g^H$. Since each factor $x+a_iy+b$ is easily seen
  to be in $\rupint{2}$, the product is in $\gsub_2$.  The rest follow
  easily from the definitions.
\end{proof}

Let's explore the definition of $\gsubpos_2$. We can write $f(x,y)$ as a
polynomial in either $x$ or $y$:

\begin{eqnarray} \label{eqn:p3-1} 
f(x,y) =& f_0(x)+f_1(x)y + \cdots + f_n(x)y^n\\
=& f^0(y)+f^1(y)x + \cdots + f^m(y)x^m \label{eqn:p3-2} \\
=& \sum_{i,j} c_{ij} x^i y^j \label{eqn:p3-3}
\end{eqnarray}

The
homogeneous part $f^H$ is $c_{0n}y^n + \cdots c_{n0}x^n$.
Consequently,  $n=m$, and  $f_n(x)$ and $f^n(x)$ are
non-zero constants. The coefficient polynomials $f_i(x)$ have
$x$-degree  $n-i$ since the total degree of $f$ is $n$ and 
the leading coefficient of $f_i$ is $c_{i,n-i}$ which is non-zero.  We
can summarize:

 \begin{quote}
   A polynomial $f\in\xsub_2$ given in
   \eqref{eqn:p3-1},\eqref{eqn:p3-2} is in $\gsubpos_2$ iff
   \begin{itemize}
   \item $n=m$
   \item The degree of $f_i$ and $f^i$ is $n-i$.
   \item The leading coefficients of $f_0,\dots,f_n$ (and
     $f^0,\dots,f^n$) are positive.
   \end{itemize}
 \end{quote}
 
 Of course the coefficient polynomials are highly interrelated. We
 will see shortly that consecutive coefficients interlace. A simple
 property is that if $f$ in \eqref{eqn:p3-1} is in $\gsubpos_2$ then
 since $f(x,1)\in\allpoly$ we have that
$$ f_0 + f_1 + \cdots +f_n\in\allpoly$$

\index{mutually interlacing!signs}
If $f_1\greateqeq \cdots \greateqeq f_n$ is a sequence of mutually
interlacing polynomials then we know that, from left to right, the
roots appear in reverse order: $f_n,\dots,f_1$. It follows that for
any $\alpha\in\reals$ there is an $i$ such that
\begin{align*}
  f_1(\alpha),\dots,f_i(\alpha) &\quad\text{all have the same sign $\epsilon$}\\
  f_{i+1}(\alpha),\dots,f_n(\alpha) &\quad\text{all have the same sign
  $-\epsilon$}.
\end{align*}

We can use this observation to construct a family of polynomials that
satisfy substitution and have all coefficients in $\allpoly$, but are
not in $\rupint{2}$. 

\begin{lemma}\label{lem:mi-2}
Suppose $f_1\greateqeq\cdots\greateqeq f_n$ and
  $g_1\greateqeq\cdots\greateqeq g_n$ are two sequences of mutually
  interlacing polynomials. If 
\[
h(x,y)  = \sum f_i(x)g_i(y) = \sum F_i(x)y^i = \sum G_i(y)x^i 
\]
then
\begin{enumerate}
\item $h(x,\alpha)$ and $h(\alpha,y)$ are in $\allpoly$ for all $\alpha\in\reals$.
\item All $F_i$ and $G_i$ are in $\allpoly$.
\item Consecutive $F_i$'s and $G_i$'s interlace.
\end{enumerate}
\end{lemma}
\begin{proof}
  Note that  $h(x,\alpha) = \sum f_i(x)g_i(\alpha)$.  The
  observation above shows that $g_i(\alpha)$ has one sign change, so
  we can apply Lemma~\ref{lem:sign-1}. This establishes the first
  part.

  Since $F_0 = h(x,0)$ we see that $F_0\in\allpoly$. Now since
  $\{g^\prime(y)\}$ is a sequence of mutually interlacing polynomials we
  know that $\sum f_i(x)g_i^\prime(y)$ satisfies the first part, and its
  constant term is $F_1$. Continuing, all $F_i$ are in $\allpoly$.

  Next, for any $\beta\in\reals$ the sequence $\{(x+\beta)g_i(x)\}$ is
  mutually interlacing. By the second part the coefficient $\beta
  F_i+F_{i-1}$ of $y^i$ is in $\allpoly$. It follows that $F_i$ and
  $F_{i-1}$ interlace.
\end{proof}

Although substitution holds for $h$, the degree condition does not. And, it is
not true that $h\in\gsubclose_2$ (take $f_i=g_i)$. However, by
Lemma~\ref{lem:mi-same-order} we know that $h$ is stable.
\index{polynomial!stable}
\index{interlacing!mutually}

\section{Inequalities for coefficients of a  quadratic in $\rupint{2}$} 

Let's look at the simplest non-trivial polynomials in $\rupint{2}$.
Consider a quadratic polynomial $f(x,y)\in\rupint{2}$ where we write the
terms in a grid

\centerline{\xymatrix{
a_{02}\, y^2 \ar@{-}[d] \ar@{-}[dr]  & & \\
a_{01}\, y \ar@{-}[dr] \ar@{-}[d] \ar@{-}[r] & \ar@{-}[dr]  \ar@{-}[d] a_{11}\, xy & \\
a_{00} \ar@{-}[r] & \ar@{-}[r] a_{10}\, x & a_{20}\, x^2
}}

We have inequalities for the coefficients on each of the three outside
lines, and the center square.

\index{quadratics in $\rupint{2}$}
\begin{lemma}\label{lem:p2-quad-ineq}
  If $f\in\rupint{2}$ is a quadratic with the above coefficients, then
  \begin{align*}
    0 & \le a_{01}^2 - 4 a_{02}a_{00} \\
    0 & \le a_{10}^2 - 4 a_{20}a_{00} \\
    0 & \le a_{11}^2 - 4 a_{02}a_{20} \\
    0 & \le a_{01}a_{10} - a_{00}a_{11}
  \end{align*}
\end{lemma}
\begin{proof}
  Since $f(x,0)\in\allpoly$, the terms on the bottom row constitute
  the coefficients of a polynomial in $\allpoly$, and the inequality
  is just Newton's inequality. Similarly for the left most terms. The
  diagonal terms are the coefficients of $f^H$, which is also a
  polynomial in $\allpoly$.
  
  If we solve the equation for $x$, then the discriminant is
  \begin{align*}
    \Delta_x &= {\left( a_{10} + a_{11}\,y \right) }^2 -
    4\,a_{20}\,\left( a_{00} + a_{01}\,y + a_{02}\,y^2 \right) \\
    \intertext{The discriminant of $\Delta_x$ as a function of $y$ is}
    \Delta_y & = 16\left(a_{02}\,{a_{10}}^2 - a_{01}\,a_{10}\,a_{11} +
    a_{00}\,{a_{11}}^2 + {a_{01}}^2\,a_{20} -
    4\,a_{00}\,a_{02}\,a_{20}\right) \\
    \intertext{Since $\Delta_x\ge0$, we know $\Delta_y\le0$. Rewriting
      $\Delta_y\le0$ yields}
&   a_{11}\left(a_{10}\,a_{01} - a_{00}\,a_{11}\right) \ge
  a_{02}\left(a_{10}^2 - 4 a_{00}\,a_{20}\right) + a_{20}\,a_{01}^2
  \end{align*}
  Now the all coefficients $a_{20},a_{11},a_{02}$ are positive
  since $f^H\in\allpolypos$, and by the above $a_{10}^2 - 4
  a_{00}\,a_{20}\ge0$, so we conclude $a_{10}\,a_{01} - a_{00}\,a_{11}\ge0$.
\end{proof}

\begin{remark}\label{rem:p2-quad-ineq}
    We can determine when each of the inequalities is strict. The
    first three are strict when the corresponding polynomials have all
    distinct roots. If $a_{10}a_{01}-a_{00}a_{11}=0$ then
    \[ 0 = a_{02}\left(a_{10}^2 - 4 a_{00}\,a_{20}\right) +
    a_{20}\,a_{01}^2 \] and both summands are zero by the lemma. Since
    $a_{20}>0$ we see $a_{01}=0$. Interchanging the roles of $x$ and
    of $y$ shows that $a_{10}=0$. Finally, $a_{10}^2 - 4
    a_{00}\,a_{22}=0$ implies $a_{00}=0$. We conclude that $f(x,y)$
    looks like 

\centerline{\xymatrix@=.4cm{
        a_{02}\, y^2 \ar@{-}[d] \ar@{-}[dr]  & & \\
        0 \ar@{-}[dr] \ar@{-}[d] \ar@{-}[r] & \ar@{-}[dr]  \ar@{-}[d] a_{11}\, xy & \\
        0 \ar@{-}[r] & \ar@{-}[r] 0 & a_{20}\, x^2 }}

  \end{remark}

\begin{remark}
    Polynomials of degree two in $\rupint{2}$ are more than three
    quadratic polynomials spliced together. Consider the two variable
    polynomial 
\[
\begin{array}{ccc}
16\,y^2 \\ 8\,y & 17\,x\,y \\ 1 & 2\,x & 1\,x^2
\end{array}
\]
The three polynomials on the boundary ($x^2+2x+1$, $x^2+17x+16$,
$x^2+8x+16$) have all real roots, but the quadrilateral inequality is
not satisfied, so the two variable polynomial is not in $\rupint{2}$.
This phenomena will be considered later \mypage{sec:hives and horns}.

It is clear from the proof that if all four conditions of the lemma
are satisfied then the polynomial is in $\rupint{2}$. 
  \end{remark}

\index{quadrilateral inequalities}
\index{rhombus inequalities}
\begin{remark}
  We will see later \mypage{prop:p2plus-inequality} that the last
  inequality (the quadrilateral or rhombus inequality) is a
  consequence of interlacing of the coefficient polynomials. It can't
  be improved.  For example, the diagram of coefficients of the
  polynomial $(x+\epsilon\, y+1)(\epsilon\, x+y+1)$ is

\centerline{\xymatrix{
        \epsilon\, y^2 \ar@{-}[d] \ar@{-}[dr]  & & \\
        (1+\epsilon{})\, y \ar@{-}[d] \ar@{-}[r] \ar@{-}[dr]& \ar@{-}[dr]
        \ar@{-}[d] (1+\epsilon^2)\, xy & \\ 
        1 \ar@{-}[r] & \ar@{-}[r] (1+\epsilon{})\, x & \epsilon{ }\,x^2 }}
    
\noindent%
    and the ratio $(a_{01}a_{10})/(a_{00}a_{11})$ equals 
\[
\frac{(1+\epsilon)^2}{1+\epsilon^2}
\]
and this goes to $1$ as $\epsilon\rightarrow0$.
If $r,s>0$ and we let 
\[
f(x,y) = (x+ \epsilon(y+1))^r \,(y+\epsilon(x+1))^s\, (\epsilon
x+y+1)(x+\epsilon y+1)
\]
then as $\epsilon\rightarrow0$ the ratio $(a_{r,s+1}a_{r+1,s})/(a_{r,s}a_{r+1,s+1})$
converges to $1$. Thus, there are no Newton inequalities with constant
greater than $1$. 

The last inequality in Lemma~\ref{lem:p2-quad-ineq} corresponds to the
two adjacent triangles with a vertex in the lower left corner. The two
other corners do not give the inequalities
    \begin{align}
      \label{eqn:quad-ineq}
      a_{10}a_{11} - a_{01}a_{20} &> 0 \quad\quad(\text{right corner})\\
      a_{01}a_{11} - a_{10}a_{02} &> 0 \quad\quad(\text{upper corner})\notag
    \end{align}
as we see by expanding $(x+y-1)^2$:

\index{quadrilateral inequalities}
\index{rhombus inequalities}

\centerline{\xymatrix{
         y^2 \ar@{-}[d] \ar@{-}[dr]  & & \\
        -2\, y \ar@{-}[d] \ar@{-}[r] \ar@{-}[dr]& \ar@{-}[dr]  \ar@{-}[d] 2\, xy & \\
        1 \ar@{-}[r] & \ar@{-}[r] {-2}\, x & x^2 }}

    We will see later (Proposition~\ref{prop:p2plus-inequality} that
    if the coefficients are positive then the inequalities
    \eqref{eqn:quad-ineq} do hold.
  \end{remark}

\begin{example} \label{ex:bad-polly}
  We can use these ideas to show that perturbations of
  products aren't necessarily in $\rupint{2}$. Let
  $$
  f = e + (x+ay+b)(x+cy+d)$$
  where $a\ne c$ are positive. The discriminant $\Delta_y$ is
  $16(a-c)^2e$. Since this is positive for positive $e$, we conclude
  that there are perturbations of products that are not in $\rupint{2}$.

  For example, if $f = (x+y)(x+2y)+e$, then $f\in\rupint{2}$ if $e\le0$,
  and $f\not\in\rupint{2}$ if $e\ge0$. 
\end{example}

  Although it is difficult  to realize a set of numbers as the
  coefficients of a polynomial in $\gsubplus_2$, the determinant is
  the only restriction for a rhombus.

  \begin{lemma}\label{lem:coef-restrict-rhomb}
    Suppose that $\smalltwo{a}{b}{c}{d}$ is a matrix with positive
  entries. The following are equivalent:
  \begin{enumerate}
  \item $ad-bc>0$
  \item $\smalltwo{a}{b}{c}{d} =
    \smalltwo{\alpha_{01}}{\alpha_{11}}{\alpha_{00}}{\alpha_{10}}$
    for some $ \sum \alpha_{ij}\,x^i\,y^j\in\gsubplus_2$.
  \end{enumerate}
  \end{lemma}
  \begin{proof}
    We may assume that $c=1$. It suffices to consider the polynomial
    \begin{gather*}
      f(x,y)=(1+r_1\,x+s_1\,y)(1+r_2\,x+s_2\,y).\\
\intertext{We want to find positive $r_1,s_1,r_2,s_2$ such that}
\begin{pmatrix}
  a & b \\ 1& d
\end{pmatrix}=
\begin{pmatrix}
  s_1+s_2 & r_1\,s_2+r_2\,s_1 \\ 1 & r_1+r_2
\end{pmatrix}.
\intertext{Now it is easy to see that}
(0,ad) = \bigl\{r_1\,s_2+r_2\,s_1\,\mid\, s_1+s_2=a, r_1+r_2=d,r_1,r_2,s_1,s_2>0\bigr\}.
    \end{gather*}
Since $b\in(0,ad)$ by the determinant hypothesis we have found $f\in\gsubplus_2$.
  \end{proof}

There are inequalities for the coefficients of $f(x,y)\in\rupint{2}$ that
come from Proposition~\ref{prop:gen-laguerre}

\index{inequalities!for $\rupint{2}$}
\begin{lemma}
  If $f(x,y)=\sum_0^n f_i(x)y^i\in\rupint{2}$ and $0<k<n$ then
\[
\sum_0^k f_i(x) \, f_{k-i}(x)(-1)^{k-i} \ge0
\]
\end{lemma}
\begin{proof}
  For any $\alpha$ we apply Proposition~\ref{prop:gen-laguerre} to $f(\alpha,y)$.
\end{proof}

The inequality isn't strict. If $k=0$ then then the sum is $f_0(x)^2$
which takes on the value $0$.

\section{Solution curves}
\label{sec:p2-sol-curves}

 We have been discussing polynomials in $\gsubpos_2$ without actually
 looking at them. The geometric perspective gives important
 information about these polynomials.  This is not surprising, since
 the condition that distinguishes $\gsubpos_2$ from $\xsub_2$ concerns
 the structure of the homogeneous part, and the homogeneous part
 constrains the geometry of the graph.

\index{graph!of $\rupint{2}$}
If $f(x,y)$ is any polynomial, then the \emph{graph} of $f$ is defined 
to be
$$ G_f = \left\{ (x,y) \,\vert\, f(x,y)=0 \right\}$$
For example, if we choose a fifth degree polynomial
\begin{equation}
    \label{eqn:f8}
f =(x + y) (1 + x + 2 y) (1 + x + 3 y) (2 + x + 3 y) (2 + x + 5 y)
\end{equation}
    then the graph of $g=4f+3\,\frac{\partial f}{\partial x}$ (see
    Figure~\ref{fig:graph-8}) is surprising linear outside of the
    central region.  We will see later that $g$ is in $\gsubpos_2$.


\begin{figure}[htbp]
\rotatebox{90}{%
  \resizebox{5cm}{5cm}{
  \begin{pspicture}(0,0)(8,10)
    \pscurve(0,2.5)(2,6)(3.5,9.5)
    \pscurve(1,0)(3,5)(4.5,9.5)
    \pscurve(.5,0)(2.5,3)(5,9.5)
    \pscurve(2,0)(4,4.5)(8,9.5)
    \pscurve(3,0)(4,4.5)(5.5,9.5)
  \end{pspicture}
}}
    \caption{The graph of a polynomial in $\gsubpos_2$}
    \label{fig:graph-8}
\end{figure}

We can decompose $G_f$ into $n$ curves, where $n$ is the degree of
$f$. Define $r_i\colon\reals\longrightarrow\reals$ by setting $r_i(y)$ to
the $i$-th largest root of $f(x,y)=0$. Since $f$ satisfies
$x$-substitution each $r_i$ is well defined. The roots are continuous
functions of the coefficients, and so each $r_i$ is continuous. We
call $r_i(x)$ a \emph{solution curve} \index{solution curves} of $f$.
These are the analogs of the zeros of a polynomial in one variable. We
can write
$$
G_f = r_1(\reals) \cup \cdots \cup r_n(\reals)$$

\begin{lemma} \label{lem:p2-rays}
  If $f\in\gsubpos_2$ then asymptotically the graph $G_f$ is
  approximately a collection of infinite rays. The curve $r_i(\reals)$
  is asymptotically a ray with slope given by the $i$-th largest root
  of $f^H(x,1)$. Consequently, for $x,y$ sufficiently large $G_f$ lies
  in the union of the upper left quadrant and the lower right quadrant.
  \end{lemma}
  \begin{proof}
    If the  degree of $f$ is $n$ the polynomial $f^{H}(1,y)$ has $n$
    roots $\beta_1,\dots,\beta_n$ by  Lemma~\ref{lem:sub-fH2}.
    Consider
    $$f(x,y) = f^{H}(x,y) + \sum_{i+j<n}c_{ij}x^iy^{j}$$
    For $x$ large  and $y$ approximately $\beta_ix$ we see that
    $$f(x,y) \approx f^{H}(x,y) + O(x^{n-1})$$
    and hence for large
    $|x|$ there is a root $y$ of $f(x,y)=0$ where $y$ is close to
    $\beta_ix$.
  
    Thus, for $|x|$ large $f(x,y)$ has $n$ roots, approximately equal
    to $\beta_1x,\dots,\beta_nx$.  Since $f^H$ has all positive
    coefficients, the roots of $f^{H}(1,y)$ are negative, and so all
    the $\beta_i$ are negative.  Each $r_i(x)$ is in the upper left
    quadrant for $x$ large, and in the lower right quadrant for $x$
    negative and $|x|$ large.
  \end{proof}

  The polynomial $y^2-x^2+1$ shows that {$x$-substitution} alone does not
  imply {$y$-substitution}. Theorem~\ref{thm:sub-xy} remedies this situation by
  giving a condition that along with $x$-substitution  implies
  $y$-substitution.  This theorem is the best way to show that a
  polynomial is in $\gsubpos_2$.

\begin{theorem} \label{thm:sub-xy}
  Suppose that $f(x,y)$ is a polynomial that satisfies
  $x$-substitution, and $f^H$ satisfies the positivity condition. Then
  $f$ satisfies $y$-substitution, and is in $\rupint{2}$.
\end{theorem}

\begin{proof}
  Since each solution curve $r_i$ is continuous and and is
  asymptotically a ray that lies in the upper left quadrant or the
  lower right quadrant it follows that the graph of each $r_i$ meets
  every horizontal line, so the equation $r_i(x)=a$ has a solution for
  every $a$.  In particular, for any $y$ there are $n$ solutions to
  $f(x,y)=0$ since $f(a,r_i(a))=0$ for each $i$.  Since $f(x,y)$ has
  degree $n$ in $x$ this implies that $f(x,y)$ satisfies
  {$x$-substitution.}
\end{proof}

\begin{remark}
  It is important in the definition of substitution that the
  homogeneous part has all positive coefficients.  For instance,
  choose $g,h\in\allpoly$, and define $f(x,y) = g(x)h(y)$.  The
  homogeneous part of $f$ is a monomial $x^ny^m$, and does not satisfy
  positivity, and so is not in $\rupint{2}$ The graph of $f$ consists of
  $n$ horizontal lines, and $m$ vertical lines. We will see later that
  if $g,h\in\allpolypos$ then $f$ is in the closure of $\rupint{2}$.
\end{remark}

Here is an important way of constructing polynomials in $\rupint{2}$. We
will later show that these two conditions are equivalent.
\index{determinants!in $\gsubpos_2$}

\begin{lemma} \label{lem:det-aic}
  If $A$ and $C$ are symmetric  matrices, and $A$ is positive definite
  then 
  $$f(x,y) = \det(xA+yI+C)\in\rupint{2}.$$
\end{lemma}
\begin{proof}
This is Example~\ref{ex:det-aic}
\end{proof}

    \begin{lemma} \label{lem:sol-curves}
\index{solution curves} 
    If a solution curve of $f\in\rupint{2}$ is horizontal or vertical at
    a point $(a,b)$, then $(a,b)$ is the intersection of at least two
    solution curves.
  \end{lemma}
  \begin{proof}
    If the solution curve is horizontal at $(a,b)$ then
    $\dfrac{dy}{dx}(a,b)=0$. Since a solution curve is implicitly
    defined by $f(x,y)=0$, we can differentiate to find $f_x +
    \frac{dy}{dx}\,f_y=0$. This shows that $f_x(a,b)=0$. Consequently,
    the point $a$ is a double point of the function $f(x,b)$, and so
    $(a,b)$ lies on two solution curves. The vertical case is similar.
  \end{proof}

  \begin{lemma}\label{lem:sol-curves-dec}
    The solution curves of a polynomial in $\rupint{2}$ are
    non-increasing. They always go down and to the right.
  \end{lemma}
  \begin{proof}
    If we ever had a solution curve that was increasing, then it would
    have a local minimum since it is asymptotic to a line of negative
    slope. A horizontal line slightly above this minimum would
    intersect it in two distinct  points, contradicting the definition
    of solution curve. Similarly, considering the curve as a function
    of $y$ we see it can never go left.
  \end{proof}

\begin{figure}[htbp]
  \begin{center}
    \leavevmode
    \epsfig{file=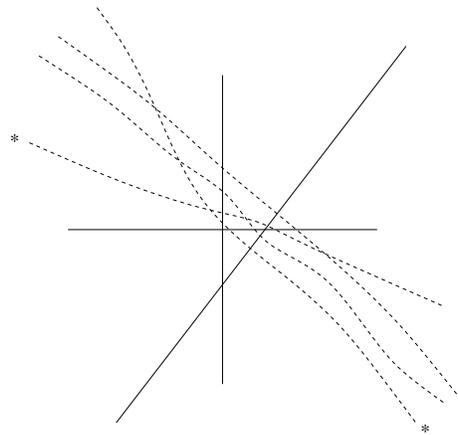,width=2.5in}
    \caption{The intersection of $f(x,y)=0$ and $y=ax+b$}
    \label{fig:fig-xy-line}
  \end{center}
\end{figure}

\begin{lemma} \label{lem:p2-xx}
  If $f$ is a polynomial in $\gsubpos_2$ then $f(x,x)\in\allpoly$.
  More generally, $f(ax+b,cx+d)$ is in $\allpoly$ for any $a,b,c,d$
  such that $a/c$ is either greater than the largest root of $f^H$ or
  less than the smallest root of $f^H$.

\end{lemma}
\begin{proof}
  Since each $r_i(\reals)$ is connected and asymptotically is a line with
  negative slope, any line of positive slope intersects $r_i(\reals)$.
  See Figure~\ref{fig:fig-xy-line}. This proves the first part.

Under the hypothesis on $a/c$ any solution curve either has asymptotes
whose slopes are both larger (or both smaller) than $a/c$. Any line
with slope $a/c$ will intersect such a curve.
\end{proof}

The hypothesis on $a/c$ are necessary. In Figure~\ref{fig:fig-xy-line}
it is easy to find lines that do not intersect the the solution curve
whose ends are marked ``*''.

 
The proof of the following theorem requires polynomials in four
variables, and is proved in Theorem~\ref{thm:sub-mx}. 

\begin{theorem}
If $f\in\gsubpos_2$  then for positive
$a,b,c,d$ we have
$$f(ax+by+u,cx+dy+v)\in\gsubpos_2$$
\end{theorem}

As an example, notice that if $f\in\gsubpos_2$ then $f(x,x+y)\in\gsubpos_2$. 

In addition to considering the intersection of the graph of $f$ with
lines we can consider the intersection of $f$ with the hyperbola
$xy=-1$. These intersections correspond to solutions of $f(x,-1/x)=0$.
Since this is not a polynomial, we multiply by $x^n$ where $n$ is the
degree of $f$. The degree of $x^nf(x,-1/x)$ is $2n$.

For example, if $f\in\allpolypos(n)$ then $F(x,y) = y^nf(x/y)$ is in
$\gsubpos_2$. Obviously $x^nF(x,-1/x) = (-1)^n f(-x^2)$. Since all roots
of $f$ are negative, each root of $f$ gives rise to $2$ roots of
$x^nF(x,-1/x)$  and so all $2n$ roots are accounted for.

The reverse of one variable requires a negative sign. 
If we use $1/x$ instead of $-1/x$ we get stable polynomials.
\seepage{lem:s-basic}

\index{reverse!in $\rupint{2}$}

\begin{figure}[htbp]
  \begin{center}
    \leavevmode
    \epsfig{file=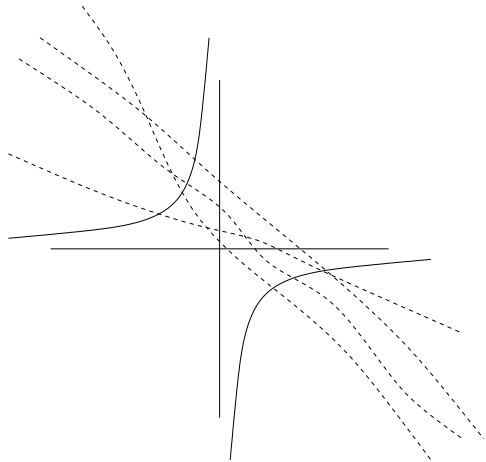,width=2.5in}
    \caption{The intersection of $f(x,y)=0$ and $xy=-1$}
    \label{fig:p2-xy}
  \end{center}
\end{figure}

\begin{lemma} \label{lem:p2-xy}
  If $f\in\gsubpos_2(n)$ then $x^nf(x,-1/x)\in\allpoly(2n)$.
\end{lemma}
\begin{proof}
  Since all asymptotes of $f$ have negative slope the graph of $xy=-1$
  meets the graph of $f$ in $n$ points in the upper left quadrant, and
  $n$ times in the lower right quadrant - see Figure~\ref{fig:p2-xy}.
  This gives $2n$ solutions to $f(x,-1/x)=0$, and since this is the
  degree of $x^nf(x,-1/x)$ the conclusion follows.
\end{proof}

\begin{remark}
  Harmonic functions are an important class of functions. A function
  is harmonic if it is the real part of an analytic function. As long
  as the degree is at least two then no polynomial in $\gsubpos_2$ is
  harmonic. A geometric explanation is that if $f$ is a harmonic
  polynomial then the real part of $f(x+iy)$ has asymptotes given by
  the $2n$ rays with angles $\pi(2k+1)/2n$, $k=1,\dots,2n$
  (\cite{gomez}). Thus, the real part has asymptotes whose slopes are
  positive and negative, and hence is not in $\gsubpos_2$.
\end{remark}

\begin{example}
    Recall (Example~\ref{ex:lax-example}) that polynomials of the form
    $f(x,y)=|I+xA+yB|$, where $A,B$ are $n$ by $n$ symmetric matrices,
    satisfy a strong line intersection property: every line through
    the origin meets the graph of $f$ in $n$ points. This fails for
    polynomials in $\rupint{2}$. For instance, if we define
$$ f(x,y)=
  \begin{vmatrix}
    4+x+2y & -3 \\ -3 & -2+x+y/2 
  \end{vmatrix} = -17 + 2\,x + x^2 - 2\,y + \frac{5\,x\,y}{2} + y^2
  $$
  then $f(x,y)\in\rupint{2}$, and it is clear from the graph of $f$
  (Figure~\ref{fig:lax-2}) that there is a range of lines through the
  origin that do not meet the graph of $f$.

\begin{figure}[htbp]
  \centering
  \includegraphics*[width=1.5in]{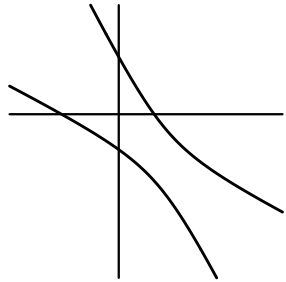}
  \caption{Satisfies substitution, fails to meet all lines through the
    origin}
  \label{fig:lax-2}
\end{figure}
  \end{example}

\section{Interlacing  in $\gsubpos_2$}
\label{sec:p2-interlacing-p2}

\index{\ zzlesslesseq@$\lesslesseq$!in $\gsubpos_2$} 
\index{interlacing!in $\gsubpos_2$}
\begin{definition}
  If $f,g\in\gsubpos_2$ then $f$ and $g$ \emph{interlace} iff
  $f+\alpha g$ is in $\pm\gsubpos_2$ for all $\alpha$. If $f$ and $g$
  have the same degree then it is possible that $(f+\alpha g)^H$ has
  negative coefficients. That's why we we require that $f+\alpha
  g\in\pm\gsubpos_2$, which we recall means that either $f+\alpha
  g\in\gsubpos_2$, or $-(f+\alpha g)\in\gsubpos_2$. If in addition the
  degree of $f$ is one more than the degree of $g$ then we say
  $f\lesslesseq g$. In this case $(f+g)^H = f^H$, and so the only
  condition we need to verify is substitution. This leads to an
  equivalent definition in terms of substitutions:
  \begin{quote}
    If $f,g\in\gsubpos_2$ then $f \lesslesseq g$ iff for all
    $a\in\reals$ we have $f(x,a)\lesslesseq g(x,a)$ ( or for all
    $a\in\reals$ we have     $f(a,y)\lesslesseq g(a,y)$).
  \end{quote}
\end{definition}

\begin{figure}[htbp]
  \begin{center}
    \leavevmode
    \epsfig{file=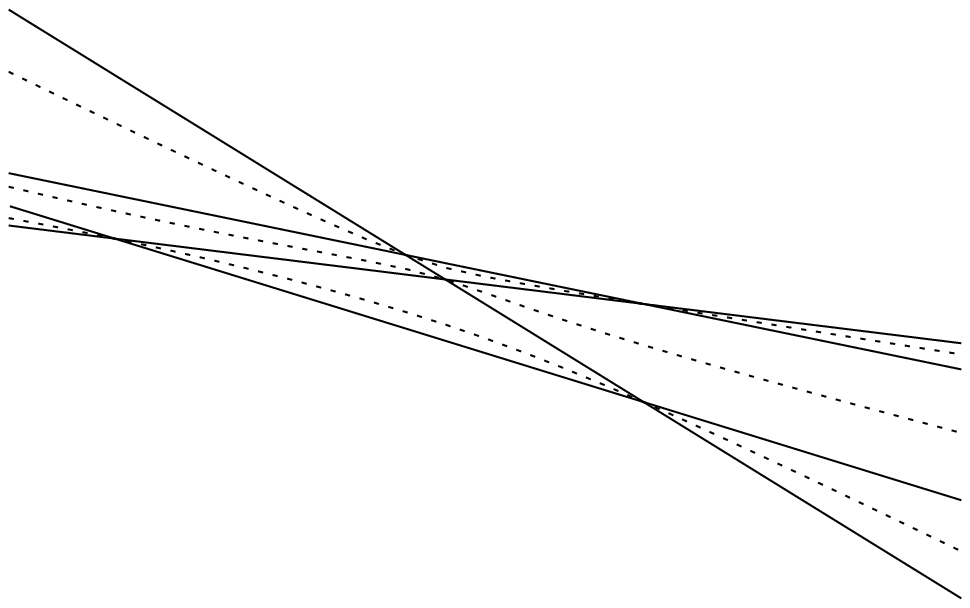,width=3in}
    \caption{The graphs of two interlacing polynomials in $\gsubpos_2$}
    \label{fig:fdxdy}
  \end{center}
\end{figure}

Figure~\ref{fig:fdxdy} shows  the graphs of the two interlacing polynomials
$$f=(x+y+1)(x+2y+3)(x+5y+4)(x+3y+2)$$ 
and $\frac{\partial f}{\partial x} + \frac{\partial f}{\partial y}. $
From this we see that there really is a geometric interpretation to
interlacing: if $f\lesslesseq g$ then the solution curves of $g$
\index{solution curves} 
interlace the solution curves of $f$.
\index{interlacing!geometric meaning in $\gsubpos_2$}

\index{interlacing!strict, in $\rupint{2}$}
We define strict interlacing $f \lessless g$ to mean that $f
\lesslesseq g$, and the graphs of $f$ and $g$ are disjoint. 
It is easy to see that this is equivalent to saying that
$f(x,\alpha)\lessless g(x,\alpha)$ for all $\alpha\in\reals$. 

The following lemma covers the fundamental properties of interlacing
in $\gsubpos_2$.

\begin{theorem}
  If $f,g,h\in\gsubpos_2$ and  $f\lesslesseq g,$  $f\lesslesseq h$,
  $\alpha,\beta$ positive then
\begin{itemize}
\item $\alpha g+\beta h\in\gsubpos_2$
\item $f \lesslesseq \alpha g+\beta h$
\end{itemize}
\end{theorem}
\begin{proof}
  Lemma~\ref{lem:sum-1} shows that $\alpha g+\beta h$ satisfies
  substitution. Since $(\alpha g+\beta h)^H = \alpha g^H + \beta h^H$
  it follows that $\alpha g+\beta h$ satisfies the positivity condition.
\end{proof}

\begin{theorem} \label{thm:gsub-diff2}
  $\gsubpos_2$ is closed under differentiation. If $f\in\gsubpos_2$
  then $f\lesslesseq \frac{\partial f}{\partial x_i} $ for $1\le i \le
  2$.  
\end{theorem}

\begin{proof}
Choose $\alpha\in\reals$, and consider $g = f + \alpha \frac{\partial
  f}{\partial x_1}$. Since $g^H = f^H$, we will show that
$g\in\rupint{2}$ by showing $g$ satisfies substitution. If we choose
$\beta\in\reals$ then
$$
g(x,\beta) = f(x,\beta) + \alpha \frac{\partial f}{\partial
  x_1}(x,\beta)$$
Since derivatives in one variable interlace, it follows
from     Theorem~\ref{thm:sub-xy} that $g\in\rupint{2}$, and and that
$f\lesslesseq \frac{\partial f}{\partial x_1}$. 
\end{proof}

One of the surprising things about this definition is that interlacing
linear polynomials are highly constrained. In $\allpoly$, any two
polynomials of degree one interlace. In $\gsubpos_2$, the only way
that degree one polynomials interlace is for them to be essentially
one-dimensional. Suppose that $ax+by+c \greateqeq s x+ t y + u$ where
$a,b,s,t$ are positive. The requirement that
$(ax+by+c)+\alpha(sx+ty+u)$ lies in $\gsubpos_2$, implies that $sx+ty$
is a multiple of $ax+bx$. If we define $f(x) = x+c$ and $g(x) =
x+t/\alpha$ then we can express these two interlacing polynomials in
terms of interlacing one variable polynomials.  Namely, $f(ax+by) =
ax+by+c$, $\alpha g(ax+by) =sx+ty+u$, and $f\greateqeq g$. Thus there
are no intrinsically two-dimensional interlacing linear polynomials.
Geometrically, all this means is that the graph of interlacing linear
polynomials must be parallel lines.

A consequence of this observation is that there appears to be no
simple creation of \emph{mutually interlacing polynomials}.
\index{mutually interlacing} \index{polynomials! mutually interlacing}
We have seen how to create mutually interlacing polynomials using
linear combinations of the products of $n-1$ factors from a polynomial
with $n$ factors. The fact that this set of polynomials is mutually
interlacing follows from the fact that any $n$ degree one polynomials
in $\allpoly$ are mutually interlacing.

Interlacing polynomials can be decomposed.  This fact follows from
from the assumption that the polynomials $f+\alpha g$ all lie in
$\pm\gsubpos_2$, and hence the coefficients of the homogeneous part all
have the same sign.

\begin{lemma}
  Suppose $f,g\in\gsubpos_2$, $f$ and $g$ have the same total degree, and
  $f,g$ interlace.  then 
  \begin{enumerate}
  \item $f^H$ and $g^H$ are scalar multiples of each
  other. 
  \item There is an $r\in\gsubpos_2$ and positive
  $\alpha$ such that 
  \begin{enumerate}
  \item $f \lesslesseq r$
   \item $g = \alpha f \pm  r$ 
  \end{enumerate}
  \end{enumerate}
\end{lemma}
\begin{proof}
First of all, we note that $(f+ \alpha g)^H = f^H +
\alpha g^H$, so $f^H$ and $g^H$ interlace.  Since $f+\alpha
g\in\gsubpos_2$ for all $\alpha$, it follows that $(f+\alpha
g)^H\in\allpolypos$ for all $\alpha$.  From
Lemma~\ref{lem:pos-neg-2} there is a $\gamma$ such that $f+\gamma g$ is not in
$\allpolypos$, so the second case of Lemma~\ref{lem:pos-neg-2} is not
possible. Thus,  $g^H$ is a constant multiple of $f^H$.

  By the previous lemma we can choose $\alpha$ so that $g^H = \alpha
  f^H$.  We define $r = \pm(f - \alpha g)$ where the sign is chosen to
  make the coefficients of $r^H$ positive. All linear combinations of
  $f$ and $g$ are in $\gsubpos_2$, so $r$ is  in $\gsubpos_2$.
  Moreover, the total degree of $r$ is less than that of $f$ since we
  removed all the highest degree terms. Also, $f$ and $r$ interlace
  since their linear span is the same as the linear span of $f$ and
  $g$.
\end{proof}

We use this definition to define $\greateqeq$ for $\gsubpos_2$. Say
that $f\greateqeq g$ iff there is an $r\in\gsubpos_2$ such that $g =
\alpha f+r$ for some positive $\alpha$, and $f \lesslesseq r$.

\section{Linear transformations on $\rupint{2}$}
\label{sec:p2-lt}

We study linear transformations on $\rupint{2}$. 
Just as in one variable, linear transformations preserving
$\gsubpos_2$ preserve interlacing.  The proof is immediate from the
definition of interlacing.

\begin{theorem} \label{thm:p2-LT}
  Suppose that $T\colon{}\allpoly\longrightarrow\gsubpos_2$ and
  $S\colon{}\gsubpos_2\longrightarrow\gsubpos_2$ are linear transformations.  If
  $f,g\in\allpoly$ interlace then $Tf,Tg$ interlace.  If
  $f,g\in\gsubpos_2$ interlace then $Sf,Sg$ interlace.
\end{theorem}

\begin{cor} \label{cor:p2-deriv}
  If $f\lesslesseq g$ are in $\gsubpos_2$, and $a,b$ are positive then
$$ a \frac{\partial f}{\partial x} + b \frac{\partial f}{\partial y} 
\lesslesseq
a \frac{\partial g}{\partial x} + b \frac{\partial g}{\partial y} 
$$
If $f\lessless g$ then $\frac{\partial f}{\partial x_i}
  \lessless \frac{\partial g}{\partial x_i} $
\end{cor}
\begin{proof}
  The map $f\mapsto a \frac{\partial f}{\partial x} + b \frac{\partial
    f}{\partial y} $ maps $\gsubpos_2$ to itself, so it preserves
  interlacing.  The second part follows form the one variable result.

\end{proof}

  \begin{lemma} \label{lem:p2-xx-int}
    If $f(x,y) \lesslesseq g(x,y)$ in $\rupint{2}$ then
    $f(x,x)\lesslesseq g(x,x)$.
  \end{lemma}
  \begin{proof}
    The map $y\mapsto x$ is a linear transformation from $\rupint{2}$ to
    $\allpoly$, and so it preserves interlacing.
  \end{proof}

\index{induced transformation}
\index{transformation!induced}

A linear transformation $\allpoly\longrightarrow\allpoly$ sometimes induces a
linear transformation $\gsubpos_2\longrightarrow\gsubpos_2$. Analyzing
the geometry of the graph of the image of a polynomial under this
induced transformation yields information about new transformations on
$\allpoly$ or $\allpolypos$.

\begin{theorem}  \label{thm:p2p2}
  Suppose that\/ $T\colon{}\allpoly\longrightarrow \allpoly$ is a linear
  transformation that preserves degree, and maps polynomials with
  positive leading coefficients to polynomials with positive leading
  coefficients.  The induced linear transformation
  $$T_\ast(x^iy^j) = T(x^i)\,y^j$$
  defines a linear transformation from $\gsubpos_2$ to itself.
\end{theorem}

\begin{proof}
  Suppose $f = \sum f_i(x)y^i$ where $f_i$ has degree $n-i$. Since $T$
  preserves degree, the degree of $T(f_i)$ is $n-i$, and so $(Tf)^H$ is
  a sum of terms $x^{n-i}($ leading coefficient of $Tf_i)$. Since $T$
  preserves the sign of the leading coefficient, $(Tf)^H$ has all
  positive terms.

By Theorem~\ref{thm:sub-xy} it suffices to show that
  $T_\ast(f)$ satisfies $x$-substitution.  If we choose
  $b\in\reals$ then
  $$
  (T_\ast\,f)(x,b) = T(f(x,b))$$
  Since $f(x,b)\in\allpoly$ we
  know $T(f(x,b))$ is in $\allpoly$, and so $T_\ast$ satisfies $x$-substitution.
\end{proof}

We will revisit induced transformations in
Chapter~\ref{cha:topology}.\ref{sec:top-induced-trans}.

\begin{remark}\label{rem:no-induced}
  It is important to see that a linear transformation
  $T\colon{}\allpoly\longrightarrow\allpoly$ does not in general induce a
  linear transformation $\rupint{2}\longrightarrow\gsub_2$. The
  assumptions of degree and positivity are essential. The proof shows
  that substitution will always be met. For instance, consider 
$$ T\colon{} g \mapsto g(\diffd)x.$$
This satisfies $T(1)=x$, $T(x)=1$ and $T(x^k)=0$ for $k>1$. Since
$$T(x+y)^2 = T(x^2+2xy+y^2) = 2y + xy^2$$
we see that $T(x+y)^2\not\in\rupint{2}$ since the homogeneous part is $xy^2$.
See \chapsec{topology}{top-induced-trans}.
\index{linear transformation!doesn't induce on $\rupint{2}$}
\end{remark}

\begin{cor} \label{cor:hihj} \index{Hermite polynomials!transformations}
  If $T(x^iy^j)=H_i(x)H_j(y)$ then $T\colon{}\gsubpos_2\longrightarrow\gsubpos_2$.
\end{cor}
\begin{proof}
  The maps $x^i\mapsto H_i(x)$ and $y^j\mapsto H_j(y)$ map $\allpoly$
  to itself by Corollary~\ref{cor:hermite}.
\end{proof}

\index{Hermite polynomials!identities}
We can make use of a special function
identity to verify Corollary~\ref{cor:hihj} in a special case.  We
have 
\begin{align*}
   T \,(x+y)^n &= \sum_{i=0}^n
   \binom{n}{k}T(x^ky^{n-k}) \\
&= \sum_{i=0}^n \binom{n}{k}
   H_{k}(x)H_{n-k}(y) \\
& = 2^{n/2} H_n\left(\dfrac{x+y}{\sqrt{2}}\right)
\end{align*}
and the last polynomial is in $\gsubpos_2$ since $H_n$ is in $\allpoly$. 

Theorem~\ref{thm:p2p2} is not true\footnote{ For example, we saw in
  \chapsec{operators}{rising-factorial} that \index{rising factorial}
  $T(x^n) = \rising{x}{n}$ is such a map. If we set $ f =
  (1+x+y)(1+2x+y)(2+3x+y)$ then the coefficients of $T_\ast f$
  interlace but $T_\ast f$ does not satisfy {$x$-substitution} since
  $(T_\ast f)(3,y)$ has imaginary roots. } if we only have that
$T\colon{}\allpolypos \longrightarrow \allpoly$.

\begin{lemma}\label{lem:simple-p2}
  Consider the maps $T\colon{}f(x)\mapsto f(x+y)$ and $S\colon{}f(x)\mapsto F(x,y)$,
  where $F(x,y)$ is the homogeneous polynomial corresponding to $f$.
  These are linear transformations
  $T\colon{}\allpoly\longrightarrow\gsubpos_2$ and
  $S\colon{}\allpolypm(n)\longrightarrow\gsubpos_2(n)$.
\end{lemma}
\begin{proof}
  The map $T$ is linear and preserves degree, and $S$ is linear if we
  restrict to polynomials of the same degree. 
\end{proof}

\begin{cor} \label{cor:exp-xy}
  The linear transformation  $T(x^iy^j) = \frac{x^iy^j}{i!j!}$ maps
  $\gsubpos_2$ to $\gsubpos_2$.
\end{cor}
\begin{proof}
  Apply Theorem~\ref{thm:p2p2} to the exponential map in each variable.
\end{proof}

There are a few instances where we can describe all the coefficients
of $y^i$ in a polynomial belonging to $\gsubpos_2$. Here are three
examples - in each case the left sides are of the form $f(x+y)$ so the
right hand sides are in $\gsubpos_2$. $H_n$ is the Hermite polynomial,
and the identity is a consequence of the Taylor series and the
identity $H_n^\prime=2nH_{n-1}$.  \index{Hermite polynomials}

\begin{align*}
  (x+y+1)^n & = \sum_{i=0}^n \binom{n}{k} (x+1)^k y^{n-k} \\
  H_n(x+y) &= \sum 2^k \binom{n}{k} H_{n-k} y^k
\end{align*}

\index{Taylor series}
If we  expand $f(x+y)$ in its Taylor series then for any
$f\in\allpolypos$
\begin{equation}
  \label{eqn:fxyp2}
  f(x+y)= \sum_i {f^{(i)}(x)} \,\frac{y^i}{i!} \in\gsubpos_2
\end{equation}

A linear transformation from $\allpoly$ to $\allpoly$ preserving
$\allpolypos$ determines a linear transformation from $\allpoly$ to
$\gsubpos_2$.

\begin{lemma} \label{lem:p2-homog}
  Suppose $T\colon{}\allpoly\longrightarrow\allpoly$ is a linear
  transformation that preserves degree, and maps $\allpolypos$ to
  itself. If $T(x^i)=f_i$ and we define $S(x^i) = y^if_i(x/y)$,
   then $S\colon{}\allpolypos\longrightarrow \gsubpos_2$.
\end{lemma}

\begin{proof}
  Suppose that $g=\sum_{i=0}^na_ix^i$ is in $\allpolypos$.  The
  homogeneous part of $S(g)$ is $a_nS(x^n)=a_ny^nf_n(x/y)$ whose
  homogenization is in $\allpolypos$ since $T$ maps $\allpolypos$ to
  $\allpolypos$. To verify substitution, choose $\alpha$ and consider
  $$
  S(g)(x,\alpha) = \sum a_iS(x^i)(x,\alpha) = \sum
  a_i\alpha^if_i(x/\alpha)$$
  This last polynomial is in $\allpolypos$, as the diagram shows

    \centerline{ \xymatrix{
 \allpolypos \ar@{->}[rrr]^{x\mapsto \alpha{ x}} \ar@{..>}[d]_{S(g)(x,\alpha)}
     &&& \allpolypos \ar@{->}[d]^T\\
 {\allpolypos} &&&
        \allpolypos \ar@{->}[lll]^{x\mapsto x/\alpha}
      }}

\end{proof}

\begin{example}
  If we choose the affine map $T(f) = f(x+1)$ then $S(x^i)=(x+y)^i$,
  and $S(f)=f(x+y)$.
\end{example}

\begin{lemma}\label{p2-pk}
  If $T\colon{}\allpoly\longrightarrow\allpoly$ and $T(x^n)=p_n$, then for
  any $\alpha$ the linear transformation $S(x^k) =
  x^{n-k}\,p_k(\alpha)$ maps $\allpolypm(n)\longrightarrow\allpoly(n)$.
\end{lemma}
\begin{proof}
  The proof follows from the diagram

\centerline{\xymatrix{
\allpolypm(n) \ar@{.>}[d]_{S} \ar@{->}[rrr]^{homogenization} &&& {\rupint{2}}
\ar@{->}[d]^{T_\ast}  \\
\allpoly(n) &&&  {\rupint{2}} \ar@{->}[lll]_{x\mapsto\alpha} && 
}}

\end{proof}

\section{Applications to linear transformations on $\allpoly$}
\label{sec:p2-lin-p}

We now apply results about $\rupint{2}$ to get information about linear
transformations on $\allpoly$. 

\begin{cor} \label{cor:p2-st}
  Assume that the linear transformations
  $T,S\colon{}\allpoly\longrightarrow\allpoly$ map $x^n$ to a polynomial of
  degree $n$ with positive leading coefficient, and maps $\allpolypos$
  to itself. . Define a transformation $V$ on $\allpoly(n)$ by $V(x^i)
  = T(x^i)\,S(x^{n-i})$.  Then $V$ maps $\allpolypos(n)$ to itself.
\end{cor}
\begin{proof}
  $V$ is the composition
 

    \centerline{ \xymatrix{
 \allpolypos(n) \ar@{->}[rrr]^{homogenize} \ar@{..>}[d]_V
        &&& \ar@{->}[d]^{T_\ast\times S_\ast} \rupint{2}(n)\\
 {\allpoly(n)} 
        \ar@{<-}[rrr]^{(x,y)\mapsto(x,x)}
        &&& {\rupint{2}(n)}}
      }

where $T_\ast$ and $S_\ast$ are defined in Theorem~\ref{thm:p2p2}.
The conclusion follows form Lemma~\ref{lem:p2-xx}. 
  \end{proof}

  \begin{cor} \label{cor:hermite-xy}
    \index{Hermite polynomials!transformations} If
    $H_i(x)$ is the Hermite polynomial then the transformation
    $x^i\mapsto H_i(x)\, H_{n-i}(x)$ maps $\allpolypos(n)$ to $\allpoly(n)$.
    If $L_i(x)$ is the Laguerre polynomial then the transformation
    $x^i\mapsto L_i(x)\, L_{n-i}(x)$ maps $\allpolypos(n)$ to $\allpolypos(n)$.
    Also, the transformation $x^i\mapsto L_i(x)\, H_{n-i}(x)$ maps
    $\allpolypos(n)$ to $\allpoly(n)$.
  \end{cor}

  \begin{cor}  \label{cor:p2-n-iT}
    Assume that the linear transformation
    $T\colon{}\allpoly\longrightarrow\allpoly$ map $x^n$ to a polynomial of
    degree $n$ with positive leading coefficient. The linear
    transformation $x^i\mapsto x^{n-i}T(x^i)$ maps $\allpolypos(n)$ to
  itself.
  \end{cor}
  \begin{proof}
 Take $S$ to be the identity in Corollary~\ref{cor:p2-st}.   
  \end{proof}

The next result  follows from the corollary, but we will refine
this result in Lemma~\ref{lem:hermite-hxn}.

  \begin{lemma} \label{cor:herm2}
     The transformation $x^i\mapsto H_i(x)x^{n-i}$ maps $\allpolypos(n)$
    to $\allpoly(n)$.
  \end{lemma}
 
  \begin{lemma}
    Suppose that $f\in\allpoly$, and
    $T\colon{}\allpoly\longrightarrow\allpoly$. Then
$$ T_\ast(f(x+y))(x,x)\in\allpoly$$
  \end{lemma}
  \begin{proof}
    We know that $f(x+y)\in\gsubpos_2$, and so
    $T_\ast(f(x+y))\in\gsubpos_2$. Apply Lemma~\ref{lem:p2-xx}.
  \end{proof}

\begin{cor} \label{cor:t1/z} \index{T@$T_{1/z}$}
  Suppose $T\colon{}\allpoly\longrightarrow\allpoly$ is a linear
  transformation that preserves degree and maps $\allpolypos$ to
  itself. The \index{\Mobius\ transformation}\Mobius\ transformation $T_{1/z}$ satisfies
  $T_{1/z}:\allpolypos\longrightarrow\allpolyneg$.
\end{cor}
\begin{proof}
  Let $T(x^i)=f_i(x)$, and recall that $T_{1/z}(x^i) = \rev{f_i}$.
  Let $S(x^i) = y^if_i(x/y)$ as in Lemma~\ref{lem:p2-homog}.  Notice that
  $S(x^n)(1,y) = y^nf_n(1/y) = T_{1/z}(y)$ and so
  $S(f)(1,y)=T_{1/z}(f)$.  Now $S$ maps $\allpolypos$ to $\gsubpos_2$
  by Lemma~\ref{lem:p2-homog}, and hence
  $T_{1/z}(f)=S(f)(1,y)\in\allpoly$.

  See \mypage{lem:p-to-p2} for a different proof.
\end{proof}

\begin{cor} \label{cor:herm-1}
  The linear transformation  $x^i\mapsto \rev{L_i}(-x)$ maps
  $\allpolypos\longrightarrow\allpolypos$. 
\end{cor}
\begin{proof}
Apply Corollary~\ref{cor:t1/z} to  Lemma~\ref{lem:laguerre-2}.
\end{proof}

The reason that we can substitute $x=y$ in a polynomial in
$\gsubpos_2$ is that the line $x=y$ must intersect every solution
curve. The line $x+y=1$ does not necessarily meet every solution
curve, but an easily met assumption about the homogeneous part will
allow us to substitute $1-x$ for $y$.

\begin{lemma}
  Suppose $T\colon{}\allpoly\longrightarrow\allpoly$ maps $\allpolypos(n)$ to
  itself, and define $S(x^i) = (1-x)^{n-i}\,T(x^i)$.  Assume that the
  leading coefficient of $T(x^i)$ is $c_i$, and write $f(x) = \sum a_i
  x^i$. If $f\in\allpolypos$ and either condition below is met then
  $S(f)\in\allpoly$.
\begin{itemize}
\item all roots of $\sum a_ic_ix^i$ are greater than $-1$.
\item all roots of $\sum a_ic_ix^i$ are less than $-1$.
\end{itemize}
\end{lemma}
\begin{proof}

$S$ is determined by the composition

    \centerline{ \xymatrix{
 \allpolypos(n) \ar@{->}[rrr]^{homogenize} \ar@{..>}[d]_S
        &&& \ar@{->}[d]^{T_\ast} \gsubpos_2(n)\\
 {\allpoly(n)} 
        \ar@{<-}[rrr]^{(x,y)\mapsto(x,1-x)}
        &&& {\gsubpos_2(n)}}
      }

If we choose $f\in\allpolypos$ as in Corollary~\ref{cor:p2-st} we know
$T_\ast(F)\in\gsubpos_2$, where $F$ is $f$ homogenized. The asymptotes
of $T_\ast(F)$ have slopes that are the roots of the homogeneous part
of $T_\ast(F)$, which is $\sum a_ic_ix^i$. If either condition is met,
then all solution curves meet the line $x+y=1$. 
\end{proof}

  The linear transformation $T\colon{}x^k\mapsto H_k(x)x^{n-k}$ acting on
  $\allpoly(n)$ becomes very simple if we make a transformation. We
  need the identity
  \begin{gather}
    \label{eqn:hermite-id}
    T(x-2+\alpha)^n = H_n(\frac{\alpha x}{2})\\
\intertext{If we substitute $\alpha=0$ then}
 T\,(x-2)^n = H_n(0) = 
\begin{cases}
  0 & \text{ if $n$ is odd}\\
  (-1)^{n/2}(n-1)! &\text{ if $n$ is even}
\end{cases}\notag
\end{gather}

This degenerate behavior at $x=2$ is reflected in the following lemma.

  \begin{lemma}\label{lem:hermite-hxn}
    The linear transformation $T\colon{}x^k\mapsto H_k(x)x^{n-k}$ acting on
    $\allpoly(n)$ maps
    $\allpolyint{(2,\infty)}\cup\allpolyint{(-\infty,2)}$ to
    $\allpoly$. 
  \end{lemma}
  \begin{proof}
    Consider the transformation $S(f) = T(f(x-2))$ acting on
    $\allpoly(n)$, and let $W(x^k) = H_k(x)$. The key observation
    (from \eqref{eqn:hermite-id}) is that
\begin{equation}\label{eqn:hxn-1}
 S(x^k) = T(x-2)^k = x^{n-k} H_k(0)
\end{equation}
Therefore, if $g=\sum a_ix^i$, and the
homogenization $G$ is $g$, then
$$ W_\ast(G) = \sum a_i\,y^{n-k}H_k(x)$$
so we get the basic relation
$$ S(g)(y) = W_\ast(G)\big\vert_{x=0}$$
If $g\in\allpolypm$, then $G$ is in $\rupint{2}$, and
$W_\ast(G)\in\rupint{2}$, so $S(g)\in\allpoly$. Therefore,
$T(g)\in\allpoly$. 

  \end{proof}

\begin{remark}
  We can use \eqref{eqn:hermite-id} to get precise information about
  the possible range of $S$ in \eqref{eqn:hxn-1}.  Suppose that
  $\diffi_n$ is the smallest interval containing all the roots of
  $H_n$. The roots of $H_n(\frac{\alpha x}{2})$ lie in the interval
  ${\frac{2}{\alpha}\diffi_n}$. If we knew that $T$ preserves
  interlacing then from Lemma~\ref{lem:find-range} we conclude from
  \eqref{eqn:hermite-id} that if $0<\alpha<1$ then $$
  T\colon{}
  \allpolyint{(2+\alpha,\infty)}(n) \longrightarrow
  \allpolyint{\frac{2}{\alpha}\diffi_n} $$
\end{remark}

\section{Properties of   $\rupint{2}$}
\label{sec:sub-derivatives2} 

The main result of this section is that the coefficients of a
polynomial in $\gsubpos_2$ are in $\allpoly$.  This  allows us
to prove that a polynomial is in $\allpoly$ by identifying it as a
coefficient of a polynomial in $\gsubpos_2$.

\begin{cor} \label{cor:p2-coef}
  If $f\in\gsubpos_2$ then all coefficients $f_i$ in \eqref{eqn:p3-1}
  and all coefficients $f^i$ in \eqref{eqn:p3-2} are in $\allpoly$.
  Moreover, $f_{i}\lesslesseq f_{i+1}$ and $f^{i}\lesslesseq f^{i+1}$ for
  $0\le i <n$.
\end{cor}
\begin{proof}
  If $f\in\gsubpos_2$ then
$ \left( \left(\frac{\partial}{\partial y}\right)^i\,f\right)\,(x,0)
= i!f_i(x)$
and hence $f_i\in\allpoly$ by Theorem~\ref{thm:gsub-diff2} and substitution.

Since $f +\alpha\frac{\partial f}{\partial y}$ is in $\gsubpos_2$ for any
$\alpha$ the coefficient of $y^i$ is in $\allpoly$.  This coefficient
is $f_i+\alpha f_{i+1}$ and so we conclude that $f_{i+1}$ and $f_i$
interlace. Since the degree of $f_i$ is greater than the degree of
$f_{i+1}$ we find $f_i \lesslesseq f_{i+1}$.
\end{proof}

  There is a simple condition that guarantees that all the
  coefficients strictly interlace.

  \begin{lemma}\label{lem:p2-strict}
    If $f(x,y)\in\rupint{2}$ and $f(x,0)$ has all distinct
    roots, then
    \begin{enumerate}
    \item All $f_i$ have distinct roots.
    \item $f_i \lessless f_{i+1}$
    \end{enumerate}
  \end{lemma}
  \begin{proof}
    Assume that $f_0(r)=f_1(r)=0$. Let $g(x,y) = f(x+r,y)$, and write
    $g = \sum g_i(x)y^i = \sum g^j(y)x^j = \sum a_{ij}x^iy^j$. Since
    $g_0(0)=g_1(0)=0$, we find $a_{0,0}=a_{0,1}=0$. Thus, $g^0$ is
    divisible by $x^2$. Since $g^1$ interlaces $g_0$, $g^1$ is
    divisible by $x$. This implies that $a_{1,0}=0$, and hence $g_0$
    has a double root at zero, but this contradicts the hypothesis
    that $f_0$ has no repeated roots since $g_0(x) = f_0(x+r)$. Since
    $f_0\lessless f_1$, it follows that $f_1$ has all distinct roots.
    Continuing, we see all $f_i$ have distinct roots.
  \end{proof}

The polynomial $f(x,y)=(1+x+y)(1+2x+y)$ shows that $f(x,0)$ can have all
distinct roots, yet $f(0,y)$ can have repeated roots.

\begin{cor} \label{cor:fct-interlace}
  If $f\in\gsubpos_2$ and $f_0\lessless f_1$ then $f_i\lessless f_{i+1}$
  for $0\le i < n$.
\end{cor}

\begin{lemma} \label{lem:p4-only-way}
  If $f$ given in \eqref{eqn:p3-1} satisfies \emph{y}-substitution then
  the sequence of coefficients is  log concave:
  $\smalltwodet{f_i}{f_{i+1}}{f_{i+1}}{f_{i+2}}\le0$ for $0\le i \le
  n-2$. If $f_0$ and $f_1$ have no roots in common then 
$\smalltwodet{f_i}{f_{i+1}}{f_{i+1}}{f_{i+2}}<0$
\end{lemma}
\begin{proof}
  \index{Newton's inequalities} 
  
  Since $f$ has all real roots for any fixed $x$, we can apply
  Newton's inequalities (Corollary~\ref{cor:newton}) to find that
  $f_j^2(x) > f_{j-1}(x)f_{j+1}(x)$ unless $f_0(r) = f_1(r) = 0$ for
  some$r$ in which case we only have $f_j^2(x) \ge
  f_{j-1}(x)f_{j+1}(x)$.
\end{proof}

\begin{cor} \label{cor:product-1}
 Set 
 $$
 f(x,y) = \prod_{i=1}^n(x+b_i+c_iy) = \sum_{i=0}^n f_j(x)y^j$$
 Assume that the $c_i$ are positive for $1\le i \le n$.  Then the
 coefficients $f_j$ are polynomials of degree $n-j$, and
$$ f_0 \lesslesseq f_1 \lesslesseq f_2 \lesslesseq \dots $$
If all $b_i$ are  distinct then the interlacings are all
strict.
\end{cor}

\begin{proof}
Lemma~\ref{lem:det-aic} shows that $f\in\gsubpos_2$. By
Corollary~\ref{cor:fct-interlace}
it remains to show that $f_0 \lessless f_1$. If we expand
$f$ we find (see Lemma~\ref{lem:product}) 
  \begin{align*}
     f_1(x) &= \sum_{j=1}^n  {c_j} \,\frac{f_0(x)}{x + b_j}
  \end{align*}
  Since all coefficients $c_j$ are positive we can apply
  Lemma~\ref{lem:sign-quant} to conclude that $f_0 \lessless f_1$ since $f_0$
  has all distinct roots.
\end{proof}

\begin{cor}
  If $f = \sum f_i(x)y^i$ and $g = \sum g_i(x)y^i$ are both in
  $\gsubpos_2(n)$ then
$$ \sum_i f_i(x)g_{n-i}(x) \in\allpoly$$
\end{cor}
\begin{proof}
  The polynomial in question is the coefficient of $y^n$ in the
  product $fg$.
\end{proof}

Note that if we take $f=\sum a_i x^{n-i}y^i$ and $g=\sum g_i
x^{n-i}y^i$ where $\sum a_ix^i$ and $\sum b_ix^i$ are both in
$\allpolypos$ then we conclude that $\sum
a_{n-i}b_ix^{2i}\in\allpoly$. This is a simple modification of the
Hadamard product. \index{Hadamard product}

\begin{remark} \label{rem:rolle-2}
  \index{Taylor series} \index{Rolle's Theorem}
  Since the first two terms of the Taylor series \eqref{eqn:fxyp2} of
  $f(x+y)$ are $f$ and $yf^\prime$ it follows that $f\lesslesseq
  f^\prime$. Thus, we can prove Rolle's theorem (Theorem~\ref{thm:rolle}) without
  analyzing the behavior of the graph of $f$ and $f^\prime$.
\end{remark}

\section{The analog of $\allpolypos$}
\label{sec:p2plus}

A polynomial in $\allpoly$ with positive leading coefficient is also
in $\allpolypos$ if and only all of its coefficients are positive. We
define $\gsubplus_2$ in the same way:

\index{\ Pgsubplus2@$\gsubplus_2$}

$$
\gsubplus_2 = \left\{ f\in\gsubpos_2\mid \text{ all coefficients of $f$
    are positive}\right\} 
$$

We start with a  useful alternative criterion.

\begin{lemma}\label{lem:in-p2plus}
  Suppose $f\in\rupint{2}$ and write $f=\sum f_i(x)y^i$. If
  $f_0\in\allpolypos$ then $f\in\gsubplus_2$.
\end{lemma}
\begin{proof}
  Since $f_i(x)\lesslesseq f_{i+1}(x)$ and since
  $f_0(x)\in\allpolypos$ it follows that all $f_i$ are in
  $\allpolypos$. Now the leading coefficients of all $f_i$'s are terms
  of $f^H$ and are therefore positive. Consequently, all $f_i$'s have
  all positive coefficients, and thus
  $f(x,y)\in\gsubplus_2(n)$.
\end{proof}

We discuss the generalization of $\allpolypos$ to more than two
variables in \chapsec{pd}{gsubpm}, and more properties of
$\gsubplus_2$ can be found there. 

Our next result is that a simple translation takes a polynomial in
$\rupint{2}$ to one in $\gsubplus_2$. 

\begin{lemma}
  If $f(x,y)\in\rupint{2}(n)$, then there is an $\alpha$ such that\\
  $f(x+\alpha,y)\in\gsubplus_2(n)$. 
\end{lemma}
\begin{proof}
  If $f(x,y) = \sum f_i(x)y^i$, then we can choose $\alpha$ such that
  ${f_0(x+\alpha)}$ is in $\allpolypos$. It follows from
  Lemma~\ref{lem:in-p2plus} that $f(x+\alpha,y)\in\gsubplus_2$.
\end{proof}

Interlacing of polynomials in $\gsubplus_2$ stays in $\gsubplus_2$.

\begin{lemma}
  If $f\in\gsubplus_2$ and $f\lesslesseq g$ then $g\in\gsubplus_2$
\end{lemma}
\begin{proof}
  If we write $f = \sum f_i(x)y^i$ and $g=\sum g_i(x)y^i$ then we know
  that $f_0 \lesslesseq g_0$. Thus $g_0\in\allpolypos$ and the
  conclusion follows from Lemma~\ref{lem:in-p2plus}.
\end{proof}

If $f$ is in $\rupint{2}$ then $f(x,\alpha x)\in\allpoly$ for positive
$\alpha$, but this can be false for negative $\alpha$. However, it is true
for $\gsubplus_2$.

\begin{prop}\label{prop:p2-xax}
  If $f(x,y)\in\gsubplus_2(n)$ then $f(x,\alpha x)\in\allpoly$ for all
  $\alpha$. 
\end{prop}
\begin{proof}
  If suffices to consider $\alpha<0$. We will show that the line
  $y=\alpha x$ meets the graph of $f$ in $n$ points. Without loss of
  generality we may assume that there are $r$ roots of $f^H(x,1)$ greater
  than $\alpha$, and $n-r$ less than $\alpha$, for some $r$, where
  $0\le r\le n$.
  
\index{solution curves} 
  Consider the upper left quadrant. Recall that the solution curves
  are asymptotic to lines whose slopes are the roots of $f^H(x,1)$.
  Thus, there are $n-r$ solution curves that are eventually above the
  line $y=\alpha x$. Since $f\in\gsubplus_2$, each of these solution
  curves meets the $x$-axis in $(-\infty,0)$ which is below the line
  $y=\alpha x$. Thus, there are $n-r$ intersection points in the upper
  left quadrant.
  
  Similarly, there are $r$ intersections in the lower right quadrant.
  We've found $n$ intersections, and the proposition is proved.
\end{proof}

Our next result is that all the
homogeneous parts of a polynomial in $\gsubplus_2$ are in
$\allpolypos$. There are simple counterexamples that show that this is
not true for arbitrary polynomials in $\gsubpos_2$. The proof uses
properties of homogeneous polynomials in three variables - see
Lemma~\ref{lem:sub-total-degree}.

\begin{lemma} \label{lem:sub-total-degree-2}
  Suppose $f\in\gsubplus_2(n)$ and $g_r$ is the polynomial consisting
  of all terms of total degree $r$. If we set $h_r(x) = g_r(x,1)$ then
  $h_r\in\allpolypos$ for $1\le r \le n$ and we have interlacings
$$ h_n \lesslesseq h_{n-1} \lesslesseq h_{n-2} \lesslesseq \cdots$$
\end{lemma}

If $f\in\gsubpos_2$ then the graphs of $f$ and $f^H$ have the same
asymptotic behavior. Although these asymptotes do not interlace in
general, they do if $f\in\gsubplus_2$:

\begin{lemma} \label{lem:p2plus-fh}
  If $f\in\gsubplus_2$ then $f(x,\alpha)\greateqeq f^H(x,\alpha)$ for $\alpha$
  sufficiently large.  
\end{lemma}
\begin{proof}
  Assume $f\in\gsubplus_2(n)$ and let $g_r$ be the homogeneous polynomial
  consisting of all terms of $f$ of degree $r$.  We can write
  \begin{equation}
    f(x,\alpha) = \sum_{i=0}^n g_i(x,\alpha)
  \end{equation}
  We know from Lemma~\ref{lem:sub-total-degree-2} that $g_i(x,\alpha)\lesslesseq
  g_{i-1}(x,\alpha)$ for $1\le i \le n$.  We show that if $\alpha$ is
  sufficiently large then $f(x,\alpha)$ sign interlaces
  $f^H(x,\alpha)$ and that the sign of $f(x,\alpha)$ on the $i$-th
  largest root of $f^H$ is $(-1)^{n+i}$.
  
  If $f^H(\beta,1)=0$ then
  $f^H(\alpha\beta,\alpha)=g_n(\alpha\beta,\alpha)=0$ since $f^H$ is
  homogeneous. Upon computing $f(\alpha\beta,\alpha)$ we find
  $$
  f(\alpha\beta,\alpha) = 0 + g_{n-1}(\alpha\beta,\alpha) +
  O(\alpha^{n-2})$$
  Now since $g_n(x,\alpha)\lesslesseq
  g_{n-1}(x,\alpha)$ and the leading coefficients of both
  $g_n(x,\alpha)$ and $g_{n-1}(x,\alpha)$ are positive we know that
  the sign of $g_{n-1}(\alpha\beta,\alpha)$ is $(-1)^{n+i}$ where
  $\beta$ is the $i$-th largest root of $f^H(x,1)$.  Since
  $g_{n-1}(\alpha\beta,a)$ is $O(\alpha^{n-1})$ it follows that for
  $\alpha$ sufficiently large
  $$
  sgn\,g_{n-1}(\alpha\beta,\alpha) = sgn\, f(\alpha\beta,\alpha) =
  (-1)^{n+i}$$
  and the conclusion follows.

\end{proof}

\section{The representation of $\rupint{2}$ by determinants}
\label{sec:p2-det-rep}

The main result of this section is that there are simple determinant
representations for polynomials in $\gsubplus_2$ and $\rupint{2}$.
Before we can prove that polynomials in $\gsubplus_2$ can be
represented by determinants, we need the following result due to
Vinnikov and Helton, see \cites{vinnikov,vinnikov-helton,lax-conj}.

\begin{theorem}\label{thm:vinnikov}
  Suppose that $f(x,y)$ is a polynomial with the property that
  $f(\alpha x,\beta x)\in\allpoly$ for all $\alpha,\beta\in\reals$. If
  $f(0,0)=1$ then there are symmetric matrices $B,C$ such that 
$$ f(x,y) = |I+xB+yC|.$$
\end{theorem}

\begin{theorem}\label{thm:det-thm}
  If $f\in\gsubplus_2$ and $f(0,0)=1$ then there are positive definite symmetric
  matrices $D_1,D_2$ so that 
$$ f(x,y) = |I + x D_1 + y D_2|.$$
\end{theorem}
\begin{proof}
  Since Proposition~\ref{prop:p2-xax} shows that $f(x,y)$
  satisfies the hypothesis of Theorem~\ref{thm:vinnikov}, Vinnikov's
  result shows that there are symmetric matrices $B,C$ so that
  $f(x,y)=|I+xB+yC|$. The result now follows from
  Corollary~\ref{cor:pos-eigen}.
\end{proof}

\begin{cor}\label{cor:p2-det}
  If $f\in\rupint{2}$ then there is a symmetric matrix $C$, a positive
  definite diagonal matrix $D$, and a constant $\beta$ so that
  $$f(x,y) =\beta |xI + y D+C|.$$
\end{cor}
\begin{proof}
  There is an $a$ so that $f(x+a,y)\in\gsubplus_2$, so we can find
  positive definite matrices $D_1$,$D_2$ and a constant
  $\alpha=f(a,0)$ so that
  \begin{align*}
    f(x+a,y) &= \alpha|I+xD_1+yD_2| \\
    f(x,y) &= \alpha|I+(x-a)D_1+yD_2| \\
\intertext{Since $D_1$ is positive definite, it has a positive
  definite square root $D_1^{1/2}$}
    f(x,y) &= \alpha|D_1|\cdot\bigl|D_1^{-1/2}(I-aD_1)D_1^{-1/2} + xI + y
    D_1^{-1/2}D_2D_1^{-1/2}\bigr|\\
\intertext{If we let $A$ be the symmetric matrix
  $D_1^{-1/2}(I-aD_1)D_1^{-1/2}$, $B$ the positive definite matrix
  $D_1^{-1/2}D_2D_1^{-1/2}$, and $\beta=\alpha|D_1|$ then}
f(x,y) &= \beta|xI+yA+B| \\
\intertext{Now let $ODO^t=A$ where $O$ is orthogonal and $D$ is diagonal}
f(x,y) &= \beta|xI+ yD + O^tBO|
  \end{align*}
Since $O^tBO$ is symmetric, the corollary is proved.
\end{proof}

\begin{remark}
  We can use the determinant representation of $\rupint{2}$ to show that
  $\rupint{2}$ is closed under differentiation without using any
  geometry.  Recall that \index{principle submatrix} the
  characteristic polynomials of a symmetric matrix and any of its
  principle submatrices interlace (see Theorem~\ref{thm:principle-1}).
    
    Suppose that $f(x,y)\in\rupint{2}(n)$, represent it by $|xI+yD+C|$ as
    above, and let $M=xI+yD+C$.  Let $\{d_i\}$ be the diagonal
    elements of $D$. Since the only occurrences of $x$ and $y$ in $M$
    are on the diagonal, it is easy to see that
\begin{gather*}
\frac{\partial}{\partial x} \ |xI+yD+C| =
|M[1]| + |M[2]| + \cdots +|M[n]| \\
\frac{\partial}{\partial y} \ |xI+yD+C| =
d_1\,|M[1]| + d_2\,|M[2]| + \cdots + d_n\,|M[n]| \\
\end{gather*}

If we substitute $\alpha$ for   $y$ we see that 
$$ |xI + \alpha D+ C| \lesslesseq |M[i]](x,\alpha)$$
since principle submatrices interlace. All the polynomials $|M[i]|$
have positive homogeneous part, so we  simply add the interlacings to
conclude that \\ $f(x,\alpha) \lesslesseq \frac{\partial}{\partial x}
f(x,\alpha)$. It follows that $f \lesslesseq \frac{\partial
  f}{\partial x}$, and in particular the derivative is in $\rupint{2}$.

\end{remark}

  \section{When some coefficients are zero}
  \label{sec:when-some-coeff}

If $f(x)$ is a polynomial in $\allpoly$, then we know
two facts about zero coefficients:
\begin{enumerate}
\item If the constant term is zero then $f$ is divisible by $x$.
\item If two consecutive coefficients are zero, then all the coefficients
  of terms of smaller degree are zero.
\end{enumerate}
We investigate the implications of zero coefficients for polynomials
in $\rupint{2}$. We will show that if a row (or column) has two
consecutive zero coefficients then there is a large triangular region
of zero coefficients. For example, the coefficient array of the
polynomial $f(x,y)=(x+y)^3(x+y+1)^2$ has a triangular block of zeros:

\[
\begin{array}{cccccc}
1 \\ 2 & 5 \\ 1 & 8 & 10 \\ 0 & 3 & 12 & 10 \\ 0 & 0 & 3 & 8 & 5 \\
0 & 0 & 0 & 1 & 2 & 1
\end{array}
\]
The polynomial $\partial f/\partial x$ has constant term zero, yet is
irreducible.

\begin{lemma} \label{lem:p2-zero}
  If $f = \sum a_{i,j}x^iy^j \in \rupint{2}$  has a row (or column) with two
  consecutive  zero coefficients, say $a_{r,s-1}=a_{r,s}=0$, then
then $a_{i,j} = 0$ for all $i+j \le r+s$.
\end{lemma}
\begin{proof}
  We use the fact that if $g(x)$ interlaces $x^k h(x)$, then $g$ is
  divisible by $x^{k-1}$. If we write $f = \sum f_i(x)y^i$ then the
  hypothesis says that $f_r$ has two consecutive zeros, so $f_r$ is
  divisible by $x^s$. We then know that $f_{r+1}$ is divisible by
  $x^{s-1}$, $f_{r+2}$ is divisible by $x^{s-2}$ and so on. 

  Next, write $f = \sum F_j(y)x^j$. The above paragraph shows that
  $a_{r+1,k}=a_{r,k}=0$ for $0\le k \le s-1$, so $F_{k}$ has two
  consecutive zeros, and hence is divisible by $x^r$. 
Continuing, $F_{k+1}$ is divisible by $x^{r-1}$, and so on. Thus,
we've found that all $a_{ij}$ are zero if $i+j\le r+s$. 
\end{proof}

If we consider polynomials in $\gsubplus_2$ that have a quadrilateral
where the strict rhombus inequality is not satisfied, then there is
also a large triangular region with zeros.

\begin{lemma}\label{lem:p2-strict-quad-1}
  If $f = \sum a_{i,j}x^iy^j\in\gsubplus_2$ and there are $r,s$ such
  that \\ $a_{r,s} a_{r+1,s+1} - a_{r+1,s}a_{r,s+1}=0$ then $a_{i,j} = 0
  $ for all $i+j\le r+s+1$.
\end{lemma}
\begin{proof}
Write $f = \sum f_i(x)y^i$.  The expression $a_{r,s} a_{r+1,s+1} -
a_{r+1,s}a_{r,s+1}$ is the determinant of consecutive coefficients of
two interlacing polynomials. If it is zero, then by
Corollary~\ref{cor:log-con-coef} we have $a_{r,s}=a_{r,s+1}=0$. The
result now follows from Lemma~\ref{lem:p2-zero}.
\end{proof}

\begin{cor}\label{cor:p2-strict-quad}
  If $f = \sum a_{i,j}x^iy^j\in\gsubplus_2$ and $f(0,0)\ne0$ then all
  rhombus inequalities are strict.
\end{cor}

As with one variable, we find that agreement at two consecutive
coefficients implies lots of agreement.

\begin{cor}
  Suppose that $f\greateqeq g$ in $\gsubplus_2$, $f=\sum
  a_{i,j}x^iy^j$, and $g = \sum b_{i,j}x^iy^j$. If there are $r,s$ so
  that $a_{r,s}=b_{r,s}$ and $a_{r+1,s} = b_{r+1,s}$ then $a_{i,j} =
  b_{i,j}$ for $i+j\le r+s+1$. 
\end{cor}
\begin{proof}
  $f-g$ is in $\gsubplus_2$, and has two consecutive zeros.
\end{proof}

We can replace the consecutive zero condition with equal minimum
degree.

\begin{lemma}
  Suppose that $f=\sum f_i(x)y^i\in\rupint{2}$, and let $d$ be the
  minimum $x$-degree. If two consecutive $f_i$ have degree $d$ then so
  do all earlier ones.
\end{lemma}
\begin{proof}
  We first differentiate $d$ times with respect to $f$, and so assume
  that two consecutive coefficients are zero. By
  Lemma~\ref{lem:newton-zeros} all earlier coefficients are zero for
  all $x$, and so are identically zero. 
\end{proof}

\section{The Hadamard product}
\label{sec:hadamard}

\index{Hadamard product} 

We can realize the Hadamard product as a coefficient in a polynomial
in two variables.  Begin with two polynomials of the same degree:

\begin{align*}
  f(x) & = a_0 + a_1x + a_2x^2 + \dots + a_nx^n \\
  g(x) & = b_0 + b_1x + b_2x^2 + \dots + b_nx^n \\
\intertext{Now reverse $g$ and homogenize it.}
  G(x,y) & = b_nx^n + b_{n-1}x^{n-1}y + b_{n-2}x^{n-2}y^2 + \dots + b_0y^n. \\
\intertext{If we multiply, and write as a series in $y$}
  f(y)G(x,y) & = a_0b_n x^n \\
  & + (a_1b_nx^n + a_0b_{n-1}x^{n-1}) y  + \dots \\
  & + (a_0b_0 + \dots + a_nb_nx^n)y^n + \dots \\
  & + a_nb_ny^{2n} 
\end{align*}
\noindent%
\index{Hadamard product!proof}
then the Hadamard product $f\ast g$ is the coefficient of $y^n$.

\begin{theorem} \label{thm:hadamard-1}
  The Hadamard product\footnote{It does not define a map
    $\allpoly\times\allpoly\longrightarrow\allpoly$. For instance
    $(x^2-1)\ast(x^2-1)=x^2+1$ is not in $\allpoly$.}  of polynomials
  defines a bilinear map
  $\allpoly(n)\times\allpolypm(n)\longrightarrow \allpoly(n)$. If
  $f_1\greateq f_2$ then $f_1\ast{}g \greateq f_2\ast{}g$.
\end{theorem}
\begin{proof}
Without loss of generality we may we assume that $g\in\allpolypos$.
Since $g\in\allpolypos$ its reversal is in $\allpolypos$, and
the homogenization $G(x,y)$ is in $\gsubpos_2$. Since $f\ast g$ is a
coefficient of $f(y)G(x,y)$ it is in $\allpoly$. Interlacing follows
since for fixed $g$ the map $f\mapsto f\ast g$ is a linear map
$\allpoly\longrightarrow\allpoly$. 
\end{proof}

Actually, $f(y)$ is not in $\gsubpos_2$, but rather in the closure of
$\gsubpos_2$ (see \chap{topology}). However, the coefficients of
$f(y)G(x,y)$ are in $\allpoly$.

  If we look more closely at the product $f(y)G(x,y)$ we see that it
  equals
\[
\cdots + (f\ast xg)\,y^{n-1} + (f\ast g)\,y^n + \frac{1}{x} (xf\ast
  g)\,y^{n+1} + \cdots 
\]
Consequently,
\begin{cor}\label{cor:hadamard-newton}
If $f,g,\in\allpolypos(n)$ then
\begin{align*}
  (1)\quad &
 (f\ast g)^2 \ge \left(\frac{n+1}{n}\right)^2\, \frac{1}{x}
(f\ast xg)(xf\ast g)\qquad\text{for $x\ne0$} \\
(2) \quad &
\begin{vmatrix}
  f\ast g & xf \ast g \\ f\ast xg & xf \ast xg
\end{vmatrix}\ge 0 \qquad\text{for $x\ge0$} \\
(3)\quad  & 
f\ast xg - \frac{1}{x} (xf \ast g) \greateqeq f\ast g
\end{align*}
  \end{cor}
  \begin{proof}
    The first part is Newton's inequality \mypage{thm:newton}. The
    second part follows from the first and the identity $xf\ast xg =
    x(f\ast g)$. Since $f\ast g$ and $\frac{1}{x}(xf\ast g)$
    interlace, and $\frac{1}{x}(xf\ast g)$ has degree $n-1$ we have
    the interlacings
\[
f\ast xg \greateqeq f\ast g \greateqeq \frac{1}{x} (xf \ast g)
\]
which implies the last part.
  \end{proof}

If we rewrite the last part in terms of coefficients we  get
  \begin{cor}
    If $f = \sum a_ix^i$, $g=\sum b_ix^i$ and $f,g\in\allpolypos$ then
\[
\sum a_i(b_{i+1}-b_i)x^i \in\allpoly.
\]
  \end{cor}

  Surprising, we can do better than Theorem~\ref{thm:hadamard-1}.  We
  use a differential operator in $\gsubpos_2$, and then evaluate it at
  $0$ to get a polynomial in $\allpoly$. See
  Corollary~\ref{cor:ifact-schur} for an alternate proof. Define
$$ x^i \,\ast'\,x^j = 
\begin{cases}
  i!\, x^i & i=j \\0 & \text{otherwise}
\end{cases}
$$
Notice that $f\ast' g = \expoper{}^{-1}\,f\ast g$. In terms of
coefficients
$$ f \ast' g = 
a_0b_0 + 1!a_1b_1x + 2!a_2b_2x^2 + \dots + n!a_nb_nx^n
$$

\index{differential operators}
\index{operator!differential}

\begin{theorem} \label{thm:hadamard-2}
    If $f\in\allpolypm(n), g\in\allpoly(n)$ then $f\ast'  g\in\allpoly(n)$. 
\end{theorem}
\begin{proof}
  Again assume that $f\in\allpolypos$.  Let $G(x,y)$ be the
  homogeneous polynomial in $\gsubpos_2$ corresponding to $g$. By
  Lemma~\ref{lem:pd-diff-prod} we know $f(\frac{\partial}{\partial
    x})G\in\gsubpos_2$, so evaluation at $x=0$ gives a polynomial in
  $\allpoly$. Consequently, the map $f\times g\mapsto
  (f(\frac{\partial}{\partial x})G)(0,y)$ defines a linear
  transformation $\allpolypm(n)\times
  \allpoly(n)\longrightarrow\allpoly(n)$.
  
  The polynomial $f\ast' g$ is the reverse of
  $(f(\frac{\partial}{\partial x})G)(0,x)$. To verify this, we check
  it for $f=x^r$ and $g=x^s:$
  $$
  \left(\frac{\partial}{\partial
      x}\right)^r\left(x^sy^{n-s}\right)(0,x) =
  \begin{cases}
    0 & r\ne s\\
    r!x^{n-r} & r=s
  \end{cases}
  $$
\end{proof}

\begin{cor} \label{cor:hadamard-2}
  If $f\in\allpolypm(n), g_1\greateqeq g_2\in\allpoly(n)$ then
  $$f\ast' g_1 \greateqeq f\ast' g_2$$
\end{cor}



The last lemma showed that
$\allpolypm(n)\ast\allpoly(n)\subset\expoper{}(\allpoly)$. If we
allow analytic functions  then we get equality.

\begin{cor}\label{cor:hadamard-f}
  $\allpolyposf \ast \allpoly = \expoper{}(\allpoly)$
\end{cor}
\begin{proof}
 By taking limits we know that $\allpolyposf \ast \allpoly \subset
 \expoper{}(\allpoly)$.  Since $\expoper{}(f) = e^x \ast f$ we see the the
 containment is an equality. 
\end{proof}

There is a partial converse to Theorem~\ref{thm:hadamard-2}.

\begin{prop}
Suppose that $f$ is a polynomial.  If for all $g\in\allpolypos$ we have
that $f\ast g\in\allpoly$, then $\expoper{f}\in\allpoly$.
If for all $g\in\allpoly$ we have
that $f\ast g\in\allpoly$, then $\expoper{f}\in\allpolypmclose$.
\end{prop}

\begin{proof}
  Since $(1+x/n)^n\in\allpolypos$ we know $f\ast(1+x/n)^n\in\allpoly$.
  Now $(1+x/n)^n$ converges uniformly to $e^x$, so $f\ast
  e^x\in\allpoly$. The first part follows since $f\ast e^x = \exp(f)$,
  so it remains to show that $f\in\allpolypm$ in the second case.
  
  If $f=\sum a_ix^i$ then $f\ast x^i(x^2-1) = x^i(a_{i+2}x^2-a_i)$ so
  it follows that  $a_i$ and $a_{i+2}$ have the same sign. If
  $a_{i+2}$ is non-zero, and $a_{i+1}=0$ then consider $f\ast
  x^i(x+1)^2 = x^i(a_{i+2}x^2+a_i)$. Since $a_{i+2}$ and $a_i$ have
  the same sign, it follows that $a_i=0$. If we consider $f\ast
  x^{i-j+2}(x+1)^j$ we see that all 
  coefficients $a_{i+1},a_i,a_{i-1},\dots$ are zero. Thus, $f=x^if_1$
  where $f_1$ has all non-zero terms. Since the signs of the
  coefficients of $f_1$ are either the same or are alternating, we
  get that $f_1\in\allpolypm$.
\end{proof}

As an application of  Theorem~\ref{thm:hadamard-2} consider
\begin{lemma}
If $f\in\allpoly$ then $ (1-x)^{-n}\,\ast\,f$ is in $\allpoly$ for $n=0,1,2,\dots$.
\end{lemma} 
\begin{proof}
  From \eqref{eqn:laguerre} we know that $\expoper\,(1-x)^{-n}$ is in
  $\allpolyposf$ since $e^x\in\allpolynegf$ and
  $L_n(-x)\in\allpolypos$.  The conclusion
  follows from Corollary~\ref{cor:hadamard-f}.
\end{proof} 

The lemma explains the curious fact that
$f\mapsto(1+x+x^2+\dots)\ast f$ sends $\allpoly$ to itself, yet
$(1+x+x^2+\dots)=1/(1-x)$ is not in $\allpolyf$.

If we re-express the lemma in terms of coefficients we find that the
map \index{rising factorial}
\begin{equation}
  \label{eqn:nii}
  x^i \mapsto \frac{\rising{n}{i}}{i!}x^i
\end{equation}
which is just Theorem~\ref{thm:multiplier}.

\section{Some differential operators}
\label{sec:p2-diff-op}

We use the fact that polynomials in $\allpoly$ factor into linear
factors to get differential operators on $\gsubpos_2$. We would like
to know that if $f\in\gsubpos_2$ then the differential operator
$f(\frac{\partial}{\partial x},\frac{\partial}{\partial y})$
maps $\gsubpos_2$ to itself, but the best we are able to prove is
when $f$ is a product of linear factors - see Lemma~\ref{lem:pd-diff-prod}. 

\begin{lemma} \label{lem:fpx}
  If $f\in\allpoly$, $g\in\gsubpos_2$ then $f(\frac{\partial}{\partial
    x})g\in\gsubpos_2$. 
\end{lemma}
\begin{proof}
  It suffices to show that $g +\alpha \frac{\partial g}{\partial x}$
  is in $\gsubpos_2$, but this is immediate from $g \lesslesseq
  \frac{\partial g}{\partial x}$.
\end{proof}

Here is another differential operator that acts on $\gsubpos_2$.

\begin{cor} \label{cor:fpp}
  If $f\in\allpolypos$ and $g\in\gsubpos_2$ then
  $f\left(-\frac{\partial^2}{\partial x\,\partial y}\right)g\in\gsubpos_2$. 
\end{cor}

\begin{proof}
  If suffices to show that $h(x,y)=\left(a-\frac{\partial^2}{\partial
      x\,\partial y}\right)g$ is in $\gsubpos_2$ for positive $a$.
  Since  homogeneity is certainly satisfied since
  $h^H=ag^H$ it suffices to check substitution. If $\alpha\in\reals$
  then
  $$
  \Bigl(a-\frac{\partial^2}{\partial x\,\partial
      y}\Bigr)\,g\,(x,\alpha) = ag(x,\alpha)-\partialx
  \left(\partialy g\right)(x,a)
$$
and the conclusion follows from $g(x,\alpha)\lesslesseq \frac{\partial
  g}{\partial y}(x,\alpha)$ and Corollary~\ref{cor:fagp}.
\end{proof}

\begin{example}
  If we take $f=(x-1)^m$ and $g=(x+y+1)^n$ then
  \begin{align*}
       \Bigl(\frac{\partial^2}{\partial x\,\partial
      y}-1\Bigr)^m\,(x+y+1)^n &=
\sum_{i=0}^n (-1)^{m-i} \binom{m}{i}    \Bigl(\frac{\partial^2}{\partial x\,\partial
      y}\Bigr)^i \,(x+y+1)^n \\
&= \sum_{2i\le n} \binom{m}{i} \falling{n}{2i}\,(x+y+1)^{n-2i} \\
\intertext{and so we find that }
&
\sum_{2i\le n} (-1)^{m-i} \frac{(x+y+1)^{n-2i}}{i!(m-i)!(n-2i)!} \in\gsubpos_2
  \end{align*}
\end{example}

 The differential operator
$g(x\frac{\partial}{\partial y})$ maps $\gsubpos_2$ to itself.

\begin{lemma} \label{lem:gxdy}
  If $g\in\allpolypos$, $f\in\gsubpos_2$ then 
$g(x\frac{\partial}{\partial y})\,f(x,y)\in\gsubpos_2$.
\end{lemma}
\begin{proof}
  Since all coefficients of $g$ are positive the homogeneous part of
  $g(x\frac{\partial}{\partial y})\,f(x,y)$ has all positive
  coefficients. To show that $g(\alpha\frac{\partial}{\partial
    y})\,f(\alpha,y)\in\allpoly$ it suffices to note that $(b+\alpha
  \frac{\partial}{\partial y})\,f(\alpha,y)\in\allpoly$.
\end{proof} 

\begin{example}
  If we take $g=(x+1)^m$ and $f=(x+y+1)^n$ then
  \begin{align*}
    (x\frac{\partial}{\partial y}+1)^m \,(x+y+1)^n &=
\sum_{i=0}^m \binom{m}{i} (x\frac{\partial}{\partial y})^i (x+y+1)^n
\\
&= \sum_{i=0}^m \binom{m}{i} x^i(x+y+1)^{n-i}\falling{n}{i} \\
\intertext{and so}
\sum_{i=0}^n \binom{m}{i} x^i(x+y+1)^{n-i}\falling{n}{i}& \,\in\gsubpos_2
  \end{align*}
\end{example}

\begin{cor} \label{cor:fgxdy}
  If $f\in\allpoly(n)$, $g\in\allpolypos(n)$, $g=\sum b_i\,x^i$ then
$$ h(x,y) = \sum b_i\, x^i\, f^{(i)}(y)\in\gsubpos_2$$
\end{cor}
\begin{proof}
The homogeneous part of $h$ is 
$$ a_n \sum b_i x^i \falling{n}{i}  y^{n-i}$$ where $a_n$ is the leading
coefficient of $f$. Since $h^H(x,1) = a_nn!\, \expoper{(g)}$
it follows that $h^H$ is in $\allpoly$. Since all coefficients $b_i$
are the same sign $h^H(x,1)$ is actually in $\allpolypos$.
We now apply Theorem~\ref{thm:hadamard-2} for fixed $y$, which shows
that substitution holds.

\end{proof}

\begin{cor} [Hermite-Poulin]
\index{Hermite-Poulin Theorem}
  If $f,g$ are given in Corollary~\ref{cor:fgxdy} then 
$$\sum b_i\, f^{(i)}(y)\in\allpolypos$$
\end{cor}

\begin{cor} [Schur-\Polya] \label{cor:f-schur}
\index{Schur-\Polya Theorem}
    If $f,g$ are given in Corollary~\ref{cor:fgxdy} then 
$$ \sum b_i\, x^i\, f^{(i)}(x) \in\allpoly$$
\end{cor}

Here is a proof of a special case of Theorem~\ref{thm:hadamard-2}. 
\index{Hadamard product!variant} 

\begin{cor} [Schur] \label{cor:ifact-schur}
  \index{Schur's Theorem} If $f,g\in\allpolypos$,  $f =
  \sum a_ix^i$ $ g = \sum b_i x^i$, then
$$\sum i!\,a_i\,b_i\,x^i\in\allpoly$$
\end{cor}
\begin{proof}
  Let $G(x,y)$ be the homogenization of $g$. Since $g\in\allpolypos$
  we know that $G(x,y)\in\gsubpos_2$.  Since $f^{(i)}(0) = i!a_i$ it
  follows that
$$ \sum i!\,a_i\,b_i\,x^i = 
f(x\frac{\partial}{\partial y})\,G(x,y)\vert_{y=0}$$ and hence is in $\allpoly$.
\end{proof}

\added{2/1/7}%
  It is possible for two different linear transformations
  $\allpoly\longrightarrow\allpoly$ to agree on all polynomials up to
  degree $n$. For instance, suppose that $f$ has degree $n+1$ and all
  distinct roots. Choose $g\in\allpoly$ such that $f-g = \alpha
  x^{n+1}$ for some non-zero $\alpha$. If we define
\[
 S(h) = f(\diffd)h\quad\text{and}\quad T(h)=g(\diffd)h
\]
then $(T-S)(h) = \alpha \diffd^{n+1}h$, and this is zero if the degree
of $h$ is at most $n$. 

Differential operators  satisfy a uniqueness condition.
\index{differential operators!uniqueness}

\begin{lemma}
  Suppose that $T(g) = f(\diffd)g$ where $f\in\allpoly$ has degree
  $n$.  If $S$ is a  linear transformation
  $\allpoly\longrightarrow\allpoly$ such that
\[
 S(x^i) = T(x^i)\quad \text{for}\ 0\le i\le n+2
\]
then $S=T$.
\end{lemma}
\begin{proof}
  Since $(x-t)^{n+3}\lesslesseq (x-t)^{n+2}$ for any real $t$ we have
\[
 S(x-t)^{n+3} \lace S(x-t)^{n+2} = T(x-t)^{n+2}
\]
$f$ has degree $n$ so  $T(x-t)^{n+2}$ has a factor of
$(x-t)^2$ and therefore $S(x-t)^{n+3}$ has $t$ as a root. Thus
\begin{align*}
  0 = S(x-t)^{n+3}(t) 
&=
S(x^{n+3})(t) + \sum_0^{n+2} \binom{n+3}{i}\,(-t)^{n+3-i}\, S(x^i)(t) \\
&= 
S(x^{n+3})(t) + \sum_0^{n+2} \binom{n+3}{i}\,(-t)^{n+3-i}\, T(x^i)(t) \\
&=  S(x^{n+3})(t) + \biggl[ T(x-t)^{n+3}(t) - T(x^{n+3})(t)\biggr] 
\end{align*}
Now  $t$ is a root of $T(x-t)^{n+3}$ so we conclude that $T(x^{n+3}) = S(x^{n+3})$.
\end{proof}
If we take $f=1$ we get
\begin{cor}
  If $T$ is a  linear transformation
  $\allpoly\longrightarrow\allpoly$ such that 
$T(x^i)=x^i$ for $i=0,1,2$  then $T$ is the identity.
\end{cor}

\section{Realizing interlacing}
\label{sec:realizing-interlacing}
\index{extension!of 2 interlacing polynomials}

There is a  close relationship between interlacing and
products:  A pair of interlacing polynomials can be realized as the
first two coefficients of a product. If the degrees are different the
representation is essentially unique.  A product is also 
determined by its first term and its homogeneous part.

\begin{lemma} \label{lem:product}
  Assume that $f_0\in\allpoly(n)$ is monic, $f_1$ has positive leading coefficient,
  and $f_0\lessless f_1$.  There exist unique $b_i,c_i$ such that
  \eqref{eqn:product-1} is satisfied, and all $c_i$ are positive. The
  product is in $\gsubpos_2$.
\begin{equation}\label{eqn:product-1}
  \prod_{i=1}^n(x+b_i+c_iy) = f_0(x) + f_1(x)y + \cdots 
\end{equation}
\end{lemma}

\begin{proof}
 Expanding the product, we find
  \begin{align*}
    f_0(x) & = \prod_{i=1}^n(x+b_i) \\
    f_1(x) & = \sum_{j=1}^n\,c_j\prod_{i\ne j}^n(x+b_i) \\
    \intertext{We can rewrite $f_1$ in terms of $f_0$}
    f_1(x) &= \sum_{j=1}^n  {c_j} \,\frac{f_0(x)}{x + b_j}
  \end{align*}
If we are given $f_0,f_1$ then we find the $b_i$ as roots of $f_0$,
and the $c_i$ as the coefficients of $f_1$ as in Lemma~\ref{lem:sign-quant}.
\index{interpolation!basis}

\end{proof}

\begin{remark}\label{rem:strict-limit}
    There is straightforward way to express $f\lesslesseq g$ as a
    limit of strict interlacings. We can find non-negative $a_i$ such
    that 
    \[ \prod(x+a_i\,y+b_i) = f + y\,g + \cdots \] It's enough to perturb
    the constants. For positive $\epsilon$ and \emph{distinct}
    $\epsilon_i$ define
\[ \prod(x+(a_i+\epsilon)\,y+(b_i+\epsilon_i)) = f_\epsilon + y\,g_\epsilon + \cdots \]
    We claim that have that
    \begin{enumerate}
    \item $f_\epsilon\lessless g_\epsilon$ for all positive $\epsilon$
      and sufficiently small $\epsilon_i$.
    \item $\displaystyle\lim_{\epsilon,\epsilon_i\rightarrow0^+} f_\epsilon = f$
    \item $\displaystyle\lim_{\epsilon,\epsilon_i\rightarrow0^+} g_\epsilon = g$
    \end{enumerate}
    We only need to check the first claim. Since the roots of
    $f_\epsilon$ are $-(b_i+\epsilon_i)$ and the $\epsilon_i$ are
    distinct, we see that if the $\epsilon_i$ are small enough then
    all roots of $f_\epsilon$ are distinct. Since $a_i+\epsilon$ is
    positive it follows that $f_\epsilon\lessless g_\epsilon$.
  \end{remark}

We can also specify the first term and the homogeneous part.

\begin{lemma}\label{lem:product-2}
  If $h\in\allpolypos$ and $f_0\in\allpoly$ are monic of the same
  degree then there is a product
$$ g(x,y) = \prod_i(a_i y + x + b_i) = f_0(x) + f_1(x)y + \cdots  $$
such that $g\in\gsubpos_2$ and $g^{H}(x,1) = h(x)$.
\end{lemma}
\begin{proof}
If we write $h(x) = \prod(x+a_i)$ and $f_0(x) = \prod(x+b_i)$ then any
ordering of the $a_i$'s and $b_i$'s gives a product.  Since all $a_i$
 have the same sign, the product is in $\gsubpos_2$.
\end{proof}

As opposed to the unique representation of $f\lessless g$ as a
product, we see from Lemma~\ref{lem:product-nl}  that $f\greateq g$ are the first two
coefficients of many products.

\begin{lemma} \label{lem:product-nl}
Given $f_1\greateq f_0$ there exist $b_i,c_i,d_i$ such
that \eqref{eqn:product-3} is satisfied, and the determinants
$\smalltwodet{1}{b_i}{c_i}{d_i}$ are negative for all $i$. The product
is in the closure of $\rupint{2}$.
\begin{equation}\label{eqn:product-3}
 F(x,y) = \prod_{i=1}^n(x+b_i+c_iy+d_ixy) = f_0(x) + f_1(x)y + \cdots 
\end{equation}
\end{lemma}

\begin{proof}
  Define $\alpha_i,\beta_i,\beta$ by
\begin{align*}
  f_0(x) &= \prod_{i=1}^n(x+\alpha_i)\\
  f_1(x) &= \beta f_0(x) + \sum_{i=1}^n \beta_i \frac{f_0(x)}{x+\alpha_i}\\
  \intertext{Choose $d_i$ so that $\sum_{i=1}^n d_i = \beta$. We see
    that $f_0,f_1$ are the leading terms of the product}
  & \prod_{i=1}^n (x +\alpha_i + d_ixy + (\beta_i+\alpha_id_i)y)\\
\intertext{since the first term is obviously $f_0$ and the second
  term is}
\sum_i \frac{f_0}{x+\alpha_i} (d_ix+\beta_i+\alpha_id_i) &=
\left(\sum_i d_i\right) f_0   + 
\sum_i \beta_i \frac{f_0}{x+\alpha_i} 
\end{align*} 
Identify $b_i=\alpha_i,c_i=\beta_i+\alpha_id_i$ and note that the
determinant is equal to $\beta_i$ which is negative since $f_1
\greateq f_0$. See see Lemma~\ref{lem:fyg} and
Corollary~\ref{cor:sub-xyc} for similar results.
\end{proof}

Another consequence of looking at coefficients of products is the
following lemma, which gives a ``canonical'' polynomial that
interlaces two polynomials.

\begin{lemma}
  If $f=\prod(x+b_i)$ and $c_i,d_i$ are non-negative, then we have a
  diagram of interlacing polynomials

\centerline{
\xymatrix{
f \ar@{<-}[rr] \ar@{<-}[d]& & \ar@{<-}[d] {\displaystyle\sum c_i \,\frac{f}{x+b_i}}  \\
{\displaystyle\sum  d_j \, \frac{f}{x+b_j}} \ar@{<-}[rr] && {\displaystyle\sum_{i\ne j} c_i d_j \,
  \frac{f}{(x+b_i)(x+b_j)}} 
}
}

\end{lemma}
\begin{proof}
  Lemma~\ref{lem:sub-prod-1} shows that Theorem~\ref{thm:product-4} holds for arbitrary
  products, so the interlacings follow from the fact that the four
  terms in the diagram are the coefficients of $1,y,z,yz$ in the
  product
$$ \prod_i (x+ c_i y + d_i z + b_i).$$
\end{proof}

Another interpretation of this lemma is that if we are given $f,g,h$
with positive leading coefficients, and $f\lessless g,h$ then there
is a canonical choice of $k$ such that $\smallsquare{f}{g}{h}{k}$ where all
interlacings are $\lessless$. In the case that $g=h=f^\prime$ the
canonical polynomial is $\frac{f^{\prime\prime}}{2}$.
Since these polynomials are the coefficients of a polynomial in
$\gsubpos_3$  we will see  that their determinant is negative. We can verify
this directly:

\begin{multline}\label{eqn:p2-canon} 
\begin{vmatrix}
  f & \displaystyle\sum_i c_i \,\frac{f}{x+b_i} \\
\displaystyle\sum_j  d_j \, \frac{f}{x+b_j} & \displaystyle\sum_{i\ne j} c_i d_j \,
  \frac{f}{(x+b_i)(x+b_j)} 
\end{vmatrix} = \\
f^2\left(\displaystyle\sum_{i\ne j} \frac{c_id_j}{(x+b_i)(x+b_j)} -
\displaystyle\sum_{i,j}\frac{c_id_j}{(x+b_i)(x+b_j)}\right)\\\notag
 =  - f^2\left(\displaystyle\sum_i \frac{c_id_i}{(x+b_i)^2}\right) \ <
0
\end{multline}

\section{Karlin's conjecture}
\label{sec:karlin}

We show that certain matrices formed from the coefficients of a
polynomial in $\gsubplus_2$ are totally positive for positive values
of $x$. A special case of this yields a solution to a version of
Karlin's problem \cites{csordas00,dimitrov98,karlin}.

A Hankel \index{Hankel matrix} matrix has equal entries on the minor
diagonals. Given a sequence $\aaa = (\alpha_1,\alpha_2,\dots)$ and a
positive integer $d$ we can form the Hankel matrix $H(\aaa;d)$:

\[
H(\aaa;d) =\begin{pmatrix}
  \alpha_1 & \alpha_2 & \alpha_3 & \hdots & \alpha_d \\
  \alpha_2 & \alpha_3 & \alpha_4 & \hdots & \alpha_{d+1} \\
  \alpha_3 & \alpha_4 & \alpha_5 & \hdots & \alpha_{d+2} \\
  \vdots   & \vdots   & \vdots   & \ddots &  \vdots \\
  \alpha_d & \alpha_{d+1} & \alpha_{d+2} & \hdots & \alpha_{2d}
\end{pmatrix}
\]

We are also interested in the Toeplitz matrix $T(\aaa;d)$ associated
to $H(\aaa;d)$. If we reverse all the columns of $H(\aaa;d)$ we get
$T(\aaa;d)$: 

\index{Toeplitz matrix}

$$
T(\aaa;d) = \begin{pmatrix}
  \alpha_d & \alpha_{d-1} & \alpha_{d-2} & \hdots & \alpha_1 \\
  \alpha_{d+1} & \alpha_{d} & \alpha_{d-1} & \hdots & \alpha_{2} \\
  \alpha_{d+2} & \alpha_{d+1} & \alpha_d & \hdots & \alpha_{3} \\
  \vdots   & \vdots   & \vdots   & \ddots &  \vdots \\
  \alpha_{2d} & \alpha_{2d-1} & \alpha_{2d-2} & \hdots & \alpha_{d}
\end{pmatrix}
$$

  We can form $T(\aaa;d)$ from $H(\aaa;d)$ by shifting the rightmost column
$d-1$ places, the previous column $d-2$ places, and so on.
Consequently, we note the useful relation
$$ det\left(T(\aaa;d)\right) =
(-1)^{\binom{d}{2}}det\left(H(\aaa;d)\right) 
$$

The matrix $\varphi(f)$ is a Toeplitz matrix. The first result picks
out a submatrix of $\varphi(f)$.

\index{Toeplitz matrix}

\begin{lemma}\label{lem:toeplitz-1}
  Suppose that $f=a_0 + a_1x + a_2x^2 + \dots\in\allpolypos$. For any
  positive integer the matrix $T(\,(a_0,a_1,\dots);d)$ is totally
  positive.   
\end{lemma}
\begin{proof} 
  The matrix $T(\,(a_0,a_1,\dots);d)$, with rows and columns reversed,
  is a submatrix of $\varphi(f)$.
\end{proof}

We have an important consequence of this simple observation.

\begin{theorem} \label{thm:karlin-1}
  Suppose that $f(x,y) = \sum f_i(x)y^i \in\gsubplus_2$. The Toeplitz matrix\\  
$T(\,(f_0(x),f_1(x),f_2(x),\dots)\,;d)$ is totally positive for all
non-negative $x$. In addition 
\[  (-1)^{\binom{d}{2}}\,det\,H(\,(f_0(x),f_1(x),f_2(x),\dots)\,;d)\ge0
\text{ for } x\ge 0
\]
\end{theorem}

\begin{cor} \label{cor:karlin-1}
  If $f\in\allpolypos$ then
  \begin{enumerate}[(1)]
  \item
    $T(\,(f(x),\frac{f^\prime(x)}{1!},\frac{f^{(2)}(x)}{2!},\frac{f^{(3)}(x)}{3!},\dots 
    )\, );d)$ is totally positive for positive $x$.
  \item
    $(-1)^{\binom{d}{2}}\, det \,H
    (\,(f(x),\frac{f^\prime(x)}{1!},\frac{f^{(2)}(x)}{2!},\frac{f^{(3)}(x)}{3!},\dots )\, ;d) \ge 0 $ for $x\ge0$.
  \end{enumerate}
\end{cor}
\begin{proof}
  Use the Taylor series for $f(x+y)$ and apply Theorem~\ref{thm:karlin-1}.
\index{Taylor series}
\end{proof}

\index{determinants!of coefficients of $\gsubplus_2$}

We can restate the last part by saying that the sign of the
determinant of the matrix whose $ij$ entry is
$\frac{f^{(i+j)}}{(i+j)!}$ is $(-1)^{\binom{d}{2}}$ whenever
  $f\in\allpolypos$, and $x\ge0$.

\index{Karlin's conjecture}

Karlin's conjecture \cite{csordas00} is that 
\begin{equation}\label{eqn:karlin-conj}
(-1)^{\binom{d}{2}}\, det \,
    H(\,(f(x),{f^\prime(x)},{f^{(2)}(x)},{f^{(3)}(x)},\dots )\, );d)
    \ge 0  \text{ for $x\ge0$} 
  \end{equation}
  It is easy to construct counterexamples to
  this.\footnote{$f=(x+1)(x+2)(x+3)$ and $d=3$ is one.}  Karlin's
  conjecture is true in special cases. Here's an
  example from  \cite{dimitrov98}.

\begin{cor}
  Define $\expoper{}_\ast^{-1}(x^iy^j) = j!\,x^iy^j$, and let
  $f\in\allpolypos$. If \\ $\expoper{}_\ast^{-1}(f(x+y))\in\gsubpos_2$ then
  Karlin's conjecture \eqref{eqn:karlin-conj} holds.
\end{cor}

\index{determinants!of translates}
We can also apply Theorem~\ref{thm:tp} to determinants of translates
of a function. Suppose
  $\affa(x)=x+1$ and let $T(y^n)=\falling{y}{n}$. Recall that
  \begin{align*}
    \sum_i i^n \frac{y^i}{i!} & = e^y T^{-1}(y^n).\\
\intertext{We want to determine the action of $e^{\affa y}$ on
  polynomials:}
e^{\affa y}x^n &=  \sum_i \frac{y^i}{i!} (x+i)^n \\
&= \sum_{r,i} \binom{n}{r} x^r i^{n-r} \frac{y^i}{i!} \\
&= \sum_r \binom{n}{r}x^r \, \sum_i i^{n-r} \frac{y^i}{i!} \\
&= e^y \sum_r \binom{n}{r}x^r\,T^{-1}y^{n-r} \\
&= e^y T_\ast^{-1} (x+y)^n\\
\intertext{where $T_\ast^{-1}$ is the induced map.  By linearity}
e^{\affa y}f(x) &= e^y T^{-1}_\ast\,f(x+y)\\
&= \sum_{k=0}^\infty \frac{\affa^k f(x)}{k!}y^k\\
  \end{align*}

Now we know that $T^{-1}:\allpolypos\longrightarrow\allpolypos$ so for any
positive value $\alpha$ of $x$  
$$\left[T^{-1}_\ast f(x+y)\right](\alpha,y) = T^{-1} f(\alpha+y)
\in\allpolypos$$
since $f(\alpha+y)\in\allpolypos$ and consequently
$e^y T^{-1}_\ast\,f(x+y)\in\allpolyposf$ for any positive $x$. We can
now apply Theorem~\ref{thm:tp} to conclude

\begin{cor}\label{cor:tp-affa}
  If $f\in\allpolypos$ then 
$T((f,\affa f,\affa^2 f/2!,\affa^3 f/3!,\dots);d)$ is totally positive for
all $x\ge0$. 
\end{cor}

\section{Determinants of coefficients}
\label{sec:coef-det}

\index{determinants!of coefficients in $\gsubpos_2$}
In Lemma~\ref{lem:log-con-coef-2} we saw that the signs of certain two by two
determinants of interlacing polynomials were always positive. Since
any interlacing $f\lessless g$ can be viewed as the first two
coefficients of a polynomial in $\gsubpos_2$ by Lemma~\ref{lem:product}, we can
generalize these results to determinants of coefficients of
a polynomial in $\gsubpos_2$. 

So, let's begin by considering the coefficients of $x^iy^j$ in a
product. 

\begin{theorem} \label{thm:product-4}
  Suppose $f\in\gsubpos_2$, $f(x,0)$ has all distinct roots, and set
\begin{equation} \label{eqn:product-4}
 f(x,y)  = \sum_{i,j} d_{ij}x^iy^j = \sum_i f_i(x)y^i
\end{equation}
For $1\le i,j\le n$ we have
\begin{align*}
  (rows) & \hspace*{1cm}\begin{vmatrix}{d_{ij}}&{d_{i+1,j}}\\
{d_{i+1,j}}&{d_{i+2,j}}\end{vmatrix} < 0 \\
  (columns) & \hspace*{1cm}\begin{vmatrix}{d_{ij}}&{d_{i,j+1}}\\
{d_{i,j+1}}&{d_{i,j+2}}\end{vmatrix} \le 0 \\
  (squares) &
  \hspace*{1cm}\begin{vmatrix}{d_{ij}}&{d_{i+1,j}}\\
{d_{i,j+1}}&{d_{i+1,j+1}}\end{vmatrix} 
  < 0
\end{align*}
If the roots of $f(x,0)$ aren't all distinct then all $<$ signs are
replaced with $\le$.  
\end{theorem}

\begin{proof}
  We know that $f_0\lessless f_1\lessless\cdots$.  The $(rows)$
  inequalities are Newton's inequalities since all $f_i\in\allpoly$.
  The $(columns)$ results follow upon writing the sum as $\sum
  g_i(y)x^i$. These inequalities are $\le$ since it necessarily true
  that $f(0,y)$ has all distinct roots - consider $(x+y+1)(x+2y+1)$.
  The $(squares)$ are negative since consecutive $f_i$'s interlace
  (Corollary~\ref{cor:log-con-coef}).
\end{proof}

\begin{remark}
  We can also derive this from Lemma~\ref{lem:p2-quad-ineq} in the same way
  that Newton's inequalities are derived from properties of the
  quadratic. 
\end{remark}

\added{9/2/7} There are polynomials in $\rupint{2}$ for
  which all terms of odd degree are zero, For example, consider the
  coefficient array of $ \bigl( (x+y)^2 -1\bigr)^3.$

\[
\begin{matrix}
-1 & . & 3 & . & -3 & . & 1 \\
. & 6 &.& -12 &.& 6 \\
3 &.& -18 &.& 15 \\
.& -12 &.& 20 \\
-3 &.& 15 \\
. & 6 \\
1
\end{matrix}
\]
In such a case we always have alternating signs. 
\index{coefficients!alternating signs}

\begin{lemma}\label{lem:p2-even-alternate}
  If $f(x,y)\in\rupint{2}$ has the property that all the coefficients of
  odd degree are zero then the signs of the non-zero coefficients alternate. 
\end{lemma}
\begin{proof}
  If $f(x,y) = \sum f_i(x)y^i$ then Lemma~\ref{lem:even-alternate}
  implies that the signs of the non-zero coefficients of each $f_i$
  alternate. Writing $f(x,y) = \sum g_i(y) x^i$ shows that the signs
  of each $g_i$ alternate. Thus all the signs in the even rows and
  columns alternate, as do the signs in the odd rows and columns. 

  If $f = \sum a_{ij}x^iy^j$ then the quadralateral inequality applied
  to the coefficients $\{a_{00},a_{10},a_{01},a_{11}\}$ yie;ds
  $a_{00}a_{11}-a_{10}a_{01}<0$. Since $a_{01}=a_{10}=0$ we conclude
  that $a_{00}$ and $a_{11}$ have opposite signs. Thus, all non-zero
signs alternate.
  
\end{proof}

\begin{example}
  When $f\in\gsubpos_2$ is given by a determinant then the positivity
  of certain determinants in Theorem~\ref{thm:product-4} is a consequence of a
  classical determinantal identity. Suppose $f=|xD+yE+C|$ where $D,E$
  are diagonal matrices with positive entries, and $C$ is symmetric.
  Assume the coefficients are given the expansion in
  \eqref{eqn:product-4}. Let $\Delta=det(C)$, $\Delta_i = \det(C[i])$,
  and $\Delta_{ij}=\det(C[i,j])$ where $C[i]$ is $C$ with the i-th row
  and column deleted. Then
  \begin{align*}
    d_{00} &= \Delta \\
    d_{01} &= \sum d_i \Delta_i\\
    d_{10} &= \sum e_i \Delta_i\\
    d_{11} &= \sum_{i\ne j} d_ie_j\Delta_{ij}\\
    \smalltwodet{d_{00}}{d_{10}}{d_{01}}{d_{11}} &=
  \sum_{i\ne j} \Delta\,\Delta_{ij} d_i e_j - \sum_{i,j}
  \Delta_i\,\Delta_j\,d_ie_j \\
&= \sum_{i\ne j}
\smalltwodet{\Delta}{\Delta_i}{\Delta_j}{\Delta_{ij}} c_id_j - \sum_i
\Delta_i^2d_ie_i\le0 
  \end{align*}
since $\smalltwodet{\Delta}{\Delta_i}{\Delta_j}{\Delta_{ij}}$ is
negative by Sylvester's identity for the minors of the adjugate.
\index{Sylvester's identity}
\end{example}

Suppose that $f = \sum a_{ij}x^iy^j$ is in $\gsubplus_2(n)$, and has
all positive coefficients. We are interested in the matrices 
\begin{equation}
  \label{eqn:p2-matrix-1}
  \begin{pmatrix}
    a_{n,0} &  \dots & a_{0,0} \\
    \vdots &  \ddots & \vdots \\
    0  & \dots  & a_{0,n}
  \end{pmatrix}
\,\text{and}\,
  \begin{pmatrix}
    a_{0,n} &  \dots & 0 \\
    \vdots &  \ddots & \vdots \\
     a_{0,0}  & \dots  & a_{n,0}
  \end{pmatrix}
\end{equation}

The determinant is $\prod a_{i,n-i}$ which is positive since all the
coefficients of $f$ are positive. In addition, Theorem~\ref{thm:product-4} implies
that all $2$ by $2$ submatrices
$\smalltwodet{a_{r+1,s}}{a_{r,s}}{a_{r+1,s+1}}{a_{r,s+1}}$ are
positive. This is not terribly exciting, but there is a more
interesting property of matrices to consider:

\index{totally positive!coefficients}
\index{polynomials!totally positive coefficients}
\begin{definition}
    A polynomial $f(x,y)$ has
  \emph{totally positive coefficients} if either  matrix in \eqref{eqn:p2-matrix-1}
  is totally positive. Note that if one of the matrices is totally
  positive, so is the other. 
\end{definition}

\begin{example}
  Consider a concrete example: choose $$f=\left( 1 + x + y \right)
  \,\left( 2 + x + y \right) \,\left( 1 + x + 2\,y \right).$$ If we
  expand $f$ we get
$$2 + 5\,x + 4\,x^2 + x^3 + 7\,y + 11\,x\,y + 4\,x^2\,y + 7\,y^2 +
5\,x\,y^2 + 2\,y^3$$
and the corresponding matrix is 
$$
\begin{pmatrix}
   2 & 7 & 7 & 2 \\ 0 & 5 & 11 & 5 \\ 0 & 0 & 4 & 4 \\ 0 & 0 & 0 & 1 \\
\end{pmatrix}
$$
The four by four determinant is 40. The matrices of the
determinants of all three by three and all two by two submatrices are
non-negative.
\end{example}

We are interested in polynomials in $\gsubplus_2$ with totally
positive coefficients.  There is a special case where we can explicitly
compute some of these determinants: $f = (x+y+1)^n$.

\begin{lemma} \label{lem:multinomial-det}
  For $d=0,\dots,n$ 
$$
\det \begin{pmatrix} \binom{n}{d,0} & \dots & \binom{n}{0,0} \\
\vdots & & \vdots \\
\binom{n}{d,d} & \dots & \binom{n}{0,d}
\end{pmatrix}
= \prod_{i=0}^d \binom{n}{i}
$$
\end{lemma}
\begin{proof}\footnote{Thanks to Ira Gessel for this argument.}
The $i,j$ entry of the matrix is
$\frac{n!}{(d-j)!\,i!\,(n-i+j-d)!}$. If $D$ is the determinant, then
we can factor $n!$ from all the entries, and 
$(d-j)!$ from all the columns so that
$$D = \frac{(n!)^{d+1}}{0!\cdots d!}
\,\det\left(\frac{1}{i!(n-i+j-d)!}\right)$$
If we multiply each column 
by $1/(n-i+j-d)!$ then
$$ D = \frac{(n!)^{d+1}}{(0!\cdots d!)((n-d)!\cdots(n)!)}
\,\det\left(\binom{n-d+j}{i}\right)$$
The last determinant has entries that are polynomials in $n$. If we
set $p_i = \binom{n-d}{i}$ then the matrix is 
$\left(p_i(n+j)\right)$.  Using a  generalization of the
Vandermonde determinant \cite{kratt} we can evaluate this determinant, 
and the result is exactly the desired formula.
\end{proof}

\index{totally positive}
\index{\ TP@$TP$}

\begin{lemma} \label{lem:tp-fg}
  If $f$ has totally positive coefficients then $\frac{\partial
    f}{\partial x}$ has totally positive coefficients. If $a$ is positive then
  $(x+a)f$ has totally positive coefficients.
\end{lemma}
\begin{proof}
  The argument is best described by an example, so assume that $f$ has
  total degree $3$.  The determinants of the submatrices of the
  derivative matrix are just constant multiples of the determinants of
  the original matrix:
\begin{gather*}
f=
\begin{pmatrix}
  a_{3,0} & a_{2,0} & a_{1,0} & a_{0,0} \\
  0 & a_{2,1} & a_{1,1} & a_{0,1} \\
  0 & 0 & a_{1,2} & a_{0,2} \\
  0 & 0 & 0 & a_{0,3} 
\end{pmatrix}\\
\ \frac{\partial f}{\partial x}=
\begin{pmatrix}
   a_{2,1} & a_{1,1} & a_{0,1} \\
   0 & 2a_{1,2} & 2a_{0,2} \\
   0 & 0 & 3a_{0,3} 
\end{pmatrix} =
\begin{pmatrix}
  1&0&0\\0&2&0\\0&0&3
\end{pmatrix}\,
\begin{pmatrix}
  a_{21} & a_{11} & a_{01} \\ 0 & a_{12} & a_{02} \\ 0 & 0 & a_{03}
\end{pmatrix}
\end{gather*}
Next, assume that $f\in\gsubplus_2$ has totally positive coefficients, and let
$\alpha$ be positive. The coefficient matrices of $f$ and
$(x+\alpha)f$ are
\begin{gather*}
\begin{pmatrix}
  a_{3,0} & a_{2,0} & a_{1,0} & a_{0,0} \\
  0 & a_{2,1} & a_{1,1} & a_{0,1} \\
  0 & 0 & a_{1,2} & a_{0,2} \\
  0 & 0 & 0 & a_{0,3} 
\end{pmatrix} \\
\begin{pmatrix}
  a_{3,0} & a_{2,0} +\alpha a_{3,0}& a_{1,0} +\alpha a_{2,0} & a_{0,0} 
  +\alpha a_{1,0} & \alpha a_{0,0}\\
  0 & a_{2,1}  & a_{1,1} +\alpha a_{2,1}& a_{0,1} + \alpha a_{1,1} &
  \alpha a_{1,0}\\
  0 & 0 & a_{1,2} & a_{0,2}+\alpha a_{1,2} & \alpha a_{0,2} \\
  0 & 0 & 0 & a_{0,3} & \alpha a_{0,3}
\end{pmatrix}
\end{gather*}

We can write the coefficient matrix of $(x+a)f$ as a product
$$ 
\begin{pmatrix}
  a_{3,0} & a_{2,0} & a_{1,0} & a_{0,0} &0\\
  0 & a_{2,1} & a_{1,1} & a_{0,1} &0\\
  0 & 0 & a_{1,2} & a_{0,2} &0\\
  0 & 0 & 0 & a_{0,3} &0 
\end{pmatrix}
\ 
\begin{pmatrix}
  1&\alpha&0&0&0\\0&1&\alpha&0&0\\0&0&1&\alpha&0\\0&0&0&1&\alpha\\
  0&0&0&0&1&  
\end{pmatrix}
$$
The first matrix of the product is totally positive by assumption,
and the second is easily seen to be totally positive.  Since the
product of totally positive matrices is totally positive \cite{karlin}
it follows that the coefficient matrix of $(x+\alpha)f$ is totally
positive.
\end{proof}

\begin{cor}\label{cor:tpc}
  If $f,g\in\allpolypos$ have  positive leading coefficients then
  $f(x)g(y)$ has totally positive coefficients.
\end{cor}
\begin{proof}
  First of all, $g(y)$ has  totally positive coefficients since the matrix of
  $g(y)=\sum_0^n b_ix^i$ is
$$
\begin{pmatrix}
b_n & \dots & b_1 & b_0 \\
0 & \dots & 0 & 0 \\
\vdots & & \vdots & \vdots \\
0 & \dots & 0 & 0
\end{pmatrix}
$$
 If we write ${f(x) =
  \alpha(x+a_1)\dots(x+a_n)}$ then we can inductively apply the lemma.
\end{proof}

\index{determinants!of coefficients in $\gsubplus_2$ !not totally positive}

We conjecture (Question~\ref{ques:tnnc}) that all polynomials that are
product of linear terms in $\gsubplus_2$ have totally positive
coefficients.  It is not the case that all polynomials in
$\gsubplus_2$ are totally positive. Consider the example where
    \begin{gather*}
      M =
      \begin{pmatrix}
        1&0&0\\0&1&0\\0&0&1
      \end{pmatrix}
+
x\,\begin{pmatrix}
  13& 9& 7  \\9 & 7& 5 \\7& 5& 4
\end{pmatrix}
+
y\,
\begin{pmatrix}
  5& 7& 8 \\7& 11& 12\\8& 12& 14
\end{pmatrix}\\
\intertext{The determinant of $M$ is in $\gsubplus_2$ and equals}
 1 + 24\,x + 16\,x^2 + 2\,x^3 + 30\,y + 164\,x\,y +
62\,x^2\,y + 22\,y^2 + 64\,x\,y^2 + 4\,y^3 \\
\intertext{with coefficient matrix}
\begin{pmatrix}
 4&0&0&0\\22&64&0&0\\30&164&62&0\\ 1&24&16&2 
\end{pmatrix}
    \end{gather*}

    The determinant of the three by three matrix in the lower left
    corner is $-1760$. Of course, all the two by two submatrices have
    positive determinant.

We can't show that $f\in\gsubplus_2$ has totally positive
coefficients, but we can show that the matrix of coefficients is
TP$_2$. This is a simple consequence of the rhombus inequalities, and
reflects the geometric fact that the the graph of the coefficients is
log-concave. 

\index{rhombus inequalities}

\begin{lemma}
  If $f\in\gsubplus_2$ then the coefficient array
  \eqref{eqn:p2-matrix-1} is TP$_2$. 
\end{lemma}
\begin{proof}
  We describe a special case; the general case is no
  different. Suppose $f = \sum a_{i,j}x^iy^j$, and we wish to show
  that $\smalltwodet{a_{0,2}}{a_{2,2}}{a_{0,0}}{a_{2,0}}$ is
  non-negative. Consider a portion of the coefficient array

\centerline{\xymatrix{
a_{0,2} \ar@{-}[r] \ar@{-}[d] \ar@{-}[dr] &
a_{1,2} \ar@{-}[r] \ar@{-}[d] \ar@{-}[dr] &
a_{2,2}           \ar@{-}[d]  \\
a_{0,1} \ar@{-}[r] \ar@{-}[d] \ar@{-}[dr] &
a_{1,1} \ar@{-}[r] \ar@{-}[d] \ar@{-}[dr] &
a_{2,1}           \ar@{-}[d]  \\
a_{0,0} \ar@{-}[r]            &
a_{1,0} \ar@{-}[r]            &
a_{2,0}
}}
From the rhombus inequalities we have
\begin{align*}
  a_{0,1}\,a_{1,0} & \ge a_{0,0}\,a_{1,1}\\
  a_{1,1}\,a_{2,0} & \ge a_{1,0}\,a_{2,1}\\
  a_{0,2}\,a_{1,1} & \ge a_{0,1}\,a_{1,2}\\
  a_{1,2}\,a_{2,1} & \ge a_{1,1}\,a_{2,2}
\intertext{Multiplying these inequalities and canceling yields the
  desired result:}
  a_{0,2}\,a_{2,0} & \ge a_{0,0}\,a_{2,2}\\
\end{align*}
In general, the various vertices contribute as follows
\begin{description}
\item[corner] The corner vertices each occur once, on the correct
  sides of the inequality.
\item[edge] The vertices on the interior of the boundary edges each
  occur twice, once on each side of the inequality.
\item[interior] The interior vertices occur twice on each side of the inequality.
\end{description}
Again, multiplying and canceling finishes the proof.
\end{proof}

\section{Generic polynomials and $\gsubsep_2$}
\label{sec:generic-p2}

\index{generic polynomials in $\rupint{2}$}

In one variable the interior of the space of polynomials consists of
polynomials whose roots are all distinct.  The analog in $\gsubpos_2$ is
polynomials whose substitutions have all distinct roots. We also have
the analog of $\allpolysep$, the set of all polynomials whose roots
are at least one apart from one another.

\begin{definition}
  $\gsubgen_2$ consists of all polynomials $f(x,y)$ in $\gsubpos_2$ such that 
  \begin{enumerate}
  \item $f(x,\alpha)$ has all distinct roots, for all $\alpha\in\reals$.
  \item $f(\alpha,x)$ has all distinct roots, for all $\alpha\in\reals$.
  \end{enumerate}

  If we define $\affa_x f(x,y)=f(x+1,y)$ and $\affa_y f(x,y) =
  f(x,y+1)$ then  $\gsubsep_2$ consists of all those
  polynomials in $\rupint{2}$ satisfying
\begin{itemize}
\item $f \greateqeq \affa_x f$
\item $f \greateqeq \affa_y f$
\end{itemize}

For any fixed $\alpha$ we see that= both
  $f(x,\alpha)$ and $f(\alpha,x)$ are in $\allpolysep$, so
  $f\in\gsubgen_2$. 

\end{definition}

The graph of a polynomial in $\gsubgen_2(n)$ is easy to describe: it's
a set of $n$ non-intersecting curves. If $f$ is in $\gsubsep_2$ then
all the curves are at least one apart on every vertical line, and they
are also at least one apart on every horizontal line. In addition, 

  \begin{lemma}
\index{solution curves} 
    If $f\in\gsubgen_2$ then the  solution curves are always decreasing.
  \end{lemma}
  \begin{proof}
    If the graph of $f(x,y)=0$ had either horizontal or vertical
    points, then by Lemma~\ref{lem:sol-curves} the graph would have double
    points.  However, the solution curves are disjoint since
    $f\in\gsubgen_2$. We conclude the the slope is always the same
    sign. Since the solution curves are asymptotic to lines with
    negative slope, the solution curves are everywhere decreasing.
  \end{proof}

  \begin{example}
    We can  use polynomials in $\allpolysep$ to construct
    polynomials in $\gsubsep_2$. If $f(x)\in\allpolysep$, then
    $f(x+y)$ is easily seen to be in $\gsubsep_2$. The graph of
    $f(x+y)$ is just a collection of parallel lines of slope $1$ that
    are at least one apart in the $x$ and $y$ directions. 
  \end{example}

We can easily determine the minimum separation of the roots when
the polynomial is a quadratic given by a determinant. We start with
one variable.  Define 

\index{quadratics in $\rupint{2}$}
\begin{equation} 
 f(x) = 
\begin{vmatrix}
  a_1x+ \alpha & \gamma \\ \gamma & a_2x+\beta
\end{vmatrix}
\end{equation}

The square of the distance between the roots of a polynomial $ax^2+bx+c$ is
$(b^2-4ac)/a^2$.  The squared difference between the roots of
$f$ is
$$\frac{1}{a_1^2a_2^2} ((\alpha + \gamma)^2 -4 a_1a_2(\alpha\beta-\gamma^2)) =
\frac{1}{a_1^2a_2^2} (4a_1a_2\gamma^2 + (\alpha a_2 - \beta a_1)^2)
$$
and hence the minimum squared distance is $4\gamma^2/(a_1a_2)$. To find the minimum
of 
\begin{equation} \label{eqn:xy-det}
 f(x,y) = 
\begin{vmatrix}
  a_1x+ b_1y+\alpha & \gamma \\ \gamma & a_2x+b_2y+\beta
\end{vmatrix}
\end{equation}
we note that if $y$ is fixed then by the above the minimum distance
is $2|\gamma|/\sqrt{a_1a_2}$. If $x$ is fixed then the minimum distance is 
$2|\gamma|/\sqrt{b_1b_2}$. It follows that

\begin{lemma}\label{lem:xy-det-2}
  The polynomial  $f(x,y)$ of \eqref{eqn:xy-det} is in $\gsubsep_2$ iff
$$
   |\gamma| \ge \frac{1}{2}\max\left\{ \sqrt{a_1a_2},\sqrt{b_1b_2} \right\}
$$
  
\end{lemma}

\begin{remark}
  The condition that the solution curves of a polynomial in
  $\gsubsep_2$ are at least one apart in the horizontal and vertical
  directions does not imply that they are one apart in all directions.
  From the above we know that $f = (x+y+1/2)(x+y-1/2)\in\gsubsep_2$. The
  line $x=y$ meets the graph, which is two parallel lines, in the
  points $\pm(1/4,1/4)$. Thus, there are two points on different
  solution curves that are only $\sqrt{2}/2=.707$ apart.
\end{remark}

  \begin{example}
    The polynomial $xy-1$ is a limit of polynomials in
    $\gsubsep_2$. We know that the limit
$$\lim_{\epsilon\rightarrow0^+} \smalltwodet{x+\epsilon
  y}{1}{1}{\epsilon x+y} = xy-1$$
expresses $xy-1$ as a limit of polynomials in $\rupint{2}$. In order to
apply the lemma we need to multiply each side by
$\smalltwodet{1}{0}{0}{\epsilon^{-1/2}}$, which yields
$$\lim_{\epsilon\rightarrow0^+} \epsilon\, \smalltwodet{x+\epsilon
  y}{\epsilon^{-1}}{\epsilon^{-1}}{ x+\epsilon^{-1}y} = xy-1$$
Since $2|c_{12}| = 2 |c_{12}|/\sqrt{d_1d_2} = 2\epsilon^{-1/2}$ we see
that the roots are at least one apart when $0<\epsilon<1/4$.

  \end{example}

Polynomials in $\gsubsep_2$ satisfy the expected properties.

\begin{lemma} \ 
  \begin{enumerate}
  \item $\gsubgen_2$ and $\gsubsep_2$ are closed under differentiation.
  \item The coefficients of any power of $y$ (or of $x$) of a polynomial in
    $\gsubsep_2$ are in $\allpolysep$. 
  \end{enumerate}
  \end{lemma}
  \begin{proof}
    The statement follows since differentiation preserves the
    properties of all distinct roots, and of interlacing.  If we
    differentiate sufficiently many times, and substitute zero we can
    conclude that all coefficients are in $\allpolysep$.
  \end{proof}
  
  We can use properties of symmetric matrices to construct polynomials
  in $\gsubgen_2$. We first show \cite{parlett}*{page 124}

\begin{lemma} \label{lem:tridiagonal}
  If $M$ is a symmetric matrix whose superdiagonal entries
  are all non-zero, then the eigenvalues of $M$ are distinct.
\end{lemma}
\begin{proof}
  If we delete the first column and the last row of $\alpha I + M$,
  then the resulting matrix has all non-zero diagonal entries, since
  these entries come from the superdiagonal of $M$. Since
  the product of the diagonal entries is the product of the
  eigenvalues, the matrix is non-singular. Thus, if $M$ is $n$ by $n$,
  then $\alpha I+M$ contains a submatrix of rank $n-1$ for any choice
  of $\alpha$. If $M$ had a repeated eigenvalue, then there would be
  an $\alpha$ for which $\alpha I + M$ had rank at most $n-2$. Since
  $\alpha I+M$ has rank at least $n-1$, this contradiction establishes
  the lemma.
\end{proof}

\begin{cor} \label{cor:tridiagonal}
  Suppose $C$ is a symmetric  matrix whose superdiagonal is
  all non-zero, and $D$ is a diagonal matrix of positive diagonal
  entries. Then,
$$ |xI + y D +C| \in \gsubgen_2$$ Moreover, 
$ |xI + y D +\alpha C| \in \gsubsep_2$ for sufficiently large $\alpha$.
\end{cor}
\begin{proof}
  Let $f(x,y)=|xI+yD+C|$. The roots of $f(x,\alpha)$ are the
  eigenvalues of $\alpha D+C$, and they are all distinct by the
  lemma. If $f(\alpha,y)=0$ then 
$$ \left| \alpha D^{-1} + yI + \sqrt{D^{-1}} C
  \sqrt{D^{-1}}\right|=0$$
Thus, the roots of $f(a,y)$ are the eigenvalues of 
$\alpha D^{-1}  + \sqrt{D^{-1}} C   \sqrt{D^{-1}}$. Since 
$\sqrt{D^{-1}} C   \sqrt{D^{-1}}$ satisfies the conditions of the
lemma, the roots of $f(\alpha,y)$ are all distinct. Thus,
$f(x,y)\in\gsubgen_2$. 

 The roots of $ |xI + y D +\alpha C|$ are the roots of
$f(\frac{x}{\alpha},\frac{y}{\alpha})$, and so the second part
follows.
\end{proof}

A similar argument shows the following.

\begin{cor}
  If $D_1$ is a diagonal matrix with distinct positive entries, and
  $D_2$ is a positive definite matrix whose superdiagonal is all
  non-zero, then
\[ |I + xD_1 + yD_2| \in\gsubgen_2\cap\gsubplus_2.\]
\end{cor}

\begin{cor}
  Every polynomial in $\rupint{2}$ is the limit of polynomials in
  $\gsubgen_2$. 
\end{cor}
\begin{proof}
  By Corollary~\ref{cor:p2-det} if $f(x,y)\in\rupint{2}$ then there is a
  diagonal matrix $D$ and a symmetric matrix $C$ such that
  $f(x,y)=|xI+yD+C|$. Let $U_\epsilon$ be the matrix that is all zero
  except for $\epsilon$'s in the sub and super diagonal. The
  polynomial $|xI + yD+C+U_\epsilon|$ is in $\gsubgen_2$ for
  $|\epsilon|$ sufficiently small, and it converges to $f(x,y)$ as
  $\epsilon\rightarrow0$. 
  
\end{proof}
\begin{example}
  Here is an example of a cubic polynomial in
  $\gsubsep_2$. Let 
\begin{gather*} f(x,y) = 
\begin{vmatrix}
  x+ y & 1 & 0 \\ 1 & x+2y & 1 \\ 0 & 1 & x+3y
\end{vmatrix}
= 
-2\,x + x^3 - 4\,y + 6\,x^2\,y + 11\,x\,y^2 + 6\,y^3
\end{gather*}

The roots of $f(x,\alpha)$ are

$$
-2\,\alpha,-2\,\alpha - {\sqrt{2 + \alpha^2}},-2\,\alpha + {\sqrt{2
    + \alpha^2}}
$$
and the minimum distance between these roots is clearly
$\sqrt{2}$. If $f(\alpha,y)=0$ then the roots are 
$$
 \frac{-\alpha}{2},\frac{-2\,\alpha - {\sqrt{6 + \alpha^2}}}{3},\frac{-2\,\alpha +
  {\sqrt{6 + \alpha^2}}}{3}
$$
\noindent%
and the minimum is seen to be $1/\sqrt{2}$. Consequently, replacing
$y$ by $y/\sqrt{2}$, $x$ by $\sqrt{2}x$, and multiplying by $\sqrt{2}$
yields

$$ -4x + 4x^3 - 4y + 12 x^2y + 11xy^2 + 3y^3 \in\gsubsep_2 
$$
and is a polynomial satisfying $\delta_x = \delta_y = 1$. 
\end{example}

The hyperbola $xy-1$ is a good example of a polynomial whose graph has
no intersections. This can be generalized.

\begin{lemma}
If $D$ is a positive definite matrix with all distinct roots, and
$g\in\rupint{2}$ then $|D-gI|\in\gsubgen_2$.     
  \end{lemma}
  \begin{proof}
    If $D$ is an $n$ by $n$ matrix then the homogeneous part of
    $|D-gI|$ is $g^n$. It suffices to show  for each $\alpha$
    that $|D-g(x,\alpha)I|$ has all distinct roots. There are $a,b,c$ such
    that $g(x,\alpha) = a(x+b)(x+c)$. After translating and scaling we
    may assume that $g(x,\alpha) = x^2-1$. Now
\[
|D - g(x,\alpha)I| = |(D+I) - x^2 I|
\]
Since $D$ has all distinct eigenvalues, $|D+I-xI|$ has all distinct
positive roots. Taking square roots shows that $|D+I-x^2I|$ has all
distinct roots.
  \end{proof}

\section{Integrating generic polynomials}
\label{sec:p2-int-generic}

We can integrate generic polynomials. To do this we first need to
introduce some quantitative measures of root separation. 
We  extend the one variable measure $\delta(f)$ of separation to
polynomials in $\rupint{2}$:
\begin{align*}
  \delta_x(f) &= \inf_{\alpha\in\reals} \ \delta(f(x,\alpha)) \\
  \delta_y(f) &= \inf_{\alpha\in\reals} \ \delta(f(\alpha,y)) \\
\end{align*}
\index{solution curves} 
If the solution curves of $f$ cross, then both $\delta_y(f)$ and
$\delta_x(f)$ are zero. Thus, $f\in\gsubgen_2$ if either
$\delta_x(f)>0$ or $\delta_y(f)>0$.  For an example, the polynomial
$$ f(x,y) = (x+ay)(x+ay+b)$$
has $\delta_x(f)=b$ and $\delta_y(f)=b/a$. This shows that 
$(\delta_x(f),\delta_y(f))$ can take any pair of positive values.

Applying \eqref{eqn:sep-prop} to $f\in\rupint{2}$ yields 
$$
0\le s \le t \le
  \delta_y(f) \implies (\forall \alpha\in\reals)\ f(\alpha,y+s)\greateqeq
  f(\alpha,y+t)
$$
This shows that $f\in\gsubsep_2$ if and only if $\delta_x(f)\ge1$ and
$\delta_y(f)\ge1$.

\begin{prop}\label{prop:interlace-gen}
  Suppose that $f\in\gsubgen_2$,  $0<\alpha\le\delta_y(f)$, and
  $w(x)$ is a non-negative integrable function on $\reals$. If we
  define
$$ g(x,y) = \int_0^{\alpha} f(x,y+t)w(t)\,dt$$
then $g(x,y)\in\gsubgen_2$, and $\delta_y(g)\ge \frac{1}{2}\alpha$.
\end{prop}
\begin{proof}
  Since $g^H = \left(\int_0^\alpha w(t)\,dt\right) f^H$, we see that $g$
  satisfies  the homogeneity  condition.
  Proposition~\ref{prop:family-int} implies that $g$ satisfies substitution, and
  Proposition~\ref{prop:family-int-2} shows that the separation number is as
  claimed. Since $\delta_y(g)>0$, we know that $g\in\gsubgen_2$. 
\end{proof}

\begin{remark}
  If $f\in\gsubsep_2$ then $\int_0^1 f(x,y+t)\,dt$ is in $\gsubgen_2$,
  but even the simplest example fails to be in $\gsubsep_2$. Consider
  \begin{align*}
    f(x,y) &= (x+y)(x+y+1) \\
    \int_0^1 f(x,y+t)\,dt &= g(x,y) = 5/6 + (x+y)(x+y+2)\\
& \delta_x(g) = \delta_y(g) = \sqrt{2/3} < 1
  \end{align*}

\end{remark}

Here's a curious corollary.

\begin{cor}
  Suppose that $f\in\gsubgen_2$, and $\delta_y(f)\ge1$. If $\frac{\partial
    }{\partial y}h = f$, then $\Delta_y(h)\in\gsubgen_2$. 
\end{cor}
\begin{proof}
  Take $w=1$, and $\alpha=1$. The corollary follows from
  Proposition~\ref{prop:interlace-gen} since 
$$ \int_0^1 f(x,y+t)dt = h(x,y+1)-h(x,y) = \Delta_y(h)$$
\end{proof}

If $f\in\allpoly$, then it follows from the definitions that
$$ \delta(f) = \delta_x\,f(x+y) = \delta_y\, f(x+y) $$
Consequently, 

\begin{cor}
  If $f\in\allpoly$, $w(t)$ positive, and $0<\alpha\le \delta(f)$ then
\[ \int_0^\alpha f(x+y+t)w(t)\,dt\in\gsubgen_2.\]
\end{cor}

\section{Recursive sequences of polynomials}
\label{sec:p2-recursive}

In this section we show how to construct polynomials in $\gsubpos_2$
that are analogs of orthogonal polynomials. The next lemma (a special case of
Lemma~\ref{lem:sub-recursive}) is the key idea, for it allows us to construct a
new polynomial in $\gsubpos_2$ from a pair of interlacing ones.

\begin{lemma}
  If $f\lesslesseq g$ are both in $\gsubpos_2$, and $a,b,d$ are positive
  then
$$ (ax+by+c)f - dg \lesslesseq f$$
\end{lemma}

Two constructions are immediate consequences of the lemma.
The first one shows that certain two-dimensional analogs of
orthogonal polynomials are in $\gsubpos_2$.
\index{orthogonal polynomials!analogs in $\gsubpos_2$}

\begin{lemma}
  Let $a_n,b_n,d_n$ be positive, and define a sequence of polynomials
  by
  \begin{align*}
    P_{-1}&=0\\
    P_0 & = 1\\
    P_{n+1} &= (a_nx+b_ny+c_n)\,P_n - d_nP_{n-1}\\
\intertext{Then all $P_n$ are in $\gsubpos_2$ and }
   P_0 & \lessgreateq P_1 \lessgreateq P_2 \lessgreateq \cdots
  \end{align*}
\end{lemma}

Note that $P_n(x,\alpha)$ and $P_n(\alpha,y)$ are orthogonal
polynomials for any choice of $\alpha$.  These polynomials also have a
nice representation as determinants that shows that they are in
$\gsubpos_2$. If we let $\ell_i = a_ix+b_iy+c_i$ then we can construct
a matrix
\index{determinants!in $\gsubpos_2$}
$$ 
\begin{pmatrix}
  \ell_1 & \sqrt{d_2} & 0 & 0 & \hdots &0 & 0\\
\sqrt{d_2} & \ell_2 & \sqrt{d_3} & 0 & \hdots &0 & 0 \\
0 & \sqrt{d_3} & \ell_3 & \sqrt{d_4} & \hdots & 0 & 0\\
0 & 0 & \sqrt{d_4} & \ell_4 & \hdots & 0 & 0\\
\vdots&\vdots &\vdots &\vdots & \ddots & \vdots & \vdots\\
0 & 0 & 0 & 0 & \hdots & \sqrt{d_n} & \ell_n
\end{pmatrix}
$$
whose determinant is $P_n$.  If we set $\ell_i = c_i$ then the matrix
is called a \emph{Jacobi matrix}, and its characteristic polynomial is
an orthogonal polynomial  in one variable. \index{Jacobi matrix}

\begin{lemma}
  Let $a_n,b_n,d_n,e_n$ be positive, and define a sequence of polynomials
  by
  \begin{align*}
    P_{-1} &=0\\
    P_0 & = 1\\
    P_{n+1} &= (a_nx+b_ny+c_n)\,P_n - \left(d_n\frac{\partial}{\partial
      y} + e_n\frac{\partial}{\partial x}\right)P_n
  \end{align*}
\noindent%
  {Then all $P_n$ are in $\gsubpos_2$ and 
$\quad P_0 \lessgreateq P_1 \lessgreateq P_2 \lessgreateq \cdots$}
\end{lemma}

\section{Bezout matrices}
\label{sec:bezout-matrices}

Any two polynomials $f,g$ of the same degree determine a matrix called the
\index{Bezout matrix}Bezout matrix. We give a simple argument that the
Bezout matrix is positive definite if and only if $f$ and $g$
interlace. We also investigate properties of the Bezout matrix for polynomials in
$\rupint{2}$.

\begin{definition}
  If $f(x),g(x)$ are polynomials of the same degree $n$ the Bezout
  matrix of $f$ and $g$ is the symmetric matrix $B(f,g)=(b_{ij})$ such
  that
\[
\frac{f(x)g(y)-f(y)g(x)}{x-y} = \sum_{i,j=1}^n b_{ij}x^{i-1}y^{j-1}
\]
\end{definition}

In order to present the key result \cite{bezout} we introduce the
\index{Vandermonde vector}Vandermonde vector.
\[
V_n(x) = (1,x,\dots,x^{n-1})^t.
\]

\begin{theorem}[Lander]
 If $f,g$ are polynomials of degree $n$ with no common
zeros, and the zeros of $f$ are $r_1,\dots,r_n$ then 
$B(f,g) = (V^t)^{-1}DV^{-1}$ where
\begin{gather*}
  V = \bigl( V_n(r_1),\dots,V_n(r_n)\bigr) \\
  D = diag\bigl({g(r_1)f'(r_1)},\dots, {g(r_n)f'(r_n)}\bigr)
\end{gather*}
\end{theorem}

\begin{cor}
  If $f\in\allpoly$ and $f,g$ have no common roots then $f$ and $g$
  interlace if and only if the Bezout matrix is positive definite.
\end{cor}
\begin{proof}
If the Bezout matrix is positive definite then $g$ and $f'$ have the
same sign at the roots of $f$. Since $f'$ sign interlaces $f$, so does
$g$. The converse is similar.
\end{proof}

This is a classical result. More interesting and recent is that using
the Bezout matrix we can parametrize the solution curves of certain
polynomials in $\rupint{2}$. We paraphrase Theorem 3.1 \cite{bezout}

\begin{theorem}
 Three polynomials $f_0,f_1,f_2$ in one variable determine a map from
 the complex line $\complexes$ to the complex plane $\complexes^2$:
\begin{equation}\label{eqn:bezout-1}
z\mapsto \left(\frac{f_0(z)}{f_1(z)},\frac{f_2(z)}{f_1(z)}\right)
\end{equation}
The image is a rational curve defined by the polynomial
\begin{equation}
  \label{eqn:bezout}
  \Delta(x,y) = det( B(f_0,f_2) + x B(f_0,f_1) + y B(f_2,f_1))
\end{equation}
  
\end{theorem}

Choose mutually interlacing polynomials $f_0,f_1,f_2$ in
$\allpoly(n)$. Each of the Bezout matrices $B(f_0,f_1)$, $B(f_0,f_2)$,
$B(f_1,f_2)$ is positive definite since each pair is interlacing. It
follows that $\Delta(x,y)$ in \eqref{eqn:bezout} belongs to
$\gsubplus_2$. The solution curves are parametrized by
\eqref{eqn:bezout-1}. 

Consider an example:

\begin{align*}
&
\begin{array}{lll}
 f_1&=(2+x) (8+x) (12+x) \\
 f_2&=(3+x) (9+x) (13+x)\\
 f_3&=(7+x) (11+x) (15+x)
\end{array} &
B(f_0,f_1) &=
\left(
\begin{array}{lll}
 12600 & 2922 & 159 \\
 2922 & 785 & 47 \\
 159 & 47 & 3
\end{array}
\right)
\\[.2cm]
B(f_0,f_2) &=
 \left(
\begin{array}{lll}
 90456 & 19074 & 963 \\
 19074 & 4109 & 211 \\
 963 & 211 & 11
\end{array}
\right)
&
B(f_1,f_2) &=
\left(
\begin{array}{lll}
 89568 & 17292 & 804 \\
 17292 & 3440 & 164 \\
 804 & 164 & 8
\end{array}
\right) 
\end{align*}

\begin{multline*}
\Delta(x,y) = 
(17325 x^3+175200 y x^2+155803 x^2+242112 y^2 x+465664 y x+191775 x \\
 +73728 y^3+240384 y^2+223136 y+61425)/17325
\end{multline*}

{ All solutions to $\Delta(x,y)=0$ are of the form}
\begin{align*}
x &= \frac{(t+2) (t+8) (t+12)}{(t+3) (t+9) (t+13)}&
y & =   \frac{(t+7) (t+11) (t+15)}{(t+3) (t+9) (t+13)}
\end{align*}

The graph of $\Delta(x,y)$ is in Figure~\ref{fig:bezout}.
The three curves are parametrized as follows:

\begin{align*}
  \text{the bottom curve}\quad& t\not\in(-13,-3) \\
  \text{middle on left, top on right}\quad& t\in(-13,-9)\\
  \text{middle on right, top on left}\quad& t\in(-9,-3)
\end{align*}

\begin{figure}
  \centering
{\psset{unit=.5cm}
  \begin{pspicture}(-4,-4)(4,4)
\pscurve(-4,-0.450303)(-3.5,-0.51786)(-3.,-0.588587)(-2.5,-0.664472)
        (-2.,-0.749517)(-1.5,-0.852927)(-1.,-1.)(-0.5,-1.2751)
        (0.,-1.89583)(0.5,-2.86313)(1.,-3.94566)
\pscurve(-4,2.66865)(-3.5,2.23475)(-3.,1.80214)(-2.5,1.37164)
        (-2.,0.944875)(-1.5,0.525434)(-1.,0.122344)(-0.5,-0.240788)
        (0.,-0.84375)(0.5,-1.33518)(1.,-1.77119)(1.5,-2.20215)
        (2.,-2.63407)(2.5,-3.06727)(3.,-3.50147)(3.5,-3.93638)
\pscurve(-2.,3.11193)(-1.5,1.99286)(-1.,0.901094)(-0.5,-0.102605)
        (0.,-0.520833)(0.5,-0.704029)(1.,-0.827423)(1.5,-0.922518)
        (2.,-1.00384)(2.5,-1.07777)(3.,-1.14735)(3.5,-1.21418)
        (4.,-1.27916)
  \end{pspicture}
}
  \caption{The graph of a parametrized $\rupint{2}$}
  \label{fig:bezout}
\end{figure}

Now suppose that $f \greateq g\in\rupint{2}$. The Bezout matrix with respect to
$x$ is

\[
B_x(f,g) = \frac{f(x,y)g(z,y)-f(z,y)g(x,y)}{x-z}
\]
 
Since $f(x,\alpha)$ and $g(x,\alpha)$ interlace for all $\alpha$ it
follows that $B_x(f,g)$ is a matrix polynomial in one variable that is
positive definite for all $y$. In particular, $v^t\, B_x(f,g)\,v>0$
for all non-zero vectors $v$. However, $v^tB_x(x,y)v$ is generally not a
stable polynomial.  \index{matrix polynomial!positive definite}

  The polynomial that determines the Bezout matrix has all
  coefficients in $\allpoly$.
\index{Bezout matrix}
\begin{lemma}
  If $f\lesslesseq g$ in $\allpoly$ and 
\[ B(x,y) =
\frac{1}{x-y}\begin{vmatrix}
  f(x) & f(y) \\ g(x) & g(y)
\end{vmatrix}
\]
then all coefficients of powers of $x$ or of $y$ are in $\allpoly$.
If $f,g\in\allpolypos$ then all coefficients of $B$ are positive, and
all coefficients are in $\allpolypos$. 
\end{lemma}
\begin{proof}
  Write $g(x) = \sum a_i \frac{f}{x-r_i}$ where $a_i\ge0$. Then
  \begin{align}\label{eqn:bezout-poly}
    B(x,y) & = \frac{1}{x-y} \biggl( f(x) \sum a_i \frac{f(y)}{y-r_i} -
 f(y) \sum a_i \frac{f(x)}{x-r_i}\biggr)\\
&= \sum a_i \frac{f(x)}{x-r_i}\frac{f(y)}{y-r_i}
  \end{align}
Since $\{f(x)/(x-r_i)\}$ is mutually interlacing the result follows
from Lemma~\ref{lem:mi-2}. If $f,g\in\allpolypos$ then $f(x)/(x-r_i)$
has all positive coefficients, so $B$ has all positive coefficients. 
\end{proof}

\begin{remark}
  It is not the case that $B(x,y)\in\gsubclose_2$. If it were then
  $B(x,x)$ is in $\allpoly$, but $B(x,x) =
  \smalltwodet{f(x)}{f'(x)}{g(x)}{g'(x)}$ which is stable by
  Lemma~\ref{lem:mi-same-order}. 
\end{remark}


\chapter{Polynomials in several variables }
\label{cha:pd}

\renewcommand{\TimeStampStart}{Monday, January 07, 2008: 12:56:15}
\mytoday

\section{The substitution property}
\label{sec:substitution-defn-pd}

Just as for one variable there is a
substitution property that all our polynomials will satisfy.

\index{substitution}
\index{substitution!x-substitution}
\begin{definition}
  If $f$ is a polynomial in variables $x_1,\dots,x_d$ then $f$
  satisfies \emph{$x_i$-substitution} if for every choice of
  $a_1,\dots,a_d$ in $\reals$ the polynomial
  $$
  f(a_1,\dots,a_{i-1},x_i,a_{i+1},\dots,a_d)$$
  has all real roots.
  We say that $f$ satisfies \emph{substitution} if it satisfies
  $x_i$-substitution for $1\le i \le d$. 

  It is clumsy to write out the variables so explicitly, so we
  introduce some notation. Let
  $\xx=(x_1,\dots,x_d)$, $\aaa = (a_1,\dots,a_d)$, and define
  $$
  \xx^\aaa_i = (a_1,\dots,a_{i-1},x_i,a_{i+1},\dots,a_d)$$
  The
  \emph{degree} of a monomial $x_1^{i_1}\cdots x_d^{i_d}$ is the sum
  $i_1+\cdots+i_d$ of the degrees. The \emph{total degree} of a
  polynomial is the largest degree of a monomial. The
  \emph{$x_i$-degree} of a polynomial is the largest exponent of
  $x_i$.
\index{index set} 
\index{\ aaaI@$\diffi$}
\index{\ aaaIs@$\vert\diffi\vert$} 
An \emph{index set} $\diffi$ is a sequence $(i_1,i_2,\dots)$ of
non-negative integers.  A \emph{monomial in $\xx$} can be written
$$
\xx^\diffi = x_1^{i_1}x_2^{i_2}\cdots$$
The size of an index set is
$|\diffi|=i_1+i_2+\cdots$. The degree of $\xx^\diffi$ is $|\diffi|$.
\end{definition}

We can write a  polynomial $f$   in terms of coefficients
of monomials. If  $\xx=(x_1,\dots,x_d)$ and
$\yy=(y_1,\dots,y_e)$ then we  write
 \begin{equation} \label{eqn:index}
f(\xx,\yy)  = \sum_{\diffi} f_{\diffi}(\xx) \yy^\diffi.
\end{equation}

We can restate the definition of $x_i$-substitution:
  \begin{itemize}
   \item $f(\xx)$ satisfies $x_i$-substitution iff
    $f(\xx^\aaa_i)\in\allpoly$ for all $\aaa\in\reals^d$.
  \end{itemize}
  All the polynomials we consider will belong to this set:
$$ \xsub_d = \left\{\text{all polynomials in $d$ variables that
    satisfy substitution} \right\}$$
\index{\ Pxsub@$\xsub_d$}

In the case $d=1$ the set $\xsub_1$ is just $\allpoly$. Consider these
examples in several variables.
\begin{example} 
  All products of linear factors such as
$$ (x_1+2x_2-x_3+1)(x_1-x_2)(x_2+x_3)$$
satisfy substitution. More generally, if the factors of a product are
linear in each variable separately then the product satisfies
substitution. An example of such a polynomial is
$$
(x_1x_2+ x_1+3)(x_1x_2x_3+2x_1x_4+1)$$
\end{example}

\begin{lemma} \label{lem:sub-closure-1}
  If $f$ is a polynomial, and there are $f_j\in\xsub_d$ such that
  $\lim f_j = f$ then $f$ satisfies substitution.
\end{lemma}
\begin{proof}
  Observe that for any $\aaa\in\reals^d$ and $1\le i \le d$ we have
$$ f(\xx^\aaa_i) = \lim_{j\rightarrow\infty}   f_j(\xx^\aaa_i)$$
The result now follows from the corresponding result
(Lemma~\ref{lem:poly-closure-1}) for one variable. 
\end{proof}

\index{factorization property}
Substitution satisfies a factorization property.
\begin{lemma} \label{lem:sub-factorization}
  If $fg\in\xsub_d$ then $f,g\in\xsub_d$.
\end{lemma}
\begin{proof}
  The result is obvious for polynomials in one variable. The result
  follows from the fact that $(fg)(\xx_i^\aaa) = f(\xx_i^\aaa)\,g(\xx_i^\aaa)$.
\end{proof}

The reversal transformation is one of the few operations that preserve
the property of substitution.

\index{reverse!in $\rupint{d}$}

\begin{prop} \label{prop:sub-reverse}
  If $f$ has $x_i$-degree $k_i$ and $f=\sum a_\sdiffi
  \xx^\diffi$ then the reverse of $f$ is defined to be
$$ rev(f) = \sum a_\sdiffi \xx^{\diffk - \diffi}$$
 where $\diffk = (k_1,\dots)$. If $f\in\xsub_d$ then $rev(f)\in\xsub_d$.
\end{prop}
\begin{proof}
  Let $\aaa=(a_1,\dots,a_d)$. If we substitute we see that
  \begin{align*}
    rev(f)(\xx^\aaa_i) &= \sum a_\sdiffi\ a_1^{k_1-i_1} \cdots x_i^{k_i-i_i}
    \cdots a_d^{k_d-i_d} \\
    &= a_1^{k_1} \cdots a_d^{k_d} \sum a_\sdiffi\ a_1^{-i_1}\cdots
    x_i^{k_i-i_i} \cdots a_d^{-i_d}      \\ 
    &= a_1^{k_1} \dots a_d^{k_d} rev(f(\xx_i^\bbb))
  \end{align*}
where $\bbb = (1/a_1,\dots,1/a_d)$. Since reversal preserves roots in
one variable, $rev(f)$ satisfies substitution.
\end{proof}

Interlacing can be easily defined for polynomials satisfying
substitution. 
\index{\ zzlesslesseq@$\lesslesseq$!in $\xsub_d$}
\index{interlacing!in $\xsub_d$}
\begin{definition}
  If $f$ and $g$ are polynomials that satisfy substitution, then we
  say that $f$ and $g$ \emph{interlace} if for every real $\alpha$ the
  polynomial $f+\alpha g$ satisfies substitution. 

\end{definition}


The relation of interlacing is preserved under limits.

\begin{lemma} \label{lem:sub-closure-2}
  Suppose $f_i,g_i$ are in $\xsub_d$, and $lim_{i\rightarrow\infty}f_i
  =f$, $\lim_{i\rightarrow\infty}g_i= g$ where $f,g$ are polynomials.
  If $f_i$ and $g_i$ interlace for all $i$ then $f$ and $g$ interlace.
\end{lemma}

\begin{proof}
  Since $f_i+\alpha g_i$ is in $\xsub_d$ for all $i$, and $f_i+\alpha
  g_i$ converges to $f+\alpha g$ it follows that $f+\alpha
  g\in\xsub_d$. Consequently $f$ and $g$ interlace.
\end{proof}

\section{The class $\rupint{d}$.}
\label{sec:sub-equal-degree}

\index{degree!total}
\index{degree!of a monomial}
\index{degree!of $x_i$}

Before we can define $\rupint{d}$ we need a few definitions that
generalize the definitions for $d=2$.

\index{\ aaahomog@$f^H$}
\begin{definition}
  If $f$ has  degree $n$, then the \emph{homogeneous part} of 
  $f = \sum a_\sdiffi x^\diffi$ is
\begin{equation} \label{eqn:homog-part}
 f^H = \sum_{|\diffi|=n} a_\sdiffi \xx^\diffi
\end{equation}
\end{definition}

The homogeneous part controls the asymptotic behavior of the graph,
just as the leading coefficient does in one variable. 

\index{positivity condition}
\index{condition!positivity}
\begin{definition} \label{defn:positivity-condition}
A polynomial $f(x_1,\dots,x_d)$ of degree $n$ satisfies the \emph{positivity
  condition} if all coefficients of monomials of degree at most $n$
are positive.
\end{definition}

We  now define 

\index{\ Pgsubdpos@$\rupint{d}$}
\begin{definition} \label{defn:pd}
  \begin{align*}
    \rupint{d} &= \text{ all $f$ in $d$ variables such that} \\
    &\bullet\quad  f \text{ satisfies substitution}\\
    &\bullet\quad f^H \text{ satisfies the positivity condition
      (Definition~\ref{defn:positivity-condition}).}
  \end{align*}

$\rupint{d}$ is the analog of ``polynomials with positive leading
sign and all real roots''.

\end{definition}

\index{\ Pgsubdn@$\rupint{d}(n)$}
We let $\rupint{d}(n)$ denote the subset of $\gsub_d$ that consists of
all polynomials of  degree $n$.  


The results in Lemma~\ref{lem:sub-fH2} also hold for polynomials in $\rupint{d}$. The
homintogeneous part is always in $\rup{d}$.

\begin{lemma} \label{lem:sub-fH}
  If $f\in\rupint{d}$ then $f^H\in\gsub_d$.   If $f\in\xsub_d$ then
  $f^H\in\xsub_d$.  
\end{lemma}

\begin{proof}
  Since $f^H = (f^H)^H$ the positivity condition is satisfied.
  It remains to show that if $f\in\xsub_d$ then $f^H$ satisfies
  substitution.  Assume that the  degree of $f$ is $n$, and write
 \begin{align}
   f & = \sum_{|\diffi|=n} c_\sdiffi x^\diffi
   + \sum_{|\diffi|<n} c_\sdiffi x^\diffi. \nonumber\\
   \intertext{Replace $\xx$ by $ \xx/\epsilon$ and multiply by
     $\epsilon^n$:} \epsilon^nf(\xx/\epsilon) & = \sum_{|\diffi|=n} c_\sdiffi x^\diffi +
   \sum_{|\diffi|<n} \epsilon^{n-|\diffi|} c_\sdiffi x^\diffi. \label{eqn:sub-fH}
 \end{align}
 It follows that $\epsilon^n f(\xx/\epsilon)$ converges to $f^H$ as $\epsilon$
 approaches $0$. Lemma~\ref{lem:sub-closure-1}  shows that $f^H$ satisfies
 substitution.
\end{proof}

As is to be expected, substituting for several variables is a good
linear transformation.

\begin{lemma} \label{lem:sub-sub-a}
  Let $\xx=(x_1,\dots,x_d)$, $\yy=(y_1,\dots,y_e)$ and choose
  $\aaa\in\reals^e$. The substitution map $\xx^\diffi \yy^\diffj
  \mapsto \xx^\diffi \aaa_\sdiffj$ is a linear transformation
  $\rupint{d+e} \longrightarrow \rup{d}$.
\end{lemma}
\begin{proof}
  Using induction it suffices to prove this in the case that $e$ is
  $1$.  Let $\xx^a = (x_1,\dots,x_d,a)$.  Choose $f\in\rupint{d}$, and
  let $g=f(\xx^a)$. Since $f$ satisfies substitution so does $g$. The
  homogeneous part $f(\xx^a)^H$ is found by deleting all the terms
  involving $x_{d+1}$ in $f^H$. Since all coefficients of $f^H$ are
  positive, so are those of $f(\xx^a)^H$.  it follows that
  $g\in\rupint{d}$.
\end{proof}

Unlike $\xsub_d$, not all products of linear terms are in $\rupint{d}$.
The lemma below shows that the signs must be coherent for this to
occur. Note that there are no constraints on the signs of the constant
terms.

\begin{lemma} \label{lem:sub-prod-1}
  Suppose 
\begin{equation} \label{eqn:sub-prod-1}
f(\xx) = \prod_{k=1}^n (b_{k} + a_{1k}x_1 + \cdots +
  a_{dk}x_d).
\end{equation}
Assume that all coefficients of $x_1$ are positive.
  Then $f\in\rupint{d}$ iff for each $1\le i \le d$ all coefficients
  of $x_i$ are positive.
\end{lemma}
\begin{proof}
  The homogeneous part of $f$ is
  $$
  f^H(\xx) = \prod_{k=1}^n ( a_{1k}x_1 + \cdots + a_{dk}x_d).$$
  If
  the signs of the coefficients of $x_i$ are positive then the
  coefficients of $f^H(x_1,\dots,x_d)$ are all positive.
  
  Conversely, the terms in $f^H$ involving only variables $x_1$ and
  $x_i$ are easily found. They factor into
  $$
  \prod_{k=1}^n ( a_{1k}x_1 + a_{ik}x_i)$$
  If $f^H$ satisfies the
  positivity property then the coefficients of this polynomial all
  have the same sign.  Since the coefficients of $x_1$ are all
  positive, the coefficients of $x_i$ are positive. This is
  Lemma~\ref{lem:all-neg-converse}.
\end{proof}

\index{\ zzlesslesseq@$\lesslesseq$!in $\rupint{d}$} 
\index{interlacing!in $\rupint{d}$}
\begin{definition}
  If $f,g\in\rupint{d}$ then $f$ and $g$ \emph{interlace} iff $f+\alpha g$ is in
  $\pm\rupint{d}$ for all $\alpha$. If in addition the degree of $f$ is one
  more than the degree of $g$ then we say $f\lesslesseq g$. In this
  case $(f+g)^H = f^H$, and so the only condition we need to verify is
  substitution. This leads to an equivalent  definition in terms of
  substitutions:
  \begin{quote}
      If $f,g\in\rupint{d}$ then 
  $f \lesslesseq g$  iff $f(\xx^\aaa_i) \lesslesseq g(\xx^\aaa_i)$ for 
  all $\aaa\in\reals^d$ and $1\le i \le d$.
  \end{quote}
\end{definition}

The following important result is an immediate consequence of the
definition. See the proof of Theorem~\ref{thm:only-roots}.

\begin{theorem} \label{thm:pdpd}
  If $T\colon{}\rupint{d}\longrightarrow \gsub_d$ is a linear transformation
  then $T$ preserves interlacing.
\end{theorem}

\begin{lemma} \label{lem:sub-LC-1}
  If $f,g\in\rupint{d}$ satisfy $f\lesslesseq g$ then $f^H\lesslesseq g^H$.
\end{lemma}
\begin{proof}
  The transformation $f\mapsto f^H$ maps $\rupint{d}$ to itself.
\end{proof}

Both $\xsub_d$ and $\rupint{d}$ satisfies a loose form of closure. We
will study the general question of closure in the next chapter.

\begin{lemma}
  Assume $f$ is a polynomial in $d$ variables that satisfies the
  homogeneity condition.  If there are polynomials $f_i\in\rupint{d}$
  such that $\lim f_i = f$ then $f$ is in $\rupint{d}$.
\end{lemma}
\begin{proof}
  Applying Lemma~\ref{lem:sub-closure-1} shows that $f$ satisfies substitution.
  By hypothesis $f$ satisfies the positivity condition, and hence
  $f\in\rupint{d}$.
\end{proof}

\subsection{Quadratics in $\rupint{d}$}

There is an effective criterion to determine if a quadratic polynomial
is in $\rupint{d}$. If we write $f = \sum
a_\sdiffi \xx^\sdiffi$ then for all values of $x_2,\dots,x_d$ the
polynomial $f$ must have all real roots as a quadratic in $x_1$. This
is true if and only if the discriminant is non-negative. If we write
out the discriminant then we have a quadratic form that must be
positive semi-definite. 

For simplicity we only consider the case $d=3$, but the general case
presents no difficulties. Suppose that

\begin{gather*}
  f(x,y,z) = 
a_{000} + z\,a_{001} + z^2\,a_{002} + y\,a_{010} + y\,z\,a_{011}
+ y^2\,a_{020} + \\  
 x\,a_{100} +   x\,z\,a_{101} +     x\,y\,a_{110} +
x^2\,a_{200} \\
\intertext{then the discriminant is}
{a_{100}}^2 + 2\,z\,a_{100}\,a_{101} + z^2\,{a_{101}}^2 +
2\,y\,a_{100}\,a_{110} + 2\,y\,z\,a_{101}\,a_{110} \\ +
y^2\,{a_{110}}^2  - 4\,a_{000}\,a_{200} -
4\,z\,a_{001}\,a_{200} - 4\,z^2\,a_{002}\,a_{200}\\ -
4\,y\,a_{010}\,a_{200}  - 4\,y\,z\,a_{011}\,a_{200} -
4\,y^2\,a_{020}\,a_{200}
\end{gather*}

Now we can homogenize this, and consider it as a quadratic
$q(y,z,w)$. Since the discriminant is non-negative, the quadratic must
be non-negative for all values of $y,z,w$. Consequently, the matrix
$Q$ of the quadratic form $q$

\begin{equation}\label{eqn:quad-in-pd}
\begin{pmatrix}
  {a_{110}}^2 - 4\,a_{020}\,a_{200} & a_{101}\,a_{110} -
  2\,a_{011}\,a_{200} & 
  a_{100}\,a_{110} - 2\,a_{010}\,a_{200} \\
   a_{101}\,a_{110} - 2\,a_{011}\,a_{200} &
  {a_{101}}^2 - 4\,a_{002}\,a_{200} & a_{100}\,a_{101} -
  2\,a_{001}\,a_{200}\\ 
  a_{100}\,a_{110} - 2\,a_{010}\,a_{200}, & a_{100}\,a_{101}
  - 2\,a_{001}\,a_{200} & {a_{100}}^2 - 4\,a_{000}\,a_{200}
\end{pmatrix}
\end{equation}
is positive semi-definite.  Conversely, if $Q$ is positive
semi-definite then $f(x,y,z)$ satisfies $x$-substitution. Note that
the diagonal elements of $Q$ are positive if and only if
$f(x,y,0),f(x,0,z),f(0,y,z)$ each have all real roots. We summarize
for $d=3$:

\begin{lemma}
  A quadratic polynomial is in $\rupint{3}$ if and only if it satisfies
  the positivity conditions and the matrix  \eqref{eqn:quad-in-pd} is
  positive semi-definite.
\end{lemma}

We define $\gsubsep_d$ to be all those polynomials in $\rupint{d}$ such
that all substitutions are in $\allpolysep$. 
It is easy to construct quadratic polynomials in $\gsubsep_d$. 

\begin{lemma}
  The polynomial \[
  (a_1x_1+\cdots+a_dx_d)(b_1x_1+\cdots+b_dx_d)-c^2\]

  is in
  $\gsubsep_d$ if $\displaystyle \quad |c| \ge \dfrac{1}{2}\max_{1\le i < d}
  \sqrt{a_ib_i}$
\end{lemma}
\begin{proof}
  Notice that Lemma~\ref{lem:xy-det-2} does not depend on constants,
  so if we substitute for all but the $i$-th variable, the minimum
  distance between roots is $2|c|/\sqrt{a_ib_i}$.
\end{proof}

\section{Properties of   $\rupint{d}$}
\label{sec:sub-derivatives} 

In this section we establish the important fact that $\rupint{d}$ is
closed under differentiation. An important consequence of this will be 
results about the coefficients of polynomials in $\rupint{d}$.

As in $\rupint{2}$ we only need to know that substitution holds in one variable to
conclude that substitution holds in all variables. This is a very
useful result.

\begin{theorem} \label{thm:sub-all-sub}
  Suppose that $f$ is a polynomial in $d$ variables, and satisfies
  (see Definition~\ref{defn:pd}) 
  \begin{itemize}
  \item the positivity condition
  \item $x_i$-substitution for some $x_i$
  \end{itemize}
then $f$ satisfies substitution, and is in $\rupint{d}$.
\end{theorem}
\begin{proof}
  It suffices to show that if $f$ satisfies $x_1$-substitution then it
  satisfies $x_2$-substitution. For any $\aaa=(a_1,\dots,a_d)$ let
  $g(x_1,x_2) = f(x_1,x_2,a_3,\dots,a_d)$. Since $g$ satisfies
  $x_1$-substitution and $g^H$ is a subset of $f^H$, we can apply
  Theorem~\ref{thm:sub-xy} to conclude that $g$ satisfies
  $x_2$-substitution.
\end{proof}

\begin{theorem} \label{thm:pd-T}
  Suppose that $T\colon{}\rupint{d}\longrightarrow \gsub_d$ is a linear
  transformation that preserves degree, and for any monomial
  $\xx^\diffi$ the coefficients of $T(\xx^\diffi)^H$ are all
  positive. If $\yy=(y_1,\dots,y_e)$ then the linear transformation
$$T_\ast(\xx^\diffi \yy^\diffj) = T(\xx^\diffi)\,\yy^\diffj$$
defines a linear transformation from $\rupint{d+e}$ to itself.
\end{theorem}

\begin{proof}
  The assumptions on the image of a monomial imply that $T_\ast(f)$
  satisfies the  positivity condition for
  $f\in\rupint{d+e}$. By Theorem~\ref{thm:sub-all-sub} it suffices to show that
  $T_\ast(f)$ satisfies $x_1$-substitution.  If we choose
  $\bbb\in\reals^e$ then
  $$
  (T_\ast\,f)(\xx,\bbb) = T(f(\xx,\bbb))$$
  Since
  $f(\xx,\bbb)\in\rupint{d}$ by Lemma~\ref{lem:sub-sub-a} we know
  $T(f(\xx,\bbb))$ is in $\rupint{d}$, and so satisfies
  $x_1$-substitution.
\end{proof}

An important special case is when $d$ is $1$.
\begin{cor} \label{cor:pd-T}
  Suppose $T\colon{}\allpoly\longrightarrow\allpoly$ is a linear
  transformation that preserves degree, and maps polynomials with
  positive leading coefficient to polynomials with positive leading
  coefficient.  We define a linear transformation on $\rupint{d}$ by
  $$
  T_i(x_1^{e_1}\cdots x_d^{e_d}) = x_1^{e_1}\cdots
  x_{i-1}^{e_{i-1}}\,T( x_i^{e_i})\, x_{i+1}^{e_{i+1}}\cdots  x_d^{e_d}
  $$
  $T_i$ maps $\rupint{d}$ to itself.
\end{cor}

\begin{theorem} \label{thm:gsub-diff}
$\rupint{d}$ is closed under differentiation. If $f\in\gsub_d$ then
$f\lesslesseq \frac{\partial f}{\partial x_i} $ for $1\le i \le d$.
\end{theorem}

\begin{proof}
The linear transformation $T$ of $\allpoly$ defined by $Tg = g+
ag^\prime$ maps $\allpoly$ to itself for any $a$. Applying  Theorem~\ref{thm:pd-T}
shows that $f + a\frac{\partial f}{\partial x_i}$ is in $\rupint{d}$ for
all $a$, and the conclusion follows.

\end{proof}

\begin{cor} \label{cor:pd-coef}
  If $f=\sum f_\diffi(\xx)\yy^\diffi$ is in $\rupint{d+e}$ and $\diffj$
  is any index set then $f_\diffj(\xx)$ is in $\rupint{d}(|\diffj|)$.
\end{cor}
\begin{proof}
  If $\diffj=(j_1,\dots,j_e)$ then we can apply Theorem~\ref{thm:gsub-diff} since
  $$
  f_{\diffj}(\xx) = \frac{1}{j_1!j_2!\cdots j_e!}
  \left(\frac{\partial^{j_1}}{\partial
      y_1}\cdots\frac{\partial^{j_e}}{\partial y_e}\right)f \biggm|_{\yy=0}
$$
\end{proof}

In $\rupint{2}$ we know that consecutive coefficients interlace. The
same is true in $\rupint{d}$, and the proof is the same.

\begin{lemma} \label{lem:pd-consecutive}
  Suppose that $f(\xx,\yy)\in\rupint{d+e}$, and we write
$$ f(\xx,\yy) = \cdots+f_\diffi(\xx)\yy^\diffi + f_{\diffi\cup
  y_1}(\xx)y_1\yy^\diffi+\cdots$$ 
then $f_\diffi \longleftarrow f_{\diffi\cup y_1}$.
\end{lemma}

\begin{theorem}
  If $f,g,h\in\rupint{d}$ and  $f\lesslesseq g,$  $f\lesslesseq h$ then
\begin{itemize}
\item $g+h\in\rupint{d}$
\item $f \lesslesseq g+h$
\item For any $f,g\in\rupint{d}$ then $fg\in\gsub_d$.
\item If $g\lesslesseq k$ then $f-k\in\rupint{d}$.
\end{itemize}
\end{theorem}
\begin{proof}
  As we saw for two variables, Lemma~\ref{lem:sum-1} shows that $g+h$ satisfies
  substitution. Since $(g+h)^H =
  g^H + h^H$ it follows that $g+h$ satisfies  and positivity
  condition.
  
  The product $fg$ certainly satisfies substitution. Since $(fg)^H =
  f^H\,g^H$ we see that the positivity condition is satisfied.  The
  last one follows from the corresponding one variable result
  Lemma~\ref{lem:add-interlace}.
\end{proof}

Two by two determinants of coefficients are non-positive. This is an
immediate consequence of Theorem~\ref{thm:product-4}.

\begin{cor}\label{cor:product-4d}
  Suppose that $f(\xx,y,z) = \sum f_{ij}(\xx) y^iz^j$ is in
  $\rupint{d+2}$. Then
$$ 
\begin{vmatrix}
  f_{i,j}(\xx) &   f_{i+1,j}(\xx)\\
  f_{i,j+1}(\xx) &   f_{i+1,j+1}(\xx)
\end{vmatrix} \le 0
$$
\end{cor}

\index{determinants!involving $f,f_x,f_y$}
\index{Taylor series}
If we apply this corollary to the Taylor series for
$f(x_1+y,x_2+z,x_3,\dots,x_d)$ we get

\begin{cor} \label{cor:product-4d-1}
  If $f(\xx)\in\rupint{d}$, $d\ge2$, then 
  $\begin{vmatrix}
    f & \frac{\partial }{\partial x_1 } f \\
\frac{\partial }{\partial x_2 } f & \frac{\partial }{\partial x_1 } 
\frac{\partial }{\partial x_2 } f 
  \end{vmatrix} \le0$.
\end{cor}
 
\begin{remark}
  If we are given a polynomial $f$ that satisfies the above inequality
  for all pairs of distinct variables, then it does not imply that
  $f\in\rupint{d}$. However, this is true for \index{multiaffine
    polynomial}multiaffine polynomials \seepage{sec:mult-polyn}.
\end{remark}

\section{The analog of $\allpolypos$}
\label{sec:gsubpm}

$\allpolypos$ is the set of all polynomials in $\allpoly$ with all
positive coefficients.  The analog $\gsubplus_d$ of $\allpolypos$ is
also defined in terms of the signs of the coefficients. $\gsubplus_d$
is closed under interlacing, and the homogeneous part of any
polynomial in $\rupint{d}$ is always in $\gsubplus_d$.

\index{\ Pgsubplus@$\gsubplus_d$}

\begin{definition}
  $\gsubplus_d(n)$ is the collection of all polynomials
  $f\in\rupint{d}(n)$ such that all coefficients of monomials of
  degree at most $n$ are positive.
  
  That is, if $f = \sum \aaa_\sdiffi \xx^\diffi$ has total degree $n$,
  then $\aaa_\sdiffi$ is positive for all $|\diffi|\le n$.
\end{definition}  
We can easily characterize products of linear terms that are in
$\gsubplus_d$. It follows from Lemma~\ref{lem:sub-prod-1} that a product
\eqref{eqn:sub-prod-1} is in $\gsubplus_d$ iff all $a_{ik}$ are
positive, and all $b_k$ are positive.

\begin{lemma} \label{lem:sub-a}\ 

  \begin{enumerate}
  \item   $\gsubplus_d$ is closed under differentiation.
  \item   If $f\in\gsubplus_d$ then for any
  positive  $a$ we have $f(a,x_2,\dots,x_d)\in\gsubplus_{d-1}$.
  \item If $f\in\gsubplus_{d+e}$ then all coefficients $f_\sdiffi$ given
  in \eqref{eqn:index} are in $\gsubplus_d$.
\end{enumerate}
\end{lemma}
\begin{proof}
  We only need to check positivity for the first statement, and this
  is immediate.  The second statement is obvious, but it is false if
  we substitute a negative value.\footnote{ $f=(1+x_1+x_2)(3+x_1+x_2)$
    is in $\gsubplus_2$, but $f(-2,x_2) = x^2_2-1$ is not in
    $\allpolypm$.  }  Finally, differentiate sufficiently many
  times until $f_\sdiffi$ is the constant term, and then substitute
  $0$ for $\yy$.
\end{proof}


We can determine if a polynomial is in $\gsubplus_d$ just by looking at
the constant terms with respect to each variable.

\begin{theorem} \label{thm:foxo}
  Suppose that $f\in\rupint{d}$ and for each $1\le i\le d$ the polynomial
  \\ $f(0,\dots,0,x_i,0,\dots,0)$ is in $\allpolypos$.  Then,
  $f$  is in $\pm\gsubplus_d$.
\end{theorem}
\begin{proof}
 We may assume that
 $$f(0,\dots,0,x_i,0,\dots,0)$$
 has all positive, or all negative,
 signs.  The proof now proceeds by induction on the number of
 variables. For $d=2$ the hypothesis implies that in $$
 f(x_1,x_2) =
 f_0(x_1) + \cdots + f_n(x_1)x^n_2$$
 the coefficient $f_0$ is in
 $\allpolypos$. We may assume that all coefficients of $f_0$ are
 positive. Since $f_0 \lesslesseq f_1$ we see that all coefficients of
 $f_1$ have the same sign. Continuing, all the coefficients of $f_i$
 have sign $\epsilon_i$. Now the hypothesis also implies that
 $f(0,x_2)$ is in $\allpolypos$. Since $f(0,x_2)$ contains the
 constant terms of all $f_i$ it follows that all $\epsilon_i$ have the
 same sign.  Thus all coefficients of all the $f_i$ are positive.
  
  In general, if $f\in\rupint{d}$ then by induction we know that all
  substitutions
  $$f(x_1,\dots,x_{i-1},0,x_{i+1},\dots,x_d)$$
  are in $\gsubplus_{d-1}$.
  We may assume that $f(0,x_2,\dots,x_d)$ has all positive terms, and
  then use overlapping polynomials as above to conclude that all terms
  are positive. 
\end{proof}

$\gsubplus_d$ is closed under interlacing.

\begin{theorem} 
  If $f\in\gsubplus_d$, $g\in\rupint{d}$ and $f\lesslesseq g$ then $g$
  is in $\gsubplus_d$.
\end{theorem}
\begin{proof} 
  The proof is by induction on the number of variables. For
  $\allpolypos$ the result is true.  Write
\begin{align*}
f(\xx) &= f_0(x_2,\dots,x_d) + f_1(x_2,\dots,x_d)x_1 + \cdots
  +f_n(x_2,\dots,x_d)x_1^n.\\
g(\xx) &= g_0(x_2,\dots,x_d) + g_1(x_2,\dots,x_d)x_1 + \cdots
  +g_{n-1}(x_2,\dots,x_d)x_1^{n-1}.
\end{align*}
Since $f_i \lesslesseq g_i$ for all $i$, and
  all $g_i$ are in $\gsubplus_{d-1}$ by induction, it follows that the
  coefficients of $g_i$ are either all positive, or all negative.
  Since $g\in\rupint{d}$ the homogeneous part $g^H$ has all positive
  coefficients, and so each $g_i$ has some positive coefficient. This
  shows that $g\in\gsubplus_d$.
\end{proof}

\index{positivity!of substitutions} The positivity of substitutions
implies that a polynomial in $\rupint{d}$ is in $\gsubplus_d$.
Note that if $f(x)$ is positive for all positive $x$ then it is not
necessarily the case that all coefficients are positive: consider
$1+(x-1)^2$. 

\begin{lemma}
  If $f(\xx)\in\rupint{d}$ and $f(\aaa)$ is positive for all
  $\aaa\in(\reals^+)^d$ then $f\in\gsubplus_d$. 
\end{lemma}
\begin{proof}
  We prove this by induction on $d$. If $d=1$ and $f(x)$ has positive
  leading coefficient then we can write $f =
  \alpha\prod(x+r_i)$. Since this is positive for all positive $x$ it
  follows that no $r_i$ can be negative, and therefore all
  coefficients are positive.

  In general, write $f(\xx,y)\in\rupint{d+1}$ as $\sum
  f_i(\xx)y^i$. Now for any $\aaa\in(\reals^+)^d$ the hypothesis
  implies that $f(\aaa,y)$ is positive for all positive $y$. Thus all
  coefficients $f_i(\aaa)$  are positive, and so by induction all
  coefficients of $f(\xx,y)$ are positive.
\end{proof}

The next two results show the relation between homogeneous polynomials and 
the positivity condition.

\begin{lemma}
  If $f$ is a homogeneous polynomial then \\ $f\in\rupint{d}$ iff \\ 
$f(1,x_2,\dots,x_d)\in\gsubplus_{d-1}$.
\end{lemma}
\begin{proof}
  The assumption that $f$ is homogeneous is equivalent to the equality 
  $f=f^H$. The result now follows from the definition of $\gsubplus_d$. 
\end{proof}

\begin{lemma}
  If $f\in\rupint{d}(n)$ let $F$ be the corresponding homogeneous
  polynomial.  We construct $F$ by replacing each monomial
  $\xx^\diffi$ by $\xx^\diffi x_{d+1}^{n-|\diffi|}$. A necessary and
  sufficient condition that $F\in\rupint{d+1}$ is that $f\in\gsubplus_d$.
\end{lemma}
\begin{proof}  All coefficients of $F^H$ are
  positive.  If $F\in\rupint{d+1}$ then $F^H = F$. Since $f$ is found
  by substituting $1$ for $x_{d+1}$ in F, all coefficients of $f$ must
  be positive. Conversely, since the coefficients of $F$ are the same
  as the coefficients of $f$ we see that the assumption on $f$ implies
  that $F$ satisfies positivity.
\end{proof}

We know that the terms of maximal total degree of a polynomial in
$\rupint{d}$ constitute a polynomial that is also in $\gsub_d$. If
$f\in\gsubplus_d$ then there is a stronger result:

\begin{lemma} \label{lem:sub-total-degree}
  Suppose $f\in\gsubplus_d(n)$ and $h_r$ is the polynomial consisting
  of all terms of total degree $r$. Then $h_r\in\gsubplus_d(r)$ for
  $1\le r \le n$.
\end{lemma}
\begin{proof}
  If we homogenize then 
$$F(\xx,z) = g_n + zg_{n-1} + \dots + z^ng_0$$
Since  $g_r$ is the coefficient of $z^{n-r}$ the conclusion follows.
\end{proof}

Here is a condition for interlacing in $\gsubplus_d$ that only
involves positive linear combinations. 

\begin{lemma} \label{lem:pdplus-int}
  If $f,g\in\gsubplus_d$ and $deg(f)>deg(g)$ then $f\lesslesseq g$ if
  and only if these two conditions hold
  \begin{align*}
    (1)\, & f + \alpha g\in\gsubplus_d \text{ for all $\alpha>0$}\\
    (2)\, & f + \alpha g x_i\in\gsubplus_d \text{ for all $\alpha>0$ and
      some $i$}\\
  \end{align*}
\end{lemma}
\begin{proof}
  Assume that $f \lesslesseq g$.
  Both left hand sides are have homogeneous part equal to $f^H$ plus
  perhaps a contribution for $g^H$, so they both satisfy the
   homogeneity conditions. We need to verify substitution. If we
  consider substituting for all but the $x_i$ variable, and restrict
  $\aaa$ to be non-negative, then both $f(\xx^\aaa)$ and $g(\xx^\aaa)$
  are in $\allpolypos$. The conclusion now follows from
  Lemma~\ref{lem:ispm-2}.

  Conversely, if $(1)$ and $(2)$ are satisfied then we can again use
  Lemma~\ref{lem:ispm-2} to see that $f\lesslesseq g$.
\end{proof}

\index{polar derivative!in $d$ variables}
Just as for one variable we can define a polar derivative. If
$f\in\gsubplus_d$ then write $f = \sum a_\sdiffi \xx^\diffi$. The
corresponding homogeneous polynomial $F = \sum a_\sdiffi \xx^\diffi
y^{n-|\diffi|}$ is in $\gsubpos_{d+1}$, and if we differentiate with
respect to $y$ and then set $y=1$ we get the polar derivative
\begin{align}
  \sum (n-|\diffi|)a_\sdiffi \xx^\diffi & \label{eqn:polar-pd} \\
\intertext{which can also be written in the more familiar fashion}
 n f - x_1\frac{\partial f}{\partial x_1} - \dots - x_n\frac{\partial
  f}{\partial x_n} & \label{eqn:polar-pd-2}
\end{align}
We will return to the polar derivative in Section~\ref{sec:xdxydy}.

Rather than forming the homogeneous polynomial with respect to all the
variables, we can  homogenize some of them.

\index{homogenizing!in $d$ variables}

\begin{lemma}
  If $f(\xx,\yy) = \sum a_\sdiffi(\xx)
  \yy^\sdiffi\in\gsubplus_{d+e}(n)$ then  \\
$g(\xx,\yy,z) = \sum a_\sdiffi(\xx)
  \yy^\sdiffi z^{n-|\sdiffi|} \in\gsubpos_{d+e+1}(n)$
\end{lemma}
\begin{proof}
  $g^H=g$, so $g$ satisfies 
  homogeneity condition. Since 
$$ g(\xx,\yy,z) = z^n f(\xx,\frac{\yy}{z})$$
it follows that $g$ satisfies substitution.
\end{proof}

We have a more general form of polar derivative:

\begin{cor}
  If $f(\xx,\yy) = \sum a_\sdiffi(\xx)
  \yy^\sdiffi\in\gsubplus_{d+e}(n)$ then  \\ 
$$ f(\xx,\yy) \lesslesseq \sum (n-|\diffi|)a_\sdiffi(\xx)\yy^\sdiffi.$$
\end{cor}

\section[Coefficient inequalities, hives, and Horn]%
{Coefficient inequalities for $\gsubplus_2$, hives, and the  Horn problem}

\label{sec:hives and horns}

Every polynomial in $\gsubplus_2$ determines a homogeneous polynomial
in $\rupint{3}$. If $f(x,y,z)\in\gsub_3$ is such a homogeneous polynomial
then the three polynomials $f(x,1,0)$, $f(0,x,1)$, $f(1,0,x)$ have all
real roots. We call these the \emph{boundary polynomials} of $f$. We
consider the converse:

\index{boundary polynomial problem}

\begin{center}\textbf{The boundary polynomial problem:}\end{center}
\begin{quote}
 Given three polynomials $f_1(x),f_2(x),f_3(x)$ with all real roots,
  when is there a polynomial
  $f(x,y,z)\in\rupint{3}$ whose boundary polynomials are $f_1,f_2,f_3$?
\end{quote}

Surprisingly, this is related to a famous problem involving
matrices.

\index{additive Horn problem}
\index{Horn problem!additive}

\begin{center}
  \textbf{The additive Horn problem:} 
\end{center}
\begin{quote}
  Fix an integer $n$, and let
  $\alpha,\beta,\gamma\in\reals^n$. Is there a triple $(A,B,C)$ of
  Hermitian matrices with $A+B+C=0$, and eigenvalues
  $\alpha,\beta,\gamma$? If so, we say \emph{Horn's additive problem
    is solvable for $\alpha,\beta,\gamma$.}\cite{speyer}
\end{quote}

We begin by looking at a particular example, and then establishing
inequalities for coefficients of polynomials in $\gsubplus_2$.
Taking logs of these inequalities gives a structure called a
\emph{hive}. We then recall the connection between hives and Horn's problem
to get information about the boundary polynomial problem.

\index{hive}
\index{Horn's conjecture}

Suppose that we are given a polynomial  of
degree $n$ in $\gsubplus_3$ that is homogeneous:
\begin{gather*}
  f(x,y,z) = (z + x + y) (4z + 2 x + y) (2z + 3 x + y) (3z + 5 x + y) \\
  =24z^4 + 112\,xz^3 + 186\,x^2z^2 + 128\,x^3z + 30\,x^4 + 50\,yz^3 + 175\,xyz^2 + 
  190\,x^2yz + \\ 61\,x^3y +  35\,y^2z^2 + 80\,xy^2z + 41\,x^2y^2 +
  10\,y^3z + 11\,xy^3 + y^4.
  \end{gather*}
 We can represent the terms of $f(x,y)$ in a triangular diagram.

\centerline{ \xymatrix@=.4cm{
    y^4 \ar@{-}[d]_{} \ar@{-}[dr] \\
    10\, y^3z \ar@{-}[dr] \ar@{-}[d]_{} \ar@{-}[r] & 11\, xy^3 \ar@{-}[d] \ar@{-}[dr]^{} \\
    35\, y^2z^2 \ar@{-}[dr] \ar@{-}[d]_{}\ar@{-}[r] & 80\,xy^2z
    \ar@{-}[dr] \ar@{-}[d] \ar@{-}[r] & 41\,  x^2y^2 \ar@{-}[d]
    \ar@{-}[dr]^{} \\ 
    50\, yz^3 \ar@{-}[d]_{} \ar@{-}[dr] \ar@{-}[r] & 175 \,xyz^2 \ar@{-}[dr]
    \ar@{-}[d] \ar@{-}[r] &
    190\,x^2yz \ar@{-}[dr] \ar@{-}[d] \ar@{-}[r]&  61\, x^3y  \ar@{-}[d] \ar@{-}[dr]^{} \\
    24\,z^4 \ar@{-}[r]_{} \ar@{-}[r] & 112\, xz^3 \ar@{-}[r]_{} & 186\,x^2z^2
    \ar@{-}[r]_{} & 128\, x^3z \ar@{-}[r]_{} & 30\,x^4 }}

The three boundary polynomials of $f$ occur on the boundary of this
diagram, and are 
\begin{align*}
  f(x,1,0) &= 30x^4+61x^3+41x^2+11x+1 \\
  f(0,x,1) &= x^4+10x^3+35x^2+50x+24 \\
  f(1,0,x) & = 24x^4+112x^3+186x^2+128x+30
\end{align*}

If we only list the coefficients then we get the triangular diagram 

\centerline{ \xymatrix@=.4cm{
    1 \ar@{-}[d]_{} \ar@{-}[dr] \\
    10\,  \ar@{-}[dr] \ar@{-}[d]_{} \ar@{-}[r] & 11\,  \ar@{-}[d] \ar@{-}[dr]^{} \\
    35 \ar@{-}[dr] \ar@{-}[d]_{}\ar@{-}[r] & 80
    \ar@{-}[dr] \ar@{-}[d] \ar@{-}[r] & 41 \ar@{-}[d]
    \ar@{-}[dr]^{} \\ 
    50 \ar@{-}[d]_{} \ar@{-}[dr] \ar@{-}[r] & 175  \ar@{-}[dr]
    \ar@{-}[d] \ar@{-}[r] &
    190 \ar@{-}[dr] \ar@{-}[d] \ar@{-}[r]&  61  \ar@{-}[d] \ar@{-}[dr]^{} \\
    24\, \ar@{-}[r]_{} \ar@{-}[r] & 112 \ar@{-}[r]_{} & 186
    \ar@{-}[r]_{} & 128 \ar@{-}[r]_{} & 30 }}

\vspace*{.1cm}

This array of numbers has the property that the coefficients in each of
the three kinds of rhombus satisfy $bc \ge ad$:

\centerline{ \xymatrix@=.4cm{
b \ar@{-}[r] \ar@{-}[d] \ar@{-}[dr]& d \ar@{-}[d] & & 
a \ar@{-}[r]\ar@{-}[dr] &  b \ar@{-}[d] \ar@{-}[dr] & & &
a \ar@{-}[d] \ar@{-}[dr] \\
a \ar@{-}[r] & c & & 
& c \ar@{-}[r] & d & &
b \ar@{-}[r] \ar@{-}[dr] & c \ar@{-}[d] & & &\\
&&&&&&&&d
}}

\index{quadrilateral inequalities}
\index{rhombus inequalities}

\begin{prop}\label{prop:p2plus-inequality}
  For any homogeneous polynomial in $\gsubplus_3$, construct the
  triangular array of coefficients. In every rhombus formed by two
  adjacent triangles the product of the two numbers on the common edge
  is at least the product of the two remaining vertices.

\end{prop}

\begin{proof}
  The inequalities corresponding to the three different kinds of
  rhombus follow by applying the inequalities for polynomials in
  $\rupint{2}$ (Theorem~\ref{thm:product-4}) to $f(x,y,1)$, $f(x,1,z)$,
  and $f(1,y,z)$.
\end{proof}

Suppose we choose $f(x,y)=\sum a_{ij}x^iy^j$ in $\gsubplus_2$. Since
all coefficients are positive, we can define real numbers
$h_{ijk}=\log(a_{ij})$ where $i+j+k=n$. If we write the inequalities
of Proposition~\ref{prop:p2plus-inequality} in terms of the $h_{ij}$ 

\begin{gather*}
  a_{i+1,j}\, a_{i,j+1}  \ge a_{i,j}\, a_{i+1,j+1} \\
  a_{i+1,j-1}\, a_{i,j-1} \ge a_{i,j}\,  a_{i+1,j-2}\\
  a_{i-1,j+1}\, a_{i-1,j} \ge a_{i,j}\,  a_{i-2,j+1}
\end{gather*}

then they translate into inequalities about the $h_{ijk}$'s:

\begin{gather*}
  h_{i+1,j,k-1} + h_{i,j+1,k-1} \ge h_{i,j,k} + h_{i+1,j+1,k-2}\\
  h_{i+1,j-1,k} + h_{i,j-1,k+1} \ge h_{i,j,k} + h_{i+1,j-2,k+1}\\
  h_{i-1,j+1,k} + h_{i-1,j,k+1} \ge h_{i,j,k} + h_{i-2,j+1,k+1}
\end{gather*}

\index{hive}

This is precisely the definition of a \emph{hive} \cite{speyer}.
There is a wonderful theorem that relates hives to the Horn problem.
The boundary of a hive consists of the three vectors  $(v_1,v_2,v_3)$
formed by taking differences of consecutive terms on the boundary. 
In terms of the $a_{ij}$'s this is

\begin{align*}
  v_1 &= (h_{n-1,1,0}-h_{n,0,0},\, h_{n-2,2,0}-h_{n-1,1,0},\cdots,h_{0,n,0}-h_{1,n-1,0})\\
  v_2 &= (h_{0,n-1,1}-h_{0,n,0},\, h_{0,n-2,2}-h_{0,n-1,1},\cdots,h_{0,0,n}-h_{0,1,n-1})\\
  v_3 &= (h_{1,0,n-1}-h_{0,0,n},\, h_{2,0,n-2}-h_{1,0,n-1},\cdots,h_{n,0,0}-h_{n-1,0,1})
\end{align*}
 then \\[.2cm]

\index{hive theorem}
\noindent\textbf{ The Hive Theorem:} Horn's  problem is solvable
for $(v_1,v_2,v_3)$ if and only if $(v_1,v_2,v_3)$  is the boundary
of a hive.\cite{speyer}\\[.2cm]

We now determine the boundary of the hive corresponding to $f$. 

\begin{align*}
  v_1 &= \biggl(
  \log\,\frac{a_{n-1,1}}{a_{n,0}},\log\,\frac{a_{n-2,2}}{a_{n-1,1}},\dots,\biggr) \\
&= \text{ the logs of consecutive coefficients of the boundary
  polynomial } f(x,1,0)\\
  v_2 &= \biggl(
  \log\,\frac{a_{0,n-1}}{a_{0,n}},\log\,\frac{a_{0,n-2}}{a_{0,n-1}},\dots,\biggr) \\
 &= \text{ the logs of consecutive coefficients of the boundary
  polynomial } f(0,x,1)\\
  v_3 &= \biggl(
  \log\,\frac{a_{1,0}}{a_{0,0}},\log\,\frac{a_{2,0}}{a_{1,0}},\dots,\biggr) \\
 &= \text{ the logs of consecutive coefficients of the boundary
  polynomial } f(1,0,x)
\end{align*}

\begin{cor}
  If $f(x,y,z)\in\rupint{3}$ is homogeneous, and $v_1,v_2,v_3$
  are the logs of the consecutive coefficients of the boundary
  polynomials, then Horn's problem is solvable for
  $v_1,v_2,v_3$.
\end{cor}

Although polynomials in $\gsubplus_2$ determine hives, the converse is
not true: here is an example of a hive that does not come from a
polynomial. We don't write the hive, but rather the coefficients of
polynomial; it's easy to check that all of the multiplicative rhombus
inequalities are satisfied. If the array were determined by a
polynomial in $\gsubplus_2$ then the polynomial corresponding to the
coefficient of $y$, $11 + 8 x + 3x^2$, would be in $\allpoly$, but it
isn't.

\centerline{
\xymatrix@=.4cm{
1 \ar@{-}[d]_{} \ar@{-}[dr] \\
6 \ar@{-}[dr] \ar@{-}[d]_{} \ar@{-}[r] & \ar@{-}[d] 3  \ar@{-}[dr]^{} \\
11 \ar@{-}[d]_{} \ar@{-}[dr]\ar@{-}[d] \ar@{-}[r]& \ar@{-}[dr]
 \ar@{-}[d] 8 \ar@{-}[r]  & 3 \ar@{-}[d]  \ar@{-}[dr]^{} \\
6 \ar@{-}[r]_{} & 11 \ar@{-}[r] & \ar@{-}[r]  6&  1
}}

\index{Horn's conjecture} Horn's problem is connected to the conjecture
(now solved - see \cite{fulton-horn}) that specifies necessary and
sufficient conditions for vectors $\alpha,\beta,\gamma$ to satisfy
Horn's problem.  We can use this to get necessary conditions for the
boundary polynomial problem.  Here is one example.  We first order
$\alpha,\beta,\gamma$:
\begin{align*}
  \alpha: & \alpha_1 \ge \alpha_2 \ge \cdots \ge \alpha_n\\
  \beta: & \beta_1 \ge \beta_2 \ge \cdots \ge \beta_n\\
  \gamma: & \gamma_1 \ge \gamma_2 \ge \cdots \ge \gamma_n
\end{align*}

Weyl's
inequalities are
$$
\alpha_i + \beta_j + \gamma_{n+2-i-j} \ge 0 \text{ whenever } i+j-1\le
n
$$ 

We now apply this to our problem. For our particular polynomial we
have that 

\begin{align*}
  v_1 &= \biggl(  \log\,\frac{61}{30},\log\,\frac{41}{61}
,\log\,\frac{11}{41},\log\,\frac{1}{11}\biggr)\\
  v_2 &= \biggl(  \log\,\frac{10}{1},\log\,\frac{35}{10}
,\log\,\frac{50}{35},\log\,\frac{24}{50}\biggr)\\
  v_3 &= \biggl(  \log\,\frac{112}{24},\log\,\frac{186}{112}
,\log\,\frac{128}{186},\log\,\frac{30}{128}\biggr)
 \end{align*}
 
 Notice that they are all decreasing - this is true in general, and is
 a consequence of Newton's inequalities.  Since these quotients are
 decreasing, they are in the correct order to apply the Weyl
 inequalities. Using the values of $v_1,v_2,v_3$ from above, the
 following is a necessary condition for the $a_{i,j}$ to be the
 coefficients of boundary polynomials.

\begin{gather*}
\log\, \frac{a_{n-i,i}}{a_{n-i+1,i-1}} +
\log\,\frac{a_{0,n-i}}{a_{0,n-i+1}} +
\log\, \frac{a_{n+2-i-j,0}}{a_{a+2-i-j+1}}  \ge 0 \\
\frac{a_{n-i,i}\,a_{0,n-i}\,a_{n+2-i-j,0}}
    {a_{n-i+1,i-1}\, a_{0,n-i+1}\,a_{a+2-i-j+1}}  \ge 1 
\end{gather*}


  \section{Log concavity in $\gsubplus_2$}
  \label{sec:log-concave-p2}

\index{log concave!in $\gsubplus_2$}

We say earlier (\mypage{sec:log-concave-p}) that the graph of the log
of the coefficients of a polynomial in $\allpolypos$ is concave. We
now prove a similar result for polynomials in $\gsubplus_2$.  Given a
polynomial $\sum a_{ij}x^iy^j$ in $\gsubplus_2$, the set of points
\[\bigl\{ (i,j,\log(a_{ij})\,)\mid a_{ij}x^iy^j\ \text{is a term of f}\ \bigr\}\]
is called the graph of the log of the coefficients of $f$. We view
this as the graph of a function defined on the points in the plane
corresponding to the exponents appearing in $f$. 

For example, the diagram below displays the log of the coefficients of
a polynomial in $\gsubplus_2$. The coordinates of the point in the lower
left corner are $(0,0)$, since $24$ is the constant term. The
coordinates of the rightmost point are $(4,0)$, and the uppermost
point's coordinates are $(0,4)$.

\centerline{ \xymatrix@=.4cm{
    \log 1 \ar@{-}[d]_{} \ar@{-}[dr] \\
    \log 10\,  \ar@{-}[dr] \ar@{-}[d]_{} \ar@{-}[r] & \log 11\,  \ar@{-}[d] \ar@{-}[dr]^{} \\
    \log 35 \ar@{-}[dr] \ar@{-}[d]_{}\ar@{-}[r] & \log 80
    \ar@{-}[dr] \ar@{-}[d] \ar@{-}[r] & \log 41 \ar@{-}[d]
    \ar@{-}[dr]^{} \\ 
    \log 50 \ar@{-}[d]_{} \ar@{-}[dr] \ar@{-}[r] & \log 175  \ar@{-}[dr]
    \ar@{-}[d] \ar@{-}[r] &
    \log 190 \ar@{-}[dr] \ar@{-}[d] \ar@{-}[r]&  \log 61  \ar@{-}[d] \ar@{-}[dr]^{} \\
    \log 24\, \ar@{-}[r]_{} \ar@{-}[r] & \log 112 \ar@{-}[r]_{} & \log
    186
    \ar@{-}[r]_{} & \log 128 \ar@{-}[r]_{} & \log30 }}

\begin{prop}
  If $f\in\gsubplus_2$ then the graph of the log of the coefficients
  is concave. Every triangle in the planar diagram determines a
  triangle in the graph, and the graph lies on one side of the plane
  containing that triangle.
\end{prop}
\begin{proof}
  We will prove the second statement, for that is what we mean by
  saying that the graph is concave. Note that it suffices to show that
  the plane through a triangle lies above all of the three adjacent
  triangles. We now show that this is a consequence of the rhombus
  inequalities.

  Consider the  triangle $T$ with solid lines 

\centerline{\xymatrix{
(r,s+1,\log a_{r,s+1}) \ar@{-}[drr] \ar@{-}[d] \ar@{..}[rr]& &
(r+1,s+1,\log a_{r+1,s+1}) \ar@{..}[d] \\
(r,s,\log a_{r,s}) \ar@{-}[rr] & & (r+1,s,\log a_{r+1,s})
}}

 The plane in $\reals^3$ containing $T$ is
\begin{multline*}
(\log  a_{r+1,s} -\log  a_{r,s})x+
(\log  a_{r,s+1} -\log  a_{r,s})y
-z -\\
\bigl(r\,\log  a_{r+1,s} + s\,\log  a_{r,s+1} -(r+s+1)\log  a_{r,s}\bigr)=0
\end{multline*}

The expression on the left side takes negative values for points
above the plane, and positive values for points below the plane. If we
evaluate it at the vertex $(r+1,s+1,a_{r+1,s+1})$ of the dotted
adjacent triangle we get 
\begin{gather*}
  \log a_{r+1,s} + \log a_{r,s+1} - \log a_{r,s} - \log a_{r+1,s+1}
  =\\
\log \frac{a_{r+1,s}\,a_{r,s+1}}{a_{r,s}\,a_{r+1,s+1}} \ge 0
\end{gather*}
where the last inequality is the rhombus inequality. A similar
computation shows that the plane lies above the other two adjacent
triangles, so the graph is concave.

\end{proof}

\begin{cor}
  Suppose that $f = \sum a_{i,j}x^iy^j\in\gsubplus_2$. For all
  non-negative  $i,j,r,s$ for which $a_{r,s}\ne0$ we have
\[
\biggl(\frac{a_{r+1,s}}{a_{r,s}}\biggr)^i\,
\biggl(\frac{a_{r,s+1}}{a_{r,s}}\biggr)^j\,
\ge\,
\frac{a_{r+i,s+j}}{a_{r,s}}
\]
\end{cor}
\begin{proof}
  Evaluate the plane through the triangle 
  \[
\{(r,s,\log a_{r,s}),(r+1,s,\log a_{r+1,s}),(r,s+1,\log  a_{r,s+1})\}
\]
 at the point $(r+i,s+j,\log a_{r+i,s+j})$. The result
  is positive, and simplifying yields the inequality.
\end{proof}

\begin{remark}
  If $a_{0,0}=1$, then the inequality shows that the coefficients of
  $f$ satisfy a simple bound:
\[ a_{1,0}^i \, a_{0,1}^j \,\ge\, a_{i,j}.\]
\end{remark}

  \section{Extending two polynomials}
  \label{sec:extension}

  The boundary polynomial problem  asks when can three
  polynomials be realized as the boundary polynomials of a polynomial
  in $\gsubplus_2$. We now ask the simpler question
  \begin{quote}
    Given two polynomials $f(x),g(x)$ in $\allpolypos$, describe the
    structure of all $F(x,y)$ in $\gsubplus_2$ such that $F(x,0)=f(x)$
    and $F(0,y) = g(y)$. We say that $F$ is an extension of $f,g$.
  \end{quote}

\index{extension!of two polynomials}

\noindent%
There is a trivial requirement that guarantees the existence of an extension.
This result is similar to Lemmas~\ref{lem:product} and \ref{lem:product-2}.

\begin{lemma}
Assume that $f,g\in\allpolypos$. There is an extension of $f,g$ if and
only if $f(0)=g(0)$. 
\end{lemma}
\begin{proof}
  If there is an extension then $F(0,0) = f(0) = g(0)$.  Conversely, we
  may assume that $f(0)=g(0)=1$.  .
    We can write $f(x) = \prod(1+ r_ix)$ and $g(x) = \prod(1+s_ix)$. The
    product below is the desired extension.
\[ \prod(1+ r_i\,x + s_i\,y)\] 

\end{proof}

Suppose that $F$ is an extension of $f,g$. We know that we can write
$F(x,y) = |I+xD_1+yD_2|$ where $D_1$ and $D_2$ are positive definite. 
Moreover, since $f(x) = |I+xD_1|$  we know the eigenvalues
of $D_1$: they are the negative reciprocals of the roots of
$f$. Similarly, we know the eigenvalues of $D_2$. We may assume that
$D_1$ is diagonal. Write $D_2 = O\Lambda O^t$ where $\Lambda$ is a
diagonal matrix consisting of the eigenvalues of $D_2$, and $O$ is an
orthogonal matrix. Consequently,

\index{Orthogonal matrices}

\begin{quote}
  There is an onto map from $O(n)$, the group of orthogonal matrices,
  to the extensions of $f,g$.
\end{quote}

Since the orthogonal matrices in $O(2)$ are easy to describe
we can completely determine the extensions for a quadratic.

\index{quadratics in $\rupint{2}$}
\begin{lemma}
  Suppose $f=(1+ r_1\,x)(1++r_2\,x)$ and $g=(1+s_1\,x)(1+s_2\,x)$.
  All extensions of $f,g$ are of the form below, where
  $W$ lies in the interval determined by $r_1s_1+r_2s_2$ and $r_1s_2+r_2s_1$.

\centerline{\xymatrix{
s_1s_2\, y^2 \ar@{-}[d] \ar@{-}[dr]  & & \\
(s_1+s_2)\, y \ar@{-}[dr] \ar@{-}[d] \ar@{-}[r] & \ar@{-}[dr]
\ar@{-}[d]
 W\, xy & \\
1 \ar@{-}[r] & \ar@{-}[r] (r_1+r_2)\, x & r_1r_2\, x^2
}}
 
\end{lemma}
\begin{proof}
If the orthogonal two by two matrix is a rotation matrix,
   the general extension of $f,g$ is given by the
  determinant of

\[ 
\begin{pmatrix} 1 & 0 \\ 0 & 1 \end{pmatrix} + x\, 
\begin{pmatrix} r_1 & 0 \\ 0 & r_2 \end{pmatrix} + y\,
\begin{pmatrix}
  \cos\theta & \sin\theta \\ -\sin\theta & \cos\theta
\end{pmatrix}
\begin{pmatrix}
  s_1 & 0 \\ 0 & s_2
\end{pmatrix}
\begin{pmatrix}
  \cos\theta & -\sin\theta \\ \sin\theta & \cos\theta
\end{pmatrix}
\]
If we simplify the coefficient of $xy$ in the determinant we find it
to be 
  \[ \frac{1}{2}\bigr[(r_1+r_2)(s_1+s_2)-(r_1-r_2)(s_1-s_2)\cos(2\theta)\bigl] \]
The two values at $\cos2\theta = \pm1$ yield the bounds of the
lemma. We get the same result if the orthogonal matrix does not
preserve orientation.
\end{proof}

The two ends of the interval for $W$ are the  values that
arise from the two different products
\begin{gather*} 
(1+r_1x+s_1y)(1+r_2x+s_2y)\\
(1+r_1x+s_2y)(1+r_2x+s_1y).
\end{gather*}
 In general, if the polynomials
are $\prod(1+r_ix)$ and $\prod(1+s_iy)$ then there are $n!$ different
products:
\[ \prod \bigl(1 + r_i\,x + s_{\sigma i}\,y\bigr)\]
where $\sigma$ is a permutation of $1,\dots,n$. If $n$ is odd then
$O(n)$ is connected, so the set of all extensions is connected. In the
case $n=2$ we saw that this set is convex, and the boundary consists
of the two different products. This may be the case in general - see
Question~\ref{ques:extension}.

It follows from the above that if the quadratic $f$ has all equal
roots then there is a unique extension. This is true for any $n$.

\begin{lemma}
  If $g\in\allpolypos(n)$ then $g(x)$ and $(x+\alpha)^n$ have a unique
  extension.
\end{lemma}
\begin{proof}
  Suppose $f(x,y)$ is an extension, so that $f(x,0)=(x+\alpha)^n$ and
  $f(0,y)=g(y)$. If we write $f(x,y) = \sum f_i(x)\,y^i$ then
  $f_0=(x+\alpha)^n$ and $f_0\lesslesseq f_1$, so there is a constant
  $a_1$ such that $f_1(x) = a_1(x+\alpha)^{n-1}$. Continuing, we find
  constants $a_i$ so that $f_i(x) = a_i(x+\alpha)^{n-i}$. Since
  $f(0,y) = \sum a_i\,\alpha^{n-i}\, y^i$ we see that $g(y)$ uniquely
  determines the $a_i$.
\end{proof}

If $f = \sum a_{ij}x^iy^j\in\gsubplus_2$ is an extension of
two polynomials, then 
$\smalltwodet{a_{10}}{a_{11}}{a_{00}}{a_{10}}$ only has one term that
depends on the extension. We assume that $f = | I+x\,D+y\,C|$ where
$D$ a diagonal matrix with entries $\{d_i\}$, and $C = (c_{ij})=
\mathcal{O}E\mathcal{O}^t$ where $E$ is diagonal with entries
$\{e_i\}$. 

Since every entry of $y\,C$ has a factor of $y$, the only way we can get
a term $xy$ in the expansion of the determinant is to take $x$ and $y$
terms from the diagonal. Thus, the coefficient $a_{11}$ of $xy$ is 
\[ 
\sum_{i\ne j} d_i\,c_{jj} = \left(\sum d_i\right)\,\left(\sum
  c_{jj}\right) - \sum d_i c_{ii}
\]

Now $a_{00} = 1$, $a_{10} = \sum d_i$ and $a_{01} = \sum c_{jj}$, so
we conclude that

\[
\begin{vmatrix}
  a_{10} & a_{11} \\ a_{00} & a_{10}
\end{vmatrix}
=
\sum d_i c_{ii}
\]

If we write $\mathcal{O} = (o_{ij})$ then $c_{ii} = \sum e_jo_{ij}^2$
and consequently 
\[
\begin{vmatrix}
  a_{10} & a_{11} \\ a_{00} & a_{10}
\end{vmatrix}
=\sum d_i e_j o_{ij}^2
\]

Now the maximum of $\sum d_i e_j o_{ij}^2$ over all orthogonal
matrices $\mathcal{O}$ is $\sum d_i'e_i'$ where $\{d_i'\}$ is
$\{d_i\}$ written in increasing order, and $\{e_i'\}$ is $\{e_i\}$
written in increasing order (Lemma~\ref{lem:bourin}).  The minimum is
$\sum d_i'e_i''$ where $\{e_i''\}$ is $\{e_i\}$ written in decreasing
order We thus can make a more precise inequality:

\[
\sum d_i'\,e_i' \ge
\begin{vmatrix}
  \  a_{10} & a_{11} \\ a_{00} & a_{10}
\end{vmatrix}
\ge \sum d_i'\,e_i'' 
\]

When $n=2$ we see that 
\[
\begin{vmatrix}
  a_{10} & a_{11} \\ a_{00} & a_{10}
\end{vmatrix}
 = \cos^2(\theta)\,(d_1e_1+d_2e_2) + \sin^2(\theta)\,(d_1e_2+d_2e_1)
\]
which shows that the determinant takes every value in the interval
determined by the largest and smallest products.

We now prove the Lemma used above; it is a generalization of the 
well-known \index{rearrangement lemma}rearrangement lemma.

\begin{lemma}\label{lem:bourin}
  If $a_1,\dots,a_n$ and $b_1,\dots,b_n$ are increasing sequences of
  positive numbers then
\[
\max_{\mathcal{O}}\, \sum_{i,j=1}^n o_{ij}^2  \,a_i\,b_j = \sum _i
a_ib_i
 \]
where the maximum is over all $n$ by $n$ orthogonal matrices $\mathcal{O}=(o_{ij})$.
\end{lemma}
\begin{proof}
  If $A = diag(a_1,\dots,a_n)$ and $B = diag(b_1,\dots,b_n)$ then 
\[
\sum o_{ij}^2 \,a_i\,b_j =  trace\,(\mathcal{O}^tA\mathcal{O}B)
\]
Since orthogonal matrices are normal, we can apply Bourin's matrix
Chebyshev inequality \cite{bourin}*{page 3} to conclude that
\[
trace \,(\mathcal{O}^tA\mathcal{O}B) \le 
trace\,(\mathcal{O}^t\mathcal{O}AB) = trace\,(AB) =\sum a_ib_i.
\]
\end{proof}

\index{matrix inequality!Chebyshev}

  \section{Symmetric Polynomials}
  \label{sec:symm-polyn}
\added{17/11/7}
\index{symmetric polynomials}

  In this section we begin an investigation of the symmetric
  polynomials in $\gsubclose_d$. We say that $f(\xx)$ is symmetric if
  $f(\xx)=f(\sigma \xx)$ for all permutations $\sigma$ of
  $x_1,\dots,x_d$. For instance, the elementary symmetric functions
  are in $\gsubclose_d$, and the product of any two symmetric
  functions in $\gsubclose_d$ is again a symmetric polynomial in $\gsubclose_d$.
  
  In two variables a homogeneous polynomial $f(x,y)$ is symmetric if and only
  if it has the form
\[ a_0y^n + a_1xy^{n-1} + a_2 x^2y^{n-2} + \cdots + a_2 x^{n-2}y^2 + a_1 x^{n-1}y + a_0x^n
\]
The roots of $f$ consist of paired roots $\{r,1/r\}$, along with some
$\pm1$'s. In degree two such a polynomial is 
\[
y^2 + \bigl(r+\frac{1}{r}\bigr)xy + x^2
\]
It's immediate that the symmetric homogeneous polynomials in
$\rupint{2}(2)$ are of the form $x^2+\alpha xy+ y^2$ where $\alpha\ge2$.
We can also characterize the symmetric homogeneous polynomials of
degree two in three variables.

  \begin{lemma}
    The symmetric homogeneous polynomials in $\rupint{3}(2)$ are of the
    form
\[ x^2+y^2+z^2 + \alpha(xy+xz+yz)\qquad \text{where $\alpha\ge2$}\]
  \end{lemma}
  \begin{proof}
    Setting $z=0$ implies $x^2+\alpha x+1\in\allpoly$, and thus
    $\alpha\ge2$. Conversely, if
    \begin{xalignat*}{2}
      D &=
      \begin{pmatrix}
        r & 0 \\ 0 & 1/r
      \end{pmatrix} &
      O &=
     \frac{1}{r+1} \begin{pmatrix}
         \sqrt{r^2+r+1} & -\sqrt{r} \\
         \sqrt{r} & \sqrt{r^2+r+1}
       \end{pmatrix}
     \end{xalignat*}
    then the determinant $f$ of $xI + y D +z ODO^t$ equals
    \[
x^2+y^2+z^2+(r+1/r)(xy+xz+yz)
\]
Since $D$ is positive definite and $O$ is orthogonal it follows that
$f\in\rupint{3}$. Finally, note that $r+1/r$ takes on all values $\ge2$.
  \end{proof}

We can show that the symmetric homogeneous polynomials in $\rupint{d}(2)$
follow the same pattern, but they don't always have a determinant
representation. 

\begin{lemma}
      The symmetric homogeneous polynomials in $\rupint{d}(2)$ are of the
    form
\[ \sum_1^d x^2_i+ \alpha \sum_{i<j} x_iy_j \qquad \text{where $\alpha\ge2$}\]

\end{lemma}
\begin{proof}
  We use the definition of $\rupint{d}$. The homogeneous part is
  certainly positive, so we check substitution. Since we are
  considering a quadratic, we only need to check that the discriminant
  of $x_d$ is non-negative. Subsituting $a_i$ for $x_i$ we get
  \begin{gather*}
    \biggl( \alpha \sum _1^{d-1} a_i\biggr)^2 - 4 \biggl(\sum_1^{d-1}
    a_i^2\biggr) \\
= (\alpha-2) \biggl( (\alpha+2) \sum_1^{d-1} a_i^2 + 2\alpha \sum_{i<j} a_ia_j\biggr)
  \end{gather*}
  The  factor $\alpha-2$  is non-negative by hypothesis, and the second factor is a quadratic
  form in $a_1,\dots,a_{d-1}$ with matrix $Q$ is
\[
\begin{pmatrix}
  \alpha+2 & \alpha & \alpha & \hdots \\
  \alpha & \alpha+2 & \alpha & \hdots \\
  \alpha & \alpha & \alpha+2 & \hdots \\
  \vdots & \vdots & \vdots & \ddots
\end{pmatrix}
= 2 I +\alpha J
\]
where $J$ is the all one matrix. Now the eigenvalues of $J$ are
$d-1,0,\dots,0$, so the eigenvalues of $Q$ are
$\alpha(d-1)+2,2,\dots,2$. Thus $Q$ is positive definite, so the
discriminant is non-negative, and  our polynomial is in $\rupint{d}$.
\end{proof}

\begin{lemma}
  The polynomial $f(x,y,z,w)$ given  below has no representation of the form $|xI + y
  D_1+zD_2+wD_3|$ where $D_1,D_2,D_3$ are positive definite two by two
  matrices and $\alpha\ge2$.
\[ f(x,y,z,w) = x^2+y^2+z^2+w^2 + \alpha (xy+xz + xw + yz + yw + zw)
\]
\end{lemma}
\begin{proof}
  Since $f(x,y,0,0) = x^2+\alpha xy+y^2$ we see that we must have
  $\alpha\ge2$, and that we can write $f(x,y,0,0) = (x+r)(x+1/r)$ where
  $r\ge1$.  We may assume that $D_1$ is diagonal, with diagonal
  $r,1/r$, and so $\alpha = r+1/r$.

  We proceed to determine $D_2$.  $f$ is symmetric for all 24
  permutations of the variables. Thus $|I + x D_1| = |I+x D_2|$,
  and therefore $D_2$ has the same eigenvalues as $D_1$. We 
   write
  \[ D_2 = 
  \begin{pmatrix}
    \cos t & \sin t \\ -\sin t & \cos t
  \end{pmatrix}
  \begin{pmatrix}
  r & 0 \\ 0 & 1/r
  \end{pmatrix}
  \begin{pmatrix}
    \cos t & -\sin t \\ \sin t & \cos t
  \end{pmatrix}
\]
Since $f(1,x,0,0) = f(0,x,1,0)$ we can equate the coefficients of $x$:
\[
r+\frac{1}{r} = 
2 \cos ^2(t)+\frac{\left(r^4+1\right) \sin ^2(t)}{r^2}
\]
Solving this equation yields
\begin{align*}
  \cos t &= \pm \frac{1}{1+r}\sqrt{1+r+r^2}\\
  \sin t &= \pm \frac{1}{1+r}\sqrt{r}
\end{align*}
Consequently there are exactly four choices for $D_1$ and for $D_2$:
\begin{align*}
D_1& =\frac{1}{r+1}
\begin{pmatrix}
  e_1 \sqrt{1+r+r^2} & e_2 \sqrt{r} \\
  -e_2 \sqrt{t} & e_1\sqrt{1+r+r^2}
\end{pmatrix}
\\
D_2& =\frac{1}{r+1}
\begin{pmatrix}
  f_1 \sqrt{1+r+r^2} & f_2 \sqrt{r} \\
  -f_2 \sqrt{t} & f_1\sqrt{1+r+r^2}
\end{pmatrix}
\end{align*}
where $e_1,e_2,f_1,f_2$ are $\pm1$. If we equate the coefficients of
$x$ in $|xD_1+D_2] = |xD_1+I|$ we get
\[
\frac{2}{(1+r)^2} \bigl(e_2f_2r +e_1f_1 (1+r+r^2)\bigr) =
\frac{2e_1}{1+r}\sqrt{1+r+r^2}
\]
It's easy to see that this equation has no solution for $r\ge1$.
\end{proof}

A symmetric homogeneous polynomial $F(x,y,z)$ in $\rupint{3}(3)$ has coefficient
array
\centerline{\xymatrix@=.2cm{
&&&1 \\
&&a&&a\\
&a&&b&&a\\
1 &&a &&a&&1
}}

We first show that for any $a\ge3$ there is at least one $b$ such that
$F(x,y,z)\in\rupint{3}$. If we take $a=\alpha+1$ and $b=3\alpha$ then
$\alpha\ge2$.  Now
\[
(x+y+z)(x^2+y^2+z^2 + \alpha(xy+xz+yz)) \in\rupint{3}
\]
and has coefficient array

\centerline{\xymatrix@=.2cm{
&&&1 \\
&&1+\alpha&&1+\alpha\\
&1+\alpha&&3\alpha  &&1+\alpha\\
1 &&1+\alpha &&1+\alpha&&1
}}

Newton's inequality applied to the first and second rows yield $a\ge3$
and $b\ge2a$. However, we can get more precise results. As above, the
roots of the bottom row are $-r,-1/r,-1$, so $a = 1+r+1/r$. We assume
$r>1$. Since the second row interlaces the first, we have that
\[
-r \le \frac{-b-\sqrt{b^2-4a^2}}{2a} \le -1 \le
\frac{-b+\sqrt{b^2-4a^2}}{2a} \le -1/r
\]
These inequalities yield
\[ a(r+1/r) \ge b \ge 2a \]
Note that if $r=1$ then $a=3$ and $b=6$. This is $(x+y+z)^3$, which is
the unique extension of the bottom row and the left hand column. These
are not the best bounds; empirically we find that the intervals for
$b$ appear to be

\begin{center}
\begin{tabular}{crrlc}
  \toprule
  r=2  &\ $b\in$& $\displaystyle\bigl(\frac{29}{4}$ &,& $\displaystyle\frac{16}{2}\bigr)$ \\[.2cm]
  r=3  &\ $b\in$& $\displaystyle\bigl(\frac{90}{9}$ &,& $\displaystyle\frac{38}{3}\bigr)$ \\[.2cm]
  r=5  &\ $b\in$& $\displaystyle\bigl(\frac{326}{25}$ &,&
  $\displaystyle\frac{142}{5}\bigr)$ \\
  \bottomrule
\end{tabular}
\end{center}

\added{10-8-07}
   When is $\sum a_i \sigma_i\in\multiaff{d}$, where $\sigma_i$ is
   the $i$'th elementary symmetric polynomial in $x_1,\dots,x_d$?
   If we set all $x_i$ equal to $x$, then $\sum a_i
   \binom{d}{i}x^i\in\allpoly$. We are not able to prove the converse;
   as usual we can't get the final factorial.

   \begin{lemma}\label{lem:ma-sym}
     If $\sum a_i x^i\in\allpoly$ then
     \begin{enumerate}
     \item $\sum a_i \sigma_i \in\multiaff{d}$.
     \item $\sum a_i i! \sigma_i \in\multiaff{d}$.
     \end{enumerate}
   \end{lemma}
   \begin{proof}
     The first sum is the coefficient of $y^n$ in
\[
\sum a_i y^i \biggl(\prod_1^d(y+x_i)\biggr) = \sum a_i y^i \cdot \sum
y^{n-i}\sigma_i
\]
while the second equals
\[
\biggl(\sum a_i \diffd_y^i\biggr) \sum y^{n-i}\sigma_i\biggr|_{y=0}
\]
   \end{proof}

   \begin{example}
     We give an example that supports the possibility  that the
     converse is true. Consider the question
\[
\text{For which $m$ is } xyz - x-y-z + m = \sigma_3-\sigma_1+m\,\sigma_0 \in\multiaff{3}?
\]

We find
\begin{align*}
  x^3-x+m & \in\allpoly & \implies && |m| &\le \sqrt{4/27} \\
  x^3/6-x+m & \in\allpoly & \implies && |m| &\le \sqrt{8/9} \\
  \binom{3}{3}x^3-\binom{3}{1}x+\binom{3}{0}m & \in\allpoly & \implies && |m| &\le 2
\end{align*}
We know from Example\ref{ex:ma-xyz} that the necessary and sufficient condition is
that $|m|\le2$.
   \end{example}

\section{Linear transformations on  $\gsubpos_d$}
\label{sec:sub-mappings}

We have looked at linear transformations that map $\gsubpos_d$ to
itself. In this section we look at mappings that either delete
variables or introduce new variables. 

The repeated application Theorem~\ref{thm:pd-T} with the linear transformation
$x\mapsto ax+b$ yields

\begin{cor}
  Choose $a_1,b_1,\dots,a_d,b_d\in\reals$ where all  of the $a_i$'s
  are positive. The linear transformation below maps $\gsubpos_d$ to itself:
  $$
  f(x_1,\dots,x_d) \mapsto f(a_1x_1+b_1,\dots,a_dx_d+b_d).$$
\end{cor}

We can allow some of the $a_i$ to be zero.

\begin{lemma} \label{lem:sub-ax}
  If $f(\xx)\in\gsubpos_d$ and exactly $r$ terms of $a_1,\dots,a_d$ are
  non-zero, and they are positive, 
  then $f(a_1x_1,\dots,a_dx_d)\in\gsubpos_r$.
\end{lemma}

We can also add new variables.

\begin{lemma} \label{lem:sub-xandy}
  If $f\in\gsubpos_d(n)$ then $g=f(x_1+y_1,x_2+y_2,\dots,x_d+y_d)$ is in
  $\gsubpos_{2d}(n)$.
\end{lemma}
\begin{proof}
  Certainly $g$ satisfies substitution.  Every term
  $x_1^{i_1}\cdots x_d^{i_d}$ in $f$ gives rise to
  $(x_1+y_1)^{i_1}\cdots (x_d+y_d)^{i_d}$ in $g$.
  This shows that
  $$
  g^H(\xx) = f^H(x_1+y_1,\dots,x_d+y_d)$$
  Since $f$ satisfies the homogeneity conditions, so does $g$.
\end{proof}

\begin{lemma} \label{lem:pd-xx}
  If $f(\xx)\in\gsubpos_d$ then the polynomial 
  $g=f(x,x,x_3,\dots,x_d)$ is in $\gsubpos_{d-1}$. Also,
  $f(x,x,\dots,x)\in\allpoly$. 
\end{lemma}
\begin{proof}
  If $d$ is two then this is Lemma~\ref{lem:p2-xx}. If $d$ is larger than two
  then $g$ satisfies $x_d$-substitution.  Since $f(x,x,x_3,\dots)^H =
  f^H(x,x,x_3,\dots)$ it follows that $g$ satisfies the 
  homogeneity condition. The second assertion follows by induction.
\end{proof}

We can replace one variable  by more complicated expressions. 

\begin{lemma} \label{lem:pd-xxg}
  If $f(\xx)\in\gsubpos_d$, $\beta\in\reals$, and $a_1,\dots,a_{d-1}$
  are non-negative then
  $$ g(x_1,\dots,x_{d-1})=
  f(x_1,x_2,\dots,x_{d-1},a_1x_1+\dots+a_{d-1}x_{d-1}+\beta)\in\gsubpos_{d-1}$$ 
\end{lemma}
\begin{proof} 
  By induction on $d$. The case $d=2$ is Lemma~\ref{lem:pd-xx}. Since $g$
  satisfies the positivity  condition it remains to show
  that $g$ satisfies substitution. If we substitute $x_i=\alpha_i$ for 
  all $x_i$ except $x_1$ and $x_d$ then
  $f(x_1,\alpha_2,\dots,\alpha_{d-1},x_d)\in\gsubpos_2$. Consequently, 
  $f(x_1,\alpha_2,\dots,\alpha_{d-1},a_1x_1+b)$ is in $\allpoly$ for any
  choice of $b$. Choosing $b =
a_2\alpha_2+\dots+a_{d-1}\alpha_{d-1}+\beta$ and observing that 
$$ g(x_1,\alpha_2,\dots,\alpha_{d-1}) =
f(x_1,\alpha_2,\dots,\alpha_{d-1},a_1x_1+a_2\alpha_2+\dots+a_{d-1}\alpha_{d-1}+\beta)$$ 
it follows that  $g$ satisfies substitution.
\end{proof}

We now use these results  to show how very general linear
combinations preserve $\gsub$.

\index{matrix!preserving interlacing}
\begin{theorem} \label{thm:sub-mx}
  Suppose $f(\xx)\in\gsubpos_d$.  If $M$ is  a $d$ by $d$ matrix  with
  positive entries  then $f(M \xx )\in\gsubpos_d$. 
\end{theorem} 
\begin{proof}
Let $M=(m_{ij})$ and define a polynomial in $d^2$ variables $x_{ij}$
$$ g = f\left(\sum_i\,m_{i1}x_{i1},\dots,\sum_i\,m_{id}x_{id}\right).$$
Repeated applications of Lemma~\ref{lem:sub-ax} and Lemma~\ref{lem:sub-xandy} show that
$g$ is in $\gsubpos_{d^2}$.  We now use Lemma~\ref{lem:pd-xx} to identify all
pairs of variables $x_{ij}$ and $x_{ik}$. The resulting polynomial is 
exactly $f(M \xx)$.
\end{proof}

\begin{lemma}
  If $f,g\in\gsubpos_d$ and $f=\sum_\sdiffi a_\sdiffi \xx^\diffi$ then
  for any index $\diffk$ 
$$ \sum_{\diffi+\diffj=\diffk} a_\sdiffi\, \frac{g^{(\diffj)}}{\diffj!}
  \,\in\gsubpos_d.$$ 
\end{lemma}
\begin{proof}
  Using the Taylor series for $g$ we get
$$ f(\yy)g(\xx+\yy) = \sum_{\diffi,\diffj} a_\sdiffi\, \xx^\diffi\,
\frac{\yy^\diffj}{\diffj!} \,\frac{\partial^\diffj}{\partial \xx^\diffj}\,
g(\xx)$$ 
The coefficient of $\yy^\diffk$ is the sum in the conclusion.
\end{proof}

\section{The graph of polynomials in $\rupint{d}$}
\label{sec:pd-graph}

\index{graph!of $\rupint{d}$}

In $\rupint{2}$ we used the geometry of the graph of $f(x,y)\in\gsub_2$
to show that $f(x,x)\in\allpoly$. Now we use the fact that if
$f(\xx)\in\rupint{d}$ then $f(x,x,\dots,x)\in\allpoly$ to get information
about the geometry of the graph of $f$. 

For clarity we'll take $d=3$. Choose $f(\xx)\in\rupint{d}(n)$. The graph
${G_f}$ of $f$ is $\{(x,y,z)\mid f(x,y,z)=0\}$. Choose any point
$v=(a,b,c)$ on the plane $x_1+x_2+x_3=0$. Since $f(\xx-v) =
f(x-a,y-b,z-c)$ is in $\rupint{d}$, the equation $f(x-a,x-b,x-c)=0$ has
$n$ solutions. Denote the $i$-th largest solution by $s_i(v)$, and let
$w_i(v) = (s_i(v)-a,s_i(v)-b,s_i(v)-c)$. By construction we know that
$f(w_i(v))=0$, so $w_i(v)$ is a point on the graph of $f$. Consider
the map $v\mapsto w_i(v)$. This is a function from the plane
$x_1+x_2+x_3=0$ to the graph of $f$ that is $1-1$, and continuous
since solutions to equations are continuous functions of coefficients.
Moreover, since the union of all lines of the form $\{\nu+t(1,1,1)\}$ is
a partition of $\reals^d$, every point on the graph of $f$ is of the
form $w_i(\nu)$ for a unique $\nu$, and perhaps several $i$'s. Thus, the
graph ${G_f}$  is the union of $n$ subsets ${G_f}_i$
of $\reals^d$, each homeomorphic to $\reals^{d-1}$. These subsets
might intersect.

In addition, since all coefficients of $f^H$ are positive, there is an
$\alpha$ such if $|\xx|>\alpha$ and all coefficients of $\xx$ are
positive then $f(\xx)>0$. Similarly, we may assume that if all
coefficients are negative then $f(\xx)\ne0$. 

Let $\reals^+$ be the set of all points with all positive coefficients,
 $\reals^-$ the set with all negative coefficients, and  $S$ be
the ball of radius $\alpha$. The previous paragraph shows that the
graph ${G_f}$ separates $\reals^+\setminus S$ from
$\reals^-\setminus S$.

We have seen that if $a,b,c$ are positive and $f\in\rupint{3}(n)$, then
$f(ax+\alpha,bx+\beta,cx+\gamma)$ is in $\allpoly(n)$. This is a
special case of the following observation:

\begin{quote}
  If a curve  goes to infinity in all coordinates in $\reals^+$ and also in
  $\reals^-$ then it meets the graph of $f$ in at least $n$ points.
\end{quote}

Here's a variation: suppose that the line $\mathcal{L}$ is parallel to
a coordinate axes. For instance, $\mathcal{L}$ could be $\{
(a,b,c)+t(1,0,0)\mid t\in\reals\}$. Then, $\mathcal{L}$ meets the
graph in $n$ points. To see this, notice that there is a sufficiently
large $T$ such that the line segment from $(a+T,b,c)$ to $(a+T,0,0)$
does not meet ${G_f}$.  Neither does the ray $\{(a+T+t,0,0)\mid
t\ge0\}$. Similarly analyzing negative values we can find a
corresponding $S$ and conclude that the piecewise linear curve
consisting of the segments

\begin{gather*}
  \{(a+T+t,0,0)\mid t\ge0\}\\
(a+T,b,c)  - (a+T,0,0) \\
(a+T,b,c)  - (a-S,b,c) \\
(a-S,b,c)  - (a-S,0,0) \\
\{(a-S-t,0,0)\mid t\ge0\}\\
\end{gather*}

\noindent%
meets   ${G_f}$ in $n$ points, and consequently the
line $\{(a+t,b,c)\mid t\in\reals\}$ meets ${G_f}$ in $n$ points.

If $f\in\gsubplus_d$ then $\xx\in\reals^+$ implies that
$f(\xx)\ne0$. Any curve from a point in $\reals^+$ that goes to
infinity in $\reals^-$ must intersect every $\mathcal{G}_i$.

\section{Differential operators}
\label{sec:sub-differential}
\index{differential operators}
\index{operator!differential}

In this section we introduce differential operators on $\gsubpos_d$ of the form
$$
f(\frac{\partial}{\partial\xx}) = f(\frac{\partial}{\partial
  x_1},\dots,\frac{\partial}{\partial x_d})$$ where $f$ is a
polynomial in $d$ variables. To simplify notation, we will sometimes
write $\partial_\xx$ for $\frac{\partial}{\partial\xx}$. In
Lemma~\ref{lem:fdxdy} we will  show that $f(\partial_\xx)$ maps
$\gsubpos_d$ to itself. This requires the consideration of exponential
functions. If we restrict ourselves to polynomials then 

\begin{lemma}  \label{lem:pd-diff-prod}
  Suppose that $f$ is a product of linear terms
  $$
  f = \prod_i(b_i + a_{1i}x_1 + \cdots a_{di}x_d) = \prod_i(b_i +
  \aaa_i\cdot\xx)$$
  where all
  the coefficients $a_{ji}$ are positive and
  $\aaa_i=(a_{1i},\dots,a_{di})$. The map
  $f(\frac{\partial}{\partial\xx})$ is a linear transformation from
  $\gsubpos_d$ to itself.
\end{lemma}
\begin{proof}
  It suffices to assume that $f$ has a single factor.  If $g$ is in
  $\gsubpos_d$, then any derivative of $g$ has total degree less than
  $g$. If the constant term of $f$ is non-zero, then
  $(f(\frac{\partial}{\partial\xx})g)^H$ is a constant multiple of
  $g^H$. If the constant term is zero then
  $(f(\frac{\partial}{\partial x})g)^H$ is a positive linear
  combination of the partial derivatives of $g^H$ which implies that
  all of its coefficients are positive. In either case 
  $(f(\frac{\partial}{\partial x})g)^H$ satisfies the 
  positivity condition.
  
  Substitution is also satisfied since since all the partial
  derivatives interlace $g$, and and are combined with the correct
  signs.  

\end{proof}

The next two lemmas are generalizations  of Corollary~\ref{cor:fpp}
and Lemma~\ref{lem:gxdy}; the proof
are the same.

\begin{lemma}\label{lem:fpp-d}
  If $g\in\allpolypos$, $f(\xx)\in\gsubpos_d$ then 
$g(-\partial{x_1}{\partial{x_2}})\,f(\xx)\in\gsubpos_d$.
\end{lemma}
  
\begin{lemma} \label{lem:gxdy-d}
  If $g\in\allpolypos$, $f(\xx)\in\gsubpos_d$ then 
$g(x_1\frac{\partial}{\partial x_2})\,f(\xx)\in\gsubpos_d$.
\end{lemma}

\begin{lemma} \label{lem:pd-diff-prod-2}
  Let $\ccc^i=(\gamma^i_1,\dots,\gamma^i_d)$ be vectors of all
  positive terms, and $\bbb^i=(\beta^i_1,\dots,\beta^i_d)$ vectors of all
  non-negative terms.  The operator
$$ \prod_i \left( \ccc^i\cdot \xx -
  \bbb^i\cdot(\frac{\partial}{\partial\xx})\right)$$ 
maps $\gsubpos_d$ to itself.
\end{lemma}
\begin{proof}
  It suffices to verify it for one factor. Choose $g\in\gsubpos_d$ and 
  set 
  $$
  h = \left( \ccc\cdot\xx - \bbb\cdot\frac{\partial}{\partial\xx}\right)\,g.$$
  Since
  $h^H = (\ccc\cdot\xx)g^H$ we see that $h$ satisfies the
  homogeneity condition. Set $k=(\bbb\cdot\frac{\partial}{\partial
    \xx})\,g$. Since derivatives interlace and all coefficients of
  $\bbb$ are non-negative, we find that
  $g\lesslesseq k$.    To verify substitution we must show that
  $\ccc\cdot(\xx_i^\aaa)g(\xx_i^\aaa)-k(\xx_i^\aaa)$ is in $\allpoly$ for 
  any $\aaa\in\reals^d$.
  This follows from Lemma~\ref{lem:sign-quant}.

\end{proof}

\begin{cor} \label{cor:pd-diff-prod-2}
  Let $\ccc^i=(\gamma^i_1,\dots,\gamma^i_d)$ be vectors of all
  positive terms, and $\bbb^i=(\beta^i_1,\dots,\beta^i_d)$ vectors of
  all non-negative terms.  The map
  $$f\times g\mapsto f(\ccc\cdot\xx - \bbb\cdot\frac{\partial}{\partial
    \xx}\,)\,g(\xx)$$
  defines a map $\allpoly\times \gsubpos_d \longrightarrow \gsubpos_d$.
\end{cor}

These last  results were generalizations to $d$ variables of the
corresponding one variable results. The next result is strictly a
property of more than one variable.

\begin{theorem} \label{thm:pd-diff-prod}
  If $f(\xx,\yy)=\sum a_\sdiffi(\xx) \yy^\diffi$ has the property that
  $f(\frac{\partial}{\partial\xx},\frac{\partial}{\partial\yy})$ maps
  $\gsubpos_{d+e}$ to itself, then all coefficients $a_\sdiffi(\xx)$
  determine operators $a_\sdiffi(\frac{\partial}{\partial\xx})$ that
  map $\gsubpos_d$ to itself.
\end{theorem}
\begin{proof}
  Choose $g\in\gsubpos_d$, a large integer $N$, and let $\yy^N =
  y_1^N\cdots y_e^N$. The action of
  $f(\frac{\partial}{\partial\xx},\frac{\partial}{\partial\yy})$ on
  $g(\xx)\yy^N$ is
  \begin{align*}
    f(\frac{\partial}{\partial\xx},\frac{\partial}{\partial\yy})g(\xx)\yy^N
    &= \left(\sum
      a_\sdiffi(\frac{\partial}{\partial\xx})\,(\frac{\partial}{\partial\yy})^\sdiffi          \right)g(\xx)\yy^N \\
    &= \sum a_\sdiffi(\frac{\partial}{\partial\xx})g(\xx) \, \left(
      \frac{\partial^\diffi}{\partial\yy^\diffi} \yy^N\right)
  \end{align*} 
  Since $N$ is large, all the $\yy$ terms
  $\frac{\partial^\diffi}{\partial\xx^\diffi} \yy^N$ are non-zero and
  distinct. Since the left hand side is in $\gsubpos_{d+e}$, all its
  coefficients are in $\gsubpos_d$, and hence
  $a_\sdiffi(\partial_\xx)g(\xx)$ is in $\gsubpos_d$. Thus
  $a_\sdiffi(\partial_\xx)$ maps $\gsubpos_d$ to itself.
\end{proof}

If we choose $f$ to be a product of linear terms, then the coefficients
are generally not products.  This shows that there are polynomials $h$
that are not products such that $h(\frac{\partial}{\partial \xx})$
maps $\gsubpos_d$ to itself.

We now generalize Lemma~\ref{cor:fgxdy}, and find a Hadamard product.
\index{Hadamard product!in $\rupint{d}$}

\begin{lemma}\label{lem:fgxdy-d}
  If $\sum f_i(\xx)y^i\in\rupint{d+1}$ and $\sum a_ix^i\in\allpolypos$
  then
\begin{equation}\label{eqn:fgxdy-d}
\sum i!\,a_i f_i(\xx)y^i\in\rupint{d+1}
\end{equation}
\end{lemma}
\begin{proof}
  From Lemma~\ref{lem:gxdy-d} we know that 
$$\sum a_i(y\partial_z)^i \sum f_j(\xx)z^j \in\gsubclose_{d+2}$$
and the coefficient of $z^0$ yields \eqref{eqn:fgxdy-d}.
\end{proof}

In the following lemma we would like to replace 
$\exp_\xx\,  f^{rev}$ with $f$, where we
define $\exp_\xx =\exp_{x_1}\cdots \exp_{x_d}.$

\begin{lemma} \label{lem:diff-op-2}
  Suppose that $f,g\in\gsubpos_d(n)$ and write $f = \sum_{|\sdiffi|\le n}
  a_\sdiffi \xx^\sdiffi$. If $f^{rev} = \sum a_\sdiffi \xx^{n-\sdiffi}$
  then
$$ \left(\exp_\xx\,
  f^{rev}\right)(\frac{\partial}{\partial \xx})\,g\in\gsubpos_d$$ 
\end{lemma}
\begin{proof}
The coefficient of $y_1^n\cdots y_d^n$ in the product of $f(\yy)$ and
the Taylor series \index{Taylor series} of $g(\xx+\yy)$
\begin{gather*}
  \left(\sum a_\sdiffi\, \yy^\sdiffi\right)\left(\sum
    g^{(\sdiffi)}(\xx)\frac{\yy^\sdiffi}{\sdiffi!}\right)  \\
\intertext{equals}
\sum \frac{a_\sdiffi}{(n-\sdiffi)!} {g^{(n-\sdiffi)}}(\xx) =
\left(\sum a_\sdiffi \frac{\xx^{n-\sdiffi}}{(n-\sdiffi)!}\right)
(\frac{\partial}{\partial\xx})\,g 
\end{gather*}
and the latter expression is exactly
$\exp_{x_1}\cdots \exp_{x_d}\,  f^{rev}(\frac{\partial}{\partial\xx})\,g .$
\end{proof}

Next we generalize Corollary~\ref{cor:fxdg} and \eqref{eqn:fxdg}. 
\begin{lemma} \label{lem:pd-fxdg}
  If $g\in\allpoly$ has no roots in $[0,n]$ and $f\in\gsubpos_d(n)$
  then $$ g(x_1 \frac{\partial}{\partial x_1})f \in\gsubpos_d$$
\end{lemma}
\begin{proof}
  It suffices to take $g=-a+x_1$ where $a\not\in[0,n]$. We need to show
  that $h= -a f + x_1\frac{\partial f}{\partial x_1}\in\gsubpos_d$.
  We first observe that Corollary~\ref{cor:agxgp} implies that $h$ satisfies
  $x_1$-substitution.  It remains
  to check the signs of the coefficients of $h^H$. If
  $\diffi=(i_1,\dots,i_d)$ and $f = \sum a_\sdiffi \xx^\sdiffi$ then
  the coefficient of $\xx^\diffi$ in $h$ is $(i_1-a)a_\sdiffi$. Since
  $0\le i_1\le n$ it follows that if $a<0$ the signs of all the
  coefficients of $h^H$ are positive, and if $a>n$ the signs are all
  negative. Thus, in either case $h\in\pm\gsubpos_d$.
\end{proof}

\begin{cor}
  Suppose that $f\in\gsubpos_d(n)$, and $g_1,\dots,g_d$ have the
  property that none of them have a root in $[0,n]$.  If $f = \sum
  a_\sdiffi\xx^\sdiffi$ then
\begin{equation} \label{eqn:pd-gi} 
\sum_{\sdiffi=(i_1,\dots,i_d)}
    g_1(i_1)\cdots g_d(i_d)\,a_\sdiffi \xx^\sdiffi\ \in\gsubpos_d
\end{equation}
In particular, if $f\in\gsubpos_d$ and $g_1,\dots g_d\in\allpolypos$
then \eqref{eqn:pd-gi} holds.

\end{cor}
\begin{proof}
  Use Lemma~\ref{lem:pd-fxdg}, equation \eqref{eqn:fxdg}, and the calculation
  \begin{align*}
    g(x_1\frac{\partial}{\partial x_1}) &= 
\sum_\sdiffi g(x_1 \frac{\partial}{\partial x_1}) \,a_\sdiffi
\xx^\sdiffi \\
&= \sum_{\sdiffi=(i_1,\dots,i_d)} g(i_1)a_\sdiffi \xx^\sdiffi
  \end{align*}
\end{proof}

For example, if we take $f=(x+y+1)^n$, $g_1(x)=(x+1)^r$ and
$g_2(y)=(y+1)^s$ then 
$$ \sum_{0\le i+j\le n} (i+1)^r (j+1)^s \binom{n}{i,j} x^i y^j \ \in\
\gsubpos_2$$

Here is a partial converse to Lemma~\ref{lem:pd-diff-prod}.
\begin{prop}
  If $f(\xx)$ is a polynomial with the property that
  $f(\frac{\partial}{\partial\xx})$ maps $\xsub_d$ to itself, then $f$
  is in $\xsub_d$.
\end{prop}
\begin{proof}
  Let $f=\sum a_\sdiffi \xx^\diffi$ have $x_i$ degree $k_i$. Consider the 
  action of $f(\frac{\partial}{\partial\xx})$ on the polynomial $g = x_1^{k_1+N}\cdots
  x_d^{k_d+N}$ where $N$ is large:
  \begin{align*}
    f(\frac{\partial}{\partial\xx})g &= \sum a_I \left( \frac{\partial}{\partial
        \xx^\diffi}\right) g \\
    &= \sum a_\sdiffi
    \falling{k_1+N}{i_1}\cdots\falling{k_d+N}{i_d}x_1^{k_1+N-i_1}\cdots
    x_d^{k_d+N-i_d} \\
    \intertext{where $\diffi = (i_1,\dots,i_d)$. If
      $f(\frac{\partial}{\partial\xx})g$ is 
      in $\xsub_d$ then so is its reversal(Proposition~\ref{prop:sub-reverse})} & \sum
    a_\sdiffi 
    \falling{k_1+N}{i_1}\cdots\falling{k_d+N}{i_d} x_1^{i_1}\cdots
    x_d^{i_d}.\\
    \intertext{If we replace $x_i$ by $x_i/N$ we get that} & \sum
    a_\sdiffi \frac{\falling{k_1+N}{i_1}}{N^{i_1}} \cdots
    \frac{\falling{k_d+N}{i_d}}{N^{i_d}}\, x_1^{i_1}\cdots x_d^{i_d}.
  \end{align*}
is also in $\xsub_d$. Since this converges to $f$ we conclude that $f\in\xsub_d$.
\end{proof}

  \section{Higher order Newton inequalities}
  \label{sec:higher-newton}

  We show how the Newton inequalities \eqref{eqn:newton-1} can be interpreted as
  statements about interlacing. This perspective leads to the Newton
  interlacing proposition.

  If $g(x) = b_0 + b_1 x + b_2 x^2 + \cdots \in\allpolypos$ then
  Newton's inequalities state that $b_k^2 \ge \frac{k+1}{k}\,
  b_{k-1}b_{k+1}$. (If we know the degree of $g$ we can be more
  precise.) We rewrite this as the sequence of inequalities
\[
\frac{1\cdot b_1}{b_0} \ge  \frac{2\cdot b_2}{b_1} \ge 
\frac{3\cdot b_3}{b_2} \ge \cdots
\]

Each fraction is the negative of the root of a linear polynomial, so 
we can reinterpret these inequalities as statements about the
interlacing of a sequence of polynomials.
\[
b_0 + 1\cdot b_1\,x \greateqeq  
b_1 + 2\cdot b_2\,x \greateqeq 
b_2 + 3\cdot b_3\,x \greateqeq \cdots
\]

These polynomials
correspond to the double lines in the diagram below, where  we write the
coefficients of $g'$ and $g$:

\centerline{\xymatrix{
1\cdot b_1 \ar@{=}[d] \ar@{.}[r] & 2\cdot b_2 \ar@{.}[r]  \ar@{=}[d]
& 3\cdot b_3 \ar@{.}[r]  \ar@{=}[d]  &  \\
b_0 \ar@{.}[r] & b_1 \ar@{.}[r] & b_2 \ar@{.}[r] & 
}}

Now in this diagram $g$ and $g'$ can be considered to be the first two
rows of a polynomial in $\gsubpos_2$, as in
Figure~\ref{fig:gen-newton}. If we simply take quadratic polynomials
instead of linear polynomials, then we find that the resulting
polynomials don't always have all real roots. The solution is to
introduce binomial coefficients.

\index{Newton's inequalities!interlacing version of}
  \begin{prop}
    Suppose that $f\in\gsubplus_2(n)$. Choose a positive integer $m$.
    If $f = \sum a_{i,j}x^iy^j$ then define the polynomials $\mathcal{F}_k$ for
    $0\le k \le n-m$ by
\[
\mathcal{F}_k(x) = \sum_{j=0}^m a_{k,j}\,\binom{k}{j}\,x^j
\]
All the $\mathcal{F}_k$ have all real roots and interlace:
\[
\mathcal{F}_0 \greateqeq \mathcal{F}_1 \greateqeq \mathcal{F}_2  \greateqeq \cdots \greateqeq
\mathcal{F}_{n-m}
\]
  \end{prop}

Consider an example where $m=2$ and $n=4$. There are three
polynomials, which are identified by the double lines in Figure~\ref{fig:gen-newton}.
\begin{align*}
    \mathcal{F}_0 &= a_{00} + 2 a_{01}x + a_{02}x^2 \\
    \mathcal{F}_1 &= a_{10} + 2 a_{11}x + a_{12}x^2 \\
    \mathcal{F}_2 &= a_{20} + 2 a_{21}x + a_{22}x^2 \\
\end{align*}
The proposition asserts that $  \mathcal{F}_0 \lesslesseq
\mathcal{F}_1 \lesslesseq   \mathcal{F} _2$. 

\begin{figure}[hb]
  
\centerline{\xymatrix{
a_{04} \ar@{.}[d] \ar@{.}[dr]& & & & \\
a_{03} \ar@{.}[d] \ar@{.}[dr] \ar@{.}[r]& a_{13} \ar@{.}[d] \ar@{.}[dr] &&&\\
a_{02} \ar@{=}[d] \ar@{.}[dr] \ar@{.}[r]& a_{12} \ar@{=}[d]
\ar@{.}[dr] \ar@{.}[r]& a_{22}
\ar@{=}[d] \ar@{.}[dr] & &\\
a_{01} \ar@{=}[d] \ar@{.}[dr] \ar@{.}[r]& a_{11} \ar@{=}[d]
\ar@{.}[dr] \ar@{.}[r]& a_{21}
\ar@{=}[d] \ar@{.}[dr] \ar@{.}[r]& a_{31}\ar@{.}[d]\ar@{.}[dr]  &\\
a_{00} \ar@{.}[r]& a_{10}\ar@{.}[r] & a_{20}\ar@{.}[r] & a_{30}\ar@{.}[r] &  a_{40}\\
}}

  \caption{A polynomial in $\rupint{2}(4)$}
  \label{fig:gen-newton}
\end{figure}

\begin{proof}(of the proposition)
  If we write 
  \begin{align*}
    f(x,y) &= f_0(y) + f_1(y)x+f_2(y)x^2 + \cdots \\
\intertext{then }
\mathcal{F}_k &= (1+y)^m \ast f_k(y) \\
  \end{align*}
  By Lemma~\ref{lem:fgxdy-d} we know that $(1+y)^m \ast f(x,y)
  \in\gsubclose_2$. The conclusion of the proposition follows since
  consecutive coefficients of polynomials in $\rupint{2}$ interlace.

\end{proof}

\section{Negative subdefinite matrices}
\label{sec:neg-subdef}
\index{quadratic forms}

When does a \index{quadratic form}quadratic form $\xx Q \xx^t$ have the property that
$\xx Q \xx^t -c^2 \in\gsubpos_d$ for all $c$? The answer is
simple: $Q$ is negative subdefinite.
\index{subdefinite matrix!negative}
\index{matrix!negative subdefinite}
\index{negative subdefinite matrix}%

\begin{definition}
  A real symmetric matrix $A$ is \emph{negative subdefinite} if it has
  exactly one positive eigenvalue, and all entries are positive.
\end{definition}

Our goal is the following theorem, which will follow from
Lemmas~\ref{lem:subdef-1} and \ref{lem:subdef-2} in the next section.

\begin{theorem}\label{thm:nsd}
  Suppose $Q$ is a real symmetric matrix.  The \index{quadratic form}quadratic form $\xx Q
  \xx^t -c^2$ is in $\gsubpos_d$ for all $c$ if and only if $Q$ is
  negative subdefinite.
\end{theorem}

In this section we establish various properties of negative
subdefinite matrices. We should first observe that among all real
symmetric matrices with all positive entries, the negative subdefinite
matrices are those with the smallest possible number of positive
eigenvalues. They can't all be negative, since Perron's theorem
guarantees at least one positive eigenvalue.

 We first note
that subdefinite matrices are hereditary in the following sense:

\index{principle submatrices}

\begin{lemma} \label{lem:nsd-1}
If $A$ is a negative subdefinite matrix then all the principle
submatrices of $A$ are also negative subdefinite.  
\end{lemma}
\begin{proof}
  It suffices to show that the matrix $T$ formed by removing the last
  row and column is negative subdefinite. The eigenvalues of $T$
  interlace the eigenvalues of $A$ by Theorem~\ref{thm:principle-1}. Thus, since
  $A$ has exactly $n-1$ non-positive eigenvalues, $T$ has at least
  $n-2$ non-positive eigenvalues. Now the sum of the eigenvalues of
  $T$ is the trace of $T$ which is positive, so $T$ must have at least
  one positive eigenvalue. Thus $T$ has exactly one positive
  eigenvalue.
\end{proof}

Here is an alternative characterization of negative subdefinite
matrices.  Recall that a symmetric matrix is positive definite if and
only if the determinants of the principle submatrices are all
positive. The $k$-th leading principle submatrix is the submatrix
formed from the elements in the first $k$ rows and columns. 

\begin{lemma} \label{lem:nsd-signs}
  Suppose that the $n$ by $n$ matrix $A$ has all positive entries, and
  all principle submatrices are invertible. Then, $A$ is negative
  subdefinite if and only if the determinant of the $k$-th leading
  principle submatrix has sign $(-1)^{k+1}$, for $1\le k \le n$. 
\end{lemma}
\begin{proof}
  If $A$ is negative subdefinite, then the $k$-th principle submatrix
  is also negative subdefinite, and so it has $k-1$ negative
  eigenvalues, and one positive eigenvalue. Since the determinant is
  the product of the eigenvalues, the result follows. 
  
  Conversely, we prove by induction that the $k$-th leading principle
  submatrix is \nsd. It is trivial for $k=1$, so assume that the
  $k$-th principle submatrix is negative subdefinite. Since the
  eigenvalues of the $(k+1)$-st leading principle submatrix $P$ are
  interlaced by the eigenvalues of the $k$-th, it follows that $P$ has
  at least $k-1$ negative eigenvalues.  Now, by hypothesis the
  determinant of $P$ has sign $(-1)^k$, and so $P$ must have exactly
  $k$ negative eigenvalues, and therefore one positive one. Thus, $P$ is
  negative subdefinite.
\end{proof}

 Note that this argument also shows that if $Q$ is
  \nsd\ and $|Q|$ is non-zero, then all leading principle
  submatrices have non-zero determinant.  Here is a useful criterion
  for extending \nsd\ matrices.
  \begin{lemma}\label{lem:extend-nsd}
    Suppose that $Q=\smalltwo{A}{u}{u^t}{c}$ where all entries are
    positive, $u$ is $1$ by $n$, $c$ is a scalar, $|Q|\ne0$, and
    $|A|\ne0$. Then, $Q$ is \nsd\ if and only if these
    two conditions are met:
    \begin{enumerate}
    \item $A$ is \nsd.
    \item $uA^{-1}u^t - c >0$
    \end{enumerate}
  \end{lemma}
  \begin{proof}
    We use the Schur complement formula \cite{horn-johnson-1}*{page
      22}:
    $$
    |Q| = |A| \,|c - uA^{-1}u^t|$$
    By hypothesis $|c-uA^{-1}u^t|$ is negative, so $|Q|$ and $|A|$ have
    opposite signs, and the conclusion now follows from
    Lemmas~\ref{lem:nsd-signs} and \ref{lem:nsd-1}.
  \end{proof}

\begin{example}
  It is not difficult to find  negative subdefinite
  matrices. Consider the matrix\footnote{In the
    statistics literature this is known, after scaling, as the
    equicorrelation matrix.} $aJ_n+bI_n$ where $J_n$ is the $n$ by $n$
  matrix of all $1$'s, and $I_n$ is the identity matrix. For example,
  $J_4-I_4$ is
$$
\begin{pmatrix} 
0 & 1 & 1 & 1 \\ 
1 & 0 & 1 & 1 \\ 
1 & 1 & 0 & 1 \\
1 & 1 & 1 & 0
\end{pmatrix}
$$
The eigenvalues of $aJ_n + bI_n$ are $na+b$ and $n-1\,b$'s.  If we
choose $b<0$, $a>0$ and $a+b>0$ then $aJ_n+bI_n$ has exactly one
positive eigenvalue and all positive entries.  The matrices $J_n-I_n$
are limits of negative subdefinite matrices.

 We can use Lemma~\ref{lem:extend-nsd} to construct \nsd\
 matrices. For example, the following matrix 

$$ \begin{pmatrix}  1&2&3&4&5 \\   2&2&3&4&5 \\ 3&3&3&4&5 \\4&4&4&4&5 \\  5&5&5&5&5
\end{pmatrix}
$$
is \nsd\ since all entries are positive, and the determinants of the
leading principle submatrices are ${1, -2,3,-4,5}$.

\end{example}

  We can embed negative subdefinite matrices in higher dimensional
  matrices.
  \begin{lemma}
    If $Q$ is a $d$ by $d$ negative subdefinite matrix, then
    $\smalltwo{Q}{0}{0}{0}$ is the limit of $d+1$ by $d+1$ negative
    subdefinite matrices.
  \end{lemma}
  \begin{proof}
    Write $Q = \smalltwo{A}{v}{v^t}{c}$ where $A$ is $d-1$ by
    $d-1$. Consider the matrix 
$$ M = 
\begin{pmatrix}
  A & v & \alpha v \\v^t & c & \alpha c \\ \alpha v^t & \alpha c & e
  
\end{pmatrix}
$$

If we subtract $\alpha$ times the $d$th row from the bottom row we see
that

$$
\begin{vmatrix}
  A & v & \alpha v \\v^t & c & \alpha c \\ \alpha v^t & \alpha c & e
  
\end{vmatrix} =
\begin{vmatrix}
  A & v & \alpha v \\v^t & c & \alpha c \\ 0  & 0 & e - \alpha^2 c
\end{vmatrix} = (e - \alpha^2 c)|Q| 
$$
We choose positive $\alpha,e$ so that $e - \alpha^2 c<0$. Now $M$ has
all positive entries, and its eigenvalues are interlaced by the
eigenvalues of $Q$. Thus, $M$ has at least one positive eigenvalue, and
at least $d-1$ negative eigenvalues. Since the determinant of $Q$ is
the opposite sign from the determinant of $M$, it follows that $M$ and
$Q$ have a different number of negative eigenvalues, so $M$ had $d$
negative eigenvalues, one positive one, and so is a negative
subdefinite matrix. Taking the limit as $e$ and $\alpha$ go to zero
gives the conclusion.
  \end{proof}

  What do cones of negative subdefinite matrices look like?
  Equivalently, when are all positive linear combinations of a
  collection of negative subdefinite matrices also negative
  subdefinite? Here is an answer for two matrices.

\begin{lemma}
  Suppose that
  \begin{enumerate}
  \item $Q_1$ is negative subdefinite and invertible.
  \item $Q_2$ has all non-negative entries.
  \item $Q_1^{-1}Q_2$ has no negative eigenvalues.
  \end{enumerate}
then
\begin{enumerate}
\item $aQ_1+bQ_2$ is negative subdefinite for all positive $a,b$.
\item $Q_2$ is the limit of negative subdefinite matrices.
\end{enumerate}
\end{lemma}
\begin{proof}
We claim that if $a,b$ are positive then $aQ_1+bQ_2$ is never
singular. If it were then there is a $\nu$ such that
$$aQ_1\nu + b Q_2\nu=0$$ 
which implies that 
$$Q_1^{-1}Q_2\nu = -\frac{a}{b}\nu$$
But $Q_1^{-1}Q_2$ has no negative eigenvalue, so this isn't possible. 
Thus, as $a,b$ vary over positive reals the eigenvalues never change
sign. Since $a=1$ and $b=0$ yields $Q_1$ which has exactly one
positive eigenvalue and no zero eigenvalues, all positive linear
combinations $aQ_1+bQ_2$ have exactly one positive eigenvalue, and are
invertible. Since all entries of $Q_1$ are positive, and of $Q_2$ are
non-negative, the entries of $aQ_1+bQ_2$ are all positive. This
establishes (1). The second part follows by taking $a\rightarrow0$
and $b=1$.
\end{proof}

\begin{cor}
  Suppose $Q_1,Q_2,\dots,Q_n$  are negative subdefinite invertible
  matrices.
  \begin{enumerate}
  \item $a Q_1 + bQ_2$ is negative subdefinite for all positive $a,b$
    iff $Q_1^{-1} Q_2$ has no negative eigenvalues.
  \item If $Q_1,\dots,Q_n$ lie in a cone of negative subdefinite
    matrices then $Q_i^{-1}Q_j$ has no negative eigenvalues for $1\le
    i,j\le n$.
  \end{enumerate}
\end{cor}

\begin{remark}
  If $Q$ is an invertible negative subdefinite matrix, then the
  eigenvalues of $Q^n$ are the $n$th power of the eigenvalues of $Q$.
  Thus, if $n$ is an integer then $Q^{2n+1}$ is also negative
  subdefinite.  Moreover, if $m$ is an integer then the eigenvalues of
  $(Q^{2m+1})^{-1}(Q^{2n+1}) = Q^{2n-2m}$ are all positive. Thus, 
$$ a Q^{2n+1} + b Q^{2m+1}$$
is negative subdefinite.

\index{subdefinite matrix!negative}
\index{matrix!negative subdefinite}
\index{negative subdefinite matrix}%

Here is an example of a cone in dimension two. Define

\begin{align*}
  Q_{a,b} &= 
  \begin{pmatrix} a &  b \\  b &  a   \end{pmatrix} 
\end{align*}

As long as $ 0<a<b$, $Q_{a,b}$ is negative subdefinite. Since
$0<a^\prime < b^\prime$ implies $0 <a + t a^\prime < b+ t b^\prime$
for any positive $t$ it follows that any positive linear combination
of these matrices is still negative subdefinite.  
\end{remark}

The following property will be used in the next section. Recall that
negative \emph{semi}definite means that there are no positive eigenvalues.

\begin{theorem}[\cite{rao}] \label{thm:rao}
Suppose that $A$ is a real symmetric matrix that is not negative
semidefinite. The following are equivalent:
\begin{enumerate}
\item $A$ is negative subdefinite.
\item For every vector $\xx$, $\xx A \xx^t >0$ implies that $A\xx$ is a
  non-positive or non-negative vector.

\end{enumerate}
\end{theorem}

The following lemma is used in the next section.

\begin{lemma} \label{lem:subdef-3}
  If $Q = \smalltwo{a}{v}{v^t}{C}$ is negative subdefinite, then
  $v^t v-aC$ is positive semidefinite.
\end{lemma}
\begin{proof}
  Assume that $v^t v-aC$ is not positive semidefinite. Then, there is
  a $z$ such that $z(v^tv-aC)z^t<0$. This implies that $ zCz^t >
  \frac{1}{a}zv^tvz^t$ . Evaluate the \index{quadratic form}quadratic form at $w=
  (\alpha,z)$ where $\alpha$ is to be determined.
  \begin{align*}
    (\alpha,z)\begin{pmatrix} a& v\\ v^t & C \end{pmatrix}
    \begin{pmatrix} \alpha \\ z^t \end{pmatrix} & = 
\alpha^2 a + \alpha(zv^t + vz^t) + zCz^t \\
&> \alpha^2 a + \alpha(zv^t + vz^t) + \frac{1}{a} zv^tvz^t \\
&= \frac{1}{a}(\alpha a + zv^t)^2
  \end{align*}
If $\alpha a + zv^t\ne0$ then $wQw^t>0$, so by Theorem~\ref{thm:rao} $Qw^t$ has
either all non-negative, or all non-positive, coefficients. But
$$ Qw^t = \left( a \alpha + v z^t, \alpha v^t + Cz^t\right).$$
Since $\alpha a+vz^t$ can be either positive or negative, as $\alpha
a+vz^t$ becomes zero, so must $\alpha v^t +Cz^t$. Thus, when $\alpha =
-\frac{1}{a}vz^t$ we have
$$ (-\frac{1}{a}vz^t) v^t + Cz^t =  -\frac{1}{a}(v^tv
  - aC) z^t = 0$$
This contradicts our hypothesis that $z(v^tv-aC)z^t<0$, and so the
lemma is proved.
\end{proof}

\section{Quadratic Forms}
\label{sec:sub-quadratic-forms}

We have seen substitution of some quadratic forms before. For
instance, if $f\in\allpolypos$ then $f(-x^2)\in\allpoly$. Also, if
$f\in\allpolypos$, then $f(-xy)\in\gsubclose_2$. We will replace
$-x^2$ and $-xy$ by quadratic forms determined by negative subdefinite
matrices. We begin with substitution into the simplest polynomial: $x+c^2$:

\begin{quote}
  When is $-\xx Q \xx^t + c^2$ in $\pm\rupint{2}$?
\end{quote}

Let's first consider the case of a \index{quadratic form}quadratic form in two variables.
Take a matrix $A=\smalltwo{a_1}{a_2}{a_2}{a_3}$ with corresponding
\index{quadratic form}quadratic form $a_1x^2 + 2a_2xy + a_3y^2$.  We are interested when the
polynomial
$$ \begin{pmatrix} x & y \end{pmatrix}
\begin{pmatrix} a_1 & a_2 \\ a_2 & a_3 \end{pmatrix}
\begin{pmatrix} x \\ y \end{pmatrix} - c^2 =
a_1 x^2 + 2a_2 xy + a_3 y^2 - c^2 
$$
is in $\gsubpos_2$ for all $c$. We must have that all $a_i$ are
positive, and the discriminant should be non-negative. Thus 
$$
4a_2^2 - 4 a_1a_3 +4c^2 = - \smalltwodet{a_1}{a_2}{a_2}{a_3}+4c^2
\ge0$$
and hence $|A|\le0$.  $A$ has a positive eigenvalue, since the
sum of the eigenvalues is the trace which is positive.  Consequently,
$\xx A \xx^t - c^2\in\gsubpos_2$ iff $A$ is negative subdefinite.

\begin{lemma} \label{lem:subdef-1}
  If $Q$ is a $d$ by $d$ symmetric matrix such that $\xx Q \xx^t
  -c^2 \in\pm\gsubpos_d$ for all $c$ then
  \begin{enumerate}
  \item $Q$ has all positive entries.
  \item $Q$ has exactly one positive eigenvalue.
  \item $Q$ is negative subdefinite.
  \end{enumerate}
\end{lemma}
\begin{proof}
  Since $\xx Q \xx^t -c^2 \in\pm\gsubpos_d$ the homogeneous part is in
  $\gsubpos_d$, and so the entries of $\xx Q \xx^t$ are all the same
  sign, and are non-zero. If they are negative then letting
  $Q=(a_{ij})$ and $\xx=(x,0,\dots,0)$ yields $\xx Q \xx^t -c^2=
  a_{11}x^2 - c^2$. This has all real roots iff $a_{11}$ is positive,
  so all entries of $Q$ are positive.
  
  For the second part we only need to assume that $\xx
  Q\xx^t\in\gsubclose_d$. 
  \index{Perron's theorem}  Recall Perron's theorem
  \cite{horn-johnson-1} that says that since all entries of $Q$ are
  positive there is a unique largest positive eigenvalue $\lambda$,
  and the corresponding eigenvector $v=(v_1,\dots,v_d)$ has all
  positive entries. Suppose that $\mu$ is an eigenvalue with
  eigenvector $u=(u_1,\dots,u_d)$. Since $Q$ is symmetric and
  $\lambda\ne\mu$ it follows that $u^t Q v=0$. If we set $f(\xx)
  = \xx Q \xx^t$ then we can replace each $x_i$ by $v_ix+u_i$ and
  the resulting polynomial is in $\allpoly$ since all $v_i$ are
  positive. Consequently,
  $$
  f(xv+u) = (xv+u) Q (xv+u)^t = \lambda x^2 |v|^2 + \mu |u|^2
  \in\allpoly$$
  The only way that this can happen is if $\mu$ is
  negative or $0$, so there is exactly one positive eigenvalue.
\end{proof}

\begin{lemma} \label{lem:subdef-2}
  If $Q$ is negative subdefinite then $\xx Q \xx^t -c^2
  \in\gsubpos_d$ for all $c$.
\end{lemma}
\begin{proof}
  The homogeneous part of $\xx Q\xx^t-c^2$ is a polynomial with all
  positive coefficients since $Q$ has all positive coefficients. Thus,
  we only have to check substitution.

  Write $Q = \smalltwo{a}{v}{v^t}{C}$ where $v$ is $1\times (d-1)$ and 
  $C$ is $(d-1)\times (d-1)$, and let $\xx = (x_1,\zz)$ where $\zz =
  (x_2,\dots,x_d)$. With this notation
  \begin{align*}
    \xx Q \xx^t &= (x_1,\zz) \begin{pmatrix} a & v \\ v^t & C
    \end{pmatrix} \begin{pmatrix} x \\ \zz^t \end{pmatrix} \\
    & = ax_1^2 + x_1(v\zz^t + \zz v^t) + \zz C \zz^t
  \end{align*}

  $\xx Q\xx^t$ satisfies $x_1$ substitution if and only if the
  discriminant is non-negative.  Now $v\zz^t=\zz v^t$ since they are
  both scalars, so the discriminant condition is that
  \begin{align*}
    0 &\le 4(v\zz^t)^2 - 4a\zz C\zz^t \,=\,  4\zz ( v^tv - aC)\zz^t
  \end{align*}
  $v^tv-aC$ is a symmetric matrix, so $\xx Q\xx^t$ satisfies
  $x_1$-substitution iff $v^tv-aC$ is positive semidefinite.  Now
  apply Lemma~\ref{lem:subdef-3}.
\end{proof}

\section{The analog of $x\diffd$}
\label{sec:xdxydy}

The differential operator $\xx\cdot\boldsymbol{\partial}=\xdxydypd$ is
the analog in $\gsubplus_d$ of the operator $x\frac{d}{dx}$ on
$\allpolypos$.  Here are some elementary properties of this operator.

\begin{enumerate}
\item If $f(\xx)$ is homogeneous of degree $n$ then $(\xdxydypd)f = n
  f$.  This is the well-known Euler identity about homogeneous
  functions.
\item If $f\in\rupint{d}(n)$ then $\left((\xdxydypd)f\right)^H = n
  f^H$. This is an immediate consequence of the previous property.
\item If $g(\xx) = \sum a_\sdiffi \xx^\diffi$ then 
  $$
  f(\xdxydypd)g = \sum a_\sdiffi f(|\diffi|) \xx^\sdiffi$$
  This is
  a simple consequence of linearity and the  calculation
  $$(\xx\cdot\boldsymbol{\partial})\xx^\diffi =
  (\xdxydypd)(x_1^{i_1}\cdots x_d^{i_d}) =
  (i_1+\cdots+i_d)x_1^{i_1}\cdots x_d^{i_d}$$
\end{enumerate}

 The following is the basic fact about $\xx\cdot\boldsymbol{\partial}$.

\begin{prop}\label{prop:xd-pdpos}
  If $f\in\gsubplus_d(n)$ then $(\xdxydypd)f\in\gsubplus_d$. Moreover,
  if $g$ is the polar derivative of $f$ then
$$ n f = (\xdxydypd)f + g.$$
  We also have interlacings 
$$
f \greateqeq (\xdxydypd)f  \lesslesseq  g 
$$

\end{prop}
\begin{proof}
  These are immediate consequences of the properties of the polar
  derivative; see \eqref{eqn:polar-pd} and \eqref{eqn:polar-pd-2}.

 In the case that $d$ is two we can derive them in a different way.   Notice that
$$ \frac{\partial}{\partial z} \, f(-xz,-yz) = -(x f_x(-xz,-yz) + y
f_y(-xz,-yz))$$
Since the derivative is in $\gsubplusclose_3$, if we substitute $z=-1$
we get a polynomial in $\gsubplusclose_2$, which shows that
 $x f_x + y f_y\in\gsubplusclose_2$. Since $x f_x + y f_y$ satisfies
 the  homogeneity condition, it follows it is in
 $\gsubplus_2$. 
 
 Finally, since the derivative of $f(-xz,-yz)$ interlaces $f(-xz,-yz)$,
 and is still true when we substitute $z=-1$, the interlacing part now
 follows. It's interesting to see  that
 $f(\xdxydy)g = f(z\partial_z)g(xz,yz)\big{|}_{z=1}.$

\end{proof}

\begin{cor}\label{cor:xd-pdpos}
  If $f\in\allpolypos$ and $g=\sum \aaa_\sdiffi \xx^\diffi\in\gsubplus_d$
  then $$ \sum_\diffi \aaa_\sdiffi f(|\diffi|)\xx^\diffi\in\gsubplus_d.$$
\end{cor}
\begin{proof}
  $\alpha g + \xx\cdot\boldsymbol{\partial} (g) \in\gsubplus_d$ by the
  above proposition, so
  $(\xx\cdot\boldsymbol{\partial}+\alpha)(g)\in\gsubplus_d$ for positive
  $\alpha$. Now factor $f$, and use induction.
\end{proof}

\section{Generalized Hadamard  products}
\label{sec:hadamard-pd}

 \index{Hadamard product!general}
A generalized Hadamard product is a mapping of the form
$$\text{monomial } \times \text{ monomial } \mapsto \text{ constant }
\times \text{ monomial}
$$
This generalizes multiplier transformations which have
  the form
$$ \text{monomial} \mapsto \text{ constant } \times \text{
  monomial}$$

We have seen two Hadamard type products 

\begin{alignat*}{2}
  Theorem~\ref{thm:hadamard-1} &\quad & x^i \ast x^j & \mapsto 
  \begin{cases}
    x^i & i=j \\ 0 & i\ne j
  \end{cases}\\
  Theorem~\ref{thm:hadamard-2} &\quad & x^i \ast^\prime x^j & \mapsto 
  \begin{cases}
    i!\,x^i & i=j \\ 0 & i\ne j
  \end{cases}
\end{alignat*}

Both of these determine bilinear maps $\allpolypos \times \allpoly
\longrightarrow \allpoly$. In Proposition~\ref{prop:hadamard-prod-gen} we will determine sufficient
conditions on constants $c_i$ for a product of the form
$$
   x^i \hadprod x^j  \mapsto 
  \begin{cases}
    c_i\,x^i & i=j \\ 0 & i\ne j
  \end{cases}
  $$
  to determine a map $\allpolypos \times \allpoly \longrightarrow
  \allpoly$.  We now study two products for $\allpoly \times
  \rupint{2} \longrightarrow \allpoly$. Define

\begin{align}
\label{eqn:hadamard-type-1}  y^i \ast x^jy^k & \mapsto 
  \begin{cases}
    x^j & i=k \\ 0 & i\ne j
  \end{cases}\\
\label{eqn:hadamard-type-2}  y^i \ast^\prime x^jy^k & \mapsto 
  \begin{cases}
    i!\,x^j & i=k \\ 0 & i\ne k
  \end{cases}
\end{align}

\begin{prop}\label{prop:hadamard-type-1}
  The linear transformations in \eqref{eqn:hadamard-type-1}
  and \eqref{eqn:hadamard-type-2} determine mappings
  $\allpoly\times\rupint{2}\longrightarrow\allpoly$.  
\end{prop}
\begin{proof}
  Choose $f(x,y)=\sum f_i(x)y^i$ in $\rupint{2}$,  and $g = \sum_0^n a_i
  y^i$ in $\allpoly$. Since $\sum a_{n-i}y^i\in\allpoly$, the product 
  $$
  \left(\sum f_i(x) y^i \right)\left(\sum a_{n-j}y^j\right)$$
  is in
  $\gsubclose_2$. The coefficient of $y^n$ in the product equals $\sum a_i
  f_i(x)$ and is in $\allpoly$.  This is exactly the product determined by
  \eqref{eqn:hadamard-type-1}.

  The second product is a strengthening of the first since we can
  apply the exponential operator to \eqref{eqn:hadamard-type-2} to
  obtain \eqref{eqn:hadamard-type-1}. Note that
  $\left(g(\frac{\partial}{\partial y})y^i\right)(0) = i!a_i$. Since
  $g(\frac{\partial}{\partial y})f$ is in $\rupint{2}$ by
  Lemma~\ref{lem:pd-diff-prod}, the conclusion follows from the observation that
$$ \left(g(\frac{\partial}{\partial y})f\right)(x,0) = \sum
\frac{f_i(x)}{i!} \left(g(\frac{\partial}{\partial y})y^i\right)(0) = 
\sum f_i(x)a_i .
$$
  
\end{proof}

The same argument shows that

\begin{cor}\label{cor:hadamard-type-2}
  If we define
$$
   y^i \ast^\prime \xx^\diffj \,y^k  \mapsto 
  \begin{cases}
    i!\,\xx^\diffj & i=k \\ 0 & i\ne k
  \end{cases}
  $$ 
  then this determines a map $\allpoly \times \rupint{d+1}
  \longrightarrow \rupint{d}$.
\end{cor}

\begin{cor}
  Suppose $f(\xx,y) = \sum f_i(\xx)y^i$ is in $\rupint{d+1}$. The
  following two sums are in $\gsubclose_d$:
$$
    \sum (-1)^i \frac{f_{2i}(\xx)}{i!} \quad\quad
    \sum (-1)^i {f_{2i}(\xx)} 
$$  
\end{cor}
\begin{proof}
  We can take limits to see that Corollary~\ref{cor:hadamard-type-2} determines a
  map $\allpolyf \times \rupint{d+1}\longrightarrow\rup{d}$. These
  results now follow if we consider $e^{-y^2}\ast f$ and
  $e^{-y^2}\ast^\prime f$. 
\end{proof}

The reverse of a polynomial in $\rupint{d}$ was defined in
\eqref{eqn:reverse-in-pd}.
\begin{cor}\label{cor:hadamard-type-3}
  If $f = \sum f_i(\xx)y^i$ and $g = \sum g_i(\xx)y^i$ are both in
  $\rupint{d+1}(n)$ then
\begin{enumerate}
\item    $ \sum f_i(\xx)g_{n-i}(\xx) \in\rupint{d} $
\item   $ \sum f_i(\xx)g_{i}^{rev}(\xx) \in\gsubclose_d $
  \end{enumerate} 
\end{cor}
\begin{proof}
  For the first one we multiply and extract a coefficient as
  usual. For the second, apply (1) to $f$ and $g^{rev} = \sum
  g_i^{rev}(\xx) y^{n-i}$.
\end{proof}

  \section{THE BMV conjecture}
  \label{sec:bmv-conjecture}

\added{6/3/7}
  The Bessis-Moussa-Villani (BMV) conjecture, as reformulated by Lieb
  and Seiringer\cite{lieb-seiringer}, is

  \begin{conj}[BMV]
    If $A,B$ are positive definite matrices then\\
    $trace(I+(A+t\,B)^n)$ has all positive coefficients for all
    positive integers $n$.
  \end{conj}


  We offer a stronger conjecture.

\begin{conj}[s-BMV]
  If $A_i,B_i$ are positive definite matrices,  $g_i\in\allpolypos$,
  and $f(x,y_1,\dots,y_e)\in\rupint{1+e}$ then
     
     \begin{enumerate}
\item $det\bigl[f(x\,I + g_1(A_1+B_1y)+\cdots +g_e(A_e+B_ey))\bigr] $ has all positive coefficients.
\item  If we write the determinant as $\sum f_i(x)y^i$ then
     all $f_i\in\allpolypos$.
   \item In addition we have $\cdots \longleftarrow f_i
     \longleftarrow f_{i+1} \longleftarrow \cdots$. 
   \end{enumerate}
 \end{conj}

 If $C$ is a $d\times d$ matrix then the trace of $C$ is the
 coefficient of $x^{d-1}$ in the polynomial $|xI+C|$. Replacing $x$ by
 $x+1$ and applying this fact shows that the strong BMV conjecture
 applied to $det(xI + I+ (A+By)^n)$ 
 implies the original BMV conjecture.

 We can take limits and let the $g_i$ belong to $\allpolyposf$. A case
 of particular interest is the exponential BMV conjecture:

\begin{conj}[e-BMV]
    If $A_i,B_i$ are positive definite matrices,  
       \begin{enumerate}
       \item  $det\bigl[xI + e^{A_1+B_1y} + \cdots + e^{A_e+B_ey}\bigr]$ has all positive coefficients.
\item  If we write the determinant as $\sum f_i(x)y^i$ then
     all $f_i\in\allpolypos$.
   \item In addition \ $\cdots \longleftarrow f_i
     \longleftarrow f_{i+1} \longleftarrow \cdots$. 
       \end{enumerate}
       
\end{conj}

It is not the case that $|xI+(A+By)^n|\in\rupint{2}$ or that
$|xI+e^{A+By}|\in\gsubf_2$. For example, if we take $A=B=I$ then
the latter determinant is $(x+e^{1+y})^d$ which is a polynomial in
$x$, but an exponential polynomial in $y$. 

The leading coefficient of 
$|xI+g_1(A_1+B_1y)+\cdots+g_e(A_e+B_ey)|$ is $1$.  We have partial
information about the constant term.

\begin{lemma}
  If $A,B$ are positive definite and $g\in\allpolypos$ then the
  constant term (with respect to $x$) of $det\bigl[xI + g(A+By)\bigr]$
  has all positive coefficients.
\end{lemma}
\begin{proof}
  The constant term is $det\bigl[g(A+By)\bigr]$. If we write $g(x) =
  \prod(x+a_i)$ where the $a_i$ are positive then
\[
det\bigl[ g(A+By)\bigl] = \prod det\bigl[ a_iI +A +By\bigr] 
\]
This has all positive coefficients since $a_iI+A$ and $B$ are positive definite.
\end{proof}

There is one simple case where s-BMV(1) and e-BMV(1) hold.

\begin{lemma}
  If $A_1,B_1,\dots$ are commuting positive definite matrices then s-BMV(1) and
  e-BMV(1) hold. 
\end{lemma}
\begin{proof}
  Let ${O}$ simultaneously diagonalize all the $A_i$'s and $B_i$'s, so that $OA_iO^t =
  diag(a_{i,j})$ and $OB_iO^t = diag(b_{i,j})$ where the $a_{i,j}$'s and $b_{i,j}$'s are
  positive. Then
\begin{align*}
|xI+g_1(A_1+B_1y)+\cdots+g_e(A_e+B_ey)| &= \prod_j \bigl(x + \sum_i g_i(a_{i,j}+b_{i,j}y)\bigr) 
\end{align*}
All coefficients are positive since $g\in\allpolypos$, and 
$a_{i,j},b_{i,j}>0$. The exponential case follows by taking limits.
\end{proof}

It is not even clear that s-BMV(2) and s-BMV(3)  hold in this simple case. It would
follow from the following  

\begin{conj}
  Suppose $g\in\allpolypos$ and $a_i,b_i$ are positive  and $n$ is a
  positive integer and write
\[ 
\prod_{i=1}^n \bigl( x + g(a_i+b_i\,y)\bigr) = \sum f_i(x)y^i.
\]
 All $f_i\in\allpolypos$ and
$f_i\longleftarrow f_{i+1}$.

\end{conj}

This is a special case of the composition conjecture, which we discuss
in the next section.

  \section{The composition conjecture}
\added{6/7/7}
  \label{sec:subst-conj}

We investigate some special cases of the \emph{composition
  conjecture}

\begin{conj}\label{conj:sub}
  Suppose that $f(\xx,y_1,\dots,y_e)\in\gsubf_{d+e}$, and choose
  $g_i(y)\in\allpolyposf$. If we write
\[
F(\xx,y) = f(\xx,g_1(y),\dots,g_e(y)) = \sum f_i(\xx)y^i
\]
then 
\begin{enumerate}
\item All $f_i$ are in $\gsubclose_d$.
\item $f_0 \longleftarrow f_1 \longleftarrow f_2 \longleftarrow \cdots $
\end{enumerate}
\end{conj}

\begin{remark}
  It is not the case that $F\in\gsubf_{d+1}$; take $f=x+y$ and
  $g(y)=e^y$.  However, if $\alpha\in\reals$ then we clearly have
  $F(\xx,\alpha)\in\gsubposf_d$, since we are merely substituting
  $g_i(\alpha)$ for $y_i$. Thus, $F(\xx,y)$ satisfies substitution for
  $y$, not substitution for $\xx$, and the coefficients of $y$
  interlace.

  If $f\in\gsubposf_{d+e}$ and the composition conjecture holds then  all coefficients of $F(\xx,y)$
  are positive.
\end{remark}

We first show if we have one  exponential function then
the composition conjecture holds. The proof is for $d=1$, but an
identical argument shows it for any $d$.

  \begin{lemma}
    Suppose $f\in\rupint{2}$, and $g\in\allpolyposf$. If the coefficient
    $h_n$ of $y^n$ in $e^{-x g(y)}$ is in $\allpoly$ then the
    coefficient $k_n$ of $y^n$ in $f(x,g(y))$ is in $\allpoly$.

    If $h_n\longleftarrow h_{n+1}$ then $k_n\longleftarrow h_{n+1}$. 
  \end{lemma}
  \begin{proof}
    Let $T_g$ be the linear transformation
\[ f \mapsto \text{coefficient of $y^n$ in $f(x,g(y))$}. \]
The generating function is
\begin{align*}
  G(x,y,u,v) &= \sum T_g(x^iy^j)\frac{(-u)^i(-v)^j}{i!j!} \\
  &= \sum \bigl[ \text{Coefficient of $y^n$ in $x^ig(y)^j$}\bigr]
  \frac{(-u)^i(-v)^j}{i!j!} \\ 
&= e^{-ux} \bigl[\text{Coefficient of $y^n$ in $e^{-vg(y)}$}\bigr].
\end{align*}
By hypothesis this is in $\gsubf_3$, so $T_g$ maps $\rupint{2}$ to
$\allpoly$.

The second part is similar, and considers the map
\[ T_{g,\alpha} \mapsto
\bigl[\text{coefficient of $y^n$ in $f(x,g(y))$}\bigr] 
+ \alpha 
\bigl[\text{coefficient of $y^{n+1}$ in $f(x,g(y))$}\bigr] 
\]
  \end{proof}

  \begin{cor}
    If $f\in\rupint{2}$ then the coefficient $k_n$  of $y^n$ in $f(x,e^y)$ is in $\allpoly$.
Moreover,   $k_n\longleftarrow k_{n+1}$. 
\end{cor}
  \begin{proof}
    Recall the identity\index{Bell polynomials}
\[ 
e^{xe^y} = e^{x} \sum_{i=0}^\infty B_i(x) \frac{y^i}{i!}
\]
where $B_i$ is the Bell polynomial. Since the Bell polynomials are in
$\allpoly$ and consecutive ones interlace, the result follows from the
preceding lemma.
  \end{proof}

The following conjecture arises naturally in the lemma following.

\begin{conj}\label{conj:lt}
  Suppose that $\alpha_1,\dots,\alpha_e$ are positive and $n$ is a
  positive integer. The mapping
\[
\xx^\sdiffi y_1^{t_1}\cdots y_e^{t_e} \mapsto 
\xx^\sdiffi (\alpha_1 t_1 + \cdots + \alpha_e t_e)^n
\]
determines a map $\gsubplus_{d} \longrightarrow \gsubplus_{d+e}$. 
\end{conj}

\begin{lemma}
  Choose positive $a_i$ and let $g_i(y) = e^{a_iy}$. If
  Conjecture~\ref{conj:lt} is true then the composition
  conjecture(part 1) holds for these choices.
\end{lemma}

\begin{proof}
  For simplicity of exposition we take $e=2$. 
  Write 
\[ f(\xx,\yy) = \sum a_{\sdiffi rs} \xx^\sdiffi y_1^r y_2^s  \]
Substituting for $g_i$ yields
\[ F(\xx,y) = \sum a_{\sdiffi rs} \xx^\sdiffi e^{a_1yr} e^{a_2ys}  
= \sum a_{\sdiffi r s}  \xx^\sdiffi e^{(a_1r+ a_2s)y}\]
We find the coefficient $g_n(\xx)$ of $y^n$ by differentiating $n$ times and
letting $y=0$
\[ 
g_n(\xx) = \sum a_{\sdiffi rs} \xx^\sdiffi (a_1r + a_2s)^n   
\]

and this is in $\gsubplus_d$ by the conjecture.
\end{proof}

Unfortunately, consider 
\begin{example}
  The operator $ a_1\,y_1\frac{\partial}{\partial y_1}  + a_2\,y_2\frac{\partial}{\partial
  y_2}$ does not necessarily map $\gsubplus_{d+e}$ to itself. Consider
\begin{align*}
  f & = (x + y + z + 2) (x + y + 2 z + 1)\\
  g=yf_y+2\,z\,f_z &= \bigl[y(2 y+2 x +3 )\bigr] + \bigl[6 x+9 y+10\bigr]z + 8z^2
\end{align*}
The constant term of $g$ is two lines that meet at
$(-3/2,0)$. However, the coefficient of $z$ is a line that does not
pass through this point, so adjacent coefficients do not intersect, and
hence $g\not\in\rupint{3}$.
\end{example}

\begin{remark}
  It is interesting to compute the generating function of the linear
  transformation
\[ T_n\colon x^iy^rz^s \mapsto x^i(\alpha r+ \beta s)^n \]
We know it does not define a map $\rupint{3}\longrightarrow\allpoly$, but
it does appear to satisfy $\gsubplus_3\longrightarrow\allpoly$.
\begin{align*}
  G(x,y,z,u,v,w) &= \sum T_n(x^iy^rz^s) \frac{u^iv^rw^s(-1)^{i+r+s}}{i!r!s!} \\
&= \sum x^i (\alpha r + \beta s)^n \frac{u^iv^rw^s(-1)^{i+r+s}}{i!r!s!} \\
&= e^{-xu} \sum _{k=0}^n \alpha^k \beta^{n-k} \binom{n}{k} 
 \sum_{r,s} \frac{r^k(-v)^r\, s^{n-k}(-w)^s}{r!s!}\\
&= e^{-xu} \sum _{k=0}^n \alpha^k \beta^{n-k} \binom{n}{k} 
\bigl( \sum \frac{r^k (-v)^r}{r!}\bigr)
\bigl( \sum \frac{s^{n-k} (-v)^r}{s!}\bigr)\\
&= e^{-xu-v-w}\sum _{k=0}^n \alpha^k \beta^{n-k} \binom{n}{k} B_k(-v)B_{n-k}(-w)
\end{align*}
where $B_k$ is the Bell polynomial. If $\alpha=\beta=1$ then this
simplifies to
\[ e^{-xu-v-w} B_{n}(-v-w)\]
which is in $\gsubf_4$, but it does not appear to be in $\gsubf_4$
if $\alpha\ne\beta$.
\end{remark}

\subsection*{The \emph{\dersub} class}
\label{subsec:exmpd-subst-class}

It appears that derivatives of compositions satisfy a 
 a much stronger condition than interlacing.

\begin{definition}
  A polynomial $f(x,y)$ satisfies \emph{\dersub} if
  for $m=0,1,\dots$ and $\alpha\in\reals$
  \begin{enumerate}
 \item $\displaystyle\frac{\partial^m f}{\partial y^m} (x,\alpha) \in\allpoly$
 \item $\displaystyle\frac{\partial^m f}{\partial y^m} (x,\alpha)$ and
   $\displaystyle\frac{\partial^{m+1} f}{\partial y^{m+1}} (x,\alpha) $ interlace.
  \end{enumerate}
\end{definition}

Clearly polynomials in $\rupint{2}$ satisify \dersub. The next lemma
gives some of the properties of this class.

\begin{lemma}
  If $f(x,y)$ satisfies \dersub\ and $f(x,y) = \sum f_i(x)y^i$ then
  for $m=0,1,\dots$ 
  \begin{enumerate}
  \item $f_m\in\allpoly$.
  \item $f_m$ and $f_{m+1}$ interlace.
  \item $f(x,y+\gamma)$ satisfies \dersub\ for all $\gamma\in\reals$.
  \end{enumerate}
\end{lemma}

\begin{proof}
  Differentiating $m$ times with respect to $y$ and setting $y=0$
  shows part $1$. From the second condition of the definition we have
  that for all $\beta\in\reals$
\[ 
   \frac{\partial^{m} f}{\partial y^{m}} (x,0) 
+ \beta   \frac{\partial^{m+1} f}{\partial y^{m+1}}
(x,0)\in\allpoly
\]
which implies that  $f_i$ and $f_{i+1}$ interlace.

The third part follows from the observation that
\[
   \frac{\partial^{m} f(x,y+\gamma)}{\partial y^{m}} (x,\alpha) 
=
   \frac{\partial^{m} f(x,y)}{\partial y^{m}} (x,\alpha+\gamma)
\]
 
\end{proof}

We have our final generalization of the BMV conjecture

\begin{conj}\label{conj:fdb-ds}\ 
  \begin{enumerate}
  \item If $f(x,y)\in\rupint{2}$ and $g\in\allpoly$ then $f(x,g(y))$
    satisfies \dersub.
  \item If $A,B$ positive definite, and $g\in\allpoly$ then $|xI +
    g(A+By)|$ satisfies \dersub.
  \end{enumerate}
  
\end{conj}

For the rest of this section we consider $F(x,y) = f(x,g(y))$. 
The \fdb\ formula is an expression for the $m$'th derivative of a
composition. In our case we can write

\[
\frac{\partial^m}{\partial y^k} f(x,g(y)) = 
\sum_k f^{(k)}(x,g(y)) \, A_{m,k}(y)
\]

Observe that this is the constant term of

\[
\left(\sum_{k=0}^m A_{m,k}(y) \,\diffd_y^k\right)\cdot
\left(\sum_{k=0}^\infty f^{(k)} (x,g(y)) \frac{z^k}{k!}\right)
\]
The right hand term is $f(x,z+g(y))$, and for $\alpha\in\reals$ we see
that $f(x,z+g(\alpha))$ is in $\rupint{2}$.  Define
\[ A_m(x,y) = \sum A_{m,k}(x) y^k \]
In order to show that
$F$ satisfies \dersub\ we only need to show that $A_m(\alpha,y)\in\allpoly$, and
that $A_m(\alpha,y)\lesslesseq A_{m+1}(\alpha,y)$. 

We can do a few special cases.

\index{Bell polynomial}
\begin{lemma}
  If $f\in\rupint{2}$ then $f(x,e^y)$ satisfies \dersub.  If
  $d$ is a positive integer  then $f(x,x^d)$
  satisfies \dersub.
\end{lemma}
\begin{proof}
These follow from Lemmas~\ref{lem:fdb-ex} and \ref{lem:fdb-xd}.
\end{proof}




\section{Recursive constructions}
\label{sec:sub-recursive}

In this section we see how to construct new polynomials in $\gsubpos_d$
from old ones. We can use these constructions to get infinite
sequences of interlacing polynomials that are analogous to sequences
of orthogonal polynomials.

\begin{lemma} \label{lem:sub-recursive}
  Let $f\lesslesseq g$ be two polynomials in $\gsubpos_d$. Choose
  $b\in\reals$, and a vector $\aaa=(a_1,\dots,a_d)$ of positive real
  numbers. With these assumptions we can conclude that
  \begin{enumerate}
  \item $(b+\aaa\xx^{t})f-g \,\in\gsubpos_d$.
  \item $(b+\aaa\xx^{t})f-g \lesslesseq f$.
  \end{enumerate}
\end{lemma}
\begin{proof}
  Let $h=(b+\aaa\xx^{t})f-g$. Since $h+\alpha f =
  (b+\alpha+\aaa\xx^{t})f-g$ it suffices to show that $h\in\gsubpos_d$.
  All terms of maximal total degree are obtained from $f$ by
  multiplying by various coordinates of $\aaa$, so they are all
  positive. If we substitute for all but one variable then we have to
  show that $(a+cx)\tilde{f}-\tilde{g}$ is in $\allpoly$, where
  $\tilde{f}$ and $\tilde{g}$ are the results for substituting in $f$
  and $g$. Since $\tilde{f}\lesslesseq \tilde{g}$ the result follows
  from Corollary~\ref{cor:xbf}.
\end{proof}

If we iterate the construction of the lemma we can get an infinite
sequence of interlacing polynomials.

\index{orthogonal polynomials!higher analogs}

\begin{cor}
  Choose constants $b_1,b_2,\dots$, positive vectors
  $\aaa_1,\aaa_2,\dots$, and positive constants $c_1,c_2,\dots$. Set
  $f_{-1}=0,f_0=1$ and define $f_n$ recursively by
$$ f_n = (b_n + \aaa_n\xx^{t})f_{n-1} - c_n f_{n-2}.$$
All $f_n$ are in $\gsubpos_d$, and we have interlacings
$$f_0 \lessgreateq f_1 \lessgreateq f_2 \lessgreateq $$
\end{cor}

Since the derivative preserves interlacing, we can use derivatives to
take the place of $f_{n-2}$.

\begin{cor}
  Choose constants $b_1,b_2,\dots$, positive vectors
  $\aaa_1,\aaa_2,\dots$, and positive vectors $\ccc_1,\ccc_2,\dots$.
  Set $f_0=1$ and define $f_n$ recursively by
  $$
  f_n = (b_n + \aaa_n\xx^{t})f_{n-1} - \ccc_n
  \,\left(\frac{\partial}{\partial\xx}f_{n-1}\right)^{t}.$$
  All
  $f_n$ are in $\gsubpos_d$, and we have interlacings
  $$f_0 \lessgreateq f_1 \lessgreateq f_2 \lessgreateq $$
\end{cor} 

\begin{example}
  As an interesting special case, consider the recurrence $f_0=1$ and
  $$
  f_{n+1} = \xx f_n - \partial\xx f_n = (x_1+\cdots+x_d)f_n -
  (\frac{\partial}{\partial x_1}+ \cdots + \frac{\partial}{\partial
    x_d})f_n$$
  If we define one variable polynomials $h_0=1$ and $h_n
  = xh_{n-1} - d h^\prime_{n-1}$ then it is easy to verify that $f_n =
  h_n(x_1+\cdots+x_d)$. The $h_n$ are orthogonal polynomials
  (essentially modified Hermite \index{Hermite polynomials}
  polynomials), and so we can view these $f_n$ as generalized Hermite
  polynomials. In \chapsec{topology}{hermite-n-dim} we will see a true
  generalization of Hermite polynomials to $n$-variables.
\end{example} 

The linear combination $(b+\aaa\xx^t)g+f$ also has interlacing
properties. This generalizes Corollary~\ref{cor:fbxcg}.

\begin{lemma} \label{lem:pd-xffp}
  If $f\greateqeq g$ in $\gsubpos_d$ and $\aaa$ is positive then
$(b+\aaa\xx^t)g+f\lesslesseq g$.
\end{lemma}
\begin{proof}
  By linearity it suffices to show that the left hand side is in
  $\gsubpos_d$. The homogeneous part is $(\aaa\xx^t) g^H$, the degrees
  are correct, and substitution follows from Corollary~\ref{cor:fbxcg}.
\end{proof}

\section{Two by two matrices preserving interlacing in $\rupint{d}$}

In this section we generalize results from
Chapter~\ref{cha:poly-matrices} to $\gsubclose_d$. Many of the
properties of one variable have immediate generalizations to $d$
variables.  However, we don't have a concept of ``mutually
interlacing'' for more than one variable, so Theorems such as
\ref{thm:tnn} are only possible for two by two matrices. 

The following Lemma generalizes Corollary~\ref{cor:lin-comb-new}. The
restrictive assumptions arise from the necessity to control the
homogeneous part.

\begin{lemma} \label{lem:lin-comb-new-pd}
  Suppose that $f_1 \greateq f_2$ (or $f_1 \lessless f_2$) in
  $\rupint{d}$ have positive leading coefficients, and
  $\smalltwo{a}{b}{c}{d}\smalltwobyone{f_1}{f_2} =
  \smalltwobyone{g_1}{g_2}$. If $a,b,c,d$ are positive and if the
  determinant $ad-bc$ is positive then $g_1\greateq g_2$, and if it is
  negative $g_2 \greateq g_1$.
\end{lemma}

\begin{proof}
  The positivity assumptions imply that $g_1,g_2$ have positive
  homogeneous part.  Corollary~\ref{cor:lin-comb-new} implies that
  $g_1,g_2\in\rupint{d}$, and that the interlacings hold as stated.
\end{proof}

As long as all coefficients are positive then the matrices
generalizing Corollary~\ref{cor:2by2} preserve interlacing.

\begin{lemma}\label{lem:2by2-cor-pd}
  Suppose that ${\alpha},{\beta},\aaa,\bbb$ are
  vectors with all positive coordinates. The following matrices
  preserve interlacing in $\rupint{d}$:

$$
\begin{pmatrix}
  \alpha & \aaa\cdot\xx \\ 0 & \beta
\end{pmatrix}\quad
\begin{pmatrix}
  \aaa\cdot\xx & 0 \\ \alpha &\bbb\cdot\xx
\end{pmatrix}\quad
\begin{pmatrix}
  \aaa\cdot\xx & -\alpha \\0 & \bbb\cdot\xx
\end{pmatrix}\quad
\begin{pmatrix}
  \aaa\cdot\xx & -\alpha \\ \beta & \bbb \cdot \xx
\end{pmatrix}
$$
\end{lemma}
\begin{proof}
  All polynomials occurring in a product have positive homogeneous
  part.  Substitution and interlacing follow from the one variable
  case.
\end{proof}
Next, we generalize Proposition~\ref{prop:2by2}.

\begin{lemma} \label{lem:2by2-pd}
  Suppose that $M=\smalltwo{f}{g}{h}{k}$ is a matrix of polynomials
  with in $\rupint{d}$.  $M$ preserves interlacing if
    \begin{enumerate}
  \item $g\longrightarrow f$, $k\longrightarrow h$, and both
    interlacings are strict.
  \item The determinant $\smalltwodet{f}{g}{h}{k}$ is never zero.
  \end{enumerate}
\end{lemma}
\begin{proof}
  If $M\smalltwobyone{r}{s}=\smalltwobyone{u}{v}$  then the
  homogeneous parts of $u$ and $v$ are positive. Interlacing and
  substitution follow from the one variable case. 
\end{proof}

How can we find such matrices? Suppose that $f\in\gsubgen_2$. Recall
that this means that the graph of $f$ has no intersection
points. Choose $g\in\gsubgen_2$ so that $f\lessless g$. For instance,
we could take $g=\frac{\partial f}{\partial x}$. Define 
$$ M = 
\begin{pmatrix}
 f &   g \\
  \frac{\partial f}{\partial x} &   \frac{\partial g}{\partial x}
\end{pmatrix}
$$
It follows from Lemma~\ref{lem:inequality-1} that the determinant of
$M$ is never zero.

\section{Polynomials from matrices}
\label{sec:matrix-2-variable}

We can construct polynomials in several variables from symmetric
matrices, and with some simple assumptions the resulting polynomial is 
in $\gsubpos_d$. We begin with a  construction for a polynomial in
$d$ variables that generalizes Lemma~\ref{lem:det-aic}.

\begin{lemma} \label{lem:pairwise}
  If $C$ is symmetric and $A_2,\dots,A_d$ are pairwise commuting
  positive definite symmetric matrices, then
  $$ f(\xx) = \det(x_1I+ x_2A_2+\dots + x_dA_d+C)$$ is in $\gsubpos_d$.
\end{lemma}
\begin{proof}
  Since $A_2,\dots,A_d$ are pairwise commuting and symmetric, they are
  simultaneously diagonalizable.  There is an orthogonal matrix $P$ such that
  $PA_kP^{t}=D_k$, where $D_k = (e^k_{ij})$ is a diagonal matrix with
  all positive entries on the diagonal. We note

\begin{align*}
 f(\xx) &= \det(x_1\,I+ x_2\,A_2+\dots x_d\,A_d+C) \\
&= \det(P^{t})\, \det(x_1\, I + x_2\,D_2 + \dots + x_d\,D_d + PCP^{t})\,
det(P)
\intertext{As before this satisfies $x_1$-substitution since $PCP^t$ is 
  symmetric, and the homogeneous
  part is}
f^H(\xx) &= \det(x_1\, I + x_2\,D_2 + \dots + x_d\,D_d ) \\
&= \prod_k( x_1 + x_2\, e^2_{kk} + \dots + x_d \,e^d_{kk}).
\end{align*}
which has all positive coefficients.
\end{proof}

\begin{remark}
  A variant of the lemma is to assume that $$f(\xx)=det(x_1 D_1 +
  \cdots + x_d D_d + C)$$
  where $C$ is symmetric, all $D_i$ are
  diagonal, and $D_1$ has all positive diagonal entries. If
  $\sqrt{D_1^{-1}}$ is the diagonal matrix whose entries are the
  inverses of the positive square roots of the diagonal entries of
  $D_1$, then
$$ f(\xx) = det(D_1)\ det\left(x_1I + x_2 D_2D_1^{-1} + \cdots + x_d D_d
D_1^{-1} + \sqrt{D_1^{-1}} \,C\,\sqrt{D_1^{-1}}\right)
$$
We can now apply Lemma~\ref{lem:pairwise} to see that $f(\xx)\in\gsubpos_d$.
\end{remark}

How can we find such commuting families of positive definite matrices?
If $f$ is a polynomial, and $\lambda$ is an eigenvalue of a positive
definite matrix $A$, then $f(\lambda)$ is an eigenvalue of $f(A)$. In
particular, if $f\in\allpolypos$ then $f(A)$ is positive definite.
Since $e^x$ is a limit of polynomials in $\allpolypos$ and any two
polynomials in $A$ commute, we have

\begin{cor}
  Suppose that $f,f_2,\dots,f_d$ are polynomials in $\allpolypos,$  $A$ is
  positive definite and  $\alpha,\alpha_2,\dots,\alpha_d$ are
  positive. The following polynomials are in $\gsubpos_d$.
  \begin{align*}  
    & \det(x_1I+ x_2f_2(A)+\dots + x_df_d(A)+f(A)) \\
    & \det(x_1I+ x_2e^{\alpha_2A} +\dots + x_de^{\alpha_dA}+e^{\alpha
      A})
\end{align*}
\end{cor}

The polynomials that arise from these constructions can interlace. 

\begin{lemma} \label{lem:det-int}
  If $D$ is a diagonal matrix with all positive diagonal entries, $C$
  is symmetric, and $P$ is a principle submatrix of $xI+yD+C$, then
$$|xI+yD+C| \lesslesseq |P| \text{ in } \gsubpos_2$$ 
\index{principle submatrix}
\end{lemma}
\begin{proof}
  We may assume $P$ is the result of deleting the first row and first
  column of $|xI+yD+C|$. If $D=diag(d_i)$, $C=(c_{ij})$,
  $v=(c_{12},\dots,c_{1n})$ then for any $\alpha\in\reals$ we can write
  \begin{align*}
    g_\alpha(x,y) &= 
    \begin{vmatrix}
      x+d_1y+c_{11} +\alpha & v \\ v^t & P
    \end{vmatrix} \\[.2cm]
&= 
\begin{vmatrix}
      x+d_1y+c_{11}  & v \\ v^t & P
\end{vmatrix} +
\begin{vmatrix}
      \alpha & 0 \\ 0 & P
\end{vmatrix} \\[.2cm]
&= |xI+yD+C| + \alpha|P|
  \end{align*}
  Since $g_\alpha\in\gsubpos_2$ for any $\alpha\in\reals$ it follows
  that the desired polynomials interlace.
\end{proof}

\begin{remark}
  If we take $y=0$ then we recover the classical result that principle
  submatrices interlace. A similar result holds for $\rupint{d}$. See
  Theorem~\ref{thm:principle-1}.

If $C$ is the all zero matrix and $y=1$ then $f(x)=|xI+yD+C|=\prod(x+d_i)$ and
$|P|=\frac{f(x)}{x+d_1}$. These polynomials are the ones occurring in
quantitative sign interlacing. 
\end{remark}

  \section{Extending three interlacing polynomials}
  \label{sec:extend-three-interl}

Consider a polynomial in $\rupint{3}$

\centerline{
\xymatrix{
\ar@{.>}[d] & \ar@{.>}[d] \\
f_{01}(x) & f_{11}(x)\ar@{.>}[l] \ar@{.>}[d] & \ar@{..>}[l]\\
f_{00}(x) \ar@{<.}[u] \ar@{<.}[r] & f_{10}(x) & \ar@{.>}[l]
}
}

If we are only given $f_{00}\lesslesseq f_{01},f_{10}$ then we show
that we can add terms to $f_{00} + y f_{10} + zf_{01}$ so that the
resulting polynomial $F$ is in $\rupint{3}$. \seepage{lem:product} We say
that $F$ is an \emph{extension} of $f_{00},f_{10},f_{01}$. Given one
extension, we show how to find many more with the same initial data.

The term $f_{11}$ is not arbitrary, since it interlaces $f_{01}$ and
$f_{10}$, and satisfies $f_{01}f_{10}-f_{00}f_{11}\ge0$. We show that
if $f_{00}$ has degree two then $f_{10}f_{01}-f_{00}W\ge 0$ if and only
if there is an extension of $f_{00},f_{10},f_{01}$ such that
$W=f_{11}$. We conjecture that this is true in general
(Question~\ref{ques:p3-ext}).

  \begin{lemma}\label{lem:p3-ext-1}
    Suppose that $f\lesslesseq g$ and $f\lesslesseq h$, where $f,g,h$
    have positive leading coefficients. Then there are $F\in\rupint{3}$
    such that if we write $F(x,y,z)=\sum f_{ij}(x)\,y^iz^j$ then
\[
f_{00} = f,\quad f_{10} = g,\quad f_{01}=h
\]
  \end{lemma}
  \begin{proof}
    If we take three sequences $a_i,b_i,c_i$, a constant
    $\alpha$, and form the product
    \begin{align*}
      F &= \alpha\prod(1+a_ix+b_iy+c_iz) \\
\intertext{then}
F_{00} &= \alpha \prod(1+a_i x) \\
F_{10} &= \sum b_i \frac{F_{00}}{1+a_ix} \\
F_{01} &= \sum c_i \frac{F_{00}}{1+a_ix} 
    \end{align*}
    Given $f,g,h$ we can, following the argument of
    Lemma~\ref{lem:product}, find $\alpha$ and the sequences
    $a_i,bi,c_i$ so that $F_{00}=f$, $F_{10}=g$, $F_{01}=h$.
  \end{proof}

\begin{remark}
    We can use a determinant identity to derive a different extension of three
    polynomials. Recall \eqref{eqn:det-psd-2}. 
%
If we  let $D = diag(d_i)$, $f = |I + xD|$, $V=(V_i)$ and $W = (W_i)$ then
\begin{gather*}
  |I + x D + y VV^t + z W W^t| =
 \\
f +  y\sum V_i^2 \frac{f}{1+d_ix} + z\sum W_i^2 \frac{f}{1+d_ix} + 
yz\sum_{i<j}
\frac{f}{{(1+xd_i)(1+xd_j)}}
\begin{vmatrix}
  V_i&V_j\\W_i&W_j
\end{vmatrix}^2
\end{gather*}
Since $D,VV^t,WW^t$ are positive semi-definite, the right hand side is
in $\gsubclose_3$. (See Lemma~\ref{lem:rank}.) It follows that we have
an interlacing square

  \centerline{\xymatrix{
      {\displaystyle\sum V_i^2 \frac{f}{1+d_ix}}
      \ar@{->}[d]_{ }           
      \ar@{<-}[rrr]^{ }         
      &&&
{\displaystyle      \sum_{i<j}
\frac{f}{(1+xd_i)(1+xd_j)}
\begin{vmatrix}
  V_i&V_j\\W_i&W_j
\end{vmatrix}^2
}
      \ar@{->}[d]^{ } \\        
      f
      \ar@{<-}[rrr]^{ }         
      &&&
      {\displaystyle\sum W_i^2 \frac{f}{1+d_ix} }
}}

\end{remark}

There are many determinants that will give the same initial
data. This follows from the more general result that shows that the
initial data only constrains the diagonals.

\begin{lemma}\label{lem:p3-ext-2}
  Suppose that $D$ is a diagonal matrix, $D_1$ and $E_1$ have
  the same diagonal, as do $D_2$ and $E_2$.  If we write
  \begin{align*}
    \bigl| I + x D + y D_1 + z D_2\bigr| &= \sum U_{ij}(x)y^iz^j \\
    \bigl| I + x D + y E_1 + z E_2\bigr| &= \sum V_{ij}(x)y^iz^j \\
\text{then}\qquad
U_{00} = V_{00},\quad U_{10} = V_{10},\quad U_{01} = V_{01}
  \end{align*}
\end{lemma}
\begin{proof}
  Clearly $U_{00} = V_{00}$. Let $D_1=(d_{ij})$ and
  $E_1=(e_{ij})$. Then
  \begin{align*}
    U_{10} &= \frac{\partial}{\partial y}\bigl|I + x D + y D_1 + 0
    D_2\bigr|_{y=0} \\
&= \sum_{i,j} d_{ij} \bigl|(I+xD)[i,j]\bigr| \\
&= \sum_{i} d_{ii} \bigl|(I+xD)[i,i]\bigr|\qquad\text{since $I+xD$ is diagonal}  \\
&= \sum_{i} e_{ii} \bigl|(I+xD)[i,i]\bigr|\qquad\text{since $e_{ii}=d_{ii}$} \\
&= V_{10}
  \end{align*}
\end{proof}

If we are given $f_{00}\longleftarrow f_{10},f_{01}$ in $\allpolypos$, then by the
first lemma we can realize $f_{00},f_{10},f_{01}$ as the coefficients of
$|I+xD+yD_1+zD_2|$ where $D,D_1,D_2$ are positive definite diagonal
matrices. If we can find symmetric matrices $S_1,S_2$ with zero
diagonal such that $D_1+S_1$ and $D_2+S_2$ are positive definite, then
$f_{00},f_{10},f_{01}$ are also realized as initial coefficients of
$|I+xD+y(D_1+S_1)+z(D_2+S_2)|$.

\begin{example}
  Consider the case $f_{00}\in\allpolypos(2)$. Given the initial data
  $f_{00}\lesslesseq f_{10},f_{01}$ we first determine all $\alpha$
  that arise from Lemma~\ref{lem:p3-ext-2}. Consider
\[
g(B,C) = \left|
  \begin{pmatrix}1&0\\0&1  \end{pmatrix} + 
x   \begin{pmatrix}a_1&0\\0&a_2  \end{pmatrix} +
y   \begin{pmatrix}b_1&B\\B&b_2  \end{pmatrix} +
z  \begin{pmatrix}c_1&C\\C&c_1  \end{pmatrix}
\right|
\]
where we choose $a_1,a_2,b_1,b_2,c_1,c_2$ so that
\[
g(0,0) = f_{00} + y f_{10} + zf_{01} + \alpha yz.
\]
This representation is possible by Lemma~\ref{lem:p3-ext-1}. Recall
\begin{quote}
  If $r,s$ are positive then $\smalltwodet{r}{t}{t}{s}$ is positive
  definite if and only if $rs>t^2$.
\end{quote}

Thus, $g(B,C)\in\rupint{3}(2)$ if and only if $b_1b_2>B^2$ and
$c_1c_2>C^2$. Now $\alpha$ is the coefficient of $yz$ which is 
\[
b_1c_2+b_2c_1+2BC.\]
Using the constraints above, we find
\begin{equation}
  \label{eqn:p3-ext}
  \bigl(\sqrt{b_1c_2}-\sqrt{b_2c_1}\bigr)^2 \le \alpha \le
\bigl(\sqrt{b_1c_2}+\sqrt{b_2c_1}\bigr)^2
\end{equation}
All $f_{00},f_{01},f_{10}$ and $\alpha$ satisfying this inequality
have extensions to $\rupint{3}$. 

Conversely, assume that $f_{00},f_{10},f_{01}$ are determined by
$g(0,0)$. If we solve the inequality 
\[
f_{10}f_{01}-\alpha f_{00}\ge0
\]
for $\alpha$ we get \eqref{eqn:p3-ext}.

\end{example}

\section{Subsets of $\rupint{d}$}
\label{sec:sub-pd}

In one variable we had two different ways to talk about $\allpolypos$.
We could describe it in terms of coefficients (all the same sign), or
in terms of roots (all negative). However, $\allpolyint{(-1,1)}$
has no definition in terms of coefficients, but is given only in terms
of the location of its roots. In higher dimensions the graph
corresponds to zeros in one dimension.  We define

\begin{definition}
  If $\diffk\subset\reals^d$ then
  \begin{align*}
    \gsubint{\diffk}_{d} &= \left\{ f\in\rupint{d}\mid \text{ the graph of
        $f$ is disjoint from $\diffk$}\right\} \\
&=   \left\{ f\in\rupint{d}\mid f(\xx)=0 \implies \xx\not\in\diffk\right\}
  \end{align*}
\end{definition}

Consider some examples and properties.

\begin{enumerate}
\item If $\diffk=\emptyset$ then $\gsubint{\diffk}_d = \rupint{d}$, since
  there are no restrictions.
\item If $\diffk\subset\diffk'$ then
  $\gsubint{\diffk'}_d\subset\gsubint{\diffk}_d$. 
\item Let $v,w$ be two points of $\reals_d$, and let 
  $$
  \diffk = \{x\in\reals^d\mid x\ge w\} \cup \{x\in\reals^d\mid x\le
  v\}$$
  Every solution variety of $\gsubint{\diffk}_d$ must meet the
  segment $\overline{vw}$. This is an analog of $\allpolyint{(-1,1)}$.
\item Let $\diffk = \{(x,x,\dots,x)\mid x\in\reals\}$. Since every
  $f\in\rupint{d}$ satisfies $f(x\dots,x)\in\allpoly$, it follows
  $\gsubint{\diffk}_d=\emptyset$. 
\item Suppose $\diffk=\{\xx\in\reals^2\mid \xx\ge 0 \text{ or }
  \xx\le0\}\setminus (0,0)$. Then, $\gsubint{\diffk}_2$ consists of
  polynomials whose graphs are straight lines passing through the
  origin. 
\item Suppose that $\diffk = \{x\in\reals^n\mid \xx\ge0\}$. We will
  show that $\gsubint{\diffk}_d=\gsubplus_d$. Since all coefficients of
  $\gsubplus_d$ are positive, it follows that $f\in\gsubplus_d$ has no
  root in $\diffk$. We prove the converse by induction on $d$. For
  $d=1$ we know that all the coefficients are the same sign if and
  only if no root is non-negative. Suppose that $f\in\rupint{d}$ has no
  root $\xx$ with $\xx\ge0$. Write
  $f(\xx)=\sum f_i(x_1,\dots,x_{d-1})x_d^i$. If we substitute $x_d=0$
  then we see that $f_0(x_1,\dots,x_{d-1})$ has no root with all
  coordinates non-negative, so by induction all coefficients of
  $f_0(x_1,\dots,x_{d-1})$ have the same sign. Since all coefficients of
  $f^H$ are positive, all coefficients of $f_0(x_1,\dots,x_{d0-1})$
  are positive.
  
  Next, if $f\in\gsubint{\diffk}_d$, then $\frac{\partial f}{\partial
    x_d}\in\gsubint{\diffk}_d$ since the roots of $\frac{\partial
    f}{\partial x_d}$ lie between roots of $f$. Consequently, we
  conclude as above that $f_1(x_1,\dots,x_{d-1})$ has all positive
    coefficients. Continuing, we see $f$ has all positive
    coefficients. 
\end{enumerate}

Many different $\diffk$ can determine the same set of polynomials. For
example, the set $\diffk = \{\xx\in\reals^d\mid \xx\ge0\}$ and $\diffk'
= \diffk\setminus p$ where $p>0$ satisfy
$\gsubint{\diffk}_d=\gsubint{\diffk'}_d$. However, if we choose $p$ to
be the origin, then $\gsubint{\diffk}_d \subsetneq
\gsubint{\diffk'}_d$ since $x_1+\cdots+x_d$ is in 
$\gsubint{\diffk'}_d$ but is not in $\gsubint{\diffk}_d$. We define

\begin{definition}
$$  \overline{\diffk} = \bigcup\left\{\diffl\subset \reals^d\mid \gsubint{K}_d
  = \gsubint{\diffl}_d\right\}
$$
\end{definition}

In some cases we can determine $\overline{\diffk}$. 

\index{convex! set in $\reals^d$}
\begin{lemma}
  Suppose that $\diffk$ is a closed convex region of $\reals^d$ whose
  supporting hyperplanes have an equation in $\rupint{d}$. Then,
  $\diffk=\overline{\diffk}$. 
\end{lemma}
\begin{proof}
  Since $K\subset\overline{\diffk}$ it suffices that any point $v$ not
  in $\diffk$ is not in $\overline{\diffk}$. Pick a supporting
  hyperplane $H$ such that $\diffk$ is on one side of it, and $v$ is
  in the interior of the other side.  Choose a hyperplane $H'$
  parallel to $H$, and containing $v$. It follows that $v$ lies on the
  graph of a polynomial that is disjoint from $\diffk$, and hence
  $v\not\in\overline{\diffk}.$
\end{proof}

\begin{remark}
  In case $d=2$ we see that if $\diffk$ is a closed convex region
  whose boundary is a non-increasing curve then $\diffk =
  \overline{\diffk}$. 
\end{remark}

\index{convex!in $i$-th coordinate}
We now ask when $\gsubint{\diffk}_d$ is closed under
differentiation. We say a set $S$ is convex in the $i$-th coordinate
if whenever $v,w$ agree in all coordinates except the $i$-th, then the
segment $\overline{vw}$ is in $S$. 

\begin{lemma}
  If $\reals^d\setminus\diffk$ is convex in the $i$-th coordinate then
$\gsubint{\diffk}_d$ is closed under $f\mapsto \frac{\partial
  f}{\partial x_i}$. 
\end{lemma}
\begin{proof}
To show that  $\frac{\partial f}{\partial x_i}\in\gsubint{\diffk}_d$
we fix all but the $i$-th coordinate and consider the roots of the
resulting polynomial. Since the roots of $\frac{\partial f}{\partial x_i}$
  lie between the roots of $f$ the convexity in the $i$-coordinate
  imply that they lie in $\reals^d\setminus\diffk$.
\end{proof}

\section{Higher dimensional recurrences}
\label{sec:n-dim-recur}
\index{recurrence!higher dimensional}

In this section we investigate polynomials defined by two and three
dimensional recurrences. Some of these polynomials are in one
variable, and some in two. One of the interesting features is that we
do not always have that consecutive terms interlace. We begin by
considering polynomials $f_{i,j}$ satisfying the similar recurrences

\begin{align}
f_{i,j} &  =   (x-a_i)f_{i-1,j} + b_{ij}f_{i,j-1}   \label{eqn:nrec-1}\\
f_{i,j} &  =   (x-a_{ij})f_{i-1,j} + b_{ij} f_{i,j-1} \label{eqn:nrec-2}\\
f_{i,j} &  =   (x-a_i)f_{i-1,j} + b_{i}f_{i-1,j-1} \label{eqn:nrec-3}
\end{align}

\noindent%
We can visualize these recursions using  lattice diagrams:\\[.2cm]

\centerline{
\xymatrix{
f_{i-1,j} \ar@{-}[r]^{x-a_i} &    \ar@{-}[d]^{b_{ij}}    f_{i,j} & 
f_{i-1,j} \ar@{-}[r]^{x-a_{ij}} & \ar@{-}[d]^{b_{ij}} f_{i,j} & 
f_{i-1,j} \ar@{-}[r]^{x-a_i} &    \ar@{-}[dl]^{b_i}   f_{i,j} \\
& f_{i,j-1} &  & f_{i,j-1} & f_{i-1,j-1} & & 
}}

We will see that recurrence \eqref{eqn:nrec-1} is a disguised product,
that recurrence \eqref{eqn:nrec-2} fails unless the coefficients
satisfy certain restrictions, and
recurrence \eqref{eqn:nrec-3} comes from coefficients of products in
$\gsubpos_2$.

\begin{lemma} \label{lem:nrec-1}
  Suppose $f_{0,0}=1$ and $f$ satisfies recurrence
  \eqref{eqn:nrec-1}. Then $f_{i,j}$ is a constant multiple of\ 
  $\prod_{k=1}^i(x-a_k)$. In particular, all $f_{i,j}$ are in
  $\allpoly$, and consecutive $f_{i,j}$ interlace.
\end{lemma}
\begin{proof}
  Assume by induction of $i,j$ that 
$$f_{r,s} = \alpha_{r,s} \prod_{k=1}^r (x-a_k)$$
for all $(r,s) < (i,j)$. We then have

\begin{align*}
  f_{i,j} &= (x-a_i)f_{i-1,j} + b_{i,j}f_{i,j-1} \\ 
  &= (x-a_i) \alpha_{i-1,j}\prod_1^{i-1} (x-a_k) +
  b_{i,j}\alpha_{i,j-1}\prod_1^i(x-a_k)\\
&= (\alpha_{i-1,j} + b_{i,j}\alpha_{i,j-1})\prod_1^i(x-a_k)
\end{align*}

\end{proof}

Here is a simple recurrence whose proof that all members are in
$\allpoly$ only depends on the additivity of interlacing.

\begin{lemma} \label{lem:rec-7}
Let $h_n = \prod_1^n(x+a_i)$, and let $h_n\longleftarrow g_n$ where all
$g_n$ have positive leading coefficients. If we define $f_0=1$ and
$$ f_n = (x+a_n)f_{n-1} + g_n$$
then $h_n \greateqeq f_n$.In particular, all $f_i$ have all real roots.
\end{lemma}
\begin{proof}
  We assume by induction that $h_{n-1}\greateqeq f_{n-1}$ Multiplying
  each side by by $(x-a_n)$ yields $h_n \greateqeq
  (x-a_n)f_{n-1}$. Since we also have $h_n\longleftarrow g_n$ we can
  add these interlacings to conclude that
$$ h_n \greateqeq (x-a_n)f_{n-1}+g_n = f_n$$
\end{proof}

We can apply this lemma to construct a particular sequence of
polynomials whose members all have all real roots, but consecutive
ones do not generally interlace. 

\begin{cor}
  If $f_0=1$ and $f_n$ satisfies the recursion
$$ f_n = (x+n)f_{n-1} + (x+n-2)\prod_1^{n-1}(x+i)$$ then all $f_n$ are
in $\allpoly$.
\end{cor}
\begin{proof}
  Take $h_n = \prod_1^n(x+i)$, and $g_n=(x+n-2)h_{n-1}$.
\end{proof}

\begin{lemma} \label{lem:nrec-3}
Suppose that $f_{0,0}=1$ and $f_{i,j}$ satisfies
  \eqref{eqn:nrec-3} where all $a_i$ have the same sign. Then all
  $f_{i,j}$ are in $\allpoly$.
\end{lemma}
\begin{proof}
  We claim that the $f_{i,j}$ are the coefficients of the products
\begin{align*} \prod_{j=1}^i(x+b_jy  -a_j) &= \sum_{j=0}^i f_{i,j}(x)\,y^j\\
\intertext{This is certainly true when $i=0$. We only need to check that if we
write}
 \prod_{j=1}^i(x+b_jy  -a_j) &= \sum_{j=0}^i
  g_{i,j}(x)\,y^j \\
\intertext{then the $g_{i,j}$ satisfy the same recurrences as do the $f_{i,j}$:}
  \prod_1^{i}(x+b_jy-a_j) &= (x+b_{i}y-a_{i}) \sum_{j=0}^{i-1}
  g_{i-1,j}(x)\,y^j \\
  & = \sum ( (x-a_i)g_{i-1,j}(x) -b_{i}g_{i-1,j-1}(x))y^j
\end{align*}
and so we are done by induction.
\end{proof}

We can generalize this argument to products in more variables - this
leads to polynomials in $x$ with more indices. For example, consider
the recurrence

\begin{align}
  \label{eqn:nrec-5}
  f_{i,j,k}(x) = (x-a_i)f_{i-1,j,k} + b_if_{i,j-1,k} + c_if_{i,j,k-1}
\end{align}

\begin{lemma} \label{lem:nrec-5}
  If $f_{0,0,0}=1$ and $f_{i,j,k}$ satisfies \eqref{eqn:nrec-5}, all
  $b_i$ have the same sign, and all $c_i$ have the same sign then all
  $f_{i,j,k}$ are in $\allpoly$.
\end{lemma}
\begin{proof}
  Consider the product below, and proceed as in Lemma~\ref{lem:nrec-3}.
$$ \prod_{h=1}^i (x-a_h+b_hy+c_hz) = \sum f_{i,j,k}(x)\,y^jz^k$$
\end{proof}

As another product example, consider the recurrence

\begin{align}
  \label{eqn:nrec-6}
  f_{i,j} = (a_ix+b_i)f_{i,j-1} + (c_i x+ d_i)f_{i-1,j}
\end{align}

whose diagram is

\centerline{
\xymatrix{
f_{i-1,j} \ar@{-}[r]^{c_ix+d_i} &    \ar@{-}[d]^{a_ix+b_1}    f_{i,j} \\ 
 & f_{i,j-1}
}}

\begin{lemma} \label{lem:recur-6}
  If $f_{0,0}=1$, $f_i,j$ satisfies \eqref{eqn:nrec-6}, and $a_id_i -
  c_ib_i>0$ for all  $i$  then $f_{i,j}\in\allpoly$.
\end{lemma}
\begin{proof}
  The assumptions imply that all the factors in the product below are
  in $\gsubposclose_2$,  so we can proceed as above
  
$$ \prod_{j=1}^i((a_jx+b_j)+ (c_j+d_jx)y) = \sum_j f_{i,j}(x)y^j$$ 
\end{proof}

In fact, we can take any polynomial in $\gsubposclose_2$ that has some
free parameters, take their product, and the coefficients $f_{i,j}$ of
$y^j$ will satisfy a recurrence, and all the $f_{i,j}$ belong to
$\allpoly$, and consequently will interlace.

Here is another disguised product 

\begin{align}
  f_{i,j} &= (x-b_i)f_{i-1,j} + (y-d_j)f_{i,j-1} \label{eqn:nrec-4}
\end{align}

with diagram 

\centerline{
\xymatrix{
f_{i-1,j} \ar@{-}[r]^{x-b_i} &    \ar@{-}[d]^{y-d_j}    f_{i,j} \\ 
 & f_{i,j-1}
}}

  The proof is a simple modification of the proof of Lemma~\ref{lem:nrec-1}.

\begin{lemma}
  Suppose that $f_{0,0}=1$ and $f_{i,j}$ satisfies
  \eqref{eqn:nrec-4}. Then $f_{i,j}$ is a constant multiple of 
  $\prod_{k=1}^i(x-b_k) \ \cdot \ \prod_{h=1}^j (y-d_h)$.
\end{lemma}

We can get many recurrences from products. Here is a general
approach. Let 
\begin{align*}
  g(x,y,z) &= \sum g_r(x,z)y^r \\
\intertext{where $g\in\gsubpos_3$, choose $t_j\in\reals$, and consider the product}
\prod_{j=1}^i g(x,y,t_j) &= \sum_{i,j} f_{i,j}(x)y^j \\
\intertext{We know that all $f_{i,j}\in\allpoly$, and consecutive ones
  interlace. The recurrence is easily found:}
\prod_{j=1}^{i+1} g(x,y,t_j) &= g(x,y,t_{i+1})\ \sum_{i,s} f_{i,s}(x)y^s
\\
&= \left(\sum g_r(x,t_{i+1}) y^r\right) \ \sum_{i,s} f_{i,s}(x)y^s \\
\intertext{and comparing coefficients of $y^j$ yields}
f_{i+1,j} &= \sum_{r+s=j}g_r(x,t_{i+1})f_{i,s}(x)
\end{align*}


\chapter{The polynomial closure of $\gsubpos_{d}$}
\label{cha:topology}

\renewcommand{\TimeStampStart}{Thursday, January 17, 2008: 19:40:23}
\mytoday    

In this chapter we study the polynomial closure of $\gsubpos_d$. For
polynomials in one variable we know that $\allpoly(n)$ is exactly the
set of polynomials in the closure of $\allpoly(n)$, so the closure of
$\allpoly$ leads to analytic functions. However, there are polynomials
in the closure of $\gsubpos_d(n)$ that are not in $\gsubpos_d(n)$, so we
have two closures to consider: new polynomials and analytic functions.
As expected, most properties of $\gsubpos_d$ extend to these closures. We
consider the analytic closure in the next chapter.

\section{The polynomial closure of $\rupint{2}$}
\label{sec:topology-p2}
In this section we investigate properties of the polynomial closure of
$\rupint{2}$.  We begin with the general definitions.

\begin{definition}
  $\gsubclose_d(n)$ is the set of all polynomials that are limits of
  polynomials in $\rupint{d}(n)$.  $\gsubplusclose_d(n)$ is defined
  similarly.  We set $$
  \gsubclose_d = \gsubclose_d(1) \cup
  \gsubclose_d(2) \cup \dots.$$

We say that $f,g\in\gsubclose_d$ \emph{interlace} if $f+\alpha
g\in\gsubclose_d$  for all real $\alpha$.
\index{interlacing!in $\gsubclose_d$}

\end{definition}

Although most properties of $\rupint{2}$ extend to $\gsubclose_2$, there
are some slight differences.  The polynomial $x^ny^m$ is in
$\gsubclose_2$ since 
$$\lim_{\epsilon\rightarrow0^+} (x+\epsilon y)^n(y+\epsilon x)^m =
x^ny^m$$
and $(x+\epsilon y)^n(y+\epsilon x)^m$ is in $\rupint{2}$.
This shows that there are polynomials in $\gsubclose_2$ where the
$x$-degree is $n,$ the $y$-degree is $m,$ and the total degree is
$n+m$.  Consider this example that shows how the coefficients can
interlace in complicated ways.

\begin{example}
  The definition of $f(x,y)$ as a limit of polynomials in $\rupint{2}$
  \begin{align*}
    f(x,y) &= \lim_{\epsilon\rightarrow0^+}
(\epsilon x+y-1)(\epsilon x+y-2)(\epsilon x+y+3)(x+y+1)(x+2y+1) \\
& = (y-1)(y-2)(y+3)(x+y+1)(x+2y+1) \\
\intertext{shows that $f\in\gsubclose_2$. If we expand $f$ we get}
&= (6 + 12x + 6x^2) + ( 11 + 4x - 7x^2)y + ( -9 - 21x)y^2 +\\
&\quad\quad 
    (-13 + 2x + x^2)y^3 + (3 + 3x)y^4 + 2y^5 \\
&= f_0 + f_1\,y + f_2\,y^2 + f_3\,y^3 + f_4\,y^4 + f_5\,y^5 
  \end{align*}
  These coefficients interlace in ways that are not possible for
  polynomials in $\rupint{2},$ namely
$$ f_0 \lesseqeq f_1 \lessless f_2 \lessgreat f_3 \lessless f_4.$$

We can use the idea of this construction to show that if
$g(x)\in\allpoly$ then $g(x)\in\gsubclose_2$. Simply note that
$g(x+\epsilon y+\epsilon)$ is in $\rupint{2}$ and converges to $g(x)$ as
positive $\epsilon$ goes to $0$. Similarly, if we let $\epsilon$ go to
zero through positive values we see that 
$$ \allpolypos \subset \gsubplusclose_2$$
The mapping $f\times g\mapsto f(x)g(y)$ defines a map
$\allpoly\times\allpoly\longrightarrow\gsubclose_2$, and \\
$\allpolypos\times\allpolypos\longrightarrow\gsubplusclose_2$.

The two polynomials $f=x^ny^m$ and $g= x^{n-1}y^{m-1}$ do not
interlace. However, they do satisfy $f\plesslesseq g$ in the sense
that $f + \alpha g\in\gsub$ for all \emph{positive} $\alpha$. 
\index{positive interlacing}

In order to see that polynomials in $\gsubclose_2$ interlace we can
express them as limits of interlacings. For instance, the interlacing
$xy-1 \lesslesseq x+y$ follows from
\[
(xy-1)+\alpha(x+y) = 
\lim_{\epsilon\rightarrow0^+}\,\begin{vmatrix}
  \epsilon\,x+y+\alpha & \sqrt{1+\alpha^2} \\
\sqrt{1+\alpha^2} & x+\epsilon\,y+\alpha
\end{vmatrix}
\]
\end{example}

\begin{example}\label{rem:neg-h}
  Unlike in $\rupint{2}$, the signs of a polynomial in $\gsubclose_2$ can
  be all negative. This follows from the observation that if
  $f\in\gsubclose_2$ then $-f\in\gsubclose_2$. This is a consequence
  of the limit 

$$ \lim_{\epsilon\rightarrow0^+} \ (\epsilon x + \epsilon y -1)f = -f$$

\end{example}

If $f\in\gsubclose_2(n)$ then we define the homogeneous part $f^H$ to
be all terms of total degree $n$.  \index{homogeneous part} Unlike
$\rupint{2}$, the example above shows that all the coefficients of the
homogeneous part can be negative. There are no other possibilities -
see Lemma~\ref{lem:homog-in-p2close}.

We can apply Lemma~\ref{lem:p2-xx} to deduce

\begin{lemma} \label{lem:sub-fxx}
  If $f\in\gsubclosepos_2(n)$ then $f(x,x)$ and $x^n f(x,-1/x)$ are in
  $\allpoly$.
\end{lemma}

\begin{figure}[htbp]
  \begin{center}
    \leavevmode
    \epsfig{file=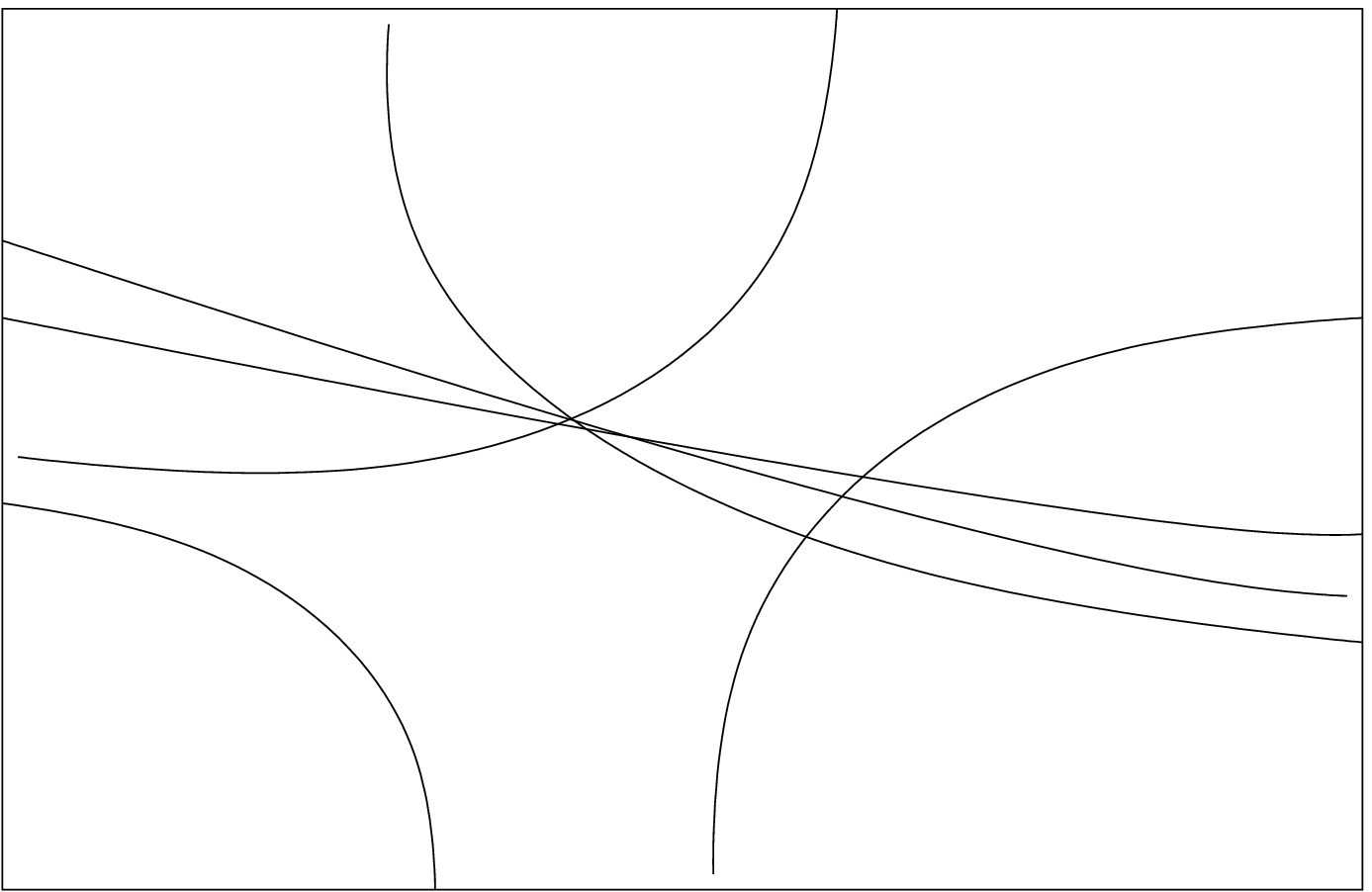,width=3in}
    \caption{The graph of of a polynomial in $\gsubclose_2$}
    \label{fig:p2close}
  \end{center}
\end{figure}

This shows that $xy+1$ isn't in $\gsubclose_2$ even though it
satisfies substitution.  The graphs of polynomials in $\gsubclose_2$
are more complicated than those in $\rupint{2}$. Consider
$$ f = (x+4y+ x y + 1)(x + y + 3x y +2)(x+5y+.5)(x+2y+1). $$ The graph 
of $f$ in Figure~\ref{fig:p2close} has two vertical asymptotes at $x=-4$ and
$x=-3$. None the less, every vertical line except these two meets the
graph in exactly four points.

\index{determinants!and realizing $\rupint{2}$}
We can use determinants to  realize members of $\gsubclosepos_2$. 
\begin{lemma} \label{lem:top-det-2}
  Suppose that $C$ is a symmetric  matrix, and $D$ is a
  diagonal matrix where the $i$-th diagonal entry is $ d_{i1}x_1 + d_{i2}x_2$
  where $d_{ij}$ is non-negative. The polynomial $det(C+D)$ is in
  $\gsubclosepos_2$. 
\end{lemma}
\begin{proof}
  Choose $\epsilon>0$ and let $D_\epsilon$ have diagonal entries
  $\sum (d_{ij}+\epsilon)x_i$.  A argument similar to Lemma~\ref{lem:pairwise}
  shows that $det(C+D_\epsilon)\in\rupint{2}$.  Taking limits gives the
  result.
\end{proof}

We now find some simple polynomials in $\gsubclose_2$ that are not in $\rupint{2}$.

\begin{lemma} \label{lem:fxy}
  If $f(x)\in\allpolypos$ then $f(-xy)\in\pm\gsubclosepos_2$. 
\end{lemma}

\begin{proof}
Since $f$ is a product of linear factors, it suffices to show that
$-xy+a^2\in\gsubclose_2$. The determinant
$$  
\begin{vmatrix}
  x+\epsilon y & a \\ a & y+\epsilon x
\end{vmatrix}
$$
is in $\rupint{2}$, and converges to $xy-a^2$. 
\end{proof}

\begin{cor} \label{cor:sub-xyc} Suppose $a>0$.
 $axy+bx+cy+d$ is in $\pm\gsubclosepos_2$ iff
  $\smalltwodet{a}{b}{c}{d}\le0$.

\end{cor}

\begin{proof}
  Let $f(x,y) =
  axy+bx+cy+d$.  Since $(1/a)f(x-c/a,y-b/a) = xy - (bc-ad)/a^2$ we can
  scale and apply Lemma~\ref{lem:fxy}.
\end{proof}

The  lemma can be derived from Theorem~\ref{thm:nsd}. Also see
Question~\ref{ques:det-ad}. It is also a special case of
Lemma~\ref{lem:fyg}. The necessity is the rhombus inequality
(Proposition \ref{prop:p2plus-inequality}). 

\index{rhombus inequalities}

\begin{lemma}
  If $xy+bx+cy+d\in\gsubclose_2$ then there are
  $\alpha,\beta,\gamma$ such that
$$ f = \begin{vmatrix} \alpha+x & \beta \\ \beta & \gamma+y
\end{vmatrix}$$
\end{lemma}

\begin{proof}
  If we choose $\gamma=b$ and $\alpha=c$ then we find that $\beta^2 =
  bc-d$. This expression is positive by the previous lemma.
\end{proof}

\begin{example}
  Since $\smalltwodet{x+1}{3}{3}{y-2}= -(2x+11) +
  y(x+1)\in\gsubclose_2$ we see that a multiaffine polynomial $f+yg
  \in\gsubclose_2$ can have $f,g$ with different leading signs. Of
  course we still have $f\greateqeq g$.
\end{example}

\begin{example} \label{ex:not-xxy}
  There are simple polynomials that are not in $\gsubclose_2$. For
  instance, $x^2\pm y$ is not in $\gsubclose_2$. Indeed,
  $f(x,y)=x^2-y$ doesn't satisfy substitution because
  $f(x,-1)\not\in\allpoly$. Similarly $x^2+y$ doesn't satisfy
  substitution.
\end{example}

We can extend Lemma~\ref{lem:fxy} by replacing $x$ by $xy+bx+cy+d$ instead of 
by $xy$. The next result expresses this extension, but uses a more
convenient replacement.

\begin{cor} \label{cor:sub-xy}
  If $f\in\allpoly$ then the polynomial
  $f((x+\alpha)(y+\beta)+\gamma)$ is in $\gsubclose_2$ if all roots of
  $f$ are either greater than $\gamma$, or all roots of $f$ are less
  than $\gamma$.
\end{cor}


If $f\in\gsubplus_2$ then $f(x,0)$ and $f(0,x)$ have all positive
coefficients, and all negative roots. Thus the graph of $f$ meets the
$x$ and $y$ axis in negative values.  If we intersect the graph with
the line $x=y$ then the intersection points are all negative.
With this in mind, we establish

\begin{lemma} \label{lem:p2plusclose}
  If $f\in\gsubcloseplus_2$ then $f(-x^2,-x^2)\in\allpoly.$
\end{lemma}
\begin{proof}
  By continuity of roots it suffices to assume that
  $f\in\gsubplus_2$.  The roots of $f(-x,-x)$ are all positive,
  so we can take their square roots.
\end{proof}

This is a special case of Lemma~\ref{lem:sub-quad-form}.

\begin{remark}
    If $f(x,y)$ satisfies $x$ substitution, and $f^H$ is nearly
    homogeneous, $f$ is not necessarily in $\rupint{2}$. For instance, if
    $f(x,y) = 1+x(x+y)$ then $f$ satisfies $x$ substitution, and the
    homogeneous part has the form \[x\times(\text{homogeneous polynomial with
      positive coefficients}).\] However, $f\not\in\gsubclose_2$ since
    $f(x,x)\not\in\allpoly$.
  \end{remark}


Here's a simple fact whose proof requires several closure properties
of $\gsubclose_2$.

\begin{lemma}
  If $ax+by+c\in\gsubclose_2$ then $a$ and $b$ have the same sign.
\end{lemma}
\begin{proof}
  After scaling and translating, it suffices to show that that $x-y$
  is not in $\gsubclose_2$. If it were the case that
  $x-y\in\gsubclose_2$ then we also know that $x-y+3\in\gsubclose_2$.
  Since $\gsubclose_2$ is closed under multiplication and
  differentiation, it follows that
$$ f(x,y) = \frac{\partial}{\partial x} (x-y)(x-y+3)(x+y) =6\,x + 3\,x^2 -
2\,x\,y - y^2 \in\gsubclose_2$$
However, $f(-1,y)=-3 + 2\,y - y^2 $ has two complex roots.
  
\end{proof}


The homogeneous part of a polynomial in $\gsubclose_2$ is not too
different from the homogeneous part of a polynomial in $\rupint{2}$.

  \begin{lemma}\label{lem:homog-in-p2close}
    If $f\in\gsubclose_2$  then all non-zero
    coefficients of $f^H$ have the same sign. There are non-negative
    integers $n,m$ and $g\in\allpolypos$ so that $f^H =\pm x^ny^m G$
    where $G$ is the homogenization of $g$.
  \end{lemma}
  \begin{proof}
    
    If $f^H$ only has one term then the result is true. If not all
    coefficients are the same sign then we can write $f^H = \cdots + r
    x^ay^{n-a} + s x^by^{n-b} + \cdots$ where $a<b$, $rs<0$, and all
    coefficients of $x^ry^{n-r}$ are zero for $a<r<b$. Let $g =
    (\partial_x)^a (\partial_y)^{n-b}f$ so that $g^H(x,y) = c y^{b-a} +
    d x^{b-a}$ where $cd<0$.
    
    Substituting $y=1$ shows that $b-a$ is either $1$ or $2$. If
    $b-a=1$ then $g = cy + dx+e$, but this is not possible by the
    previous lemma.  If $b-a=2$ then we can write $g = \alpha x^2 + \beta xy +
    \gamma y^2 + \cdots$ where $\alpha$ and $\gamma$ have opposite
    signs. After scaling $x$ and $y$ and perhaps multiplying by $-1$
    we can assume that $g = x^2 + \beta xy - y^2 + \cdots$. If
    $\beta>0$ then $\partial_y g = \beta x - 2y + (constant)$ which is
    impossible by the lemma. If $\beta<0$ then $\partial_x g
    = 2x + \beta y + (constant)$ which is also impossible.

The second statement follows immediately from the first since $f^H\in\allpoly$. 
  \end{proof}
  
  It is a consequence of the lemma that the polynomial $f^H(x,1)$ has
  all non-positive roots if $f\in\gsubclosepos_2$.

If we take a product of linear factors where the
coefficients of $x$ and $y$ don't all have the same sign, then the
resulting product isn't in $\rupint{2}$. However, if the constant terms
are all the same sign then the coefficients are in
$\allpoly$. 

\begin{cor}\label{cor:prod-lin-fac}
  Suppose that $b_i\in\reals$, and all $d_i$ have the same sign. If we
  write
$$ \prod_{i=1}^n(y+b_ix+d_i) = \sum f_i(x)y^i$$ then
\begin{enumerate}[(1)]
\item $f_i(x)\in\allpoly$ for $i=0,\dots,n-1$.
\item $f_i(x)$ and $xf_{i+1}(x)$ interlace, for $i=0,\dots,n-2$.
\end{enumerate}
\end{cor}
\begin{proof}
  By switching signs of $x$ and $y$ if necessary we may assume that
  all the $d_i$ are negative. Replace $y$ by $xy$:
$$ \prod_{i=1}^n(xy+b_ix+d_i) = \sum f_i(x)x^iy^i$$ 
Since the determinant of each factor is
$\smalltwodet{1}{b_i}{0}{d_i}=d_i$ which is negative, we can apply
Corollary~\ref{cor:sub-xyc} to conclude that the product is in
$\gsubclose_2$. Consequently, adjacent coefficients interlace, which
finishes the proof. 
\end{proof}

If the roots of a polynomial are all positive or all negative, then
the coefficients are all non-zero. We use this simple idea to
constrain where the roots of coefficients lie.

\begin{cor}\label{ref:coef-interval}
  Suppose that $f,g\in\allpolypos$ and $s<r$. If we define
$$ f(xy-ry)g(xy-sy) = \sum h_i(x)y^i$$
then $h_i(x)\in\allpolyint{[s,r]}$.
\end{cor}
\begin{proof}
  We can apply  Corollary~\ref{cor:sub-xyc} to see that
  $h_i(x)\in\allpoly$. If $\alpha>r$ then all coefficients of
  $f((\alpha-r)y)$ and $g((\alpha-s)y)$ are positive. Thus, all
  coefficients of $ (f(xy-ry)g(xy-sy))(\alpha,y)$ are positive, and
  hence $h_i(\alpha)$ is not zero.  
  
  If $\alpha<s$ then the roots of $f((\alpha-r)y)$ and
  $g((\alpha-s)y)$ are positive. Again, all coefficients of $
  (f(xy-ry)g(xy-sy))(\alpha,y)$ are non-zero, and hence $h_i(\alpha)$ is
  not zero.
\end{proof}

If we define $f_0(x) = J_0(2\sqrt{-x})$, where $J_0$ is the
\index{Bessel function}Bessel function, then we have the identity
\cite{iliev}*{(2.1.A)}:
\begin{equation}\label{eqn:legendre-1}
 f_0(xy+y)f_0(xy-y) = \sum_{n=0}^\infty P_n(x)\frac{(2y)^n}{n!n!}
\end{equation}
where $P_n(x)$ is the \index{Legendre polynomial}Legendre polynomial. 
Since $f_0$ has all positive coefficients, it follows that the
Legendre polynomials have roots in $[-1,1]$, as is well known.
We will also use this identity to determine mapping properties of a
Legendre transformation - see Lemma~\ref{lem:legendre2}.

Any polynomial in $\allpoly$ determines a mutually interlacing
sequence (see Example~\ref{ex:interpolating-mutual}). The following
lemma shows that we can find a polynomial in $\gsubclose_2$ that
interpolates two of these mutually interlacing sequences.  We will
show that given polynomials $f(x),g(x)$ with roots $\{r_i\},\{s_i\}$,
we can find $F(x,y)$ in $\gsubclose_2(n)$ that satisfies
\begin{align*}
  F(x,s_i) &= g'(s_i)\, f_{n+1-i}(x)\\
  F(r_i,y) &= f'(r_i)\, g_{n+1-i}(y)
\end{align*}

\index{interpolation!mutually interlacing}
  \begin{lemma}\label{lem:interp-mutual}
    Suppose that $f(x) = \prod(x-r_i)$, $g(x) = \prod(x-s_i)$ where
    $r_1\le\cdots\le r_n$ and $s_1\le\cdots\le s_n$. Then
    \begin{align}
      \label{eqn:interpolate-1}
      F(x,y) &= \sum_{i=1}^n
      \frac{f(x)}{x-r_i}\,\frac{g(y)}{y-s_{n+1-i}}\in\gsubclose_2\\
      \label{eqn:interpolate-2}
      G(x,y) &= \sum_{i=1}^n
      \frac{f(x+y)}{x+y-r_i}\,\frac{g(x+y)}{x+y-s_{n+1-i}}\in\rupint{2}
          \end{align}
  \end{lemma}
  \begin{proof}
    If we replace $f(x,y)$ by $f(x,\epsilon y)$ and $g(x,y)$ by
    $g(\epsilon x,y)$ where $\epsilon$ is positive, then
    \eqref{eqn:interpolate-2} implies that 
    $$\sum_{i=1}^n \frac{f(x+\epsilon y)}{x+\epsilon
      y-r_i}\,\frac{g(y+\epsilon x)}{\epsilon x+y-s_{n+1-i}}
$$
and letting $\epsilon$ go to zero establishes
\eqref{eqn:interpolate-1}. The homogeneous part of $G(x,y)$ is
$n(x+y)^{2n}$. To check substitution, we choose $\alpha\in\reals$ and
let $f_\alpha(x) = f(x+\alpha)$, $g_\alpha(x) = g(x+\alpha)$. The
result now follows from \ref{lem:convolution}.
  \end{proof}

If $\sum f_i(x)y^i\in\rupint{2}$ then $f_i$ and $f_{i+1}$
interlace. This is half true for $f_i$ and $f_{i+2}$.
\index{positive interlacing}

\begin{lemma}
If  $\sum f_i(x)y^i\in\rupint{2}$ then $f_i \plesslesseq -f_{i+2}$. 
\end{lemma}
\begin{proof}
  Since $y^2-\alpha\in\gsubclose_2$ for positive $\alpha$, the product
\[
(y^2-\alpha)\sum f_i(x)y^i = \sum (f_{i+2}(x)-\alpha f_i(x))y^i
\]
is in $\gsubclose_2$. Thus $f_{i+2}(x)-\alpha f_i(x)\in\allpoly$ for
all positive $\alpha$. 
\end{proof}

\added{12/17/07}
A Hadamard product with a polynomial in $\gsubplus_2$ is effectively
one with a polynomial in $\allpolypos$.

\index{Hadamard product!in $\gsubplus_2$}

\begin{lemma}
  If $f = \sum_0^n f_i(x)y^i\in\gsubplus_1$ has the property that
  $\sum_0^n g_i(x)f_i(x)y^i\in\gsubplus_2$ for all $\sum
  g_i(x)y^i\in\gsubplus_2$ then
$f = F(x)\cdot \sum_0^n a_i y^i$ where $F(x)\in\allpolypos$ and $\sum
a_i\frac{y^i}{i!(n-i)!}\in\allpolypos$. 
\end{lemma}
\begin{proof}
  If we let $g = (y+1)^n$ then 
\[
\sum f_i g_i y^i = \sum \binom{n}{i} f_i(x)y^i\in\gsubplus_2 \]
Applying $f$ again yields 
\[ \sum \binom{n}{i} f_i^2(x) y^i\in\gsubplus_2\]
Since adjacent coefficients interlace, Lemma~\ref{lem:int-squares}
implies that $f_i = F(x) a_i$ for some constant $a_i$. Thus
\[
\sum \binom{n}{i} f_i^2(x) y^i = n! F(x) \sum a_i
\frac{y^i}{i!(n-i)!}\in\gsubplus_2\]
which proves the theorem.
\end{proof}

  \section{Reversal and cancellation in $\rupint{2}$}
\label{sec:top-rev}

In this section we first consider two properties that are trivial for
polynomials in one variable. First, the reversal of a polynomial in
$\allpoly$ is in $\allpoly$, and second, if $xf\in\allpoly$ then
$f\in\allpoly$. Neither is immediate in $\rupint{2}$, and we do not know
the latter to be true if there are more than two variables. Finally,
we discuss differential operators that preserve $\rupint{d}$.

\index{reverse!in $\rupint{2}$}

The key to proving that $xy-1\in\gsubclose$ was the identity 
\[
xy - 1 =
\begin{vmatrix}
x & 1 \\ 1 & y  
\end{vmatrix}
\]
which can be modified to give a sequence of polynomials in $\rupint{2}$
converging to $xy-1$. Now $xy-1$ is the reversal of $x+y$; we have a
similar matrix argument that shows that the reversal is in
$\gsubclose_2$. If we use $1/x$ instead of $-1/x$ \seepage{lem:p2-xy} we get stable
polynomials \seepage{lem:s-basic}.  The key is the following well
known matrix identity for the determinant of a partitioned matrix:

\index{matrix identity}

\begin{equation}\label{eqn:det-iden}
  \begin{vmatrix}
    A_{11}& A_{12} \\ A_{21} & A_{22}
  \end{vmatrix}
=
\bigl| A_{11}\bigr|\,
\begin{vmatrix}
  A_{22} - A_{21}A_{11}^{-1} A_{12}
\end{vmatrix}
\end{equation}
 It follows easily that for matrices $A,B$ of the same size 
\[
| xB-A^2| =
\begin{vmatrix}
  xI & A \\ A & B
\end{vmatrix}
\]

\begin{lemma}\label{lem:reverse-p2-1}
If $f(x,y)\in\rupint{2}(n)$, then $y^n f(x,-1/y)\in\gsubclose_2$.  
Equivalently, if $f=\sum f_i(x)y^i\in\rupint{2}(n)$, then $\sum
f_{i}(x)(-y)^{n-i}\in\gsubclose_2$. 
\end{lemma}
\begin{proof}
If we  write $f(x,y) = |I + xD_1 + yD_2|$ where $D_1,D_2$ are
  positive definite then
\begin{align*}
x^nf(-1/x,y) &= | xI - D_1 + xy D_2| \\
&= |x(I + y D_2) - D_1| \\
& =
\begin{vmatrix}
  x I & \sqrt{D_1} \\ \sqrt{D_1} & I + y D_2
\end{vmatrix}\\
&=
\biggl|
\begin{pmatrix}
  0 & 0 \\0 & I
\end{pmatrix}+
x \begin{pmatrix}
  I & 0 \\0 & 0
\end{pmatrix}+
y\begin{pmatrix}
 0 & 0 \\ 0 & D_2 
\end{pmatrix}+
\begin{pmatrix}
  0 & \sqrt{D_1} \\\sqrt{D_1} & 0
\end{pmatrix}
\biggr|
\end{align*}
This represents $f$ as a determinant of a matrix of size $2n$ that is
clearly the limit of determinants of matrices that are in $\rupint{2}$.
\end{proof}

\begin{remark}

  Here is an alternative geometric argument.
  Define \[f_\epsilon(x,y) = (y+\epsilon x)^n f(x+\epsilon
  y,-1/(y+\epsilon x)),\]
  where $\epsilon$ is positive. Since
  $\lim_{\epsilon\rightarrow0} f_\epsilon = f$, we need to show that
  $f_\epsilon\in\rupint{2}(2n)$. The homogeneous part of $f_\epsilon$ is
  clearly $(x+\epsilon y)^n(y+\epsilon x)^n$, and so all coefficients
  of $f_\epsilon^H$ are positive. If we fix $y$, then the locus of $\{(x+\epsilon
  y,-1/(y+\epsilon x))\}$ as $x$ varies consists of two hyperbolas that
  open in the second and fourth quadrants. Consequently, every
  solution curve of $f$ meets each of these hyperbolas - see
  Figure~\ref{fig:p2-xy} . Thus, for any $y$, there are $2n$
  solutions, and substitution is satisfied.
\end{remark}

\begin{prop} \label{prop:reverse-p2} The two variable reversal of a
  polynomial in $\rupint{2}(n)$ is in $\gsubclose_2(2n)$.  The two
  variable reversal of $f$ is $\rev{f} = x^n y^n f(1/x,1/y)$.
\end{prop}
\begin{proof}
  If we apply the lemma twice we see that
  $x^ny^nf(-1/x,-1/y)\in\gsubclose_2$. Since the substitution of
  $(-x,-y)$ for $(x,y)$ preserves $\rupint{2}$, the result follows.
\end{proof}



Note that the
reversal of $x+y$ is $xy-1$, whereas the two variable reversal of $x+y$ is
$x+y$.

The  reverse of a polynomial in $\rupint{2}$ 
satisfies the assumptions of the following lemma.

  \begin{lemma}\label{lem:p2-fxdy}
  Suppose that $f(x,y)=\sum f_i(x)y^i\in\gsubf_2$ or $\gsubclose_2$, and each $f_i(x)$
  is a polynomial of  degree $i$. 
 If the leading coefficients of the $f_i$ alternate in  sign then
 $$
 g(x,y)\in\rupint{2}(n) \implies
 f(x,-\partial_y)\,g(x,y)\in\rupint{2}(n)$$
  \end{lemma}
  \begin{proof}
    If we choose $\alpha\in\reals$ then $f(\alpha,y)$ is in
    $\allpolyf$. It follows from Lemma~\ref{lem:fpx} that
    $f(\alpha,-\partial_y)\,g(\alpha,y)\in\allpoly$, and so
    $f(x,-\partial_y)\,g(x,y)$ satisfies $y$-substitution.
    
    Suppose the leading coefficient of $f_i(x)$ is $c_i$.  Since the
    sign of $c_i$ alternates the homogeneous part of
    $f_i(x)(-\partial_y)^i \,g(x,y)$ has degree $n$, and equals $c_i
    (-x)^i (\partial^i_y\,g)^H$.  The homogeneous part of
    $f(x,-\partial_y)\,g(x,y)$ has all positive signs (or all
    negative), and so the first part is proved. The second part is
    similar.
  \end{proof}

\index{cancellation}
\index{truncation}

We need a lemma about truncating polynomials in $\allpoly$ before we
can proceed to canceling $x$ in $\gsubclose_2$.
\begin{lemma}
  Suppose that $f\in\allpoly(m-1)$ does not have multiple roots. If
  $f_n\longrightarrow xf$ where all $f_n\in\allpoly(m)$ then
  $f_n(x) - f_n(0)\in\allpoly$ for $n$ sufficiently large.
\end{lemma}
\begin{proof}
  If $r_1,\dots,r_{m-1}$ are the roots of $f^\prime$, then none of
  $f(r_1),\dots,f(r_{m-1})$ are zero since $f$ has no multiple roots.
  Let $\alpha$ be the minimum value of $|f(r_1)|$,\dots,$|f(r_{m-1})|$.
  Since $f_n(0)$ converges to $0$ there is an $N$ such that if $n>N$
  we have $|f_n(0)|< \alpha/2$. And, if we choose $N$ sufficiently
  large we can require that the minimum of $\{|f_n(r) \mid
  f^\prime(r)=0\}$ is at least $\alpha/2$. Notice that the graph of
  $f_n(x)-f_n(0)$ is just the graph of $f_n(x)$ shifted by less than 
  $\alpha/2$. Since all the critical points of $f_n(x)$ are more than
  $\alpha/2$ away from the $x$ axis, $f_n(x)$ has all real roots.
\end{proof}

\begin{prop}\label{prop:cancel-2}
  If $xf\in\gsubclose_2$ then $f\in\gsubclose_2$.
\end{prop}

\begin{proof}
  Since $xf\in\gsubclose_2(m)$ there are polynomials
  $h_n(x,y)\in\rupint{2}(m)$ such that $h_n\longrightarrow xf$. We can
  write $h_n(x,y) = x\,g_n(x,y)+k_n(y)$. Since $g_n\longrightarrow f$
  it suffices to show that $g_n\in\rupint{2}$. Now $h_n$ satisfies the
  homogeneity condition, and therefore so does $g_n$. We need to show
  that $g_n$ satisfies substitution. If we choose $\alpha\in\reals$
  then $x g_n(y,\alpha)+k_n(\alpha)\longrightarrow x f(y,\alpha)$. The
  lemma above implies that $g_n(y,\alpha)\in\allpoly$ for $n\ge
  N_\alpha$. Now we can apply Lemma~\ref{lem:show-in-p2} to conclude
  that there is an $N$ (the maximum of all the $N_\alpha$'s
  corresponding to roots of the resultant) such that $n> N$ implies
  that $g_n(y,\alpha)\in\allpoly$ for all $\alpha\in\reals$.
  
\end{proof}

We would like to prove the converse of this.

\begin{lemma} \label{lem:converse-p2-diff}
  If $f(\partial_x,\partial_y)g$ is in $\gsubclose_2$ for all
  $g\in\gsubclose_2$ then  $f$ is in $\gsubclose_2$. 
\end{lemma}
\begin{proof}
Consider the calculation
  \begin{align*}
    f(x,y) &= \sum_{0\le i,j\le n} a_{ij}x^iy^j \\
    f(\partial_x,\partial_y)x^Ry^R &= \sum_{0\le i,j\le n}
    a_{ij}\falling{R}{i} \falling{R}{j}  x^{R-i}y^{R-j}\\
&= x^{R-n}y^{R-n} \sum_{0\le i,j\le n}
    a_{ij}\falling{R}{i} \falling{R}{j}  x^{n-i}y^{n-j}\\
\intertext{Replacing $x$ by $x/R$ and $y$ by $y/R$ yields}
& x^{R-n}y^{R-n}\sum_{0\le i,j\le n}
a_{ij}\frac{\falling{R}{i} }{R^i}\frac{\falling{R}{j} }{R^j}x^{n-i}y^{n-j} \in\gsubclose_2\\
\intertext{and so we can divide out the
  initial factors to conclude that}
& \sum_{0\le i,j\le n}
a_{ij}\frac{\falling{R}{i} }{R^i}\frac{\falling{R}{j} }{R^j}x^{n-i}y^{n-j} \in\gsubclose_2\\
\intertext{Now letting $R\rightarrow\infty$ shows that}
& \sum_{0\le i,j\le n} a_{ij}x^{n-i}y^{n-j} \in\gsubclose_2\\
  \end{align*}
  Taking the reversal of the last expression and applying
  Proposition~\ref{prop:reverse-p2} finishes the proof.
\end{proof}

\section{Polynomials linear in $y$}
\label{sec:topology-linear-y}

The polynomials of the form $f(x)+yg(x)$ are well understood. First of
all, since $f+\alpha g$ is in $\allpoly$ for all $\alpha$, it follows
that $f$ and $g$ interlace, and their degrees differ by at most $1$.
If $f$ and $g$ have the same degree then we can write $f = \beta g+h$
where $g\lesslesseq h$, so that $f+y g = h + (y+\beta)g$. Substituting
$y = y+\beta$ shows that we can always assume that the degrees differ
by one.  The next two lemmas show that the condition of interlacing is
all we need to belong to $\gsubclose_2$.

\begin{lemma} \label{lem:fyg}
  If $f\longleftarrow g\in\allpoly$ have positive leading
  coefficients then the polynomial $f+yg$ lies in $\gsubclosepos_2$.
\end{lemma}
\begin{proof}
If $f\longleftarrow g$ then we can write
$$
f = (ax+b)g - \sum_1^n c_i^2 \,\frac{g}{x+a_i}$$
where $g =
\prod(x+a_i)$ and $a$ is non-negative. Since $f+yg$ can be expressed
as the determinant
$$
\begin{vmatrix}
  y+ax+b & c_1 & c_2 & \hdots & c_n \\
c_1 & x+a_1 & 0& \hdots& 0 \\
c_2 & 0 & x+a_2  && 0 \\
\vdots & \vdots & & \ddots & \vdots\\
c_n & 0 & \hdots & & x+a_n
\end{vmatrix}
$$
it follows that $f+yg\in\gsubclosepos_2$.
\end{proof}

\begin{lemma} \label{lem:yfg} 
  If $f\longleftarrow g\in\allpoly$ have positive leading
  coefficients then the polynomial $yf-g$ lies in $\gsubclosepos_2$.
  \end{lemma}

\begin{proof}
The reversal of $f+yg$ is $yg-f$, so
$yg-f\in\gsubclose_2$. 

\end{proof}

\begin{remark}
  If we choose $f_0\lesslesseq f_1$ and $g_0\lesslesseq g_1$ then both
  $f_0+yf_1$ and $g_0+yg_1$ are in $\gsubclose_2$. Consequently
  their product
$$ (f_0+yf_1)(g_0+yg_1) = f_0g_0 + (f_0g_1+f_1g_0)y + f_1g_1y^2$$ 
is in $\gsubclose_2,$ and so the coefficients of the product
interlace. This explains the  Leibnitz property (Lemma~\ref{lem:leibnitz}):
\index{Leibnitz}
$$   f_0g_0  \lesslesseq f_0g_1+f_1g_0 \lesslesseq  f_1g_1.$$ 
\end{remark}

The following lemma generalizes the Leibnitz property.
\begin{lemma} \label{lem:top-leibnitz}
  Suppose $f_i\longleftarrow g_i$ have positive leading coefficients
  for $1\le i \le n$ and define
$$ h_k = \sum_{|\diffi|=k} f_\diffi\ g_{\{1,\dots,n\}\setminus\diffi}$$
Then 
$$ h_1 \longleftarrow h_2 \longleftarrow \cdots \longleftarrow h_n$$
\end{lemma}
\begin{proof}
  The $h_k$'s are consecutive coefficients of $\prod_i(f_i+y g_i)$.
\end{proof}

\begin{example}
  In the case $n=3$ the interlacings are
$$ f_1f_2f_3 \longleftarrow f_1f_2g_3+f_1g_2f_3+g_1f_2f_3
\longleftarrow
f_1g_2g_3+g_1f_2g_3+g_1g_2f_3 \longleftarrow g_1g_2g_3.
$$
The first and last interlacings are easy, but the middle one is not so
obvious. 
\end{example}

Here is a similar result for common interlacings.
\index{common interlacing}
  \begin{lemma}
    If $f_1,\dots,f_n$ are in $\allpoly$ and have a common interlacing
    then \[ \sum_{i_1<\cdots<i_s} f_{i_1}f_{i_2}\cdots f_{i_s}\in\allpoly.\]
  \end{lemma}
  \begin{proof}
    If $g\lesslesseq f_i$ then $g + y f_i\in\gsubclose_2$ and
    therefore
$\displaystyle  \prod_{i=1}^n (g + y f_i) \in\gsubclose_2. $
The coefficient of $y^s$ is in $\allpoly$ and is $g^{n-s}$ times the
desired polynomial.
  \end{proof}

\begin{cor}
  Suppose $f = f_1(x)+yf_2(x)$ and $g = g_1(x)+yg_2(x)$ are in
  $\gsubclose_2$ where $f_1,f_2,g_1,g_2$ have positive leading
  coefficients. Then $f \lessless g$ iff
  \begin{enumerate}
  \item $f_1 \lessless g_1$
  \item $f_2 \lessless g_2$
  \item $\smalltwodet{f_1}{f_2}{g_1}{g_2} < 0$
  \end{enumerate}
\end{cor}
\begin{proof}
This follows from Lemma~\ref{lem:inequality-4} and Lemma~\ref{lem:inequality-4b}.
\end{proof}

There is no known criterion for a polynomial that is quadratic in $y$
to belong to $\gsubclose_2$. The following lemma gives a general
construction for such polynomials, but it is not known if there are
quadratics in $y$ in $\gsubclose_2$ that are not given by this form.

\begin{lemma} \label{lem:det-quad}
  Suppose that $f = (x+c_1)\cdots(x+c_n)$, $f_i = \dfrac{f}{x+c_i}$,
  and $f_{ij} = \dfrac{f}{(x+c_i)(x+c_j)}$. For any choice of
  $a_i,b_i$ the following polynomial is in $\gsubclose_2$:
$$ y^2 f - y\left(\sum_i (a_i^2 + b_i^2)f_i\right) 
+ \sum_{i<j} \smalltwodet{a_i}{a_j}{b_i}{b_j}^2f_{ij}$$
\end{lemma}

\begin{proof}
  It suffices to show that the polynomial in question is the value of
  the determinant
\index{determinants!and polynomials in $\rupint{2}$}
\begin{equation} \label{eqn:det-quad}
\begin{vmatrix}
y & 0 & a_1 & \hdots & a_n \\
0 & y & b_1 & \hdots & b_n \\
a_1 & b_1 & x+c_1 & & 0 \\
\vdots & \vdots & &\ddots & \vdots\\
a_n & b_n & 0 & \hdots & x+c_n
\end{vmatrix}
\end{equation}
    
The coefficient of $y^2$ is clearly $f$. To find the coefficient of
$y$, consider the minor corresponding to the first row and column:

$$
\begin{vmatrix}
  y & b_1 & \hdots & b_n \\
b_1 & x+c_1 & & 0 \\
\vdots & & \ddots & \vdots \\
b_n & 0 & \hdots & x+c_n
\end{vmatrix}
$$

We have seen that this determinant has value $yf -\sum b_i^2 f_i$. The
$yf$ term comes from the $y^2f$ term in \eqref{eqn:det-quad}, and the
remaining part is part of the coefficient of $y$. Similarly, the value
of the minor corresponding to the $y$ in the second row and column of
\eqref{eqn:det-quad} is $yf - \sum a_i^2f_i$. Combining these gives the
coefficient of $y$.

The constant term of the determinant \eqref{eqn:det-quad} is the value
of the determinant
\[ 
\begin{vmatrix}
0 & 0 & a_1 & \hdots & a_n \\
0 & 0 & b_1 & \hdots & b_n \\
a_1 & b_1 & x+c_1 & & 0 \\
\vdots & \vdots & &\ddots & \vdots\\
a_n & b_n & 0 & \hdots & x+c_n
\end{vmatrix}
\]

Expanding using the first two rows (and some patience) gives the final
result.

\end{proof}

Here is another interpretation of the last result.
\begin{cor}
  Suppose $f\in\allpoly(n)$, $f\lesslesseq g$, $f\lesslesseq h$, and the leading
  coefficients of $f,g,h$ are all positive. There is a $k\in\allpoly(n-2)$
  so that $$k - y (g+h) + y^2 f \in\gsubclose_2$$
\end{cor}
\begin{proof}
  Use Lemma~\ref{lem:sign-quant} to write $g = \sum a_i^2 f_i$ and $h = \sum
  b_i^2 f_i$, and now apply  Lemma~\ref{lem:det-quad}. 
\end{proof}

We can generalize this corollary to three variables.

\begin{lemma}
  Suppose $f\in\allpoly(n)$, $f\lesslesseq g$, $f\lesslesseq h$, and the leading
  coefficients of $f,g,h$ are all positive. There is a $k\in\allpoly(n-2)$
  so that $$k - y\,g + z\, h + yz\, f \in\gsubclose_3$$
\end{lemma}
\begin{proof}
  Following the evaluation in Lemma~\ref{lem:det-quad}, we find that 
  \begin{gather*}
   F(x,y,z) =  \begin{vmatrix}
z & 0 & a_1 & \hdots & a_n \\
0 & y & b_1 & \hdots & b_n \\
a_1 & b_1 & x+c_1 & & 0 \\
\vdots & \vdots & &\ddots & \vdots\\
a_n & b_n & 0 & \hdots & x+c_n
\end{vmatrix} \\
= yz\,f - z\,h-y\,g + \sum_{i<j} \smalltwodet{a_i}{a_j}{b_i}{b_j}^2f_{ij}
  \end{gather*}
where $f,g,h$ are as in the proof of Lemma~\ref{lem:det-quad}. The polynomial
$k=\sum_{i<j} \smalltwodet{a_i}{a_j}{b_i}{b_j}^2f_{ij}$ is in $\allpoly$
since it is just $F(x,0,0)$. 

\end{proof}

If we take $g=h$ then $k$ is zero, and so $yz\,f -
(y+z)g\in\gsubclose_2$.   If we compute the determinant
\[
 \begin{vmatrix}
z & 0 & 0 & a_1 & \hdots & a_n \\
0 & y & 0 &  a_1 & \hdots & a_n \\
0 & 0 & w &  a_1 & \hdots & a_n \\
a_1 & a_1 & a_1 &   x+c_1 & & 0 \\
\vdots & \vdots & \vdots &  &\ddots & \vdots & \\
a_n & a_n & a_n & 0 &  \hdots & x+c_n
\end{vmatrix} 
\]
then we get $yzw\,f - (yz-yw-zw)g$ which is therefore in
$\gsubclose_4$. More generally, if we compute the determinant of 
\[
\begin{vmatrix}
  diag(y_i) & (A,\dots,A)^t \\
(A^t,\dots,A^t) & diag(x+c_i)
\end{vmatrix}
\] 
where $A=(a_1,\dots,a_n)$ then we find that 
\[
\left(y_1\cdots y_n\right)f - g\left(\sum_j \frac{y_1\cdots y_n
    }{y_j}\right)\in\gsubclose_{n+1}.
\]

  One way to create polynomials of a fixed degree in $y$ is to use
  positive semi-definite matrices.

\index{rank of a matrix}
\index{positive semi-definite matrix}
  \begin{lemma}\label{lem:rank}
    If $D_1$ is an $n$ by $n$ positive definite matrix, and $D_2$ is
    positive semi-definite of rank $r$, then $|I+x
    D_1+yD_2|\in\gsubclose_2$, and has degree $r$ in $y$.
  \end{lemma}
  \begin{proof}
    If we diagonalize $D_2$, then there are exactly $r$ non-zero terms
    on the diagonal, so the maximum degree is $r$. The coefficient of
    $y^r$ is the determinant of the principle submatrix of $I + x D_1$
    corresponding to the zero diagonal entries of $D_2$, and is a
    polynomial in $\allpolypos$ of degree $n-r$.
  \end{proof}

  \begin{example}
    If $v_1,\dots,v_d$ are vectors then
$\displaystyle
\biggl|\sum_i^n x_i\, v_i^t\,v_i\,\biggr|
$
is linear in each variable since $v_i^t\,v_i$ is a semi-definite rank $1$ matrix.
  \end{example}

\index{determinants!and low rank matrices}
  \begin{example}
    There is a simple formula if $D_2$ has rank $1$.   If $v$ is a
    vector of length $n$, $A$ is $n$ by $n$, and $f = |I + xD|$ then
    \begin{align}\label{eqn:low-rank}
      |I + vv^t| &= 1 + v^tv \\
      |I + y vv^t| &= 1 + y vv^t \notag\\
      |A + yvv^t| &= |A| (1+ yv^tA^{-1}v)\notag\\
      |I + xD+yvv^t| &= f(1+yv^t(I+xD)^{-1}v)\notag\\
\intertext{If $D =diag(d_i)$ then}
      |I+xD+yvv^t| &= f + y \sum v_i^2 \frac{f}{1+xd_i}.\notag
    \end{align}
    Since the coefficient of $y$ interlaces the coefficient of $f$, we
    have another proof of part of Lemma~\ref{lem:sign-quant}. If
    $f\lessless g$ then $g = \sum a_i f/(1+xd_i)$ where $a_i>0$. If we
    set $w = (\sqrt{a_i})$ then 
\[
|I + xD_1 + yww^t| = f+ y g\]
This is a different determinant realization of $f+yg$ \seepage{lem:fyg}.
  \end{example}

\section{Constrained polynomials with zero coefficients} \added{12/2/7}
\label{sec:constr-polyn}

We look at some conditions on polynomials in $\gsubclose_2$ that imply
they have a simple form. First of all, there are no non-trivial ways
of adding one variable polynomials to get a polynomial in
$\gsubclose_2$.

  \begin{lemma}\label{lem:p2-f+g}
    If $f(x) + g(y)\in\gsubclose_2$ then at least one of the following holds
    \begin{enumerate}
    \item $f$ or $g$ is constant.
    \item $f$ and $g$ have degree $1$.
    \end{enumerate}
  \end{lemma}
  \begin{proof}
    Assume that neither $f$ nor $g$ is constant.  If $f$ is linear
    then it takes arbitrarily large positive and negative values, so
    by Lemma~\ref{lem:f+c} $g$ is linear. If both $f$ and $g$ have
    degree at least $2$ then the lemma implies that they have degree
    exactly two. Moreover, if $f$ has positive leading coefficient
    then $g$ has negative leading coefficient, but this is impossible
    for then $(f(x) + g(y))^H$ does not have all coefficients of the
    same sign.
  \end{proof}

Linear combinations of the form $f(x) + g(x)h(y)$ are trivial if $deg(h)>1$.

\index{Linear combinations!$f(x)+g(y)$}

  \begin{lemma}\label{lem:p2-f+gh}
    If $f(x) + g(x)h(y)\in\gsubclose_2$ where $deg(h)>1$ then $g$ is a constant
    multiple of $f$.
  \end{lemma}
  \begin{proof}
 If $deg(h)>2$ the set $\bigl\{r\in\reals\,\mid\,
    h(y)+r\in\allpoly\bigr\}$ is bounded. Since $h(y) +
    f(r)/g(r)\in\allpoly$ for all $r$ such that $g(r)\ne0$ we have
    that
\[ 
\biggl\{ \frac{f(r)}{g(r)}\,\mid\, g(r)\ne0\biggr\}
\] 
is bounded. Thus $f/g$ extends to a bounded rational function on
 $\reals$ which is constant by Liouville's theorem.

\index{Liouville's theorem}

If the degree of $h$ is $2$ then we may assume that $f$ and $g$ have no
common factors. If $g$ has a root then $f/g$ takes on arbitrarily
large positive and negative values, so by Lemma~\ref{lem:f+c} $h$ is
linear. Thus $g$ is constant, but this contradicts
Lemma~\ref{lem:p2-f+g} since $h$ has degree $2$. 
  \end{proof}

\index{Newton's inequalities!and zero coefficients in $\gsubclose_2$}

  If $\cdots + a\,x^i + 0\,x^{i+1} + b \,x^{i+2}+\cdots\in\allpoly$ then
  Newton's inequality implies that $ab\le0$.  The analogous result in
  two variables is much stronger. 

\index{zero coefficients!in $\rupint{2}$}

  \begin{lemma}\label{lem:p2-iso-zero}
    If $f(x,y) = \cdots + f_i(x)\,y^i + 0\,y^{i+1} +
    f_{i+2}(x)\,y^{i+2}+\cdots \in\gsubclose_2$ and $f_{i}$ and
    $f_{i+2}$ are not zero then there are $g,h\in\allpoly$ so that
    $f(x,y) = g(x)h(y)$.
  \end{lemma}
  \begin{proof}
    If we differentiate, reverse and differentiate with respect to $y$
    sufficiently many times we find that there are  constants
    $\alpha,\beta$ so that $\alpha\,f_i(x) + \beta\,
    f_{i+2}(x)\,y^2\in\rupint{2}$.  By Lemma~\ref{lem:p2-f+gh} we see that 
    $f_{i+2}$ is a constant multiple of $f_i$.

    We can thus write 
\[ 
f(x,y) = \cdots + f_{i-1}(x)\,y^{i-1} + f_i(x)\,(a + b\, y^2) y^i + \cdots 
\]
for constants $a,b$. Differentiating and reversing yields constants
$c_i$ so that
\[
c_{i-1}\, f_{i-1}(x)\,y^{i-1} + f_i(x)\,(c_i + c_{i+1}\, y^2)\, y^i
\in\gsubclose_2
\]
Applying the lemma again shows that $f_{i-1}$ is a multiple of
$f_i$. Continuing, we see that we can write
\[
f(x,y) = f_i(x)\, H(y) + f_{i+1}(x)\,y^{i+1} + \cdots \]
Reversing, we apply the same argument to finish the proof.
  \end{proof}

  \begin{cor}
    If $f\in\gsubclose_2$ has all exponents of even degree, then
    $f(x,y)=g(x)h(y)$ where $g,h\in\allpoly$.
  \end{cor}
  \begin{proof}
    If we write $f(x,y) = \sum f_i(x)y^i$ then if $i$ is odd the
    hypotheses imply that $f_i(x)=0$. The conclusion now  follows from the lemma.
  \end{proof}

\section{The polynomial closure of $\rupint{d}$}
\label{sec:sub-closure}

If $V$ is the vector space of a polynomials in $x_1,\dots,x_d$ of
total degree at most $n$, then $\rupint{d}(n)$ is contained in $V$.
The closure of $\rupint{d}(n)$ in $V$ is denoted $\gsubclose_d(n)$, and
the closure of $\gsubplus_d$ is $\gsubplusclose_d$.  We define
$$ \gsubclose_d = \gsubclose_d(1) \cup \gsubclose_d(2) \cup \dots.$$
In this section we show that most properties of $\rupint{d}$ apply to its 
closure. We also identify some specific polynomials in the closure
that are not in $\rupint{d}$.

\index{\ Pgsubclose@$\gsubclose_d$} 
\index{\ Pgsubclosen@$\gsubclose_d(n)$}

First of all, all monomials are in the closure.

\begin{lemma}
  All monomials of at most $d$ variables are in $\gsubclose_d$.
\end{lemma}
\begin{proof}
  The monomial $x_1^{i_1}\cdots x_d^{i_d}$ is the limit of these polynomials
in $\rupint{d}:$
$$ (x_1+ \epsilon (x_2+\cdots+x_d))^{i_1}\cdots(
\epsilon(x_1+\cdots+x_{d-1})+x_d)^{i_d}.$$ 
\end{proof}

Next, we have some important containments.

\begin{lemma}
  If $d<e$ then 
  $\gsubclosepos_d \subset \gsubclosepos_e$. In the case $d=1$ we have
  $\allpoly\subset \gsubclosepos_2$, $\allpolypos\subset
  \gsubplusclose_2$.
\end{lemma}
\begin{proof}
  It suffices to assume $e=d+1$. Choose $f\in\gsubclose_d(n)$. Since
  $f$ is in the closure of $\rupint{d}(n)$ there is a sequence
  $f_1,f_2,\dots$ of polynomials in $\rupint{d}(n)$ that converge to $f$.
  For any positive $\epsilon$ we know that $f_i(x_1+\epsilon
  x_{d+1},x_2,\dots,x_d)$ is in $\rupint{d+1}$. Since $\lim_{\epsilon\rightarrow0}
  f_i(x_1+\epsilon x_{d+1},x_2,\dots,x_d) = f_i$ we see that $f$ is the
  limit of polynomials in $\rupint{d+1}$. It is clear that if the
  homogeneous part of $f$ has all positive coefficients, then so does
  the homogeneous part of \\ 
  $f(x_1+\epsilon x_{d+1},x_2,\dots,x_d)$.
  This implies that $\gsubclosepos_d \subset \gsubclosepos_e$.

  If $f\in\allpoly(n)$ then $f(x+\epsilon y)$ has homogeneous part
  $c(x+\epsilon y)^n$ and so all $f(x+\epsilon y)$ are in $\rupint{2}$. 
  If $f\in\allpolypos$ then all its coefficients are positive, and so
  the same is true of $f(x+\epsilon y)$. 
\end{proof}

We can use determinants to  realize members of $\gsubclosepos_d$.  The
proof is the same as Lemma~\ref{lem:top-det-2}.
\begin{lemma} \label{lem:top-det}
  Suppose that $C$ is a symmetric matrix, and $D$ is a
  diagonal matrix where the $i$-th diagonal entry is $\sum d_{ij}x_i$
  where $d_{ij}$ is non-negative. The polynomial $det(C+D)$ is in
  $\gsubclosepos_d$. 
\end{lemma}

Before we go too far, let's see that there are polynomials that are
\emph{not} in $\gsubclose_d$. Consider this  simple consequence of  Lemma~\ref{lem:pd-xxg} :

\begin{lemma}
  If  $f(x,\dots,x)$ does not have all real roots then $f\not\in\gsubclose_d$.
\end{lemma}

This implies that $xyz+1$ is not in $\gsubclose_3$, because $\pm x^3+1$
does not have all real roots. We can slightly generalize the lemma to
say that if we replace some of the variables by $\pm x$ and the
resulting polynomial is not in $\gsubclose_d$ then neither is
$f$. This shows that $xy+z$ is not in $\gsubclose_3$ since $\pm x^2+z$ 
is not in $\gsubclose_2$ by Example~\ref{ex:not-xxy}. 

There are simple polynomials in $\gsubclose_3$. Since 
$$\begin{vmatrix}
  x&1&0\\1&y&1\\0&1&z
\end{vmatrix} = xyz-x-z
$$
it follows from Lemma~\ref{lem:top-det} that $xyz-x-z\in\gsubclose_3$.

\begin{lemma} \label{lem:elem-sym}
  The elementary symmetric functions are in $\gsubclosepos$.
\end{lemma}
\begin{proof}
  The $i$-th elementary symmetric function of $x_1,\dots,x_d$ is the
  coefficient of $y^i$ in the product
  $$
  (y+x_1)(y+x_2)\cdots(y+x_d)$$
  Since each of the factors is in
  $\gsubclosepos_{d+1}$, the product is in $\gsubclosepos_{d+1}$, and
  so the coefficients of $y$ are in $\gsubclosepos_d$.
\end{proof}

As long as there is no obvious restriction,  the interior of
$\xsub_d$ is contained in $\gsubclose_d$.

\begin{lemma}
  If $f\in\, int\,\xsub_d$ and $f^H$ has all non-negative terms then
  $f\in\gsubclose_d$. 
\end{lemma}
\begin{proof}
  With these assumptions on $f$ the reason that $f$ might not be in
  $\rupint{d}$ is that $f^H$ might have some zero coefficients. If $f$
  has total degree $n$ then we write
  \begin{xalignat*}{2}
    f & = \sum_\diffi c_\diffi \xx^\diffi &
    g & = \sum_{|\diffi|=n,c_\diffi=0}  \xx^\diffi 
  \end{xalignat*}
  For any $\alpha>0$ the polynomial $f+\alpha g$ satisfies the
  homogeneity condition. Let $\mathcal{O}$ be an open neighborhood in
  $\xsub_d(n)$ that contains $f$. For sufficiently small $\alpha$ we
  have that $f+\alpha g$ is in $\mathcal{O}$, so for these $\alpha$ it
  follows that $f+\alpha g$ satisfies substitution.  Since the
  homogeneous part of $f+\alpha g$ has all positive coefficients
  $f+\alpha g$ is in $\gsubpos_d(n)$, and converge to $f$.
\end{proof}

Recall the definition of interlacing:
\begin{quote}
  $f,g\in\gsubposclose_d$ \emph{interlace} iff $f+\alpha g$ is in
  $\gsubposclose_d$ for every $\alpha$.
\end{quote}

\index{interlacing!in $\gsubposclose_d$}

The following are all consequences of the corresponding results for
$\gsubpos_d$.

\begin{theorem} \label{thm:gclose}
 Suppose $f\in\gsubclosepos_d$.
  \begin{enumerate}
  \item If  $g\in\gsubclosepos_e$  then
    $fg\in\gsubclosepos_{d+e}$. 
  \item If  $f\lesslesseq g,h$ and
    $g,h\in\gsubclosepos_d$ then $f \lesslesseq g+h$ and
    $g+h\in\gsubclosepos_{d}$.
  \item $f$ and $\frac{\partial f}{\partial x_i}$ interlace.  
  \item Linear transformations extend from
    $\allpoly\longrightarrow\allpoly$ to
    $\gsubclosepos_d\longrightarrow\gsubclosepos_d$.
  \item Polynomial limits of polynomials in $\gsubclosepos_d$ are in
    $\gsubclosepos_d$.
  \item If $f = \sum a_\diffi(\xx)\yy$ then
    $a_\diffi(\xx)\in\gsubclosepos_d$
  \end{enumerate}
\end{theorem}



We can generalize Lemma~\ref{lem:fyg} to $\gsubpos_d$, but we no longer have a nice
determinant representation.

\begin{lemma}\label{lem:fyg-d}
  If $f\lesslesseq g\in\gsubpos_d$ then $f+x_{d+1}g\in\gsubposclose_{d+1}$.
\end{lemma}
\begin{proof}
  To simplify notation we assume that $d=3$. Define a polynomial in
  four variables
  $$
  F_\epsilon = f(x+\epsilon u,y+\epsilon u,z+\epsilon u) +
  ug(x+\epsilon u,y+\epsilon u,z+\epsilon u). $$
  Since
  $\lim_{\epsilon\rightarrow0} F_\epsilon = f+ug$ it suffices to
  verify that $F_\epsilon\in\gsubpos_4$.  The homogeneous part
  satisfies
  $$
  F_\epsilon^H = f^H(x+\epsilon u,y+\epsilon u,z+\epsilon u) +
  ug^H(x+\epsilon u,y+\epsilon u,z+\epsilon u). $$
  Since both $f^H$
  and $g^H$ have all positive terms, it follows that $F_\epsilon^H$
  also does.  It remains to check substitution.  Setting
  $y=\alpha,z=\beta,u=\gamma$ we need to show that
$$
F_\epsilon(x,\alpha,\beta,\gamma) = f(x+\epsilon
\gamma,\alpha+\epsilon \gamma,\beta+\epsilon \gamma) + \gamma
g(x+\epsilon \gamma,\alpha+\epsilon \gamma,\beta+\epsilon \gamma) $$
is in $\allpoly$. Since $f\lesslesseq g$ in $\gsubpos_3$ it follows
that  
$$ f(w,\alpha+\epsilon \gamma,\beta+\epsilon \gamma) \lesslesseq 
g(w,\alpha+\epsilon \gamma,\beta+\epsilon \gamma)
$$
Since any linear combination of these two polynomials is in $\allpoly$
we add $\gamma$ times the second one to the first and replace $w$ by
$x+\epsilon\gamma$ to finish the proof.
\end{proof}

If we don't introduce a new variable then the result is easier.
\begin{cor} \label{cor:xffp}
  If $f\greateqeq g$ in $\gsubpos_d$ then $x_i g+f\in\gsubposclose_d$
  for $1\le i\le d$.
\end{cor}
\begin{proof}
  Use Lemma~\ref{lem:pd-xffp}
\end{proof}

The polynomial $f+yg+zh+yzk$ is in $\gsubclose_3$ if $f,g,h,k$ satisfy
certain interlacing assumptions.  Here's a more general result.

\begin{lemma}
  If $\yy=(y_1,\dots,y_d)$ and  $F=\sum_\sdiffi f_\sdiffi(x)g_\sdiffi(\yy)$ 
   satisfies
  \begin{enumerate}
  \item   $x$-substitution
  \item $deg\,( f_{00\dots0}(x)) \ge deg\,(f_\sdiffi(x))+ deg \,(g_\sdiffi(\yy))$
    for all $\diffi$
  \item all $f_\sdiffi$ and $g_\sdiffi$  have all non-negative coefficients
  \end{enumerate}
then $F\in\gsubclosepos_{d+1}$
\end{lemma}
\begin{proof}
  The notation $f_{0\dots0}$ refers to the terms that contain no
  $y_i$'s. Let $G_\epsilon = F(x+\epsilon y,\yy)$ where $y=y_1+\dots+
  y_d$ and $\epsilon>0$. If $n= deg(f_{0\dots0}(x))$ then condition two
  implies that $n\ge deg\,(f_\sdiffi(x+\epsilon y)g_\sdiffi(\yy))$ so
  the degree of $G_\epsilon$ is $n$. The homogeneous part of
  $G_\epsilon$ includes a term $(x+\epsilon y)^n$, and the third
  condition guarantees that all other contributions are non-negative,
  so $G_\epsilon$ has all positive terms. If we let
  $\aaa=(a_1,\dots,a_d)$ then $G_\epsilon(x,\aaa) = F(x+\epsilon
  a,\aaa)$ where $a=a_1+\dots+a_d$.  Now $F(x,\aaa)$ has all real
  roots by the first condition, and so the translate by $\epsilon a$ is also in
  $\allpoly$, and so $G_\epsilon(x,\aaa)$ satisfies $x$-substitution.
  Since $\lim _{\epsilon\rightarrow0}G_\epsilon = F$, it follows that
  $F\in\gsubclosepos_{d+1}$.
\end{proof}

\begin{cor} \label{cor:f00}
  Suppose that $$F(x,y,z)=f_{00}(x)+yf_{10}(x)+zf_{01}(x)+yzf_{11}(x)$$
  satisfies
  \begin{enumerate}
  \item $F(x,a,b)\in\allpoly$ for all $a,b\in\reals$.
  \item $f_{00}\in\allpolypos$.
  \item $deg(f_{00}) \ge deg(f_{10})+1$
  \item $deg(f_{00}) \ge deg(f_{01})+1$
  \item $deg(f_{00}) \ge deg(f_{11})+2$.
  \item All $f_{ij}$ have non-negative coefficients
  \end{enumerate}
then $F\in\gsubclosepos_3$.
\end{cor}

\begin{remark}
  We can multiply the Taylor series \index{Taylor series} and extract
  coefficients to derive some properties  of derivatives
  (Theorem~\ref{thm:gsub-diff2}, Corollary~\ref{cor:fpp}). Starting
  with $f\in\rupint{2}$ and writing
  \begin{align*}
    f(x+u,y+v) &= f(x,y)+ f_x(x,y)u + f_y(x,y)v\, + \\
&    f_{xx}(x,y)\frac{u^2}{2} + f_{xy}(x,y)uv +
    f_{yy}(x,y)\frac{y^2}{2} +\dots
  \end{align*}
we first multiply by $u+\alpha$. The coefficient of $u$ is $\alpha f_x
+ f$. Since this is in $\gsubposclose_2$ for all $\alpha$ we conclude
$f\lesslesseq f_x$.

Next, choose $\alpha>0$ and multiply by
$u^2-2\alpha\in\gsubclose_4$. The coefficient of $u^2$ is $\alpha
f_{xx} -f$, so we conclude that $f-\alpha f_{xx}\in\gsubpos_2$ for all
positive $\alpha$. If instead we multiply by $1-\alpha uv\in\gsubpos_2$ we
find that the coefficient $f- \alpha f_{xy}$ of $uv$ is in $\gsubpos_2$. 

Since $(u+1)(v+1)-1 = uv+u+v$ is in $\gsubposclose_2$, the coefficient
$f_x+f_y+f_{xy}$ of $uv$ is in $\gsubpos_2$.

In general, if $g(u,v) = \sum a_{rs} u^rv^s\in\gsubposclose_2$ then
\begin{align*}
  f(x+u,y+v) &= \sum \diffd_x^i \diffd_y^j f\, \frac{u^i}{i!}
  \frac{v^j}{j!}\\
f(x+u,y+v)g(u,v) &= 
\sum a_{rs} \diffd_x^i \diffd_y^j f\, \frac{u^{i+r}}{i!}
  \frac{v^{j+s}}{j!}\\
\intertext{The coefficient of $v^nu^m$ is}
& \sum \frac{a_{n-i,m-j}}{i!j!} \diffd_x^i \diffd_y^j f
\end{align*}
and is in $\gsubpos_2$.
\end{remark}

\begin{lemma}\label{lem:pd-limit}\added{11/17/07}
    Suppose that $f(x_1,\dots,x_d)=\sum \aaa_\sdiffi \xx^\sdiffi$ is a non-zero
    polynomial, and define
\[
g(\xx) = f\bigl(\epsilon_{10}+\sum \epsilon_{1j}x_j,\dots,\epsilon_{d0}+\sum\epsilon_{dj}x_j\bigr)
\]
If $\diffi$ is an index for which $\aaa_\sdiffi=0$ then the
coefficient of $\xx^\sdiffi$ in $g(\xx)$ is a non-zero polynomial in
the $\epsilon_{ij}$.
  \end{lemma}
  \begin{proof}
    A monomial $\aaa_\sdiffj\xx^\sdiffj$ in $f$ contributes a product 
\begin{equation}\label{eqn:pd-limit}
\aaa_\sdiffj \prod_{k=1}^d \left( \epsilon_{k0}+ \sum_i \epsilon_{ki}x_i\right)^{j_k}
\end{equation}
to $g(\xx)$. Every term in this expansion has the form
$
\epsilon_{1?}^{j_1}\cdots\epsilon_{d?}^{j_d}.
$
where each ``?'' represents an index in $0,\dots,d$.
Consequently the coefficient of $\xx^\sdiffi$ in $g(\xx)$ is a sum of
terms of the form 
\[
\text{(coefficient of $\xx^\sdiffj$)} \times\text{(at most one term arising
  from $\xx^\sdiffj$)}
\]
and so there is no cancellation. Note that \eqref{eqn:pd-limit} has
terms of every index of degree at most $deg(f)$.
  \end{proof}

  \begin{lemma}\label{lem:p-limit}
    If $f(\xx)\in\gsubclose_d$ has non-negative coefficients then it is a
    limit of polynomials in $\gsubplus_d$.
  \end{lemma}
  \begin{proof}
    Define
    \[ f_\epsilon(\xx) = f\bigl(\epsilon_{10}+\sum
    \epsilon_{1j}x_j,\dots,\epsilon_{d0}+
    \sum\epsilon_{dj}x_j\bigr).\] By Lemma~\ref{lem:pd-limit} the
    terms in $f_\epsilon$ with coefficient zero in $f$ are a sum of
    polynomials in the $\epsilon_{ij}$ with positive coefficients in
    $g(\xx)$. Observe that $f$ is the limit of $f_\epsilon$ through
    positive values of $\epsilon_{ii}\rightarrow1$ and $\epsilon_{i\ne
      j}\rightarrow0$. All the $f_\epsilon$ are in $\rupint{d}$, have all
    positive coefficients for all indexes of degree at most $deg(g)$,
    and so $f$ is a limit of polynomials in $\gsubplus_d$.
  \end{proof}

Interlacings in $\gsubclose_d$ can be approximated by interlacings in $\rupint{d}$.
  \begin{cor}
    If $f\longleftarrow g$ in $\gsubclose_d$ then there are
    $f_n,g_n$  such that
    \begin{itemize}
    \item $\lim f_n = f$
    \item $\lim g_n = g$
    \item     $f_n\longleftarrow g_n$
    \end{itemize}
  \end{cor}
  \begin{proof}
    Since $f+yg\in\gsubclose_{d+1}$ we let $h_\epsilon$ be the
    polynomial determined by the lemma above that converges to
    $f+yg$. If $h_\epsilon = f_\epsilon(\xx) + y g_\epsilon(\xx) +
    \cdots$ then $f_\epsilon,g_\epsilon$ are the desired polynomials.
  \end{proof}

  \section{Multiaffine polynomials}
  \label{sec:mult-polyn}

  A \emph{multiaffine polynomial}\index{multiaffine polynomial} is a
  polynomial that is linear in each variable. We look at a few
  properties of such polynomials; for a more detailed discussion see
  [bbs]. 

  \begin{definition}
    \[ \multiaff{d} = \bigl\{ f\in\gsubclose_d \mid 
    \text{all variables have degree 1}\bigr\}
\]
  \end{definition}

We can write a multiaffine polynomial as
\[ \sum_{\sdiffi\subset\{1,\dots,d\}} \aaa_\sdiffi \xx_\sdiffi \] A
multiaffine polynomial in $\multiaff{d}$ can be visualized as an
assignment of real numbers to the $2^d$ vertices of the $d$-cube. 

There is no simple way of determining if a polynomial $f(\xx)$ is
multiaffine. Br\"and\'en\cite{branden-hpp} proved that a necessary and sufficient condition
is that 
\[ 
\frac{\partial f}{\partial x_i}\, \frac{\partial f}{\partial x_j} -
f\,\frac{\partial^2 f}{\partial x_i\partial x_j} \ge 0 \quad\text{ for all $1\le i,j\le n$.}
\]
In the case $d=2$ this is  the simple criterion 
\[
a + bx + cy + d xy \in\multiaff{2} \Leftrightarrow ad-bc \le 0
\]

\begin{example}\label{ex:ma-xyz}
  Here is an application for a polynomial of total degree
  $3$. Consider the question:
\[
\text{For which $m$ is } f(x,y,z)=xyz - x-y-z + m  \in\multiaff{3}?
\]
Using the criterion we find
\[
\frac{\partial f}{\partial x} \frac{\partial f}{\partial y} - 
f\cdot \frac{\partial^2 f}{\partial x \partial y}
=
z^2 - zm +1
\]
This is non-negative if and only if $|m|\le 2$. A similar calculation
shows

\[
xyz - ax-by-cz + m  \in\multiaff{3} \text{ iff } a,b,c\ge0 \text{ and
} m^2 \le 4abc
\]

\end{example}

Here are three examples of general multiaffine polynomials in
$\multiaff{d}$; the first two are trivial to show multiaffine, but the
last is the Grace polynomial \index{Grace polynomial} $\grace{d}$ and does not
have an elementary proof that it is in $\multiaff{d}$.

\begin{align*}
  \prod_1^d (x_i + y_i) \\
  \prod_1^d (x_iy_i-1) \\
  \sum_{\sigma\in sym(d)} \prod_1^d (x_i + y_{\sigma i})
\end{align*}

It's easy to construct multiaffine polynomials using this lemma:

\begin{lemma}
  If $f\in\gsubclose_d$ then the multiaffine part of $f$ is in $\multiaff{d}$.
\end{lemma}
\begin{proof}
  The usual reversal, differentiation, reversal shows that if $\sum
  f_i(\xx)y^i\in\gsubclose_d$ then
  $f_0(\xx)+f_1(\xx)y\in\gsubclose_d$. Repeating this for each
  variable shows that the multiaffine part of $f$ is in $\gsubclose_d$.
\end{proof}

\begin{example}
  If we let $f = (1+x+y+z)^n$ then the multiaffine part is
\[
1 + n(x+y+z) + n(n-1)(xy+xz+yz) + n(n-1)(n-2)xyz 
\]
and is in $\multiaff{3}$.
\end{example}

This last example is an example of a multiaffine polynomial that is
symmetric in its variables. See Lemma~\ref{lem:ma-sym}.

Multiaffine polynomials are closed under multiplication in the following sense

\begin{cor}
  If $f,g\in\multiaff{d}$ the the multiaffine part of $fg$ is also in $\multiaff{d}$.
\end{cor}

We can also construct multiaffine polynomials using determinants. If
$S$ is symmetric and $D_1,\dots,D_d$ are positive definite then the
multiaffine part of 
\[ |S + x_1D_1 + \dots + x_dD_d| \]
is multiaffine. We don't always have to take the multiaffine part.
If $v_1,\cdots,v_d$ are vectors then
\[ |S + x_1 v_1^tv_1 + \cdots + x_d v_d^tv_d |\]
is multiaffine.

There are no interesting multiaffine bilinear forms.

\begin{lemma}\label{lem:multi-bilinear}
  If $f = \sum_1^n a_{ij}x_i y_j\in\multiaff{2d}$ then we can
  write
\[ f = \bigl(\sum_1^n \alpha_i x_i\bigr)\cdot\bigl(\sum_1^n
\beta_i y_i\bigr)
\]
where all $\alpha_i$ have the same sign, and all $\beta_i$ have the
same sign.
\end{lemma}
\begin{proof}
  Write $(a_{ij}) = (v_1,\dots,v_d)^t$ so that
\[
f = (x_1,\dots,x_d)\begin{pmatrix}v_1\\\vdots\\ v_d\end{pmatrix}
\begin{pmatrix}y_1\\\vdots\\ y_d\end{pmatrix}
\]
We see $\frac{\partial^2f}{\partial x_1 \partial x_j}=0$, and hence
$\frac{\partial f}{\partial x_1}\frac{\partial f}{\partial x_j}\ge0$.
If $\yy = (y_1,\dots,y_d)^t$ then $\frac{\partial f}{\partial x_1} =
v_1\cdot y$ and $\frac{\partial f}{\partial x_j} =
v_j\cdot y$. Consequently,
\[ (v_1\cdot \yy)(v_j\cdot \yy) \ge0
\]
If the product of two linear forms is never negative then they must be
multiples of one another. Thus we can write $v_j = \alpha_j v_1$ for
some constant $\alpha_j$ which leads to the the desired form
\[ f = \bigl[ \alpha_1x_1 + \cdots + \alpha_d x_d\bigr] \cdot \bigl[ 
v_1\yy\bigr]
\]
 Each of the factors is in $\multiaff{d}$ so all the
coefficients of a factor have the same sign.
\end{proof}

The following might be true for arbitrary multiaffine polynomials with
positive coefficients.

\index{Hadamard product}

\begin{lemma}
  Polynomials with positive coefficients in $\multiaff{2}$ are closed
  under Hadamard product.
\end{lemma}
\begin{proof}
  We know $a+bx + cy + dxy$ is in $\multiaff{2}$ iff
  $\smalltwodet{b}{a}{d}{c}\ge0$. If $a,b,c,d$ are positive then
  $\smalltwodet{b}{a}{d}{c}$ is positive semi-definite. Since positive
  semi-definite matrices are closed under Hadamard product, the
  conclusion follows.
\end{proof}

  \section{Multiaffine polynomials with complex coefficients}
  \label{sec:mult-polyn-with}

  Multiaffine polynomials with complex coefficients are much more
  complicated than those with real coefficients. A multiaffine
  polynomial with real coefficents determines polynomials with complex
  coefficients, but not all of them.  We give this simple
  construction, and then restrict ourselves to multiaffine polynomials
  with two variables, where we will see an equivalence with certain
  \Mobius\ transformations.

  Suppose that $f(x_1,\dots,x_d)\in\multiaff{d}$. If
  $\sigma_1,\dots,\sigma_d$ lie in the upper half plane then 
  $f(x_1+\sigma_1,\dots,x_d+\sigma_d)$ is multiaffine, has complex
  coefficients, and is non-vanishing for $\xx\in\uhp$.

  So, assume that $f(x,y) = \alpha x + \beta  + \gamma y x + \delta
  y$. Solve $f(x,y)=0$ for $y$:
\[
y = - \frac{\alpha x+ \beta}{\gamma x+ \delta}.
\]
If $M$ is the \Mobius\ transformation with matrix $\smalltwo
{\alpha}{\beta}{\gamma}{\delta}$ then $y = - M(x)$. If $x\in\uhp$ and
$f\in\hb{2}$ then $y\in\complexes-\uhp$, for otherwise we have a
solution with both $x,y$ in the upper half plane. Thus, $x\in\uhp$
if and only if  $M(x)\in\uhpc$. Consequently,
\[
f \text{ is complex multiaffine } \Longleftrightarrow
M:\uhp\longrightarrow\uhpc.
\]

For the rest of this section we determine conditions on a \Mobius\
transformation to map $\uhp\longrightarrow\uhpc$.  

If $M$ maps the upper half plane to itself then the image of the upper
half plane is either a half plane parallel to the real line, or the
interior of a circle lying in the upper half plane.

Suppose that $\smalltwo{\alpha}{\beta}{\gamma}{\delta}$ is a matrix
with complex entries that maps the upper half plane to itself.  Since
our entries are complex, we can divide by the square root of the
determinant, and so without loss of generality we may assume that
\textbf{the determinant is one}.  In the case the image of the upper
half plane is a half plane contained in the upper half plane we have
\[
\begin{pmatrix}
  \alpha&\beta\\\gamma&\delta
\end{pmatrix} = 
\begin{pmatrix}
  e & f  \\ g& h
\end{pmatrix}+
\imag\,t\,\begin{pmatrix}
  g & f  \\ 0 & 0
\end{pmatrix}
\]
where $\smalltwo{e}{f}{g}{h}$ is a matrix with real entries and
determinant $1$, and $t\ge0$. 

The more interesting case is when the image of the real line is a
circle constrained in the upper half plane.  The image of the upper half
plane is either the interior or the exterior of the circle. We first
find the center of the circle\footnote{Thanks to David Wright for this
  information.}

\[
  \text{center}  =
  \frac{\beta\overline{\gamma}-\alpha\overline{\delta}}{\delta\overline{\gamma}-\gamma\overline{\delta}} 
\]

\noindent
The image of $\infty$ is $\alpha/\gamma$, so the radius is the distance from the
center to $\alpha/\gamma$:
\[
  \text{radius} = \frac{1}{|2 \Im(\gamma\overline{\delta})|}
\]

The image of the circle is entirely in the upper half plane if the
center lies in the upper half plane, and the radius is less than the
imaginary part of the center.  This gives the middle condition
below. If this holds, then the image of the upper half plane lies in
the circle if and only if $0$ is not the image of a point in the upper
half plane. Since $-\beta/\alpha$ maps to zero, we have the three
conditions
\begin{align*}
  \Im(\text{center}) &>0 &
\Re\left(
  \begin{vmatrix}
    \overline{\alpha} &\beta\\\overline{\gamma}&\delta
  \end{vmatrix}\right) & \ge1 &
  \Im(\beta/\alpha) &\ge 0
\end{align*}

It is more interesting to write these conditions in terms of the real
and imaginary parts. Write
\begin{align*}
  \alpha &= a_1 + \imag a_2 &
  \beta &= b_1 + \imag b_2  \\
  \gamma &= c_1 + \imag c_2 &
  \delta &= d_1 + \imag d_2 
\end{align*}

Note that if we want to show $\Im(\sigma/\tau)>0$ it is enough to show
$\Im(\sigma\overline{\tau})>0$. The three conditions are

\begin{align*}
  2 \left(b_1 c_1+b_2 c_2-a_1 d_1-a_2 d_2\right) \left(-c_1 d_2+c_2
    d_1\right)&  >0 & \text{(center)} \\
  -b_1 c_1-b_2 c_2+a_1 d_1+a_2 d_2 & \ge 1 &\text{(radius)} \\
  a_1 b_2-a_2 b_1&  >0 & \text{(interior)}
\end{align*}

Now the determinant of $M$ is one, which gives us two equations

\begin{align*}
  -b_2 c_1-b_1 c_2+a_2 d_1+a_1 d_2 & = 0 & \text{(imaginary  part)} \\
  -b_1 c_1+b_2 c_2+a_1 d_1-a_2 d_2 & = 1 & \text{(real part)}
\end{align*}

If we reexpress these in terms of determinants we get
\begin{align*}
  \biggl(
  \begin{vmatrix} a_1 & b_1 \\ c_1 & d_1 \end{vmatrix}
  +
  \begin{vmatrix} a_2 & b_2 \\ c_2 & d_2 \end{vmatrix}
  \biggr)
  \begin{vmatrix} c_1 & c_2 \\ d_1 & d_2 \end{vmatrix} 
  & > 0 & \text{(center)}\\
  \begin{vmatrix} a_1 & b_1 \\ c_1 & d_1 \end{vmatrix}
  +
  \begin{vmatrix} a_2 & b_2 \\ c_2 & d_2 \end{vmatrix}
  & \ge 1 &  \text{(radius)} \\
  \begin{vmatrix} a_1 & b_1 \\ c_1 & d_1 \end{vmatrix}
  -
  \begin{vmatrix} a_2 & b_2 \\ c_2 & d_2 \end{vmatrix}
  & = 1 &  \text{(real part)} 
\end{align*}

and these simplify to four conditions

\begin{align*}
  \begin{vmatrix} a_1 & b_1 \\ c_1 & d_1 \end{vmatrix} & > 1 & 
  \begin{vmatrix} a_2 & b_2 \\ c_2 & d_2 \end{vmatrix} & \le 0 & 
  \begin{vmatrix} a_1 & a_2 \\ b_1 & b_2 \end{vmatrix} & > 0 & 
  \begin{vmatrix} c_1 & c_2 \\ d_1 & d_2 \end{vmatrix} & > 0 & 
\end{align*}

Expressed in terms of the original coefficients this is

\begin{align*}
  |\Re(M)| &\ge 1 & |\Im(M)| & \le 0 & \Im(\alpha/\beta)& \le0 &
  \Im(\gamma/\delta)& \le0
\end{align*}

\begin{example}
  The \Mobius\ transformation with matrix
  $M=\smalltwo{1+\imag}{\imag}{2+\imag}{1+\imag}$ maps the real line to
  the circle with center $1+\imag/2$ and radius $1/2$. Since the 
  solution to $Mx=0$ is $\imag/(1+\imag)= -1/2(1+\imag)$, the upper
  half plane maps to the interior of the circle. The corresponding
  complex multiaffine polynomial is
\[
\imag + (1+\imag)x + (1+\imag)y + (2+\imag)xy
\]
The four conditions are
\[
|\Re(M)| = 1\qquad |\Im(M)|=0 \qquad Im(\alpha/\beta) = -1\qquad
\Im(\gamma/\delta)=-1/2. \]
\end{example}

\section{Nearly multiaffine polynomials}
\label{sec:nearly-mult-polyn}

A \emph{nearly multiaffine polynomial} is a polynomial that is degree
one in all but one variable.

\begin{lemma}
  If $yf(x) + z g(x)\in\gsubclose_3$ then $f$ is a constant multiple of $g$.
\end{lemma}
\begin{proof}
  We first remove all common factors of $f$ and $g$ so that it
  suffices to show that $f$ and $g$ are constant. Since
  $yf+zg\in\gsubclose_3$ it follows that $f$ and $g$ interlace. We
  consider the possibilities:

  First assume $f\lessless g$; the case $g\lessless f$ is the same. If
  we let $y=x$ then $(xf) + z g\in\gsubclose_2$. But this is
  impossible, since this implies that $g$ and $xf$ interlace yet
  $deg(xf) = 2+deg(g)$.

  We next assume $f\greateq g$. All roots of $g$ are smaller than the
  largest root $s$ of $f$. If $r>s$ then $(x-r)f+zg\in\gsubclose_2$,
  and hence $(x-r)f\lesslesseq g$. But this implies that $g$ has a
  root in $[s,r]$, which is a contradiction.

  Thus, $f$ and $g$ must be constant, which proves the lemma.
\end{proof}

\begin{cor}
  If $ f(x,\yy)=\sum_1^d y_i f_i(x)\in\gsubclose_{d+1}$ then there are positive
  $\alpha_i$ and $g\in\allpoly$ such that
\[ f(x,\yy) = g(x)\, (\alpha_1y_1 + \cdots \alpha_d y_d) \]
\end{cor}

We next show that there are no interesting matrices $M(x)$ such that $\yy
M \zz^t\in\gsubclose_d$ for vectors $\yy,\zz$ of variables. The only
such $M$ have the form $g(x)\, v^t\cdot w$ where $v,w$ are vectors of
positive constants.

\begin{lemma}
  If $f(x,\yy,\zz)=\sum_1^d y_i\,z_j f_{ij}(x)\in\gsubclose_{2d+1}$ then there are
  $g\in\allpoly$ and positive $\alpha_i,\beta_i$ so that
\[ 
f(x,\yy,\zz) = g(x) \bigl(\sum \alpha_i y_i\bigr)\bigl(\sum \beta_i
z_i\bigr).
\]
\end{lemma}
\begin{proof}
  Consider the matrix $M(x) = (f_{ij}(x))$. For any $\alpha\in\reals$
  we know from Lemma~\ref{lem:multi-bilinear} that $M(\alpha)$ has
  rank $1$, so for all distinct $i,j,k,l$
\[ f_{ij}(\alpha)f_{kl}(\alpha)= f_{il}(\alpha)f_{kj}(\alpha)
\]
and hence $f_{ij}f_{kl}= f_{il}f_{kj}$ as polynomials. Now
\[
f_{11}\cdot M = (f_{11}f_{ij}) =
\begin{pmatrix}
  f_{11}\\ \vdots \\ f_{1d}
\end{pmatrix}
\begin{pmatrix}
  f_{11} & \dots & f_{d1}
\end{pmatrix}
\]
since $f_{11}f_{ij} = f_{1j}f_{j1}$. It follows that
\[
f_{11}\cdot f(x,\yy,\zz) = \bigl( \sum y_if_{1i}\bigr)\bigl( \sum
z_if_{i1}\bigr).
\]
Since each factor is in $\gsubclose_{2d+1}$ we can use the corollary
to write
\[
f_{11}\cdot f(x,\yy,\zz) = g(x) \bigl( \sum \alpha_i y_i\bigr)\bigl( \sum
\beta_iz_i\bigr). 
\]
Clearly $f_{11}$ divides $g$, giving the result.
\end{proof}

We can apply the characterization of complex multiaffine
polynomials to get properties of nearly multiaffine polynomials.

\begin{lemma}
  Suppose $f(x) + g(x)y + h(x)z + k(x) yz\in\gsubclose_3$. If
  $\sigma\in\uhp$, $\delta^2 =
  \smalltwodet{f(\sigma)}{g(\sigma)}{h(\sigma)}{k(\sigma)}$ and 
$\smalltwodet{\alpha}{\beta}{\gamma}{\delta} = 
 \smalltwodet{f(\sigma)/\delta}{g(\sigma)/\delta}{h(\sigma)/\delta}{k(\sigma)/\delta}$
 then
\[ 
\begin{vmatrix}
  \Re(\alpha) & \Re(\beta) \\ \Re(\gamma) & \Re(\delta)
\end{vmatrix}\ge1 \qquad 
\begin{vmatrix}
  \Im(\alpha) & \Im(\beta) \\ \Im(\gamma) & \Im(\delta)
\end{vmatrix}\le0  
\]
\end{lemma}

  \section{The polar derivative and the Schur-\Szego\ map}
  \label{sec:polar-deriv-schur}

\added{1/12/6}

\index{polar derivative}

  The polar derivative was defined for polynomials in $\gsubplus_d$
  \seepage{eqn:polar-pd} but its properties still hold for
  $\gsubclose_d$ as well if we restrict ourselves to a single variable. If
  $f(\xx,y)\in\rupint{d+1}$ then the polar derivative with respect to
  $y$, which we shall denote as $\partial_y^{polar}$, is defined as the
  composition
\[ . \xrightarrow{\text{reverse y}}
 \,.\,
\xrightarrow{\text{differentiate}}
\,.\,
 \xrightarrow{\text{reverse y}} \,.
\]
$\partial_y^{polar}$ defines a map from $\gsubclose_{d+1}$ to itself
since reversal maps $\gsubclose_{d+1}$ to itself. As in one variable
we have the formula
\[
\partial_y^{polar} f = n f - y \partial_y\,f
\]
where $n$ is the degree of $y$ in $f$. The key fact about polar
derivatives is

\begin{lemma}\label{lem:polar-in-pd}
  If $f(\xx,y)\in\gsubclose_{d+1}$ then
  \begin{enumerate}
  \item  $\partial_y $ and  $\partial_y^{polar}$ commute
  \item  $\partial_y f$ and  $\partial_y^{polar} f$ interlace.
\end{enumerate}
\end{lemma}
\begin{proof}
  Commutativity is the calculation
\[
\partial_y \partial_y^{polar} f = n f' - f'-yf'' = (n-1)f'-yf''=
 \partial_y^{polar}\partial_y
\]
  Choose $\alpha\in\reals$ and let $g(\xx,y) = f(\xx,y-\alpha)$. Since
  $g\in\gsubclose_{d+1}$ we know that its polar derivative is in
  $\gsubclose_{d+1}$ and therefore 
$ n g - y \partial_y g\in\gsubclose_{d+1}$. Substituting $y+\alpha$ for
$y$ yields 
\[ nf(\xx,y) - (y+\alpha)f(\xx,y) = \partial_y^{polar}f - \alpha
\partial_y f \in\gsubclose_{d+1}
\]
which proves the lemma.
\end{proof}

We can apply this to derive an important theorem. We call the
transformation the \emph{Schur-\Szego\ map}.
\index{Schur-\Szego map}

\begin{theorem}
  If $f = \sum_0^na_iy^i\in\allpoly(n)$ and $\sum
  g_i(\xx)y^i\in\gsubclose_{d+1}(n)$ then
\[
\sum a_i\, i!\,(n-i)!\, g_i(\xx) \in\gsubclose_d
\]

Equivalently, the transformation $y^i\times f(\xx)y^j\mapsto
\begin{cases}
  0 & i\ne j\\ i!\,(n-i)!\,f(\xx) & i=j
\end{cases}$

 defines a map $\allpoly(n)\times \gsubclose_{d+1}(n)\longrightarrow \gsubclose_d(n)$.
\end{theorem}
\begin{proof}
  By the lemma we see that if $f(x)\in\allpoly$ and $F(x,y)$ is the
  homogenization of $f$ then 
\[ F(\partial_y,\partial_y^{polar}) \colon \gsubclose_{d+1}
\longrightarrow\gsubclose_{d+1}\]
since $\partial_y+\alpha
\partial_y^{polar}:\gsubclose_{d+1}\longrightarrow\gsubclose_{d+1}$
and the regular and polar derivatives commute.  It's easy to see that
\[
\partial_y^i \bigl(\partial_y^{polar}\bigr)^{n-i}y^j =
\begin{cases}
  i!(n-i)! & i=j \\ 0 & i\ne j
\end{cases}
\]
so we have
\begin{multline}
  F(\partial_y,\partial_y^{polar}) \sum_i g_i(\xx)y^i = 
\sum_{i,j} a_i  \partial_y^i
\bigl(\partial_y^{polar}\bigr)^{n-i}y^j g_j(\xx)\,\\
=   \sum a_i\, i!\,(n-i)!\, g_i(\xx) \in\gsubclose_d
\end{multline}
\end{proof}

There are two important consequences that previously were proved using
Grace's theorem. \index{Grace's theorem}

\begin{cor}[Schur-\Szego]
  If $\sum a_i x^i\in\allpoly$ and $\sum b_ix^i\in\allpolypos$  then
\[ \sum a_i b_i\, i!\, (n-i)!\, x^i \in\allpoly. \]
\end{cor}
\begin{proof}
  Since $\sum b_ix^i\in\allpolypos$ we know that $\sum b_i\,
  (-xy)^i\in\gsubclose_2$. Now apply the theorem.
\end{proof}

\begin{cor}[\bbs]
  If $\sum a_i x^i\in\allpoly$ and $g(x)\in\allpoly$ then
\[
\sum (n-i)!\, a_{n-i}\, g^{(i)}(x) \in\allpoly.
\]
\end{cor}
\begin{proof}
  Apply the theorem to 
\[
g(x+y) = \sum g^{(i)}(x) \frac{y^i}{i!} \in\rupint{2}.
\]
\end{proof}

More generally, if $T\colon \allpoly\longrightarrow\rupint{2}$ is a map
\seepage{lem:p-to-p2} then we get a map 

\centerline{\xymatrix{
\allpoly\times\allpoly \ar@{-->}[rr] \ar@{->}[dr]_{id\times T} &&
\allpoly \\
& \allpoly\times\rupint{2} \ar@{->}[ur]_{\text{Schur-\Szego}}
}}

  \section{Determining the direction of interlacing}
  \label{sec:determ-direct-interl}

  If all  linear combinations of $f$ and $g$ are in $\gsubclose_d$
  then either $f\longleftarrow g$ or $g\longleftarrow f$. Determining
  the direction of interlacing is easy if the degrees are different;
  here are some ways of determining the direction when the degrees
  are the same. First of all,

  \begin{lemma}
    If $f\lace g$ in $\rupint{d}$ have the same degree then either
    $f\greateqeq g$ or $g \greateqeq f$.
  \end{lemma}
  \begin{proof}
    The homogeneous part of $\alpha f + g$ is positive for large
    $\alpha$, and negative for very negative $\alpha$.  Thus, there is
    an $\alpha$ for which one of the coefficients of $(\alpha f+g)^H$
    is zero, but this is only possible if $(\alpha f+g)^H=0$. For this
    $\alpha$ define $r = \alpha f +g$. Now $f\lace r$, and since $r$
    has smaller degree than $f$ we have $f\lesslesseq r$. Thus $g =
    \alpha f + r$, from which the conclusion follows.
  \end{proof}

  \begin{lemma}
    If $f,g\in\rupint{d}$, $f\greateqeq g$, and $g \greateqeq f$ then $g$
    is a scalar multiple of $f$.
  \end{lemma}
  \begin{proof}
    $f\greateqeq g$ implies that there are $\alpha>0$ and $r$ such
    that $f\lesslesseq r$, $f = \alpha g + r$ and $r$ has positive
    leading coefficient. But $g\greateqeq f$ implies $f = (1/\alpha) g
    - (1/\alpha)r$ where $-r$ has positive leading coefficient. It
    follows that $r=0$. 
  \end{proof}

  \begin{lemma}
    Suppose that $f \lace g$ in $\rupint{d}$. If either 
\begin{align*}
  \forall \ a_2,\dots,a_d\in\reals & \quad f(x,a_2,\dots,a_d)\longleftarrow
  g(x,a_2,\dots,a_d) \\
\intertext{or}
\exists a_2,\dots,a_d\in\reals  & \quad f(x,a_2,\dots,a_d)\greateq g(x,a_2,\dots,a_d) 
\end{align*}
then $f \longleftarrow g$.
  \end{lemma}
  \begin{proof}
    We may assume that both $f$ and $g$ have degree $n$. If
    $g\greateqeq f$ then we know that 
\[
 g(x,a_2,\dots,a_d) \greateqeq f(x,a_2,\dots,a_d)
\]
The hypothesis gives the other direction, so 
\[
g(x,a_2,\dots,a_d) = \alpha f(x,\alpha_2,\dots,\alpha_d)
\]
Since this holds for all $x$ and $\alpha_i\in\reals$ this implies
that $g = \alpha f$, so $f\longleftarrow g$ still holds.

If we have strict interlacing for one choice of $a_i$ then we can't
also have $g(x,a_2,\dots,a_d)\greateqeq f(x,a_2,\dots,a_d)$. 
  \end{proof}

We only need interlacing in one variable to determine interlacing 
for $\rupint{2}$.

\begin{lemma}
  If $f,g\in\rupint{2}$ and $f(x,a)\greateqeq g(x,a)$ for all
  $a\in\reals$ then $f\greateqeq g$.
\end{lemma}
\begin{proof}
  It suffices to show that $f^H$ and $g^H$ are scalar
  multiples. Geometrically the hypotheses say that the solution curves
  for $g$ interlace the solution curves for $f$. For large positive
  $x$ we therefore have $f^H(x) \greateqeq g^H(x)$ and for small
  negative $x$ we have $f^H(-x) \greateqeq g^H(-x)$ which implies $f^H$ is a
  scalar multiple of $g^H$. 
\end{proof}

Substitution for $\rupint{d}$  preserves the \index{direction of interlacing}
direction of interlacing  but this is not true for $\gsubclose_d$.  For example
\[
\begin{vmatrix}
  x & 1 & 1 \\ 1 & y & 1 \\ 1 & 1 & z
\end{vmatrix}
=
2-x-y + z(xy-1)
\]
so $2-x-y\greateqeq xy-1$. However, substitution of $y=2$ and $y=-2$
gives interlacings in different directions:
\[
\begin{array}{rlcrl}
-x & = (2-x-y)(x,2) & \lesseqeq & (xy-1)(x,2) & = 2x-1\\
-x+4 & = (2-x-y)(x,-2) & \greateqeq & (xy-1)(x,-2) & = -2x-1  
\end{array}
\]



In the following we let $\aaa=(a_1,\dots,a_d)\in\reals^d$ and
$\bbb=(b_1,\dots,b_d)\in(\reals^+)^d$. From Lemma~\ref{lem:pd-xxg} we
know that if $f(\xx)\in\gsubclose_d$ then $f(\aaa+\bbb t)\in\allpoly$.
The converse is also true, see \cite{bbs}.

\begin{lemma}\label{lem:abt-pd}
  If $f\lace g$
  in $\gsubclose_d$ and $f(T)\longleftarrow g(T)$ for all
  $\aaa\in\reals^d$ and $\bbb\in(\reals^+)^d$ where $T=\aaa+\bbb t$
    then $f\longleftarrow g$.
\end{lemma}
\begin{proof}
  We first assume $f,g\in\rupint{d}(n)$. If we take
  $\aaa=(0,a_2,\dots,a_d)$ and $\bbb_\epsilon =
  (1,\epsilon,\dots,\epsilon)$  then $f(\aaa+t\bbb_\epsilon)
  \greateqeq g(\aaa+t\bbb_\epsilon )$. Since $f(\aaa+t\bbb_\epsilon)$
  and $g(\aaa+t\bbb_\epsilon)$ also have degree $n$ it follows that
  $f(\aaa+t\bbb)\greateqeq g(\aaa+t\bbb)$. Thus
\[
\forall a_i\in\reals\qquad f(t,a_2,\dots,a_d) \greateqeq g(t,a_2,\dots,a_d)
\]
so $f\greateqeq g$.

Now assume that $f,g\in\gsubclose_d$ have the same degree. Let
$X=x_1+\cdots+x_d$, $f_\epsilon(\xx)=f(\xx+\epsilon X)$, and
$g_\epsilon(\xx)=g(\xx+\epsilon X)$. We claim
\begin{enumerate}
\item $f_\epsilon,g_\epsilon\in\rupint{d}$
\item $\lim_{\epsilon\rightarrow0} f_\epsilon = f$, $\lim_{\epsilon\rightarrow0} g_\epsilon = g$.
\item $f_\epsilon(T)\longleftarrow g_\epsilon(T)$.
\end{enumerate}
Given 1,2,3 then by the first part we have that
$f_\epsilon(\xx)\longleftarrow g_\epsilon(\xx)$. All terms have the same
degree, so using Lemma~\ref{lem:lace-limit} $f(\xx)\greateqeq g(\xx)$. 

Now 1) follows from Lemma~\ref{lem:pd-xxg}, 2) is clear, and 3) follows from the
first part as follows. If 
\[
S = \aaa+t\bbb + \epsilon\biggl[\sum a_i + t \sum b_i\biggr](1,1,\dots,1)
\]
then $f_\epsilon(T) = f(S)$ and $g_\epsilon(T) = g(S)$
\end{proof}

  Interlacing in $\gsubclose_d$ is the limit of interlacing in $\rupint{d}$.
  \begin{cor}\label{cor:int-limits-d}
    If $f\longleftarrow g$ in $\gsubclose_d$ then there are
    $f_n,g_n$ such that 
    \begin{enumerate}
    \item $f_n,g_n\in\rupint{d}$.
    \item $f_n \longleftarrow g_n$.
    \item $\displaystyle\lim_{n\rightarrow\infty} f_n = f$ and
      $\displaystyle\lim_{n\rightarrow\infty} g_n = g$. 
    \end{enumerate}
  \end{cor}
  \begin{proof}
    First assume $f\lesslesseq g$. 
    Since $f + x_{d+1}\in\gsubclose_{d+1}$ there are
    $F_n\in\rupint{d+1}$ such that $\lim F_n = f+ x_{d+1}g$. If we write
    $F_n = f_n + x_{d+1}g_n + \cdots$ then $f_n,g_n$ satisfy the
    three conclusions.

    If $f\longleftarrow g$ then we can write $g = \alpha f + r$ where
    $f\lesslesseq r$. Now apply the first paragraph to $f$ and $r$.
  \end{proof}

  \begin{cor}
    If $f\greateqeq g$ in $\gsubclose_d$ then
$
\begin{vmatrix}
  f & g \\ \frac{\partial f}{\partial x_i} & \frac{\partial g}{\partial x_i}
\end{vmatrix} \le 0
$  for $i=1,\dots,d$.
  \end{cor}
  \begin{proof}
    We approximate and reduce the problem to the corresponding result
    for $\rupint{d}$. If suffices to assume that $f\lesslesseq g$ and
    $i=1$. Since $f+x_{d+1}g\in\gsubclose_{d+1}$ there are $F_n = f_n+
    x_{d+1}g_n + \cdots \in\rupint{d+1}$ such that $\lim_{n\rightarrow\infty} F_n =
    f+x_{d+1}g$. Now by \index{Taylor's theorem}Taylor's theorem
\[
F_n(x_1+t,x_2,\dots,x_{d+1})=
f_n + x_{d+1}g + t\biggl( \frac{\partial f_n}{\partial x_1} +
x_{d+1}\frac{\partial g_n}{\partial x_1}\biggr) + \cdots
\]
so Corollary~\ref{cor:product-4d} tells us that
\[
\begin{vmatrix}
  f_n & g_n  \\
\frac{\partial f_n}{\partial x_1} & 
\frac{\partial g_n}{\partial x_1}
\end{vmatrix}\le0
\]
and taking limits finishes the proof.
  \end{proof}

\begin{cor}
  If $f\lace g\in\gsubclose_d$ and $  \begin{vmatrix}
    f & g \\ \frac{\partial f}{\partial x_i} & \frac{\partial g}{\partial x_i}
  \end{vmatrix} < 0$  for some $\xx\in\reals^d$ and some $i$ where $1\le i \le d$
  then $f\greateqeq g$.
\end{cor}
\begin{proof}
  We know $f\greateqeq g$ or $g\greateqeq f$. If the latter holds then 
\[
  \begin{vmatrix}
    f & g \\ \frac{\partial f}{\partial x_i} & \frac{\partial g}{\partial x_i}
  \end{vmatrix} \ge 0 
\]
but this contradicts the hypothesis.
\end{proof}

\section{Reversal in $\rupint{d}$}
\label{sec:reversal-pd}

We use the properties of the graph of a polynomial in $\rupint{d}$ to
prove that the reversal is in $\gsubclose_d$. 

  \begin{lemma}\label{lem:reverse-in-pd}
    If $f(\xx)\in\rupint{d}(n)$ then the reversal $f^{rev}$ is in
    $\gsubclose_d(nd)$, where
\begin{equation}\label{eqn:reverse-in-pd}
 f^{rev} (\xx) = (x_1\cdots x_d)^n
f\left(\frac{1}{x_1},\dots,\frac{1}{x_d}\right)
\end{equation}
\end{lemma}
  \begin{proof}
Define
$$ g_\epsilon(\xx) = \left(\prod_{i=1}^d (x_i+\epsilon \sum_{i\ne j}
  x_j)\right)^n \
f\left(\frac{1}{x_1+\epsilon\sum_{1\ne j} x_j},\dots,
\frac{1}{x_d+\epsilon\sum_{d\ne j} x_j}\right)
$$
Since $\lim_{\epsilon\rightarrow0} g_\epsilon(\xx)=f^{rev}(\xx)$,
it suffices to show that $g_\epsilon\in\rupint{d}(nd)$.  Now $g_\epsilon$
clearly satisfies  homogeneity, so we just
have to check substitution. Let $x_1=t,x_2=a_2,\dots,x_d=a_d$, and set
$a= a_2+\cdots a_d$. It suffices to show that the second factor, which
is

$$
f\left(\frac{1}{t+\epsilon a},\frac{1}{\epsilon t+ \epsilon
    a+(1-\epsilon)a_2},\dots,\frac{1}{\epsilon t+\epsilon a +
    (1-\epsilon)a_d}\right),$$
has exactly $nd$ roots.  So, consider
the curve
$$\mathcal{C} = \left(\frac{1}{t+\epsilon a},\frac{1}{\epsilon t+ \epsilon
    a+(1-\epsilon)a_2},\dots,\frac{1}{\epsilon t+\epsilon a +
    (1-\epsilon)a_d}\right)$$

The curve $\mathcal{C}$ has singularities at
$$
-\epsilon a, -a - \frac{1-\epsilon}{\epsilon}a_2,\dots, -a -
\frac{1-\epsilon}{\epsilon}a_d $$
As $t$ goes to $\pm\infty$ the curve
goes to $0$, so there are $d$ components in the closure of the image
of $\mathcal{C}$. A component has two coordinates that are unbounded;
the other coordinates only vary by a bounded amount. The curve goes to
$+\infty$ in one coordinate, and is bounded in the others. In the
other direction, it goes to $-\infty$ in the other coordinate, and is
bounded in the remaining ones. Consequently, the curve is eventually
above the graph in one direction, and below in the other, so each
component meets the graph in $n$ points. Since there are $d$
components of $\mathcal{C}$ we have $nd$ intersection points, and so
substitution is satisfied.
  \end{proof}
 
  Using reversal, an argument entirely similar to
  Lemma~\ref{lem:converse-p2-diff} shows that

\begin{lemma} \label{lem:converse-pd-diff}
  If $f(\partial_{\xx})g$ is in $\gsubclose_d$ for all
  $g\in\gsubclose_d$ then  $f$ is in $\gsubclose_d$. 
\end{lemma}

  \begin{lemma}\label{lem:pd-fxdy}
    Suppose that $f(\xx,y)=\sum f_i(\xx)y^i\in\gsubf_{d+1}$ or
    $\gsubclose_{d+1}$, and the degree of $f_i(\xx)$ is $i$.

 If the homogeneous parts of the $f_i$ alternate in  sign then
    $$
    g(\xx,y,\zz)\in\rupint{d+e+1}(n) \implies
    f(\xx,-\partial_y)\,g(\xx,y,\zz)\in\rupint{d+e+1}(n).$$

  \end{lemma}
  \begin{proof}
    If $e>0$ so that there are some $\zz$ variables, then we replace
    $f(\xx,y)$ by $f(\xx+\epsilon \zz,y)$, which still satisfies the
    degree hypothesis. We now follow the proof of Lemma~\ref{lem:p2-fxdy}
  \end{proof}

  \section{Induced transformations}
\label{sec:top-induced-trans}

\index{induced transformation}
\index{transformation!induced}

\index{transformation!induced}  A linear transformation $T$ on
polynomials in one variable determines  an induced transformation
defined on polynomials in two variables:
$$ T_\ast(x^iy^j) = T(x^i)\, y^j.$$

Conversely, we will see that polynomials in $\gsubclose_2$ can
determine linear transformations on $\allpoly$.  We have seen 
examples of induced transformations \seepage{thm:pd-T} in $\rupint{d}$
where the results are easier.

We first make the important observation that if
$T\colon{}\allpoly\longrightarrow\allpoly$ then $T$ does not necessarily
determine a linear transformation
$T_\ast:\gsubclose_2\longrightarrow\gsubclose_2$. As in
Remark~\ref{rem:no-induced}, we take $T(g) = g(\diffd)x$ and note that
(but see Lemma~\ref{lem:p2-induction-3})
$$ T(xy-1) = y-x\not\in\gsubclose_2.$$
\index{linear transformation!doesn't induce on $\rupint{2}$}

However, if $T\colon{}\allpoly\longrightarrow\allpoly$ then the induced
transformation always satisfies $x$-substitution, as the calculation
below shows:

$$
(T_\ast f)(x,\alpha) = T(f(x,\alpha)) \in\allpoly\text{ since }
f(x,\alpha)\in\allpoly.$$

\begin{definition}
  If $T$ is a linear transformation on polynomials in one variable
  then we say that \emph{$T$ satisfies induction} if the induced
  transformation $T_\ast$ determines a map
  $T_\ast:\rupint{2}\longrightarrow\gsubclose_2$. If
  $T_\ast:\gsubplus_2\longrightarrow\gsubclose_2$ then we say that $T$
  \emph{satisfies induction on $\gsubplus_2$}.

\index{satisfies induction} 
\index{linear transformation!satisfies induction}
\index{satisfies induction!on $\gsubplus_2$} 
\index{linear transformation!satisfies induction on $\gsubplus_2$}
\end{definition}

Note that if $T$ satisfies induction then $T$ maps $\allpoly$ to
itself, since $f\in\allpoly$ implies $T(f) =
T_\ast(f)\in\gsubclose_2$. Also, if $T$ and $S$ satisfy induction,
then so does their composition $ST$.  Here are two assumptions that
guarantee that induction is satisfied.

\begin{lemma}\label{lem:good-induced}
  Suppose that $T\colon{}\allpoly\longrightarrow\allpoly$. If either of the
  these two conditions hold then $T$ satisfies induction.
  \begin{enumerate}
  \item There is an integer $r$ such that if $T(x^i)\ne0$ then the
    degree of $T(x^i)$ is $i+r$, and $T(x^i)$ has positive leading coefficients. 
  \item  $T(g) = f\ast g$ where $f\in\allpolypos$. 
  \end{enumerate}
\end{lemma}
\begin{proof}
  The first part is a slightly more general statement that
  Theorem~\ref{thm:p2p2}; the proof is similar. In this case we have 
  $T_\ast:\rupint{2}\longrightarrow\gsub_2$. 
  
  In the second case, suppose that $h(x,y)\in\rupint{2}(n)$, and
  $f\in\allpolypos(r)$. If $r$ is at least $n$ then the homogeneous
  part of $T(h)$ has all positive terms. However, if $r<n$ then
  $T(h)^H$ has some zero coefficients.  Replace $f$ by
  $f_\epsilon=f(x)(1+\epsilon x)^n$ where $\epsilon$ is positive. The
  corresponding transformation $S(g) = f_\epsilon(x)\ast g$ satisfies
  $S_\ast(h)\in\rupint{2}$ since $f_\epsilon\in\allpolypos(n+r)$. Taking
  limits as $\epsilon$ goes to zero shows that
  $T_\ast(h)\in\gsubclose_2$.

\end{proof}

We can resolve the problem with the transformation $g\mapsto
g(\diffd)f$ by introducing a negative factor

\begin{cor}\label{cor:sat-ind-gdf}
  If $f(x)$ is in $\allpoly$ then the linear transformation 
  $T\colon{}g\mapsto g(-\diffd)f$  satisfies induction. 
\end{cor}
\begin{proof}
  We calculate
$$ \sum T(x^i) \frac{(-y)^i}{i!} = \sum (\diffd^i f)\frac{y^i}{i!} = f(x+y).$$
Since $f(x+y)\in\rupint{2}$, we can apply Lemma~\ref{lem:p2-induction-3}.
\end{proof}

  \begin{lemma}\label{lem:sat-induction}
    If $F(x,y)= \sum f_i(x)(-y)^i$ is in $\gsubclose_2$  then
    \begin{enumerate}
    \item The linear transformation $T\colon{}x^i\mapsto f_i$
      maps $\allpoly\longrightarrow\allpoly$.
    \item $T$ satisfies induction.
    \end{enumerate}

  \end{lemma}
  \begin{proof}
By Lemma~\ref{lem:reverse-p2-1} the reverse of $F$ satisfies
    \begin{align*}
      F_{rev}(x,y) &= \sum f_i(x) y^{n-i} \in\gsubclose_2\\
\text{Choose}\quad
      h(x,y) &= \sum h_i(x)\,y^i\in\rupint{2}\\
      T(h_i(x)) &= \text{coefficient of $y^n$ in } h_i(y)F_{rev}(x,y) \\
      T_\ast(h(x,y)) &= \text{coefficient of $z^n$ in } h(z,y)F_{rev}(x,z) 
    \end{align*}
Since $h(z,y)$ and $F_{rev}(x,z)\in\gsubclose_3$, it follows that
$T_\ast(h)\in\gsubclose_2$. Finally, (2) holds since $T$ satisfies induction.
  \end{proof}

As is often the case, if we  replace multiplication by a
differential operator and add a constraint, we can introduce a factorial.

  \begin{lemma}\label{lem:p2-induction-3}
    Suppose that $f(x,y) = \sum \frac{f_i(x)}{i!}{(-y)^i}$ is in 
    $\gsubclose_2$, and each $f_i$ is a polynomial of degree $i$ with
    positive leading coefficient.
    \begin{enumerate}
    \item The linear transformation $T\colon{}x^i\mapsto f_i$
      maps $\allpoly\longrightarrow\allpoly$.
    \item $T$ satisfies induction.
    \end{enumerate}
      \end{lemma}
  
\begin{proof}
      Choose $g(x,y)=\sum g_i(y)x^i$ in $\rupint{2}$.  By  Lemma~\ref{lem:pd-fxdy}
      we know that $  f(x,-\partial_z)g(y,z) \in\gsubclose_3$.
 The coefficient of $z^0$ is in $\gsubclose_2$ and equals
$$
\sum \frac{f_i(x)}{i!} (\partial_z)^i g_i(y) z^i = \sum  f_i(x) g_i(y) = (T_\ast g)
$$

\end{proof}

We can show a linear transformation preserves all real roots by
applying it to test functions.

\begin{lemma}\label{lem:1-xy-testfct}
If $T$ is a linear transformation on polynomials and
$T_\ast(1-xy)^n\in\gsubclose_2$ for $n=1,2,\dots$ then
$T:\allpoly\longrightarrow\allpoly$. 
\end{lemma}
\begin{proof}
  Let $g= \sum_0^r a_ix^i\,\in\allpoly$ and  calculate
  \begin{align*}
    g(-\diffd_y) T_\ast\biggl(1-\frac{xy}{n}\biggr)^n\biggr|_{y=0} 
&=
\sum_{i=0}^r a_i (-1)^i T_\ast \sum_{i=0}^n
\biggl(\frac{-xy}{n}\biggr)^k \binom{n}{k}\biggl|_{y=0} \\
&=
\sum_{i=0}^r a_i \frac{\falling{n}{i}}{n^i} T(x^i)
  \end{align*}

  Since $g(-\diffd_y)$ maps $\gsubclose_2$ to itself, and evaluation yields a
  polynomial in $\allpoly$, we see the above polynomial is in
  $\allpoly$. Finally, taking the limit of polynomials yields
\[
\lim_{n\rightarrow\infty} \sum_{i=0}^r a_i \frac{\falling{n}{i}}{n^i}
T(x^i)
= \sum a_i T(x^i) = T(g(x))
\]
\end{proof}

\section{Linear Transformations}
\label{sec:application-bilinear}

One of the ways that we use to show that a polynomial has all real
roots is to identify it as a coefficient of a polynomial in
$\gsubclose_2$. We use this idea to show that some linear
transformations preserve $\allpoly$.

\begin{theorem} \label{thm:xfi}
  Suppose $f$ is a polynomial in $\gsubpos_2$ of degree $n$ with
  coefficients given in \eqref{eqn:p3-1}.  The linear transformation
  defined by $T\colon{}x^i \mapsto f_{n-i}$ maps $\allpoly(n)$ to itself.
  Similarly the transformation $x^i\mapsto f_i$ maps $\allpoly$ to itself.
\end{theorem}
\begin{proof}
  Choose a polynomial $g\in\allpoly$. Since $g(y)\in\gsubclosepos_2$ 
  the result of multiplying $g(y)f(x,y)$ is in
  $\gsubclosepos_2$, and hence the coefficient of $y^n$ is in
  $\allpoly$. The coefficient of $y^n$ is exactly $Tg$ since
  \begin{align*}
    g(y)f(x,y) &= \left(\sum a_j y^j\right) \ \left(\sum
      f_j(x)y^j\right) \\
&= \sum_k \left(y^k \sum_{i+j=k} a_j f_{i}(x)\right)
  \end{align*}
If we replace $g$ by its reverse then the transformation is
$x^i\mapsto f_i$.
\end{proof}

\begin{cor}
  If $f\in\allpoly$ then the map $T\colon{}x^i\mapsto \frac{f^{(i)}}{i!}$
  preserves roots and interlacing.
\end{cor}
\begin{proof}
\index{Taylor series}
  Apply the Theorem to the Taylor series
$$ f(x+y) = f(x) + f^\prime(x)y + \cdots $$
\end{proof}

An alternative argument follows from the observation that $T(g) =
g(\diffd{})\expoper{}(f)$.  Or, we could apply $\expoper{}$ to
Theorem~\ref{thm:xfi}.

\begin{lemma}
  If $F(x,y)$ is the homogenization of $f$ then define a map
  $T\colon{}f\times g\mapsto F(-\frac{\partial}{\partial x},y)\,g$. This
    defines a linear transformation  $$
    \allpolypos\times\,\gsubclosepos_2 \longrightarrow \gsubclosepos_2.$$
\end{lemma}
\begin{proof}
  Since $f\in\allpolypos$ we can write $F=\prod(\alpha_ix+y)$ where
  all $\alpha_i$ are positive. From Lemma~\ref{lem:pd-diff-prod-2} we know that 
 $$ \prod_i(\epsilon x + y - \alpha_i \frac{\partial}{\partial x})$$
maps $\gsubpos_2$ to itself. The result follows by letting $\epsilon$
go to zero.
\end{proof}

A similar application of Lemma~\ref{lem:pd-diff-prod-2} shows

\begin{lemma}
  If $f\in\allpoly$ and $g\in\gsubposclose_2$ then $f(y-\diffd_x)g(x,y)
  \in\gsubclosepos_2$.
\end{lemma}
\begin{proof}
  If suffices to show that $(y+\epsilon x -\diffd_x)h(x,y)\in\gsubpos_2$ if
  $h\in\gsubclosepos_2$ and $\epsilon$ is positive. This follows from
  Lemma~\ref{lem:pd-diff-prod-2}.
\end{proof}

There is an interesting identity related to this last result:
\begin{align*}
f(y+\diffd_x) &= \sum_i f^{(i)}(\diffd_x) g \, \frac{y^i}{i!}\\
\intertext{Since this is linear in $f$ only need to verify it for
$f=x^k$}
\sum_i (x^k)^{(i)} (\diffd_x) g\, \frac{y^i}{i!} &=
\sum_i \frac{\falling{k}{i} }{i!} \diffd_x^{k-i}g\, y^i \\
&= \sum_i \binom{k}{i} \diffd_x^{k-i} y^i\, g \\
&= (y+\diffd_x)^k \,g
\end{align*}

The next result constructs a transformation by looking at coefficients.

\begin{theorem} \label{thm:pd-prod-index}
  Choose $f(\xx,\yy)\in\gsubpos_{2d},$ write $f$ in terms of
  coefficients $f = \sum_\diffi a_\diffi(\xx)\yy^\diffi$, and let
  $\diffk$ be an index set. The linear transformation defined by 
$$ \xx^\diffj \mapsto a_{\diffk-\diffj}(\xx)$$
maps $\gsubposclose_d$ to itself.
\end{theorem}
\begin{proof}
  If $g(\xx)\in\gsubpos_d$ then $g(\xx)\in\gsubposclose_{2d}$, and
  hence $g(\yy)f(\xx,\yy)\in\gsubposclose_{2d}$.  If $g(\yy)=
  \sum_\diffj b_\diffj \yy^\diffj$ then we can write
$$ g(\yy)f(\xx,\yy) = \sum_{\diffi,\diffj} b_\diffj\,a_\diffi(\xx)
\yy^{\diffi+\diffj}. $$
The coefficient of $\yy^\diffk$ in this product is in
$\gsubclosepos_d,$ and equals $$\sum_{\diffi+\diffj=\diffk} \,
b_\diffj\,a_\diffi(\xx)$$
 This is exactly $T(g)$ since 
$$ T(g(\xx)) = T\left( \sum_\diffj b_\diffj\, \xx^\diffj\right) =
\sum_\diffj b_\diffj a_{\diffk-\diffj} (\xx)
$$
\end{proof}

\begin{cor}
  If $f\in\gsubpos_d$ and $\diffk$ is an index set then the map
$$\xx^{\diffk-\diffi} \mapsto \frac{1}{\diffi!}
\frac{\partial^\diffi}{\partial \xx^\diffi}\, f(\xx)$$
defines a map $\gsubpos_d\mapsto \gsubpos_d$.
\end{cor}
\begin{proof}
  Since $f(\xx)\in\gsubpos_d$ we know that $f(\xx+\yy)\in\gsubpos_{2d}$.
  The Taylor series of $f(\xx+\yy)$ is
\begin{equation}\label{eqn:multinomial-taylor}
  f(\xx+\yy) = \sum_\sdiffi \frac{\yy^\diffi}{\diffi!}
  \frac{\partial^\diffi}{\partial \xx^\diffi}\, f(\xx).
\end{equation}
\index{Taylor series}
  The corollary now follows from the proof of Theorem~\ref{thm:pd-prod-index}.
\end{proof}

\section{The diamond product}
\label{sec:bilinear-diamond}

\index{diamond product!general}
\index{diamond product!original}
\index{diamond product!original signed}
\index{product!diamond}
Recall (\chapsec{linear}{bilinear-diamond-1}) the diamond product 
$\binom{x}{r}\mydiamond{} \binom{x}{s}  = \binom{x}{r+s}.$
We  show that $f\Diamond g$ has all real roots for certain sets of
$f$ and $g$ by identifying it as a diagonal of a polynomial in
$\gsubclose_3$.  In order to do this we need an identity \cite{wagner}:
$$
f \Diamond g = \sum_i \frac{(x+1)^i x^i}{i!i!}\,
\diffd{}^if\,\diffd{}^ig$$

It turns out that the signed version of this identity also gives a
bilinear product that often preserves roots, so we will begin with it.

\begin{figure}[htbp]
  \begin{center}
    \leavevmode
    \epsfig{file=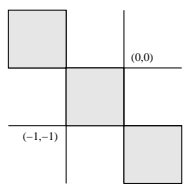,width=1in}
    \caption{ Domains of the signed diamond product}
    \label{fig:diamond-2}
  \end{center}
\end{figure}

\begin{lemma} \label{lem:signed-diamond}
  Define a bilinear map 
$$ S(f,g) = \sum_i (-1)^i\frac{(x+1)^i x^i}{i!i!}\,
\diffd{}^if\,\diffd{}^ig$$
If $(f,g)$  has roots in any one of the four quarter planes below
\footnote{ There are examples where $S(f,g)\not\in\allpoly$ if $f,g$
  have roots in any one of the shaded regions of Figure~\ref{fig:diamond-2}} (see
Figure~\ref{fig:diamond-2}) then $S(f,g)\in\allpoly$.
\begin{enumerate}[(1)]
\item $(-\infty,0)\times(-\infty,-1)$
\item $(-\infty,-1)\times(-\infty,0)$
\item $(0,\infty)\times(-1,\infty)$
\item $(-1,\infty)\times(0,\infty)$
\end{enumerate}
\end{lemma}

\begin{proof}
 
We begin with two consequences of Taylor's theorem
  \begin{align}
f(x+xy) & = \sum_i\frac{x^i\diffd{}^if}{i!}\,y^i \notag \\
g(x+z+xz) &= \sum_j\frac{(x+1)^j\diffd{}^jg}{j!}\,z^j \notag \\
\intertext{  The signed diamond product is the signed diagonal in the product}
    f(x+xy)g(x+z+xz) & = \left(
      \sum_i\frac{x^i\diffd{}^if}{i!}\,y^i\right)\cdot \left(
      \sum_j\frac{(x+1)^j\diffd{}^jg}{j!}\,z^j\right) \label{eqn:signed-diamond-2}
\end{align}

For easy of exposition we define $f(x,y)\in\gsubcloseneg_2$ if and
only if $f(x,-y)\in\gsubclose_2$.  The four domains of the lemma lead
to the four cases: $$
  \begin{matrix}
    (a) & f(\,x(y+1)\,)\in\gsubcloseneg_2 & g(\,(x+1)(y+1)-1\,)\in\gsubcloseneg_2\\
    (b) & g(\,x(y+1)\,)\in\gsubcloseneg_2 & f(\,(x+1)(y+1)-1\,)\in\gsubcloseneg_2\\
    (c) & f(\,x(y+1)\,)\in\gsubposclose_2 & g(\,(x+1)(y+1)-1\,)\in\gsubposclose_2\\
    (d) & g(\,x(y+1)\,)\in\gsubposclose_2 & f(\,(x+1)(y+1)-1\,)\in\gsubposclose_2\\
  \end{matrix}
$$

For instance, suppose that (1) holds. If
$f\in\allpolyint{(-\infty,0)}$ then $f(x+xy)=f(x(y+1))
\in\gsubcloseneg_2$. If $g(x)\in\allpolyint{(-\infty,-1)}$ then
$g(x-1)\in\allpolyint{(-\infty,0)}$, and hence
$g(x+z+xz)=g((x+1)(z+1)-1)\in\gsubcloseneg_2$.  Consequently, case (a)
above holds.  Both of the factors of the left hand side of
\eqref{eqn:signed-diamond-2} are in $\gsubcloseneg_2$, so the product
is as well. The signed diagonal of the right hand side of
\eqref{eqn:signed-diamond-2} is in $\allpoly$ by
Lemma~\ref{lem:diagonal-signed}.

\end{proof}

\begin{prop} \label{prop:diamond-2}
  If $f\times g\in \left(\allpolyint{(-1,\infty)}\times\allpolypos\right) \cup
\left(  \allpolypos\times\allpolyint{(-1,\infty)}\right)$ \\ then $f\Diamond
  g\in\allpoly$. 
\end{prop}

\begin{proof}
The cases are similar to Lemma~\ref{lem:signed-diamond}.
Suppose that $f\in\allpolypos$ and
$g\in\allpolyint{(-1,\infty)}$. Then we
see that $f(-x(y+1))\in\gsubclosepos_2$. As before
$g(x+z+xz)\in\gsubclosepos_2$, so 
$$ h(x,y,z) = f(-x(y+1))\, g(x+z+xz) \in\gsubclosepos_2.$$
The signed diagonal of $h(x,y,z)$ is the diamond product of $f$ and $g$.
\end{proof}

We now consider a different map
$m:\allpoly\times\allpoly\longrightarrow\allpoly$.

\index{diamond product!Hermite}
\index{Hermite polynomials!diamond product}
\begin{prop} \label{prop:diamond-3}
  Suppose $m(f,g)=f(\diffd)g$, and $T(x^n) = H_n$. Then
  $f\mydiamondA{T}{m}g:\allpoly\times\allpoly\longrightarrow\allpoly$. 
\end{prop}
\begin{proof}
  We will show that $f\mydiamondA{T}{m}g = (Tf)(2\diffd)g(x)$ which
  implies the conclusion using Corollary~\ref{prop:fofd} and Corollary~\ref{cor:hermite}. First
  of all, from $(H_n)^\prime = 2nH_{n-1}$ it follows that
  \begin{align*}
    \diffd^k\, H_n &= 2^k\falling{n}{k} H_{n-k} \\
\intertext{and consequently}
T^{-1}(\diffd^k\,H_n) &= 2^k\falling{n}{k} x^{n-k} = 2^k \diffd^k x^n\\
\intertext{By linearity in $H_n$ and $x^n$}
T^{-1}(\diffd^k Tg) &= (2\diffd)^k g \\
\intertext{By linearity in $\diffd^k$ it follows that for any
  polynomial $h$}
T^{-1}(h(\diffd)Tg) &= h(2\diffd)g\\
\intertext{and hence choosing $h=Tf$ yields}
T^{-1}(Tf(\diffd)Tg) &= (Tf)(2\diffd)g.
  \end{align*}
\end{proof}

If we use the difference operator $\Delta(f)=f(x+1)-f(x)$ instead of
the derivative, then we have a similar result, except that the Hermite
polynomials are replaced with the Charlier polynomials.
\index{Charlier polynomials!diamond product}
\index{diamond product!Charlier polynomials}

\begin{prop}
  Suppose $m(f,g)=f(\Delta)g$, and $T(x^n) = C^\alpha_n$. Then
  $f\mydiamondA{T}{m}g:\allpolyalt\times\allpoly\longrightarrow\allpoly$. 
\end{prop}
\begin{proof}
  We will show that $f\mydiamondA{T}{m}g = (Tf)(\Delta)g(x)$ which
  implies the conclusion using Corollary~\ref{prop:fofd} and Corollary~\ref{cor:charlier-x}. Since
 $\Delta(C^\alpha_n) = nC^\alpha_{n-1}$ it follows that
  \begin{align*}
    \Delta^k\, C^\alpha_n &= \falling{n}{k} C^\alpha_{n-k} \\
\intertext{and consequently}
T^{-1}(\Delta^k\,C^\alpha_n) &= \falling{n}{k} x^{n-k} =  \diffd^k x^n\\
\intertext{By linearity}
T^{-1}(h(\Delta)Tg) &= h(\diffd)g\\
\intertext{and hence choosing $h=Tf$ yields}
T^{-1}Tf(\Delta)Tg &= (Tf)(\diffd)g.
  \end{align*}
\end{proof}

\index{falling factorial!diamond product}%
\index{rising factorial!diamond product}%
A diamond-type product associated with falling and rising factorials
has elementary proofs of its properties. (But see
Question~\ref{ques:diamond-falling}.)

\begin{lemma} \label{lem:diamond-falling}
  If $T\colon{}x^n\mapsto  \falling{x}{n}$ and $S\colon{}x^n\mapsto\rising{x}{n}$ then \\
  $\mydiamond{TS}:\allpolyalt\times\allpolyalt\longrightarrow\allpolyalt$
  where $\mydiamond{TS}(f\times g)= S^{-1}(Tf\times Tg)$.

\end{lemma}
\begin{proof}
  The proof follows from the diagram

\centerline{
\xymatrix{
\allpolyalt\times\allpolyalt \ar@{->}[rr]^{T\times T} \ar@{.>}[d]_{\mydiamond{TS}} && 
\allpolyalt\times\allpolyalt  \ar@{->}[d]^{multiplication}\\
\allpolyalt \ar@{<-}[rr]^{S^{-1}} && \allpolyalt
}}
\end{proof}

\section{Substituting into quadratic forms}
\label{sec:nsd-sub}

\index{subdefinite matrix!negative}
\index{matrix!negative subdefinite}
\index{negative subdefinite matrix}%

We have seen that if $f\in\allpolypos$ then
$f(-xy)\in\pm\gsubposclose_2$. Now $xy$ is the \index{quadratic
  form}quadratic form corresponding to the matrix
$\smalltwo{0}{1}{1}{0}$ which is a limit of negative subdefinite
matrices. Here is a generalization:

\begin{lemma}
  Suppose that the $d$ by $d$ matrix $Q$ is the limit of negative
  subdefinite matrices. If $f\in\allpolypos$ then $f(-\xx
  Q\xx^t)\in\pm\gsubclose_d$.
\end{lemma}
\begin{proof}
  If we write $f=\prod(x+c_i)$ where all $c_i$ are positive, then each
  factor $-\xx Q \xx^t + c_i$ of $f(-\xx Q\xx^t)$ is in
  $\pm\gsubclose_d$, so their product is in it as well.
\end{proof}

For instance, if we take $Q=
\left(\begin{smallmatrix}
  0&1&1 \\ 1&0&1 \\1&1&0
\end{smallmatrix}\right)
$ with associated \index{quadratic form}quadratic form $xy+xy+zy$ then $Q$ is a limit of
negative subdefinite matrices since $Q=J_3-I_3$. Consequently, if
$f(x)$ has all negative roots then $f(-(xy+xy+yz))\in\pm\gsubclose_3$.


We have seen that if $f(x)\in\allpolypos$ then $f(-x^2)\in\allpoly$.
Here's a different generalization to $\gsubpos_d$.

\begin{lemma} \label{lem:sub-quad-form}
  Suppose that $Q$ is an $e$ by $e$ negative subdefinite matrix.
  If $f(\xx)\in\gsubplus_d(n)$ then $ f(-\xx Q\xx^t,\dots,-\xx
  Q\xx^t)\in\pm\gsubpos_e(2n)$. 
\end{lemma}
\begin{proof}
  If $f = \sum a_\sdiffi\,\xx^\sdiffi$ then
  $$g(\xx) = f(-\xx Q\xx^t,\dots,-\xx Q\xx^t) = \sum a_\sdiffi(-\xx
  Q\xx^t)^{|\sdiffi|}.$$
  Since $f\in\gsubplus_d$ we know that
  $f(x,\dots,x)$ has all negative roots. For every $c0$ the
  polynomial $-\xx Q\xx^t +c^2 $ satisfies substitution. Consequently,
  for every root of $f(x,\dots,x)=0$ we substitute for all but one
  variable, and find two roots. This accounts for all $2n$ roots of
  $g$.  Since the degree and positivity conditions are clearly
  satisfied for $(-1)^ng$, we see $(-1)^n g\in\gsubpos_e(2n)$.
\end{proof}

  The only bilinear forms in $\gsubclose_{2d}$ are products.
\index{bilinear forms in $\gsubclose_{2d}$}
  \begin{lemma}
    If $f = \sum_{i,j=1}^d\, a_{ij}\,x_i\,y_j\in\gsubclose_{2d}$ then
    there are non-negative $b_i,c_i$ so that
\[
  \sum_{i,j=1}^d\, a_{ij}\,x_i\,y_j = \pm(b_1\,x_1 + \cdots + b_d\,x_d)
(c_1\,y_1 + \cdots + c_d\,y_d)
\]
  \end{lemma}
  \begin{proof}
    We will show that all two by determinants
    $\smalltwodet{a_{ik}}{a_{il}}{a_{jk}}{a_{jl}}$ are zero. This
    implies that all rows are multiples of one another, and the
    conclusion follows.  It suffices to show that
    $\smalltwodet{a_{11}}{a_{12}}{a_{21}}{a_{22}}=0$. If we set
    $x_j=y_j=0$ for $j>2$ then
\[
a_{11}\,x_1\,y_1 + a_{12}\,x_1\,y_2 + a_{21}\,x_2\,y_1 + a_{22}\,x_2\,y_2
\in\gsubclose_4.
\]
If all of $a_{11},a_{12},a_{21},a_{22}$ are zero then the determinant
is zero, so we may assume that $a_{11}>0$. By
Lemma~\ref{lem:homog-in-p2close} we see that $a_{12},a_{21},a_{22}$ are
non-negative. Substituting $x_2=y_2=1$ and applying
Corollary~\ref{cor:sub-xyc} shows that $a_{11}a_{22}-a_{12}a_{21}\le0$. If
$a_{12}$ is non-zero then substituting $x_2=y_1=1$ and using the lemma
we conclude that $a_{21}a_{12}-a_{11}a_{22}\le0$. We get the same
inequality if $a_{21}$ is non-zero. Consequently, the
determinant is zero.
  \end{proof}

\section{Simultaneous negative pairs of matrices}
\label{sec:simultaneous}

In Lemma~\ref{lem:sub-quad-form} we substituted the \emph{same} quadratic from
into a polynomial.  When can we substitute two different quadratic
forms into a polynomial in $\rupint{2}$ and have the resulting polynomial
be in $\rupint{d}$? We need a condition on pairs of quadratic forms.

\index{simultaneous negative pair!of polynomials}
\index{simultaneous negative pair!of matrices}

\begin{definition}
  If $p(x),q(x)$ are polynomials with positive leading coefficients,
  then we say that they are a \emph{simultaneous negative pair} if
  there is some value $x$ such that $p(x)\le 0$ and $q(x)\le 0$. If
  $Q_1$ and $Q_2$ are two matrices, we say that they are a
  \emph{simultaneous negative pair of matrices} if for all vectors of
  constants $\aaa$ the polynomials $\xx_1^\aaa Q_1 \xx_1^{\aaa\,t} $
  and $\xx_1^\aaa Q_2 \xx_1^{\aaa\,t}$ are a simultaneous negative
  pair.

\end{definition}

The first question to address is when does $\{Q,Q\}$ form a
simultaneous negative pair?

\begin{lemma}
  Suppose that $Q$ has all positive entries.   $\{Q,Q\}$ is a
  simultaneous negative pair iff $Q$ is \nsd.
\end{lemma}
\begin{proof}
  Assume that $\{Q,Q\}$ is a simultaneous negative pair.  We show that
  $\xx Q\xx'-c^2\in\rupint{d}$, which implies that $Q$ is \nsd. Since
  $\xx Q\xx'-c^2$ has positive homogeneous part, we need to verify
  substitution.  Choose a vector $\aaa$ and consider $\xx^\aaa Q
  {\xx^\aaa} '$.  For large values of $x$ this goes to infinity. Since
  there is a value of $x$ making $\xx^\aaa Q {\xx^\aaa} '$ non-positive, it
  follows that for fixed $\aaa$, $\xx^\aaa Q {\xx^\aaa} '$ takes on all
  positive values. Consequently $\xx^\aaa Q {\xx^\aaa} '-c^2$ has two
  real zeros, and so substitution is satisfied. 

  The converse is similar.
\end{proof}

  \begin{prop}
    Suppose that $f\in\gsubplus_2(n)$. If $Q_1$ and $Q_2$ form a
    simultaneous negative pair of matrices then $g(\xx) = f(-\xx
    Q_1\xx^t,-\xx Q_2\xx^t)$ is in $\pm\rupint{d}(2n)$.
  \end{prop}

  \begin{proof} 

    Choose a vector of constants $\aaa$ and consider the parametrized
    curve in the plane given by
$$ \mathcal{C} = \left\{ (-\xx_1^\aaa Q_1 \xx_1^{\aaa\,t},-\xx_1^\aaa Q_2
  \xx_1^{\aaa\,t})\mid x_1\in\reals\right\}.$$
Since $Q_i$ has all positive coefficients, the limit of $-\xx_1^\aaa Q_i
\xx_1^{\aaa\,t}$ as $x_1\longrightarrow\infty$ is $-\infty$. Thus, the
curve $\mathcal{C}$ is eventually in the lower left quadrant. Since
$Q_1,Q_2$ are a simultaneous negative pair of matrices we know that 
 $\mathcal{C}$ meets the upper right quadrant. It then
follows that $\mathcal{C}$ meets each of the solution curves of $f$,
since the graph of $f$ meets the $x$ axis to the left of the origin,
\index{solution curves} 
and the $y$ axis below the origin. Each solution curve  yields two
intersections, so we find all $2n$ solutions. As before, $\pm g(\xx)$
satisfies the positivity condition, so $\pm
g\in\rupint{d}(2n)$.
\end{proof}

\begin{example}
  What are conditions on a pair of matrices that make them a
  simultaneous negative pair?  Suppose $Q_1 = \smalltwo{c}{v}{v^t}{C}$
  and $Q_2 = \smalltwo{b}{w}{w^t}{B}$ are  negative
  subdefinite $d$ by $d$ matrices.

  Choose a vector $\aaa$ of constants.  We will find the value $V_1$
  of $\xx_1^\aaa Q_2 \xx_1^{\aaa\,t}$ evaluated at the smallest
  solution to $\xx_1^\aaa Q_1 \xx_1^{\aaa\,t}=0$, and the value $V_2$
  of $\xx_1^\aaa Q_1 \xx_1^{\aaa\,t}$ evaluated at the smallest
  solution to $\xx_1^\aaa Q_2 \xx_1^{\aaa\,t}=0$. If for every $\aaa$
  either $V_1$ or $V_2$ is negative then $Q_1,Q_2$  form a
  simultaneous negative pair. This will be satisfied if we can show that
  $V_1 + V_2$ or $V_1V_2$ is never positive.

Without loss of generality we may scale by
  positive values, and thus assume that $b=c=1$.

\begin{align*}
  \xx_1^\aaa Q_1 \xx_1^{\aaa\,t} & =
  (x_1,\aaa)\smalltwo{1}{v}{v^t}{C}\smalltwobyone{x_1}{\aaa^t} \\
  &= x_1^2 + 2\aaa v^t x_1 + \aaa C\aaa^t\\
  \xx_1^\aaa Q_2 \xx_1^{\aaa\,t} & =
  (x_1,\aaa)\smalltwo{1}{w}{v^t}{B}\smalltwobyone{x_1}{\aaa^t} \\
  &= x_1^2 + 2\aaa w^t x_1 + \aaa B\aaa^t\\
  \intertext{To find $V_1$, we solve for $x_1$ and take the smaller
    root}
  x_1 &= -\aaa v^t - \sqrt{(\aaa v^t)^2 - \aaa C \aaa^t}\\
  \intertext{Since $Q_1$ is negative subdefinite we know $x_1$ is
    real. We then substitute into $\xx_1^\aaa Q_2 \xx_1^{\aaa\,t}$. We
    then do the same for $V_2$. Simplification occurs, and } V_1+V_2
  &= 2(\aaa v^t - \aaa w^t)\left( (\aaa v^t - \aaa w^t) + \sqrt{\aaa
      (v v^t - C)\aaa^t} - \sqrt{\aaa(w w^t - B)\aaa^t}\right)
\end{align*}

Notice that if $v=w$ then $V_1+V_2=0$, so $Q_1,Q_2$ form a
simultaneous negative pair in this case. Reconsidering our earlier
normalization, this shows that if two negative subdefinite matrices
have an equal row then they form a simultaneous negative pair.

\end{example}

\begin{example} \label{ex:xy-and-yz}
In this example we will show that $xy$ and $yz$ are limits of
simultaneous negative pairs. 
Let $Q_1$ and $Q_2$ be the quadratic forms corresponding to the
matrices $M_1$ and $M_2$:
\begin{align*}
  M_1 &= \begin{pmatrix} e & 1 & e \cr 1 & e & e^2 \cr e & e^2 & e^4 \end{pmatrix}\\
  M_2 &= \begin{pmatrix} e^4 & e^2 & e \cr e^2 & e & 1 \cr e & 1 & e
  \end{pmatrix}
\end{align*}
We will show that $M_1$ and $M_2$ are a simultaneous negative pair for
$0<e<1$. When $e=0$ the quadratic forms are $2xy$ and $2yz$. The
determinant of $M_1$ and of $M_2$ is $(e-1)^2e^3(1+e)$, which is
positive for $0<e<1$. Since the eigenvalues for $e=.1$ are $(1.1,
-0.9, -0.0008)$, we see that $M_i$ is negative subdefinite
for $0<e<1$. Notice that $Q_1(x,y,z) = Q_2(z,y,x)$.

If we solve for
$y$ in $Q_1$, and substitute into $Q_2$ we get

\begin{align*}
 V_1(x,z) &= \frac{ (1 - e) ( x - z)}{e} \times \\
&  \bigg[ \left( 2 +
      2\,e - e^2 - e^3 - e^4 \right) \,x +  
      \left( e^2 + e^3 - e^4 \right) \,z\, + \\ & 2\,\left( 1 + e \right)
      \,{\sqrt{(1- e^2)x^2 + (e^4 - e^5)z^2}} \bigg] 
  \end{align*}
  
  We next determine the sign of $V_1$ in the $xz$ plane. First, note
  that $V_1(x,z)=0$ has solution $x=z$. The remaining factor has the
  form $ ax + bz + (cx^2+dz^2)^{1/2}=0$, where $a,b,c,d$ are positive.
  The solution to such an equation is found to be of the form
  $x=\alpha z$ and $x=\beta z$, . In our case, it is found that
  $\alpha$ and $\beta$ are positive, $\alpha>\beta$, and
  $\alpha\beta=1$. The solution is only valid for negative $z$. Thus,
  the sign of $V_1(x,z)$ is as follows. Consider the four rays
  in clockwise order around the origin:
\begin{enumerate}[A:]
\item  $(0,0) - (\infty,\infty)$
\item  $(0,0) - (-\infty,-\alpha\infty)$
\item  $(0,0) - (-\infty,-\infty)$
\item  $(0,0) - (-\infty,-\beta\infty)$
\end{enumerate}

Then $V_1$ is positive in the sectors AB and CD, and negative in the
sectors BC and DA. By symmetry, $V_2(x,z) = V_1(z,x)$, and since
$\alpha\beta=1$ the roots of the complicated factor yield lines with
the same slope. Thus $V_2(x,z)$ is negative in the sectors AB and CD,
and positive in the sectors BC and DA. Consequently, $V_1V_2$ is never
positive.

We conclude that for every value of $x,z$ we can find a value of $y$
for which both quadratic forms are non-positive, so $M_1$ and $M_2$
form a negative simultaneous pair.

\end{example}

\begin{cor} \label{cor:-xy-yz}
If $f(x,y)\in\gsubplus_2(n)$ then $f(-xz,-yz)\in\pm\gsubclose_3(2n)$.
\end{cor}

  \begin{cor}
    Suppose that the linear transformation  $x^i\mapsto \alpha_i\,
    x^i$ maps $\allpolypos$ to itself.
    The linear transformation
$$ x^iy^j\mapsto \alpha_{i+j} x^iy^j$$
maps $\gsubplus_2\longrightarrow\rupint{2}$.
  \end{cor}

\begin{proof}
  If $f(x,y) = \sum a_{ij}x^iy^j$ then by Corollary~\ref{cor:-xy-yz} we know 
\begin{align*}
\sum  a_{ij} (-1)^{i+j}x^iy^j z^{i+j} &\in\pm\gsubclose_3\\
\intertext{ Applying the induced  transformation $z^k\mapsto \alpha_k\,z^k$
  yields that}
 \sum a_{ij} (-1)^{i+j}x^iy^j \alpha_{i+j}z^{i+j} &\in\pm\gsubclose_3\\
\intertext{Substituting $z=-1$ shows that}
\sum a_{ij}  x^iy^j \alpha_{i+j} &\in\rupint{2}
\end{align*}
\end{proof}

\begin{cor} \label{cor:ij-pd}
  The linear transformation $x^iy^j\mapsto \dfrac{x^iy^j}{(i+j)!}$
  maps $\gsubplus_2$ to itself.
\end{cor}

Denote by $T$ the transformation of Corollary~\ref{cor:ij-pd}.
Here are some consequences:
\begin{enumerate}
\item If we take $f(x,y) = (x+y+1)^n$ then applying $T$ yields
$$ \sum_{0\le i,j\le n} \frac{n!}{i!j!(i+j)!} x^i y^j \in\gsubplus_2$$
\item If we take $g\in\allpolypos$ then $g(-xy)\in\pm\gsubplus_2$. If
  $g(-x) = \sum b_i x^i$ then
$$ \sum b_i \frac{x^i y^i}{(2i)!} \in\gsubplus_2$$
\end{enumerate}

\section{The interior of $\gsubpos_d$}
\label{sec:sub-interior}

We have seen (\chapsec{polynomials}{topology}) that the
polynomials in one variable with all distinct roots can be
characterized as the interior of $\allpoly$. 
We can extend these ideas to $\gsubpos_d$.

\index{\ Pgsubd@$\gsubpos_d$!interior} 
\index{\ zzlessless@$\lessless$!in  $\gsubpos_d$}%
 The set of all polynomials of degree at most $n$ in
$\gsubpos_d$ is a vector space $V$ of finite dimension, and so has the
standard topology on it. A polynomial $f$ of degree $n$ is in the
interior of $\gsubpos_d$, written $f\in int\,\gsubpos_d$, if there is an
open neighborhood $\mathcal{O}$ of $f$ in $V$ such that all
polynomials in $\mathcal{O}$ are also in $\gsubpos_d$. The boundary of
$\gsubpos_d$ is defined to be $\gsubpos_d \setminus int\,\gsubpos_d$.

We say that $f \lessless g$ if there are open neighborhoods of $f,g$
such that if $\tilde{f}$ is in the open neighborhood of $f$, and
$\tilde{g}$ is in the open neighborhood of $g$ then $\tilde{f}
\lesslesseq \tilde{g}$. Both $f$ and $g$ are necessarily in the
interior of $\gsubpos_d$.

Surprisingly,  products of linear factors are on the boundary.
\begin{lemma}
  If $n>1$ and $f = \displaystyle\prod_{k=1}^n (b_k + x_1 + a_{2k}x_2
  + \cdots + a_{dk}x_d)$ is in $\gsubpos_d$ then $f$ is in the boundary
  of $\gsubpos_d$.
\end{lemma}

\begin{proof}
  We first make a small perturbation so that all the coefficients
  $a_{ij}$ are distinct. Consider the perturbation
  $$
  h = (\epsilon+(b_1+x_1+\dots)(b_2+x_1+\dots))\ 
  \displaystyle\prod_{k=3}^n (b_k + x_1 + a_{2k}x_2 + \cdots +
  a_{dk}x_d)
$$
If we set all but one of $x_2,\dots,x_d$ to zero, then
Example~\ref{ex:bad-polly} shows that we can choose $\epsilon$ such that 
the first factor does not satisfy substitution.
\end{proof}

There is an effective way of determining if a polynomial is in the
interior of $\gsubpos_2$. It is only practical for polynomials of
small degree. In order to do this, we must introduce the
resultant.  \index{resultant}

\begin{definition}
  If 
  \begin{align*}
    f(x) &= a_0 + a_1x + \cdots + a_nx^n \\
    g(x) &= b_0 + b_1x + \cdots + b_mx^m
  \end{align*}
  then the resultant of $f$ and $g$ is defined to be the determinant
  of the $n+m$ by $n+m$ matrix 
\index{determinants!and the resultant}
$$
\begin{vmatrix}
  a_0 & a_1 & \hdots &    a_n &     &    & \\
      & a_0 & a_1    & \hdots & a_n &    & \\
      & \hdotsfor{4} & \\
      &              &    & a_0  & \hdots    & a_n\\
  b_0 & b_1 & \hdots &    b_m &     &    & \\
      & b_0 & b_1    & \hdots & b_m &    & \\
      & \hdotsfor{4} & \\
      &              &  &   b_0 & \hdots   & b_m
\end{vmatrix}
$$
The resultant $R_f$ of $f(x,y)$ and $\frac{\partial f}{\partial
  y}(x,y)$ is a polynomial in one variable.  The resultant has the
property that $R_f(a)=0$ if and only if $f(a,y)$ has a double root.
\index{resultant}
\end{definition}

We can use the resultant to show that a polynomial satisfies
substitution.

\begin{lemma}
  If $f(x,y)$ satisfies the following two conditions then $f$ satisfies
  y-substitution:
  \begin{itemize}
  \item $f(\alpha,y)$ is in $\allpoly$ for some $\alpha$.
  \item $R_f$ has no real roots.
  \end{itemize}
\end{lemma}
\begin{proof}
  If there is an $a$ for which $f(a,y)$ does not have all real roots,
  then as $t$ goes from $\alpha$ to $a$ there is a largest value of $t$
  for which $f(t,y)$ has all real roots. For this value we must have
  that $f(t,y)$ has a double root, but this is not possible since the
  resultant $R_f$ has no real roots.
\end{proof}

We can extend this lemma to the case where $R_f$ has some real roots. 
The proof is similar.

\begin{lemma}\label{lem:show-in-p2}
  Suppose $f(x,y)$ has resultant $R_f$, and that the distinct roots of $R_f$
  are $r_1<r_2 < \cdots < r_s$. If $f$ satisfies the condition below
  then $f$ satisfies y-substitution.
  \begin{itemize}
  \item Set $r_0=-\infty$ and $r_{s+1}=\infty$. For every $0\le i \le
    s$ there is an $\alpha_i\in(r_i,r_{i+1})$ such that
    $f(\alpha_i,y)\in\allpoly$.
  \end{itemize}
\end{lemma}

\begin{lemma}
  If $f\in\gsubpos_2$  has the property that the
  resultant $Res(f,f_y)$ has only simple roots then $f$ is in the
  interior of $\gsubpos_2$.
\index{resultant}
\end{lemma}

\begin{proof}
  The resultant is a continuous function of the coefficients, so if
  $g$ is close to $f$ then $R_g$ is close to $R_f$. Consequently we
  can find an open neighborhood $\mathcal{O}$ of $f$ so that all
  polynomials in $\mathcal{O}$ have a resultant with simple roots, and
  exactly as many roots as $R_f$.  Moreover, we may choose
  $\mathcal{O}$ small enough that each interval satisfies the
  condition of the previous lemma. We now apply the lemma to conclude
  that every function in $\mathcal{O}$ satisfies y-substitution. We can
  also choose $\mathcal{O}$ small enough to insure that all its members
  satisfy the homogeneity conditions, so that
  $\mathcal{O}\subset\gsubpos_2$.
\end{proof}

\begin{example}
Of course, the resultant must fail to show that linear products are in the 
interior, since they are on the boundary. The resultant of 
\[
\left( x + {a_1} + y\,{b_1} \right) \,\left( x + {a_2} + y\,{b_2}
\right) \,\left( x + {a_3} + y\,{b_3} \right) 
\]
 is 
$$
- {b_1}\,{b_2}\,\,{b_3}\,{\left( {a_1} - {a_2} + y\,{b_1} - y\,{b_2} \right) }^2\,
    {\left( {a_1} - {a_3} + y\,{b_1} - y\,{b_3} \right) }^2\,{\left(
        {a_2} - {a_3} + y\,{b_2} - y\,{b_3} \right) }^2  
$$
and in general 
$$ Resultant\left(\prod_{i=1}^n x+a_i + b_iy \right) =
b_1\cdots b_n \ \prod_{i<j} \left(a_j-a_i + y(b_j-b_i)\right)^2
$$
There are many values of $y$ for which the resultant is $0$, and they
are all double roots, so the previous lemma doesn't apply.

\end{example}

\section{The cone of interlacing}
\label{sec:cone}

\index{quantitative sign interlacing}

Lemma~\ref{lem:sign-quant} in \chap{polynomials} explicitly describes all the
  polynomials that interlace a given polynomial of one variable. In
  this section we study the collection of all the polynomials that
  interlace a given polynomial in $\gsubclose_d$. This set forms a cone,
  and in some cases we are able to determine its dimension.

\index{interlacing!cone}
\index{cone of interlacing}

\begin{definition}
  Suppose that $f\in\gsubclosepos_d$. The
  \emph{interlacing cone of $f$} is defined to be
$$ \mycone{f} =\left\{ g\in\gsubclosepos_d\,\vert\,  f\lesslesseq g
\right\}.$$
\end{definition}

We should first note that if $g_1,g_2\in\mycone{f}$ and $a,b>0$ then
$g=ag_1+bg_2$ satisfies $f\lesslesseq g$ and $g\in \gsubclose_d.$, so
$\mycone{f}$ is closed under positive linear combinations. Thus,
$\mycone{f}$ is indeed a cone. In some cases we can explicitly
determine the cone of interlacing. See Corollary~\ref{cor:cone-1}.

The dimension of a cone is the least number of elements whose
non-negative  linear combinations span the cone. From Lemma~\ref{lem:cone-1} we
see that $\text{dim }\mycone{x^ny^m} = 1$.  

The first lemma is an  analog of the  one variable fact that all
polynomials that interlace $x^n$ are of the form $cx^m$.

\begin{lemma} \label{lem:cone-1}
Suppose that $f\in\gsubclose_2$ has the property that $f$ and $x^ny^m$ 
interlace. If $f$ has $x$-degree at most $n$, and
$y$-degree at most $m$ then there are $a,b,c,d$ such that
$$ f = x^{n-1}y^{m-1}(a+bx+cy+dxy).$$
\end{lemma}

\begin{proof}
For any $\alpha\in\reals$ we substitute to find that $x^n\alpha^m$ and 
$f(x,\alpha)$ interlace. Since the $x$-degree of $f$ is at most $n$ we 
know that there are constants $r(\alpha)$ and $s(\alpha)$  depending
on $\alpha$ such that 
$$
f(x,\alpha) = r(\alpha)x^n + s(\alpha)x^{n-1}.$$ However, we do not
yet know that $r(\alpha)$ and $s(\alpha)$ are polynomials in $\alpha$.
If we solve $f(1,\alpha)=f(2,\alpha)=0$ for $r(\alpha)$ and
$s(\alpha)$ then we see that $r(\alpha)$ and $s(\alpha)$ are rational
functions of $\alpha$. Thus, $f(x,y)=x^{n-1}(r(y)+s(y)x)$ and so we
can conclude that $x^{n-1}$ divides $f$. Similarly, $y^{m-1}$ divides
$f$.  Since $x^{n-1}y^{m-1}$ divides $f$, the conclusion follows from
the degree assumption.
\end{proof}

\begin{cor} \label{cor:cone-1}
$$ \mycone{x^ny^m} = \left\{ a x^{n-1}y^{m-1} \text{where $a>0$}\right\}.$$
\end{cor}

\begin{proof}
  The $x$ and $y$ degrees of a polynomial $g$ such that $x^ny^m\lesslesseq
  g$ must be $n-1$ and $m-1$ respectively. Apply Lemma~\ref{lem:cone-1}.
\end{proof}

The next result is obvious in $\allpoly$. It's worth
noting  however that the positivity condition is not true
for polynomials in one variable that do not have all real roots.\footnote%
{  The product of two polynomials with positive and negative
  coefficients can have all positive coefficients:
$ (2x^2+3x+2)(x^2-x+2) = 2x^4 + x^3 + 3x^2 + 4x+4.$
}

\begin{theorem} \label{thm:p2-unique-prod}
  If $fg\in\gsubpos_2$ then either $f,g\in\gsubpos_2$ or $-f,-g\in\gsubpos_2$.
\end{theorem}
\begin{proof}
  From Lemma~\ref{lem:sub-factorization}  we know
  that $f$ and $g$ satisfy substitution. It
  remains to show that $\pm f^H$ has all positive coefficients. We know
  that $(fg)^H = f^H\,g^H,$ and since $fg\in\gsubpos_2$ we can factor 
  $$ f^H\,g^H = (fg)^H = \prod_i (a_ix + b_iy)$$
  where all $a_i,b_i$ are positive. It follows that $f^H$ and $g^H$ are
  products of some of the terms $a_ix+b_iy,$  and the theorem is proved.
\end{proof}

The following lemma is harder than it seems. See the short proof in
\cite{carroll}.
\begin{lemma} \label{lem:is-a-poly}
  If $f(x,y)$ is a continuous function defined for all $x,y$ with the
  property that for any $\alpha\in\reals$ both $f(x,\alpha)$ and $f(\alpha,y)$
  are polynomials then $f$ is itself a polynomial.
\end{lemma}

We use the last lemma and the factorization property
(Theorem~\ref{thm:p2-unique-prod}) to determine what polynomials interlace $f^n$.

\begin{lemma}
  If $f,g\in\gsubpos_2$ and $f^n\lesslesseq g$ then there is a
  $g_1\in\gsubpos_2$ such $f\lesslesseq g_1$ and $g= f^{n-1}g_1$.
\end{lemma}

\index{quantitative sign interlacing}
\begin{proof}
If we substitute $y=\alpha$ then $f(x,\alpha)^n\lesslesseq
g(x,\alpha)$ so from Lemma~\ref{lem:sign-quant} we conclude
\begin{align*} g(x,\alpha) &= f(x,\alpha)^{n-1}g_1(x,\alpha) \\
\intertext{where we only know that $g(x,y)$ is a polynomial in $x$ for
  any value of $y$. Similarly we get }
g(\beta,y) & = f(\beta,y)^{n-1}g_2(\beta,y) \\
\intertext{where $g_2(x,y)$ is a polynomial in $y$ for any fixed $x$.
  Now}
\frac{g(x,y)}{f(x,y)^{n-1}} & = g_1(x,y) = g_2(x,y) 
\end{align*}
so $g_1$ satisfies the hypothesis of Lemma~\ref{lem:is-a-poly} and hence is a
polynomial in $x$ and $y$. Finally, to see that $f\lesslesseq g_1$
note that since $f^n\lesslesseq f^{n-1}g_1$ we know that
$$ f^{n-1}\,(f+\alpha g_1) \in\gsubpos_2$$
We now apply Theorem~\ref{thm:p2-unique-prod} to conclude that $f+\alpha
g_1\in\gsubpos_2,$ and so $f\lesslesseq g_1$.
\end{proof}

\begin{cor}
  $\mycone{(x+y)^n} = \bigl\{ \alpha (x+y)^{n-1}\,\bigl|\,\alpha\in\reals\bigr\}$
\end{cor}

\begin{cor}
  If $f\in\gsubpos_2$ then $\text{dim }\mycone{f} = \text{dim
    }\mycone{f^k}$ for any $k\ge 1$.
\end{cor}

Next, we  have an analog of quantitative sign interlacing. 
\index{quantitative sign interlacing}

\begin{theorem} \label{thm:cone-2}
  If $h\in\gsubclosepos_2$, $f(x)=(x-a_1)\dots(x-a_n)$ and
  $g(y)=(y-b_1)\dots(y-b_m)$ satisfy $f(x)g(y)\lesslesseq h$ then
  \begin{equation}
    \label{eqn:cone-1}
    h(x,y) = \sum_{i,j} c_{ij}\, \frac{f(x)}{x-a_i}\,\frac{g(y)}{y-b_j}
  \end{equation}
  where all $c_{ij}$ are non-negative.  Conversely, if
  $h$ is given by \eqref{eqn:cone-1} where all $c_{ij}$ have the same
  sign or are zero then $f(x)g(y)\lesslesseq h(x,y)$.
\end{theorem}

\index{quantitative sign interlacing}
\begin{proof}
  If we fix $y$ then from quantitative sign interlacing (Lemma~\ref{lem:sign-quant}) we
  can write
  \begin{align*}
    h(x,y) &= \sum \alpha_i(y) \,\frac{f(x)}{x-a_i} \\
\intertext{and by holding $x$ fixed }
    h(x,y) &= \sum \beta_j(x) \,\frac{g(y)}{y-b_j} \\
\intertext{Substituting $y=b_j$ yields}
h(x,b_j) &= \beta_j(x)g^\prime(b_j) \\
&= \sum \alpha_i(b_j) \, \frac{f(x)}{x-a_i} \\
\intertext{so $\beta_j(x)$ is a polynomial in $x$. Next substituting $x=a_i$
  and then $y=b_j$ gives }
h(a_i,y) &= \alpha_i(y)f^\prime(a_i) \\
h(a_i,b_j) &= \alpha_i(b_j)f^\prime(a_i) \\
\intertext{and so}
h(x,y) &= \sum_j \frac{h(x,b_j)}{g^\prime(b_j)}\, \frac{g(y)}{y-b_j}\\
       &= \sum_{i,j} \frac{\alpha_i(b_j)}{g^\prime(b_j)}\,
       \frac{f(x)}{x-a_i} \frac{g(y)}{y-b_j}\\
       &= \sum_{i,j} \frac{h(a_i,b_j)}{f^\prime(a_i)g^\prime(b_j)}\,
       \frac{f(x)}{x-a_i} \frac{g(y)}{y-b_j}
  \end{align*}
  
  We now determine the sign of $c_{ij}$. The sign of $h(a_i,b_j)$
  isn't arbitrary, for we must have that $f(x)\lesslesseq h(x,\alpha)$
  for any $\alpha$, and thus all the coefficients 
\begin{align*}
  s_i(\alpha) &= \sum_j\frac{h(a_i,b_j)}{f^\prime(a_i)g^\prime(b_j)}
  \,\frac{g(\alpha)}{\alpha-b_j} \\
  \intertext{of $\frac{f(x)}{x-a_i}$ must have the same sign. If we
    evaluate $s_i(x)$ at a root $b_k$ of $g$ we get} s_i(b_k) &=
  \frac{h(a_i,b_k)}{f^\prime(a_i)g^\prime(b_k)} \, g^\prime(b_k) \, =
  \,\frac{h(a_i,b_k)}{f^\prime(a_i)}
\end{align*}
Consequently $sgn\, s_i(b_k) = (-1)^{n+i} sgn(h(a_i,b_k))$ since the
sign of $f^\prime(a_i)$ is $(-1)^{n+i}$. Thus the sign of $h(a_i,b_k)$ is
$c_k(-1)^i$ where $c_k$ depends on $k$. A similar argument shows that 
the sign is also $d_i(-1)^k$, so we conclude that the sign is
$e(-1)^{i+j}$ where $e$ is a constant. This implies that the sign of 
$$\frac{h(a_i,b_j)}{f^\prime(a_i)g^\prime(b_j)}$$
is constant, as desired. The converse follows easily.
\end{proof}

\begin{cor}
  If all $a_i,b_j$ are distinct, $f(x)=(x-a_1)\dots(x-a_n),$ and
  $g(y)=(y-b_1)\dots(y-b_m)$ then
$\text{dim }\mycone{f(x)g(y)}\le nm$.
\end{cor}
\begin{proof}
  The basis for the cone is all the products 
$\frac{f(x)}{x-a_i}\,\frac{g(y)}{y-b_j}$.
\end{proof}

\section{Products of linear factors}
\label{sec:topology-products-1}

In this section we look at the interlacing properties of products of
linear factors. We construct polynomials in $\gsubplusclose_2$ whose
coefficients are mutually interlacing polynomials.

  \begin{theorem} \label{thm:product-3} Let $p_1,\dots,p_n$ and
    $q_1,\dots,q_n$ be two sequences of polynomials with positive
    leading coefficients, and define $$F =
    \prod_1^n(p_i+q_i\,y) = f_0(x)+f_1(x)y +\cdots+f_n(x)y^n.$$
\begin{enumerate}
\item If $p_i\lesslesseq q_i$ for $i=1,\dots,n$ then $f_i\lesslesseq
  f_{i+1}$ for $i=1,\dots,n$.
\item If $p_i\greateqeq q_i$ in $\allpolypos$ for $i=1,\dots,n$ then
  $f_i\greateqeq f_{i+1}$ for $i=1,\dots,n$.
\end{enumerate}
\end{theorem}
\begin{proof}
If  $p_i\lesslesseq  q_i$ then $p_i + yq_i$ is in $\gsubposclose_2$,
and so the product $F$ is in $\gsubposclose_2$. The coefficients of
$F$ are in $\allpoly$, and hence adjacent coefficients interlace.

In the second case, we know that consecutive coefficients interlace -
the problem is determining which direction.  We prove this by
induction. We need only to look at the two leading coefficients, so
let
  \begin{align*}
    p &= c x^r + dx^{r-1} + \cdots \\
    q &= \alpha x^r + \beta x^{r-1} + \cdots \\
    f_i &= a_i x^n + b_i x^{n-1} + \cdots \\
    f_{i+1} &= a_{i+1} x^n + b_{i+1} x^{n-1} + \cdots \\
    \intertext{Consecutive terms of the product $(p+yq)F$ are}
 &
    \left(\alpha \,a_{i-1} + c\,a_i \right)x^{n+r} + \left(\beta
      \,a_{i-1} + d\,a_i + \alpha \,b_{i-1} + c\,b_{i}
    \right)x^{n+r+1}+\cdots\\
&
    \left(\alpha \,a_{i} + c\,a_{i+1} \right)x^{n+r} + \left(\beta
      \,a_{i} + d\,a_{i+1} + \alpha \,b_{i} + c\,b_{i+1}
    \right)x^{n+r+1}+\cdots\\
\intertext{The determinant of the first two coefficients is}
&
\begin{vmatrix}
  \alpha a_{i-1} + c\,a_i & \beta a_{i-1} + d\,a_i + \alpha \,b_{i-1} +
  c\,b_{i}\\
\alpha \,a_{i} + c\,a_{i+1} & 
\beta a_{i} + d\,a_{i+1} + \alpha \,b_{i} + c\,b_{i+1}
\end{vmatrix}
\intertext{which can be written as}
& 
\begin{vmatrix}
  c & \alpha \\ d & \beta
\end{vmatrix}
\begin{vmatrix}
  a_i & a_{i+1} \\ a_{i-1} & a_i
\end{vmatrix}
+ \alpha^2
\begin{vmatrix}
  a_{i-1} & a_i \\ b_{i-1} & b_i
\end{vmatrix}
+
c^2
\begin{vmatrix}
  a_i & a_{i+1} \\ b_i & b_{i+1} 
\end{vmatrix}
+ \alpha c
\begin{vmatrix}
  a_{i-1}&a_{i+1} \\ b_{i-1} & b_{i+1}
\end{vmatrix}
  \end{align*}
The first determinant is positive since $p\greateqeq q$, the second is
Newton's inequality, and the last three follow from $f_i\greateqeq
f_{i+1}$. 
\end{proof}

\begin{cor}\label{cor:product-2}
  Suppose that constants $a_i,b_i,c_i,d_i$ have positive determinants
  $\smalltwodet{b_i}{a_i}{c_i}{d_i}$ and that
  $a_id_i$ is positive for $1\le i \le n$. If
\[
 \prod_{i=1}^n(a_ix+b_i +(d_ix +c_i)y  ) = f_0 + f_1y + \dots f_ny^n
\]
then 
\[
 f_0 \greateqeq f_1 \greateqeq \cdots \greateqeq f_n.
\]
\end{cor}

\begin{proof}
  If we set $p_i = d_ix+c_i$ and $q_i = a_ix+b_i$ then the
  relationship $p_i\greateqeq q_i$ is equivalent to $-c_i/d_i >
  -b_i/a_i$. Since $a_id_i$ is positive the condition is equivalent to 
  the positivity of the  determinant.
  Now apply Theorem~\ref{thm:product-3}.
\end{proof}

\begin{cor}
  Suppose that constants $c_i$ and $d_i$ satisfy 
\[ 0 < d_1 < c_1 < d_2 < c_2 < \cdots < d_n < c_n \]
If 
\[
 f = \prod_i^n(xy+x+c_iy+d_i)= f_0(x)+\cdots+f_n(x)y^n
\]
 then $(f_0,\dots,f_n)$
is mutually interlacing. \index{mutually interlacing}
\end{cor}
\begin{proof}
  The determinants of the factors are
  $\smalltwodet{1}{1}{c_i}{d_i}=d_i-c_i>0$, so we can apply
  Corollary~\ref{cor:product-2} to conclude that $$
  f_0 \greateqeq f_1 \greateqeq
  \cdots \greateqeq f_n.$$
  The hypothesis are that $f_0=\prod(x+d_i)
  \greateqeq \prod(x+c_i)=f_n$. It follows that this sequence of
  polynomials is mutually interlacing.
\end{proof}

\begin{remark}
  The interlacings in the last three results are strict if we assume
  that the interlacings in the hypothesis are strict.
\end{remark}

What happens in Corollary~\ref{cor:product-1} if some of the coefficients of $x$ are
zero?  Theorem~\ref{thm:product-2} shows that we still have interlacing, but some 
interlacings change from $\greateq$ to $\lessless$.

\begin{theorem} \label{thm:product-2}
  Suppose that $f\in\allpoly$ is a polynomial of degree $r$ and let 
$$ f(y)\prod_{i=1}^n (y+a_ix+b_i) = \sum_{i=0}^{n+r} p_i(x)y^i$$
Assume that $\prod(a_ix+b_i)$ has all distinct roots, all $a_i$ have
the same sign $\epsilon$, and write $f(y)=\sum c_iy^i$.
\begin{itemize}
\item If $c_0c_1\epsilon>0$ then
$$p_0 \greateq p_1 \greateq \dots \greateq p_n \lessless p_{n+1}
\lessless \dots \lessless p_{n+r-1}$$
\item If $c_0c_1\epsilon<0$ then
$$p_0 \lesseq p_1 \lesseq \dots \lesseq p_n \lessless p_{n+1}
\lessless \dots \lessless p_{n+r-1}$$
\end{itemize}
\end{theorem}

\begin{proof}
  Assume $\epsilon>0$. Both $f(y)$ and the product are in
  $\gsubclose_2$, so all the coefficients interlace. It remains to see
  how they interlace. The first two coefficients are
\begin{align*}
  p_0 & = c_0 \prod(a_ix+b_i) \\
p_1 & = c_1p_0 + c_0\sum_j\prod_{i\ne j}(a_ix+b_i) \\
&= c_1p_0 + c_0\sum_j\frac{p_0}{a_ix+b} \\
&= c_1p_0 + c_0\sum_j\frac{1}{a_i}\frac{p_0}{x+b/a_i}.
\end{align*}
Assume that $c_0c_1\epsilon>0$, so $c_0,c_1\ne0$.  By hypothesis $p_0$
has all distinct roots.  The coefficient $c_0/a_i$ has the same sign
as $c_1$, so $p_1 \lesseq p_0$.  Now apply Corollary~\ref{cor:fct-interlace}.  The
remaining cases are similar.
\end{proof}

\section{The interlacing of products}
\label{sec:topology-products}

In one variable, every product $\prod(x+a_i)$ is interlaced by many
different products. In two variables, the only products
$\prod(x+a_iy+b_i)$ that are interlaced by products are images of one
variable polynomials. In particular, most products in two variables
are not interlaced by any other products.  This is actually a result
about lines in the plane, and the proof requires the fundamental
result due to Sylvester and Gallai \cite{aigner-ziegler}:
\index{Sylvester-Gallai}

\begin{quote}
  In any configuration of $n$ points in the plane, not all on a line,
  there is a line which contains exactly two of them.
\end{quote}

Since this is a property of points, lines, and incidence, we can
dualize:

\begin{quote}
  In any configuration of $n$ lines in the plane, not all parallel nor
  all coincident, there is a point which lies on exactly two of them.
\end{quote}

\begin{theorem} \label{thm:int-prod}
  Suppose that $f,g\in\gsubpos_2(n)$ are products of linear factors, $fg$
  has no repeated factors, and $f,g$ interlace.  Then there are
  interlacing polynomials $h,k\in\allpoly$ with homogenizations $H,K$,
  and constants $\alpha,\beta$ such that either
\begin{xalignat*}{2}
  f(x,y) & = h(x+\alpha y+\beta) & g(x,y) & = k(x+\alpha y+\beta)
  \text{ or} \\
  f(x,y) & = H(x,\alpha y+\beta) & g(x,y) & = K(x,\alpha y+\beta)
\end{xalignat*}
\end{theorem}
\begin{proof}
  The graph of $f$ consists of $n$ lines of negative slope in the
  plane, as does the graph of $g$. By assumption, all these lines are
  distinct. If the lines are all parallel then we get the first case
  of the conclusion, and if they are all coincident then we have the
  second case. 
  
  If we substitute any value for $x$ then the resulting polynomials
  interlace, so we see that the intersections of the graph $G$ of $fg$
  with any vertical line almost always consist of $2n$ points, and
  these points alternate between being roots of $f$, or roots of $g$.
  The same is true for horizontal lines.
  
  Consider a point $(u,v)$ of the plane that lies on an even number of
  lines of $G$. By Sylvester-Gallai we know that there is at least one
  such point. Slightly to the left of $v$ a horizontal line meets the
  lines through $(u,v)$ in points that (from the left) are in
  $f,g,f,g,\dots$ (say), and slightly to the right of $v$ a horizontal
  line meets the lines through $(u,v)$ in $g,f,g,f,\dots$. Since a
  horizontal line meets all the lines of $G$, we get the important
  conclusion that on the horizontal line through $(u,v)$ all the lines of
  $G$ must intersect, and all such intersections have even degree. A
  similar statement holds for the vertical
  lines. Figure~\ref{fig:lines} shows an example of such
  intersections, where solid lines are factors of $f$, and dashed
  lines are factors of $g$.

\begin{figure}[htbp]
  \begin{center}
    \leavevmode
    \epsfig{file=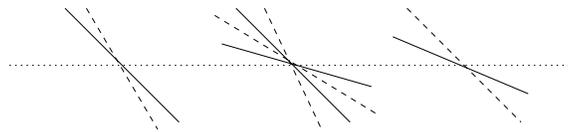,width=3in}
    \caption{The graph of $fg$ meeting a horizontal line}
    \label{fig:lines}
  \end{center}
\end{figure}

Now, take the leftmost such vertical line: this is the line through
the intersection point of $G$ that has an even number of lines through
it, and has smallest $x$ coordinate.  On this line, take the point of
intersecting lines with the largest $y$ coordinate. Label the
intersection $(s,t)$. A line of $G$ has negative slope, and can not
meet the line $y=t$ to the right of $(s,t)$, and also the line $x=s$
below $(s,t)$. This implies that all lines pass through $(s,t)$ and
this contradiction establishes the theorem.
\end{proof}

\begin{remark}
  Although the lemma shows that there are no interesting interlacing
  products of linear terms, there are interesting \emph{positive}
  interlacing products. Suppose that $f = \prod_1^n(x+a_iy+b_i)$ is a
  product in $\rupint{2}$. We claim that
\[
\frac{f(x,y)}{x+a_1+b_1} \poslace \frac{f(x,y)}{x+a_2+b_2} \poslace 
\cdots \poslace \frac{f(x,y)}{x+a_n+b_n}.
\]
Indeed, each of them interlaces $f$, so any positive linear
combination does as well. More generally, if $M = I+xD_1+yD_2$ where
$D_1,D_2$ are positive definite and $M[i]$ is the $i$th principle
submatrix then
\[
|M[1]| \poslace |M[2]|  \poslace \cdots \poslace |M[n]|.
\]

\end{remark}
\section{Characterizing transformations}
\label{sec:sub-characterize}

Just as in the one variable case, we can characterize linear
transformations  $T$ of  $\gsubpos_2$ that satisfy $f\lesslesseq Tf$. As
expected, $T$ is a linear combination of derivatives.

\changed{20/12/6}
    \begin{lemma}
      If $T:\rupint{2}\longrightarrow\gsub_2$ satisfies $f\lesslesseq Tf$
      for all $f\in\rupint{2}$ then there are non-negative $c,d$
      such that $Tf = \bigl(c \frac{\partial}{\partial x}  + d
      \frac{\partial}{\partial y }\bigr) f$.
    \end{lemma}

    \begin{proof}
      It suffices to show it holds for polynomials of degree $n$. Let
      $t_1,t_2$ be positive. Since
      \begin{gather*}
        (t_1x+t_2y)^n \lesslesseq T(t_1x+t_2y)^n
      \end{gather*}
      there is a constant $\alpha_t$ depending on $t_1,t_2$ such that
      \begin{equation}\label{eqn:p2-ch-1}
        T(t_1x+t_2y)^n = \alpha_t (t_1x+t_2y)^{n-1}. 
      \end{equation}
      We will determine $T$ by equating coefficients, but first we need to
find $\alpha_t$. Since $T$ decreases degree we can write
\begin{gather*}
  T(x^iy^{n-i}) = a_i x^iy^{n-1-i} + b_i x^{i-1}y^{n-i}
\end{gather*}
Thus
\begin{align}\label{eqn:p2-ch-2}
  T(t_1x+t_2y)^n &= t_1^nT(x^n) + n t_1^{n-1}t_2\, T(x^{n-1}y) + \cdots
  \\ \notag
&=
t_1\,b_n x^{n-1}  + n t_1^{n-1}t_2\bigl( a_{n-1} x^{n-1} + b_{n-1} x^{n-2}y\bigr) +
\cdots
\\ \notag &=
\alpha_t\bigl( t_1^{n-1} x^{n-1} + (n-1) t_1^{n-2}t_2
\,x^{n-2}y + \cdots\bigr)
\end{align}
Equating coefficients of $x^{n-1}$ in \eqref{eqn:p2-ch-1} yields
\begin{multline*}
\sum\binom{n}{i}t_1^it_2^{n-i}T(x^iy^{n-i})  \\ =
(t_1b_n + n b_{n-1}t_2)
\sum\binom{n-1}{j}t_1^jt_2^{n-1-j}T(x^jy^{n-1-j}) 
\end{multline*}
and conbined with \eqref{eqn:p2-ch-2} yields
\begin{gather*}
  t_1 b_n + n b_{n-1}t_2 = \alpha_t
\end{gather*}
Substituting this into \eqref{eqn:p2-ch-1}  and equating coefficients of
$t_1^it_2^{n-i}$ yields
\begin{align*}
  T(x^iy^{n-i}) &= b_n \,i\, x^{i-1}y^{n-i} + a_{n-1}\,(n-i) x^iy^{n-1-i}
  \\ &=
  \bigl( b_n \frac{\partial}{\partial x} + 
 a_{n-1} \frac{\partial}{\partial y}\bigr)\,x^iy^{n-i}
\end{align*}
The conclusion now follows by linearity.
    \end{proof}

This holds more generally; we sketch the proof.

\begin{lemma}
      If $T:\rupint{d}\longrightarrow\gsub_d$ satisfies $f\lesslesseq Tf$
      for all $f\in\rupint{d}$ then there are non-negative $c_i$
      such that $Tf = \bigl(\sum c_i \, \frac{\partial}{\partial x_i}\bigr)\,f$
    \end{lemma}
    \begin{proof}
      As above,
\begin{align*}
\left(\sum t_i x_i\right)^n &\lesslesseq T\, \left(\sum t_i x_i\right)^n \\
T\, \left(\sum t_i x_i\right)^n  &= \alpha_t \left( \sum t_ix_i\right)^{n-1} \\
\intertext{Write}
T(x_1^{i_1}\cdots x_d^{i_d}) & = \sum a_{i_1,\dots,i_d;j}x_1^{i_1}\cdots
x_j^{i_j-1}\cdots x_d^{i_d} \\
\intertext{Equating coefficients of $x_1^{n-1}$ yields}
\alpha_t &= t_1a_{n0\cdots0 }+ n   \sum_{j>1} t_j
  a_{n-10\cdots1\cdots0}
\end{align*}
Finally, equating coefficients uniquely determines $T$ on monomials,
which completes the proof.
\end{proof}

\section{Orthogonal polynomials in $\gsubpos_d$}
\label{sec:ortho-polys-d}

\index{Appell polynomials}\index{polynomials!Appell}
\index{orthogonal polynomials!Appell}

There are many different definitions for orthogonal polynomials in
several variables (see \cites{suetin,dunkl}). We
know two infinite families of orthogonal polynomials in more than one
dimension that are also in $\gsubpos_d$: the Hermite polynomials
(\chapsec{new-analytic}{hermite-n-dim}) and Appell polynomials (this section).

Knowing that an orthogonal polynomial is in $\gsubpos_2$ allows us to
recover some of the one dimensional properties concerning
roots. Although such a polynomial does not factor into linear factors,
all of its coefficients (of  powers of either $x$ or of $y$) factor
into linear factors, and adjacent coefficients interlace.

\index{Rodrigues' formula!orthogonal polynomials in $\rupint{d}$}
We begin by studying Rodrigues' formulas in $\gsubpos_d$.
\index{Rodrigues' formula} In one variable a Rodrigues' formula
represents a family of polynomials by an expression of the form $ch^{-1}
D^n(h p^n)$ where $h$ is a function, $c$ is a constant, and $p$ is a
polynomial. Here are three classic examples:

\index{Rodrigues' formula!Laguerre}
\index{Laguerre polynomials!Rodrigues' formula}
\index{Rodrigues' formula!Hermite}
\index{Hermite polynomials!Rodrigues' formula}
\index{Rodrigues' formula!Jacobi}
\index{Jacobi polynomials!Rodrigues' formula}

\begin{tabular}{rrl}
Laguerre & $L_n^{(\alpha)}(x)$ &= $\dfrac{1}{n!}
\dfrac{1}{  e^{-x}x^\alpha} \ \left(\dfrac{d}{dx}\right)^n\ \left(
  e^{-x} x^\alpha x^n\right)$\\
Hermite & $H_n(x)$ &= $(-1)^n \dfrac{1}{ e^{-x^2}} 
\ \left(\dfrac{d}{dx}\right)^n\ \left(
  e^{-x^2} \right)$\\
Jacobi & $P_n^{\alpha,\beta}(x)$ &= $\dfrac{(-1)^n}{2^nn!}
\dfrac{1}{ (1-x)^\alpha(1+x)^\beta }\quad \times $\\[.3cm]
&&\ $\left(\dfrac{d}{dx}\right)^n\ \left(
   (1-x)^\alpha(1+x)^\beta (1-x^2)^n \right)$
\end{tabular}

The next proposition generalizes the Rodrigues' formula for Laguerre
polynomials to define polynomials in $d$ variables.

\index{Rodrigues' formula!Laguerre in $\rupint{d}$}
\index{Laguerre polynomials!Rodrigues' formula for $\rupint{d}$}
\begin{prop} \label{prop:rod-laguerre}
  Suppose that $\bbb=(\beta_1,\dots,\beta_d)$ and
  $\diffi=(i_1,\dots,i_d)$ are all non-negative, and all entries of
  $\aaa=(\alpha_1,\dots,\alpha_d)$ are at least $-1$. If
  $f\in\gsubposclose_d$ then the following is a polynomial in
  $\gsubposclose_d$:
$$ \dfrac{1}{e^{-\xx\cdot\bbb}\xx^\aaa}
 \ \dfrac{\partial^{|\diffi|}}{\partial \xx^\diffi}
\ \left( e^{-\xx\cdot\bbb} \xx^\aaa f^{}\right)
$$
\end{prop}

\begin{proof}
  To reduce the necessary notation, let's assume $d=2$.  We can not
  just apply results about differentiation of polynomials in
  $\gsubclosepos_2$ since the objects we are differentiating are not
  polynomials when the $\alpha_i$ are not integers and the $\beta_i$
  aren't all zero.  So, we enlarge our class of objects:
$$ P(r,s) =
\left\{ e^{-\beta_1 x - \beta_2 y} x^{\alpha_1} x^r y^{\alpha_2} y^s
 f\,\mid\,f\in\gsubposclose_2\right\}
$$
We then are going to show that 
\begin{align} 
  \frac{\partial }{\partial x}&:
  P(r,s)\longrightarrow
  P(r-1,s)\label{eqn:appell-1} \\
  \frac{\partial }{\partial y}&:
  P(r,s)\longrightarrow
  P(r,s-1) \label{eqn:appell-2}
\end{align}
Since
 $$
 e^{-\beta_1 x - \beta_2 y} x^{\alpha_1} x^r y^{\alpha_2} y^s f
 \in  P(r,s)$$
it will follow that 
$$\dfrac{\partial^{n+m}}{\partial x^n \partial y^m}
\left( e^{-\beta_1 x - \beta_2 y} x^{\alpha_1} x^r y^{\alpha_2} y^s
 f\right)
\in P(0,0)$$
which implies the conclusion. 

It suffices only to prove  \eqref{eqn:appell-1}; the proof
of \eqref{eqn:appell-2} is similar. Differentiating,

\begin{multline*}
\frac{\partial }{\partial x} \left( e^{-\beta_1 x - \beta_2 y}
  x^{\alpha_1} x^r y^{\alpha_2} y^s f\right) = e^{-\beta_1 x - \beta_2
  y} x^{\alpha_1} x^{r-1} y^{\alpha_2} y^s \quad \times \\ \left( -\beta_1 x f +
  (\alpha_1+r)f+x\frac{\partial f}{\partial x}\right)
\end{multline*}
Since $r$ is at least $1$, $\alpha_1+r$ is non-negative, so we know
that that $-\beta_1xf+(\alpha_1+r)f+x\frac{\partial f}{\partial x}
\in\gsubpos_2$, and hence \eqref{eqn:appell-1} holds.
\end{proof}

There are several varieties of Appell polynomials. We  begin with
the simple, no parameter version, and then introduce more complexity.
The version below is in \cite{suetin} for $d=2$. If
$\diffi=(i_1,\dots,i_d)$ and $\xx=(x_1,\dots,x_d)$ then define

\begin{align*}
A_\diffi(\xx) & =
\dfrac{1}{\diffi!}\dfrac{\partial^{|\diffi|}}{\partial\xx^\diffi}
\left(\ \left(1-\Sigma \xx\right)^{|\diffi|}\xx^\diffi\ \right)
\\
&= \dfrac{1}{i_1!\cdots i_d!} \dfrac{\partial^{i_1+\cdots+i_d}}%
{\partial x_1^{i_1}\cdots \partial x_d^{i_d}}%
\left(\ \left(1-x_1-\cdots-x_d\right)^{i_1+\cdots+i_d}
x_1^{i_1}\cdots x_d^{i_d}\ \right)
\end{align*}

\begin{prop} \label{prop:appell-1}
  $A_\diffi(\xx)$ is in $\gsubposclose_d$.
\end{prop} 
\begin{proof}
    Since $(1-\Sigma\,\xx)^{|\diffi|}$ is in $\gsubpos_d$, the product
    with $x_1^{i_1}\cdots x_d^{i_d}$ is in $\gsubposclose_d$. Since
    $\gsubposclose_d$ is closed under differentiation,
    $A_\diffi(\xx)\in\gsubposclose_d$. 
\end{proof}

\begin{prop} \label{prop:appell-2}
  Adjacent Appell polynomials interlace.
\end{prop}
\begin{proof}
Without loss of generality,   we have to show that
\[    A_{i_1,\dots,i_d}(\xx) \lesslesseq A_{i_1,\dots,i_d-1}(\xx).\]
Since we can ignore constant factors, we need to show that
\begin{multline*}
 \dfrac{\partial^{|\diffi|}}{\partial x_1^{i_1}\cdots \partial x_d^{i_d}}%
\left(\ \left(1-\Sigma \xx\right)^{|\diffi|}x_1^{i_1}\cdots
  x_d^{i_d}\ \right) \lesslesseq \\
 \dfrac{\partial^{|\diffi|-1}}{\partial x_1^{i_1}\cdots \partial x_d^{i_d-1}}%
\left(\ \left(1-\Sigma \xx\right)^{|\diffi|-1}x_1^{i_1}\cdots
  x_d^{i_d-1}\ \right) 
\end{multline*}
{Since differentiation preserves interlacing, we need to
  show that}
  \begin{gather*}
 \dfrac{\partial}{\partial x_d}
\left(\ \left(1-\Sigma\xx\right)^{|\diffi|}x_1^{i_1}\cdots
  x_d^{i_d}\ \right) \lesslesseq 
\left(\ \left(1-\Sigma\xx\right)^{|\diffi|-1}x_1^{i_1}\cdots
  x_d^{i_d-1}\ \right) \\
\end{gather*}
\noindent{Differentiating and factoring the left side yields}
\begin{multline*}
(1-\Sigma\xx)^{|\diffi|-1}x_1^{i_1}\dots x_d^{i_d-1}
\left(-|\diffi|x_d + i_d(1-\Sigma\xx)\right) \lesslesseq\\
\left(\ \left(1-\Sigma\xx\right)^{|\diffi|-1}x_1^{i_1}\cdots
  x_d^{i_d-1}\ \right) 
\end{multline*}
Since the left side is a multiple of the right side by a linear term
that is in $\gsubpos_d$, the last interlacing holds, and so the
proposition is proved.
\end{proof}

The polynomials for $d=1,2$ are
\begin{align*}
A_n(x) &=  \dfrac{1}{n!} \dfrac{d^n}{dx^n}\left(\ (1-x)^n x^n\
\right)\\
A_{n,m}(x_1,x_2) &= \dfrac{1}{n!m!} \dfrac{\partial^{n+m}}%
{\partial x_1^{n} \partial x_2^{m}}%
\left(\ \left(1-x_1-x_2\right)^{n+m}
x_1^{n} x_2^{m}\ \right)
\end{align*}

The one dimensional polynomials are a modified Legendre polynomial.
The first few of the one and two dimensional polynomials are listed in
Chapter~\ref{cha:tables} The Appell polynomials of different
dimensions are related since $A_{i_1,\dots,i_d}(x_1,\dots,x_d) =
A_{i_1,\dots,i_d,0}(x_1,\dots,x_d,0).$

 Next we consider Appell polynomials with extra parameters
(\cite{suetin}). Let $\aaa=(\alpha_1,\dots,\alpha_d)$, and
define

\begin{align*}
A^{\aaa}_\diffi(\xx) & =
\dfrac{\xx^{-\aaa}}{\diffi!}\dfrac{\partial^{|\diffi|}}{\partial\xx^\diffi}
\left(\ \left(1-\Sigma \xx\right)^{|\diffi|}\xx^{\diffi+\aaa}\ \right)
\\
&= \dfrac{x_1^{-\alpha_1}\cdots x_d^{-\alpha_d}}{i_1!\cdots i_d!}
\dfrac{\partial^{i_1+\cdots+i_d}}%
{\partial x_1^{i_1}\cdots \partial x_d^{i_d}} \\
&\hspace*{1in}
\left(\ \left(1-x_1-\cdots-x_d\right)^{i_1+\cdots+i_d} x_1^{i_1+\alpha_1}\cdots
  x_d^{i_d+\alpha_d}\ \right)
\end{align*}

\begin{prop} \label{prop:appell-3} If  $\alpha_i\ge-1$ for $1\le i \le
  d$ then 
  $A^{\aaa}_\diffi(\xx)\in\gsubpos_d$
\end{prop}

\begin{proof}
Since $(1-\Sigma\diffi)^{|\diffi|}\xx^\diffi\in\gsubposclose_d$, the
result follows from Proposition~\ref{prop:rod-laguerre} with $b=0$.
\end{proof}

If we attempt to generalize the Rodrigues' formula for the Laguerre
polynomials with the definition
$$ \dfrac{1}{e^{-\xx\cdot\bbb}\xx^\aaa}
 \ \dfrac{\partial^{|\diffi|}}{\partial \xx^\diffi}
\ \left( e^{-\bbb\cdot\xx} \xx^\aaa \xx^\diffi\right)
$$
then it is easy to see that this factors into a product of one
dimensional Laguerre polynomials.  We use $\Sigma \xx$ to force the
$x_i$'s to be dependent. The following polynomials are in $\gsubposclose_d$.

$$ \dfrac{1}{e^{-\xx\cdot\bbb}\xx^\aaa}
 \ \dfrac{\partial^{|\diffi|}}{\partial \xx^\diffi}
\ \left( e^{-\bbb\cdot\xx} \xx^\aaa (\Sigma\xx)^{|\diffi|}\right)
$$

Here is a two dimensional example that is $(28)$ from
\cite{suetin}*{Page 159}. Since $xy-1\in\gsubposclose_2$, the
following is in $\gsubposclose_2$ for $b>0$:

$$
\dfrac{1}{e^{-b(x+y)}}  \ \dfrac{\partial^{n+m}}{\partial x^n\,\partial
  y^m}
\, \left(e^{-b(x+y)}(xy-1)^{n+m}\right)
$$

\index{Legendre polynomials!generalized}\index{polynomials!Legendre!generalized}%
\index{Rodrigues' formula!generalized Legendre}
The Rodrigues' formula for the Legendre polynomials is
$$ P_n(x) = \frac{1}{2^nn!} \left(\frac{d}{dx}\right)^n (x^2-1)^n$$
We can give a Rodrigues' formula for an analog of Legendre polynomials, but
the resulting family does not appear to have any recursive
properties. Let $Q$ be a negative subdefinite symmetric $d$ by $d$
matrix, and let $\xx=(x_1,\dots,x_d)$. For any index set $\diffi$ we
define

$$ P_\sdiffi(\xx) = \dfrac{\partial^{|\diffi|}}{\partial \xx^\diffi} \
(\xx Q \xx^t -1)^{|\sdiffi|} $$

Since $Q$ is negative subdefinite, the \index{quadratic form}quadratic form $\xx Q\xx^t-1$
is in $\gsubpos_d$, and since multiplication and differentiation
preserve $\gsubpos_d$, it follows that $P_\sdiffi(\xx)\in\gsubpos_d$.

\begin{remark}
  There is a different definition of Appell polynomials in \cite{dunkl} that does
  not lead to polynomials in $\gsubposclose_d$. In the case that $d=2$ these
  polynomials are, up to a constant factor, defined to be
$$ U_{n,m}(x,y) = (1-x^2-y^2)^{-\mu+1/2}\, \dfrac{\partial^{n+m}}{\partial
  x^n\,\partial y^m} \left(1-x^2-y^2\right)^{n+m+\mu-1/2}
$$
where $\mu$ is a parameter. Since
$$U_{2,0}(x,y) = -\left( \left( 3 + 2\,\mu  \right) \,\left( 1 -
    (2+ 2\,\mu)\,x^2 - y^2   \right)  \right)$$
it is clear that $U_{2,0}(x,y)$ does not satisfy substitution, and so
is not in $\gsubpos_2$.
\end{remark}

When we generalize the classical orthogonal polynomials, there are
many choices and no one way of doing it. We can notice that these
generalized Legendre polynomials are similar to the definition of
Appell polynomial in \cite{dunkl}, except that he uses a positive
definite matrix (the identity) leading to the \index{quadratic form}quadratic form (in two
variables) $x^2+y^2-1$. Replacing the identity by a negative
subdefinite matrix leads to a definition of generalized Appell
polynomials as
$$
\frac{1}{(\xx Q\xx^t -1)^\mu} \partial^\sdiffi \ (\xx Q \xx^t
-1)^{|\sdiffi|+\mu}$$
where $\mu$ is a parameter. That this polynomial
is in $\gsubpos_d$ follows from the following theorem, whose proof is
similar to Proposition~\ref{prop:rod-laguerre} and omitted.

\begin{prop}
  If $f\in\gsubpos_d$, $\mu$ any parameter, $\diffi$ any index set,
  then $$
  \frac{1}{f^\mu}\ \dfrac{\partial^{|\diffi|}}{\partial
    \xx^\diffi} \ f^{|\sdiffi|+\mu} \in\gsubpos_d$$.
\end{prop}

\section{Arrays of interlacing polynomials}
\label{sec:top:arrays}

The polynomials determined by Rodrigues' type formulas form interlacing
arrays. For instance, consider the polynomials
$$
f_{n,m} = \dfrac{\partial^{n+m}}{\partial x^n\,\partial y^m} \,
g^{n+m}
$$
where $g\in\rupint{2}$ is of degree $2$. Since $\frac{\partial
  g}{\partial x}$ is a polynomial of degree $1$ in $\rupint{2}$ we know
that
\begin{align*}
  g^{n+m-1}\,  \frac{\partial g}{\partial x} & \lesslesseq g^{n+m-1}\\
\intertext{and therefore }
  \dfrac{\partial^{n+m-1}}{\partial x^{n-1}\,\partial y^m}\,
\frac{\partial}{\partial x}   g^{n+m} &
 \lesslesseq \dfrac{\partial^{n+m-1}}{\partial x^{n-1}\,\partial y^m}  g^{n+m-1}\\
\end{align*}
which implies that $f_{n,m} \lesslesseq f_{n-1,m}$. These interlacing
take the following form, where $h\longleftarrow k$ stands for $h
\lesslesseq k$:

\centerline{
\xymatrix{
f_{0,2} &            f_{1,2}  \ar@{<-}[l]           & f_{2,2} \ar@{<-}[l]\\ 
f_{0,1} \ar@{->}[u] &f_{1,1} \ar@{->}[u] \ar@{<-}[l]& f_{2,1} \ar@{<-}[l]\ar@{->}[u] \\
f_{0,0} \ar@{->}[u] &f_{1,0} \ar@{->}[u] \ar@{<-}[l]& f_{2,0} \ar@{<-}[l]\ar@{->}[u] \\
}}

We consider two examples for $g$ that are of the form $\xx Q\xx^t-1$
where $Q$ is negative subdefinite.  

\begin{example}
First, take $g = xy-1$. 
In this case we have an explicit formula for $f_{n,m}$:

\begin{align*}
  f_{n,m} &= \dfrac{\partial^{n+m}}{\partial x^n\,\partial y^m} \,
(xy-1)^{n+m} \\
&= \sum \dfrac{\partial^{n+m}}{\partial x^n\,\partial y^m} \,
\binom{n+m}{i} x^i y^i (-1)^{n-i} \\
&= \sum_{i=\max(n,m) }^{n+m} 
\binom{n+m}{i}\falling{n}{i} \falling{m}{i} x^{i-n} y^{i-m}  \\
&= (n+m)! \sum_{i=\max(n,m) }^{n+m} \binom{i}{i-n,i-m} x^{i-n}y^{i-m}
\end{align*}

Using this formula we can verify that $f_{n,m}$ satisfies a simple
two term recurrence relation
\index{recurrence!two dimensional}

\begin{equation*}
  f_{n,m} = \frac{1}{n}\left( \,(n+m)^2 \,y\,f_{n-1,m} -
    m(n+m)(n+m-1)f_{n-1,m-1}\, \right)
\end{equation*}

If we substitute $y=1$ we get an array of one variable polynomials
such that every polynomial interlaces its neighbors.

\hspace*{-.6cm}%
$
\begin{array}{ccc}
   x & 2\,x^2 & 6\,x^3 \\
  -2 + 4\,x & -12\,x + 18\,x^2 & -72\,x^2 + 96\,x^3 \\
   -12 + 18\,x &24 - 144\,x + 144\,x^2&360\,x - 1440\,x^2 +
   1200\,x^3 \\
   -72 + 96\,x&360 - 1440\,x + 1200\,x^2&-720 + 8640\,x - 21600\,x^2 + 14400\,x^3\\
\end{array}
$

\end{example}

\begin{example}
  Next, we consider the \index{quadratic form}quadratic form corresponding to the negative
  subdefinite matrix $\smalltwo{1}{2}{2}{1}$. If we define 
$$
f_{n,m} = \dfrac{\partial^{n+m}}{\partial x^n\,\partial y^m} \,
(x^2 + y^2 + 4xy -1)^{n+m}
$$
then we again get an array of interlacing polynomials.  The degree of
$f_{n,m}$ is $n+m$. In this case if we substitute $y=0$ then we get an
array of interlacing one variable polynomials, but there does not
appear to be any nice recurrence relation.

\hspace*{-.5cm}%
$
\begin{array}{cccc}
 1 & 4\,x & -4 + 36\,x^2 \\
  2\,x & -8 + 24\,x^2 & -216\,x + 408\,x^3 \\ 
   -4 + 12\,x^2 & -144\,x + 240\,x^3 & 432 - 4896\,x^2 + 6000\,x^4 \\
  -72\,x + 120\,x^3 &288 - 2880\,x^2 + 3360\,x^4& 24480\,x - 120000\,x^3 + 110880\,x^5
\end{array}
$

\end{example}

\chapter{Polynomials satisfying partial substitution}
\label{cha:partial}

\renewcommand{\TimeStampStart}{Monday, December 17, 2007: 17:22:09}
\mytoday    

In this chapter we are interested in polynomials in two kinds of
variables that satisfy substitution only for positive values of the
variables of the second kind.  Such polynomials naturally arise from
transformations that map $\allpolypos$ to $\allpoly$.

\section{Polynomials satisfying partial substitution in two variables}
\label{sec:partial-class}

We introduce the class $\partialpoly{1,1}$ of polynomials $f(x,y)$ that
satisfy substitution for positive $y$.

\begin{definition} \label{defn:partialpoly}
  The class $\partialpoly{1,1}$ consists of all polynomials $f(x,y)$
  satisfying
  \begin{enumerate}
  \item The homogeneous part $f^H$ is in $\allpolypos$.
  \item $f(x,\beta)\in\allpoly$ for all $\beta \ge 0$.
  \end{enumerate}

If all the coefficients are non-zero, then $f\in\partialpolypos{1,1}$.
The polynomials in the closure of $\partialpoly{1,1}$ are
$\partialpolyclose{1,1}$, and those in the closure of
$\partialpolypos{1,1}$ are $\partialpolyposclose{1,1}$.
\end{definition}

We have the inclusion $ \rupint{2} \subsetneq \partialpoly{1,1}$.  The
next example shows that $\partialpoly{1,1}$ is not closed under
differentiation, and so the  containment is strict.

\begin{example}

  $\partialpoly{1,1}$ is not closed under differentiation with respect
  to $y$. However, if $f\in\partialpoly{1,1}(n)$ then there is an
  $\alpha\ge0$ such that $\frac{\partial f}{\partial
    y}(x,y+\alpha)\in\partialpoly{1,1}(n)$.

  The calculations below show that $f\in\partialpoly{1,1}$, but that 
$\frac{\partial f}{\partial  y}\not\in\partialpoly{1,1}$. The
calculations do show that 
$\frac{\partial f}{\partial  y}(x,y+\alpha)\in\partialpoly{1,1}$ for
$\alpha\ge 1/8$.
  \begin{align*}
    f  &={\left( -3 + x + y \right) }^2\,\left( 18 - 19\,x + 5\,x^2 -
      22\,y + 10\,x\,y + 5\,y^2 \right)\\
    f_y &= 2\,\left( -3 + x + y \right) \,\left( 51 - 45\,x + 10\,x^2 - 48\,y + 20\,x\,y + 10\,y^2 \right)\\
\text{roots of } f: \quad&\quad
 \left\{ 3 - y,3 - y,\frac{19 - 10\,y \pm {\sqrt{1 + 60\,y}}}{10}\right\}\\
\text{roots of } f_y: \quad
&\quad \left\{ 3 - y,\frac{45 - 20\,y \pm {\sqrt{15}}\,{\sqrt{-1 +
8\,y}}}{20}\right\}
  \end{align*}
\end{example}

Every vertical line $x=\beta$ to the left of the smallest root  of $f(x,0)$ has
$n-1$ roots of $\frac{\partial f}{\partial  y}(\beta,y)$. Since these
$n-1$ solution curves are asymptotic to lines whose slopes lie between
the slopes of the asymptotes of $f$, every sufficiently positive
horizontal line will meet all $n-1$ curves.

What do the graphs of polynomials in $\partialpoly{1,1}$ look like?
The solution curves do not
necessarily go from the upper left quadrant to the lower right
quadrant as is the case for polynomials is $\gsubpos_2$, but they can
loop back to the upper left quadrant.

  \begin{figure}[htbp]
    \centering
    \includegraphics*[height=3in]{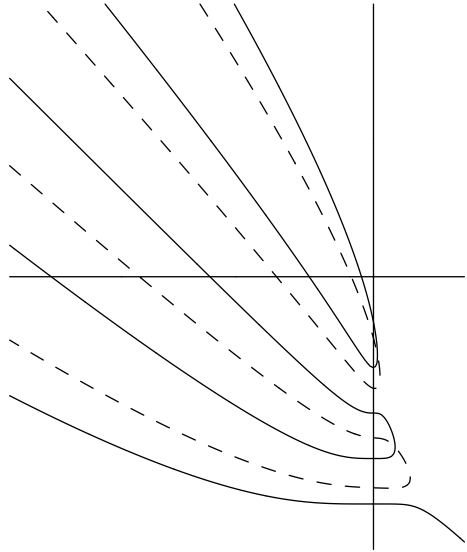}
    \caption{$L_5(x;y)$ and its derivative (dashed) with respect to $y$}
    \label{fig:laguerre2}
  \end{figure}

In general, since $f^H$ has all negative roots, the asymptotes of $f$
are lines in the upper left quadrant with negative slopes. As
\index{solution curves}
$x\rightarrow-\infty$, the solution curves eventually have positive
$y$-coordinates. Thus, there is an $\alpha_1$ such that
$f(\alpha,y)\in\allpolyalt$ for all $\alpha<\alpha_1$, and also an
$\alpha_0$ such that $f(\alpha,y)\in\allpoly$ for all
$\alpha<\alpha_0$.  In Figure~\ref{fig:laguerre2} $\alpha_0$ is
about $.1$, and $\alpha_1$ is about $-12.6$. A formal proof of this
general fact follows the argument of   Theorem~\ref{thm:sub-xy}.

\section{Examples}
\label{sec:partial-two-examples}

Polynomials satisfying partial substitution arise naturally in several
ways.

\begin{example}\index{recurrence!Laguerre}
  The generalized Laguerre polynomials $L_n(x;y)$ are in
  $\partialpoly{1,1}$ since their homogeneous part is
  $\frac{1}{n!}(x+y)^n$, and they they satisfy the three term
  recurrence \cite{szego}*{(5.1.10)}
\begin{gather*}
n L_n(x;y) = (-x+2n+y-1)L_{n-1}(x;y) - (n+y-1)L_{n-2}(x;y)\\
L_0(x;y)=1\quad L_1(x;y) = x+y+1
\end{gather*} 
\end{example}

The following lemma generalizes the fact that three term
recurrences determine orthogonal polynomials.

\begin{lemma}\label{lem:3-term-partial}
  Suppose that $a_i,b_i,d_i,e_i$ are positive for $i\ge0$. If
  \begin{gather*}
    f_n(x;y) = (a_n x + b_n y + c_n)f_{n-1}(x;y) - (d_n y +
    e_n)f_{n-2}(x;y) \\
f_0(x;y) = 1 \quad f_{-1}(x;y) = 0
  \end{gather*}
then $f_n(x;y)\in\partialpoly{1,1}$ for $n=0,1,\dots$.
\end{lemma}
\begin{proof}
  The construction yields polynomials that satisfy the degree
  requirements, and the homogeneous part is $\prod_i (a_i x + b_i y)$.
  If $y$ is positive then the recurrence defines a three term
  recurrence for orthogonal polynomials, so they are in $\allpoly$.
\end{proof}

Note that if $\alpha$ is sufficiently negative, then $f_n(x;\alpha)$
satisfies a recurrence that does not define a sequence of orthogonal
polynomials, so we do not expect any polynomials from this
construction to lie in $\rupint{2}$.

\begin{example}\index{Charlier polynomials!two variable}
  We can construct two variable Charlier polynomials from the one variable
  polynomials by
\[ C_n(x;y) = (-1)^nC_n^y(-x).\]
They also satisfy  a recurrence relation $C_{-1}(x;y)=0$, $C_0(x;y)=1$
and
\[ C_{n+1}(x;y) = (x+y+n)\,C_n(x;y) - n\, y\,C_{n-1}(x;y). \]
It follows from the Lemma that $C_n(x;y)\in\partialpoly{1,1}$.
Figure~\ref{fig:charlier} shows $C_5(x;y)$.

\begin{figure}[htbp]
  \centering
  \includegraphics*[width=3in]{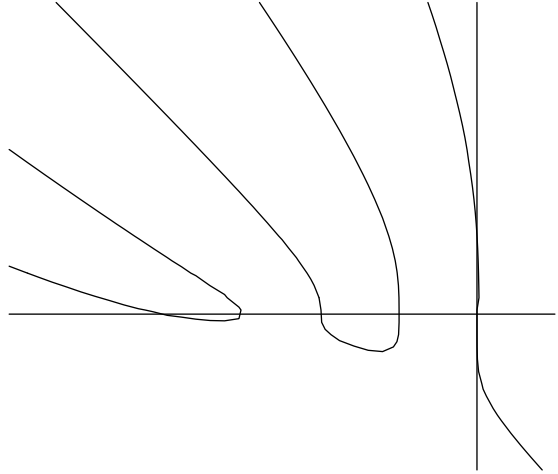}
  \caption{The Charlier polynomial $C_5(x;y)$}
  \label{fig:charlier}
\end{figure}
\end{example}

 If $f\in\allpolypos$ then
$f(-xy)\in\gsubclose_2$, and consequently is in
$\partialpolyclose{1,1}$. However, more is true:
\begin{lemma}\label{lem:partial-xy}
  If $f\in\allpoly$, then $f(xy)\in\partialpolyclose{1,1}$.
  If $f\in\allpolypos$, then $f(xy)\in\partialpolyposclose{1,1}$.
\end{lemma}
\begin{proof}
  Since $\partialpoly{1,1}$ is closed under multiplication, the
  comment above shows that it suffices to show that
  $1+xy\in\partialpolyclose{1,1}$. Consider $f_\epsilon (x,y) =
  1+(x+\epsilon y)(y+\epsilon x)$. For positive $\epsilon$ the graph
  of $f_\epsilon$ is a hyperbola, and the minimum positive value is
  found to be $d_\epsilon=\frac{2\sqrt{\epsilon}}{|1-\epsilon^2|}$.
  Since the asymptotes have slope $-\epsilon$ and $-1/\epsilon$, the
  polynomial $f_\epsilon-d_\epsilon$  meets every horizontal
  line above the $x$-axis in two points, and the degree and
  homogeneity conditions are met, so
  $f_\epsilon-d_e\in\partialpoly{1,1}$. Since
  $\lim_{\epsilon\rightarrow0} (f_\epsilon-d_\epsilon) = f$ it follows
  that $f\in\partialpolyclose{1,1}$.
\end{proof}

  \begin{example}
    We can construct polynomials  that satisfy partial substitution
    but not the homogeneity condition in a manner that is similar to
    the construction of polynomials in $\gsubplus_2$. Suppose that $A$
    is a positive definite $n$ by $n$ matrix, and $B$ is a symmetric
    matrix that has positive and negative eigenvalues. Define $f(x,y)
    = |I+y A + x B|$. If $y$ is positive then $I+yA$ is positive
    definite, and so all roots of $f(x,y)$ are real. The homogeneous
    part of $f(x,y)$ is $|yA+xB|$.  Since $A$ is positive definite,
    the roots of $f^H(x,1)$ are the roots of $|I+x
    A^{-1/2}BA^{-1/2}|$. Since $B$ has positive and negative
    eigenvalues, Sylvester's law of inertia shows that
    $A^{-1/2}BA^{-1/2}$ does also. Consequently, the roots of $f^H$
    are not all the same sign, and thus $f\not\in\partialpoly{1,1}$.
    
    Consider a concrete example. $A$ was chosen to be a random
    positive definite matrix, $C$ was a random symmetric matrix, $D$
    was the diagonal matrix with diagonal $(1,1,1,-1)$, and $B = CDC$.
    If $f(x,y) = |I +y A + x B|$, then Figure~\ref{fig:ranpartialpoly}
    shows that $f(x,y)\not\in\rupint{2}$. (The dashed line is the
    $x$-axis.) 

\begin{align*} f(x,y) & = 1 + 11.42\,x + 4.84\,x^2 - 7.13\,x^3 -
  3.67\,x^4 + 15.48\,y \\
&+ 42.56\,x\,y - 9.60\,x^2\,y - 21.88\,x^3\,y + 
  12.05\,y^2 + 13.80\,x\,y^2  \\ &- 5.94\,x^2\,y^2 + 1.99\,y^3 +
  0.96\,x\,y^3 + 0.07\,y^4
\end{align*}

  \begin{figure}[htbp]
    \centering
    \includegraphics*[width=2in]{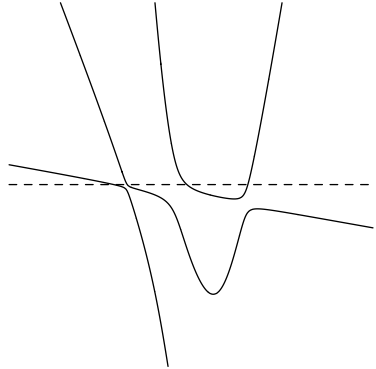}
    \caption{A polynomial of the form $|I+yA+xB|$ }
    \label{fig:ranpartialpoly}
  \end{figure}
\end{example}

We summarize these observations in a lemma. The last two parts are
proved in the same way as the results for $\rupint{2}$.

\begin{lemma}\label{lem:det-partial-poly}
  If $A$ is positive definite, and $B$ is symmetric with positive and
  negative eigenvalues then $f(x,y) = |I+yA+xB|$ satisfies
  \begin{enumerate}
  \item $f(x,y)$ has asymptotes of positive and negative slope.
  \item $f(x,\alpha)\in\allpoly$ for $\alpha\ge0$.
  \item $f(\alpha,y)\in\allpoly$ for all $\alpha$.
  \item $f(\alpha x,x)\in\allpoly$ for all $\alpha$.
  \end{enumerate}
\end{lemma}

\section{Induced transformations}
\label{sec:partialpoly-induced}

We now look at induced linear transformations on $\allpolypos$.

\begin{lemma} \label{lem:partial-trans}
  If $T\colon{}\allpolypos\longrightarrow\allpolypos$, $T$ preserves degree
  and the sign of the leading coefficient, then
  $T_\ast:\partialpolypos{1,1}\longrightarrow\partialpolypos{1,1}$.
  If $T\colon{}\allpolypos\longrightarrow\allpoly$ then
  $T_\ast:\partialpolypos{1,1}\longrightarrow\partialpoly{1,1}$.
\end{lemma}
\begin{proof}
  The hypotheses imply that $T_\ast f$ satisfies 
  homogeneity assumption.  We next verify substitution for $y$.  For
  $\beta\ge0$ we have $$(T_\ast f(x,y))(x,\beta) =
  T(f(x,\beta))\in\allpolypos$$
  since $f(x,\beta)\in\allpolypos$, and
  consequently $T_\ast(f)\in\partialpolypos{1,1}$. The  second part is
  similar.
\end{proof}

We can easily determine the slopes of the asymptotes of $T_\ast
f$. Suppose that $f=\sum a_{i,j}x^iy^j$, and $c_i$ is the leading
coefficient of $T(x^i)$. Then
$$ (T_\ast f)^H(x,1) = \left(\sum a_{i,j}T(x^i)y^j\right)^H (x,1) =
  \sum a_{i,n-i}\,c_ix^i$$
The slopes of the asymptotes are the roots of this last polynomial.

\begin{cor}
  Suppose $T\colon{}\allpolypos\longrightarrow\allpolyneg$, $T$ preserves
  degree and the sign of the leading coefficient.  If
  $f\in\allpolypos$ then $T_\ast f(x+y)\in\partialpolypos{1,1}$.
\end{cor}
\begin{proof}
  Since $f(x+y)\in\partialpolypos{1,1}$ we can apply
  Lemma~\ref{lem:partial-trans}.
\end{proof}

 

\begin{lemma}\label{lem:partial-homog-1} 
  Choose a positive integer $n$, and define the linear transformation
  $S\colon{}x^i\mapsto T(x^i)y^{n-i}$. If
  $T\colon{}\allpolypos\longrightarrow\allpolypos$ preserves degree and the
  sign of the leading coefficient then
  $S\colon{}\allpolypos\longrightarrow\partialpolypos{1,1}$.
\end{lemma}
\begin{proof}
  If $f\in\allpolypos$ then the homogenization $F$ of $f$ is in
  $\partialpoly{1,1}$. Since $S(f) = T_\ast(F)$ we can apply
  Lemma~\ref{lem:partial-trans} to conclude that
  $T_\ast(F)\in\partialpolypos{1,1}$.
\end{proof}

\index{rising factorial}
\begin{cor} \label{cor:partial-fall}
  If $f\in\allpolypos(n)$  and $S(x^i) =
  \rising{x}{i}y^{n-i}$ then $S(f)\in\partialpolypos{1,1}$. 
\end{cor}
\begin{proof}
  The map $x^i\mapsto \rising{x}{i}$ satisfies the hypotheses of
  Lemma~\ref{lem:partial-homog}.
\end{proof}

\begin{example}
  If we choose $f=(x-1)^n$ then from the corollary
$$
    S(\,(x+1)^n\,) =  \sum_{i=0}^n \rising{x}{i} y^{n-i}\binom{n}{i}
    = n! \sum_{i=0}^n \binom{x+i-1}{i}\frac{y^{n-i}}{(n-i)!}
$$
The last polynomial is in $\partialpolypos{1,1}$.
\end{example}

  \begin{lemma}\label{lem:partial-homog}
    If $T\colon{}\allpolypos\longrightarrow\allpolypos$ and we define a
    linear transformation $S(x^n) = y^n\bigl[ T(x^n)(x/y)\bigr]$ then
    $S\colon{}\allpolypos\longrightarrow\partialpoly{1,1}$.
  \end{lemma}
  \begin{proof}
    if $f(x) = \sum a_i x^i$ is in $\allpolypos$ then
    \begin{align*}
      S(f)(x,\alpha) &= \sum_{k=0}^n a_i \alpha^i \bigl[
      T(x^i)(x/\alpha)\bigr] \\
      &= T(f(\alpha x))(x/\alpha)
    \end{align*}
    The last expression is in $\allpolypos$ for positive
    $\alpha$. Finally, $S(f)^H$ is the homogenization of $T(x^n)$, and
    has all positive coefficients since $T$ maps to $\allpolypos$. 
  \end{proof}

\section{Applications to transformations}
\label{sec:partial-one-var}

The graph of a polynomial in $\partialpoly{1,1}$ goes up and to the left.
If we intersect it with an appropriate curve we will find sufficiently
many intersections that determine a polynomial in $\allpoly$.
Using these techniques we can establish mapping properties of
transformations. 

\begin{lemma} \label{lem:partial-x-x}
  If $f\in\partialpoly{1,1}$, $f(\alpha,y)\in\allpoly$ for $\alpha\le0$,
  and $f(0,y)\in\allpolypos(n)$ then $f(x,x)\in\allpolypos$.
\end{lemma}
\begin{proof}
  We will show that the graph of $f$ intersects the line $y=x$ in $n$
  points that have negative $x$ coordinate. We know that the graph of
  $f$ has asymptotes that have negative slope, so $f$ is eventually in
  the third quadrant. However, the assumption that
\index{solution curves} 
  $f(x,0)\in\allpolypos$ means that all solution curves begin on the
  negative $y$ axis. Consequently, all the curves must intersect the
  line $y=x$. See Figure~\ref{fig:partial-x-x}.

 \begin{figure}[htbp]
  \begin{center}
    \reflectbox{\includegraphics*[width=1.5in]{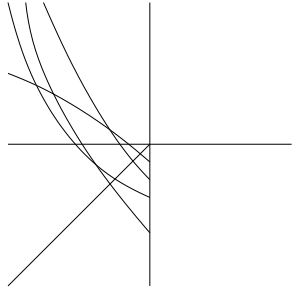}}
    \caption{The intersection of $f(x,y)=0$ and $y=x$. }
    \label{fig:partial-x-x}
  \end{center}
\end{figure}
\end{proof}

\begin{lemma} \label{lem:partial-T-xx}
  Suppose that $T\colon{}\allpolypos\longrightarrow\allpolypos$, where $T$
  preserves degree, and $T(x^i)$ has leading coefficient $c_i$. Choose
  $f\in\allpolypos(n)$, and define  $S(x^i) = T(x^i)x^{n-i}$. If $\sum
  c_ia_i x^i\in\allpolyint{(-\infty,-1)}$ where $f=\sum a_ix^i$ then
  $S(f)\in\allpolypos$. 
\end{lemma}
\begin{proof}
  Since $f\in\allpolypos$ we know that the homogenization $F$ of $f$
  is in $\partialpoly{1,1}$.  From Lemma~\ref{lem:partial-trans} the induced map
  $T_\ast$ satisfies $T_\ast(F)\in\partialpoly{1,1}$. By hypothesis, the
  asymptotes all have slope less than $-1$.  Each solution curve
  meets the $x$-axis in $(-\infty,0)$ since
  $(T_\ast(F))(x,0)=T(x^n)\in\allpolypos$. Consequently, the graph of
  $T_\ast F$ meets the line $x=-y$ in $n$ points, and these are the
  roots of $S(f)$. See Figure~\ref{fig:partial-T-xx}.

 \begin{figure}[htbp]
  \begin{center}
    \reflectbox{\includegraphics*[width=1.5in]{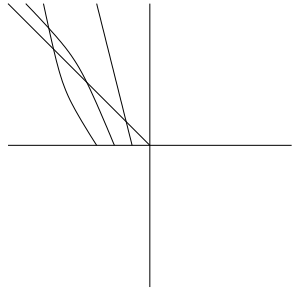}}
    \caption{The intersection of $T_\ast(F)$ and $y=-x$. }
    \label{fig:partial-T-xx}
  \end{center}
\end{figure}

\end{proof}

\begin{figure}[htbp]
  \begin{center}
    \leavevmode
    \reflectbox{\includegraphics*[width=1.5in]{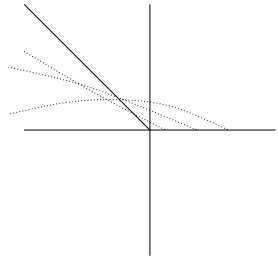}}
    \caption{The intersection of $T_\ast(F)$ and $y=-x$. }
    \label{fig:rising}
  \end{center}
\end{figure}

\begin{cor} \label{cor:half-ff}
  The linear transformation 
, \index{rising factorial} $V:x^i\mapsto
  \rising{x}{i}x^{n-i}$ maps $\allpolyint{(0,1)}$ to $\allpolypos$.
\index{falling factorial} Similarly, $U:x^i\mapsto
  \falling{x}{i}x^{n-i}$ maps $\allpolyint{(0,1)}(n)$ to $\allpolyalt$.
\end{cor}
\begin{proof}
  If we define $T(x^i) = \rising{x}{i}$ and $F$ is the homogenization of $f$
  then $T_\ast(F)(x,x) = V(f)$.
  If $f\in\allpolyint{(0,1)}(n)$ then $F\in\gsubpos_2(n)$ and so from
  Corollary~\ref{cor:partial-fall} we know that $T_\ast(F)\in\partialpoly{1,1}$. Since
  $(T_\ast(F))(x,0)=\rising{x}{n}$ the graph of $T_\ast(F)$ meets the $y$ axis
  below the line $x=-y$. Since the leading coefficient of $\rising{x}{i}$ is
  $1$ the asymptotes of $T_\ast(F)$ are the roots of $f$, which are
  between $0$ and $1$. Thus $T_\ast(F)$ is asymptotic to lines that
  lie below $x=-y$, and so every solution curve of $T_\ast(F)$ meets
  $x=-y$ above the $x$ axis.  See Figure~\ref{fig:rising}.
\end{proof}

\begin{lemma}
  If $T\colon{}\allpolypos\longrightarrow\allpoly$ preserves degree, $T(x^n)$
  is monic, $a>0$, and we define
  $V(x^n)=x^n\bigl[ T(x^n)(-a/x)\bigr]$
  then $V:\allpolyint{(-a,0)}\longrightarrow\allpoly$. 
\end{lemma}
\begin{proof}
  If we define $S(x^n) = y^n\,\bigl[(T\,x^n)(x/y)\bigr]$ then $V(f) =
  (Sf)(-a,x)$. Since $T(x^n)$ is monic we find that $(Sf)(x,0) =
  f(x)$, so the roots of $(Sf)(x,0)$ are exactly the roots of $f$.
  Consequently, the vertical line through $-a$ is to the left of the
  intersection of the graph of $Tf$ with the $x$-axis, and so meets
  all the solution curves since $Sf\in\partialpoly{1,1}$ by
  Lemma~\ref{lem:partial-homog}.
\end{proof}

  \begin{cor}\label{cor:falling-x}
    If $V:x^n\mapsto -(x-1)(2x-1)\cdots((n-1)x-1)$ then\\
$V:\allpolyint{(-1,0)}\longrightarrow\allpolyint{(0,1)}$.

The transformation $x^n\mapsto(x+1)(2x+1)\cdots((n-1)x+1)$ maps\\
$\allpolyint{(0,1)}\longrightarrow \allpolyint{(-1,0)}$.
  \end{cor}
  \begin{proof}
    We know $V:\allpolyint{(-1,0)}\longrightarrow\allpoly$ by applying
    the Lemma to $T\colon{}x^n\mapsto \rising{x}{n}$.  To see that $V$ has
    range $\allpolyint{(0,1)}$ we consider the intersection of the
    graph with the line $x=-y$. If $f =\sum a_i x^i$ and $S$ is as in
    the Lemma then \[S(x^n) = (x+y)(x+2y)\cdots(x+ny)\] so $(Sf)(x,-x)
    = a_0 + a_1x$. Since $f\in\allpolypos$, $a_0$ and $a_1$ are
    positive, and so the graph of $Vf$ does  meets the line $x=-y$ once
    in the upper left quadrant. 

    The homogeneous part of $S(f)$ is $x(x+y)(x+2y)\cdots(x+(n-1)y)$,
    so the graph of $S(f)$ has asymptotes with slopes
    $1,1/2,\cdots,1/(n-1)$, and one vertical asymptote. For positive
    $y$ the graph of $S(f)$ consists of $n$ curves starting at the
    roots of $f$. The vertical asymptote begins below $y=-x$ and so
    meets the line $y=-x$. Thus, the other curves lie entirely below
    $y=-x$, and consequently the line $x=-1$ meets these solution
    curves beneath the line $y=-x$, and hence all the roots of $V$ are
    in $(0,1)$.

    Replace $x$ by $-x$ for the second statement.
  \end{proof}

\begin{figure}[htbp]
    \centering
    \includegraphics*[width=2in]{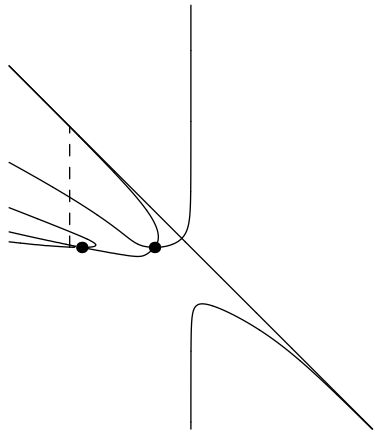}
    \caption{An example for Corollary~\ref{cor:falling-x}}
    \label{fig:partial-falling-x}
  \end{figure}

  For example, consider Figure~\ref{fig:partial-falling-x}, the graph
  of $Sf$, where $f(x)$ is the polynomial $(x+.3)^3(x+.9)^3$. The two
  dots are are at $(-.9,0)$ and $(-.3,0)$ which correspond to the
  roots of $f$. Except for the vertical asymptote, the graph of $Sf$
  lies under the line $x=-y$. The dashed line $x=-1$ meets the graph
  in points whose $y$ coordinates lie in $(0,1)$.

\section{\Mobius\ transformations}
\label{sec:partial-mobius}

We use properties of the graph of $T_\ast f(xy)$ to establish
properties of \Mobius\ transformations acting on linear
transformations.  We know from Corollary~\ref{cor:t1/z}
that if $T$ maps $\allpolypos$ to itself then so does $T_{1/z}$. It is
\index{T@$T_{1/z}$}%
more complicated when the domain is not $\allpolypos$.

\index{M\"{o}bius transformation} 

\begin{figure}[htbp]
  \begin{center}
    \leavevmode
    \epsfig{file=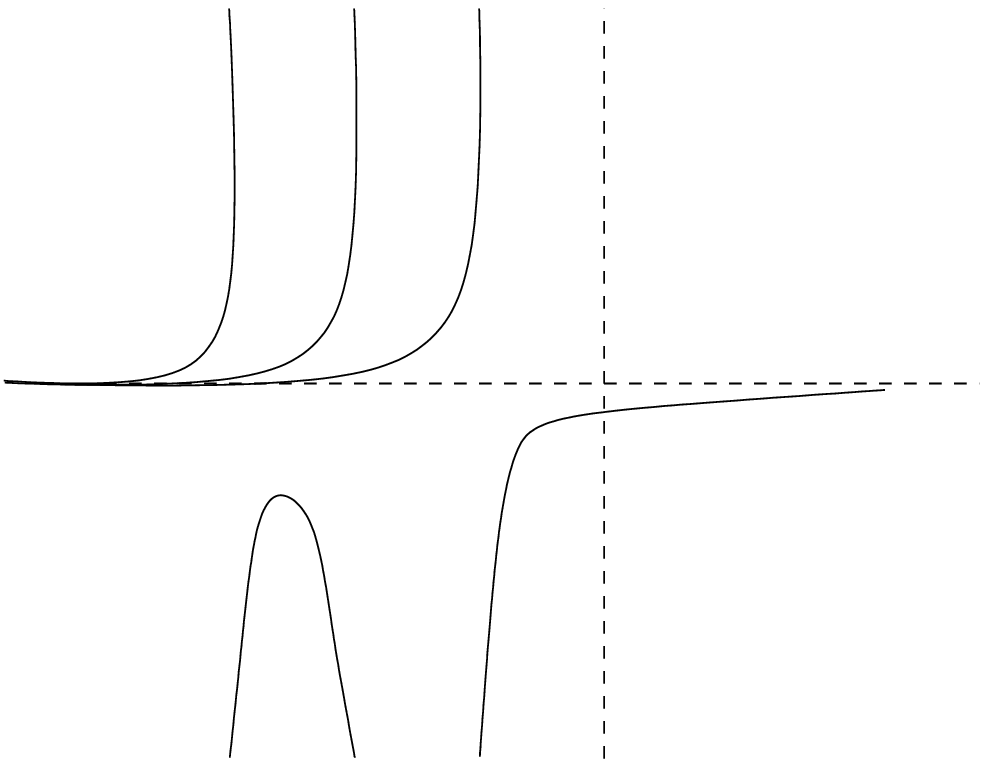,width=2.5in}
    \caption{A graph of $T_\ast(f(xy))$ where $f\in\allpolypos$}
    \label{fig:Tfxy}
  \end{center}
\end{figure}

The graph of $T_\ast f(xy)$ is different from the graphs of
Figure~\ref{fig:laguerre2} and Figure~\ref{fig:partial-x-x} since there
are vertical and horizontal asymptotes. Figure~\ref{fig:Tfxy} is the
graph of $T_\ast(f(xy))$ where \index{rising factorial}
$T(x^i)=\rising{x+1}{i}$ and $f\in\allpolypos$ has degree $3$.  

\begin{lemma} \label{lem:partial-fxy}
  Let $T\colon{}\allpolypos\longrightarrow\allpolypos$ where $T$ preserves
  degree, and $c_i$ is the leading coefficient of $T(x^i)$.  For any
  positive integer $n$ and $f=\sum a_ix^i$ in $\allpolypos(n)$
  \begin{enumerate}
  \item $(T_\ast f(xy))(x,\beta) \in\allpolypos$ for all $\beta>0$.
  \item $T_\ast f(xy)\in\partialpolyclose{1,1}$ 
  \item The vertical asymptotes of \/$T_\ast f(xy)$ occur at the roots 
    of $T(x^n)$.
  \item The horizontal asymptotes are of the form $xy=s$, where $s$ is 
    a root of $\sum a_ic_ix^i$.

  \end{enumerate}
\end{lemma}
\begin{proof}
  The first part follows from
  $$
  T_\ast(f(xy))(x,\beta) = \sum a_i T(x^i)\beta^i = \sum a_i
  T(\beta^ix^i) = T(f(\beta x)).$$
  
  As the coefficient of $y^n$ is $T(x^n)$ there are vertical
  asymptotes at the roots of $T(x^n)$.  We know that
  $f(xy)\in\partialpolypos{1,1}$ so
  $T_\ast f(x,y)\in\partialpolypos{1,1}$. The horizontal asymptotes
  are governed by the terms with equal $x$ and $y$ degrees. The
  polynomial so determined is $\sum a_i c_i x^i y^i$, which
  establishes the last part. Lemma~\ref{lem:signs-T} shows that all the $c_i$
  have the same sign, so $s$ is negative.
\end{proof}

\begin{figure}[htbp]
  \begin{center}
    \leavevmode
    \epsfig{file=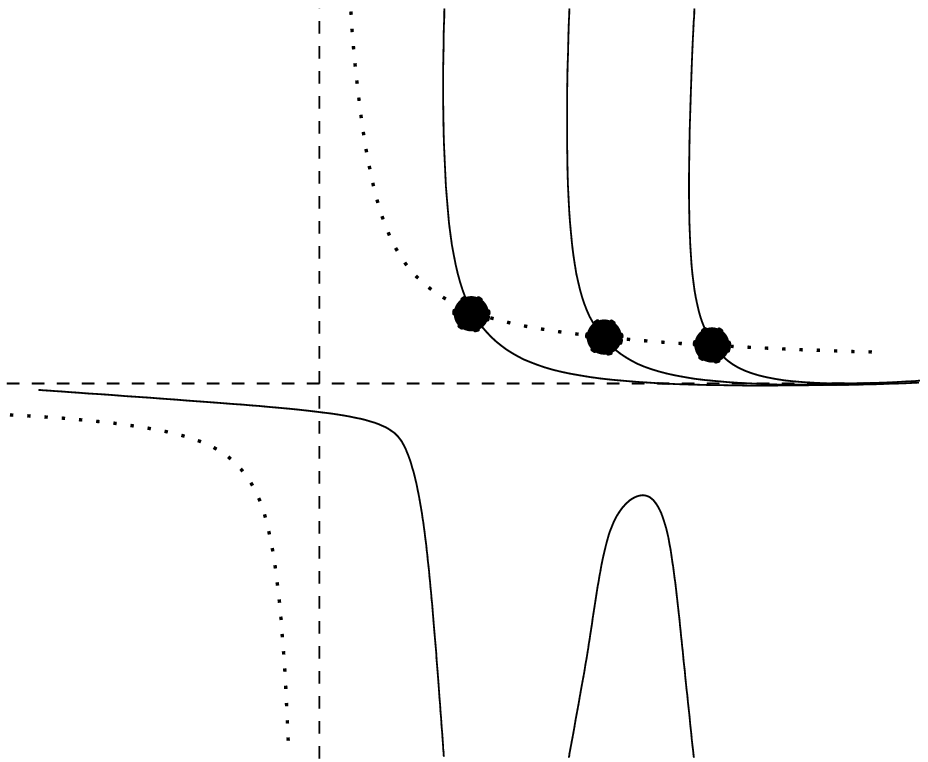,width=2.5in}
    \caption{The roots of $T_{1/z}(f)$ for $f\in\allpolyalt$}
    \label{fig:Tfxy2}
  \end{center}
\end{figure}

If we assume that $f\in\allpolyalt$ then Figure~\ref{fig:Tfxy} is
reflected around the $y$ axis. Figure~\ref{fig:Tfxy2} shows the
hyperbola and the intersection points with the graph that give the
roots in the next lemma. 

\begin{cor} \label{cor:partial-T-1/z2}
  Let $T\colon{}\allpolyalt\longrightarrow\allpolyalt$ where $T$ preserves
  degree, and $c_i$ is the leading coefficient of $T(x^i)$.  For any
  $f=\sum a_ix^i$ in $\allpolyalt$ such that $\sum c_i a_i x^i$ is in
  $\allpolyint{(0,1)}$ the polynomial $T_{1/z}(f)$ is in
\index{T@$T_{1/z}$}%
  $\allpolyalt$.
 \end{cor}

 \begin{proof}
   As usual let $T_\ast(x^iy^j) = T(x^i)y^j$. If we define
   $T(x^i)=p_i$ then 
   \[
   T_{1/z}(f) = \sum a_i\, x^i \,p_i(1/{x}) =
   \sum a_i \bigl[ T_\ast(x^iy^i)(1/x,x)\bigr] =
   T_\ast(f(xy))\,({1}/{x},{x})\]
   Thus, the roots of $T_{1/z}$ are
   found at the intersection of $T_\ast(f(xy))$ and the hyperbola
   $xy=1$. As long as the solution curves of $T_\ast(f(xy))$ go to
   zero more rapidly than $1/x$ then $xy=1$ will intersect every
   solution curve.  By Lemma~\ref{lem:partial-fxy} the horizontal asymptotes
   are of the form $xy=s$ where $s$ is a root of $\sum c_ia_ix^i$. By
   assumption $0<s<1$, and hence $T_{1/z}(f)$ has all real roots.
 \end{proof}
 
Similarly, we have
\begin{cor} \label{cor:partial-T-1/z}
  Let $T\colon{}\allpolypos\longrightarrow\allpolyneg$ where $T$ preserves
  degree, and $c_i$ is the leading coefficient of $T(x^i)$.  For any
  $f=\sum a_ix^i$ in $\allpolypos$ such that $\sum c_i a_i x^i$ is in
  $\allpolyint{(-1,0)}$ the polynomial $T_{-1/z}(f)$ is in
  $\allpolypos$.
 \end{cor}


Here is another consequence of changing signs.

\begin{cor} \label{cor:partial-T1/z2}
  Let $T\colon{}\allpolypos\longrightarrow\allpolyalt$ where $T$ preserves
  degree, and $c_i$ is the leading coefficient of $T(x^i)$.  For any
  $f=\sum a_ix^i$ in $\allpolyalt$ such that $\sum c_i a_i (-x)^i$ is in
  $\allpolyint{(0,1)}$ the polynomial $T_{1/z}(f)$ is in
\index{T@$T_{1/z}$}%
  $\allpolyalt$.
 \end{cor}
 \begin{proof}
   If we let $S(x) = T(-x)$ then $S$ satisfies the hypothesis of
   Corollary~\ref{cor:partial-T-1/z}. The conclusion follows since $S_{-1/z} = T_{1/z}$.
 \end{proof}

 We can apply Corollary~\ref{cor:partial-T-1/z2} to the monic Laguerre
 and Charlier polynomials (Corollary~\ref{cor:laguerre-x} and
 Corollary~\ref{cor:charlier-x}) since all the $c_i$ are $1$.

 \begin{cor} \label{cor:mobius-laguerre} If $\tilde{L}_n^{(\alpha)}$
   is the monic Laguerre polynomial then the transformation $x^i
   \mapsto \rev{\tilde{L}_i(\alpha;x)}$ maps $\allpolyint{(0,1)}$ to
   $\allpoly$.
 \end{cor}

 \begin{cor} \label{cor:mobius-charlier}
   If $C_n^\alpha$ is the Charlier polynomial then the transformation
   $x^i \mapsto \rev{C_i^\alpha}$ maps $\allpolyint{(0,1)}$ to
   $\allpoly$.
 \end{cor}

In the falling  factorial the constant terms of the factors
changed. Now we let the coefficient of $x$ change. Consider the linear
transformation 

\begin{equation}
  \label{eqn:factorial-variant}
  T(x^n) = (1-x)(1-2x)(1-3x)\cdots(1-nx)
\end{equation}

\begin{lemma} \label{lem:factorial-variant}
  The linear transformation \eqref{eqn:factorial-variant} maps
  $\allpolyint{(0,1)}$ to $\allpoly$.
\end{lemma}
\begin{proof}
  If $S(x^n) = \falling{x-1}{n}$ then $T= S_{1/z}$. The result now
  follows from Corollary~\ref{cor:falling-x}.
\end{proof}

Next we consider the \index{\Mobius\ transformation}\Mobius\
transformation $z\longrightarrow (1-z)$.

\begin{lemma} \label{lem:homog-xy-1}
Suppose that $T\colon{}\allpolypos\longrightarrow\allpoly$, and let
$f(x)=\sum a_ix^i\in\allpolypos(n)$. Define $c_i$ to be the leading
coefficient of $T(x^i)$. We can conclude that 
$\sum a_i (1-x)^{n-i} T(x^i)\in\allpoly$ if any of the following hold:
\begin{enumerate}
\item $T(x^n)\in\allpolyint{(-\infty,1)}$ and $\sum
  a_ic_ix^i\in\allpolyint{(-\infty,-1)}$.
\item $T(x^n)\in\allpolyint{(1,\infty)}$ and $\sum
  a_ic_ix^i\in\allpolyint{(-1,0)}$.
\item $T(x^n)\in\allpolyint{(-\infty,1)}$ and $\sum
  a_ic_ix^i\in\allpolyalt$.
\end{enumerate}
\end{lemma}
  
\begin{figure}[htbp]
  \begin{center}
    \leavevmode
    \epsfig{file=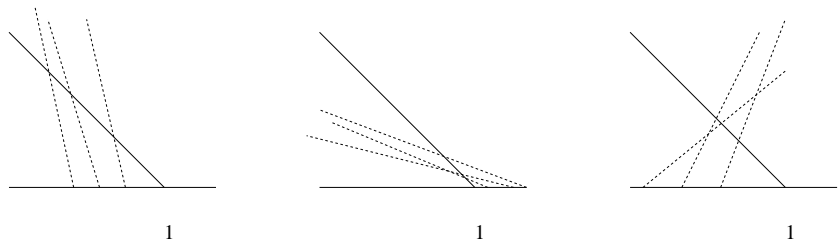,width=4.5in}
    \caption{The geometry of the solution curves of Lemma~\ref{lem:homog-xy-1} }
    \label{fig:three-xy}
  \end{center}
\end{figure}

\begin{proof}
  \index{solution curves} We need to show that each solution curve of
  $T_\ast(f)$ meets the line $x+y=1$. Since
  $T\colon{}\allpolypos\longrightarrow\allpoly$ there are solution curves for
  all positive $y$. The asymptotes of $T_\ast(f)$ are the roots of
  $\sum a_ic_ix^i$, and $T_\ast(f)(x,0)=T(x^n)$ so the geometry of the
  solution curves is given in Figure~\ref{fig:three-xy}. It is clear
  that each solution curve must meet the line $x+y=1$.
\end{proof}

\section{Partial substitution in several variables}
\label{sec:partial-d}
  
  We now generalize $\partialpoly{1,1}$ to more variables.  In short,
  we consider polynomials in two kinds of variables such that if we
  substitute arbitrary values for all but one of the $x$ variables,
  and only non-negative values for all the $y$ variables, then the
  resulting polynomial is in $\allpoly$. In more detail,

\begin{definition} \label{defn:partialpoly-d}
  The class $\partialpoly{d,e}$ consists of all polynomials
  $f(\xx;\yy)$ where $\xx=(x_1,\dots,x_d)$ and $\yy=(y_1\dots,y_e)$
  satisfying
  \begin{enumerate}
  \item The total degree equals the degree of each variable.
  \item The homogeneous part $f^H$ is in $\gsubpos_{d+e}$.
  \item $f(\xx;\beta) \in\rupint{d}$ where all
    $\beta=(\beta_1,\dots,\beta_e)$ are non-negative.
  \end{enumerate}
  We define $\partialpolyclose{d,e}$ to be the set of all polynomials
  that are the limit of polynomials in $\partialpoly{d,e}$.  The set
  of all polynomials in $\partialpoly{d,e}$ with all positive
  coefficients is $\partialpolypos{d,e}$.
 \end{definition}

 Notice that $\partialpolypos{d,0}$ consists of polynomials in
 $\rupint{d}$ with all non-negative coefficients, and is exactly
 $\gsubplus_{d+e}$.  We have the inclusions
 
 $$
 \gsubplus_{d+e} = \partialpolypos{d+e,0} \subseteq
 \partialpolypos{d+e-1,1}\subseteq \cdots \subseteq
 \partialpolypos{1,d+e-1}$$.

If $d,e>1$ then substitution of positive values determines two maps

\begin{align*}
  \partialpolypos{d,e} & \longrightarrow \partialpolypos{d-1,e} \text{
    (substituting for an $x$ value)}\\
  \partialpolypos{d,e} & \longrightarrow \partialpolypos{d,e-1} \text{
    (substituting for a $y$ value)}
  \end{align*}
  
  $\partialpoly{d,e}$ and $\partialpolypos{d,e}$ are closed under
  multiplication. We can not define interlacing in
  $\partialpolypos{d,e}$ as closed under all linear combinations,
  since we don't allow negative coefficients. If
  $f,g\in\partialpolypos{d,e}$ we define $f\plesslesseq g$ to mean that
  $f+\alpha g\in\partialpolypos{d,e}$ for all non-negative $\alpha$. For
  example, if $\ell$ is linear and in $\partialpolypos{d,e}$, and
  $f\in\partialpolypos{d,e}$, then $\ell f \plesslesseq f$.
\index{positive interlacing}

A linear polynomial $c+\sum a_i x_i + \sum b_i y_i$ is in
$\partialpoly{d,e}$ if and only if it is in $\rupint{d+e}$. The
quadratic case is much more interesting, and is considered in
Section~\ref{sec:partial-quad}.

It is easy to construct polynomials in $\partialpoly{d,e}$ by
generalizing the inductive construction of Lemma~\ref{lem:3-term-partial}.
The proof is straightforward and omitted.

\begin{lemma}\label{lem:3-term-partial-de}
  Suppose that $\aaa_i,\bbb_i,\ddd_i,\eee_i$ are vectors of positive
  constants for $i\ge0$, and that $\ccc_i$ can be be positive or
  negative. If
\begin{gather*}
    f_n(\xx;\yy) = (\aaa_n\cdot \xx + \bbb_n\cdot\yy +
    \ccc_n)f_{n-1}(\xx;\yy) - (\ddd_n\cdot \yy +
    \eee_n)f_{n-2}(\xx;\yy) \\
f_0(\xx;\yy) = 1 \quad f_1(\xx;\yy) = \aaa_1\cdot\xx + \bbb_1\cdot \yy + \ccc_1
\end{gather*}
then $f_n(\xx;\yy)\in\partialpoly{d,e}$ for $n=0,1,\dots$.
\end{lemma}

\section{Quadratic forms}
\label{sec:partial-quad}

Our goal in this section is to determine conditions for the
\index{quadratic form}quadratic form $ \begin{pmatrix}\xx &
  \yy\end{pmatrix} Q
\begin{pmatrix}\xx & \yy\end{pmatrix}^t$ to belong to
$\partialpoly{d,e}$. We first consider the case when $d=1$, and then
when $e=1$.

A matrix $Q$ is called \emph{copositive} if $xQx^t\ge0$ for every
vector $x$ with all non-negative entries. Unlike positive definite
matrices, there is no good characterization of copositive matrices.

  \begin{lemma}
    Suppose that $Q = \smalltwo{a}{v}{v^t}{C}$, where $C$ is symmetric
    and $a$ is positive. The \index{quadratic form}quadratic form $f(x;\yy) =
    (x,\yy)Q(x,\yy)^t$ has all real roots for all non-negative $\yy$
    if the matrix $v^tv-aC$ is copositive.
  \end{lemma}
  \begin{proof}
    Expanding the \index{quadratic form}quadratic form yields
  \begin{align*}
    (x,\yy)\begin{pmatrix} a& v\\ v^t & C \end{pmatrix}
    \begin{pmatrix} x \\ \yy^t \end{pmatrix} & = 
a x^2  + x(\yy v^t + v\yy^t) + \yy C\yy^t 
  \end{align*}
Since $\yy v^t = v\yy^t$, the discriminant of the equation is $\yy(v^tv-aC)\yy^t$.

  \end{proof}

\begin{lemma}\label{lem:quad-nsd}
Suppose that $Q=\begin{pmatrix}  A &  u \\ u^t  & c \end{pmatrix}$
where all entries are positive, $A$  is a $d$ by $d$ symmetric matrix,
and $c$ is a scalar. The \index{quadratic form}quadratic form
$ f(\xx;y)=\begin{pmatrix}\xx & y\end{pmatrix} Q
\begin{pmatrix}\xx &  y\end{pmatrix}^t$
is in $\partialpoly{d,1}$ iff $Q$ is \nsd.
\end{lemma}
\begin{proof}
  If $Q$ is \nsd\ then $f(\xx,y)\in\rupint{d+1}$, and since all terms
  are positive it is in $\partialpoly{d,1}$.
  
  For the converse, we only need to show that
  $f(\xx;y)\in\rupint{d+1}$, since this implies that $Q$ is \nsd. Since
  all entries of $Q$ are positive, $f(\xx,y)$ satisfies the degree
  condition. Since $f(\xx,y)^H=f(\xx,y)$, it suffices to show that
  $f(\xx,y)$ satisfies substitution for any choice of values for $\xx$
  and $y$. By hypothesis, $f(\xx,1)\in\rupint{d}$. Now since
  $f(\xx;z) = z^2 f(\frac{\xx}{z},1)$ and $f(\frac{\xx}{z};1)$
  satisfies substitution,  it follows that $f(\xx;z)$
  satisfies substitution.
\end{proof}

  \begin{prop}\label{prop:quad-nsd}
Suppose that $Q=\begin{pmatrix}  A &  u \\ u^t  & C \end{pmatrix}$
where all entries are positive, $A$  is a $d$ by $d$ symmetric matrix,
and $C$ is an $e$ by $e$ symmetric matrix. The \index{quadratic form}quadratic form
$$ \begin{pmatrix}\xx & \yy\end{pmatrix} Q
\begin{pmatrix}\xx^t \\ \yy^t\end{pmatrix}$$
is in $\partialpoly{d,e}$ iff the two conditions are met:
\begin{enumerate}
\item $A$ is \nsd.
\item $u A^{-1} u^t - C$ is copositive.
\end{enumerate}    
  \end{prop}
  \begin{proof}
If the \index{quadratic form}quadratic form is in $\partialpoly{d,e}$ then substituting
$\yy=0$ shows that $\xx A\xx^t\in\rupint{d}$, so $A$ is negative
subdefinite. 

If we substitute non-negative values $\beta$ for $\yy$, expand and
regroup the \index{quadratic form}quadratic form, we get
\begin{align}
\begin{pmatrix}\xx & \beta\end{pmatrix} 
\begin{pmatrix}
  A &  u \\ u^t  & C
\end{pmatrix}
\begin{pmatrix}\xx \\ \beta\end{pmatrix}
&=
\xx A \xx^t + \beta u^t \xx^t + \xx u \beta^t + \beta C \beta^t \notag\\
&=
\begin{pmatrix}\xx & 1 \end{pmatrix} 
\begin{pmatrix}
  A &  u\beta^t \\ \beta u^t & \beta C\beta^t
\end{pmatrix}
\begin{pmatrix}\xx \\ 1\end{pmatrix} \label{eqn:qf-proof}
\end{align}
By Lemma~\ref{lem:quad-nsd}, $f(\xx,z)$ is in $\rupint{d+1}$ iff 
$\smalltwo{A}{u\beta^t}{\beta u^t}{\beta C\beta^t}$ is \nsd. By
Lemma~\ref{lem:extend-nsd} this is the case iff 

$$ 0 \le (u\beta^t) (A^{-1})(\beta u^t) - ( \beta C \beta^t) = \beta(
uA^{-1}u^t-C)\beta^t$$
 which finishes the proof.

  \end{proof}

\begin{remark}
    Here is an alternative derivation in the case that $A$ is $2$ by
    $2$. We only need to determine when the \index{quadratic form}quadratic form satisfies
    substitution. This will be the case iff the discriminant, with
    respect to the variable $x_1$, is positive. This discriminant is a
    \index{quadratic form}quadratic form, and its discriminant, with respect to $x_2$,
    equals
    $$
    a\,|A|\, \yy  (uA^{-1}u^t - C)\yy$$
    The condition that this is
    non-positive is equivalent to the positivity condition since $A$
    is \nsd, and therefore $|A|$ is negative.
  \end{remark}

\begin{remark}
  Let's work out an example in light of this theory. We start
  with
  \begin{align*}
    g(x;y) &= x^2+2xy+y^2 + 4x+3y+2 \\
\intertext{and the homogenized form is}
    G(x;y,z) &= x^2+2xy+y^2 + 4xz+3yz+2z^2 \\
    &= (x,y,z)\left(\begin{smallmatrix}1&1&2\\1&1 & 3/2
        \\2&3/2&2\end{smallmatrix}\right)(x,y,z)^t\\
\intertext{To show that this is in $\partialpoly{1,2}$ we note $A  =
  (1)$, $u = (1,2)$, $C = \smalltwo{1}{3/2}{3/2}{2}$}
uA^{-1}u^t-C &= \smalltwo{1}{2}{2}{4}-\smalltwo{1}{3/2}{3/2}{2}\,=\, \smalltwo{0}{1/2}{1/2}{2}
  \end{align*}
which  is copositive. Substituting $z=1$ shows
$g(x;y)\in\partialpoly{1,1}$. 
\end{remark}

\begin{remark}
  A \index{quadratic form}quadratic form in $\partialpoly{1,1}$ is in $\rupint{2}$. To see
  this, let $Q=\smalltwo{a}{b}{b}{c}$, and assume
  $(x\,y)Q(x\,y)^t\in\partialpoly{1,1}$. Applying Proposition~\ref{prop:quad-nsd} we
  see that $b^2/a-c>0$ or equivalently $\smalltwodet{a}{b}{b}{c}<0$,
  and hence $Q$ is \nsd.  The example in the previous remark is a
  quadratic polynomial in $\partialpoly{1,1}\setminus\rupint{2}$, but it
  comes from a quadratic from  three variables.
\end{remark}

\begin{lemma}
  $x_1y_1+\cdots x_dy_d\in\partialpolyclose{d,d}$
\end{lemma}
\begin{proof}
  We will give the computations for $d=2$, and describe the general
  construction. If we define
  \begin{align*}
    Q_\epsilon &=
    \begin{pmatrix}
      \epsilon^3 & \epsilon & 1 & \epsilon \\ \epsilon & \epsilon^3 & \epsilon & 1 \\
      1 & \epsilon & \alpha & \alpha \\
      \epsilon & 1 & \alpha & \alpha
    \end{pmatrix}
    \intertext{then $(x_1,x_2,y_1,y_2)Q_\epsilon(x_1,x_2,y_1,y_2)^t$
      converges to $x_1y_1+x_2y_2$.  It suffices to show that
      $u_\epsilon A_\epsilon u_\epsilon$ has positive entries, for we
      just choose $\alpha$ small enough so that $u_\epsilon A_\epsilon
      u_\epsilon-C_\alpha$ is positive.}  u_\epsilon
    A_e^{-1}u_\epsilon &=
\begin{pmatrix}
  1&\epsilon \\\epsilon & 1
\end{pmatrix}
\begin{pmatrix} \epsilon^3 & \epsilon \\ \epsilon & \epsilon^3
\end{pmatrix}^{-1}
\begin{pmatrix}
  1&\epsilon \\ \epsilon & 1
\end{pmatrix} \\
&=\frac{1}{\epsilon(1+\epsilon)(1+\epsilon^2)} 
\begin{pmatrix}
  2+\epsilon+\epsilon^2 & 1+\epsilon+2\epsilon^2 \\
1+\epsilon+2\epsilon^2 &  2+\epsilon+\epsilon^2
\end{pmatrix}
  \end{align*}
In general, we define $J(a,b)$ to be the $d$ by $d$ matrix whose
diagonal is $a$, and remaining elements are $b$. We let 
$A_\epsilon = J(\epsilon,\epsilon^3)$, $U_\epsilon=J(1,\epsilon)$ and
$C_\epsilon=J(\alpha,\alpha)$. With some work it can be verified that
$U_\epsilon A_e^{-1}U_\epsilon$ has all positive entries.
\end{proof}

Note that this gives an element of $\partialpolyclose{2,2}$ that is
not in $\gsubclose_4$. If it were, we could substitute $x_2=y_2=1$, but
$x_1y_1+1$ is not in $\gsubclose_2$.

\begin{remark}
  We can construct polynomials of higher degree in $\partialpoly{1,1}$
  by multiplying quadratics. For instance, $\left(\begin{smallmatrix}
      1&1&1\\1&1&\epsilon\\1&\epsilon&c\end{smallmatrix}\right)$
  satisfies the conditions of Proposition~\ref{prop:quad-nsd} for
  $0<c<1$. If we choose $0<c_i<1$ and let
  $\epsilon\longrightarrow0^+$, then the following polynomial is in
  $\partialpolyclose{1,1}$:
$$ \prod_i(x^2+2xy+y^2+2x+c_i)$$
\end{remark}


\chapter{The analytic closure of $\gsubpos_d$}
\label{cha:new-analytic}
 
\renewcommand{\TimeStampStart}{Thursday, January 17, 2008: 19:23:27}
\mytoday

$\allpolyf_d$ is the closure of $\gsubpos_d$ under uniform convergence
on compact domains. The uniform closure of $\gsubplus_d$ is
$\allpolyposf_d$.  Many results about these closures follow by taking
limits of results about polynomials. For instance, we can conclude
from Lemma~\ref{lem:p2-xx} that if $f(x,y)\in\allpolyf_2$ then
$f(x,x)\in\allpolyf_2$. Similarly, we can also extend definitions from
polynomials to these closures. As an example, we say that
$f,g\in\gsubf_d$ \emph{interlace} if $f+\alpha g\in\gsubf_d$ for
all real $\alpha$.

\section{The analytic closure}
\label{sec:analytic-closure}

The properties of $\allpolyf_d$ are similar to those of $\allpolyf$,
except that we do not have a characterization of the members as we do
in one variable. To start off,

\begin{lemma}
$\allpolyf_d$ consists of analytic functions whose domain is $\complexes^d$. 
\end{lemma}
\begin{proof}
  If $f(\xx)$ is in $\allpolyf_d$ then it is a limit of polynomials
  $f_n(\xx)$ in $\gsubpos_d$. If we choose $\aaa\in\reals^d$ then
  $f(\xx^\aaa_i)$ is the uniform limit of polynomials
  $f_n(\xx^\aaa_i)$ in $\allpoly$ and hence $f(\xx^\aaa_i)$ is
  analytic. Consequently, $f$ is analytic in each variable separately,
  and so we can apply Hartog's Theorem \cite{hormander} to conclude that
  $f(\xx)$ is analytic.
\index{Hartog's Theorem}
\end{proof}

From Section~\ref{sec:sub-closure} we know that
$\allpolyf_e\subset\allpolyf_d$ if $e<d$, and $\allpolyposf\subset
\gsubposf_2$. The most important non-polynomial function in
$\allpolyf$ is $e^x$. Thus $e^x$ is in all $\gsubf_d$. There is
another exponential function that is in $\gsubf_d$.

\begin{lemma} \label{lem:exy}
  If $d\ge2$ then $e^{-xy}$ is in $\gsubf_d$. Also,
  $e^{xy}\not\in\allpolyf_d$. 
\end{lemma}
\begin{proof}
From Corollary~\ref{cor:sub-xyc} we know that $n - xy$ is in $\gsubposclose_2$.
  Writing   $e^{-xy}$ as a limit 
\begin{equation} \label{eqn:exy}
 e^{-xy} = \lim_{n\rightarrow\infty} n^{-n} \,(n - xy)^n
\end{equation}
it follows that all the factors of the product are in
$\gsubposclose_2,$ and hence $e^{-xy}$ is in $\gsubf_2$.

If $e^{xy}\in\gsubf_2$ then substituting $x=y$ yields
$e^{x^2}\in\allpolyf$ which is false. 
\end{proof}

 Since $f(x,y)=e^{-xy}$
is in $\gsubposf_2$ we can apply Lemma~\ref{lem:sub-fxx} to find that $f(x,x)
= e^{-x^2}$ is in $\allpolyf,$ as we saw earlier.  Since $\allpolyf
\subset \gsubposf_d$ we see $e^{-x^2}$ is in $\gsubposf_d$. We
can combine these members of $\gsubposf_d$ to conclude

\begin{prop} \label{prop:e-xhx}
  If $H$ is a $d$ by $d$ matrix with non-negative entries then
  $e^{-\xx H\xx^t}$ is in $\gsubposf_d$. 
\end{prop}
\begin{proof}
  It suffices to show that $e^{-\sum a_{ij}x_ix_j}$ is in
  $\gsubposf_d$ for all non-negative $a_{ij}$. Since
$$ e^{-\sum a_{ij}x_ix_j} = \prod e^{-a_{ij}x_ix_j}$$
and all the factors of the product are in $\gsubposf_d$ it follows
that the product itself is in $\gsubposf_d,$ and so is the left
hand side.
\end{proof}

The non-negativity of the coefficients of $H$ is important. For
instance, if $H=(-1)$ then $e^{-xHx} = e^{x^2}$ and this is not in
$\allpolyf$  by \eqref{eqn:type-2}.

\begin{remark} 
  We can easily construct functions in $\allpolyf_2$.  If $\sum |a_i|
  < \infty$ then the infinite product $\prod_{i=1}^\infty(1+a_i z)$ is
  an entire function. Suppose we choose sequences of positive real
  numbers $\{a_i\},\{b_i\}$ so that $\sum a_i$ and $\sum b_i$ are
  finite.  The product
  $$f(x,y)=\displaystyle \prod_1^\infty (1+a_i x+ b_i y)$$
  is an
  entire function since
$$ \sum |a_i x + b_i y| \le |x|\,\sum |a_i| + |y|\,\sum|b_i| <
\infty$$
Since $f$ is the limit of the partial products which are in
$\gsubpos_2$, $f(x,y)$ is in $\allpolyf_2$.
\end{remark}

\index{sin and cos}
\begin{example} \label{ex:sincos}
 We will show that  $\sin(x) + y\cos(x)\in\allpolyf_2$,
    which is an example of a member of $\allpolyf_2$ that is not a
    product.
    We start with the two product approximations to $\sin$ and $\cos$:

\begin{align*}
  f_n(x) &= x\,\prod_{k=1}^n
\left(1 - \frac{x^2}{k^2\pi^2}\right) \\
g_n(x) & = \prod_{k=0}^n
\left(1 - \frac{4x^2}{(2k+1)^2\pi^2}\right) 
\end{align*}
\noindent%
and notice that $f_{n+1}\lesslesseq g_n$. Lemma~\ref{lem:fyg} implies
that $f_{n+1}(x) + yg_n(x)\in\gsubclose_2$. Taking the limit, we find
that $\sin(x) + y \cos(x)\in\allpolyf_2$.
  \end{example}
  
  The same limit argument used in Lemma~\ref{lem:exy} can be applied
  to Corollary~\ref{cor:fpp}:

\begin{lemma} \label{lem:efpp}
  If $f\in\allpolyposf$ and $g\in\gsubposclose_d$ then 
$f(-\partial_x\,\partial_y)g\in\gsubposclose_d$.
\end{lemma}

This has a very important corollary:

\begin{cor}\label{cor:p2-exy}
  If $f\in\gsubclose_d$ then $e^{-\partial x_1\,\partial x_2}\, f
  \in\gsubclose_d$. More generally,\\
  $e^{-\partial\xx\cdot\partial\yy}:\gsubclose_{2d}\longrightarrow\gsubclose_{2d}$. 
\end{cor}
\begin{proof}
  Since $e^{-xy}\in\gsubf_2$ we can apply limits to Lemma~\ref{lem:efpp}. The
  second part follows from
\[
e^{-\partial\xx\cdot\partial\yy} = e^{-\partial\xx_1\partial\yy_1}
\cdots e^{-\partial\xx_d\partial\yy_d} 
\]
\end{proof}

\begin{cor}
  The map $x^k\mapsto L_k(xy)$ maps $\allpolypos\longrightarrow\gsubclose_2$.
\end{cor}
\begin{proof}
  From \eqref{eqn:laguerre-quad-diff} the following diagram commutes

\centerline{\xymatrix{
x^n \ar@{->}[d]_{exp} \ar@{.>}[rr] && L_n(xy)  \ar@{<-}[d]^{e^{-\partial_x\partial_y}}\\
\frac{x^n}{n!} \ar@{->}[rr]_{x\mapsto -xy} && \frac{x^n(-y)^n}{n!}
}}

and consequently this commutes:

\centerline{\xymatrix{
    \allpolypos \ar@{->}[d]_{exp} \ar@{.>}[rrr]^{x^n\mapsto L_n(xy)}
    &&& \gsubclose_2  \ar@{<-}[d]^{e^{-\partial_x\partial_y}}\\ 
    \allpolypos \ar@{->}[rrr]_{x\mapsto -xy} &&& \gsubclose_2 }}
\end{proof}

Lemma~\ref{lem:efpp} has another important consequence:

\index{\ diag}
\index{\ diag$_1$}

\index{diagonal polynomials}
\begin{theorem} \label{thm:diagonal}
  If $f=\sum a_{ij}x^iy^j$ in $\gsubplus_2$ then  the diagonal
  polynomials $\diag(f) =  \sum a_{ii}x^i \quad$ and
  $\quad \diag_1(f)=\sum a_{ii}\,i!\,x^i$ are in
  $\allpolypos$.
\end{theorem}
\begin{proof}
  Since $f$ is in $\gsubplus_2(n)$ the homogenization $F = \sum
  a_{ij}x^iy^jz^{n-i-j}$ is in $\gsubpos_3$. The linear transformation
  $g\mapsto e^{-\partial_x \partial_y} g$ maps $\gsubposclose_2$ to
  itself, and by Theorem~\ref{thm:pd-T} it extends to a linear
  transformation mapping $\gsubposclose_3$ to itself. We compute
  \begin{align}
     e^{-\partial_x \partial_y} F &=  e^{-\partial_x \partial_y} \sum_{i,j}
     a_{ij}x^iy^jz^{n-i-j} \notag \\
&= \sum_{i,j,k} a_{ij} \frac{(-1)^k}{k!} \left(\frac{\partial}{\partial
  x}\right)^k \left(\frac{\partial}{\partial y}\right)^k x^i y^j z^{n-i-j}
 \notag \\
&= \sum_{i,j,k} a_{ij} \frac{(-1)^k}{k!} \falling{i}{k}  \falling{j}{k}   x^{i-k} y^{j-k} z^{n-i-j}
\label{eqn:diag-1}\\
\intertext{Substituting $x=y=0$ gives a polynomial in $\allpoly$}
     e^{-\partial_x \partial_y} F \big\vert_{x=y=0} & =
 \sum_k a_{kk} (-1)^k k! z^{n-2k} \notag
\end{align}
The reverse of this last polynomial is in $\allpoly$, so $\sum a_{kk}
(-1)^k k! z^{2k}$ is in $\allpoly$. Now all $a_{kk}$ are positive
since $f\in\gsubplus_2$, and so $\sum a_{kk}(-1)^k k! x^k$ is in
$\allpoly$. Replacing $x$ by $-x$, and an application of the
exponential map shows that $\sum a_{ii} x^i$ is also in $\allpoly$.
\end{proof}

A similar proof establishes a similar result in higher dimensions, but
we are not able to eliminate the negative factors.

\begin{lemma} \label{lem:diagonal-signed}
  If $f(\xx,y,z) = \sum f_{ij}(\xx) y^iz^j$ is in $\gsubpos_{d+2}(n)$
  then\\ $\sum (-1)^k k!\,f_{kk}(\xx)$ and $\sum (-1)^k\,f_{kk}(\xx)$
  are in $\gsubpos_d$.
\end{lemma}
\begin{proof}
  We follow the proof above and apply the operator
  $e^{-\partial_y\partial_z}$, but we do not need to homogenize. For
  the second part, apply the first part to $exp(f)$ where  the
  exponential map operates on the $y$ variable.
\end{proof}

  \begin{remark}
    Here is an alternative proof that the Hadamard product maps
    $\allpolypos\times\allpolypos\longrightarrow\allpolypos$. It does
    not show that the domain could be enlarged to
    $\allpoly\times\allpolypos$.  

    If $f,g\in\allpolypos$ then
    $f(x)g(y)\in\gsubplusclose_2$. Applying Theorem~\ref{thm:diagonal}
    shows that the diagonal of $f(x)g(y)$ is in $\allpolypos$. But the
    diagonal is  Hadamard product: if $f(x) = \sum a_i x^i$, $g(y) =
    \sum b_iy^i$ then the diagonal of $\sum a_i b_j x^iy^j$ is $\sum
    a_ib_ix^i$. 
  \end{remark}
\index{Hadamard product!in $\rupint{d}$}

\begin{cor}
  If $\sum f_i(\xx)y^i$ and $\sum g_i(\xx)y^i$ are both in
  $\gsubpos_{d+1}$ then 
$$ \sum (-1)^if_i(\xx)g_i(\xx) \in\gsubpos_d $$
\end{cor}
\begin{proof}
  The product $\left(\sum f_i(\xx)y^i\right)\,\left(\sum
    g_i(\xx)z^i\right)$ is in $\gsubposclose_{d+2}$. Now apply
  Lemma~\ref{lem:diagonal-signed}.
\end{proof}

This is false without the factor $(-1)^i$. Just take $f=g$; the sum
$\sum f_\sdiffi(\xx)g_\sdiffi(\xx)$ is positive if the constant term
of $f$ is non-zero.  We note some similar results:

\begin{cor}
  If $f = \sum a_{ij}x^iy^j$ is in $\gsubplus_2$ then for any
  non-negative integers $r,s$  the polynomial 
$$ \sum_i a_{i+r,i+s} \frac{(i+r)!(j+s)!}{r!s!i!}x^i \quad\in\allpoly$$
\end{cor}
\begin{proof}
  Instead of substituting $x=y=0$ in \eqref{eqn:diag-1} take the
  coefficient of $x^ry^s$ and apply a similar argument.
\end{proof}

\begin{cor}
  If $f=\sum a_{ij}x^iy^j$ is in $\gsubplus_2(n)$ then define $$s_k(x)
  = \sum_{i-j=k} a_{i,j} x^i$$
  For $-n<k<n$ we have that
  $s_k\in\allpoly$, and $s_k,s_{k-1}$ have a common interlacing.  The
  map $x^k\mapsto s_k$ sends $\allpolypos(n)$ to itself.
\end{cor}

\begin{proof}
  For any positive $\alpha$ the polynomial $g=(x^k+\alpha x^{k-1})f$
  is in the closure of $\gsubplus_2$, and $\diag g = s_k + \alpha
  s_{k-1}$. By Proposition~\ref{prop:pattern} $s_k$ and $s_{k-1}$ have a common
  interlacing. If $h = \sum b_k y^k$ is in $\allpolypos(n)$ then
  $h(y)$ is in the closure of $\gsubplus_2(n)$, and the following
  calculations show that $T(h)$ is in $\allpoly$.
  \begin{align*}
    \diag h(y) f(x,y) &= \diag \sum_{i,j,k} b_k a_{ij}x^iy^{j+k} \\
&= \sum_k b_k \diag \sum_{i,j} a_{ij} x^i y^{j+k} \\
&= \sum_k b_k \sum_{i=j+k} a_{ij} x^i \\
& = \sum b_k s_k \quad=\quad T(h)
  \end{align*}
\end{proof}

\begin{cor} \label{cor:diag-binom}
  For any $n$ these polynomials  have all real roots:
$$ \sum_i \binom{n}{i,i}x^i \quad\quad \sum_i \frac{x^i}{i! (n-2i)!}$$
\end{cor}
\begin{proof}
  The two polynomials in question are the diagonal polynomials of a polynomial
  in $\gsubpos_2$, namely $(x+y+1)^n$.
\end{proof}

\begin{remark}
  \index{even part} 
  
  Here's a variation on the even part (\chapsec{linear}{hurwitz}). If
  $f=\sum a_ix^i$ is in $\allpolyalt$ then $f(x+y)$ is in
  $\gsubplus_2$. Now $\diag(x+y)^k$ is either $0$ if $k$ is odd, or
  $\binom{2m}{m} x^{m}$ if $k=2m$ is even.  Consequently, $ \diag_1
  f(x+y) = \sum \frac{(2i)!}{i!}\, a_{2i} x^i$ is in $\allpoly$.
\end{remark}

  In the remark we considered $f(x+y)$. If we consider
  $f(\,(x+1)(y+1)\,)$ then we get a quite different transformation:

\begin{cor} \label{cor:fx1y1}
  The linear transformation $x^r \mapsto \sum_i \binom{r}{i}^2 x^i$
  maps $\allpolyalt$ to itself.
\end{cor}
\begin{proof}
If $f\in\allpolyalt$ then $f(xy)\in\gsubclosepos_2$, and so
after substituting we have that $f(\,(x+1)(y+1)\,)\in \gsubpos_2$. If we
write $f=\sum a_ix^i$ then
\begin{align*}
  f(\,(x+1)(y+1)\,) &= \sum_i a_i (x+1)^i(y+1)^i \\
\diag\, f(\,(x+1)(y+1)\,) &= \sum_i a_i \sum_j \binom{i}{j}^2\,x^j
\end{align*}
which establishes the result.
\end{proof}

\index{partial sums}
In \chapsec{polynomials}{poly-no-roots} we saw that the partial sums
of the function $e^x$ in $\allpolyf$ do not  have all real roots, but
adding another factorial lands in $\allpoly$.

\begin{lemma} \label{lem:sumsofF}
  If $f = \sum \displaystyle{a_i \frac{x^i}{i!}}$ is in $\allpolyf$ then
  for any $n$ the polynomial  \\
$\displaystyle{ \sum_{i=0}^n \binom{n}{i} a_i\, x^i}
=\displaystyle{ n!\,\sum_{i=0}^n \frac{a_i}{i!(n-i)!}\, x^i}$ has all real
roots. 
\end{lemma}
\begin{proof}
  The operator $f(\partial_x)$ maps $\gsubpos_2$ to itself, so
  consider
  \begin{align*}
    f(\partial_x)(x+y)^n\Big\vert_{x=0} & = \sum_{i,j} a_i y^{n-j}
    \binom{n}{j} \diffd^i x^j\Big\vert_{x=0} \\
& = \sum _k a_k\, y^{n-k} \binom{n}{k}
  \end{align*}
Taking the reverse establishes the lemma.
\end{proof}

\index{Hadamard product!in $\allpolyf$}
We can take limits of the Hadamard product to conclude
from Theorem~\ref{thm:hadamard-2}
\begin{cor} \label{cor:hadamard-3}
If $f\in\allpolyposf$ and $g\in\allpolyf$ then $\expoper{}^{-1} f\ast
g$ is in $\allpolyf$.  
\end{cor}

Here's a representation of $\allpolyf$ as the image of one exponential.

\begin{lemma}\label{lem:represent-f}
  \begin{align}
    \label{eqn:rep-f-1}
    \allpolyf &= \left\{ T(e^x)\mid T\colon{}\allpolypos\longrightarrow\allpoly\right\}\\
    \allpolyposf &= \left\{ T(e^x)\mid T\colon{}\allpolypos\longrightarrow\allpolypos\right\}
  \end{align}
\end{lemma}
\begin{proof}
  Since $T\colon{}\allpolypos\longrightarrow\allpoly$ extends to
  $T\colon{}\allpolyposf\longrightarrow\allpolyf$, we know that
  $T(e^x)\in\allpolyf$. Conversely, choose $f\in\allpolyf$ and define
$$ T_f(g) = f \ast' g = f \ast \expoper{}^{-1}(g)$$
We know that $T_f:\allpolyposf\longrightarrow\allpolyf$. Since
$T_f(e^x)=f$, we have equality. The second case is similar.
\end{proof}

Unlike $\allpoly$, factors of elements in $\allpolyf_2$ aren't
necessarily in $\allpolyf_2$. A simple example is $x =
\left(e^{xy}\right)\left(x e^{-xy}\right)$. Note that 
$x e^{-xy}\in\allpolyf_2$, yet $e^{xy}\not\in\allpolyf_2$. However, we
do have

\begin{lemma}\label{lem:jensen1}
  If $f(x,y)e^{-xy}\in\allpolyf_2$, then $f(x,\alpha)$ and
  $f(\alpha,x)$ are in $\allpolyf$ for any choice of
  $\alpha\in\reals$. 
\end{lemma}
\begin{proof}
  If $f(x,y)e^{-xy}\in\allpolyf_2$, then $f(x,\alpha)e^{-\alpha
    x}\in\allpolyf$. Since $e^{\alpha x}\in\allpolyf$ it follows that
  $f(x,\alpha)\in\allpolyf$. 
\end{proof}

  \section{Differential operators}
  \label{sec:diff-oper}

  If we substitute derivatives for some variables, the result maps
  $\gsubclose_d$ to itself.  This generalization of the corresponding
  property for $f(\diffd)$ acting on $\allpoly$ has many important uses.

  \begin{lemma}\label{lem:fdxdy}
    If $f(\xx)\in\gsubclose_d$ then
    $f(-\partial_\xx)\colon\gsubclose_d\longrightarrow\gsubclose_d$. 
  \end{lemma}
  \begin{proof}
    Since $e^{-\partial_\xx\partial_\yy}$ maps $\gsubclose_{2d}$ to
    itself the lemma follows from the identity
\[
e^{-\partial_\xx\partial_\yy}\,g(\xx)\,f(\yy)\biggl|_{\yy=0} \,=\,
f(-\partial_\xx)g(\xx).
\]
By linearity we only need to check it for monomials:
\[
e^{-\partial_\xx\partial_\yy} \xx^\sdiffi \yy^\sdiffj\biggl|_{\yy=0} \, =\,
\frac{(-\partial_\xx)^\sdiffj}{\diffj!} \xx^\sdiffi \,\cdot\, (\partial_\yy)^\sdiffj
\yy^\sdiffj\biggl|_{\yy=0}\, =\,
(-\partial_\xx)^\sdiffj \xx^\sdiffi
\]
  \end{proof}

  We use the curious equality below to derive a corollary
  due to Sokal and Lieb \cite{lieb-sokal}.
  \begin{lemma}\label{lem:expdxdy}
    If $f(x,y)$ is any polynomial then
\[
f(x,y-\partial_x) = e^{-\partial x\partial y} f(x,y)
\]
  \end{lemma}
  \begin{proof}
    By linearity we may assume that $f(x,y) = x^sy^r$. The left hand
    side is
    \begin{align*}
      (y-\partial_x)^rx^s &=
\sum_{i=0}^r y^{r-i}\binom{r}{i}(-1)^i \partial_x^i x^s \\
&= \sum_{i=0}^r y^{r-i}x^{s-i} \frac{(-1)^i}{i!}
\frac{r!}{(r-i)!}\frac{s!}{(s-i)!} \\
&= \sum_{i=0}^\infty \frac{(-\partial_x\partial_y)^i}{i!} x^s y^r
\end{align*}
{which is exactly the right hand side.}
  \end{proof}

  \begin{cor}\label{cor:sokol-lieb}
    If $f(\xx,\yy)\in\rupint{2d}$ then
    \begin{enumerate}
    \item $f(\xx,\yy-\partial_\xx)\in\gsubclose_{2d}$.
    \item $f(\xx,-\partial_\xx)\in\gsubclose_{d}$.
    \end{enumerate}
  \end{cor}
  \begin{proof}
    If we iterate Lemma~\ref{lem:expdxdy} we get
\[
f(x_1,y_1-\partial x_1,\dots,x_d,y_d-\partial x_d) =
e^{-\partial x_1\partial y_2}\cdots e^{-\partial x_d\partial y_d}
f(\xx,\yy)
\]
The first part now follows from Corollary~\ref{cor:p2-exy}. If we set
$\yy=0$ we get the second part.
  \end{proof}

If we substitute $\partial x_i$ for $x_i$ then the result is much
easier to prove.

  \begin{lemma}\label{lem:diff-d-one}
    If $f = \sum f(x_1,\dots,x_{d-1})x_d^k$ is in $\rupint{d}$ and $f(0)
    \ne0$ then the
    linear transformation $g\mapsto \sum f_k(x_1,\dots,x_{d-1})
    \frac{\partial^k}{\partial x_d^k}(g)$ maps $\rupint{d}$ to itself.
  \end{lemma}
  \begin{proof}
    The homogeneous part of $T(g)$ is $f(0)g^H$, so $T$ satisfies the
    homogeneity condition. If we choose $a_1,\dots,a_{d-1}\in\reals$
    then
\[
      T(g)(a_1,\dots,a_{d-1},x_d) =
      f(a_1,\dots,a_{d-1},\partial_{x_d}) g(a_1,\dots,a_{d-1},x_d)
\]
and this is in $\allpoly$ since $f(a_1,\dots,a_{d-1},x)$ and
$g(a_1,\dots,a_{d-1},x)$ are in $\allpoly$.
\end{proof}

\index{Hadamard product!in $\rupint{d}$}
Next is a Hadamard product result.

\begin{lemma}\label{lem:hp-in-pd}
  If $\sum f_\sdiffi(\xx)\yy^\sdiffi$ and   $\sum
  g_\sdiffi(\xx)\yy^\sdiffi$ are in $\gsubclose_{2d}$ then
\[
\sum _\sdiffi (-1)^\sdiffi\, \sdiffi!\,  f_\sdiffi(\xx)
\,g_\sdiffi(\xx) \in\gsubclose_d
\]
\end{lemma}

\begin{proof}
  A calculation shows:
\[
  \biggl(\sum f_\sdiffi(\xx)(-\partial\yy)^\sdiffi\biggr)
  \biggl(\sum g_\sdiffi(\xx)\yy^\sdiffi\biggr) \biggr|_{\yy=0}
=
\sum _\sdiffi (-1)^\sdiffi\, \sdiffi!\,  f_\sdiffi(\xx)
\,g_\sdiffi(\xx) 
\]
\end{proof}

\begin{lemma}
  If $Q=(q_{ij})$ is a symmetric $d$ by $d$ matrix, then $e^{-\xx
    Q\xx^t}\in\allpolyf_d$ if and only if all elements of $Q$ are non-negtive.
\end{lemma}
\begin{proof}
  If all elements are non-negative then it is easy to see that $e^{-\xx
    Q\xx^t}\in\allpolyf_d$. 
  Since $\alpha-xy\in\gsubclose_2$ for $\alpha>0$,  compute
  \[e^{-(ax^2+2bxy+cy^2)}(\partial_x,\partial_y)\,(\alpha-xy) = \alpha-xy +2b\]
  Since this is in
  $\gsubclose_2$, it follows that $\alpha+2b\ge0$ and so $b$ is non-negative.
\end{proof}

  \section{Limits and transformations}
  \label{sec:limits-transf}

  If $f$ and $g$ are entire functions and $T$ is a linear
  transformation then we define $f(T)g$ to be the \emph{formal} series
  \begin{equation}
    \label{eqn:ftg}
    \sum_{i=0}^\infty \sum_{j=0}^\infty a_i b_j T^i(x^j)
  \end{equation}
where $f = \sum_0^\infty a_ix^i$ and $g = \sum_0^\infty b_ix^i$. This
is only a formal sum since there is no apriori reason why $f(T)g$ should
exist. For a full discussion of the complexities when $T = a\diffd + x
\diffd^2$ see \cite{wolowski-convex}.

We are interested in the case where $f,g\in\allpolyf$ and
$T\colon\allpoly\longrightarrow\allpoly$. We assume that
$f_n\rightarrow f$ and $g_n\rightarrow g$ where convergence is uniform
convergence on compact subsets. The problem is then to determine when
\[
\lim_{n\rightarrow\infty} f_n(T)g_n\ \text{exists and equals\ 
  $f(T)g$}.
\]

The simplest case is when all limits involve polynomials.

\begin{lemma}\label{lem:limit-1}
  If \begin{enumerate}
\item $T$ decreases degree.
\item $h(T)\colon\allpoly\longrightarrow\allpoly$ for all $h\in\allpoly$.
\item $f\in\allpolyf$
\end{enumerate}
 then $f(T)\colon\allpoly\longrightarrow\allpoly$.
\end{lemma}
\begin{proof}
  If $f_n\rightarrow g$ and $g\in\allpoly(m)$ then both $f_n(T)g$ and
  $f(T)g$ are polynomials since $T^k(g)=0$ for $k\ge m$. Since the
  coefficients of $f_n$ converge to the coefficients of $f$ we see
  that $f_n(T)g\rightarrow f(T)g$. By hypothesis all
  $f_n(T)g\in\allpoly$, so $f(T)g\in\allpoly$.
\end{proof}

\begin{example}
  For example, since $f(\diffd)$ maps $\allpoly$ to itself it follows
  that $e^{\diffd}:\allpoly\longrightarrow\allpoly$. Of course, this
  is trivial since $e^\diffd f(x) = f(x+1)$, but the same argument
  applies to $e^{-\partial_x\partial_y}$ (see Corollary~\ref{cor:p2-exy}).
\end{example}

\index{Montel's theorem}
In order to work with entire functions we recall Montel's theorem:
\begin{quote}
  Any locally bounded sequence of holomorphic functions $f_n$ defined on
  an open set $D$ has a subsequence which converges uniformly on compact subsets to
  a holomorphic function $f$.
\end{quote}

It follows that if $T$ is an operator that takes bounded sequences to
bounded sequences and if $f_n\longrightarrow f$ then
$Tf_n\longrightarrow Tf$. For instance, the $k$'th derivative is a bounded
operator by Cauchy's integral formula:
\[ 
f^{(k)}(a) = 
\frac{k!}{2\pi\imag}\oint_C 
\frac{f(z)}{(z-a)^{k+1}} \,dz
\]
where $C$ is a small circle centered at $a$.  More generally, if
$f(x,y) = \sum f_i(x)y^i$ is a polynomial in two variables then 
the operator
\[
T\colon g \mapsto \sum f_i(x) D^i g
\]
is a bounded operator since it is a sum of  products of  locally bounded
functions $f_i$ times bounded operators. We therefore have

\begin{lemma}\label{lem:limit-diff}
  If $T$ is a differential operator as above then
  $\lim_{n\rightarrow\infty} T\bigl(1-\frac{xy}{n}\bigr)^n$ exists,
  equals $T(e^{-xy})$, and is an entire function.
\end{lemma}

In order to prove that sequences are uniformly bounded we use the bounds
on the size of a polynomial.

\begin{lemma}\label{lem:limit-tpp}
  Suppose that  $T$ is a linear transformation and
  \begin{enumerate}
  \item  $T_\ast\colon\gsubclose_2\longrightarrow\gsubclose_2$
  \item $T(1)\in\reals\setminus0$ and  $f(0)\ne0$.
  \item $f_n\rightarrow f$ in $\gsubf_2$ where $f_n\in\gsubclose_2$.
  \end{enumerate}
 then $T_\ast (f_n) \longrightarrow T_\ast(f)$. 
\end{lemma}
\begin{proof}
  Write $f_n = \sum_i F_{n,i}(x)y^i$ and $f = \sum_i F_i(x)y^i$. Now
  $Tf_n= \sum T(F_{n,i})y^i$ is in $ \gsubclose_2$ by hypothesis, so
  for any $\alpha\in\reals$ we have $(T\,f_n)(\alpha,y)\in\allpoly$ and
  therefore we can apply Lemma~\ref{lem:abs-bound-2}:
\[
    \sup_{|y|\le r} (T\,f_n)(\alpha,y) \le
T(F_{n,0})
\exp\biggl( r \frac{T\,F_{n,1}(\alpha)}{T(F_{n,0})} + 
2 r^2 \frac{(T\,F_{n,1}(\alpha)\,)^2}{(TF_{n,0})^2} + 
r^2 \frac{T\,F_{n,2}(\alpha)}{T(F_{n,0})}  \biggr)
\]
If we define
\begin{align*}
  a = &\inf_n T(F_{n,0})  &   A =& \sup_n T(F_{n,0}) \\
  B_r =&  \sup_{n,|\alpha|\le r} T(F_{n,1}(\alpha)) &
  C_r =&  \sup_{n,|\alpha|\le r} T(F_{n,2}(\alpha)) 
\end{align*}
then $A,B_r,C_r$ are finite and $a$ is non-zero by hypothesis so
\[
    \sup_{|\alpha|\le r,|y|\le r} (T\,f_n)(\alpha,y) \le
A
\exp\biggl( r \frac{B}{a} + 
2r^2 \frac{B^2}{a^2} + 
r^2 \frac{C}{a}  \biggr)
\]
Thus $T(f_n)$ is uniformly bounded and therefore converges to $T(f)$.
\end{proof}

\section{Totally positive functions}
\label{sec:tot-pos-fct}

We raise some questions about totally positive functions and
$\gsubplus_2$. 
\index{totally positive!function}
\index{totally positive!function on $\reals^2_+$}

\begin{definition}
  A function $f(x,y)$ is \emph{strictly totally positive} if for every positive
  integer $n$ and sequences $x_1<x_2<\cdots<x_n$ and
  $y_1<y_2<\cdots<y_n$ the determinant
  \begin{equation}
    \label{eqn:tot-pos-fct}
    \begin{vmatrix}
      f(x_1,y_1) & f(x_1,y_2) & \hdots & f(x_1,y_n) \\
      f(x_2,y_1) &  &  & \vdots   \\
      \vdots & & &  \\
      f(x_n,y_1) & \hdots &   & f(x_n,y_n)
    \end{vmatrix}
  \end{equation}
is positive. If the sequences $\{x_i\}$ and $\{y_i\}$ are restricted
to positive values then we say that $f(x,y)$ is a strictly totally positive
function on $\reals^2_+$. 
\end{definition}

If we take $n=1$ we see that a necessary condition for $f(x,y)$ to be
strictly totally positive on $\reals_+^2$ is that $f(x,y)$ is positive for $x$
and $y$ positive.  A natural question (\ref{ques:tot-pos-fct}) is to
determine conditions on polynomials $g(x,y)$ so that $1/g$ is
strictly totally positive on $\reals_+^2$. 

We know three functions for which this holds. It is well
known that $e^{xy}$ is strictly totally positive, and $e^{xy}$ is of the form
$1/g$ where $g=e^{-xy}\in\gsubf_2$. The simplest function in
$\gsubplus_2$ is $x+y$; this is strictly totally positive on $\reals_+^2$. That
this is so follows from a famous evaluation. Recall Cauchy's double
alternant:
\index{Cauchy's identity}
$$
\operatornamewithlimits{det}_{1\le i,j\le n}
\left(\frac{1}{x_i+y_j}\right) =
\frac{\displaystyle \prod_{1\le i< j\le n} (x_i-x_j)(y_i-y_j)}%
{\displaystyle \prod_{1\le i,j\le n} (x_i+y_j)}
$$
The numerator is positive for any pair of strictly increasing
sequences; the denominator is positive since all $x_i$ and $y_j$ are
positive. 

For the last example, there is a quadratic analog
\eqref{eqn:borchardt} of Cauchy's formula due to Borchardt\cite{singer}
where $perm$ is the permanent. Since the permanent of a matrix with
all positive entries is positive, it follows that $(x+y)^{-2}$ is
strictly totally positive on $\reals_+^2$.

\index{Borchardt's identity}

\begin{equation}\label{eqn:borchardt}
\operatornamewithlimits{det}_{1\le i,j\le n} 
\left(\frac{1}{(x_i+y_j)^2}\right) = 
\operatornamewithlimits{det}_{1\le i,j\le n} 
\left(\frac{1}{x_i+y_j}\right)  
\operatornamewithlimits{perm}_{1\le i,j\le n} 
\left(\frac{1}{x_i+y_j}\right) 
\end{equation}

If $h(x)$ is any function that is positive for positive $x$, and
$f(x,y)$ is totally positive on $\reals_+^2$ then 
$h(x)f(x,y)$ is strictly totally positive on $\reals_+^2$. This follows from
simple properties of the determinant:

$$  \operatornamewithlimits{det}_{1\le i,j\le n}
\bigl(h(x_i)f(x_i,y_j)\bigr)  = h(x_1)\cdots h(x_n)\ 
 \operatornamewithlimits{det}_{1\le i,j\le n} \bigl(f(x_i,y_j\bigr)>0$$
 
 Consequently, if $g\in\gsubf_2$ satisfies the property that $1/g$ is
 strictly totally positive on $\reals_+^2$ then so does $e^{-x^2}g$.

It is not true that $1/g$ is positive semi-definite if
$g\in\gsubplus_2$ - see \cite{bhatia}. 

 We need a few properties of positive definite matrices
 \cite{johnson}.

\index{positive definite!ordering of matrices}

\begin{definition}
  Let $A,B$ be \index{Hermitian matrices}Hermitian matrices.  We write
  $B \prec A$ if the matrix $A-B$ is positive definite.  If $A-B$ is
  positive semidefinite we write $B\preceq A$. The relation $\preceq$
  makes the set of all Hermitian matrices into a partially ordered
  set.
\end{definition}

Here are some properties of $\prec$ that we will use.
\begin{enumerate}
\item $A\prec B$ implies $det(A) < det(B)$.
\item $A\prec B$ implies $B^{-1} \prec A^{-1}$.
\item $A\prec B$ and $C$ positive definite implies $A+C\prec B+C$.
\item $A\prec B$ and $C$ positive definite implies $CAC\prec CBC$.
\end{enumerate}

If $f(x)$ is a positive definite matrix for all $x$ in an interval on
the real line, then we say that $f$ is \emph{increasing} if $a<b$
implies $f(a)\prec f(b)$. $f$ is \emph{decreasing} if $-f$ is
increasing. 

  \begin{lemma}\label{lem:positive-def-2}
If $A$ and $B$ are positive definite matrices, and 
    $0<x_1<x_2$, $0<y_1<y_2$ then
$$
\begin{vmatrix}
\displaystyle  \frac{1}{|I+x_1A + y_1B|} &
\displaystyle  \frac{1}{|I+x_2A + y_1B|}\\
\displaystyle  \frac{1}{|I+x_1A + y_2B|}&
\displaystyle  \frac{1}{|I+x_2A + y_2B|}
\end{vmatrix} >0
$$
  \end{lemma}
  \begin{proof}
    If we let $C = x_2A$, $D=y_2B$, $x=x_1/x_2$, $y = y_1/y_2$ then $0<
    x,y<1$ and $C,D$ are positive definite. We need to show that
$$
\begin{vmatrix}
\displaystyle  \frac{1}{|I+xC +yD|} &
\displaystyle  \frac{1}{|I+C+ yD|}\\
\displaystyle  \frac{1}{|I+xC + D|}&
\displaystyle  \frac{1}{|I+C + D|}
\end{vmatrix} >0
$$
Since all the matrices are positive definite the determinants of
positive linear combinations are positive. Thus  we need to show that
\footnote{Thanks to Han Engler for this argument.}
\begin{gather*}
{|I+xC +D|} {|I+C + yD|} >
{|I+xC+ yD|}
{|I+C + D|}\\
\intertext{which is equivalent to}
\displaystyle \frac{|I+xC +yD|}{|I+xC+ D|} < 
\frac{|I+C + yD|}{|I+C + D|} 
\end{gather*}
If we denote the left hand side by $F(x)$ then we will show that
$F(x)$ is an increasing function. 

\begin{align*}
  I+xC+D & \text{ is increasing}\\
  (I+xC+D)^{-1} & \text{ is decreasing}\\
  D^{1/2}(I+xC+D)^{-1}D^{1/2} & \text{ is decreasing}\\
  (y-1)D^{1/2}(I+xC+D)^{-1}D^{1/2} & \text{ is increasing}\\
 I+ (y-1)D^{1/2}(I+xC+D)^{-1}D^{1/2} & \text{ is increasing}\\
\intertext{The determinant of the last expression is $F(x)$ - this
  follows from the identity}
|(\alpha A+B)B^{-1}| &= |I+\alpha A^{1/2}B^{-1}A^{1/2}|
\end{align*}

  \end{proof}

\section{Hermite polynomials}
\label{sec:hermite-n-dim}
\index{Hermite polynomials!in $\gsubpos_d$} 

We use generating functions to define Hermite polynomials in $d$
variables. The main result is that they are in $\gsubpos_d$.

Let $S$ be a $n\times n$ Hermitian matrix, $\by=(y_1,\dots,y_n)$ and
$\bx=(x_1,\dots,x_n)$. The Hermite polynomials determined by $S$ are
defined to be the coefficients of the exponential generating function

\begin{equation} \label{eqn:hermite-d}
 f_S(\xx,\yy) \,=\, \exp\left( - \by S\by^\ast - 2 \by S\bx^\ast\right) = \sum_\diffi 
H_\diffi(\bx) \frac{(-\by)^\diffi}{\diffi!}
\end{equation}
where  $\diffi = (i_1,\dots,i_n), $  $\diffi! = i_1!\cdots i_n!,$ 
$\by^\diffi = y_1^{i_1}\cdots y_n^{i_n},$ and $\zz^\ast$ is the
conjugate transpose of $\zz$. The factor $2 \by S\bx^\ast$ has a
negative sign, so we can  use Proposition~\ref{prop:e-xhx} to
conclude that $f_S\in\allpolyf_{2d}$.

If we take $n=1$ and $S=(1)$ then the exponential is
$e^{-y_1^2-2x_1y_1}$ which is exactly the generating function of the
one variable Hermite polynomials.

If we  multiply \eqref{eqn:hermite-d} by $\exp(\,-\xx S\xx^t\,)$ then
$$  \exp(\,-\xx S\xx^t\,) \ \sum H_\diffi \frac{(-\yy)^\diffi}{\diffi!} =
\exp(\, -(\yy+\xx)S(\,\yy+\xx)^t\,).
$$
\index{Rodrigues' formula!Hermite}
We use this to derive the Rodrigues' formula for the Hermite
polynomials.  The Taylor series of  $g(\xx)$  asserts
\begin{align*}
\index{Taylor series}
  g(\xx+\yy) &= \sum_\diffi 
  \left(\frac{\partial}{\partial\xx}\right)^\diffi g(\xx)  \,
\frac{\yy^\diffi}{\diffi!}\\ 
\intertext{so if we choose $g(\xx) = \exp(-\xx S\xx^t)$ then}
\exp(-(\yy+\xx)S(\yy+\xx)^t ) &= \sum 
  \left(\frac{\partial}{\partial\xx}\right)^\diffi \,
e^{-\xx S\xx^t}\,
\frac{\yy^\diffi}{\diffi!}\\
&= \sum e^{-\xx S\xx^t} H_\diffi(\xx) \frac{(-\yy)^\diffi}{\diffi!}.\\
\intertext{Equating coefficients of $\yy$ yields the Rodrigues' formula}
\index{Rodrigues' formula}
H_\diffi(\xx) &= (-1)^{\vert\diffi\vert}\,e^{\xx S\xx^t} \,
  \left(\frac{\partial}{\partial\xx}\right)^\diffi \,
e^{-\xx S\xx^t}.\\
\intertext{If we rewrite this as}
e^{-\xx S\xx^t} \,H_\diffi(\xx) &= (-1)^{\vert\diffi\vert}
  \left(\frac{\partial}{\partial\xx}\right)^\diffi \,
e^{-\xx S\xx^t}
\end{align*}
then it follows that $H_\diffi(\xx)$ satisfies substitution.

In order to show that $H_\diffi$ is in $\gsubpos_n$ it suffices to
determine the homogeneous part of $H_\diffi$, and we use the Rodrigues
formula for that. We first note that if $e_i$ is the coordinate vector
with a $1$ in the $i$-th place, and $0$ elsewhere then
$$
  \frac{\partial}{\partial x_i} \, \xx S \xx^t = ( e_i S + S e_i^t)
  \xx. 
$$
Consequently, if we differentiate $g(\xx)\exp(-\xx S\xx^t)$ where $g$ is
  some polynomial then 
\begin{align*}
   \frac{\partial}{\partial x_i} g(\xx) \exp(-\xx
      S\xx^t) &= k(\xx) \exp(-\xx S\xx^t)\\
\intertext{where $k(\xx)$ is a polynomial, and moreover }
k(\xx)^H &= (\,(e_i S + S e_i^t)   \xx \,) g(\xx)^H. \\
\intertext{ Upon writing  }
  \left( \frac{\partial}{\partial x_i}\right)^\diffi  \exp(-\xx
      S\xx^t) &= j(\xx) \exp(-\xx S\xx^t)\\
\intertext{we conclude that  the homogeneous part of $j$ is}
j(\xx)^H & = \prod \,(\,(e_i S + S e_i^t)   \xx \,).
\end{align*}

By assumption all coefficients of $S$ are positive, so all the
coefficients of $j^H$ are positive, and by the above it follows that
the homogeneous part of $H_\diffi$ has all positive coefficients. This
proves

\begin{theorem} \label{thm:hermite-pd}
  If $S$ is a symmetric matrix with all positive coefficients then the Hermite
  polynomials determined by \eqref{eqn:hermite-d} are in $\gsubpos_d$. 
\end{theorem}

\section{Hermite polynomials in two variables}
\label{sec:hermite-2-var}

The results of this section exist simply to explain properties that
can be observed in a particular graph of a Hermite polynomial in two
variables. If we take
$Q=\smalltwo{1}{3}{3}{2}$ and $\diffi=(2,3)$ then 
\begin{multline*}
  H_{2,3} = -16\,\bigl( 108\,x - 207\,x^3 + 54\,x^5 + 114\,y -
  876\,x^2\,y + 432\,x^4\,y - 990\,x\,y^2 \\ + 1206\,x^3\,y^2 - 332\,y^3
  + 1420\,x^2\,y^3 + 744\,x\,y^4 + 144\,y^5 \bigr)
\end{multline*}
The homogeneous part factors into $-32\,{\left( 3 + x \right)
}^2\,{\left( 2 + 3\,x \right) }^3$, which explains the asymptotes in 
Figure~\ref{fig:h23}. The rest of the graph  of $H_{2,3}$ is quite surprising:

\begin{figure}[htbp]
  \centering
  \includegraphics*[width=2in]{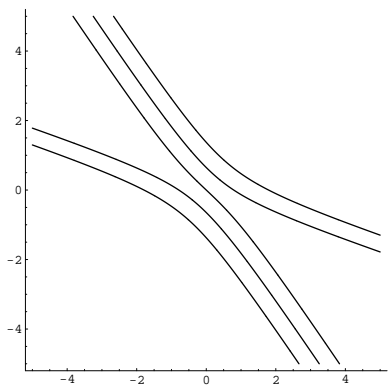}
  \caption{The graph of $H_{2,3}$}
  \label{fig:h23}
\end{figure}

The following result shows that all the $H_{i,j}$ are in $\gsubgen_2$,
but doesn't explain the striking regularity.

\begin{lemma}
  If $S$ is not a diagonal matrix then  $H_{i,j}\in\gsubgen_2$.
\end{lemma}
\begin{proof}
  Let $S=\smalltwo{a}{b}{b}{c}$ where $b\ne0$, and set $Q =
  ax^2+2bxy+cy^2$. Suppose that $H_{i,j}$ has intersecting solution
  curves. If so, then $H_{i,j}$ and $\frac{\partial}{\partial
    x}H_{i,j}$ have a common zero. Now since
  \begin{align*}
    \frac{\partial}{\partial x}H_{ij} &=
    \frac{\partial}{\partial x}\left( e^Q\, \partial_x^i\,\partial_y^j\,
      e^{-Q}\right)\\
&= e^Q\,\partial_x^{i+1}\,\partial_y^j \,e^{-Q} + Q_x
e^Q\,\partial_x^i\,\partial_y^j\,e^{_Q} \\&= H_{i+1,j} + Q_x H_{ij}
  \end{align*}
it follows that $H_{i,j}$ and $H_{i+1,j}$ have a common zero. Since
$e^Q$ has no zeros, both 
$\partial_x^i\,\partial_y^j\,e^{_Q}$ and
$\partial_x^{i+1}\,\partial_y^j\,e^{_Q}$ have a common root. Since
$\partial_x^{i}\,\partial_y^j\,e^{_Q}\lesslesseq
\partial_x^{i+1}\,\partial_y^j\,e^{_Q}$, we conclude from
Corollary~\ref{cor:p2-deriv} that 
$\partial_x^{i-1}\,\partial_y^j\,e^{_Q}$ and
$\partial_x^{i}\,\partial_y^j\,e^{_Q}$ have a common root. Continuing,
we see
$\partial_y^j\,e^{_Q}$ and
$\partial_x\,\partial_y^j\,e^{_Q}$ have a common root. Eliminating
derivatives of $y$, we conclude that 
$\partial_y\,e^{_Q}$ and
$\partial_x\,\partial_y\,e^{_Q}$ have a common root. Because $H_{1,0}
= -2 b x - 2 c y$, we can solve for $x$ and substitute into 
$$H_{1,1} = 2\,\left( -b + 2\,a\,b\,x^2 + 2\,b^2\,x\,y + 2\,a\,c\,x\,y
  + 2\,b\,c\,y^2 \right) $$
which yields $-2b$. Since $b\ne0$, there are no common roots.
\end{proof}

We follow Nunemacher in finding the exact equations of the
asymptotes. Let
  \begin{align*}
    Q &= a x^2    + 2 b x y + cy^2 \\
    H_{n,m}(x,y) &=     -e^{ax^2 + 2b x y + c y^2}\left(\frac{\partial}{\partial x}\right)^n
\left(\frac{\partial}{\partial y}\right)^m e^{-({ax^2 + 2b x y + c
    y^2})}\\
H_{n,m+1} &= -2(bx+cy) H_{n,m} + \partial_y H_{n,m} \\
H_{n+1,m} &= -2(ax+by) H_{n,m} + \partial_x H_{n,m} \\
\end{align*}

From the recursions it easily follows that the homogeneous part of
$H_{n,m}$ is $(-2)^{n+m}(ax+by)^n(bx+cy)^m$. Consequently, there are
$n$ asymptotes in one direction, and $m$ in another. We proceed to
find the precise equations of these asymptotes. 
It's easy to see that we can write 
\begin{align*}
  H_{n,m} &= (ax+by)^n Q_{n,m}(x,y) +(ax+by)^{n-2}Q_{n-2,m}(x,y) + \cdots \\
  H_{n,m+1} &= (ax+by)^n Q_{n,m+1}(x,y) +(ax+by)^{n-2}Q_{n-2,m+1}(x,y) + \cdots \\
 &= -2(bx+cy)\left[(ax+by)^n Q_{n,m}(x,y) +(ax+by)^{n-2}Q_{n-2,m}(x,y) + \cdots \right]\\
& + \partial_y\left[(ax+by)^n Q_{n,m}(x,y) +(ax+by)^{n-2}Q_{n-2,m}(x,y) + \cdots \right]\\
\end{align*}
{Now when we substitute $x=b,y=-a$, all terms with a factor
  of $(ax+by)$ vanish, so equating terms of like degree yields}
\begin{align*}
Q_{n,m+1}(b,-a) &= -2(bx+cy) Q_{n,m}(b,-a) \\
Q_{n-2,m+1}(b,-a) &= -2(bx+cy) Q_{n-2,m}(b,-a) \\
\dots \\
\intertext{Thus we see that}
Q_{n-2r,m}(b,-a) &= (-2(bx+cy))^m Q_{n-2r,0}(b,-a)\\
\intertext{The asymptotes \cite{nunemacher} are $ax+by=\alpha$ where
  $\alpha$ is the root of}
& t^nQ_{n,m}(b,-a) + t^{n-2}Q_{n-2,m}(b,-a) + \cdots \\
&= (-2(bx+cy))^m \left[t^nQ_{n,0}(b,-a) + t^{n-2}Q_{n-2,0}(b,-a) + \cdots \right]\\
\intertext{and so the asymptotes of $H_{n,m}$ and $H_{n,0}$ are the
    same. But, }
H_{n,0}(x,y) &= 
-e^{ax^2 + 2b x y }\left(\frac{\partial}{\partial x}\right)^n
 e^{-({ax^2 + 2b x y })}\\
\intertext{consists of parallel lines, so we can substitute $y=0$ to
  see that the asymptotes $ax+by=\alpha$ to $H_{n,m}$ satisfy
  $H_{n,0}(\alpha,0)=0$. Since}
H_{n,0}(x,0) &= e^{ax^2} \left(\frac{\partial}{\partial x}\right)^n e^{-ax^2}
\end{align*}

is the Hermite polynomial $H_n(x/\sqrt{a})$ we conclude:

\begin{prop}
  The asymptotes of $H_{n,m}$ are of the form
$$ ax+by = \alpha \quad\quad bx+cy=\beta$$
where $H_n(\alpha/\sqrt{a})=0$ and $H_m(\beta/\sqrt{c})=0$.
\end{prop}

  The coefficients of some Hermite polynomials obey a recursion
  similar to the three term recursion for orthogonal polynomials.
  Choose a positive symmetric matrix $S$, fix a positive integer $n$, and define

\begin{align*} H_n(\xx)  & = e^{\xx S \xx'} \left(\frac{\partial}{\partial
    x_1}\right)^n e^{-\xx S \xx'}\\
&= \sum_{i=0}^n f_i(x_1,\dots,x_{d-1}) x_d^i\\
\intertext{If we let $v=(v_1,\dots,v_d)$ be the first row of $S$, then
  it easy to see that}
H_{n+1} &= -2(v_1x_1 + \cdots + v_dx_d) H_{n} - 2n v_1\, H_{n-1}\\
\intertext{We can use this recursion to prove by induction on $n$
  that}
f_k &= \left(\sum_{i=1}^{d-1}\frac{n+1}{n-k}
  \,\frac{v_i}{v_d}\,x_i\right)\,f_{k+1} -
\left(\frac{(n+1)(n+2)}{2(n-k)}\,\frac{v_1}{v_d^2}\right) \,f_{k+2}\\
\intertext{If $d=2$, and $S =\smalltwo{a}{b}{b}{c}$ then}
f_k(x) &= \left(\frac{k+1}{n-k}\,\frac{a}{b}\right)x  f_{k+1}(x) -
\left(\frac{(k+1)(k+2)}{2(n-k)}\,\frac{a}{b^2}\right) f_{k+2}(x)
\end{align*}

We get a three term recursion if we consider $H_{n,1}$ instead of 
$H_{n,0}$, but the coefficients of $H_{n,2}$ have a five term recursion
that includes $x^2f_{n+2}$ and $x f_{n+3}$.


\chapter{Extending $\rupint{d}$}
\label{cha:extending}

\renewcommand{\TimeStampStart}{Monday, December 17, 2007: 17:22:08}
\mytoday

The goal in this Chapter is to investigate
polynomials that only satisfy substitution for some values of $y$. 
There are several different approaches.

\section{Taking the image of $\gsubpos_2$}

We define the set $\image{2}$ of polynomials as images of polynomials
in $\gsubpos_2$. Lemma~\ref{lem:represent-f}  provides the motivation for our analogues, as it
shows that we can define $\allpolyf$  in terms of linear
transformations. Recall its results:

  \begin{gather*}
  \allpolyf = \left\{ T(e^x)\mid T\colon{}\allpolypos\longrightarrow\allpoly\right\}    \\
  \allpolyposf = \left\{ T(e^x)\mid    T\colon{}\allpolypos\longrightarrow\allpolypos\right\}    
  \end{gather*}

We now  consider two variables. 

  \begin{definition}
    \begin{gather*}
      \image{2} = \left\{ T_\ast \,(f)\mid f\in\gsubplus_2 \wedge
        T\colon{}\allpolypos\longrightarrow\allpoly \wedge \text{$T$
          satisfies induction} \right\} \\      
      \imagepos{2} = \left\{ T_\ast \,(f)\mid f\in\gsubplus_2 \wedge
        T\colon{}\allpolypos\longrightarrow\allpolypos \wedge \text{$T$
          satisfies induction}  \right\} \\      
      \imagef{2} = \text{ closure of } \image{2}\\
      \imageposf{2} = \text{ closure of } \imagef{2}
    \end{gather*}

  \end{definition}

The following lemma enumerates basic properties of $\image{2}$.

\begin{lemma}
  Suppose $g = T_\ast f$ is in $ \image{2}$. Then
  \begin{enumerate}
  \item $\gsubplus_2\subsetneq \imagepos{2} \subsetneq \image{2}
    \subsetneq \partialpoly{1,1}$.
  \item $g^H(x,1)\in\allpoly$
  \item $g(x,\alpha)\in\allpolypos$ for $\alpha\ge0$.
  \item $\frac{\partial g}{\partial x}\in\image{2}$
  \item $\frac{\partial g}{\partial y}\in\image{2}$
  \item The coefficient of any power of $y$ is in $\allpolypos$.
  \item If $h(x)\in\allpolypos(n)$ the the homogenization $H=y^n
    h(x/y)$ is in $\image{2}$. 
  \item If $S\colon{}\allpolypos\longrightarrow\allpolypos$ then $S_\ast:\image{2}\longrightarrow\image{2}$.
  \end{enumerate}
  If we assume that $g\in\imagepos{2}$ then   parts $(4),(5)$, and $(8)$
  hold with $\image{2}$ replaced by $\imagepos{2}$ and  in $(2)$ we have
  that $g^H(x,1)\in\allpolypos$.
\end{lemma}
\begin{proof}
  If we take $T$ to be the identity then $\gsubplus_2\subseteq
  \image{2}$. The example below shows that this is a proper containment.
  Since $g^H = T_\ast(f^H)$ we see $g^H(x,1)\in\allpoly$. If
  $\alpha\ge0$ then $f(x,\alpha)\in\allpolypos$, and so $g(x,\alpha) =
  T_\ast f(x,\alpha) \in\allpolypos$.  If we define $S(h) =
  \frac{\partial}{\partial x}\,T(h)$ then
  $S\colon{}\allpolypos\longrightarrow\allpolypos$, and so $\frac{\partial
    g}{\partial x}\in\image{2}$.  Next, $\frac{\partial g}{\partial y} =
  T_\ast(\frac{\partial f}{\partial y})$, and so $\frac{\partial
    g}{\partial y}\in\image{2}$ since $\frac{\partial f}{\partial
    y}\in\gsubplus_2$. Since we can differentiate with respect to $y$,
  and substitute $y=0$, we see that all coefficients of powers of $y$
  are in $\allpolypos$. The coefficients of $H$ are all positive, and
  substitution is clear. The last one is immediate from the definition:
  $S_\ast(T_\ast f) = (ST)_\ast f$.

\end{proof}

\begin{example}
 Members of $\imagepos{2}$ do not necessarily satisfy substitution for
 negative $\alpha$.  Consider $T\colon{}x^n\mapsto \rising{x}{n}$ and
 $f=(x+y+1)^2\in\gsubplus_2$. If $g(x,y) = T_\ast \,f $ then
 $g\in\image{2}$, yet $g(x,-2)
 = x^2-x+1$ has complex roots.
\end{example}

\subsection{Constructing elements in $\image{2}$}

We can construct elements of $\image{2}$ by explicitly exhibiting them
in the form $T_\ast f$ as in the definition, or by a  general
construction. 

\begin{example}
  The generalized Laguerre polynomials $L_n(x;y)$ are in $\imagepos{2}$.
  This is a consequence of the surprising identity
  \eqref{eqn:laguerre-surprise}. If we define $T\colon{}x^n\mapsto
  \falling{x}{n}$ then
  \begin{equation}
    \label{eqn:laguerre-surprise}
    T_\ast^{-1}\,\rising{x+y+1}{n} = {n!} \,  L_n(x;y)
  \end{equation}
  Now $T^{-1}$ maps $\allpolypos$ to itself and
  $\rising{x+y+1}{n}\in\gsubplus_2$.  It follows from
  \eqref{eqn:laguerre-surprise} that $ L_n(x;y)\in\imagepos{2}$. I don't
  know a good proof for \eqref{eqn:laguerre-surprise}; there is an
  uninteresting inductive proof.

\end{example}

\begin{example} \index{Hermite polynomials}
  We claim $xy-1\in\image{2}$. First of all, the transformation
  $T\colon{}x^n\mapsto H_n(x/2)$ maps $\allpoly$ to itself, and satisfies
  $T(1)=1$, $T(x)=x$, $T(x^2)=x^2-2$. Now choose a positive integer
  $m$, and define $S(f) = T^m(f)$. Then
  $S\colon{}\allpoly\longrightarrow\allpoly$ and 
$$ S(1) = 1 \quad\quad S(x) = x \quad\quad S(x^2) = x^2 - 2m$$
Set $\epsilon = 1/(2m)$, and  define 
\begin{align*}
f_m &= S_\ast(x+ \epsilon y)(y+\epsilon x)\\
&= xy(1+\epsilon^2) + \epsilon y^2 +  \epsilon x^2 -1  
\intertext{then}
\lim_{m\rightarrow\infty} f_m &= xy-1
\end{align*}
A similar argument shows that $xy-\alpha\in\image{2}$ for any positive
$\alpha$. Notice that we are using a limit of transformations to show
that an element is in $\image{2}$, not a sequence of elements. 
\end{example}

\begin{example}
  The map $T\colon{}x^i\mapsto \rising{x}{i}$ satisfies
  $\allpolypos\longrightarrow\allpolypos$. If $f(x)\in\allpolypos$,
  then the homogenization of $f(x)$ is in $\gsubpos_2$. Consequently,
  if we define $S(x^i) = \rising{x}{i}y^{n-i}$ then
  $S(f)\in\imagepos{2}$.
\end{example}

\begin{example}

We can use transformations to get elements of $\image{2}$ or $\imagepos{2}$ from
elements in $\allpolypos$. Suppose that $T$ is a linear transformation
satisfying $T\colon{}\allpolypos\longrightarrow\allpolypos$, and $S$ is a
linear transformation then the composition below gives elements of
$\image{2}$ from elements of $\allpolypos$.

\centerline{\xymatrix{
    {\allpolypos}
      \ar@{->}[r]^{{ S  }}         
      &
      {\gsubpos_2} 
      \ar@{->}[r]^{{ T_\ast   }}         
      & {\image{2}}
}}

Useful choices of $S$ are 
\begin{align*}
  f(x) &\mapsto f(x+y) \\
  f(x) & \mapsto \text{ homogenization of } f(x)
\end{align*}
In the latter case the composition is $x^i\mapsto T(x^i)y^{n-i}$ where
$n$ is the degree of $f$. 

\end{example}

\subsection{The graphs of polynomials in $\image{2}$}

The graphs of polynomials in $\image{2}$ resemble polynomials in
$\rupint{2}$ in the upper left quadrant, since they satisfy substitution
for positive $y$, and their homogeneous part has all negative roots.
Here are two examples.  Figure~\ref{fig:laguerre2} shows the graph of
$L_5(x;y)$. In general, the Laguerre polynomials $L_n(x;y)$ have a
series expansion \eqref{eqn:laguerre-gen} \index{Laguerre polynomials}
that shows that $L_n(x;y)$ is a polynomial of degree $n$ in both $x$
and $y$, with all positive terms.  The leading term of the coefficient
of $x^k$ is
$$\frac{y^{n-k}}{(n-k)!k!}$$
and hence the homogeneous part of
$L_n(x;y)$ is $\frac{1}{n!}(x+y)^n$.  Consequently, the solution
curves are all asymptotic to lines with slope $-1$.  Also, $L_n(0;y) =
\binom{n+y}{n}$ implies that the graph of $ L_n(x;y)$ meets the
$y$-axis in $-1,-2,\dots,-n$.  All lines $y=\alpha x$ where $\alpha$
is positive will meet the curve in negative $y$ values, so
$L_n(x;\alpha x)\in\allpolypos$.
$\image{2}$ is closed under differentiation with respect to $y$, as
can also be seen in Figure~\ref{fig:laguerre2}.

\begin{figure}[htbp] \label{fig:tstar}
  \begin{center}
    \includegraphics*[width=2in]{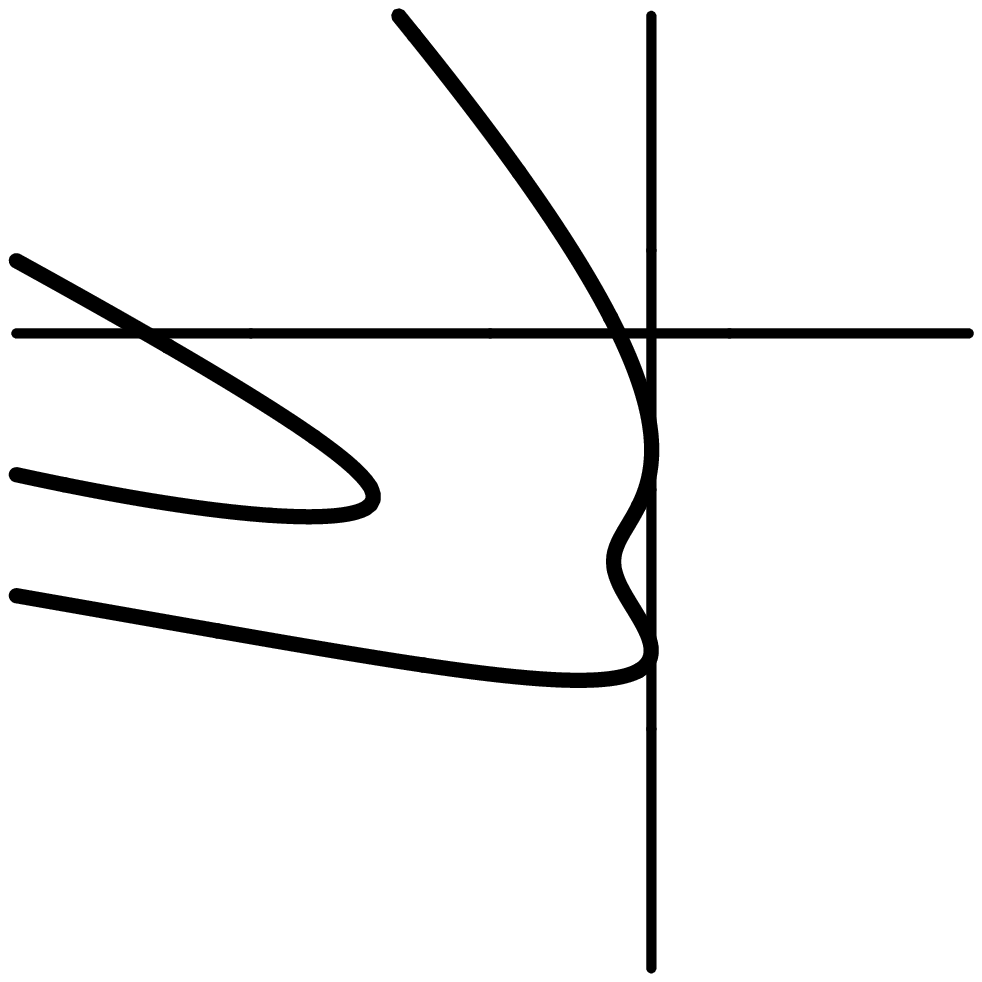}
    \caption{The graph of $T_\ast f$}
  \end{center}
\end{figure}

For the next example, we take $T\colon{}x^n\mapsto \rising{x}{n}$, and pick
(at random)
$$f=(x + y + 1)(x + 2y + 2)(x + 3y + 1)(x + 5y + 2).$$
Figure~\ref{fig:tstar} shows the graph of $T_\ast f$.  In
Figure~\ref{fig:laguerre2} the solution curves all meet the $y$
axis, but this is not true in Figure~\ref{fig:tstar}.
Also, the line $y=x$ does not intersect all the solution curves, and
so $f(x,x)\not\in\allpoly$. 


If we pick our slopes appropriately, then we can intersect all
solution curves.

\begin{lemma}
  Suppose that $f(x,y)\in\imagepos{2}$, and that $\alpha$ is the largest
  root of $f^H$. If $0>\beta>\alpha$ then $f(x,\beta x)\in\allpolypos$.
\end{lemma}
\begin{proof}
  The line $y=\beta x$ meets the $x$-axis to the right of intersection
  points of the solution curves, and is eventually less than every
  solution curve, so it must meet every one.
\end{proof}

\subsection{A partial \Polya-Schur result}

Next, we look at two actions on $\image{2}$ by polynomials.

\begin{lemma}
  Suppose that $g\in\image{2}$ and $h\in\allpolypos$. Then
  \begin{enumerate}
  \item $h(y)g(x,y)\in\image{2}$.\\[.2cm]
  \item $h(\frac{\partial}{\partial x})g(x,y)\in\image{2}$.\\[.2cm] 
  \item $h(\frac{\partial}{\partial y})g(x,y)\in\image{2}$. 
  \end{enumerate}
If $g\in\imagepos{2}$ then all the statements are true with
$\image{2}$ replaced by $\imagepos{2}$. 
\end{lemma}
\begin{proof}
  Assume that $g=T_\ast f$ as in the definition. If
  $h(y)\in\allpolypos$, then $h(y)\in\gsubposf_2$, so
  $h(y)f(x,y)\in\allpolyposf_2$. The first part now follows from
$$ T_\ast h(y)f(x,y) = h(y) T_\ast f(x,y) \ = \  h(y) g(x,y)$$
For the second part, it suffices to show that $g+\alpha g'\in\image{2}$
if $\alpha\ge0$. This follows from the fact that $k\mapsto k+\alpha k'$
maps $\allpolypos$ to itself. To verify the third part, if suffices to
show that $g+\alpha g_y\in\image{2}$ for positive $\alpha$. But
$T_\ast(f+\alpha f_y) = g+ \alpha g_y$ since $T_\ast$ commutes with
$\frac{\partial}{\partial y}$, and $f+\alpha f_y\in\gsubpos_2$. 
\end{proof}

We can use this result to go from an element of $\image{2}$ to a
transformation.

\begin{prop}
  If $T_\ast(e^{xy})\in\imageposf{2}$ then
  $T\colon{}\allpolypos\longrightarrow\allpolypos$. In other words, 
  suppose $g(x,y)=\sum g_i(x)\frac{y^i}{i!}\in\imageposf{2}$. The linear transformation $T\colon{}x^i\mapsto g_i$
  defines a map $\allpolypos\longrightarrow\allpolypos.$
\end{prop}
\begin{proof}
  If we choose $h(x) = \sum_0^n a_i x^i\in\allpolypos$ then  The coefficient of
  $y^n$ in $h^{rev}(\frac{\partial}{\partial y})g(x,y)$ is
$$ \sum a_i g_i(x) = T(h)$$
Since $h^{rev}(\frac{\partial}{\partial y})g(x,y)\in\imagef{2}$, we know that all coefficients are
in $\allpolypos$. 
\end{proof}

\section{Polynomials in several variables with  positive coefficients}
\label{sec:wagner}

In this section we look at the set $\pospoly_{d,e}$ of certain
polynomials in $d+e$ variables that have non-negative coefficients.
The motivation for this set comes from considering these three properties of
$\rupint{2}$. Assume that $f(x,y)\in\gsub_2$.
\begin{enumerate}
\item $f(x,\alpha)\in\allpoly$ for all $\alpha$.
\item $f(\alpha,x)\in\allpoly$ for all $\alpha$.
\item $f(x,\alpha\, x)\in\allpoly$ for all positive $\alpha$.
\end{enumerate}


$\pospoly_{d,e}$ is constructed by starting with $\gsubplus_d$ and
adding $e$ new variables so that the resulting polynomials satisfy
certain substitution conditions.  We next define interlacing, and
verify that it obeys the expected properties. The graphs of
polynomials in $\pospoly_{0,2}$ are more complicated than those in
$\gsubplus_2$, but their behavior in each quadrant is explained by
their defining conditions.  We find some general constructions for
polynomials in $\pospoly_{d,e}$ that are not necessarily in
$\gsubplus_{d+e}$.  We then look at a special class of polynomials in
$\pospoly_{d,e}$ - those that are linear in the last $e$ variables.
They can be identified with collections of polynomials that are
defined on the corners of a cube.  The three variable case of these
polynomials leads to matrices that preserve interlacing.

\subsection{Definitions }
\label{sec:wag-defn}

In order to better motivate our definition, we begin by restating the
definition for interlacing mentioned above.  Let $F(x,y) = f(x) + y
g(x)$. If $F(x,\alpha)$ and $F(x,\alpha x)$ are in $\allpolypos$ for
all non-negative $\alpha$ then $f$ and $g$ interlace. This suggests our
first definition:

\begin{definition}
  $\pospoly_{1,1}$ is the set of all polynomials $F$ in $x,y$ with non-negative
  coefficients that satisfy the two conditions

  \begin{align*}
    F(x,\alpha) \in \allpolypos & \text{ for all non-negative } \alpha \\
    F(x,\alpha x) \in \allpolypos & \text{ for all non-negative } \alpha 
  \end{align*}
\end{definition}

The general definition extends polynomials in $\gsubplus_d$ rather
than polynomials in $\allpolypos$. We let 
$\xx=(x_1,\dots,x_d)$, $\yy=(y_1,\dots,y_e)$, and 
$\yy^\aaa_\bbb = (\alpha_1x_{b_1}^{c_1},\dots,\alpha_e x_{b_e}^{c_e})$. 

\begin{definition}
  $\pospoly_{d,e}$ is the set of all polynomials $f(\xx,\yy)$ in $d+e$
  variables with non-negative coefficients such that for all non-negative
  $\alpha_1,\dots,\alpha_e$, all $b_i$ which satisfy $b_i\in\{1,2,\dots,d\}$
  and all $c_1,\dots,c_e$ which are $0$ or $1$  we have
  that$ f(\xx,\yy^\aaa_\bbb)\in\gsubplus_d.$
  
  In other words, if we substitute either a non-negative constant or a
  non-negative constant multiple of some $x$ variable for each $y$
  variable then the resulting polynomial is in $\gsubplus_d$. If $d$
  is zero, then we substitute either a non-negative constant or a
  non-negative multiple of $x$ for each $y$ variable.

\end{definition}

Notice that the conditions for a polynomial to be in $\pospoly_{1,2}$
are a subset of those that define $\pospoly_{0,3}$. Consequently, 
$\pospoly_{0,3} \subset \pospoly_{1,2}$. More generally we have

\begin{lemma}
  If $d\ge2$ and $e\ge0$ then
  $\pospoly_{d,e}\subset\pospoly_{d-1,e+1}$ and $\pospoly_{0,e+1}\subset\pospoly_{1,e}$.
\end{lemma}
\begin{proof}
  Since a polynomial in $\pospoly_{0,e+1}$ satisfies all the
  conditions of $\pospoly_{1,e}$ plus a few more involving
  substitutions for the first variable, we have
  $\pospoly_{0,e+1}\subset\pospoly_{1,e}$.
  If $d\ge2$ then choose $f\in\pospoly_{d,e}$ and consider the
  substitution conditions that must be met for $f$ to be in
  $\pospoly_{d-1,e+1}$:
  \begin{align}
    \label{eqn:wag-eq-1}
    f(x_1,\dots,x_{d-1},\alpha_1,\alpha_2x_{b_2}^{c_2},\dots,\alpha_{e+1}x_{b_{e+1}}^{c_{e+1}})&
 \in\gsubplus_{d-1}\\   
   f(x_1,\dots,x_{d-1},\alpha_1x_{b_1},\alpha_2x_{b_2}^{c_2},\dots,\alpha_{e+1}x_{b_{e+1}}^{c_{e+1}})&
   \in\gsubplus_{d-1}\\
\intertext{\eqref{eqn:wag-eq-1} holds since }
   f(x_1,\dots,x_{d-1},x_d,\alpha_2x_{b_2}^{c_2},\dots,\alpha_{e+1}x_{b_{e+1}}^{c_{e+1}})&
   \in\gsubplus_{d}\notag
  \end{align}
  and substitution of a constant for $x_d$ gives a polynomial in
  $\gsubplus_{d-1}$. \eqref{eqn:wag-2} holds because we can substitute
  any positive multiple of a variable for a variable in $\gsubplus_{d}$, and the
  result is in $\gsubplus_{d-1}$.
\end{proof}

$\pospoly_{d,e}$ is not empty since
$\gsubplus_{d+e}\subset\pospoly_{d,e}$.  We allow some coefficients to
be zero, so  $\pospoly_{d,e} \subset \pospoly_{d,e+1}$, and therefore
we have the containments

\xymatrix@-1pc{
&&& & &{\allpolypos= \mathcal{Q}}_{1,0} \ar@{}[d]_\subset  & = & {\pospoly_{0,1}} \ar@{}[d]_\subset \\
&& & {\gsubplus_2=\pospoly_{2,0}} \ar@{}[d]_\subset & \subset & {\pospoly_{1,1}}
\ar@{}[d]_\subset & \supset & {\pospoly_{0,2}}  \ar@{}[d]_\subset \\
& {\gsubplus_3=\pospoly_{3,0}} & \subset & {\pospoly_{2,1}} & \subset & {\pospoly_{1,2}}& \supset & {\pospoly_{0,3}}
}

\begin{example}
If we choose $f_i\lesslesseq g_i\in\allpolypos$ then it follows from
the results in the next section that
$$\prod_{i=1}^n (f_i+yg_i) \in\pospoly_{1,1}$$

 The polynomial $f(x,y)=x+y^2$ is in $\pospoly_{1,1}$ but is not in
 $\pospoly_{0,2}$ since $f(1,y)\not\in\allpoly$.  The polynomial 
\begin{equation}\label{eqn:wag-poly-bad}
f = 2x + xy + (1+3x)z + yz
\end{equation}
is in $\pospoly_{1,2}$ but is not in $\pospoly_{0,3}$.  To verify this
we just make the four substitutions that define $\pospoly_{1,2}$, and
check that in each case the polynomials have all real roots.  It is
not in $\pospoly_{0,3}$ since $f(.01,x,.01x)$ has imaginary roots.


\end{example}

\begin{example}\label{ex:wag-3}
  Suppose that $f\lesslesseq g$ where $f,g\in\allpolypos$ have
  positive leading coefficients. We know that $(xg,f,g)$ is totally
  interlacing. If $\alpha$ is positive then the matrix below on the
  left is totally non-negative
$$ 
\begin{pmatrix}
  1 & \alpha & 0 \\ 0 & 1 & \alpha
\end{pmatrix}
\begin{pmatrix}
  xg \\ f \\g
\end{pmatrix}
=
\begin{pmatrix}
  xg+\alpha f\\ f+\alpha g
\end{pmatrix}
$$
and hence
\begin{equation}
  \label{eqn:wag-3}
  xg+\alpha f \longleftarrow f+\alpha g
\end{equation}
If we define
\begin{equation}
  \label{eqn:wag-4}
  h(x,y,z) = xg+zf+y(f+zg)
\end{equation}
then \eqref{eqn:wag-3} shows that $h(x,y,\alpha)\in\pospoly_{1,1}$ for
positive $\alpha$. Also, replacing $\alpha$ by $1/\alpha$ and multiplying
by $\alpha$ yields $\alpha xg+f\longleftarrow \alpha f+g$. Multiplying
the right by $x$ shows that
$$ xg+(\alpha x)f \longleftarrow f + (\alpha x)g$$
which implies $h(x,y,\alpha x)\in\pospoly_{1,1}$. Consequently,
$h(x,y,z)\in\pospoly_{1,2}$. 
\end{example}

\subsection{Elementary properties}
\label{sec:wag-elem}

We define interlacing, and then show that 
$\pospoly_{d,e}$  satisfies the expected closure properties.
If we consider the definition of $\pospoly_{1,1}$, we see that we can
define interlacing in $\allpolypos$ as follows:
\begin{quote}
  $f\lesslesseq g$ in $\allpolypos$ if and only if $deg(f)>deg(g)$ and
\\  $f + y g \in\pospoly_{1,1}$.
\end{quote}

If $f,g\in\pospoly_{d,e}$ then we say that $f \lesslesseq g$ if and
only if the total $x$-degree of $f$ is greater than that of $g$, and
$f + y_{e+1} g\in\pospoly_{d,e+1}$. All this means is that
$f\lesslesseq g$ if and only if
\begin{align*}
  \text{total $x$-degree}(f) & >   \text{total $x$-degree}(g)\\
  f(\xx,\yy) + \alpha g(\xx,\yy) & \in\pospoly_{d,e}\\ 
  f(\xx,\yy) + \alpha x_i g(\xx,\yy) & \in\pospoly_{d,e} \text{ (for any
  $1\le i \le d$) }
\end{align*}

We use $\lesslesseq$ to define $\longleftarrow$:
\begin{quote}
  $f\longleftarrow g$ iff we can write $g=\alpha f+h$ where
  $\alpha\ge0$ and $f\lesslesseq h$.
\end{quote}
In order to show that $f \lesslesseq g$ all we need to do is to show
that all substitutions interlace.
\begin{lemma} \label{lem:wag-basic-int}
  Suppose that $f,g\in\pospoly_{d,e}$, and that the $x$-degree of $f$
  is greater than the $x$-degree of $g$. The following are equivalent
\begin{enumerate}
\item $f(\xx,\yy^\aaa_\bbb) \lesslesseq g(\xx,\yy^\aaa_\bbb) \text{ in
    $\gsubplus_d$ for all appropriate $\aaa$ and $\bbb$}$
\item    $f(\xx,\yy) + y_{e+1}g(\xx,\yy)  \in \pospoly_{d,e+1}$
\item $    f(\xx,\yy)  \lesslesseq g(\xx,\yy)$ in $\pospoly_{d,e}$.
  \end{enumerate}
\end{lemma}
\begin{proof}
  The second and third are equivalent since they are the definition of
  interlacing. To show that $(3)$ implies $(1)$ note that interlacing
  implies that $f+y_{e+1}g\in\pospoly_{d,e+1}$. We can substitute for
  $y_{e+1}$ and conclude that
  \begin{align*} 
    f(\xx,\yy^\aaa_\bbb) + \alpha \, g(\xx,\yy^\aaa_\bbb) & \in\gsubplus_d \\
f(\xx,\yy^\aaa_\bbb) + \alpha x_i \, g(\xx,\yy^\aaa_\bbb) & \in\gsubplus_d \\
  \end{align*} 
and the conclusion follows from  Lemma~\ref{lem:pdplus-int}. That $(1)$ implies
$(3)$ is similar.
\end{proof}

\begin{lemma} \label{lem:wag-closure}
  Assume that $f,g,h\in\pospoly_{d,e}$.
  \begin{enumerate}
  \item  $fg\in\pospoly_{d,e}$.
  \item If $f\longleftarrow g$ then $hf \longleftarrow hg$.
  \item If $f\longleftarrow g$ and $f\longleftarrow h$ then
    $g+h\in\pospoly_{d,e}$ and $f \longleftarrow g+h$.
  \item If $g\longleftarrow f$ and $h\longleftarrow f$ then
    $g+h\in\pospoly_{d,e}$ and $g+h \longleftarrow f$.
  \item If $\alpha_1,\dots,\alpha_d$ and $\beta_1,\dots,\beta_e$ are
    non-negative then \\
    $f(\alpha_1 x_1,\dots,\alpha_d x_d,\beta_1 y_1,\dots,\beta_e
    y_e)\in\pospoly_{d,e}$.
  \end{enumerate}
\end{lemma}

\begin{proof}
  We first consider the case where all ``$\longleftarrow$'' are
  ``$\lesslesseq$''.  The first one is immediate from the observation
  that $ fg(\xx,\yy^\aaa_\bbb) = f(\xx,\yy^\aaa_\bbb) \, g(\xx,\yy^\aaa_\bbb)$ The
  second follows from the first since $hf + y_{e+1}hg =
  h(f+y_{e+1}g)$, and $f+y_{e+1}g\in\pospoly_{d,e+1}$ by definition of
  interlacing.

For assertion $(3)$, in order to verify that $f +
y_{e+1}(g+h)\in\pospoly_{d,e+1}$ we use the properties of interlacing
in $\gsubplus_d$. Since $f\lesslesseq g$, we know that 
by Lemma~\ref{lem:wag-basic-int} 

\begin{align*}
  f(\xx,\yy^\aaa_\bbb) & \lesslesseq g(\xx,\yy^\aaa_\bbb) \\
  f(\xx,\yy^\aaa_\bbb) & \lesslesseq h(\xx,\yy^\aaa_\bbb) \\
\intertext{ and adding these  interlacings together gives} 
f(\xx,\yy^\aaa_\bbb) & \lesslesseq  g(\xx,\yy^\aaa_\bbb) + h(\xx,\yy^\aaa_\bbb)
\end{align*}

It follows that $f+ y_{e+1}(g+h)\in\pospoly_{d,e+1}$.  The proof of the
next statement is similar. The last statement is immediate from the
definition.

In the general case $f\longleftarrow g$ we use the representation
$g=\alpha f+k$ where $\alpha\ge0$ and $f\lesslesseq h$. For instance,
to establish $(2)$, we know that $f\lesslesseq k$ implies
$hf\lesslesseq hk$, and therefore $hg=\alpha hf + hk$ implies
$hf\longleftarrow hg$.
\end{proof}


Factors of polynomials in $\pospoly_{d,e}$ are generally in $\pospoly_{d,e}$.

\begin{lemma}
  If $fg\in\pospoly_{1,e}$, and $f$ has only non-negative coefficients
  then $f\in\pospoly_{1,e}$.  If $fg\in\pospoly_{d,e}$, and $f$ has
  only non-negative coefficients and $f(\xx,\yy^\aaa_\bbb)\in\gsubplus_d$
  for all appropriate $\yy^\aaa_\bbb$ then $f\in\pospoly_{d,e}$.

\end{lemma}
\begin{proof}
  If some substitution resulted in a complex root for $f$, then it
  would be a complex root for $fg$ as well. Since $f$ has non-negative
  coefficients, $f$ is in $\pospoly_{d,e}$.
\end{proof}

\subsection{Constructing polynomials in $\pospoly_{1,2}$}
\label{sec:wag-construct}

We can easily construct polynomials in $\pospoly_{1,2}$ using
products. If $f\lesslesseq g$, $h\lesslesseq k$ are polynomials in
$\allpolypos$ then $f+ yg\in\pospoly_{1,2}$, $h+zk\in \pospoly_{1,2}$,
and therefore $(f+yg)(h+zk)\in\pospoly_{1,2}$. However, these
polynomials are also in $\gsubcloseplus_3$.  We do not know if
$\pospoly_{1,2}$ is different from $\gsubcloseplus4_3$.

\begin{lemma}
  \label{lem:wag-construct}
  Suppose that $f,g,h,k\in\allpolypos$ have all positive coefficients.
  If $\smallsquare{f}{g}{h}{k}$ where all interlacings are
  ``$\lessless$'' or all are ``$\lessless$'' and the determinant
  $\smalltwodet{f}{g}{h}{k}$ has no negative roots then
$$ f + y g + z h + yz k \in\pospoly_{1,2}$$
\end{lemma}
\begin{proof}
  It suffices to show that $f+\alpha g \lessless k+\alpha h$ and
  $f+\alpha x g \lessless k+\alpha x h$ for all positive $\alpha$. The
  interlacing assumptions imply that all four linear combinations are
  in $\allpolypos$. We now follow the proof of Lemma~\ref{lem:inequality-4b} to
  show that the first interlacing holds since the determinant in the
  conclusion is never zero for negative $x$. Again by following the
  proof of Lemma~\ref{lem:inequality-4b}, the second holds since
  $\smalltwodet{f}{xg}{h}{xk} = x \smalltwodet{f}{g}{h}{k}$.
\end{proof}

Note that Example~\ref{ex:wag-3} is a special case of the Lemma. Take 
$\smallsquare{xg}{f}{f}{g}$. Since
$\smalltwodet{xg}{f}{g}{f}=xg^2-f^2$ clearly has no negative roots,
the Lemma applies.

It is easy to find polynomials in $\pospoly_{1,2}$ using
Lemma~\ref{lem:wag-construct}.  This is a special case of a more
general construction in the next section.

  \begin{lemma}
    Choose $w(x,y,z)\in\gsubplus_{3}$, and write
  $$
  w(x,y,z) = f(x) + g(x)y + h(x)z + k(x)yz + \cdots $$ 
Then $f(x) + g(x)y + h(x)z + k(x)yz \in \pospoly_{1,2}$.
   \end{lemma}
   \begin{proof}
     $f,g,h,k$ are all in $\rupint{2}$, and Theorem~\ref{thm:product-4}
     shows that the determinant is never zero.
   \end{proof}

  For instance, if $v(x)\in\allpolypos$, then the Taylor series
  expansion of $v(x+y+z)$ shows that
$$ v(x)+ v'(x)\cdot(y+z) + \frac{v''(x)}{2}yz\in\pospoly_{1,2}$$

\subsection{Polynomials linear in the $y$ variables}
\label{sec:wag-linear}

An interesting class of polynomials in $\pospoly_{d,e}$ consists of
those polynomials that are linear in the $y$ variables. In general, if
$f$ is a polynomial in $\pospoly_{d,e}$ that is linear in
$y_1,\dots,y_d$ then we can represent $f$ by labeling the vertices of
a $d$-cube with appropriate coefficients of $f$.

We can find such polynomials in $\pospoly_{d,e}$ by taking initial terms of
polynomials in $\gsubplus_{d+e}$.

\begin{lemma}
  Suppose $f(\xx,\yy)\in\gsubplus_{d+e}$, and let $g(\xx,\yy)$ consist
  of all terms of $f$ that have degree at most one in every $y$
  variable. Then $g\in\pospoly_{d,e}$.
\end{lemma}
\begin{proof}
  For simplicity we illustrate the proof with $e=2$; the general proof
  is the same. If we write
\begin{align*}
    f(\xx,y,z) &= f_{00}(\xx) + f_{10}(\xx)y + f_{01}(x)z +
    f_{11}(x)yz + \cdots \\ 
    g(\xx,y,z) &= f_{00}(\xx) + f_{10}(\xx)y + f_{01}(x)z +
    f_{11}(x)yz  \\ 
\intertext{then the coefficient of $yz$ in $(y+x_i^r\alpha)(z+x_j^s\beta)f(\xx,y,z)$ is}
g(\xx,x_i^r\alpha,x_j^s\beta)&= f_{00}(\xx) + f_{10}(\xx) x_i^r\alpha + f_{01}(x)x_j^s\beta +
    f_{11}(x)x_i^rx_j^s\alpha\beta  
  \end{align*}
Consequently, we see that since $\alpha,\beta$ are positive,
$$g(\xx,\alpha,\beta), g(\xx,x_i\alpha,\beta), g(\xx,\alpha,x_j\beta),
g(\xx,x_i\alpha,x_j\beta)$$ are all in $\gsubplus_d$, and so
$g(\xx,y,z)\in\pospoly_{d,e}$. 
\end{proof}

\begin{lemma} \label{lem:wag-reverse}
  Suppose $f(x,\yy)\in\pospoly_{1,e}$ and has $x$ degree $r$ and $y_i$
  degree $s_i$. The reverse polynomial below is in $\pospoly_{1,e}$:
  $$
  x^r y_1^{s_1}\cdots y_e^{s_e}\,f(
  \frac{1}{x},\frac{1}{y_1},\dots,\frac{1}{y_e})
$$
\end{lemma}

\begin{proof}
  Any substitution for $y_i$ by a constant yields a polynomial in some
  $\pospoly_{1,e^\prime}$ where $e^\prime<e$. We can also scale each
  $y_i$, so we only need to consider the substitution $y_i\rightarrow
  x$. Thus we need to verify that
\begin{equation}\label{eqn:wag-2}
 x^{r+s_1+\cdots+s_e} f(\frac{1}{x},\dots,\frac{1}{x})\in\allpolypos
\end{equation}
If a monomial in $f(x,\yy)$ is $a_\diffi x^i\yy^\diffi$ with total
degree $i+|\diffi|$, then in \eqref{eqn:wag-2} the degree is
$n-(i+|\diffi|)$ where $n = r+s_1+\cdots+s_e$. Consequently,
\eqref{eqn:wag-2} is the usual reverse of $f(x,x,\dots,x)$ which we
know to be in $\allpolypos$.
\end{proof}

For instance, if $d=1$ then Lemma~\ref{lem:wag-reverse} shows that if
\begin{align*}
   f_{11}(x) +y\, f_{21}(x) + z\,f_{12}(x) + yz\, f_{22}(x)
   &\in\pospoly_{1,2}\\
\intertext{then the reverse is also in $\pospoly_{1,2}$}
   yz\,f_{11}(x) + z\, f_{21}(x) + y\,f_{12}(x) + \, f_{22}(x) &\in\pospoly_{1,2}
\end{align*}

\subsection{Matrices preserving interlacing}
\label{sec:wag-preserve}

Matrices determined by polynomials linear in $y$ and $z$ are closed
under multiplication. We can use this fact to show that these matrices
also preserve interlacing. This is a special case of
Theorem~\ref{thm:multiply-vectors}.  
\index{matrix!preserving interlacing}

\begin{lemma}
  Suppose that $f,g\in\pospoly_{d,2}$  are linear in $y$
  and $z$ and write
  \begin{align*}
    f & = f_{11} + y f_{12} + z f_{21} + yz f_{22} \\
    g & = g_{11} + y g_{12} + z g_{21} + yz g_{22} \\
    \intertext{Represent $f,g$ by matrices in the following way}
    F &= \begin{pmatrix}{f_{12}}&{f_{22}}\\{f_{11}}&{f_{21}}\end{pmatrix}\\
    G &= \begin{pmatrix}{g_{12}}&{g_{22}}\\{g_{11}}&{g_{21}}\end{pmatrix}\\
\intertext{If the product is $ FG = \begin{pmatrix}{h_{12}}&{h_{22}}\\{h_{11}}&{h_{21}}\end{pmatrix}$}
\intertext{then  $h_{11} + y h_{12} + z h_{21} + yz h_{22} \in\pospoly_{d,2}$}
  \end{align*}
\end{lemma}

\begin{proof}
  The proof is an application of the Leibnitz rule \index{Leibnitz} in
  disguise.
We know that
\begin{align*}
  f_{11} + y f_{12} & \lesslesseq f_{21} + y f_{22} \\
  g_{11} + z g_{21} & \lesslesseq g_{12} + z g_{22} \\
\intertext{and multiplying by the right hand sides }
(f_{11} + y f_{12})(g_{12} + z g_{22}) & \lesslesseq (f_{21} + y
f_{22})(g_{12} + z g_{22}) \\
(f_{21} + y f_{22})( g_{11} + z g_{21})  & \lesslesseq (f_{21} + y
f_{22})( g_{12} + z g_{22}) \\
\intertext{and adding them together yields}
(f_{11} + y f_{12})(g_{12} + z g_{22}) + (f_{21} + y f_{22})( g_{11} +
z g_{21}) & \lesslesseq (f_{21} + y  
f_{22})(g_{12} + z g_{22}) \\
\intertext{The matrix representation of the left hand polynomial is exactly $FG$.}
\end{align*}
\end{proof}

\begin{cor}
  If $f\in\pospoly_{1,2}$ is  given as in the lemma, $g_1\longleftarrow
  g_2$ and we set 
  $$
  \begin{pmatrix}{f_{12}}&{f_{22}}\\{f_{11}}&{f_{21}}\end{pmatrix}\,
  \begin{pmatrix}{g_2}\\{g_1}\end{pmatrix} = 
  \begin{pmatrix}{h_2}\\{h_1}\end{pmatrix}
$$
  then
  $h_1\longleftarrow h_2$.
\end{cor}

\begin{proof}
  
  We check that $g = (g_1+ y g_2)(1+z)$ is in $\pospoly_{1,2}$. Since
  $g_1\longleftarrow g_2$ we know that $g_1+yg_2\in\pospoly_{1,1}$.
  Multiplying by $1+z$ shows $g\in\pospoly_{1,2}$.
  
  The matrix of $g$ is $\smalltwo{g_2}{g_2}{g_1}{g_1}.$ Applying the
  lemma and considering only the first column of the product gives the
  conclusion.
\end{proof}

We can use the determinant condition of Lemma~\ref{lem:wag-construct} to find
such matrices.

\begin{cor}
  Suppose that $f,g,h,k\in\allpolypos$ have all positive coefficients.
  If $\smallsquare{f}{g}{h}{k}$ where all interlacings are $\lessless$
  and the determinant $\smalltwodet{f}{g}{h}{k}$ has no negative roots
  then the matrix $\smalltwo{h}{k}{f}{g}$ preserves interlacing for
polynomials in $\allpolypos$.
\end{cor}

\begin{example}
  The determinant of the polynomials in the interlacing square below
  has only positive roots, so the corresponding matrix preserves
  interlacing. If we multiply this matrix by $\smalltwo{0}{x}{1}{0}$
  which preserves interlacing, we find that the matrix of polynomials
  below preserves interlacing for polynomials in $\allpolypos$.

\xymatrix@-1pc{
&&&& 4+88x \ar@{->}[rrr]^{\lessgreat} &&& 8 + 260x + 216x^2 \\
&&&& 32+64x \ar@{->}[rrr]_{\lesseq} \ar@{->}[u]^{\greateq} &&&
 64 + 416 x \ar@{->}[u]_{\lessless}
}

\vspace*{.2cm}
 
$$
\begin{pmatrix}
32x+64x^2 &  64x + 416 x^2 \\
4+88x     &  8 + 260x + 216x^2
\end{pmatrix}
$$

\end{example}

We can not give a characterization of the matrices that preserve
interlacing, but we can find a few restrictions. We begin with 
the linear case.

\begin{lemma}
  Suppose $f,g\in\pospoly_{1,2}$ have the property that for all
  $p\longleftarrow q\in\pospoly_{1,2}$ it holds that
  $fp+gq\in\pospoly_{1,2}$. Then $g\longleftarrow f$.
\end{lemma}
\begin{proof}
  If we choose $p = \alpha q$ for positive $\alpha$ then
  $p\lesslesseq q$, and so $p(\alpha f + g)\in\pospoly_{1,2}$. Since
   factors of polynomials in $\pospoly_{1,2}$ are in $\pospoly_{1,2}$, we
  conclude that $\alpha f + g\in\pospoly_{1,2}$. Next, if we choose $p =
  \alpha x q$ then again $p\lesslesseq q$, and so $ x f +
  g\in\pospoly$. It follows that $g\lesslesseq f$.
\end{proof}

\begin{lemma}
  Suppose that 
  $$
  \begin{pmatrix}{f_{12}}&{f_{22}}\\{f_{11}}&{f_{21}}\end{pmatrix}\,
  \begin{pmatrix}{g_2}\\{g_1}\end{pmatrix} = 
  \begin{pmatrix}{h_2}\\{h_1}\end{pmatrix}
  $$
  and the $f_{ij}$ have the property that whenever
  $g_1\lesslesseq g_2$ in $\pospoly_{1,2}$ then it holds that
  $h_1\lesslesseq h_2$ in $\pospoly_{1,2}$. Then, the $f_{ij}$ are in
  $\pospoly_{1,2}$ and form an
  interlacing square
\index{interlacing square} \index{square!interlacing}

\raisebox{0.5\depth}{\xymatrix@-1pc{
f_{12} \ar@{<-}[r]&  f_{22} \ar@{->}[d] \\
f_{11} \ar@{<-}[u]&  f_{21} \ar@{->}[l] 
}}

\end{lemma}
\begin{proof}
  The horizontal arrows follow from the previous lemma. If we take a
  limit, we may assume that $g_1$ is zero. Thus, we have that
  $f_{11}g_2 \lesslesseq f_{12}g_2$ which implies that
  $f_{11}\lesslesseq  f_{12}$ and that $f_{11},f_{12}$ are in
  $\pospoly_{1,2}$. Taking $g_2$ to be zero gives the properties of the
  second column.
\end{proof}

\section{Epsilon substitution}
  We have looked at properties of polynomials $f(x,y)$ where
  $f(x,\alpha)$ is guaranteed to be in $\allpoly$ for only positive
  $\alpha$. We now reduce the range of acceptable $\alpha$.

\newcommand{\eppoly}[1]{\mywp^{\epsilon}_{#1}}

  \begin{definition}
$\eppoly{d}$ consists of those polynomials $f(x;y_1,\dots,y_d)$ with
the property that there is an $\epsilon$, depending on $f$, such that 
$$ |\alpha_1|\le\epsilon,\dots,|\alpha_d|\le\epsilon\quad
\implies\quad f(x;\alpha_1,\dots,\alpha_d)\in\allpoly
$$
  \end{definition}

\begin{figure}[htbp]
  \centering
 \includegraphics*[width=2in]{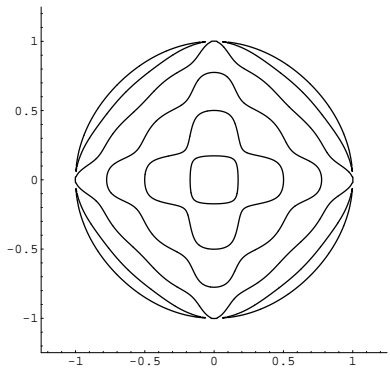}
  \caption{The graph of $p_{6,6,12}$}
  \label{fig:h6612}
\end{figure}

It's easy to tell from the graph if a polynomial is in
$\eppoly{1}$. For instance,  Figure~\ref{fig:h6612} is the graph of a
polynomial of degree $12$. We can see that
every horizontal line $y=\alpha$ intersects the graph in $12$ points
if $|\alpha|$ is sufficiently small, and consequently
$f\in\eppoly{1}$. 

Here are a few simple facts that need no proof.

\begin{enumerate}
\item If $f,g\in\eppoly{d}$ then $fg\in\eppoly{d}$.
\item If $f\in\eppoly{d}$ then $\frac{\partial f}{\partial x}
  \in\eppoly{d}$.
\item If $f\in\eppoly{d}$ then $f(x;0,\dots,0)\in\allpoly$.
\end{enumerate}

Let's consider two examples.
\begin{example} 
 We will show that
    $\sum_{i=0}^n
    P_i(x)\binom{n}{i}y^{n-i}\in\eppoly{1}$ where $P_i$
    is the Legendre polynomial. 
We make use of the formula \cite{szego}
\begin{equation} \label{eqn:legendre-oval}
 \sum_{i=0}^n\binom{n}{i}P_i(x)y^{n-i} =
(1+2xy+y^2)^{n/2}P_n((x+y)(1+2xy+y^2)^{-1/2})
\end{equation}


We show that we can take $\epsilon=1$. 
If $P_n(\beta)=0$, then we know that $|\beta|<1$. Solving
$$\beta = (x+y)(1+2xy+y^2)^{-1/2}$$
for $x$ yields
$$x = -y + y\,{\beta }^2 - {\sqrt{{\beta }^2 - y^2\,{\beta }^2 +
    y^2\,{\beta }^4}}$$
Now if $|y|<1$, then $|x|<1$.  In this case not only is
\eqref{eqn:legendre-oval} in $\eppoly{1}$, but we have the stronger
result that
$$ |\alpha|\le 1 \quad\implies\quad \sum_{i=0}^n
    P_i(x)\binom{n}{i}\alpha^{n-i}\in\allpolyint{(-1,1)}
$$

\end{example}

\begin{example}

Define the polynomials of degree $n$ 
\begin{gather*}
  f_n(x,y) = \frac{\partial^n}{\partial y^n}\left(x^2+y^2-1\right)^n\\
  h_n(y) = \left(\frac{d}{dy}\right)^n (y^2-1)^n
  \intertext{Figure~\ref{fig:leg66a} is the graph of $f_6$. We have
    the relationship} f_n(x,y) = 0 \text{ iff }
  h_n\left(\frac{y}{\sqrt{1-x^2}}\right)=0
\end{gather*}

Note that $h$ is, up a factor, the Legendre polynomial, and so has
the same roots.  Since $f_n(0,y)=h_n(y)$ we see that the intersection
of the graph of $f_n$ with the $y$ axis is at the roots of the
Legendre polynomial. As $n\rightarrow\infty$, these roots become dense
in $(-1,1)$.

The relationship between $f_n$ and $h_n$ shows that the graph of $f_n$
consists of ovals through the roots of $h_n$ on the $y$ axis that all
meet at $(1,0)$ and $(-1,0)$. Consequently, any horizontal line close
enough to the $x$-axis will meet the graph in $n$ points, and so
$f_n\in\eppoly{1}$. 

Since the smallest positive root of the Legendre polynomial becomes
arbitrarily small, we see that the $\epsilon$ for which the interval
$-\epsilon\le \alpha\le \epsilon$  implies $f_n(x,\alpha)\in\allpoly$
becomes arbitrarily small. None the less, all $f_n$ are in
$\eppoly{1}$. 

   \begin{figure}[htbp]
  \begin{center}
    \leavevmode
    \includegraphics*[width=2in]{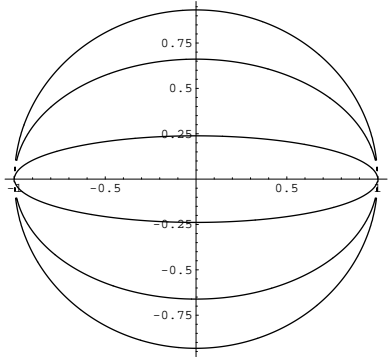}
    \caption{The graph of $\partial_y^6(x^2+y^2-1)^6=0$}
    \label{fig:leg66a}
  \end{center}
\end{figure}

\end{example}

In the case $d=1$, this last observation says that the constant term
(with respect to $y$) of a polynomial in $\eppoly{1}$ is in
$\allpoly$. If all derivatives with respect to $y$ are in $\eppoly{1}$
then all the coefficients are in $\allpoly$. 

We now show that $2$-variable analogs \eqref{eqn:legendre-2} of the 
Legendre polynomials are in $\eppoly{1}$, and consequently all their
coefficients are in $\allpoly$. Figure~\ref{fig:h6612} is the graph of
$p_{6,6,12}$. 

\begin{equation}
  \label{eqn:legendre-2}
 p_{n,m,r}(x,y)=\left(\frac{\partial}{\partial x}\right)^n 
       \left(\frac{\partial}{\partial y}\right)^m \,
       \left(x^2+y^2-1\right)^{r} 
\end{equation}

\begin{figure}[htbp]
  \begin{center}
    \leavevmode
    \includegraphics*[width=2in]{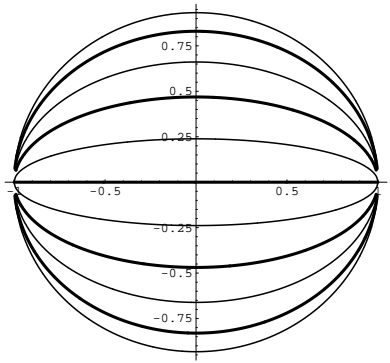}
    \caption{The graph of $\partial_y^6(x^2+y^2-1)^6=0$ and its $y$ derivative}
    \label{fig:leg66}
  \end{center}
\end{figure}

Since differentiation with respect to $x$ preserves $\eppoly{1}$, we
may assume that $n=0$. We need to consider how differentiation with
respect to $y$ affects the graph of $p_{0,m,r}$. In
Figure~\ref{fig:leg66} we have the graph of $p_{0,6,6}$ and
$p_{0,7,6}$. The light lines are in $p_{0,6,6}$ (see
Figure~\ref{fig:leg66a}) and the dark lines are in $p_{0,7,6}$. 
It is easy to see that 
$$ \text{the $x$ degree of }p_{0,m,r} =
\begin{cases}
  2r-n & \text{ $n$ even}\\
  2r-n-1 & \text{ $n$ odd}
\end{cases}
$$

We determine a point $(a,b)$ in $p_{0,7,6}$ by considering the line $x=a$,
and finding the intersection of the line with the graph. If $-1<a<1$
then there are $6$ intersection points. If $(a,b)$ is in the
  graph of $p_{0,7,6}$ then $b$ is a root of $\frac{\partial}{\partial
    y}p_{0,6,6}(a,y)$, and so there are five possible $b$'s, and they
  interlace the six intersection points. 

Thus, the graph consists of five curves joining $(-1,0)$ and $(1,0)$,
each one passing between a different pair of adjacent curves of
$p_{0,6,6}$, as can be seen in the figure.

This is the general situation: $h_{0,0,r}$ consists of $r$ circles
$x^2+y^2=1$. Each successive $y$ derivative introduces a curve between
existing curves, each joining $(-1)$ and $(1,0)$.

\begin{prop}\label{prop:pnmr}
  All $p_{n,m,r}$ in \eqref{eqn:legendre-2} lie in $\eppoly{1}$. All
  coefficients of powers of $y$ in $p_{n,m,r}$ lie in $\allpolyint{[-1,1]}$.
\end{prop}
\begin{proof}
  The only part that remains to be proved is that the roots of the
  coefficients lie in $[-1,1]$, but this is clear from the graph.
\end{proof}

We  can show that $p_{0,m,r}$ only meets the axis in $\pm1$.  We claim that
\begin{equation}\label{eqn:pkl-legendre}
\left(\frac{\partial}{\partial y}\right)^n \,
       \left(x^2+y^2-1\right)^{r}\,(x,0)
=
\begin{cases}
  n! \binom{r}{n/2}\,(x^2-1)^{r-n/2} & r \text{ even}, r\ge n \\
0 & \text{otherwise}
\end{cases}
\end{equation}
This is an easy consequence of Taylor series:
$$ 
\sum_{n=0}^\infty \left[\left(\frac{\partial}{\partial y}\right)^n \,
       \left(x^2+y^2-1\right)^{r}\right]\,(x,0) \, \frac{z^n}{n!} =
(x^2+z^2-1)^r
$$ so the left hand side of \eqref{eqn:pkl-legendre} is the
coefficient of $z^n$ in $(x^2+z^2-1)^n$ which yields the right
hand side of \eqref{eqn:pkl-legendre}.


\chapter[Generating functions][Generating functions]%
{Generating functions of linear transformations}
\label{cha:generate}

\renewcommand{\TimeStampStart}{Friday, January 18, 2008: 10:30:32}
\mytoday

The goal of this chapter is to explain the use of generating functions
to prove that linear transformations preserve particular spaces of polynomials.

\section{The idea of generating function}

Suppose $T$ is a linear transformation and $S$ is a set of
polynomials. Generating functions are a method for proving that $T$
maps $S$ to $S$. At its simplest, there is a function $G$, the
\emph{base generating function}, and a set of functions $S'$ such that
\begin{equation}\label{eqn:gf-1}
T(G)\in S' \implies T:S\longrightarrow S.
\end{equation}
We say that we have a \emph{\Polya-Schur} result if in addition
\begin{equation}\label{eqn:gf-2}
T:S\longrightarrow S \implies T(G)\in S'.
\end{equation}

The general method relies on expressing $T(f)$ in terms of $T(G)$,
where $f\in S$.  $T(G)$ is the induced map on $G$, defined by
$T(\xx^\sdiffi\,\yy^\sdiffj) = T(\xx^\sdiffi)\,\yy^\sdiffj$.  We need
the following data, where
$\xx=(x_1,\dots,x_d)$ and $\yy=(y_1,\dots,y_d)$:\\[.2cm]

\begin{tabular}{lp{3.5in}}
  \toprule
  $S\qquad$ & A set of polynomials in $\xx$.\\
  $S'$ & A set of functions in $\xx$ and $\yy$ satisfying\\
&  $f(\xx,\yy)\in S' \implies f(\xx,0)\in S$.\\
  $G$ & The base generating function. \\
  $W$ & A map from $S$ to maps  $S'\rightarrow S'$.\\
  \bottomrule
\end{tabular}\\[.2cm]

These combine in the fundamental identity for $f\in S$:

\begin{equation}
  \label{eqn:fund-iden}
   (Wf) \,T(G) \bigr|_{\yy=0} = T(f)
\end{equation}

Let's verify that the assumptions above imply \eqref{eqn:gf-1}. If
$f\in S$ and $T(G)\in S'$ then $W(f) T(G)$ lies in $S'$. Evaluation at
$\yy=0$ gives an element of $S$, so $T(f)\in S$.

\begin{table}[h]
  \centering
  \begin{tabular}{lllrc}
    \toprule
    $S$ & $S'$ & $G$ & $U$ & Comments \\
    \midrule
    $\allpoly$ & $\gsubf_2$ & $e^{-xy}$ & $-\diffd_y$\\
    $\rupint{d}$ & $\gsubf_{2d}$ & $e^{-\xx\cdot\yy}$ &
    $-\diffd_\yy$\\
    $\stabled{1}$ & $\stabledf{2}$ & $e^{xy}$ & $\diffd_y$ \\
    $\multiaff{d}$ & $\multiaff{2d}$ & $\prod_1^d (x_i+y_i)$ &
    $-\diffd_\yy$\\
    $\allpolysep$ & ? & $(1+x)^y$ & $\Delta$ \\
    $\allpoly$ & $\allpolyposf$ & $e^{x}$ & $\diffd_y$ & multiplier transformation \\
    $\stabled{1}$ & $\stabledf{1}$ & $e^{x}$ & $\diffd_y$ & multiplier transformation \\
    $\allpoly(n)$ & $\rupint{2}$ & $(x+y)^n$ & $-\diffd_y$ & bounded
    degree\\
\bottomrule    
  \end{tabular}
  \caption{Generating function data}
  \label{tab:gf}
\end{table}

How do we find $G$ and $W$? All we need is a dimension decreasing
transformation $U^\sdiffi$. We define $W(f) = f(U)$, construct dual
polynomials $p_\sdiffi(\yy)$ satisfying

\[U^\sdiffi p_\sdiffj(\yy)\bigl|_{\yy=0}=
  \begin{cases}
    0 & \diffi \ne \diffj \\
    1 & \diffi = \diffj
  \end{cases}
\] 

\noindent%
and define 
\[ G(\xx,\yy) = \sum _\sdiffi \xx^\sdiffi p_\sdiffi(\yy). \]

The fundamental identity is formally trivial, since we only need to
check on monomials:

\[
U^\sdiffi\, T(G)\bigl|_{\yy=0} =
 U^\sdiffi \sum T(\xx^\sdiffj)\,p_\sdiffj(\yy)\bigl|_{\yy=0} =
 \sum T(\xx^\sdiffj)\, U^\sdiffi p_\sdiffj(\yy)\bigl|_{\yy=0} =
 T(\xx^\sdiffi).
\]

The only difficult part is the verification of the mapping requirement

\begin{equation}\label{eqn:gf-map}
  f\times g\mapsto f(U)g \quad \text{satisfies}\quad S\times S'
  \longrightarrow S'.
\end{equation}

  In order to use the method of generating functions effectively we
  must be able to compute $T(G)$, and have some knowledge of the
  members of $S'$. This isn't always easy.

Proving a \Polya-Schur result is harder, and often requires some
analytic information in order to show that $T(G)\in S'$. For instance,
if we take $T$ to be the identity it isn't even clear why $T(G)=G$ is
in $S'$.

We now discuss, with simple examples,  each of the cases in Table~\ref{tab:gf}. 

\begin{example}
  \emph{Transformations  $\allpoly\longrightarrow\allpoly$}

  We take $S=\allpoly$, $S'=\gsubf_2$, and $U = -\diffd$, with dual
  polynomials $p_i(y) = (-y)^i/i!$. The base generating function is
\[
 e^{-xy} =\sum_0^\infty x^i \frac{(-y)^i}{i!}
\]
and is in $\gsubf_2$. We need to verify the mapping condition
\eqref{eqn:gf-map}. Let $f(x)\in\allpoly$ and
$g(x,y)\in\gsubclose_2$. We know  that
$f(-\diffd)g\in\gsubf_2$, and therefore we can apply the method of
generating functions to prove that transformations preserve all real roots.

\begin{lemma}\label{lem:gf-had-p}
  The Hadamard product $x^i\hadprod x^j\mapsto
  \begin{cases}
    0 & i\ne j \\ x^i & i=j
  \end{cases}
$\\
 determines a map $\allpolypos\times\allpoly\longrightarrow \allpoly$.
\end{lemma}
\begin{proof}
Choose $f\in\allpolypos$ and define $T_f(g) = f\hadprod g$. If $f=x^n$
then it's easy to see that $T_f(e^{-xy}) = (-xy)^n$. It follows that
in general $T(e^{-xy}) = f(-xy)$. Since $f\in\allpolypos$ we know that
$f(-xy)\in\gsubclose_2$, which is contained in $\gsubf_2$. 
\end{proof}

It is interesting that the base generating function is equivalent to
base generating functions for each positive integer.

\begin{lemma}

  If $T_\ast(1-xy)^n\in\gsubclose_2$ for $n=1,2,\dots$ then
  $T\colon\allpoly\longrightarrow\allpoly$. 
\end{lemma}
\begin{proof}
  Since $f(-\diffd_y)$ preserves $\gsubclose_2$ we know that
  $f(-\diffd_y)T_\ast(1-xy/n)^n$ is in $\gsubclose_2$. The result follows
  from the claim that
\[
\lim_{n\rightarrow\infty}
f(-\diffd_y)T\biggl(1-\frac{xy}{n}\biggr)^n\biggr|_{y=0} = T(f)
\]
By linearity it suffices to show this for $f(x) = x^k$ in which case
\begin{align*}
(-\diffd_y)^kT_\ast\biggl(1-\frac{xy}{n}\biggr)^n\biggr|_{y=0} 
&= 
\diffd^k \sum \binom{n}{i}  \biggl(\frac{y}{n}\biggr)^i T(x^i)\biggr|_{y=0} \\
&=
\binom{n}{k}n^{-k} T(x^k) \longrightarrow T(x^k).
\end{align*}

Alternatively, see Corollary~\ref{cor:sat-ind-gdf}
\end{proof}

Of course, the two approaches are equivalent:

\begin{lemma}

  If $T_\ast(e^{-xy})\in\gsubclose_2$  then
  $T_\ast(1-xy)^n)\in\gsubclose_2$ for $n=1,2,\dots$.

  Conversely, if $T(1)\ne0$ and $T_\ast(1-xy)^n\in\gsubclose_2$ for
  all $n$ then $T_\ast(e^{-xy})\in\gsubf_2$.

\end{lemma}

\begin{proof}
   Since $1-yz\in\gsubclose_2$ the operator $1+z\partial_y$ preserves
   $\gsubf_3$. We have
\[
(1+z\partial y)^n T_\ast(e^{-xy})\biggr|_{y=0} =
\sum (-1)^k z^k T(x^k)\binom{n}{k} = T_\ast(1-xz)^n
\]

  Since $(1-xy/n)^n\rightarrow e^{-xy}$ the second part follows from
  Lemma~\ref{lem:limit-tpp}. 
\end{proof}

\end{example}

\begin{example} \emph{Transformation $\gsubclose_d\longrightarrow\gsubclose_d$}
  
The argument for $\gsubclose_d$ is exactly the same as for
$\allpoly$. If $\diffi=(i_1,\dots,i_d)$ we let
\[
U^\sdiffi = \left(\frac{\partial}{\partial y_1}\right)^{i_1} \cdots
\left(\frac{\partial}{\partial y_d}\right)^{i_d}
\]
The dual polynomials are
\[
p_\sdiffi = \frac{y_1^{i_1}\cdots y_d^{i_d}
  (-1)^{i_1+\cdots+i_d}}{i_1!\cdots i_d!}
\]
For an example, here is a Hadamard product result.

\begin{lemma}
  Suppose that $f_1,\dots,f_d$ are in $\allpolypos$. The linear
  transformation \\ $g\mapsto \bigl( f_1(x_1)\cdots
  f_d(x_d)\bigr)\hadprod g$ defines a map $\gsubclose_d\longrightarrow
  \gsubclose_d$.
\end{lemma}
\begin{proof}
 The generating function of the map is $f_1(-x_1y_1)\cdots
  f_d(-x_dy_d)$ which is in $\gsubclose_{2d}$ since all $f_i\in\allpolypos$.
\end{proof}

\end{example}

\begin{example}\emph{Transformations $\stabled{1}\longrightarrow\stabled{1}$}

Since $f\times g\mapsto f(\diffd)g$ maps
$\stabled{1}\times\stabled{1}\longrightarrow\stabled{1}$ we know that
$T(e^{xy})\in\stabledf{2}$ implies $T:\stabled{1}\longrightarrow\stabled{1}$.

We start wwith a Hadamard product result.
\begin{lemma}
  The map $f\times g\mapsto f\ast g$ determines a
  map $\allpolypos\times \stabled{1}\longrightarrow\stabled{1}$.
\end{lemma}
\begin{proof}
  If we fix $f\in\allpolypos$ then the generating function of
  $g\mapsto g\ast f$ is $\exp f(xy)$. Since $\exp f\in\allpolypos$ we
  know that $\exp f(xy)\in\stabled{2}$.
\end{proof}
\end{example}

\begin{remark}
  If $T$ does not map $\stabled{1}\longrightarrow\stabled{1}$ but rather maps
  $\stabled{1}\longrightarrow\stabled{1} \cup 0$ then there might not be a
  \Polya-Schur type result. 

  Consider the Hadamard product $T\colon f\mapsto (x^2+1)\hadprod
  f$. In coordinates this is
\[
T(\sum a_i x^i) = a_0 + a_2 x^2. \]
Note that $T(x)=0$, and since all $a_i$ are non-negative,
$T:\stabled{1}\longrightarrow \stabled{1}\cup 0$. 

The generating function of $T$ is $(xy)^2/2 + 1$ which is \emph{not}
in $\stabled{2}$, so there is no \Polya-Schur type result. In addition
the induced map fails to map $\stabled{2}\longrightarrow\stabled{2}$:
\[
T(x+y)^2 = x^2 + y^2 \not\in\stabled{2}
\]

In general, the map $g\mapsto f\times g$ has generating function $\exp
f(xy)$, and this is in $\stabled{2}$ if and only if $\exp
f\in\allpolypos$. 
\end{remark}

\begin{cor}
  The following linear transformations map $\stabled{1}$ to $\stabled{1}$.
  \begin{align*}
    Exponential\quad x^k\mapsto & \frac{x^k}{k!}x^k \\
    Binomial\quad x^k\mapsto & \binom{n}{k}x^k \\
    Laguerre\quad x^k\mapsto & L_k(-x) \\
    Laguerre^{REV}\quad x^k\mapsto & L_k^{REV}(-x) \\
  \end{align*}
\end{cor}
\begin{proof}
  See Table~\ref{tab:gen-fct-1}.
\end{proof}

Next we have an inner product result. If $f = \sum f_i(x)y^i$ and $g =
\sum g_i(x)y^i$ then we define
\[
{\bigl<f,g\bigr>} = \sum i! f_i g_i
\]
This is clearly an inner product. The generating function of the map
\[
T_f: g\mapsto {\bigl<f,g\bigr>}
\]
is
\[
\sum T_f(x^iy^j) \frac{u^iv^j}{i!j!} = 
\sum \bigl<f,x^iy^j\bigr> \frac{u^iv^j}{i!j!} = 
\sum j! f_j x^i \frac{u^iv^j}{i!j!} = 
e^{xu} f(x,v)
\]
Since this is in $\stabled{2}$ if $f\in\stabled{2}$ we conclude
\begin{lemma} \label{lem:stable:ip}
  If $f\in\stabled{2}$ then the map $g\mapsto
  {\bigl<f,g\bigr>}$ satisfies
  $\stabled{2}\times\stabled{2}\longrightarrow \stabled{1}$.

\end{lemma}

Sums of squares of polynomials in $\allpolypos$ are not necessarily
stable, but the sum of the squares of the coefficients of a polynomial
in $\gsubplus_2$ is stable. Since $\gsubplus_2\subset\stabled{2}$ we can
apply the inner product.

\begin{cor}
  If $\sum f_i(x)y^i\in\gsubplus_2$ then $\sum f_i^2$ is in $\stabled{1}$.
\end{cor}

For example, using the \index{Taylor series}Taylor series $\sum
\frac{f^{(i)}}{i!}y^i$ yields that if $f\in\allpolypos$ then
\[
\sum \biggl(\frac{f^{(i)}}{i!}\biggr)^2 \in\stabled{1}
\]

\begin{cor}
  If $T\colon\allpoly\longrightarrow\allpoly$, 
  $T\colon\allpolypos\longrightarrow\allpolypos$, and $T$ preserves degree then 
$\displaystyle{
\sum \left(\frac{T(x^i)}{i!}\right)^2 \in\stabledf{1}
}.
$
\end{cor}

\index{Laguerre polynomials}
For instance, if $L_n$ is the Laguerre polynomial then the
transformation $x^n\mapsto L_n(-x)$ satisfies the hypotheses of the
corollary, so 
\[
\sum_n \frac{L_n(-x)^2}{n!n!}\in\stabledf{1}.
\]

\begin{example}\emph{Multiaffine polynomials}

  The unusual aspect of multiaffine polynomials is that the generating
  function is a polynomial. There are two choices for $U$, and they
  lead to base generating functions that are the reverse of one
  another. In each case the mapping properties follow from the
  properties for $\gsubclose_d$. 

  If we take $U^\sdiffi = (-1)^{|\sdiffi|}\diffd_{\yy_\sdiffi}$ then
  the dual polynomials are $(-1)^{|\sdiffi|}\yy_\sdiffi$ and the
  generating function is
  \[ 
\sum_\sdiffi \xx_\sdiffi (-1)^{|\sdiffi|}\yy_\sdiffi =
\prod_1^d(1-x_iy_1)
\]

We can also take $U^\sdiffi =
\diffd_{\yy_{\{1,\dots,d\}\setminus\sdiffi}}$ with dual polynomials 
$\yy_{\{1,\dots,d\}\setminus\sdiffi}$. The generating function is
\[
\sum_\sdiffi \xx_\sdiffi \yy_{\{1,\dots,d\}\setminus\sdiffi} =
\prod_1^d(x_i+y_i)
\] 
\end{example}

\begin{example}\emph{Transformations $\allpolysep\longrightarrow\allpolysep$}
  
We do not yet understand this case. We have $S=\allpolysep$,
$U=\Delta_y$, with dual polynomials $\falling{y}{i}/i!$. The base
generating function is
\[
\sum_{i=0}^\infty x^i \frac{\falling{y}{i}}{i!} = (1+x)^y
\]
Note that the base generating function is not entire, but is analytic
for $|x|<1$ and all $y$. It is not clear what $S'$ is, and we do not
know how to prove the mapping property.
\end{example}

\begin{example}\emph{Polynomials with separated roots}
  
We start with the fact that  $\Delta$ maps $\allpolysep$ to
itself. The dual polynomials are $\falling{y}{n}/n!$. Thus, the base
generating function is
\[
\sum x^i \frac{\falling{y}{i}}{i!} = (1+x)^y
\]
However, in this case we do not know the correct space of functions
$S'$ that contain the base generating function.
\end{example}

\begin{example}\emph{Multiplier transformations for $\allpoly$}
  
Let $T:x^i\mapsto a_ix^i$ be a multiplier transformation. If we treat
this a map on polynomials then we know that if $T(e^{-xy})\in\gsubf_2$
then $T:\allpoly \longrightarrow \allpoly$. However, since $T$ is a
multiplier transformation
$T(e^{-xy}) = T(e^x)(-xy)$. It follows that if $T(e^x)\in\allpolyposf$
then $T(e^{-xy})\in\gsubf_2$, and hence
$T:\allpoly\longrightarrow\allpoly$. 
\end{example}

\begin{example}\emph{Multiplier transformations for $\stabled{1}$}
  
Let $T:x^i\mapsto a_ix^i$ be a multiplier transformation. As above we
know that if $T(e^{xy})\in\stabledf{2}$ then
$T:\stabled{1}\longrightarrow\stabled{1}$.
Now $T(e^{xy}) = T(e^x)(xy)$, and the only way for $f(xy)$ to belong
to $\stabled{2}$ is if $f\in\allpolyposf$. 

\begin{lemma}
  If the multiplier transformation $T$ satisfies
  $T(e^x)\in\allpolyposf$ then
  \begin{enumerate}
  \item $T:\allpoly\longrightarrow\allpoly$
  \item $T:\stabled{1}\longrightarrow\stabled{1}$
  \end{enumerate}
\end{lemma}
\end{example}

\begin{example}\emph{Polynomials of bounded degree}

In this case we have that $T(G)\in S'$ implies that $V:S\longrightarrow
S$ where $V$ is a transformation depending on $T$.

We take $U=-\diffd_y$ and $G=(x+y)^n$. 
\[
U^k T(G)\bigl|_{y=0} = \sum (-\diffd_y)^k T(x^{n-k}) y^k
\binom{n}{k}\bigl|_{y=0} = (-1)^k T(x^{n-k})\falling{n}{k}
\]
If we define $V(x^k) = T(x^k)\falling{n}{k}$ then by taking the revese
and negating $x$ we see that
$V\colon\allpoly(n)\longrightarrow\allpoly(n)$. 

In this case even the identity yields an interesting result:
\begin{lemma}
  The mapping $x^k\mapsto x^k \falling{n}{k}$ maps $\allpoly(n)\longrightarrow\allpoly(n)$.
\end{lemma}

\end{example}

\section{Multiplier  transformations}
\label{sec:anal-gen-fct}

\index{generating function} 

The main results of this section could be derived as special cases of
the previous section, but the proofs are easy.  We are concerned with
the generating functions of multiplier transformations $x^i\mapsto
a_ix^i$. If $T(x^i)=a_ix^i$ then the generating function of $T$ is

\begin{equation}
  \label{eqn:anal-gen-fct}
T_\ast(e^{-xy})   = \sum_{i=0}^\infty T(x^i)\, \frac{(-y)^i}{i!} = 
\sum_{i=0}^\infty a_i \ \frac{(-xy)^i}{i!} = F(-xy)
\end{equation}
where 
\[
F(x) = \sum a_i \frac{x^i}{i!}
\]

We will consider the generating
function to be a function  $T(e^{-x}) = F(x)$ of one variable. 
Our goal is to prove that $(A)\implies(C)$ for multiplier
transformations. 

\begin{theorem}[\Polya-Schur] \label{thm:polya-schur}\ 

\begin{enumerate}
\item   $\allpolyposf$ is precisely the set of generating functions of
  linear transformations $x^i\mapsto a_ix^i$ that map $\allpolypos$ to
  $\allpolypos$.
\item   $\allpolyf$ is precisely the set of generating functions of
  linear transformations $x^i\mapsto a_ix^i$ that map $\allpolypos$ to
  $\allpoly$.
\end{enumerate}

\end{theorem}

\begin{proof}
  If $T\colon{}\allpolypos\longrightarrow\allpolyneg$ then 
  $$T(e^x) = \lim_{n\rightarrow\infty}
  T\,\left(1+\frac{x}{n}\right)^n.$$
  Sine $(1+x/n)^n\in\allpolypos$ we
  conclude from Lemma~\ref{lem:abs-bound-2} that the polynomials in the limit are in $\allpolypos$ by
  assumption on $T$.    Thus, the generating function is in
  $\allpolyposf$ by Lemma~\ref{lem:limit-tpp}. A similar observation shows that the generating
  function in the second case is in $\allpolyf$.
  
  Conversely, assume that the generating function $f$ is in
  $\allpolyposf$ and choose $g(x)$ in $\allpolypos$.  By
  Theorem~\ref{thm:hadamard-2} the Hadamard-type product $\ast^\prime$ satisfies
\begin{equation}\label{eqn:hadamard-5}
  T(g) = F \ast^\prime g \in\allpoly
\end{equation}
and therefore $T(g)\in\allpolypos$ since
  all coefficients are positive.  In the second case $f\in\allpolyf$
  and $g\in\allpolypos$, so we only get that $T(g)\in\allpoly$.
\end{proof}

\begin{remark}
  This result is known as the \Polya-Schur Theorem \cite{polya-schur}.
  \index{multiplier~transformations} Part (1) is sometimes stated in
  terms of maps $\allpoly\longrightarrow\allpoly$, but this can be
  deduced from the following lemma.
\end{remark}
 
\begin{lemma}
  If $T\colon{}\allpolypos\longrightarrow\allpolypos$ is a multiplier
  transformation, then $T$ extends to a linear transformation
  $\allpoly\longrightarrow\allpoly$.
\end{lemma}
\begin{proof}
  Since the generating function $F$ of $T$ is in $\allpolyposf$, the
  result follows from \eqref{eqn:hadamard-5} and Theorem~\ref{thm:hadamard-2}.
\end{proof}

Here's another consequence that can be found in \cite{polya-schur}.

\begin{lemma}\label{lem:f-pos}
  If $f\in\allpolyf$ has all positive coefficients then
  $f\in\allpolyposf$. 
\end{lemma}
\begin{proof}
 If $f\in\allpolyf$ then consider the linear transformation given 
by \eqref{eqn:hadamard-5}. If $g\in\allpolypos$ then $Tg\in\allpolypos$
since $f$ has all positive coefficients. Thus, $f\in\allpolyposf$. 

\end{proof}

{\reversemarginpar\marginpar[(B)$\implies$(C)]{}}
\begin{lemma} \label{lem:jensen}
\index{multiplier~transformations}
\index{Taa@$T(1+x)^n$}
  If $T(x^i) = a_ix^i$ and $T(\,(1+x)^n\,)\in\allpoly$ for 
  infinitely many $n$ then $T(e^x)\in\allpolyf$. 
\end{lemma}
\begin{proof}
  If we substitute $x/n$ for $x$ then we see that
$$ \lim_{n\rightarrow\infty}  T\left(1+\frac{x}{n}\right)^n\ =
T(e^x) \,\in\allpolyf$$
By Theorem~\ref{thm:polya-schur} we are done. 
\end{proof}

If $T$ is as in the last lemma, then since
$(x+1)^n\lesslesseq(x+1)^{n-1}$, we know that $T(x+1)^n\lesslesseq
T(x+1)^{n-1}$. In \cite{csordas75} it is shown that the interlacing is
strict.


\section{Linear transformations $\allpoly\longrightarrow\allpoly$} 
\label{sec:top-root-pres}

Table~\ref{tab:gen-fct-1} lists some examples of generating functions
of linear transformations that map $\allpoly$ to $\allpoly$. The first
four are elementary, and the rest are standard formulas (e.g. \cite{grads}). In the
Hadamard product and $f(x\diffd)$ we must choose
$f\in\allpolypm$.\index{Hadamard product} In the table $I_0(z)$ is
the modified Bessel function of the first kind \index{Bessel function}
and $J_0(z)$ is the Bessel function.  
In the \emph{Multiplier} entry
we are given a linear transformation of the form $T(x^i)=a_ix^i$ which
maps $\allpoly$ to itself. The generating function of $T$ in the sense
of the previous section is $f(x)=\sum a^i\frac{x^i}{i!}$ and is in
$\allpolyf$.

\index{generating function}

\begin{table}[htbp]
$$
\begin{array}{lrrclll}
\toprule
\text{Name} &&&&& \multicolumn{1}{r}{T_\ast(e^{-xy})\qquad\qquad} \\
\midrule
\displaystyle
\genfct{Affine}{(ax+b)^i}{e^{-(ax+b)y}}{}
\genfct{Derivative}{ix^{i-1}}{-ye^{-xy}}{}
\genfct{Derivative}{f(\diffd)x^{i}}{f(-y)e^{-xy}}{}
\genfct{Derivative}{f(x\diffd)x^{i}}{f(-xy)e^{-xy}}{}
\genfct{Exponential}{\displaystyle\frac{x^i}{i!}}{J_0(2\sqrt{-xy})}{}
\genfct{Hadamard}{f(x)\ast x^i}{(\expoper{}f)(-xy)}{}
\index{Hadamard product!generating function}
\genfct{Hadamard$^\prime$}{[x^r\ast x^r=r!x^r]\ast x^i}{ f(-xy)}{}
\index{Hadamard product!variant!generating function}
\genfct{Hermite}{H_i}{e^{-2xy-y^2}}{}
\index{Hermite polynomials!generating function}
\genfct{Laguerre}{L_i(-x)}{e^{-y}J_0(2\sqrt{xy})}{}
\index{Laguerre polynomials!generating function}
\genfct{Laguerre$\rev{}$}{\rev{L_i}(-x)}{e^{-xy}J_0(2\sqrt{y})}{}
\genfct{Laguerre}{L_n^i}{e^{-y}L_n(x+y)}{}
\genfct{Multiplier}{a_ix^i}{f(-xy)}{}
\bottomrule
\end{array}
$$
  \caption{Generating Functions for transformations $\allpoly\longrightarrow\allpoly$}
  \label{tab:gen-fct-1}
\end{table}

\begin{example}
  If $T\colon{}\allpoly\longrightarrow\allpoly$ doesn't preserve degree and
  the sign of the leading coefficient then the generating function
  might not be in $\gsubf_2$. Consider $T(g) = G(\diffd)f(x)$. We have
  seen that this linear transformation doesn't extend to a map
  $\gsubclose_2\longrightarrow\gsubclose_2$ when $f(x)=x$. An easy
  calculation shows that the generating function is $(x-y)e^{-xy}$, which
  isn't in $\gsubf_2$. 
\end{example}

\begin{remark} \label{rem:rolle-3}
  \index{Rolle's Theorem}
  The generating function of the linear transformation $f\mapsto
  f+\alpha f^\prime$ is $(1+\alpha y)e^{-xy}$. Since this is in $\allpolyf_2$
  it follows that $f+\alpha f^\prime$ is in $\allpoly$ for all
  $f\in\allpoly$. This is another proof of Rolle's theorem (Theorem~\ref{thm:rolle}).
\end{remark}

\begin{example} \index{Hermite polynomials}
  What happens if $T$ doesn't map $\allpoly$ to itself? The generating
  functions should be functions that are not in $\allpolyf_2$. Let's
  look at some examples. The Hermite transformation $T(x^i) = H_i$
  maps $\allpoly$ to itself. Since this is not an affine map we know
  that the inverse does not map $\allpoly$ to itself. Using identities
  in \cite{roman} we can find the generating function  of
  $T^{-1}$:
\begin{align*}
  T^{-1}(x^n) &= \left(\frac{i}{2}\right)^n
  H_n\left(\frac{-ix}{2}\right)\\
 T_\ast(e^{-xy}) &= \sum_{n=0}^\infty \left(\frac{i}{2}\right)^n
  H_n\left(\frac{-ix}{2}\right) \frac{(-y)^n}{n!} \\
&= \sum_{n=0}^\infty 
  H_n\left(\frac{-ix}{2}\right) \frac{\left(\frac{(-iy)}{2}\right)^n}{n!}
  \\
  &= e^\frac{-2xy+y^2}{4}
\end{align*}
This is not in $\allpolyf_2$ since substituting $x=0$ gives
$e^{y^2/4}$ which is not in $\allpolyf$.


Next, consider \index{rising factorial} $T(x^i)=\rising{x}{i}$ which
maps $\allpolypos$ to itself.  We can compute

$$T_\ast(e^{-xy}) = \sum_{i=0}^\infty T(x^i)\frac{(-y)^i}{i!} =
\sum_{i=0}^\infty x(x+1)\cdots(x+i-1)\frac{(-y)^i}{i!} = (1+y)^{-x}
$$
Although the generating function  has a simple form
it isn't in $\allpolyf_2$. This is because it is not an entire
function: $T_\ast(e^{-xy})(1,y)=(1+y)^{-1}$ has a pole at $y=-1$.

An example of a different sort is given by the generating function of
the $q$-derivative. Recall if $\affa(x) = qx$ then the $q$-derivative
is $\dfrac{f(qx) - f(x)}{x(q-1)}$. The generating function is 
$$ \frac{e^{-qxy} - e^{-xy}}{x(q-1)}$$
which is just the $q$-derivative of $e^{-xy}$. It's not in $\allpolyf$
since substituting $y=1$ gives a function that has complex roots.

\end{example}

\begin{example}
  If $T(x^i)=H_i$, then using the generating function we can show that
  $T(\allpoly)\not\subset \expoper{}(\allpoly)$ without exhibiting any
  particular polynomial in $\expoper{}(\allpoly)$ that is not in
  $T(\allpoly)$. If $T(\allpoly)\subset \expoper{}(\allpoly)$ then
  $\expoper{}^{-1}T(\allpoly)\subset \allpoly$. We will show that the
  generating function of $\expoper{}^{-1}T$ isn't in $\allpolyf_2$.
  We compute
  \begin{align*}
    \sum_i \expoper{}^{-1}\,T(x^i) \frac{(-y)^i}{i!} & =
    \expoper{}_\ast^{-1}\, \left(\sum_i T(x^i)\frac{(-y)^i}{i!}\right) \\
& = \expoper{}_\ast^{-1}\, e^{-2xy-y^2} \\
&= e^{-y^2} \sum i! \frac{(-2xy)^i}{i!} \\
&= \frac{e^{-y^2}}{1+2xy}
  \end{align*}
This last function is not entire.
\end{example}

\section{Linear transformations $\allpoly\longrightarrow\allpoly$ -
  general properties}
\label{sec:gen-fct-gen-prop}

In this section we establish general properties of generating
functions. In the next section we will apply some of the results to
show that various linear transformations preserve $\allpoly$. 

If a linear transformation $T$ has generating function $F(x,y)$, then
it is easy to compute the generating function of the \index{\Mobius\ 
  transformation}\Mobius\ transformation $T_M$
(\chapsec{linear}{mobius}), where $M:z\mapsto \frac{az+b}{cz+d}$.  To
begin with, there is a very useful relationship between the generating
functions of a linear transformation $T$ and the transformation
$T_{1/z}$. Recall
\index{T@$T_{1/z}$}%
(\chapsec{linear}{mobius}) that if $T(x^n)=p_n(x)$ then $T_{1/z}(x^n)
= \rev{p_n}$. The generating function of $T_{1/z}$ is
$$ \sum_{n=0}^\infty T_{1/z}(x^n) \frac{(-y)^n}{n!} =
\sum_{n=0}^\infty x^np_n(1/x) \frac{(-y)^n}{n!} = 
\sum_{n=0}^\infty p_n(1/x) \frac{(-xy)^n}{n!} = F(1/x,xy)$$

More generally, if $M$ is a \Mobius\ transformation
$M(z)=\frac{az+b}{cz+d}$ then we define $T_M(x^n) = (cx+d)^n
T(x^n)\left(\frac{ax+b}{cx+d}\right)$, and the generating function is
$$ \sum_i (cx+d)^i T\left(\frac{ax+b}{cx+d}\right)\,\frac{(-y)^i}{i!} =
F\left(\frac{ax+b}{cx+d},(cx+d)y\right)$$


\index{T@$T_{1/z}$}

If we first compose with an affine transformation then the generating
function has an exponential factor.

\begin{lemma}\label{lem:affine-gen-fct}
  Suppose that $T$ is a linear transformation with generating function
  $F(x,y)$. If $S(f) = T(f(ax+b))$ then the generating function of $S$
  is $e^{-by}F(x,ay)$. 
\end{lemma}
\begin{proof}
  A computation:
  \begin{align*}
    \sum_n S(x^n) \frac{(-y)^n}{n!} &= 
    \sum_n T((ax+b)^n) \frac{(-y)^n}{n!} \\
&=     \sum_n \sum_k a^k T(x^k) b^{n-k}\binom{n}{k} \frac{(-y)^n}{n!} \\
&=     \sum_k  a^k T(x^k)  \frac{(-y)^k}{k!} \sum_{n\ge k}
\frac{(by)^{n-k}}{(n-k)!}\\ 
&=  \sum_k   T(x^k) \frac{(-ay)^k}{k!} e^{-by}\\
& = e^{-by}F(x,ay)
  \end{align*}
\end{proof}

The next  result shows that the coefficients of the polynomials defining a 
root preserving transformation are highly constrained. 

\begin{lemma}
  Suppose that $T\colon{}\allpoly\longrightarrow\allpoly$.
  Choose a non-negative integer $r$ and let $d_i$ be the coefficient
  of $x^r$ in $T(x^i)$.  Then, the series $ \sum_{i=0}^\infty d_i
  \frac{x^i}{i!}$ is in $\allpolyf$, and $x^i\mapsto d_ix^i$ maps
$\allpolypos\longrightarrow\allpoly$.

 In addition, for any $n$ the polynomial $\sum
  \binom{n}{i}d_ix^i$ is in $\allpoly$.
\end{lemma}
\begin{proof}
  The first sum is the coefficient of $x^r$ in  the generating
  function of $T$.
  The second follows from Lemma~\ref{lem:sumsofF}.
\end{proof}

In particular, the  generating function of the constant
terms of a root preserving linear transformation is in $\allpolyf$.

The leading coefficients are similarly constrained.

\index{coefficients!leading}

\begin{lemma} \label{lem:pd-lead-coef}
  Suppose $T\colon{}\allpoly\longrightarrow\allpoly$,
  $T\colon{}\allpolypos\longrightarrow\allpolyneg$, $T$ preserves degree, and
  the leading coefficient of $T(x^i)$ is $c_i$.Then $ \sum c_i
  \frac{x^i}{i!}\in\allpolyf$, and $x^i\mapsto c_ix^i$ maps
  $\allpolypos$ to itself.
\end{lemma}
\begin{proof}
  Corollary~\ref{cor:t1/z} shows that 
  $T_{1/z}:\allpoly\longrightarrow\allpoly$. The
\index{T@$T_{1/z}$}%
  constant term of $T_{1/z}(x^i)$ is $c_i$, so we can apply the
  previous lemma.
\end{proof}

We have seen (Corollary~\ref{cor:leading-coef}) that if $T\colon{}x^n\mapsto f_n(x)$
 maps $\allpoly\longrightarrow\allpoly$, then the
map $x^n\mapsto c_nx^n$, where $c_n$ is the leading coefficient of
$f_n$, also maps $\allpoly$ to itself. Evaluation has similar
properties.

  \begin{cor}\label{cor:eval-f2}
    Suppose that $T\colon{}x^n\mapsto f_n(x)$ preserves degree and maps
  $\allpoly\longrightarrow\allpoly$, then for any
  $\alpha\in\reals$ the multiplier transformation \\ $x^n\mapsto
  f_n(\alpha)x^n$ maps $\allpolypos\longrightarrow\allpoly$.
  \end{cor}
  \begin{proof}
    If the generating function $T_\ast(e^{-xy}\in\allpolyf_2$, then evaluation
    is in $\allpolyf$, and so
$$\sum f_n(\alpha) \frac{(-y)^n}{n!}\in\allpolyf$$
from which the conclusion follows.
  \end{proof}

If we have two linear transformations $S,T$ with generating functions
in $\gsubf_2$ then we can multiply their generating functions by
$\sum a_iy^{n-i}z^{n-i}$:
$$ \left(\sum_{i=0}^\infty S(x^i)\frac{(-y)^i}{i!}\right)
\left(\sum_{i=0}^\infty T(x^i)\frac{(-z)^i}{i!}\right)
\left(\sum_{i=0}^n  a_iy^{n-i}z^{n-i}\right)
$$
If $\sum_0^n a_ix^i\in\allpolypos$ then the third factor is in
$\gsubposf_2$ so all the factors are in $\gsubf_3$. The
coefficient of $y^nz^n$ is 
$\sum \frac{a_i\,S(x^i)T(x^i)}{i!i!}$. Consequently, the linear
transformation $x^n \mapsto \sum a_i\,\frac{S(x^i)T(x^i)}{i!i!}$
maps $\allpolypos$ to $\allpoly$. We can remove the factorials.

\begin{lemma}\label{lem:gen-fct-ST}
  Suppose $S,T$ are linear transformations with generating functions
in $\gsubf_2$. The map $x^i \mapsto (-1)^iS(x^i)T(x^i)$ maps
$\allpolypos\longrightarrow \allpoly$.
\end{lemma}
\begin{proof}
  If $\sum a_ix^i\in\allpolyalt$ then we can apply Lemma~\ref{lem:fpp-d}.
  Simply observe (as above) that since
\begin{multline*}
\left(\sum_{i=0}^n  a_i(-1)^i\partial_y^i\partial_z^i\right)
 \left(\sum_{i=0}^\infty S(x^i)\frac{(-y)^i}{i!}\right)
\left(\sum_{i=0}^\infty T(x^i)\frac{(-z)^i}{i!}\right)\bigg\vert_{y=z=0}
\\ = \sum_{i=0}^n a_i(-1)^i S(x^i)T(x^i)
\end{multline*}
and the differential operator maps $\gsubf_3$ to itself so the right
hand side is in $\allpoly$.
\end{proof}

We can use the relationship between linear transformations
$\allpoly\longrightarrow\allpoly$ and elements of $\allpolyf_2$ to
construct new linear transformations from old. We start with a linear
transformation, determine the corresponding element of $\allpolyf_2$,
manipulate it to create a new element of $\allpolyf_2$, and then convert
back to a linear transformation. Unfortunately, we don't know many ways
 to get new elements of $\allpolyf_2$ from old ones. Differentiation
 in $\allpolyf_2$ yields some simple transformations. Suppose that

 \begin{align*}
 f(x,y) & = \sum T(x_i)\frac{(-y)^i}{i!}\\
   \frac{\partial f}{\partial x} & = \sum \frac{d}{dx} T(x^i)\frac{(-y)^i}{i!}\\
   & \text{ and so $\frac{\partial }{\partial x}$ corresponds to
       $g\mapsto \frac{d}{dx}(Tg)$}\\
   \frac{\partial f}{\partial y} & =- \sum  T(x^i)\frac{(-y)^{i-1}}{(i-1)}\\
   & \text{ and  $\frac{\partial }{\partial y}$ corresponds to
       $g\mapsto -T(xg)$}
 \end{align*}

Multiplication in $\allpolyf_2$ gives a  convolution of linear
transformations.

\begin{lemma}
  Let $S,T\colon{}\allpoly\longrightarrow\allpoly$ preserve degree and the
  sign of the leading coefficient. Then the linear
  transformation below also maps $\allpoly\longrightarrow\allpoly$.
$$x^n \mapsto \sum_{i=0}^n\binom{n}{i}T(x^i)S(x^{n-i})$$
\end{lemma}

\begin{proof}
  Multiplying two elements in $\gsubf_2$ yields an element of
  $\gsubf_2$:
$$ \left(\sum T(x^i)\frac{(-y)^i}{i!}\right)
\left(\sum S(x^j)\frac{(-y)^j}{j!}\right) =
\sum_{n=0}^\infty\left(\sum_{i=0}^n\binom{n}{i}T(x^i)S(x^{n-i})\right)\frac{(-y)^n}{n!}
$$
which establishes the lemma.
\end{proof}

\begin{remark}
  \index{convolution} This really is a property of linear
  transformations that satisfy induction. Indeed, if $S,T$ satisfy
  induction, then since $W(f) = (S_\ast T_\ast f(x+y))(x,x)$ we see
  that $W(f)\in\allpoly$. The example $S(g) = g(\diffd)x$, $T(f)=f$
  yields $W(1)=1+x^2$. This is not a counterexample, since $S$ does
  not satisfy induction. 
\end{remark}


\begin{remark}
  We revisit the characterization results of
  \chapsec{linear}{characterize}. If $T$ is a degree preserving linear
  transformation such that $Tf$ and $f$ interlace for all
  $f\in\allpoly$, then the linear transformation $T_\alpha(f) = T(f) +
  \alpha f$ maps $\allpoly$ to itself. Thus, for all $\alpha$, the
  generating function of $T_\alpha$ is in $\gsubf_2$. If $G$ is the
  generating function of $T$, then the generating function of
  $T_\alpha$ is $G + \alpha e^{-xy}$. Consequently, $G$ and $e^{-xy}$
  interlace.

  We know that such linear transformations are given by $Tf = axf + bf'
  + cf$ where $a$ and $b$ have the same sign. The generating function
  of this transformation is $(ax+by+c)e^{-xy}$. This makes it clear
  why we must have the sign condition. See
  Question~\ref{ques:interlace-exy}. 

  Since we do not have a characterization of functions in $\gsubf_2$,
  this point of view does not give a \emph{proof} of the
  characterization theorems, just an understanding of them.
\end{remark}

\begin{lemma}
  If $T$ is an onto linear transformation defined on polynomials in one
  variable and the linear transformation $S\colon{}T(x^i)y^j \mapsto
  T(x^j)y^i$ maps $\gsubpos_2$ to itself, then $T$ is an affine
  transformation. 
\end{lemma}
\begin{proof}
Since $S(T(x^i)y^j) = T(x^j) y^i$ by linearity
we see $S(T(f(x))y^j) = T(x^j)f(y)$ for any polynomial $f$. Choosing
$f=T^{-1}(x^i)$ shows that $S(x^iy^j) = T(x^j)T^{-1}(y^i)$. 
The compositions

\centerline{ \xymatrix{ \allpoly \ar@{->}[r] & \gsubclose_2
    \ar@{->}[r]^S & \gsubclose_2
    \ar@{->}[r] & \allpoly \\
    f(x) \ar@{|->}[r]^{x\mapsto x} & f(x) \ar@{|->}[r] & T(f)
    \ar@{|->}[r]^{y=1} & T(f) \\
    f(x) \ar@{|->}[r]^{x\mapsto y} & f(y) \ar@{|->}[r] & T^{-1}(f)
    \ar@{|->}[r]^{x=1} & T^{-1}(f) }} 
\noindent%
show that both $T$ and $T^{-1}$
map $\allpoly$ to itself, and consequently $T$ is affine.
\end{proof}

\section{Applications: Linear transformations $\allpoly\longrightarrow\allpoly$ }
\label{sec:top-gen-fct}

We now have a powerful tool for showing that a linear transformation
maps $\allpoly$ to itself. We use the generating function to show that
many linear transformations preserve roots. 

\begin{example} \label{ex:fi}
If we start with an element of $\allpolyf_2$ then we get a linear
transformation. For instance, if $f\in\allpoly$, then
$f(x+y)\in\gsubpos_2\subset\allpolyf_2$. The Taylor series 
\index{Taylor series} of $f$ is
\[
 f(x+y) = \sum f^{(i)}(x)\frac{y^i}{i!} = \sum (-1)^if^{(i)}(x)\frac{(-y)^i}{i!} 
\]
and consequently the linear transformation $x^i\mapsto (-1)^if^{(i)}$ is a
  map $\allpoly\longrightarrow\allpoly$. Precomposing with $x\mapsto
  -x$ shows that $x^i\mapsto f^{(i)}$ maps $\allpoly\longrightarrow\allpoly$.
\end{example}

In Corollary~\ref{cor:laguerre-x} we showed that the Laguerre transformation maps
$\allpolyalt$ to itself. Now we show that it actually maps $\allpoly$
to itself.

\begin{lemma} \label{lem:laguerre-2}
\index{Laguerre polynomials} \index{polynomials!Laguerre}
  The mapping $x^n \mapsto L_n(-x)$ (the $n$-th Laguerre polynomial)
  maps $\allpoly$ to itself.
\end{lemma}
\begin{proof}
  It suffices to know that the generating function of the
  transformation is in $\allpolyf_2$. From Table~\ref{tab:gen-fct-1}
  the generating function is $J_0(2\sqrt{xy})$ where the Bessel
  \index{Bessel function} function $J_0(z)$ is given by the series
\begin{align*}
  J_0(z) &= \sum_{k=0}^{\infty} \frac{(-1)^k}{2^{2k}k!k!} z^{2k}\\
  J_0(2\sqrt{xy}) &= \sum_{k=0}^{\infty} \frac{(-1)^k}{k!k!} x^ky^k\\
  \intertext{The \index{Bessel function} Bessel function has a product
    formula}
  J_0(z) &= \prod_{k=1}^\infty \left( 1 - \frac{z^2}{r^2_i}\right)\\
  \intertext{and hence}
  J_0(2\sqrt{xy}) &= \prod_{k=1}^\infty \left( 1 - \frac{2xy}{r^2_i}\right)\\
\end{align*}
Since the factors of this last product are all in $\gsubposclose_2$,
it follows that $J_0(2\sqrt{xy})$ is  in $\gsubposf_2$, and
hence the generating function given of $x^n\mapsto L_n(-x)$  is in
$\gsubposf_2$.
\end{proof}

Since the linear transformation $T\colon{}x^n\mapsto L_n$ maps $\allpoly$ to
$\allpoly$, we know that $T^{-1}$ does not map $\allpoly$ to
itself. It's easy to verify that $$ T^{-1} = \expoper{}^{-1}\circ T
\circ \expoper{}^{-1}.$$ Since $\expoper{}^{-1}$ doesn't map
$\allpoly$ to itself, it isn't surprising that $T^{-1}$ doesn't either.

\begin{lemma} \label{lem:laguerre-3}
  If $n$ is a positive integer, then the linear transformation
  $x^i\mapsto L_n^i(-x)$ maps $\allpoly$ to itself.
\end{lemma}

\begin{proof}
  From Table \ref{tab:gen-fct-1} the generating function $F(x,y)$ of
  the transformation is \mbox{$e^yL_n(x+y)$}. Since
  $F(x,y)\in\gsubf_2$, the linear transformation maps $\allpoly$ to
  itself.
\end{proof}

We can apply \Mobius\ transformations to show that the reverse of the
Laguerre transformation also maps $\allpoly$ to itself.

\begin{lemma} \label{lem:laguerre-4}
  The linear transformation $S(x^n) = \rev{L_n}(-x)$ maps $\allpoly$ to
  $\allpoly$. 
\end{lemma}
\begin{proof}
  The generating function of $T(x^n)=L_n(-x)$ is $e^yJ_0(2\sqrt{xy})$.
  Since $S = T_{1/z}$ we know the generating function $G(x,y)$ of $S$
  is $ F(1/x,xy)
  = e^{-xy}J_0(2\sqrt{y})$. Now $e^{-xy}\in\gsubf_2$, and
  $J_0(2\sqrt{y})\in\allpolyf$, and hence 
  $G(x,y)\in\gsubf_2$.
\end{proof}

\index{T@$T_{1/z}$}%

The next result is an  immediate consequence of Lemma~\ref{lem:gen-fct-ST}.

\index{Laguerre polynomials} \index{polynomials!Laguerre}
\index{Hermite polynomials} \index{polynomials!Hermite}
\begin{cor} \label{cor:hnln}
  If $H_n$ is the Hermite polynomial, and $L_n$ is the Laguerre
  polynomial, then the following linear transformations map
  $\allpolypos$ to $\allpoly$.
  \begin{itemize}
  \item $x^n\mapsto H_n(x)^2$
  \item $x^n\mapsto L_n(-x)^2$
  \item $x^n\mapsto H_n(x)L_n(-x)$
  \end{itemize}
\end{cor}

Lemma~\ref{lem:tprod} characterizes the linear transformations of the form \\
$T(x_n)=\prod_{i=1}^n (x+a_i)$ that map $\allpoly$ to itself. We
conclude
\begin{lemma}
  The only choice of constants $a_i$ satisfying
$$ \sum_{n=0}^\infty\left( \prod_{i=0}^n (x+a_i)\ \frac{y^n}{n!}\right) \in\allpolyf_2$$ 
is $a_1=a_2=a_3\cdots$. 
\end{lemma}

We can use generating functions to get results converting complex
roots by applying appropriate affine transformations. See \cite{bump}.
  \begin{lemma}\label{lem:riemann-1}
    Consider the  linear transformations 

\begin{align*}
T\colon{}f(x)\mapsto &     x\,f(x+1)+(x-1)\,f(x-1)\\
S\colon{}f(x)\mapsto &  (\imag x+ \frac{1}{2})f(x-\imag) + (\imag
    x-\frac{1}{2})f(x+\imag)
  \end{align*}
The following diagram commutes

  \centerline{\xymatrix{
\allpoly
      \ar@{<-}[d]_{{x\mapsto \imag(\frac{1}{2}-x)    }}           
      \ar@{->}[rrr]^{{ S  }}         
      &&&
\allpoly
      \ar@{->}[d]^{{ x\mapsto \imag x + \frac{1}{2}   }} \\        
\allpolyint{\frac{1}{2}+\imag\reals}
      \ar@{->}[rrr]^{{ T   }}         
      &&&
\allpolyint{\frac{1}{2}+\imag\reals}
}}

  \end{lemma}
  \begin{proof}
    It is easy to verify that the diagram commutes at the element
    level, so it suffices to show that $S$ maps $\allpoly$ to itself.
    We compute the generating function of $S$:

    \begin{gather*}
      \sum_{n=0}^\infty \left( (\imag x+\frac{1}{2}) (x-\imag)^n + 
                               (\imag x-\frac{1}{2})
                               (x+\imag)^n\right)
                               \frac{(-y)^n}{n!} \\
=\frac{1}{2}\left(
-e^{\imag \,\left( -1 + \imag \,x \right) \,y} + 
  e^{\imag \,\left( 1 + \imag \,x \right) \,y} + 
  \left( 2\,\imag  \right) \,e^{\imag \,\left( -1 + \imag \,x \right) \,y}\,
   x + \left( 2\,\imag  \right) \,
   e^{\imag \,\left( 1 + \imag \,x \right) \,y}\,x\right) \\
 = \imag \,\left( 2\,x\,\cos (y) + \sin (y) \right) {e^{-x\,y}}
\end{gather*}
Now Example~\ref{ex:sincos} shows that
$2x\cos(y)+\sin(y)\in\allpolyf_2$, so all factors are in
$\allpolyf_2$, and hence $S$ maps $\allpoly$ to itself.

  \end{proof}

We can  determine more transformations that are sums of shifted arguments
using Lemma~\ref{lem:f-of-sin}.
\index{sin and cos}
\begin{lemma}
  If $f(x)=\sum a_i x^i$ is in $\allpolyint{(-1,1)}$ and we define
$$ T(g) = \sum_k a_k \,(g(x +k\imag)+ g(x-k\imag)\,)$$
then $T$ maps $\allpoly$ to itself.
\end{lemma}
\begin{proof}
  It suffices to compute the generating function of $T$. Notice that
  the generating function of the transformation $x^n\mapsto
  (x+\alpha)^n$ is $e^{-xy-\alpha y}$. Consequently, the generating
  function of $T$ is
\[ 2 e^{-xy} \, \sum_k a_k\left( \frac{e^{-k y \imag} + e^{k y
      \imag}}{2}\right)  =
2 e^{-xy} \, \sum_k a_k \cos(ky)  
\]
Since $f(x)\in\allpolyint{\Delta}$ we know that $\sum a_k
\cos(ky)\in\allpolyf$ - see the proof of Lemma~\ref{lem:cheby}.
Thus the generating function is in $\gsubf_2$.
\end{proof}

For example, we can take $f=(x+1)^n$. The lemma shows that if
$g\in\allpoly$ then 
$$ \sum_{k=0}^n \binom{n}{k}\left(g(x+k\imag) +
  g(x-k\imag)\right)\in\allpoly.$$ 

Now we have  properties of the Legendre and Jacobi polynomials.

\begin{lemma}\label{lem:legendre2}
  The linear transformation $x^k\mapsto P_k(x)/k!$ defines a map
  $\allpoly\longrightarrow\allpoly$. 
\end{lemma}
\begin{proof}
  Since $J_0(x)\in\allpolyposf$, it follows from
  \eqref{eqn:legendre-1} that $f_0(xy+y)f_0(xy-y)\in\rupint{2}$.
  Consequently the generating function for the transformation is in
  $\gsubf_2$.
\end{proof}

We can generalize this to Jacobi polynomials. We use the
\emph{Mathematica} definition of $P_n^{\lambda,\mu}(x)$.  
\index{Jacobi polynomials}

\begin{lemma}
  For $\lambda,\mu>-1$ the linear transformation $$x^n\mapsto
  \dfrac{n!}{\rising{\lambda+1}{n}\rising{\mu+1}{n}}\,P_n^{\lambda,\mu}(x)$$
    defines a map $\allpoly\longrightarrow\allpoly$.
\end{lemma}
\begin{proof}
  From \cite{iliev}*{(2.1.A)} we let
  \begin{align*}
    f_\lambda(z) &= \sum_{k=0}^\infty \frac{1}{k!}
    \frac{\Gamma(\lambda+1)} {\Gamma(\lambda+k+1)} \,z^k \\
  \intertext{then we have that}
f_\lambda(xz-z)f_\mu(xz+z) &= \sum_{n=0}^\infty
\frac{2^n}{\rising{\lambda+1}{n}\rising{\mu+1}{n}}\,P_n^{\lambda,\mu}(x)\,z^n
    \end{align*}
Now $f_\lambda$ is a modified Bessel function and is in
$\allpolyposf$, so we follow the argument of the previous lemma.
\end{proof}

\section{Applications: the diamond product}
\label{sec:hermite-diamond}

We can show that a bilinear mapping $T$ satisfies
$T\colon\allpolyint{\diffi}\times\allpoly\longrightarrow\allpoly$ by showing
that for every $f\in\allpolyint{\diffi}$ the linear transformation
$T_f(g)=T(f,g)$ maps $\allpoly$ to itself. 

We now show that the diamond product based on the Hermite polynomials
maps $\allpoly\times\allpoly\longrightarrow\allpoly.$  The proof shows
that the generating function of the transformation is in
$\gsubaltf_2$ by using some special function identities for Hermite
polynomials.

Recall (\chapsec{operators}{bilinear-diamond}) that if $T$ is a linear
transformation then the {diamond product} is given by $f\mydiamond{T} g
= T^{-1}(T(f)\, T(g))$.

\index{diamond product!Hermite}
\index{product!diamond}
\index{Hermite polynomials}
\index{polynomials!Hermite}

\begin{lemma} \label{lem:diamond-hermite}
  Suppose $T(H_n)=x^n$ where $H_n$ is the Hermite polynomial. The
  diamond product $\mydiamond{T}$ defines a mapping
  $\allpoly\times\allpoly\longrightarrow\allpoly$. 
\end{lemma}
\begin{proof}
  We fix $f\in\allpoly$, and show that the generating function of
  $g\mapsto f\mydiamond{T} g$ is in $\gsubf_2$.  We recall two
  identities that can be found in \cite{roman}
  \begin{align}
    \label{eqn:hermite-inv}
  T(x^n) &= \left(\frac{i}{2}\right)^n\,
  H_n\left(\frac{-ix}{2}\right)\\
H_n(x) \,H_k(x) &= \sum_{j=0}^{\min(n,k)}
2^j\binom{n}{j}\binom{k}{j}j! \,H_{n+k-2j}(x) \notag
  \end{align}

We simplify matters by initially setting $f = x^n$. The generating
function of $g\mapsto x^n\mydiamond{T} g$ is
  \begin{align*}
    F(x,y) &= \sum_{k=0}^\infty T^{-1}(T(x^n)\, T(x^k))\,\frac{(-y)^k}{k!}
    \\
&=\sum_{k=0}^\infty T^{-1}\left(
\left(\frac{i}{2}\right)^nH_n\left(\frac{-ix}{2}\right) \, 
\left(\frac{i}{2}\right)^kH_k\left(\frac{-ix}{2}\right) 
\right)
\,\frac{(-y)^k}{k!}\\
&=\sum_{k=0}^\infty T^{-1}\left(
(H_n\,H_k)\left(\frac{-ix}{2}\right) 
\right)
\,\left(\frac{i}{2}\right)^{n+k}\frac{(-y)^k}{k!}\\
&=\sum_{k=0}^\infty T^{-1}
\left(
\sum_{j=0}^{\min(n,k)}
2^j\binom{n}{j}\binom{k}{j}j! \,H_{n+k-2j}\left(\frac{-ix}{2}\right)
\right)
\,\left(\frac{i}{2}\right)^{n+k}\frac{(-y)^k}{k!}\\
\intertext{Applying the identity \eqref{eqn:hermite-inv} yields}
H_{n+k-2j}\left(\frac{-ix}{2}\right) &=
\left(\frac{i}{2}\right)^{-n-k+2j}T(x^{n+k-2j})\quad\quad\text{and so}\\
F(x,y)&=\sum_{k=0}^\infty 
\sum_{j=0}^{\min(n,k)}
2^j\binom{n}{j}\binom{k}{j}j! \, \left(\frac{i}{2}\right)^{-n-k+2j}x^{n+k-2j}
\,\left(\frac{i}{2}\right)^{n+k}\frac{(-y)^k}{k!}\\
&= e^{-xy} \left(\frac{2x+y}{2}\right)^n\\
\intertext{By linearity, we conclude that the generating function
  $G(x,y)$ for $g\mapsto f\mydiamond{T} g$ is}
& e^{-xy}f\left(\frac{2x+y}{2}\right)\\
\end{align*}
Since $G(x,y)\in\gsubf_2$, and so the diamond product maps
$\allpoly\times\allpoly\longrightarrow\allpoly$.

\end{proof}

The Hermite diamond product has a simple description using the Hermite
basis, and using this description we can find some interesting linear
transformations. The diamond product is simply $H_n\mydiamond{T}H_m =
H_{n+m}$. Upon fixing $m=1$   the linear transformation
$f\mapsto f\mydiamond{T}H_1$ defines a map
$\allpoly\longrightarrow\allpoly$. Using the definition of the diamond
product above, this shows that the linear transformation  $H_n\mapsto
H_{n+1}$ determines a map $\allpoly\longrightarrow\allpoly$.
\label{hermite-increment}
Since $H_n\mydiamond{T}4x^2 = H_n\mydiamond{T}(H_2+2H_0)$, the map
$H_i\mapsto H_{i+2}+2H_i$ maps $\allpoly\longrightarrow\allpoly$.

\section{Applications: generalized Hadamard products}
\label{sec:gen-fct-hadamard}

In this section we characterize several generalized Hadamard products.
\index{Hadamard product!general} We have seen that the two Hadamard
products $x^i\ast x^i = x^i$ and $x^i \ast^\prime x^i = i!x^i$
map $\allpoly\times\allpolypos\longrightarrow\allpoly$. The next
result characterizes such general Hadamard products.

\begin{prop}\label{prop:hadamard-prod-gen}
  Suppose that $g(x) = \displaystyle \sum_{i=0}^\infty a_i \frac{x^i}{i!i!}$. The
  generalized Hadamard product
$$ x^i \hadprod x^j = \begin{cases} a_i x^i & i=j \\ 0 & \text{otherwise}
\end{cases}
$$
defines a map
$\allpoly\hadprod\allpolypos\longrightarrow\allpoly$ if and
only if $g(x)\in\allpolyposf$. 
\end{prop}
\begin{proof}
  Choose $f=\sum b_i x^i\in\allpolypos$, and consider the map
  $T_f:g\mapsto f\hadprod g$. We will show that the
  generating function of $T_f$ is in $\allpolyf_2$. Compute:
  \begin{align*}
    \sum T_f(x^n) \frac{(-y)^n}{n!} &= \sum (x^n \hadprod f)
    \frac{(-y)^n}{n!} \\
&= \sum b_n a_n \frac{x^n(-y)^n}{n!} \\
&= (f\ast^\prime g)(-xy)
  \end{align*}
  Since $f\in\allpolypos$, and $g\in\allpolyposf$, we know that
  $f\ast^\prime g\in\allpolypos$, and consequently $(f\ast^\prime
  g)(-xy)\in\gsubposf_2$.

  Conversely, since the Hadamard product extends to $\allpolyf$ the
  following is in $\allpolyposf$: $ e^x \hadprod e^x = \displaystyle
  \sum a_i \frac{x^i}{i!i!}$
\end{proof}


We can  generalize this to products defined by
$$
x^{i_1}\times x^{i_2} \times \cdots \times x^{i_d} \mapsto 
\begin{cases}
a_i x^i & i=i_1=i_2\cdots= i_d\\
0 &\text{otherwise}
\end{cases}
$$
where all $a_i$ are positive.  A similar argument shows that a
necessary condition that this
product determines a linear transformation
$\left(\allpolypos\right)^d\longrightarrow\allpolypos$ is
$$\sum_{i=0}^\infty \frac{a_i}{(i!)^d}x^i\in\allpolyposf$$

Next we consider maps $\gsubplus_2\longrightarrow\allpolypos$ that
have the form
\begin{equation}
  \label{eqn:had-prod-diag}
  x^iy^j \mapsto 
  \begin{cases}
    a_i x^i & i=j \\ 0 & \text{otherwise}
  \end{cases}
\end{equation}
If $a_i=1$ this was called the diagonal map $diag$, and if $a_i=i!$ this
was called $diag_1$. See Theorem~\ref{thm:diagonal}.

\begin{lemma}\label{lem:had-prod-diag}
  The map \eqref{eqn:had-prod-diag} defines a map
  $T\colon{}\gsubplus_2\longrightarrow\allpolypos$ iff $G=\sum a_i
  \frac{x^i}{i!i!}\in\allpolyposf$. 
\end{lemma}
\begin{proof}
  Since $T(e^xe^y) = \sum a_i \frac{x^i}{i!i!}$ the condition
  $G\in\allpolyposf$ is necessary.  From Theorem~\ref{thm:diagonal} we
  see that $a_i=i!$ determines a transformation that satisfies the
  conclusions of the theorem.  Consequently, if $G\in\allpolyposf$
  then we can express $T$ as a composition

\centerline{\xymatrix{
{\gsubplus_2} \ar@{->}[rr]^T \ar@{->}[d]_{a_i=i!} && \allpolypos \\
\allpolypos \ar@{->}[urr]_{x^i\mapsto \frac{a_i}{i!}x^i}
}}
The map $x^i\mapsto \frac{a_i}{i!}x^i$ defines a map
$\allpolypos\longrightarrow\allpolypos$ since
$\sum\frac{a_i}{i!}\,\frac{x^i}{i!}\in\allpolyf$ by hypothesis.
\end{proof}

Note that $diag_1(f(x)g(y)) = f\ast' g$.  
Consequently,  Theorem~\ref{thm:hadamard-2}, in the case both polynomials
are in $\allpolypos$, follows from this lemma, since we can use the
embedding
$\allpolypos\times\allpolypos\longrightarrow\gsubplusclose_2$ given by
$f\times g\mapsto f(x)g(y)$.

We next show that the only products of the form 
\begin{equation}\label{eqn:prod-factors}
x^i \times x^j \mapsto a_{i+j}x^{i+j}
\end{equation}
that determine maps $\allpoly\times\allpoly\longrightarrow\allpoly$
are compositions.

\begin{lemma}
A product \eqref{eqn:prod-factors} determines a map 
$T\colon{}\allpoly\times\allpoly\longrightarrow\allpoly$ iff $\sum a_n
\frac{x^n}{n!}\in\allpolyf$. This is equivalent to saying that the
product factors through $\allpoly$:

\centerline{\xymatrix{
\allpoly\times\allpoly \ar@{->}[rrr]^T
\ar@{->}[d]_{\text{multiplication}} &&& \allpoly\\
\allpoly \ar@{->}[urrr]_{\quad\quad\quad\text{multiplier transformation}}
}}
\end{lemma}
\begin{proof}
It suffices to evaluate $T(e^x,e^x)$:
\begin{align*}
  T(x^x,e^x) &=
\sum_{i,j=0}^\infty \frac{a_{i+j}}{i!j!}x^{i+j} \\
&= \sum_{n=0}^\infty a_n\frac{x^n}{n!}\sum_{i+j=n} \frac{n!}{i!j!} \\
&= \sum_{n=0}^\infty a_n\frac{(2x)^n}{n!}
\end{align*}
\end{proof}

We can only state necessary conditions in the more general cases.
The proofs are similar to the above, and are omitted. 
\begin{lemma}
  If the map \eqref{eqn:gen-had-1} induces
  $\rupint{2d}\longrightarrow\rup{d}$ then $\sum \aaa_\sdiffi
  \frac{\xx^\sdiffi}{(\diffi!)^2}\in\allpolyf_d$.
  If the map \eqref{eqn:gen-had-2} induces
  $\rupint{2d}\longrightarrow\rup{d}$ then $\sum \aaa_\sdiffi
  \frac{\xx^\sdiffi}{\diffi!}\in\allpolyf_d$.
  \begin{align}
    \label{eqn:gen-had-1}
    \xx^\diffi \yy^\diffj & \mapsto 
    \begin{cases}
      \aaa_\sdiffi \xx^\diffi & \diffi=\diffj \\
      0 & \text{otherwise}
    \end{cases}\\
    \label{eqn:gen-had-2}
    \xx^\diffi \yy^\diffj & \mapsto 
      \aaa_{\sdiffi+\sdiffj} \xx^{\diffi+\diffj} 
    \end{align}
  \end{lemma}
\section{Generating functions on $\allpoly(n)$}
\label{sec:gen-fct-pn}

The utility of generating functions is that if a certain function (the
generating function) computed from a linear transformation is in some
space, then the linear transformation has some nice mapping
properties. In this section we look at linear transformations such as
$x^r \mapsto x^{n-r}$ that are only defined on $\allpoly(n)$. The
following proposition does not characterize such generating functions,
since there is a necessary factorial.


\begin{prop} \label{prop:gen-fct-pn}
  Suppose that $T$ is a linear transformation of polynomials of degree
  at most $n$ such that $T_\ast(x+y)^n\in\gsubclose_2$. Then, the linear
  transformation $x^i \mapsto \frac{T(x^i)}{i!}$ maps $\allpoly(n)$ to
  itself.
\end{prop}
\begin{proof}
Choose $g(x) = \sum a_i x^i\in\allpoly(n)$. Since
$g^{rev}(\frac{\partial}{\partial y})$ maps $\gsubclose_2$ to itself, we
see that
\begin{align*}
g^{rev}(\frac{\partial}{\partial y})T_\ast(x+y)^n\bigg{|}_{y=0} &=
g^{rev}(\frac{\partial}{\partial
  y})\sum_{i=0}^n\binom{n}{i}T(x^i)y^{n-i}\bigg{|}_{y=0}\\
&= \sum_{i=0}^n\binom{n}{i}T(x^i)\left(g^{rev}(\frac{\partial}{\partial
  y})y^{n-i}\bigg{|}_{y=0}\right)\\
&= \sum_{i=0}^n \binom{n}{i} T(x^i) (n-i)! a_i\\
&= n!\, T\circ \expoper{} (g)
\end{align*}
and hence $T\circ \expoper{}(g)\in\allpoly(n)$. 

\end{proof}

The converse needs an extra hypothesis.

\begin{prop}
  Suppose $T\colon{}\allpoly(n)\longrightarrow\allpoly(n)$. If 
  $T_\ast(x+y)^n$ satisfies the  homogeneity condition then
   $T_\ast(x+y)^n\in\rupint{2}(n)$.
\end{prop}
\begin{proof}
  We only need to check substitution, and if we substitute $\alpha$
  for $y$ then $T_\ast(x+y)^n(x,\alpha) = T(x+\alpha)^n\in\allpoly$.
\end{proof}

\begin{example}
Here are a few examples.
\begin{enumerate}
\item \Mobius\ transformations. Suppose $T(x^i) =
  (ax+b)^i(cx+d)^{n-i}$. Then
$$ T_\ast(x+y)^n = \sum_{i=0}^n \binom{n}{i} (ax+b)^i(cx+d)^{n-i}y^{n-i}
  = (\,(ax+b) + y(cx+d)\,)^n$$
\item Polar derivative. If $T(x^i) = (n-i)x^i$ then
$$ T_\ast(x+y)^n = \sum_{i=0}^n \binom{n}{i} (n-i)x^iy^{n-i} =
ny(x+y)^n$$
\item Reversal (with a negative sign). If $T(x^i) = (-x)^{n-i}$ then
$$ T_\ast(x+y)^n = \sum_{i=0}^n \binom{n}{i} (-1)^{n-i}x^{n-i}y^{n-i} = (1-xy)^n$$

\item Hermite. \index{Hermite polynomials!transformations} Let $T(x^n) =
  H_n$. From~\eqref{eqn:iden-diag-7} we know that 
  $$
  T_\ast(x+y)^n = T(x^n)(x+2y) = H_n(x+2y) \in\rupint{2}.$$
  This
  implies that $x^n\mapsto H_n/n!$ maps $\allpoly(n)$ to itself, but
  in this case we know $T$ maps $\allpoly(n)$ to itself.

\item  In Lemma~\ref{lem:binomial-xdk} we saw that the linear
  transformation $T\colon{}x^k\mapsto \falling{x+n-k}{n}$ maps
  $\allpolyalt(n)\longrightarrow\allpoly(n)$. The generating function

$$ T_\ast(x+y)^n = \sum_{i=0}^n \binom{n}{i} \falling{x+n-k}{n}y^{n-i}
= (y+1)^n \rising{x}{n}
$$
 is in $\gsubclose_2(2n)$.  Proposition~\ref{prop:gen-fct-pn}
only allows us to conclude that
\begin{equation}
\label{eqn:xnk} x^k \mapsto \frac{\falling{x+n-k}{n}}{k!} \text{ maps }
\allpoly(n)\longrightarrow \allpoly(n).
\end{equation}
Since it is not true that $ x^k \mapsto {\falling{x+n-k}{n}}$ maps
$\allpoly(n)$ to itself, we see that the factorial in
Proposition~\ref{prop:gen-fct-pn} is necessary.

\end{enumerate}

\end{example}

\section{Higher dimensional generating functions}
\label{sec:new-analytic-higher}

If $T\colon{}\gsubpos_d \longrightarrow
\gsubpos_d$ is a linear transformation, then its generating function
is a function of $2d$ variables:
\[
 T_\ast(e^{-\xx\yy})= \sum_\sdiffi T(\xx^\sdiffi) \frac{(-\yy)^\sdiffi}{\diffi!}
\]

Here are some examples.

  The generating function of the identity transformation is  $e^{-\xx\cdot\yy}$.

Suppose that $T\colon{}\allpoly\longrightarrow\allpoly$ and let $F(x,y)$ be the
generating function of $T$. If we define $T_\ast(x^iy^j) = T(x^i)y^j$
then the generating function $F_\ast$ of $T_\ast$ is
\begin{align*} F_\ast(x,y,u,v) &= \sum_{i,j=0}^\infty T(x^i) y^j \frac{(-u)^i}{i!}
\frac{(-v)^j}{j!}  \\ &= \left(\sum_{i=0}^\infty T(x^i) \frac{(-u)^i}{i!}\right)
\ \left(\sum_{j=0}^\infty y^j \frac{(-v)^j}{j!}\right) = F(x,u)e^{-yv}
\end{align*}
\noindent%
More generally, if we have another linear transformation
$S\colon{}\allpoly\longrightarrow\allpoly$ with generating function $G(x,y)$
then  the generating function of $x^iy^j \mapsto T(x^i)S(y^j)$ 
is  $F(x,u)G(y,v)$.

If $T\colon{}\gsubpos_2 \longrightarrow \gsubpos_2$ is a linear
transformation and $a,b$ are positive, then $S(f) = a
\frac{\partial}{\partial x}T(f) + b \frac{\partial}{\partial x}T(f)$
satisfies $T(f)\lesslesseq S(f)$. If $T$ has generating function $F$,
then the generating function of $S$ is $ a F_x + b F_y$.

The generating function of the linear transformation
$g\mapsto f(\partial_x,\partial_y)g$  has generating
function $f(u,v)e^{-xu-yv}$. This follows using linearity from the calculation
\begin{align*} \sum_{i,j} \partial_x^r \partial_y^s (x^i y^j)
\frac{(-u)^i}{i!}\frac{(-v)^j}{j!} & = 
\left(\sum_i\partial_x^r x^i \frac{(-u)^i}{i!}\right)\,
\left(\sum_j\partial_y^r y^j \frac{(-v)^j}{j!}\right) \\ &=
\left(u^r e^{-xu}\right)\, 
\left(v^s e^{-yv}\right) 
\end{align*}
In general, the generating function of $g\mapsto f(\partial_\xx)g$
is $f(\xx)e^{-\xx\cdot\yy}$. 
Similarly, the generating function of $f\mapsto
f(\partial_\xx)g$ is $g(\xx+\yy)$.



\begin{prop}\label{prop:ps-2d}
If $T$ is a linear transformation on $\gsubclose_{d}$ and
$T_\ast(e^{-\xx\cdot\yy})\in\gsubf_{2d}$ then
  $T\colon\gsubclose_d\longrightarrow\gsubclose_d$. 
\end{prop}
\begin{proof}
Since $f(\partial_\xx)$ maps $\gsubf_{2d}$ to itself the proof is the
same as the case $d=1$. 
\end{proof}

  \begin{cor}\index{Hermite polynomials!in $\gsubpos_d$} 

    If $S$ is a d by d positive symmetric matrix, and
    $H_{\sdiffi}(\xx)$ is the corresponding Hermite polynomial then the
    linear transformation
$ \xx^\sdiffi \mapsto H_\sdiffi(\xx)$
maps $\gsubclose_d$ to $\gsubclose_d$.
  \end{cor}
  \begin{proof}
    The generating function of the Hermite polynomials
    \eqref{eqn:hermite-d} is $\exp(-\yy S\yy^* -2\yy S\xx^*)$ which is in
    $\gsubf_{2d}$. We can now apply the Proposition.
  \end{proof}

  The proposition implies that differential operators preserve
  $\rupint{d}$.

\begin{lemma}\label{lem:p2-diff}
  If $f\in\rupint{d}$ then the linear transformation
$$ T\colon{}g \mapsto f(\partial_\xx) g$$
maps $\rupint{d}$ to itself.
\end{lemma}
\begin{proof}
  The generating function of $T$ is $f(\yy)e^{-\xx\cdot\yy}$ which is in
  $\gsubf_{2d}$. 
\end{proof}

Table~\ref{tab:gen-fct-pd} lists some generating functions in higher
dimensions. Note that these transformations do not all map $\rupint{d}$ to
itself since not all of the generating functions are in $\gsubf_{2d}$.
The operator $\expoper{}_\xx$ is
$\expoper{}_{x_1}\cdots\expoper{}_{x_d}$, and the generalized Hurwitz
transformation is defined as

\begin{equation}\label{eqn:hurwitz-in-pd}
 T(\xx^\sdiffi) = \begin{cases}
\xx^{\sdiffi/2} & \text{ if all coordinates of $\diffi$ are even} \\
0 & \text{otherwise}
\end{cases}
\end{equation}

\begin{table}[htbp]
$$
\begin{array}{lrrclll}
\toprule
\text{Name} &&&&& \multicolumn{1}{l}{\qquad\qquad T_\ast(e^{-xy})} \\
\midrule
\displaystyle 
\genfctpd{differentiation}{f(\partial_\xx)\xx^\sdiffi}{f(\yy)
  e^{-\xx\cdot\yy}}{}
\genfctpd{Hermite}{H_\diffi(\xx)}{e^{-\yy S \yy^\prime - 2\yy S
    \xx^\prime}}{}
\genfctpd{Hadamard}{f(\xx)\ast
  \xx^\sdiffi}{\expoper{}_\xx\,f(-x_1y_1,\dots,-x_d y_d)}{}
\genfctpd{Hurwitz}{\xx^{\sdiffi/2}}{ \cosh(y_1\sqrt{x_1})\cdots\cosh(y_d\sqrt{x_d})}{}
\genfctpd{identity}{\xx^\sdiffi}{e^{-\xx\cdot\yy}}{}
\bottomrule
\end{array}
$$
  \caption{Generating Functions in higher dimensions}
  \label{tab:gen-fct-pd}

\end{table}

\section{Higher dimensional multiplier transformations}
\label{sec:sec:p2-ps}

The \Polya-Schur theorem does not generalize to two variables. 
Although we will see that all multiplier transformations are products
of one dimensional multipliers, there are examples of transformations
that do not map $\rupint{2}$ to itself, yet have a two variable
generating function in $\gsubclose_2$.

\cite{bbs} observed that all higher dimensional multiplier
transformations are just products of one dimensional maps:

\begin{lemma}\label{lem:d-multiplier}
  The following are equivalent:
  \begin{enumerate}
  \item $T:\xx^\diffi \longrightarrow \alpha_\sdiffi \xx^\diffi$ maps
    $\gsubclose_d$ to itself.
  \item $T$ is a product of one-dimensional transformations.
  \end{enumerate}
\end{lemma}
\begin{proof}
  We first assume that $d=2$, so assume that
  $T(x^iy^j)=\alpha_{ij}$. Following \cite{bbs} we apply $T$ to two test functions
  \begin{align*}
    Tx^iy^j(1+x)(1+y) &= x^iy^j \bigl(\alpha_{i,j} + \alpha_{i+1,j}x +
    \alpha_{i,j+1}y + \alpha_{i+1,j+1}xy\bigr)\\
    Tx^iy^j(1+x)(1-y) &= x^iy^j \bigl(\alpha_{i,j} + \alpha_{i+1,j}x -
    \alpha_{i,j+1}y - \alpha_{i+1,j+1}xy\bigr).
  \end{align*}
  By  Proposition~\ref{prop:cancel-2} or  Lemma~\ref{lem:hb-basic} 
  both of the factors are in $\gsubclose_2$, and by
  Theorem~\ref{thm:product-4}
  \begin{align*}
    \alpha_{i,j}    \alpha_{i+1,j+1}  -  \alpha_{i+1,j} 
    \alpha_{i,j+1} & \ge0 \\
    -\alpha_{i,j}    \alpha_{i+1,j+1}  +  \alpha_{i+1,j} 
    \alpha_{i,j+1} & \ge0 
  \end{align*}
and therefore 
\begin{equation}
  \label{eqn:uhpp-test}
      \alpha_{i,j}    \alpha_{i+1,j+1}  =  \alpha_{i+1,j} 
    \alpha_{i,j+1} 
\end{equation}

  We now consider the \emph{support} $S$ of $T$ - that is, the set of
all $(i,j)$ such that $\alpha_{i,j}\ne0$.  If $r$ is a non-negative
integer then the composition 

\centerline{\xymatrix{ \allpoly \ar@{->}[rr]^{{x^n\mapsto x^ny^r} } &&
    \gsubclose_2 \ar@{->}[r]^{T } & \gsubclose_2
    \ar@{->}[rr]^{{x^ny^r\mapsto x^n} } && \allpoly }} 
\noindent%
determines a map
$\allpoly\longrightarrow\allpoly$.  Since this map is a multiplier
transformation in one variable we know that the non-zero coefficients
have no gaps, and so form an interval. Thus, all
intersections of the support $S$ with a horizontal or vertical line
are intervals. The identity \eqref{eqn:uhpp-test} shows that
we can not have the configurations
\[
\begin{array}{c|ccc|c}
  \ne0 & 0 &\qquad \qquad\qquad& \ne0 & \ne 0 \\
  \cline{1-2}  \cline{4-5}
  \ne0 & \ne0 & & 0 & \ne 0
\end{array}
\]
in the support, so it follows that the support is a rectangle.

We may assume that $\alpha_{0,0}=1$. Using the recursion
\eqref{eqn:uhpp-test} and induction yields $\alpha_{ij} =
\alpha_{i,0}\,\alpha_{0,j}$.  It follows that $T$ is a product of the
one dimensional transformations
\begin{align*}
  T_r:&x^n\mapsto \alpha_{n,0}x^n &
  T_c:&y^n\mapsto \alpha_{0,n}y^n
\end{align*}

The case for general $d$ is no different; we use the fact that the
intersection with every $d-1$
dimensional face is a product to conclude the support is a product. 
\end{proof}


\begin{example}
  Here is an example of a multiplier transformation that does not map
  $\rupint{2}\longrightarrow\gsub_2$ but whose two variable generating
  function is in $\gsubf_2$. The two variable generating function of
  $T\colon{}x^iy^j\mapsto (i+j)^2 x^iy^j$ is

$$
  \sum_{i,j=0}^\infty
  \frac{(i+j)^2}{i!j!}x^iy^j = e^{x+y}(x+y)(x+y+1)$$
  which is in
  $\gsubf_2$. We know $T$ maps $\gsubplus_2$ to itself, but $T$ does not
  map $\rupint{2}$ to itself.  If
 
$$ f = -193 + 43\,x + 21\,x^2 + x^3 + 262\,y + 196\,x\,y + 14\,x^2\,y
+ 448\,y^2 + 64\,x\,y^2 + 96\,y^3$$
\noindent%
then $Tf(-2,y)$ has complex roots. Note that the usual generating
function is

$$ \sum_{i,j=0}^\infty (i+j)^2  x^iy^j \frac{(-u)^i(-v)^j}{i!\,j!} = 
e^{-ux-vy}(ux+vy)(ux+vy-1)
$$
and the latter expression is not in $\gsubf_4$.

\end{example}

\begin{remark}
  In Corollary~\ref{cor:ij-pd} we saw  that the linear transformation
  $x^i y^j \mapsto \dfrac{x^iy^j}{(i+j)!}$ maps $\gsubplus_2$ to
  itself, yet it does not map $\gsubpos_2$ to itself. The generating
  function is
$$ \sum_{i,j} \frac{x^i y^j}{(i+j)!} = \frac{xe^x - ye^y}{x-y}. $$
This can be seen as follows. If the generating function is $s$, then
the terms of $\frac{x}{y}(s-1)$ largely cancel, leaving only
$\frac{x}{y}e^x - e^y$. This generating function is not in
$\allpolyf_2$ since it is not an entire function.
\end{remark}

\index{coefficients!leading}

Just as in one variable (Lemma~\ref{lem:leading-coef}) a linear transformation
$T\colon{}\gsubposclose_2\longrightarrow\gsubposclose_2$ determines a simpler
linear transformation.

\begin{lemma} \label{p2:leading-coef}
  If $T\colon{}\gsubposclose_2\longrightarrow\gsubposclose_2$, and we can
  write
$$ T(x^iy^j) = g_{ij}(x,y) = c_{ij} x^iy^j + \text{ terms of lower
  degree}$$
then the linear transformation $S(x^iy^j)= c_{ij}x^iy^j$ also maps
$\gsubposclose_2\longrightarrow\gsubposclose_2$.
\end{lemma}
\begin{proof}
  Choose $f = \sum a_{ij}x^iy^j$ in $\gsubposclose_2$. Substitute
  $x/\alpha$ for $x$, $y/\beta$ for $y$ and apply $T$. Next, substitute
  $\alpha x$ for $x$ and $\beta y$ for $y$. The result is that the
  following polynomial is in $\gsubposclose_2$:
  \begin{align*}
    \sum_{i,j} a_{ij} \alpha^i \beta^j g_{ij} ( x/\alpha,y/\beta) & \\
\intertext{Next, using the fact that}
\lim_{\alpha\rightarrow0^+}\lim_{\beta\rightarrow0^+} \alpha^i \beta^j
g_{ij}(x/\alpha,y/\beta) &= c_{ij}x^iy^j\\
\intertext{we see that }
\lim_{\alpha\rightarrow0^+}\lim_{\beta\rightarrow0^+} \sum_{i,j}
a_{ij}\alpha^i \beta^j
g_{ij}(x/\alpha,y/\beta) &= \sum_{i,j} a_{ij}c_{ij}x^iy^j  = 
S(f) \in\gsubposclose_2\\
  \end{align*}
\end{proof}

\begin{example}
  The two variable generating function of $x^ry^s\mapsto f(r+s)x^ry^s$
  can be easily found.
\begin{align*} \sum (i+j)^n \frac{x^iy^j}{i!j!} &= e^{x+y}
T_{\ast\ast}^{-1}(x+y)^n \\
\intertext{where $T\colon{}x^n=\falling{x}{n}$, and $T_{\ast\ast}^{-1}(x^ry^s) =
  T^{-1}(x^r)T^{-1}(y^s)$. By linearity we find the generating function}
\sum f(x+y) \frac{x^iy^j}{i!j!} &= e^{x+y}
T_{\ast\ast}^{-1}f(x+y) \,=\, (e^xT^{-1}f)(x+y)
\end{align*}
where the last equality follows from Example~\ref{ex:bin-type}. It
follows that if $f\in\allpolypos$ then this generating function is in
$\gsubf_2$. We saw that the linear transformation corresponding to
$f(x)=x^2$ does not map $\rupint{2}$ to itself.

\end{example}

\section{Differential operators}
\label{sec:gen-fct-diff-ops}

\index{differential operators} \index{operator!differential}

We continue the investigation of differential operators. We first
compute the generating functions of differential operators acting on
$\allpoly$, extract some properties, show some eigenpolynomials are in
$\allpoly$, and then extend some of the results to $\rupint{d}$.

Consider the linear transformation $T$ defined on polynomials of one
variable

$$
  T(g) =  \sum_{i=0}^n f_i(x) g^{(i)}(x)
$$
We define $f(x,y) = \sum f_i(x)y^i$, and write $T(g) =
f(x,\diffd)g$. The generating function of $T$ is

\begin{align*}
  \sum_{j=0}^\infty T(x^j) \frac{(-y)^j}{j!} &=
  \sum_{j=0}^\infty \sum_{i=0}^n f_i(x)\diffd^i\,x^j \frac{(-y)^j}{j!} \\
&= \sum_{i=0}^n f_i(x)\diffd^i  \sum_{j=0}^\infty \,x^j \frac{(-y)^j}{j!}
\\
&= \sum_{i=0}^n f_i(x)\diffd^i  e^{-xy} \\
&=  \sum_{i=0}^n f_i(x)(-y)^i  e^{-xy} \\
&= f(x,-y)e^{-xy}
\end{align*}

Thus, we have

\begin{prop} \label{prop:diff-op-poly}
If $f\in\gsubclose_2$ then $f(x,-\diffd)$ maps $\allpoly$ to
$\allpoly$. Conversely, if

\begin{enumerate}
\item $f(x,-\diffd)$ maps $\allpoly$ to $\allpoly$.
\item The coefficients of the homogeneous part of $f(x,y)$ are all positive.
\end{enumerate}
then $f\in\rupint{2}$.  
\end{prop}

\begin{proof}
  The generating function is $f(x,y)e^{-xy}\in\gsubf_2$. In
  Lemma~\ref{lem:diff-op-a} we saw that $f(x,\alpha)\in\allpoly$ for
  all $\alpha\in\reals$, so $f\in\rupint{2}$.
\end{proof}

\begin{remark}
  For example, if we take $F=x+y$ then the linear transformation is
  $f\mapsto xf-f'$. We know that this is in $\allpoly$ since
  $xf\lesslesseq f\lesslesseq f'$ implies $xf-f'\lesslesseq f$.
\end{remark}

\begin{cor}
  Suppose that $T\colon{}\allpoly\longrightarrow\allpoly$. The linear
  transformation
\[ g\mapsto \sum T(x^k) \frac{g^{(k)}(x)}{k!}(-1)^k\]
maps $\allpoly\longrightarrow\allpoly$.
\end{cor}
\begin{proof}
  If $f(x,y)$ is the   the generating function of $T$, then
  $f(x,y)\in\rupint{2}$. The map is $g(x) \mapsto f(x,-\diffd)g(x)$. 
\end{proof}

  In \cite{tanja01} they show that for $n$ sufficiently large there is
  an eigenpolynomial $f(x,\diffd)$ of degree $n$.  We now show that
  the eigenpolynomials  are in $\allpoly$ if
  $f\in\gsubclose_2$, and $f$ satisfies a degree condition.
  \begin{lemma}\label{lem:eigenpoly}
    Choose $f(x,y)\in\gsubclose_2(d)$
    $$
    f(x,y) = f_0(x) + f_1(x)y + f_2(x)y^2 + \cdots + f_d(x)y^d$$
    where the
    degree of $f_i$ is $i$. If $n$ is sufficiently large then there is
    a polynomial $p\in\allpoly(n)$ and constant $\lambda$ such
    that $f(x,\diffd)p=\lambda p$.
  \end{lemma}
  \begin{proof}
    Since $f(x,\diffd)$ maps $\allpoly$ to itself, and preserves
    degree, we will apply Lemma~\ref{lem:eigen-1}. It suffices to show
    that there is a dominant eigenvalue, which is the same as finding
    a dominant leading coefficient. 

    

    The $r$-th diagonal element of $M$ is the coefficient
    $c_r$ of $x^r$ in $f(x,\diffd)x^r$.  Denoting the leading
    coefficient of $f_i$ by $C_i$,

\index{eigenpolynomial}

    $$
    c_d = \sum_{i=0}^d c_i\, \falling{r}{i} = C_d \falling{r}{d} +
    O(r^{d-1}).$$
    For $n$ sufficiently large, $c_n$ is the largest
    eigenvalue of $M$. 

  \end{proof}

  \begin{remark}
    The identity $e^{-(\partial_x+\partial_y)^2}(x+y)^n = (x+y)^n$
    shows that the operator $e^{-(\partial_x+\partial_y)^2}$ has
    polynomial eigenvalues of every degree.
  \end{remark}

We can generalize some of the above to more variables. Suppose that
$f(\xx,\yy)\in\rupint{d+e}$ and we define a differential operator

\begin{equation}\label{eqn:diff-op-d}
  T(g(\xx)) =  \sum_\sdiffi f_\sdiffi(\xx) \diffd^\sdiffi g(\xx)
\end{equation}

It is easy to see that the generating function of $T$ is simply
$f(\xx,\yy)e^{-\xx\cdot\yy}$.

\begin{prop}
Suppose
\begin{enumerate}
\item $f(\xx,-\diffd)$ maps $\rupint{d}$ to $\gsub_d$.
\item The coefficients of the homogeneous part of $f(\xx,\yy)$ are all positive.
\end{enumerate}
then $f\in\rupint{d+e}$.  

Conversely, if $f(\xx,\yy)\in\gsubclose_{2d}$ then $f(\xx,-\diffd)$
maps $\gsubclose_d$ to itself.
\end{prop}
\begin{proof}
  If suffices to show that $f(\xx,\yy)$ satisfies substitution. The
  proof is the same as Lemma~\ref{lem:diff-op-a}: apply $f(\xx,-\diffd)$ to
  $e^{-\aaa\cdot\xx}$ to conclude that
  $f(\xx,\aaa)e^{-\aaa\cdot\xx}\in\rupint{d}$. Multiplying by
  $e^{\aaa\cdot\xx}$ finishes the proof.

  The second part follows form the proposition.
\end{proof}

\index{f(x,D)}
  Which polynomials $f(x,y)$ determine operators $f(x,\diffd)$ that
  map $\allpolypos\longrightarrow\allpoly$, and don't map
  $\allpoly\longrightarrow\allpoly$? Here's a necessary condition. 

  \begin{lemma}
    Suppose that  $f(x,y)$ determine an operator $f(x,\diffd)$ that
  maps $\allpolypos\longrightarrow\allpoly$. If $f(x,-y)$ satisfies degree
  and positivity then $f(x,-y)\in\partialpoly{1,1}$.
  \end{lemma}
  \begin{proof}
    Let $f(x,y) = \sum a_{ij}x^iy^j$.  Since $e^{-\alpha
      x}\in\allpolypos$ for positive $\alpha$, we see that
$$ \left(\sum a_{ij} x^i\diffd^j\right)e^{-\alpha x} =
\left(\sum a_{ij} x^i(-\alpha)^j\right)e^{-\alpha x} 
$$
is in $\allpolyf$, and hence $f(x,-\alpha)\in\allpoly$ for all positive
$\alpha$. Since $f(x,-y)$ satisfies degree and positivity,
$f(x,-y)\in\partialpoly{1,1}$. 
  \end{proof}

\section{Generating functions for maps $\gsubpos_2\longleftrightarrow\allpoly$}
\label{sec:gen-fct-p2-p}

If $T\colon{}\gsubpos_2\longrightarrow\allpoly$ then the generating
function of $T$ lies in $\allpolyf_3$.  We compute the generating
functions for many of these maps.

The most basic map of all is evaluation. The generating function of
$f(x,y)\mapsto f(z,\alpha)$ is

$$ \sum T(x^iy^j)\frac{(-u)^i(-v)^j}{i!j!} = 
\sum z^i \alpha^j\frac{(-u)^i(-v)^j}{i!j!} = e^{-uz - \alpha v}$$

Next, consider the diagonal map $f(x,y)\mapsto f(z,z)$. The generating
function is

$$ \sum T(x^iy^j)\frac{(-u)^i(-v)^j}{i!j!} = 
\sum z^{i+j}\frac{(-u)^i(-v)^j}{i!j!} = e^{-(u + v)z}$$

The map that extracts the coefficient of a fixed monomial also maps
$\rupint{2}\longrightarrow\allpoly$. If $T(x^iy^j) = z^i$ if $j=k$ and
$0$ otherwise then the generating function is
$$ \sum T(x^iy^j)\frac{(-u)^i(-v)^j}{i!j!} = 
\sum_i z^i  \frac{(-u)^i(-v)^k}{i!k!} = \frac{(-v)^k}{k!} e^{-uz} $$

Finally, we consider the map $f(x,y)\mapsto f(\diffd,x)$. More
precisely, this is the map $x^i y^j \mapsto \frac{d^i}{dz^i}z^j =
\falling{j}{i}z^{j-i}$. 

\begin{align*}
  \sum_{i\le j} \falling{j}{i} z^{j-i} \frac{(-u)^i(-v)^j}{i!j!} &=
\sum_{j=0}^\infty \frac{(-zv)^j}{j!}  \sum_{i=0}^j
\binom{j}{i}\left(\frac{-u}{z}\right)^i \\ &=   
\sum_{j=0}^\infty \frac{(-zv)^j}{j!} \left(1+\frac{-u}{z}\right)^j = e^{-zv + uv}
\end{align*}

We summarize these few results in Table~\ref{tab:gen-fct-p2-p}

\begin{table}[htbp]
$$
\begin{array}{lrrclll}
\toprule
\text{Name} &&&& \multicolumn{1}{l}{\qquad T_\ast(e^{-\xx\yy})} \\
\midrule
\displaystyle 
\genfctpp{evaluation at $\alpha$}{f(z,\alpha)}{  e^{-uz-\alpha v}}{}
\genfctpp{diagonal}{f(z,z)}{e^{-z(u+v)}}{}
\genfctpp{coefficient}{\text{coef. of $y^k$}}{\frac{v^k}{k!}e^{-uz}}{}
\genfctpp{derivative}{f(-\diffd,z)}{e^{-v(u+z)}}{}
\bottomrule
\end{array}
$$
  \caption{Generating Functions for  $\gsubpos_2\longrightarrow\allpoly$}
  \label{tab:gen-fct-p2-p}

\end{table}

We can construct transformations
  $\allpoly\longrightarrow\rupint{2}$ by homogenizing  transformations
  $\allpoly\longrightarrow\allpoly$.

  \begin{lemma}\label{lem:p-to-p2}
    Suppose $T\colon\allpoly\longrightarrow\allpoly$ and
    $T(x^n)\in\allpolypos(n)$. If $T(x^i)=f_i(x)$ and
    $T_H(x^i)=f_i(x/y)y^i$ then
    $T_H:\allpoly\longrightarrow\rupint{2}$. 
  \end{lemma}
  \begin{proof}
    If $g=\sum a_ix^i\in\allpoly(n)$ then $ T_H(g) = \sum a_i f_i(x/y)
    y^i.$  The homogeneous part of $T_H(g)$ is $f_n(x/y)y^n$ which has
    all positive terms since $T(x^n)\in\allpolypos$. Substituting
    $\alpha$ for $y$ yields
\begin{multline*}
  T_H(g)(x,\alpha) = \sum a_i f_i(x/\alpha)\alpha^i = \sum (\alpha^i
  a_i) f_i(x/\alpha) = T(g(\alpha x))(x/\alpha)
\end{multline*}
which is in $\allpoly$.
  \end{proof}

If the generating function of $T$ is $G(x,y)$ then the generating
function of $T_H$ is $G(x/y,yz)$.

Note that evaluation at $x=1$ shows that $T\colon x^i\longrightarrow
f_i^{rev}$ maps $\allpoly\longrightarrow\allpoly$. This gives another
proof of                                     Corollary~\ref{cor:t1/z}.

  \begin{example}\index{Laguerre polynomials} 
    If $T$ is the affine transformation $x\mapsto x+1$ then $T_H(g) =
    g(x+y)$. For a more complicated example, if $T(x^i) = L_i(-x)$ is
    the Laguerre transformation then 
\[ T_H(x^n) = \sum_0^n \frac{1}{k!} \binom{n}{k}x^ky^{n-k} \]
  \end{example}

\section{Linear transformations $\allpolypos\longrightarrow\allpoly$}
\label{sec:gen-fct-partial}

In this section we consider maps
$\allpolypos\longrightarrow\allpolypos$.  If we have a linear transformation
$T\colon{}\allpolyalt\longrightarrow\allpolyalt$ then we can construct a map
$S\colon{}\allpolypos\longrightarrow\allpolypos$ by $S(f)(x) = T(f(-x))(-x)$.
The two generating functions satisfy
    \begin{align*}
      f(x,y) &= \sum T(x^n) \frac{(-y)^n}{n!} \\
      g(x,y) &= \sum S(x^n) \frac{(-y)^n}{n!}\,=\, \sum (-1)^nT(x^n)(-x) \frac{y^n}{n!} \\
    \end{align*}
and so $g(x,y) = f(-x,-y)$.
    
\begin{theorem} \label{thm:partial-gen-fct}
  Suppose $T\colon{}\allpolypos\longrightarrow\allpolypos$, $T$ preserves
  degree and the sign of the leading coefficient. If
  \begin{equation}
    \label{eqn:partial-gen-fct}
 F(x,y) =    T_\ast(e^{-xy}) = \sum_{i=0}^\infty T(x^i) \frac{(-y)^i}{i!} 
  \end{equation}
then $F(x,y)\in\partialpolyf{1,1}$. If
$T\colon{}\allpolyalt\longrightarrow\allpolyalt$ then $F(x,-y)\in\partialpolyf{1,1}$.
\end{theorem}
\begin{proof}
It suffices to assume that $T\colon{}\allpolypos\longrightarrow\allpolypos.$
Since $e^{-xy}\in\gsubf_2\subset\partialpolyf{1,1}$, we know that
$F(x,y) = T_\ast(e^{-xy})\in\partialpolyf{1,1}$.
\end{proof}

\begin{theorem}
  Suppose that $F(x,y)=\sum f_i(x)\dfrac{(-y)^i}{i!}$ is in
  $\partialpolyf{1,1}$. If each $f_i$ is a polynomial and we define
  $T(x^i) = f_i/i!$ then $T\colon{}\allpolypos\longrightarrow\allpoly$.
\end{theorem}

\begin{proof}
  We follow the proof for $\gsub_d$, but 
  since $\partialpoly{1,1}$ is not closed under differentiation we
  need to multiply. This is why we have the factorial appearing.
\end{proof}

Since $\partialpoly{1,1}$ is closed under multiplication, so is
$\partialpolyf{1,1}$. Look at Table~\ref{tab:pos-gen-fct}.
We know that $e^{-\alpha x}\in\gsubf_2$, and
$(1+y)^x\in\partialpolyf{1,1}$ since the latter is the generating
function of a linear transformation
$\allpolyalt\longrightarrow\allpolyalt$. Their product is in
$\partialpolyf{1,1}$, and is the generating function for the Charlier
transformation.

Table \ref{tab:pos-gen-fct} lists the generating functions of some
linear transformations $\allpolyalt\longrightarrow\allpolyalt$.
Table~\ref{tab:pos-gen-fct-2} lists the generating functions of linear
transformations $\allpolypos\longrightarrow\allpoly$ whose image is not
$\allpolypos$. 
Notice that most of these generating functions are not in
$\allpolyf_2$ - indeed, they are not analytic for all real values of
the parameters, but only positive $y$.

\begin{table}[htbp]
  \centering
$$
\begin{array}{lrrclll}
\toprule
\text{Name} &&&&& \multicolumn{1}{l}{\qquad\qquad T_\ast(e^{-xy})} \\
\midrule
\genfct{Binomial}{\displaystyle\binom{x}{i}}{\,_1F_1(-x,1,-y)}{}
\genfct{Charlier}{C_n^\alpha}{e^{-\alpha y}(1+y)^x}{}
\genfct{Falling Factorial}{\falling{x}{i}}{(1-y)^{x}}{}
\text{Hurwitz}&\makebox[.1in]{} & \displaystyle \sum_{i=0}^\infty &  x^i  &%
\displaystyle\frac{y^{2i}}{(2i)!} &\quad=\quad \cosh(\sqrt{xy^2})&{}\\
\text{Hurwitz}&\makebox[.1in]{} & \displaystyle \sum_{i=0}^\infty &  x^i  &%
\displaystyle\frac{y^{4i}}{(4i)!} &\quad=\quad
\frac{1}{2}\left(\cos(x^{1/4}y)+\cosh(x^{1/4}y)\right){}\\ 
\genfct{q-series}{(-1)^{\binom{i}{2}}x^i}{\sqrt{2} \cos (xy -
  \pi/4)}{}
\genfct{Rising Factorial}{\left(\rising{x}{i}\mapsto x^i\right)}{\displaystyle e^{x(1-e^{-y})}}{}
\bottomrule
\end{array}
$$

  \caption{Generating functions for
    $\allpolyalt\longrightarrow\allpolyalt$}
  \label{tab:pos-gen-fct}
\end{table}

\begin{table}[htbp]\label{tab:pos-gen-fct-2}
  \centering
$$
\begin{array}{lrrclll}
\toprule
\text{Name} &&&&& \multicolumn{1}{l}{\qquad\qquad T_\ast(e^{-xy})} \\
\midrule
\genfct{Hermite}{H_n^{rev}}{e^{2y-x^2y^2}}{}
\genfct{Hermite}{x^nH_n}{e^{2x^2y-x^2y^2}}{}
\genfct{Laguerre}{L_n^{rev}}{e^{xy}J_0(2\sqrt{y})}{}
\bottomrule
\end{array}
$$
\caption{Generating functions for $\allpolypos\longrightarrow\allpoly$}
\end{table}

\begin{example}
  We can determine the closed formula for the generating function
  $F(x,y)$ of $T\colon{}x^i\mapsto (-1)^{\binom{i}{2}} x^i$ given in
  Table~\ref{tab:pos-gen-fct} by considering four sums determined by
  the index mod $4$. The result is
\begin{align*} F(x,y) = \sum_{i=0}^\infty (-1)^{\binom{i}{2}}\, x^i\, \frac{y^i}{i!}
  & =\sqrt{2} \cos (xy - \pi/4). \\
  \intertext{The cosine has a representation as an infinite product}
  \cos(x) &= \prod_{k=0}^\infty \left( 1 -
    \frac{4x^2}{(2k+1)^2\pi^2}\right) \intertext{so the generating
    function can be expressed as the infinite product} F(x,y) &=
  \prod_{k=0}^\infty \left( 1 -
    \frac{4(xy-\frac{\pi}{4})^2}{(2k+1)^2\pi^2}\right).
  \intertext{The product is not in $\allpolyf_2$ as  the factorization
    shows}
    \left( 1 -
    \frac{4(xy-\frac{\pi}{4})^2}{(2k+1)^2\pi^2}\right) &= \left(1-
    \frac{2(xy-\frac{\pi}{4})}{(2k+1)\pi}\right) \left(1+
    \frac{2(xy-\frac{\pi}{4})}{(2k+1)\pi}\right)\\ &= F_1(x,y)F_2(x,y)
\end{align*}
$F_1(x,y)\in\allpolyf_2$, but $F_2(x,-y)\in\allpolyf_2$. 
This reflects the fact that the linear transformation
$x^i\mapsto (-1)^{\binom{i}{2}}x^i$ maps $\allpolypos$ to
$\allpoly$, and does not map $\allpoly\longrightarrow\allpoly$.
\end{example}

\begin{example}
  One way to find elements in $\partialpoly{1,1}$ is to
  compute the generating function of  simple linear
  transformations. The generating function of $U(x^n)=0$ if $n$ odd, and $x^n$
  if $n$ even is
\begin{align*}
    \sum_{n=0}^\infty x^{2n}\frac{y^{2n}}{(2n)!} &= \cosh(xy) \\
\intertext{while the generating function of $V(x^n)=0$ if $n$ even,
  and $x^n$ if $n$ odd is}
    \sum_{n=0}^\infty x^{2n+1}\frac{y^{2n+1}}{(2n)!} &= \sinh(xy) \\
  \end{align*}
  Let $T_e(x^n) = 0$ if $n$ is odd, and $x^{n/2}$ if
  $n$ is even. We know that
  $T_e:\allpolypos\longrightarrow\allpolyneg$.  Consequently, the map
  $S(f)=T_e(f)(1-x^2)$ maps $\allpolypos\longrightarrow\allpoly$. The
  generating function of $T_e$ is for positive $x$
  \begin{align*}
   \sum_{n=0}^\infty x^n \frac{y^{2n}}{(2n)!} 
&= \cosh(y\sqrt{x})\\
\intertext{and the generating function of $S$ is
}
   \sum_{n=0}^\infty (1-x^2)^n \frac{y^{2n}}{(2n)!} 
&= \cos(y\sqrt{1-x^2}) \\
\intertext{Finally, let $T_o(x)=0$ if $n$ even, and $x^{(n-1)/2}$ if
  $n$ odd. The generating function of $T_o$ for positive $x$ is}
   \sum_{n=0}^\infty x^n \frac{y^{2n+1}}{(2n+1)!} 
&= \sqrt{x}\sinh(y\sqrt{x})
  \end{align*}
\end{example}

\section{More generating functions}
\label{sec:gen-finite}

In this section we list in Table~\ref{tab:gen-finite} a few
generating functions for linear transformations that map
$\allpolyint{finite}\longrightarrow\allpolyint{\diffj}$. 
These functions are not entire - in order to converge we must assume
that $|x|,|y|<1$.

\begin{table}[htbp] \label{tab:gen-finite}
  \begin{center}
$$
    \begin{array}{lrrclll}
\toprule
\text{Name} &&&&& \multicolumn{1}{l}{\qquad T_\ast(e^{-xy})} \\
\midrule
\genfct{Chebyshev}{T_i}{e^{xy}cosh(\sqrt{y^2(x^2-1)})}{}
\genfct{Legendre}{P_i}{e^{xy}{\,_0F_1(;1;\frac{1}{4}(y^2-1)x^2)}}{}
\genfct{Euler}{A_i}{\displaystyle \ \frac{e^{xy} - x
  e^{xy}}{e^{xy}-xe^y}}{}\\
\genfct{Factorial}{\prod_{k=1}^i(1-kx)}{(1+xy)^{-1+1/x}}{}
\bottomrule
    \end{array}
$$
    \caption{Miscellaneous generating functions}
    
  \end{center}
\end{table}

Since the method of finding bivariate generating functions may not be
familiar, I'll describe how to find the generating function for the
Euler transformation.
We compute the generating function of $x^n\mapsto A_n$ by first
computing the generating function of $x^n \mapsto B_n =
x^nA_n(\frac{x+1}{x})$. We use the recurrence for $B_n$ to derive a
partial differential equation for the generating function, which we
use the method of characteristics to solve.  So recall
(Lemma~\ref{lem:euler-mobius}) that $B_n$ satisfies the recursion $B_{n+1} =
(x+1)(B_n+xB_n^\prime)$. Substituting this into the generating
function for $T$
\begin{align}
  F(x,y) &= \sum_{n=0}^\infty B_n\ \frac{y^n}{n!}\notag\\
&= \sum_{n=0}^\infty (x+1)(B_n+xB_n^\prime) \ \frac{y^n}{n!}\notag\\
\intertext{with some  manipulation yields the partial differential equation}
(x+1)F(x,y) &= (-x^2-x)F_x(x,y) + F_y(x,y)\label{eqn:pde}
\end{align}

The method of characteristics assumes that $x,y$ are functions of
$s,t$, where we consider $t$ a variable, and $s$ a parameter. We 
choose these functions to make each side of \eqref{eqn:pde} an
exact differential. Thus, we have to solve the three equations
\begin{align}
  \frac{d}{dt} x(s,t) &= - x^2 -x \label{eqn:pde-1}\\
  \frac{d}{dt} y(s,t) &= 1 \label{eqn:pde-2}\\
  \frac{d}{dt} F(s,t) &= (x+1)F(s,t) \label{eqn:pde-3}
\end{align}
with the initial conditions
\begin{align*}
  x(s,0) &= s\\
  y(s,0) &= 0\\
  F(s,0) &= 1
\end{align*}
We solve \eqref{eqn:pde-1} and \eqref{eqn:pde-2} for $x,y$, and then
use these two solutions to express  \eqref{eqn:pde-3} entirely in
terms of $s,t$. We solve that equation, solve for $s,t$ in terms of
$x,y$, substitute into our last solution, and we are done! Of course,
a computer algebra system is indispensable for these calculations.

The solution is $F(x,y)=\displaystyle\frac{e^y}{1+x-xe^y}$. We then
apply the first \Mobius\ transformation $z\mapsto 1/z$ which yields
the generating function $F(1/x,xy)$, and then $z\mapsto z-1$ which
gives $F(1/(x-1),(x-1)y)$ which simplifies to our result in the
table.


\chapter{Recurrences with polynomial coefficients }
\label{cha:recur}

\index{recurrence!with polynomial coefficients}

How can we construct sequences of real-rooted polynomials $\{f_i\}_{i=0}^\infty$
such that for all positive $n$ the $f_i$ satisfy the recurrence
\begin{equation}\label{eqn:rec-seq}
f_n(x) g_0(x) + f_{n+1}(x)g_1(x) + \cdots + f_{n+d}(x) g_d(x) = 0
\end{equation}
where the $g_i$ are polynomials in $x$? We say that the sequence
$\{f_i\}$ satisfies a \emph{recursion with polynomial coefficients.}
We will show several ways of doing constructing such
sequences. Unfortunately, we don't know how to reverse the
process. That is, given a sequence that satisfies a recursion with
polynomial coefficients, prove that all the terms have all real roots.

\begin{example}
Here are two simple examples. If $f_i = x^i$ then we have the
recurrence

\[
f_n(x)\cdot x - f_{n+1}(x) = 0.
\] 

  \index{Chebyshev polynomials}

A less trivial example is given by the Chebyshev polynomials $U_k$ of the second kind.
They satisfy the recurrence

\[
U_n(x) - U_{n+1}(x) \cdot 2x + U_{n+2} = 0
\]

In each of these cases the consecutive terms interlace, but
this will not be true in general.

The \index{Hermite polynomial}Hermite polynomial $H_n$ satisfies the
recursion 
\[ H_{n+1} = 2x H_n - 2n H_{n-1}
\]
This is \emph{not} a recursion with polynomial coefficients since the
factor $(2n)$ depends on the index. However, they do satisfy a
differential recursion with polynomial coefficients 
\[
H_n = 2x H_{n-1} - H_{n-1}'.
\]
\end{example}

Every sequence that satisfies a recursion with polynomial coefficients
also satisfies a determinantal identity. This has nothing to do with
the property of having real roots or not. Given a sequence as in
\eqref{eqn:rec-seq} we have the identity
\[
\begin{pmatrix}
  g_0 & g_1 & \cdots & g_d
\end{pmatrix}
\begin{pmatrix}
  f_n & f_{n+1} & \cdots & f_{n+d}\\
  f_{n+1} & f_{n+2} & \cdots & f_{n+d+1}\\
\vdots & & & \vdots \\
  f_{n+d} & f_{n+d+1} & \cdots & f_{n+2d}\\
\end{pmatrix}
=0
\]
If  the $g_i$ are not identically zero the determinant is zero:
\[
\begin{vmatrix}
  f_n & f_{n+1} & \cdots & f_{n+d}\\
  f_{n+1} & f_{n+2} & \cdots & f_{n+d+1}\\
\vdots & & & \vdots \\
  f_{n+d} & f_{n+d+1} & \cdots & f_{n+2d}\\
\end{vmatrix}
=0
\]

\section{The general construction}
\label{sec:general-construction}

Many of our recursions arise from repeatedly multiplying by a
polynomial. Here's a trivial example:
\begin{align*}
  f_0 &= f\\
  f_n &= g\, f_{n-1}
\end{align*}
where $f,g\in\allpoly$. More generally, we turn polynomial
multiplication into matrix multiplication, and the characteristic
polynomial becomes the recursion.

So, our construction requires
\begin{enumerate}
\item A $d$ by $d$ matrix $M$ that might have polynomial
entries, or perhaps even linear transformations as entries. 
\item A vector $v$ of polynomials.
\item For $i=0,1,\dots,$ all the entries of $M^i\,v$ have all real roots.
\end{enumerate}
If $\sum_0^d a_i\lambda^i$ is the characteristic polynomial of $M$ and
$n=0,1,\dots,$ then
\[
\sum_{i=0}^d  a_i \,(M^{i+n}v) = 0.
\]
Thus if we write 
\[
M^i\,v = (p_{i,1},\dots,p_{i,d})^t
\]
then each sequence of polynomials $\{p_{i,k}\}_{i=1}^\infty$ has all
real roots, and satisfies the same recursion
\[
\sum_{i=0}^d a_i p_{i+n,k} =0.
\]

The vector $v$ of polynomials can be
\begin{enumerate}
\item A collection of mutually interlacing polynomials.
\item All coefficients of a polynomial in $\gsubclose_2$.
\item The initial coefficients of a function in $\gsubf_2$.
\end{enumerate}

\section{Recursions from  $f(y)$ and $g(x,y)$}
\label{sec:rati-funct-from-p2}

In our first construction we always get the same recursion; this is
due to the fact that our matrix is lower triangular. We start with the
following data
\begin{align*}
  f(y) &= \sum_0^d a_i y^i & \in\allpoly \\
  g(x,y) &= \sum g_i(x)y^i & \in\gsubclose_2\\
  f(y)g(x,y) &= \sum h_i(x) y^i
\end{align*}
and we have the relation
\[
\begin{pmatrix}
  a_0 & 0 & 0 & \dots & 0 \\
  a_1 & a_0 & 0 & \dots & 0 \\
  a_2 & a_1 & a_0  & \dots & 0 \\
\vdots & & & \ddots & \\
  a_d & a_{d-1} & a_{d-2} & \dots & a_0
\end{pmatrix}
\begin{pmatrix}
  g_0 \\ g_1 \\ g_2 \\ \vdots \\ g_d
\end{pmatrix}
=
\begin{pmatrix}
  h_0 \\ h_1 \\h_2 \\ \vdots \\ h_d
\end{pmatrix}
\]

Since $f(y)g(x,y)\in\gsubclose_2$ it follows that all the coefficients
$h_0,\dots,h_d$ have all real roots. We  follow the outline in the
previous section. The characteristic polynomial is $(\lambda-a_0)^d$. 
If we write 
\[
M^i\,v = (p_{i,1},\dots,p_{i,d})^t
\]
then all $\{p_{i,k}\}$ satisfy
\[
p_{n,k} - \binom{d}{1}a_0 p_{n-1,k} + \cdots + (-1)^d a_0^d p_{n-d,k}=0
\]

Note that the recursion is effectively independent of $f$ if
$a_0\ne0$, since we can rescale to make $a_0=1$.  
Here's an example where  $d=3$:
\begin{align*}
  f &= (y+1)^3 & g &= (x+2y+1)(x+3y+2) \\
  M &=
  \begin{pmatrix}
    1 & 0 & 0 & 0 \\
    3 & 1 & 0 & 0 \\
    3 & 3& 1 & 0\\
    1 & 3 & 3 & 1
  \end{pmatrix} &
\end{align*}
\begin{align*}
  v &=
\begin{pmatrix}
   x^2+3 x+2 \\ 5 x+7 \\6\\0 
  \end{pmatrix}
& 
Mv &=
\begin{pmatrix}
x^2+3 x+2 \\ 3 x^2+14 x+13 \\ 3 x^2+24 x+33 \\ x^2+18 x+41
\end{pmatrix}\\
M^2 v &= 
\begin{pmatrix}
  x^2+3 x+2\\6 x^2+23 x+19\\15 x^2+75 x+78\\20 x^2+135 x+181
\end{pmatrix}
&
M^3v &= 
\begin{pmatrix}
  x^2+3 x+2\\9 x^2+32 x+25\\36 x^2+153 x+141\\84 x^2+432 x+474
\end{pmatrix}
\end{align*}

The four sequences
\[
\begin{matrix}
   x^2+3 x+2 & x^2+3 x+2 &   x^2+3 x+2 &   x^2+3 x+2 & \dots\\
 5 x+7 &  3 x^2+14 x+13 & 6 x^2+23 x+19 & 9 x^2+32 x+25 & \dots \\
6 &  3 x^2+24 x+33 & 15 x^2+75 x+78 & 36 x^2+153 x+141& \dots \\
0 &  x^2+18 x+41& 20 x^2+135 x+181 & 84 x^2+432 x+474& \dots
\end{matrix}
\]
have all real roots and  satisfy the recurrence
\[
p_n - 3 p_{n-1} + 3 p_{n-2} - p_{n-3}=0.
\]

\section{Recursions from  $f(y,z)$ and $g(x,y,z)$}
\label{sec:rati-funct-from-p3}

Given a polynomial in $\rupint{3}$ and $v$ the vector of the first $k$
coefficients of a polynomial in $\rupint{2}$   we can
construct a matrix that preserves such initial sequences. This gives us
recursions of real rooted polynomials. We 
give the general construction, and then give some examples.

\begin{construction}\index{recurrences!from $\rupint{3}$}
  We are given 
  $f = \sum f_{i,j}(x)\,y^iz^j$ in $\gsubclose_3$, and $g(x,y) =
  \sum_0^r g_i(x)y^i$ is the first $r+1$ terms of a function in
  $\gsubf_2$. .

We  define
  \begin{align*}
    M &= \left(f_{r-i,j}\right)_{ 0\le i,j\le r} \\
    v_0 &= (g_0,\dots,g_r)
  \end{align*}
  
By Theorem~\ref{thm:multiply-vectors} we know that all polynomials
\[ M^k v_0 = v_k = (v_{0,k},\dots,v_{r,k}) \]
have all real roots. Suppose that $M$ satisfies
\[
 \sum_0^s \alpha_i M^i = 0\]
Then  the sequences
\begin{align*}
    g_0 & = v_{0,0},v_{0,1}, v_{0,2}, \dots \\
    g_1& = v_{1,0},v_{1,1}, v_{1,2}, \dots \\
   & \dots \\
    g_r&= v_{r,0},v_{r,1}, v_{r,2}, \dots \\
 \end{align*}
 
all satisfy the same recurrence
   \[ \sum_{i=0}^s \alpha_i p_i = 0
\]
\end{construction}

\begin{example}
  We start with 
\[
\begin{vmatrix}
  x+1 & 1 & 1 \\ 1 & y & 1 \\ 1 & 1 & z
\end{vmatrix}=
y z (1+x) -y-z+(1-x) \in\gsubclose_3
\]
The matrix 
\[ M =
\begin{pmatrix}
  -1 & 1-x \\ 1+x & -1
\end{pmatrix}
\]
preserves initial sequences of length two. 
$M$ satisfies its
characteristic polynomial
\[ \lambda^2 + 2\lambda+x^2 = 0
\]
If  $g_0(x) + y \,g_1(x)\in \gsubclose_2$ then
the  polynomial sequences 
\[
\begin{pmatrix}
  -1 & 1-x \\ 1+x & -1
\end{pmatrix}^n
\begin{pmatrix}
  g_0 \\f_0
\end{pmatrix}
=
\begin{pmatrix}
  g_n \\f_n
\end{pmatrix}
\]
are in $\allpoly$ and satisfy 
\begin{align*}
  g_n + 2\, g_{n-1} + x^2 \,g_{n-2} &= 0 \\
  f_n + 2\, f_{n-1} + x^2 \,f_{n-2} &= 0 
\end{align*}

Now we choose $g_0=-1$ and $g_1=x$. Since $g_0+yg_1=xy-1\in\gsubclose_2$ we can
apply the above construction
\begin{align*}
  v & = \begin{pmatrix}-1 \\ x\end{pmatrix} 
&
 Mv & = \begin{pmatrix} -x^2+x+1 \\ -2 x-1  \end{pmatrix} 
& 
 M^2v & = \begin{pmatrix} 3 x^2-2 x-2 \\ -x^3+4 x+2  \end{pmatrix} 
& \dots
\end{align*}
giving us the two sequences of polynomials
\[
\begin{matrix}
  g: &-1&-x^2+x+1 & 3 x^2-2 x-2&x^4-x^3-7 x^2+4 x+4,  & \dots  \\
  f: & x&-2 x-1&-x^3+4 x+2&4 x^3+x^2-8 x-4, & \dots
\end{matrix}
\]
 Consecutive $g_i$'s or $f_i$'s 
  do not necessarily interlace, but they all are in
 $\allpoly$. However, $g_i$ and $f_i$ do interlace.

\end{example}

\begin{example}
  The polynomial $xy + xz + yz -1$ is in $\gsubclose_3$ since 
$
\left(\begin{smallmatrix}
  0 & 1 & 1 \\1 & 0 & 1 \\ 1 & 1 & 0
\end{smallmatrix}
\right)
$ is \nsd. For a particular choice of $v$ we have  nice formulas for
the two polynomial sequences.
\begin{align*}
  M &= \begin{pmatrix} x & -1 \\ 1 & x \end{pmatrix} \\
  \charpoly{M} &= y^2 - 2xy + x^2+1 \\
  p_n &= 2x p_{n-1} - (x^2+1)p_{n-2} \\
  v &= (x,1) \\
  M^kv &= (f_k,g_k)\\
  f_k &= \sum_i x^{n-2i}\binom{n}{2i}(-1)^i \\
  g_k &= \sum_i x^{n-2i+1}\binom{n}{2i-1} (-1)^{i+1}
\end{align*}
We know that $f_k$ and $g_k$ have all real roots since they are
derived from the even and odd parts of $(x-1)^n$

\end{example}

\begin{example}
  In this example we construct a sequence of polynomials in two
  variables that satisfy a three term recursion.  We begin with the
  \index{Grace polynomial}Grace polynomial that is in $\gsubclose_4$
\[ 
\grace{2} = 2x_1x_2 + (x_1+x_2)(y_1+y_2) + 2 y_1y_2
\index{recurrence relations!and rational functions}
\]

Assume that $f = f_0(x_1,x_2) +
f_1(x_1,x_2)y_1\in\gsubclose_3$. Consider the coefficients of $y_1$ in
$f\cdot\grace{2}$

\begin{align*}
  \text{Coefficient of $y_1^0$}:\quad & [2x_1x_2+y_2(x_1+x_2)]f_0 \\
  \text{Coefficient of $y_1^1$}:\quad & [2x_1x_2+y_2(x_1+x_2)]f_1+[2y_2+x_1+x_2]f_0\\
  \text{Coefficient of $y_1^2$}:\quad & [2y_2+x_1+x_2]f_1 \\
\end{align*}

The coefficient of $y_1^1$ is a polynomial in $y_2$, so the
coefficients of $y_2^0$ and $y_2^1$ interlace:

\begin{align*}
  \text{Coefficient of $y_2^0$}:\quad & (x_1+x_2)f_0 + 2x_1x_2 f_1  \\
  \text{Coefficient of $y_2^1$}:\quad & 2 f_0 + (x_1+x_2)f_1 
\end{align*}

It follows that the matrix 
\[ M =
\begin{pmatrix}
  x_1+x_2 & 2x_1x_2 \\ 2 & x_1+x_2
\end{pmatrix}
\]
preserves interlacing pairs of polynomials in two variables. As expected, the
determinant, $(x_1-x_2)^2$, is non-negative.

The sequence of vectors
\[
\begin{pmatrix}  f_0 \\ f_1\end{pmatrix},\,
M
\begin{pmatrix}  f_0 \\ f_1\end{pmatrix},\,
M^2
\begin{pmatrix}  f_0 \\ f_1\end{pmatrix},\cdots
\]
consists of polynomials in $\gsubclose_2$. Write $M^i\begin{pmatrix}
  f_0 \\ f_1\end{pmatrix} = \begin{pmatrix}
  g_i\\h_i\end{pmatrix}$. The characteristic polynomial of $M$ is
\[
 z^2 - 2(x_1+x_2)z + (x_1-x_2)^2
\]
and since $M$ satisfies its characteristic polynomial we have
\[
M^{n+2}\cdot  v - 2(x_1+x_2)M^{n+1}\cdot v + (x_1-x_2)^2 M^n\cdot v =
0
\]
for any vector $v$. 
Therefore, we have the recurrences
\begin{align*}
  g_{n+2} - 2(x_1+x_2)g_{n+1} + (x_1-x_2)^2g_n &= 0 \\
  h_{n+2} - 2(x_1+x_2)h_{n+1} + (x_1-x_2)^2h_n &= 0 \\
\end{align*}

If we want to construct examples then we must begin with 
$g_0\longleftarrow h_0$, and not with $g_0,g_1$. For example, if we
take a degenerate case $g_0=1$, $h_0=0$ then the first few pairs are

\[
\begin{pmatrix} 1 \\ 0 \end{pmatrix},\quad 
\begin{pmatrix} x_1+x_2 \\ 2 \end{pmatrix},\quad 
\begin{pmatrix} x_1^2 + 6 x_1 x_2 + x_2^2 \\ 4x_1 + 4 x_2 \end{pmatrix},\quad 
\begin{pmatrix} x_1^3 + 15 x_1^2 x_2 + 15 x_1 x_2^2 + x_2^3 \\ 6 x_1^2 + 20 x_1 x_2 + 6 x_2^2 \end{pmatrix},\quad 
\]
It's not hard to see that we have an explicit formula 
\[ 
g_n(x_1,x_2) = \sum_{i=0}^n \binom{2n}{2i}x_1^i x_2^{n-i}
\]
and the $g_i$ satisfy the recurrence above with initial conditions
$g_0=1$ and  $g_1 = x_1+x_2$.

\end{example}

\section{Recursions from $f(x,-\diffd_y)$ and $g(x,y)$}
\label{sec:recursions-from-fx}

We know that if $f(x,y)\in\gsubclose_2$ then $f(x,-\diffd_y)$ maps
$\gsubclose_2\longrightarrow\gsubclose_2$. This guarantees that our
recursions will consist of polynomials with all real roots.

We assume
\begin{align*}
  f(x,y) &= \sum_i f_i(x)y^i & \in\gsubf_2 \\
  f(0,0) &\ne 0 \\
  g(x,y) &= \sum_{i=0}^d g_i(x)y^i & \in\rupint{2}
\end{align*}
The action of $f(x,-\diffd_y)$ never decreases $y$-degree, so it
determines a $d+1$ by $d+1$ matrix acting on the coefficients of
polynomials in $\rupint{2}(d)$.

\begin{example}
  We assume
  \begin{align*}
    f(x,y) &= (x+y)^3 & g(x,y) &= g_0(x) + g_1(x)y + g_2(x)y^2
  \end{align*}
The action of $f(x,-\diffd_y)$ on $g$ is given by the matrix
\[
\begin{pmatrix}
  x^3 & -3 x^2 & 6x \\ 0 & x^3 & -6x^2\\  0 & 0 & x^3
\end{pmatrix}
\]
This gives recurrences of the form
\[
p_n - 3 x^3p_{n-1} + 3 x^6 p_{n-2} - x^9 p_{n-3} = 0
\]
\end{example}

\begin{example}
  In this case we let $f = e^{-xy}$, and $g = \sum_0^4 g_i(x)y^i$. 
  The matrix is
\[
\begin{pmatrix}
  1 & -x & x^2 & -x^3 & x^4 \\
  0 & 1 & -2x & 3 x^2 & -4 x^3 \\
  0 & 0 & 1 & -3x & 6x^2 \\
  0 & 0 & 0 & 1 & -4x \\
  0 & 0 & 0 & 0 & 1
\end{pmatrix}
\]
The characteristic polynomial is $(\lambda-1)^5$, so these recurrences are of the form
\[
p_n - 5p_{n-1} + 10  p_{n-2} -  10p_{n-3} + 5 p_{n-4} - p_{n-5} = 0.
\]

\end{example}

\section{Recurrences from $f(-\diffd_x,y)$ and $g(x,y)$}
\label{sec:recurences-from-f}

If we use $f(-\diffd_x,y)$ instead of $f(x,-\diffd_y)$ then we get
recurrences involving derivatives.
We assume
\begin{align*}
  f(x,y) &= \sum_i f_i(x)y^i & \in\gsubf_2 \\
  g(x,y) &= \sum_{i=0}^d g_i(x)y^i & \in\rupint{2}
\end{align*}
We give a very simple example - note that the degree of the recurrence
comes from the number of coefficients we consider, and not from the
degree of $f$.

\begin{example}
  We let 

\begin{align*}
f(x,y) &= y - \diffd_x & g &= \sum g_i(x)y^i.
\end{align*}

We consider only the first four coefficients of $g$, so the action
of $y-\diffd_x$ is
\[
-g_0' + y(g_0-g_1') + y^2( g_1-g_2') + y^3(g_2-g_3').
\]
As a matrix this is
\[
\begin{pmatrix}
  -\diffd_x & 0 & 0 & 0 \\
1 &   -\diffd_x & 0 & 0  \\
0 & 1 &  -\diffd_x & 0  \\
0 & 0 & 1 &  -\diffd_x  \\
\end{pmatrix}
\begin{pmatrix}
  g_0 \\ g_1 \\ g_2 \\ g_3
\end{pmatrix}
\]
The characteristic polynomial is $(\lambda+\diffd_x)^4$, so the
recurrence is
\[
p_n + 4 p_{n-1}^{'} + 6 p_{n-2}^{(2)}+ 4 p_{n-3}^{(3)}+  p_{n-4}^{(4)}=0.
\]
\end{example}

\section{Recurrences from mutually interlacing polynomials}
\label{sec:rati-funct-from-mi}

If a matrix preserves \index{polynomials!mutually interlacing}mutually
interlacing polynomials then we get recursive sequences of polynomials
in $\allpoly$.

\begin{example}
  The matrix and initial vector
  \begin{align*}
 M &=    \begin{pmatrix}
      1 & x+1 & x+1 \\ 1 & 1 & x+1 \\ 1 & 1 & 1
    \end{pmatrix}
& v &= 
\begin{pmatrix}
  x-1\\x\\x+1
\end{pmatrix}
  \end{align*}
determine three polynomial sequences $M^iv = (f_i,g_i,h_i)$ that
satisfy the recurrence
\[
p_n = 3\,p_{n-1} + 3\,x\,p_{n-2} + x^2\,p_{n-3}.
\]
To find the initial terms we need to compute $Mv$ and $M^2v$. 
\begin{align*}
  Mv &= 
\begin{pmatrix}2 x^2+4 x \\x^2+4 x \\3 x \end{pmatrix}
&
M^2v &= 
\begin{pmatrix}x^3+10 x^2+11 x \\ 6 x^2+11 x \\3 x^2+11 x \end{pmatrix}
\end{align*}
It follows that the following sequence of polynomials has all real
roots. It also appears that $f_i\lesslesseq f_{i-1}$ but that doesn't
follow from the general construction, and I don't know how to prove it.
\begin{align*}
  h_0 &= x+1 \\
  h_1 &= 3x \\
  h_2 &= 3x^2+11x\\
  h_n &= 3h_{n-1} + 3x h_{n-2} + x^2 h_{n-3}
\end{align*}
\end{example}

\begin{example}
  In this example we start with the recurrence, and find a matrix that
  preserves mutually interlacing polynomials. Consider the sequence 
  \begin{gather*}
q_0 = q_1 = q_2 = \cdots = q_{d} = 1 \\
q_k = q_{k-1} + x q_{k-d} \text{ for } k> d \\
q_{d+1} = 1+x,\quad q_{d+2}=1+2x,\quad\cdots\quad q_{2d+1}=1+(d+1)x\\
\end{gather*}

Rewrite the recurrence in matrix terms:
$$
\begin{pmatrix}
  1 & 0 & 0 & 0 & \hdots & x \\
  1 & 0 & 0 & 0 & \hdots & 0 \\
  0 & 1 & 0 & 0 & \hdots & 0 \\
  \vdots &  & \ddots & \ddots &   & \vdots \\
  0 & 0 & \hdots &  &  1 & 0
\end{pmatrix}
\begin{pmatrix}
  q_{k-1} \\q_{k-2} \\ \vdots \\ \\ q_{k-d}
\end{pmatrix}
=\begin{pmatrix}
  q_{k} \\q_{k-1} \\ \vdots \\  \\ q_{k-d+1}
\end{pmatrix}
$$

The $d+1$ by $d+1$ matrix preserves mutually interlacing polynomials.  Since the
vector $(q_{d+1},\dots,q_{2d+1})$ is mutually interlacing, we can
repeatedly apply the Lemma to conclude that $(q_{k},\dots,q_{k+d})$ is
a vector of mutually interlacing vectors. In particular, all $q_n$
have all real roots.

\end{example}

If a two by two matrix preserves interlacing polynomials then we get a
recurrence. This is really just a special case of a general matrix
preserving the first $k$ coefficients of a polynomial, and also a
special case of a matrix preserving mutual interlacing.

\begin{example}
  The matrix $
  \begin{pmatrix}
    1 & x \\ 1 &0
  \end{pmatrix}$ preserves interlacing. If we let
\[
  \begin{pmatrix}
    1 & x \\ 1 &0
  \end{pmatrix}^n
  \begin{pmatrix}
    1 \\1
  \end{pmatrix}=
  \begin{pmatrix}
    p_n \\p_{n-1}
  \end{pmatrix}
\]
then all $p_n\in\allpoly$. The characteristic polynomial is
$\lambda^2-\lambda-x$, so we have the recurrence
\[ p_n = p_{n-1} + x p_{n-2}
\]

\end{example}

\section{Alternative approaches}

There are many specialized ways to prove facts about
recurrences. \cite{ryavec-redmond} has some very interesting results
about recurrences with polynomial coefficients. 

Of course, there is a 1-1 correspondence between recurrences with
polynomial coefficients and rational functions. Here we start with a
rational function.

If we  write
\[
\frac{1}{1-(x+y)^n} = \sum_{j=0}^\infty
\frac{h_{n,j}(x)}{(1-x^n)^{j+1}} y^j
\]
then we will show that the roots of the coefficients $h_{n,j}$ lie on $n$
equally spaced rays through the origin.  
For example, if $n=2$ we have 

\[
\frac{1}{1-(x+y)^2} = \frac{1}{1-x^2} +
 \frac{2x}{\left(1-x^2\right)^2}y +
\frac{3x^2+1}{\left(1-x^2\right)^3}y^2 + 
\frac{4\left(x^3+x\right)}{\left(1-x^2\right)^4} y^3 + \cdots
\]
and all the roots lie on the imaginary axis.
There is a simple expression for these
coefficients
\[
  \frac{1}{1-(x+y)^n} = \sum_i (x+y)^{ni} 
  = \sum_j y^j \cdot \sum_{i=0}^\infty \binom{ni}{j}x^{ni-j}.
\]
 The coefficient of $y^j$ is $x^{-j}g(x^n)$ where
\[ 
g_{n,j}(x) = \sum_{i=0}^\infty \binom{ni}{j}x^i = \frac{f_{n,j}(x)}{(1-x)^{j+1}}.
\]

The series converges for $|x|<1$. If $x$ is positive then the series
is positive or divergent. Thus, $g_{n,j}$ has no positive roots.
Here are some values of the sum

\[
\begin{array}{r|ccc}
  &n=1 & n=2 & n=3 \\[.2cm]\hline
j=0 & \frac{1}{1-x} & \frac{1}{1-x} & \frac{1}{1-x} \\[.2cm]
j=1 &  \frac{x}{(1-x)^2} &  \frac{2x}{(1-x)^2} &  \frac{3x}{(1-x)^2}
\\
j=2 &  \frac{x^2}{(1-x)^3}  &  \frac{3x^2+x}{(1-x)^3}  &  \frac{6x^2+3x}{(1-x)^3} 
\end{array}
\]

We will show that $f_{n,j}\in\allpoly$, which establishes the claim. 
We write $g_{nj}(x)$ in terms of a differential operator. First, let
$\binom{nx}{j} = \sum a_k x^k$.

\begin{multline*}
g_{n,j}(x) 
= \sum_i \binom{ni}{j}x^i 
= \sum_{i,k} a_k i^k x^i  \\
= \sum_{i,j} a_k  (x\diffd)^k x^i 
= \bigl(\sum_i (x\diffd)^k a_k\bigr) \frac{1}{1-x}
=\binom{nx\diffd}{j}\frac{1}{1-x}.
\end{multline*}

It simplifies the argument if we  eliminate the denominator, so
we will show that $F_j\in\allpoly$ where
\[
(nx\diffd)(nx\diffd-1)\cdots(nx\diffd-j+1)\frac{1}{1-x} =
\frac{F_j}{(1-x)^{j+1}}
\]

This yields the recurrence

\[
F_{n,j+1} = \bigl[\bigl( (j+1)n +j\bigr)x + j\bigr] F_{n,j} -
x(x-1)F_{n,j}'
\] 

This recurrence is of the form
\[
p_n = (ax+b)p_{n-1} - x(x-1)p_{n-1}'
\]
Since $x(x-1)$ is positive at the roots of $p_n$ it follows that all
$f_j$ have all real roots and all positive coefficients.


\chapter{Matrices}
\label{cha:matrices}
 
\renewcommand{\TimeStampStart}{Tuesday, August 01, 2006: 21:20:32}
\mytoday

The theme of this chapter is that there are several ways to represent
sequences of interlacing polynomials by Hermitian matrices.  There is
a 1-1 correspondence between Hermitian matrices and interlacing
sequences of monic polynomials.  Since orthogonal polynomials can be
realized as the characteristic polynomials of certain Hermitian
matrices, many results about orthogonal polynomials can be  generalized
to Hermitian matrices.

\section{Basics}
\label{sec:matrix-basics}

We first recall some basic facts about matrices.  $A^\ast$ denotes
the conjugate transpose of $A$. The matrix $A$ is \emph{Hermitian} if
$A=A^\ast$.  Now we are only considering real polynomials, so we could
restrict ourselves to symmetric matrices, but it is no harder to
consider Hermitian matrices.  A matrix $U$ is unitary if $UU^\ast=I$,
where $I$ is the identity matrix. A Hermitian matrix has all real
eigenvalues.  The spectral theorem says that there is a unitary matrix
$U$ and a diagonal matrix $\Lambda$ consisting of the eigenvalues of $A$
such that $A = U\Lambda U^\ast$.  If all entries of $A$ are real, then
$U$ can be chosen to be real, so that it is an orthogonal matrix. 

\index{characteristic polynomial}
\index{\ cp@\charpoly{}}
\index{\ zzlesseq@$\lesseq$!defined for matrices}

The characteristic polynomial of $A$ is $\vert xI - A|$, and is
written $\charpoly{A}$.  The set of roots of the characteristic
polynomial of $A$ is written $\lambda(A)$. We say that $A\greateq B$
if $\charpoly{A} \greateq \charpoly{B}$, $A\lessless B$ iff
$\charpoly{A}\lessless\charpoly{B}$, and so on.  
\index{positive definite} A Hermitian matrix is \emph{positive definite} if all of
its eigenvalues are positive, and is \emph{positive semidefinite} if
the eigenvalues are all non-negative.

If the matrix $A$ is invertible then interlacing properties of $A$
determine interlacing properties of $A^{-1}$.  

\begin{itemize}
\item If $A$ is invertible and $\lambda(A) = (a_1,\dots,a_n)$, then
$\lambda(A^{-1}) = (\dfrac{1}{a_1}.\dots,\dfrac{1}{a_n})$.
\item If $A$ is positive definite and $A \greateq B$ then $B^{-1} \greateq A^{-1}$.
\item If $A$ is positive definite and $A \lessless B$ then $A^{-1} \lessless B^{-1}$.

\index{reverse!of a polynomial}

\item  The characteristic polynomial of $A^{-1}$ is 
  $$
  \vert x I - A^{-1} \vert = \vert - A^{-1}\vert \,\vert I - xA
  \vert = \vert - A^{-1}\vert \, \vert (1/x)I - A \vert x^n $$
  where
  $A$ is $n$ by $n$. Since the reverse of a polynomial $f$ is $\rev{f}
  = x^n f(1/x)$, we can write
$$ \charpoly{A^{-1}} = \frac{1}{\vert -A\vert} \rev{\left(\charpoly{A}\right)}$$
\end{itemize}

\section{Principle submatrices}
\label{sec:matrix-submatrix}

In this section we establish  the important result that Hermitian matrices
correspond to sequences of interlacing polynomials.

\index{principle submatrix}

\begin{definition}
  For $n$ by $n$ matrix $M$ , and index sets
  $\alpha,\beta\subset\{1,\dots,n\}$, the (sub)matrix that lies in the
  rows indexed by $\alpha$ and columns indexed by $\beta$ is denoted
  $M\{\alpha,\beta\}$. If $\alpha=\beta$ then we write $M\{\alpha\}$; such
  a submatrix is called a \emph{principle} submatrix
  (\chapsec{pd}{sub-quadratic-forms}).  Sometimes we want to refer to
  a submatrix via the deletion of rows; in this case we let
  $M[\alpha]$ be the submatrix resulting from the deletion of
  the rows and columns not listed in $\alpha$.  The size
  $\vert\alpha\vert$ of an index set $\alpha$ is its cardinality.
\end{definition}

We know that if $A$ is a principle submatrix of $B$ then $B\lesslesseq
A$.

Suppose that $f$ is a polynomial with roots
$\roots(f)=(a_1,\dots,a_n)$.  Set $\Lambda$ equal to the diagonal
matrix whose elements are $(a_1,\dots,a_n)$, and let $y^\ast =
(y_1,\dots,y_n)$. The following calculation shows the relation between
the characteristic polynomial of $\Lambda$ and an extension of
$\Lambda$. Suppose that

\begin{align*}
A & = 
\begin{pmatrix}
\Lambda & \vdots & y \\
\hdotsfor{3}\\
y^\ast & \vdots & a
\end{pmatrix} \\
\intertext{The characteristic polynomial of $A$ is}
det( xI-A) & =
det 
\begin{pmatrix}
xI-\Lambda & \vdots & -y \\
\hdotsfor{3}\\
-y^\ast & \vdots & x-a
\end{pmatrix} \\
&=
det
\begin{pmatrix}
  I & \vdots & 0 \\
\hdotsfor{3}\\
\left[(xI-\Lambda)^{-1}y\right]^\ast & \vdots & 1
\end{pmatrix}
\begin{pmatrix}
  xI-\Lambda&\vdots & -y \\
\hdotsfor{3}\\
-y^\ast & \vdots & x-a
\end{pmatrix} \times \\
& \quad\quad
\begin{pmatrix}
  I & \vdots & (xI-\Lambda)^{-1}y \\
\hdotsfor{3} \\
0 & \vdots & I
\end{pmatrix} \\
&=
det
\begin{pmatrix}
  xI-\Lambda & \vdots & 0 \\
\hdotsfor{3}\\
0 & \vdots & (x-a) - y^\ast(xI-\Lambda)^{-1}y
\end{pmatrix}\\
&=
\left[(x-a)-y^\ast(xI-\Lambda)^{-1}y\right]\,det(xI-\Lambda)\\
&=
(x-a)f(x) - \sum_{i=1}^n y^2_i \dfrac{f(x)}{x-a_i}
\end{align*}

If $p_i$ is the characteristic polynomial of $A\{1,\dots,i\}$, then
from repeated applications of Theorem~\ref{thm:principle-1} we have a sequence of
interlacing polynomials
$$
p_n \lesslesseq p_{i-1} \lesslesseq \dots \lesslesseq p_0$$
Conversely, we can use the above calculation to show that given the
$p_i$'s, we can find a Hermitian matrix determining them.

\begin{theorem} \label{thm:principle-2}
  Suppose that $p_0,p_1,\dots,p_n$ is a sequence of polynomials such
  that the degree of $p_i$ is $i$, and $p_i \lesslesseq p_{i-1}$ for
  $1\le i \le n$.  There is a Hermitian matrix $A$ such that the
  characteristic polynomial of $A\{1,\dots,i\}$ is a 
  multiple of $p_i$, for $0\le i \le n$.
\end{theorem}

\begin{proof}
  We prove this by induction on $n$.  The case $n=1$ is clear, so
  assume that $A$ is an $n$ by $n$ Hermitian matrix such that the
  characteristic polynomial of $A\{1,\dots,i\}$ is a multiple of
  $p_i$, for $0\le i \le n$.  Choose an orthogonal matrix $U$ so that
  $UAU^\ast$ is the diagonal matrix $\Lambda$ whose diagonal consists
  of the roots of $p_n$. We need to find a number $b$ and vector
  $z^\ast = (z_1,\dots,z_n)$ so that the characteristic polynomial of
\begin{equation} \label{eqn:mat-1}
\begin{pmatrix}
  A&\vdots & z \\
\hdotsfor{3}\\
z^\ast & \vdots & b
\end{pmatrix}
\end{equation}
  
is equal to a multiple of $p_{n+1}$.  Write

\begin{equation} \label{eqn:mat-2}
B = 
\partmatrix{U}{0}{0}{1}
\partmatrix{\Lambda}{z}{z^\ast}{b}
\partmatrix{U^\ast}{0}{0}{1}
=
\partmatrix{\Lambda}{y}{y^\ast}{a}
\end{equation} 
  
The characteristic polynomial of $A$ is $cp_n$ for some constant $c$.
The characteristic polynomial of \eqref{eqn:mat-1} is the same as the
characteristic polynomial $B$ in \eqref{eqn:mat-2}.  Since $p_{n+1}
\lesslesseq p_n$, we can apply Lemma~\ref{lem:sign-quant} and find an $a$ and
$y_i^2$ such that $$p_{n+1} = (x-a)p_n - \sum y_i^2 \frac{p_n(x)}{x-a_i}$$
where $\roots(p_n) = (a_1,\dots,a_n)$. The desired matrix is $U^\ast BU$.
\end{proof}

  If $p_{i+1} \lessless p_i$ for all $i$ then all of the coefficients
  $y_i^2$ are non-zero.  In this case there are $2^n$ choices of
  $y_i$'s, and so there are $2^n$ different $A$'s.

\begin{example}
  The construction of a matrix from a sequence of characteristic
  polynomials involves solving many equations.  Here are a few
  examples of $n$ by $n$ matrices $M$ where we are given a polynomial $f$
  and the determinant of the first $i$ rows and columns is a constant
  multiple of $f^{(i)}$.
  \begin{enumerate}
  \item If $f=(x-a)(x-b)$ then $M = 
    \begin{pmatrix}
      \frac{a+b}{2} & \frac{a-b}{2} \\
      \frac{a-b}{2} & \frac{a+b}{2} 
    \end{pmatrix}
$
\item If $f = (x-1)(x-2)(x-3)$  then  $M = 
  \begin{pmatrix}
    2 & \frac{-1}{\sqrt{3}} & 0 \\
\frac{-1}{\sqrt{3}} & 2 & \frac{-2}{\sqrt{3}} \\
0 & \frac{-2}{\sqrt{3}} & 2
  \end{pmatrix}
$
\item If $f = x^2(x-1)$ then $M = 
  \frac{1}{3} \begin{pmatrix}
    1 & -1 & -1 \\ -1 & 1 & 1 \\ -1 & 1 & 1
  \end{pmatrix}
$
  \end{enumerate}
\end{example}

\begin{example}
  If a Hermitian matrix $H$ is partitioned into diagonal blocks, and the
  polynomials determined by each of these blocks strictly interlace,
  we do not have strict interlacing for the polynomials determined by
  $H$.  Assume that $A$ is an $r$ by $r$ Hermitian matrix determining
  the polynomial sequence $p_r \lessless \cdots \lessless p_1$, and 
  $B$ is an $s$ by $s$ Hermitian matrix determining
  the polynomial sequence $q_s \lessless \cdots \lessless q_1$. The
  partitioned matrix
$$
\partmatrix{A}{0}{0}{B}
$$
determines the polynomial sequence
$$ p_1 \lessgreat p_2 \lessgreat \cdots p_r 
\lessgreateq p_r\,q_1 \lessgreateq p_r\,q_2 \lessgreateq \cdots
\lessgreateq p_r\,q_s
$$
\end{example}

If a Hermitian matrix $M$ is diagonal with diagonal entries
$(a_1,\dots,a_n)$, then the characteristic polynomial is $f(x) =
(x-a_1)\cdots(x-a_n)$. The characteristic polynomial of a single
row/column deletion $M[i]$ equals $f(x)/(x-a_i)$.  Since
these are the polynomials occurring in quantitative sign interlacing
it is not surprising that the characteristic polynomials of single
deletions are also important. For instance, an immediate consequence
of Theorem~\ref{thm:principle-1} is that for any choice of positive $b_i$  the
polynomials
$$ \sum_{i=1}^n b_i\, \charpoly{M[i]}$$
have all real roots.

Analogous to the formula 
$$ \frac{d}{dx} f(x) = \sum_{i=1}^n \frac{f(x)}{x-a_i}$$
we have

$$ \frac{d}{dx} \charpoly{M} = 
\sum_{i=1}^n \charpoly{M[i]}$$

This can be established by considering the derivative of only the
terms in the expansion of the determinant that have diagonal entries
in them.

\begin{example}
  It is not the case that the characteristic polynomials of the single
  deletion submatrices of an $n$ by $n$ matrix span an $n$-dimensional
  space.  The matrix $M=\smalltwo{0}{1}{1}{0}$ has eigenvalues $1,-1$
  but $M[1] = M[2] = (0)$, so they have the same characteristic
  polynomial, $x$.  More generally, for any symmetric $C$ with
  distinct eigenvalues none of which are $1$ or $-1$, the matrix $M$
  below has $M[1] = M[2]$, so the characteristic polynomials do not
  span an $n$-dimensional space.
$$ 
\begin{pmatrix}
  0 & 1 & 0 \\ 1 & 0 & 0 \\ 0 & 0 & C
\end{pmatrix}
$$
\end{example}

An $n$ by $n$ matrix $M$ has $2^n$ principle submatrices. Their
characteristic polynomials can be  combined to form  polynomials that
are in $\gsubclosepos_{n+1}$.

\begin{lemma}
  Let $M$ be an $n$ by $n$ symmetric matrix. For each
  $\alpha\subset\{1,\dots,n\}$ let $\xx^\alpha = \prod_{i\in\alpha}
  x_i$. Then the polynomials below are in $\gsubclosepos_{n+1}$.
  \begin{align*}
    & \sum_{\alpha\subset\{1,\dots,n\}}
    \charpoly{M\{\alpha\}}\,\xx^\alpha \\
    & \sum_{\alpha\subset\{1,\dots,n\}}
    det\,{M\{\alpha\}}\,\xx^\alpha \\
  \end{align*}
\end{lemma}
\begin{proof}
  The polynomial $f(x,\xx)$ in question is nothing other
  than the characteristic polynomial of $M+D$ where $D$ is the
  diagonal matrix whose diagonal is $(x_1,\dots,x_n)$. It remains to
  see that $f\in\gsubclosepos_{n+1}$. Consider $g_\epsilon =
  \charpoly{M+D+\epsilon zI}$ where $z=x_1+\cdots+x_n$.  The
  homogeneous part of $g_\epsilon$ is $\prod_i
  (x+x_i+\epsilon(x_1+\cdots+x_n)$  and for $\epsilon>0$ this is in
  $\gsubpos_{n+1}$.  If we substitute values for $x_1,\dots,x_n$ then
  the resulting polynomial is the characteristic polynomial of a
  symmetric matrix and so is in $\allpoly$.  Since
  $\lim_{\epsilon\rightarrow0^+} g_\epsilon = f$ it follows that $f$
  is in the closure of $\gsubpos_{n+1}$.

  The second assertion follows from the first by setting $x=0$.

\end{proof}

For example, if $M= \smalltwo{a}{b}{b}{c}$ then the two polynomials
are
\begin{align*}
  det\smalltwodet{a-x-x_1}{b}{b}{c-x-x_2}
  &   = (x-a)(x-c)-b^2 + (x-a)x_2 + (x-c)x_1 + x_1x_2 \\ 
  det\smalltwodet{a-x_1}{b}{b}{c-x_2}& =(ac-b^2) -a x_2 -cx_1 + x_1x_2
\end{align*}

\section{The Schur complement}
\label{sec:matrix-schur}

The Schur complement of a matrix plays a role that is similar to
contraction for a graph.  Although the initial definition is
unintuitive, there is a simple interpretation in terms of deletion in
the inverse that is easy to understand.

\index{Schur complement}

\begin{definition}
  For an $n$ by $n$ matrix $M$ , and index sets
  $\alpha,\beta\subset\{1,\dots,n\}$ where $a$ and $\beta$ have the
  same size,  the Schur complement of the (sub)matrix that lies in the
  rows indexed by $\alpha$ and columns indexed by $\beta$ is
$$ A/[\alpha,\beta] = A\{\alpha^\prime,\beta^\prime\} - 
A\{\alpha^\prime,\beta\}  A\{\alpha,\beta\}^{-1}  A\{\alpha,\beta^\prime\}
$$
 If $\alpha=\beta$ then we write $M/[\alpha]$. Notice that if $\alpha$
 has size $k$, then the Schur complement has dimension $n-k$. 
\end{definition}

\begin{example}
  If $M$ is diagonal with diagonal entries $(a_1,\dots,a_n)$, and
  $\alpha\subset \{1,\dots,n\}$, then the Schur complement
  $M/[\alpha]$ is the diagonal matrix $M\{\alpha^\prime\}$.
\end{example}

The Schur complement may be computed by taking the inverse, deleting
appropriate rows and columns, and taking the inverse of the resulting
matrix.

\begin{theorem}[\cite{ando}] \label{thm:ando-1}
  If $A$ is invertible, then
$ A/[\alpha,\beta] =
 \left( A^{-1}\{\beta^\prime,\alpha^\prime\} \right)^{-1}$
\end{theorem}

An interesting property is  that two disjoint deletions commute.

\begin{theorem}[\cite{ando}] \label{thm:ando-2}
If $\alpha_1,\alpha_2,\beta_1,\beta_2$ are index sets with
$\vert\alpha_1\vert = \vert\beta_1\vert, \vert\alpha_2\vert =
\vert\beta_2\vert$, and
$\alpha_1\cap\alpha_2=\beta_1\cap\beta_2=\emptyset$ then
$$ \left(M/[\alpha_1,\beta_1]\right)/[\alpha_2,\beta_2] = 
\left(M/[\alpha_2,\beta_2]\right)/[\alpha_1,\beta_1] = 
M/[\alpha_1\cup\alpha_2,\beta_1\cup\beta_2]
$$
\end{theorem}

We can now show that Schur complements preserve interlacing.
\begin{theorem}[\cite{smith}] \label{thm:schur-1}
  If $A$ is a   definite (positive or negative) Hermitian  matrix, and
  if $|\alpha| = |\beta|=1$, then
$ A \lesslesseq A/[\alpha,\beta]$.
\end{theorem}

\begin{proof}
  Since $A^{-1}$ is Hermitian and its principle submatrices interlace,
  we get $A^{-1} \lesslesseq A^{-1}\{\beta^\prime,\alpha^\prime\}$.  Since
  all eigenvalues of $A$ have the same sign, when we take inverses
  interlacing is preserved:
    $$ A   = \left(A^{-1}\right)^{-1} \lesslesseq
    \left(A^{-1}\{\beta^\prime,\alpha^\prime\}\right)^{-1} 
= A/[\alpha,\beta]$$
\end{proof}

\begin{cor} \label{cor:schur-1}
If $M$ is non-singular definite Hermitian, and we set $p_k =
\charpoly{M/[\{1..k\}]}$ then
$$ p_n \lesslesseq p_{k-1} \lesslesseq \dots \lesslesseq  p_1$$
\end{cor}

\begin{proof}
  Since $(M/[\{1..k\}])/[\{k+1\}] = M/[\{1..k+1\}]$ by Theorem~\ref{thm:ando-2}, we
  can apply Theorem~\ref{thm:schur-1}.
\end{proof}

Another simple consequence is that we can realize polynomial sequences
by characteristic polynomials of Schur complements.

\index{reverse!of a polynomial}
\begin{theorem}
  Suppose that $p_0,p_1,\dots,p_n$ is a sequence of polynomials such
  that the degree of $p_i$ is $i$, and $p_i \lesslesseq p_{i-1}$ for
  $1\le i \le n$.  There is a Hermitian matrix $A$ such that the
  characteristic polynomial of $A/[\{1,\dots,i\}]$ is a constant
  multiple of $p_i$, for $0\le i \le n$.
\end{theorem}

\begin{proof}
Let $B$ be a Hermitian matrix such that $\charpoly{B\{\{1..k\}\}} =
\rev{(p_i)}$. $B^{-1}$ is the desired matrix.

\end{proof}

\section{Families of matrices}
\label{sec:matrix-family}

\index{locally interlacing family!of matrices}
\index{Feynman-Hellmann theorem}
\index{family!of interlacing polynomials}
\index{family!of locally interlacing polynomials}
\index{family!of  polynomials}
\index{locally interlacing family}

We consider properties of the eigenvalues of continuous families of
Hermitian matrices. We begin with a definition of \emph{a locally
interlacing family} of matrices.

\begin{definition}
If $\{H_t\}$ is a continuous family of Hermitian matrices such that
the family of their characteristic polynomials is a locally interlacing
family, then we say that $\{H_t\}$ is a \emph{locally interlacing family}.
\end{definition}

If the entries of a Hermitian matrix are  functions of a parameter,
then we can determine the derivatives of the eigenvalues with respect
to this parameter.  This result is known as the Hellmann-Feynman
theorem:
\begin{theorem}
  Let $H(x)$ be an $n$ by $n$ Hermitian matrix whose entries have
  continuous first derivatives with respect to a parameter $x$ for
  $x\in(a,b)$. Let $\lambda(H(x)) = (\lambda_1(x),\dots,\lambda_n(x))$
  be the eigenvalues of $H(x)$, with corresponding eigenvectors
  $u_1,\dots,u_n$. If $H^\prime$ is the matrix
  formed by the derivatives of the entries of $H(x)$, then for
  $j=1,\dots,n$ 
$$ \frac{d\lambda_i}{dx} = \frac{u_i^\ast\, H^\prime\,u_i}{u_i^\ast u_i}$$
\end{theorem}

\begin{proof}
  More generally, we consider an inner product for which $(H_t u,v) =
  (u,H_t v)$, $(u,v) = (v,u)$, and $(v,v)>0$ if $v\ne0$.  From the
  fact that the $u_i$ are eigenvectors, we find that for any $s\ne t$
  \begin{align*}
    (H_t u_i(t),u_i(s)) & = \lambda_i(t)\, (u_i(t),u_i(s)) \\
    (H_s u_i(t),u_i(s)) & =     (H_s u_i(s),u_i(t)) \\
                        & = \lambda_i(s)\, (u_i(t),u_i(s)) \\
\intertext{Subtracting and dividing by $t-s$ gives the difference
    quotient}
  \left( \frac{H_t - H_s}{t-s} u_i(t),u_i(s)\right) & = 
\frac{\lambda_t - \lambda_s}{t-s} (u_i(t),u_i(s))
  \end{align*}
The result follows by taking the limit as $s$ goes to $t$.
\end{proof}

\begin{cor}  \label{cor:matrix-pd}
  If $H(x)$  is an $n$ by $n$ Hermitian matrix whose entries have
  continuous first derivatives with respect to a parameter $x$ for
  $x\in(a,b)$, and $H^\prime(x)$ is positive definite for all $x$,
  then $\{H_t\}$ is a locally
  interlacing family.
\end{cor}

\begin{example}
  Suppose that $H$ is Hermitian. Consider the family
  $\{e^{tH}\}$. The derivative of this family is
  $$
  \frac{d}{dt} \left(e^{tH}\right) = He^{tH}$$
  If $v$ is an
  eigenvector of $H$ with eigenvalue $a$, then $He^{tH}v = H(e^{at}v)
  = ae^{at}v$ and so $v$ is an eigenvector of $H^\prime$ with
  eigenvalue $ae^{at}$.  These eigenvalues are increasing if and only
  if $a$ is positive, so we conclude that $\{e^{tH}\}$ is locally
  interlacing if and only if $H$ is positive definite.
\end{example}

\begin{cor} \label{cor:matrix-tp}
  If $\{H_t\}$ is  a family of Hermitian matrices such that
  $\frac{d}{dt}H_t$ is strictly totally positive, then the families of 
  characteristic polynomials of any principle submatrix are locally
  interlacing. 
\end{cor}

\begin{proof}
  Since every principle submatrix of $H^\prime$ is positive definite,
  the result follows from Corollary~\ref{cor:matrix-pd}.
\end{proof}

The next result is known as the \emph{monotonicity theorem}, since it
states that adding a positive definite matrix increases all the
eigenvalues.
\index{monotonicity theorem}

\begin{cor}
  If $A,B$ are Hermitian matrices, and $B$ is positive definite, then
  the family $\{A+tB\}$ is locally interlacing.  
\end{cor}

\begin{proof}
  The derivative $\frac{d}{dt}(A+tB)$ is simply $B$.
\end{proof}

\begin{example}
  Consider the family $\{H+tI\}$, where $H$ is Hermitian. Since $I$ is
  positive definite, it is an locally interlacing family. If the
  characteristic polynomial of $H$ is $f(x)$, then 
$$ \vert| xI - (H+tI)\vert = \vert (x-t)I - H\vert = f(x-t)$$

  Another elementary family is given by $\{tH\}$.  Since
$\vert xI - tH\vert = f(x/t)$ it is locally interlacing.
\end{example}

We can relate interlacing of polynomial sequences with their Hermitian
matrices, but first we need to recall the \emph{positive definite
ordering} (p. \pageref{lem:positive-def-2}) of Hermitian matrices.

\index{positive definite!ordering of matrices}

\begin{lemma}
  Suppose that $A,B$ are Hermitian.  They determine polynomial
  sequences $p_i = \vert|xI - A\{1..i\}\vert$ and 
  $q_i = \vert|xI - B\{1..i\}\vert$.  If $B \prec A$ then $p_i
  \prec q_i$ for $1\le i \le n$.
\end{lemma}

\begin{proof}
  Consider the family $M_t=tA + (1-t)B$.  The derivative is $A-B$ which is
  positive definite, so it is a locally interlacing family.  Since
  $M_0=B$ and $M_1=A$ it follows that  the roots of $M_t$ are
  increasing as a function of $t$, so $p_i \prec q_i$.

\end{proof}

\section{Permanents}
\label{sec:permanents}

\index{permanent}
If $A=(a_{ij})$ is a matrix then the \emph{permanent} of $A$, written
$per(A)$, is the determinant without alternating signs:
$$ per(A) = \sum_{\sigma\in S_n}\ \prod_{i=1}^n\ a_{i,\sigma(i)}$$
Unlike the determinant, polynomials such as $per(x I + C)$ where $C$
is symmetric almost never have all real roots. If we define $J$ to be
the all $1$ matrix, then there is a beautiful conjecture
\cite{haglund97}

\begin{conj}
  If $A=(a_{ij})$ is a real $n$ by $n$ matrix with non-negative
  entries which are weakly increasing down columns then the permanent of
  $A+xJ$ has all real roots.
\end{conj}

The permanents of principle submatrices of $A+xJ$ do not interlace the
permanents of $A+xJ$, but we do have common interlacing:

\begin{lemma}
  Suppose that $A$ satisfies the conditions above, and the conjecture
  is true. Then $per(A+xJ)$ and $per( (A+xJ)[1])$ have a common interlacing.
\end{lemma}
\begin{proof}
  If we let $M=A+xJ$ so that $M = M_0$ where
  \begin{align*}
    M_\alpha &= 
    \begin{vmatrix}
      a_{11}+x + \alpha & v+x \\ w+x & A[1] + x J[1]
    \end{vmatrix} \\
\intertext{By linearity of the permanent we have}
per(M_\alpha) &=
per
\begin{vmatrix}
  a_{11}+x & v+x \\w+x & A[1]+xJ[1]
\end{vmatrix} + 
\alpha\ per(A[1] + x J[1])
  \end{align*}
If we choose $\alpha$ to be negative then the conjecture implies that
$per(M_\alpha)$ has all real roots. Consequently, $per(M) + \alpha\ 
per(M[1])$ has all real roots for all negative $\alpha$. This implies
that $per(M)$ and $per(M[1])$ have a common interlacing. 
\end{proof}

\begin{remark}
  The proof can be easily modified to show that if we remove the first
  row and any column from $A$ then the permanents have a common
  interlacing. Empirically, it appears that all the submatrices have a
  common interlacing. 
\end{remark}

\section{Matrices from polynomials}
\label{sec:pd-mat-from-poly}

If $A$ is a constant matrix and $f(x)$ is a polynomial then $f(xA)$ is
a matrix whose entries are polynomials in $x$. We look for conditions
on $A$ or $f$ that imply that all entries of $f(xA)$ are in
$\allpoly$.

A simple case is when $A$ is a diagonal matrix $diag(a_1,\dots,a_n)$.
The resulting matrix $f(xA)$ is a diagonal matrix. To check this,
write $f(x)=\sum b_ix^i$:

\begin{align*}
  f(xA) &= \sum b_i x^i diag(a_1^i,\dots,a_n^i)) \\
&= diag\left(\sum b_i(xa_1)^i,\dots,\sum b_i(xa_n)^i\right) \\
&= diag(f(a_1x),\dots,f(a_nx))
\end{align*}

If $A$ is diagonalizable then we can write $A=S\Lambda S^{-1}$ where
$\Lambda = diag(d_1,\dots,d_n)$.  The entries of $f(xA)$ are linear
combinations of $f(d_1x),\dots,f(d_nx)$ since

\begin{align*}
  f(xA) &= \sum b_i x^i S \Lambda^i S^{-1} \\
&= S\left( \sum b_i x^i \Lambda^i\right) S^{-1} \\
&= S\,diag(f(d_1x),\dots,f(d_nx))\,S^{-1}
\end{align*}

For example, if $A=\smalltwo{1}{1}{0}{a}$ then $A$ is diagonalizable,
and $A^n = \smalltwo{1}{\frac{a^n-1}{a-1}}{0}{a^n}$ so that
$$ f(xA) = 
\begin{pmatrix}
  f(x) & \frac{1}{a-1}(f(ax)-f(x)) \\ 0 & f(ax)
\end{pmatrix}
$$

It is easy to find polynomials $f$ for which the entries of the above
matrix are not in $\allpoly$. 

A different kind of example is given by $A=\smalltwo{1}{0}{1}{1}$
which is not diagonalizable, yet $f(xA)$ has an especially simple
form. Since $A^n = \smalltwo{1}{0}{n}{1}$ it follows that
$$ f(xA) = \sum b_ix^i\smalltwo{1}{0}{i}{1} = \smalltwo{f}{0}{xf^\prime}{f}$$

This example can be generalized to the matrices whose entries above
the diagonal are $0$, and the rest are $1$. For instance, when $n=4$

$$
\begin{pmatrix}
  1&0&0&0 \\ 1 &1 & 0 & 0 \\ 1 & 1 & 1 & 0 \\ 1 & 1 & 1 & 1
\end{pmatrix}^n =
\begin{pmatrix}
  1 & 0 & 0 & 0 \\
\binom{n}{1} & 1 & 0 & 0 \\
\binom{n+1}{2} & \binom{n}{1} & 1 & 0 \\
\binom{n+2}{3} & \binom{n+1}{2} & \binom{n}{1} & 1
\end{pmatrix}
$$

It follows easily that

$$
f(xA) = 
\begin{pmatrix}
  f & 0 & 0 & 0 \\
 x\diffd(f) & f & 0 & 0 \\
x \diffd^2(xf) & x\diffd(f) & f & 0 \\
x \diffd^3(x^2f) & x\diffd^2(xf) & x\diffd(f) & f
\end{pmatrix}
$$

Thus for these matrices $f(xA)$ has all its entries in $\allpoly$ for
any choice of $f$.  It is easy to characterize those polynomials for
which $f(xA)$ is a matrix whose entries are all in $\allpoly$ for
every choice of $A$.

\begin{lemma} \label{lem:f(xA)}
  If $f$ is a polynomial then the following are equivalent:
  \begin{enumerate}
  \item $f(xA)$ has all entries in $\allpoly$ for every matrix $A$.
  \item $f(x) = cx^n(ax+b)$ for constants $a,b,c$.
  \end{enumerate}
\end{lemma}
\begin{proof}
  If $(2)$ holds then 
  $$f(xA) = cx^n(axA^{n+1} + bA^n).$$
Consequently, all entries of $f(xA)$ are of the form $x^n$ times a
linear term, and so all entries are in $\allpoly$.

Conversely, assume $f(xA)$ has all entries in $\allpoly$ for any choice
of $A$. Consider the matrix 
$A = \smalltwo{0}{1}{-1}{0}$
and its powers:

\begin{xalignat*}{2}
  \smalltwo{0}{1}{-1}{0}^{4k}   & =  \smalltwo{1}{0}{0}{1} &  
  \smalltwo{0}{1}{-1}{0}^{4k+1} & =  \smalltwo{0}{1}{-1}{0} \\
  \smalltwo{0}{1}{-1}{0}^{4k+2} & =  \smalltwo{-1}{0}{0}{-1} & 
  \smalltwo{0}{1}{-1}{0}^{4k+3} & =  \smalltwo{0}{-1}{1}{0}
\end{xalignat*}

If $f(x) = \sum_0^n a_ix^i$ then
\begin{align*}
  f(xA) &= 
  \begin{pmatrix}
  \sum_i a_{4i}\,x^{4i} - a_{4i+2}\,x^{4i+2}  & 
  \sum_i a_{4i+1}\,x^{4i+1} - a_{4i+3}\,x^{4i+3} \\
 -\left(\sum_i a_{4i+1}\,x^{4i+1} - a_{4i+3}\,x^{4i+3}\right) &
  \sum_i a_{4i}\,x^{4i} - a_{4i+2}\,x^{4i+2} 
  \end{pmatrix} \\
&=
\begin{pmatrix}
  f_e(-x^2) & xf_o(-x^2) \\
  -xf_o(-x^2) & f_e(-x^2)
\end{pmatrix}
\end{align*}
where $f_e$ and $f_o$ are the even \index{even part} and odd parts of
$f$. If we consider $B = \smalltwo{0}{1}{1}{0}$ then a similar
computation shows that
$$ f(xB) = \begin{pmatrix}
  f_e(x^2) & xf_o(x^2) \\
  -xf_o(x^2) & f_e(x^2)
\end{pmatrix}
$$

Now $f_e(x^2)\in\allpoly$ iff $f_e\in\allpolyaltclose$, and 
$f_e(-x^2)\in\allpoly$ iff $f_e\in\allpolyposclose$, and hence $f_e(x)
= ax^r$ for some $a,r$. Similarly $f_o(x) = bx^s$ for some $b,s$, and
thus $ f(x) = a x^{2r} + b x^{2s+1}$.  If we choose $A=(1)$ we see
that $f\in\allpoly$, and hence $2r$ and $2s+1$ must be consecutive
integers. This concludes the proof of the lemma.
\end{proof}

\chapter{Matrix Polynomials }
\label{cha:mat-poly}

\renewcommand{\TimeStampStart}{Monday, December 17, 2007: 17:22:07}
\mytoday

\index{matrix polynomial}
\index{hyperbolic matrix polynomials}

  A \emph{matrix polynomial} is a polynomial whose coefficients are
  matrices. Equivalently, it is a matrix whose entries are
  polynomials. We assume all matrices are $\mu$ by $\mu$. We write 
\[
f(x) = A_0 + \cdots + A_n x^n
\]
The leading coefficient of $f(x)$ is the matrix $A_n$. We are
interested in the class of \emph{hyperbolic} matrix polynomials that
are a nice generalization of $\allpoly$. 

\begin{definition}
  The \emph{hyperbolic matrix polynomials in $d$ variables of degree
    $n$}, denoted $\hyper{d}(n)$,  consists of all matrix polynomials \\
  $f(\xx)=\sum A_\diffi \xx^\diffi$ such that
  \begin{enumerate}
  \item $<f\cdot v,v>=v^tf(\xx)v\in\rupint{d}$ for all $v\in\reals^\mu$.
  \item $f(\xx)$ has degree $n$.
  \item If $|\diffi|=n$ then $A_\diffi$ is positive definite.
  \end{enumerate}
  $\hyperpos{d}(n)$ is the subset of $\hyper{d}(n)$ where all coefficients
  are positive definite. 
If we don't wish to specify the degree we write $\hyper{d}$ and $\hyperpos{d}$.
\end{definition}

A hyperbolic matrix polynomial of degree zero is a constant matrix,
and must be positive definite.

If $g\in\rupint{d}$, and $P$ is a positive definite matrix then
$g\,P\in\hyper{d}$. Although this is a trivial construction, the
embedding $\rupint{d}\longrightarrow \hyper{d}$ is important.

  \section{Introduction}
  \label{sec:introduction}
  In this section we establish a few simple properties about
  hyperbolic matrix polynomials. The role of positive coefficients in
  $\allpoly$ is played by positive definite polynomials in
  $\hyper{d}$. For instance, a matrix polynomial of degree $1$ in
  $\hyper{d}$ has the form
\[
A_1x_1+\cdots +A_dx_d + S
\]
where the $A_i$ are positive definite, and $S$ is symmetric. If $S$ is
positive definite then the polynomial is in $\hyperpos{d}$. 

The analog of multiplying by a positive constant is replacing $f$ by
$A^tfA$. The leading coefficients are positive definite, and if $w=Av$
then
\[ <A^tfA\cdot v,v> = v^tA^tfAv = <f\cdot w,w>\in\rupint{d} .\]

In the case that $d$ is $1$ then we can multiply on the left and right
by the square root of the inverse of the leading coefficient, and so
we may assume that the coefficient of $x^n$ is the identity matrix.
\begin{lemma}
  Any principle submatrix of a matrix in $\hyper{d}$ is
  also in $\hyper{d}$. In particular, all diagonal entries are in
  $\rupint{d}$. 
\end{lemma}
\begin{proof}
  This is immediate, for we just need to set the $i$-th coordinate of
  $v$ equal to zero to show the hyperbolic property for the 
  $i$-th principle submatrix. The leading coefficients are positive
  definite since all principle submatrices of a positive definite
  matrix are positive definite. 
\end{proof}

  \begin{lemma}\label{lem:hyp-diag}
    Suppose that $f(x)$ is a diagonal matrix
    $diag(f_1(x),\dots,f_\mu(x))$. The following are equivalent:
    \begin{enumerate}
    \item $f\in\hyper{1}$.
    \item $f_1,\dots,f_\mu$ have a common interlacing.
    \end{enumerate}
  \end{lemma}
  \begin{proof}
    If $\mu=2$ then 
\[
(a,b)\smalltwo{f_1}{0}{0}{f_2}(a,b)^t = a^2 f_1 + b^2 f_2 \in\allpoly 
\]
Thus all positive linear combinations of $f_1$ and $f_2$ lie in
$\allpoly$, and so $f_1$ and $f_2$ have a common interlacing. It
follows that any two of the $f_i$ have a common interlacing, and so
they all have a
common interlacing.

Conversely, if they have a common interlacing then all positive linear
combinations are in $\allpoly$, and so
\[
(a_1,\dots,a_\mu)diag(f_1,\dots,f_\mu)(a_1,\dots,a_\mu)^t = a_1^2 f_1
+ \cdots a_\mu^2 f_\mu \in\allpoly.
\]

  \end{proof}

  \begin{remark}
\index{tensor product}
\index{Kronecker product}
    We can use this result to show that the tensor (Kronecker) product
    does not preserve hyperbolic matrix polynomials. If $f =
    \smalltwo{x}{0}{0}{x-1}$ then $f\otimes f\otimes f$ is an $8$ by
    $8$ diagonal matrix whose diagonal contains $x^3$ and
    $(x-1)^3$. Since these two do not have a common interlacing,
    $f\otimes f\otimes f$ is not in $\hyper{1}$.
  \end{remark}

Replacing $x$ by $Ax$ in a polynomial with all real roots does not usually
yield a hyperbolic matrix polynomial.

  \begin{cor}
    Suppose that $A$ is a positive definite matrix with eigenvalues
    $d_1,\dots,d_\mu$, and $f\in\allpoly$. Then $f(xA)\in\hyper{1}$ if
    and only if $\{f(d_1x),\dots,f(d_\mu x)\}$ have a common interlacing.
  \end{cor}
  \begin{proof}
    Let $O$ be orthogonal, and $D$ diagonal so that $A = O^tDO$. Then 
\[ f(xA) = O^t\,f(xD)\,O = O^t\,diag(f(xd_i))\,O \]
Thus, $f(xA)\in\hyper{1}$ if and only if
$diag(f(xd_i))\in\hyper{1}$. Lemma~\ref{lem:hyp-diag} completes the proof.
  \end{proof}

\begin{remark}
    There is a simple criterion for $\{f(d_ix\}_i$ to satisfy the
    corollary. Suppose that $\roots(f) = (r_i)$, and $\alpha>1$. The
    roots of $f(\alpha x)$ are $r_1/a<\cdots < r_n/\alpha$. If they
    have a common interlacing then
\[
r_1/\alpha < r_1 < r_2/\alpha < r_2 \cdots 
\]
and thus $\alpha< \min_{i>j} r_i/r_j$. If $0<\alpha<\beta$ then
$f(\alpha x)$ and $f(\beta x)$ have a common interlacing if and only
if $f(x)$ and $f(\beta x/\alpha)$ do. We conclude that $\{f(d_ix)\}$
has a common interlacing if and only if
\[
\max _{i, j} \frac{d_i}{d_j} \le \min_{i>j} \frac{r_i}{r_j}
\]
Note that if this condition is satisfied then they are actually
mutually interlacing.
\index{interlacing!common}\index{common!interlacing}
\index{mutually interlacing}
  \end{remark}

If we add off-diagonal entries that are all equal, then we can
determine when the matrix polynomial is in $\hyper{1}$. Define
\[
\sigma_n = \max_{a_1,\dots,a_n} 
\frac{a_1a_2+a_2a_3+a_3a_4 + \cdots +a_{n-1}a_n}%
{a_1^2+\cdots +a_n^2}
\]
Some values of $\sigma_n$ are
\[
\sigma_2=\frac{1}{2} \qquad\sigma_3=\frac{\sqrt{2}}{2}
\qquad\sigma_4=\frac{1+\sqrt{5}}{4} \qquad\sigma_5=\frac{\sqrt{3}}{2}
\]

\begin{lemma}\label{lem:hyp-tridiag}
The $n$ by $n$ tridiagonal matrix below is in $\hyper{1}$ if and only if
$f+2\alpha g\in\allpoly$ for $|\alpha|\le \sigma_n$.
\[
m(x)=\begin{pmatrix}
  f & g & 0 & \hdots \\
  g & f & g & 0 & \hdots \\
  0 & g & f & g & 0 & \hdots\\
  0 & 0 & \ddots & \ddots & \ddots & \ddots \\
\end{pmatrix}
\]
\end{lemma}
\begin{proof}
  If $v=(a_1,\dots,a_n)^t$, then
\[
v^tmv = (a_1^2+\cdots+a_n^2)f + 2(a_1a_2+a_2a_3+a_3a_4 + \cdots
+a_{n-1}a_n)g
\]
Dividing by the first factor, we see that this is in $\allpoly$ if and
only if $f+2\alpha g\in\allpoly$ for $|\alpha|\le \sigma_n$.
\end{proof}

The off-diagonal element can even be a constant if it is small enough.

\begin{cor}
  Suppose that $f\in\allpoly$. The following are equivalent:
  \begin{enumerate}
  \item $\smalltwo{f}{c}{c}{f}\in\hyper{1}$
  \item $|c|\le \displaystyle\min_{f'(\alpha)=0}|f(\alpha)|$
  \end{enumerate}
\end{cor}
\begin{proof}
  Geometrically, the second condition means that the largest value of
  $|c|$ satisfies  $f+b\in\allpoly$ for all
  $|b|\le|c|$. Since $2\sigma_2=1$, the result now follows from the
  lemma.
\end{proof}

The off diagonal entries are not necessarily in $\rupint{d}$. This can be
deduced from Lemma~\ref{lem:hyp-tridiag} by taking $g$  small
relative to $f$. For example, 
\[
\begin{pmatrix}
  x^2-1 & \epsilon (x^2+1) \\ \epsilon(x^2+1) & x^2-1
\end{pmatrix}\in\hyper{1} \Longleftrightarrow |\epsilon|\le 1
\]

The reverse preserves $\hyper{1}$.
\index{reverse!in $\hyper{1}$}

\begin{lemma}\label{lem:hyp-reverse}
  If $f(x)=\sum_{i=0}^n A_ix^i\in\hyper{1}$, and $A_0$ is positive
  definite, then the reverse $f^{rev}(x)=\sum_{i=0}^n A_ix^{n-i}$ is
  also in $\hyper{1}$. 
\end{lemma}
\begin{proof}
Since $A_0$ is positive definite, the leading coefficient of $f^{rev}$
is positive definite. The conclusion follows from the calculation
\[
<vf^{rev},v> = <vx^nf(1/x),v> = x^n<vf(1/x),v> = <vf,v>^{rev}.
\]
\end{proof}

 The product of hyperbolic matrix polynomials is generally
  not hyperbolic since the product of the leading coefficients, which
  are positive definite matrices, is not usually positive
  definite. We do have
  \begin{lemma}
    If $f,g\in\hyper{1}$ then $f^t\,g\,f\in\hyper{1}$.
  \end{lemma}
  \begin{proof}
    If $A$ is the leading coefficient of $f$ and $B$ is the leading
    coefficient of $g$, then the leading coefficient of $f^t\,g\,f$ is
    $A^tbA$ which is positive definite. If $v\in\reals^\mu$ and $w =
    fv$ then
\[
v^t\,f^t\,g\,f\,m\,v = w^t\,g\,w \in\allpoly
\]
  \end{proof}

  \section{Interlacing}
  \label{sec:interlacing}
\index{interlacing!hyperbolic matrix polynomials}
\index{hyperbolic matrix polynomials!interlacing}

There are two equivalent ways to define interlacing in $\hyper{d}$:
the usual linearity definition, or reduction to one variable. We start
with the latter, since it is more precise.

\begin{definition}
  Suppose that $f,g\in\hyper{d}$ have degree $n$ and $n-1$
  respectively. We define
  \begin{align*}
    f & \lesslesseq g & \text{if and only if} &\quad v^tfv \lesslesseq
    v^tgv\qquad \text{for all }v\in\reals^\mu\setminus 0\\
    f & \lessless g & \text{if and only if} & \quad v^tfv \lessless
    v^tgv\qquad \text{for all }v\in\reals^\mu\setminus 0\\
  \end{align*}
  $\greateq$ and $\greateqeq$ are defined similarly.
\end{definition}

\begin{lemma}
  Suppose that $f,g\in\hyper{d}$ have degree $n$ and $n-1$
  respectively. The following are equivalent:
  \begin{enumerate}
  \item $f+\alpha g\in\hyper{d}$ for all $\alpha\in\reals$.
  \item $f\lesslesseq g$
  \end{enumerate}
\end{lemma}
\begin{proof}
  By hypothesis we know that 
  \begin{align*}
  v^t(f+\alpha g)v =   v^tfv + \alpha v^t g v & \in \rupint{d}\qquad\text{for all }\alpha\in\reals
  \end{align*}
  and consequently we have that $ v^tfv  \lesslesseq v^tgv$. The
  converse is similar.  
\end{proof}

We now have some easy consequences of the definition:

\begin{lemma}
  Assume that $f,g,h\in\hyper{d}$. 
  \begin{enumerate}
  \item If $f\lesslesseq g$ and $f\lesslesseq h$ then $f\lesslesseq
    g+h$. In particular, $g+h\in\hyper{d}$.
  \item If $f\lesslesseq g\lesslesseq h$ then $f-h\lesslesseq g$.
  \item If $k\in\rupint{d}$ then $fk\in\hyper{d}$.
  \item If all $a_i$ are positive and $f\in\hyper{d}$ then
\[
(a_1x_1+\cdots a_dx_d+b)f \lesslesseq f
\]
\item If the $i$-th diagonal element of $f$ is $f_i$,
  and of $g$ is $g_i$ then $f\greateqeq g$ implies $f_i\greateqeq g_i$.
\item If $F\lesslesseq G$ in $\rupint{d}$ and $M_1,M_2$ are positive
  definite matrices then $F\,M_1\lesslesseq G\,M_2$ in $\hyper{d}$.
\item If $g_1,\dots,g_r$ have a common interlacing in $\rupint{d}$, and
  $M_1,\dots,M_r$ are positive definite matrices then
$ g_1\,M_1+\cdots+ g_r\,M_r\in\hyper{d}.$
  \end{enumerate}
\end{lemma}
\begin{proof}
  By hypothesis we know that 
$    v^tfv  \lesslesseq v^tgv$ and 
$ v^tfv  \lesslesseq v^thv$.
  Since $f,g,h\in\hyper{d}$ we know that the leading coefficients of
  these three terms are positive, so we can add interlacings in
  $\rupint{d}$:
\[
v^tfv \lesslesseq v^tgv+v^thv = v^t(g+h)v
\]
The proof of the second one is similar. The third one follows since
$\rupint{d}$ is closed under multiplication, and $k$ is a scalar:
\[
<fk\cdot v,v> = k\,<f\cdot v,v> \in\rupint{d}
\]
The fourth follows from the third one. For the next one, note that if
$e_i$ is the vector with $1$ in the $i$-th place and zeros elsewhere,
then $e^t_ife_i$ is the $i$-th diagonal element of $f$. 

If $F\lesslesseq G$ then $v^tM_1v$ and $v^tM_2v$ are positive so
$v^tFM_1v\lesslesseq v^tGM_2v$. If $G\lesslesseq g_i$ for $1\le i \le
r$ then $G\,I\lesslesseq g_i \,M_i$ by the previous result, and the
conclusion follows by adding interlacings.
\end{proof}

  The interlacing of even linear polynomials leads to some
  complexity. 
  \begin{lemma} Suppose that $A,C$ are positive definite, and $B,D$
    are symmetric.
    \begin{enumerate}
    \item $xI-B \greateqeq xI-D$ iff $B\ge D$ (i.e. $B-D$ is positive definite).
    \item $xA-B \greateqeq xC-D$ iff 
      \begin{equation}
        \label{eqn:hyp-linear}
        \frac{v^tBv}{v^tAv} \ge 
        \frac{v^tDv}{v^tCv}\quad\text{for all $v$}
      \end{equation}
    \end{enumerate}
  \end{lemma}    
  \begin{proof}
    If $xI-B \greateqeq xI-D$ then for all non-zero $v$ we have
  \begin{align*}
    v^t(xI-B)v & \greateqeq v^t(xI-D)v \\
    (v^tv)x - v^tBv & \greateqeq (v^tv)x-v^tDv     
  \end{align*}
  These are linear polynomials with positive leading coefficients, and
  they interlace iff $v^tBv \ge v^tDv$. Thus $v^t(B-D)v\ge0$, so $B\ge
  D$.

  In the second case 
\[
v^tAv\,x- v^tBv \greateqeq v^tCv\,x - v^tDv 
\]
for all non-zero $v$. Solving each equation for $x$ and using the
first part yields the conclusion.
    
  \end{proof}

\noproblem{give an example of degree 2 $<$ degree 1 that is not derivative}
  The condition \eqref{eqn:hyp-linear} is not well understood.  Let
  $\lambda_{min}(W)$ and $\lambda_{max}(W)$ be the smallest and
  largest eigenvalues of a matrix $W$. Suppose that matrices $A,B,C,D$
  are as in the lemma and \eqref{eqn:hyp-linear} is satisfied for all
  non-zero $v$.  If $W$ is any symmetric matrix then the Rayleigh-Ritz
  inequality \cite{horn-johnson-1} says that
\[
\lambda_{min}(W) \le \frac{v^tWv}{v^tv} \le \lambda_{max}(W) \quad \text{for
  all }v\ne 0
\]
We may assume that $v$ has norm $1$. The values of $v^t\,W\,v$ lie
in the interval $(\lambda_{min}(W),\lambda_{max}(W))$. If
\eqref{eqn:hyp-linear} holds it follows that

\begin{gather}
  \label{eqn:hyp-ineq-1}
        (v^tBv)\,(v^tDv) \ge 
(v^tAv)\,(v^tCv)\quad\text{for all $v$}\notag\\
  \lambda_{max}(B)  \lambda_{max}(D) \le   \lambda_{min}(A)  \lambda_{min}(C)
\end{gather}

The values of $\frac{v^t\,B\,v}{v^t\,Av}$ lie in the interval
$[\lambda_{min}(BA^{-1}), \lambda_{max}(BA^{-1}]$ and therefore 

\begin{gather}
  \label{eqn:hyp-ineq-2}
  \lambda_{min}(BA^{-1})  \le  \lambda_{max}(DC^{-1}) \qquad\text{and}\\
  \lambda_{min}(BA^{-1})   \lambda_{min}(CD^{-1}) \le 1 \notag
\end{gather}

Conversely, if we know that the stronger condition 
\[
  \lambda_{max}(BA^{-1})  <  \lambda_{min}(CD^{-1}) 
\]
holds then we can conclude that \eqref{eqn:hyp-linear} holds.

$\hyper{d}$ behaves well under induced linear transformations.

\begin{lemma}
  Suppose that $T\colon{}\rupint{d}\longrightarrow\gsub_d$, and if $d$ is $1$
  then $T$ preserves the sign of the leading coefficient. Define a map
  on $\hyper{d}$ by $T(\sum A_\diffi\xx^\diffi) = \sum A_\diffi
  T(\xx^\diffi)$. Then $T\colon{}\hyper{d}\longrightarrow\hyper{d}$.
\end{lemma}
\begin{proof}
  If $f\in\hyper{d}$ then the leading coefficient of $T(f)$ is
    positive definite. The result now follows from
\[
T \,< f\cdot v,v> = <T(f)\cdot v,v>
\]

\end{proof}

The derivative behaves as expected. The proofs follow from the
previous lemma, and the arguments for $\allpoly$. 

\begin{lemma}
  Suppose that $f,g\in\hyper{d}$ and $h\in\allpolypos$.
  \begin{enumerate}
  \item $\frac{\partial f}{\partial x_i}\in\hyper{d}$.
  \item $f \lesslesseq \frac{\partial f}{\partial  x_i}$
  \item If $f\lesslesseq g$ then $\frac{\partial f}{\partial  x_i}
    \lesslesseq \frac{\partial g}{\partial  x_i}$
  \item $h(\partial/\partial x_i) f\in\hyper{d}$
  \end{enumerate}
\end{lemma}

\begin{example}
  What does the graph of a diagonal matrix polynomial look like?
  Consider $f = \smalltwo{(x-1)(x-3)}{0}{0}{(x-2)(x-4)}$. We have the
  interlacings
\[
(x-1)(x-2)I \greateqeq f \greateqeq (x-3)(x-4)I
\]
If we graph $(1,v)^tf(1,v)$ then Figure~\ref{fig:hyper-3} shows that
it stays in the regions determined by $(x-1)(x-2)=0$ and
$(x-3)(x-4)=0$. The small dashed lines are the graph of
$(x-1)(x-3)=0$, and the large dashed lines are the roots of
$(x-2)(x-4)=0$. This is a general phenomenon -- see
page~\pageref{sec:root-zones}. 
\begin{figure}[htbp]
  \centering
  \includegraphics*[height=2in]{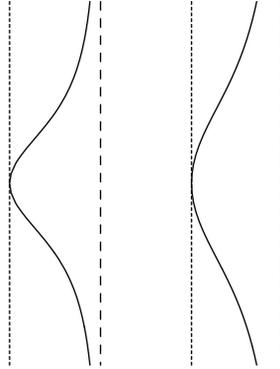}
  \caption{The graph of a diagonal matrix polynomial}
  \label{fig:hyper-3}
\end{figure}
\end{example}

\section{Common interlacings}
\label{sec:common-interlacings}
\index{interlacing!hyperbolic matrix polynomials}
\index{hyperbolic matrix polynomials!interlacing}
\index{common interlacing!hyperbolic matrix polynomials}
\index{hyperbolic matrix polynomials!common interlacing}

If $f\in\hyper{1}$ and $g\in\allpoly$, then $f\lesslesseq gI$ if and
only if $v^tfv \lesslesseq g$ for all $v\in\reals^\mu$. If either
condition holds we say that \emph{$f$ has a common interlacing}.
\index{common interlacing!in $\hyper{1}$}

The usual definition\cite{markus} of hyperbolic polynomials allows the
vectors $v$ to be complex. This is equivalent to the existence of a
common interlacing. Recall that $w^\ast$ is the conjugate transpose.

\begin{lemma}\label{lem:root-zone-c}
  Suppose that $f$ is a matrix polynomial whose leading coefficient is
  the identity. The following are equivalent.
  \begin{description}
  \item[1] $w^\ast fw\in\allpoly$ for all $w\in\complexes^\mu$.
  \item[2a] $\{v^tfv \mid v\in\reals^\mu\}$ has a common interlacing.
  \item[2b] $f$ has a common interlacing.
  \end{description}
\end{lemma}

\begin{proof}
  $2a$ and $2b$ are equivalent by definition, so assume that $1$
  holds. Write $w=u+\imag v$ where $u,v$ are real vectors. Then
  \begin{align*}
    w^\ast f w &= u^\ast fu + u^\ast f(\imag v) + (\imag v)^\ast fu +
    (\imag v)f\imag v \\
&= u^t fu + v^t f v
  \end{align*}
Replacing $v$ by $\alpha v$ shows that $u^t fu +\alpha^2 v^t f
v\in\allpoly$ for all $\alpha$, so $u^tf u$ and $v^tfv$ have a common
interlacing. Conversely, if all $u^tf u$ and $v^tfv$ have a common
interlacing then $w^\ast fw\in\allpoly$. 
\end{proof}

Do all hyperbolic polynomials in $\hyper{1}$ have a common interlacing?
This is an unsolved question, but it is true for two by two matrices.

\begin{lemma}\label{lem:hyper-2x2}
Any two by two matrix   $F\in\hyper{1}$ has a common interlacing.
\end{lemma}
\begin{proof}
  If we write $F=\smalltwo{f}{h}{h}{g}$ and let $v=(\beta,1)$ then
\[
v^tfv = \beta^2\,f + 2\beta\, h + g
\]
Lemma~\ref{lem:mat-poly-2} guarantees the existence of a common interlacing.
\end{proof}

There is a simple condition that implies a common interlacing.
\begin{lemma}
  Assume $f(x) = x^2 I + Ax+B$ where $A,B$ are symmetric. If there is
  an $\alpha$ such that $f(\alpha)$ is negative definite then $f$ has
  a common interlacing.
\end{lemma}
\begin{proof}
  Since $v^tf(\alpha)v<0$ and $v^tf(\beta)v>0$ for $\beta$
  sufficiently large, it follows that $f \lessless (x-\alpha)I$. 
\end{proof}

We have seen that the polynomials of a diagonal hyperbolic matrix
polynomial have a common interlacing. This is true in general.

\begin{cor}
  The diagonal elements of a polynomial in $\hyper{1}$ have a common
  interlacing. 
\end{cor}
\begin{proof}
  Any two by two principle submatrix is in $\hyper{1}$, and by
  Lemma~\ref{lem:hyper-2x2} the diagonal elements have a common
  interlacing. Since any two diagonal elements have a common
  interlacing, they all have a common interlacing.
\end{proof}

\begin{cor}
  If $f\in\hyper{1}$ has a common interlacing then
  the determinant and trace of $f$ are in $\allpoly$. 
\end{cor}
\begin{proof}
  If $det(f(\lambda))=0$ then there is a $w\in\complexes$ so that
  $f(\lambda)w=0$, and thus $w^\ast f(\lambda) w=0$. By
  Lemma~\ref{lem:root-zone-c} we know that $w^\ast fw\in\allpoly$ for
  all complex $w$, so $\lambda$ is real.

  Since the diagonal elements have a common interlacing, their sum is
  in $\allpoly$.
\end{proof}

\begin{lemma}
  If $f\lesslesseq g$ in $\hyper{1}$ then $trace(f)\lesslesseq
  trace(g)$ in $\allpoly$.
\end{lemma}
\begin{proof}
  For any $\alpha\in\reals$ we know that $f+\alpha g\in\hyper{1}$. Now
  the trace is a linear map $\hyper{1}\longrightarrow\allpoly$,
  so $trace(f) + \alpha\, trace(g)\in\allpoly$, which finishes the proof.
\end{proof}

\begin{example}
  Here is a simple example of a polynomial with a common
  interlacing. If $v = (a,b)$ and
  \begin{align*}
    f &= x \begin{pmatrix}1 &0\\0&1\end{pmatrix} + 
\begin{pmatrix}0&1\\1&0\end{pmatrix}\\
\intertext{then}
v^tfv &= (a^2+b^2)x + 2ab 
  \end{align*}
Consequently, the roots of $v^tfv$ lie in the closed interval $[-1,1]$
and so \[ (x-1)(x+1)I \lesslesseq f.\]
The determinant of $f$ is $x^2-1$, which is in $\allpoly$, as the
corollary tells us.
\end{example}

\index{adjugate}

  Empirical investigations suggest that the adjugate of a matrix that
  determines a polynomial in $\gsubplus_d$ is in $\hyper{d}$. We prove
  this for $d=1$, and show that there is a common interlacing for all $d$.

  Assume that $D_1,\dots,D_d$ are positive definite $\mu$ by $\mu$
  matrices, and consider the matrix
\[
M = Id + x_1D_1 + \cdots + x_d D_d.
\]
Recall that $M[i,j]$ is the matrix formed by removing the $i$'th row
and $j$'th column of $M$. The \index{adjugate}adjugate of $M$ is
\[
Adj[M] = \left( \left| (-1)^{i+j}M[i,j]\right|\right).
\]
The adjugate satisfies
\[
Adj[M]\cdot M = |M| Id
\]
and our conjecture is that $Adj(M)\in\hyper{d}$. First of all, we show
that $Adj(M)$ has a common interlacing.
\begin{lemma}
  If $M$ is as above and $v\in\reals^\mu$ then $|M|\lesslesseq v^t\,
  Adj(M)\,v$. In particular, $v^t\,Adj(M)\,v\in\gsubcloseplus_d$.
\end{lemma}
\begin{proof}
  Choose $v=(v_i)\in\reals^\mu$ and define $u =
  ((-1)^iv_i)$. Introducing a new variable $y$, define
\[
N = M + y\cdot uu^t.
\]
Now $N\in\gsubcloseplus_{d+1}$ since $uu^t$ is positive semi-definite.
Also, $uu^t$ has rank one, so $N$ is linear in $y$. Expanding $|N|$
yields
\begin{align*}
|N| &= |M| + y \sum_{i,j} (-1)^{i+j}\,v_iv_j\,|M[i,j]|\\
&= |M| + y\, v^t \,Adj(M)\,v
\end{align*}
We conclude that $|M|\lesslesseq v^t\,Adj(M)\,v$ since consecutive
coefficients interlace.
\end{proof}

In order to show that $Adj(M)\in\hyper{d}$we need to show that the
coefficients of the homogeneous part are positive definite. It's easy
for $d=1$.

\begin{lemma}
 If $D$ is positive definite then $Adj(Id+xD)\in\hyper{1}$. 
\end{lemma}
\begin{proof} 
    Let $M= Id + x D$. The degree of $|M[i,j]|$ is at most $\mu-1$, so
    the degree of $Adj(M)$ is at most $(\mu-1)\mu$. The coefficient of
    $x^{\mu(\mu-1)}$ is $((-1)^{i+j}\,|D[i,j]|) = |D|D^{-1}$. Since
    $D^{-1}$ is positive definite, we see that $Adj(M)\in\hyper{1}$.
  \end{proof}

\index{norm!of a matrix}
  We can sometimes substitute matrices for constants, and we get
  hyperbolic matrices with a common interlacing.  Recall that the norm
  $\norm{M}$ is the maximum eigenvalue of a positive definite matrix
  $M$ and satisfies
\[ \norm{M} = \sup_{v\ne0} \frac{v^t\,M\,v}{v^t\,v}. \]

\begin{lemma}
  Suppose that $M$ is positive definite and $deg(f)>deg(g)$. 
  \begin{enumerate}
  \item If $f+\alpha g\in\allpoly$ for $0\le \alpha\le \norm{M}$ then $fI+gM\in\hyper{1}$.
  \item $fI \lesslesseq fI+gM$
  \end{enumerate}
\end{lemma}
\begin{proof}
  The leading coefficient of $fI+gM$ is $I$. To check the inner product
  condition, compute 
\[
v^t\left( fI+gM\right)v = {v^tv}\left(f + \frac{v^tMv}{v^tv}g\right)
\]
Since $\frac{v^tMv}{v^tv}\le \norm{M}$ the inner product condition
holds. Next,
\[
f\,I + \alpha(fI+gM) = (1+\alpha)\,\left( fI +
  \frac{\alpha}{1+\alpha}gM\right)\]
and the latter polynomial is in $\hyper{1}$ since
$\norm{\frac{\alpha}{1+\alpha}M}\le \norm{M}$ for all positive $\alpha$.
\end{proof}

  \begin{lemma}
    If $f\in\hyper{d}$ has a common interlacing then $|f|\in\rupint{d}$.
  \end{lemma}
  \begin{proof}
    If we substitute for all but one of the variables then we get a
    polynomial in $\hyper{1}$ with a common interlacing. We know that
    its determinant is in $\allpoly$, so substitution is satisfied.

    Assume that $f$ has degree $n$.  If we write $f = \sum_{\sdiffi}
    \aaa_\sdiffi \xx^\sdiffi$ then the homogeneous part of $|f|$ is
    $\vert\sum_{|\sdiffi|=n} \aaa_\sdiffi \xx^\sdiffi\vert$. If we
    replace each monomial $\xx^\sdiffi$ by a new variable
    $\yy_\sdiffi$ then we know that 
    $|\sum_{|\sdiffi|=n} \aaa_\sdiffi\yy_\sdiffi|$ has all positive
    coefficients since all the $\aaa_\sdiffi$ are positive
    definite. Replacing the $\yy_\sdiffi$ by $\yy^\sdiffi$ shows that
    $|f|^H$ has all positive coefficients.

  \end{proof}

  \begin{lemma}
    If $f\lesslesseq g$ in $\hyper{1}$ then $f + y g\in\hyperclose{2}$.
  \end{lemma}
  \begin{proof}
    It suffices to show that 
\[
F_\epsilon(x,y) = f(x+\epsilon y) + (y+\epsilon x)g(x+\epsilon y)
\in\hyper{2}
\]
for positive $\epsilon$. Now the homogeneous part is
\[
F_\epsilon^H = 
f^H(x+\epsilon y) + (y+\epsilon x)g^H(x+\epsilon y)
\]
and $f^H,g^H$ are matrix polynomials with positive definite
coefficients it follows that $F_\epsilon^H$ is a sum of polynomials
with positive definite coefficients, and so all its coefficients are
positive definite.

Next, we verify substitution. Fix $v\in\reals^n$ and let $G(x,y) =
v^tfv + y v^t g v$. By hypothesis $v^tfv\lesslesseq v^tgv$, so we know
that $G(x+\epsilon y,y+\epsilon x)\in\rupint{2}$. Thus, $F_\epsilon\in\hyper{2}$.
  \end{proof}

  \begin{lemma}
    If $gI\lessless f$ in $\hyper{1}$ then
    \begin{enumerate} 
    \item $|gI+yf|\in\gsubclose_2$
    \item If we write $|I + yf| = \sum f_i(x)y^i$ then $f_i\in\allpoly$.
    \item $f_i \lesslesseq g f_{i+1}$.
    \end{enumerate}
  \end{lemma}
  \begin{proof}
    The first part follows from the previous lemma. If we expand the
    determinant we get
\[
|gI+yf| = \sum g^{\mu-i} f_i(x)y^i
\]
Since this is in $\rupint{2}$ we see that $g^{\mu-i}f_i \lessless
g^{\mu-(i-1)}f_{i+1}$ which implies that $f_i\lesslesseq g f_{i+1}$. 
  \end{proof}

The next corollary generalizes the fact the the determinant and trace
are in $\allpoly$.

\begin{cor}
  If $f\in\hyper{1}$ has a common interlacing then
$\displaystyle\sum_{|S|=k} |f[S]| \in\allpoly$.
\end{cor}
\begin{proof}
  We have $S\subset\{1,\dots,n\}$ and $f[S]$ is the submatrix of $f$
only using the rows and columns in $S$. The sum in question is the
coefficient of $y^{n-k}$ in $|ygI+f|$, and this is in $\allpoly$ by
the last lemma.
\end{proof}

\begin{remark}
    \added{7/4/7}
    If we drop the assumption that matrices are positive definite then
    we can find a $2\times2$ matrix $M$ that satisfies
    $v^tMv\in\allpoly$ for all $v$, yet has no common
    interlacing. Choose
\[
F(x,y,z) = f + y\, g + z\, h + yz\,k \in\gsubclose_3
\]
Choose $\alpha,\beta\in\reals$ and $\alpha\ne0$.
Since
\[
\begin{pmatrix}
\alpha & \beta  
\end{pmatrix}
\begin{pmatrix}
  f & g \\ h & k
\end{pmatrix}
\begin{pmatrix}
  \alpha\\\beta
\end{pmatrix}
=
\alpha^2\bigl[ f + \frac{\beta}{\alpha}(g+h) +
\bigl(\frac{\beta}{\alpha}\bigr)^2 k \bigr]
\]
we see it is in $\allpoly$ since it equals
$F(x,\beta/\alpha,\beta/\alpha)$. For a particular example, consider
\begin{multline*}
  (x + y + z + 1) (x + 2 y + 3 z + 5) \\=
5 + 6 x + x^2 + y(7 + 3 x) + z(8 + 4 x) + 5 yz
+ 2y^2 + 3 z^2
\end{multline*}
Taking coefficients yields the matrix
\[
M =
\begin{pmatrix}
 5 + 6 x + x^2 &  8 + 4 x \\7 + 3x & 5
\end{pmatrix}
\]
The table below shows that there is no common interlacing.

\begin{center}
\begin{tabular}{rcc}
  \toprule
  b & \multicolumn{2}{c}{\text{roots of $(1\, b)M(1\,b)^t$}} \\
  \midrule
  -5 & -35.1 & -5.8 \\
  5 & 2.0 &  26.9\\
  \bottomrule
\end{tabular}
\end{center}

\end{remark}

  \section{Orthogonal-type recurrences}
  \label{sec:orthogonal-sequences}
\index{orthogonal polynomials!hyperbolic matrix polynomials}
\index{hyperbolic matrix polynomials!orthogonal sequences}

\label{sec:orth-type-recurr}

The fact that
$\quad
f\lesslesseq g \lesslesseq h \implies f-h\lessless g
\quad$
is all we need to establish orthogonal type recurrences in
$\hyper{d}$. 

\begin{lemma}
  Define a sequence of matrix polynomials by 
  \begin{align*}
    f_0 &= U \\
    f_1 &= A_1x_1 + \cdots +A_dx_d+S \\
    f_k &= (a_{k,1}x_1+\cdots +a_{k,d}x_d+b_k) f_{k-1} - c_k f_{k-2}
  \end{align*}
  where $U$ and the $A_i$ are positive definite,  $S$ is symmetric,
  all $c_k$ are positive, and all $a_{i,j}$ are positive.

  Then all $f_k$ are in $\hyper{d}$ and $f_k \lesslesseq f_{k-1}$.
\end{lemma}
\begin{proof}
  We note that $f_0$ and $f_1$ are in $\hyper{d}$. By induction we
  have that
\[
(a_{k,1}x_1+\cdots a_{k,d}x_d+b_k) f_{k-1} \lesslesseq f_{k-1}
\lesslesseq f_{k-2}
\]
and the observation above finishes the proof.
\end{proof}

The off-diagonal entries of a matrix polynomial in $\hyper{1}$ might
not have all real roots. Here's an explicit example. The
coefficient of $x$ in $f_1$ below is positive definite, as is $f_0$,
so $f_2\in\hyper{1}$. However, the off diagonal entry of $f_2$ has two
complex roots.  The diagonal elements have roots $(-2.4,.4)$ and
$(-2.5,1)$, so they have a common interlacing.

In general, all the diagonal elements of the sequence $f_0,f_1,\dots$
form an orthogonal polynomial sequence. The off diagonal elements
obey the three term recursion, but the degree $0$ term might not be
positive, and so they might not have all real roots.

\begin{align*}
  f_0 &=
  \begin{pmatrix}
    2 & 3 \\ 3 & 6
  \end{pmatrix} \\
f_1 &= 
  \begin{pmatrix}
    1 & -1 \\ -1 & 2
  \end{pmatrix}x + I \\
f_2 &= (x+1)f_1 - f_0 \\ 
&= \left(
\begin{array}{ll}
 x^2+2 x-1 & -x^2-x-3 \\
 -x^2-x-3 & 2 x^2+3 x-5
\end{array}
\right)
\end{align*}

\begin{example}\label{ex:hyp-p2}
Next, an example in two variables. 
The coefficients of $x$ and $y$ in $f_1$ below are positive
definite, but not all entries of $f_2$ are in $\rupint{d}$ since the
homogeneous part has positive and negative coefficients.

\begin{align*}
  f_0 &=
  \begin{pmatrix}
    2 & 3 \\ 3 & 6
  \end{pmatrix} \\
f_1 &= 
  \begin{pmatrix}
    1 & -1 \\ -1 & 2
  \end{pmatrix}x +
  \begin{pmatrix}
    3 & 1 \\ 1 & 2
  \end{pmatrix}y + I  \\
f_2 &= (x+2y+1)f_1 - f_0  \\ 
&= \left(
\begin{array}{ll}
 x^2+5 y x+2 x+6 y^2+5 y-1 & -x^2-y x-x+2 y^2+y-3 \\
 -x^2-y x-x+2 y^2+y-3 & 2 x^2+6 y x+3 x+4 y^2+4 y-5
\end{array}
\right)
\end{align*}
\end{example}

  \section{Problems with products}
  \label{sec:probl-with-prod}

It's reasonable to think that a product such as 

\[ (xI-S_1)(xI-S_2)\cdots(xI-S_n) \] 
should be in $\hyper{1}$, but this fails, even in the diagonal case.
If we take each $S_i$ to be a diagonal matrix, then the product is a
diagonal matrix, and we know that a necessary condition is that there
is a common interlacing.  To try to satisfy this condition, we can
require that all eigenvalues of $S_1$ are less than all eigenvalues of
$S_2$, which are less than all eigenvalues of $S_3$, and so on.
However, this does not always produce hyperbolic matrix polynomials.
If we perturb a diagonal matrix that gives a hyperbolic polynomial by
conjugating some of the terms, then we can lose hyperbolicity. 

\begin{figure}[htbp]
  \centering
  \includegraphics*[height=1.5in]{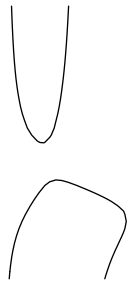}
  \caption{The graph of a  matrix product not in $\hyper{1}$}
  \label{fig:hyp-not}
\end{figure}

If
$e^{\imag\,\tau} = \smalltwo{\cos\tau}{\sin
  \tau}{-\sin\tau}{\cos\tau}$ then the graph of the inner product of
the polynomial
\[
f(x) = \bigl(xI-e^{-2\imag}\smalltwo{1}{0}{0}{2}e^{2\imag}\bigr)
\bigl(xI-e^{-\imag}\smalltwo{3}{0}{0}{4}e^{\imag}\bigr)
\]
is given in Figure~\ref{fig:hyp-not}, and is clearly not in
$\hyper{1}$.

\begin{figure}[htbp]
  \centering
  \includegraphics*[width=3in,height=2in]{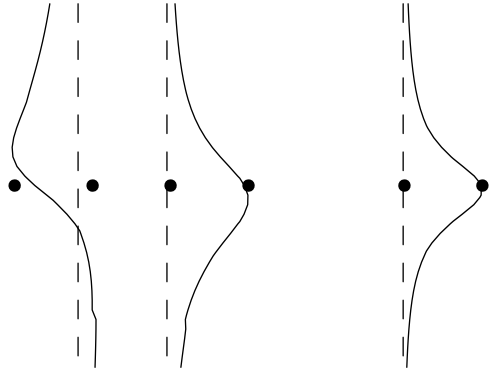}
  \caption{The graph of a hyperbolic matrix polynomial}
  \label{fig:hyp-mat-poly}
\end{figure}

Here is an example of a product that is in $\hyper{1}$. 
\[
\biggl(xI-
\begin{pmatrix}
 1.24 & -0.42 \\
 -0.42 & 1.75
\end{pmatrix}
\biggr)
\left(xI-
\begin{pmatrix}
 3.97 & -0.15 \\
 -0.15 & 3.02
\end{pmatrix}
\right)
\left(xI-
\begin{pmatrix}
 6.99 & -0.04 \\
 -0.04 & 6.00
\end{pmatrix}
\right)
\]
where the eigenvalues of the matrices are $(1,2)$, $(3,4)$,
$(6,7)$.  The graph of $(1,y)^tf(1,y)$ is in
Figure~\ref{fig:hyp-mat-poly}. The three dashed lines are the
asymptotes, and the dots are the roots of the determinant, drawn on the
$x$-axis. The graph shows that $f\in\hyper{1}$, since each horizontal
line meets the graph in exactly three points.

  Some constructions that you would expect to give hyperbolic
  polynomials give stable matrix polynomials. For instance, if $A$ is
  positive definite then $(xI+A)^2$ is generally not hyperbolic, but
  is a stable matrix polynomial. \seepage{lem:stable-power}

  \section{Sign Interlacing}
  \label{sec:sign-interlacing}
\index{sign interlacing!hyperbolic matrix polynomials}
\index{hyperbolic matrix polynomials!sign interlacing}

\noproblem{ I use $gI<f$ and $f<gI$ - are these different?}

  Recall that $f$ sign interlaces $g$ if the sign alternates as we
  evaluate $f$ at the roots of $g$. There is a similar definition for
  $f\in\hyper{1}$, where the role of alternating signs is replaced by
  alternating positive and negative definite matrices.

\begin{lemma}\label{lem:hyp-sign}
  Suppose that $g\in\allpoly$ has roots $r_1<\cdots<r_n$.\\  If
  $f\in\hyper{1}$ satisfies $gI \lessless f$ then
\[
(-1)^{n+i} f(r_i)\ \text{is positive definite for }\ 1\le i \le n
\]
\end{lemma}
\begin{proof}
The hypotheses imply that $  v^t(gI)v = (v^tv)g \lessless v^tfv$, so
$v^tfv$ sign interlaces $g$. Thus, $(-1)^{n+i}v^tf(r_i)v>0$  for all
non-zero $v$, and hence $(-1)^{n+i}f(r_i)$ is positive definite.
\end{proof}

The converse can be used to show that a matrix polynomial is
in $\hyper{1}$.

\begin{lemma}\label{lem:hyp-si}
  Suppose that $g\in\allpoly$ has roots $r_1<\cdots<r_{n+1}$.  If
\begin{enumerate}
\item $(-1)^{n+i} f(r_i)\ \text{is positive definite for }\ 1\le i \le n+1$
\item $f$ has degree $n$.
\item The leading coefficient of $f$ is positive definite.
\end{enumerate}
then $f\in\hyper{1}$ and $gI\lessless f$.
\end{lemma}
\begin{proof}
From the first condition we get that $v^tfv$ sign interlaces
$g$. Since the degree of $g$ is one more than the degree of $v^tfv$,
we conclude that $g \lessless v^tfv$. This implies that
$v^tfv\in\allpoly$, and that $gI\lessless f$.
\end{proof}

\begin{construction}
We can use sign interlacing to construct elements of $\hyper{1}$.
We start with the data
\begin{align*}
  \text{real numbers} & \qquad p_1 < p_2 <\cdots< p_{n+1}\\
\text{positive definite matrices} &\qquad  m_1,\cdots,m_{n+1}
\end{align*}
We claim that the following is in $\hyper{1}$:
\begin{equation}
  \label{eqn:hyper-sign}
 f(x)= \sum_{k=1}^{n+1} \, m_{k} \, 
\prod_{i\ne k}(x-p_i)
\end{equation}
In order to verify that $f\in\hyper{1}$ we see that
\[
  f(p_k) =  \prod_{i\ne k}(p_k-p_i)\, m_k
\]
and the sign of the coefficient of $m_k$ is $(-1)^{n+k+1}$ since the
$p_i$ are increasing.  Since the leading coefficient of $f(x)$ is $
\sum_k m_k $ which is positive definite we can now apply the lemma.

Note that if $\mu=1$ this is exactly Lemma~\ref{lem:sign-quant}
\end{construction}

\begin{construction}
    There is a similar construction for polynomials in
    $\hyper{d}$. The idea is that we construct a polynomial by
    adding together lots of interlacing polynomials, multiplied by
    positive definite matrices. So, we start with
    \begin{align*}
      \text{positive definite $d$ by $d$ matrices}\quad &
      D_1,\dots,D_d \\ 
      \text{positive definite $d$ by $d$ matrices}\quad &
      M_1,\dots,M_d  
    \end{align*}
    Let
    \begin{align*}
      A &= I + x_1D_1 + \cdots + x_d D_d \\
      f(\xx) &= |A|\\
      f_i(\xx) &= \text{determinant of $i$'th principle submatrix of
        $A$}\\
      g(\xx) &= f_1(\xx)M_1 + \cdots +f_d(\xx)M_d
    \end{align*}
    Since $f\lesslesseq f_i$ it follows that $fI\lesslesseq
    f_iM_i$. Addition yields our conclusion: $fI\lessless g$ in $\hyper{d}$.
  \end{construction}

\index{Hadamard product!matrix polynomials} The Hadamard product of
two matrix polynomials $f=\sum A_ix^i$ and $g=\sum B_ix^i$ is $f\ast
g= \sum A_i\ast B_i x^i$ where $A_i\ast B_i$ is the usual Hadamard
product of matrices.  In certain circumstances $\hyper{1}$ is closed
under the Hadamard product $\ast$ with a positive definite matrix.

\begin{lemma}
  Suppose that $f\lessless gI$ where $f\in\hyper{1}$ and
  $g\in\allpoly$. If $P$ is a positive definite matrix then 
  $f\ast P\lessless gI$. In particular, $f\ast P\in\hyper{1}$.
\end{lemma}
\begin{proof}
  If the roots of $g$ are $r_1<\cdots<r_n$ then we know that
  $(-1)^{n+i}f(r_i)$ is positive definite. Now the Hadamard product of
  two positive definite matrices is again positive definite
  \cite{horn-johnson-1}, so 
\[(-1)^{n+i}(f\ast P)(r_i) = (-1)^{n+i}f(r_i)\ast P \] is positive
  definite, as is the leading coefficient of $f\ast P$. The conclusion
  follows.
  \end{proof}

The Hadamard product of three hyperbolic matrix polynomials can be in
$\hyper{1}$. 

  \begin{lemma}
    If $f,g,h\in\hyper{1}$ have the same common interlacing then \\ $f\ast
    g\ast h\in\hyper{1}$.  
  \end{lemma}
  \begin{proof}
    Let $r_1,\dots,r_{n+1}$ be the roots of the common interlacing. We
    know that 
\[
(-1)^{n+i}f(r_i),\quad (-1)^{n+i}g(r_i),\quad(-1)^{n+i}h(r_i)
\]
are all positive definite, and so $(-1)^{n+i}(f\ast g\ast h)(r_i)$ is
positive definite. We now apply Lemma~\ref{lem:hyp-si}.
  \end{proof}

  \section{Root zones}
  \label{sec:root-zones}
  \index{root zones} \index{hyperbolic matrix polynomials!root zones}

Hyperbolic matrix polynomials in one variable are generalizations of
polynomials with all real roots. There are three different notions of
the zeros of a hyperbolic matrix polynomial $f(x)$, and all reduce to the
same idea if the matrix dimension is $1$:

\begin{enumerate}
\item The roots of the diagonal elements.
\item The roots of the determinant of $f$.
\item The roots of all the quadratic forms $v^tfv$.
\end{enumerate}

\noproblem{add diagonal root comment, and nicer determinant description}

  What does the graph of a hyperbolic matrix polynomial $f\in\hyper{1}$
  look like? When $\mu=2$ we can graph $(1,y)f(1,y)^t$ --
  Figure~\ref{fig:hyp-mat-poly} is a typical picture.

  The first point to note is that the graph does not contain large $x$
  values. This can be seen as follows. It suffices to consider $v$
  satisfying $|v|=1$, and that the leading coefficient of $f$ is the
  identity. In this case $v^tfv$ is monic and the coefficients of
  $v^tfv$ are all bounded (by compactness), and so the roots are
  bounded.

  If a polynomial $f\in\hyper{1}$ has a common interlacing $g$, then if we
  graph $f$ and the vertical lines corresponding to the roots of $g$,
  then the strip between any two consecutive such vertical lines
  contains exactly one solution curve of $f$.

The usual terminology is as follows. If $f\in\hyper{1}$ define
$\Delta_k$ to be the set of all $x$ such that there is a $v$ where
$x$ is the $k$'th largest root of $v^tfv$. $\Delta_k$ is called the
$k$'th \emph{root zone}\index{root zone}.

Interlacing of root zones implies interlacing.

\begin{lemma}
  If the root zones of two hyperbolic matrix polynomials alternate
  then the polynomials interlace.
\end{lemma}
\begin{proof}
  Let $r_1,\dots,r_n$ be the root zones of $r$, and $s_1,\dots,s_n$
  the root zones of $s$. The hypotheses mean that the zones are
  arranged as in the diagram

\[
 \underbrace{\hspace*{1cm}}_{r_1}\,f_1\,\underbrace{\hspace*{1cm}}_{s_1}g_1
 \underbrace{\hspace*{1cm}}_{r_2}\,f_2\,\underbrace{\hspace*{1cm}}_{s_2}g_2
\cdots\cdots
g_{n-1} \underbrace{\hspace*{1cm}}_{r_n}\,f_n\,\underbrace{\hspace*{1cm}}_{s_n}g_n
\]
The $f_i$'s and the $g_i$'s are values that separate the various root
zones. Let $f=\prod(x-f_i)$ and $g = \prod(x-g_i)$. The interlacings
\[
r\greateqeq fI\greateqeq s\qquad r,s\greateqeq gI
\]
easily imply that $r\greateqeq s$.
\end{proof}

\begin{example}
  The converse of this lemma is not true -- it fails even for
  diagonal hyperbolic matrix polynomials. Let
\[
f =
\begin{pmatrix}
  (x+32)(x+36) & 0 \\ 0 & (x+14)(x+35)
\end{pmatrix}
\qquad f' =
\begin{pmatrix}
  2(x+34) & 0 \\ 0 & 2(x+24.5)
\end{pmatrix}
\]
We know that $f\lesslesseq f'$ yet their root zones overlap (Figure~\ref{fig:overlap}).

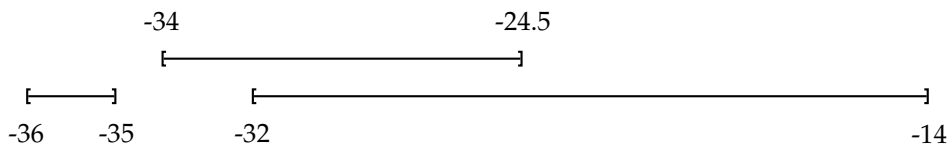
\begin{figure}[h]
  \centering

\begin{pspicture}(0,0)(22,2.5)
\psset{xunit=.6cm,yunit=.5cm}
\psline{[-]}(0,2)(2,2)
\psline{[-]}(5,2)(20,2) 
\psline{[-]}(3,3)(11,3) 
\rput(0,1){-36}
\rput(2,1){-35}
\rput(5,1){-32}
\rput(20,1){-14}
\rput(3,4){-34}
\rput(11,4){-24.5}
\end{pspicture}
  
  \caption{Overlapping root zones}
  \label{fig:overlap}
\end{figure}

\end{example}
\begin{lemma}
The diagonal elements of a hyperbolic polynomial have a common
interlacing.
\end{lemma}
\begin{proof}
  If $f\in\hyper{1}$ then each two by two principle submatrix is in
  $\hyper{1}$, and the two diagonal elements have a common interlacing.
  Consequently, every two diagonal polynomials have a common
  interlacing, so they all have one.
\end{proof}

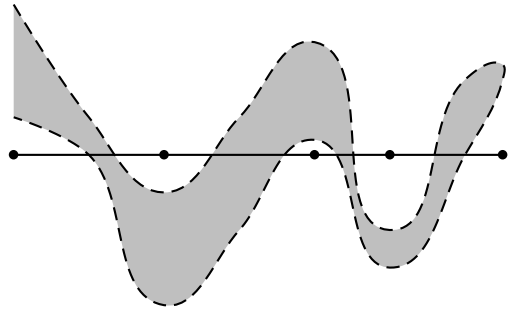
\begin{figure}[htbp]
    \centering
  \begin{pspicture}(0,-2)(7,2)
\pscurve[fillstyle=solid,fillcolor=lightgray,linestyle=dashed]
(0,2)(1,.5)(2,-.5)(3,.5)(4,1.5)(5,-1)(6,1)(6.5,1.2)%
(6.5,1)(6,0)(5,-1.5)(4,.2)(3,-,5)(2,-2)(1,0)(0,.5)
\psline(0,0)(6.5,0)
\psdots(0,0)(2,0)(4,0)(5,0)(6.5,0)
  \end{pspicture}
    \caption{The graphs of $v^tfv$ as $v$ varies}
    \label{fig:vfv}
  \end{figure}

If the degree of a matrix polynomial is greater than two then we do
not have nice pictures such as Figure~\ref{fig:hyp-mat-poly}.
We can make another graph by drawing all of the graphs $v^tfv$ for all
non-zero $v$. The result looks like Figure~\ref{fig:vfv}. 
The points on the $x$ axis show that  all the curves $v^tfv$ have a common
interlacing. If $\lambda_{min}(M)$ and $\lambda_{max}(M)$ are the
smallest and largest eigenvalues of a symmetric matrix $M$ then we
have the inequality

\[
\lambda_{min}(f(x)) \le \frac{v^tfv}{v^tv} \le \lambda_{max}(f(x))
\]
The dashed boundary curves in the figure are the graphs of the maximum
and minimum eigenvalues of $f(x)$ as a function of $x$.

\index{Hadamard product!matrix polynomials}
The Hadamard product preserves $\hyper{1}$ if one of the terms is a
diagonal polynomial. 

\begin{lemma}
  If $f,g\in\hyper{1}$ and $f$ is  diagonal then $f\ast g\in\hyper{1}$.
\end{lemma}
\begin{proof}
  Since the diagonal elements of $f$ and of $g$ have common
  interlacings $F$ and $G$, the result follows from the fact that
  $f_{ii}\lesslesseq F$ and $g_{ii}\lesslesseq G$ implies
  (Lemma~\ref{lem:hadamard-interlace})
\[
f_{ii}\ast g_{ii} \lesslesseq F\ast G
\]
\end{proof}

The roots of the diagonal fall in the root zones. Indeed, if the
diagonal elements of $f$ are $f_i$, and $e_i$ is a coordinate vector,
then $e_i^tfe_i = f_i$. 

The location of the roots of the determinant of a hyperbolic matrix
polynomial fit nicely into the root zone picture.

\begin{lemma}
  Suppose that $f\in\hyper{1}(n)$ has a common interlacing. Then,
  $det(f)$ has $n\mu$ roots, all real. In the interior of each root
  zone all $f(x)$ are positive definite or negative definite. There
  are $\mu$ zeros of the determinant in each root zone.
\end{lemma}
\begin{proof}
  This follows from Lemma~\ref{lem:hyp-sign}, since we know how many
  positive eigenvalues ($\mu$ or $0$) there are in each region between
  the root zones. As $\alpha$ procedes from one root of the common
  interlacing to the next, there is a change from all eigenvalues of
  one sign to all eigenvalues of the opposite sign. Each time an
  eigenvalue changes sign, the determinant must be zero.
\end{proof}

\begin{example}
  Suppose that $\mu=n=3$, and that the roots of $det(f)$ are
  $r_1,\dots,r_9$. Then
\[
  \begin{array}{rcl}
\text{interval } & \text{number of positive roots} & \\
(-\infty,r_1) & 0 \\
(r_1,r_2) &  1 \\
(r_2,r_3) & 2\\
(r_3,r_4) & 3\\
(r_4,r_5) & 2\\
(r_5,r_6) & 1\\
(r_6,r_7) & 0\\
(r_7,r_8) & 1\\
(r_8,r_9) & 2\\
(r_9,\infty) & 3
  \end{array}
\]

We can show that a polynomial is in $\hyper{1}$ by checking if all
roots of the determinant are real, and that the number of positive
eigenvalues in each interval follows the patterns described above.
\end{example}

\section{Quadratics and Newton inequalities}
\label{sec:quadratics}

Quadratic matrix polynomials have been widely studied \cite{lancaster}. We derive a
simple inequality for eigenvalues of the quadratic, and use it in the
usual way to get Newton inequalities for hyperbolic matrix polynomials.

\index{Newton's inequalities!hyperbolic matrix polynomials}

There are necessary and sufficient conditions for a quadratic matrix
polynomials to belong to $\hyper{1}$, but they are not as elementary
as for polynomials.

\begin{lemma}
  Suppose $f(x) = Ax^2 + Bx+C$ where $A$ and $C$ are non-negative
  definite. 
   \begin{enumerate}
   \item If $f\in\hyper{1}$ then $\lambda_{max}(BA^{-1})\,\lambda_{max}(BC^{-1})\ge4$.
   \item If  $\lambda_{min}(BA^{-1})\,\lambda_{min}(BC^{-1})\ge4$ then
     $f\in\hyper{1}$.
   \end{enumerate}
\end{lemma}
\begin{proof}
  If $f\in\hyper{1}$ then 
\[
v^tAv\, x^2 + v^tBv\,x + v^tCv\in\allpoly
\]
Thus the discriminant must be non-negative, and so
\[
(v^tBv)^2 \ge 4 (v^tAv)(v^tCv)
\]
Now $v^tAv$ and $v^tCv$ are non-negative, so this is equivalent to
\[
\frac{v^tBv}{v^tAv}\ge 4\frac{v^tCv}{v^tBv} 
\]
We now apply \eqref{eqn:hyp-ineq-2}.
\end{proof}

\begin{cor}
  Suppose that $f=\sum A_ix^i\in\hyperpos{1}(n)$. For {k=1,\dots,n-1}
  the following inequality holds:
\[
\lambda_{max}(A_kA_{k+1}^{-1})\, \lambda_{max}(A_kA_{k-1}^{-1}) \ge
(1+\frac{1}{k})(1+\frac{1}{n-k}) 
\]
\end{cor}
\begin{proof}
  Since $\hyperpos{1}$ is closed under differentiation and reversal,
  we can follow the usual proof of Newton's inequalities.
\end{proof}

The coefficients of the associated quadratic form are not arbitrary.
Suppose that $f(x)=\sum A_ix^i$ is a matrix polynomial in
$\hyperpos{1}(n)$.  We know that for all $v\ne0$ and $0\le i \le n$
\[
\lambda_i^{min} \le \frac{v^tAv}{v^tv} \le\lambda_i^{max}
\]
Define the two vectors
\begin{align*}
  \Lambda_{min} & = (\lambda_0^{min},\dots,\lambda_n^{min})\\
  \Lambda_{max} & = (\lambda_0^{max},\dots,\lambda_n^{max})
\end{align*}
The inequality above shows that for $|v|=1$ all the polynomials $v^tfv$ are
polynomials of the form $\ppoly(\ccc)$ where $0\le\Lambda_{min}\le \ccc\le
\Lambda_{max}$.  \seepage{lem:polytope-2}.

\section{Newton inequalities in $\hyper{2}$}
\label{sec:coeff-ineq-hyper2}
\index{Newton's inequalities!matrix polynomials}
\index{rhombus inequalities!matrix polynomials}

The coefficients of a quadratic in $\hyper{2}$ satisfy inequalities, and
they lead to inequalities for general polynomials in $\hyper{2}$. We
begin with an example. Consider the polynomial $f(x,y)$ \mypage{ex:hyp-p2}:

\centerline{\xymatrix{
         {\smalltwo{6}{2}{2}{4}}y^2 \ar@{-}[d] \ar@{-}[dr]  & & \\
        {\smalltwo{5}{1}{1}{4}} y \ar@{-}[d] \ar@{-}[r] \ar@{-}[dr]& 
      \ar@{-}[dr]  \ar@{-}[d] {\smalltwo{5}{-1}{-1}{6}} xy & \\
        {\smalltwo{-1}{-3}{-3}{-5}} \ar@{-}[r] & 
        \ar@{-}[r] {\smalltwo{2}{-1}{-1}{3}} x &
        {\smalltwo{1}{-1}{-1}{2}}x^2 }}

If we compute $v^tfv$ then we get a polynomial in $\rupint{2}$. For
instance, if we take $v=(2,3)$ we get

\centerline{\xymatrix{
         84\,y^2 \ar@{-}[d] \ar@{-}[dr]  & & \\
        68\, y \ar@{-}[d] \ar@{-}[r] \ar@{-}[dr]& 
      \ar@{-}[dr]  \ar@{-}[d] 62\, xy & \\
        -85 \ar@{-}[r] & 
        \ar@{-}[r] 23\, x &
        10\,x^2 }}
This polynomial satisfies the rhombus inequalities, so we can use this
information to get inequalities for the original polynomial.
In general, we start with $f(x,y) = \sum A_{i,j}x^iy^j$ in
$\hyperpos{2}$. If $v$ is any vector then
\[
v^tfv = \sum (v^tA_{i,j}v)\,x^iy^j \in\gsubpos_2
\]
The rhombus inequalities (Proposition~\ref{prop:p2plus-inequality})
imply that the following holds for all vectors $v$

\[
(v^t\,A_{i,j+1}\,v)\,(v^t\,A_{i+1,j}\,v) >
(v^t\,A_{i,j}\,v)\,(v^t\,A_{i+1,j+1}\,v)
\]
Using \eqref{eqn:hyp-ineq-2} this implies rhombus inequalities in $\hyper{2}$:
\index{rhombus inequalities!for $\hyper{2}$}

\[
\lambda_{max}(A_{i,j+1}A_{i,j}^{-1}) \ge
\lambda_{min}(A_{i+1,j+1}A_{i+1,j}^{-1})
\]

\section{Singular polynomials}
\label{sec:singular-polynomials}

We now look at some examples of matrix polynomials $f$ for which the
leading coefficient of a matrix polynomial is positive semi-definite
and singular, and $v^tfv\in\allpoly$ for all $v\in\reals^\mu$. Unlike
hyperbolic matrix polynomials, the root zones can be infinite, or all
of them can  even be equal.

\begin{figure}[htbp]
  \centering
  \includegraphics*[width=2in]{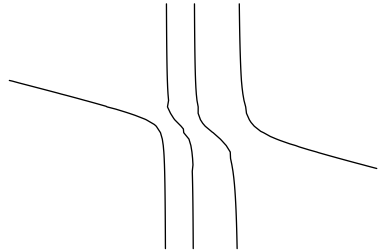}
  \caption{Positive semi-definite leading coefficient}
  \label{fig:hyp-psd}
\end{figure}

First, we consider the matrix $F = \smalltwo{f}{h}{h}{0}$ where
$deg(f)> deg(h)$. If $f$ is monic then the leading coefficient of $F$
is $\smalltwo{1}{0}{0}{0}$, which has eigenvalues $0$ and $1$, and is
positive semi-definite.  The quadratic form is $\alpha^2f+2\alpha\beta
h$ which is in $\allpoly$ for all $\alpha,\beta$ if and only if
$f\lesslesseq h$. Figure~\ref{fig:hyp-psd} is the graph of such a
polynomial.
 
Assume that $f\lesslesseq h$.  If $\epsilon>0$ and $F_\epsilon =
\smalltwo{f}{h}{h}{\epsilon f}$ then $F_\epsilon\in\hyper{1}$ since
its leading coefficient is $\smalltwo{1}{0}{0}{\epsilon}$, and the
quadratic form is a linear combination of $f$ and $h$. Consequently,
$F$ is the limit of polynomials in $\hyper{1}$. Note that $F$ has a
common interlacing, but  the root zones can be unbounded. 

Next, we have an example with no common interlacing. Let
\[
F =
\begin{pmatrix}
  (x+1)(x+3) & 3x+5 \\ 3x+5 & 8
\end{pmatrix}.
\]
The leading coefficient is again $\smalltwo{1}{0}{0}{0}$ and the
quadratic form factors:
\[
(\alpha,\beta)\,F\,(\alpha,\beta)^t =
(\alpha(x+1)+2\beta)(\alpha(x+3)+4\beta).
\]
Consequently, for every $x_0$ there are two $v$ such
that the quadratic form $v^tF(x_0)v$ is zero. Thus, the two root zones
are both equal to $\reals$, there is no common interlacing, and so $F$
is not a limit of polynomials in $\hyper{1}$.

\section{Matrix polynomials with determinants in $\gsubclose_d$}
\added{7/8/7}
\label{sec:matrix-polynomials-1}

In this last section we consider the properties the set of matrix
polynomials with determinants in $\gsubclose_d$. Define

\[
\calM = \bigl\{ \text{matrix polynomial M}\,\mid\, |\text{M}|\in\gsubclose_d
\text{ for some $d$}\bigr\}
\]

If $y$ is a variable not occuring in $A$ or $B$ then we say that
$A\matint B$ if and only if $A + y B\in\calM$.
  
If the matrices are $1$ by $1$ then this is just the usual definition
of interlacing. One interesting question is what properties of
$\gsubclose_d$ extend to $\calM$? Clearly $\calM$ is closed under
substitution of real values for variables, and of $x_i+a$ for $x_i$.

Here are some basic interlacing properties:

\begin{lemma}\label{lem:mp-1}
  If $A\matint B$ then
  \begin{enumerate}
  \item $B\matint -A$
  \item $A + \alpha\,B\matint B$ for $\alpha\in\reals$.
  \item $A \matint B + \alpha\, A$ for $\alpha\in\reals$.
  \item If $C\in\calM$ then $AC\matint BC$ and $CA\matint CB$. 
  \end{enumerate}
\end{lemma}

\begin{proof}
If $|A+yB|\in\gsubclose_d$ then reversing just $y$ yields
$|yA-B|\in\gsubclose_d$. Next, replacing $y$ by $y+\alpha$ shows that
$|A+(y+\alpha)B|\in\gsubclose_d$. From $B\matint -A$ we get that
$B+\alpha(-A) \matint -A$, and applying the first part again gives $A \matint B +
(-\alpha) A$. For the next one, note that $|AC+yBC| = |C||A+yB|$.  
\end{proof}

 Some constructions:

\begin{lemma} \label{lem:mp-2} \ 
\begin{enumerate}
\item  If $A$ is symmetric, all $D_i$ and $B$ are positive semi-definite, then \[A+\sum x_iD_i
  \matint B.\]
\item If $Z = diag(z_1,\dots,z_d)$ and $A$ is positive semi-definite then
  \begin{enumerate}
  \item $Z+A \in\calM$
  \item $Z-A \in\calM$
  \item $AZ+I\in\calM$
  \item $ZA + I\in\calM$
  \item If $W=diag(w_i)$ and $O$ is orthogonal then $OZ+WO\in\calM$,
  \end{enumerate}
\end{enumerate}
\end{lemma}

\begin{proof}
  The first one is the statement that $|A + y B +\sum x_i D_i
  |\in\gsubclose_d$. For the second, reverse each variable separately,
  and then reverse all at once. The next follows from the fact that $A+Z = A +
  \sum z_iD_i$ where $D_i$ is the positive semi-definite matrix that
  is all zero except for a $1$ in location $(i,i)$. If we reverse all
  variables we have
\[
z_1\cdots z_d |A + Z^{-1}| = |AZ+I| = |ZA+I|.
\]
For the last one, $Z + O^{-1}WO$ is a sum of variables times positive
semi-definite matrices.
\end{proof}

If $a<b$ then $x+a \greateq x+b$. The corresponding result for $\calM$
is

\begin{lemma} If $A<B$ are positive definite and $Z = diag(z_i)$ then
  \[ A+Z \matint B+Z. 
\]
  
\end{lemma}
\begin{proof}
  Since $B-A$ is positive definite, Lemma~\ref{lem:mp-2} implies that
  \begin{align*}
    A+Z & \matint B-A. \\
    \intertext{Applying Lemma~\ref{lem:mp-1} yields the conclusion}
    A+Z & \matint B-A + (A+Z) = B+Z.
  \end{align*}
\end{proof}

\begin{cor}
  With the same hypotheses, $BZ+I \matint AZ+I$.
\end{cor}
\begin{proof}
  Reverse all variables.
\end{proof}

As \cite{bbl} observed, it is an interesting question if $A+Z\matint
B+Z$ implies $A\le B$.



\part{Polynomials with complex roots}

\chapter{Polynomials non-vanishing in a region}
\label{cha:nv}
In this chapter we introduce polynomials that are non-vanishing (NV) in a
region $\calD$. We show how geometric properties of $\calD$ imply
properties of the sets of polynomials such as closure under
differentiation, or under reversal. 

  \begin{definition}
    Suppose that $\mathcal{D}\subsetneq\complexes$.  Define
    
    \begin{align*}
        \nv{d}{\calD} &= \bigl\{ f\,\mid\, \text{$f$ has complex coefficients
          and } f(\sigma)\ne0 \text{ for all }
        \sigma\in \calD\bigr\}\\
        \nvr{d}{\calD} &= \bigl\{ f\,\mid\, \text{$f$ has real coefficients
          and } f(\sigma)\ne0 \text{ for all }
        \sigma\in \calD\bigr\}
    \end{align*}
  \end{definition}

Here are some examples of the spaces we will consider.

\index{Hurwitz polynomial}
\index{right half plane polynomial}
\index{upper half plane polynomial}
\index{real upper half plane polynomial}

\begin{example}
If $\calD$ is the product of $d$ open upper half planes, then
a polynomial in $\nv{d}{\calD}$ is called a \emph{upper half plane polynomial}.
It is a  \emph{real upper half-plane polynomial} if it is in
$\nvr{d}{\calD}$. 
The spaces of such polynomials are denoted $\hb{d}$ and $\rhb{d}$.
Many results about these polynomials can be found in
\cite{branden-hpp}, \cite{bbs} and Chapter~\ref{cha:complex-coef}.
\end{example}

\begin{example}

  If $\calD$ is a product of $d$ open right half planes then a
  polynomial in $\nv{d}{\calD}$ is called a \emph{right half-plane
    polynomial}, or a \emph{Hurwitz polynomial}. $\nv{d}{\calD}$ is denoted
  $\stabledc{d}$ and  $\rnv{d}{\calD}$ is denoted $\stabled{d}$. In the
  $\stabledc{1}$ case, a right half-plane polynomial has all roots in
  the closed left half plane. Such polynomials are usually called
  stable, or Hurwitz-stable. More information about these polynomials
  can be found in \cite{sokal-wagner} and Chapter~\ref{cha:stable}.
\end{example}

\begin{example}
  If $\calS$ is a vertical strip such as $\{z\,\mid\, -1< \Re(z)<1\}$
  then we consider $\calD=\complexes\setminus\calS$. The set
  $\nv{1}{\complexes\setminus\calS}$ consists of polynomials whose zeros
  all lie in the strip $\calS$. More information about them can be
  found in Section~\ref{sec:prop-polyn-with}.
\end{example}

\begin{example}
  The open unit disk $\Delta$ in $\complexes$ has two generalizations to
  $\complexes^d$. The first is the  open ball  $\Delta_d$ of radius
  $1$, and the second is the product $\Delta\times \cdots \times
  \Delta$. The spaces are considered in 
  Section~\ref{sec:open-unit-ball}.  
\end{example}

\begin{example}\label{ex:cone-12}
  Define the sector $\sector{n} =\left\{ z\mid -\pi/n < \arg(z) <
    \pi/n\right\}$.  For instance, $\sector{1}$ is
  $\complexes\setminus(-\infty,0)$, and $\sector{2}$ is the open right
  half plane. We consider properties of product of sectors
  $\nv{d}{\prod_i(\sector{n_i})}$ in  various examples in this
  chapter and Section~\ref{sec:polyn-non-vanish}.
\end{example}

\section{Fundamental properties}
\label{sec:fund-prop}

In order to derive interesting properties we we usually assume that
our domain is an open cone $\calC$ in $\complexes^d$. This means $\calC$ is
closed under addition and multiplication by positive constants. We
have a few properties that hold for general domains. The proofs are
immediate from the definition.

\begin{lemma}\ 

  \begin{enumerate}
  \item   $f(\xx)g(\xx)\in\nv{d}{\calD}$ if and only if
  $f(\xx)\in\nv{d}{\calD}$ and $g(\xx)\in\nv{d}{\calD}$.
\item Suppose $f(\xx,\yy)\in\nv{d+e}{\calD\times\calE}$ where
  $\calD\subset\complexes^d$ and $\calE\subset\complexes^e$. If
  $\alpha\in\calE$ then $f(\xx,\alpha)\in\nv{d}{\calD}$.
  \end{enumerate}
\end{lemma}

The following classical theorem of Hurwitz also holds for analytic
functions in $d$ variables.

\begin{theorem}[Hurwitz]
  Let $(f_n)$ be a sequence of functions which are all analytic and
  without zeros in a region $\Omega$. Suppose, in addition, that
  $f_n(z)$ tends to $f(z)$, uniformly on every compact subset of
  $\Omega$. Then $f(z)$ is either identically zero or never equal to
  zero in $\Omega$.
\end{theorem}

The next lemma is an immediate consequence of the general Hurwitz's theorem:

\begin{lemma}[Polynomial closure]\label{lem:nv-hur} If $\calD$ is
  open and $\{f_i\}$ is a sequence of polynomials in $\nv{d}{\calD}$
  that converge uniformly on compact subsets to a polynomial $g$ then
  either $g$ is zero or is in $\nv{d}{\calD}$.
\end{lemma}

The next lemma gives several simple ways to construct new polynomials
in $\nv{d}{\cone}$ from old ones. The proofs easily follow from the
definition. 

\begin{lemma}
  Assume that $\cone$ is a cone.
  \begin{enumerate}

\item If $f(\xx)\in\nv{d}{\cone}$ then
  $f(\xx+\yy)\in\nv{2d}{\cone\times\cone}$ where
  $\xx=(x_1,\dots,x_d)$ and $\yy=(y_1,\dots,y_d)$.
\item If $a$ is a positive constant and
  $f(\xx)\in\nv{d}{\cone}$ then $f(a\,x_1,\dots,a\,x_d)\in \nv{d}{\cone}$.
\item If $a_i$ are  positive constants, $\cone_i$ are cones, and 
  $f(\xx)\in\nv{d}{\prod\cone_i}$ then $f(a_1\,x_1,\dots,a_d\,x_d)\in \nv{d}{\prod\cone_i}$.
  \end{enumerate}
  
\end{lemma}

The homogeneous part $f^H$ of a polynomial consists of  the
terms of maximum total degree.

\begin{cor}\label{cor:hbc-homog}
If $f\in\hbc{d}{\cone}$ then $f^H\in\nv{d}{\cone}$.   
\end{cor}
\begin{proof}
  This follows from the above, the fact that $f^H$ is never zero
  and
\[f^H(x) = \lim_{r\rightarrow\infty} \frac{f(r\xx)}{r^n}
\]
where $n$ is the degree of $f^H$.

\end{proof}

The next result refines the earlier substitution result.

  \begin{lemma}[Substitution]\label{lem:hbc-1a}
    Suppose that $\calD$ is open and $\overline{\calD}$ is the closure of $\calD$.  If
    $f(\xx,\yy)\in\nv{d+e}{\calD\times\calE}$ and
    $\alpha\in\overline{\calD}$ then $f(\alpha,\yy)\in\nv{e}{\calE}$
    or $f(\alpha,\yy)=0$.
  \end{lemma}
  \begin{proof}
If $\alpha\in\calD$ then $f(\alpha,\yy)\in\nv{d}{\calE}$. If
$\alpha\in\overline{\calD}$ then we approach $f(\alpha,\yy)$ by sequences
in $\nv{d}{\calD\times\calE}$, so we can apply Lemma~\ref{lem:nv-hur}.
  \end{proof}

The constant in the next lemma only depends on $\cone$ and is called
the \emph{reversal constant} \index{reversal constant}.

\begin{lemma}[Reversal]\label{lem:hbc-9}Assume
  $\cone\subset\complexes$ is a cone and   suppose that $\sum x^i
  f_i(\yy)\in\nv{d+1}{\cone\times\calD}$. There is a unique
  $\rho\in\complexes$ with $|\rho|=1$ such that
\[
 \sum (\rho x)^{n-i} \,f_i(\yy) \in\nv{d+1}{\cone\times\calD}.
\]
\end{lemma}
\begin{proof}
  If $\cone^{-1} = \{z^{-1}\,\mid\, z\in\cone_1\}$ then
  $\cone^{-1}$ is also a cone, whose elements are the conjugates of
  $\cone$. Since $\cone$ and $\cone^{-1}$ only differ by a
  rotation there is a unique rotation $\rho\in\complexes$ so that
  $\rho\,\cone^{-1}=\cone$. 
If $(\sigma_1,\dots,\sigma_{d+1})\in\cone\times\calD$ then 
\[
 \sum (\rho \sigma_1)^{n-i} \,f_i(\sigma_2,\dots,\sigma_{d+1}) =
(\rho\sigma_1)^n \sum \bigl(\frac{\rho}{\sigma_1}\bigr)^i
\,f_i(\sigma_2,\dots,\sigma_{d+1}) 
\]
and this is non-zero since $\rho/\sigma_1\in\cone$. See Table~\ref{tab:cones}. 
\end{proof}

\begin{table}
  \centering
  \begin{tabular}{rrrc}
    \toprule
    $\cone$ & $\cone^{-1}$ & $\rho$ & $\cone/\cone$ \\
    \midrule
    $\uhp$ & $-\uhp$ & $-1$ & $\complexes\setminus(-\infty,0)$\\
    $\rhp$ & $\rhp$ & $1$ & $\complexes\setminus(-\infty,0)$\\
    $\quada$ & $\quadd$ & $\imag$ & $\rhp$\\
    $\sector{n}$ & $\sector{n}$ & 1 & $\sector{2n}$\\
    \bottomrule
  \end{tabular}
  \caption{Properties of cones}
  \label{tab:cones}
\end{table}

We now consider when a linear polynomial is in $\nv{d}{\calD}$.
The following are elementary.

\begin{lemma}
  The following are equivalent, where $\cone_i\subset\complexes$ are cones.
  \begin{enumerate}
  \item $\sum_1^d a_i\,x_i \in \nv{d}{\prod \cone_i}$.
  \item $0 \notin \sum_1^d a_i \cone_i$
  \item $\sum_1^d a_i\,\cone_i$ is contained in a half plane.
  \end{enumerate}
\end{lemma}

\begin{example}
  Suppose $\cone = \sector{2}\times\sector{4}$ and consider $x
  + \alpha y$. In order for $\sector{2}+ \alpha \sector{4}$ to
  remain in a half plane, we must have that $\alpha
  \sector{4}\subset\sector{2}$. Thus we conclude that
\[ x + \alpha y \in\nv{2}{\sector{2}\times\sector{4}} \Longleftrightarrow -\pi/4\le \alpha
\le \pi/4
\]

  If $r>0$ and $\cone\subset\complexes$ then $x+ry\in\nv{2}{\cone^2}$, and
  hence $\rho xy + r$ and $xy + r/\rho$ are in $\nv{2}{\cone^2}$.

\end{example}

There's a simple condition if all the planes are half planes.

\begin{lemma}
If all cones $\cone_i$ are the same half plane $\cone$, then the following are
equivalent: 
  \begin{enumerate}
  \item $\sum_1^d a_i\,x_i \in \nv{d}{\cone^d}$.
  \item All non-zero $a_i$ have the same argument.
  \end{enumerate}
\end{lemma}
\begin{proof}
If $\sum \alpha_i x_i\in\nv{d}{\cone^d}$ then we may assume that
$\alpha_1$ is non-zero, so we may divide through and assume
$\alpha_1=1$. We now need to show that all non-zero $\alpha_i$ are
positive. If some $\alpha_i$ is not positive, then $\cone$ and
$\alpha\cone$ are different, and so their sum equals
$\complexes$. The converse is immediate.
\end{proof}

The corollary is simple and the proof is omitted.

\begin{cor}
  If $\cone$ is a half plane then
  $\beta+x_1+\alpha_2x_2+\cdots+\alpha_dx_d\in\nv{d}{\cone^d}$ if and
  only if 
  \begin{enumerate}
  \item All $\alpha_i$ are non-negative.
  \item $\beta$ is in the closure of $\cone$. 
  \end{enumerate}
\end{cor}

Following \cites{sokal-wagner,bbs} we have

\begin{lemma}\label{lem:nv-14}
  Suppose $\cone\subset\complexes$ and $f(\xx)$ is homogeneous. The following
  are equivalent
  \begin{enumerate}
  \item $f(\xx)\in \nv{d}{\cone^d}$
  \item $f(\aaa+t\bbb)\in\nv{1}{\cone/\cone}$ where $\aaa,\bbb>0$.
  \end{enumerate}
\end{lemma}
\begin{proof}
  $\cone/\cone$ denotes the set $\{\alpha/\beta\,\mid\,
  \alpha,\beta\in\cone\}$.  If $\alpha/\beta\in \cone/\cone$ where
  $\alpha,\beta\in\cone$ and $f$ has degree $n$ then by homogeneity
\[
f(\aaa+(\alpha/\beta)\bbb) = \beta^n \, f(\beta\,\aaa + \alpha\,\bbb)
\ne 0
\]
since $\beta\,\aaa + \alpha\,\bbb\in\cone^d$. Conversely, let
$\sigma_1,\dots,\sigma_d\in\cone$. Since $\cone$ is a cone there are
$\alpha,\beta\in \cone$ and $\aaa,\bbb>0$ such that 
$\beta\,\aaa + \alpha\,\bbb = (\sigma_1,\dots,\sigma_d)$. Thus
\[
f(\sigma_1,\dots,\sigma_d) = f(\beta\,\aaa + \alpha\,\bbb) =
\beta^n\,f(\aaa+\bbb(\alpha/\beta))\ne0
\].
\end{proof}

\section{Differentiation}
\label{sec:nv-diff}

We can not differentiate with respect to a variable and stay in
the same class, but we can if the corresponding factor is a half plane.

  \begin{lemma}\label{lem:nv-15}
    If $f(x,\yy)\in\nv{d+1}{\cone\times\calD}$ and $\cone$ is a  half plane then
    $\frac{\partial}{\partial x}f(x,\yy)\in\nv{d+1}{\cone\times\calD}$.
  \end{lemma}
  \begin{proof}
Suppose that $(\alpha_1,\dots,\alpha_{d+1})\in\cone\times\calD$. Since
\[
\bigl(\frac{\partial f}{\partial x}\bigr)
(\alpha_1,\dots,\alpha_{d+1}) =
\biggl[ \frac{\partial}{\partial x}\,
(x,\alpha_2,\dots,\alpha_{d+1})\biggr](\alpha_1)
\]
it suffices to assume that $d=1$. So, if $g\in\nv{1}{\cone}$ then
  all roots lie in the complement of $\cone$. Now $\cone$ is a
  half plane by hypothesis, so its complement is  a closed half plane,
  and in particular is convex. By Gauss-Lucas the roots of $g'$ lie in
  the convex hull of the roots of $g$, and hence lie in the complement
  of $\cone$. Consequently, $g'\in\nv{1}{\cone}$.
  \end{proof}

  \begin{lemma}\label{lem:nv-16}If $\cone$ is a  half plane,
    $f(x,\yy) = \sum x^i f_i(\yy)\in\nv{d+1}{\cone\times\calD}$ then
    either $f_i(y)=0$ or
\[
f_i(\yy) \in\nv{d}{\calD}\quad
\text{for $i=0,1,\dots$}
\]
 \end{lemma}
 \begin{proof}
     Differentiate with respect to $x$ until the desired coefficient is the constant term,
  and then substitute zero for $x$.
 \end{proof}

 \begin{example}
   The assumption that $\cone$ is a half plane is essential. We give
   an example of a polynomial in $\nv{2}{\sector{4}}$ whose
   coefficients are not all in $\nv{1}{\sector{4}}$.  Let
\[ f(x) = (x-1)^2 + 4  = (x-1-2\imag)(x-1+2\imag). \]
Clearly $f(x)\in\nv{1}{\sector{4}}$, but $f'(x) =
2(x-1)\not\in\nv{1}{\sector{4}}$. Now consider 
\[ f(x+y) = (x-1)^2 + 4 + 2y(x-1) + y^2 \]
We know $f(x+y)\in\nv{2}{\sector{4}}$, and the coefficient of $y^0$
in $f(x+y)$ is in $\nv{1}{\sector{4}}$ but the coefficient of $y^1$
is not.
 \end{example}

\begin{lemma}\label{lem:nv-18}
  If $\sum x^i f_i(\yy) \in \nv{d+1}{\cone\times\calD}$ and $\cone$ is a half
  plane then
\[
f_i(\yy) + x f_{i+1}(\yy) \in\nv{d+1}{\cone\times\calD}\quad\text{for $i=0,1,\dots$}
\]
\end{lemma}
\begin{proof}
Because $\cone$ is a half plane we can differentiate and
reverse, so all the polynomials below are in $\nv{d+1}{\cone\times\calD}$.
\begin{align*}
  \text{reversing}\quad & \sum (\rho x)^{n-i}\,f_i(\yy) \\
  \text{ $n-k-1$ differentiations}\quad &
\sum (\rho )^{n-i}\, x^{k+1-i} \falling{n-i}{n-k-1}\,f_i(\yy) \\
\text{reverse} \quad &
\sum (\rho )^{n-i}\, (\rho x)^{i} \falling{n-i}{n-k-1}\,f_i(\yy) \\
\text{ $k$ differentiations}\quad &
\sum (\rho )^{n-i}\, (\rho x)^{i-k}
\falling{n-i}{n-k-1}\falling{i}{k}\,f_i(\yy) \\
\text{which equals}\quad  &
   \rho^n\bigl( f_k(\yy) (n-k)!k! + x \,f_{i+1}(\yy) (n-k-1)!(k+1)!\bigr)
\end{align*}
Dividing and rescaling yields the result.
\end{proof}

\begin{cor}\label{cor:nv-19}
  If $f(\xx)\in\nv{d+1}{\cone\times\calD}$ and $\cone$ is a half plane then
\begin{gather*}
f(\xx) + y \frac{\partial f}{\partial x_1} \in
\nv{d+2}{\cone\times\cone\times\calD}\\
\rho\,y\,f(\xx) +  \frac{\partial f}{\partial x_1} \in
\nv{d+2}{\cone\times\cone\times\calD}
\end{gather*}
where $\rho$ is the reversal constant.
\end{cor}
\begin{proof}
  Since $f(x_1+y,x_2,\dots,x_d)\in\nv{d+2}{\cone\times\cone\times\calD}$ we
  can apply the lemma to the Taylor series
\[
f(x_1+y,\cdots,x_d) = f(\xx) + \frac{\partial f}{\partial x_1} y +
\cdots
\]
The second statement follows by reversing the first.
\end{proof}

We have seen that coefficients of the homogeneous part of a polynomial
in $\gsubclose_d$ all have the same sign. A similar result holds for
$\nv{d}{\cone}$: if we multiply by the appropriate scalar then all
the coefficients are positive:

\begin{lemma}\label{lem:nv-fh}
  If $\cone$ is a half plane and $f\in\nv{d}{\cone^d}$ then all
  non-zero coefficients of $f^H$ have the same argument.
\end{lemma}

\begin{proof}  
  We first assume that $f$ is homogeneous and use induction on $d$.
  For $d=2$ we note that $x+a y\in\nv{2}{\cone^2}$ if and only if
  $a\ge0$. Thus if $f\in\nv{2}{\cone^2}$ is homogeneous then
\[ f = \alpha \prod_i (x+a_i y)\ \text{where all $a_i\ge0$} \]
so all non-zero coefficients have argument $\arg(\alpha)$.

If $f(\xx,y)\in\nv{d+1}{\cone^{d+1}}$ is homogeneous then we can write
\[
f(\xx,y) = \sum_i y^i\, f_i(\xx)
\]
where all $f_i\in\nv{d}{\cone^d}$ and are homogeneous.
By induction all coefficients of $f_i$ have the same argument.

By Lemma~\ref{lem:nv-14} with $\aaa=(1,\dots,1,0)$ and
$\bbb=(0,\dots,0,1)$ we have that 
\[
f(1\dots,1,t) = \sum t^i\, f_i(1,\dots,1) \in \allpolypos
\]
since $\cone/\cone=\complexes\setminus(-\infty,0)$.
Now the argument of $f_i(1,\dots,1)$ equals the argument of the
coefficients of $f_i$, and since $f(1,\dots,1,t)\in\allpolypos$ all
coefficients have the same argument.

If $f$ is not homogeneous and has degree $n$ then the conclusion
follows from  Corollary~\ref{cor:hbc-homog}.
\end{proof}

\section{Interlacing}
\label{sec:nv-interlacing}

Interlacing was initially defined by the positions of roots, and then
by linearity. Now we define it algebraically.

\begin{definition} Suppose $f,g\in\nv{d}{\calD}$. 
  We say that \emph{$g$ interlaces $f$ in $\cone$}, written
  $f\nvlace g$, if and only if 
\[f(\xx) + y g(\xx)\in\nv{d+1}{\calD\times\cone}\]
\end{definition}

\begin{remark}
  Note that if $\cone'\subset\cone$ then interlacing in $\cone$
  implies interlacing in $\cone'$, but the converse is false. For
  example, we claim that
\[ x^2 + y(ax+b) \in\nv{2}{\sector{4}\times\sector{4}}\quad\text{if
  $a,b>0$}
\]
This is clear, since both the image of $\sector{4}$ by $x^2$ and
$y(ax+b)$ lies in $\sector{2}$. However, we also have that $f(x,y) =
x^2+y(x+1)$ is not in $\nv{2}{\sector{2}\times\sector{4}}$. To see
this, we need to exhibit one zero in $\sector{2}\times\sector{4}$:
\[
f(\,(1+2\imag)/8\,,\, (19-42\imag)/680\,) = 0
\]
However, it's not had to check that
$x^2+y(x+\imag)\in\nv{2}{\sector{2}\times\sector{4}}$, and therefore is also
  in $\nv{2}{\sector{4}\times\sector{4}}$
\end{remark}

We first have some properties of interlacing that do not require one
of the cones to be a half plane.

\begin{lemma} Assume that $f,g,h\in\nv{d}{\calD}$, and $\cone$ is a
  cone.

  \begin{enumerate}
  \item $fg \nvlace fh$ in $\cone$ iff $g\nvlace h$ in
    $\cone$.
  \item If $-1\not\in\cone$ then $f\nvlace f$ in $\cone$.
  \item If $\rho$ is the reversing constant for $\cone$ then
$f\nvlace g$ in $\cone$ implies $\rho g \nvlace f$ in $\cone$.
\item If $0\not\in\cone_1+\cone_2$ and
  $f\in\nv{2}{\cone_1\times\cone_2}$ then $xf\nvlace f$ in
    $\cone_2$. 
  \end{enumerate}
  
\end{lemma}
\begin{proof}
The first is trivial. To check if $f\nvlace f$ we need to see
that $f(\xx)+y\,f(\xx) = f(\xx)\,(y+1)$ is non-zero. Now $f(\xx)\ne0$,
and $y+1\ne0$ by hypothesis. For the next one, $f+y
g\in\nv{d+1}{\cone\times\cone_0}$ implies that  $\rho g+
f\in\nv{d+1}{\cone\times\cone_0}$. 

Finally, to check $xf\nvlace f$ we check that $xf+yf \in\nv{2}{\cone_1\times\cone_2}$.
Simply note that $(x+y)f\ne0$ by hypothesis.
\end{proof}

Most interesting results about interlacing require a half plane.  As
always, we expect adjacent coefficients as well as the derivative to
interlace. The corollary is a consequence of Lemma~\ref{lem:nv-18} and
Corollary~\ref{cor:nv-19}.

\begin{cor} 
  Suppose $\cone$ is a half plane and $f = \sum x^i f_i(\yy) \in\nv{d}{\cone\times\calD}$.
  \begin{enumerate}
  \item   $f_i \nvlace f_{i+1}$\ for $i=0,1,\dots$.
  \item   $f \nvlace \frac{\partial f}{\partial x}$.
  \end{enumerate}
\end{cor}

\begin{lemma}\label{lem:nv-24}
  If $\cone$ is a half plane and $f,g,h,k\in\nv{d}{\calD}$ then (all
  interlacings in $\cone$)
  \begin{enumerate}
  \item $f\nvlace g$ and $f\nvlace h$ implies
    $f\nvlace ag+bh$ for $a,b>0$.
  \item $f\nvlace g$ and $h\nvlace k$ implies
    $fh\nvlace fk+gh \nvlace gk$.
  \item $g\nvlace f$ and $h\nvlace f$ implies
    $ag+bh\nvlace f$ for $a,b>0$.
  \item If $f\nvlace g \nvlace h$ then $f+ \rho h
    \nvlace g$.
  \end{enumerate}
\end{lemma}
\begin{proof}
All these results follow using the interlacing of adjacent
coefficients. For the first we have
\[ (f+yg)(f+yh) = f^2 + yf(g+h) + y^2 gh\]
so $f^2 \nvlace f(g+h)$ and therefore $f\nvlace
g+h$. The second is similar. For the third we use
\[
(g+y\rho f)(g+y\rho h) = g^2 + gy\rho (f+h)
\]
so $g\nvlace \rho(f+h)$ and reversing again yields
$f+h\nvlace g$. The last follows from the third since $ \rho
h\nvlace g$. 
\end{proof}

Next is the recurrence for orthogonal-type polynomials.

\begin{lemma}
Assume $\cone$ is a half plane and $-1\not\in\cone$. If
$f_0,f_1\in\nv{d}{\cone\times\calD}$ then
\begin{align*}
  f_1 & \nvlace f_0 \\
  f_n &= (a_nx_1+b_n) f_{n-1} + \rho\, c_n \,f_{n-2} \quad\text{where $a_n,b_n,c_n>0$.}
\end{align*}
{then} $f_n  \nvlace f_{n-1}$ in $\cone$.
\end{lemma}
\begin{proof}
  By induction we may assume that $f_{n-1}\nvlace f_{n-2}$.
  Since $-1\not\in\cone$ we know $b_n f_{n-1}\nvlace
  f_{n-1}$. Since $\cone+\cone$ doesn't contain zero we know $x_1
  f_{n-1}\nvlace f_{n-1}$, and thus $ (a_n
  x_1+b_n)f_{n-1}\nvlace f_{n-1}$. By reversal of
  $f_{n-1}\nvlace f_{n-2}$ we have $\rho f_{n-2}\nvlace
  f_{n-1}$, so addition gives the conclusion.
\end{proof}

\begin{example}
  We can use differentiation to derive properties of $\nv{1}{\sector{4}}$ even
  though $\nv{1}{\sector{4}}$ isn't closed under differentiation. Since
  $x+y^2$ and $ x+(y+1)^2$ are in
  $\nv{2}{\sector{2}\times\sector{4}}$ the coefficients of powers of
  $x$ in $\bigl[x+y^2\bigr]^2\, \bigl[x+(y+1)^2\bigr]^2$ interlace.
  Thus all the following polynomials are in $\nv{1}{\sector{4}}$ and
  interlace in $\sector{2}$.
\begin{gather*}
    y^4 (y+1)^4  \nvlace  2 y^2 (y+1)^2 \left(2 y^2+2
      y+1\right)  \nvlace \\
    6 y^4+12 y^3+10 y^2+4 y+1  \nvlace  2 \left(2 y^2+2
   y+1\right)  \nvlace 1
\end{gather*}

\end{example}

We can use the example above to give a general construction.

\begin{lemma}
  If $f,g\in\nv{1}{\sector{4}}$ and $f\nvlace g$ in
  $\sector{2}$ then
\[ y^2f(y) \nvlace f(y) + y^2g(y) \nvlace g(y)
\quad\text{in $\sector{2}$}
\]
\end{lemma}
\begin{proof}
  Since $x+y^2$ and $f(y) + x g(y)$ are in
  $\nv{2}{\sector{2}\times\sector{4}}$ their product is also. Thus
\[
y^2 f(y) + x\bigl[ f(y) +y^2 g(y)\bigr] + x^2 \big[ g(y)\bigr]
\quad\in\nv{2}{\sector{2}\times\sector{4}}. 
\]
\end{proof}

\begin{remark}
  It follows from the lemma that the  two matrices
\[
\begin{pmatrix}
  1 & y^2 \\ 0 & 1
\end{pmatrix}\qquad\qquad
\begin{pmatrix}
  y^2 & 0 \\ 1 & y^2
\end{pmatrix}
\]
map pairs of polynomials in $\nv{1}{\sector{4}}$ that interlace in
$\sector{2}$ to another pair of polynomials in $\nv{1}{\sector{4}}$
that interlace in $\sector{2}$
\end{remark}

\section{Some determinant properties}
\label{sec:some-matr-prop}
In this section we establish some properties of determinants we need to
construct polynomials. 

If $A$ and $B$
 are anti-diagonal matrices with anti-diagonal
entries $(a_i)$ and $(b_i)$ then 

  \begin{equation}
    \label{eqn:skew-det-iden}
    \begin{vmatrix}
      A & I \\ -I & B
    \end{vmatrix} 
= (a_1b_n+1)\cdots(a_nb_1+1)
  \end{equation}
The proof is an easy induction.
For instance, if $n=4$ then 
\[
\begin{vmatrix}
 . & . & . & a_1 & 1 & . & . & . \\
 . & . & a_2 & . & . & 1 & . & . \\
 . &  a_3 & . & . & . & . & 1 & . \\
 a_4 & . & . & . & . & . & . & 1 \\
 -1 & . & . & . & . & . & . & b_1 \\
 . & -1 & . & . & . & . & b_2 & . \\
 . & . & -1 & . & . & b_3 & . & . \\
 . & . & . & -1 & b_4 & . & . & .
\end{vmatrix}
=(a_1b_4+1)(a_2b_3+1)(a_3b_2+1)(a_4b_1+1)
\]
where the dots "." are zeros.  A similar argument shows that if $I$ is
the $n$ by $n$ identity matrix then $
\begin{vmatrix}
  I & \imag I \\ \imag I & I
\end{vmatrix}=2^n
$.

\begin{lemma}\label{lem:skew-non-zero}
  If $A$ is a skew symmetric matrix, $B$ is positive definite, and $S$
  is symmetric then
\[ | A+B+\imag S|\ne0\]
\end{lemma}
\begin{proof}
  $A+B+\imag S\ne0$ since $B$ is positive definite. 
  If $|A+B+\imag S|=0$ then there are non-zero real $u,v$ such that 
  \begin{align*}
    (A+B+\imag S)(u+\imag v)& =0 \\
\intertext{Now $A,B,S,u,v$ are real, so we can separate real and
  imaginary parts:}
  Au + Bu - Sv &= 0 \\
  Av+Bv+Su &= 0\\
  \intertext{The inner product $<r,s>=r^ts$ satisfies $<Ar,r>=0$ since
    $A$ is skew-symmetric. Thus} 
  <Bu,u> - <Sv,u> &= 0 \\
  <Bv,v> + <Su,v>&=0 \\
\intertext{$S$ is symmetric, so adding the last two equations yields}
  <Bu,u> + <Bv,v> &= 0
  \end{align*}
However, $B$ is positive definite and $<Br,r>$ is positive  unless $r=0$. Thus
$u=v=0$, a contradiction. 
\end{proof}

\begin{cor}\label{lem:nv-not-zero}
  If $S,B$ are symmetric and either one is positive definite then
  $|S+\imag B|$ is not zero.
\end{cor}

\begin{proof}
  Take $S=0$. If necessary, multiply by $\imag$. 
\end{proof}

\begin{lemma}\label{lem:abba}
   \[
\begin{vmatrix}
  A & B  \\ -B & A
\end{vmatrix} = |A+\imag B|\cdot|A-\imag B|
\]
If $A$ is a skew-symmetric matrix then 
$\begin{vmatrix}
  -A & B  \\ -B & -A
\end{vmatrix} =
\begin{vmatrix}
  A & B  \\ -B & A
\end{vmatrix} 
$
\end{lemma}
\begin{proof}
Assume all matrices are $n$ by $n$. We compute
\begin{gather*}
2^{2n}
\begin{vmatrix}
  A &  B\\  -B & A 
\end{vmatrix}
=
\left|
  \begin{pmatrix}
    I & - \imag I \\ - \imag I & I
  \end{pmatrix}
\begin{pmatrix}
  A &  B\\  -B & A 
\end{pmatrix}
  \begin{pmatrix}
    I &  \imag I \\  \imag I & I
  \end{pmatrix}
\right|  \\ =
\begin{vmatrix}
  2(A+\imag B) & 0 \\ 0 & 2(A-\imag B)
\end{vmatrix}
\end{gather*}

The second part follows from the first.
\end{proof}

\begin{cor}
   If $A$ is a skew-symmetric matrix, $S$ is symmetric, and
  $B$ is positive definite then 
\[
\begin{vmatrix}
  A & S + \imag B \\ -S-\imag B & A
\end{vmatrix}\ne0
\]
\end{cor}
\begin{proof}
The determinant equals $|  2(A-B+\imag S)|\cdot|2(A+B-\imag S)|$, and
is non-zero by Lemma\ref{lem:skew-non-zero}.
\end{proof}

In the next result we determine properties of eigenvalues.

  \begin{lemma}\label{lem:pos-def-cpx}
    Suppose that $A,B$ are positive definite and $\alpha+\beta\imag$ is an
    eigenvalue of $A+B\imag$.
    \begin{enumerate}
    \item  $0<\alpha<\rho(A)$, $0<\beta<\rho(B)$, $\beta/\alpha \le \rho(BA^{-1})$.
    \item If $A$ is only symmetric then $\beta$ is positive.
    \item If $B$ is only symmetric then $\alpha$ is positive.
    \end{enumerate}
  \end{lemma}
  \begin{proof}
    Let $u+ \imag v$ be an eigenvector of $\alpha+\beta\imag$, where
    $u,v$ are real. Then%
\footnote{Thanks to math.sci.research and Tony O'Conner for this argument.}
    \begin{align*}
      (A+B\imag) (u+\imag v)&= (\alpha+\beta\imag)(u+\imag v)\\
      \intertext{Separating real and imaginary parts yields}
      Au-Bv &= \alpha u - \beta v\\
      Bu+Av &= \alpha v + \beta u\\
      \intertext{$<x,y> = x^ty$ is the usual inner product.
         The symmetry of $A$ and $B$ yields:}
      <u,Au> - <u,Bv> &= \alpha\,<u,u> -\,\beta\,<u,v> \\
      <u,Bv> + <v,Av> &= \alpha\,<v,v> +\, \beta\,<u,v>
    \end{align*}
\begin{align*}    
\intertext{We can now solve for $\alpha$, and a similar
        calculation yields $\beta$}
      \alpha &= \frac{<u,Au> + <v,Av>}{<u,u>+<v,v>}  \le \rho(A) \\
      \beta &= \frac{<u,Bu> + <v,Bv> }{<u,u> +<v,v> } \le \rho(B)\\
      \frac{\beta}{\alpha} &= \frac{<u,Bu> + <v,Bv> }{<u,Au> + <v,Av> } \le
      \rho(BA^{-1})
    \end{align*}
    We use these formulas to establish the various parts of the lemma.
  \end{proof}

  \begin{figure}
    \centering
  \begin{pspicture}(-1,-.5)(2,2.5)
      \psline(-.5,1)(2.5,1)
      \psline(1,-.5)(1,2.5)
      \pspolygon[fillstyle=solid,fillcolor=yellow](1,1)(2,2)(2,1)
      \pscurve{<-}(1.5,1.3)(2,1.4)(2.5,1.3)
      \rput[l](2.6,1.3){eigenvalues of $A+B\imag$}
    \end{pspicture}
    
    \caption{Location of roots of $|-xI+A+B\imag|$ if $B>A$.}
    \label{fig:ang-32}
  \end{figure}

\section{Determinants}
\label{sec:nv-determinants}

We can construct elements of various $\nv{d}{\cone}$ using
determinants.

\begin{cor}\label{cor:nv-det}
  If $\gamma\in\complexes\setminus0$, $\plane$ is the upper half
  plane, $D_i$ are $n\times n$ positive definite matrices and $S$ is
  symmetric then
\[ \bigl| \gamma\, S + \sum_1^d x_k\,D_k\bigr|
\in\nv{d}{(\gamma\plane)^d} \]
\end{cor}
\begin{proof}
  If $\gamma\sigma_1,\dots,\gamma\sigma_d\in\gamma\plane$ where
  $\sigma_k = \alpha_k+\imag \beta_k$ then
\begin{gather*} 
\biggl| \gamma\, S + \sum_1^d x_k\,D_k\biggr|(\gamma\sigma_1,\dots,\gamma\sigma_d) 
= \gamma^n \bigl| S + \sum \sigma_k\,D_k\bigr| \\
= \gamma^n \bigl| \bigl(S + \sum \alpha_k\,D_k\bigr) + \imag\bigl(
\sum \beta_k\,D_k\bigr)\bigr| 
\end{gather*}
and this is non-zero by the lemma since $\beta_k>0$ and so $\sum
\beta_k D_k$ is positive definite. 
\end{proof}

Next we construct polynomials in $\nv{d}{\cone^d}$ where $\cone$ is
not a half plane.

\begin{lemma}\label{lem:nv-det-1}
  If $A_k,B_k$ are positive definite and $\cone$ is the first quadrant
 then
\[ \biggl| \sum_1^d x_k (A_k + \imag B_k)\biggl| \in\nv{d}{\cone^d}
\]
\end{lemma}
\begin{proof}
  If $\alpha_k+\imag \beta_k\in\cone$ then $\alpha_k$ and $\beta_k$
  are positive, and 
  \begin{equation}\label{eqn:nv-cone}
    \biggl| \sum_1^d (\alpha_k+\imag\beta_k) (A_k + \imag B_k)\biggl| =
    \biggl| \sum_1^d (\alpha_kA_k -\beta_kB_k) + \imag
\sum_1^d (\alpha_kB_k +\beta_kA_k)\biggr| 
  \end{equation}
Since  $\sum_1^d (\alpha_kB_k +\beta_kA_k)$ is positive definite
the result follows from Lemma~\ref{lem:nv-not-zero}.
\end{proof}

\begin{lemma}\label{lem:nv-det-2}
  If $A_k$ is positive definite, $A_k > B_k $ and $\cone
  =\{a+b\imag\,|\,a>b>0\}$ then
\[ \biggl| \sum_1^d x_k (A_k + \imag B_k)\biggl| \in\nv{d}{\cone^d}
\]
\end{lemma}
\begin{proof}
  The hypothesis imply that
\[
\alpha_kA_k -\beta_kB_k = A_k(\alpha_k-\beta_k) + \beta_k(A_k-B_k) \]
is positive definite, so the result follows from
Lemma~\ref{lem:nv-not-zero} and \eqref{eqn:nv-cone}.
\end{proof}

  \begin{lemma}\label{lem:poly-from-det}
 Suppose $A$ and $B$ are symmetric and $f(x) = |-xI+A+B\imag|$.
    \item If $B>A>0$ then $f\in\nv{1}{\calD}$ where
      $\calD=\{a+b\imag\mid a>b>0\}$. 
  \end{lemma}
  \begin{proof}
    Since $B>A$ we know that $\rho(BA^{-1})>1$. From the lemma the
    eigenvalues $\alpha+\beta\imag$ satisfy $\beta>\alpha>0$. The
    roots of $f$ are the negative eigenvalues of $A+B\imag$, so they
    lie in the region of the first quadrant bounded by $x=y$ and
    $x=0$.  See Figure~\ref{fig:ang-32}.

  \end{proof}

\section{Linear transformations}
\label{sec:nv-lt}
\index{Hadamard product! for cones}

Properties of the Hadamard product in $\allpoly$ depended on the
identification of the Hadamard product as a coefficient:
\begin{quote}
  If $f(x)\in\allpolypos$ and $g\in\allpoly$ then $f(-xy)$ and 
  $g(y)\in\gsubclose_2$, and the Hadamard product is a coefficient of
  $f(-xy)g(y)$. 
\end{quote}
The reason why $f(-xy)\in\gsubclose_2$ is that $\alpha-xy\in\gsubclose_2$
if and only if $\alpha>0$, and this follows by reversing $x+\alpha
y$.  A similar argument works for half planes.

\begin{lemma}
  Suppose that $\cone$ is a half plane with reversal constant $\rho$,
  and $\ell$ is the ray $\rho^{-1}\reals^+$. The Hadamard product 
\[ \sum a_ix^i \times \sum g_i(\yy)x^i \mapsto \sum a_i g_{n-i}(\yy)x^i \]
 determines a map
\[ \allpolyint{\ell} \times\nv{d}{\cone^d} \longrightarrow
\nv{d}{\cone^d}
\]

\end{lemma}

\begin{proof}
  If $r>0$ then $x+r y\in\nv{2}{\cone^2}$ so $\rho xy +
  r\in\nv{2}{\cone^2}$ and hence $xy-(-r/\rho)\in\nv{2}{\cone^2}$.
  Thus, if $f(x)\in\allpolyint{\ell}(n)$ then
  $f(xy)\in\nv{2}{\cone^2}$. Write $f = \sum a_i x^i$ and $g = \sum
  g_i(\yy) x^i$. The coefficient of $z^n$ in 
\[ f( x z) g(z) =  \sum_{i,j} a_i\,  x^i\, z^i\, g_j(\yy) \,z^{j}
\]
 is in $\nv{d}{\cone}$ and equals
$\displaystyle \sum a_i g_{n-i}(\yy) x^i.$
\end{proof} 

If we consider $f(x+y)$ we get
\begin{lemma}
If $\cone$ is a half plane then  $f(x)\times \sum b_iy^i\mapsto
\sum f^{(i)}(x)\, b_{n-i}/i!$ \\  determines a map
\[
\nv{2}{\cone^2}\times \nv{1}{\cone} \longrightarrow \nv{2}{\cone^2}
\]
\end{lemma}
\begin{proof}
  Write $f(x+y) = \sum f^{(i)}(x)\, y^i/i!$. The coefficient of $y^n$ in
  $f(x+y)\sum b_iy^i$ is
\[ \sum \frac{f^{(i)}(x)}{i!}\,b_{n-i}.\]
\end{proof}

\begin{lemma}\label{lem:nv-x1x2}
  Suppose $\cone$ is a half plane with reversal constant $\rho$.
  \begin{enumerate}
  \item If $f\in\nv{d}{\cone^d}$ then $f + \rho\, \partial_{x_1}
    \,\partial_{x_2} f\in\nv{d}{\cone^d}$.
\item The linear transformation $ g\times f\mapsto g(\rho\,
  \partial_{x_1}\partial_{x_2})f$ determines a map
\[ \allpolypos \times \nv{d}{\cone^d} \longrightarrow \nv{d}{\cone^d} \]
  \end{enumerate}
\end{lemma}
\begin{proof}
  The first part follows from Lemma~\ref{lem:nv-24} and $f \nvlace
  \partial_{x_2}f \nvlace
  \partial_{x_1}\bigl(\partial_{x_2} f\bigr)$. The second part follows
  by induction.
\end{proof}

\section{Analytic functions}
\label{sec:nv-anal}

We define $\nvf{d}{\cone}$ to be the uniform closure on compact
subsets of $\nv{d}{\cone}$. Here are some exponential functions in
$\nv{d}{\cone}$. 

\begin{lemma}
  Suppose $\cone\subset\complexes$ is a cone.
  \begin{enumerate}
  \item If   $\beta\in\overline{\cone}\setminus0$ then
    $e^{x/\beta}\in\nvf{1}{\cone}$.
  \item If $\rho$ is the reversal constant then
$e^{\rho xy} \in\nvf{2}{\cone}$.
  \end{enumerate}
\end{lemma}
\begin{proof}
  Since $\beta+x\in\nv{1}{\cone}$ we have
  $\biggl(1+\frac{x/\beta}{n}\biggr)^n$ is in $\nv{1}{\cone}$, and
  hence $e^{x/\beta}\in\nv{1}{\cone}$.  For the second part we know
  $x+y\in\nv{2}{\cone^2}$, so $1+\rho xy\in\nv{2}{\cone^2}$. Taking
  limits and following the first part finishes the proof.
\end{proof}

\begin{prop}
  If $\cone$ is a half plane with reversal constant $\rho$ then
  \begin{enumerate}
  \item $f\mapsto e^{\rho \partial_{\xx}\cdot \partial_{\yy}}f$ determines a map
    \[\nv{d}{\cone^d}\longrightarrow \nv{d}{\cone^d}\]
  \item $f(\xx)\times g(\xx)\mapsto f(\rho\partial \xx)g(\xx)$ determines
    a map  \[\nv{d}{\cone^d}\times \nv{d}{\cone^d}\longrightarrow \nv{d}{\cone^d}.\]
  \end{enumerate}
\end{prop}
\begin{proof}
  We know that $g + \rho
  \partial_{x_1}\partial_{y_1}g\in\nv{d}{\cone^d}$. Thus 
  $g\mapsto \bigl(1 + \frac{\rho
    \partial_{x_1}\partial_{y_1}}{n}\bigr)^n\,g$ maps
  $\nv{d}{\cone^d}$ to itself, since the homogeneous part of $\bigl(1
  + \frac{\rho
    \partial_{x_1}\partial_{y_1}}{n}\bigr)^n\,g$ equals
  $g^H$. Consequently, taking limits shows that the homogeneous part
  of $e^{\rho \partial_{x_1}\partial_{y_1}}g$ equals $g^H$, and in
  particular is not zero. The result now follows from
  Lemma~\ref{lem:nv-hur} and the observation that
  \[ \rho \,\partial_{\xx}\cdot\partial_{\yy} = \rho
  \partial_{x_1}\partial_{y_1} + \cdots + \rho
  \partial_{x_d}\partial_{y_d}\]

The second part relies on the identity
\[
e^{\rho\partial_\xx\partial_\yy}\,g(\xx)\,f(\yy)\biggl|_{\yy=0} \,=\,
f(\rho\partial_\xx)g(\xx).
\]
By linearity we only need to check it for monomials:
\[
e^{\rho\partial_\xx\partial_\yy} \xx^\sdiffi \yy^\sdiffj\biggl|_{\yy=0} \, =\,
\frac{(\rho\partial_\xx)^\sdiffj}{\diffj!} \xx^\sdiffi \,\cdot\, (\partial_\yy)^\sdiffj
\yy^\sdiffj\biggl|_{\yy=0}\, =\,
(\rho\partial_\xx)^\sdiffj \xx^\sdiffi
\]

\end{proof}

\begin{example}
  The infinite product
$
\cos(x)  = \prod_{k=0}^\infty
\left(1 - \frac{4x^2}{(2k+1)^2\pi^2}\right) \\
$
shows that $\cos(\gamma x)\in\nv{1}{\gamma \plane}$. Taking derivatives
shows 
\[\cos(\gamma x) \nvlace -\gamma \sin(\gamma x)
\quad\text{and}\quad
\sin(\gamma x) \nvlace \gamma \cos(\gamma x)
\]
This gives us elements in $\nv{2}{\gamma\plane}$:
\[\cos(\gamma x)  -\gamma\, y\,\sin(\gamma x)
\quad\text{and}\quad
\sin(\gamma x) + \gamma\, y\,\cos(\gamma x)
\]

\end{example}

\section{Homogeneous polynomials}

Define $\homog{d}$ to be the set of homogeneous polynomials contained
in $\hb{d}$. If we don't want to specify $d$ we write $\homog{\ast}$.
By Lemma~\ref{lem:nv-fh} we can write $f(\xx)\in\homog{d}$ as $\alpha\cdot
g(\xx)$ where $g(\xx)$ has all positive coefficients.  If $d=2$ then
$f(x,y)\in\homog{d}$ is of the form $\alpha\cdot G(x,y)$ where $G$ is
the homogenization of a polynomial in $\allpolypos$. In general,
$f\in\homog{d}$ is a multiple of the homogenization of a polynomial in
$\gsubplus_{d-1}$, since $\gsubclose_{d-1}=\rhb{d-1}$
(Theorem~\ref{thm:rhb=p}). 

The first lemma shows that there is nothing special about the upper
half plane in the definition of $\homog{d}$; any half plane through
the origin will do.

\begin{lemma}
  $\homog{d} = \bigcap_{\gamma\in\complexes\setminus0} \nv{d}{(\gamma\plane)^d}$
\end{lemma}
\begin{proof}
  If $f\in\homog{d}$ has degree $n$ and $\sigma\in\plane^d$ then
  $f(\gamma\sigma) = \gamma^n\,f(\sigma)\ne0$ so
  $f\in\nv{d}{(\gamma\plane)^d}$.

  Conversely, if $f(\xx)$ has degree $n$ and is in the intersection,
  $\sigma\in\plane^d$ and $\gamma\ne0$ then $f(\gamma\sigma)\ne0$. But
  $f(\gamma\sigma)$ is a polynomial in $\gamma$, so the solutions to
  $f(\gamma\sigma)=0$ must all be $\gamma=0$. Thus $f(\gamma\sigma)$
  is divisible by $\gamma^n$, and $f(\gamma\sigma) =
  \gamma^n\,f(\sigma)$. Since this holds for all $\sigma\in\plane^d$
  we see that $f(\gamma\xx) = \gamma^n\,f(\xx)$, so $f$ is homogeneous.
\end{proof}

There is a simple determinant construction for polynomials in
$\homog{d}$, but not all such polynomials can be so realized.

\begin{cor}\label{cor:nv-det-2}
  If all $D_i$ are positive definite then $\displaystyle \bigl| \sum_1^d
  x_k\,D_k\bigr| \in \homog{d}.$
\end{cor}
\begin{proof}
  This follows from Corollary~\ref{cor:nv-det}.
\end{proof}

Next, we construct homogeneous polynomials from elementary symmetric
functions. Fix a positive integer $n$, and let
$\sigma_i(y_1,\dots,y_n)$ be the $i$'th elementary
symmetric function of $y_1,\dots,y_n$. If $\diffi$ is an index set then
\[
\sigma_\sdiffi(\xx) = \prod_{i\in\sdiffi} \sigma_i(x_{i,1},\dots,x_{i,n})
\]

\begin{lemma}
  For positive integers $n$ and $m$ we have that $\displaystyle \sum_{|\diffi|=m}
  \sigma_\sdiffi(\xx) \in\homog{*}$.
\end{lemma}
\begin{proof}
Since $uz-1\in\hb{2}$ 
\[
\prod_{i=1}^n (x_{i,1}z-1) = \sum _k
\sigma_k(x_{i,1},\dots,x_{i,n})(-z)^{n-k} \in\hb{*}
\]
and thus
\[
\prod_{j=1}^d
\biggl(\sum_k\sigma_k(x_{j,1},\dots,x_{j,n})(-z)^{n-k}\biggr)\in\hb{*}
\]
The coefficient of $z^m$ is $\pm
\sum_{|\sdiffi|=m}\sigma_{\sdiffi}(\xx)$. This is in $\homog{*}$ since
all monomials have degree $m$.
\end{proof}

Now we have the Schur-\Szego\ theorem for homogeneous polynomials.
\index{Schur-\Szego!theorem}

\begin{lemma}
If $f = \sum a_\sdiffi\,\xx^\sdiffi$,    $g = \sum
b_\sdiffi\,\xx^\sdiffi$,   $deg(f)=deg(g)$, and $f,g\in\homog{d}$ then
\[
 \sum a_\sdiffi\,b\sdiffi\, \diffi!\,\xx^\sdiffi \in\homog{d}
\]
\end{lemma}
\begin{proof}
  If $f$ has degree $n$ then $f(\partial\xx)g(\xx+\yy)$ is in
  $\homog{*}$ and equals
\[
\sum_{|\sdiffi|=|\sdiffj|=n} a_\sdiffi\, \partial\xx^\sdiffi\,
b_\sdiffj\, (\xx+\yy)^\sdiffj =
 \sum a_\sdiffi\,b\sdiffi\, \diffi!\,\xx^\sdiffi
\]
since 
$(\partial\xx^\sdiffi)(\xx+\yy)^\sdiffj =
\begin{cases}
  0 & \diffi\ne\diffj \\ \diffi!\,\yy^\sdiffi & \diffi=\diffj
\end{cases}
$
\end{proof}

\begin{lemma}
  If $\sum a_\sdiffi \,\xx^\sdiffi \in\hb{d}$ then $\sum a_\sdiffi
  \sigma_\sdiffi(\xx)\, \diffi!\in\hb{d}$.
\end{lemma}
\begin{proof}
  Note that
\[
(\partial x)^i \bigl[ (x+x_1)\cdots(x+x_n)\bigr]\bigr|_{x=0} =
\sigma_i(\xx)\,i!
\]
and thus
\[
\bigl(\sum a_\sdiffi\, \partial
\xx^\sdiffi\bigr)\,\biggl(\prod_{i=1}^d\,\prod_{k=1}^n \, (x_i +
x_{i,k}) \biggr)\biggr|_{\xx=0} = \sum a_\sdiffi \sigma_\sdiffi(\xx)\,
\diffi!.
\]
\end{proof}

\begin{cor}
  The linear transformation $\xx^\sdiffi \mapsto \sigma_\sdiffi(\xx)\,
  \diffi!$ maps $\hb{*}\longrightarrow\hb{*}$ and
  $\homog{*}\longrightarrow\homog{*}$. 
\end{cor}

\section{\Mobius\ transformations}

Suppose $Mz = \frac{az+b}{cz+d}$ is a \Mobius\ transformation where
$|M|=ad-bc\ne0$. If $f(x_1,\dots,x_r)$ is homogeneous of degree $n$
then we say that $f$ is \emph{M-invariant}\index{M-invariant} if
$(Mf)(\xx) = |M|^nf(\xx)$ for all \Mobius\ transformations $M$ where
$|M|\ne0$ and 
\[
(Mf)(x_1,\dots,x_r) = \bigl[(cx_1+d)\cdots(cx_r+d)\bigr]^n\cdot
f\bigl(\frac{ax_1+b}{cx_1+d},\dots, \frac{ax_r+b}{cx_r+d}\bigr)
\]

We will see (Lemma~\ref{lem:nv-mob-1}) that there are no M-invariant
polynomials in $\homog{d}$, but we can find polynomials $f(\xx,\yy)$
such that
\begin{enumerate}
\item $f$ is M-invariant.
\item $f(\xx,-\yy)\in\homog{d}$
\end{enumerate}
We will call such polynomials \emph{\index{Ruelle polynomials}Ruelle
  polynomials} \cite{ruelle}. It is easy to show that a polynomial is
M-invariant; the difficulty is showing that it is  $\homog{d}$. 

\begin{lemma}\label{lem:nv-mob-1}
  If $f\in\hb{n}$ then there is a \Mobius\ transformation $M$ so that
  $Mf\not\in\hb{n}$. 
\end{lemma}
\begin{proof}
  Choose $\alpha_i$ so that $f(\alpha_1,\dots,\alpha_n)=0$, and choose
  $\beta$ so that $\beta+\alpha_i$ is in the upper half plane for all
  $i$. If $Mz=z-\beta$ then the roots of 
\[
Mf(x_1,\dots,x_n) = f(x_1-\beta,\dots,x_n-\beta)
\] are $\beta+\alpha_1,\dots,\beta+\alpha_n$ which are all in the
upper half plane, so $Mf\not\in\hb{n}$.
\end{proof}

We now give several constructions of M-invariant polynomials. They all
depend on the fact that $x-y$ is M-invariant:

\[
M(x-y) = (cx+d)(cy+d)\biggl[\frac{ax+b}{cx+d}-\frac{ay+b}{cy+d}\biggr]
= |M| (x-y)
\]

\begin{construction}
  If $(a_{ij})=A$ is an $r$ by $r$ matrix and
  $X=\diag(x_1,\dots,x_r)$, $Y=diag(y_1,\dots,y_r)$ then we define
\[
f(\xx,\yy) = det(AX-YA) = det(a_{ij}(x_i-y_j)).
\]
$f(\xx,\yy)$ is M-invariant since
\[
Mf = det(a_{ij}M(x_i-y_j)) = det(|M|a_{ij}(x_i-y_j)) = |M|^r f(X,Y).
\]
\end{construction}

\begin{construction}
  With the same setup as above we let
\[
g(\xx,\yy) = per(AX-YA)
\]
\index{permanent}%
where $per$ is the permanent. The same argument shows $g(\xx,\yy)$ is M-invariant.
\end{construction}

If we expand the determinant or the permanent we find that these two
constructions are special cases of

\begin{construction}
If $\alpha_\sigma\in\complexes$ for each permutation $\sigma$ then define
\[ h(\xx,\yy) = \sum_{\sigma\in sym(r)}   \alpha_\sigma \, \prod_{i=1}^r
(x_i-y_{\sigma i}).
\]
This is clearly M-invariant. Ruelle \cite{ruelle} proved that all
Ruelle polynomials can be written in this form, but the
$\alpha_\sigma$ are not unique.
\end{construction}

\begin{example}
  If $J$ is the $r$ by $r$ all $1$ matrix then Cauchy's determinant
  formula \cite{kratt} gives
\begin{align*}
f(\xx,\yy) & = |JX-YJ| = |(x_i-y_j)| = \prod_{i<j}(x_i-x_j)\, \cdot\,
\prod_{i<j}(y_i-y_j) \\
&= \sum_{\sigma\in sym(r)} sign(\sigma) \prod_1^r(x_i-y_{\sigma i})
\end{align*}
Note that $f(X,-Y)$ is not in $\hb{2r}$ since it contains $x_i-x_j$
factors. However, we will see later 
that
\[
per(XJ-JY) = \sum_{\sigma\in sym(r)}  \prod_1^r(x_i-y_{\sigma i})
\] 
is a Ruelle polynomial. We call this the \index{Grace polynomial}Grace
polynomial, denoted $\grace{r}$. For instance, 
\begin{align*}
  \grace{2} &= (x_1-y_1)(x_2-y_2) + (x_1-y_2)(x_2-y_1) \\
  &= 2x_1x_2 - (x_1+x_2)(y_1+y_2)+ 2 y_1y_2 \\
  &= 2\sigma_2(x_1,x_2)\sigma_0(y_1,y_2) -
\sigma_1(x_1,x_2)\sigma_1(y_1,y_2) + 2\sigma_0(x_1,x_2)\sigma_2(y_1,y_2)
\end{align*}

\end{example}

The next lemma is Ruelle's original definition of a Ruelle polynomial.

\begin{lemma}[Ruelle]
  $f(\xx,\yy)$ is a Ruelle polynomial iff whenever there is a circle
  separating $\sigma=(\sigma_1,\dots,\sigma_n)$ from
  $\tau=(\tau_1,\dots,\tau_m)$ then $f(\sigma,\tau)\ne0$.
\end{lemma}
\begin{proof}
  Choose a \Mobius\ transformation such that $M\sigma$ is in the upper
  half plane and $M\tau$ is in the lower half plane. Since
  $f(M\xx,-M\yy)$ is in $\homog{d}$ we have that
  $f(M\sigma,-M\tau)\ne0$. The converse is harder, and can be found in
  \cite{ruelle}. 
\end{proof}

\begin{lemma}[Ruelle] If $A$ is unitary then $f(\xx,\yy)=|AX-YA|$ is a
  Ruelle polynomial.  
\end{lemma}
\begin{proof}
  We have seen that $f$ is M-invariant. If $A$ is unitary then
  $A^{-1}=A^*$ so
\[
f(\xx,\yy) = |AX-YA| = |A|\, |X-A^{\ast}YA| 
\]
If the rows of $A$ are $A_1,\dots,A_r$ then
\[
f(\xx,\yy) = |A| \,det\biggl( X + \sum y_i\, A^*_i\,A_i\biggr)
\]
Since all $A_i^* A_i$ are positive semi-definite the determinant is in
$\hb{2r}$. 
\end{proof}

\begin{example}
  In this example we give the Ruelle polynomial that is constructed
  using the general three by three orthogonal matrix \cite{bernstein}
  where $\beta ={a^2+b^2+c^2+d^2}$.

\[
A=\frac{1}{\beta}\left(
\begin{array}{lll}
 {a^2+b^2-c^2-d^2} & {2 (b c+a d)} &
 {2 (b d-a   c)} \\
 {2 (b c-a d)} & {a^2-b^2+c^2-d^2} &
 {2 (a b+c   d)} \\
 {2 (a c+b d)} & {2 (c d-a b)} &
   {a^2-b^2-c^2+d^2}
\end{array}
\right)
\]

If we expand the determinant  we get
the following permutation representation:

\[
\begin{array}{rl}

 &  (y_3+x_1) (y_2+x_2) (y_1+x_3) + \\
 4 (a c+b d)^2
& (y_2+x_1) (y_3+x_2)(y_1+x_3)+\\
4 (a c-b d)^2 
&  (y_3+x_1) (y_1+x_2) (y_2+x_3) + \\
4 (b^2-c^2) (a^2-d^2)
&(y_1+x_1)   (y_3+x_2) (y_2+x_3) +\\
4 (a^2-b^2) (d^2-c^2)
&(y_2+x_1) (y_1+x_2) (y_3+x_3) +\\
\left(a^2-b^2+c^2-d^2\right)^2
&   (y_1+x_1) (y_2+x_2) (y_3+x_3)
\end{array}
\]
\end{example}

In order to show that the Grace polynomial is a Ruelle polynomial we
need to recall Grace's theorem:

\begin{theorem}[Grace]
  If $f= \sum_0^n a_ix^i$, $g = \sum_0^n b_ix^i$, and the roots of $f$
  are separated by a circle from the roots of $g$ then
\[
\sum_0^n a_i\, b_{n-i}\,(-1)^i\, i!\, (n-i)! \ne0
\]
\end{theorem}

\begin{cor}
  The Grace polynomial is a Ruelle polynomial.
\end{cor}

\begin{proof}
  We know the Grace polynomial is M-invariant; we now show it
  satisfies the second condition.
  Suppose that $x_1,\dots,x_n$ and $y_1,\dots,y_n$ are in the upper
  half plane. Grace's theorem applied to $\prod(x-x_i)$ and
  $\prod(x+y_i)$ shows that
\begin{gather*}
0 \ne \sum (-1)^i \sigma_i(x_1,\dots,x_n)\cdot
\sigma_{n-i}(y_1,\dots,y_n) (-1)^i\,i!\,(n-i)! \\
= \sum \sigma_i(x_1,\dots,x_n)\sigma_{n-i}(y_1,\dots,y_n)\,i!\,(n-i)!
\end{gather*}

Comparing monomials shows that this is exactly $\grace{r}(\xx,-\yy)$,
and it is non-zero for substitutions in the upper half plane.
\end{proof}

Here's a simple one-variable corollary.

\added{7/15/7}
\begin{cor}
  If $a_i,b_i\in\reals$ then
\[
\sum_{\sigma\in sym(n)} \prod_i(x+a_i + b_{\sigma i}) \in\allpoly
\]
\end{cor}
\begin{proof}
  Replace $y_i$ by $b_i$ and $x_i$ by $x+a_i$ in the Grace polynomial.
\end{proof}

The conjecture below generalizes Grace's theorem:

\begin{conj}
  If  $n$ is a positive integer, $\sigma_i(x_k)$ the $i$'th
  elementary symmetric polynomial on $x_{k,1},\dots,x_{k,n}$ and
  $\sigma_{\sdiffi}(\xx) = \prod_{k} \sigma_{i_k}(x_k)$ then
\[ \sum_{|\sdiffi|=n} \sigma_\sdiffi(\xx)\,\diffi! \in\homog{*}
\]
\end{conj}

\section{Properties of polynomials with zeros in a strip}
\label{sec:prop-polyn-with}

Assume that $\calS$ is the strip $\{z\,\mid\, -1 < Re(z) < 1\}$. 
We will consider properties of $\nv{1}{\complexes\setminus\calS}$. If
$f\in\nv{1}{\complexes\setminus\calS}$ then all roots of $f$ lie in
$\calS$. Since $\calS$ is convex we have that
\begin{quote}
  $\nv{1}{\complexes\setminus\calS}$ is closed under differentiation.
\end{quote}
The reverse of $x-\alpha$ has root $1/\alpha$, so
\begin{quote}
  $\nv{1}{\complexes\setminus\calS}$ is not closed under reversal.
\end{quote}
A non-trivial property is this interlacing result:

\begin{lemma}
  If $\calS$ is a vertical strip and
  $f\in\nv{1}{\complexes\setminus\calS}$ then for all real $\alpha$
\[
f + \imag \,\alpha\, f' \in \nv{1}{\complexes\setminus\calS}
\]
\end{lemma}

The proof follows from the next lemma.

\begin{lemma}
  Suppose that $f = \prod(x-r_k)$, and define 
\[
 g = f + \imag \alpha \sum_k a_k \frac{f(x)}{x-r_k}
\]
where $\alpha\in\reals$ and all $a_k$ are non-negative. If $a+b\imag$
is a root of $g$ then
\[
 \min_k \Re(r_k) \le a \le \max_k \Re(r_k)
\]
\end{lemma}
\begin{proof}
  Let $r_k = s_k + \imag t_k$. Dividing by $f$ yields 
  \begin{align*}
    1 + \imag \,\alpha\, \sum_k a_k \frac{1}{a+b\imag - s_k - \imag t_k} 
&= 0 \\
    \sum_k a_k \frac{a-\imag b - s_k + \imag t_k}{|a+b\imag - s_k -\imag t_k|^2} 
&= \imag/\alpha\\
\intertext{Taking the real part yields}
    \sum_k a_k \frac{a-  s_k}{|a+b\imag - s_k - \imag t_k|^2 }
&= 0
  \end{align*}
Thus $a$ can't be less than all $s_k$, nor greater than all $s_k$. 
\end{proof}

\begin{cor}
  If $f\in\allpoly$ then the map $g\mapsto f(\imag\diffd)g$ determines
  a linear transformation
\[ \nv{1}{\complexes\setminus\calS}\longrightarrow\nv{1}{\complexes\setminus\calS}\]
\end{cor}

\begin{cor}
  If $g\in\nv{1}{\complexes\setminus\calS}$ then
\[ g + g'' + g^{(4)}/2! + g^{(6)}/3! + \cdots \in
\nv{1}{\complexes\setminus\calS}
\]
\end{cor}
\begin{proof}
  Apply the previous corollary to $f=e^{-x^2}$.
\end{proof}

The polynomial $xy-1$ is in $\nv{2}{\complexes\setminus\calS}$ since
$|x|>1$ if $x\not\in\calS$. If $x,y\notin\calS$ then $|xy|>1$, and
hence $xy-1\ne0$. This observation will be generalized in the next section

\section{Polynomials non-vanishing on the open unit ball}
\label{sec:open-unit-ball}

Suppose that $\Delta_d$ is the open unit ball in $\complexes^d$. That
is,
\[
\Delta_d = \{z\in\complexes^d\,\mid\, |z|<1 \}
\]
We will construct some simple non-trivial elements of
$\nv{d}{\calD}$ where $\calD\subset\complexes$, and use them to show that $\nv{d}{\Delta^d}$ has
non-trivial elements. 
Recall that $\nv{d}{\calD}$ contains trivial products of the form
\[
\prod_{i=1}^d \prod_j(x_i-\alpha_{ij})
\]
where $\alpha_{ij}\not\in\calD$. The next lemma constructs elements
that are not products.

\begin{lemma}\label{lem:not-prods}
  Suppose that the ball $B = \{z\,\mid\, |z-\sigma|<r\}$ is contained
  in $\complexes^d\setminus\calD$. If $\sigma=(\sigma_1,\dots,\sigma_d)$
  and $0<b<r^d$ then
\[
(x-\sigma_1)\cdots(x-\sigma_d) - b \in\nv{d}{\calD}
\]
\end{lemma}
\begin{proof}
  If $(\tau_1,\dots,\tau_d)\not\in B$ then $|\tau_i-\sigma_i|\ge r$,
  so $|(\tau_1-\sigma_1)\cdots(\tau_d-\sigma_d)| \ge r^d$. Since
  $0<b<r^d$ it follows that
\[ (\tau_1-\sigma_1)\cdots(\tau_d-\sigma_d)-b \in\nv{d}{\calD}
\]
\end{proof}

\begin{example}
  If $\plane$ is the upper half plane, then the ball
  $\{z\,\mid\,|z+\imag|<1\}$ is contained in
  $\complexes^d\setminus\plane^d$. It follows that
\[ (x_1+\imag)\cdots(x_d+\imag)-1 \in\hb{d}
\] 
A similar argument shows that
\[
 (x_1+1)\cdots(x_d+1)-1 \in\stabled{d}
\]
\end{example}

Here's a small extension to $\nv{d}{\Delta_d}$.

  \begin{lemma}
    If  $f(x)$ has all its zeros in the open unit disk then
\[ 
f(x_1\cdots x_d)\in\nv{d}{\Delta_d}
\]
  \end{lemma}
  \begin{proof}
    Use the lemma with $B = \Delta_d$, and multiply many of the
    polynomials together.
  \end{proof}

It is easy to see that $xy-1/2 \in\nv{2}{\complexes^2\setminus\Delta\times\Delta}$ since if
neither $x$ nor $y$ is in $\Delta\times\Delta$ then both have absolute
value at least one, so their product has absolute value at least
one. There is a simple geometric condition for a multiaffine
polynomial to belong to $\nv{2}{\complexes^2\setminus \Delta^2}$.

\index{multiaffine polynomial}

\begin{lemma}
Suppose that$f(x,y) = a + bx+cy +d xy$, and let
$M=\smalltwo{a}{c}{b}{d}$. 
\begin{enumerate}
\item $f\in \nv{2}{\complexes^2\setminus\Delta^2} \Longleftrightarrow
  M\colon \complexes^2\setminus\Delta^2 \longrightarrow \Delta^2$
\item $f\in \nv{2}{\complexes^2\setminus\overline{\Delta}^2} \Longleftrightarrow
  M\colon \complexes^2\setminus\overline{\Delta}^2 \longrightarrow \overline{\Delta}^2$
\end{enumerate}
\end{lemma}
\begin{proof}
  Assume that $|y|\ge1$. If $f(x,y)=0$ we must show that
  $|x|<1$. Solving for $x$
\[
x = - \frac{a+cy}{b+dy} = - M(y)
\]
Thus $|x|<1$ if and only if $|M(y)|<1$. The second case is similar.
\end{proof}

\index{Kakeya}
Recall the result \mypage{thm:kakeya}  that if $0<a_1<\cdots
<a_n$ then $\sum a_ix^i$ has all its roots in the open unit
ball. Here's a conjecture in two variables 
\begin{conj}
  If $0<a_1 < \cdots < a_n$ then
\[ \sum (a_i + a_j) x^iy^j \in\nv{2}{\complexes^2\setminus \overline{\Delta}^2}
\]
\end{conj}

We can establish the multiaffine case of this conjecture.

\begin{lemma}
  If $0<a<b$ then $2a + (a+b)(x+y) + 2b xy \in\nv{2}{\complexes^2\setminus\overline{\Delta}^2}$.
\end{lemma}

\begin{proof}
  This is a question about \Mobius\ transformations. We take three
  points on the unit circle, and find the unique circle containing
  their image. We check that $M(\infty)$ lies in this circle, and that
  this circle lies in the unit circle. 
  \begin{align*}
    M(1) &= \frac{3a+b}{3b+a} & M(-1) & = 1 \\
    M(\imag) &=\frac{(2+i) a+i b}{a+(1+2 i) b} & M(\infty) &=
    \frac{a+b}{2b}\\
    \text{center} &= \frac{2 (a+b)}{a+3 b} & \text{radius} &= \frac{b-a}{a+3 b}
  \end{align*}
Since the center is on the real line, it follows that the image is the
circle with diameter $M(1)\,M(-1)$. Thus the image lies in the closed unit circle.
\end{proof}

{\psset{unit=2.5cm}
\begin{figure}
  \begin{pspicture}(-1.5,-1.2)(1.5,1.2)

    \pscircle(0,0){1}
    \pscircle(.666,0){.333}
    \psline(-1,0)(1,0)
    \psline(0,-1)(0,1)
    \psline{->}(1,1)(.666,.05)
    \psline{->}(1,-1)(.5,-.05)
    \rput[l](1,-1){$\frac{a+b}{2b}$}
    \rput[l](1,1){$\frac{2a+2b}{a+3b}$}
    \psdots(.333,0)(.666,0)(.5,0)(1,0)
    \rput[l](1.05,.1){$M(-1)$}
    \rput[r](.33,.1){$M(1)$}
  \end{pspicture}
  
  \caption{The image of the circle under a \Mobius\ transformation}
  \label{fig:mobius-circle}
\end{figure}
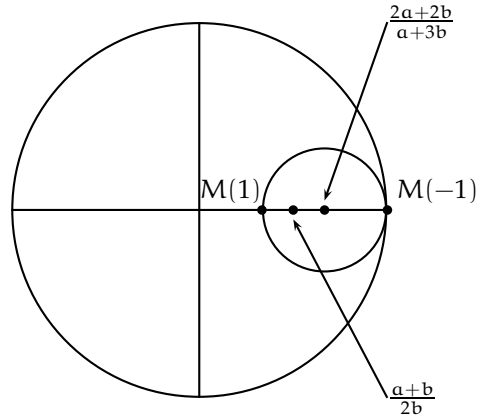
}

The following lemma constructs a stable polynomial from a polynomial
in $\nv{1}{\complexes\setminus\Delta}$.

\begin{lemma}
  If $0<a_0<\cdots < a_n$ then $\sum_0^n a_i (x+1)^i(1-x)^{n-i}\in\stabled{1}$.
\end{lemma}
\begin{proof}
  The \Mobius\ transformation $z\mapsto \frac{1+z}{1-z}$ maps the
  right half plane to $\complexes\setminus\Delta$. It follows that if
  $f(x) \in\nv{1}{\complexes\setminus\Delta}$ then
  $f((1+x)(1-x))\in\stabled{1}$. If we let $f(x) = \sum a_ix^i$ then
  $f\in\nv{1}{\complexes\setminus\Delta}$.  Multiplying by $(1-x)^n$
  yields the result.
\end{proof}

\section{Polynomials non-vanishing in a sector}
\label{sec:polyn-non-vanish}

We begin with a simple construction of polynomials with no roots in a sector.
We claim that if $f$ and $g$ have
  positive leading coefficient then
\[ 
f\in\allpolypos(n) \And g\in\allpolypos(m) \implies f(x) + g(y)
\in\nvr{2}{\sector{2n}\times\sector{2m}} 
\]
To see this, we first claim that the image of $\sector{2n}$ under $f$
is the open right half plane. If $f=\prod(x+r_i)$ where $r_i\ge0$, and
$\sigma\in\sector{n}$, then $r_i+\sigma\in\sector{2n}$, so $f(\sigma)$
is the product of $n$ points whose arguments are all less than
$\pi/(2n)$, Similarly, if $\tau\in\sector{2m}$ then $g(\tau)$ is also
in the right half plane, so their sum is as well. In particular, their
sum is not zero.

For example,
\[
 x^n + y^m \in\nvr{2}{\sector{2n}\times\sector{2m}}
\]

This construction  easily generalizes. 

\begin{lemma}
  Suppose that $n_1,\dots,n_d$ are positive integers. If all
  $f_i(x)\in\allpolyposclose(n_i)$ have positive leading coefficient then
\[
 f_1(x_1) + \cdots + f_d(x_d) \in
\nvr{d}{\sector{2n_1}\times\cdots\times\sector{2n_d}} 
\]
\end{lemma}


\chapter{Upper half plane polynomials}

\label{cha:complex-coef}
 
\renewcommand{\TimeStampStart}{Sunday, February 10, 2008: 11:14:48}

There is a natural subset of polynomials with complex coefficients
that has properties similar to $\allpoly$: those whose roots lie in
the lower half plane. However, we call them \emph{upper half plane
  polynomials} because the generalization to more variables uses the
fact that they are non-zero in the upper half plane. For a quick
summary with proofs see \cite{upper}.

We first recount properties that follow from the general theory of
Chapter~\ref{cha:nv}, and then investigate  polynomials
in one variable.

\section{Properties of  upper half plane polynomials}

We apply the general results of Chapter~\ref{cha:nv} to derive results
about polynomials non-vanishing in the upper half plane.  If $\uhp$
is the upper half plane then $\up{d} = \nv{d}{\uhp^d}$.  The
reversal constant is $-1$, and
$\uhp/\uhp=\complexes\setminus(-\infty,0)$. If $f,g\in\up{d}$ then
we write $f\ulace g$ if $f + z g\in\up{d+1}$.  From
Chapter~\ref{cha:nv} we have

\begin{lemma}\label{lem:hb-basic} \

  \begin{enumerate}
  \item $\up{d}$ is closed under differentiation.
  \item If $f\in\up{d}$ then $f\ulace \partial f/\partial x_i$
  \item If $\sum x^if_i(\yy)\in\up{d+1}$ then $f_i(y)\in\up{d}$ and
    $f_i\ulace f_{i+1}$.
  \item  $fg\in\up{d}$ if and only if $f\in\up{d}$ and $g\in\up{d}$.
  \item If $f\ulace g$ then $g\ulace -f$.
  \item If $f\ulace g\ulace h$ then $f-h\ulace g$.
  \item If $S$ is symmetric and all $D_i$ are positive definite then 
$ \bigl| S + \sum_1^d x_i D_i\bigr| \in\up{d} $.

\item The Hadamard product $\sum a_ix^i\times\sum f_i(\yy)x^i\mapsto \sum
  a_if_{n-i}(\yy)x^i$ determines a map $\allpolypos\times
  \up{d}\longrightarrow \up{d}$.

  \item $f(\xx)\times g(\xx)\mapsto f(-\partial \xx)g(\xx)$ determines
    a map  \[\up{d}\times \up{d}\longrightarrow \up{d}.\]

\item The following are equivalent
  \begin{enumerate}
  \item $f(\xx)\in\up{d}$
  \item $f(\aaa+t\bbb)\in\up{1}$ for $\aaa\in\reals^d,\bbb>0$.
  \end{enumerate}
  \end{enumerate}

\end{lemma}
\begin{proof}
  The only modification is that we do not need homogeneity in the last
  one. Following the proof of Lemma~\ref{lem:nv-14} we see (1)
  implies (2). Conversely, consider $\sigma_1,\dots,\sigma_d$ in $\uhp$.
  Since $1$ is in the closure of $\uhp$ we can find an $\alpha$
  and positive $\aaa,\bbb$ such $(\sigma_1,\dots,\sigma_d) = \aaa +
  \alpha\bbb$. Thus,
  $f(\sigma_1,\dots,\sigma_d)=f(\aaa+\bbb\alpha)\ne0$.
\end{proof}

There are two reversals in $\up{d}$, a single variable reversal with a
minus sign, and a full reversal with no signs.

\begin{lemma}\label{lem:hb-reversal}\ 
  \begin{enumerate}
    \item If $\sum x^if_i(\yy)\in\up{d}$ has degree $n$ then $\sum
    (-x)^{n-i}\,f_i(y)\in\up{d}$.
\item If $x_i$ has degree $e_i$ then $x_1^{e_1}\cdots x_d^{e_d}
  f(1/x_1,\dots,1/x_d)\in\up{d}$. 
  \end{enumerate}
\end{lemma}
\begin{proof}
  Only the second one requires proof. If $\sigma_i\in\uhp$ then
  $1/\overline{\sigma_i}\in\uhp$. Thus 
\[ f(1/\sigma_1,\dots,1/\sigma_d) =
\overline{f(1/\overline{\sigma_1},\dots, 1/\overline{\sigma_d})} \ne 0.
\]
\end{proof}

\mynote{need real and im parts in Pd}

\begin{theorem}\label{thm:rhb=p}
  $\rup{d} = \gsubclose_d$.
\end{theorem}
\begin{proof}
  We first show $\rup{d}\subset\gsubclose_d$. 
  If $f(\xx)\in\rup{d}$ then define $g_\epsilon(\xx)=f(\xx+\epsilon
  (x_1+\cdots+x_d))$ where $\epsilon>0$. We may assume that $f^H$ has
  all non-negative terms, so $g_\epsilon^H$ has all positive terms.
  Define $\aaa=(0,r_2,\dots,r_d)$ and $\bbb=(t,0,\dots,0)$.  Since
  $\rup{d}$ is closed $g_\epsilon(\aaa+\bbb t)\in\rup{1}$.  Since
  $\rup{1}=\allpoly$, we see that $g_\epsilon$ satisfies substitution,
  and therefore $g_\epsilon\in\rup{d}$. The conclusion follows since
  $\lim_{\epsilon\rightarrow0^+} g_\epsilon = f$.

 Next we show $\rup{d}\subset\rup{d}$ which implies
 $\gsubclose_d=\rup{d}$. If $f\in \rup{d}$ and is homogeneous then
 $f(\aaa+\bbb t)\in\allpolypos$ for all $\aaa,\bbb>0$, so
 $f\in\rup{d}$  by Lemma~\ref{lem:hb-basic}. 
 If $f$ is not homogeneous then choose $\ccc\in\reals^d$ so that
 $f(\xx+\ccc)$ has all positive coefficients. Since
 $f(\xx+\ccc)\in\gsubpos_d$ we can homogenize it  $y^nf(\xx/y) =
 F(\xx,y)$. Since $F(\xx,y)\in\rup{d+1}$ and is homogeneous we know
 from the previous paragraph that $F(\xx,y)\in\rup{d+1}$, and thus
 $f(\xx+\ccc)=F(\xx,1)\in\rup{d}$. Now $\ccc$ is real so it follows
 that $f(\xx)\in\rup{d}$. 
\end{proof}

  \begin{lemma}\label{lem:p-square}
    If $f(\xx,y^2)\in\rup{d}$ then $f(\xx,y)\in\gsub_d$.
  \end{lemma}
  \begin{proof}
    If $\sigma,\sigma_i\in\uhp$ then we can write $\sigma=\tau^2$
    where $\tau\in\uhp$. Thus
\[
f(\sigma_1,\dots,\sigma_d,\sigma) =
f(\sigma_1,\dots,\sigma_d,\tau^2)\ne0
\]
which proves the lemma.
  \end{proof}
Of course the converse is false -- consider $x+1$ and $x^2+1$.

\section{Constructions from determinants}
\label{sec:constr-from-determ}

  We can  construct polynomials in $\hb{d}$ using
  skew-symmetric matrices. 
\begin{lemma}\label{lem:p2-skew}
If $A$ is skew symmetric and $B$ is positive definite then 
\[ 
\begin{vmatrix}
  x\,A & B \\ -B & y\,A
\end{vmatrix} \in\gsubclose_2
\]
\end{lemma}
\begin{proof}
If $b = B^{-1/2}=b^t$ then
\[
\begin{pmatrix}
  b&0\\  0&b
\end{pmatrix}
\begin{pmatrix}
  x\,A & B \\ -B & y\,A
\end{pmatrix} 
\begin{pmatrix}
  b&0\\  0&b
\end{pmatrix}
=
\begin{pmatrix}
  x\,bA b&I\\
  -I&y\, bAb
\end{pmatrix}
=
\begin{pmatrix}
  x\,C & I \\ -I & y\, C
\end{pmatrix}
\]
where $bAb=C$ is skew symmetric.  We can find an orthogonal matrix
$\mathcal{O}$ and  $D$ an anti-diagonal matrix with positive anti-diagonal
$(d_i)$ such that $\mathcal{O}C\mathcal{O}^t=
\begin{pmatrix}
  0 & D \\ -D & 0
\end{pmatrix}
$. Clearly
\[
\begin{pmatrix}
  \mathcal{O}&.\\
  .&\mathcal{O}
\end{pmatrix}
\begin{pmatrix}
  x\,C & I \\ -I & x\,C
\end{pmatrix} 
\begin{pmatrix}
  \mathcal{O}^t&.\\
  .&\mathcal{O}^t
\end{pmatrix}
=
\begin{pmatrix}
  . & xD & I & . \\
-xD & . & . & I \\
-I & . & . & yD \\
. & -I & -yD & .
\end{pmatrix}
\]
where the dot is a zero matrix. By \eqref{eqn:skew-det-iden} this
determinant is $\prod_1^n (d_id_{n-i}xy-1)$ which is in $\gsubclose_2$.

\begin{example}
  If $A$ is the skew-symmetric $n$ by $n$ matrix that is all $1$ above
  the diagonal, and all $-1$ below, then
\[
\begin{vmatrix}
  x\,A & I \\ -I & A
\end{vmatrix} = \left[\sum (-x)^k \binom{n}{2k}\right]^2
\]
We already knew that this is in $\allpoly$ since it's just the even
part of $(x-1)^n$.
\end{example}

\end{proof}

The following corollary was proved by a different method in
\cite{frenkel}. 

\begin{cor}
  Suppose $A$ is skew-symmetric and $B$ is positive definite. The
  determinant and pfaffian of $
  \begin{pmatrix}
    xA & B \\ -B & A
  \end{pmatrix}
$ have alternating coefficients.
\end{cor}
\begin{proof}
  Substituting $y=1$ in Lemma~\ref{lem:p2-skew} shows that all the
  roots are positive, and hence the coefficients alternate.  Since the
  pfaffian is the square root of the determinant all the roots are
  negative, and so the coefficients again alternate.
\end{proof}

We can construct polynomials in $d$ variables of even total degree.

\begin{cor}\label{hbc-alt}
  If $A$ is skew-symmetric, and $B_1,\dots,B_d$ are positive definite
  then
\[ 
f(\xx)=\left|
  \begin{pmatrix} A & 0 \\ 0 & A  \end{pmatrix}+
 x_1 \begin{pmatrix} 0 & B_1 \\ -B_1 & 0  \end{pmatrix}+\cdots+
 x_d \begin{pmatrix} 0 & B_d \\ -B_d & 0  \end{pmatrix}
\right|\in\hb{d}
\]
$f(\xx)$ is a square, and all non-zero monomials of $f(\xx)$ have even total degree.
\end{cor}
\begin{proof}
  If we substitute $\alpha_k+\beta_k\imag$ for $x_k$ where all
  $\beta_k$ are positive, then the determinant equals
\[
\left|
  \begin{pmatrix} A & \sum \alpha_kB_k \\ -\sum \alpha_k B_k & A  \end{pmatrix}
  + \imag\,
  \begin{pmatrix}
    0 & \sum \beta_kB_k \\ -\sum \beta_k B_k &0
  \end{pmatrix}
\right|
\]
and this matrix is non-zero by Lemma~\ref{lem:skew-non-zero} since
$\sum \beta_k B_k$ is positive definite. Since the matrix is
anti-symmetric $f(\xx)$ is a square.  Finally, using
Lemma~\ref{lem:abba},
\[
f(-\xx) =   
\begin{pmatrix}
    A & -\sum x_kB_k \\ \sum x_k B_k & A
  \end{pmatrix}
=
\begin{pmatrix}
    -A & -\sum x_kB_k \\ \sum x_k B_k & -A
  \end{pmatrix}
= f(\xx)
\]

\end{proof}

\begin{example}
  If $d=1$, $B_1$ is the identity, and $A$ is the matrix of the last
  example then
\begin{gather*}
\left|
  \begin{pmatrix}
    A & 0 \\ 0 & A 
  \end{pmatrix}
+
x\begin{pmatrix}
 0 & I \\ -I & 0 
\end{pmatrix}
\right|
=
\left[\sum_{i=0}^n (-1)^i x^{n-2i}\binom{n}{2i}\right]^2 \\
|A + \imag x I| = \sum_{i=0}^n (-1)^i x^{n-2i}\binom{n}{2i} \\
\end{gather*}

\end{example}

In the last corollary the matrix had real coefficients and the
polynomial was a square. If we allow complex coefficients then we can
take the square root:

\begin{cor}
  If $A$ is skew-symmetric and all $B_k$ are positive definite then 
\[ g(\xx) = \bigl| A +\imag\, \sum_1^d x_i\,B_i\bigr| \in\rhb{d}\]
and all monomials have even degree.
\end{cor}

\begin{proof} The determinant in Corollary~\ref{hbc-alt} equals 
  \[
   \bigl| A +\imag\, \sum_1^d x_i\,B_i\bigr| 
 \bigl| A   -\imag\, \sum_1^d x_i\,B_i\bigr| \]
If $C$ denotes the sum then 
\[
|A+\imag C| = |(A+\imag C)^t| = |-A+\imag C| = |A-\imag C|
\]
Thus both factors are equal and since the conjugate of $|A+\imag C|$
is 
$|A-\imag C|$ they both have real coefficients. This also shows that
$g(-\xx) = g(\xx)$, and there fore all terms have even degree.
  
\end{proof}

\begin{lemma}
  Suppose that $D_1,D_2,D_3$ are positive definite matrices, and $S$
  is symmetric. If 
\[
f(x,y) = |x\,D_1 + y\,D_2 + S + \imag\, D_3|
\]
then $f(x,y)\in\rup{2}$.
\end{lemma}
\begin{proof}
  The homogeneous part of $f$ is $|xD_1+yD_2|$ which has all positive
  coefficients by Lemma~\ref{lem:pos-def-sum}. For any
  $\alpha\in\reals$ 
\[
f(x,\alpha) = |x\,D_1 + (\alpha\, D_2+S) + \imag\, D_3|
\]
and by  this is in $\polycpx$ since
$\alpha\, D_2+S$ is symmetric.
\end{proof}

\section{$\sltwo$}
\label{sec:sltwo}

In this section we characterize the set of matrices of $\sltwo$ whose entries
interlace. We define

\begin{definition}
  \[
\sltwop = \bigl\{\smalltwo{g}{f}{k}{h}\,\mid\, gh-fk=1 \text{ and }
  f\ulace g\ulace k, f\ulace h\ulace k \bigr\}
\]
\end{definition}

We begin with some elementary properties.

\begin{lemma} Suppose $\smalltwo{g}{f}{k}{h}\in\sltwop$ is not constant.
  \begin{enumerate}
  \item $deg(g)+deg(h) = deg(f) + deg(k)$.
  \item If $sign(r)$ is the sign of the leading coefficient of $r$
    then $sign(g)\cdot sign(h) = sign(f) \cdot sign(k)$.
  \item $\smalltwo{g}{f}{k}{h}^{-1}\not\in\sltwop$.
  \end{enumerate}
\end{lemma}
\begin{proof}
  The first two follow from the fact that $gh-fk$ is a constant, so
  the leading coefficients of $gh$ and $fk$ must cancel. The
  determinant of the inverse is $1$, but all the interlacings go in the
  opposite direction.
\end{proof}

\begin{remark}
  Consider some necessary conditions for  the matrix $M=\smalltwo{g}{f}{k}{h}$ 
 to belong to $\sltwop$.
\begin{enumerate}
\item $M$ has determinant one.
\item Degrees of adjacent entries differ by at most one.
\item All entries are in $\allpoly$.
\item Either $f\ulace g$ or $f\ulace h$ or $g\ulace k$ or $h\ulace k$.
\end{enumerate}

The matrix $\smalltwo{1}{x}{x}{x^2-1}$ satisfies all conditions except
(1).  If (2) fails then not all adjacent polynomials can interlace.
It is possible for only condition (3) to fail -- consider
$\smalltwo{1}{x}{x}{x^2+1}$. In addition, only condition (4) can fail
\[
\left(
\begin{array}{ll}
 x^2+3 x & x^2+3 x+2 \\
 x^2+3 x-\frac{1}{2} & x^2+3 x+\frac{3}{2}
\end{array}
\right).
\]
All the entries are in $\allpoly$, but no two  entries interlace.
\end{remark}

\begin{lemma}\label{lem:sl-1}
  If $M=\smalltwo{g}{f}{k}{h}$ has determinant $1$ and three out four
  interlacings then $M\in\sltwop$.
\end{lemma}
\begin{proof}
  Using the first part of the next lemma we can multiply by certain of
  the first three matrices of the lemma to assure that all leading
  coefficients are positive. The interlacings were preserved, so we
  know that the degrees must be $\smalltwo{n+b}{n+a+b}{n}{n+a}$ where $a,b\ge0$. 

  Suppose that we do not know that $g\longrightarrow f$. Since
  $h\longrightarrow f$ we know that $h$ sign interlaces $f$. Now 
  $fk-gh=1$, so $g$ has the same sign at the roots of $f$ as $h$ does,
  so $g$ sign interlaces $f$. Since $deg(g)\le deg(f)$ it follows that
  $g\longrightarrow f$. The remaining cases are similar.
\end{proof}

We now show that certain matrices preserve $\sltwop$.

\begin{lemma}\label{lem:sl-2}
  Suppose $a>0$. Multiplying on the right or left by any of these
  matrices maps $\sltwop$ to itself:
\[
\begin{matrix}
\smalltwo{0}{-1}{1}{0} & \smalltwo{0}{1}{-1}{0} &
\smalltwo{-1}{0}{0}{-1} &
\smalltwo{1}{ax+b}{0}{1} & \smalltwo{1}{0}{-ax-b}{1} 

  \end{matrix}
\]
\end{lemma}
\begin{proof}
  The first three follow easily from the fact that
\[ f\ulace g \quad \Leftrightarrow \quad  -f \ulace -g  \quad \Leftrightarrow \quad  -g \ulace f
 \quad \Leftrightarrow \quad  g \ulace -f.
\]
If $\smalltwo{g}{f}{k}{h}\in\sltwop$ then 
\[ 
\begin{pmatrix}
  g & f \\ k & h
\end{pmatrix}
\begin{pmatrix}
  1 & ax+b \\ 0 & 1
\end{pmatrix}
=
\begin{pmatrix}
  g & (ax+b)g+f \\ k & (ax+b)k + h
\end{pmatrix}
\]
Now if we have $r\ulace s$ then $r + ys \in\rup{2}$. Since $a\ge0$ we
may substitute 
\[ r+ (y+ ax+b) s = r+ (ax+b) s + y s \quad \implies\quad
r+(ax+b)s\ulace s.
\]
This shows that we have three out of the four interlacings, and since
the determinant is one, the conclusion follows from
Lemma~\ref{lem:sl-2}. Multiplying on the other side is similar.

The last one follows from
\[
\begin{pmatrix}
  1&0\\-ax-b&1
\end{pmatrix}=
\begin{pmatrix}
  0&-1\\1&0
\end{pmatrix}
\begin{pmatrix}
  1&ax+b\\0&1
\end{pmatrix}
\begin{pmatrix}
  0&1\\-1&0
\end{pmatrix}
\]

\end{proof}

\begin{prop}
  $\sltwop$ is generated by 
\[
\begin{pmatrix}
{0} & {-1}\\{1} & {0}
\end{pmatrix}\qquad
\begin{pmatrix}
{c}&{0}\\{0}&{1/c}   
\end{pmatrix}\qquad
\begin{pmatrix}
{1}&{ax+b}\\{0}&{1} 
\end{pmatrix}
\]
where $a\ge0, b\in\reals,c\ne0$.
\end{prop}
\begin{proof}
  The degree of $M=\smalltwo{g}{f}{k}{h}$ is $deg(g)+deg(h)$. If $M$
  has degree at least two then we show that there are matrices
  $A,B,M_1$ such that $M=AM_1B$ where $A,B$ are certain products of
  the first two generators, and $M_1\in\sltwop$ has lower degree.

  We can multiply on the left or right by the first two generators so
  that the leading coefficients in the top row are positive. If the
  matrix is $\smalltwo{g}{f}{k}{h}$ then $f\ulace g$, and since they
  have positive leading coefficients $f\longleftarrow g$. Thus we can
  write $f = (ax+b)g -s$ where $a\ge0$, $g\lessless s$, and $s$ has positive leading
  coefficient. We write
\[
\begin{pmatrix}
  g&f\\k&h
\end{pmatrix}
=
\begin{pmatrix}
  g & (ax+b)g-s \\ k & h 
\end{pmatrix}
=
\begin{pmatrix}
  g&-s \\k & -h_1
\end{pmatrix}
\begin{pmatrix}
  1& ax+b\\0 & 1
\end{pmatrix}
\]
where $h_1$ has positive leading coefficient, and $h = (ax+b)k -
h_1$. We show that $\smalltwo{g}{-s}{k}{-h_1}\in\sltwop$,
and by Lemma~\ref{lem:sl-1} we need to show that $k\lessless h_1$.

Now $deg(g)+deg(h_1) = deg(k)+deg(s)$, and  $deg(g)=deg(s)+1$ so we
have $deg(h_1) = deg(k)-1$. From $-gh_1 + sk = 1$ we see that $h_1$
and $g$ have the same signs on the roots of $k$. Since $g$ alternates sign
on the roots of $k$ so does $h_1$. From $deg(h_1)< deg(k)$ it follows
that $k\lessless h_1$.

If we reduce a matrix to one of degree one, then after multiplying by
appropriate generators  it has the form
$\smalltwo{c}{ax+b}{0}{1/c}$, which is easily seen to be a product of
the generators. 

\end{proof}

\begin{cor}
  $\sltwop$ is closed under multiplication.
\end{cor}

\begin{cor}\label{cor:sl-and-ma}
  If $\smalltwo{g}{f}{k}{h}\in\sltwop$ then $f + y\,g+z\,h+yz\,k\in\gsubclose_3$.
\end{cor}
\begin{proof}
  It is clear that each of the first two generators of $\sltwop$ determines a
  polynomial in $\gsubclose_3$. The last generator determines the
  polynomial $ax+b + y + z$ which in is $\gsubclose_3$ since $a$ is
  non-negative.

  To complete the proof we verify that multiplication by each of the
  generators preserves the property of being in $\gsubclose_3$. We
  start with a polynomial $F=f+yg+zh+yzk\in\gsubclose_3$. There
  are three cases:
  \begin{description}
  \item[$\smalltwo{0}{-1}{1}{0}$] This matrix transforms
    $\smalltwo{g}{f}{k}{h}$ to $\smalltwo{-k}{-h}{g}{f}$ which
    corresponds to 
    \[ -k - yh + zf+yzg = - ( k + yh + f(-z) + y(-z)g \]
    which equals $-zF(x,y,-1/z)\in\gsubclose_3$.
  \item[$\smalltwo{c}{0}{0}{1/c}$] This transforms $F$ to $c(f+yg +
    (z/c^2)h + y(z/c^2)k)=F(x,y,z/c^2)$ which is in $\gsubclose_3$.
  \item[$\smalltwo{1}{ax+b}{0}{1}$] This transforms $F$ to
\[ (ax+b)g+f + y g + z[ (ax+b)k+h] + yz\,k =
F(x,ax+b+y,z)
\]
which is in $\gsubclose_3$.
  \end{description}
\end{proof}

\begin{cor}\label{cor:sl-on-ma}
  $\sltwop$ acts on polynomials $f(x) +
  g(x)y+h(x)z+k(x)yz\in\gsubclose_3$ by matrix multiplication.
\end{cor}

  \section{Matrices of nearly quadratic polynomials}
  \label{sec:matr-init-sequ}

  In this section we consider matrices formed from polynomials in
  $\rup{3}$ that have degree at most $2$ in $y$ and $z$. We show how to
  construct such matrices of arbitrarily large degree and 
  constant determinant. 

  \begin{definition}
    $M_3 = \bigl\{ (f_{ij})_{0\le i,j\le 2} \,\mid\,
    \sum_{i,j=0}^{2} f_{ij}(x)\,y^iz^j \in\rup{3}\bigr\}$
  \end{definition}

Many operations on $F(x,y,z)=\sum f_{ij}y^iz^j$ corespond to matrix
multiplication. We only consider $\rup{3}$, but the arguments apply to all
$\rup{k}$. Let
\[
M =
\begin{pmatrix}
  f_{00} & f_{10} & f_{20} \\
  f_{01} & f_{11} & f_{21} \\
  f_{02} & f_{12} & f_{22} 
\end{pmatrix}
\]

\begin{example}
  Suppose we replace $y$ by $cy$, where $c>0$. Since
  $F(x,cy,z)\in\rup{3}$ the new matrix is in $M_3$ and equals

\[
\begin{pmatrix}
  f_{00} & c\, f_{10} & c^2\,f_{20} \\
  f_{01} &c\, f_{11} & c^2\,f_{21} \\
  f_{02} & c\,f_{12} & c^2\,f_{22} 
\end{pmatrix} 
= 
M\,
\begin{pmatrix}
  1 & . & . \\ . & c & . \\ . & . & c^2
\end{pmatrix}
\]
If we replace $z$ by $cz$ then we multiply on the left

\[
\begin{pmatrix}
  f_{00} &  f_{10} & f_{20} \\
c\,  f_{01} &c\, f_{11} & c\,f_{21} \\
c^2\,  f_{02} & c^2\,f_{12} & c^2\,f_{22} 
\end{pmatrix} 
= 
\begin{pmatrix}
  1 & . & . \\ . & c & . \\ . & . & c^2
\end{pmatrix}\,M
\]
\end{example}

\begin{example}
  Next we consider $F(x,y+a,z)$ which is in $\rup{3}$ for all $a$. The
  corresponding matrix equals

\[
\begin{pmatrix}
  f_{00} + a f_{10} + a^2 f_{20} & f_{10} + 2a\,f_{20} & f_{20}\\
  f_{01} + a f_{11} + a^2 f_{21} & f_{11} + 2a\,f_{21} & f_{21}\\
  f_{02} + a f_{12} + a^2 f_{22} & f_{12} + 2a\,f_{22} & f_{22}
\end{pmatrix}
=
M\, 
\begin{pmatrix}
  1 & . & . \\ a & 1 & . \\a^2 & 2a & 1
\end{pmatrix}
\]
If we consider $F(x,y,z+a)$ then we multiply on the left by:
\[
\begin{pmatrix}
 1 & a & a^2 \\ . & 1 & 2a \\ . & . & 1
\end{pmatrix}\, M
\]

\end{example}

\begin{example}
  If $a$ is positive then $F(x,y+ax,z)\in\rup{3}$ and
  $F(x,y,z+ax)\in\rup{3}$. These correspond to multiplying on the
  right  by
\[
\begin{pmatrix}
  1 & . & . \\ ax & 1 & . \\a^2x^2 & 2ax & 1
\end{pmatrix}
\]
or on the left by 
\[
\begin{pmatrix}
1 & ax & a^2x^2 \\ . & 1 & 2ax \\ . & . & 1
\end{pmatrix}
\]

\end{example}

\begin{example}
  If we reverse with respect to $y$ then we get
\[
\begin{pmatrix}
  f_{20} & - f_{10} & f_{00} \\
  f_{21} & - f_{11} & f_{01} \\
  f_{22} & - f_{12} & f_{01} 
\end{pmatrix}
=
M\,\begin{pmatrix}
  .&.&1 \\ .&-1&. \\ 1 & . & .
\end{pmatrix}
\]
If we reverse with respect to $z$ then we multiply on the left.
\end{example}

\begin{example}
  In order to get started we need some simple matrices. Since
  $(y+z)^2\in \rup{3}$ the following matrix is in $M_3$
\[
\begin{pmatrix}
  . & . & 1 \\ . & 2 & . \\ 1 & . & .
\end{pmatrix}
\]
and has determinant $-2$.
\end{example}

\begin{example}
  Finally, we show how to combine these matrices to make matrices in
  $M_3$ with determinant $-2$. Let
\[
C = 
\begin{pmatrix}
  1 & . & . \\ x & 1 & . \\x^2 & 2x & 1 
\end{pmatrix}
\qquad
A = \begin{pmatrix}
  .&.&1 \\ .&-1&. \\ 1 & . & .
\end{pmatrix}
\qquad
B=\begin{pmatrix}
  . & . & 1 \\ . & 2 & . \\ 1 & . & .
\end{pmatrix}
\]
The matrix $BC(AC)^n$ has degree $2n+2$. For instance, $BC(AC)^1$ equals
\[\left(
\begin{array}{lll}
 x^4-2 x^2+1 & 2 x^3-2 x & x^2 \\
 2 x^3-2 x & 4 x^2-2 & 2 x \\
 x^2 & 2 x & 1
\end{array}
\right)
\]
\end{example}

\section{$\polycpx$, a subset of $\up{1}$}
\label{sec:poly-cpx}

We now consider polynomials in $\up{1}$ with positive leading
coefficients and strict interlacing.

\begin{definition}
  \begin{gather*} 
\polycpx = \left\{ f+ g\imag \mid \text{$f,g\in\allpoly$ have positive leading
      coefficients and } f\lessless g\right\} \\
\polycpxpos = \left\{ f+ g\imag \mid \text{$f,g\in\allpolypos$ have positive leading
      coefficients and } f\lessless g\right\} \\
\polycpxclose = \left\{ f+ g\imag \mid \text{$f,g\in\allpoly$ have positive leading
      coefficients and } f\lesslesseq g\right\} \\
\polycpxf = \text{ the uniform closure of } \polycpx
\end{gather*}
\end{definition}
 
If $h(x) = \sum a_i x^i$ then we define $h_\Re = \sum \Re(a_i)x^i$
and $h_\Im = \sum \Im(a_i)x^i$. With this notation, the interlacing
condition of the definition of $\polycpx$ is that 
$h_\Re\lessless h_\Im$  where $h = f + g \imag$.

\begin{remark}
  There are two reasons we only consider $\lessless$, and not also
  $\greateq$. The first is that if we consider a product
  $\prod(x+\sigma_i) = f(x) + \imag g(x)$ where $\sigma_i\in\complexes$
  then the degree of $f(x)$ is greater than the degree of $g(x)$. 
  
  Second, if the degrees of the real and imaginary polynomials are
  equal then the leading coefficient isn't positive nor even real;
  it's complex. However, this is not  a serious problem, for we can
  multiply it to have the correct form:
  \begin{lemma}\label{lem:fix-i}
    If $f\greateq g$, $a>0$ is the leading coefficient of $f$ and
    $b>0$ is the leading coefficient of $g$ then
\[ (a-b\imag)(f+\imag g)\in\polycpx
\]
  \end{lemma}
  \begin{proof}
    The leading coefficient of $(a-b\imag)(f+\imag g)$ is
    $a^2+b^2$. Now
    \begin{gather*}
      (a-b\imag)(f+\imag g) = af+bg + \imag(ag - bf) \\
      \begin{pmatrix}
        a & b \\ -b & a
      \end{pmatrix}
      \begin{pmatrix}
        f\\g
      \end{pmatrix}=
      \begin{pmatrix}
        af+bg \\ ag - bf
      \end{pmatrix}
    \end{gather*}
The determinant is $a^2+b^2$, and $a$ is positive, so the interlacing
on the right hand side is strict   (Corollary~\ref{cor:lin-comb-new}).

It remains to determine the sign of the leading coefficient of
$ag-bf$. Write
\begin{align*}
  f &= c_0 + \cdots + c_{n-1}x^{n-1} + a\,x^n \\
  g &= d_0 + \cdots + d_{n-1}x^{n-1} + b\,x^n 
\end{align*}
Since $f\greateq g$ we know that
$\smalltwodet{d_{n-1}}{b}{c_{n-1}}{a}>0$ and this is the leading
coefficient of $ag-bf$.
  \end{proof}
  
\end{remark}

\begin{remark}
  If $F = f+ \imag g$ then $F$ has no real roots. If $r\in\reals$ were
  a root of $F$ then $f(r)+\imag g(r)=0$, which implies that
  $f(r)=g(r)=0$. This contradicts the hypothesis that $f$ and $g$ have
  no roots in common.
\end{remark}

If $f\in\allpolypos$ then there are further
restrictions on the location of the roots.

  \begin{lemma}
    If $f,g\in\allpolypos$, $f\lessless g$, and $f,g$ have positive
    leading coefficients then the roots of $f+\imag g$ have negative
    real part and negative imaginary part.
  \end{lemma}
  \begin{proof}
    We know the imaginary part is negative. If $f = \prod_k(x-r_k)$ and
    $f(\alpha) + \imag g(\alpha)=0$ then there are positive $a_k$ such
    that
    \begin{align*}
     \sum_j a_j\left( f/(x-r_j)\right) & = g(x) \\
      \prod_k(\alpha-r_k) + \imag \sum_j a_j\frac{1}{\alpha-r_j}\prod_k(\alpha-r_k) &= 0
      \\
      1 + \imag \sum_j a_j /(\alpha-r_j) &= 0\\
\intertext{Taking the imaginary part yields}
\Re(\alpha) \sum_j  \frac{a_j}{|\alpha-r_j|^2} &= \sum_j a_j \frac{r_j}{|\alpha-r_j|^2}
    \end{align*}
Since the $a_j$'s are positive and the $r_j$'s are negative, it
follows that the real part is negative.
  \end{proof}

It's easy to construct polynomials in $\polycpx$ from polynomials in $\allpoly$:

  \begin{lemma}\label{lem:polycpx-basic}
    Suppose that $f\in\allpoly$, $\sigma$ is in the upper half plane
    and $\alpha$ is positive. Then
    \begin{enumerate}
    \item $f(x+\sigma)\in\polycpx$.
    \item $\displaystyle\int_0^1 f(x+\alpha\imag t)\,dt \in\polycpx$.
    \item $\displaystyle\int_0^1 f(x+\alpha\imag t)\,dt \in\polycpx$.
    \item $\imag^{n}\,f(-\imag x)\in\polycpx$ if $f\in\allpolypos(n)$.
    \item If $g(x)\in\polycpx$ then $g(\alpha x)\in\polycpx$ for
      positive $\alpha$.
    \item If $f\in\allpoly$ has all distinct roots then $f(x)+ \imag
      f'(x)\in\polycpx$. 
    \end{enumerate}
  \end{lemma}
  \begin{proof}
    The roots of $f(x+\sigma)$ all have imaginary part $-\Im(\sigma)$
    which is negative.  Corollary~\ref{cor:int-of-i-2} shows that all
    the roots of the integral have imaginary part equal to
    $-\alpha/2$, and so the integral is in $\polycpx$.

If $f = x^n + a x^{n-1} + \cdots$ then 
$ \imag^{n}\,f(-\imag x) = x^n + \imag a x^{n-1} + \cdots$
Since $f(x)\in\allpolypos$, all the roots of $f(-\imag x)$ lie in the
lower half plane. Thus the roots are in the correct location, and the
leading coefficients are positive. 

The remaining ones are obvious.
  \end{proof}

A   useful property of $\polycpx$ is the following, which follows
from the fact that $\polycpx\subset\up{1}$, and the leading
coefficient of the product is positive.

\begin{lemma}\label{lem:cpxpoly-1}
  $\polycpx$ and $\polycpxf$ are  closed under multiplication.
\end{lemma}

\begin{remark} \label{rem:i-mult}
  We can also prove that $\polycpx$ is closed under multiplication
  using properties of matrices.  It suffices to show that if $f+g
  \imag $ and $h+k\imag$ are in $\polycpx$, then so is their product
  $(fh-gk) + (fk+gh)\imag$.  Notice that both $fh-gk$ and $fk+gh$ have
  positive leading coefficients, and that the degree of $fh-gk$ is one
  more than the degree of $fk+gh$. If we write these terms as a matrix
  product
$$
\begin{pmatrix}
  fh-gk \\ fk+gh
\end{pmatrix} =
\begin{pmatrix}
  h & -k \\ k & h
\end{pmatrix}
\begin{pmatrix}
  f \\ g
\end{pmatrix}
$$
then the conclusion follows from Lemma~\ref{cor:2by2-a}.
  
\end{remark}

Here is a simple consequence.

\begin{cor}\label{cor:fg-imag}
  If $f\lessless g$ and $n$ is a positive integer then
  \begin{align*}
    \text{If}\ n=2m \hspace*{1in} & \\
    \sum_{k=0}^m f^{2k}g^{2m-2k}\binom{2m}{2k}(-1)^k &
    \lessless &
\sum_{k=0}^m f^{2k+1} g^{2m-2k-1} \binom{2m}{2k+1} (-1)^k \\
    \text{If}\     n=2m+1 \hspace*{1in} &\\ \sum_{k=0}^m f^{2k}g^{2m+1-2k}\binom{2m+1}{2k}(-1)^k &
    \lessgreat &
\sum_{k=0}^m f^{2k+1} g^{2m-2k} \binom{2m+1}{2k+1} (-1)^k \\
  \end{align*}
\end{cor}
\begin{proof}
  Expand $(f+\imag g)^n$.
\end{proof}

The following corollary is useful.

\begin{cor}\label{cor:cpxclose-lt}
  If $T\colon\allpoly\longrightarrow\allpoly$ preserves the degree
  and sign of the leading coefficient then
  $T\colon{}\polycpxclose\longrightarrow\polycpxclose$.  If $T$ also
  preserves strict interlacing then
  $T\colon{}\polycpx\longrightarrow\polycpx$.
\end{cor}
\begin{proof}
  Since $T$ is a transformation that only involves real coefficients,
  we know that
\[ \Re(T(f)) = T(\Re(f)) \text{ and }\Im(T(f)) = T(\Im(f)). \]

\end{proof}

\begin{cor}\label{cor:cpxpoly-fp}
  These linear transformations map $\polycpx$ to itself.
  \begin{enumerate}
  \item $f\mapsto f+ \alpha f'$ where $\alpha\in\reals$.
  \item $f\mapsto g(\diffd) f$ where  $g(x)\in\allpoly$
  \item $f\mapsto \expoper{} (f)$.
  \item $f \mapsto H_n$.
  \end{enumerate}
\end{cor}
\begin{proof}
  The only observation required is that the second one follows from
  the first.
\end{proof}

\begin{remark}
  Note that the first transformation is equivalent to the statement
  that if $g \lessless h$ then $g + \alpha g' \lessless h +
  \alpha h'$. This is Lemma~\ref{lem:inequality-1}.
\end{remark}

\begin{cor}
If $f\in\polycpx$ then
\begin{enumerate}
\item $f(x+\imag) + f(x-\imag)\in\polycpxclose$
\item $(1/\imag)\left(f(x+\imag) - f(x-\imag)\right)\in\polycpxclose$
\end{enumerate}
\end{cor}
\begin{proof}
  The first one is $2\cos(\diffd)f$, and the second is $2\sin(\diffd)f$.
\end{proof}

  \begin{lemma}
    Suppose that $S_1,S_2:\allpoly\longrightarrow\allpoly$ are linear
    transformations that preserve $\lesslesseq$ and the sign of the
    leading coefficient, and satisfy $S_1(f)\lesslesseq S_2(f)$ for
    all $f\in\allpoly$. Define $T(f) = S_1(f) + \imag S_2(g)$. Then,
    \begin{itemize}
    \item $T\colon{}\allpoly\longrightarrow\polycpxclose$
    \item $T\colon{}\polycpxclose\longrightarrow\polycpxclose$
    \end{itemize}
  \end{lemma}

  \begin{proof}
    The first is immediate from the hypotheses of the Lemma. 
    Since $T(f+\imag g) = S_1(f)- S_2(g) + \imag(S_1(g) + S_2(f))$
    the degrees and the signs of the leading coefficients are
    correct. We have
    \begin{align*}
      S_1(f) - S_2(g) & \lesslesseq S_1(g) \text{ since } S_1(f)\lesslesseq
      S_1(g)\lesslesseq S_2(g)\\
       & \lesslesseq S_2(f) \text{ since } S_1(f)\lesslesseq
      S_2(f)\lesslesseq S_2(g)\\ 
    \end{align*}
      Adding these two interlacings gives the desired interlacing.
  \end{proof}

    Here are some examples of linear transformations satisfying the
    hypothesis of the lemma. Suppose that
    $U:\allpoly\longrightarrow\allpoly$ preserves $\lesslesseq$ and
    the sign of the leading coefficient. 
    \begin{alignat*}{2}
      S_1(f) &= U(xf)\quad\quad & S_2(f) &= U(f)\\
      S_1(f) &= U(f) & S_2(f) &= U(f)'\\
      S_1(f) &= U(f) & S_2(f) &= U(f')
    \end{alignat*}

\begin{remark}
  Here is a different proof of part of the Hermite-Biehler theorem.
  We  prove that if $f\lessless g$ and the leading coefficients of
  $f$ and $g$ are positive, then $h(x)=f(x)+\imag g(x)$ has all roots
  in $-\plane$. First of all, notice that $h(z)$ can not have any real
  roots, since if $h(\alpha)=0$ and $\alpha$ is real then
  $f(\alpha)=g(\alpha)=0$ which contradicts $f\lessless g$. Next, it
  follows from from Example~\ref{ex:powerofi} that $
  (x+\imag)^n=f_\Re+g_\Im\imag$ is in $\polycpx$ since $f_\Re\lessless
  g_\Im$, and $ (x+\imag)^n\in\up{1}$. Since any pair of
  strictly interlacing polynomials can be reached by a path of
  strictly interlacing polynomials, it follows that any pair of
  strictly interlacing polynomials with positive leading coefficients
  must have all roots in the lower half plane, since the roots can
  never be real.

\end{remark}

  The real parts of the roots of a polynomial in $\polycpx$ are
  constrained by the location of the roots of the real part of the
  polynomial. See   Question~\ref{ques:affine-and-i}. 

\begin{lemma}
  Suppose that $f = f_\Re(x)+\imag f_\Im(x)\in\polycpx$, and assume that all the
  roots of $f_\Re$ lie in an interval $I$. Then, the real parts of the
  roots of  $f$ also lie in $I$.
\end{lemma}
\begin{proof}
  We may assume $f_\Re$ is monic, and write $f_\Re(x) = \prod(x-r_k)$,
  $f_\Im(x) = \sum a_k f_\Re/(x-r_k)$ where the $a_k$ are non-negative. We
  show that if $\Re(\sigma)$ is greater than all $r_k$ then
  $f(\sigma)\ne0$. The case where $\Re\sigma$ is less than all the roots
  is similar. Dividing by $f_\Re$
\[
0=1 + \imag \sum_k a_k \frac{1}{\sigma-r_k} 
\]
Since $\Re\sigma>r_k$, all terms in the sum have positive real part, and
so $\sigma$ is not a zero of $f$.
\end{proof}

\begin{figure}
    \centering
{
        \psset{xunit=.5cm,yunit=.5cm}
    \begin{pspicture}(0,0)(12,6)
      \psline(0,5)(10,5)
      \pscurve(1,5)(1.5,4)(2,5)
      \pscurve(3,5)(3.5,3)(4,5)
      \pscurve(5,5)(5.5,2)(6,5)
      \pscurve(7,5)(6.8,3)(5.8,1.8)(5.5,0)
      \pscurve(8,5)(8.5,3)(9,5)
      \rput(12,5){Real axis}
      \rput(1,5.5){f}
      \rput(2,5.5){g}
      \rput(3,5.5){f}
      \rput(4,5.5){g}
      \rput(5,5.5){f}
      \rput(6,5.5){g}
      \rput(7,5.5){f}
      \rput(8,5.5){g}
      \rput(9,5.5){f}
    \end{pspicture}
}
    \caption{Trajectories for $f+t\,\imag\,g$}
    \label{fig:i-trajectory}
  \end{figure}
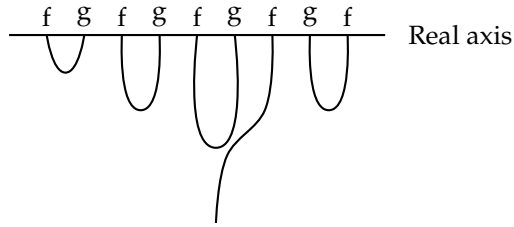

Suppose that $f\lessless g$ and consider the plot
(Figure~\ref{fig:i-trajectory}) of the roots of $F_t=(1-t)f + \imag\,t\, g$
for $0\le t \le 1$. From the lemma we know that the real part of
the roots of $F_t$ lies in the interval determined by the roots of $f$. At $t=0$
the roots are the roots of $f$, and as $t$ increases the roots move to
the roots of $g$, except for one root whose imaginary part goes to $-\infty$.

Here's a variation of the lemma where we start with polynomials of
equal degree.
  
  \begin{lemma}
    Suppose that $g\weaki f$ where $f$ and $g$ have the same degree.
    If $f+\imag g = \prod(x-r_k-\imag s_k)$, $f$ has roots
    $a_1\le\cdots \le a_n$, $g$ has roots $b_1\le\cdots\le b_n$ then
\[
r_k \in\bigl[ \frac{a_1+b_1}{2},\frac{a_n+b_n}{2}\bigr] 
\]
  \end{lemma}
  \begin{proof}
    Suppose that $\alpha+\imag \beta$ is a root of $f+\imag g$. Then
    \begin{gather*}
      \prod(\alpha+\imag \beta - a_k) + \imag \prod(\alpha+\imag\beta
      - b_k) = 0 \\
      \prod\frac{\alpha+\imag\beta - b_k}{\alpha+\imag\beta - a_k} =
        \imag \\
      \prod\frac{(\alpha-b_k)^2 + \beta^2}{(\alpha-a_k)^2 + \beta^2} =
        1 
    \end{gather*}
    Since $g\weaki f$ we have that $b_k\le a_k$. Thus, if $b_k<a_k$
    and $\alpha < (a_k+b_k)/2$ then the $k$th term is less than one,
    and if $\alpha > (a_k+b_k)/2$ then the $k$th term is greater than
    $1$.  This implies the result.
  \end{proof}

\section{Some simple consequences}
\label{sec:some-simple-cons}
  
  \begin{cor} \label{cor:imag-1}
    If $f(x,y)\in\rup{2}$  and $\sigma$ has positive imaginary part
     then $f(x,\sigma)\in\polycpx$.
  \end{cor}
  \begin{proof}
Since $\rup{2}\subset\up{2}$ substitution of $\sigma$ yields a
polynomial in $\up{1}$, and this is in $\polycpx$ since the
coefficient of $x^n$ is positive.
   \end{proof}

  \index{even and odd parts!in $\rup{2}$}
  \begin{cor}
    If $f(x,y)\in\rup{2}$ and 
$$
f(x,y) = f_0(x) + f_1(x)y+f_2(x)y^2 + \cdots$$
then
$$
f_0 - f_2 + f_4 - \cdots \lesslesseq f_1 -f_3+f_5-\cdots
$$
  \end{cor}
  \begin{proof}
    This is just a restatement of Corollary~\ref{cor:imag-1}, where
    $f_\Re$ is the left hand side, and $f_\Im$ is the right hand side.
  \end{proof}

\begin{lemma}
  If $f=\displaystyle\sum_{ij} a_{ij}x^iy^j\in\gsubpos_2(2n)$ then
  $\displaystyle\sum_{i\equiv j\pmod{2}}(-1)^{(n+i+j)/2} a_{ij}x^iy^j$
   is in $\gsubpos_2$
\end{lemma}
\begin{proof}
  Since $g(x,y)=\imag^n\,f(-\imag x,-\imag
  y)\in\up{2}$,  the expression of the conclusion is the real part
  of $g$.
\end{proof}

\begin{cor}
  If $f\lessless g$ in $\rup{2}$ then the matrix
  $\smalltwo{f}{-g}{g}{f}$ maps a pair of interlacing polynomials in
  $\rup{2}$ to a pair of interlacing polynomials in $\gsub_2$. 
\end{cor}
\begin{proof}
  The action of the matrix on $\smalltwobyone{h}{k}$ is the same as the
  multiplication of $h+\imag k$ by $f+\imag g$. 
\end{proof}

\begin{cor}
Suppose that $f(x,y)\in\rup{2}$. Then
\begin{enumerate}
\item $\frac{1}{\imag}\left( f(x,y+\imag) - f(x,y-\imag)\right) \in \up{2}$
\item $ f(x,y+\imag) + f(x,y-\imag) \in \up{2}$
\item $ f(x+\imag,y-\imag) + f(x-\imag,y+\imag) - f(x+\imag,y+\imag) -
  f(x-\imag,y-\imag)  \in \up{2}$
\end{enumerate}
\end{cor}
\begin{proof}
 The first one is $\sin(\diffd)$ and the second is $\cos(\diffd)$
 applied to the $y$ variable. If we apply $\sin$ to (1)
 using differentiation with respect to $x$, then we get (3).
\end{proof}

\begin{lemma}\label{lem:gsub-sub-i}
  If $f(x,y)\in\rup{2}$ then 
 $\displaystyle\int_0^1f(x,y+\imag t)\,dt\in\rhb{2}$.
\end{lemma}
\begin{proof}
  If we write the integral as $g+\imag h$ where $g,h$ have real
  coefficients, then 
\[\int_0^1f(\alpha,y+\imag t)\,dt = g(\alpha,y) + \imag h(\alpha,y)\]
Since the left hand side is in $\polycpx$, we see that $g(\alpha,y)
\lessless h(\alpha,y)$ for all $y$. It's clear that the
homogeneous parts of $g$ and $h$ have all positive coefficients, and
so they are in $\rup{2}$ and interlace.
\end{proof}

\begin{remark}
    If $f \lesslesseq g$ in $\polycpx$ then we can write $g$ in terms
    of the roots of $f$, but the coefficients
    might be complex. For instance, suppose
    \begin{align*}
      D_1 & = \left(
\begin{array}{llll}
 95 & 40 & 78 & 61 \\
 40 & 35 & 39 & 48 \\
 78 & 39 & 70 & 60 \\
 61 & 48 & 60 & 70
\end{array}
\right) &
D_2 &= \left(
\begin{array}{llll}
 90 & 65 & 73 & 77 \\
 65 & 57 & 65 & 68 \\
 73 & 65 & 83 & 73 \\
 77 & 68 & 73 & 90
\end{array}
\right) \\
S &=\left(
\begin{array}{llll}
 3 & 6 & 3 & 3 \\
 6 & 1 & 2 & 7 \\
 3 & 2 & 6 & 5 \\
 3 & 7 & 5 & 6
\end{array}
\right) 
& f &= |x I + y D_1 + S + \imag D_2| 
    \end{align*}
Since $D_1,D_2$ are positive definite, and $S$ is symmetric, we know
that $f \in\rhb{2}$. If  $f = f_0(x) + f_1(x)y+\cdots $
then $f_0 \lesslesseq f_1$. If $\roots{(f_0)}=(r_i)$ then we can write
\begin{align*}
  f_1 &= \sum_1^4 a_i\,\frac{f_0}{x-r_i} & \text{where}\\
a_1 &= -17.-292. \imag &
a_2 &= -1.38-13.2 \imag \\
a_3 &= -0.775-12.3 \imag &
a_4 &= 3.17-2.45 \imag
\end{align*}
  \end{remark}

\section{Interlacing in $\polycpx$}

We make the usual definition of interlacing in $\polycpx$ (closure
under linear combinations), and it turns out that we have already seen
the definition! We then gather some properties of interlacing in
$\polycpx$. We begin with the case of unequal degrees.

\begin{definition}
  If $f,g\in\polycpx$, then $f \lessless g$ if and only if $f +
  \alpha g\in\polycpx$ for all $\alpha\in\reals$.
\end{definition}

Suppose that $f = f_\Re + \imag f_\Im$ and $G = g_\Re + \imag g_\Im$.
$f\lessless g$ is equivalent to
$f_\Re + \imag f_\Im + \alpha(g_\Re + \imag g_\Im)\in\polycpx$, which,
if we express this in terms of the definition of $\polycpx$ is
$$
f_\Re + \alpha g_\Re < f_\Im + \alpha g_\Im.$$

\noindent
It follows from
Lemmas~\ref{lem:inequality-4} and \ref{lem:inequality-4b} that we have

\begin{lemma}\label{lem:i-interlace}
Suppose that  $f = f_\Re + \imag f_\Im$ and $g = g_\Re +
  \imag g_\Im$ are in $\polycpx$. Then 
  $f\lessless g$ if and only if
  \begin{enumerate}
  \item $f_\Re \lessless g_\Re $
  \item $f_\Im \lessless g_\Im $
   \item $\smalltwodet{f_\Re}{f_\Im}{g_\Re}{g_\Im}<0$
  \end{enumerate}
\end{lemma}

Such polynomials exist. For instance, 
\[
    (8 + 14x + 7x^2 + x^3)+\imag(22 + 25x + 6x^2)  \lessless
    (16 + 17x + 4x^2)+\imag(29 + 15x) 
\]

These polynomials arise from a general construction - see
Lemma~\ref{lem:make-i}. 

Here are a few simple properties of $\lessless$ for $\polycpx$. These
are the $\polycpx$-analogs of the usual interlacing properties in
$\allpoly$. Note that multiplication preserves strict
interlacing. Many of these can also be easily proved using the
non-vanishing definition of $\polycpx$, but it interesting to see that
we can prove them just using properties of interlacing polynomials.

\begin{lemma} \label{lem:i-analog}
Assume $f,g,h\in\polycpx$.

  \begin{enumerate}
  \item  $f\lessless f'$.
  \item If $f\lessless g$ and $f\lessless h$ then $f\lessless \alpha g +
    \beta h$ for positive $\alpha,\beta$. In particular, $\alpha f +
    \beta g\in\polycpx$.
  \item If $f \lessless g \lessless h$ then $f - h \lessless g$.
  \item If $f\lessless g$ and $\Im(\sigma)>0$ then $f+\sigma g\in\polycpx$.
  \item If $f\lessless g$   then $fh\lessless gh$.
  \item If $f\in\polycpxpos$ then $f(-\diffd)g\in\polycpx$.
  \item If $f\in\polycpxpos$ and $f^{rev} = f_\Re^{rev}
    +\imag\,f_\Im^{rev}$ then $f^{rev}\in\polycpxpos$.
  \item Suppose that $f$ factors as $\prod(x-\sigma_k)$.
    \begin{enumerate}
    \item $\dfrac{f}{x-\sigma_k}\in\polycpx$
    \item $f \lessless \sum a_k \dfrac{f}{x-\sigma_k}$
      for any non-negative $a_k$.
    \end{enumerate}
  \end{enumerate}
\end{lemma}

\begin{proof}
  If $f = f_\Re + \imag f_\Im$ then $f'  = f_\Re' + \imag
  f_\Im'$. Since $f_\Re \lessless f_\Im$, the interlacing conditions
  hold. The determinant condition is Lemma~\ref{lem:inequality-1}. 
  
  In order to see that $f\lessless \alpha g + \beta h$ we have to
  check the three conditions, and that $\alpha g+\beta h\in\polycpx$.
  The first two follow from the additivity of interlacing, and the
  third one follows from the linearity of the determinant:

$$ 
\begin{vmatrix}
  f_\Re & f_\Im \\ \alpha g_\Re+\beta h_\Re & \alpha g_\Im + \beta h_\Im
\end{vmatrix}
=
\alpha
\begin{vmatrix}
  f_\Re & f_\Im \\  g_\Re &  g_\Im 
\end{vmatrix}
+
\beta
\begin{vmatrix}
  f_\Re & f_\Im \\  h_\Re &  h_\Im
\end{vmatrix}
$$

Since the real parts of $f$ and $\alpha g+\beta h$ interlace, as do
the imaginary parts, and the determinant is negative, we can apply
Lemma~\ref{lem:square-1}. This shows that the real and imaginary parts
of $\alpha g+\beta h$ interlace, and hence it is is in $\polycpx$. 

If $f\lessless g\lessless h$ then the argument is nearly the same as
the previous one.

If we write $\sigma = a +b\imag$ then $f+\sigma g = (f+a g)+ \imag
(bg)$. Since $f+a g\lessless bg$ it suffices to consider $f+\imag g$. 
Expanding into real and imaginary parts
\begin{gather*}
  f+\imag g = (f_\Re - g_\Im) + \imag(f_\Im + g_\Re)\\
\intertext{and using the interlacings 
$f_\Re \lessless f_\Im + g_\Re \lessless g_\Im$ yields}
\Re(f+\imag g) = f_\Re - g_\Im \lessless f_\Im + g_\Re = \Im(f+\imag g)
\end{gather*}

Since $f\in\polycpxpos$ both $f_\Re$ and $f_\Im$ have all negative
roots, and consequently $f_\Re(-x)f_\Im(-x)$ has negative leading
coefficient, and positive constant term. Consequently,
Proposition~\ref{prop:fofd-direction} implies that we have the
interlacing square 

  \centerline{\xymatrix{
      f_\Im(-\diffd)g_\Re 
      \ar@{->}[d]_{{    }}           
      \ar@{<-}[rrr]^{{   }}         
      &&&
      f_\Im(-\diffd)g_\Im 
      \ar@{->}[d]^{{    }} \\        
      f_\Re(-\diffd)g_\Re 
      \ar@{<-}[rrr]^{{    }}         
      &&&
      f_\Re(-\diffd)g_\Im
}}

It follows that 
\[
      f_\Re(-\diffd)g_\Re - f_\Im(-\diffd)g_\Im \lessless
      f_\Re(-\diffd)g_\Im + f_\Im(-\diffd)g_\Re
\]

Since all the roots of the real and
imaginary part are positive, reversal preserves interlacing.

Since $f<g$ we know that $f+\alpha g\in\polycpx$ for all real $\alpha$.
Thus $hf + \alpha hg = h(f+\alpha g)\in\polycpx$ since $\polycpx$ is
closed under multiplication.

Since interlacing is preserved by addition it will suffice to show
that $f\lessless f/(x-\sigma_k)$. If $\sigma_k = r_k + \imag s_k$ then
we compute
$$ f + \alpha\frac{f}{x-r_k-\imag s_k} =
\frac{f}{x-r_k-\imag s_k}\,\cdot\,(x-r_k-\imag s_k + \alpha).
$$
Since both factors are in $\polycpx$, so is their product, and thus we
have shown interlacing.
\end{proof}

There are two choices for the Hadamard product of polynomials in
$\polycpx$. Suppose that $ f = \sum \alpha_ix^i$ and $g = \sum
\beta_ix^i$.
  \begin{enumerate}
  \item Same definition as before: $f\ast g = \sum \alpha_i\beta_i x^i$.
  \item Separate the real and complex parts: $f\ast^{\imag} g = f_\Re\ast g_\Re +
    \imag\,f_\Im\ast g_\Im$.
  \end{enumerate}
Surprisingly, each one preserves $\polycpx$.

\index{Hadamard product!in $\polycpx$}

  \begin{cor}\label{cor:had-prod-i}
    If $f\in\polycpx$, $g\in\polycpxpos$ then $f\ast^{\imag} g\in\polycpxpos$.
  \end{cor}
  \begin{proof}
    This is Lemma~\ref{lem:hadamard-interlace}.
  \end{proof}

  \begin{cor}\label{cor:had-prod-i-2}
    If $f,g\in\polycpxpos$ then $f\ast
    g\in\polycpx$.\footnote{Unlike $\allpoly$, it's
      not true that the Hadamard product $\ast$  maps
      $\polycpx\times\polycpxpos\longrightarrow\polycpx$} 
  \end{cor}
  \begin{proof}
Write  in terms of real and imaginary parts:
\[
f\ast g = \bigl(f_\Re\ast g_\Re - f_\Im\ast g_\Im\bigr) +
\imag\bigl(f_\Re\ast g_\Im + f_\Im \ast g_\Re\bigr)
\]
Since Hadamard product preserves interlacing in $\allpolypos$, the
interlacing of the two sides follows from the interlacings
\[
f_\Re\ast g_\Re \lessless
f_\Re\ast g_\Im + f_\Im \ast g_\Re \greateq
 f_\Im\ast g_\Im
\]
The result is only in $\polycpx$ since $f_\Re\ast g_\Re - f_\Im\ast
g_\Im$ might have negative signs.
  \end{proof}

We can construct interlacing polynomials in $\polycpx$ from
polynomials in $\gsubplus_2$ and $\gsubplus_3$.

\begin{lemma}\label{lem:make-i}  \ 
\begin{enumerate}
\item   If $f_0(x)+f_1(x)y+ \cdots \in\gsubplus_2$ and $f_0$ has all
  distinct roots then 
$ f_i + \imag f_{i+1} \lessless f_{i+1} + \imag f_{i+2}$.\\
\item If $\sum f_{i,j}(x)y^iz^j\in\gsubplus_3$ and $f_{0,0}$ has all
  distinct roots  then 
  \begin{align*}
    f_{i,j} + \imag\, f_{i+1,j} & \lessless f_{i,j+1} + \imag\,
    f_{i+1,j+1} \\
    f_{i,j} + \imag\, f_{i,j+1} & \lessless f_{i+1,j} + \imag\,
    f_{i+1,j+1} 
  \end{align*}
\end{enumerate}
\end{lemma}
\begin{proof}
  Since $f_0 \lessless f_1$ by Lemma~\ref{lem:p2-strict}, we know from
  Lemma~\ref{lem:p4-only-way} that
  $\smalltwodet{f_i}{f_{i+1}}{f_{i+1}}{f_{i+2}}<0$.  From
  Corollary~\ref{cor:fct-interlace} we see that $f_i\lessless f_{i+1}
  \lessless f_{i+2}$. This establishes the first part.

For the second part, since $f_{0,0}$ has all distinct roots, it follows
that all interlacings are strict. Thus, in order to verify the
interlacing in $\polycpx$ we need to show that the determinant 
$\smalltwodet{f_{i,j}}{f_{i+1,j}}{f_{i,j+1}}{f_{i+1,j+1}}$ is
positive for all $x$. If we substitute $\alpha$ for $x$, then
Corollary~\ref{cor:p2-strict-quad} says that all quadrilateral
inequalities are strict if $f(\alpha,0,0)\ne 0$. Moreover, if they
aren't strict, then $f(\alpha,y,z)$ satisfies the hypothesis of
Lemma~\ref{lem:p2-strict-quad-1}, and so there is a triangular region
of zeros. It follows that $f_{0,0}(x)$ and $f_{1,0}(x)$ have a common
root, which is a contradiction.

\end{proof}
\begin{example}
Consider an example. Write \[(1 + x + y + z) (2 + x + 2 y + z) (4
+ x + 3 y + 2 z)\] as any array  with the usual
horizontal and vertical interlacings. 

 \centerline{ \xymatrix@=.4cm{
      2\ar@{->}[d]_{}  \\
    10+5x  \ar@{->}[d]_{} \ar@{<-}[r] &
    9 \ar@{->}[d]  \\
    16 + 17\,x + 4\,x^2   \ar@{->}[d]_{}\ar@{<-}[r] & 
    29 + 15\,x      \ar@{->}[d] \ar@{<-}[r] & 
    13 \ar@{->}[d]  \\ 
    8 + 14\,x + 7\,x^2 + x^3  \ar@{<-}[r] & 
    22 + 25\,x + 6\,x^2  \ar@{<-}[r] &
    20+11x  \ar@{<-}[r]&  
    6  \\
  }}
\ \\[.1cm]
\noindent%
Each vertical and horizontal arrow determines a polynomial
in $\polycpx$. If we only consider the horizontal arrows then we have
the following interlacings in $\polycpx$:

 \centerline{ \xymatrix@=.4cm{
    (10+5x)+\imag(9)  \ar@{->}[d]_{}   \\
    (16 + 17x + 4x^2)+\imag(29 + 15x) \ar@{->}[d]_{}\ar@{<-}[r] & 
    (29 + 15x)+\imag(13)      \ar@{->}[d]   \\ 
    (8 + 14x + 7x^2 + x^3)+\imag(22 + 25x + 6x^2)  \ar@{<-}[r] & 
    (22 + 25x + 6x^2)+\imag(20+11x) \\
    &  \ar@{->}[u]     (20+11x)+\imag(6) 
      \\
  }}
\end{example}

  \begin{lemma}
    \label{lem:i-integral}
    If $f(x)\in\polycpx$ then $\displaystyle \int_{0}^1 f(x+\imag
    t)\,dt\in\polycpxclose$. 
  \end{lemma}
  \begin{proof}
    If we write $f = f_\Re + \imag f_\Im$ then we know
\begin{align*}
F_\Re &= \int_{0}^1 f_\Re(x+\imag t)\,dt \in\polycpx & 
F_\Im &= \int_{0}^1 f_\Im(x+\imag t)\,dt \in\polycpx 
\end{align*}
Since $f_\Re+\alpha f_\Im\in\allpoly$ for all $\alpha$ it follows
that $F_\Re+ \alpha F_\Im\in\polycpx$ for all $\alpha\in\reals$. Thus,
$F_\Re\lesslesseq F_\Im$ in $\polycpxclose$, and therefore
$\int_{0}^1f(x+\imag t)=F_\Re + \imag
F_\Im\in\polycpxclose$. 
\end{proof}



  \section{Particular interlacings}

It is easy to describe all the polynomials that interlace a given
polynomial in $\allpoly$. This is an unsolved problem in $\polycpx$. 
We  have a number of examples and simple observations.

\begin{example}
  In $\allpoly$, $(x-a)^2\lesslesseq x-b$ if and only if $a=b$.  The
  first surprise is that there are many polynomials interlacing
  $(x+\imag)^2$.  We show that $\sigma$ satisfies
  $(x+\imag)^2\lessless x-\sigma$ in $\polycpx$ if and only if
  $\sigma$ lies in the disk of radius $1$ centered at $-\imag$.

To see this, write $\sigma = u+\imag v$. The interlacing requirements
 are met if $|u|<1$ and $v>0$. The determinant requirement
is that
\[
0 > \begin{vmatrix}
  x^2 -1 & 2x \\ x+ u & v
\end{vmatrix} = x^2(v-2) -2xu-v 
\]
This means that the discriminant is non-positive, and thus
$u^2+(v+1)^2\le 1$. The general quadratic appears to also have a
simply stated answer - see Question~\ref{ques:i-interlace}.
\end{example}

\begin{example}
  Another difference from interlacing in $\allpoly$ is that
  \[
(x+\imag)^2\greateq (x+\imag+\alpha)^2 \quad\text{if and only if}\quad
  0<\alpha<2.
\] 
Note that
\[
(x+\imag+\alpha)^2 = (x+\imag)^2 + 2\alpha(x+\alpha/2+\imag)
\]
In order to have interlacing we need that $\alpha>0$ and that
$(x+\imag)^2 \lessless (x+\alpha/2+\imag)$. The latter happens exactly
when $|\alpha|<2$ by the first example

More generally, if $x+\imag)^2 \greateq (x+\alpha + \beta \imag)^2$,
and $\beta$ is not zero, then 
\[
(x+\imag+\alpha)^2 - (x+\imag)^2 = 2(\alpha \imag -\imag +\beta)x +
\cdots
\]
If $\alpha\ne1$ then the difference is not in $\polycpx$.
\end{example}

\begin{example}
  If we replace $x$ by $\beta x$ in the first two examples, and
  observe that interlacing is preserved under $x\mapsto x+r$ where
  $r\in\reals$,  then we see that
  \begin{enumerate}
  \item $(x+\sigma)^2 \lessless (x+\tau)$ iff $|\sigma-\tau| < 2\,\Im(\sigma)$.
\item $(x+\sigma)^2 \greateq (x+\sigma+\alpha)^2$ iff $0<\alpha<2\,\Im(\sigma)$.
  \end{enumerate}
  As $\sigma$ approaches the real line, both interlacings converge to
  the usual interlacings in $\allpoly$. 
\end{example}

\begin{example}
  Here are some of the properties of the polynomials interlacing  $(x+\imag)^n$.
  If $(x+\imag)^n\lessless f(x)$ and we differentiate $n-2$ times then
  $(x+\imag)^2 < f^{(n-2)}$.  Now the single root of $f^{(n-2)}$ is
  the average of the roots of $f$.  Consequently we conclude:
  \begin{quote}
    The average of the roots of $f$ lies in the unit disk centered at $\imag$.
  \end{quote}

\noindent%
Next, we claim that if $(x+\imag)^n \lessless (x+r \imag)^{n-1}$ where
$r$ is real, then $0<r<\frac{n}{n-1}$.  To see this, let
$f=(x+\imag)^n$ and $g= (x+r \imag)^{n-1}$, and note that
\begin{align*}
  f_\Re &= x^n + \cdots  & f_\Im &= n x^{n-1} + \cdots \\
  g_\Re &= x^{n-1} + \cdots  & f_\Im &= (n-1) x^{n-2}r +\cdots\\
\end{align*}
so  coefficient of $x^{2n-2}$ in $f_\Re g_\Im - f_\Im g_\Re$ is
$(n-1)r-n$. Since the determinant is negative, the leading coefficient
must be negative, and so $r \le n/(n-1)$. 

If $(x+\sigma)^n \greateq (x+\tau)^n$ then differentiating $n-2$ times
shows that $(x+\sigma)^2\greateq(x+\tau)^2$, and so $0<\tau-\sigma<2$. 

\end{example}

\begin{example}
  Given $n$, we can determine an $r$ so that if $s\in(r,n/(n-1))$ then
  $(x+\imag)^n\lessless (x+s\, \imag)^{n-1}$. To do this, we compute the
  discriminant of the determinant. This is a polynomial in $s$. The
  intervals between consecutive distinct roots of the discriminant have
  the property that the determinant has no multiple roots for $s$ in
  the interval. If all roots of the determinant are complex, and the
  leading coefficient is negative then the determinant is everywhere
  negative. 

  Once we have these intervals then it is easy to verify that the real
  and complex parts interlace for values of $s$ in these intervals.
  Table~\ref{tab:intervals-i} lists a few of these intervals.

\begin{table}[htbp]
\[
  \begin{array}{rll}
    \toprule
    n & \multicolumn{2}{c}{\text{interval}}\\
    \midrule
    2 & (0,&2) \\
    3 & (1/9,&3/2) \\
    4 & (.23,&4/3)\\
    5 & (.34,&5/4) \\
    6 & (.42,&6/5)
\end{array}
\]
  \caption{Intervals where $(x+\imag)^n\lessless (x+ s\, \imag)^{n-1}$}
  \label{tab:intervals-i}
\end{table}
  
\end{example}

  \section{The geometry of interlacing in $\polycpx$}
  \label{sec:geom-interl-polycpx}
There is a geometric condition for the interlacing of polynomials in
$\polycpx$. Although it is in general hard to verify it does help to
understand why polynomials in $\polycpx$ fail to interlace.

Suppose that $f,g\in\polycpx$ satisfy $deg(f)=deg(g)+1$ and have no
common factors. If they don't interlace then there is an
$\alpha\in\reals$ and $\sigma$ in the closed upper half plane such
that $(f + \alpha g)(\sigma)=0$. Now if $\alpha=0$ all the roots are
in the lower half plane, and the roots are continuous functions of
$\alpha$, so we see that
\begin{quote}
  If $f $ and $g$ don't interlace then there are $\alpha,t\in\reals$
  such that \\  $(f+\alpha g)(t)=0$. If they do interlace then $(f+\alpha
  g)(t)$ is never zero.
\end{quote}

Now if $f(t) + \alpha g(t)=0$ then $\frac{f}{g}(t)\in\reals$, so $\arg(f/g)$
is a multiple of $\pi$. This gives us a way to check interlacing.
Let
\begin{align*}
  f &= \prod(x-r_i) & g = \prod(x-s_i)
\end{align*}
\index{arg sum}
Thus, if $deg(f)=n$ then $f$ and $g$ interlace iff the \emph{arg sum}
\begin{equation}
  \label{eqn:arg-sum}
  \arg\frac{f}{g} = \sum_1^n \arg(x-r_i) - \sum_1^{n-1} \arg(x-s_i)
\end{equation}
is not a multiple of $\pi$ for any real $x$. Sometimes we can use
geometric arguments to show that this is indeed the case.

\begin{example}
  Suppose that the roots of $f$ and $g$ all have the same imaginary
  part, and their real parts alternate. (Figure~\ref{fig:arg-1})
  Now
\[
   \sum_1^n \arg(x-r_i) - \sum_1^{n-1} \arg(x-s_i) =
\arg(x-r_n) + \sum_1^{n-1} \bigl( \arg(x-r_i) - \arg(x-s_i)\bigr).
\]
Since the roots alternate all of these differences have the same sign.
The arg sum is the sum of the shaded angles in the figure.
It's clear that their sum is between $\pi$ and $0$ for all $x$, and so they interlace.

\begin{figure}[h]
  \centering
{\psset{yunit=1.5cm,xunit=1.5cm}
\begin{pspicture}(0,.5)(5,2.2)
  \psline(0,2)(5.2,2)
  \psline(0,1)(5.2,1)
  \rput(6,2){$\Im(z)=0$}
  \psline(1,1)(3.5,2)
  \psline(2,1)(3.5,2)
  \psline(3,1)(3.5,2)
  \psline(4,1)(3.5,2)
  \psline(5,1)(3.5,2)
  \rput(1,.6){$f_1$}
  \rput(2,.6){$g_1$}
  \rput(3,.6){$f_2$}
  \rput(4,.6){$g_2$}
  \rput(5,.6){$f_3$}
  \rput(3.5,2.2){$x$}
  \pswedge[fillcolor=red,fillstyle=solid](3.5,2){1.5}{-33}{0}
  \pswedge[fillcolor=red,fillstyle=solid](3.5,2){1}{-116}{-63}
  \pswedge[fillcolor=red,fillstyle=solid](3.5,2){1.5}{-158}{-146}
\end{pspicture}
  }
  \caption{Interlacing example}
  \label{fig:arg-1}
\end{figure}
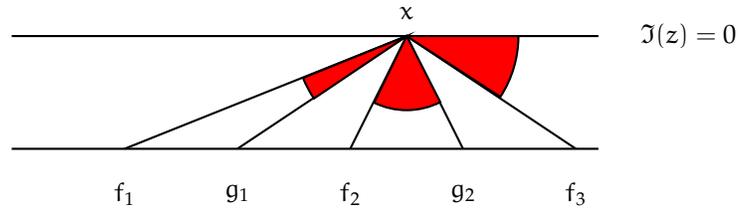
\end{example}

\begin{example}
  Now we look at the behavior at infinity. If $deg(f)=deg(g)+1$ and
  $x$ goes to positive infinity then the angle between any  of the
  $2n-1$ points and $x$ is close to $0$, so the arg sum converges to
  $0$. As $x$ goes to $-\infty$ the arg sum goes to $\pi$. See
  Figure~\ref{fig:arg-plot-3} for an example of a plot of the arg sum
  of two interlacing polynomials.

\begin{figure}[h]
  \centering
{
\psset{xunit=.1cm}
\begin{pspicture}(-35,0.105)(15,3.2)
\psline[linestyle=dashed,linecolor=blue](-35,3.14)(15,3.14)
\psline[linestyle=dashed,linecolor=blue](-35,0)(15,0)
\rput(17,3.14){$\pi$}
\rput(17,0){$0$}
 \pscurve(-35,3.07)(-34,3.07)(-33,3.07)(-32,3.06)
        (-31,3.06)(-30,3.06)(-29,3.05)(-28,3.05)
        (-27,3.04)(-26,3.03)(-25,3.03)(-24,3.02)
        (-23,3.01)(-22,2.99)(-21,2.97)(-20,2.94)
        (-19,2.89)(-18,2.79)(-17,2.54)(-16,1.88)
        (-15,1.38)(-14,1.64)(-13,2.13)(-12,2.35)
        (-11,2.4)(-10,2.34)(-9,2.2)(-8,1.98)
        (-7,1.69)(-6,1.47)(-5,1.46)(-4,1.58)
        (-3,1.6)(-2,1.4)(-1,1.08)(0,0.771)
        (1,0.556)(2,0.423)(3,0.34)(4,0.284)
        (5,0.244)(6,0.214)(7,0.191)(8,0.173)
        (9,0.158)(10,0.145)(11,0.135)(12,0.126)
        (13,0.118)(14,0.111)(15,0.105)
\end{pspicture}
}
  
  \caption{An argument plot for interlacing polynomials}
  \label{fig:arg-plot-3}
\end{figure}
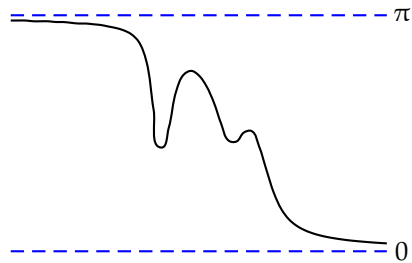

\end{example}

\begin{example}
  The polynomials whose roots are given in Figure~\ref{fig:arg-2}
  don't interlace. This can be seen from the plot of the arg sum in
  Figure~\ref{fig:arg-plot-4} since the curve crosses the line $y=0$.
  The failure to interlace can also be checked using
  Lemma~\ref{lem:i-interlace}.

  \begin{figure}[h]
    \centering
\begin{pspicture}(-3,-2)(4,1.2)
  \psline(-3,0)(3,0)
  \rput(4,0){$\Im(z)=0$}
  \psline(1,-1)(0,0)
  \psline(2,-1)(0,0)
  \psline(0,-1)(0,0)
  \psline(-1,-1)(0,0)
  \psline(-2,-1)(0,0)
  \rput(0,-1.4){$f_3$}
  \rput(-1,-1.4){$f_2$}
  \rput(-2,-1.4){$f_1$}
  \rput(1,-1.4){$g_1$}
  \rput(2,-1.4){$g_2$}
  \rput[l](2,-2){$2-\imag$}
  \rput(1,-2){$1-\imag$}
  \rput(0,-2){$-\imag$}
  \rput(-1,-2){$-1-\imag$}
  \rput[r](-2,-2){$-2-\imag$}
  \rput(0,.2){$0$}
\end{pspicture}
    
    \caption{Non-interlacing example}
    \label{fig:arg-2}
  \end{figure}
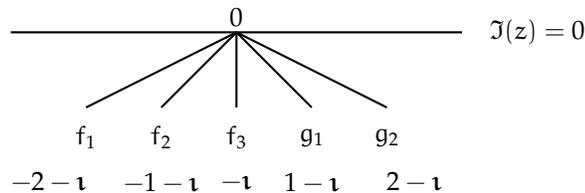

  \begin{figure}[h]
    \centering
{
\psset{xunit=.3cm,yunit=.5cm}
    \begin{pspicture}(-10,-2.36)(5,3)
\psline[linestyle=dashed,linecolor=blue](-10,0)(5,0) 
\psline[linestyle=dashed,linecolor=blue](-10,3.1)(5,3.1)
\rput(6,0){$0$} 
\rput(6,3.1){$\pi$} 
 \pscurve(-10,2.98)(-9,2.96)(-8,2.92)(-7,2.87)
        (-6,2.8)(-5,2.68)(-4,2.47)(-3,2.01)
        (-2,0.888)(-1,-0.785)(0,-2.21)(1,-2.36)
        (2,-1.33)(3,-0.485)(4,-0.178)(5,-0.0623)
\end{pspicture}
}
    \caption{An argument plot for non-interlacing polynomials}
    \label{fig:arg-plot-4}
  \end{figure}
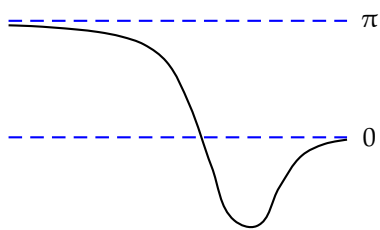

\end{example}

\begin{example}
  In this example we give four roots, for which no choice of a fifth
  root leads to interlacing. Consider Figure~\ref{fig:arg-plot-5}. The
  arg sum is 
\[
\angle f_10g_1 + \angle f_20g_2 + \angle 10f_2 = {\pi}/{2} + {\pi}/{2} + \angle 10f_3 >
\pi
\]
so no choice of $f_3$ can give interlacing polynomials.

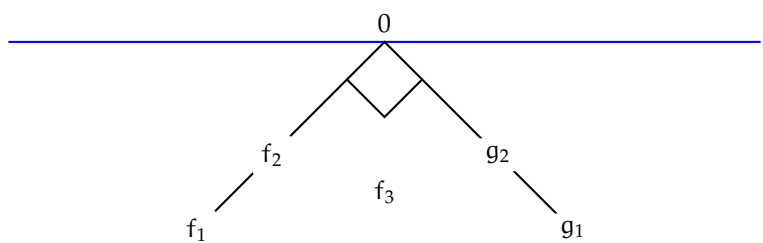
\begin{figure}[h]
  \centering
{
\psset{unit=.5cm}
  \begin{pspicture}(-10,-5)(10,.5)
    \psline[linecolor=blue](-10,0)(10,0)
    \psline(-5,-5)(0,0)(5,-5)
    \psline(-1,-1)(0,-2)(1,-1)
    \rput*(-5,-5){$f_1$}
    \rput*(5,-5){$g_1$}
    \rput*(-3,-3){$f_2$}
    \rput*(3,-3){$g_2$}
    \rput(0,-4){$f_3$}
    \rput(0,.5){$0$}
  \end{pspicture}
}
  \caption{Roots of non-interlacing polynomials}
  \label{fig:arg-plot-5}
\end{figure}
\end{example}



\begin{example}
  We revisit the problem of finding all $\beta$ for which
  $(x+\imag)^n$ and $(x+\beta\imag)^{n-1}$ interlace. The arg sum is
\[
A(x,\beta) = n \arctan(1/x) - (n-1) \arctan(\beta/x)
\]
The unique solution to $A(x,\beta)=0$ is $\beta=n/(n-1)$ as we saw
before. The solution to $A(x,\beta)=\pi$ determines the lower
bound. We find the smallest possible $\beta$ by solving
\begin{align*}
  A(x,\beta) &= 0 \\
  \frac{\partial}{\partial x} A(x,\beta) &= 0
\end{align*}
Using the addition formula for arc tangents
\[
\arctan x + \arctan y = \arctan\frac{x+y}{1-x y}
\]
we see that these are polynomial equations. The first few solutions
are in Table~\ref{tab:intervals-i}. 
\end{example}

\section{Orthogonal type recurrences}
\label{sec:orthog-i}

We can define orthogonal type recurrences whose members all interlace
in $\polycpx$; the only new condition is that the imaginary part of
the constant term is positive.

\begin{lemma}\label{lem:orthog-i}
    Choose constants $a_k>0$, $\Im(b_k)>0$, $c_k>0$.  If
    $p_{-1}=0$, $p_0=1$ and
\[
p_{n+1} = (a_n\,x+b_n)\,p_n - c_n\,p_{n-1}
\]
then $p_n\in\polycpx$ and we have interlacings $p_1\lessgreat p_2
\lessgreat p_3 \lessgreat \cdots$.
  \end{lemma}
  \begin{proof}
    We prove that $p_{n+1}\lessless p_n$ by induction on $n$. We
    first note that $p_1\in\polycpx$ since its only root has negative
    imaginary part by the hypotheses on $a_0$ and $b_0$. Thus,
    $p_1\lessless p_0$ since all linear combinations are in
    $\polycpx$.

    Next, assume that $p_n\lessless p_{n-1}$. All the leading
    coefficients of the polynomials in the interlacing
\[
(a_nx+b_n)\,p_n\lessless p_n \lessless c_n p_{n-1}
\]
are positive, and so it follows from Lemma~\ref{lem:i-analog} that
\[
(a_nx+b_n)\,p_n - c_n p_{n-1} \lessless  p_n
\]
which finishes the proof. 

  \end{proof}

\begin{example}
We get two interrelated recurrences if we separate into real and
imaginary parts. Consider the simple example

\[
p_{n+1} = (x+\imag)\,p_n - p_{n-1}
\]

If we write $p_n = r_n + \imag s_n$ then

\begin{gather*}
r_{n+1} = x\,r_n - r_{n-1}-s_{n} \\
s_{n+1} = x\,s_n - s_{n-1}+ r_n
\end{gather*}

The  recurrences for just the $r$'s and $s$'s are the same:

\begin{gather*}
  r_{n+1} = 2x\,r_n - (3+x^2)\,r_{n-1} + 2x\, r_{n-2} - r_{n-3}\\
  s_{n+1} = 2x\,s_n - (3+x^2)\,s_{n-1} + 2x\, s_{n-2} - s_{n-3}\\
\end{gather*}

It is hard to see how one would ever show that the $r$'s and $s$'s
interlace if we only had these recurrences.
\end{example}

  \begin{example}
    Consider the Chebyshev recurrence $p_{n+1} = 2x\,p_{n}-p_{n-1}$,
    except we start with $p_1=x+\imag$. An easy induction shows that
    \[ p_n = T_n + \imag U_{n-1} \] We can see directly that this is
    in $\polycpx$ since the Chebyshev polynomials $T_n$ and $U_n$ have
    positive leading coefficients and satisfy $T_n\lessless U_{n-1}$.
  \end{example}

  Here's a simple generalization of that fact.

  \begin{lemma}
    Suppose orthogonal polynomial sequences $\{p_i\}$ and $\{q_i\}$
    both satisfy the recurrence
\[ r_{n+1}(x) = (a_nx+b_n)r_n(x) - c_nr_{n-1}(x)\]
where all $a_n$ and $c_n$ are positive. If \[ 
p_0=1,\ p_1=x,\ q_0=0,\ q_1=1 \]
then
\[
p_n\lessless q_n \qquad\text{and}\qquad
\begin{vmatrix}
  p_n& p_{n+1}\\ q_n & q_{n+1}
\end{vmatrix}>0.
\]
  \end{lemma}
  \begin{proof}
    Let $w_i$ satisfy the same recurrence with $w_0=1$ and
    $w_1=x+\imag$. By induction $w_n = p_n + \imag q_n$, and $w_n$
    satisfies the conditions of Lemma~\ref{lem:orthog-i}. Thus,
    $w_n\lessless w_{n+1}$ which gives the conclusions of the lemma.
  \end{proof}

  If we have a sequence of interlacing polynomials $f_k\lesslesseq
  f_{k-1}$ and we let $p_n = f_n + \imag f_{n-1}$ then all $p_n$ are
  in $\polycpxclose$. Table~\ref{tab:rec-rel} lists recurrence
  relations for the $p_n$ for several polynomial sequences.

  \begin{table}\label{tab:rec-rel}
    \centering
\[
\begin{array}{rrl}
  \toprule
  \multicolumn{1}{l}{\text{Rising Factorial}} & p_n = & \rising{x}{n} + \imag\, \rising{x}{n-1}\\
 & p_n  = & \bigl[ n-1-\imag + x\bigr] p_{n-1} + \imag\bigl[
n-3+x\bigr] p_{n-2} \\
\midrule
  \multicolumn{1}{l}{\text{Hermite}} & p_n = & H_n + \imag H_{n-1}\\
& p_n = & 2x p_{n-1} - p_{n-1}' \\
\midrule  
  \multicolumn{1}{l}{\text{Bell}} & p_n =& B_n + \imag B_{n-1}\\
& p_n = &   x( p_{n-1} + p_{n-1}')\\
\midrule
  \multicolumn{1}{l}{\text{Euler}} & p_n = & B_n + \imag B_{n-1}\\
& p_n =&  (nx-\imag) p_{n-1} + (n-2)\,\imag\,x\, p_{n-2} + (x-x^2)
(p_{n-1}' + \imag p_{n-2}')\\
\midrule  
  \multicolumn{1}{l}{\text{Laguerre}} & p_n = & L_n(-x) + \imag L_{n-1}(-x)\\
&n\, p_n =&  \bigl[ x + n - \imag\,(n-2)\bigr] p_{n-1} + \imag\,(n-2)p_{n-2} +
(1+\imag) p_{n-1}'\\
\midrule  
  \multicolumn{1}{l}{\text{Legendre}} & p_n = & P_n + \imag P_{n-1}\\
& p_n = & (x+\imag)p_{n-1} - (n-2)p_{n-2} + 2x^2 p_{n-1}' - 2x p_{n-2}'\\
\midrule  
  \multicolumn{1}{l}{\text{Gegenbauer}} & p_n = & G^{(2)}_n + \imag G^{(2)}_{n-1}\\
& n(n-1)\,p_n = & 2(n+2)(n-1)\,x\,p_{n-1} - (n+1)(n+2)p_{n-2} +
2(1-x^2)p_{n-1}'\\
\midrule  
\end{array}
\]

\caption{Recurrences for $f_n + \imag f_{n-1}$}
  \end{table}

\section{The transformation $x\mapsto \exp(-\imag x)$}

We have seen (Lemma~\ref{lem:f-of-sin}) that
$\sin(x)-\alpha\in\allpolyf$ when $-1\le \alpha\le1$. The next lemma
extends this to complex exponentials. We then use this result to
determine mapping properties of the Chebyshev polynomials.

\begin{lemma}\label{lem:cpxpoly-2}
  If $|\sigma|\le1$  then $e^{-\imag x}+\sigma\in\polycpxf$.
\end{lemma}
\begin{proof}
If we write $x = a + b\imag$ then 
\[
e^{-\imag x}+\sigma = e^{-a\imag}\,e^{b} + \sigma
\]
Now $|e^{-a\imag} |=1$, and $e^b>1$ since $b>0$, so the sum is
non-zero since $|\sigma|\le1$.

\end{proof}

\begin{prop}\label{prop:polycpx-3}
  If $f\in\allpolyint{\Delta}$ then $f(e^{-\imag x})\in\polycpxf$.
\end{prop}
\begin{proof}
  This is an immediate consequence of  Lemmas~\ref{lem:cpxpoly-1} and
  \ref{lem:cpxpoly-2}. This is essentially due to \Polya: see \cite{iliev}*{page 50}.
\end{proof}

\begin{cor}\label{lem:cheby}
  If $T_n$ and $U_n$ are the Chebyshev polynomials then $x^n\mapsto
  T_n(x)$ and $x^n\mapsto U_n(x)$ map
  $\allpolyint{(-1,1)}\longrightarrow\allpolyint{(-1,1)}$.
\end{cor}
\begin{proof}
  We only consider $T_n$; the case of $U_n$ is similar. $T_n(x)$ is defined by
  $T_n(\cos x)= \cos(nx)$, so the diagram below commutes at
  the element level. Proposition~\ref{prop:polycpx-3} implies that the
  top arrow maps as claimed, and hence $x^n\mapsto \cos nx$ maps as in
  the diagram. Since $g(\cos x)\in\allpolyf$ if and only if
  $g(x)\in\allpolyint{(-1,1)}$, the proof is complete.

  \centerline{\xymatrix{
\allpolyint{(-1,1)}
\ar@{-->}[drrr]^{{\ \ x^n\mapsto\cos nx}}
      \ar@{.>}[d]_{{x^n\mapsto T_n(x)    }}           
      \ar@{->}[rrr]^{{x\mapsto e^{-\imag nx}   }}         
      &&&
      {{\polycpxf}}
      \ar@{->}[d]^{{ f\mapsto \Re e(f)    }} \\        
\allpolyint{(-1,1)}
      \ar@{<-}[rrr]_{{\cos(x)\mapsto x    }}         
      &&&
      {\allpolyf}
}}

\end{proof}

\begin{example}
 $\dfrac{1-e^{- \imag y}}{y}$ is another function   in
  $\polycpxf$. To see that this is so, note that
       \begin{align*}
         \frac{1-e^{-y\imag}}{y} &= \frac{1-\cos y}{y} + \imag
         \frac{\sin y}{y} \\
         & =  \frac{2\,\sin^2(y/2)}{y} + \imag \frac{\sin y}{y}
       \end{align*}
       The infinite product expansion shows that $\sin^2(y/2)/y$ and
       $\sin(y)/y$ are in $\allpolyf$. Moreover, if we approximate them
       by their partial products then they interlace ($\lesslesseq$).
       For instance, for zeros up to $\pm2\pi$ the roots are

\begin{align*}
  \sin^2(y/2)/y &:\quad\quad -2\pi,-2\pi,0,0,2\pi,2\pi \\
\sin(y)/y &:\quad\quad -2\pi,-\pi,0,\pi,2\pi
\end{align*}
\end{example}

If we know that a polynomial is in $\allpolyint{\Delta}$ then 
we can use Proposition~\ref{prop:polycpx-3} to get information about
trigonometric and exponential polynomials.  We use a result due to
Enestr\"{o}m-Kakeya \cite{borwein} to create a polynomial in
$\allpolyint{\Delta}$. 

\index{Kakeya}
\begin{theorem}\label{thm:kakeya}
  If $f(x) = a_0 + a_1 x+ \cdots + a_n x^n $ and 
\begin{equation}\label{eqn:kakeya} 
a_n > a_{n-1} >  \cdots    > a_0 > 0
\end{equation}
then  $f\in\allpolyint{\Delta}$. 
\end{theorem}
 
\begin{cor}
If the real numbers $\{a_k\}$ satisfy \eqref{eqn:kakeya} then 
\begin{enumerate}
\item $\sum_{k=0}^n a_k \sin(kx)$ has all real roots.
\item $\sum_{k=0}^n a_k \cos(kx)$ has all real roots.
\item The roots of $\sum_{k=0}^n a_k \sin(kx)$ and  $\sum_{k=0}^n a_k
  \cos(kx)$ interlace.
\end{enumerate}
\end{cor}
\begin{proof}
  Since $f(x) = a_0 + a_1 x+ \cdots + a_n x^n $ is in
  $\allpolyint{\Delta}$, we can apply Proposition~\ref{prop:polycpx-3}
  to conclude that $ \sum a_k e^{-\imag k x} \in\polycpxf$. This
  implies (1), (2), and (3).
\end{proof}

If we take limits we get \cite{polya-szego1}*{III 205}

\begin{cor}\label{cor:phi-int}
  If $\phi(x)$ is a positive increasing function on $(0,1)$ with
  $\int_0^1\phi(x)\,dx$ finite then
$$ \int_0^1 \phi(x) e^{-\imag t x}\,dx \in\polycpxf.$$
\end{cor}
\begin{proof}
  Since the integral of $\phi$ is finite we have 
$$ \int_0^1 \phi(x) e^{-\imag t x}\,dx = \lim_{n\rightarrow\infty}
\sum_{k=1}^n \frac{1}{n} \phi(k/n) e^{-\imag t k /n}
$$
Since $\phi$ is increasing the coefficients $(1/n)\phi(k/n)$ are positive and
increasing. It follows from Proposition~\ref{prop:polycpx-3} that the
sum is in $\polycpxf$, and hence the limit is as well.
\end{proof}

\section{Newton inequalities}
\label{sec:newton-inequalities}

We have seen  mutually interlacing polynomials  evaluated at points in the
upper half plane (Figure~\ref{fig:mutual}). We now give an example of the
ratio of the coefficients of a polynomial in $\gsubplus_2$ evaluated at a point in
$\quada$. If $f=\sum f_i(x)y^i$ then Figure~\ref{fig:quada} shows the
ratios of a fourth degree polynomial evaluated at a point $\sigma\in\quada$.

\begin{figure}
  \centering
{
\newpsobject{showgrid}{psgrid}{subgriddiv=1,griddots=10,gridlabels=6pt}
\psset{xunit=1cm,yunit=.6cm}
\begin{pspicture}(0,0)(7,7)
\psline(6,0)(0,0)(0,6)
\psdots(1,1)(2,3)(4,4)(5,6)
\rput[l](1,1){$\quad \frac{f_0(\sigma)}{f_1(\sigma)}$}
\rput[l](2,3){$\quad \frac{f_1(\sigma)}{f_2(\sigma)}$}
\rput[l](4,4){$\quad \frac{f_2(\sigma)}{f_3(\sigma)}$}
\rput[l](5,6){$\quad \frac{f_3(\sigma)}{f_4(\sigma)}$}
\end{pspicture}
}
  \caption{ Coefficient ratios of $f\in\gsubplus_2$ and $\sigma\in\quada$}
  \label{fig:quada}
\end{figure}
 There are two striking facts in this picture:
 \begin{enumerate}
 \item The ratios are all in the first quadrant.
 \item The real and imaginary parts are increasing.
 \end{enumerate}

 The first part is easy; the second part is contained in the following
 lemma. We can differentiate, so we only need to establish it for
 $f_0$ and $f_1$. Since $f_0\lesslesseq f_1$, we can write
\[
\frac{f_0}{f_1} = \biggl(\sum \frac{a_i}{x+r_i}\biggr)^{-1}
\]
where the $a_i$ and the $r_i$ are positive. If $\sigma\in\quada$ then
$\sigma+r_i\in\quada$, and so $\sum (\sigma+r_i\bigr)^{-1}\in\quadd$,
which establishes the first observation above.

The  lemma below is another generalization of Newton's
inequality. If we let $f = \sum a_ix^i\in\allpolypos$ then
$f(x+y)\in\gsubplus_2$. Choosing $\sigma=0$ in the Proposition gives
that
\[
0 \le \Re \frac{a_i}{a_{i+1}} \le \Re \frac{a_{i+1}}{a_{i+2}} 
\]
since $f^{(i)}(0)/i!=a_i$. The imaginary part is vacuous.

\begin{lemma}\label{lem:newton-p2-3}
  Suppose that $\sum {f_i(x)}\,y^i\in\gsubplus_2(n)$. For all
  $\sigma\in\quada$ and $0\le i\le n-2$ 

\[
\begin{array}{ccccc}
 0 & \le & \Re\ \frac{f_i(\sigma)}{f_{i+1}(\sigma)} &  \le & \Re\
 \frac{f_{i+1}(\sigma)}{f_{i+2}(\sigma)} \\[.4cm]
 0 & \le & \Im\ \frac{f_i(\sigma)}{f_{i+1}(\sigma)}   &\le & \Im\
 \frac{f_{i+1}(\sigma)}{f_{i+2}(\sigma)} \\
  \end{array}
\]
\end{lemma}
\begin{proof}
Since $\sigma\in\quada$ we know that $f(\sigma,y)\in\polycpx$. The
inequality now follows from Lemma~\ref{lem:newton-i-2}. 
\end{proof}

  If all the roots of a polynomial lie in a quadrant then there are
  Newton inequalities for the real parts, imaginary parts, and magnitudes of the
  coefficients.

  \begin{lemma}\label{lem:newton-i-2}
    Suppose that $\sum_0^n a_i x^i$ has all roots in a quadrant. The table
    below summarizes the properties of the sequence $
    \frac{a_0}{a_1},\,\frac{a_1}{a_2},\,\frac{a_2}{a_3},\dots$ where
    ``$neg,\downarrow$'' means it is a decreasing sequence of negative
    terms, and so on.
    \\[.2cm]

\centering{
\begin{tabular}{crrrr}
  \toprule
  & $\quada$ & $\quadb$ & $\quadc$ & $\quadd$ \\\midrule
  Real parts & $neg,\downarrow$ & $pos,\uparrow$ & $pos,\uparrow$ & $neg,\downarrow$ \\[.2cm]
Imaginary parts &
$neg,\downarrow$ & $neg,\downarrow$ & $pos,\uparrow$ & $pos,\uparrow$ \\[.2cm]
Magnitude & $\uparrow$& $\uparrow$& $\uparrow$& $\uparrow$ \\\bottomrule
\end{tabular}
}
  \end{lemma}
  \begin{proof}
    As usual, we  differentiate, reverse and differentiate; the result
    is a quadratic with all roots in the original quadrant. We will
    work with the third quadrant since we can use properties of
    $\polycpx$. The upper two quadrants follow using conjugation, and
    the fourth quadrant is similar.

    Thus, we may assume that our polynomial is 
    \[ (x+a)(x+ b) + \imag c(x+d) \qquad 0< a < d < b, \quad 0<c \] We
    simply need to compute the real parts, imaginary parts and
    magnitudes, and show that the inequalities are satisfied.  If we
    write the polynomial as $a_0+a_1x+a_2x^2$ then
\begin{align*}
  a_0 &= ab+\imag cd &
  a_1 &= a+b+\imag c \\
  a_2 &=1 \\
  \Re\frac{a_0}{a_1} &= \frac{b a^2+b^2 a+c^2 d}{a^2+2 b a+b^2+c^2}&
  \Re\frac{a_1}{a_2} &= a+b\\
  \Im\frac{a_0}{a_1} &= \frac{c (-a b+d b+a d)}{a^2+2 b a+b^2+c^2}&
  \Im\frac{a_1}{a_2} &= c \\
  \biggl|\frac{a_0}{a_1}\biggr|^2  &=  \frac{a^2 b^2+c^2 d^2}{(a+b)^2+c^2} &
  \biggl|\frac{a_1}{a_2}\biggr|^2  &=  (a+b)^2+c^2
\end{align*}
All six terms are clearly positive, and it is easy to check that all
three inequalities are satisfied.

  \end{proof}

The bound for absolute values is  stronger that the bound we
get from considering sectors - see~\pageref{eqn:newton-cpx-2}.

\chapter{Stable polynomials}
  \label{cha:stable}

\renewcommand{\TimeStampStart}{Monday, March 10, 2008: 10:27:45}
\mytoday

  A one variable polynomial is \emph{stable} (sometimes called \emph{Hurwitz
    stable}) if all its roots lie in the closed left
  half plane. This is a well studied class of polynomials with many
  applications - e.g. \cite{MR1787335}.  There are two classes of
  stable polynomials in one variable:
\begin{align*}
\stabled{1} &= \text{ real coefficients and roots in
  the closed left half plane}\\
\stabledc{1} &= \text{ complex coefficients and roots in
  the closed left half plane}\\
\end{align*}

Since  polynomials whose roots are in the closed left half plane are
non-zero in the right half plane, we can apply the results of
Chapter~\ref{cha:nv} to get properties of stable polynomials in $d$
variables. 

\begin{definition}\ 

  $\stabledc{d} = \left\{\text{\parbox{3.5in}{All polynomials
      $f(x_1,\dots,x_d)$ with complex coefficients such that
      $f(\sigma_1,\dots,\sigma_d)\ne0$ for all
      $\sigma_1,\dots,\sigma_d$ in the right half plane. If we don't
      need to specify $d$ we simply write $\stabledc{}$.}}\right.$ We call
such polynomials  \emph{stable polynomials}. In one variable they are
often called Hurwitz stable.

\end{definition}

After establishing the basic properties of stable polynomials in $d$
variables, we focus our attention on stable polynomials with real
coefficients in one variable. We will see that whenever a polynomial
is shown to be strictly positive, it's likely that it's stable. There
are three kinds of interlacing for stable polynomials, and
\emph{positive interlacing}  is especially useful for
showing that a polynomial is stable.  \index{positive interlacing}
\index{interlacing!positive}

For a short summary of the properties of stable polynomials, see
\cite{aspects}.

\section{Stable polynomials in $d$ variables}
\label{sec:stable-polynomials-d}

We apply the general results of Chapter~\ref{cha:nv} to derive results
about polynomials non-vanishing in the open right half plane ($\rhp$).    The
reversal constant is $1$, and
$\rhp/\rhp=\complexes\setminus(-\infty,0)$. If $f,g\in\stabledc{d}$ then
we write $f\hlace g$ if $f + z g\in\stabledc{d+1}$. From
Chapter~\ref{cha:nv} we have

\begin{lemma}\label{lem:s-basic} Suppose that $f\in\stabledc{d}$.

  \begin{enumerate}
  \item If $\sigma_i\in\rhp$ and $r_i>0$ then
    $f(r_1x_1+\sigma_1,\dots,r_dx_d+\sigma_d)\in\stabledc{d}$.
  \item $f(x_1+y,x_2,\dots,x_d)\in\stabledc{d+1}$.
  \item If $\sigma\in\reals$ then
    $f(\imag\sigma,x_2,\dots,x_d)\in\stabledc{d-1}\cup\{0\}$. 
  \item If $f\hlace g$ then $g\hlace f$.
  \item If $f\hlace g\hlace h$ then $f+h\hlace g$.

  \item $\stabledc{d}$ is closed under differentiation. That is, if
    $g\in\stabledc{d}$ then $\frac{\partial}{\partial
      x_i}f\in\stabledc{d}\cup\{0\}$. 
  \item $f\hlace \partial f/\partial x_i$
  \item If $\sum x^if_i(\yy)\in\stabledc{d+1}$ then $f_i(y)\in\stabledc{d}\cup\{0\}$ and
    $f_i\hlace f_{i+1}$.

  \item $f(\xx)\times g(\xx)\mapsto f(\partial \xx)g(\xx)$ determines
    a map  \[\stabledc{d}\times \stabledc{d}\longrightarrow \stabledc{d}.\]
    \item If $\sum x^if_i(\yy)\in\stabledc{d}$ has degree $n$ then $\sum
    x^{n-i}\,f_i(y)\in\stabledc{d}$.
\item If $x_i$ has degree $e_i$ then $x_1^{e_1}\cdots x_d^{e_d}
  f(1/x_1,\dots,1/x_d)\in\stabledc{d}$.

\item The Hadamard product $\sum a_ix^i\times\sum f_i(\yy)x^i\mapsto \sum
  a_if_{n-i}(\yy)x^i$ determines a map $\allpolypos\times
  \stabledc{d}\longrightarrow \stabledc{d}$.
  \item If $S$ is skew-symmetric and all $D_i$ are positive definite then \\
$ \bigl| S + \sum_1^d x_i D_i\bigr| \in\stabled{d} $.

\item The following are equivalent
  \begin{enumerate}
  \item $f(\xx)\in\stabled{d}$
  \item $f(\aaa+t\bbb)\in\stabled{1}$ for $\aaa,\bbb>0$.
  \end{enumerate}
  \end{enumerate}

\end{lemma}
\begin{proof}
 Part (1) shows a) implies b).
 Conversely, consider $\sigma_1,\dots,\sigma_d$ in $\rhp$.
  Since $1$ is in the closure of $\rhp$ we can find an $\alpha$
  and positive $\aaa,\bbb$ such $(\sigma_1,\dots,\sigma_d) = \aaa +
  \alpha\bbb$. Thus,
  $f(\sigma_1,\dots,\sigma_d)=f(\aaa+\bbb\alpha)\ne0$.
\end{proof}

\section{Real stable polynomials}
\label{sec:real-stable-polyn-1}

We now consider properties that depend on real coefficients.
It is immediate that $\allpolypos\subset\stabled{1}$; this is still true
for more variables.  $\gsubplus_d$ was defined in terms of $\rup{d}$
and the positivity of the coefficients. We can also 
express this in terms of non-vanishing.

  \begin{lemma}\label{lem:pdpos-up}
    $f(\xx)\in\gsubplus_d$ if and only if
    $f(\frac{\sigma_1}{\sigma_0},\dots,\frac{\sigma_d}{\sigma_0})\ne0$
    for all $\sigma_i$ in the upper half plane.
  \end{lemma}
  \begin{proof}
    We know that $f(\xx)\in\gsubplus_d$ if and only if its
    homogenization is in $\rup{d+1}$. Now the homogenization equals
    $x^nf(\frac{x_1}{x},\cdots,\frac{x_d}{x})$ where $n$ is the degree
    of $f$, so the lemma is now clear.
  \end{proof}
  \begin{lemma}
    $\gsubplusclose_d = \stabled{d}\cap\gsubclose_d$
  \end{lemma}
  \begin{proof}
    We first show that $\gsubplus_d \subset \stabled{d}\cap\rup{d}$.
    Suppose that $f(\xx)\in\gsubplus_d$, and $\tau_1,\dots,\tau_d$ are
    in the right half plane. Now $\imag \tau_1,\dots,\imag\tau_d$ lie in
    the upper half plane, and so by Lemma~\ref{lem:pdpos-up}
\[ 0 \ne f\bigl(\frac{\imag\tau_1}{\imag},\dots,\frac{\imag\tau_d}{\imag}\bigr)
= f(\tau_1,\dots,\tau_d)
\]
Thus, $f(\xx)\in\stabled{d}$.  Next, if $f\in \stabled{d}\cap\gsubclose_d$
then $f\in\gsubclose_d$, and $f$ has all non-negative coefficients, so
$f\in\gsubplusclose_d$ by Lemma~\ref{lem:p-limit}
  \end{proof}

\begin{cor} 
 The unsigned reversal of a polynomial in $\gsubplus_d$ is in $\stabled{d}$.
\end{cor}
\begin{proof}
  The only point is that $\gsubplus_d\subset\stabled{d}$.
\end{proof}

\begin{lemma}\label{lem:h-is-pos}
  If $f\in\stabled{d}$ then all coefficients have the same sign.
\end{lemma}
\begin{proof}
  By Lemma~\ref{lem:pd-limit} and following the proof of
  Lemma~\ref{lem:p-limit} we can write $f$ as a limit of
  $f_\epsilon\in\stabled{d}$ where all coefficients of $f_\epsilon$
  are non-zero. Since the coefficients of a polynomial in $\stabled{1}$
  with all non-zero coefficients have the same sign, an easy induction
  shows that all the coefficients of $f_\epsilon$ have the same
  sign. Thus, the non-zero coefficients have the same sign.
\end{proof}



\begin{lemma}
  If $f(\xx)\in\stabled{d}$ and the homogenization $F$ of $f$ is in
  $\stabled{d+1}$ then $ f\in\gsubplusclose_d$.
\end{lemma}
\begin{proof}
  If $F$ is in $\stabled{d+1}$ then for $\tau_0,\dots,\tau_d\in\rhp$ we know
  that
\[
f\bigl(\frac{\tau_1}{\tau_0},\dots,\frac{\tau_d}{\tau_0}\bigr)\ne0
\]
If $\sigma_1,\dots,\sigma_d\in\uhp$ then we can choose $\tau_0\in\rhp$
so that $\tau_0^{-1}$ rotates all of the $\sigma_i$ to the right half
plane. Thus $\sigma_i = \tau_i/\tau_0$ where $\tau_i\in\rhp$, and so
$f$ is non-vanishing on the upper half plane. Thus $f\in\gsubclose_d$,
and by Lemma~\ref{lem:h-is-pos} all coefficients of $f$ have the same
sign. Consequently, $f$ is in $\gsubplusclose_d$.
\end{proof}

\begin{remark}
    Here is a particular example of a polynomial in $\stabled{2}$ whose homogenization
    is not in $\stabled{3}$. Consider 
    $f(x,y)=1+(x+y)+(x+y)^2$ with homogenization $F(x,y,z)$. $F$ has a
    root in the right half plane: $F(1,-.23-1.86\imag,2\imag)=0$.
  \end{remark}


Multiplying each variable by a new variable is similar to
homogenizing.

\index{homogenizing!in $d$ variables!multiply by $y$}
\begin{lemma}
  $f(y\,\xx)\in\stabled{d+1}$ if and only if $f(\xx)\in\polypos{d}$.
\end{lemma}
\begin{proof}
  If $f(y\xx)\in\stabled{d+1}$ then substituting $1$ shows
  $f(\xx)\in\stabled{d}$. If $\sigma_1,\dots,\sigma_d\in\uhp$ then
  choose $\tau\in\rhp$ so that $\tau\sigma_1,\dots,\tau\sigma_d\in\rhp$.
  Since $\tau^{-1}\in\rhp$
\[ f(\sigma_1,\dots,\sigma_d) = f(\tau^{-1}(\tau\sigma_1,\dots,\tau\sigma_d))
\ne0
\]

Thus $f\in\rup{d}$, so $f\in\polypos{d}$.  Conversely, assume that
$f\in\polypos{d}$. We know that the homogenization of $f$ is in
$\polypos{d+1}$. Writing $f = \sum \aaa_\sdiffi \xx^\sdiffi$ and
taking the reverse yields
\[
\sum \aaa_\sdiffi \xx^\sdiffi \, y^{n-|\sdiffi|} \implies 
\sum \aaa_\sdiffi \xx^\sdiffi \, y^{|\sdiffi|} = f(y\cdot \xx) \in\stabled{d+1}
\]
so the proof is complete.
\end{proof}

  \begin{lemma}\label{lem:h-square}
    If $f(\xx,y^2)\in\stabled{d}$ then $f(\xx,y)\in\gsubplus{d}$.
  \end{lemma}
  \begin{proof}
    See the proof of the corresponding result for $\rup{d}$, Lemma~\ref{lem:p-square}.
  \end{proof}

  \begin{lemma}
    Suppose that $f(x) = \sum_0^n a_ix^i\in\stabled{1}$. If some internal
    coefficient is zero then $f(x) = x^r g(x^2)$ where $g\in\allpolypos$.
  \end{lemma}
  \begin{proof}
 We can write
\[
f(x) = x^r \times \prod_1^s(x+a_i) \times \prod_1^t((x+b_i)^2+c_i)
\]
where $a_i$ is positive, $b_i,c_i\ge0$ and $b_ic_i\ne0$. If the second
factor is non-empty then there are no internal zeros. If any $b_i$ is
non-zero there are no internal zeros, so $f(x) = x^r \prod(x^2+c_i)$
from which the conclusion follows.
  \end{proof}

\index{Hermite-Biehler}
  \begin{lemma}[Hermite-Biehler]
      Suppose that $f(\xx)$ is a polynomial, and
  write $f(\xx) = f_E(\xx)+f_O(\xx)$ where $f_E(\xx)$
  (resp. $f_O(\xx)$) consists of all terms of even (resp. odd)
  degree. Then
  \begin{quote}
    $f\in\stabledc{}$ if and only if $f_E\hlace f_O$.
  \end{quote}
  \end{lemma}
  \begin{proof}
    Since $f(\xx) = f_E(\xx)+f_O(\xx)\in\stabled{}$ we know $f(-\imag
    \xx) = f_E(-\imag \xx) + f_O(-\imag x)\in\up{}$. Now $f_E(-\imag
    \xx)$ has all real coefficients, and $f_O(-\imag \xx)$ has all
    imaginary coefficients, so by Hermite-Biehler for $\up{}$ we know
    that 
\[
f_E(-\imag\xx) + (y/\imag)f_O(-\imag \xx)\in\up{}
\]
Returning to $\stabledc{}$ yields $f_E(\xx) + yf_O(\xx)\in\stabledc{}$.
  \end{proof}

\section{Polynomials (not) in $\stabled{d}$}
\label{sec:polyn-not-stabl}

We give some examples of polynomials  that are (not)
in $\stabled{d}$.

  \begin{example}
    Since the map $z\mapsto z^2$ maps $\rhp\longrightarrow
    \complexes\setminus(-\infty,0)$ it follows that the equation
    $x^2+y^2=0$ has solutions in the right half plane, so $x^2+y^2\not\in\stabled{2}$.
  \end{example}

  \begin{example}
    We claim $xy+a^2\in\stabled{2}$ if $a\in\reals$. As above, the image
    of $\rhp$ under $xy$ is $\complexes\setminus(-\infty,0)$, so adding
    a positive amount will not make $xy$ equal to zero. 
  \end{example}

  \begin{example}
    If $\sigma$ is complex (and not negative real) then
    $x^2-\sigma\not\in\stabled{1}$. Indeed, one of $\sigma^{1/2}$,
    $-\sigma^{1/2}$ lies in the right half plane. Consequently,
    $xy-\sigma\not\in\stabled{2}$. 
  \end{example}

  These last two examples imply that 

  \begin{lemma}If $f$ is a polynomial with positive coefficients then
    $f(xy)\in\stabled{2}$ if and only if $f\in\allpolypos$.
  \end{lemma}

  The stable polynomials of degree $1$ in $\stabled{2}$ are clearly
  $ax+by+c$ where $a,b$ are positive and $c$ is non-negative. If
  $a,b,c$ are positive then we know that
\[
ax^2+bx+c \in\begin{cases}
\allpolypos & \text{if}\ b^2-4ac\ge0 \\
\stabled{1} & \text{if}\ b^2-4ac\le0 
\end{cases}
\]

\begin{lemma}
  If $a,b,c,d$ are positive then
$axy+bx + cy + d \in \stabled{2}$.
\end{lemma}
\begin{proof}
If $a+bx+cy+dxy=0$ then $y = -\frac{a+bx}{c+dx}$. By the following
lemma the \Mobius\ transformation $z\mapsto \frac{a+bz}{c+dz}$ maps
the right hand plane to itself. Thus, if $x\in\rhp$ then
$y\not\in\rhp$, and so $a+bx+cy+dxy\in\stabled{2}$.

This can also be derived from the corresponding criterion for $\up{2}$.
\end{proof}

\begin{remark}
  It is not the case that a multiaffine polynomial with all positive
  coefficients is in $\stabled{3}$. If 
\[f(x,y,z) = (135/4) + x+y+z+xy+xz+yz+24xyz\]
then $f(1/2+\imag,1/2+\imag,1/2+\imag)=0$, so $f\not\in\stabled{3}$.
\end{remark}

\begin{lemma}
  If the  \Mobius\ transformation $M$ satisfies $Mz= \frac{a+bz}{c+dz}$
  then the following are equivalent:
  \begin{enumerate}
  \item $M$ maps the  right half plane to itself.
  \item At least one of the following holds
    \begin{enumerate}
    \item $a,b,c,d$ have the same  sign (or some are zero).
    \item $b,d$ have the same sign and $ad=bc$. In this case the image
      is the single point $b/d$.
    \end{enumerate}
  \end{enumerate}
\end{lemma}
\begin{proof}
  If $r+\imag s\in\rhp$ then we need to show that $M(r+\imag
  s)\in\rhp$. Equivalently, we show that if $r>0$ then $\Re M(r+\imag
  s)>0$. Now
\[
M(r+\imag s) = 
\frac{b d s^2+(a+b r) (c+d r)}{(c+d r)^2+d^2 s^2}
+ \imag 
\frac{(b c-a d) s}{(c+d r)^2+d^2 s^2}
\]
and we can ignore the positive denominator, so we need to show that
\begin{equation}
  \label{eqn:mobius-h}
b d s^2+(a+b r) (c+d r) > 0
\end{equation}
but this is clear if $a,b,c,d,r$ are positive. If $ad=bc$ then
$Mz=b/d$ which is in the right half plane since $b/d>0$. 

Conversely, assume \eqref{eqn:mobius-h} holds. Assume that $abcd\ne0$;
the remaining cases are easily taken care of. Taking $s$ large shows
that $bd$ is positive, so $b$ and $d$ have the same sign. Taking $s=0$
and $r$ close to zero shows that $ac$ is positive, so $a$ and $c$ have
the same sign. Take $s=0$. If $ab<0$ then the factor $(a+br)$ has a
positive root, and so $(a+b r) (c+d r)$ will be negative for some $r$
close to that root, unless the root is also a root of the other
factor. In this case $a+bz$ is a multiple of $c+dz$, and so $ad=bc$.
\end{proof}

There are many quadratic forms in $\stabled{d}$. The next lemma
generalizes the fact that $x^2-1\in\allpoly$, and $x^2+1\in\stabled{1}$.

    \begin{lemma}\label{lem:stable-qf}
      If $Q$ is a $d$ by $d$ \nsd\ matrix then
\[
\xx^tQ\xx + \alpha \in
\begin{cases}
  \stabled{d} & \text{if $\alpha\ge0$}\\
  \gsubclose_d &\text{if $\alpha\le0$}
\end{cases}
\]
    \end{lemma}
    \begin{proof}
      If $\alpha<0$ then $\xx^tQ\xx +\alpha\in\gsubclose_d$ by
      Lemma~\ref{lem:subdef-2}. Now if
      $f(x_1,\dots,x_d)\in\hb{d}$ then $f(\imag x_1,\dots,\imag
      x_d)\in\stabled{d}$.  Since $\gsubclose_d\subset\hb{d}$
      it follows that
\[
-\bigl( \xx^tQ\xx-\alpha\bigr) = (\imag \xx)^t Q (\imag\xx)+\alpha \in
\stabled{d}.
\]
    \end{proof}

\begin{figure}[ht]
    \centering
{
\psset{xunit=.5cm,yunit=.5cm}
    \begin{pspicture}(0,0)(10,10)
       \pscurve(0,9)(3,5)(4,0)
       \pscurve(5,0)(5,3)(7,0)
       \pscurve(8,0)(7,4)(10,3)
       \pscurve(1,10)(4,8)(3,10)
       \pscurve(0,10)(4,6)(5,5)(6,6)(10,6)
  \end{pspicture}
}
    
    \caption{The graph of a stable polynomial in two variables}
    \label{fig:stable}
  \end{figure}
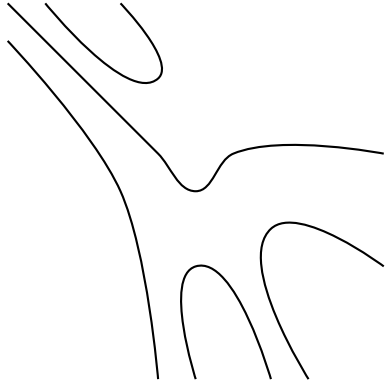

  The graphs of stable polynomials (Figure~\ref{fig:stable} is a
  sketch of one) are different from the graphs of polynomials in
  $\rup{2}$.

The following famous theorem of Heilman-Lieb determines a multiaffine
stable polynomial.

\begin{theorem}
Let $G = (V,E)$ be a graph, $V = \{1, \dots, n\}$. To
each edge $e = ij \in E$ assign a non-negative real number
$\lambda_{ij}$. Then the polynomial
\[
M_G(z) = \sum_{\text{M is a matching}} \prod_{ij\in M}
\lambda_{ij}z_iz_j
\]
is stable.
\end{theorem}

\section{Stable polynomials with complex coefficients}
\label{sec:stable-polyn-with}

In this section we consider polynomials with complex coefficients in
one variable.
Since a $90^o$ rotation moves the right half plane to the upper half
plane, we can easily transfer properties from $\polycpx$ to $\stabledc{1}$.  If $f$
has degree $n$ define $\phi(f) = \imag^n f(-\imag x)$. Note that the
factor of $\imag^n$ implies that $\phi$ is \emph{not} a linear
transformation.

\begin{lemma}\label{lem:polycpx-stable}
  If $f$ has positive leading coefficient then $f\in\stabledc{1}$ iff
  $\phi(f)\in\polycpx$.
\end{lemma}
\begin{proof}
  We may assume $f$ is monic.  All roots of $f(-\imag x)$ lie in the
  lower half plane, and the leading coefficient is $(-\imag)^n$. Thus,
  $\phi(f)$ is monic, and all roots are in the lower half plane.
\end{proof}

Interlacing in $\stabledc{1}$ is different form $\allpoly$. Any definition of
interlacing should satisfy $f \lesslesseq f'$, and this should mean
that all linear combinations of $f$ and $f'$ are stable. However, this
is false. If $f=x+1$ then $f - 2f' = x-1$ is not stable.

$\polycpx$ is invariant under the transformation $x\mapsto x+\alpha$,
where $\alpha\in\reals$, but $\stabledc{1}$ is preserved by $x\mapsto
x+\imag \alpha$.  In addition, $\stabledc{1}$ is preserved by $x\mapsto
x+\alpha$ if $\alpha>0$. Consequently, there are two new kinds of
interlacing along with $\hlace$.
\begin{align*}
\intertext{Assume that $f,g\in\stabledc{1}$}
f\iposlace g &\quad \text{iff}\quad &
  \forall\alpha\in\reals & \quad f+\alpha\,\imag\ g\in\stabledc{1}\\
f\hposlace g &\quad \text{iff}\quad &
  \forall\alpha\in\reals>0& \quad f+\alpha\, g\in\stabledc{1}\\ 
  f\hlace g &\quad \text{iff} & 
\forall \alpha\in\rhp & \quad f+\alpha g\in\stabledc{1}
  \end{align*}

We will see later that $f\iposlace f'$ and $f \hposlace f''$.
Taking limits shows that
\begin{lemma}
  If $f,g\in\stabledc{1}$ then $f\hlace g\implies f\iposlace g
  \ \text{and}\ 
  f\hposlace g$
\end{lemma}

\begin{remark}
  It is not true that   $f\hposlace g \implies f \iposlace g$. If
  $f=x^2$ and $g=1$ then clearly $f\hposlace g$, but $x^2+\imag$ has
  roots with positive and negative real parts.
\end{remark}
  The next result allows us to transfer interlacing results from
  $\polycpx$ to $\stabledc{1}$.

\begin{lemma}
  $f\hlace g$ in $\stabledc{1}$ if and only if
  $\phi(f)\iposlace\phi(g)$ in $\polycpx$.
\end{lemma}
\begin{proof}
Since $\phi$ is a bijection between $\stabledc{1}$ and $\polycpx$, the
computation below establishes the lemma.
  \begin{align*}
    \phi(f+\alpha\,\imag\,g) &= \imag^n(\, f(-\imag x) +
    \alpha\,\imag\,g(-\imag x)\,) \\
    &= \phi(f) - \alpha \,\imag^{n-1}g(-\imag x) \\
    &= \phi(f) - \alpha\,\phi(g)
  \end{align*}
\end{proof}

Most properties of interlacing in $\allpoly$ carry over to $\iposlace$
in $\stabledc{1}$, except for some occasional sign changes.

\begin{lemma} \label{lem:h-analog}
Assume $f,g,h\in\stabledc{1}$.

  \begin{enumerate}
  \item If $f\iposlace g$   then $fh\iposlace gh$.
    \item If $\sigma\in\rhp$ then $(x+\sigma)f \iposlace f$.
  \item  $f\iposlace f'$.
  \item If $f\in\allpolypos$ then $f(\imag\,\diffd)g\in\stabledc{1}$.
  \item If $f\iposlace g$ and $f\iposlace h$ then $f\iposlace \alpha g +
    \beta h$ for positive $\alpha,\beta$. In particular, $\alpha f +
    \beta g\in\stabledc{1}$.
  \item If $f\iposlace g \iposlace h$ in $\stabledc{1}$ then $f+h \iposlace g$.
  \item $x-\sigma\in\stabledc{1}$ iff $\sigma\in\rhp$.
  \item $x\hlace1$
  \end{enumerate}
\end{lemma}

\begin{proof}
The first is immediate. The second is a consequence of 
\[
(x+\sigma)f + \alpha\imag f = (x+\sigma+\alpha\imag)f
\]
and $x+\sigma+\alpha\imag\in\stabledc{1}$.  The chain rule shows that
$\phi(f') = \phi(f)'$. Since $f\in\stabledc{1}$ implies that
$\phi(f)\in\polycpx$, we know that $\phi(f)\lesslesseq \phi(f)' =
\phi(f')$.

We only need to check the next one for $f = x+a$ where $a>0$. We see
that $f(\imag\, \diffd)g = \imag g' + a g$ which is in $\stabledc{1}$ by the
previous part.

The next part follows from the corresponding result
(Lemma~\ref{lem:i-analog}) for $\polycpx$.  We have
interlacings
\[ \imag^n f(\imag x) \lesslesseq \imag^{n-1} g(\imag x)
\lesslesseq  \imag^{n-2} h(\imag x)
\]
in $\polycpx$ where $f$ has degree $n$. In $\polycpx$ we have
\[
\imag^n(f+h)(\imag x) = \imag^n f(\imag x) -\imag^{n-2} h(\imag x)
\lesslesseq \imag^{n-1} g(\imag x) 
\]
and the result follows.
The last two are trivial.
\end{proof}

We have orthogonal type recurrences in $\stabledc{1}$, except the last sign
is reversed.

\begin{cor}\label{cor:orthog-h}
If $a_i,c_i>0$, $\Re(\beta_i)\ge0$ then define
\begin{align*}
  p_{-1} &=0 \\
  p_0 &=1 \\
  p_n &= (a_n x + \beta_n)\,p_{n-1} + c_n\,p_{n-2}
\end{align*}
 All $p_n$ are stable, and
\[
\cdots p_3 \hlace p_2 \hlace p_1
\]
\end{cor}
\begin{proof}
  Using induction, these follow from the last part of
  Lemma~\ref{lem:h-analog} and the interlacings
\[
(x+ \beta_n)p_{n-1} \hlace p_{n-1} \hlace p_{n-2}
\]
\end{proof}

\begin{lemma} \label{lem:stable-1}
  Suppose that $f(x)=\prod(x-\sigma_k)$ is in $\stabledc{1}$.
 If all $\alpha_k$ are positive and
\begin{equation}\label{eqn:stable-1}
 g(x) = \sum_k \alpha_k \frac{f(x)}{x-\sigma_k} 
\end{equation}
then $g\in\stabledc{1}$ and $f\hlace g$.
\end{lemma}
\begin{proof}
  By additivity of interlacing for $\stabledc{1}$ it suffices to show that 
$f \hlace f/(x-\sigma_k)$, or $(x-\sigma_k)h\hlace h$ where $h =
f/(x-\sigma_k)$. This interlacing holds since $x-\sigma_k\hlace1$.
\end{proof}

Here are some lemmas that use the geometry of $\complexes$ and are
useful in proving that polynomials are stable.  The first lemma gives
useful information about the mappings determined by quotients.

  \begin{lemma}\label{lem:stable-fractions}
    Assume that $r,s>0$ and $\sigma\in\quada$.
    \begin{align*}
      \frac{1}{r+\sigma} & \in\quadd &
      \frac{1}{(r+\sigma)(s+\sigma )} & \in \text{lower half plane} \\
      \frac{\sigma}{r+\sigma} & \in\quada &
      \frac{\sigma}{(r+\sigma)(s+\sigma )} & \in \text{right half plane}
    \end{align*}
  \end{lemma}
  \begin{proof}
    Suppose $\sigma = a+b\imag$ where $a,b>0$. The left two
    follow from 
    \begin{align*}
      \frac{1}{r+\sigma} & = \frac{a+r-\imag b}{|r+\sigma|^2} & 
      \frac{\sigma}{r+\sigma} & = \frac{a^2+b^2+ar + \imag b r}{|r+\sigma|^2}\\
    \end{align*}
    Since $1/(r+\sigma)$ and $1/(s+\sigma)$ are in the fourth
    quadrant, their product is in the lower half plane. Since
    $\sigma/(r+\sigma)$ is in the first quadrant, and $1/(s+\sigma)$
    is in the fourth quadrant, their product is in the right half plane.
    
  \end{proof}

In $\allpoly$ we know that $f$ and $f''$ don't interlace, but
$f-\alpha f''\in\allpoly$ for all positive $\alpha$. There is a
similar result for $\stabledc{1}$. 

\begin{lemma}\label{lem:stable-2}
  Suppose that $f = \prod(x-\sigma_i)$ is in $\stabledc{1}$. If all
  $\alpha_{ij}$ are non-negative then
\begin{equation}\label{eqn:stable-2}
f  \hposlace \sum_{i\ne j} \alpha_{ij}\frac{f}{(x-\sigma_i)(x-\sigma_j)}
\end{equation}
\end{lemma}
\begin{proof}
  If $\tau$ is a root of the sum of the right and left hand sides of
  \eqref{eqn:stable-2} then, after dividing by $f(\tau)$, we know that
\[
1 + \sum_{i\ne j} \alpha_{ij}\frac{1}{(\tau-\sigma_i)(\tau-\sigma_j)} = 0
\]
By Lemma~\ref{lem:stable-fractions} the sum lies in the lower half plane.
Thus $\tau$ does not satisfy \eqref{eqn:stable-2}, so
$\tau$ must lie in the open left half plane. Similarly for the lower
right quadrant.
\end{proof}


If we restrict all the coefficients to be equal, the lemma gives an
alternate proof that
that $f\hposlace f''$.

\begin{lemma}
  If $f\in\stabledc{1}$, $\Re(\alpha)\ge0$, $\beta\ge0$ then 
\begin{align*}
\alpha f+\beta f' & \in\stabledc{1}.\\
\alpha f + \beta f'' & \in\stabledc{1}
\end{align*}
\end{lemma}
\begin{proof}
  Assume that $f = \prod(x-\sigma_i)$ where $\Re\sigma_i<0$. 
  Factoring out $f(x)$, let $\tau$ be a  solution of
  \[ \alpha +\beta \sum\frac{1}{x-\sigma_i} =0
\]
If $\Re(\tau)\ge0$ then $\Re 1/(\tau-\sigma_i)>0$, so the sum 
has positive real part, and can't equal $-\alpha/\beta$.  The second
follows from Lemma~\ref{lem:stable-2}.

\end{proof}

We know the following corollary holds more genearally in
$\stabledc{d}$, but in $\stabledc{1}$ we can prove it by induction.

\begin{cor}\label{cor:stable-fofd}
  If $f,g\in\stabledc{1}$ then $    f(\diffd)g  \in\stabledc{1}$.
\end{cor}
\begin{proof}
  Factoring $f$, it suffices to show that $g+\alpha g'\in\stabledc{1}$ for
  $\Re(\alpha)>0$. This follows from the lemma.
\end{proof}

Linear transformations on $\allpoly$ induce linear transformations on
$\stabledc{1}$. 

\begin{lemma}
  Suppose $T\colon{}\allpoly\longrightarrow\allpoly$ is a linear
  transformation that preserves degree, the sign of the leading
  coefficient, and strict interlacing. The linear transformation
\[ x^k\mapsto \imag^k\, T(x^k)(-\imag\,x) \]
maps $\stabledc{1}$ to itself.
\end{lemma}
\begin{proof}
  Since $T$ induces a linear transformation in $\polycpx$
  \mypage{cor:cpxclose-lt}, the conclusion follows from the diagram

\centerline{\xymatrix{
{\sum a_kx^k} \ar@{-->}[rr] \ar@{->}[d]_{\phi} &&
{\sum a_k (-\imag)^k\,T(x^k)(\imag\,x)} \\
{\imag^n\,\sum a_k(-\imag\,x)^k} \ar@{->}[rr] &&
{\imag^n\sum a_k(-\imag)^kT(x^k)} \ar@{->}[u]_{\phi^{-1}}
}}
\end{proof}

\index{Hermite polynomials!mapping in $\stabled{1}$}
\index{Laguerre polynomials!mapping in $\stabled{1}$}
\begin{cor}\label{cor:h-hermite}
  The transformations $x^k\mapsto \imag^kH_k(-\imag x)$ and $x^k\mapsto
  \imag^kL_k(\imag x)$
  map $\stabledc{1}$ to $\stabledc{1}$, and $\stabled{1}$ to $\stabled{1}$.
\end{cor}
\begin{proof}
  The maps $x^n\mapsto H_n(x)$ and $x^n\mapsto L_n(-x)$ satisfy the
  hypotheses of the lemma.  It follows from the above that they map
  $\stabledc{1}$ to itself. Since $\imag^kH_k(-\imag x)$ and $
  \imag^kL_k(\imag x)$ are real polynomials, the maps also send
  $\stabled{1}$ to itself.
\end{proof}

\index{Hadamard product}
  If all roots of two polynomials lie in given quadrants then we know
  the location of the roots of their Hadamard product. In the table
  below the notation $-\stabledc{1}$ means all roots in the right half
  plane, and $-\polycpx$ means all roots in the upper half plane.
  Recall that $\sum a_i x^i \ast \sum b_i x^i = \sum a_ib_ix^i$. 

  \begin{lemma}\label{lem:had-table} (The location of $f\ast g$)
\[    \begin{array}{lc|c|c|c|c|}
\multicolumn{2}{c}{}      &\multicolumn{4}{c}{Quadrant}\\
\multicolumn{1}{c}{} & \multicolumn{1}{c}{} & \multicolumn{1}{c}{1} & 
\multicolumn{1}{c}{2} & \multicolumn{1}{c}{3} & \multicolumn{1}{c}{4}
\\\cline{3-6}
&&&&&\\[-.1cm]
& 1 & \polycpx & -\stabledc{1} & -\polycpx & \stabledc{1} \\[.2cm]
& 2 & -\stabledc{1} & -\polycpx & \stabledc{1} & \polycpx \\[.2cm]
Quadrant\ & 3 & -\polycpx & \stabledc{1} & \polycpx & -\stabledc{1} \\[.2cm]
& 4 & \stabledc{1} & \polycpx & -\stabledc{1} & -\polycpx \\[.1cm]\cline{3-6}
\end{array}
\]
\end{lemma}
\begin{proof}
  If $f$ and $g$ are in the third quadrant then they are in
  $\polycpxpos$, and so their Hadamard product is in $\polycpx$.
  The Lemma now follows from the identity
\[
f(\alpha x)\ast g(\beta x) = \bigl[ f\ast g\bigr](\alpha\beta x)
\]
which holds for all complex $\alpha,\beta$.
\end{proof}

Note that this shows that the Hadamard product of stable polynomials
with complex coefficients is not necessarily stable. For example
\[
(x+1+\imag)\ast(x+1+2\imag) = x-1+3\imag.
\]

\section{The analytic closure}
\label{sec:analytic-closure-1}

We define the analytic closure of $\stabled{1}$ and $\stabledc{1}$ as usual:

\begin{align*}
\stabledf{d} & = \text{ the uniform closure of } \stabled{d}\\
\stabledcf{d} & = \text{ the uniform closure of } \stabledc{d}
\end{align*}

We can determine the exponentials that are in $\stabledc{1}$,

\begin{lemma}\ 

\begin{enumerate}
\item  
$    e^{\alpha x}, e^{\alpha x^2}$ are in $\stabledf{1}$ for all positive $\alpha$.
\item
$    e^{\alpha\,\imag\, x}$ is in $\stabledc{1}$ for all
$\alpha\in\reals$.
\item $\sinh(\alpha x)/x$ and $\cosh(\alpha x)$ are in $\stabledf{1}$ for
  $\alpha\in\reals$. 
\item
$    e^{- x},    e^{- x^2}$  are not in $\stabledc{1}$.
\end{enumerate}
\end{lemma}
\begin{proof}
  Consider the limits
  \begin{align*}
    e^{\alpha x} & = \lim_{n\rightarrow\infty} \bigl(1+ \frac{\alpha
      x}{n}\bigr)^n \\
    e^{\alpha x^2} & = \lim_{n\rightarrow\infty} \bigl(1+ \frac{\alpha
      (x+1/n)^2}{n}\bigr)^n \\
    e^{\alpha\, \imag\,x} & = \lim_{n\rightarrow\infty} \bigl(1+ \frac{\alpha
      \imag x}{n}\bigr)^n \\
  \end{align*}
In the first two cases  the part being exponentiated is in $\stabled{1}$,
and in $\stabledc{1}$ in the third, so the limits are as claimed.
The infinite product formulas 
\begin{align*}
  \sinh(x)/x &= \prod_{n=1}^\infty \left(1 + \frac{x^2}{n^2\pi^2}\right)
  &
\cosh(x) &= \prod_{n=1}^\infty \left(1 + \frac{4x^2}{(2n-1)^2\pi^2}\right)
\end{align*}
show that $\sinh(x)/x$ and $\cosh(x)$ are in $\stabledf{1}$.

To show non-existence we use Corollary~\ref{cor:stable-fofd}. If
$e^{-x}\in\stabledf{1}$ then $e^{-\diffd}{f} \in\stabled{1}$, but this
fails for $f = 2x+1$. If $e^{-x^2}\in\stabledf{1}$ then
$e^{-\diffd^2}f\in\stabled{1}$, and this fails for $f = (2x+1)^2+1$.
\end{proof}

\section{Real stable polynomials in one variable}
\label{sec:real-stable-polyn}

The class $\stabled{1}$ shares many properties with $\allpolypos$, as the
first result shows.

\begin{lemma}
  \
  \begin{enumerate}
  \item If $f,g$ have real coefficients then $fg\in\stabled{1}$ iff $f$
    and $g$ are in $\stabled{1}$.
  \item If $f\in\stabled{1}$ then all coefficients of $f$ have the same sign.
  \item If $f\in\stabled{1}$ then $f^{rev}\in\stabled{1}$.
  \end{enumerate}
\end{lemma}
\begin{proof}
  The first is obvious.  Factor monic $f\in\stabled{1}$ as
\[\prod(x-r_j)\,\cdot\,\prod (x-\sigma_k)(s-\overline{\sigma}_k)\]
where the $r_j$ are negative, and the $\sigma_k$ have negative real
part. Expanding $(x-\sigma_k)(s-\overline{\sigma}_k)=x^2 -
2\Re(\sigma_k)x+|\sigma_k|^2$ shows that all the factors have
positive coefficients, and so the coefficients of $f$ are all
positive.

The roots of $f^{rev}$ are the inverses of the roots of $f$. Since the
roots of $f$ lie in the open left half plane, so do the roots of
$f^{rev}$. 
\end{proof}

The converse of $2$ is false:
\[
(x-1-2 \imag) (x-1+2 \imag) (x+2-\imag) (x+2+\imag)
=
x^4 + 2 x^3 + 2 x^2 + 10 x + 25
\]
\index{positive coefficients!not stable}

There are simple conditions for the stability of polynomials of small
degree. 
\begin{lemma} \label{lem:stable-reasons}
Assume all these polynomials have real coefficients.
  \begin{enumerate}
    \item A quadratic $ax^2+bx+c$  is  stable
    if and only if all coefficients have the same sign.
  \item A cubic $x^3+ax^2+bx+c$ is stable if and only if all
    coefficients are positive and $ab>c$.
  \item A quartic $x^4+ax^3+bx^2+cx+d$ is stable if and only if all
    coefficients are positive and $abc> c^2+a^2d$.

  \end{enumerate}
  
\end{lemma}
\begin{proof}
  We know all coefficients must be positive. The first one is trivial.
  For the remaining two, we convert to $\polycpx$, and the interlacing
  condition becomes the inequalities.

  \begin{gather*}
    x^3+ax^2+bx+c \quad\text{is stable iff}\\
    x^3-bx + \imag(ax^2-c)\quad\text{is in $\polycpx$ iff}\\
    x^3 - bx \lesslesseq ax^2 -c \quad \text{iff}\\
    \sqrt{b} > \sqrt{c/a} \quad\text{which is equivalent to $ba>c$}
  \end{gather*}

  \begin{gather*}
    x^4+ax^3+bx^2+cx+d \quad\text{is stable iff}\\
    x^4 - bx^2 + d + \imag(ax^3-cx) \quad\text{is in $\polycpx$ iff}\\
    x^4 - bx^2 + d \lessless ax^3-cx \quad\text{iff}\\
    \frac{1}{2}\left(b-\sqrt{b^2-4d}\right)\le \frac{c}{a} \le
        \frac{1}{2}\left(b+\sqrt{b^2-4d}\right) \quad\text{iff}\\
        c^2+ a^2d < abc
  \end{gather*}
\end{proof}

In ove variable $F_E$ is a sum of even powers, and $f_O$ is a sum of
odd powers so we can refine the previous result.
Recall that  $f_e$ and $f_o$ are the even and odd parts of $f$.

\begin{lemma}\label{lem:stable-even}
  Assume that $f$ is a polynomial of degree $n$ with all positive
  coefficients. 
\[
f\in\stabled{1} \ \text{if and only if}\
\begin{cases}
  f_e\lessless f_o & n\ \text{even}\\
  f_e\greateq f_o & n\ \text{odd}
\end{cases}
\]
\end{lemma}
\begin{proof}
  We first assume that the degree $n$ of $f\in\stabled{1}$ is even. If we
  write $f(x) = f_e(x^2) + x f_o(x^2)$ then
  \[ \phi(f) = \imag^n\bigl[ f_e(-x^2) -\imag x\,f_o(-x^2)\bigr] \]
  Since $n$ is even $f_e(-x^2) \lessless x\,f_o(-x^2)$ and hence
  $f_e(x) \lessless f_o(x)$, which implies $\psi(f)\in\polycpxpos$.
  The converse direction follows by reversing the argument.

  If $n=2m+1$ then
\[\phi(f) = (-1)^m \left(x f_o(-x^2) + \imag f_e(-x^2)\right)\]
which implies that $xf_o(-x^2)\lessless f_e(-x^2)$, and hence
$f_e\greateq f_o$. 
\end{proof}

It follows from the lemma that if the degree of $f\in\stabled{1}$ is even
then $f_e+\imag f_o\in\polycpx$. If the degree is odd, then 
$(1-\imag)(f_e + \imag f_o)\in\polycpx$.

If $\tau$ is a root of $f\in\stabled{1}$ then $f/(x-\tau)$ is not necessarily in
$\stabled{1}$. If we pair the complex roots then we get interlacing:

\begin{lemma}\label{lem:stable-int}
  Suppose that $f(x)\in\stabled{1}$, and write 
  \[f(x) = \prod_{j=1}^n (x-r_j)\,\cdot\,\prod_{k=1}^m
  (x-\sigma_k)(s-\overline{\sigma}_k)\] If all $\alpha_k,\beta_j$ are
  non-negative and $\Re(\sigma_k)<0$ then
\[
f\hlace
\sum_{j=1}^n \beta_j \frac{f(x)}{x-r_j} + \sum_{k=1}^m \alpha_k 
\frac{x - \Re(\sigma_k)}{(x-\sigma_k)(x-\overline{\sigma}_k)}\,f(x)
\]
\end{lemma}
\begin{proof}
  The conclusion follows from Lemma~\ref{lem:stable-1} and the observation
  that
\[
\frac{x - \Re(\sigma_k)}{(x-\sigma_k)(x-\overline{\sigma}_k)} =
\frac{1}{2}\left(\frac{1}{x-\sigma_k} + 
\frac{1}{x-\overline{\sigma}_k} \right)
\]
\end{proof}

\begin{example}
  It is not easy to determine all polynomials interlacing a given
  polynomial. For example, we will find all linear polynomials
  interlacing $(x+1)^2$ by translating to $\polycpx$.  Assume that
  $(x+1)^2 \lesslesseq \sigma x + \tau$. Then
\[
(x+1)^2 + \alpha \imag (\sigma x+ \tau) \in\stabledc{1}
\]for all real $\alpha$. Since the sum of the roots is
$-2-\alpha\imag\sigma$ and lies in the left half plane for all
$\alpha$ it follows that $\Im(\sigma)=0$. After scaling, we may assume
that
\begin{align*}
(x+1)^2 & \lesslesseq x+a + b \imag.\\
\intertext{ We now translate this interlacing to $\polycpx$:}
x^2-1 + 2\imag x & \lesslesseq x + a\imag -b
\end{align*}
The determinant condition is
\[
0 > \begin{vmatrix}
  x^2-1 & 2x \\ x-b & a
\end{vmatrix} =
(a-2) x^2 + 2xb -a
\]
Solving the equations yields
\begin{enumerate}
\item $(x+1)^2 \lesslesseq  x + \tau$ in $\stabledc{1}$ if and only
  if $\tau$ lies in a circle of radius one centered at $1$.
\item $(x+1)^2 \lesslesseq  x + t$ in $\stabled{1}$ if and only
$0<t<2$.
\end{enumerate}

If $0<a<1$ then we can give an explicit polynomial in $\stabled{2}$
that shows that $(x+1)^2\lessless x+a$.

\begin{gather*}
\begin{vmatrix}
 y+x \left(\alpha ^2+\sqrt{\alpha ^2+1} \alpha +1\right)+1 & \alpha  \\
 -\alpha  & y+x \left(\alpha ^2-\sqrt{\alpha ^2+1} \alpha +1\right)+1
\end{vmatrix} \\
 = (1+\alpha^2)(1+x)^2 + 2(1+\alpha^2)\biggl(x +
 \frac{1}{1+\alpha^2}\biggr)y + y^2
\end{gather*}

We have seen the first kind of region in $\polycpx$. Restricted to
real polynomials, we get an interval  in $\stabled{1}$, but just a single
point  in $\allpolypos$. 

A similar calculation shows the following:
\[
(x+1)(x+2)\lesslesseq (x+1)+a(x+2)\ \text{in}\ \stabled{1}\ \text{iff}\ a> -1/2
\]

\end{example}

  If $f=\prod(x-w_i) \in\allpoly$ then in Lemma~\ref{lem:sign-quant}
  and Lemma~\ref{lem:sign-1} we gave conditions on constants $c_i$
  such that $\sum c_if/(x-w_i)$ is in $\allpoly$.  If
  $f\in\allpolypos$ we now ask what choice of constants determines a
  stable interlacing polynomial. The example above shows that some
  constants can be negative.

  \begin{lemma}
    Suppose that $f = (x+r_1)\cdots(x+r_n)$ where
    $0<r_1<\cdots<r_n$. If
\[
a_k+a_{k+1}+\cdots +a_n\ge0 \ \text{for k=1,\dots,n}
\]
and at least one sum is positive, then $f \hlace \displaystyle \sum
a_k \frac{f}{x+r_k}$ in $\stabled{1}$.
  \end{lemma}
  \begin{proof}
    We need to verify that the quotient 
\[
\frac{1}{f} \sum a_k \frac{f}{x+r_k} = \sum \frac{a_k}{x+r_k}
\]
maps $\Re(\sigma)>0$ to the right half plane. Let
$\sigma=\alpha+\beta\imag$ with $\alpha>0$. Then,
\begin{align*}
\Re\ \sum \frac{a_k}{\sigma+r_k} &= \Re\ \sum a_k
\frac{r_k+\alpha-\beta\imag}{(r_k+\alpha)^2+\beta^2} \\
& = \sum a_k \frac{r_k+\alpha}{(r_k+\alpha)^2+\beta^2} \\
&\ge
\frac{1}{(r_n+\alpha)^2+\beta^2}\sum a_k(\alpha+r_k)
\end{align*}
This is positive since
\begin{align*}
  \sum a_k(\alpha+r_k) &= 
   (\alpha+r_1)(a_1+ \cdots + a_n) + (r_2-r_1) (a_2+\cdots + a_n) + \\
  & (r_3-r_2)(a_3+\cdots+a_n ) + \cdots + (r_n-r_{n-1})a_n >0
\end{align*}
by hypothesis.
\end{proof}
\begin{example}
  If $0<r_1<r_2$ and we choose $a_1+a_2>0, a_2>0$ then
  $f=(x+r_1)(x+r_2)\lessless a_1(x+r_2) + a_2(x+r_1)$ in $\stabled{1}$.
  Scaling so that $a_1+a_2=1$, this is $f \lesslesseq x+
  r_1+(r_2-r_1)a_1$. Since $1\ge a_1$, it follows that
\[ (x+r_1)(x+r_2) \lessless x+t \ \text{in $\stabled{1}$ if}\ 0<t<r_2
\]

We can determine all $t\in\reals$ for which $(x+r_1)(x+r_2)\lessless
x+t\ \text{in $\stabled{1}$}$. Converting to $\polycpx$ yields
$x^2+(r_1+r_2)x\imag - r_1r_2 \lessless x+\imag t$. The interlacing
conditions are satisfied, and the determinant condition yields
\[ (x+r_1)(x+r_2) \lessless x+t \ \text{in $\stabled{1}$ iff}\ 0<t<r_1+r_2
\]

\end{example}

Although the example shows that the hypotheses of the lemma do not
determine all interlacing polynomials, there is a sense in which they
are the best possible.

\begin{lemma}
  If $a_1,\dots,a_n$ have the property that for all $0<r_1<\cdots<r_n$ we
  have 
\[
f \lessless \sum_{k=1}^n a_k
\frac{f}{x+r_k}
\ \text{in $\stabled{1}$}
\]
where $f=(x+r_1)\cdots(x+r_n)$ then
\[
a_k+a_{k+1}+\cdots +a_n\ge0 \ \text{for k=1,\dots,n}
\]
\end{lemma}
\begin{proof}
  Assume that $a_k+\cdots+a_n<0$. Take $r=r_1=\cdots=r_{k-1}$ and
  $s=r_k=\dots=r_n$. Then the real part of the quotient, where
  $\sigma=\alpha+\beta\imag$,  is equal to 
  \begin{gather*}
    \sum_{i=1}^{k-1} a_i \frac{r+\alpha}{|\sigma+r|^2} + \sum_{i=k}^na_i
    \frac{s+\alpha}{|\sigma+s|^2} =\\
\frac{r+\alpha}{|\sigma+r|^2}\sum_{i=1}^{k-1}a_i +
\frac{s+\alpha}{|\sigma+s|^2}\sum _{i=k}^n a_i
  \end{gather*}
  Since the second sum is negative this expression will be negative if
  we can choose constants so that
\[
 \frac{r+\alpha}{|\sigma+r|^2} \biggl/ \frac{s+\alpha}{|\sigma+s|^2} < \epsilon
\]
for appropriately small positive $\epsilon$. 
This is the case if we pick $|\sigma| \gg s \gg r \gg \alpha$.
\end{proof}

The next lemma follows from the corresponding facts for $\stabledc{1}$. 

\begin{lemma}
  If $f,g\in\stabled{1}$ then $f(\diffd)g  \in\stabled{1}$. If
  $f\in\allpolypos$ then $f(\diffd^2)g\in\stabled{1}$. 
\end{lemma}

Combining the positive and negative cases, we get
\begin{cor}\label{cor:stable-fpp}
  If $f\in\allpolypos$ then
\[
f+ \alpha \,f'' \in
\begin{cases}
  \allpoly & \alpha\le 0\\
  \stabled{1} & \alpha>0
\end{cases}
\]
\end{cor}

  We describe some interactions among $\allpolypos$, $\stabled{1}$,
  $\stabledc{1}$ and $\polycpx$.

  \begin{lemma}
    If $F\in\polycpx$ has all coefficients in the first quadrant then
    $F\in\stabledc{1}$.
  \end{lemma}
  \begin{proof}
    Write $F=f+\imag g$. If $f = \prod(x-r_k)$ and $g/f = \sum
    \frac{a_k}{x-r_k}$ where all $a_k$ are positive then $(f+\imag
    g)(\sigma)=0$ implies
\[
1 + \imag\, \sum \frac{a_i}{|\sigma-r_i|^2} (\overline{\sigma}-r_i)=0
\]
It follows that the real part of $\sigma$ lies in the interval $(\min
r_k,\max r_k)$. If all the coefficients are positive then all roots of
$f$ are negative, so the real parts of $\sigma$ are negative, and
$F\in\stabledc{1}$. 
  \end{proof}

The intersection of $\stabledc{1}$ and $\polycpx$ contains polynomials
that do not have all positive coefficients. If $f+\imag g$ is 
\begin{gather*}
   (x-.08) (x+.7) (x+1) (x+1.7) (x+2) +\imag (x+.3)
  (x+.9) (x+1.7) (x+1.8)\\
  \intertext{then the roots are all are in the third quadrant}
 (-1.9-0.1 \imag,-1.7-0.003 \imag,-0.9-0.04 \imag,-0.5-0.5 \imag,-0.2-0.2 \imag)
\end{gather*}
\index{interlacing!in $\allpolypos$ and $\stabled{1}$}
\begin{lemma}
    Assume that $f,g\in\allpolypos$. 
    If $f\place g$ in $\allpoly$ then $f\hlace g$ in $\stabled{1}$. 
  \end{lemma}
  \begin{proof}
    If $f\place g$ then $f+yg\in\gsubplus_2\subset\stabled{2}$, so
    $f\hlace g$. 
  \end{proof}

The converse is false.

\begin{remark}
    The image of $\stabled{1}$ in $\polycpx$ is easily described. If $f =
    \sum^n a_kx^k$ then
\begin{multline*}
 \phi(f) = \imag^n \sum^n x^k a_k (-\imag x)^k \\
= \sum_{k\equiv n\pmod{2}} x^k a^k(-1)^{(n-k)/2} +
 \imag\,\sum_{k\not\equiv n\pmod{2}} x^k a^k(-1)^{(n-k-1)/2} 
\end{multline*}
Thus the real part consists of the terms of even degree with
alternating signs and the odd part consists of terms of odd degree
with alternating signs, or vice versa. This is equivalent to
$\overline{g(x)} = g(-x)$. Since $\phi$ is just a rotation, the roots
of $\phi(f)$ are all in the lower half plane, and are symmetric around
the imaginary axis. Note that the image of $\stabled{1}$ isn't $\polycpxpos$.
\end{remark}

\index{truncation}
  \begin{lemma}\label{lem:trunc}
    If $f\in\allpolypos$, and $T(x^i)=x^{i-1}$ for $i\ge1$ and $T(1)=0$
    then $T\colon\allpolypos\longrightarrow\stabled{1}$.
  \end{lemma}
  \begin{proof}
    Note that $T(f) = (f(x)-f(0))/x$. If $f = \prod(x+a_i)$ where all
    $a_i$ are positive then $\alpha$ is a root of $T(f)$ implies that
    $f(\alpha)=f(0)$, and therefore
\[
\prod\biggl(\frac{x}{a_i}+1\biggr) = 1
\]
If $x$ has positive real part then $x/a_i$ also has positive real
part, so $|1+x/a_i|>1$. Consequently the product above can't be $1$,
and so no roots are in the right half plane.
  \end{proof}

\section{Coefficient inequalities in $\stabled{1}$}
\label{sec:coeff-ineq-stable}

Newton's inequalities do not hold for $\stabled{1}$. This can be seen
easily in the quadratic case. If $f = a_0+a_1x+a_2x^2\in\stabled{1}$ then

\[
\frac{a_1^2}{a_0a_2}\quad
\begin{cases}
  \quad\ge 4 & \text{if the roots are real}\\[.3cm]
  \quad\le 4 & \text{if the roots are complex}
\end{cases}
\]
Moreover, the quotient can take on any value in $(0,\infty)$.

However, there are many inequalities involving the coefficients. For
instance, if $ f = \sum a_i x^i\in\stabled{1}$ then we know that
\[
a_0 - a_2 x + \cdots \lessless a_1 - a_3 x + \cdots 
\]
and therefore $\smalltwodet{a_1}{a_3}{a_0}{a_2}>0$ which is equivalent
to
\begin{equation}\label{eqn:stable-det-1}
\frac{a_1a_2}{a_0a_3} > 1
\end{equation}
A much more general result is that the Hurwitz matrix is totally
positive.\seepage{rem:hurwitz-mat} This means that 
\[
\begin{pmatrix}
  a_1 & a_3 & a_5 & \dots & \dots \\
a_0 & a_2 & a_4 & \dots & \dots \\
0&  a_1 & a_3 & a_5 & \dots  \\
0& a_0 & a_2 & a_4 & \dots  \\
\vdots & \vdots & \vdots & \vdots & \ddots
\end{pmatrix}
\]
is totally positive. For example, if we choose $f\in\stabled{1}(4)$ then 
\[
0 \le
\begin{vmatrix}
  a_1 & a_3 & 0 \\ a_0 & a_2 & a_4 \\ 0 & a_1 & a_3
\end{vmatrix}
= a_1a_2a_3 - a_1^2a_4- a_3^2a_0
\]
This inequality is the necessary and sufficient condition
(Lemma~\ref{lem:stable-reasons}) for a positive quartic to be stable.
\index{Hurwitz matrix}

  Here's a simple yet surprising conclusion about initial segments of
  stable polynomials.
  \begin{lemma}\label{lem:unexpected}
    If $\sum a_ix^i\in\allpolypos$ then
\[
\sum_{k=0}^1 a_ix^i \poslace \sum_{k=0}^2 a_ix^i \poslace 
\sum_{k=0}^3 a_ix^i \poslace \sum_{k=0}^4 a_ix^i 
\quad\text{in $\stabled{1}$}
\]
In particular, each of these partial sums is stable. 
  \end{lemma}
  \begin{proof} We begin with  the middle interlacing; the first is trivial.
 We need to verify that
\[ \alpha(a_0+a_1x+a_2x^2) + (a_0+a_1x+a_2x^2+a_3x^3) =
(1+\alpha)(a_0+a_1x+a_2x^2) + a_3x^3
\] 
is stable. From \eqref{eqn:stable-det-1} we know that $a_1a_2 > a_0a_3$, so
\[
(1+\alpha)a_1\,(1+\alpha)a_2 > (1+\alpha)a_0\,a_3
\]
and therefore we have positive interlacing. For the last interlacing,
we need to show that 
\[ (1+\alpha)(a_0+a_1x+a_2x^2+a_3x^3) + a_4x^4
\] is stable for positive $\alpha$. To do this we need to show that
\[
(1+\alpha)^3a_1a_2a_3 > a_4(1+\alpha)^2a_1^2 + (1+\alpha)^2a_3^2a_0
\]
The Newton inequalities give
\begin{align*}
  a_1^2 & \ge \frac{8}{3} a_0a_2 &
  a_2^2 & \ge \frac{9}{4} a_1a_3 &
  a_3^2 & \ge \frac{8}{3} a_2a_4 &
\end{align*}
and as a consequence
\begin{align*}
  a_2a_3 & \ge 6 a_1a_4 & a_1a_2 &\ge 6 a_0a_3
\end{align*}
and so
\begin{align*}
  (1+\alpha)^3 a_1a_2a_3 &\ge 6(1+\alpha)^2a_1^2a_4 \\
  (1+\alpha)^3 a_1a_2a_3 &\ge 6(1+\alpha)^2a_3^2a_0 \\
\end{align*}
Addition yields the desired result.
  \end{proof}

Note that the first two only require that $f$ is stable, but the last
interlacing uses $f\in\allpolypos$. 

  \begin{remark}
    It is not true that the higher partial sums are stable. For
    $k\ge5$ the partial sum $1+x+\cdots+x^k/k!$ of the exponential
    series is not stable \cite{zemyan}. The sum of the first $k+1$
    coefficients of $(1+x/n)^n$ converges to $1+x+\cdots+x^k/k!$, and
    $(1+x/n)^n$ is in $\allpolypos$. Thus, for $n$ sufficiently large
    the partial sums aren't stable. For instance, the sum of the first
    six terms of $(1+x/n)^n$ is stable for for $n\le23$ and not stable
    for $n>23$.
  \end{remark}

\section{Positive interlacing }
\label{sec:posit-interl}

We study  positive interlacings $f\hposlace g$. If $f,g$ in
$\allpoly$, then one way of showing $f\lessless g$ is sign
interlacing. This can be viewed as describing the behavior of the map
\[
\frac{g}{f}\colon \reals \longrightarrow \reals\setminus\infty.
\]
In order to prove that $f\hposlace g$ in $\stabled{1}$, we look at the
image of the first quadrant rather than the real line. The next lemma
gives criteria for $\poslace$ and $\hposlace$. The proof is omitted.

\begin{lemma}\label{lem:quot-geom}\label{lem:quot-pos}\ 
  \begin{enumerate}
  \item $f\poslace g$ iff $\frac{f}{g}\colon \quada\longrightarrow \rhp$.
  \item $f\hposlace g$ iff $\frac{f}{g}\colon \quada\longrightarrow \complexes\setminus(-\infty,0)$.
  \end{enumerate}
  
\end{lemma}

If $f,g\in\allpoly$ then the implication $ f\lesslesseq g \implies
f\hposlace g $ is trivial. It's also easy in $\stabled{1}$.

\begin{lemma}\label{lem:stable-less-pless}
  If $f,g\in\stabled{1}$ and $f\lesslesseq g$ implies $f\hposlace g$.
\end{lemma}
\begin{proof}
  We know that the image of the closed right plane
  under $x\mapsto \frac{f}{g}(x)$ lies in the open right half
  plane. Thus, $f\hposlace g$.
\end{proof}

We have seen \seepage{lem:pos-1} that if $f\hposlace g,h$ in
$\allpolypos$ then $g+h$ might not have all real roots. However,
\begin{lemma}\label{lem:pos-1-stab}
  If $f \hposlace g,h$ are in $\allpolypos$ then $g+h$ is stable.
\end{lemma}
\begin{proof}
From Lemma~\ref{lem:quot-geom} we know that evaluation of $g/f$ and
$h/f$ at a point in the closed first quadrant gives a point in the
right half plane, and so their sum lies there too.
\end{proof}

  \begin{lemma}
    If $f\hlace g$  then $f \hposlace x\,g$ 
  \end{lemma}
  \begin{proof}
    We assume $f,g\in\stabled{1}$; the case $f,g\in\stabledc{1}$ is similar.
    Since $f\lessless g$ we know that
    $f/g\colon\overline{\quada}\longrightarrow \rhp$. 
    To show $f \hposlace xg$ we choose $\alpha\ge0$ and prove that
    $f+\alpha xg\in\stabled{1}$. We do this by establishing that $f+\alpha
    xg \hposlace f$. If $\sigma\in\overline{\quada}$ then
\[
\frac{f+\alpha xg}{f} (\sigma) = 1 + \alpha \sigma
\frac{f(\sigma)}{g(\sigma)}
\]
Now $(g/f)(\sigma)\in\rhp$ since $f\lessless g$ and
$\alpha\sigma\in\overline{\rhp}$ so $\frac{f+\alpha xg}{f}(\sigma)$
    is in $\rhp+\overline{\rhp}\times\rhp = \complexes\setminus(0,\infty)$.
  \end{proof}

Next is a Leibniz-like rule for second derivatives.

\begin{lemma}
  If $f,g\in\allpolypos$ then $fg'' + f''g\in\stabled{1}$.
\end{lemma}
\begin{proof}
  Recall Lemma~\ref{lem:quot-geom}.  If we evaluate $f/f''$ and
  $g/g''$ at a point in the closed first  quadrant then we get
  a point in the lower  half plane. Thus $f/f'' + g/g''$
  is stable.
\end{proof}

Note that $fg''+f''g$ might not lie in $\allpolypos$. For example, if
$f=(x+1)^2$ and $g = (x+2)^2$ then $fg''+f''g$ is always positive, and
so has no real zeros.

Positive interlacings can be multiplied.

\begin{lemma}\label{lem:mult-lace}
  If $f\lessless g$ and $h\lessless k$ in $\allpolypos$ 
  \begin{enumerate}
  \item  $fh  \hposlace gk$ in $\stabled{1}$.
  \item  $fk  \hposlace gh$ in $\stabled{1}$.
  \item  $f^2 \hposlace g^2$ in  $\stabled{1}$. 
  \end{enumerate} 
\end{lemma}
\begin{proof}
  We know that $f/g$ and $h/k$ both map the first quadrant to itself,
  so $(fh)/(gk)$ maps the first quadrant to the upper half plane, so
  $fh\hposlace gk$. The last is a special case of the first. Since
  $k/h$ maps $\quada$ to $\quadd$, we see that $(fh)/(gk)$ maps to
  $\quadd\cup\overline{\quada}\cup\quadb\setminus0$, from which the
  second statement follows.
 \end{proof}

\begin{cor}\label{cor:orthog-h-2}
If $a_i,b_i,c_i>0$,  then the sequence defined by
\begin{align*}
  p_{-1} &=0 \\
  p_0 &=1 \\
  p_n &= (a_n x + b_n)\,p_{n-1} + c_n\,p_{n-2}
\end{align*}
then all $p_n$ are stable, and
\[
\dots p_3 \hposlace p_2 \hposlace p_1
\]
\end{cor}

\begin{proof}
  From Corollary~\ref{cor:orthog-h} we know that $p_n \lesslesseq
  p_{n-1}$, and so $p_n\hposlace p_{n-1}$ by
  Lemma~\ref{lem:stable-less-pless}.
\end{proof}

The preceding lemma can be considered as a construction of stable
orthogonal polynomials. \index{orthogonal polynomials!stable}
There is a simple case where the roots are all known explicitly
\cite{hoggatt}. If $f_1=1$, $f_2=x$ and 
\[  f_{n+1}(x) = x\,f_n(x) + f_{n-1}(x)\]
then the polynomials are known as the \index{Fibonacci polynomials} 
Fibonacci polynomials  and have purely imaginary roots:
\[2\imag\cos(k\pi/n),\quad 1\le k \le n-1. \]

\begin{remark}
    If we take two stable polynomials
    \begin{align*}
f &=  x^4+6 x^3+14 x^2+10 x+4 & g & =        x^2+2 x+2
    \end{align*}
    then how do we go about showing that $f \hposlace g$?  If we use
    the definition then we need to show that $f+t\,g\in\stabled{1}$ for
    all positive $t$. Since this is a quartic we could apply the
    criterion (Lemma~\ref{lem:stable-reasons}) for a quartic to be
    stable. A more general technique  uses the fact
    (Lemma~\ref{lem:quot-pos}) that $f/g$ must map the closed first
    quadrant to $\complexes\setminus(-\infty,0]$. If $\sigma\in\overline{\quada}$
    then this would follow from $\Im(f(\sigma)/g(\sigma))>0$. To show
    this we just need to show that the polynomial
    $\Im(f(a+b\imag)g(a-b\imag))$ is  positive for $a,b\ge0$. If we are lucky,
    this polynomial will have no minus signs, as in this case:
\begin{multline*}
{f(a+b\imag)}g(a-b\imag) = 
2 b \left(a^2+b^2\right)^3 \times
\\
(a^5+6 a^4+2 b^2 a^3+16 a^3+8 b^2
  a^2+27 a^2+b^4 a+8 b^2 a+24 a+2 b^4+3 b^2+6) 
\end{multline*}
This approach doesn't quite work here, since $a=b=0$ has imaginary
part $0$. However, $f(0)/g(0)=2$ which is not on the negative axis, so
we conclude that $f\hposlace g$.

 \end{remark}

The \index{Bezout polynomial}Bezout polynomial is stable.

\begin{lemma}
  The Bezout polynomial
  $\frac{1}{x-y}\smalltwodet{f(x)}{f(y)}{g(x)}{g(y)}$ is stable if
  $f\place g\in\allpolypos$. 
\end{lemma}
\begin{proof}
  Recall \eqref{eqn:bezout-poly}
  \begin{align*}
    B(x,y) &= \sum a_i \frac{f(x)}{x-r_i}\frac{f(y)}{y-r_i} \\
    \frac{B(x,y)}{f(x)f(y)} &= \sum a_ i\frac{1}{x-r_i}\frac{1}{y-r_i} 
  \end{align*}
  If $f,g\in\allpolypos$,  $\sigma\in Q_1$, $\tau\in Q_1$ then 
$\frac{1}{\sigma-r_i}\frac{1}{\tau-r_i}$ is in the lower half plane,
and so the sum misses $(-\infty,0)$. If $\sigma\in Q_1$, $\tau\in Q_4$ then 
$\frac{1}{\sigma-r_i}\frac{1}{\tau-r_i}$ is in the right half plane,
and so the sum again misses $(-\infty,0)$. The remaining cases are
similar, so we conclude $f(x)f(y) \hposlace B(x,y)$.

\end{proof}

When a polynomial is always positive then it is not necessarily
stable. Here's a familiar case where it is. 

\begin{lemma}\label{lem:stable-square-1}
  If $f\lessless g$ in $\allpolypos$ then $\smalltwodet{f}{g}{f'}{g'}$
  is stable. If $f$ has all distinct roots then 
 $\smalltwodet{f}{f'}{f'}{f''}$ is  stable.
\end{lemma}
\begin{proof}
  If we write $f=\prod(x+r_i)$ and $g = \sum a_i f/(x+r_i)$ then
  \[ f'g-fg' = f^2\,\sum \frac{a_i}{(x+r_i )^2} \] and
  Lemma~\ref{lem:quot-geom} applies. Use $f\lessless f'$ for the
  second one.
\end{proof}

\begin{example}
  Here's a nice example of the lemma. Assume that $f_0,f_1,f_2,\dots$
  is an orthogonal polynomial sequence with strict
  interlacing: $f_{k+1}\lessless f_k$ for $k=1,\dots$. The
  Christoffel-Darboux formula \cite{szego} states that

\index{Christoffel-Darboux formula}

\begin{gather*}
  f_0(x)f_0(y) + f_1(x)f_1(y)+\cdots+f_n(x)f_n(y) =
  \frac{k_n}{k_{n+1}}\frac{1}{x-y}
\begin{vmatrix} f_{n}(y) & f_{n+1}(y) \\ f_{n}(x) &
  f_{n+1}(x)\end{vmatrix}\\
\intertext{for positive constants $k_n$. If we let $x=y$ then }
 f_0(x)^2 + f_1(x)^2+\cdots+f_n(x)^2 =
  \frac{k_n}{k_{n+1}}
\begin{vmatrix} f_{n}(x) & f_{n+1}(x) \\ f_{n}'(x) &
  f_{n+1}'(x)\end{vmatrix}\\
\end{gather*}
It follows from the Lemma that

\[
 f_0(x)^2 + f_1(x)^2+\cdots+f_n(x)^2 \quad \text{is stable}
\]

\end{example}

We can also state Lemma~\ref{lem:stable-square-1} in terms of
Wronskians: if $f\lessless g$ in $\allpolypos$ then the Wronskian
$W(f,g)=fg'-f'g$ is stable.  \index{Wronskian}

More generally we have
  \begin{lemma}\label{lem:stable-square-3}
    The square determined by $|I + xD_1+yD_2+zD_3| $, where all $D_i$
    are positive definite, $D_1$ is diagonal, and $D_2,D_3$ have all
    positive entries, yields a stable polynomial.
  \end{lemma}
  \begin{proof}
    The assumptions mean that we are to show that
    $\smalltwodet{f}{g}{h}{k}$ is stable, where
\[
f+gy+hz+kyz + \cdots = | I + x D_1 + yD_2 + z D_3|
\]
and all matrices are positive definite.  Assume that $D_1$ is
diagonal $(r_i)$, $D_2 = (a_{ij})$, $D_3=(b_{ij})$. We have
\begin{align*}
  f &= \prod(1+r_ix)\\ 
  g &= \sum_i a_{ii} f/(1+r_ix)\\
  h &= \sum_j b_{jj} f/(1+r_jx)\\ 
  k &= \sum_{i\ne j} \begin{vmatrix} a_{ii}&b_{ij}\\
    a_{ij}&b_{jj}\end{vmatrix}  \frac{f}{(1+r_ix)(1+r_jx)}\\ 
gh - fk &= f^2\left(
 \sum_i \frac{a_{ii}b_{ii}}{(1+r_ix)^2} +\sum_{i\ne j}
\frac{a_{ij}b_{ij}}{(1+r_ix)(1+r_jx)}\right)
\end{align*}
Since $D_1,D_2,D_3$ are positive definite, their diagonals are
positive, so $(gh-fk)/f^2$ lies in the lower half plane by
Lemma~\ref{lem:stable-fractions}. Lemma~\ref{lem:quot-pos} implies
that $f^2\hposlace gh-fk$, and therefore $gh-fk$ is stable.

  \end{proof}

\begin{cor}
  If $f,g,h,k$ are as above, and $0<\alpha\le1$ then
  $\smalltwodet{\alpha f}{g}{h}{k}$ is stable.
\end{cor}
\begin{proof}
  The above computation shows that
\[
gh - \alpha f k = gh - fk  + (1-\alpha)fk = 
f^2\left( \sum \frac{a_ib_i}{(1+r_ix)^2} + (1-\alpha)\frac{k}{f}\right)
\]
From Lemma~\ref{lem:quot-geom} we know that if we evaluate at
a point in the first quadrant then all the terms lie in the lower
half plane.
\end{proof}

We can slightly improve Lemma~\ref{lem:stable-square-1}.

  \begin{cor}\label{cor:stable-sq-2}
    If $f\in\allpolypos$ and $0<\alpha\le1$ then $\smalltwodet{\alpha
      f}{f'}{f'}{f''}$ is stable.
  \end{cor}
  \begin{proof}
Applying  the calculations of the lemma to
\[
f(x+y+z) = f(x) + f'(x)y + f'(x)z + f''(x)yz/2 + \cdots \in\gsubplus_3
\]
shows that
$(\alpha/2) ff''- f'f'$ is stable.
\end{proof}

\index{Newton's inequalities!and stability}
Newton's inequality says that a certain determinant is positive. If we
introduce a second variable, this becomes a statement that a polynomial
is stable. The first part follows from
Lemma~\ref{lem:stable-square-3}. 

\begin{lemma}\label{lem:stable-ff}
  If $f = \sum f_i(x)y^i$ is in $\gsubplus_2$ then 
\[
\begin{matrix}
f_0^2 & \hposlace & f_1^2 - f_0f_2 & \ \text{in $\stabled{1}$}\\
f_0\,f_1 & \hposlace & f_1^2 - f_0f_2 & \ \text{in $\stabled{1}$}\\
f_0^2 & \hposlace &x(f_1^2 - f_0f_2) & \ \text{in $\stabled{1}$}\\
f_0\,f_1 & \hposlace & x(f_1^2 - f_0f_2) & \ \text{in $\stabled{1}$}\\
\end{matrix}
\]
In particular, $f_1^2 - f_0f_2$ is (weakly) stable, and  all  coefficients are non-negative.
\end{lemma}
\begin{proof}
  We compute $f_1^2 - f_0f_2$ using the determinant representation for
  $f$. So, assume that $f = |I + xD_1 + yD_2|$. We may assume that
  $D_1$ is diagonal, $D_1 = diag(r_i)$, and that $D_2 =
  (d_{ij})$. Then
  \begin{align*}
    f_0(x) &= \prod_i (1+ r_i x)\\
f_1(x)  &= \sum_i d_{ii} \frac{f_0}{1+r_ix} \\
f_2(x) &= \sum_{i<j} \begin{vmatrix}d_{ii}&d_{ij}\\ d_{ij}&d_{jj}
\end{vmatrix} \frac{f_0}{(1+r_ix)(1+r_jx)}\\
f_1^2 - f_0f_2 &= \sum_i d_{ii}^2 \frac{f_0^2}{(1+xr_i)^2} +
\sum_{i<j} (d_{ij}^2 + d_{ii}d_{jj}) \frac{f_0^2}{(1+xr_i)(1+xr_j)}\\
\frac{f_1^2 - f_0f_2}{f_0^2} &= \sum_i \frac{d_{ii}^2 }{(1+xr_i)^2} +
\sum_{i<j} \frac{(d_{ij}^2 + d_{ii}d_{jj})}{(1+xr_i)(1+xr_j)}
  \end{align*}
Now the diagonal entries of a positive definite matrix are positive,
so all the coefficients of the expansion of $f_1^2 - f_0f_2$ are
positive. Consequently, evaluation at a point in the closed first
quadrant yields a point in the open lower half plane, which
establishes the first part. 

For the second interlacing we have
\[
\frac{f_1^2 - f_0f_2}{f_0\,f_1} = \frac{f_0}{f_1}\left(\sum_i \frac{d_{ii}^2 }{(1+xr_i)^2} +
\sum_{i<j} \frac{(d_{ij}^2 + d_{ii}d_{jj})}{(1+xr_i)(1+xr_j)}\right)
\]

When we evaluate a  point in the closed first quadrant, $f_0/f_1$ is
in the first quadrant, and the part in the parentheses is in the lower
half plane. Consequently, their product lies in an open half plane $P$
(the rotation of the lower half plane by the argument of the point)
which misses the negative real axis.

Since $f_1^2 - f_0f_2$ is stable, all of its coefficients have the
same sign. If we let $x=0$ then this is positive, by Newton's
inequality for $f(0,y)$. Since the constant term is positive, all
terms are positive.

When  we introduce the factor of $x$, then the sum in parentheses for
$x(f_1^2-f_0f_2)/f_0^2$ lies in the right half plane by
lemma~\ref{lem:stable-fractions}. The case where the denominator is
$f_0f_1$ yields a point in the product of the first quadrant and right
half plane, and so again misses the negative real axis.

\end{proof}

If we choose $\beta>0$ then $f_0(\beta),f_1(\beta),f_2(\beta)$ are the
first three terms of a polynomial in $\allpolypos$, so we know that 
$f_1(\beta)^2 > 2 f_0(\beta)\,f_2(\beta)$. This suggests that we can
improve the result above.

\begin{cor}
  If $0<\alpha<2$ and $f,f_i$ are as above then 
\begin{enumerate}
\item $ f_1^2 - \alpha f_0f_2$ is stable.
\item $\smalltwodet{f_{k+1}}{f_{k+2}}{f_{k}}{f_{k+1}}$ is stable.
\end{enumerate}
\end{cor}
\begin{proof}
  For the second, we differentiate $k$ times and apply the first part.
  Using the notation of the lemma, the first part follows from the
  identity 
\begin{multline*}
\frac{f_1^2-\alpha f_0f_2}{f_0^2} =\\
\sum_i \frac{d_{ii}^2}{(1+x r_i)^2} + 
\alpha\sum_{i<j} \frac{d_{ij}^2}{(1+r_i x)(1+r_jx)} +
(2-\alpha)\sum_{i<j} \frac{d_{ii}d_{jj}}{(1+r_ix)(1+r_jx)}
\end{multline*}

\end{proof}

We can use the lemma to construct stable polynomials from polynomials
in $\allpolypos$.

\begin{cor}\label{cor:stable-deriv}
  If $f\in\allpolypos$ and $a,b,c>0$ then $af^2 + b ff' + cx(ff''-(f')^2)$ is stable.
\end{cor}
\begin{proof}
  From Lemma~\ref{lem:stable-ff} we know that each of the three terms
  positive interlace $f^2$, so any positive combination is stable
\end{proof}

The next lemma arises from considerations involving the Euler
polynomials. Unlike the last lemma, we must look carefully at the
terms of the quotient.

  \begin{lemma}\label{lem:euler-stable}
    If $f\in\allpolypos(n)$ then 
\[ Q=x\biggl(nf^2 + (x-1)^2\bigl[ f f' + x(ff'' -
f'^2)\bigr]\biggr)\in\stabled{1}
\]
  \end{lemma}
  \begin{proof}
    We show that  $Q$ and $f^2$ interlace positively. First, 
\[
Q/f^2 = x\biggl( n + (x-1)(f'/f) + x\bigl[ff''-f'^2\bigr]/f^2\biggr).
\]
If  $ f = \prod (x+r_i)$ where all $r_i$ are positive then
\begin{align*}
  \frac{f'}{f} &= \sum_{i=1}^n \frac{1}{x+r_i} & 
  \frac{f'^2-ff''}{f^2} &= \sum_{i=1}^n \frac{1}{(x+r_i)^2}
\end{align*}
Thus the quotient satisfies
\begin{align*}
  Q/f^2 &= x\biggl( n + (x-1)\sum \frac{1}{x+r_i} - x \sum \frac{1}{(x+r_i)^2}\biggr)
  \\
&= x \sum \biggl( 1 + \frac{x-1}{x+r_i} - \frac{x}{(x+r_i)^2}\biggr)\\
& = \sum \frac{x r_i(1+r_i+x)}{(x+r_i)^2}
\end{align*}
If we substitute $\alpha+\beta\imag$ for $x$ then
\begin{multline*}
\frac{x r_i(1+r_i+x)}{(x+r_i)^2} =\\
\frac{r_i \left(\beta ^4+\left(2 \alpha ^2+3 r_i \alpha +\alpha +r_i (r_i+2)\right)
\beta ^2+\alpha  (r_i+\alpha )^2 (r_i+\alpha +1)\right)}{|r_i+\alpha|^2}
\end{multline*}
Thus, the sum is in the first quadrant if $\alpha+\beta\imag$ is in
the first quadrant.
  \end{proof}

\section{Linear transformations}
\label{sec:line-transf-1}

Linear transformations preserve linearity, so the following is
immediate.

\begin{lemma}
  If $T\colon{}\stabled{1}\longrightarrow\stabled{1}$ preserves degree and the sign
  of the leading coefficient then $T$ preserves interlacing and
  positive interlacing.
\end{lemma}

\index{Hadamard product!for $\stabled{1}$}

\begin{lemma}
    If $f,g\in\stabled{1}$  then $f\ast g\in\stabled{1}$.
\end{lemma}
\begin{proof}
If $f\in\allpolypos$ then the method of generating functions shows
that the Hadamard product maps
$\allpolypos\times\stabled{1}\longrightarrow\stabled{1}$. However, we need a
different approach when $f\not\in\allpolypos$.
There does not appear to be a direct way to do this. We follow
\cite{garloff-wagner}: convert to $\rup{1}$, take Hadamard products
there, and use the fact that $f\lesslesseq g$ and $h\lesslesseq k$
implies $f\ast h\lesslesseq g\ast k$.

\end{proof}

\begin{lemma}\label{lem:h-rising}
  If $T\colon x^n\longrightarrow \rising{x}{n}$ then $T$ maps
  $\stabled{1}$ to $\stabled{1}$ and  $\stabledc{1}$ to $\stabledc{1}$.
\end{lemma}
\begin{proof}
  It suffices to show that $T$ maps $\stabledc{1}$ to itself; we proceed
  by induction. So assume that $f\in\stabledc{1}$ and
  $T(f)\in\stabledc{1}$. If $a\in\rhp$ then will show that $T(x+a)f
  \hposlace T(f)$ which implies that $T(x+a)f\in\stabledc{1}$. From
  \eqref{eqn:recur-5a} we know that
\[ 
T(x+a)f = (x+a)\,T(f) + T(xf').
\]
If $\sigma\in\overline{\quada}$ then
\[
\frac{T(x+a)f}{T(f)}(\sigma) = (\sigma + a) +
\frac{T(xf')}{T(f)}(\sigma)
\]
Now $T(xf')\hposlace T(f)$ so  $\frac{T(xf')}{T(f)}(\sigma)$ is in
$\complexes\setminus(0,\infty)$. Thus $\frac{T(x+a)f}{T(f)}(\sigma)$ is in
$\complexes\setminus(0,\infty)$  and the lemma is proved.
\end{proof}

\begin{lemma}\label{lem:h-falling}
  If $T\colon x^n\longrightarrow \falling{x}{n}$ then $T^{-1}$ maps
  $\stabled{1}$ to $\stabled{1}$ and  $\stabledc{1}$ to $\stabledc{1}$.
\end{lemma}
\begin{proof}
  We proceed as above. Let $S = T^{-1}$. From \eqref{eqn:recur-5}  we have
\[
S(x+a)f = (x+a)\,S(f) + x \bigl( S(f)\bigr)'
\]
To show $S(x+a)f \hposlace S(f)$ we choose
$\sigma\in\overline{\quada}$ and then
\begin{equation}
  \label{eqn:h-ff}
  \frac{S(x+a)f}{S(f)} (\sigma) = (\sigma+a) + \sigma\, \frac{
    \bigl(S(f)\bigr)'}{S(f)}(\sigma)
\end{equation}
Since $S(f)\in\stabledc{1}$ by induction and $S(f) \lessless
\bigl(S(f)\bigr)'$ we have 
$\bigl(S(f)\bigr)'/S(f)\,(\sigma)\in\rhp$ and therefor \eqref{eqn:h-ff}
$\in\rhp+\rhp\times\overline{\rhp}=\complexes\setminus(0,\infty)$.
\end{proof}

  \begin{lemma}\label{lem:h-even}
    If $T\colon x^n\mapsto
    \begin{cases}
      x^{n/2}& \text{$n$ even} \\0 & \text{$n$ odd}
    \end{cases}$ then $T:\stabled{1}\longrightarrow\allpolypos$.    
  \end{lemma}
  \begin{proof}
    If $f\in\stabled{1}$ has degree $n$ and we write $f(x) = f_e(x^2) +
    xf_o(x^2)$ then 
\[ \imag^n( f_e(-x^2) + (-\imag x) f_o(-x^2)) \in\polycpx
\]
Thus, $f_e$ and $f_o$ are in $\allpoly$, and so $T(f)=f_e$ since all
coefficients are positive.
  \end{proof}

\index{Charlier polynomials!and stable polynomials}
The Charlier polynomials are defined by
\begin{align*}
C_{-1}^\alpha& =0 \qquad C_0^\alpha=1 \\
  C_{n+1}^{\alpha}(x) &= (x-n-\alpha)C_{n}^{\alpha}(x) - \alpha n
  C_{n-1}^{\alpha}(x)
\end{align*}

  \begin{lemma}\label{lem:h-charlier}
    If $\alpha\in\rhp$ and $T\colon x^n\mapsto C_n^\alpha$ then 
$ T^{-1}$ maps  $\stabledc{1}$ to $\stabledc{1}$, and $\stabled{1}$ to $\stabled{1}$.
  \end{lemma}
  \begin{proof}
    We proceed as with the rising and falling factorials. We show by
    induction that for $a\in\rhp$ we have that $T^{-1}(x+a)f \hposlace
    T^{-1}(f)$. From the recurrence
    \begin{align*}
      T^{-1}(x+a)f &= (x+\alpha+a)T^{-1}f + (x+\alpha)
      \bigl(T^{-1}(f)\bigr)\\
\intertext{we choose $\sigma\in\rhp$ and consider the quotient}
\frac{T^{-1}(x+a)f}{T^{-1}(f)}(\sigma) &=
(\sigma+\alpha+a) + (\sigma+a)  \frac{\bigl(T^{-1}(f)\bigr)'}{T^{-1}(f)}(\sigma)
    \end{align*}
We see that the quotient lies in $\rhp + \rhp\times\rhp \subset
C\setminus(0,\infty)$ since $a,\alpha$ and $\sigma$ are all in
  $\rhp$. 
  \end{proof}

If we have two mutually interlacing sequences then the sum of products
in reverse order is in $\allpoly$ \seepage{lem:convolution}. If we use the
original order we get a stable polynomial.

\begin{lemma}\label{lem:mi-same-order}
  Suppose that $f_1,f_2, \cdots, f_n$ and $g_1, g_2, \cdots, g_n$ are
  two sequences of mutually interlacing polynomials in $\allpolypos$. Then
\[
f_1g_1+f_2g_2+ \cdots+ f_n g_n\ \in\ \stabled{1}
\]
\end{lemma}
\begin{proof}
  If $1\le k \le n$ then $f_1\greateqeq f_k$ and $g_1\greateqeq
  g_k$. If $\sigma$ is in the closed first quadrant then
  $\frac{f_k}{f_1}(\sigma)$ and $\frac{g_k}{g_1}(\sigma)$ are in the
  open fourth quadrant. Thus $\frac{f_kg_k}{f_1g_1}(\sigma)$ are in
  the open lower half plane. Addition shows that
\[
\frac{1}{f_1g_1} \sum_{k=1}^n f_kg_k \quad \text{is in the open lower half
plane}
\]
which finishes the proof.
\end{proof}

We can apply this lemma to integrate families of interlacing
interlacing polynomials. We  also view these results as linear
transformations on $\allpolysep$. 

\begin{lemma}
  If $f,g\in\allpolysep$ then $\displaystyle\int_0^1 f(x+t)g(x+t)\,dt$ is stable.
\end{lemma}
\begin{proof}
  The proof is the same as Lemma~\ref{lem:convolution-int}, except we
  use the preceding lemma instead of Lemma~\ref{lem:convolution}.
\end{proof}

\begin{cor}
  If $f\in\allpolysep$ then $\displaystyle\int_0^1 e^t\,f(x+t)\,dt$ is
  stable.
\end{cor}
\begin{proof}
  Apply the lemma to $e^x$ and $f$, and then factor out $e^x$ from the
  result.
\end{proof}

\begin{example}
  If we take $f\in\allpolysep$ then we see that $\int_0^1
  f(x+t)^2\,dt$ is stable. In particular, 
\[
\int _0^1 \bigl(\rising{x}{n}\bigr)^2 \,dt\quad \text{ is stable}.
\]
\end{example}

  \begin{lemma}\label{lem:h-integral}
    If $f\in\stabled{1}$ then $\displaystyle \int_0^1 f(x+t)\,dt\in\stabled{1}$.
  \end{lemma}
  \begin{proof}
    Suppose $deg(f)=n$, let $g(x) = \imag^n f(-\imag x)$, and $h(x) =
    \int_0^1 f(x+t)\,dt$. Since $g(x)\in\polycpx$ we have from
    Lemma~\ref{lem:i-integral} 
    \begin{align*}
      \int_0^1 g(x+\imag t)\,dt&\in\polycpx & \text{and therefore}\qquad &
      \int_0^1 \imag^nf (-\imag x+  t)\,dt \in\polycpx. 
    \end{align*}
    The last integral is $\imag^n h(-\imag x)$, so $h\in\stabled{1}$.
  \end{proof}

It is surprising that we do not need to make any other assumptions
about $f$ in order that the difference  is stable.

\index{difference operator}
  \begin{cor}\label{cor:stable-difference}
    If $f\in\stabled{1}$ then $\Delta(f) = f(x+1) - f(x) \in\stabled{1}$.
  \end{cor}
  \begin{proof}
    Since $f'\in\stabled{1}$ we have
\[
f(x+1)-f(x) = \int_0^1 f(x+t)'\,dt \in\stabled{1}.
\]
  \end{proof}

  \begin{example}
    Here's a simple example. We know that $x^n\in\stabled{1}$, and
    the roots of $\Delta x^n = (x+1)^n-x^n$ are 
\[ -\frac{1}{2} -\imag\, \frac{\sin \frac{2\pi k}{n}}{-2 +
  2\cos \frac{2\pi k}{n}}\qquad k=1,\dots,n
\]
They all lie in the open left half plane, so $\Delta x^n\in\stabled{1}$.
  \end{example}

  \section{The totally stable matrix conjecture}
  \label{sec:totally-stable-matr}

  Recall that if $f(x)=\sum a_ix^i$ is in $\allpolypos$ then the matrix

\begin{equation}\label{eqn:tp-stab}
\begin{pmatrix}
  a_0 & a_1 & a_2 & a_3 & a_4 & a_5 & a_6 & \hdots \\
0&  a_0 & a_1 & a_2 & a_3 & a_4 & a_5 &  \hdots \\
0&0&  a_0 & a_1 & a_2 & a_3 & a_4 &   \hdots \\
0&0&0&  a_0 & a_1 & a_2 & a_3 &    \hdots \\
\vdots&\vdots&&&   & \ddots & \ddots &  \ddots \\
\end{pmatrix}
\end{equation}
 is totally positive. This means that all minors are
 non-negative. Now, assume that
$f(x,y)=\sum f_i(x)y^i$ and form the matrix

\begin{equation}\label{eqn:tp-stab-2}
\varphi(f)(x,y)=
\begin{pmatrix}
  f_0 & f_1 & f_2 & f_3 & f_4 & f_5 & f_6 & \hdots \\
0&  f_0 & f_1 & f_2 & f_3 & f_4 & f_5 &  \hdots \\
0&0&  f_0 & f_1 & f_2 & f_3 & f_4 &   \hdots \\
0&0&0&  f_0 & f_1 & f_2 & f_3 &    \hdots \\
\vdots&\vdots&&&   & \ddots & \ddots &  \ddots \\
\end{pmatrix}.
\end{equation}

If $f\in\gsubplus_2$ and $\alpha>0$ we know that
$f(\alpha,y)\in\allpolypos$, so all minors of $\varphi(f)(\alpha,y)$
are positive. Now, whenever we have a polynomial that is positive for
$x>0$  we should consider whether it is stable. Empirical
evidence suggests the following:

\begin{conj}[The totally stable matrix conjecture]

  If $f(x,y)\in\gsubplus_2$ then all minors of the matrix
  \eqref{eqn:tp-stab-2} are stable.
\end{conj}

Here's a simple consequence:

\begin{lemma}
  If $f\in\gsubplus_2$ and all minors of $\varphi(f)$ are stable, then
  all minors have non-negative coefficients.
\end{lemma}
\begin{proof}
  If $g(x)$ is a minor then since $g$ is stable we know that all
  coefficients have the same sign. If we set $x=\alpha\ge0$ then the
  matrix $\varphi(f(\alpha,y))$ is just the matrix of a one variable
  polynomial in $\allpolypos$ and we know that all of its minors are
  non-negative. Thus, $g(\alpha)$ is non-negative for all positive
  $\alpha$.  Since all coefficients are the same sign they are
  positive.
\end{proof}

\begin{example}
For a simple example, consider $\varphi(1+x+y)^4 =$

\[
\begin{pmatrix}
(1+x)^4    & 4(1+x)^3  & 6(1+x)^2  & 4(1+x)  & 1  & 0    & \hdots \\
0&   (1+x)^4 & 4(1+x)^3   & 6(1+x)^2  & 4(1+x)  & 1  &  \hdots \\
0&0&    (1+x)^4 & 4(1+x)^3  & 6(1+x)^2  & 4(1+x)   &    \hdots \\
0&0&0&  (1+x)^4  & 4(1+x)^3  & 6(1+x)^2  &     \hdots \\
\vdots&\vdots&&   & \ddots & \ddots &  \ddots \\
\end{pmatrix}.
\]
Every minor is a power of $(1+x)$ times a minor of the coefficients,
which we know to be positive, so all minors are stable.
\end{example}

\begin{example}
  Lemma~\ref{lem:stable-ff} shows that the $2$ by $2$ determinants
  $\smalltwodet{f_k}{f_{k+1}}{f_{k+1}}{f_{k+2}}$ are stable.
\end{example}

\begin{example}
  If $f(x)\in\allpolypos$ then $f(x+y)\in\gsubplus_2$, and the
 entries of the matrix come from the Taylor series.
\[
\varphi(f(x+a))=
\begin{pmatrix}
  f & f^{(1)} & f^{(2)}/2! & f^{(3)}/3! & f^{(4)}/4! & f^{(5)}/5! & f^{(6)}/6! & \hdots \\
0&  f & f^{(1)}/1! & f^{(2)}/2! & f^{(3)}/3! & f^{(4)}/4! & f^{(5)}/5! &  \hdots \\
0&0&  f & f^{(1)}/1! & f^{(2)}/2! & f^{(3)}/3! & f^{(4)}/4! &   \hdots \\
0&0&0&  f & f^{(1)}/1! & f^{(2)}/2! & f^{(3)}/3! &    \hdots \\
\vdots&\vdots&&&   & \ddots & \ddots &  \ddots \\
\end{pmatrix}.
\]
We know that some of the two by two determinants are stable. For example,
\[
\begin{vmatrix}
 f^{(3)}/3! & f^{(4)}/4!\\ 
 f^{(4)}/4! & f^{(5)}/5! 
\end{vmatrix}=
\frac{1}{4!4!} \bigl(\frac{4}{5}gg'' - g'g'\bigr)
\]
where $g=f^{(3)}$. Since $4/5<1$ Corollary~\ref{cor:stable-sq-2}
applies to show that it is stable.

\end{example}

\begin{remark}\label{rem:hurwitz-mat}
  If $f = \sum a_i x^i$ is a stable polynomial  then
  $f_e\lesslesseq f_o$. 
  Proposition~\ref{prop:hurwitz-totally-pos} therefore implies that
\[
\begin{pmatrix}
  a_0 & a_2 & a_4 & \cdots \\
0 & a_1 & a_3 & a_5 & \cdots \\
0 & 0 & a_0 & a_2 & a_4 & \cdots\\
0 & 0 & 0 & a_1 & a_3 & a_5 & \cdots \\
&&&\ddots &\ddots & \ddots

\end{pmatrix}
\]
is totally positive.  This matrix is known as the Hurwitz matrix
\cite{asner}.
\index{Hurwitz matrix}
\end{remark}

  Earlier \seepage{cor:tpc} we saw that the matrix of coefficients of a
  polynomial in $\gsubplus_2$ is totally positive$_2$, but there are
  examples where some three by three determinants are
  negative. However, there is a surprising conjecture that is not
  fazed by this fact.

  \begin{conj}
    If $D_1,D_2,D_3$ are positive definite matrices and $M$ is the
    matrix of coefficients of $y,z$ of $ |I + x\,D_1+y\,D_2+z\,D_3|$
   then all minors of $M$ are stable polynomials in $x$. 
  \end{conj}

  For example, if we use the same example with a particular choice for
  $D_1$

\[
      M =
      \begin{pmatrix}
        1&0&0\\0&1&0\\0&0&1
      \end{pmatrix}
+
     x\, \begin{pmatrix}
        0&0&0\\0&1&0\\0&0&0
      \end{pmatrix}
+
y\,\begin{pmatrix}
  13& 9& 7  \\9 & 7& 5 \\7& 5& 4
\end{pmatrix}
+
z\,
\begin{pmatrix}
  5& 7& 8 \\7& 11& 12\\8& 12& 14
\end{pmatrix}
\]
then the coefficient matrix is a $4$ by $4$ matrix, and the
non-constant $3$ by $3$ minors are
\[
\left(
\begin{array}{lll}
 -1620 x^3-7560 x^2-4884 x+1760 & -1080 x^2-3512 x-4544 & -768 x-4096  \\
 -1080 x^2-3496 x-3376 & -720 x-1312 &   \\
 -744 x-2728 &  &   
\end{array}
\right)
\]
All entries are stable polynomials. Note that if we substitute zero
for $x$ we get the example on page~\pageref{cor:tpc}.

  \section{Determinants of Hankel matrices }
  \label{sec:determ-hank-matr}

If $p_1,p_2,\dots$ is a sequence of polynomials then the $d$ by $d$  Hankel
matrix \seepage{sec:karlin} of this sequence is

\[
H(\{p_i\};d) =\begin{pmatrix}
  p_1 & p_2 & p_3 & \hdots & p_d \\
  p_2 & p_3 & p_4 & \hdots & p_{d+1} \\
  p_3 & p_4 & p_5 & \hdots & p_{d+2} \\
  \vdots   & \vdots   & \vdots   & \ddots &  \vdots \\
  p_d & p_{d+1} & p_{d+2} & \hdots & p_{2d}
\end{pmatrix}
\]
In this section we study properties of Hankel matrices. This includes
the kind of polynomial (all real roots, stable, etc.),  the sign of
the coefficients, and interlacing properties.

A polynomial sequence $\{p_n(x)\}_{n\ge0}$ is called
\emph{$q$-log-concave} if all coefficients of 
\index{q-log-concave}
\[
 p_n\,p_{n+2} - p_{n+1}^2
=
\begin{vmatrix}
  p_n & p_{n+1} \\ p_{n+1} & p_{n+2}
\end{vmatrix}
= H[\{p_n,p_{n+1},p_{n+2}\};2]
\]
are positive for $n\ge0$. If all coefficients are negative it is
\emph{$q$-log-convex} \cite{liu-wang-log}.  \index{q-log-convex} A
stronger condition is that the polynomial $p_{n+1}^2- p_n\,p_{n+2}$
is stable. For instance,
Lemma~\ref{lem:stable-ff} shows that the coefficients of a polynomial
in $\gsubplus_2$ are $q$-log-convex.

\begin{lemma}\index{Bell polynomials}
  If $B_n$ is the Bell polynomial then $B_{n+2}B_{n} - B_{n+1}^2$ is
  weakly stable with positive leading coefficient. In particular, it
  is $q$-log-convex.
\end{lemma}

\begin{proof}
  The Bell polynomials satisfy the recurrence
$ B_{n+1} = x(B_n + B_n')$, and therefore
\[
B_{n+2}B_{n} - B_{n+1}^2 = x\left( B_n^2 + B_nB_n' +
  x(B_nB_n''-B_n^{'2})\right)
\] Since $B_n\in\allpolyposclose$ it follows from
Corollary~\ref{cor:stable-deriv} that this is in $\stabled{1}$.
\end{proof}

In addition, we have shown that $B_n^2 \hposlace
B_nB_{n+2}-B_{n+1}^2$. 

  \begin{lemma}
    If $A_n$ is the Euler polynomial then $A_{n+2}A_n-A_{n+1}^2$ is
    weakly stable with positive  coefficients. 
  \end{lemma}
  \begin{proof}
    The Euler polynomial satisfies the recurrence
\[
A_{n+1} = x((n+1)A_{n} + (1-x)A_{n}') 
\]
and so substitution yields
\begin{align*}
A_{n+2}A_n-A_{n+1}^2 &=  
x\biggl( (1+n)A_n^2 + (x-1)^2 \bigl[ A_nA_n' + x(A_nA_n'' - (A_n')^2\bigr]\biggr)
\end{align*}
This sum is the same as the sum in Lemma~\ref{lem:euler-stable},
except that it has $n+1$ in place of $n$. This only adds a point in
the first quadrant to the quotient, so we conclude that
$A_{n+2}A_n-A_{n+1}^2$ is stable. The value at $0$ is known to be
positive \cite{liu-wang-log}, so all terms are positive. Since $A_n$
has 0 for a root, we only get weakly stable.
\end{proof}

The following general result is easy, since everything factors. The
work is finding the constant.

\begin{lemma}\label{lem:hankel}
  If $n$ and $d$ are positive then  $|H[\{\rising{x}{n},\rising{x}{n+1},\dots\};d]|$ is in
  $\allpoly$. 
\end{lemma}

\begin{proof}

We first prove the result for $d=3$. 
  If $n=1$ and $d=3$ then the matrix is
\[
\left(
\begin{array}{lll}
 x  & x (x+1)  & x (x+1) (x+2)  \\
 x (x+1)  & x (x+1) (x+2)  & x (x+1) (x+2) (x+3)  \\
 x (x+1) (x+2)  & x (x+1) (x+2) (x+3)  & x (x+1) (x+2) (x+3) (x+4) 
\end{array}
\right)
\]
with determinant 
\[
2 x^3 (x+1)^2 (x+2)
\]
If we let $f_n(x)$ be the determinant of
$H[\{\rising{x}{n},\rising{x}{n+1},\dots\};3]$ then it is easy to see
that 
\begin{align*}
f_n(x) & = \rising{x}{n} \, \rising{x}{n+1}\, \rising{x}{n+2}\,\times \\
& \begin{vmatrix}
  1 & 1 & 1 \\ x+n & x+n+1 & x+n+2 \\
(x+n)(x+n+1) & (x+n+1)(x+n+2) & (x+n+2)(x+n+3) \\
\end{vmatrix} \\
&  = 2 \rising{x}{n} \, \rising{x}{n+1}\, \rising{x}{n+2}
\end{align*}

In general, the answer is
\[
(d-1)!(d-2)!\cdots 1!\quad
\rising{x}{n}\rising{x}{n+1}\cdots\rising{x}{n+d-1}
\]
It's easy to find the polynomial factors. If we remove them and
replace $x+n$ by $y$ we are left with
\[
\begin{vmatrix}
  1 & \dots & 1 \\
  \rising{y}{1} & & \rising{y+d-1}{1} \\
\vdots && \vdots \\
\rising{y}{d-1} & \dots & \rising{y+d-1}{d-1}
\end{vmatrix}
 \]
From \cite{kratt} we see that this equals the Vandermonde determinant
\[
\begin{vmatrix}
  1 & \dots & 1 \\
y & & y+d-1 \\
\vdots & & \vdots \\
y^d & \dots & (y+d-1)^{d-1}
\end{vmatrix}
\]
which equals
\[
\prod_{0<i<j<d} \bigl([y+j] - [y+i]\bigr) =
\prod_{0<i<j<d}(j-i)
\]
and this is precisely $(d-1)!\cdots1!$.
\end{proof}

\begin{example}
  The Chebyshev T polynomials have trivial Hankel determinants. First
  of all,
\[
\begin{vmatrix}
  T_{n} & T_{n+1}\\
  T_{n+1} & T_{n+2}
\end{vmatrix} =
\begin{vmatrix}
  T_{n} & T_{n+1}\\
  2xT_{n}-T_{n-1} & 2xT_{n+1}-T_{n}
\end{vmatrix} =
\begin{vmatrix}
  T_{n-1} & T_{n}\\
  T_{n} & T_{n+1}
\end{vmatrix}
\]
Continuing, we find they are all equal to $x^2-1$. Thus, all the
Hankel matrices $H[\{T_n\};2]$ are equal to $x^2-1$, so they all have
all real roots. 

Since $T_{n+1} = 2xT_{n} -T_{n-1}$ has coefficients that do not depend
on $n$, all the higher Hankel determinants are zero.

\end{example}

We end this section with some empirical results about the Hankel determinants of
the form $H(\{p_n\};d]$. ``cpx'' means roots in all quadrants, and
\checkmark means it has been proved.

\index{polynomials!Euler}\index{Euler polynomials}
\index{polynomials!Bell}\index{Bell polynomials}
\index{polynomials!Narayana}\index{Narayana polynomials}
\index{polynomials!Laguerre}\index{Laguerre polynomials}

\begin{tabular}{llcc}
  \toprule
  Family & $d$ & Kind of polynomial & Sign of coefficients\\ \midrule
  Bell &$2$ & stable \checkmark & positive \checkmark \\
  & $\ge3$ &  cpx & positive \\
  \midrule
  Euler & $2$ & stable \checkmark & positive \checkmark \\
  & $\ge3$ & cpx & positive \\
  \midrule
  Narayana & $2$ & stable & positive \\
  & $3$ & stable & positive \\
  & $\ge3$ & cpx & positive \\
  \midrule
  Laguerre & $2$ & stable & negative \\
   & $3$ & stable & negative \\
   & $4$ & stable & positive \\
   & $n$ & stable & $(-1)^{\binom{n}{2}}$\\
   \bottomrule
\end{tabular}

\section{Constructing stable polynomials}
\label{sec:stable-two-variables}

We  construct matrices whose determinant is a stable
polynomial. Recall
that a matrix $A$ is called \emph{skew-symmetric}
\index{skew-symmetric} if $A^T = -A$. If $A$ is skew-symmetric then
\begin{enumerate}
\item The eigenvalues of $A$ are purely imaginary.
\item If $D$ is positive definite, then the eigenvalues of $A-D$ have
  negative real part.
\end{enumerate}

Using the latter property, and the fact that $f(x)\in\stabled{1}$ if and only if
$\imag^nf(-\imag x)\in\polycpx(n)$ we have 
 we have
\begin{lemma}
  If $D_1,\dots$ are positive definite, $S$ is symmetric, and $A$ is
  skew-symmetric then
  \begin{enumerate}
\item $    \bigl| I + x_1\,D_1 + \cdots + x_d\,D_d + \imag S\bigr |
  \in\stabledc{1} $\\
\item $    \bigl| I + x\,D_1 + \cdots + x_d\,D_d +  A\bigr |  \in\stabled{d} $\\
  \end{enumerate}
\end{lemma}

For example, the determinant of 
\begin{gather*}
\begin{pmatrix} 1&0\\0&1\end{pmatrix} + 
x \begin{pmatrix} 1 &0 \\0 & 2\end{pmatrix}
+y\begin{pmatrix} 3 & 1\\1 & 2\end{pmatrix}
+\begin{pmatrix} 0&1\\-1&0\end{pmatrix}\\
\intertext{is in $\stabled{2}$ and equals}
2 x^2+8 y x+3 x+5 y^2+5 y+2.
\end{gather*}

\begin{remark}
  It's easy to construct a determinant representation $\alpha|I+xD+A|$
  of $f\in\stabled{1}$, where $\alpha>0$. Factoring $f$, it suffices to
  show this for polynomials of degree one and two - the general case
  follows by taking direct sums.
  \begin{align*}
    x+a &= a\left| (1) + x\,(1/a) + (0)\right|\\
    (x+a)^2+b^2 &= a^2\left|
      \begin{pmatrix}1&0\\0&1 \end{pmatrix} +
      {x}\begin{pmatrix}1/a&0\\0&1/a \end{pmatrix} +
      \begin{pmatrix} 0 & {b}/{a} \\ -{b}/{a} & 0      \end{pmatrix}
      \right|
    \end{align*}
  \end{remark}

\index{Jacobi matrix}
\begin{example}
A special case of the matrix $ xI + D + A$ is the Jacobi-like matrix
\[  
\left(
\begin{array}{lllllll}
 x+a_1 & b_1 & 0 & 0 & 0 & 0 & \dots\\
 -b_1 & x+a_1 & s_1 & 0 & 0 & 0 & \dots \\
 0 & -s_1 & x+a_2 & b_2 & 0 & 0 &\dots\\
 0 & 0 & -b_2 & x+a_2 & s_2 & 0 &\dots\\
 0 & 0 & 0 & -s_2 & x+a_3 & b_3 \\
 0 & 0 & 0 & 0 & -b_3 & x+a_3 \\
\vdots & \vdots &\vdots &\vdots&& \ddots & \ddots
\end{array}
\right)
\]
whose determinant is stable. 
The determinant of the submatrix consisting of
the first $n$ rows and columns is an $n$'th degree polynomial $p_n$.
The recurrence relations for $p_n$ are
\begin{align*}
  p_{2k} &= (x+a_k)\,p_{2k-1} + b_k^2\,p_{2k-2} \\
  p_{2k+1} &= (x+a_{k+1})\,p_{2k} + s_k^2\,p_{2k-1} \\
\end{align*}
If all the $a_i$ are positive then we know from the Lemma above that
all $p_k$ are stable ; this also follows from the recurrence of
Corollary \ref{cor:orthog-h}.
Moreover, the corollary shows that they interlace.

\end{example}

If a polynomial $f(x)$ is in $\allpolypos$ then all its roots are
negative, so $f(x^2)$ has  purely imaginary roots. There is a similar
result for $\gsubplus_2$.

\begin{lemma}\label{lem:stable-fxx}
If $f(x,y)\in\gsubplus_2$ then $f(x,x^2)\in\stabled{1}$. 
\end{lemma}
\begin{proof}
  Since $f$ has real coefficients it suffices to show that
  $f(\sigma,\sigma^2)\ne0$ where $\sigma$ is in the first quadrant. If
  so, then $\sigma^2$ is in the upper half plane, and hence
  $f(\sigma,\sigma^2)\ne0$ since $f\in\rup{2}$.

  Alternatively, this is a consequence of the next lemma.
\end{proof}

  \begin{lemma}
If $D_1$ and $D_2$ are positive definite then $|I + x D_1+x^2D_2|$ is
stable. 
\end{lemma}
\begin{proof}
If $\sigma=a+b\imag$ is in the first quadrant then
\[
\bigl|I + \sigma D_1+\sigma^2D_2\bigr| =
\bigl|(I + a D_1 +(a^2-b^2)D_2) + \imag(bD_1 + 2abD_2)\bigr|
\]
The coefficient of $\imag$ is positive definite, so the determinant is non-zero.


\end{proof}

  \section{Stable matrix polynomials}
  \label{sec:stable-matr-polyn}

\index{stable matrix polynomials}
\index{matrix polynomials!stable}

  A hyperbolic matrix polynomial $M$ is determined by the properties of the
  family of polynomials $v^tMv$. If we posit that all the quadratic
  forms are stable we get \smp{s}.

  \begin{definition}
    An $n$ by $n$  matrix polynomial $M$ is \emph{stable} if $v^tMv$
    is stable for all non-zero $v\in\reals^n$, and the leading
    coefficient is positive definite.
  \end{definition}

  Most properties of hyperbolic matrix polynomials carry over to
  \smp{s}. For instance, we define

\[
\begin{array}{cccccccc}
  M_1 &  \hlace& M_2 & iff & v^t\,M_1\,v &\hlace& v^tM_2v &
  \forall v\in\reals^n\setminus0 \\[.2cm]
  M_1 &  \hposlace & M_2 & iff & v^t\,M_1\,v &\hposlace& v^tM_2v &
  \forall v\in\reals^n\setminus0 \\
\end{array}
\]

The following lemma is easily proved using properties of stable
polynomials.

\begin{lemma} Suppose that $f,g,h$ are \smp{s}, and
  $k\in\stabled{1}$. 
  \begin{enumerate}
  \item All elements of $\hyperpos{1}$ are \smp{s}.
  \item $f\hlace g$ iff $f + \sigma
     g$ is a \smp\ for all $\sigma\in\rhp$.
        \item If $f\hlace g$ and $f\hlace h$ then $f\hlace
    g+h$. In particular, $g+h$ is a \smp.
  \item If $f\hlace g\hlace h$ then $f+h\hlace g$.
  \item $fk$ is a \smp.
  \item If all $a_i$ are positive  then
$
(a_1x_1+\cdots a_dx_d+b)f \hlace f
.$
\item If the $i$-th diagonal element of $f$ is $f_i$,
  and of $g$ is $g_i$ then $f\greateqeq g$ implies $f_i\greateqeq g_i$.
   \item $\frac{ df}{dx}$ is a \smp.
  \item $f \hlace \frac{ df}{dx}$
  \item If $f\hlace g$ then $\frac{df}{dx}
    \hlace \frac{dg}{dx}$
  \end{enumerate}

\end{lemma}

\index{positive interlacing}
\index{interlacing!positive}

We can construct stable matrix polynomials using positive interlacing.
We start with the data
\begin{align*}
  \text{positive real numbers} & \qquad a_1,b_1,\dots,a_n,b_n\\
\text{positive definite matrices} &\qquad  m_1,\cdots,m_{n}
\end{align*}
We claim that 
\[
 M(x)= \sum_{k=1}^{n} \, m_{k} \, 
\prod_{i\ne k}(x+a_i)^2+b_i^2
\]
 is a stable matrix polynomial. In order to verify this  we first define
\[
g(x) = \prod_{i}(x+a_i)^2+b_i^2
\]
and  we know that 
\[
g(x) \hlace  \frac{g(x)}{(x+a_k)^2+b_k^2} = \prod_{i\ne k}(x+a_i)^2+b_i^2.
\]
Next,  
\[
v^t\,M\,v = \sum_k (v^t\,m_k\,v) \frac{g(x)}{(x+a_k)^2+b_k^2}
\]
Since $v^tm_kv\,$ is positive, it follows
that $g\hlace v^t\,M\,v\,$ since each of the terms positively
interlaces $g(x)$, and positive interlacing is closed under positive
linear combinations.

More generally, we can also construct stable matrix polynomials using
common interlacing.

\begin{lemma}
  If $f_1,\dots,f_r$ have a common interlacing in $\stabled{1}$ and
  $m_1,\dots,m_r$ are positive definite matrices then
\[ f_1\,m_1 + \cdots + f_r\,m_r\ \text{is a stable matrix
  polynomial}.\]
\end{lemma}
\begin{proof}
  The leading coefficient is a positive linear combination of positive
  definite matrices and so is positive definite. If $g\hlace f_i$
  then 
\[v^t\,g\,v =(v^tv) g \hlace (v^t\,m_i\,v)f_i = v^t\, f_im_i\,v\]
 and therefore
  \[ v^t\,g\,v \hlace v^t\,(f_1m_1+\cdots+f_rm_r)\,v.\]
\end{proof}

  A quadratic polynomial is stable if all its coefficients are
  positive. The analogous result holds for matrix polynomials.

  \begin{lemma}
    If $M_1,M_2,M_3$ are positive definite matrices then the matrix
    polynomial $M_1+x\,M_2 + x^2\,M_3$ is a stable matrix polynomial,
    and its determinant is stable.
  \end{lemma}
  \begin{proof}
    We must verify that $v^t(M_1+x\,M_2 + x^2\,M_3)v$ is stable for
    all non-zero $v\in\reals^\mu$. Now all the coefficients of $(v^tM_1v) + x
    (v^tM_2v) + x^2(v^tM_3v)$  are positive since the $M_i$
    are positive definite, and , so this is a stable
    polynomial.

    If $\lambda$ is a root of $\bigl|M_1+x\,M_2 + x^2\,M_3\bigr|$ then
    $\lambda\ne0$ since $M_1$ has  positive determinant. There is a
    non-zero vector $v$ in $\complexes^\mu$ such that
    \begin{gather*}
      (M_1+\lambda\,M_2 + \lambda^2\,M_3)v=0 \\
\intertext{Recall that if $v^\ast$ is the conjugate transpose, then $v^\ast\,M\,v>0$ for
    non-zero $v$ and positive definite $M$. 
Multiplying on the left by $v^\ast$ yields}
      (v^\ast\, M_1\,v) + \lambda (v^\ast\, M_2\,v) + \lambda^2(v^\ast\, M_3\,v)=0
    \end{gather*}
    Since this is a quadratic with positive coefficients, it is a
    stable polynomial, and so $\lambda$ lies in the left half plane.
  \end{proof}

Surprisingly, powers such as $(xI+A)^d$ are generally not stable matrix
polynomials. We do have

\begin{lemma}\label{lem:stable-power}
  If $A$ is positive definite, and $B$ is positive semi-definite then
  $(xI+A)^2+B$ is a stable matrix polynomial, and generally not hyperbolic.
\end{lemma}
\begin{proof}
  First take $B=0$. If $v\in\reals^\mu$ is non-zero then all coefficients of 
\[
x^2 (v^t\,v) + x( 2 v^t\,A\,v) + (v^t\,A^2\,v)
\]
are positive, so $(xI+A)^2$ is a stable matrix polynomial.
Next, we show that  the
discriminant $4((v^t\,A\,v)^2-(v^t\,v)(v^t\,A^2\,v))$ is never
positive. Define an inner
product $<v,w> = v^t\,w$. Then
\[
(v^t\,A\,v)^2-(v^t\,v)(v^t\,A^2\,v) = <v,Av>^2 - <v,v><Av,Av>\ \le 0
\]
by the \index{Cauchy-Schwartz inequality}Cauchy-Schwartz inequality. 
If $B$ is positive semi-definite then the discriminant is
\[
4((v^t\,A\,v)^2-(v^t\,v)(v^t\,A^2\,v)) - 4(v^t\,v)(v^t\,B\,v)
\]
which is still non-positive.
\end{proof}

There is a similar result for cubics, but we need a condition on the
eigenvalues to get stability.

\begin{lemma}
  Suppose $A$ is a positive definite matrix with minimum and maximum
  eigenvalues $\lambda_{min}$ and $\lambda_{max}$. If
  $\frac{\lambda_{max}}{\lambda_{min}}\le 9$ then $(xI+A)^3$ is a
  stable matrix polynomial.
\end{lemma}
\begin{proof}
  We need to show that if $v$ is any non-zero vector then 
\[
x^3(v^t\,v) + x^2(3v^t\,A\,v)+ x(3v^t\,A^2\,v)+(v^t\,A^3\,v)
\]
is stable. Since all the coefficients are positive, by
Lemma~\ref{lem:stable-reasons} this is the case precisely when
\[
\frac{(3v^t\,A\,v)(3v^t\,A^2\,v)}{(v^t\,v)(v^t\,A^3\,v)}\ge1
\]
Now this is equivalent to
\[
\frac{(v^t\,A\,v)}{(v^t\,v)}\times\frac{(v^t\,A^2\,v)}{(v^t\,A^3\,v)}\ge
\frac{1}{9}
\]
The first factor lies in $[\lambda_{min},\lambda_{max}]$, and
the second lies in $[1/\lambda_{max},1/\lambda_{min}]$, so
their product lies in 
$\left[\frac{\lambda_{min}}{\lambda_{max}},\frac{\lambda_{max}}{\lambda_{min}}\right]$.
By hypothesis $\lambda_{min}/\lambda_{max}\ge1/9$, so the proof is complete.
\end{proof}

The ratio $\lambda_{max}/\lambda_{min}$ is known as the
\index{condition number}\emph{condition number} of $A$, so we could restate
the hypothesis to say that the condition number is at most $9$.

\chapter{Transformations in the complex plane} 

\label{cha:complex}

\renewcommand{\TimeStampStart}{Friday, January 11, 2008: 12:02:21}

In this chapter we are interested in polynomials whose roots might be
complex. If $\cals$ is a region in the complex plane, then
$\allpolyint{\cals}$ means all polynomials whose roots all lie in
$\cals$. We will give some examples of linear transformations $T$ and
regions $\diffi,\diffj$ such that
$T\colon{}\allpolyint{\diffi}\longrightarrow\allpolyint{\diffj}$.

It is difficult to describe the behavior of a transformation on all
polynomials, so we sometimes restrict ourselves to a special subset of
polynomials.  Define the function
$\phi:\complexes\mapsto\allpolyint{\complexes}(n)$ by $z\mapsto
(x-z)^n$. Given a transformation $T$ we are interested in determining
the image  $T\,\phi$. In the best case we can find a a simple
function $\phi'$ so that the diagram commutes:

\begin{equation}\label{eqn:diagram-phi}
  {\xymatrix{
{\complexes}
      \ar@{.>}[drrr]_{\phi' }           
      \ar@{->}[rrr]^{{\phi:\,z\mapsto(x-z)^n} }         
      &&&
\allpolyint{\complexes}(n)
      \ar@{->}[d]^{T } \\        
      &&&
\allpolyint{\complexes}(n)
}}
\end{equation}

Given a set $\mathcal{P}$ of polynomials we are also interested in
all the roots of polynomials in $\mathcal{P}$. We extend the meaning
of $\roots(f)$ to a set of polynomials. Define
\[
\roots(\mathcal{P}) = \{\zeta\in\complexes\,\mid\,\exists
f\in\mathcal{P}\wedge f(\zeta)=0\}
\]
If $\cals$ is a region in
$\complexes$, and $T$ is a transformation, then since $\phi
\cals \subset \allpolyint{\cals}$ we know that 
\begin{equation}\label{eqn:lambda}
 \roots( \,T (\phi
\cals)\,) \subseteq \roots(\, T(\allpolyint{\cals})\,).
\end{equation}
We are interested in situations when \eqref{eqn:lambda} is an
equality. See Question~\ref{ques:cpx-equal}.
   We generalize the Hermite-Biehler theorem that characterizes
  polynomials whose roots are in the lower half plane in terms of
  interlacing properties of their real and imaginary parts.

  The closed disk of radius $r$ centered at $\sigma$ is
  $\Delta(\sigma,r)$; the unit disk is $\Delta$, and the boundary of
  the unit disk is the unit circle $\partial\Delta$.

\section{The derivative and interlacing}
\label{sec:cpx-derivative}

We recall the Gauss-Lucas theorem and prove a converse.

The fundamental result about the derivative is the  the
\index{Gauss-Lucas theorem}Gauss-Lucas 
theorem which states that the roots of the derivative of $f$ lie in
the convex hull of the roots of $f$. It follows that

\begin{theorem}
  If $\cals$ is a convex subset of $\complexes$ then the linear transformation $f\mapsto
  f'$ maps $\allpolyint{\cals}$ to itself.
\end{theorem}

The converse is also true. 

\begin{prop}\label{prop-lucas-gauss}
  If $T$ is a linear transformation with the property that
  $T\colon{}\allpolyint{\cals}\longrightarrow\allpolyint{\cals}$ for all
  convex regions $\cals$, then $Tf$ is a multiple of some derivative.
\end{prop}
\begin{proof}
  We first choose $\cals$ to be the single point $0$. It follows that
  $T(x^n) = a_n x^{e_n}$, since $T(x^n)$ can only have $0$ as a root.
  Ignoring the trivial case, there is some non-zero $a_k$.  Let $s$ be
  the smallest integer so that $a_s\ne 0$. Next, consider
  $\cals=\{-1\}$. Since $(x+1)^s\in\allpolyint{\cals}$ and $T(x+1)^s =
  a_s x^{e_s}$, it follows that $e_s=0$.

Consider $\cals=\{-\lambda\}$. The only members of
  $\allpolyint{\cals}$ are multiples of  powers of $(x+\lambda)$. We compute
  $T(x+\lambda)^r$ in two ways:
  \begin{alignat*}{2}
    T(x+\lambda)^r &= \beta_r (x+\lambda)^{m_r } \quad & =\quad \beta_r
    \sum_{j=0}^{m_r} \binom{m_r}{j} x^j \lambda^{m_r-j} \\
    &=  \sum_{k=0}^{r} \binom{r}{k}T( x^k)\lambda^{r-k} \quad&=\quad
    \sum_{k=0}^{r} \binom{r}{k} a_k x^{e_k}\lambda^{r-k} 
  \end{alignat*}

Since these are  identities in $\lambda$ that hold for infinitely
many values of $\lambda$, the coefficients of powers of $\lambda$ must
be equal.  From consideration of the  coefficient of $\lambda^{r-k}$
we find

\begin{equation}\label{eqn:cpx-3}
 \binom{r}{k}a_k x^{e_k} = \beta_r \binom{m_r}{m_r-r+k} x^{m_r-r+k}
\end{equation}

Choosing $k=s$ shows that $\beta_r\ne 0$ for $r\ge s$. Consequently,
$a_k\ne 0$ for $k\ge s$. Substituting $r=k$ shows that $\beta_r=a_r$.
If we take $r\ge s$ then comparing exponents shows that $e_k =
m_r-r+k$. It follows from this that $e_k = k-c$ for some constant $c$,
and $k\ge s$.  Since we know that $e_s=0$, it follows that  $T(x^k) =
a_k x^{k-s}$ for $k\ge s$ and is $0$ for $k < s$. In addition, $m_r =
r-s$. Substituting this into \eqref{eqn:cpx-3} yields

 $$a_r = \begin{cases} \frac{a_s}{s!} \frac{r!}{(r-s)!} & r\ge
  s\\
0 & r< s
\end{cases}$$
\noindent%
which  implies that $T(x^r) = a_s \diffd^s x^r$, as promised.

\end{proof}

 Every point of $\cals$ is a root of a derivative of a polynomial
supported on the boundary.

  \begin{lemma}
    If $\cals$ is a convex region of $\complexes$, and $\diffd$ is the
    derivative, then
$$ \roots(\, \diffd(\, \allpolyint{\partial\cals}\,)\,) = \cals.$$
  \end{lemma}
  \begin{proof}
If $\sigma,\tau\in\partial\cals$ then $f(x) =
(x-\sigma)(x-\tau)\in\allpolyint{\partial\cals}$. Since $f'
=2x-(\sigma+\tau)$, we see that
$\frac{1}{2}(\sigma+\tau)\in\roots(\diffd\allpolyint{\partial\cals})$.
Since $\cals$ is convex, the set of all midpoints of pairs of vertices on the
boundary equals all of $\cals$.

  \end{proof}

\section{Properties of algebraic Interlacing}  

In this section we consider some questions about the location of the
roots of algebraic interlacing polynomials.

  \begin{definition}
    Suppose $f(z)=\sigma (z-r_1)\cdots(z-r_n)$. We say that $f\cpxlace
    g$ if and only if we can write
$$ g(z) = f(z)\left(\frac{\rho_1}{x-r_1} + \cdots +
  \frac{\rho_1}{x-r_1}\right)$$
where $\rho_1,\dots,\rho_n$ are all real and non-negative.
  \end{definition}

It is immediate from this definition that if $f,g\in\allpoly$ then $f\lesslesseq g$ 
iff $f\cpxlace g$. For any $f$ we have that $f\cpxlace f'$, since we
know
$$ f'(x) = f(x)\,\left( \frac{1}{x-r_1}+\cdots+\frac{1}{x-r_n}\right).$$
The Gauss-Lucas property holds for this definition of interlacing, and
the proof is nearly the same, see \cite{marden}*{Theorem 6.1'}.

\begin{lemma}\label{lem:gauss-lucas-1}
Suppose that $\cals$ is a convex region, and
$f\in\allpolyint{\partial \cals }$.
 If $f\cpxlace g$ then
$g\in\allpolyint{\cals}$. 
\end{lemma}

We now extend the Gauss-Lucas Theorem by showing that every possible
point is a root of some interlacing polynomial.  Let
$Hull(\zeta_1,\dots,\zeta_n)$ denote the convex hull of the points
$(\zeta_1,\dots,\zeta_n)$. If these points are the roots of $f$, then
we write $Hull(f) = Hull(\zeta_1,\dots,\zeta_n)$.

\begin{prop}
  If $f$ is any polynomial, then every point in $Hull(f)$ is a root of
  some $g$ where $f \cpxlace g$.  Consequently, 
  $$ Hull(f) = \roots(\,\{ g \mid f\cpxlace g\}\,)$$
\end{prop}
\begin{proof}
  lemma~\ref{lem:gauss-lucas-1} shows that
  $  \roots(\,\{ g \mid f\cpxlace g\}\,) \subseteq Hull(f)$.
  Let $\zeta=\{\zeta_1,\dots,\zeta_n\}$ be the roots of $f$. 
  We first show that if $0\in Hull(f)$ then $0\in\roots(\,\{ g
  \mid f\cpxlace g\}\,)$. For $\beta=(\beta_1,\dots,\beta_n)$ we
  define
\begin{gather}
 f_\beta = \sum_{k=1}^n \beta_k \frac{f}{x-\zeta_k}\label{eqn:fbeta}\\
\intertext{The requirement that $f_\beta(0)=0$ implies that}
 0 = \pm \zeta_1\cdots\zeta_n\, \left(\frac{\beta_1}{\zeta_1} + \cdots
   +   \frac{\beta_n}{\zeta_n}\right)\notag
\end{gather}
Thus, it suffices to see that $0$ is in the convex hull of
$\{1/\zeta_1,\dots,1/\zeta_n\}$, which is clear, since $0$ is in the
convex hull of $\{\zeta_1,\dots,\zeta_n\}$.

Choose $\alpha$ in the convex hull of $f$. 
Define the Blaschke  (see last section) transformation
$$ T\colon{} z \mapsto \frac{z-\alpha}{1-\overline{\alpha}z}$$
Consider $f_\beta(z)$ and the transformed polynomial $g_\beta(z) =
(1-\overline{\alpha}z)^nf_\beta(T\,z)$. If $\beta$ is chosen so that
$g_\beta(0)=0$ then $f_\beta(\alpha)=0$. So, we need to know that if
$\alpha$ is in the convex hull of the $\zeta_i$ then $0$ is in the
convex hull of $T\zeta_i$. But since $\alpha\in Hull(f)$ there are
$\gamma_i$ so that
\begin{gather*}
  \alpha = \sum \gamma_i \zeta_i \\
\intertext{and so}
0 = T(\alpha) = \sum \gamma_i \,T(z_i)
\end{gather*}
Thus, $0$ is in the hull, and we are done.
\end{proof}

This last proposition showed that there are no restrictions on the
locations of roots of polynomials interlacing a given polynomial.
However, there are restrictions on the simultaneous location of
\emph{all} the roots of a given polynomial interlacing $f$. Write
$f(x)=(x-r_1)\cdots(x-r_n)$, and define
$$ f_\beta(z) = \sum \beta_i \frac{f(z)}{z-r_i}$$
Letting $s_1,\dots,s_{n-1}$ be the roots of $f_\beta$, and
substituting $z=r_k$ we have
\begin{gather*}
 \prod_i (r_k - s_i) =  \beta_k \prod_{j\ne k} (r_k-r_j)  \\
\intertext{and since $\beta_k$ is positive we have}
 \sum_i \arg(r_k-s_i) = \sum_{j\ne k} \arg(r_k-r_j) 
\end{gather*}
This says that the sum of the angles (mod $2\pi$) formed at $r_k$ by
the roots of $f_\beta$ is a fixed quantity. We can therefore derive
restrictions such as
\begin{quote}
  The roots of $f_\beta$ can not satisfy $$ 0\le \arg(r_k-s_i) <
  \frac{1}{n-1}\sum_{j\ne k}\arg(r_j-r_k) \quad\quad\text{for } 1\le
  i\le n-1$$
\end{quote}
For example, if $n=3$ this  restriction becomes:
\begin{quote}
  The  angles at a root of $f$ formed by the two roots can not both be
  less than half the angle at the vertex. 
\end{quote}
In the picture below, it follows that $f_\beta$ has one root in each
triangle marked $a$ (or one each in $b$, or one each in $c$).
The interior lines are the angle bisectors.

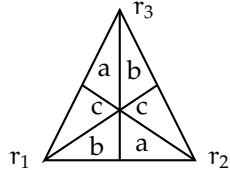
\begin{figure}[htbp]
  \centering
\begin{pspicture}(0,0)(3,2)
\psline(1,0)(3,0)(2,2)(1,0)(2.5,1)
\psline(3,0)(1.5,1)
\psline(2,2)(2,0)
\uput[l](1,0){$r_1$}
\uput[r](3,0){$r_2$}
\uput[r](2,2){$r_3$}
\rput(1.7,.2){b}
\rput(1.7,.7){c}
\rput(1.8,1.2){a}
\rput(2.3,.2){a}
\rput(2.3,.7){c}
\rput(2.2,1.2){b}
\end{pspicture}
  
  \caption{Where the roots are.}
  \label{fig:arg-pic}
\end{figure}

\section{The map $\phi$}

When we have a formula for $T(x-z)^n$ we can use this information to
determine a map $\phi'$ so that \eqref{eqn:diagram-phi} is
satisfied. We begin with multiplier transformations.

\begin{example}

Multiplier transformations have a very simple action in the complex
plane.  Suppose that $T$ is a multiplier transformation. We have a
fundamental commutativity diagram that is immediate from the
definition of multiplier:

\centerline{\xymatrix{
    .
      \ar@{->}[d]_{x\mapsto \sigma x}
      \ar@{->}[rrr]^{T} 
      &&&
     .
      \ar@{->}[d]^{x\mapsto \sigma x} \\
      .
      \ar@{->}[rrr]^{T}
      &&&
      .
}}

\begin{lemma}
  Suppose that $T$ is a multiplier transformation $x^k\mapsto a_kx^k$
  and $r_1,\dots,r_n$ are the roots of $T(x-1)^n$. The definition of
  $\phi'$ below makes \eqref{eqn:diagram-phi} commute.
  $$\phi': z \mapsto \prod_k(x-r_kz)$$
\end{lemma}

\begin{proof}
  If $z\in \complexes$, then $T(x-z)^n(x)=0$ implies that
  $T(\frac{x}{z}-1)^n(x)=0$. Since $T$ is a multiplier transformation
  this is equivalent to $T(x-1)^n(\frac{x}{z})=0$, which implies the
  result. 
\end{proof}

\end{example}

\begin{example}

  If we use the formula \eqref{eqn:laguerre-addition} for Laguerre
  polynomials with the transformation $T(x^k)=L_k(x)$ we get
\begin{align}
    T(x+y)^n &= (y+1)^n T(x^n)(\frac{x}{y+1})\label{eqn:laguerre-new-1}\\
    T\left(\frac{x+y}{y+1}\right)^n &= T(x^n)(\frac{x}{y+1})\notag\\
\intertext{and by linearity yields}
    T\left(f\left(\frac{x+y}{y+1}\right)\right) &= T(f)(\frac{x}{y+1})\notag
  \end{align}

If we substitute $\alpha(y+1)$ for $x$ in \eqref{eqn:laguerre-new-1}
then we see that
$$  T(x+y)^n\,(\alpha(y+1)) = (y+1)^n T(x^n)(\alpha) = (y+1)^n L_n(\alpha)$$ and that
therefore the roots of $T((x-y)^n)$ are given by $\alpha(1-y)$, where
$\alpha$ is a root of $L_n$. In particular,  it follows that

\begin{lemma} \label{lem:laguerre-phi}
Suppose $T(x^n) = L_n(x)$ has roots $r_1,\dots,r_n$. The following
definition of $\phi'$ makes \eqref{eqn:diagram-phi} commute:

$$ \phi': z\mapsto \prod_{k=1}^n (\,x- r_i(1+z)\,)$$
\end{lemma}

For instance, if $\alpha$ is positive then the roots of
$T(x+\alpha\imag)^n$ lie on vertical rays emanating from the roots of
$L_n$. 

\end{example}

\begin{example}
We use the addition formula \eqref{eqn:her-tf-1} for Hermite
polynomials to determine $\phi'$. 
\begin{lemma} 
Suppose $T(x^n) = H_n(x)$ has roots $r_1,\dots,r_n$. The following
definition of $\phi'$ makes \eqref{eqn:diagram-phi} commute:

$$ \phi': z\mapsto \prod_{k=1}^n (\,x- (r_i + (1/2)( 1+z))\,)$$

\end{lemma}
\begin{proof} 
  From the addition formula \eqref{eqn:her-tf-1}  with $f=(x+1)^n$ we
  have
  \begin{gather*}
    T(x+1)^n\bigl|_{x\mapsto x+y} = T_\ast(x+2y+1)^n
\intertext{If we let $z=2y+1$ then}
    T(x+1)^n\,(x-1/2+z/2) = T_\ast (x+z)^n
  \end{gather*}
If $r$ is a root of $T(x+1)^n$, then $x-1/2+z/2=r$ yields a
root of the right hand side. Thus $x = r+1/2 + z/2$.    

\end{proof}

\end{example}

\begin{example}
  Suppose that $T(x^n) = x^n H_n(1/x) = H_n^{rev}(x)$. We have the
  following identity that can be derived from the similar identity
  for the Hermite polynomial transformation.

  \begin{equation}
    \label{eqn:hermite-rev}
    T(x+2u-2)^n = u^n\,T(x^n)\,(\frac{x}{u})
  \end{equation}

  \begin{lemma}
    Suppose $T(x^n) = x^n H_n(1/x)$ has roots $r_1,\dots,r_n$. The
    following definition $\phi'$ makes \eqref{eqn:diagram-phi}
    commute: 
$$ \phi': z\mapsto \prod_k\,\left( x - r_k(1+z/2)\,\right)$$

  \end{lemma}

  \begin{proof}
    We must show that $T(x-z)^n = \prod_k\,\left( x -
      r_k(1+z/2)\,\right)$. The conclusion follows as before.
  \end{proof}
\end{example}
\section{Lines, rays, and sectors}

We show that some transformations preserve lines or rays or sectors in
the complex plane.

\begin{lemma}
  Suppose $T\colon{}x^n\mapsto a_nx^n$ is a multiplier transformation that
  maps $\allpolypos$ to itself.  For any $\sigma\in\complexes$, we
  have $T\colon{}\allpolyint{\sigma\reals}\longrightarrow\allpolyint{\sigma\reals}$.
\end{lemma}
\begin{proof}
Suppose $f(x) = \prod_1^n(x-\sigma r_i)$ where $r_i\in\reals$, and let
$g(x) = \prod(x-r_i)$. Then since $T$ is a multiplier transformation
\begin{gather*}
  T\,f(x) = T(f(\sigma x))(x/\sigma) = T(  \sigma^n g(x))(x/\sigma) =
  \sigma^n T(g(x))(x/\sigma)
\end{gather*}
Since $T\colon{}\allpoly\longrightarrow\allpoly$ we know that
$T(g(x))\in\allpoly$, and so $TF\in\allpolyint{\sigma\reals}$.
\end{proof}

Similarly, the roots of $(x-\imag-\alpha)^n$ lie on the horizontal
lines $\{\imag r_k +\reals\}$.  



  Here is a very general way of finding linear transformations that
  preserve lines.

  \begin{prop}
    Suppose $f\in\allpolyf$. The linear transformation \\ $g\mapsto
    f(\sigma\diffd)g$ maps $\allpolyint{\xi+\sigma\reals}$ to
    itself for any $\xi\in\reals$ and $\sigma\in\complexes$.
  \end{prop}
  \begin{proof}
    Since $f$ is a limit of polynomials in $\allpoly$, it suffices to
    prove the result for {linear} polynomials. So, choose
    $g\in\allpolyint{\xi+\sigma\reals}$; it suffices to show that
    $(\sigma\diffd+\beta)g = \sigma g' + \beta g$ is also in
    $\allpolyint{\xi+\sigma\reals}$ for any
    $\beta\in\reals$. Since $\frac{d}{dx}g( x+\gamma) =
    g'( x+\gamma)$,
    we take $\xi=0$. 
    
    If $g\in\allpolyint{\cals}$ and we write $g(x) =
    \prod_1^n(x-\sigma r_k)$ with $r_i\in\reals$, then $$g(\sigma x) =
    \sigma^n\prod(x-r_k).$$
    If we define $h(x) = \sigma^{-n}g(\sigma
    x)$ then $h\in\allpoly$. Note also that $h'(x) = \sigma^{-n+1}
    g'(\sigma x)$.  When $ k(x) = \sigma g'(x) + \beta g(x)$ then
$$
k(\sigma x) = \sigma g'(\sigma x) + \beta g(\sigma x) 
= \sigma^{n}( h'(x) + \beta h(x)).
$$
 Since $h\in\allpoly$, it follows that $\sigma^{n}k(\sigma
x)\in\allpoly$, and hence the roots of $k(x)$ are real multiples of $\sigma$.
  \end{proof}

  \begin{cor}
    If $f(x)\in\allpolyint{\alpha+\imath\reals}$ then 
     $f(x+1)+f(x-1)\in\allpolyint{\alpha+\imath\reals}$.
  \end{cor}
  \begin{proof}
    Simply note that 
$$ f(x+1)+f(x-1) = (e^\diffd + e^{-\diffd})f = 2 \cos(\imath\diffd)f$$
and $\cos(x)\in\allpolyf$. 
  \end{proof}

  \begin{remark}
    This transformation appears in \cite{bump}, where it is used to
    show that certain polynomials have roots whose real part is
    $1/2$. 
  \end{remark}

In the the next two examples the regions are rays and horizontal lines.

\begin{lemma}
  Consider $T\colon{}x^k\mapsto H_k(x)x^{n-k}$ acting on $\allpoly(n)$. Every
  polynomial with roots on a ray in the complex plane based at $2$ is
  mapped by $T$ to a line through the origin with the same slope.
\end{lemma}
\begin{proof}
 Assume that $f$ has all its roots on the line $2+\alpha \sigma$,
where $\sigma\in\complexes$ is fixed, and $\alpha\ge0$. If we define
$S(f)=T(f(x-2))$ , then it suffices to determine the action of $S$ on
polynomials with roots on the ray through $\sigma$. Since
$\sigma^{-n}f(\sigma x))$ has all positive roots, the result
follows from Lemma~\ref{lem:hermite-hxn} 
\end{proof}

\begin{lemma}\label{lem:hermite-on-cpx-line}
  If $T(x^k)= H_k(x)$ then $\allpolyint{\reals+\alpha{\imath}}
  \longrightarrow \allpolyint{\reals+\alpha{\imath}/2}$.  In other
  words, if $f(x)$ is a polynomial whose roots all have imaginary part
  $\alpha\imath$, then the roots of $T(f)$ all have imaginary part
  $\frac{\alpha\imath}{2}$.
\end{lemma}

\begin{proof}

The  addition formula for Hermite polynomials
\eqref{eqn:hermite-addition}  implies that
\begin{align}
T(x^n)(x+y) &= T_\ast(x+2y)^n\notag\\
\intertext{and by linearity }
T(f)(x+y) &= T_\ast f(x+2y)
\intertext{Replacing $x+y$ by $x$ in this equation yields}
T(f) &= (T_\ast(f(x+2y))\,(x-y)\label{eqn:her-tf-1}
  \end{align}

 Equation \eqref{eqn:her-tf-1} implies the following diagram commutes:

\centerline{\xymatrix{
\allpolyint{\reals+\alpha{\imath}} \ar@{.>}[d]_T \ar@{->}[rr]^{x\mapsto
  x+\alpha{\imath}} && \allpoly \ar@{->}[d]^T  \\
\allpolyint{\reals+\alpha{\imath}/2} &&  \ar@{->}[ll]^{x\mapsto x+\alpha{\imath/2}}
\allpoly
}} 
\end{proof}

Multiplier transformations preserve sectors because rotation commutes
with such transformations.

\begin{lemma}
  Suppose $T\colon x^n\mapsto a_nx^n$ maps $\allpoly$ to itself and has all
  positive coefficients. If $\mathcal{S}$ is a sector of angle at most
  $\pi$ then $T$ maps $\allpolyint{\mathcal{S}}$ to itself. 
\end{lemma}
\begin{proof}
  Since $T$ preserves interlacing and the sign of the leading
  coefficient it maps $f+\imag g\in\polycpx$ to $Tf+\imag
  Tg\in\polycpx$. If $|\sigma|=1$ is a rotation then $T$ maps
  $\allpolyint{-\sigma\plane}$ to itself. Consequently.
\[
T\colon{}\allpolyint{-\sigma\plane\cap-\plane}=\allpolyint{-\sigma\plane}\cap\allpolyint{-\plane}
\longrightarrow
\allpolyint{-\sigma\plane}\cap\allpolyint{-\plane}=\allpolyint{-\sigma\plane\cap-\plane}
\]
We can choose $\sigma$ so that $-\sigma\plane\cap-\plane$ is a sector
of any desired opening of at most $\pi$. Another rotation to
$\mathcal{S}$ finishes the proof.

\end{proof}

\section{The set of roots }

In this section we look at properties of the set of roots. 
If we are in the fortunate situation of \eqref{eqn:diagram-phi} then
for any subset $\cals$ of $\complexes$ we know that

$$ \roots(\, T \phi \cals\,) = \roots(\, \phi' \cals\,)$$

In general we can not describe the behavior so well. 
We next  note that we can refine Lemma~\ref{lem:find-range} as follows:

\begin{lemma}
  Suppose that $T\colon{}\allpolyint{(a,b)}\longrightarrow\allpoly$ preserves
  interlacing. Then
$$\roots(\, T(\phi\,(a,b))\,) = \roots(\, T(\allpolyint{(a,b)})\,).$$
\end{lemma}
\begin{proof}
  The proof follows from the proof of Lemma~\ref{lem:find-range} and
  using the continuity of roots.
\end{proof}

\begin{figure}[htbp]
  \resizebox{2in}{2in}{%
  \begin{pspicture}(0,-5)(9,5)
    \pscircle(.1,0){.1}
    \pscircle(.5,0){.5}
    \pscircle(1.25,0){1.25}
    \pscircle(2.5,0){2.5}
    \pscircle(4.5,0){4.5}
    \psline(0,0)(9,0)
  \end{pspicture}
  }
    \caption{Image of the unit circle under Laguerre transformation}
    \label{fig:laguerre-circle}
\end{figure}
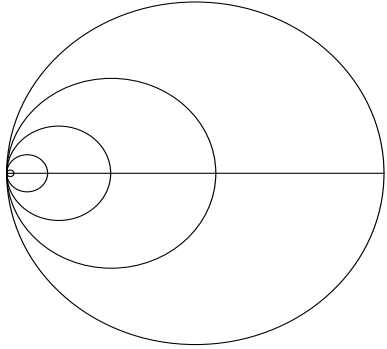

If $\ell_n$ is
the largest root of $L_n$, then from Lemma~\ref{lem:laguerre-phi}
$$
\roots(\, T(\phi \Delta)\,) = \bigcup \{ \text{roots of } T(x-z)^n\, \mid
\, |z|\le 1\}$$
consists of all points in the circle of radius
$\ell_n$ centered at $\ell_n$. See Figure~\ref{fig:laguerre-circle}
which shows $\roots(\,T(\phi (\partial\Delta))\,)$ when $n$ is $5$.  Since
$\ell_n$ goes to infinity as $n$ goes to infinity, the union of all
these circles covers the right half plane. This proves

\begin{prop}\label{prop:laguerre-image}
  If $T(x^k)=L_k$, then $\{\Re(\zeta)\ge0\} \subseteq \roots(\, \,T(\phi \Delta)\,).$
\end{prop}

In other words, if $\Re(\xi)\ge0$ then there is a $z$ and positive
integer $n$ where $|z|\le1$ such that $\xi$ is a root of $T(x-z)^n$.


\section{Complex Transformations involving $f(x+\imag)$}
\label{sec:linear-complex}

\index{linear transformation!complex}

We introduce some transformations that involve complex numbers, and
show that they preserve $\allpoly$ because they are other known
transformations in disguise. 

To begin, define the transformation $T(f) = f(x+i) + f(x-i)$.  If $f$
has all real coefficients, then $T(f)$ has all real coefficients,
since $Tf$ is equal to its conjugate.

\begin{lemma}
  The transformation $T\colon{}f(x)\mapsto f(x+i)+f(x-i)$ preserves roots. 
\end{lemma}

\begin{proof}
  Since $e^{a\diffd} = f(x+a)$ it follows that $ \cos(\diffd)f = \frac{1}{2}Tf$.
  Since $\cos(x)$ is in $\allpolyf$, $\cos(\diffd)$ preserves roots, and thus so
  does $T$.
\end{proof}

\begin{remark}
  More generally, if we use the fact that $\cos (ax+b)\in\allpolyf$ we
  see that 
\[
\cos(a \diffd + b) f(x) = \frac{1}{2}\bigl( 
e^{\imag b} f(x+\imag a)+ e^{-\imag b} f(x-\imag a)\bigr)
\]
is in $\allpoly$.   \cite{cardon04}
\end{remark}

We can iterate the construction of the following lemma. If we let the number of
terms go to zero then the sum becomes a convolution. See
\cite{cardon}.
\index{convolution}

\begin{lemma}
  If $a_1,\dots,a_r,b_1,\dots,b_r$ are real numbers, and $f\in\allpoly$, then
$$\frac{1}{2^r} \sum f(x\pm a_1i \pm a_2i \pm\cdots\pm
a_ri)e^{\imag(\pm b_1 \pm b_2\cdots\pm b_r)} \in\allpoly$$
where the sum is over all $2^r$ choices of $\pm$ signs. 
\end{lemma}
\begin{proof}
  This is $\cos(a_1\diffd+b_1)\cdots\cos(a_r\diffd+b_r)f$.
\end{proof}

We get a similar result using $\sin$ instead of $\cos$. The proof is omitted.

\begin{lemma}
  The map $f(x) \mapsto \frac{ f(x+i) - f(x-i)}{2i}$ maps $\allpoly$
  to itself.
\end{lemma}

Next we look at a variation of the above.

\begin{lemma}\label{lem:sum-of-imag}
  If $n$ is a positive integer and $f(x)\in\allpoly$ then 
$$  f(x) + \sum_{k=-n}^n f(x+k\imag) \quad \in\allpoly$$
\end{lemma}
\begin{proof}
  We compute that
  \begin{align*}
    f(x) + \sum_{k=-n}^n f(x+k\imag) &= 2(1 + \cos(\diffd) +
    \cos(2\diffd) +\cdots +     \cos(n\diffd)) \,f(x)  \\
    &= 2 \left[\frac{\sin((n+1)x/2)}{\sin(x/2)} \, \cos(nx/2)\right](\diffd)\,f(x)
  \end{align*}
Now $\cos(nx/2)$ is in $\allpolyf$, and so is the first factor, since
we can write an infinite product expansion for it, and the factors on
the bottom are canceled out by terms on the top.
\end{proof}

\begin{cor}\label{cor:int-of-i}
  If $f(x)\in\allpoly$ then $\displaystyle\int_{-1}^1 f(x+\imag t)\,dt\in\allpoly$.
\end{cor}
\begin{proof}
 Approximate the integral by the sum
$$ \frac{1}{2n} \left( f(x) + \sum_{k=-n}^n f(x+\imag k/n)\right).$$
This sum is in $\allpoly$, for we can scale $f(x)$ in the Lemma.

Here is an alternative proof not using sums. If $T(f)$ is the integral
in the Corollary, then
\begin{gather*}
  T(x^n) = \int_{-1}^1 (x+\imag t)^n\,dt = \frac{1}{\imag} \left(
  \frac{(x+\imag)^{n+1}}{n+1} -\frac{(x-\imag)^{n+1}}{n+1} \right) \\
= 2 \frac{e^{\imag\diffd} -e^{-\imag\diffd}}{2\imag} \,
  \frac{x^{n+1}}{n+1} = 2\sin(\diffd) \intoper{} x^n
\end{gather*}
Thus $T(f) = 2\sin(\diffd)\intoper{} f$. This in in $\allpoly$ since 
\begin{equation}\label{eqn:sinD}
 \sin(\diffd) \intoper{}f = \left( \frac{\sin(x)}{x}(\diffd)\right)
\diffd (\intoper{}f) = \left( \frac{\sin(x)}{x}(\diffd)\right) f
\end{equation}
and $\frac{\sin(x)}{x}$ is in $\allpolyf$.
\end{proof}

Translating and scaling this result gives the following corollary,
which will prove useful when we look at polynomials with complex
coefficients (Chapter~\ref{cha:complex-coef}).
\begin{cor}\label{cor:int-of-i-2}
  If $\alpha>0$ and $f(x)\in\allpoly$ then the all the roots of
  $\displaystyle\int_0^1 f(x+\alpha\imag t)\,dt$ have imaginary part $-\alpha/2$.
\end{cor}
\begin{proof}
  Observe that
$$ \int_0^1 f(x+\alpha\imag t)\,dt =
\frac{1}{2}\int_0^2 f(f+\alpha\imag t/2)\,dt =
\frac{1}{2}\int_{-1}^1f([x-\alpha\imag/2] + \alpha\imag t/2)\,dt
$$
Since the last integral has all real roots, the roots of the first
integral have imaginary part $-\alpha/2$.
\end{proof}


\section{The Hermite-Biehler problem}

  This section is concerned with determining how many roots of a
  polynomial with complex coefficients lie in the upper half plane.
  If a polynomial $h$ has $\alpha$ roots in the $\mathcal{U}$pper half
  plane, and $\beta$ roots in the $\mathcal{L}$ower half plane, then
  we define $\updn(h)=(\alpha,\beta).$

\begin{quote}
  \textbf{ The Hermite-Biehler problem:} Suppose that $f,g$ have all real
  coefficients and no common real roots.  Express $\updn(f+\imag g)$
  in terms of the arrangement of the roots of $f$ and $g$.
\end{quote}

\index{Hermite-Biehler problem}

Suppose $f$ and $g$ are two polynomials with real coefficients, with
perhaps some complex roots. The polynomial $f+\imag g$ has the
following fundamental property which follows from the fact that $f$ and
$g$ have all real coefficients:
\begin{quote}
  If $f+\imag g$ has a real root, then $f$ and $g$ have that root in common.
\end{quote}
Thus, if $f$ and $g$ have no real roots in common, then all
the roots of $f+\imag g$ are complex.  Our approach is based on this
simple consequence:

\begin{quote}
If $f_t$ and $g_t$ are a continuous family of real polynomials where
$f_t$ and $g_t$ have no common root for all $0\le t\le1$, then
$\updn(f_0+\imag g_0) = \updn(f_1+\imag g_1)$ since the roots can
never cross the real axis. 
\end{quote}

We begin by describing three transformations of $f+\imag g$ that
preserve the number of roots in each half plane. We assume that $f$
and $g$ have no common real roots. If $\updn(f+\imag g) =
(\alpha,\beta)$ then $ \updn(-f-\imag g) = (\alpha,\beta)$. Taking
conjugates shows that
$$\updn(f-\imag g) = \updn(-f+\imag g) = (\beta,\alpha)$$
Consequently, we will always assume that $f$ and $g$ have positive
leading coefficients. 

  \begin{enumerate}
  \item \textbf{ Join adjacent roots}
    
    Suppose that $u,v$ are adjacent roots of $g$ that are not
    separated by a root of $f$, and let $u_t = tu+(1-t)v$. The
    transformation below doesn't ever create a $g_t$ with a root in
    common with $f$ since $u_t$ lies between $u$ and $v$. Thus
    $\updn(h) = \updn(h_1)$ since $h_0=h$.
    
    $$ h_t = f + \imag \frac{g}{x-v}\cdot(x-u_t)  $$

  \item \textbf{ Replace double roots with $x^2+1$}
    
    Suppose that $u$ is a double root of $g$. The following
    transformation moves a pair of roots of $g$ into the complex
    plane, missing all real roots of $f$, and lands at $\pm\imag$.
    Consequently, we get a factor of $x^2+1$. Again, $\updn(h) =
    \updn(h_1)$.

$$ h_t = f + \imag \frac{g}{(x-u)^2}\, \left((x-(1-t)u)^2+t\right)$$

\item \textbf{ Replace complex roots with $x^2+1$}

  If $a+b\imag$ is a root of $g$ then since $g$ has all real
  coefficients $g$ has a factor of $(x-a)^2+b^2$. The transformation
  below converts this factor to $x^2+1$.

$$ h_t = f+\imag\frac{g}{(x-a)^2+b^2} \, \left((x-(1-t)a)^2+ (1-t)b^2+t\right)$$

\end{enumerate}

We can also do the above constructions to $f$. If we do them as many
times as possible, then we have found a polynomial that we call the
\textbf{HB-simplification} of $h$:

$$ H = (x^2+1)^rF + \imag (x^2+1)^sG$$

The constructions above guarantee that $\updn(h)=\updn(H)$. In
addition, $F$ and $G$ interlace, since we have removed all consecutive
roots belonging to only one of $f,g$. We can now state our result:

  \begin{prop}\label{prop:hb}
    Suppose $f$ and $g$ are polynomials with real coefficients,
    positive leading coefficients, and no common real roots. Let
    $H=(x^2+1)^rF + \imag (x^2+1)^s G$ be the HB-simplification of
    $h=f+\imag g$, and set $m = \max(r,s)$. Then
$$ \updn(h) = 
\begin{cases}
  (m,m)+(0,deg(F)) & \text{if }F\longleftarrow G\\
  (m,m)+(deg(G),0) & \text{if }G\longleftarrow F\\
\end{cases}
$$
  \end{prop}

  \begin{proof}
    
    If $F$ and $G$ have the same degree then without loss of
    generality we may assume that $F\greateq G$.  Suppose that $r>s$.
    Consider
    
    $$
    h_t = (x^2+1)^r\,f + \imag(tx^2+1)(x^2+1)^s\,g$$
    Since
    $deg(h_t) = deg(h)$, $\updn(h_1) = \updn(h)$. We can do a similar
    transformation if $r<s$, so we may assume that both exponents are
    equal to $m=\max(r,s)$.  If write $(x^2+1)^mH_1=H$ then we see
    that $\updn(h) = \updn(h_1) + (m,m)$.  Since $H_1 = F+\imag G$, we
    can apply the classical Hermite-Biehler theorem to conclude that
    $\updn(H_1) = (0,deg(F))$. If $F\lessless G$ then we consider
    $G_\epsilon=\epsilon F+G$ and let $\epsilon$ go to zero. The
    proposition is now proved.

For completeness, here's a quick proof of the
Hermite-Biehler theorem.

Suppose $h=f+\imag g$ with $f\greateq g$, and the roots of $f$ are
$v_1,\dots,v_n$ then there is a positive $\sigma$ and positive $a_i$
so that
\begin{gather*}
  g = \sigma f + \sum a_i \frac{f}{x-v_i} \\
\intertext{A root $\rho$ of $f+\imag g$ satisfies}
  \sigma + \sum \frac{a_i}{\rho-v_i}= \imag  \\
  \sigma + \sum \frac{a_i(\overline{\rho}-v_i)}{|\rho-v_i|^2} = \imag 
\end{gather*}
from which it follows that $\rho$ has negative imaginary part.

  \end{proof}

  \begin{example}
    Here is an example. Suppose that 
    \begin{align*}
      f &= (x^2+2)(x^2+3)(x-2)(x-4)(x-5)(x-8)(x-9)(x-10)\\ 
      g &= (x^2+5)(x-1)(x-3)(x-6)(x-7) \\
\intertext{We write the roots of $f$ and $g$ in order, and successively
  remove pairs of adjacent roots of $f$ (or of $g$)}
        & g f g f f g g f f f \\
        & g f g g g f f f \\
        & g f g f f f \\
        & g f g f
    \end{align*}

Thus $F\lessless G$. Since $f$ has degree $10$ and $g$ has degree $6$,
we see that $r=(10-2)/2=4$, and $s = (6-2)/2=2$. Thus $m=4$ and
$$ \updn(f+\imag g) = (4,6)$$
  \end{example}
  
The image of the upper half plane $\plane$ under a counterclockwise rotation of
angle $\alpha$ is $e^{\imag \alpha}\plane$. A polynomial $f(x)$ has all its
roots in this half plane if and only if 
$f(e^{\alpha\imag})\in\allpolyint{\plane}$, and this is the case
exactly when
$$ \Re(\,f(e^{\alpha\imag})\,)\lessless \Im(\,f(e^{\alpha\imag})\,)$$
and the leading coefficients have opposite signs.
However, it is difficult in general to express the real and imaginary
parts of $f(e^{\imag \alpha})$ in terms of coefficients of $f$. When
$e^{\imag\alpha}$ is a root of unity then there is an explicit
formula. We describe the general case, and then look at three special
cases. 

Choose a positive integer $n$ and let $\zeta=e^{2\pi\imag/n}$. Write
$f= \sum a_i x^i$ and define
\begin{align*}
  f_k(x) &= \sum_{s\equiv k\pmod{n}} a_s x^s \\
\intertext{ We can express $f(\zeta x)$ in terms of the $f_k$:}
f(\zeta x) &= \sum_{k=0}^{n-1} f_k(\zeta x) \,=\, \sum_{k=0}^{n-1} \zeta^k f_k(x) \\
\intertext{If we define}
F_0(x) &= \sum_{k=0}^{n-1} \cos(2k\pi/n)\,\Re(f_k(x)) - \sin(2k\pi/n)\,\Im(f_k(x)) \\
F_1(x) &= \sum_{k=0}^{n-1} \sin(2k\pi/n)\,\Re(f_k(x)) + \cos(2k\pi/n)\,\Im(f_k(x)) \\
\intertext{then since $\zeta = \cos(2k\pi/n) +\imag \sin(2k\pi/n)$}
f(\zeta x) &= F_0(x) + \imag F_1(x)
\end{align*}
and $F_0$ and $F_1$ have all real coefficients. Consequently,
$f\in\allpolyint{\zeta \plane}$ iff $F_0\lessless F_1$ and the leading
coefficients have opposite signs. 

If $f(x)$ has all real coefficients then the imaginary parts are zero
and we have
\begin{align*}
F_0(x) &= \sum_{k=0}^{n-1} \cos(2k\pi/n)\,f_k(x)\\
F_1(x) &= \sum_{k=0}^{n-1} \sin(2k\pi/n)\,f_k(x)
\end{align*}

We now look at three special cases.

\begin{example}
  If $n=3$ then $\zeta = \frac{-1}{2} + \frac{\sqrt{3}}{2}\imag$ , and
    $\zeta^2 =  \frac{-1}{2} - \frac{\sqrt{3}}{2}\imag$.  We define
\begin{align}\label{eqn:f0f1}
  F_0(x) &= \Re(f_0) + \frac{-1}{2}\Re(f_1) -
  \frac{\sqrt{3}}{2}\Im(f_1) + \frac{1}{2}\Re(f_2) -
  \frac{\sqrt{3}}{2}\Im(f_2) \\
  F_1(x) &= \Im(f_0) + \frac{\sqrt{3}}{2}\Re(f_1) +
  \frac{-1}{2}\Im(f_1) + \frac{\sqrt{3}}{2}\Re(f_2) -
  \frac{1}{2}\Im(f_2) \label{eqn:f0f1-2}
\intertext{If the coefficients of $f$ are real then}
  F_0(x) &= f_0 + \frac{-1}{2}f_1  + \frac{1}{2}f_2  \notag \\
  F_1(x) &=  \frac{\sqrt{3}}{2}f_1 + \frac{\sqrt{3}}{2}f_2  \notag
\end{align}

 We conclude:

 \begin{lemma}
   All roots of $f$ have argument between $2\pi/3$ and $5\pi/3$ if and
   only if $F_0$ and $F_1$ in \eqref{eqn:f0f1},\eqref{eqn:f0f1-2}
   satisfy $F_0\lessless F_1$ and their leading coefficients have
   opposite signs. If $f$ has all real coefficients then the
   interlacing criterion is
$$ 2f_0 - f_1+f_2 \lessless f_1+f_2.$$
 \end{lemma}
\end{example}

\begin{example}
  We now take $n=4$. Since $e^{2\pi\imag/4}=\imag$, if we follow 
 the calculation above, then
  \begin{align*}
    F_0 &= \Re(f_0) - \Im(f_1) - \Re(f_2) + \Im(f_3) \\
    F_1 &= \Im(f_0)  +\Re(f_1) -\Im(f_2) - \Re(f_3)
  \end{align*}
  and consequently $f$ has roots with only negative real parts if and
  only if $F_0\lessless F_1$. This can be further simplified. If we
  assume that $f(x)$ has all real coefficients then 
  \begin{align*}
    F_0 &= f_0  - f_2 \\
    F_1 &= f_1  - f_3
  \end{align*}
  If we write $f(x) = p(x^2) + x q(x^2)$ then $F_0(x) = p(-x^2)$ and
  $F_1(x) = xq(-x^2)$.  In this case, the result is also known as the
  Hermite-Biehler \cite{holtz} theorem:

\begin{theorem}[Hermite-Biehler]
  Let $f(x) = p(x^2) + x q(x^2)$ have all real coefficients. The
  following are equivalent:
  \begin{enumerate}
  \item All roots of $f$  have negative real part.
  \item  $p(-x^2)$ and $xq(-x^2)$ have simple interlacing roots, and
    the leading coefficients are the same sign.
  \end{enumerate}
\end{theorem}
\index{Hermite-Biehler}

\end{example}

\begin{example}
  We next take $n=8$, where $\zeta = \frac{1+\imag}{\sqrt{2}}$. We
  assume that all coefficients of $f$ are real. Define
  \begin{align}
    F_0(x) &= 
{f_0} + \frac{{f_1}}{{\sqrt{2}}} - \frac{{f_3}}{{\sqrt{2}}} - {f_4} -
\frac{{f_5}}{{\sqrt{2}}} + \frac{{f_7}}{{\sqrt{2}}}
\label{eqn:8-1} \\
F_1(x) &= 
\frac{{f_1}}{{\sqrt{2}}} + {f_2} + \frac{{f_3}}{{\sqrt{2}}} -
\frac{{f_5}}{{\sqrt{2}}} - {f_6} - \frac{{f_7}}{{\sqrt{2}}}
 \label{eqn:8-2}
  \end{align}

  \begin{lemma}
    If $f(x)$ has all real coefficients, then $\Re(r)< \Im(r)$ for
    every root $r$ of $f(x)$ if and only if $F_0,F_1$ in
    \eqref{eqn:8-1},\eqref{eqn:8-2} have leading coefficients of
    opposite signs, and satisfy $F_0 \lessless F_1$.
  \end{lemma}
\end{example}

Suppose we are given a line $\ell$ in the complex plane such that the
angle the line makes with the real line is a rational multiple of
$2\pi$.  We can find an interlacing restriction such that a polynomial
has all its roots in a particular half plane determined by $\ell$ if
and only if the interlacing condition is met, along with agreement or
disagreement of the signs of the leading coefficients.

We can combine these interlacing and sign conditions, provided the
relevant angles are rational. We can summarize this as

\begin{lemma}
  Suppose that $\cals$ is a convex region determined by $n$ lines, and
  each line makes a rational angle with the real axis.  Then, we can
  find $n$ interlacing and sign restrictions such that a polynomial
  has all its roots in $\cals$ if and only if these $n$ restrictions
  are satisfied.
\end{lemma}

\section{The B\^ocher-Walsh problem}
\label{sec:bocher-walsh-problem}

\index{B\^ocher-Walsh theorem}

The essence of the B\^ocher-Walsh theorem can be stated
\cite{sheil-small} in the following form
\begin{quote}
  If all roots of $f$ lie in the upper half plane and all roots of $g$
  lie in the lower half plane then $\smalltwodet{f}{f'}{g}{g'}$ has no
  real zeros.
\end{quote}

\index{Gauss-Lucas theorem}

This is reminiscent of the Gauss-Lucas theorem, which in general
follows from the fact that if $f$ has all its roots in the upper half
plane then so does $f'$. Unlike the Hermite-Biehler theorem, it's easy
to determine the location of the roots of the determinant.

\begin{lemma}
  Suppose that $f \algless f_1$ has all roots in the upper half plane,
  and $g\algless g_1$ has all its roots in the lower half plane. If
  $deg(f) = n > m = deg(g)$ then
\[
\begin{vmatrix}
  f & f_1 \\ g & g_1
\end{vmatrix}
\ \text{has}\ 
\begin{cases}
$n-1$ &  \text{roots in the upper half plane}\\
  $m$ &   \text{roots in the lower half plane}
\end{cases}
\]
\end{lemma}
\begin{proof}
  If the determinant is zero then $fg_1-gf_1=0$, so $f_1/f =
  g_1/g$. Because of algebraic interlacing, we can find non-negative
  $a_i,b_i$ so that
\[
\sum a_i\, \frac{1}{x-r_i} =\sum b_i\, \frac{1}{x-s_i} 
\]
where $\roots{(f)} = (r_i)$ and $\roots{(g)}=(b_i)$. If we evaluate
the equation at any real number the left hand side is in the upper half
plane, and the right hand side is in the lower half plane, so there
are no real roots.

Thus, as $f$ and $g$ vary among polynomials satisfying the hypothesis
the roots of the determinant never cross the real axis. We may
therefore choose $f = (x-\imag)^n$, $f_1=f'$, $g = (x+\imag)^m$ and
$g_1=g'$, in which
case
\[
\begin{vmatrix}
  f & f_1 \\ g & g_1 
\end{vmatrix}
=
\begin{vmatrix}
  (x-\imag)^n &  n(x-\imag)^{n-1} \\
  (x-\imag)^m &  m(x-\imag)^{m-1} 
\end{vmatrix}
= (x-\imag)^{n-1}(x+\imag)^{m-1}((n-m)x+(n+m))
\]
This determinant has $n-1$ roots in the upper half plane, and $m$ in
the lower half plane, so the lemma is proved.
\end{proof}

If $n=m$ then the proof shows that there are $n-1$ roots in each half plane.

  \section{Newton's inequalities in the complex plane}
  \label{sec:newton-cpx}

  If we constrain the location of the roots of a polynomial to a
  region $\mathcal{S}$ of the complex plane, then geometric properties
  of $\mathcal{S}$ lead to restrictions on the Newton quotients. It
  turns out that the inequalities only depend on the behavior of the
  Newton quotients on quadratics, so we define
\[
NQ(\mathcal{S}) = \inf_{f\in\allpolyint{\mathcal{S}}(2)} \ 
\left|\frac{a_1^2}{a_0a_2}\right|
\]
where $f(x) = a_0+a_1x+a_2x^2$. Consider some examples. The usual
argument shows that $NQ(\reals^+)=4$. Similarly, $NQ(\sigma
\reals^+)=4$ for any $\sigma\in\complexes\setminus0$. Since
$(x-\sigma)(x+\sigma)$ has Newton quotient zero, we see that
$NQ(\sigma\reals)=0$ for any $\sigma\in\complexes\setminus0$. The
Newton quotient of $(x-\sigma)^2$ is $4$, and therefore
$NQ(\mathcal{S})\le4$ for any $\mathcal{S}$. 

\begin{prop}
  Suppose that $\mathcal{S}$ is a region of $\complexes$ such that
$\mathcal{S}$ and $\mathcal{S}^{-1}$ are convex. If
  $f(x)=\sum a_ix_i$ is in $\allpolyint{\mathcal{S}}(n)$ then
  \begin{equation}
    \label{eqn:newton-cpx-1}
    \left|\frac{a_k^2}{a_{k-1}a_{k+1}}\right| \ge
    \frac{NQ(\mathcal{S})}{4}\,\left(\frac{k+1}{k}\,\cdot\,\frac{n-k+1}{n-k}\right)
  \end{equation}
\end{prop}
\begin{proof}
  We follow the usual proof of Newton's inequalities.  Write
\[ f(x) = \sum_{k=0}^n b_k \binom{n}{k}\,x^k \]
where $\binom{n}{i}b_i=a_i$. The derivative of $f$
satisfies
\[ f'(x) = n\sum_{k=1}^{n} b_k \binom{n-1}{k-1} \,x^{k-1}\]
By Gauss-Lucas, the derivative is in $\allpolyint{\mathcal{S}}$ since
$\mathcal{S}$ is convex. If we form the reverse of the polynomial,
$x^nf(1/x) = \sum \binom{n}{i}b_{n-i}x^i$ then the roots of the reverse are the
inverses of the roots of $f$, and thus lie in
$\mathcal{S}^{-1}$. Since $\mathcal{S}^{-1}$ is convex, any derivatives
of the reverse are in $\allpolyint{\mathcal{S}^{-1}}$.
Consequently, if we apply the derivative $k$ times, reverse, apply
the derivative $n-k-2$ times, and reverse again, the resulting
polynomial is in $\allpolyint{\mathcal{S}}$:
\begin{equation}\label{eqn:cpx-newton-1}
 n\cdots (n-k+1)\biggl( b_k \binom{2}{0} + b_{k+1}\binom{2}{1}x
+ b_{k+2}\binom{2}{2}x^2\biggr) \in\allpolyint{\mathcal{S}}
\end{equation}
The Newton quotient of this quadratic equation satisfies
\[
\frac{4\,b_{k+1}^2}{b_k\,b_{k+2}} \ge NQ(\mathcal{S})
\]
Substituting $b_k = a_k /\binom{n}{k}$ and simplifying yields the
conclusion.
  
\end{proof}

In general it is hard to determine $NQ(\mathcal{S})$. The next lemma
gives a good lower bound for the Newton quotient that is often enough
to determine $NQ(\mathcal{S})$.

\begin{lemma}
  If $0<\alpha<\beta<2\pi$ and $r,s>0$ then the Newton quotient for
  $(x-re^{\imag \alpha})(x-se^{\imag\beta})$ is at least
  $4\cos^2(\frac{\alpha-\beta}{2})$. 
\end{lemma}
\begin{proof}
  The Newton quotient is invariant under multiplication by any complex
  number, so we multiply by $(1/r)e^{-(\alpha+\beta)/2}$. Thus we may
  assume that $f=(x-e^{\imag \gamma})(x-te^{-\imag\gamma})$ where
  $\gamma = (\beta-\alpha)/2$. If $t>0$ then the Newton quotient is
\[
\left|\frac{{\left( -e^{\imag \,\gamma } - {t\,e^{-\imag \,\gamma }}
    \right) }^2}{t}\right| =  \frac{1}{t}+t +
2\cos(2\gamma)
\]
This is minimized for $t=1$, and has value $4\cos^2(\gamma)$.
\end{proof}

Sectors satisfy the conditions of the proposition.

\begin{cor}
  If $0<\alpha<\beta<\pi$ and $\mathcal{S}$ is the sector $\alpha \le
  \arg z \le \beta$ then
\[
\left|\frac{a_k^2}{a_{k-1}a_{k+1}}\right| \ge
\cos^2(\frac{\alpha-\beta}{2})\,\left(\frac{k+1}{k}\,\cdot\,\frac{n-k+1}{n-k}\right)
\]
\end{cor}
\begin{proof}
  If $u,v$ are any two points in the sector, then the difference of
  their arguments is at most $\beta-\alpha$. From the lemma the
  Newton quotient is minimized when $\alpha$ and $\beta$ are as
  different as possible, so the minimum is attained at $(\alpha-\beta)/2$.
\end{proof}

\begin{cor}\label{cor:newton-quadrant}
  If $f(x) = \sum a_i x^i$ is a polynomial of degree $n$ satisfying
  the condition that the real parts of all the roots have the same
  sign, and the imaginary parts of the roots all have the same sign
  then 
\[
   \left|\frac{a_k^2}{a_{k-1}a_{k+1}}\right| \ge
    \frac{1}{2}\,\left(\frac{k+1}{k}\,\cdot\,\frac{n-k+1}{n-k}\right)    
\]
\end{cor}
\begin{proof}
  The region $\mathcal{S}$ is a quadrant. Thus $\alpha-\beta=\pi/2$,
  and $4\cos^2(\frac{\alpha-\beta}{2}) =2$.

\end{proof}

\begin{cor}
  If $f(x) = \sum a_i x^i$ is a polynomial of degree $n$ whose roots
  all have the same argument then
\[
   \left|\frac{a_k^2}{a_{k-1}a_{k+1}}\right| \ge
    \left(\frac{k+1}{k}\,\cdot\,\frac{n-k+1}{n-k}\right)
\]
\end{cor}

We can give a lower bound for more general regions.

\begin{cor}
  If $\mathcal{S}$ and $\mathcal{S}^{-1}$ are convex regions, and if
  $\mathcal{S}$ is contained within a sector of angle width $\alpha$
  then for $f\in\allpolyint{\mathcal{S}}$ 
\[
   \left|\frac{a_k^2}{a_{k-1}a_{k+1}}\right| \ge
    \cos^2(\alpha/2)\,\left(\frac{k+1}{k}\,\cdot\,\frac{n-k+1}{n-k}\right)
\]

\end{cor}

For an example of such a region, consider a circle of radius 1 around
$\sigma$ where $|\sigma|>1$. Simple trigonometry shows that the cosine of one
half the angle between the two tangents to the circle is
$\sqrt{|\sigma|^2-1}/|\sigma|$. We conclude that for all $f$ with roots in this
circle we have
\[
    \left|\frac{a_k^2}{a_{k-1}a_{k+1}}\right| \ge
    \frac{|\sigma|^2-1}{|\sigma|^2}
    \, \left(\frac{k+1}{k}\,\frac{n-k+1}{n-k}\right)
\]

If we restrict $f$ to have all real coefficients then the results are
simpler. 
\begin{cor}
    Suppose that $\mathcal{S}$ and $\mathcal{S}^{-1}$ are convex
    regions of $\complexes$. If 
  $f(x)=\sum a_ix_i$ is in $\allpolyint{\mathcal{S}}(n)$ then
  \begin{equation}
    \label{eqn:newton-cpx-2}
    \left|\frac{a_k^2}{a_{k-1}a_{k+1}}\right| \ge
        \inf_{z\in\mathcal{S}\cap\overline{\mathcal{S}}}
    \cos^2(\arg(z))\, \left(\frac{k+1}{k}\,\frac{n-k+1}{n-k}\right)
  \end{equation}

\end{cor}

\begin{cor}\label{cor:newton-cpx-3}
  If $0\le\alpha\le\pi$ and $\mathcal{S}$ is one of the two sectors
  \begin{enumerate}
  \item $-\alpha \le \arg z \le \alpha$ 
  \item $\pi-\alpha \le \arg z \le \pi+\alpha$ 
  \end{enumerate}
  and $f$ is a polynomial with all real coefficients in
  $\allpolyint{\mathcal{S}}(n)$ then
\[
    \label{eqn:newton-cpx-3}
    \left|\frac{a_k^2}{a_{k-1}a_{k+1}}\right| \ge
    \cos^2(\alpha)\,\left(\frac{k+1}{k}\,\frac{n-k+1}{n-k}\right)    
\]
\end{cor}


\chapter{Equivariant Polynomials}
\label{cha:epoly}

\renewcommand{\TimeStampStart}{Friday, February 24, 2006: 08:32:23}
\mytoday  

In this chapter we  consider
polynomials that are fixed under the actions of a cyclic group.

\section{Basic definitions and properties}
\label{sec:epoly-basic}

\index{equivariant}
\index{group}
\index{cyclic}

A polynomial is \emph{equivariant} if it is fixed under the action of
a group. For polynomials in one variable the group we consider is
$\intmod{n}$, the integers mod $n$, for which we take the
multiplicative generator $\zeta = e^{2\pi i/n}$. The action of a group
element $\zeta^k$ takes a polynomial $f(z)$ to $f(\zeta^k z)$. Since
$\intmod{n}$ is cyclic we only need to verify that $f(z) =
f(\zeta z)$ in order to know that $f$ is equivariant.

When $n$ is $1$ then all polynomials are equivariant and we are left
with the theory of polynomials in $\allpoly$ and $\allpolypm$. In
general, we want to consider equivariant polynomials that have the
largest possible number of real roots. Unfortunately, this is not
quite the definition because we must treat zero roots specially: 

\begin{definition}
  The set $\epoly{n}$ consists of all polynomials with real
  coefficients fixed under the action of $\intmod{n}$ such that
  if we write $f(z)= z^s g(z)$ where $g(0)\ne0$ then $g$ has degree
  $nk$ and $k$ real roots. We also define (where $g$ is as in the
  preceding sentence)
  \begin{align*}
    \epolypos{n} & = \{ f\in\epoly{n}\vert \text{ all real roots of
      $g$ are positive}\} \\
    \epolyneg{n} & = \{ f\in\epoly{n}\vert \text{ all real roots of
      $g$ are negative}\} \\
  \end{align*}
\end{definition}

We should make a few observations about this definition. First, if
$f(z)$ is an equivariant polynomial with a term $cz^k$ then $f(z) =
f(\zeta z)$ implies that $cz^k = c(\zeta z)^k$ which requires that
$\zeta^k=1$. Since $\zeta$ is a primitive $n$-th root of unity, the
non-zero terms of an equivariant polynomial have degree divisible by
$n$. This shows that $s$ in the definition is also divisible by $n$.
Moreover, if $g$ in the definition has a real root $r$ then it also
has roots $r\zeta, r\zeta^2,\dots,r\zeta^{n-1}$. If $n$ is odd then
these are all complex.  Thus when $n$ is odd the number of non-zero
real roots of any equivariant polynomial is at most the degree of $g$
divided by $n$. If $n$ is even then the number of \emph{positive}
roots is at most the degree divided by $n$. In addition, if $n$ is
even then $g$  has positive and negative roots, so both 
$\epolypos{n}$ and $\epolyneg{n}$ are empty.

From the remarks above we see that if monic $f\in\epoly{n}$ has real roots
(or positive roots if $n$ is even) $r_1,\dots,r_k$ then we can write
\begin{align*} f(z) &= z^{ns} \prod_k (z-r_k)(z-\zeta r_k)\cdots
  (z-\zeta^{n-1}r_k)\\ 
& =  \tilde{f}(z^n) \\
\intertext{where}
\tilde{f}(x) &= z^s\prod_k (x - r_k^n)
\end{align*}

This gives us a useful representation for polynomials in $\epoly{n}$.
A polynomial $f(z)$ is in $\epoly{n}$ iff there is a polynomial
$\tilde{f}\in\allpoly$ (or $\allpolyaltclose$ if $n$ even) such that
$f(z) = \tilde{f}(z^n)$. Since this correspondence is linear, we have a 
linear transformation
$$ \widetilde{}:\epoly{n}\longrightarrow
\begin{cases}
\allpolyaltclose & \text{ $n$  even} \\
\allpoly & \text{ $n$ odd} 
\end{cases}
$$
\noindent%
that is a bijection. 

\begin{example}
\index{Hermite polynomials}
\index{polynomials!Hermite}
  Consider the case $n=2$. Since $\epoly{2}=\{g(x^2)\mid
  g\in\allpolyaltclose\}$ if $g$ has non-zero roots $r_1,\dots,r_m$
  then $g(x^2)$ has non-zero roots $\pm\sqrt{r_i}$ and hence
  $\epoly{2}\subset\allpoly$. The Hermite polynomials $H_n$ are
  invariant under $Z_2$. Corollary~\ref{cor:hermite} shows that the map
  $x^k\mapsto H_k(x)$ defines a linear transformation
  $\epoly{2}\longrightarrow\epoly{2}$. \index{Hermite polynomials} 
  
  If $n>2$ then $\epoly{n}$ is not contained in $\allpoly$.
\end{example}
  
  It is at first surprising that the derivative does not preserve
  $\epoly{n}$:

\begin{lemma} \label{lem:epoly-deriv}
  If $n$ is odd and $f\in\epoly{n}$ then $f^\prime\not\in\epoly{n}$ but
  $f^{(n)}\in\epoly{n}$. For any $k\ge0$ we  have that $z^k
  f^{(k)}(z)$ is in $\epoly{n}$.
\end{lemma}

\begin{proof}
  From the definition
  $$
  f^\prime(\zeta z) = \lim_{\zeta\rightarrow0} \frac{f(\zeta (z+h)) -
    f(\zeta z)}{\zeta h} = \zeta^{-1} f^\prime(z)$$ This shows that
  $f(z)$ is not equivariant, but that $zf^\prime(z)$ is equivariant.
  Similarly $f^{(k)}(\zeta z) = \zeta^{-k} f^{(k)}(z)$ so that $z^k
  f^{(k)}(z)$ is equivariant. It remains to show that $f^\prime$ has
  the right number of real roots.
  
  \index{Rolle's Theorem} From Rolle's theorem, if a polynomial $p$
  has $s$ real roots, then the derivative has at least $s-1$ real
  roots, and they interlace the real roots of $p$. It might be the
  case that $p^\prime$ has more roots than $p$, but we always have
  that if $p$ has $a$ roots, and $p^\prime$ has $b$ roots then
  $a\equiv b+1\pmod{2}$.
  
  Without loss of generality we may assume that $\tilde{f}$ has all
  distinct roots, all non-zero. Assume that $n$ is odd.  If $f$ has
  degree $nm$ then $\tilde{f}$ and $f$ have $m$ non-zero real roots.
  The derivative $f^\prime$ has either $m-1$ or $m-2$ non-zero real
  roots, and an $n-1$ fold root at $0$. The total number of roots of
  $f^\prime$ has opposite parity from $m$, so $f^\prime$ has exactly
  $m-1$ non-zero roots. Since $zf^\prime$ is equivariant, we can write
  $zf^\prime(z) = g(z^n)$ where $g\in\allpoly$ and hence
  $zf^\prime(z)\in\epoly{n}$.
  
  Next, $f^{(2)}$ has at least $m-2$ non-zero roots, and an $n-2$ fold
  root at $0$. The same argument as above shows that $z^2 f^{(2)}$ has
  an $n$ fold root at $0$, and exactly $m-1$ non-zero roots.
  Consequently $z^2f^{(2)}$ is in $\epoly{n}$. Continuing, $f^{(i)}$
  has an $n-i$ fold root at $0$ and exactly $m-1$ non-zero roots. As
  above, $z^i f^{(i)}$ is in $\epoly{n}$.

Since $z^nf^{(n)}$ is in $\epoly{n}$ we can divide by $z^n$ and thus
$f^{(n)}$ has exactly $m-1$ non-zero roots, and is also in
$\epoly{n}$.
\end{proof}

Notice that although there is a  bijection between $\epoly{n}$ and
$\allpoly$ some results are more natural in $\epoly{n}$. For example, 
if $f(x)=g(x^5)$ then $f^{(5)}$ is in $\epoly{5}$. In terms of $g$
this says that if  $g\in\allpoly$ then 
$$
120\,g^\prime(x) + 15000\,x\,g^{\prime\prime}(x) + 
  45000\,{x^2}\,g^{(3)}(x) + 
  25000\,{x^3}\,g^{(4)}(x) + 
  3125\,{x^4}\,g^{(5)}(x)
$$
is also in $\allpoly$.

There is a general principle at work here: properties 
of $\allpoly$ extend to $\epoly{n}$, but we often have to replace $x$
with $z^n$, or $f^\prime$ with $f^{(n)}$.

\begin{cor} \label{cor:epoly-binom}
  If $n$ is odd and $k$ is a positive integer then the linear
  transformation $T\colon{}x^j\mapsto \falling{jn}{k}x^j$ maps $\allpoly$ to itself.
  In addition, $x^j\mapsto \binom{jn}{k}x^j$ also maps $\allpoly$ to
  itself. If $n$ is even the domain of the maps is $\allpolyalt$.
\end{cor}
\begin{proof} 
  The map $T$ is the composition

\centerline{
\xymatrix{
\allpoly \ar@{->}[r]^{\tilde{}\star} & \epoly{n} \ar@{->}[rr]^{z^kf^{(k)}} 
&& \epoly{n} \ar@{->}[r]^{\tilde{}} & \allpoly
}}
\noindent%
where $\tilde{}\star$ is the map $g(z)\mapsto g(z^n)$. To verify this
claim, we calculate the image of a monomial $z^j$:

\centerline{
\xymatrix{
z^j \ar@{|->}[r]^{\tilde{}\star} & z^{jn} \ar@{|->}[rr]^{z^kf^{(k)}} 
&& z^k\falling{jn}{k} z^{jn-k} \ar@{|->}[r]^{\tilde{}} & \falling{jn}{k} z^j
}}

To derive the second part, divide by the constant $k!$. 
\end{proof} 

\section{Interlacing}
\label{sec:epoly-interlacing}

We make the usual definition of linear interlacing:
\begin{definition}
  Two polynomials $f,g\in\epoly{n}$ \emph{interlace} if and only if $f+\alpha 
  g$ is in $\epoly{n}$ for all real $\alpha$. 
\end{definition}

It is very easy to determine interlacing in $\epoly{n}$. 

\begin{lemma}
  Suppose $f,g$ are in $\epoly{n}$. If $f$ and $g$ interlace then
  $\tilde{f}$ and $\tilde{g}$ interlace. The converse is true if $n$
  is odd.
\end{lemma}
\begin{proof}
  Since the map $f\mapsto \tilde{f}$ is linear it follows that if
  $f,g$ interlace then $\tilde{f}$, $\tilde{g}$ interlace. Conversely,
  if $\tilde{f}$ and $\tilde{g}$ interlace then
  $\tilde{f}+\alpha\tilde{g}$ is in $\allpoly$, and hence $f+\alpha g
  = (\tilde{f}+\alpha\tilde{g})(z^n)$ is in $\epoly{n}$ if $n$ is odd.
  \end{proof}
  
  If $n$ is even then the converse is not true. For instance, if
  $\tilde{f}=x-1$ and $\tilde{g}=x-2$ then $\tilde{f}(x^2)$ and
  $\tilde{g}(x^2)$ do not interlace. However, the \emph{positive} roots of
  $\tilde{f}(x^2)$ and $\tilde{g}(x^2)$ interlace.

\begin{cor}
  If $f,g\in\epoly{n}$ and $f,g$ interlace then the degrees of $f$ and
  $g$ differ by $-n$, $0$, or $n$.
\end{cor}
\begin{proof}
  The degrees of $\tilde{f}$ and $\tilde{g}$ differ by at most one.
\end{proof}

\begin{cor}
  If $f\in\epolypos{n}$  then 
$ f \greateqeq zf^\prime \greateqeq z^2 f^{(2)} \greateqeq \cdots $.
\end{cor}
\begin{proof}
  If we write $z^i f^{(i)} = g_i(z^n)$ then the proof of
  Lemma~\ref{lem:epoly-deriv} shows that $g_i \greateqeq g_{i+1}$ and the
  Corollary follows.
\end{proof}

Every linear transformation on $\allpoly$ (or $\allpolyalt$ is $n$
even) induces a linear transformation on $\epoly{n}$, and vice versa.
So suppose that $T\colon{}\allpoly \longrightarrow \allpoly$. We define a
linear transformation $T_\ast$ by the composition

\centerline{
\xymatrix{
\epoly{n} \ar@{->}[r]^{\tilde{}}  & \allpoly \ar@{->}[d]^T \\
\epoly{n} \ar@{.>}[u]^{T_\ast}     &  \allpoly \ar@{->}[l]^{\tilde{}_\ast}
}}

We can also use the diagram to define $T$ if we are given $T_\ast$, so 
linear transformations on $\epoly{n}$ are not very interesting.  For
instance, we have 

\begin{lemma}
  If $f,g\in\epoly{n}$ then the Hadamard product $f\ast g$ is in $\epoly{n}$.
\end{lemma}
\begin{proof}
  Observe that $(f\ast g)(z) = (\tilde{f} \ast \tilde{g})(z^n)$.
\end{proof}

\section{Creating equivariant polynomials and the Hurwitz theorem}
\label{sec:epoly-average}

\index{Hurwitz's theorem!for $\intmod{d}$}

Whenever a group acts on functions we can construct equivariant
functions by averaging over the group.  In our case, we define the
average of a polynomial over $\intmod{n}$ to be 

\begin{equation}
  \label{eqn:average}
  \hat{f}(z) = \frac{1}{n} \sum _{k=0}^{n-1} f(\zeta^k z)
\end{equation}

The average is obviously equivariant, and has a very simple
interpretation in terms of coefficients. If $f = \sum a_k x^k$ then
\begin{align}
  \hat{f}(z) &= \frac{1}{n} \sum_k a_k \left( \sum_{j=0}^{n-1}
    \zeta^{jk}\right)z^k \notag \\
\intertext{and since $\zeta$ is a root of unity}
\sum_{j=0}^{n-1} \zeta^{jk} &= \begin{cases} 0 & k \nmid n \\ n & k
  \mid n 
\end{cases} \notag \\
\intertext{This gives us our fundamental formula:}
\hat{f}(z) &= \sum_{k\equiv0\pmod{n}} a_k x^k \label{eqn:average-coef}
\end{align}

\index{even part} Since the even part of a polynomial is just the
average of $f$ over $\intmod{2}$ the Hurwitz theorem Theorem~\ref{thm:hurwitz}
\index{Hurwitz's theorem!for $\intmod{2}$} can be restated:
\begin{quote}
  The average over $\intmod{2}$ of a polynomial in $\allpolypos$ is
  in $\allpolyalt$.
\end{quote}
More generally we have

\begin{theorem} \label{thm:epoly-hurwitz}
  If $n$ is odd and $f\in\allpolyalt$ then
  $\hat{f}\in\epolypos{n}$. 
\end{theorem}
\begin{proof}
  A stronger result has been proved earlier. Here we give an alternate
  analytic proof of a part of Theorem~\ref{thm:hurwitz-gen}.  The key idea is to
  apply \eqref{eqn:fxdg} since we can also express $\hat{f}$ as
  $g(x\diffd)f$ where
$$g(z) = \frac{1}{n}\frac{\sin (\pi z)}{\sin (\pi z/n)}.$$
In order to
see this, notice that the zeros of $g$ are at the integers not
divisible by $n$, and if $k$ is an integer then $g(kn)=1$. (If $n$ is
even then $g(kn)=(-1)^k$.) 

If $f\in\allpolyalt(m)$ then $g(z)$ has $m - \lfloor m/n\rfloor$ zeros
in $(0,m)$. By Corollary~\ref{cor:agxgp} it follows that $g(x\diffd)f$ has at least
$\lfloor m/n\rfloor$ zeros. However, if $f = \sum a_i x^i$ then
$$ g(x\diffd)f = \sum g(i)a_ix^i = \sum_{n\mid i} a_i x^i = \hat{f}$$
Now $\hat{f}$ has degree $n\lfloor m/n \rfloor$ and hence has at most
$\lfloor m/n\rfloor$ real roots. Consequently $\hat{f}$ has exactly
$\lfloor m/n\rfloor$ real roots, and is in $\epolypos{n}$.
\end{proof}

\begin{cor}
  If $f=\sum a_kx^k$ is in $\allpolyalt$ (or $\allpolyaltf$) then $\sum 
  a_{nk}x^k$ is in $\allpolyalt$ (or $\allpolyaltf$).
\end{cor} 

\begin{cor}
  If $n$ is odd then $\sum_{k=0}^\infty \frac{z^k}{(nk)!}$ is in
  $\allpolyaltf$. 
\end{cor}

\index{positive interlacing}
\begin{cor}
  If $f\lessless g$ in $\allpoly$ then $\hat{f}\plesslesseq \hat{g}$.
\end{cor}
\begin{proof}
  The average is a linear transformation
  $\allpolyalt\longrightarrow\epoly{n}$. 
\end{proof}

\section{Hermite polynomials}
\label{sec:epoly-hermite}
\index{Hermite polynomials}
\index{Hermite polynomials!equivariant}
\index{recurrence!Hermite}

The Hermite polynomials have analogs in $\epoly{n}$.  If we replace
$x^2$ by $x^n$ in \eqref{eqn:hermite-4} we get polynomials $H_k^n$ that are
nearly equivariant. 

\begin{align*}
 H_k^n(z) &= (-1)^k e^{z^n} \left( \frac{d}{dz}\right)^k e^{-z^n} 
\end{align*}

From the definition we find the recurrence

\begin{align*}
  H_{k+1}^n(z)  &= nz^{n-1}H_k^n(z) - H_k^n(z)^\prime \\
  \intertext{which shows that the degree of $H_k^n$ is $k(n-1)$.  If we write
    $z^kH_k^n(z) = g_k^n(z^n)$ then the degree of $g^n_k$ is $k$ and } 
g_{k+1}^n(x) &= (nx+k)g_k^n(x) -  nxg^\prime(x)
\end{align*}

From the recurrence we see that $H_k^n$ satisfies $H_k^n(\zeta z) =
\zeta^{-k}H_k^n(z)$. The recurrence also shows that all $g_k^n$ are in
$\allpolyalt$, and hence

\begin{lemma}
  The  polynomials $z^kH_k^n$ are in $\epolypos{n}$ (or
  $\epoly{n}$ if $n$ is even). In particular, $H_{nk}^n(z)$ is in
  $\epoly{n}$. 
\end{lemma}

We can also express  $H_k^n$ as a composition:

\begin{align*}
  H_k^n(z) &= (nz^{n-1} - \diffd)^k\,(1) 
\end{align*}


\chapter{Polynomials modulo an ideal }
\label{cha:ideal}

\renewcommand{\TimeStampStart}{Monday, January 07, 2008: 11:44:59}
\mytoday

  If $\mathcal{S}$ is a set of polynomials in variables $\xx$, and
  $\mathcal{I}$ is an ideal in $\reals[\xx]$ (or $\complexes[\xx]$)
  then in this chapter we study $\init{\mathcal{S}}{\mathcal{I}}$. There are
  two cases that we are particularly interested in:
  \begin{enumerate}
  \item $\init{\allpoly}{\{x^n\}}$ 
  \item $\init{\rup{2}\,}{\{y^n\}}$ 
  \end{enumerate}
We want to find intrinsic descriptions of
$\init{\mathcal{S}}{\mathcal{I}}$ as well as to resolve the usual questions of
interlacing, coefficient inequalities, and linear transformations.

\section{Initial segments of $\allpoly$}
\label{sec:init-sequ-allp}

An element of $\init{\allpoly}{\{x^n\}}$ can be viewed as a polynomial
$f = \sum_{k=0}^{n-1} a_ix^i$ such that there is an extension to
$\allpoly$. By this we mean a polynomial $\tilde{f}$ in $\allpoly$
such that the first $n$ coefficients of $\tilde{f}$ and $f$ are equal.
For instance, $1+7x+28x^2 \in\init{\allpoly}{\{x^3\}}$ since there is
a polynomial, namely $(1+x)^7$, that begins $1+7x + 28x^2 + \cdots$.
We sometimes say that elements of $\init{\allpoly}{\{x^{n+1}\}}$ are
$n$-initial segments of $\allpoly$. Note that
$\allpoly(n-1)\subset\init{\allpoly}{\{x^n\}}$ since
$f\in\allpoly(n-1)$ is its own extension.

Although $\init{\allpoly}{\{x^2\}} = \{a+b x\mid a,b\in\reals\}$ is
trivial, even the quadratic case for $\allpolypos$ takes some effort.

\begin{lemma}
  $\init{\allpolypos}{\{x^3\}} = \left\{a+bx+cx^2\,\mid\, a,b,c>0 \And b^2> 2ac\right\}$
\end{lemma}
\begin{proof}
  If $f=\sum a_ix^i\in\allpolypos(n)$ then Newton's inequalities yield
\[   \frac{a_1^2 }{a_0a_2} \ge \frac{2}{1}\frac{n}{n-1} > 2
\]
Conversely, assume that $b^2>2ac$. If we can find
$f=a+bx+cx^2+\cdots\in\allpolypos$ then if we  divide by $a$, and
replace $x$ by $xb/a$ we can find $1+x + Cx^2+\dots\in\allpoly$ where
$0<C<1/2$. It therefore suffices to show that  if $0<C<1/2$ then we
can find $f=1+x+Cx^2+\cdots \in\allpolypos$.

The quadratic equation shows that
\[
1 + x + tx^2 \in\allpoly \text{ for all } 0<t<1/4
\]
Squaring yields that 
\[
(1+x+tx^2)^2 = 1+2x + (2t+1)x^2 + \cdots \in\allpoly \text{ for all }
0<t<1/4 \]
and replacing $2x$ by $x$ yields
\[
1+x + (1/4+ t/2)x^2 + \cdots \in\allpoly \text{ for all }
0<t<1/4 \]
If we let $s = 1/4+t/2$ we have
\[
1+x + sx^2 + \cdots \in\allpoly \text{ for all }
0<s<3/8 \]
Continuing by induction, if we expand $(1+x+tx^2)^{2^n}$, we get that

\[
1+x + sx^2+ \cdots \in\allpoly \text{ for all }
0<s<1/2 - 2^{-n} \] 
Consequently, if we choose $n$ so that $2^{-n}<1/2-C$ then we can realize
$1+x+Cx^2$ as a $2$-initial segment.
\end{proof}

The argument in the lemma shows that powers of $1+x+ax^2$ converge,
after renormalization, to $1+x+x^2/2$. Here's a more general result
that shows why the exponential function appears.

\begin{lemma}\label{lem:init-exp}
  If $f = 1+ a_1x+\cdots+a_{n-1}x^{n-1} + x^n$ then
\[
\lim_{n\rightarrow\infty} f(x/n)^n = e^{a_1x}
\]
\end{lemma}
\begin{proof}
  If we write $f = \prod(x+r_i)$ then $\prod r_i=1 $ and $\sum
  1/r_i=a_1$. 
  \begin{align*}
    \lim_{m\rightarrow\infty} f(x/m)^m &= \lim_{x\rightarrow\infty}
    \prod(x/m+r_i)^m \\
    & = \prod_1^n
    \lim_{m\rightarrow\infty}\left(\frac{x}{mr_i}+1\right)^m 
   = \prod_1^n e^{x/r_i} = e^{a_1x}
  \end{align*}
\end{proof}

\index{exponential function}
Not only do the powers converge to an exponential function, but the
coefficients are bounded by the coefficients of $e^x$. In addition,
Lemma~\ref{lem:init-exp} shows that there are polynomials in
$\init{\allpolypos}{\{x^{n+1}\}}$ whose coefficients are arbitrarily
close to $1/i!$.

\begin{lemma}
  If $1+x+a_2x+\cdots+a_nx^n\in\allpolypos$ then $0<a_i<1/i!$.
\end{lemma}
\begin{proof}
  Newton's inequalities give that $1/2 > a_2$. Assuming that the Lemma
  is true up to $n-1$, we see that
  \begin{align*}
    a_{n-1}^2 & > \frac{n}{n-1}a_{n-2}a_n &\text{Newton}\\
    a_n & < \frac{1}{(n-1)!^2}\frac{n-1}{n}\frac{1}{a_{n-2}}  &
      \text{induction}\\
      a_n & < 1/n! & \text{induction}
  \end{align*}
\end{proof}

If $f = b_0 + b_1x + \cdots \in\init{\allpolypos}{\{x^{n+1}\}}$
then
\[
\frac{1}{b_0}\,f\,\bigl(\frac{b_0\,x}{b_1}\bigr) = 1 + x + \cdots
\]
is also in $\init{\allpolypos}{\{x^{n+1}\}}$ Thus, in order to
describe $\init{\allpolypos}{\{x^{n+1}\}}$ it suffices to determine
which polynomials $1+x+a_2x^2+\cdots+a_nx^n$ are in
$\init{\allpolypos}{\{x^{n+1}\}}$. The lemma shows that the set of all
such $(a_2,\cdots,a_n)$ is bounded. 

If $1+x+a_2x^2+a_3x^3\in\init{\allpolypos}{\{x^{4}\}}$ then we know
several bounds:
\begin{align*}
  1/2 &> a_2 > 0 \\
  1/6 &> a_3 > 0 \\
  a_2^2 & > (3/2) a_3 & Newton\\
  1+a_3 & \ge 2a_2 & \eqref{eqn:init-newton}
\end{align*}
Empirical evidence suggests that these are also sufficient.

\section{Initial segments of $\rupint{2}$}
\label{sec:init-segm-rupint{2}}

We now introduce $\init{\gsubplus_2\,}{\{y^{n+1}\}}$. Its elements can be
considered to be sums $f_0(x)+\cdots + f_n(x)y^n$ arising from a
polynomial in $\rupint{2}$. 
Here are two examples of elements in $\init{\rupint{2}\,}{\{y^n\}}$.

\begin{itemize}
\item Let $f\in\allpolypos$.
  The Taylor series \index{Taylor series} for $f(x+y)$ shows
\[
f+ y\,\frac{f^\prime}{1!}+y^2\,\frac{f^{\prime\prime}}{2!}+\dots +y^{n}\frac{f^{(n)}}{n!}
  \quad\in \init{\rupint{2}\,}{\{y^{n+1}\}}
\]
\item If we choose $f = \prod(x+a_i+b_iy)$ and set $g=\prod(x+a_i)$ then 
\[  g +y\, \sum b_j \frac{g}{x+a_j} 
+y^2\,\sum_{j\ne k} b_j b_k \frac{g}{(x+a_j)(x+a_k)}
\quad\in \init{\rupint{2}\,}{\{y^3\}}
\]
\end{itemize}

The case $n=2$ is very familiar.
\begin{align*}
  \init{\rupint{2}\,}{\{y^2\}} &= \left\{ f+y\,g\,\mid\, f\lesslesseq g \text{
      in } \allpoly\right\}\\
  \init{\gsubplus_2\,}{\{y^2\}} &= \left\{ f+y\,g \,\mid\, f\lesslesseq g \text{
      in } \allpolypos\right\}
\end{align*}

We think that the characterization for $\init{\gsubplus_2\,}{\{y^3\}}$
is similar to the result for $\init{\allpolypos}{\{x^3\}}$. Namely,
$f_0+y\,f_1+y^2\,f_2$ is a $2$-initial segment if and only if
$f_0\lesslesseq f_1 \lesslesseq f_2$ and $f_1^2-2f_0f_2\ge0$ for
positive $x$.  Necessity is easy since we just substitute for $\alpha$
and apply the result for $\init{\allpolypos}{\{x\}}$. The problem is
existence.

\section{Coefficient inequalities}
\label{sec:coeff-ineq-init}

There are Newton's inequalities for $\init{\allpolypos}{\{x^n\}}$. If
we choose $f=\sum_0^{n-1}a_ix^i$ in $\init{\allpolypos}{\{x^n\}}$ and $f$ has an
extension of degree $r$ then Newton's inequalities tell us that
\begin{align*}
  \frac{a_{k+1}^2}{a_ka_{k+2}} & \ge \frac{k+1}{k}\cdot\frac{r-k+1}{r-k} &
  0 \le k < n-2\\
\intertext{Since we don't know $r$ all we can say is that}
  \frac{a_{k+1}^2}{a_ka_{k+2}} & > \frac{k+1}{k} &
  0 \le k < n-2\\
\end{align*}

In addition to these Newton inequalities, we have determinant
inequalities for $\init{\allpolypos}{\{x^{n}\}}$. Any subdeterminant of
\eqref{eqn:tp} that only involves $a_0$ up to $a_{n-1}$ must be
non-negative. For instance, if
$a_0+a_1x+a_2x^2+a_3x^3\in\init{\allpolypos}{\{x^4\}}$  then 
\begin{align}
\begin{vmatrix}
  a_1 & a_2 & a_3 \\ a_0 & a_1 & a_2 \\ 0 & a_0 & a_1
\end{vmatrix} & \ge 0 \notag\\
\intertext{which gives the inequality}
a_1^3 + a_0^2a_3 &\ge 2a_0a_1a_2 \label{eqn:init-newton}
\end{align}

If we are given $a_0+\cdots +a_nx^n\in\init{\allpolypos}{\{x^{n+1}\}}$ then there
is an $f\in\allpolypos$ with $f=\sum a_ix_i$, but we don't know the
degree $m$ of $f$. We can also use Newton's inequalities to get lower bounds on
$m$. First of all, 
\begin{align*}
  \frac{a_k^2}{a_{k-1}a_{k+1}} & \ge
  \frac{k+1}{k}\left(1+\frac{1}{m-k}\right) & \text{for $1\le k < n$}\\
\intertext{Solving for $m$ yields}
m & \ge k +\frac{m+1}{m}\left(\frac{a_k^2}{a_{k-1}a_{k+1}}-\frac{k+1}{k}\right)^{-1}
\end{align*}
These inequalities hold for $1\le k <n$. For example, if $k=1$ the
inequality is
\begin{equation}\label{eqn:ideal-bound}
m \ge 1+ 2\left(\frac{a_1^2}{a_0a_2}-2\right)^{-1}
\end{equation}
We see $m$ increases as the Newton quotient gets closer to its minimum
value. It's not surprising that the exponential series is the unique
series where the Newton quotient equals its lower bound for all
indices. 

For an example, we claim that the smallest extension of $1+nx +
\binom{x}{2}x^2$ has degree $n$. Indeed, if we compute the bound
\eqref{eqn:ideal-bound} we get exactly $n$.

\begin{remark}
  $\init{\allpoly}{\{x^n\}}$ isn't closed. If we have a sequence
  $f_k\in \init{\allpoly}{\{x^n\}}$ such that $f_k\rightarrow f$ then
  the extensions might not converge to a polynomial. For instance,
  \eqref{eqn:ideal-bound} shows that
$1+x+x^2/2\not\in\init{\allpoly}{\{x^3\}}$. It is, however, the limit
of $(1+x/n)^n$ in $\init{\allpoly}{\{x^3\}}$. 
\end{remark}

\section{There are many extensions }
\label{sec:there-are-many}

We  show that polynomials in $\allpolypos$ and $\gsubplus_2$ have
extensions of arbitrary degree.

\begin{lemma}
  If $f\in\init{\allpolypos}{\{x^n\}}$ has an extension of degree $m$,
  then $f$ has  extensions of degree $m+k$ for $k=1,2,\dots$.
\end{lemma}
\begin{proof}
  We observed earlier \mypage{sec:rapid} that if $f\in\allpoly(n)$
  then there is a positive $\alpha$ such that $f + \alpha
  x^{n+1}\in\allpoly$. Using this fact, we establish the lemma by
  induction. 
\end{proof}
\begin{lemma}
  If $f\in\init{\gsubplus_2\,}{\{y^n\}}$ has an extension of degree $m$,
  then $f$ has extensions of degree $m+k$ for $k=1,2,\dots$. 
\end{lemma}
\begin{proof}
  If $f\in\gsubplus_2(n)$, and $g(x,y)$ consists of the terms of
  degree $n$ then we will show that there is an $\epsilon$ such that
\[
 F_\epsilon(x,y) = f(x,y) + \epsilon(x+y)g(x,y) 
\]
is in $\gsubplus_2(n+1)$. Since the asymptotes of the graph of $f$ are
also asymptotes of the graph of $F_\epsilon$ for any $\epsilon$, we
see that there is an $N$ such that if $|\alpha|>N$ then
$F_\epsilon(\alpha,y)\in\allpoly$. Also, for any $\alpha$ we can, by
the observation in the proof above, find an $\epsilon$ such that
$F_\epsilon(\alpha,y)\in\allpoly$. Since $[-N,N]$ is compact, we can
find an $\epsilon$ that works for all points in this interval, and
hence it works for all of $\reals$. 
\end{proof}

A polynomial in $\allpoly(n)$ determines an element in
$\init{\allpoly}{\{x^r\}}$ for $1\le r <n$. However, it is not true
that an element of $\allpolyf$ determines a member of
$\init{\allpoly}{\{x^r\}}$. For instance, $e^x\in\allpolyf$, but it is
known \cite{zemyan} that no $k$-initial segment is in $\allpoly$
except for $k=1$. Corollary~\ref{cor:kurtz-4} is an example of the
other extreme. The polynomial
\[\sum_1^\infty x^k 2^{-{k^2}} \]
has the property that \emph{all} initial segments are in
$\allpolypos$.

\section{Interlacing}
\label{sec:interlacing-1}

\index{interlacing!for ideals}
We can define interlacing as usual. We say that $f$ and $g$ interlace,
written $f\thicksim g$, in $\init{\allpoly}{\mathcal{I}}$ if $f+\alpha
g\in\init{\allpoly}{\mathcal{I}}$ for all $\alpha\in\reals$.  If we
restrict $\alpha$ to positive values we write $f\poslace g$.

Here is a simple result.

\begin{lemma}
  Suppose $ a+bx+cx^2$ and $A+Bx+Cx^2$ are in $\init{\allpolypos}{\{x^3\}}$.
  Then,  $ a+bx+cx^2\poslace A+Bx+Cx^2$ if and only if $bB \ge aC+cA$.
\end{lemma}
\begin{proof}
  $ a+bx + cx^2+ \alpha (A+Bx+Cx^2)$ is in $\init{\allpolypos}{\{x^3\}}$ if and
  only if 
\[
(b+\alpha B)^2 > 2(a+\alpha A)(c+\alpha C)
\]
and this holds for all positive $\alpha$ if and only if $bB\ge aC+Ac$ holds.
\end{proof}

It is not easy to show that two polynomials do or don't interlace.
For example we show that $x$ and $(x+1)^2$ do not interlace in
$\init{\allpoly}{\{x^3\}}$. If they did then $(x+1)^2 + \alpha x\in
\init{\allpoly}{\{x^3\}}$ for all $\alpha$. Newton's inequality says
that
\[ (2+\alpha)^2 \ge 2 \]
and this doesn't hold for all $\alpha$, so they don't interlace. It
\emph{is} the case that for any positive $t$
\[
x \poslace (x+t)^2\ \text{in}\ \init{\allpoly}{\{x^3\}}
\]
since $(2t+\alpha)^2>2t^2$ for $\alpha\ge0$.

\index{Newton's inequalities}
Since there is no concept of degree in $\init{\allpoly}{\{x^{n+1}\}}$
there is no idea of the direction of interlacing, so we can not expect
simple inequalities between interlacing polynomials. We do have
\begin{lemma}
  Suppose $\sum a_ix^i$ and $\sum b_ix^i$ are interlacing polynomials
  in $\init{\allpoly}{\{x^{n+1}\}}$. The following inequality holds
  for $0\le k <n-1$:
\begin{multline*}
\biggl( (k+2)a_kb_{k+2} + (k+2)a_{k+2}b_k - 2 (k+1)\bigr)
a_{k+1}b_{k+1}\biggr)^2 \\
\le 4(k+1)^2
\biggl(a_{k+1}^2 - \frac{k+2}{k+1}a_ka_{k+2}\biggr)
\biggl(b_{k+1}^2 - \frac{k+2}{k+1}b_kb_{k+2}\biggr)
\end{multline*}
\end{lemma}
\begin{proof}
  Apply Newton's inequality to $\sum (a_k+\alpha b_k)x^k$. The result is
  a quadratic in $\alpha$; the inequality is the statement that the
  discriminant is non-positive.
\end{proof}

\section{Linear transformations}
\label{sec:line-transf}

A linear transformation $\allpoly\longrightarrow\allpoly$ does not in
general determine a linear transformation 
$\init{\allpoly}{\{x^n\}}\longrightarrow\init{\allpoly}{\{x^n\}}$
since $T$ applied to different extensions can give different
results. We will see that Hadamard products do preserve
$\init{\allpoly}{\{x^n\}}$, and there are some other interesting
preserving transformations arising from $\gsubplus_3$. 

We recall an elementary result about ideals.

\begin{lemma}
  If $\mathcal{I},\mathcal{J}$ are ideals of $\reals[x]$ then a linear
  transformation $T$ on $\reals[x]$ such that $T(\mathcal{I})\subset\mathcal{J}$
  determines a map $T^\ast:
  \init{\allpoly}{\mathcal{I}}\longrightarrow
  \init{\allpoly}{\mathcal{J}}$. 
\end{lemma}
\begin{proof}
  The point is that the map involves a choice, and the hypothesis
  makes the map well defined. We define $T^\ast$ by

  \centerline{\xymatrix{
      \allpoly
      \ar@{<-}[d]_{\text{choose extension} }           
      \ar@{->}[rrr]^{T }         
      &&&
      \allpoly
      \ar@{->}[d]^{\text{project} } \\        
            {\init{\allpoly}{\mathcal{I}}}
      \ar@{.>}[rrr]^{T^\ast }         
      &&&
            {\init{\allpoly}{\mathcal{J}}}
}}

\end{proof}

For example, the derivative map $\diffd:f\mapsto f'$ maps the ideal generated
by $x^n$ to the ideal generated by $x^{n-1}$. We therefore have
\[
\diffd\colon \init{\allpoly}{\{x^n\}}\longrightarrow
\init{\allpoly}{\{x^{n-1}\}}
\]
More generally, we have
\begin{lemma}
  If $f\in\allpoly(r)$ then
$
f(\diffd)\colon \init{\allpoly}{\{x^n\}}\longrightarrow
\init{\allpoly}{\{x^{n-r}\}}
$.
\end{lemma}

$\allpoly$ is closed under products, and 
$\mathcal{I}\times\mathcal{I} \subset\mathcal{I}$  so the product induces a
linear transformation 
\[
\init{\allpoly}{\mathcal{I}}\times\init{\allpoly}{\mathcal{I}}
\longrightarrow \init{\allpoly}{\mathcal{I}}
\]

If $f = \sum
a_ix^i$ and $g = \sum b_i x^i$ then the matrix representation of 
\[
\init{\allpoly}{\{x^n\}}\times\init{\allpoly}{\{x^n\}}
\longrightarrow \init{\allpoly}{\{x^n\}}
\] is given by

\[
\begin{pmatrix}
  b_0\\b_1 \\\vdots\\b_{n-1}
\end{pmatrix} \mapsto
\begin{pmatrix}
  a_0 \\ a_1 & a_0 \\ \vdots & \vdots & \ddots \\
a_{n-1} & a_{k-2} & \dots & a_0
\end{pmatrix}
\begin{pmatrix}
  b_0\\b_1 \\\vdots\\b_{n-1}
\end{pmatrix} 
\]

We get more interesting matrices if we use two variables. Consider
\[
\gsubplusclose_2\times\allpoly \xrightarrow{\text{ multiplication }} \gsubplusclose_2
\xrightarrow{\text{ extract coefficient }} \allpoly
\]
If we take the coefficient of $y^n$ then we get a map
\[
\gsubplusclose_2\times\init{\allpoly}{\{x^{n+1}\}}
\xrightarrow{\text{ multiplication }}
\init{\gsubplusclose_2\,}{\{x^{n+1}\}}
\xrightarrow{\text{ coefficient of $y^n$ }}
\init{\allpoly}{\{x^{n+1}\}}
\]

For example, take
$g(y)=\sum b_iy^i\in\allpoly$ and $f(x,y) = \sum
a_{ij}x^iy^j\in\gsubclose_2$. Then $fg\in\gsubclose_2$, and the
coefficient of  $y^n$, is in $\allpoly$. The
coefficient of $x^jy^n$ is $\sum_j b_ia_{j,n-i}$. This gives a matrix
that preserves $\init{\allpoly}{\{x^{n+1}\}}$. For instance, if we consider
$2$-initial segments and take $n=2$ then the mapping
\[
\begin{pmatrix}
  b_0\\b_1 \\b_2
\end{pmatrix} \mapsto
\begin{pmatrix}
a_{02} & a_{01} & a_{00} \\
a_{12} & a_{11} & a_{10} \\
a_{22} & a_{21} & a_{20} \\
\end{pmatrix}
\begin{pmatrix}
  b_0\\b_1 \\b_2
\end{pmatrix} 
\]
maps $\init{\allpoly}{\{x^3\}}$ to itself.  Similarly, we can  find
matrices that preserve $\init{\rupint{2}\,}{\{y^n\}}$.  If we take
$f(x,y,z)\in\rupint{3}$ and $g(x,z)\in\gsub_2$ then

\index{matrix!preserving initial segments}

\begin{theorem} \label{thm:multiply-vectors} Suppose
  $f\in\gsubposclose_3$ where  $f = \sum f_{i,j}(x)\,y^iz^j$. If we
  define
  \begin{align*}
    M &= \left(f_{n-i-1,j}\right)_{ 0\le i,j\le n-1}
  \end{align*}
  then $M$ maps $\init{\rupint{2}}{\{y^n\}}$ to itself.
\end{theorem}

Consider an example. 
  Choose $g\in\allpolypos$, and consider
  $g(x+y+z)$ in Theorem~\ref{thm:multiply-vectors}. The matrix  $M$ is
 
\[ \begin{pmatrix}
   \dfrac{g^{(2)}}{2!0!} & \dfrac{g^{\prime}}{1!0!} & \dfrac{g}{0!0!} \\[.4cm]
   \dfrac{g^{(3)}}{2!1!} & \dfrac{g^{(2)}}{1!1!} & \dfrac{g^{\prime}}{0!1!} \\[.4cm]
   \dfrac{g^{(4)}}{2!2!} & \dfrac{g^{(3)}}{1!2!} & \dfrac{g^{(2)}}{0!2!} 
  \end{pmatrix}
\]

We can use the theorem to find matrices $M$ and vectors $v$ such that
$M^kv$ consists of interlacing polynomials for $k=1,2,\dots$. For
example, take
\begin{align*}
  f &= (x + y + 1) (x + 2  y + 2)(2  x + y +  1) \\
&= (2 x^3+7 x^2+7 x+2) + (7 x^2+14 x+6)y+(7 x+6)y^2+2y^3 \\
v &= \bigl(2 x^3+7 x^2+7 x+2, 7 x^2+14 x+6,7 x+6,2\bigr) \\
  g &=(x + y + z + 1)^2(x + 2  y + 3  z + 4)^2 \\
  M &= 
\left(
\begin{array}{llll}
 6 (4 x+7) & 22 x^2+86 x+73 & 2 (x+1) (x+4) (4 x+7) & (x+1)^2 (x+4)^2 \\
 30 & 2 (29 x+53) & 2 \left(17 x^2+70 x+62\right) & 6 (x+1) (x+2) (x+4) \\
 0 & 37 & 2 (23 x+44) & 13 x^2+56 x+52 \\
 0 & 0 & 20 & 12 (x+2)
\end{array}
\right)
\end{align*}
Since $M^kv\in\init{\rupint{2}\,}{\{y^4\}}$ for $k=1,\dots$, it follows
that all vectors $M^kv$ consist of interlacing polynomials.
If we take $g = (x+y+z+1)^6$ then the matrix $M$ is composed of
multinomial coefficients:

\[
\left(
\begin{array}{llll}
\binom{6}{3,0}\, (x+1)^3 & \binom{6}{2,0}\, (x+1)^4 & \binom{6}{1,0}\, (x+1)^5 & \binom{6}{0,0}\, (x+1)^6 \\[.2cm]
\binom{6}{3,1}\, (x+1)^2 & \binom{6}{2,1}\, (x+1)^3 & \binom{6}{1,1}\, (x+1)^4 & \binom{6}{0,1}\, (x+1)^5 \\[.2cm]
\binom{6}{3,2}\, (x+1) & \binom{6}{2,2}\, (x+1)^2 & \binom{6}{1,2}\, (x+1)^3 & \binom{6}{0,2} \, (x+1)^4 \\[.2cm]
\binom{6}{3,3} & \binom{6}{2,3}\, (x+1) & \binom{6}{1,3}\, (x+1)^2 & \binom{6}{0,3} \, (x+1)^3
\end{array}
\right)
\]

\index{Hadamard product!for ideals}

The Hadamard product defines a map
$\allpolypos\times\allpoly\longrightarrow\allpoly$. Since
$f\times\cdot$ maps $\{x^n\}$ to itself, we have a map
\[
\init{\allpolypos}{\{x^n\}}\times
\init{\allpoly}{\{x^n\}}\longrightarrow
\init{\allpoly}{\{x^n\}}
\]

The converse is slightly more complicated.
\begin{lemma}
  Suppose that the Hadamard product with $f$ maps the closure of
  $\init{\allpolypos}{\{x^n\}}$ to itself. Then $\exp(f)\in\init{\allpolypos}{\{x^n\}}$.
\end{lemma}
\begin{proof}
  $e^x$ is in the closure of $\init{\allpolypos}{\{x^n\}}$.
\end{proof}

  Here is an interesting construction of matrices that preserve
  interlacing in $\allpolypos$. Recall $\charpoly{M}$ is the
  characteristic polynomial of $M$, and $M[1,2]$ is $M$ with the
  first two rows and columns removed.

  \begin{lemma}\label{lem:pres-int-det}
    If $M$ is positive definite then the matrix below preserves
    interlacing in $\allpolypos$.
\[
\begin{pmatrix}
  \charpoly{M[1]} &   \charpoly{M} \\
  \charpoly{M[1,2]} &   \charpoly{M[2]}
\end{pmatrix}
\]
  \end{lemma}
  \begin{proof}
     The determinant of 
\begin{gather*}
-xI + M + 
\left(\begin{smallmatrix} y & 0 & \hdots \\ 0 & z & 0 & \hdots \\
\vdots & \vdots & & \text{\huge 0}
\end{smallmatrix}\right)\\
\intertext{equals}
\charpoly{M} + y\, \charpoly{M[1]}+ z\, \charpoly{M[2]} + y\,z\, \charpoly{M[1,2]}
    \end{gather*}
Now apply Theorem~\ref{thm:multiply-vectors}.
  \end{proof}

\section{Stable polynomials}
\label{sec:stable-polynomials-1}
\index{polynomial!stable}

Some stable cases are easy. 
\begin{lemma}\ 
  \begin{enumerate}
  \item $\init{\stabled{1}}{\{x^2\}}= \left\{ a_0+a_1x+a_2x^2\,\mid\, a_0,a_1,a_2>0\,
    \right\}$
  \item $\init{\stabled{1}}{\{x^4\}}= \left\{  a_0+a_1x+a_2x^2+a_3x^3 \,\mid\, a_0,a_1,a_2,a_3>0
      \And a_1a_2> a_0a_3\, \right\}$
  \item $\init{\stabled{1}}{\{x^5\}}\supset \left\{ 
      a_0+\cdots+a_4x^4\,\mid\,  a_0,\dots,a_4\ge0 
      \And  a_1a_2a_3 \ge a_1^2a_4+a_3^2a_0\, \right\}$
  \end{enumerate}
  
\end{lemma}

\begin{proof}
  Lemma~\ref{lem:stable-reasons} shows that if the inequalities are
  strict then the polynomials are already in $\stabled{1}$. The quadratic
  is trivial, and we know from \eqref{eqn:stable-det-1} that
  $a_1a_2>a_0a_3$ characterizes cubic stable polynomials. If $f\in\stabled{1}$ then by
  Lemma~\ref{lem:stable-reasons} we know that its $4$-initial segment is
  in $\stabled{1}$.
\end{proof}

\begin{example}

  We conjecture that
\[
\init{\stabled{2}}{\{y^2}\} \eqques \left\{ f\lesslesseq g\mid f,g\in\stabled{1}\right\}
\]
The left hand side is certainly contained in the right hand side; we
give an example to show that certain elements of the right hand side
are in the left.

We know that $(x+1)^2\lesslesseq x+a$ in $\stabled{1}$ if and only if
$0<a<2$.  Choose positive $\alpha$ and consider the matrix
\begin{gather*}
M = 
\begin{pmatrix} 1 & 0 \\ 0 & 1\end{pmatrix}+
y
\begin{pmatrix} e_1 & 0  \\ 0 & e_2\end{pmatrix}+
\begin{pmatrix} 0 & -\alpha  \\ \alpha & 0\end{pmatrix} + 
\\
x\,\begin{pmatrix}
1+\alpha^2 - \alpha\sqrt{1+\alpha^2} & 0 \\
0 & 1+\alpha^2 + \alpha\sqrt{1+\alpha^2} 
 \end{pmatrix}
\end{gather*}
The determinant of M is
\[
(1+\alpha^2)(1+x)^2 + y\biggl[ (e_1+e_2) + x\biggl((1+\alpha^2)(e_1+e_2)
+ (e_1-e_2)\alpha\sqrt{1+\alpha^2}\biggr)\biggr]+\cdots
\]
Thus, the first term is a constant multiple of $(x+1)^2$, so we need
to determine what the coefficient of $y$ is. If we divide the
coefficient of $x$ in this coefficient by the
constant term $e_1+e_2$ we get
\[
1+\alpha^2 + \frac{e_1-e_2}{e_1+e_2}\alpha\sqrt{1+\alpha^2}
\]
Now as $\alpha,e_1,e_2$ range over the positive reals, the expression
above takes all values in $(1/2,\infty)$, and therefore all linear
terms $x+a$ where $0<a<2$ are realized by a determinant. Thus,
$(x+1)^2 + y(x+a)$  lies in $\init{\stabled{1}}{\{y^2\}}$.

\end{example}

The next result is the $2\times2$ analog of
Theorem~\ref{thm:multiply-vectors}; the general case is similar.
  \begin{lemma}
    Suppose $h(x,y)=\sum h_i(x)y^i$ is in $\stabled{2}$.  For
    appropriate $n$ the matrix
    $\smalltwo{h_{n-1}}{h_{n-2}}{h_{n}}{h_{n-1}}$ maps
    $\init{\stabled{2}}{\{y^2\}}$ to itself.
  \end{lemma}
  \begin{proof}
    We can write $F = f + g y+\cdots$ where
    $F\in\stabled{2}$.  The coefficients of $y^n$ and $y^{n-1}$ in $hF$
    are $h_n f + h_{n-1}g$ and $h_{n-1}f +  h_{n-2}g$. The result
    follows since adjacent coefficients interlace in $\stabled{2}$.
  \end{proof}

\section{Other ideals}
\label{sec:other-ideals}

A representative  element of $\init{\gsubplus_2\,}{\{x^2,y^2\}}$ is a
polynomial $a+bx+cy + dxy$. The characterization is very simple.

\begin{lemma}
  $\init{\gsubplus_2\,}{\{x^2,y^2\}} = \left\{ a+bx + cy + dxy\,\bigl|\,
    a,b,c,d>0 \And \smalltwodet{a}{b}{c}{d}<0\right\}$
\end{lemma}
\begin{proof}
  If $f = a+bx+cy+dxy+\cdots$ is in $\rupint{2}$ then
  Theorem\ref{thm:product-4} says that $\smalltwodet{a}{b}{c}{d}<0$.
  Conversely, assume the determinant is negative. By dividing by $a$
  and scaling $x$ and $y$ it suffices to show that $1+x+y+Dxy\in
  \init{\gsubplus_2\,}{\{x^2,y^2\}}$ where $0<D<1$. Consider the product
  \begin{multline*}
    (1+\alpha x+\beta y)(1+\,(1-\alpha)x+\,(1-\beta)y) \\
    = 1 + x + y + (\alpha(1-\beta)+\beta(1-\alpha))xy + \cdots
  \end{multline*}
Clearly since $0<D<1$ we can find $)<\alpha,\beta<1$ so that
$(\alpha(1-\beta)+\beta(1-\alpha))=D$.
\end{proof}

Although $2xy-1\not\in\rupint{2}$, the product
\[
(x+y+1)(x+y-1) = -1 + 2xy + \cdots
\]
shows that $2xy-1\in\init{\rupint{2}\,}{\{x^2,y^2\}}$. A similar argument
shows that

\begin{lemma}
  $\init{\rupint{2}\,}{\{x^2,y^2\}} = \left\{ a+bx + cy + dxy\,\bigl|\,
    \smalltwodet{a}{b}{c}{d}\le0\right\}$
\end{lemma}

The elements of $\init{\gsubplus_3\,}{\{y^2,z^2\}}$ are interlacing
squares with non-negative determinant
\[ \arbupsquare{f_{00}}{f_{10}}{f_{01}}{f_{11}}{\lesslesseq} \]
arising from $\sum f_{ij}(x)y^iz^j\in\gsubplus_3$. See the next
section.

\begin{lemma}\label{lem:init-yz}
  If $f\lesslesseq g,h$ then $f+yg+zh\in\init{\rupint{3}\,}{\{y^2,yz,z^2\}}$.
\end{lemma}
\begin{proof}
  The argument similar to Lemma~\ref{lem:product}. If we write
  \begin{gather*}
    f = \prod(1+r_i\,x)\\
    g = \sum a_i\,f/(1+r_i\,x)\\
    h = \sum b_i\,f/(1+r_i\,x)\\
\intertext{then}
\prod(1+r_i\,x+a_i\,y+b_i\,z) = f + g\,y + h\,z + \sum_{i\ne j}
a_ib_j\frac{f}{(1+r_ix)(1+r_jx)}yz+ \cdots
\end{gather*}

\end{proof}

\begin{lemma}
   $\init{\rupint{2}\,}{\{x-y\}} \approx \allpoly$
\end{lemma}
\begin{proof}
  A representative of $\init{\rupint{2}\,}{\{x-y\}}$ is a polynomial of
  the form $f(x,x)$ where $f(x,y)\in\rupint{2}$. Now we know that
  $f(x,x)\in\allpoly$ so this gives a map
  $\init{\rupint{2}\,}{\{x-y\}}\longrightarrow \allpoly$. Conversely, if
  $g\in\allpoly$ has positive leading coefficients and we define
  $f(x,y) = g((x+y)/2)$ then $f(x,y)\in\rupint{2}$ and $f(x,x)=g$. 
\end{proof}

If we make a divisibility assumption then we get induced maps.
  \begin{lemma}\label{lem:ideal-falling}
    Suppose that the linear transformation $T$ satisfies
    \begin{enumerate}
    \item $T(x^n)$ divides $T(x^{n+1})$ for $n\ge1$.
    \item $T\colon\allpolyint{\mathcal{R}}\longrightarrow\allpolyint{\mathcal{R}}$
    \item $T^{-1}\colon\allpolyint{\mathcal{S}}\longrightarrow\allpolyint{\mathcal{S}}$
    \end{enumerate}
then
\[
\begin{matrix}
T&\colon&\ \init{\allpolyint{\mathcal{R}}}{\{x^n\}}& \longrightarrow&
\init{\allpolyint{\mathcal{R}}}{\{T(x^n)\}}\\[.2cm]
T^{-1}&\colon&\ \init{\allpolyint{\mathcal{S}}}{\{T(x^n)\}}& \longrightarrow&
\init{\allpolyint{\mathcal{S}}}{\{x^n\}}
\end{matrix}
\]
  \end{lemma}
  \begin{proof}
    If $g = \sum a_i x^i$ is any polynomial then we claim that there
    exist $\alpha_i\in\reals$ such that
    \begin{align}
      T(g\,x^n) & \in\ \{T(x^n)\} \label{eqn:ideal-fact-1}\\
      g\,T(x^n) &=\ \sum \alpha_i T(x^{n+i})\label{eqn:ideal-fact-2}
    \end{align}
    For \eqref{eqn:ideal-fact-1} we have
\[
T(g\,x^n) = T(x^n) \sum a_i \frac{T(x^{n+i})}{T(x^n)}
\]
which is of the form $h(x)\,T(x^n)$ by the first hypothesis. For
\eqref{eqn:ideal-fact-2}, using induction it suffices to show that it
is true for $g=x$. If we write $T(x^{n+1}) = (\alpha+\beta x)T(x^n)$
then
\[
x\,T(x^n) = \frac{1}{\beta} T(x^{n+1}) - \frac{\alpha}{\beta} T(x^n)
\]
If $f_0\in \init{\allpolyint{\mathcal{R}}}{\{x^n\}}$ then write
$  f_0 = f + g\,x^n$
where $f\in\allpolyint{\mathcal{R}}$. Applying $T$
\begin{align*}
  T(f_0) &= T(f) + T(g\,x^n)
\end{align*}
and using the facts that $T(f)\in\allpolyint{\mathcal{R}}$ and 
$T(g\,x^n)\in\{T(x^n)\}$ shows that 
$T(f_0)$ is in $\init{\allpolyint{\mathcal{R}}}{\{T(x^n)\}}$.
If $f_0\in \init{\allpolyint{\mathcal{S}}}{\{T(x^n)\}}$ then write
$  f_0 = f + g\,T(x^n)$
where $f\in\allpolyint{\mathcal{S}}$. Applying $T^{-1}$ yields
\begin{align*}
  T^{-1}(f_0) &= T^{-1}(f) + T^{-1}(g\,T(x^n))
\end{align*}
Since $T^{-1}(g\,T(x^n))= x^n\sum \alpha_ix^i$ and
$T^{-1}(f)\in\allpolyint{\mathcal{S}}$ 
we see that
$T^{-1}(f_0)$ is in $\init{\allpolyint{\mathcal{S}}}{\{x^n\}}$.
  \end{proof}

  \begin{cor}
    Let $T(x^n) = \falling{x}{n}$. Then
\[
\begin{matrix}
T&\colon&\ \init{\allpolyalt}{\{x^n\}}& \longrightarrow&
\init{\allpolyalt}{\{T(x^n)\}}\\[.2cm]
T^{-1}&\colon&\ \init{\allpolypos}{\{T(x^n)\}}& \longrightarrow&
\init{\allpolypos}{\{x^n\}}
\end{matrix}
\]

  \end{cor}

\index{falling factorial}

  \begin{cor}
    If $a+bx+cx^2\in\init{\allpolypos}{\{x(x-1)(x-2)\}}$ then
    $(b+c)^2> 2ac$.
  \end{cor}
  \begin{proof}
    If $a+bx+cx^2\in\init{\allpolypos}{\{x(x-1)(x-2)\}}$ then 
and $T(x^n) = \falling{x}{n}$ then
\[
T^{-1}(a+bx+cx^2) = a + bx + c(x^2+x) \in\init{\allpolypos}{\{x^3\}}
\]
This implies that $(b+c)^2> 2ac$.
  \end{proof}

  \section{Coefficients of extensions}
  \label{sec:coeff-extens}

  In this section we choose  $f\in\mathcal{S}/\mathcal{I}$ and a
  monomial $\xx^\diffi$, and consider properties of the set
\[
\biggl\{ \alpha\,\bigl|\, \text{$\alpha$ is the coefficient of
$\xx^\diffi$ in $F\in \mathcal{S}$ where $F\equiv f
\pmod{\mathcal{I}}$} \biggr\}
\]

Our first elementary example is a linear polynomial. We take
$\mathcal{S}=\allpolypos$, $\mathcal{I} = \{x^3\}$, $\xx^\diffi =
x^3$, and $f = a+bx$. 

\begin{lemma}
  If $a,b>0$ then
\[
\biggl\{\alpha\,\bigl|\, a+bx+\alpha
x^2\in\init{\allpolypos}{\{x^3\}}\biggr\} 
= \bigl(\,0,b^2/(2a)\,\bigr)
\]
\end{lemma}
\begin{proof}
  Use the fact that $a+bx+\alpha x^2\in\init{\allpolypos}{\{x^3\}}$\/ if
  and only if $b^2> 2a\alpha$. 
\end{proof}

The corresponding problem for $\gsubplus_2$ is unsolved. We believe
that if $f+yg\in\init{\gsubplus_2\,}{\{y^2\}}$ (equivalently
$f\lesslesseq g$) then
\[
\biggl\{h\,\bigl|\,  f+ yg+hy^2\in\init{\gsubplus_2\,}{\{y^3\}}\biggr\} 
\eqques
\biggl\{h\,\bigl|\, g\lesslesseq h \And \smalltwodet{f}{g}{g}{h}\le0
\biggr\}
\]

The conjectured solution for $\gsubplus_3$ is similar. If
$f+gy+hz\in\init{\gsubplus_3\,}{\{y^2,z^2\}}$ then
\begin{equation}\label{eqn:init-xyz}
\biggl\{k\,\bigl|\,  f+ yg+hz + kyz\in\init{\gsubplus_3\,}{\{y^2,z^2\}}\biggr\} 
\eqques
\biggl\{k\,\bigl|\, g,h\lesslesseq k \And \smalltwodet{f}{g}{h}{k}\le0
\biggr\}
\end{equation}

We don't know how to approach these problems, but we can establish
properties of the interlacing squares of the last question. Using
these, we can establish \eqref{eqn:init-xyz} when the degree of $f$ is
two. So, we consider the new question

\begin{quote}  What are the extensions of three interlacing polynomials to an
  interlacing square with a positive determinant? 
\end{quote}

We are given three
  polynomials $f,g,h$ in $\allpolypos$ with positive leading
  coefficients that satisfy $f\lessless g$, $f\lessless h$, and consider the set
  $\mathcal{K}$ of all polynomials $k$ with positive leading
  coefficient satisfying

  \index{extension!of 3 polynomials}

\begin{gather} \label{eqn:ext-sq-1}
\xymatrix@=8pt{
f \ar@{}[r]|{\lessless} & g \\
h \ar@{}[r]|{\lessless} \ar@{}[u]|{\rotatebox{270}{$\lessless$}} 
& k \ar@{}[u]|{\rotatebox{270}{$\lessless$}}  
}\\\label{eqn:ext-sq-2}
\smalltwodet{f}{g}{h}{k}>0 \qquad\text{for $x\in\reals$}
\end{gather}

First of all,
\begin{lemma}\label{lem:ext-sq}
  $\mathcal{K}$ is a non-empty bounded convex set of polynomials.
\end{lemma}
\begin{proof}
  Lemma~\ref{lem:init-yz} shows $\mathcal{K}$ is non-empty.
  Suppose that $k_1$ and
  $k_2$  satisfy \eqref{eqn:ext-sq-1} and \eqref{eqn:ext-sq-2}. If
  $0<\alpha<1$ then $\alpha$ and $1-\alpha$ are positive, so
 
\centerline{\xymatrix@=8pt{
f \ar@{}[rrr]|{\lessless} &&& g \\
h  \ar@{}[u]|{\rotatebox{270}{$\lessless$}} 
&\lessless&& \alpha\, k_1 + (1-\alpha)\,k_2\ar@{}[u]|{\rotatebox{270}{$\lessless$}}  
}}
From the determinant hypothesis:
\begin{align*}
  gh - fk_1 & > 0 \\
  gh - fk_2 & > 0. \\
\intertext{Multiplying and adding yields}
  gh - f(\alpha k_1 + (1-\alpha)k_2) & > 0 \\
\end{align*}
Thus $\alpha k_1+(1-\alpha)k_2$ satisfies the two conditions, and is
in $\mathcal{K}$.

Since $f(1)>0$ we see that $g(1)h(1)/f(1) > k(1)$. Consequently, all
coefficients of $k$ are bounded by $g(1)h(1)/f(1)$. 
\end{proof}

\begin{example}
  If $f = (x+1)^n$, $g=h=(x+1)^{n-1}$ then the interlacing requirement
  implies that $k = \alpha(x+1)^{n-2}$. The determinant requirement
  implies that $0<\alpha<1$, and so
\[
\mathcal{K} \subseteq \bigl\{\alpha(x+1)^{n-2}\,\mid\,0<\alpha<1\bigr\}
\]
In order to see that this is an equality, we proceed as follows. If
$a_i,b_i$ are non-negative then  $0\le\left(\sum a_ib_i\right) \le \bigl(\sum
a_i\bigr)\bigl(\sum b_j\bigr)$. 
So, if $0\le \alpha\le 1$ choose $a_i,b_i$ so that $\sum a_i=\sum
b_i=1$ and $\sum a_ib_i = 1-\alpha$. We have
\begin{multline*}
  \prod_1^n (x+1+a_iy + b_iz) \\=
(x+1)^n + y\,(x+1)^{n-1}+ z\,(x+1)^{n-1} +
yz(x+1)^{n-2}\biggl(\sum_{i\ne j}a_ib_j\biggr) 
\end{multline*}
This is the desired representation since 
\begin{gather*}
  \sum_{i\ne j} a_ib_j = \biggl(\sum a_i\biggr)\biggl(\sum b_j\biggr) -
  \sum a_ib_i = 1-(1-\alpha)=\alpha.
\end{gather*}
\end{example}

\begin{example}
  If $F(x,y,z)\in\gsubplus_3$ and $F = \sum f_{ij}(x)y^iz^j$ then
  $f_{11}$ is an extension of $f_{00},f_{10},f_{01}$. Thus, the
  polynomial $f_{11}\in\mathcal{K}$. In particular, if
  $f(x)\in\allpolypos$, then $f(x+y+z)\in\gsubplus_3$, so
  $f''(x)/2$ is an extension of $f,f',f'$. 
\end{example}
In the case that $f$ has degree $2$, $\mathcal{K}$ is an interval of
positive numbers. We may assume that the polynomials are as given
below, where all constants are positive:

\begin{equation}\label{eqn:ext-sq-3}
\xymatrix@=8pt{
(x+r_1)(x+r_2) 
\ar@{}[rrr]|{\lessless} 
&&& t_1(x+r_1) + t_2(x+r_2) 
\\
s_1(x+r_1) + s_2(x+r_2)
\ar@{}[rrr]|{\lessless} \ar@{}[u]|{\rotatebox{270}{$\lessless$}} 
&&& k \ar@{}[u]|{\rotatebox{270}{$\lessless$}}  
}
\end{equation}

The construction of the lemma suggests we consider the polynomial
\[
(x+r_1+s_1y+t_1z)(x+r_2+s_2y+t_2z)
\]
It provides an extension where $k$ is $s_1t_2+s_2t_1$. Note that this
value does not depend on the roots of $f$. Using these polynomials, we
have

\begin{lemma}
  The set $\mathcal{K}$ of extensions of \eqref{eqn:ext-sq-3} is the interval 
\[( s_1t_2+s_2t_1 - 2 \sqrt{s_1s_2t_1t_2},s_1t_2+s_2t_1 + 2
\sqrt{s_1s_2t_1t_2}).\] 
\end{lemma}
\begin{proof}
  The polynomial $gh-kf$ is a polynomial in $x$, and for it to be
  positive its discriminant must be negative. The discriminant is a
  quadratic in $k$ whose roots are $s_1t_2+s_2t_1 \pm 2
  \sqrt{s_1s_2t_1t_2}$.
\end{proof}

An explanation of these conditions is that when $n=2$
\eqref{eqn:ext-sq-1} and \eqref{eqn:ext-sq-2} are exactly what we need
to extend to a polynomial in $\gsubplus_3$.

\begin{lemma}
  Assume that $deg(f)=2$. The polynomials $f,g,h,k$ in $\allpolypos$
  satisfy the conditions \eqref{eqn:ext-sq-1} and
  \eqref{eqn:ext-sq-2} above if and only if 

\[ f + gy + hz + yzk \in\init{\gsubplus_3\,}{\{y^2,z^2\}}
\]
Equivalently, if $f+gy+hz\in\init{\gsubplus_3\,}{\{y^2,z^2\}}$ and $f$
has degree $2$  then \eqref{eqn:init-xyz} holds.
\end{lemma}
\begin{proof}
  We only need to show that if we are given $f,g,h,k$ satisfying the
  conditions then there is such a polynomial in $\gsubpos_3$. We start
  with the matrix
\[M=
\begin{pmatrix}
  (r_1 + x + s_2\, y + t_2\, z)/r_1 &
 d_2\, y + e_2\, z \\
 d_2\, y + e_2\, z &
 (r_2 + x + s_1\, y + t_1\, z)/r_2
\end{pmatrix}
\]
and notice that $r_1r_2\,|M|$ is in $\rupint{3}$ and equals
\[
f + g \,y + h\,z + (-2 d_2 e_2 r_1 r_2 + s_2 t_1 + s_1 t_2)yz + \cdots
\]

From the Lemma above we know that we can write
\[
k = s_1t_2+s_2t_1 + 2\beta \qquad\text{where}\ |\beta| \le
\sqrt{s_1s_2t_1t_2}
\]
Since the matrix $M$ is of the form $I+xD_1+yD_2+zD_3$ we know
that we must choose $d_2,e_2$ so that $D_1,D_2,D_3$ are positive
definite matrices.
Now $D_2 = \smalltwo{s_2/r_1}{d_2}{d_2}{s_1/r_2}$ and 
$D_3 = \smalltwo{t_2/r_1}{e_2}{e_2}{t_1/r_2}$.
the determinant must be  positive $d_2$ and $e_2$ must satisfy
\begin{align*}
  s_1s_2 / (r_1r_2) & > d_2^2 \\
  t_1t_2 / (r_1r_2) & > e_2^2 
\end{align*}
Notice that if $d_2,e_2$ satisfy these inequalities then
\[
\sqrt{s_1s_2t_1t_2} \ge r_1r_2d_2e_2
\]
Consequently we can choose $d_2,e_2$ satisfying $d_2 e_2 r_1 r_2 =
\beta$. 
\end{proof}

\begin{cor}
  If $f,g,h,k$ satisfy \eqref{eqn:ext-sq-1} and \eqref{eqn:ext-sq-2}
  and $f$ has degree $2$ then $\smalltwodet{f}{g}{h}{k}$ is a stable
  polynomial. 
\end{cor}
\begin{proof}
  This follows from the lemma and Lemma~\ref{lem:stable-square-3}.
\end{proof}

\begin{remark}
  If the determinant is always negative (rather than positive) then
  the corollary fails.  The following polynomials form an interlacing
  square with everywhere negative determinant that is not stable.\\[.2cm]
\index{polynomial!stable}

\centerline{\xymatrix@=8pt{
(x+0.9) (x+0.97) (x+9.5)
\ar@{}[rrr]|{\lessless} 
&&&  (x+0.95) (x+2.9)
\\
(x+0.94) (x+2.6)
\ar@{}[rrr]|{\lessless} \ar@{}[u]|{\rotatebox{270}{$\lessless$}} 
&&& x+1.5 \ar@{}[u]|{\rotatebox{270}{$\lessless$}}  
}}

\end{remark}


\backmatter

\part{Appendix}

\chapterstyle{section}

\chapter{Glossary of polynomials, spaces, and interlacing}

\renewcommand{\TimeStampStart}{Thursday, January 17, 2008: 19:50:26}
\mytoday  

\newcommand{\gloss}[2]{{\makebox[1.5cm][l]{#1}} & \parbox[t]{4in}{#2}\\[.2cm]}

\begin{center}
  \textbf{\large Polynomials}
\end{center}

\begin{supertabular}{ll}
\gloss{$A_n$}{Euler polynomial}
\gloss{$C_n^\alpha$}{Charlier polynomial with parameter $\alpha$}
\gloss{$C_n(x;y)$}{$(-1)^nC_n^y(-x) = (-1)^n n! L_n^{-x-n}(y)=(-1)^nn!L_n(-y;-x-n)$}
\gloss{$G_n^{(\alpha)}$}{Gegenbauer polynomial with parameter $\alpha$}
\gloss{$H_i$}{Hermite polynomial}
\gloss{$H_\sdiffi$}{Hermite polynomial in $d$ variables}
\gloss{$L_n$}{Laguerre polynomial}
\gloss{$\tilde{L}_n$}{monic Laguerre polynomial: $n! L_n$}
\gloss{$L_n^\alpha$}{Laguerre polynomial with parameter $\alpha$}
\gloss{$L_n(x;y)$}{$L_n^y(-x)$}
\gloss{$P_n$}{Legendre polynomial}
\gloss{$Q_n^{\alpha,\beta}$}{$\left(\ (x-\alpha)\affa^{-1} - \beta\ \right)^n(1)$}
\gloss{$T_n$}{Chebyshev polynomial}
\gloss{$\falling{x}{n}$}{falling factorial $x(x-1)\cdots(x-n+1)$}
\gloss{$\rising{x}{n}$}{rising factorial $x(x+1)\cdots(x+n-1)$}
\gloss{$(x;q)_n$}{$(1-x)(1-qx)\cdots(1-q^{n-1}x)$}
\end{supertabular}

\pagebreak

\begin{center}
  \textbf{\large Sets of polynomials in one variable}
\end{center}

\begin{supertabular}{ll}
\gloss{$\allpoly$}{ All polynomials in one variable with all real  roots.}
\gloss{$\allpolyf$}{    The analytic closure of $\allpoly$.}
\gloss{$\allpolypos$}{ All polynomials in $\allpoly$ with all negative  roots and  all positive (or all negative) coefficients: $\allpolyint{(-\infty,0)}$}
\gloss{$\allpolyposf$}{ The analytic closure of $\allpolypos$.}
\gloss{$\allpolyalt$}{ All polynomials in $\allpoly$ with all positive  roots and alternating coefficients: $\allpolyint{(0,\infty)}$\\}
\gloss{$\allpolyaltf$}{ The analytic closure of $\allpolyalt$.}
\gloss{$\allpolypm$}{$\allpolyalt\cup\allpolypos$}
\gloss{$\allpolyint{I}$}{ All polynomials in $\allpoly$ with all roots
  in an interval $I$.}
\gloss{$\allpolysep$}{Polynomials in $\allpoly$ whose roots are all at least one  apart.}
\gloss{$\allpolysepf$}{ The analytic closure of $\allpolysep$.}

\gloss{$\polycpx$}{Polynomials $f+\imag g$ where $f \lessless g$.}
\gloss{$\polycpxclose$}{Polynomials $f+\imag g$ where $f \lesslesseq g$.}
\gloss{$\polycpxf$}{Analytic closure of $\polycpx$.}
\gloss{$\stabled{1}$}{Stable  polynomials with real coefficients}
\gloss{$\stabledopen{}$}{Interior of $\stabled{1}$}
\gloss{$\stabledf{}$}{Analytic closure of $\stabled{1}$}
\gloss{$\stabledc{1}$}{Stable  polynomials with complex coefficients}
\gloss{$\stabledcf{1}$}{Analytic closure of $\stabledc{1}$}

\end{supertabular}

\vfill
\pagebreak

\begin{center}
  \textbf{\large Sets of polynomials in several variables}
\end{center}

\begin{supertabular}{ll}
\gloss{$\rup{d}$}{Polynomials in $d$ variables non-vanishing in the
  upper half plane}
\gloss{$\up{d}$}{Polynomials in $d$ variables and complex
  coefficients non-vanishing in the   upper half plane}
\gloss{$\stabled{d}$}{Polynomials in $d$ variables non-vanishing in the
  right half plane}
\gloss{$\stabledc{d}$}{Polynomials in $d$ variables with complex
  coefficients non-vanishing in the  right half plane}
\gloss{$\polypos{d}$}{$\up{d}\cap\stabled{d}$}
\gloss{$\rupint{d}$}{Polynomials in $d$ variables satisfying degree
  conditions, substitution, and homogeneous part in $\gsubpos_{d-1}$}
\gloss{$\gsubclose_d$}{Polynomials in the closure of $\rupint{d}$}
\gloss{$\allpolyf_d$}{The analytic closure of $\rupint{d}$}

\gloss{$\gsubplus_d$}{Polynomials in $\gsubpos_d$ with all positive
  coefficients.}
\gloss{$\gsubplusclose_d$}{Polynomials in the closure of $\gsubplus_d$}
\gloss{$\gsubposf_d$}{The analytic closure of $\gsubplus_d$}
\gloss{$\partialpoly{d,e}$}{Polynomials satisfying
  substitution for $\xx\in\reals$, $\yy\ge0$, the degree condition,
  and homogeneous part in $\allpoly_{d+e}$.}
\gloss{$\partialpolyclose{d,e}$}{Polynomials in the closure of $\partialpoly{d,e}$.}
\gloss{$\partialpolyf{d,e}$}{The analytic closure of $\partialpoly{d,e}$.}
\gloss{$\pospoly_{d,e}$}{Extension of $\gsubplus_d$ by $e$ variables
  subject to restricted substitution conditions.}
\gloss{$\stabled{d}$}{stable  polynomials with real coefficients in $d$ variables}
\gloss{$\stabledf{d}$}{analytic closure of $\stabled{d}$ in $d$ variables}
\gloss{$\gsubgen_2$}{Polynomials in $\rupint{2}$ whose solution curves
  are all distinct}
\gloss{$\gsubsep_2$}{Polynomials in $\rupint{2}$ whose substitutions are
  all in $\allpolysep.$}

\end{supertabular}

\pagebreak

\begin{center}
\textbf{Kinds of interlacing}  
\end{center}

\begin{supertabular}{ll}
\gloss{$f\ulace g$}{$f + y g \in\up{d}$}
\gloss{$f\hlace g$}{$f + y g \in\stabledc{d}$}
\gloss{$f\place g$}{$f + y g \in\polypos{d}$}
\gloss{$f\poslace g$}{$f +\alpha g\in\allpoly$ for $\alpha\ge0$}
\gloss{$f\hposlace g$}{$f +\alpha g\in\stabled{d}$ for $\alpha\ge0$}
\gloss{$f\pposlace g$}{$f +\alpha g\in\polypos{d}$ for $\alpha\ge0$}
\gloss{$f\lesslesseq g$}{The roots of $f$ and $g$ alternate, and $deg(g)+1=deg(f)$}
\gloss{$f\greateqeq g$}{The roots of $f$ and $g$ alternate,
  $deg(g)=deg(f)$, and the largest root belongs to $f$}
\gloss{$f\longleftarrow g$}{$f\lesslesseq g$ or $g\greateq g$ in $\allpoly$}
\end{supertabular}


\chapter{Tables of transformations}
 
\renewcommand{\TimeStampStart}{Monday, January 07, 2008: 11:38:37}
\mytoday  

\newcommand{\ctrans}[7]{#1 & $#2$ & $#3$ & $#4$ & $#5$ & #7 & $#6$\\[.1cm]}
\newcommand{\btrans}[7]{#1 & $#2$ & $#3$ & $#4$ & $#5$ & #7 & $#6$\\[.1cm]}
\newcommand{\trans}[7]{ #1 & $#2$ & $#3$ & $#4$ & $#5$ & #7 & $#6$\\[.1cm]}

\newcommand{\nocol}{\multicolumn{1}{c}{}}
\newcommand{\mycaption}[1]{ &&&&&& \\ \multicolumn{7}{l}{\textbf{ #1}} \\
  &&&&&& \\}
\newcommand{\mybreaktables}{%
\end{tabular}

\noindent%
\hspace*{-1cm}
\begin{tabular}{lccllll}
\nocol & \textbf{ f} & \textbf{ T(f)} & \textbf{ Domain} & \textbf{ Range}  & & \\ \hline
}

\noindent%
\hspace*{-1cm}
\begin{tabular}{lccllll}
\nocol & \textbf{ f} & \textbf{ T(f)} & \textbf{ Domain} & \textbf{ Range}  & & \\ \hline
%
%
\mycaption{Multiplier transformations $\allpoly \longrightarrow\allpoly$}
\trans{Binomial}{x^n}{\binom{nd}{k}x^n}{\allpoly}{\allpoly}{d \text{ odd}, k=1,2\dots}{\scorr{epoly-binom}}
\trans{Binomial}{x^n}{\falling{nd}{k}  x^n}{\allpoly}{\allpoly}{d \text{ odd}, k=1,2\dots}{\scorr{epoly-binom}}
\trans{Derivative}{x^n}{nx^{n-1}}{\allpoly}{\allpoly}{}{\sthm{rolle-2}}
\trans{Exponential}{x^n}{x^n/n!}{\allpoly}{\allpoly}{}{\slem{integral}}
\trans{Exponential}{x^n}{x^n n!/(kn)!}{\allpoly}{\allpoly}{k=1,2,\dots}{\slem{exp-ki}}
\trans{Exponential}{x^n}{x^n/(kn)!}{\allpoly}{\allpoly}{k=1,2,\dots}{\slem{exp-ki}}
\trans{$q$-Exponential}{x^n}{x^n/[n]!}{\allpoly}{\allpoly}{q>1}{\sthm{multiplier}}
\trans{$q$-series}{x^n}{q^{\binom{n}{2}}x^n}{\allpoly}{\allpoly}{|q|<1}{\sthm{multiplier}}
\trans{Rising factorial}{x^n}{\frac{\rising{n}{i}}{i!}x^i}{\allpoly}{\allpoly}{}{\pageref{eqn:nii}}
%
%
\mycaption{Transformations $\allpoly \longrightarrow\allpoly$}
\trans{Derivative}{f}{f'}{\allpoly}{\allpoly}{}{\pageref{ex:fi}}
\trans{Falling Factorial}{ \falling{x}{n}}{\falling{x}{n+1}}{\allpoly}{\allpoly}{}{\slem{increment}}
\trans{Hermite}{x^n}{H_n}{\allpoly}{\allpoly}{}{\scorr{hermite}}
\trans{Hermite}{H_n}{H_{n+1}}{\allpoly}{\allpoly}{}{\scorr{hermite}}
\trans{Laguerre}{x^n}{{L}_n}{\allpoly}{\allpoly}{}{\scorr{laguerre-x}}
\trans{Laguerre}{x^n}{\rev{{L}}_n}{\allpoly}{\allpoly}{}{\slem{laguerre-4}}
\trans{Laguerre}{x^n}{{L}_k^n}{\allpoly}{\allpoly}{k=1,2\dots}{\slem{laguerre-3}}
\trans{Legendre}{x^n}{P_n/n!}{\allpoly}{\allpoly}{}{\slem{legendre2}}
\trans{Rising Factorial}{\rising{x}{n}}{\rising{x}{n+1}}{\allpoly}{\allpoly}{}{\slem{increment}}
%
%
\mycaption{Transformations $\allpolypos \longrightarrow \allpoly$}
\trans{Affine}{x^n}{(-1)^{\binom{n}{2}}x^n}{\allpolypos}{\allpoly}{}{\slem{affine-x}}
\trans{Chebyshev}{x^n}{\rev{T_n}}{\allpolypos}{\allpoly}{}{\slem{cheby-rev}}
\trans{Laguerre}{x^n}{\rev{L_n}(-x)}{\allpolypos}{\allpolypos}{}{\scorr{herm-1}}

%
%
\mycaption{Transformations $\allpolypos \longrightarrow \allpolypos$}
\trans{Charlier}{C_n^\alpha}{x^n}{\allpolypos}{\allpolyneg}{?}{\scorr{charlier-x}}
\trans{Falling factorial}{ \falling{x}{n}}{x^n}{\allpolypos}{\allpolyneg}{}{\spropp{recur-5}}
\trans{Falling factorial}{ \falling{x}{n}}{\rising{x}{n}}{\allpolypos}{\allpolyneg}{}{\scorr{rising-factorial}}
\trans{Hurwitz}{x^n}{\begin{cases} x^{n/2} & n \text{ even} \\ 0 & n
    \text{ odd}\end{cases}}{\allpolypos}{\allpolyneg}{}{\sthm{hurwitz}}
\trans{Rising factorial}{x^n}{\rising{x}{n}}{\allpolypos}{\allpolypos}{}{\scorr{rising-factorial}} 
\trans{Rising Factorial}{x^n}{\rising{x}{n}/n!}{\allpolypos}{\allpolypos}{}{\slem{rising-fact}}
\mybreaktables 
%
%
\mycaption{Transformations $\allpolyalt \longrightarrow \allpoly$}
\trans{Affine}{x^n}{(-1)^{\binom{n}{2}}x^n}{\allpolyalt}{\allpoly}{}{\slem{affine-x}}
\trans{Chebyshev}{x^n}{\rev{T_n}}{\allpolyalt}{\allpoly}{}{\slem{cheby-rev}}
\trans{Hermite}{x^n}{H_nL_n}{\allpolyalt}{\allpoly}{}{\scorr{hnln}}
\trans{Hermite}{x^n}{(H_n)^2}{\allpolyalt}{\allpoly}{}{\scorr{hnln}}
\trans{Laguerre}{x^n}{(L_n)^2}{\allpolyalt}{\allpoly}{}{\scorr{hnln}}
%
%
\mycaption{Transformations $\allpolyalt \longrightarrow  \allpolyalt$}
\trans{Binomial}{x^n}{\binom{x}{n}}{\allpolyalt}{\allpolyalt}{}{\scorr{binomial}}
\trans{Binomial}{x^n}{\binom{nd}{k}x^n}{\allpolyalt}{\allpolyalt}{d \text{ even}, k=1,2\dots}{\scorr{epoly-binom}}
\trans{Binomial}{x^n}{\falling{nd}{k} x^n}{\allpolyalt}{\allpolyalt}{d \text{ even}, k=1,2\dots}{\scorr{epoly-binom}}
\trans{Charlier}{x^n}{C_n^\alpha}{\allpolyalt}{\allpolyalt}{?}{\scorr{charlier-x}}
\trans{Falling factorial}{x^n}{ \falling{x}{n}}{\allpolyalt}{\allpolyalt}{}{\spropp{recur-5}}
\trans{Legendre}{x^n}{P_M(x)/n!}{\allpolyalt}{\allpolyalt}{}{\slem{legendre}}
\trans{Rising factorial}{\rising{x}{n}}{x^n}{\allpolyalt}{\allpolyalt}{}{\scorr{rising-factorial}} 
%
%
\mycaption{Transformations $\allpolyint{I}\longrightarrow\allpolyint{I}$}
\ctrans{Charlier}{x^n}{\rev{}\,C_n^\alpha}{\allpolyint{(-1,0)}}{\allpoly}{}{\scorr{mobius-charlier}}
\ctrans{Chebyshev}{x^n}{T_n}{\allpolyint{(-1,1)}}{\allpolyint{(-1,1)}}{}{\slem{cheby}}
\ctrans{Euler}{x^n}{A_n}{\allpolyint{(-1,0)}}{\allpolyalt}{}{\scorr{euler}}
\ctrans{Euler-Frobenius}{x^n}{P_n}{\allpolyint{(-1,1)}}{\allpolyint{(-1,1)}}{}{\slem{euler-frob}}
\ctrans{Factorial}{x^n}{\prod_1^n(1-kx)}{\allpolyint{(0,1)}}{\allpolyint{(0,1)}}{}{\slem{factorial-variant}}
\trans{Falling Factorial}{x^i}{
  \falling{x}{i}^{Rev}}{\allpolyint{(0,1)}}{
  \allpolyint{(0,1)}}{}{\slem{factorial-variant}}
\ctrans{Laguerre}{x^n}{\rev{}\,{\tilde{{L}}}^{(\alpha)}_n}{\allpolyint{(0,1)}}{\allpoly}{\alpha   \ge-1,}{\scorr{mobius-laguerre}}
\ctrans{Q-series}{x^n}{(x;q)_n}{{\allpolyint{\reals\setminus(0,1)}}}{\allpolyint{\reals\setminus(0,1)}}{0<q<1,}{\sthm{qseries}}
\ctrans{Q-series}{x^n}{(x;q)_n}{\allpolyint{(0,1)}}{\allpolyint{(0,1)}}{1<q,}{\sthm{qseries}}
\ctrans{Q-series}{(x;q)_n}{x^n}{{\allpolyint{\reals\setminus(0,1)}}}{\allpolyint{\reals\setminus(0,1)}}{1<q,}{\sthm{qseries}}
\ctrans{Q-series}{(x;q)_n}{x^n}{\allpolyint{(0,1)}}{\allpolyint{(0,1)}}{0<q<1,}{\sthm{qseries}}
\trans{Rising Factorial}{x^n}{\rising{x}{n}/n!}{\allpolyint{(1,\infty)}}{\allpolyint{(1,\infty)}}{}{\slem{rising-fact}}
%
%
\mycaption{Transformations $\allpoly(n)\longrightarrow\allpoly(n)$}
\trans{Polar derivative}{x^i}{(n-i)x^{i}}{\allpoly(n)}{\allpoly(n-1)}{}{\pageref{sec:polar}}
\trans{Reversal}{x^i}{x^{n-i}}{\allpoly(n)}{\allpoly(n)}{}{\slem{reverse}}
\trans{\Mobius}{x^i}{x^i(1-x)^{n-i}}{\allpoly(n)}{\allpoly(n)}{}{\slem{axb}}
\mybreaktables 
%
%
\mycaption{Transformations  $\allpolyint{I}(n)\longrightarrow\allpolyint{I}(n)$}
%
\trans{Binomial}{x^i}{\binom{\alpha}{i}x^i}{\allpolypos(n)}{\allpolyneg(n)}{\alpha>n-2\ge 0}{\slem{falling-non-int}}    
\trans{Binomial}{x^i}{\binom{x+n-i}{n}}{\allpolyint{(0,1)}(n)}{\allpolyalt(n)}{}{\scorr{Txdi}}
\trans{Chebyshev}{x^i}{(T_r)^i}{\allpolyint{(-1,1)}(n)}{\allpolyint{(-1,1)}(nr)}{}{\scorr{cheby-2}}
\trans{Factorials}{ \falling{x}{k}}{\rising{x}{n-k}}{\allpolypos}{\allpolyneg}{}{\slem{fall-rise}}\trans{Falling factorial}{x^i}{\falling{\alpha}{i}x^i}{\allpolypos(n)}{\allpolyneg(n)}{\alpha>n-2\ge 0}{\slem{falling-non-int}}    
\trans{Falling factorial}{x^i}{\falling{x}{i}x^{n-i}}{\allpolypos(n)}{\allpolyneg(n)}{}{\slem{fall-rise}} 
\trans{Falling factorial}{x^i}{\falling{x+n-i}{n}\frac{x^{i}}{i!}}{\allpolypos(n)}{\allpolypos(n)}{}{\slem{fall-rise}} 
\trans{Hermite}{x^i}{H_i\  H_{n-i}}{\allpolypos(n)}{\allpoly(n)}{}{\scorr{hermite-xy}}
\trans{Hermite}{x^i}{x^{n-i}H_i}{\allpolyint{(2,\infty{)}}(n)}{\allpoly(n)}{}{\scorr{herm2}}
\trans{Hermite}{x^i}{x^{n-i}H_i}{\allpolypos(n)}{\allpoly(n)}{}{\scorr{p2-n-iT}}
\trans{Laguerre}{x^i}{x^{n-i}L_i}{\allpolypos(n)}{\allpoly(n)}{}{\scorr{p2-n-iT}}
\trans{Laguerre}{x^i}{L_i\  L_{n-i}}{\allpolypos(n)}{\allpoly(n)}{}{\scorr{hermite-xy}}
\trans{Rising Factorial}{x^i}{\rising{x}{i}x^{n-i}}{\allpolyint{(0,1)}(n)}{\allpolypos(n)}{}{\scorr{half-ff}}
\trans{Rising factorial}{x^i}{\rising{x}{n-i}}{\allpolypos(n)}{\allpolyneg(n)}{}{\slem{fall-rise}} 
%
%
\mycaption{Affine transformations with $\affa x = x+1$}
\trans{Binomial}{x^n}{\falling{x+d-n}{d}}{\allpolyalt}{\allpolysep}{d=1,2,\dots}{\pageref{lem:binomial-xdk}}
\trans{Charlier}{x^n}{C_n^\alpha}{\allpolyint{(-\alpha,\infty)}}{\allpolyint{(-\alpha,\infty)}\cap\allpolysep}{}{\scorr{affine-trans}}
\trans{Derivative}{f}{f^\prime}{\allpolysep}{\allpolysep}{}{\slem{deriv-diff}}
\trans{Difference}{f}{\Delta{}f}{\allpolysep}{\allpolysep}{}{\slem{deriv-diff}}
\trans{Falling Factorial}{x^n}{ \falling{x}{n}}{\allpolyalt}{\allpolyalt\cap\allpolysep}{}{\scorr{affine-trans}}
\trans{Rising Factorial}{x^n}{\rising{x}{n}}{\allpolypos}{\allpolyneg\cap\allpolysep}{}{\scorr{affine-trans}}
%
%
\mycaption{Transformations of stable polynomials}
\trans{Charlier}{x^n}{C_n^\alpha}{\stabledc{1}}{\stabledc{1}}{\alpha{\in\rhp}}{\pageref{lem:h-charlier}}
\trans{Charlier}{x^n}{C_n^\alpha}{\stabled{1}}{\stabled{1}}{\alpha{\in\rhp}}{\pageref{lem:h-charlier}}
\trans{Difference}{f}{f(x+1)-f(x)}{\stabled{1}}{\stabled{1}}{}{\pageref{cor:stable-difference}}
\trans{Hermite}{x^n}{\imag^nH_n(-\imag  x)}{\stabled{1}}{\stabled{1}}{}{\pageref{cor:h-hermite}}
\trans{Laguerre}{x^n}{\imag^n L_n(\imag x)}{\stabled{1}}{\stabled{1}}{}{\pageref{cor:h-hermite}}
\trans{Rising  factorial}{x^n}{\rising{x}{n}}{\stabledc{1}}{\stabledc{1}}{}{\pageref{lem:h-rising}}
\trans{Rising factorial}{x^n}{\rising{x}{n}}{\stabled{1}}{\stabled{1}}{}{\pageref{lem:h-rising}}
\trans{Falling  factorial}{\falling{x}{n}}{x^n}{\stabledc{1}}{\stabledc{1}}{}{\pageref{lem:h-falling}}
\trans{Falling factorial}{\falling{x}{n}}{x^n}{\stabled{1}}{\stabled{1}}{}{\pageref{lem:h-falling}}
\trans{Hurwitz}{x^n}{\begin{cases} x^{n/2} & n \text{ even} \\ 0 & n
    \text{ odd}\end{cases}}{\stabled{1}}{\allpolypos}{}{\pageref{lem:h-even}}

%
%
\mycaption{Transformations in two or more variables}
\trans{Exponential}{x^iy^j}{\frac{x^iy^j}{(i+j)!}}{\gsubplus_2}{\gsubplus_2}{}{\scorr{ij-pd}}
\trans{Evaluation}{f\times\xx^\diffi}{f(|\diffi|)\xx^\diffi}{\allpolypos\times\gsubplus_d}{\gsubplus_d}{}{\scorr{xd-pdpos}}
\trans{Polar   derivative}{f}{(\xx\cdot\boldsymbol{\partial})f}{\gsubplus_d}{\gsubplus_d}{}{\spropp{xd-pdpos}}

%
\end{tabular}

                                %




\chapter{Empirical tables of transformations}
\label{chap:mma}

\renewcommand{\TimeStampStart}{Monday, January 07, 2008: 11:26:28}
\mytoday

In the table on the next page we list the results of some empirical
testing. We constructed ten transformations for each  polynomial
family. We use the abbreviations 

\begin{tabular}{lcp{4in}}
. & \hspace{.5cm} & The transformation did not appear to map
any interval to $\allpoly$.\\
$X$ && The transformation appears to map
  $\allpolyint{X}\longrightarrow \allpoly$.\\
$\reals$ && The transformation appears to map
  $\allpolyint{(-\infty,\infty)}\longrightarrow \allpoly$.\\
$\reals^+$ && The transformation appears to map
  $\allpolyint{(0,\infty)}\longrightarrow \allpoly$.\\
$\reals^-$ && The transformation appears to map
  $\allpolyint{(-\infty,0)}\longrightarrow \allpoly$.\\
${\mathbb I}$ && The transformation appears to map
  $\allpolyint{(-1,1)}\longrightarrow \allpoly$.\\
${\mathbb I}^+$ && The transformation appears to map
  $\allpolyint{(0,1)}\longrightarrow \allpoly$.\\
${\mathbb I}^-$ && The transformation appears to map
  $\allpolyint{(-1,0)}\longrightarrow \allpoly$.
\end{tabular}

An underlined entry such as $\isknown{(-1,1)}$ means
that the transformation is known to map
$\allpolyint{(-1,1)}\longrightarrow\allpoly$. The other entries are
unsolved questions.

The polynomial families in the first column are in the glossary,
except for the last three. The Laurent-X polynomials are defined by
$p_{-1}=0$, $p_0=1$ and\\[.2cm]

\begin{tabular}{lcll}
Laurent-1 &\hspace*{1cm}&$p_{n+1}$ &= $x(p_n-p_{n-1})$ \\
Laurent-2 &\hspace*{1cm}&$p_{n+1}$ &= $x(p_n-n\,p_{n-1})$ \\
Laurent-3 &\hspace*{1cm}&$p_{n+1}$ &= $(x-n)p_n- n x\,p_{n-1}$ \\
\end{tabular}

There are some cavaets on a few entries.

$
\begin{array}{ccccl}
P_n & x^n & \mapsto & p_nx^{d-n} & \text{proved  for $(-\infty,0)$}\\
\binom{x}{n} & x^n & \mapsto & p_n & \text{only proved for each  half} \\
(x;1/2) & x^n &\mapsto & p_n/n! &  \text{perhaps it's  $(-\infty,0)\cup(0,\infty)$} 
\end{array}
$

The second table lists the actions of just a  single transformation
$x^n\mapsto \text{polynomial}$ on various subsets of the complex plane.

\pagebreak

\noindent
\hspace*{-1.5in}%
$
\begin{array}{c|cccccccccc|}
\hline
& x^n & p_n & x^n & x^n & x^n & p_n  & x^n & x^n & x^n & x^n\\
&\downarrow&\downarrow&\downarrow&
\downarrow&\downarrow&\downarrow&\downarrow&\downarrow&\downarrow&\downarrow\\ 
& p_n & x^n & p_n^{rev} & p_n^2 & p_n x^n & p_{d-n}  & p_n p_{d-n} &
p_nx^{d-n} & p_n/n! & n!p_n\\[.1cm]\hline
&&&&&&&&&&\\
H_n &
 \isknown{\reals} & . & {(-\infty,1)} &
\isknown{\reals^+} & \reals^+ & \reals^+  &
\isknown{\reals^\pm} & \isknown{(-\infty,2)} & \isknown{\reals} & {\mathbb I}\\[.3cm]
L_n & 
 \isknown{\reals} & . & \isknown{\reals} &
\isknown{\reals^+} &\reals^- & .  &
\isknown{\reals^\pm} & \isknown{\reals^-}\cup\reals^+& \isknown{\reals} & (-\infty,1)\\[.3cm]
P_n & 
{\mathbb I} & . & \reals^\pm &
{\mathbb I}^+& {\mathbb I}^+ & \reals\setminus{\mathbb I}  & \reals^\pm
&\isknown{\reals^\pm}&\isknown{\reals} & {\mathbb I}\\[.3cm] 
T_n &
\isknown{{\mathbb I}} & . & \reals^\pm & {\mathbb I}^+ & {\mathbb I}^+ & \reals^+\cup\reals^-  & 
\reals^\pm &\reals^\pm &\isknown{{\mathbb I}} & {\mathbb I}\\[.3cm]
A_n & 
\isknown{{\mathbb I}^-} &(1,\infty) &{\mathbb I} &. &. & .  &\reals & {\mathbb I} & \reals^- &
{\mathbb I}^- \\[.3cm] 
C_n^0 &
\isknown{\reals^+} & \isknown{\reals^-} &
\isknown{{\mathbb I}^+} & . & .& (1,\infty)  &
\reals^{+}  & {\mathbb I}^+ &
(-\infty,-1) & {\mathbb I}^+\\
&&&&&&&&& \cup\, \reals^+ & \\[.3cm] 
G_n^{(2)} &
{\mathbb I} & . & (-\infty,1) & {\mathbb I}^+ & {\mathbb I}^+ & \reals  &
\reals^- & (-\infty,1) 
& \reals & {\mathbb I} \\[.3cm]
\falling{x}{n} & 
\isknown{\reals^+} & \isknown{\reals^- } & \isknown{{{\mathbb I}^+}}
& . & . & (1,\infty)  &\reals^+& \isknown{{{\mathbb I}^+}} &
\isknown{(-\infty,-1)} & {\mathbb I}^+\\  
&&&&&&&&&\cup\,\isknown{\reals^+} & \\[.3cm]
\rising{x}{n} & 
\isknown{\reals^-} & \isknown{\reals^+ } & {{\mathbb I}^+}
& . & . & (-\infty,-1)  &\reals^+& \isknown{{\mathbb I}^+} & \isknown{\reals^-} & {\mathbb I}^-\\
&&&&&&&&&\cup\,\isknown{(1,\infty)} & \\[.3cm]
\binom{x}{n} &
\reals\setminus{\mathbb I}^-&{\mathbb I}^- & .&
(1,\infty) &. & {\mathbb I}^-  &
\reals^\pm & . & \isknown{(-\infty,-1)} & \isknown{\reals^+}\\[.3cm]
&&&&&&&&&\cup\,\isknown{\reals^+} & \\[.3cm]
(x;2) &
\isknown{{\mathbb I}^+} & \isknown{\reals\setminus{\mathbb I}^+ } & \reals^+ & . & . &
\reals^-  & \reals^+ & \reals^+ 
& . & {\mathbb I}^+\\[.3cm]
(x;1/2) & 
\isknown{\reals\setminus{\mathbb I}^+ } & \isknown{{\mathbb I}^+} & . & (1,\infty) & . & .
& \reals^\pm & 
. &\reals-{\mathbb I}^+ &
\reals^-\\[.3cm] 
Exp &
\isknown{\reals } & . & . & \isknown{\reals^+ } & \isknown{\reals^+
} & .  &.&.&\isknown{\reals } & \isknown{\reals}\\[.3cm]
2^{\binom{n}{2}}x^n &
. & \isknown{\reals} & . & . & . & \reals  & . & . & . & .  \\[.3cm]
(1/4)^{\binom{n}{2}}x^n &
\isknown{\reals} & . & . & \reals^+ & \reals^+ & .  & . & . & \isknown{\reals} & \isknown{\reals}  \\[.3cm]
(1/2)^{\binom{n}{2}}x^n &
\isknown{\reals} & . & . & \isknown{\reals^+} & \isknown{\reals^+} & .  & . & . &
\isknown{\reals} & .\\[.3cm] 
(-1)^{\binom{n}{2}}x^n & 
\isknown{\reals^\pm } & \isknown{\reals^\pm} & . & \isknown{\reals^+} &
. & \reals  & . & . &\isknown{\reals^\pm } & .\\[.1cm]
\text{Laurent-1} &
\mathbb{I}^+ & . & \reals-\mathbb{I}^+ & . & . & . &
\reals^+\cup\reals^- & \reals-\mathbb{I}^+ & \reals^+ & .\\[.3cm]
\text{Laurent-2} &
\isknown{\reals^+} & . & \isknown{\reals^+}\cup\reals^- & . & . & . &
\reals^+\cup\reals^- & \reals^+\cup\reals^- & \reals^+ & .\\[.3cm]
\text{Laurent-3} &
\isknown{\reals^+} & . & \mathbb{I} & . & . & \mathbb{I}^- &
\reals^+\cup\reals^- & \mathbb{I} & \reals^+ & . \mathbb{I}^+ \\[.1cm]\hline
\end{array}
$

\pagebreak

\[
  \begin{array}{c|ccccccccc|}
    \toprule
x^n\mapsto & \quada & \quadb  & \quadc  & \quadd & \stabled{1} & \stabledc{1} & \polycpx &
\allpoly & \imag\allpoly  \\
\midrule
H_n &-\polycpx&-\polycpx&\polycpx&\polycpx&\allpoly&.&\polycpx&\isknown{\allpoly}&.
\\[.3cm]
L_n &{\quadb}&{\quadb}&{\quadc}&{\quadc}&{\stabled{1}}&{\stabledc{1}}&{\quadc}&\isknown{\allpoly}&{\stabledc{1}}\\[.3cm]
P_n &-\polycpx&-\polycpx&\polycpx&\polycpx&.&.&\polycpx&.&.
\\[.3cm]
T_n &-\polycpx&-\polycpx&\polycpx&\polycpx&.&.&\polycpx&.&.\\[.3cm]
A_n &.&.&.&.&\stabled{1}&.&.&.&.
\\[.3cm]
%
%
G_n^{(2)} &-\polycpx&-\polycpx&\polycpx&\polycpx&.&.&\polycpx&.&.
\\[.3cm]
\falling{x}{n} &.&.&.&.&.&.&.&.&.\\[.3cm]
\rising{x}{n} &.&.&.&.&\stabled{1}&.&.&.&.
\\[.3cm]
\binom{x}{n} &.&.&.&.&.&.&.&.&.
\\[.3cm]
(x;2) &.&.&.&.&-\stabled{1}&.&.&.&.
\\[.3cm]
(x;1/2) &.&-\stabledc{1}&-\stabledc{1}&.&-\stabled{1}&-\stabledc{1}&.&.&-\stabledc{1}
\\[.3cm]
Exp &\isknown{\quada}&\isknown{\quadb}&\isknown{\quadc}&\isknown{\quadd}&\isknown{\stabled{1}}&\isknown{\stabledc{1}}&\isknown{\polycpx}&\isknown{\allpoly}&\isknown{\imag\allpoly} 
\\[.3cm]
2^{\binom{n}{2}}x^n &.&.&.&.&.&.&.&.&.
\\[.3cm]
(1/4)^{\binom{n}{2}}x^n &\quada&\quadb&\quadc&\quadd&\stabled{1}&\stabledc{1}&\polycpx&\allpoly&\imag\allpoly
\\[.3cm]
(-1)^{\binom{n}{2}}x^n &.&.&.&.&\allpoly&.&.&.&.
\\[.3cm]
Laurent-1 &-\polycpx&.&.&\polycpx&.&.&.&.&.
\\[.3cm]
Laurent-2 &-\polycpx&.&.&\polycpx&.&.&.&.&. 
\\[.3cm]
Laurent-3 &.&.&.&.&.&.&.&.&.
\\[.3cm]
\bottomrule
  \end{array}
\]


\chapter{Tables of determinants and integrals }
  
\renewcommand{\TimeStampStart}{Monday, January 07, 2008: 11:32:15}
\mytoday  

We assume that all polynomials in one variable have positive leading
coefficients. 

  $$
  \begin{array}{ccp{3in}c}
\twodet{a_k}{a_{k+1}}{a_{k-1}}{a_k} & < 0 & where $ \sum_0^n
    a_i x_i \in\allpoly$  has all distinct roots. If they are not all
    distinct, then it is still $<$ unless $a_k=0$ and
    $a_{k-1}a_{k+1}=0$. & \pageref{eqn:newton-1}\\[.4cm]
\twodet{a_i}{a_{i+1}}{b_i}{b_{i+1}} & < 0 & where $f = \sum a_i
x^i$, $g = \sum b_i x^i$,and $f \greateq g$.  & \pageref{eqn:log-con-coef}\\[.4cm]
\twodet{f}{g}{f'}{g'} & <0 & where $f\lessless g$ or $f\greateq g$. &
\pageref{lem:inequality-1}\\[.4cm]
\twodet{f}{g}{h}{k} & < 0 & where $deg(f) =
  deg(g)+1=deg(h)+1=deg(k)+2$ and for all all  $\alpha,\beta$, not both
  zero, we have that $\alpha f+\beta g \lessless \alpha h+\beta k$ &
  \pageref{lem:inequality-4} \\[.4cm]
\twodet{f}{\affa f}{\affa f}{\affa^2 f} & < 0 & for positive $x$
if $\affa x > x$ and $f\in\allpolypos\cap\allpolyaffine$ &
\pageref{lem:affine-det} \\[.4cm]
\twodet{f}{\affa f}{\affa f}{\affa^2 f} & > 0 & for positive $x$
if $\affa x < x$ and $f\in\allpolypos\cap\allpolyaffine$ &
\pageref{lem:affine-det} \\[.4cm]
\twodet{f}{\Delta f}{\Delta f}{\Delta^2 f} & < 0 & for positive $x$
if $f\in\allpolypos\cap\allpolysep$ &
\pageref{lem:affine-det} \\[.4cm]
\twodet{f(x)}{f(x+1)}{f(x+1)}{f(x+2)} & <0 &  for positive $x$
and  $f\in\allpolypos\cap\allpolysep$ &
\pageref{lem:affine-det-2} \\[.4cm]
\twodet{d_{ij}}{d_{i+1,j}}{d_{i,j+1}}{d_{i+1,j+1}} &  <0 & $f(x,y)
= \sum d_{ij}x^iy^j$ is in $\rupint{2}$ and $f(x,0)$ has all distinct
roots. & \pageref{eqn:product-4} \\[.4cm]
\twodet{f}{f_x}{f_y}{f_{xy}} & \le 0 & $f\in\rupint{d}$, $d\ge2$. & 
\pageref{cor:product-4d-1} \\[.4cm]
\twodet{f_{00}}{f_{10}}{f_{01}}{f_{11}} & \le 0 & $f\in\rupint{d+2}$
and $f(\xx,y,z) = f_{00}(\xx) + f_{10}(\xx)y + f_{01}(\xx)z +
f_{11}(\xx)yz + \cdots $ & 
\pageref{cor:product-4d} \\[.4cm]
\twodet{1/f(x_1,y_1)}{1/f(x_2,y_1}{1/f(x_2,y_1)}{1/f(x_2,y_2)} & > 0
& $f\in\gsubpos_2$ &\pageref{lem:positive-def-2} \\[.4cm]
  \end{array}
$$

Next, we have examples of three by three determinants, all of which
are special cases of $n$ by $n$ determinants. All these matrices are
totally positive for non-negative $x$.

\begin{xalignat*}{2}
  \begin{vmatrix}
a_2 & a_1 & a_0 \\[.2cm]a_3 & a_2 & a_1 \\[.2cm]a_4 & a_3 & a_2
  \end{vmatrix}\ 
 &&\quad \sum a_i x^i \in \allpolypos 
&\ \quad\pageref{thm:tp}
\\[.4cm]
  \begin{vmatrix}
f_2 & f_1 & f_0 \\[.2cm]f_3 & f_2 & f_1 \\[.2cm] f_4 & f_3 & f_2
  \end{vmatrix}\  
 &&\quad \sum f_i(x)y^i \in \gsubpos_2 
&\ \quad\pageref{cor:karlin-1}
\\[.4cm]
  \begin{vmatrix}
   \frac{1}{2!} f'' & f' & f \\[.2cm]
    \frac{1}{3!}f^{(3)} & \frac{1}{2!}f'' & f' \\[.2cm]
    \frac{1}{4!}f^{(4)} & \frac{1}{3!}f^{(3)} & \frac{1}{2!}f'' 
  \end{vmatrix}\  
&&\quad f \in \allpolypos 
&\ \quad\pageref{cor:karlin-1}
\\[.4cm]
  \begin{vmatrix}
    \frac{1}{2!}f(x+2) & f(x+1) & f(x) \\[.2cm]
    \frac{1}{3!}f(x+3) & \frac{1}{2!}f(x+2) & f(x+1) \\[.2cm]
    \frac{1}{4!}f(x+4) & \frac{1}{3!}f(x+3) & \frac{1}{2!}f(x+2)
  \end{vmatrix}\  
&&\quad f \in \allpolypos 
&\ \quad\pageref{cor:tp-2}
\\[.4cm]
  \begin{vmatrix}
   \frac{1}{2!} (x+2)^n & (x+1)^n & x^n \\[.2cm]
   \frac{1}{3!} (x+3)^n & \frac{1}{2!}(x+2)^n & (x+1)^n \\[.2cm]
   \frac{1}{4!} (x+4)^n & \frac{1}{3!}(x+3)^n & \frac{1}{2!}(x+2)^n
  \end{vmatrix}\  
&& 
&\ \quad\pageref{cor:tp-2}
\\
\end{xalignat*}

\vfill

\pagebreak

\noindent%
\hspace*{-2cm}%
\begin{tabular}{lp{2.5in}l}
$\displaystyle\int_{-1}^1 f(x+\imag  t)\,dt\in\allpoly$. &
 $f\in\allpoly$ &
  \pageref{cor:int-of-i}\\[.6cm]
 $\displaystyle\int_0^1 f(x+\alpha\imag t)\,dt \in\polycpx$ &
 $f\in\allpoly$ and $\alpha>0$ &
\pageref{lem:polycpx-basic} \\[.6cm]
 $\displaystyle\int_0^1 f(x+\imag t)\,dt \in\polycpx$ &
 $f\in\allpoly$  &
\pageref{lem:i-integral} \\[.6cm]
 $\displaystyle\int_0^1 f(x+ t)\,dt \in\stabled{1}$ &
 $f\in\stabled{1}$  &
\pageref{lem:h-integral} \\[.6cm]
$ \displaystyle\int_0^1 \phi(x) e^{\imag t x}\,dx \in\polycpxf.$ &
$\phi(x)\uparrow$  and positive  on $(0,1)$ &
 \pageref{cor:phi-int}\\[.6cm]
 $\displaystyle\int_0^1f(x,y+\imag t)\,dt\in\polycpx_2$ &
 $f\in\rupint{2}$ &
 \pageref{lem:gsub-sub-i}\\[.6cm]
$\displaystyle\int_0^1  \,f(x,t)\,dt\in\allpoly$ &
$f\in\gsubsep_2$ &
\pageref{prop:family-int}\\[.6cm]
$\displaystyle \int_0^1 f(x,t)g(x,-t)\, dt\, \in\allpoly$ &
  $f,g\in\allpolysep_2$   &
\pageref{lem:convolution-int} \\[.6cm]
$\displaystyle\int_0^1 \,\int_0^{1-t} \,f(x,t)g(x,s)\,ds\,dt\in\allpoly$ &
$f,g\in\gsubsep_2$ &
\pageref{cor:int-fam-3}\\[.6cm]
\end{tabular}


 \chapter{Tables of Polynomials}
\label{cha:tables}
 
\renewcommand{\TimeStampStart}{Tuesday, June 05, 2007: 14:45:23}
\mytoday

\onepoly{Appell polynomials in one variable}{$A_n(x)$}
{A_0} (x)&= 1\\
 {A_1} (x)&= 1 - 2\, x\\
 {A_2} (x)&=     1 - 6\, x + 6\, x^2\\
{A_3} (x)&= 1 - 12\, x + 30\, x^2 - 20\,     x^3 \\
{A_4} (x)&= 1 - 20\, x + 90\, x^2 - 140\, x^3 + 70\, x^4\\ 
{A_5} (x)&= 1 - 30\, x + 210\, x^2 - 560\, x^3 + 630\,     x^4 - 252\, x^5
\end{array}
$

\onepoly{Appell polynomials in two variables}{$A_{n,m}(x)$}
A_{0,0} (x,y)&= 1\\
 A_{1,0} (x,y)&= 1 - 2\, x - y  \\
 A_{0,1} (x,y)&= 1 - x - 2\, y\\
 A_{2,0} (x,y)&= 1 - 6\, x + 6\, x^2 - 2\, y + 6\, x\, y + y^2\\
 A_{1,1} (x,y)&= 1 - 4\, x + 3\,  x^2 - 4\, y + 8\, x\, y + 3\, y^2\\
 A_{0,2} (x,y) & = 1 - 2\, x + x^2 - 6\, y + 6\, x\, y + 6\, y^2
\end{array}
$

\onepoly{Charlier with parameter 1}{$C_n^1(x)$}
C_0 (x)  =& 1 \\
 C_1 (x)  =& -1 + x \\
 C_2 (x)  =&     1 - 3 x + x^2 \\
 C_3 (x)  =& -1 + 8 x - 6 x^2 + x^3 \\
 C_4 (x)  =& 1 - 24 x + 29 x^2 - 10 x^3 + x^4 \\
 C_5 (x)  =& -1 + 89 x - 145 x^2 + 75 x^3 - 15  x^4 + x^5 
\end{array}
$

\onepoly{Charlier in two variables}{$C_n(x;y) = (-1)^n C_n^y(-x)$}
C_0 (x;y)  =&  1\\ 
 C_1 (x;y)  =&  x + y \\
 C_2 (x;y)  =&  x + x^2 + 2 x y + y^2 \\
     C_3 (x;y)  =&  2 x + 3 x^2 + x^3 + 3 x y + 3 x^2 
    y + 3 x y^2 + y^3 \\
 C_4 (x;y)  =&  6 x + 11 x^2 + 6 
    x^3 + x^4 + 8 x y + 12 x^2 y + 4 x^3 y + 6 x y^2 + 6 
    x^2 y^2 + 4 x y^3 + y^4\\ 
 C_5 (x;y)  =&  24 x + 50 
    x^2 + 35 x^3 + 10 x^4 + x^5 + 30 x y + 55 x^2 y + 30 x^3 
    y + 5 x^4 y + 20 x y^2 + 30 x^2 y^2 + 10 x^3 y^2 + 10 
    x y^3 + 10 x^2 y^3 + 5 x y^4 + y^5 
\end{array}
$

\onepoly{Chebyshev}{$T_n(x)$}
 T _  0 (x) =&  1   \\
 T _  1 (x) =&  x    \\
 T _  2 (x) =&  -1 + 2 x^2    \\
 T _  3 (x) =&  -3 x + 4 x^3    \\
 T _  4 (x) =&  1 - 8 x^2 + 8 x^4    \\
 T _  5 (x) =&  5 x - 20 x^3 + 16 x^5  
\end{array}
$

\onepoly{Euler}{$A_n(x)$}
A_0(x) =& 1\\
A_1(x) =& x\\
A_2(x) =& x+x^2\\
A_3(x) =& x + 4 x^2 + x^3\\
A_4(x) =&  x + 11 x^2 + 11 x^3 + x^4 \\
A_5(x) =&     x + 26 x^2 + 66 x^3 + 26 x^4 + x^5
\end{array}
$

\onepoly{Falling factorial}{$\falling{x}{n}$}
\falling{x}{0}   &= 1 \\
 \falling{x}{1}   &= x     \\
 \falling{x}{2}   &= -x + x^2 \\
 \falling{x}{3}   &=     2\ x - 3\ x^2 + x^3 \\
 \falling{x}{4}   &= -6\ x + 11\     x^2 - 6\ x^3 + x^4 \\
 \falling{x}{5}   &= 24\ x - 50\     x^2 + 35\ x^3 - 10\ x^4 + x^5 \\
\end{array}
$

\onepoly{Falling factorial inverse (Exponential or Bell polynomials)}{$B_n(x)$}
 B_0(x) &= 1 \\
 B_1(x) &= x \\
 B_2(x) &= x + x^2 \\
 B_3(x) &=     x + 3\ x^2 + x^3 \\
 B_4(x) &= x + 7\ x^2 + 6\     x^3 + x^4 \\
 B_5(x) &=  x + 15\ x^2 + 25\ x^3 + 10\     x^4 + x^5 \end{array}
$

\onepoly{Gegenbauer}{$G_n(x)$}
G_0(x) &= 1 \\
G_1(x) &= 4 x \\
G_2(X) &= 12 x^2-2\\
G_3(x) &= 32 x^3-12 x\\
G_4(X) &= 80 x^4-48 x^2+3 \\
G_5(x) &= 192 x^5-160 x^3+24 x
\end{array}
$

\onepoly{Hermite}{$H_n(x)$}
 H _  0 (x) =&  1 \\
  H _ 1 (x) =&  2 x   \\
      H _  2 (x) =&  -2 + 4 x^2   \\
      H _  3 (x) =&  -12 x + 8 x^3   \\
      H _  4 (x) =&  12 - 48 x^2 + 16 x^4   \\
      H _  5 (x) =&  120 x - 160 x^3 + 32 x^5 
\end{array}
$

\onepoly{Laguerre}{$L_n(x)$}
 L _  0 (x)=&  1\\
 L _ 1  (x) =&  1 - x   \\
 L _  2(x) =&  \frac{1}{2}\left(   {2 - 4x + x^2}\right)\\[.1cm]
 L _  3(x) =&  \frac{1}{6}\left(   {6 - 18x + 9x^2 - x^3}\right)\\[.1cm]
 L _  4(x) =&  \frac{1}{24}\left(   {24 - 96x + 72x^2 - 16  x^3 + x^4}\right)\\[.1cm]
 L _  5(x) =&  \frac{1}{120}\left( {120 - 600x + 600x^2 - 200x^3 + 25  x^4 - x^5}\right)
\end{array}
$

\onepoly{$Q$-Laguerre}{$L_n^\affa(x)$}
 L^\affa _  0(x)=&  1\\
 L ^\affa_  1(x) =&   - x +1   \\
 L ^\affa_  2(x) =&  x^2 - 2[2]x+[2]\\
 L ^\affa_  3(x) =&  -x^3 + 3[3]x^2 - 3[2][3]x + [2][3]\\
 L ^\affa_  4(x) =&  x^4 -4[4] x^3 + 6[3][4]x^2 - 4[2][3][4]x + [2][3][4]\\
 \end{array}
$

\onepoly{Laguerre in two variables}{$L_n(x;y) = L_n^y(-x)$}

 L _ 0(x) =& 1\\
  L _ 1 (x; y) =& 1 + x + y  \\
 L _ 2\, (x; y) =& \frac{1}{2}\bigl(2 + 4\, x + x^2 + 3\, y + 2\, x\, y + y^2\bigr)\\[.1cm]
 L _ 3\, (x; y) =& \frac{1}{6}\bigl(6 + 18\, x + 9\, x^2 + x^3 + 11\, y + 15\, x\, 
 y + 3\, x^2\, y + 6\, y^2 + 3\, x\, y^2 + y^3\bigr)\\[.1cm]
 L _ 4\, (x; y) =& \frac{1}{24}\bigl(24 + 96\, x + 72\, x^2 + 16\, x^3 + x^4 + 50\, 
 y + 104\, x\, y + 42\, x^2\, y  + 4\, x^3\, y + 35\, y^2 + 36\, 
 x\, y^2 + 6\, x^2\, y^2 + 10\, y^3 + 4\, x\, 
 y^3 + y^4\bigr)\\[.1cm]
 L _ 5\, (x; y) =& \frac{1}{120}\bigl(120 + 600\, x + 600\, x^2 + 200\, x^3 + 25\, 
 x^4 + x^5 + 274\, y + 770\, x\, y  + 470\, x^2\, y + 90\, x^3\, 
 y + 5\, x^4\, y + 225\, y^2 + 355\, x\, y^2 + 120\, x^2\, 
 y^2 + 10\, x^3\, y^2 + 85\, y^3 + 70\, x\, y^3 + 10\, x^2\, 
y^3 + 15\, y^4 + 5\, x\, y^4 + y^5\bigr)
\end{array}
$

\onepoly{Legendre}{$P_n(x)$}
P_0 (x)  =& 1 \\[.1cm]
 P_1 (x)  =& x     \\[.1cm]
 P_2 (x)  =& \frac{1}{2} \bigl(-1  + 3 x^2\bigr) \\[.1cm]
 P_3 (x)  =& \frac{1}{2}\bigl(-3 x + 5\ x^3\bigr)     \\[.1cm]
 P_4 (x)  =& \frac{1}{8}\bigl(3 - 30\ x^2 + 35x^4\bigr)\\[.1cm]
 P_5 (x)  =& \frac{1}{8}\bigl(15x-70x^3+63x^5\bigr)
\end{array}
$

\onepoly{Narayana}{$N_n(x)$}
N_0(x) =&1 \\[.1cm]
N_1(x) = & x\\[.1cm]
N_2(x) = &x^2+x\\[.1cm]
N_3(x) = &x^3+3 x^2+x\\[.1cm]
N_4(x) = &x^4+6 x^3+6 x^2+x\\[.1cm]
N_5(x) = &x^5+10 x^4+20 x^3+10 x^2+x
\end{array}
$

\onepoly{Rising factorial}{$\rising{x}{n}$}
\rising{x}{0}   &= 1 \\
 \rising{x}{1}   &= x      \\
 \rising{x}{2}   &= x + x^2 \\
 \rising{x}{3}   &= 2\     x + 3\ x^2 + x^3 \\
 \rising{x}{4}   &= 6\ x + 11\ x^2 + 6\     x^3 + x^4 \\
 \rising{x}{5}   &= 24\ x + 50\ x^2 + 35\     x^3 + 10\ x^4 + x^5 \\
\end{array}
$

\onepoly{Rising factorial inverse}{$\rising{x}{n}^{-1}$}
\rising{x}{0}^{-1}   &= 1 \\
 \rising{x}{1}^{-1}   &= x     \\
 \rising{x}{2}^{-1}   &= -x + x^2 \\
 \rising{x}{3}^{-1}   &=     x - 3\ x^2 + x^3 \\
 \rising{x}{4}^{-1}   &= -x + 7\ x^2 - 6\     x^3 + x^4 \\
 \rising{x}{5}^{-1}   &= x - 15\ x^2 + 25\ x^3 - 10\     x^4 + x^5 \\
\end{array}
$


\chapter{Questions}
There are many unsolved problems about  polynomials.  
  
\renewcommand{\TimeStampStart}{Friday, January 18, 2008: 10:05:00}
\mytoday

\subsection*{Polynomials with all real roots}


\begin{question}
    Is the minimum of the roots of the Hadamard product of polynomials
    in $\allpolyint{(0,1)}(n)$ attained at $f=g=(x+1)^n$?
  \end{question}

\begin{question} \label{ques:qxn}
  Show that for all $n$ the falling factorial
  \index{falling factorial!question} $ \falling{x}{n}$ is the derivative of a polynomial
  with all real roots. (See \ref{ques:qxn2})
\end{question}


\begin{question}
  Find all polynomials $h(x,y)$ such that for all $f\greateqeq g$ we
  have $h(f,g)\in\allpoly$.
\end{question}

\begin{question}
  Show that $f_n,g_n$ defined below satisfy $f_n \longleftarrow
  f_{n-1}$ and that $g_n \longleftarrow g_{n-1}$, where
  $\binom{n}{i,i}$ and $\binom{n}{i,i,i}$ are multinomial
  coefficients. In addition, $f_n\lessless f_{n-2}$ and $g_n \lessless
  g_{n-3}$. (There is the obvious generalization.)  Corollary~\ref{cor:diag-binom}
  shows that $f_n\in\allpoly$.
  \begin{xalignat*}{2}
    f_n &= \sum_{i=0}^n \binom{n}{i,i} x^i &
    g_n &= \sum_{i=0}^n \binom{n}{i,i,i} x^i 
  \end{xalignat*}
\end{question}

\begin{question}
  Suppose $n$ is a positive integer. Show the following interlacings,
  where $B,L,N,E$ are the Bell, Laguerre, Narayana and Euler polynomials.
  \begin{align*}
    \exp^{-1} B_{n+1} & \longleftarrow \exp^{-1} B_{n} \\
    \exp^{-1} L_{n+1} & \longleftarrow \exp^{-1} L_{n} \\
    \exp^{-1} N_{n+1} & \longleftarrow \exp^{-1} N_{n} \\
    \exp^{-1} E_{n+1} & \longleftarrow \exp^{-1} E_{n} 
  \end{align*}
\end{question}


\begin{question}
  Under what conditions does $0\le\aaa\le\bbb$ imply there is a $g$
  such that $\ppoly(\ccc)\le g$ for all $\aaa\le\ccc\le\bbb$?  
\end{question}

\begin{question}
  (C. Johnson) Given $f,g\in\allpoly$, not necessarily interlacing,
  what is the maximum number of complex roots of linear combinations
  of $f$ and $g$? Can this number be computed from the arrangement of
  the roots of $f$ and $g$?
\end{question}

\begin{question}[The \fdb\ problem - \seepage{sec:fdb-problem}]

Suppose that $g\in\allpolyf$. For any $\alpha\in\reals$ and
$m=0,1,\dots,$
\begin{enumerate}
\item  $F_m(\alpha,y)$ has all real roots.
\item $F_{m+1}(\alpha,y)\lesslesseq F_m(\alpha,y)$.
\end{enumerate}
  
\end{question}

\subsection*{Polynomials with all positive coefficients}

\begin{question}  \label{ques:signs}
  If $T\colon{}\allpolypos\longrightarrow\allpoly$, then what constraints are
  there on the signs of $T(x^i)$? The example $x^i\mapsto
  (-1)^{\binom{i}{2}}x^i$ shows that the pattern $--++$ is possible.
  In addition to all signs the same, and all signs alternating, are there
  any other possibilities?
\end{question}

\begin{question}
  If $\sum a_ix^i\in\allpolypos$ and $\alpha>1$ show that $\sum
  a_i^\alpha x^i\in\allpolypos$. See \cite{gregor98} for a
  questionable argument. Compare Question~\ref{ques:aei}.
\end{question}


\begin{question}
  If $f\in\allpolypos$,   consider the set $S=\{g\mid f\plesslesseq
  g\}$. Unlike the case where $\plesslesseq$ is replaced by
  $\lesslesseq$, the set $S$ is not a cone. However, it is a union of a
  collection of cones. What is a usefful description of $S$?
\index{positive interlacing}
\end{question}

\subsection*{Matrices that preserve interlacing}

\begin{question}
    If $M=(f_{ij})$ is a matrix of polynomials with positive leading
    coefficients that preserves mutually interlacing polynomials then
    does $(f_{ij}^{\prime})$ also preserve mutual  interlacing? If so, the
    following would be true:
    \begin{quote}
      Write $M = \sum M_i\, x^i$ where each $M_i$ is a matrix of
      constants. Each $M_i$ is totally positive$_2$.
    \end{quote}
  \end{question}

\begin{question}
  If $(f_1,\dots,f_n)$ is a sequence of mutually interlacing polynomials,
  then is $\{(a_i)\in\reals^n\vert \sum a_i f_i\in\allpoly\}$ a convex
  set?
\end{question}

\begin{question}
  Describe all matrices $M$ of polynomials with the property that if
  $v$ is a vector of mutually interlacing polynomials then the vector
  $Mv$ also consists of mutually interlacing polynomials.
\end{question}

\begin{question} 
\index{polynomials! mutually interlacing} \index{mutually interlacing!question} 
  A question related to the previous one. If $v$ is fixed, describe
  all matrices $M$ such that $M^iv$ is mutually interlacing, for all
  positive integers $i$. 
  
  A special case. Find all $M=\smalltwo{f}{g}{h}{k}$ such that
$M^i\smalltwobyone{1}{1}$ consists of interlacing polynomials for all
positive integers $i$.
\end{question}

\begin{question}
  If $f \lessless g$ and $fF-gG=1$ where $deg(F)=deg(G)-1=deg(f)-2$
  then are $F$ and $G$ unique?
\end{question}

\begin{question}
  Find all $d$ by $d$ matrices that are $\allpoly$-positive
semi-definite. That is, if $f\in\allpoly$ is $\sum a_ix^i$, and  $A =
(a_0,\dots,a_{d-1})$, then $A Q A^t \ge 0$.

  Are all \allpolypsd\ matrices the sum of a positive semi-definite
  matrix and a \allpolypsd\ matrix whose non-zero entries are only on
  the anti-diagonal?

\index{P-positive semi-definite@$\allpoly$-positive semi-definite}
\end{question}

\subsection*{Homogeneous Polynomials}

\begin{question}
  Given $f\in\allpoly$, are there  $a,b,c\in \reals$ and $g\in\allpoly$ such that
  $\left(a+b\frac{\partial}{\partial x} + c\frac{\partial}{\partial y}\right)
g=f$? 
\end{question}

\begin{question}
  If $T\colon{}\allpolypos\longrightarrow\allpolypos$ then when is
  $\varphi(T)$ totally positive? It appears to be so for the Laguerre
  transformation. 
\end{question}

\subsection*{Analytic functions}

\begin{question}
  Which hypergeometric functions are in $\allpolyf$?
\end{question}

\begin{question}\index{sin and cos}
  Which polynomials $f(x,y)$ have the property that $f(x,\sin
  x)\in\allpolyf$? How about $f(x,e^x)\in\allpolyf$?
\end{question}


\begin{question}
  Is there a $q$-analog $\Gamma_q$ of the gamma function for which
  $\Gamma_q(z+1)/\Gamma_q(kz+1)$ is entire, and only has negative zeros?
\end{question}

\begin{question}
  If $f\in\allpolyf$, and for all $a\in\reals$ we know that
  $f+a\in\allpolyf$ then show that $f = ax+b$. This is trivial when
  $f$ is a polynomial.
\end{question}

\subsection*{Linear Transformations of polynomials}

\begin{question}
  Characterize all linear transformations $T$ (see \eqref{eqn:tprod})
  defined by $ T(x^n) = \prod_{i=1}^n (x+a_i)$ that map $\allpolypos$
  to $\allpoly$.
\end{question}

\begin{question}
  Can we characterize transformations $T$ that satisfy
  \begin{alignat*}{2}
    T\colon{}\quad& \allpolypos\quad & \longrightarrow\quad & \allpolyneg \\
    T^{-1}:\quad& \allpolyalt\quad & \longrightarrow\quad & \allpolyalt 
  \end{alignat*}
\end{question}

\begin{question}
  Can we characterize transformations $T$ that  satisfy
  $T\colon{}\allpolyalt\longrightarrow\allpolypos$ and $T^2$ is the identity?
  The Laguerre transformation is such an example.
\end{question}

\begin{question} \index{Hermite polynomials!question}
  Characterize all linear transformations
  $T\colon{}\allpoly\longrightarrow\allpoly$ that commute with
  differentiation ($\diffd T = T \diffd$). 
\end{question}

\begin{question}
  Suppose that $V$ is the infinite dimensional vector space of all
  linear transformations from polynomials to polynomials, and $W$ is
  the subset of $V$ consisting of all linear transformations that map 
  $\allpolypos\longrightarrow\allpoly$. Then $W$ contains no open sets.
\end{question}

\begin{question}
  If $T(x^n) = d_nx^n$ satisfies  
  $T(\allpolypos)\subset\allpolypos$ then does $T$ preserve interlacing? 
\end{question}

\begin{question}
  If $T\colon{}\allpoly\longrightarrow\allpoly$ is a linear transformation
  then describe the possible sets
  $S=\{deg(T(x)),deg(T(x^2)),deg(T(x^3)),\cdots\}$.  Since $T(x^n)$
  and $T(x^{n+1})$ interlace, we know that adjacent terms of $S$
  differ by at most $1$. The example $f\mapsto f(\diffd)g$ shows that
  $S$ can be bounded.
\end{question}


\begin{question}
  Can we find linear transformations $T$ such that
  \begin{enumerate}
  \item $T\colon{}\allpoly\longrightarrow\allpoly$.
  \item $T(x^i) = x^i$ for $0\le i < r$
  \item $T(x^r) \ne x^r$
  \end{enumerate}
$x^n\mapsto H_n(x)$  works for $r=2$.
\end{question}

\subsection*{Linear transformations that preserve roots}

\index{falling factorial!question} \index{rising factorial!question}
\begin{question} \label{ques:fact-1} Show that

  \begin{enumerate}
  \item $T(x^i)=\rising{x}{i}\, \rising{x}{n-i}$ maps 
    $\allpolypos  \longrightarrow \allpolyneg$.
  \item $T(x^i)=\rising{x}{i}\, \falling{x}{n-i}$ maps 
    $\allpolypos  \longrightarrow \allpoly$.
  \end{enumerate}

\end{question}

\begin{question}\label{ques:bell-polys}
  Lie and Wang \cite{liu-wang} showed that 
\[
\begin{vmatrix}
B_{n+1} & B_{n+2} \\  B_n & B_{n+1} 
\end{vmatrix}
\]
has all positive coefficients where $B_n$ are the Bell polynomials
\mypage{rem:charlier-comp}. Show that the determinant is actually stable.
\end{question}

\begin{question}\label{ques:hermite-1}
  Show that the linear transformation $x^i\mapsto H_i(x)x^{n-i}$ maps
  $\allpolyint{(-\infty,2)}(n) \longrightarrow\allpoly$.
\end{question}

\begin{question} \label{ques:herm-2}
  \index{Hermite polynomials!question} If $T(H_n) = x^n$, then $T$
  does not map $\allpoly$, $\allpolyalt$ or $\allpolypos$ to
  $\allpoly$. Show that $\expoper{}\circ T$ and $T\circ \expoper{}$
  both map $\allpoly$ to itself.
\end{question}

\begin{question}\label{ques:herm-q} \index{q-Hermite polynomials!question}
  Show that the transformation $x^k\mapsto H_n^q$ where $H_n^q$ is the
  $q$-Hermite polynomial maps
  $\allpolyint{(-1,1)}\longrightarrow\allpoly$ for $q>1$.
\end{question}

\begin{question} \label{ques:prod-1}
  Suppose $r_i> 1+ r_{i-1}$ are positive constants. Show that the
  transformation $x^n \mapsto \prod_{i=1}^n (x-r_i)$ maps
  $\allpolyint{(0,1)}\longrightarrow\allpoly$. 
\end{question}

\begin{question} \label{ques:leg-1}
\index{Legendre polynomials!question}\index{polynomials!Legendre}%
  Show that the transformation $x^i\mapsto \rev{P_i}$, the reverse of
  the Legendre polynomial, maps
  $\allpolyalt\cup\allpolypos\longrightarrow\allpolyneg$. The
  transformation satisfies
$$ T(x^{n+1}) = \frac{2n+1}{n+1}T(x^n) - \frac{x^2}{n+1}T( (x^n)')$$
\end{question}

\begin{question} \label{ques:leg-2}
\index{Legendre polynomials!question}\index{polynomials!Legendre}%
  Show that the transformation $x^i\mapsto {P_i}$, where $P_i(x)$ is
  the Legendre polynomial, maps
  $\allpolyint{(-1,1)}\longrightarrow\allpolyint{(-1,1)}$. 
\end{question}

\begin{question} \label{ques:jac-3}
\index{Jacobi polynomials!question}\index{polynomials!Jacobi}%
For which $\alpha,\beta$ is it the case that
$T\colon{}\allpoly\longrightarrow\allpoly$, where 
$T\colon{}x^i\mapsto \frac{{P^{\alpha,\beta}_i}}{i!}$, and
$P^{\alpha,\beta}_i(x)$ is the Jacobi polynomial?
\end{question}

\begin{question}
\index{Chebyshev polynomials!question}\index{polynomials!Chebyshev}%
    If $T$ is the linear transformation  $x^n\mapsto T_n$ and
  $Mz=\frac{z+1}{z-1}$ then
  $T_M:\allpolyint{(-1,1)}\rightarrow \allpolypos$. See
  Lemma~\ref{lem:chebyshev}. 
\end{question}

\begin{question} \label{ques:euler-1}%
\index{Eulerian polynomials!question}\index{polynomials!Eulerian}%
  Show that the transformation $E_i\mapsto x^i$, where $E_i$ is the
  Euler polynomial, maps
  $\allpolyint{(1,\infty)}\longrightarrow\allpolyalt$. 
\end{question}

\begin{question}\label{ques:diamond-exp}
  \index{diamond product!exponential-question}
  Show that the exponential diamond product $f\mydiamond{\expoper{}}g =
  \expoper{}^{-1}\left(\expoper{} f \times \expoper{} g\right)$ is a
  bilinear map $\allpolypos \times \allpolyalt \longrightarrow
  \allpoly$. This map is equivalent to $x^i \mydiamond{\expoper{}} x^j
  = \binom{i+j}{i} x^{i+j}.$
\end{question}

\begin{question}
  \index{diamond product!Laguerre-question} Show that if $T(x^i)=L_i(x)$, the
  Laguerre polynomial, then the Laguerre diamond product
  $f\mydiamond{T}g = T^{-1}\left(T(f) \times T(g)\right)$ is a
  bilinear map $\allpolyint{(1,\infty)} \times
  \allpolyint{(-\infty,1)} \longrightarrow \allpoly$.

\end{question}

\begin{question} \label{ques:diamond-falling}
  \index{diamond product!falling factorial-question}
  If $T\colon{}x^n\mapsto \falling{x}{n}$ then the diamond product
  $f\mydiamond{T}g = T^{-1}\left(T(f) \times T(g)\right)$ maps
  $\allpoly\times\allpolypos \longrightarrow \allpoly.$
  See Lemma~\ref{lem:diamond-falling}.
\end{question}

\begin{question}
  We can modify the definition of diamond product by replacing
  multiplication by the Hadamard product. Given a linear
  transformation $T$,  define  $S$ by

\centerline{
\xymatrix{
  f \times g \ar@{|->}[rr]^{T\times T} \ar@{|..>}[d]_S &&
Tf \times Tg \ar@{|->}[d]^{\text{Hadamard product}} \\
T^{-1}(Tf\ast Tg) && Tf\ast Tg \ar@{|->}[ll]_{T^{-1}}
}
}

\noindent%
\begin{enumerate}[(A)]
\item Let $T(x^n) =  \falling{x}{n}$. Show that
$S\colon{}\allpolypm\times\allpolypm\longrightarrow\allpoly$.
\item  Show that the above result holds if we use the modified Hadamard product
$x^r \ast^\prime x^r = r!\, x^r$.
\item If $T(H_i)=x^i$, \index{Hermite polynomials!question} 
then show that 
$$ S\colon{}\left(\allpolyint{(\infty,-1)}\cup\allpolyint{(1,\infty)}\right)
  \times
\left(\allpolyint{(\infty,-1)}\cup\allpolyint{(1,\infty)}\right)
  \longrightarrow \allpoly$$
\end{enumerate}
\end{question}

\begin{question} \label{ques:jacobi}
  If $T(x^n) = P_m^{\alpha,\beta}$, the $n$-th Jacobi polynomial, for
  which $\alpha,\beta$ does $T$ map $\allpolyint{-1,1}$ to itself?
  Some computations show that there are examples that do not map
  $\allpolyint{-1,1}$ to itself, but they are all where
  $\alpha\ne\beta$ and $-1< \alpha,\beta<0$.  
\index{Jacobi polynomials!question}
\end{question}

\begin{question} \label{ques:fact-q}
Show that the transformation $x^{i}\mapsto \prod(1+(-1)^i
x)=(x;-1)_i$ maps $\allpolypos\cup \allpolyalt\longrightarrow
\allpoly$. 
\end{question}

\begin{question} \label{ques:fact-q2}
  For which $a,b$ does the transformation $(x;a)\mapsto(x;b)$
  map $\allpoly\longrightarrow\allpoly$?
\end{question}


\begin{question}

What assumptions on $S,T$ do we need to make the following true?

  If $S,T$ are linear transformations that map $\allpolypos$ to
  itself, then define the transformation $U$ on $\allpolypos(n)$ by
  $x^i\mapsto T(x^i) S(x^{n-i})$.  If $S,T$ map $\allpolyalt$ to
  itself then $U$ maps $\allpolyalt(n)$ to itself and preserves
  interlacing.

\end{question}

\begin{question}
  Determine all orthogonal polynomial families $\{p_n\}$ so that
  $x^n\mapsto p_n$ maps $\allpolyint{\diffi}$ to $\allpoly$ for some
  interval $\diffi$.
\end{question}

\begin{question}
  Identify a polynomial $f$ in $\allpoly(n)$ with 
  $\alpha f$ in projective $n+1$ space. What is 
  $$\frac{ Vol\,\diffd{}\allpoly({n+1})}{Vol\, \allpoly(n)}$$
 \end{question}

 \begin{question}
   If $T\colon{}\allpoly\longrightarrow\allpoly$, then is it possible to
   compute
$\frac{Vol\,(T(\allpoly))}{Vol(\allpoly)}$? What are upper and lower
bounds for this ratio? Should we restrict ourself to $\allpoly(n)$?
 \end{question}

 \begin{question}
Let $all$ denote all polynomials.
   If $T\colon{}\allpoly\longrightarrow\allpoly$ and 
   $\frac{Vol\,(T(\allpoly))}{Vol(\allpoly)}$ is small
then $\frac{Vol\,(T(all)\cap\allpoly)}{Vol(\allpoly)}$ is large.
 \end{question}

\begin{question}
  Suppose that $T\colon{}\allpoly(n)\longrightarrow\allpoly(n)$. What are the
  obstructions to extending $T$ to a linear transformation that maps
  $\allpoly(n+1)$ to itself?
\end{question}

\begin{question}
  Show that if $T(x^i) = \falling{\alpha}{i}x^i$ where $\alpha<n-2$ is not an
  integer then $T((x+1)^n)\not\in\allpoly$.  See
  Lemma~\ref{lem:falling-non-int}.
\end{question}

\begin{question}
  When is it true that if $f\in\allpolypos$, then for all positive
  integers $i$
  $$f(i)^2> f(i-1)f(i+1).$$
  Corollary~\ref{cor:fi} shows that the inequality with
  $\ge$ in place of $>$ is true.
\end{question}

\begin{question}
  Suppose that $p_i$ is a sequence of polynomials defined by 
  \begin{align*}
    p_0 &= 1\\
    p_1 &= x\\
    p_{n+1} &= x p_n - c_n p_{n-1}
  \end{align*}
where the $c_n$ are constants. If the linear transformation
$x^n\mapsto p_n$ maps $\allpoly$ to itself  then are the $p_i$
essentially the Hermite polynomials? 
\end{question}

\begin{question}\label{ques:risingxk}
  If  $T(x^k) = \frac{\rising{x}{k}}{k!}$ then show that $T\colon{}
  \allpolyint{\reals\setminus(0,1)}\longrightarrow   
  \allpolyint{\reals\setminus(0,1)}$.
\end{question}

\begin{question}
  If  $T(x^k) = \frac{\rising{x}{k}}{k!}$ then show that for any
  $\alpha\in(0,1)$ there is an $n$ such that $T(x-\alpha)^n$ does not
  have all real roots. 

\end{question}

  \begin{question}
    If $T$ is a linear transformation, and $T(x-r)^n\in\allpoly$ for
    all positive integers $n$ and all $r$ satisfying $a\le r\le b$, then
    is it true that $T\colon{}\allpolyint{(a,b)}\longrightarrow\allpoly$?
  \end{question}

\begin{question}
Find examples of transformations $T$ and constant $a$ 
  for which $T(x+a)^n$ is a scalar, for all $n$.
This is satisfied by  $T(x^i) = \falling{x+n-i}{n}$
\end{question}

\begin{question}
  If $T\colon{}x^n \mapsto \rising{x}{n}/n!$ then show that $T$ maps
  polynomials whose roots all have real part $1/2$ to  polynomials with
  real part $1/2$. \cite{redmond}
\end{question}

\begin{question}
\index{Commutator}
In Example~\ref{ex:comm-narayan} we saw that the commutator of the two
transformations $x^n\mapsto N_n$ and $x_n\mapsto \diffd x^n$ maps
the sequence $\{x^n\}$ to an interlacing sequence. It appears that
more is true. An entry in  Table~\ref{tab:comm}, such as $\allpoly$ at the
intersection of the row $A_n$ and column $B_n$ means that empirical
evidence suggests that 
\[
T_1\colon x^n\mapsto A_n\qquad T_2\colon x^n\mapsto B_n\qquad \implies\quad 
T_1T_2-T_2T_1\colon\allpolypos\longrightarrow \allpoly
\]
\begin{table}\label{tab:comm}
  \centering
  \[
\begin{array}{cllllll}
    \toprule
    & N_n & A_n & B_n & \diffd & L_n & \rising{x}{n}  \\
    \midrule
    N_n & & \allpoly & \allpoly& \allpolypos & \allpolypos & . \\
    A_n & & & \allpoly & \allpolypos & \allpolypos & . \\
    B_n & & & & \stabled{} & \allpolypos & \allpoly \\
    \diffd & & & & & \allpolyint{(-\infty,0]} & \allpolypos \\
    L_n &&&&& & \allpolypos\\
    \bottomrule

  \end{array}
\]

  \caption{Purported range of $\allpolypos$ under the commutator map}
\end{table}
\end{question}

\subsection*{Affine Transformations}

\begin{question}
  Show that if $f(x)\greateqeq f(x+2)$ then 
     \[\smalltwodet{f(x+1)}{f(x)}{ f(x+2)}{f(x+1)} > 0\]
\end{question}


\begin{question}
  If $T(x^n)$ is $H_n(x),x^n/n!$, or $L_n(x)$ then show that the
  minimum separation is achieved at the falling factorial. That is, if
  $f\in\allpolysep$ show that
$$ \delta\, T(f) \ge\delta\, T\falling{x}{n}.$$

\end{question}

\begin{question}
  Is $\falling{x}{i}\mapsto \falling{x}{n-i}$ a bijection on $\allpolysep(n)$?
\end{question}

\begin{question}
  Show that $x^n\mapsto L_n$ ($L_n$ is the Laguerre polynomial) maps
  $\allpolysep$ to itself. \index{Laguerre polynomials!question}
\end{question}

\begin{question}
  Show that \index{rising factorial!question} $\rising{x}{n}\mapsto x^n$ maps
  $\allpolyalt$ to $\allpolysep$.
\end{question}

\begin{question}
  Suppose all $r_i$ are positive and define $$T(x^n) = \prod_{i=1}^n
  (x-r_i)$$
If $T(x^n) \in\allpolysep$ then 
  $T\colon{}\allpolyint{(0,1)}\longrightarrow\allpolyint{(1,\infty)}\cap\allpolysep$.
\end{question}

\begin{question}
  If $T\colon{}\allpolysep\longrightarrow\allpolysep$ satisfies $f\lesslesseq 
  Tf$ then is $Tf = \alpha f^\prime + \beta\Delta(f)$?
\end{question}

\begin{question}
  If $\delta(f)>1$ then is there a $g$ with all real roots such that
  $\Delta(g)=f$? 
\end{question}

\begin{question}
  Since $\prod_{i=0}^{2n-1}(x-2i)$ is a polynomial of degree $2n$ in
  $\allpolysep$, we know that
$$ p_n = \Delta^n \prod_{i=0}^{2n-1}(x-2i) \in\allpolysep(n)$$
Although it is not the case that $p_n \lesslesseq p_{n-1}$, it appears
to be nearly true. Show that if $g_n = p_n(x)/(x-\alpha)$ where
$\alpha$ is the largest root of $p_n$ then $p_{n-1} \greateqeq g_n$.
\end{question}

\begin{question}
  Which series $f(x) = \sum_0^\infty a_ix^i$ have the property that
  $f\ast g$ is in $\allpolyaffine$ if $g\in\allpolyaffine$? This is an
  affine analog of Theorem~\ref{thm:polya-schur}.
\end{question}

\begin{question}
  If $f$ has all its roots in $(-1,0)$ and $T\colon{}x^n\longrightarrow
  \rising{x}{n}$ then the roots of $Tf$ satisfy: the largest is
  larger than all roots of $f$, and the rest are all less than -1.
\end{question}

\begin{question} \label{ques:sep-bijection}
Show that the map $x^n\longrightarrow\rising{x}{n}$
is a bijection between $\allpolypos$ and  $\allpolysep\cap\allpolypos$.
\end{question}

\begin{question}
  For any $\affa$ is there an infinite family of polynomials $h_n$ so that
  \begin{itemize}
  \item $\polard{\affa}h_n = (const) h_{n-1}$
  \item $h_n\in\allpolyaffine$
  \end{itemize}
\end{question}

\begin{question}
  If $\{a_i\}$ is a sequence such that $g \lessless \sum a_ig^{(i)}$
  whenever $g$ is \index{falling factorial!question} $ \falling{x}{n}$ then do we have
  interlacing for all $g\in\allpolysep$?
\end{question}

\begin{question}
  If $f$ is a polynomial such that $f(\diffd_\affa)\allpolyaffine\subset
  \allpolyaffine$ then is $f\in\allpolyaffine$?
\end{question}

\begin{question}
  Does $Tf = x\,f(x/2) - \Delta(f)$ preserve roots in some appropriate 
  space?
\end{question}

\begin{question}
  If $f(x)\lesseq g(x)$, $f(x)\lesseq f(x+1)$, $g(x)\lesseq g(x+1)$
  then is $\smalltwodet{f(x)}{g(x)}{f(x+1)}{g(x+1)}\ne0$?
\end{question}

\begin{question}
  What are  conditions on $f$  that guarantee that  $\smalltwodet{f}{\Delta
  f}{\Delta f}{\Delta^2f}$ is never zero? This determinant is also
equal to $\smalltwodet{f(x)}{f(x+1)}{f(x+1)}{f(x+2)}$.
\end{question}

\begin{question}\label{ques:q-hermite}
Show that the positive part of the $q$-Hermite polynomials is in
$\allpolyaffine$. 
\end{question}

\begin{question}
\index{q-Hermite polynomials!question}
\index{polynomials!q-Hermite}

For which polynomials does the $q$-Hermite transformation $x^n\mapsto
H^q_n$ (see \chapsec{affine}{affine-q}) preserve real roots?
\end{question}

\begin{question}
 Recall $\delta(f)$ is the minimum distance between roots of $f$.
  What is $ \inf_f \frac{\delta(exp(f))}{\delta(f)}$?
It is 2 for quadratic and 1 for linear. 
\end{question}

\begin{question}
  What is $ \inf_f \frac{\delta(f\ast g)}{\delta(f)\delta(g)}$?
\end{question}

\begin{question}
  For which polynomials $F(x,z)$ is
 $F(x,\Delta):\allpolysep\longrightarrow\allpolysep$?  
How about $F(x,\affa):\allpolysep\longrightarrow\allpolysep$?  

\end{question}

\begin{question}
  Consider the recurrence
  $$
  p_{n+2} = (a_nx+b_n + (c_nx+d_n)\affa)p_{n+1} - e_n\affa p_n$$
  What conditions on the coefficients do we need so that all
  $p_n\in\allpolysep$, and $p_{n+1}\lesslesseq p_n$ for all $n$? 
\end{question}

\begin{question}
  Show that the real parts of the roots of $\Delta^k
  [\rising{x}{n}\rising{x}{n}]$ satisfy these two properties.
  \begin{enumerate}
  \item The real parts are all integers, or integers plus $1/2$.
  \item If $k=n+r$ where $r\ge -1$ then the real parts are all equal,
    and have value $-(2n+r-1)/2$.
  \end{enumerate}
\index{rising factorial}
\end{question}

\begin{question} \label{ques:sep-alt}
Show the following:

  The map $f\times g\mapsto f\septimes g$ satisfies 
$\allpolysep\cap\allpolyalt \times \allpolysep\cap\allpolyalt
\longrightarrow \allpolysep\cap\allpolyalt$.
\seepage{sec:affine-diamond}
\end{question}

\begin{question}
Find an explicit formula for the $q$-Laguerre
  polynomials. It is probably   similar to the definition of the
  Laguerre polynomials:
  \begin{align*}
    L_n(x) &= \frac{1}{n!}\sum_{i=0}^n
    \binom{n}{i}\frac{n!}{i!}\,(-x)^n \\
  \end{align*}

\end{question}

\begin{question}
  Suppose that $f_t,g_t,h_t$ are locally interlacing families with
  $\rho(f)\ge1$, $\rho(g)\ge1$, $\rho(h)\ge1$. Show that the triple
  convolution below is in $\allpoly$.

  $$\int_0^1 \int_0^{1-t} \, f_t(x) g_s(x) h_{1-s-t}(x)\,ds\,dt$$
\end{question}

\begin{question}
 Suppose that $f\in\allpolyaffine$ where $\affa x = qx$ with
 $q>0$. Show that if $f = \sum a_i x^i$ then (see \pageref{prop:newton-q})
\[
\frac{a_{k+1}^2}{a_ka_{k+2}} \ge \frac{1}{q}\,\frac{[k+2][n-k]}{[k+1][n-k-1]}.
\]

\end{question}

\begin{question}
  Suppose that $\affa x = qx$ where $q>1$, $\affb x = q^{1/2}x$ and
  $\daffine^{sym}f(x) = (\affa f - \affa^{-1}f)/(\affa x -
  \affa^{-1}x)$. Show that if $f\in\allpolyaffine$ then
  $\daffine^{sym} f\in\allpolyint{\affb}$.
\end{question}

\begin{question}\label{ques:affine-and-i}
  Suppose that $f \greateqeq \affa^2f$ and the  polynomial
  $f+\imag\affa f$ has roots $r_1+\imag s_1,\dots,r_n+\imag s_n$, and
  $g = (x-r_1)\cdots(x-r_n)$. Show that
\[ f \greateqeq g\greateqeq\affa f.
\]
(It then follows that $\affa f\greateqeq \affa g\greateqeq
\affa^2f$. Since $f\greateqeq \affa^2f$ the sequence $f,g,\affa f,\affa
g,\affa^2f$ is mutually interlacing, and hence $g\in\allpolyaffine$.)
\end{question}

\begin{question}
    Define the q-Charlier polynomials \cite{stanton-white} by 
\[
C_{n+1}(x,a;q) = (x-aq^n-[n])C_n(x,a;q) - a[n]q^{n-1}C_{n-1}(x,a;q),
\]
where $C_{-1}(x,a;q)=0$ and $C_0(x,a;q)=1$. If $q>1$ and $a>1$ then
show that $C_n(x,a;q)$ is in $\allpolyaffine$, where $\affa x = qx$.
  \end{question}

\subsection*{Polynomials in $\rupint{2}$}

  

\begin{question}
  If $f\in\gsubsep_2$ then show that $\Delta_x (f) = f(x+1,y)-f(x,y)$ is in
  $\gsubsep_2$. 
\end{question}

\begin{question}
  Show that the linear transformation $x^iy^j \mapsto
  \rising{x}{i}\rising{y}{j}$ maps $\rupint{2}$ to itself.
\end{question}


\begin{question}
Suppose that we have three polynomials $f,g,h\in\rupint{2}$ where
$af+bg+ch\in\pm\rupint{2}$ for all $a,b,c\in\reals$. Are $f,g,h$
necessarily linearly dependent?
\end{question}

\begin{question}
  If $f(x,y)$ is a product of linear factors, then for which linear
  transformations $T$ is $Tf$ ever a product of linear factors?
\end{question}

\begin{question}
  Suppose that $f\in\gsubpos_2$ and write
$$f = f_0 + f_1 y + f_2 y^2 + \dots + f_n y^n. $$  Consider the $d\times d$
  determinant

$$
det_d =
  \begin{vmatrix}
    f_0 & f_1 & \dots & f_{d-1} \\
    f_1 & f_2 & \dots & f_d \\
    f_2 & f_3 & \dots & f_{d+1} \\
    \vdots & & & \\
    f_{d-1} & f_d & \dots & f_{2d-2}
  \end{vmatrix}
$$

Empirical evidence suggests that
\begin{itemize}
\item If $d=2r$ is even and $d<n$ then $det_d$ always has the same
  sign $(-1)^r$.
\item If $d$ is odd then the determinant has positive and negative values.
\end{itemize}
This is similar to results for matrices of orthogonal polynomials
found in \cite{karlin-szego}. The case where $f$ is in $\gsubplus_2$
is discussed in \chapsec{homog}{poly-det}.
\end{question}

\begin{question} \textbf{(The total positivity
    conjecture)}\label{ques:tnnc} \index{polynomials!totally positive
    coefficients} \index{totally positive!question}
  Suppose that $f\in\gsubplus_2(n)$ is a product of linear terms, and
  write $f = \sum a_{ij}x^iy^j$. When is the matrix below  totally
  positive? (Lemma~\ref{lem:tp-fg} is a special case.)

$$
\begin{vmatrix}
  a_{0d} & \dots & a_{00} \\
\vdots & & \vdots \\
a_{dd} & \dots & a_{d0}
\end{vmatrix}
$$
 \end{question}

 \begin{question}
   We can form polynomials from the coefficients of a polynomial in
   $\gsubpos_2$.  Suppose that $f = \sum a_{i,j}x^iy^j$. Define
   $M(i,j;d)$ to be the determinant of the $d$ by $d$ matrix $(
   a_{i+r,j+s})$, where $0\le r,s<d$. For example, $M(i,j,1) =
   a_{i,j}$, and $M(i,j;2) = a_{i,j}a_{i+1,j+1} -a_{i,j+1}a_{i+1,j}$.

   When is $\sum M(i,j;d)x^iy^j\in\pm\gsubpos_2$?
 \end{question}

 \begin{question}
   Show that the product of two polynomials in $\gsubplus_2$ with
   totally positive coefficients has totally positive coefficients.
 \end{question}

 \begin{question}
 Suppose $f = \sum a_ix^i$ where
   $f\in\allpolypos(n)$. The polynomial $f(x+y)$ is in $\gsubpos_2$
   and
   \begin{align*}
     f(x+y) &= \sum_0^n a_i(x+y)^i \\
&= \sum_{i=0}^n \sum _{j=0}^i\binom{i}{j}a_ix^{i-j}y^j
   \end{align*}
Show that the following matrix is totally positive
\index{totally positive!question} 
\index{matrix!totally positive} 
$$
\begin{pmatrix}
  \dots & a_3 & a_2 & a_1 & a_0 \\
& \dots &3a_3 & 2a_2 & a_1 \\
& & \dots & 3a_3 & a_2 \\
&&&\dots & a_3 \\
&&&&\dots 
\end{pmatrix}
$$
This is a stronger statement than Theorem~\ref{thm:tp}.
 \end{question}


\begin{question}\label{ques:tot-pos-fct}
  Suppose that $f(x,y)\in\gsubf_2$ is positive whenever $x$ and $y$
  are positive. This is the case if $f\in\gsubplus_2$; it also holds
  for $e^{-xy-x^2}$ . When is $1/f$ is totally positive on
  $\reals_+^2$?
\end{question}

\begin{question}\index{permanent!question}
  The conjecture in \cite{haglund97} is that if $A=(a_{ij})$ is a real
  $n$ by $n$ matrix with non-negative entries which are weakly
  increasing down columns then the permanent of $A+xJ$ has all real
  roots. 

  Show that the permanents of all the submatrices of $A+xJ$ obtained
  by deleting any row and column have a common interlacing with $per(A+xJ)$.

\end{question}





\begin{question}
  If a linear transformation is defined on $\gsubpos_2$, and $T$ maps
  all products $\prod(y+a_ix+b_i)$ to $\gsubpos_2$, then does $T$ map 
  $\gsubpos_2$ to itself?
\end{question}

\begin{question}
  Which products have coefficients satisfying a recurrence?
That is, given a product $\prod(a_i + b_iy +c_i) = \sum f_i(x)y^i$ for 
which $a_i,b_i,c_i$ do we have $f_i = x\alpha_i 
f_{i-1} + \beta_i f_{i-2}$ for constants $\alpha_i$, $\beta_i$.
\end{question}


\begin{question}
  Given $f\in\allpoly(n)$  for which $r$ can we find a
  product of $n+r$ terms $ax +by+c$ so that the coefficient of $y^r$ is $f$?
  We can always do it for $r=1$, since all we need do is choose
  $g\lessless f$.
\end{question}

\begin{question}
  If $f_0\lessless f_1 \lessless f_2$ then is there a polynomial
$f_0 + f_1(x) y + f_2(x)y^2 + \cdots + f_n y^n$ that satisfies
{$x$-substitution}?

\end{question}

\begin{question}
  If $f\in\gsubpos_2$ and $\prod(a_ix+b_iy+c_i)\lesslesseq f$ then is \\ $f
  = \sum d_i \frac{f}{a_x+b_iy+c_i}$?
\end{question}

\begin{question}
  \index{Legendre polynomials!question} The Legendre polynomials can be defined
  by $P_n(x) = \frac{(-1)^n}{2^nn!}\diffd^n(x^2-1)^n$. This leads to
  the guess that the map
  $$
  x^iy^j \mapsto \alpha_{n,m}\diffd^i(x^2-1)^j$$ might map $\gsubpos_2$
  to $\allpoly$ for some constants $\alpha_{n,m}$.  This is not true
  when all $\alpha_{n,m}$ are equal to $1$. For instance, the image of
  \mbox{$(x+y+16)(x+4y+25)$} under this map has complex roots. Is
  there a choice of constants for which it is true?
\end{question}

  \begin{question}
    What are the different isotopy types  of graphs of polynomials in $\rupint{2}(n)$?.
  \end{question}

\begin{question}
  If $f\in\rupint{2}$  then let $p_n$ be an eigenpolynomial. This means
  that there exists a $\lambda_n$ such that
  $f(x,\diffd)p_n = \lambda_n p_n$. See Lemma~\ref{lem:eigenpoly}.
  \begin{enumerate}
  \item When do they exist?
  \item When is $p_n\in\allpoly$?
  \item When does $p_{n+1}\lesslesseq p_n$?
  \end{enumerate}
\index{eigenpolynomial}
\end{question}

\begin{question}
  For which diagonal matrices $D$ and symmetric matrices $C$ is the
  integral below in $\gsubpos_2$?
$$ \int_0^1 \vert xI+yD+e^{tC}\vert\,dt$$
\end{question}

\begin{question}
  Let $C$ be a function from the reals to symmetric matrices. Define
  $$
  f(x,t) = \vert xI + C(t)\vert.$$
  For every $t$ we kare given that
  $f(x,t)\in\allpoly$. What assumptions on $C$ must we make so that
  $f\in\gsubf_2$? Consider $C(t)=e^{tE}$, where $E$ is symmetric.  The
  roots of $f(x,t)$ are the negative eigenvalues of $e^{tE}$, and they
  are $\{e^{t\lambda_i}\}$ where the eigenvalues of $E$ are
  $\{\lambda_i\}$. Is this in $\gsubf_2$?
\end{question}

\begin{question}\label{ques:p2-imag}
  If $f(x,y)\in\rupint{2}$, and $\sigma$ has positive real and imaginary
  parts then we know that the roots of $f(x,\sigma)$ all lie in the
  lower half plane. Show that the argument of a root $\rho$ of
  $f(x,\sigma)$ satisfies
    $$ |\pi - \arg(\rho)| \le 2 \arg(\sigma).$$
  \end{question}

  \begin{question}
    If $f\in\gsubplus_2$, $f = \sum f_i(x)y^i$, $f_0 = \sum_k
    f_{2k}(x)$, $f_1(x) = \sum_k f_{2k+1}(x)$ then all the roots of
    $f_0$ and of $f_1$ have negative real part (they may be complex). 
  \end{question}

  \begin{question}
Can we remove the factorials from Corollary~\ref{cor:tp-2}? 
  Suppose $f\in\allpolypos$, $d$ is a positive integer, and
  $\alpha>0$. Then show that 
$$
\begin{vmatrix}
   f(\alpha+d) & f(\alpha+d-1) & \hdots &
  f(\alpha)\\
\vdots & & & \vdots \\
\vdots & & & \vdots \\
   f(\alpha+2d) & f(\alpha+2d-1) & \hdots &
  f(\alpha+d)\\
\end{vmatrix}\ge0
$$

  \end{question}

\begin{question}
  If $d$ is odd show that the derivative of $N(f)$ (see
  \eqref{eqn:tp-newt}) is positive for all $x$
\end{question}
\begin{question}
  If $d$ is even show that   $N(f)$ (see
  \eqref{eqn:tp-newt}) is positive for all $x$. (Not just for positive
  $x$.)
\end{question}

\begin{question}\label{ques:extension}
  Is the set of all extensions of two polynomials in $\allpolypos$ a
  convex set? More precisely, is it the convex hull of the $n!$
  polynomials arising from products? \mypage{sec:extension}
\end{question}

\begin{question}
  \emph{(The analog of common interlacing for $\rupint{2}$})
  If $f,g\in\rupint{2}$ and $\alpha f+\beta g\in\gsub_2$ for all
  non-negative $\alpha,\beta$ then is there an $h\in\rupint{2}$ such that
  $f\lesslesseq h$ and $g\lesslesseq h$?

  Using the one variable result, it is easy to see that there is a
  continuous function $h$ such that the graph of $h$ interlaces the
  graphs of $f$ and $g$. Can we choose such an $h$ to be in
  $\rupint{2}$?
\end{question}

\begin{question}
  If $f = \sum a_{i,j}x^iy^j\in\gsubplus_2$ then 
  the $n-2$ by $n-2$ matrix
\[
\begin{pmatrix}
  \frac{a_{i,j}\,a_{i+1,j+1}}{a_{i+1,j}\,a_{i,j+1}}
\end{pmatrix}_{0\le i,j\le n-2}
\]
has $\lfloor\frac{n-1}{2}\rfloor$  negative eigenvalues
and  $\lfloor\frac{n}{2}\rfloor$  positive eigenvalues.
\end{question}

\begin{question}
\added{6/3/7}
  (A generalization of Kurtz's theorem\mypage{thm:kurtz}) Suppose
  that $f = \sum_{i+j\le n}a_{i,j}x^iy^j$ has all positive
  coefficients. If for all relevant $i,j$ the following inequalities 
  hold then $f\in\rupint{2}$. 

\[
\frac{a_{i,j}a_{i,j+1}}{a_{i-1,j}a_{i+1,j}}\ge2,\quad
\frac{a_{i,j+1}a_{i+1,j}}{a_{i,j}a_{i+1,j+1}}\ge2,\quad
\frac{a_{i,j}a_{i+1,j+1}}{a_{i,j+1}a_{i+1,j+1}}\ge2
\]

In other words, if all rhombus inequalities hold with constant $2$
then $f\in\rupint{2}$.  It is a consequence of these inequalities that
$a_{i,j}^2\ge4 a_{i,j-1}a_{i,j+1}$ and so by Kurtz's theorem we know
that all $f_i\in\allpolypos$ where $f = \sum f_i(x)y_i$.
\end{question}

\subsection*{Polynomials in several variables }

\begin{question}
  Given $f\lessless g$ in $\allpolypos$ is there $\sum
  h_{ij}(x)y^iz^j\in\gsubplus_3$ such that
\[
h_{00} = f \qquad
h_{10} = g \qquad
\begin{vmatrix}
  h_{00} & h_{10}\\h_{01} & h_{11}
\end{vmatrix}=1
\]
\end{question}

\begin{question}
  For a fixed polynomial $f$ in $\gsubpos_d$ what are the possible
  dimensions of the cone of interlacing of $f$?
  \end{question}

\begin{question}
Can every polynomial is $\gsubpos_d$ be realized as a coefficient of a
product of linear terms with  a large number of variables?
\end{question}


\begin{question} \label{ques:hadamard}
  If $f=\sum a_\sdiffi(\xx) \yy^\diffi$, $g=\sum b_\sdiffi(\xx)
  \yy^\diffi$, $f\in\gsubplus_{d+e}(n)$ and $g\in\gsubplus_{d+e}(n)$,
  then when is the ``Hadamard Product'' below  in
  $\gsubpos_{d+e}(n)$?
$$ \sum a_\sdiffi(\xx) b_\sdiffi(\xx) \yy^\diffi $$

It appears to hold for $d=0,e=2$, or $d=e=1$. If it held
for $d=e=1$ then we could easily derive a consequence due to
\cite{garloff-wagner}: If $f_0\lessless f_1$ and $g_0\lessless g_1$
are all in $\allpolypos$ then $f_0\ast g_0 \lessless f_1\ast g_1$.
The proof goes like this: We can find $f,g\in\gsubplus_2$ such that
 \begin{align*}
   f (x,y) &= f_0(x) + f_1(x)y + \cdots \\
   g (x,y) &= g_0(x) + g_1(x)y + \cdots \\
\intertext{By Question~\ref{ques:hadamard} }
& f_0\ast g_0 + f_1\ast g_1 y + \cdots  \, \in\gsubpos_2 
 \end{align*}
Since the coefficients interlace, $f_0\ast g_o \lesslesseq f_1\ast g_1$.

\end{question}

\begin{question}\label{ques:legendre}
  The ``generalized Legendre'' polynomials are of the form
$$ \left(\frac{\partial}{\partial x_1}\right)^{n_1}\dots
\left(\frac{\partial}{\partial x_r}\right)^{n_r} 
(x_1^2 +x_2^2 + \cdots+ x_r^2 -1)^{n_1+\cdots+n_r}
$$
They do not satisfy substitution for $x_i's$ in the unit ball. 
\begin{enumerate}
\item Show that the one variable polynomials that are coefficients of
  monomials of the form $x_2^{e_2}x_3^{e_3}\cdots x_r^{e_r}$ are in
  $\allpoly$, where the $e_i$ are positive integers.
\item If the case $r=2$ show that if $$f_{2i,2j}(x^2,y^2) =
\left(\frac{\partial}{\partial x}\right)^{2i}
\left(\frac{\partial}{\partial y}\right)^{2j}
(x^2+y^2-1)^{2i+2j}$$ then $f_{2i,2j}(x,y)\in\gsubpos_2$. 
\end{enumerate}

\end{question}



\begin{question} \label{ques:aei}
  If $e>1$ and $\sum a_\sdiffi \xx^\sdiffi\in\gsubplus_d$, then
$\sum a_\sdiffi^e \xx^\sdiffi\in\gsubplus_d.$ 
\end{question}

\begin{question}
  Suppose that $Q_1,\dots,Q_d$ are negative subdefinite matrices with
  the property that for all positive $\alpha_1,\dots,\alpha_d$ the
  matrix $\alpha_1 Q_1 + \cdots + \alpha_d Q_d$ is negative
  subdefinite. (Find a general construction for such matrices.) Show
  that if $f(\xx)\in\gsubplus_d(n)$ then
$$(-1)^n f(-\xx Q_1\xx^t,\dots,-\xx Q_d \xx^t) \in\gsubpos_d(2n)$$
\end{question}

\begin{question}
  Suppose that $T\colon{}\rupint{d}\longrightarrow\allpoly$ is a linear
  transformation with the property that monomials map to monomials.
  Can we characterize such $T$? Are they constructed out of maps from
  $\rupint{2}\longrightarrow\allpoly$? And, are these of the form
  $x^iy^j\mapsto z^{i+j}$, $x^iy^j\mapsto \alpha^i z^j$,  coefficient
  extraction, and $x^iy^j\mapsto \falling{j}{i} z^{j-i}$?   For instance,
  if it's true that $x^iy^j\mapsto \falling{j}{i} z^{j-i}$ maps
  $\rupint{2}\longrightarrow\allpoly$, then we have a map

\centerline{
\xymatrix{
{\rupint{3}} \ar@{->}[rrr]^{(x,y,v)\mapsto(x,x,v)} &&& {\gsub_2}
\ar@{->}[rrr]^{f(x,\diffd)} &&& \allpoly \\
x^i y^j v^k \ar@{|->}[rrr] &&& x^{i+j}v^k \ar@{|->}[rrr]  &&&  \falling{i+j}{k}
z^{k-i-j}
}}

More generally, can we describe all \emph{multiplier transformations}?
These are linear transformations $\rupint{d}\longrightarrow\gsub_e$ that
map monomials to monomials.
\index{multiplier~transformations}
\end{question}

\begin{question}\label{ques:subdef}
Prove the following. Note that the only difficulty is to replace the
zero entries of $Q$ with small positive values so that there is still
exactly one positive eigenvalue. There is only a problem if the
determinant is zero.

\begin{lemma}
  If $Q$ is a matrix with
  \begin{enumerate}
  \item all non-negative entries
  \item exactly one positive eigenvalue
  \end{enumerate}
then $Q$ is a limit of negative subdefinite matrices.
\end{lemma}
\end{question}

\begin{question}
  Show that the Hurwitz transformation \eqref{eqn:hurwitz-in-pd} maps
  $\gsubplus_d$ to itself.
\end{question}

\begin{question} \label{ques:exp-2d}
  Show that the linear transformation $x_1^{i_1}\cdots x_d^{i_d}
  \mapsto \dfrac{x_1^{i_1}\cdots x_d^{i_d}}{(i_1+\cdots i_d)!}$ maps
  $\gsubplus_d$ to itself. (It does not map $\gsubpos_2$ to itself.)
\end{question}

\begin{question}
 \label{ques:nrec-2}
  Suppose that $f_{0,0}=1$ and $f_{i,j}$ satisfies \eqref{eqn:nrec-2}
  where the constants satisfy $b_{i,j}\ge0$ and
$$ a_{0,0} \le a_{0,j} \le a_{1,0} \le a_{1,j} \le a_{2,0} \le a_{2,j}
\le a_{3,0} \cdots \text{ for } j=1,2,\dots$$
  Then we have $f_{i,j} \greateqeq f_{i,0}$ for all $i,j$. Consecutive
  $f_{i,j}$ do not necessarily interlace.

\end{question}

\begin{question}
  For which regions $D\subset\reals^2$ is it the case that for all
  $f\in\gsubsep_3$ the integral below is in $\allpoly$?
$$ \int_D f(x,s,t)\,ds\,dt$$

\end{question}

\begin{question}
  For which regions $D\subset\reals^2$ is it the case that for all
  $f,g\in\gsubsep_2$ the integral below is in $\allpoly$?
$$ \int_D f(x,s)g(x,t)\,ds\,dt$$
In this case, the conclusion holds if $D$ is the product
$(0,1)\times(0,1)$ or if $D$ is the segment $0\le s,t$ and $s+t=1$ or
if $D$ is the triangle $0\le s,t$ and $s+t\le1$.
\end{question}

\begin{question}\label{ques:p3-ext}
  Suppose that $f_{00}\lesslesseq f_{01},f_{10}$. Show that
  $f_{10}f_{01}-f_{00}W\ge 0$ if and only if there is an extension of
  $f_{00},f_{10},f_{01}$ such that $W=f_{11}$. \seepage{lem:p3-ext-1}
\end{question}


\begin{question}
  If $f(x,y)$ has the property that $f(\alpha,y)\greateq
  g(\alpha,y+1)$ and $f(x,\alpha)\greateq f(x+1,\alpha)$ for all real
  $\alpha$ does it follow that $f\in\gsubsep$?
\end{question}

\subsection*{The polynomial closure of $\gsubpos_d$}








\begin{question}
  We know that $a+bx+cy+dxy$ is in $\gsubclose_2$ for any choice of
  $a,b,c,d$ satisfying $\smalltwodet{a}{b}{c}{d}<0$.  For which
  choices of constants $a_\sdiffi$ do we have
  $$
  \sum_\sdiffi a_\sdiffi \xx^\sdiffi \in \gsubclose_d$$
  where the
  sum is over all $2^d$ monomials $\diffi$?  For instance, if $d=3$
   a determinant such as
$$ 
\begin{vmatrix}
  x+a & d_1 & d_2 \\ d_1 & y+b & d_3 \\ d_2 & d_3 & z+c
\end{vmatrix}
$$
has the desired form. Br\"and\'en's criterion generally leads to
infeasible problems
\end{question}

\begin{question}
  Suppose that $L$ is a lattice in $\mathbb{Z}^d$ with basis vectors
  $v=\{v_1,\dots,v_e\}$. If $\diffi=(i_1,\dots,i_e)$ is an index set
  then set $\diffi\cdot v= i_1v_1+\dots+i_ev_e$. When is it true that
\begin{quote}
  If $f(\xx)=\sum a_\diffi \xx^\diffi$ is in $\gsubplus_d$ then
$ \sum   a_{\diffi\cdot v} \xx^\diffi \in\gsubcloseplus_e$
\end{quote}

Notice that by taking $v_1=(2)$ the Hurwitz theorem on even and odd
parts (Theorem~\ref{thm:hurwitz}) is a consequence of this question. 
Empirically, the question holds for $v_1=(3,3)$ or $(4,4)$ or
$(2,4)$, and fails for $v_1=(3,2)$ or $(2,1)$. The question fails
if we replace $\gsubplus_2$ by $\gsubpos_2$.

 The diagonal (Theorem~\ref{thm:diagonal}) of a polynomial in $\gsubpos_2$ is
found by taking $v_1=(1,1)$.  For another example, take $f =
(x+y+z+1)^n$, $v_1=(1,1,0)$ and $v_2=(0,0,1)$. The resulting
polynomial is in $\gsubclose_2$, and the coefficient of $x_2^k$ is
$\sum \binom{k}{i,i}x_1^i$. 
\end{question}

\begin{question} \label{ques:det-ad}
  If $f\in\gsubposclose_d$ has degree $1$ in each variable then is
  there a $d$ by $d$ matrix $A$ and a diagonal matrix of the form
  $(d_ix_i)$ such that $det(A+D)=f$? In particular, can $xyz-x-y-z+1$
  be so represented?
\end{question}

\begin{question}
  Suppose $f=\sum a_{i,j}x^iy^j\in\gsubpos_2$, and define $f_{rs} =
  \sum a_{2i+r,2j+s}x^iy^j$ where $r,s\in\{0,1\}$. Show that
$$ f_{00} + u f_{10} + u f_{01} + uv f_{11}\in\gsubposclose_4.$$
Is this true in the special case that $f$ is a product of linear factors?
\end{question}

\begin{question}
  \index{Hurwitz's theorem! for $\gsubcloseplus_2$!question} Suppose that
  $g\in\gsubcloseplus_2(n)$ and write \\ $ \displaystyle g(x,y) =
  \sum_{i=0}^n g_i(x)y^i$. Define $$
  f(x) = \sum_{i=0}^n g_i(x^{n+1})
  x^i = g(x^{n+1},x)$$
  then the absolute value of the argument of any
  root of $f$ is at least $\frac{\pi}{n+1}$.

  Note that if $n=1$ then this states that $g_0(x^2)+xg_1(x^2)$ has
  negative real part for $g_0\longleftarrow g_1$ - the classic Hurwitz 
  stability result.
\end{question}


\begin{question}
  Suppose that $f_1,\dots,f_n$ are mutually interlacing. Are there
  positive $a_i$ so that
$$ \sum_{i=1}^n a_i f_i(x)\,y^i \in\gsubposclose_2?$$

A special case of this: suppose the roots of $f$ are $r_1\le\cdots\le
r_n$. Are there positive $a_i$ so that $$\sum_{i=1}^n a_i \frac{f(x)}{x-r_i}\,
y^i\in\gsubposclose_2?$$ 
\end{question}

\begin{question} \label{ques:legendre-2}
  Let $T(x^n) = P_n$, the Legendre polynomial. Show that if
  $f\in\gsubposclose_2$, then $T_\ast(f)(\alpha,y)\in\allpoly$ for
  $|\alpha|<1.$  The case $f=(-1+xy)^n$ is \cite{szego}*{page 387,
    problem 69} and \eqref{eqn:legendre-oval}. 
\end{question}

\begin{question}
  Recall the Hermite polynomials $H_{r,s}$ in
  Theorem~\ref{thm:hermite-pd}. We know that we can represent
  $H_{r,s}$ as the determinant of a matrix.  Show that we have the
  following representation of $H_{r,1}$ as the determinant of the 
  tridiagonal matrix

\begin{gather*}
 \hspace*{-4cm} (-1)^{r+1} e^{ax^2+2bx+cy^2}\frac{\partial^{r+1}}{\partial x^r
    \partial y} \, e^{-(ax^2+2bxy+cy^2)} =\\[.5cm]
 \hspace*{-2cm} \begin{vmatrix}
    2(ax+by) & \sqrt{2a} & 0 & 0 & \hdots & \hdots & 0\\
\sqrt{2a} & 2(ax+by) & \sqrt{4a} & 0 & \hdots & \hdots & 0\\
0 & \sqrt{4a} & 2(ax+by) & \sqrt{6a} & \hdots & \hdots & 0\\
    & & \ddots & \ddots\\
  & & &  \\
0    &\hdots &\hdots &  \sqrt{2(r-2)a} & 2(ax+by) & \sqrt{2(r-1)a} & 0 \\
0    &\hdots &\hdots &0 & \sqrt{2(r-1)a} & 2(ax+by) & {\sqrt{2rb}} \\
0    &\hdots &\hdots &0 & 0 & {\sqrt{2rb}} & 2{(bx+cy)}
  \end{vmatrix}
\end{gather*}

\end{question}

\begin{question}
      For which symmetric matrices $Q$ is $xQy^t\in\gsubclose_{2d}$?
      Setting $\xx=\yy$ shows that $\xx Q\xx^t\in\gsubclose_d$, so $Q$
      is \nsd.
  \end{question}

\begin{question} Suppose that $f,g,h,k\in\allpoly$ satisfy
  \begin{enumerate}
  \item Interlacing: $\smallsquare{f}{g}{h}{k}$
  \item Determinant: $\smalltwodet{f}{g}{h}{k}\le0$
  \end{enumerate}
Show that $f+ y\, g + z \,h + yz\,k \in\gsubclose_3$.
\end{question}



\subsection*{The analytic closure of $\gsubpos_d$}

\begin{question}
  If $L_n$ is the Legendre polynomial then \cite{szego}*{page 69}
$$ f(x,y) = \sum_{i=0}^\infty L_n(x)y^n = (1-2xy+y^2)^{-1/2}$$
If $g\in\allpolyint{(-\infty,1)}(n)$ then show that the coefficients of
$y^n,\cdots,y^{2n}$ in $g(y)f(x,y)$ are in $\allpoly$. The
coefficients of $y^2,\dots,y^{n-1}$ are not necessarily in $\allpoly$.
\end{question}

\begin{question}
  When do all the coefficients of $y^i$ in a product of terms of the
  form $(a+y)^{bx+c}(d+y)^{rx+s}$ have all real roots? The coefficient
  of $y^3$ in $(1+y)^x(2+y)^x$ is not in $\allpoly$. It seems that all
  coefficients are in $\allpoly$ for $(1+y)^x(2-y)^x$.
\end{question}

\begin{question}
  Show that the linear transformation $\xx^\sdiffi \mapsto H_\sdiffi$
  where $H_\sdiffi$ is the $d$-dimensional Hermite polynomial,
\index{Hermite polynomials!in $\gsubpos_d$!question} maps $\gsubpos_d$ to itself.
\end{question}

\begin{question}
  Suppose that $\{p_n\}$ is an orthogonal polynomial system. Can we
  find constants $a_i$ so that
$$ \sum_{i=0}^\infty a_i p_i(x)y^i \in\allpolyf_2$$.
\end{question}

\subsection*{Extending $\rupint{d}$}

\begin{question}
Some questions about $\image{2}$.

\begin{enumerate}
\item Is $\image{2}$ closed under multiplication?
  \item If $f\in\allpolypos$ then is $f(-xy)\in\image{2}$?
  \item   Is $e^{-xy}\in\image{2}$? This would follow from
the above question.
\item Does Lemma~\ref{lem:3-term-partial} hold for $\image{2}$.
\end{enumerate}
\end{question}

\begin{question}
  Consider the map $x^n\mapsto L_n(x;y)$. Show that it maps
  $\allpoly\longrightarrow \partialpoly{1,1}$. 
\end{question}

\begin{question}

  Is $y\sqrt{x}$ in the closure of $\partialpoly{1,1}$?
\end{question}

\begin{question}
  Which functions analytic in $Re(z)>0$ are uniform limits (in
  $Re(z)>0$) of polynomials with all real roots?
\end{question}

\begin{question}
  If $f\in\partialpoly{1,1}$ then show that $f(x,y)$ determine an
  operator $f(x,\diffd)$ that maps
  $\allpolypos\longrightarrow\allpoly$

\end{question}

  \begin{question}
    Is $\pospoly_{d,e}$ closed under differentiation?  Are
    coefficients of $\pospoly_{d,e}$ well-behaved??  Are there
    interesting transformations that preserve $\pospoly_{d,e}$?

  \end{question}

\begin{question}
  Suppose $f,g\in\pospoly_{d,2}$ have the property that for all
  $p\longleftarrow q\in\pospoly_{d,2}$ it holds that
  $fp+gq\in\pospoly_{d,2}$. Then $g\longleftarrow f$.
\end{question}
\subsection*{Generating Functions}

\begin{question}  \label{ques:qxn2}
  If \index{falling factorial!question} $\frac{d}{dx} f_i(x) =  \falling{x}{i}$ then the
  generating function satisfies
  \begin{align*}
    \frac{\partial}{\partial x}\left( \sum f_i(x)
      \frac{y^i}{i!}\right) &= \sum  \falling{x}{i}
    \frac{y^i}{i!} \\
    &= (1+y)^x \\
\intertext{Integration with respect to $x$ shows}
\sum f_i(x) \frac{y^i}{i!} &= \frac{(1+y)^x}{\log(1+y)} + c(y) \\
\intertext{If we choose $c(y)=-1/y$ and write}
\sum g_i(x) \frac{y^i}{i!} &= \frac{(1+y)^x}{\log(1+y)} - \frac{1}{y} 
  \end{align*}
  then show that $g_i(x)$ interlaces $ \falling{x}{i}$ and hence is in
  $\allpoly$. This solves Question~\ref{ques:qxn}.

\end{question}


\begin{question}
  Can we use the identity 
$$ \sum (\Delta^if)\, \frac{\falling{y}{i} }{i!} = f(x+y)$$
to get information about linear combinations of $\Delta^i f$?
\end{question}

\begin{question}
  Is the following true?
\begin{quote}
  Suppose that $a_1\dots,a_r$ and $b_1,\dots,b_s$ are positive and
  satisfy $(\sum b_i) - (\sum a_i)> -1$. Then the Hadamard-type
  product
$$ x^m \ast^{\prime\prime} x^n \mapsto
\begin{cases}
  \frac{\displaystyle\prod_{i=1}^r (a_in)!}{\displaystyle\prod_{j=1}^s
    (b_jn)!}x^n & m=n \\ 0 &\text{otherwise}
\end{cases}
$$
maps $\allpoly\times\allpolypos\longrightarrow\allpoly$. The
assumptions on the $a$'s and $b$'s guarantee that the genus is
$0$. The problem is to show that all the roots are real.
\end{quote}

\end{question}

\begin{question}
  In Lemma~\ref{lem:eigenpoly} we saw that $f(x,\partial)$ had an
  eigenpolynomial $p_n$ of degree $n$ for $n$ sufficiently large. Show
  that there is an eigenpolynomial for all $n$, and that $p_n\lessless
  p_{n-1}$. 
\end{question}

\begin{question}\label{ques:interlace-exy}
  Show that if $G\in\gsubf_2$ interlaces $e^{-xy}$ then $G =
  (ax+by+c)e^{-xy}$ where $a$ and $b$ have the same sign. 
  
  More generally, if $f\in\gsubclose_2$, determine all $G$ such that $G\in\gsubf_2$
  interlaces $f(x,y)e^{-xy}$.

\end{question}

\subsection*{Recurrences  with polynomial coefficients}

\begin{question}
  What are necessary and sufficient conditions for
  $g_0(x),\dots.g_d(x)$ so that there exists a sequence
  $f_0,f_1,\dots$ satisfying
  \begin{enumerate}
  \item $f_i\in\allpoly$, $i=0,1,\dots$
  \item $\sum_{i=0}^d g_i f_{i+n}=0$ for $n=0,1,\dots$.
  \end{enumerate}

\end{question}

\subsection*{Matrices}

\begin{question}
  If $f = \prod(x+b_iy+c_i) = \sum p_i(x)y^i$ then can we explicitly
  construct a matrix $A$ such that $|A[1,\dots,i]|$ is a multiple of
$p_i$? 
\end{question}




\begin{question}
  If $A$ is a $d$ by $d$ matrix of non-negative constants then the
  sequence
  \begin{align*}
    h_0 &= I \\
    h_{n+1} & = x A h_n - h_n^\prime
  \end{align*}
defines a sequence of $d$ by $d$ matrices whose entries are all in $\allpoly$.
\end{question}

\begin{question}
  Suppose that $A,A^\prime$ are diagonal matrices with positive
  diagonal entries, and $C,C^\prime$ are symmetric matrices. If
  $|xI+yA+C| = |xI+yA^\prime+C^\prime|$ then what is the relationship
  between $A,A^\prime$, and $C,C^\prime$?
\end{question}

\begin{question}
  Interlacing was defined for matrices in terms of their
  characteristic polynomials. Is there a notion of ``interlacing'' for
  general linear transformations?  For example, if $S$ and $T$ are
  compact self-adjoint linear operators on a Hilbert space, then it is
  clear what it means to say that the discrete spectrums of $S$ and
  $T$ interlace. However, what does it mean to say that the continuous
  spectrums interlace? 
\end{question}

\begin{question}
Suppose $M=xI+yD+C$ where $C$ is symmetric, and $D=(d_{ij})$
is a diagonal matrix with positive diagonal. If $f(x,y) = |M|$, and
$g(x,y) = |M[1]|$ then define $h = f - (x+d_{11}y)g$. Show that
$h\in\gsubpos_2$ and $ f \lesslesseq h$. 
\end{question}

\begin{question}
  Suppose that $P_0=0$, $P_1$ is the $n$-dimensional identity matrix,
  and $P_{n+1}=x\,P_n - JP_{n-1}$ where $J$ is the $n\times n$ matrix
  of all $1$'s. The diagonals of $P_n$ are all equal to a polynomial
  $f_n$, and the off diagonals are all equal to $g_n$. It isn't true
  that $f_n\in\allpoly$, but show that $g_n\in\allpoly$ and
  $g_{n+1}\lesslesseq g_n$.
\end{question}

\subsection*{Matrix polynomials}
\label{sec:matrix-polynomials}

\begin{question}
  Do all matrix polynomials in $\hyper{1}$ have a common interlacing?
\end{question}

\begin{question}
  Is $\hyper{1}$ closed under Hadamard products?
\end{question}


\begin{question}
If $f\in\hyper{1}$ has degree $n$ then does the function 
\[ x\mapsto \text{maximum eigenvalue of $f(x)$} \] 
have exactly $n$ zeros?
\end{question}

\begin{question}
    Suppose that $D_1,\dots$ and $E_1,\dots$ are positive definite
    matrices. Show that the determinant of 
\[
(D_1+x\,E_1)\ast(D_2+x\,E_2)\ast\cdots\ast(D_n+x\,E_n)
\]
is a stable polynomial, where $\ast$ is the Hadamard product.

More generally, if $A,B$ are matrices with polynomial entries, and
$|A|$ and $B|$ are stable polynomials then is it true that $|A\ast B|$
is stable?

Note that all the polynomials in question are matrix polynomials in
$\hyper{1}$. 
  \end{question}

\subsection*{Polynomials non-vanishing in a region}

\begin{question}
  Determine all the exponential functions in $\nvf{1}{\calC}$ where
  $\calC$ is a cone in $\complexes$. 
\end{question}

\subsection*{Upper half plane polynomials}
\label{sec:ques:cpx-poly}

\begin{question}
  If $A$ is skew-symmetric and $f(\xx)\in\rup{n}$ then is
  $f(A\,\xx)\in\rup{n}$?  If $n=2$ and $A=\smalltwo{0}{1}{-1}{0}$ then
  $f(A\xx) = f(y,-x)$ which is in $\rup{2}$ if $f$ is.
\end{question}

\begin{question}
  Describe all the polynomials that interlace a given polynomial in $\up{1}$.
\end{question}

\begin{question}
  If $f\lessless g$ in $\up{1}$, and all coefficients are positive
  then show that 
$ f\cdot \overline{g} =h(x) +\imag k(x)$ where $h$ and $k$ are stable.
\end{question}

\begin{question}\label{ques:i-interlace}
  If  $(x-\sigma)(x-\tau)\in\up{1}$ then show that 
  $(x-\sigma)(x-\tau) \lessless x-\nu$ if and only if $\nu$ lies in
  the ellipse whose foci are $\sigma,\tau$, and is tangent to the $x$ axis.
\end{question}



\begin{question}

 If $f\in\allpolypos$ then show that $\int_0^1 f(-\imag t x)\,dt
\in\up{1}$. 

A consequence: 

Let $h(x) = \int_0^1 f(-\imag t x)\,dt$. Then
$$\int_{-1}^1 f(\imag -t x)\,dt = \int_0^1 f(-\imag t x)\,dt + \int
_{-1}^0 f(-\imag t x)\,dt = h(x) + \overline{h(x)} = 2 \Re(h)$$
Since $h\in\up{1}$, it follows that the integral from $-1$ to $1$
is in $\allpoly$. 
\end{question}

\begin{question}
    If $f\in\gsubplus_2$, $f = \sum f_i(x)y^i$ then  all the polynomials
$
\begin{vmatrix}
  f_i(x) & f_{i+1}(x)\\ f_{i-1}(x) & f_i(x)
\end{vmatrix}
$ are stable. In particular, all their coefficients are positive.
Thus, the sequence of coefficients of $f$ are \emph{q-log-concave}.
(This is not true if $f\in\rupint{2}$.)
More generally, show that   the $d+1$ by $d+1$ determinant of the matrix
given below is stable.
\[
\begin{vmatrix}
  f_i & f_{i+1} & \dots & f_{i+d}\\
  f_{i+1} & f_{i+2} & \dots & f_{i+d+1} \\
\vdots && \vdots \\
  f_{i+d} & f_{i+d+1} & \dots & f_{2d}
\end{vmatrix}
\]

  \end{question}

\begin{question}
  Although we can't graph a polynomial in $\up{2}$ since it has
  complex roots, we can look at the plots of the real and imaginary
  parts. Since all imaginary parts are negative, the graph of the
  imaginary parts lies in the lower half plane; the graph of the real
  parts has no such restriction. The graph of the real part looks like
  the graph of a polynomial in $\rupint{2}$, while the graph of the
  imaginary part looks very different
  (Figure~\ref{fig:graph-i2}). Explain this behavior.

  \begin{figure}
    \centering
  \newpsobject{showgrid}{psgrid}{subgriddiv=1,griddots=10,gridlabels=6pt}
    \begin{pspicture}(-5,-3)(8,0)
      \psset{xunit=.7cm,yunit=3cm}
      \pscurve(-5,-0.787764)(-4,-0.802378)(-3,-0.822464)(-2,-0.852762)
        (-1,-0.903495)(0,-0.980205)(1,-0.894454)(2,-0.629592)
        (3,-0.515771)(4,-0.46471)(5,-0.436539)(6,-0.41882)
        (7,-0.406682)(8,-0.397858)
      \pscurve(-5,-0.567099)(-4,-0.566402)(-3,-0.567142)(-2,-0.570297)
        (-1,-0.527562)(0,-0.15561)(1,-0.210714)(2,-0.468607)
        (3,-0.580056)(4,-0.630273)(5,-0.658228)(6,-0.676008)
        (7,-0.688322)(8,-0.697355)
        \pscurve(-5,-0.274119)(-4,-0.260201)(-3,-0.239376)(-2,-0.205922)
        (-1,-0.197925)(0,-0.493167)(1,-0.523814)(2,-0.530783)
        (3,-0.533156)(4,-0.533999)(5,-0.534215)(6,-0.534153)
        (7,-0.533977)(8,-0.533769)
    \end{pspicture}
    \caption{The imaginary part of a polynomial in $\polycpx_2$}
    \label{fig:graph-i2}
  \end{figure}
\end{question}

\subsection*{Stable polynomials}
\label{sec:stable-polynomials}

\begin{question}\label{ques:sub-det}
    Suppose that $f = \sum f_i(x)y^i$ is in $\gsubplus_2(n)$. Define the
    polynomial
\[
p_d(x,y) = \sum_{k=0}^{n-2d} y^k\ \det\left( f_{k+i+j}\right)_{0\le i,j< d}\ 
\]
Show that
\begin{enumerate}
\item $p_d(x,\alpha)$  is a stable polynomial  for all positive $\alpha$.
\item All coefficients (the determinants) are stable polynomials.
\end{enumerate}
  \end{question}

\begin{question}
    If $f(x,y)\in\gsubplus_2$ then all minors of the matrix
  \eqref{eqn:tp-stab-2} are stable.

\end{question}

\begin{question}
    If $f,g\in\gsubplus_2$ then all entries of the Bezout matrix
    $B(f,g)$ are  stable polynomials in $y$. 
\[
B(f,g) = \frac{f(v,y)g(u,y) - f(u,y)g(v,y)}{u-v}
\]
  \end{question}

\begin{question}
  If $f = \sum f_i(x,y)z^i\in\gsubplus_3$ and
\[
M_d = 
\begin{pmatrix}
  f_0(x,y_0) & \hdots & f_0(x,y_d) \\
  \vdots & & \vdots \\
  f_d(x,y_0) & \hdots & f_d(x,y_d)
\end{pmatrix}
\]
then the coefficient of any monomial $y_0^{e_0}\cdots y_d^{e_d}$ in
$|M_d|$ is stable. If $d=1$ this is the Bezout matrix.
In addition, the determinant of $M_d$ is stable.
\index{Bezout matrix}

\end{question}

\begin{question}
  A special case of the previous question.  If $\sum
  f_i(x)y^i\in\gsubplus_2$ then the polynomial below has all positive
  coefficients. Moreover, all coefficients of $x^iy^j$ are in
  $\allpolypos$. 
\[
\frac{1}{(x-y)(x-z)(y-z)}
\begin{vmatrix}
  f_0(x) & f_1(x) & f_2(x) \\
  f_0(y) & f_1(y) & f_2(y) \\
  f_0(z) & f_1(z) & f_2(z) \\
\end{vmatrix}
\]
\end{question}

\begin{question}
    If $g_1,\dots,g_n\in\allpolypos(r)$ and $m_1,\cdots, m_n$ are positive
    definite matrices then
\[
\bigl|  g_1m_1 + \cdots + g_nm_n\bigr| \in\stabled{}
\]
  \end{question}

  \begin{question}
    Is the determinant of a stable matrix polynomial a stable
    polynomial? 
  \end{question}



\begin{question}
If $f = \sum f_i(x)y^i\in\gsubplus_2$ then is 
\[
\frac{f_1}{f_2}-\frac{f_0}{f_1}\colon \quada \longrightarrow \quada
\]
\end{question}

\begin{question}
  Suppose $f = \prod(x+r_i)\in\allpolypos$. What polynomials in
  $\allpolypos$ interlace $f$ in $\stabled{1}$?
\end{question}



\begin{question}\label{ques:stable-int}
  If $f\hlace g$  then when are there positive definite
  $D_1,D_2$,  skew-symmetric $A$, and positive $\alpha,\beta$  such that 
\[ 
\bigl| I + x D_1 + y D_2 +A\bigr| = \alpha\, f + \beta\,g\,y + \cdots
\]
\end{question}

\begin{question}
    Suppose that $f = \sum f_i(x)y^i$ is in $\gsubplus_2(n)$. Define the
    polynomial
\[
p_d(x,y) = \sum_{k=0}^{n-2d} y^k\ \det\left( f_{k+i+j}\right)_{0\le i,j< d}\ 
\]
Show that
\begin{enumerate}
\item $p_d(x,\alpha)$  is a stable polynomial  for all positive $\alpha$.
\item All coefficients (the determinants) are stable polynomials.
\end{enumerate}
  \end{question}



\begin{question}
  Define the infinite matrix
\[
M = 
\begin{pmatrix}
. & . & . & \hdots \\
  \binom{k}{k} & . & . & \hdots \\
  . &\binom{k+1}{k} &  . & \hdots \\
  . & . & \binom{k+2}{k}  & \hdots \\
  \vdots & & & \ddots
\end{pmatrix}
\]
where the periods (.) mean zero. If 
\[
\exp(M) =
\begin{pmatrix}
  . & . & . & \hdots \\
  a_{0,0} & . & . & \hdots \\
  a_{0,1} & a_{1,1} & . & \hdots \\
  a_{0,2} & a_{1,2} & a_{2.2} & \hdots\\
  \vdots & & & \ddots
\end{pmatrix}
\]
then show that the rows $ \sum a_{i,j}x^j$ are stable.

If $k=2$ and $f_n$ is the polynomial corresponding to the $n$th row then the recurrence 
\[ f_n = (x+n-1)f_{n-1} - \frac{1}{2}\binom{n}{2}f_{n-2}\]
shows that $f_n\in\allpoly$. Since all entries of the exponential are
positive we have that $f_n\in\allpolypos\subset\stabled{1}$.  
\end{question}


\begin{question}
  If $f\in\stabled{}$ and $f(x)\ge0$ for all real vectors $x$ then are
  there $g\place h$ so that $f = g^2+h^2$.
\end{question}



\subsection*{Transformations in the complex plane}

\begin{question}
  Show that the containment in Proposition~\ref{prop:laguerre-image}
  is an equality. That is, if $T(x^n)=L_n(x)$ and $|z|\le1$ then all
  roots of $T(x-z)^n$ have non-negative real part.
\end{question}

\begin{question}
  If $T\colon{}x^n\mapsto L_n(x)$ then show that
  $\Lambda(T\allpolyint{\Delta})$ is the closed right half
  plane. Equivalently, show that if $Tf$ has a root with negative real
  part, then $f$ has a root with absolute value greater than $1$.
\end{question}

\begin{question}
  Suppose $f_\beta$ is defined as in \eqref{eqn:fbeta}. Is there a
  Sendov type result for the roots of $f_\beta$? That is, if all the
  roots of $f$ are in the unit disk, then how far away from a root of
  $f$ are the roots of $f_\beta$?
\end{question}

\begin{question}
  What linear transformations have the property that they map
  polynomials with all roots on a horizontal line in the complex plane
  to polynomials with all roots on a horizontal line? The
  transformation $x^n\mapsto H_n$ is such a transformation.
\end{question}

\begin{question}\label{ques:cpx-equal}
  If $\cals$ is a region in $\complexes$ and $T$ is a linear
  transformation, then when is it the case that
  \begin{enumerate}
  \item $\Lambda(T\allpolyint{\cals}) = \Lambda(T(\phi\cals))$?
  \item $\Lambda(T\allpolyint{\cals}) \subset Hull(\,\Lambda(T(\phi\cals))\,)$?
  \end{enumerate}
\end{question}

\begin{question}
  If $\Delta$ is the unit disk in the complex plane, then what are
  examples, other than dilations or powers of the derivative, of linear
  transformations that map $\allpolyint{\Delta}$ to itself?
\end{question}

\begin{question}
  If $T(x^n) = p_n(x)$ where $p_n$ is the Legendre or Chebyshev
  polynomial, then is it the case that $$
  T(\allpolyint{\Delta})\subseteq \allpolyint{\Delta}?$$
\end{question}

\begin{question}
  What are examples, other than the derivative, of linear
  transformations that map $\allpolyint{\Delta}$ to itself?
\end{question}

\begin{question}
  Let $T(x^n) = L_n(x)$. 
    If $\Delta$ is the unit disk and $\diffr$ the right half plane,
    then 
$$    T(\allpolyint{\Delta})\subseteq \allpolyint{\diffr}$$

That is, if all roots of $f$ have magnitude at most $1$, then the roots of
$Tf$ have non-negative real part. Is it also the case that every point
in the right half plane is the root of  $Tf$, for some
$f\in\allpolyint{\Delta}$? 
\end{question}


\chapter{Notes}

\renewcommand{\TimeStampStart}{Tuesday, January 05, 1999: 06:54:04}
\mytoday

\subsection*{Polynomials}

Interlacing has several synonyms in the literature.  Gantmacher
\cite{gantm} uses the term \emph{positive pair} for interlacing
polynomials with all negative roots. In \cite{fuhrmann} a pair of
polynomials $f,g$ with all real roots satisfying $f \lessless g$ or $f
\lesseq g$ is called a {\em real pair.}\/ If they satisfy $f \lesseq g $
and all roots are negative, they are called a {\em positive pair}.
Wagner \cite{wagner} says that {\em $f$ alternates left of $g$}\/ if $f
\lesseq g$.  Bilodeau \cite{bilodeau} says that $q$ {\em separates}
$p$ if $p \lesslesseq q$. In \cite{gregor96} two polynomials with a
common interlacing are called \emph{complementary}.  In \cite{simion}
and \cite{kim-kwon} a sequence $f \lessless g \lessless h$ is called a
{\em Sturm sequence of polynomials}.  \cite{constan} says that two
polynomials interlace if the roots \emph{se s\'{e}perent} (separate themselves).

Apparently Markov was the first to prove that differentiation
preserves interlacing. See \cite{rahman}.

Lemma~\ref{lem:add-small} is in \cite{masson}.

Quantitative sign interlacing occurs in several different places.  It
is implicit in \cite{gantm}*{page 226}.  The case $f \lesseq g$ is found
in  \cite{garloff-wagner}.  Quantitative sign interlacing is used in
\cite{godsil} to prove that the characteristic polynomial of a graph
is interlaced by the characteristic polynomial of a one point
deletion. Exercise 7 on page 146 of \cite{godsil} is an explicit
statement of quantitative sign interlacing for $f \lessless g$.
\index{quantitative sign interlacing}
Quantitative sign interlacing can also be interpreted as a property of
rational functions. See \cite{sheil-small}*{page 320}.

\cite{wagner96} contains many of the sign-interlacing results of
\chapsec{polynomials}{intro-sign-interlacing}

Corollary~\ref{cor:inequality-3} is well known; it occurs for example as an
exercise in \cite{fuhrmann}*{page 269}.

Johnson \cite{johnson} shows the convexity of the set of polynomials
that interlace a given polynomial using sign interlacing.

The basis in Section~\ref{sec:sign-quant} is sometimes called the
``interpolation basis''~\cite{fuhrmann}, since they are, up to a
multiple, the polynomials in Lagrange interpolation.
Fuhrman \cite{fuhrmann} calls the basis $f_i(x) = \frac{f(x)}{x-a_i}$
the \emph{spectral basis}.

See \cite{rahman}*{page 198} for a different proof of
Proposition~\ref{prop:pattern}.  \cite{gregor96} established
Proposition~\ref{prop:pattern} in the case of positive combinations
for the special class of Hurwitz polynomials.  These are polynomials
$f(x^2)$ where $f\in\allpolyalt$ has no repeated roots.

The identity \cite{chihara}*{page 24}
$\smalltwodet{p_n}{p_{n+1}}{p_n^\prime}{p_{n+1}^\prime}>0$ where the
$p_i$ are orthogonal follows from Lemma~\ref{lem:inequality-1} using the fact
that $p_{n+1}\lessless p_n$.

A consequence of Lemma~\ref{lem:square-1} is that the associate polynomials of a
family of orthogonal polynomials interlace.  See 
\cite{chihara}*{page 86}.  Lemma~\ref{lem:square-1} is often used implicitly -
see for example \cite{locher}*{Lemma 3.1}

\cite{barrett} proved Lemma~\ref{lem:inequality-1} in the more general case of
differentiable functions.

Corollary~\ref{cor:log-con-coef} is in \cite{constan}.

 \cite{faddeev} has many results about polynomials, with no
 references. Some of them are from \cite{polya-szego2}.
This is   translation of the book that first appeared in 1954. 
Here are a few of the items in the book:

\begin{tabular}{lp{5in}}
  problem 727 & interlacing implies all linear  combinations in $\allpoly$\\
  problem 728 & all linear combinations implies interlacing \\
  problem 729 & $f,g$ interlace implies derivatives interlace.\\
  problem 730 &$\gamma f + (\lambda+x)f'\in\allpoly$ if $\gamma>0$ or $\gamma<
  -n$\\
problem 742 & if $f-a\in\allpoly$ and $f-b\in\allpoly$, then
$f-\gamma\in\allpoly$ for $a\le \gamma\le b$
\end{tabular}

\cite{nijenhuis} has a proof of an elementary  interlacing property.

Jahnke \cite{jahnke}*{page 153} observes that $af+bf^\prime$ and $cf+df^\prime$
interlace if $ad-bc$ is not zero, where $f$ is a Bessel function. 
\index{Bessel function}

Lemma~\ref{lem:polytope-2} is a version of Kharitonov's theorem
\cite{robust} for polynomials in $\allpolypos$. 

See \cite{horwitz} or \cite{peyser} for inequalities similar to
Example~\ref{ex:roots-ineq}.

Corollary~\ref{cor:orthog-quad} is in \cite{liu-wang}, along with many
nice applications.

\subsection*{Matrices preserving interlacing}
\label{sec:matr-pres-interlace}

In \cite{goodman-sun} they use the term \emph{interlacing  negative
  zeros} for a set of mutually interlacing polynomials (with negative
zeros). 

The proof of Proposition~\ref{prop:hurwitz-totally-pos} is based on
the argument in Holtz\cite{holtz}.

\subsection*{Homogeneous}

We could also derive Newton's inequalities from Corollary~\ref{cor:log-con-coef}
using $f$ and $f^\prime$.

\subsection*{Analytic}
Theorem~\ref{thm:polya-schur-2} is in \cite{bilodeau}. The class $\allpolyf$ is
also known as the Laguerre-\Polya\ class $\mathcal{L-P}$
\cite{cravenII}. The polynomials $T(1+x)^n$ are called  Jensen
polynomials. A polynomial in $\allpolyposf$ is called type I.

In \cite{dimitrov-pena} they essentially prove Kurtz' result
(Theorem~\ref{thm:kurtz}), but with a coefficient of $4.07$ instead
of $4$. 

\subsection*{Linear Transformations}

In \cite{cravenIII} they discuss various questions concerning linear
transformations, in particular what they call linear operators
operating on entire functions. For instance, they ask for a
classification of all linear transformations that take polynomials
with zeros in a set $S$ (which might be complex) to polynomials with all
zero in a set $T$. In \cite{craven85}{Corollary 2.4} they determine
all linear transformations $T$ for which $f$ and $Tf$ interlace, but
just for transformations that are multipliers. The answer is not
expressed in terms of derivatives, but rather that the multiplier
sequence is a nonconstant arithmetic sequence all of whose terms have
the same sign.

\cite{masson} proves a special case of Lemma~\ref{lem:eigen-2}.

\cite{branden} proved the identity of Example~\ref{ex:branden} by
using properties of the diamond product.

\subsection*{Examples of Linear Transformations}

\Szego \cite{szego} is the source for all the facts about orthogonal
 polynomials.

Iserles and Norsett \cite{iserles} call the inverse of the falling factorial
the {Laguerre transformation}. 

In \cite{dubeau98} they study the interlacing properties of the
polynomials sequences determined by recurrences $p_{n+1} = 2r_nxp_n -
(1-x^2)p_n^\prime$ for interlacing sequences $(r_n)$ of
coefficients. 

See \cite{bragg99} for an identity involving Hadamard products and 
Laguerre polynomials that is similar to \eqref{eqn:laguerre}. 

\index{Hermite polynomials!modified}%
\index{umbral notation}%
Some results involving linear transformations of polynomials have been
expressed in umbral notation, which, to my mind, is just confusing.
For example, Rota \cite{rota} defines modified Hermite polynomials
$\tilde{H}_n(x)$ by
$$
\sum_{n=0}^\infty \tilde{H}_n(x) \frac{y^n}{n!} = e^{-x^2+xy}$$
  If we define
  $T(x^n) = \tilde{H}_n(x)$ and $M(x) = \sqrt{2}x$ then we have an
  elegant identity $T^2 = M^{-1}TM$. However, this is written in
  \cite{rota} as 
$$ \tilde{H}_n(\mathbf{H})= 2^{n/2}\tilde{H}_n(\frac{x}{\sqrt{2}})$$

\cite{tanja01} and \cite{shapiro-01} discuss asymptotic properties of
the roots of eigenpolynomials of arbitrary differential operators of
the form \eqref{eqn:diff-op-1}.

\subsection*{Affine Transformations}
Johnson \cite{johnson} proves Lemma~\ref{lem:deriv-diff} by a messy algebraic
analysis. The idea for the $q$-Laguerre polynomials comes from
\cite{shukla}, who considers various properties of more general
polynomials. There is no discussion about the roots, however.

The affine derivative was called $L_{(q,b)}$ in \cite{hahn}.

\subsection*{Complex Transformations}
There is an analytic proof \cite{iserles-saff} that the linear
transformation $x^n\longrightarrow T_n$ maps $\allpolyint{(-1,1)}$ to
itself.

Wagner \cite{wagner-newton} refers to a folklore theorem that is
Corollary~\ref{cor:newton-cpx-3} with $\alpha=\pi/3$.

\subsection*{Polynomials in two variables}

Lemma~\ref{lem:p2-rays} is a standard result about asymptotes. See
\cite{nunemacher} Lemma~\ref{lem:tridiagonal} can be found in \cite{parlett}

\subsection*{Polynomials in several variables}

\subsection*{Matrices}

The  calculation in Section~\ref{sec:matrix-submatrix} is
from \cite{horn-johnson-1}*{page 187}.  

Smith \cites{smith,smith-hu} extends results on interlacing of Schur
complements to singular Hermitian matrices using the Moore-Penrose
inverse. 

\subsection*{Stable}
\label{sec:stable}

Lemma~\ref{lem:stable-reasons} is in \cite{zabczyk}*{page 33}, but the
proof is more involved. \cite{liu-wang-log} showed that the Bell
polynomials are $q$-log-convex.


\cite{*}

\chapter{Bibliography}
\bibliography{biblio}

\begin{bibdiv}
\begin{biblist}

\bib{aigner-ziegler}{book}{
      author={Aigner, Martin},
      author={Ziegler, G\"unter},
       title={Proofs from the book},
   publisher={Springer},
        date={1998},
}

\bib{asw}{article}{
      author={Aissen, Michael},
      author={Schoenberg, I.~J.},
      author={Whitney, A.~M.},
       title={On the generating functions of totally positive sequences. {I}},
        date={1952},
     journal={J. Analyse Math.},
      volume={2},
       pages={93\ndash 103},
      review={\MR{0053174 (14,732d)}},
}

\bib{al-salam}{article}{
      author={Al-Salam, Waleed~A.},
      author={Ismail, Mourad E.~H.},
       title={Polynomials orthogonal with respect to discrete convolution},
        date={1976},
     journal={J. Math. Anal. and Appl.},
      volume={55},
       pages={125\ndash 139},
}

\bib{anderson}{article}{
      author={Anderson, Bruce},
       title={Polynomial root dragging},
        date={1993},
     journal={American Mathematical Monthly},
       pages={864\ndash 866},
}

\bib{ando}{article}{
      author={Ando, T.},
       title={Totally positive matrices},
        date={1987},
     journal={Linear algebra and its applications},
      volume={90},
       pages={165\ndash 219},
}

\bib{askey}{book}{
      author={Andrews, George~E.},
      author={Askey, Richard},
      author={Roy, Ranjan},
       title={Special functions},
      series={Encyclopedia of Mathematics and its Applications},
   publisher={Cambridge University Press},
     address={Cambridge},
        date={1999},
      volume={71},
        ISBN={0-521-62321-9; 0-521-78988-5},
      review={\MR{2000g:33001}},
}

\bib{asner}{article}{
      author={Asner, Bernard~A., Jr.},
       title={On the total nonnegativity of the {H}urwitz matrix},
        date={1970},
        ISSN={1095-712X},
     journal={SIAM J. Appl. Math.},
      volume={18},
       pages={407\ndash 414},
      review={\MR{0260768 (41 \#5391)}},
}

\bib{barrett}{article}{
      author={Barrett, Louis~C.},
      author={Bendixen, Gerald~E.},
       title={Separation and interlacing theorems},
        date={1965},
     journal={Quart. Appl. Math.},
      volume={23},
       pages={69\ndash 78},
      review={\MR{31 \#2523}},
}

\bib{beauzamy}{article}{
      author={Beauzamy, Bernard},
      author={D{\'e}got, J{\'e}r{\^o}me},
       title={Differential identities},
        date={1995},
        ISSN={0002-9947},
     journal={Trans. Amer. Math. Soc.},
      volume={347},
      number={7},
       pages={2607\ndash 2619},
      review={\MR{1277095 (96c:05009)}},
}

\bib{tanja01}{incollection}{
      author={Bergkvist, Tanja},
      author={Rullg\o{a}rd, Hans},
       title={On polynomial eigenfunctions for a class of differential
  operators},
        date={2001},
   booktitle={Research reports in mathematics},
   publisher={Department of Mathematics, Stockholm University},
}

\bib{bernstein}{book}{
      author={Bernstein, Dennis~S.},
       title={Matrix mathematics},
   publisher={Princeton University Press},
     address={Princeton, NJ},
        date={2005},
        ISBN={0-691-11802-7},
        note={Theory, facts, and formulas with application to linear systems
  theory},
      review={\MR{2123424 (2005i:93001)}},
}

\bib{bhatia}{article}{
      author={Bhatia, Rajendra},
       title={Infinitely divisible matrices},
        date={2006},
     journal={Americal mathemtical monthly},
}

\bib{robust}{book}{
      author={Bhattacharyya, S.P.},
      author={Chapellat, H},
      author={Keel, L.H.},
       title={Robust control: The parametric approach},
   publisher={Prentice Hall},
        date={1995},
}

\bib{bilodeau}{article}{
      author={Bilodeau, G.~G.},
       title={Generating polynomials all of whose roots are real},
        date={1991},
     journal={Mathematics Magazine},
      volume={64},
       pages={263\ndash 270},
}

\bib{boas}{book}{
      author={Boas, Ralph~Philip, Jr.},
       title={Entire functions},
   publisher={Academic Press Inc.},
     address={New York},
        date={1954},
      review={\MR{16,914f}},
}

\bib{bbl}{misc}{
      author={Borcea, Julius},
      author={Br\"and\'en, Petter},
      author={Liggett, Thomas~M.},
       title={Negative dependence and the geometry of polynomials},
        note={arXiv:0707.2340},
}

\bib{bbs-johnson}{misc}{
      author={Borcea, Julius},
      author={Br\"and\'en, Petter},
      author={Shapiro, Boris},
       title={Applications of stable polynomials to mixed determinants:
  Johnson's conjectures, unimodality and symmetrized fischer products},
        note={math.SP/0607755},
}

\bib{bbs}{misc}{
      author={Borcea, Julius},
      author={Br\"and\'en, Petter},
      author={Shapiro, Boris},
       title={{Classification of hyperbolicity and stability preservers: the
  multivariate Weyl algebra case}},
        note={arXiv:math.CA/0606360},
}

\bib{bb2}{misc}{
      author={Borcea, Julius},
      author={Br\:and\`en, Petter},
      author={Shapiro, Boris},
       title={Polya-schur master theorems for circular domains and their
  boundaries},
        note={arXiv:math/0607416},
}

\bib{borwein}{book}{
      author={Borwein, Peter},
      author={Erd{\'e}lyi, Tam{\'a}s},
       title={Polynomials and polynomial inequalities},
      series={Graduate Texts in Mathematics},
   publisher={Springer-Verlag},
     address={New York},
        date={1995},
      volume={161},
        ISBN={0-387-94509-1},
      review={\MR{97e:41001}},
}

\bib{bourin}{misc}{
      author={Bourin, Jean-Christophe},
       title={{Matrix versions of some classical inequalities}},
        note={arXiv:math.OA/0601543},
}

\bib{bragg99}{article}{
      author={Bragg, L.~R.},
       title={Trigonometric integrals and {H}adamard products},
        date={1999},
     journal={Amer. Math. Monthly},
      volume={106},
      number={1},
       pages={36\ndash 42},
}

\bib{branden-hpp}{misc}{
      author={Br\"and\'en, P.},
       title={{Polynomials with the half-plane property and matroid theory}},
        note={arXiv:math.CO/0605678},
}

\bib{branden}{misc}{
      author={Br\"and\'en, Petter},
       title={{On linear transformations preserving the P\'olya Frequency
  property}},
        note={arXiv:math.CO/0403364},
}

\bib{brenti}{book}{
      author={Brenti, Francesco},
       title={Unimodal, log-concave and {P}\'olya frequency sequences in
  combinatorics},
      series={Memoirs},
   publisher={American Mathematical Society},
        date={1989},
      volume={81},
      number={413},
}

\bib{bump}{article}{
      author={Bump, Daniel},
      author={Choi, Kwok-Kwong},
      author={Kurlberg, P{\"a}r},
      author={Vaaler, Jeffrey},
       title={A local {R}iemann hypothesis. {I}},
        date={2000},
        ISSN={0025-5874},
     journal={Math. Z.},
      volume={233},
      number={1},
       pages={1\ndash 19},
      review={\MR{2001g:11073a}},
}

\bib{cardon04}{article}{
      author={Cardon, David},
       title={Sums of exponential functions having only real zeros},
        date={2004},
     journal={Manuscripta Math.},
      volume={113},
       pages={307\ndash 317},
}

\bib{cardon-2002}{incollection}{
      author={Cardon, David~A.},
      author={Nielsen, Pace~P.},
       title={Convolution operators and entire functions with simple zeros},
        date={2002},
   booktitle={Number theory for the millennium, i (urbana, il, 2000)},
   publisher={A K Peters},
     address={Natick, MA},
       pages={183\ndash 196},
      review={\MR{1956225 (2003m:30012)}},
}

\bib{cardon}{article}{
      author={Cardon, David~A.},
       title={Convolution operators and zeros of entire functions},
        date={2002},
     journal={Proc. of the Amer. Math. Soc.},
       pages={1725\ndash 1734},
}

\bib{pinkus-pena}{article}{
      author={Carnicer, J.~M.},
      author={Pe{\~n}a, J.~M.},
      author={Pinkus, A.},
       title={On some zero-increasing operators},
        date={2002},
        ISSN={0236-5294},
     journal={Acta Math. Hungar.},
      volume={94},
      number={3},
       pages={173\ndash 190},
      review={\MR{1905723 (2003c:26020)}},
}

\bib{carroll}{article}{
      author={Carroll, F.~W.},
       title={A polynomial in each variable separately is a polynomial},
        date={1961},
     journal={Amer. Math. Monthly},
      volume={68},
       pages={42},
}

\bib{chihara}{book}{
      author={Chihara, T.~S.},
       title={Orthogonal polynomials},
   publisher={Gordon-Breach},
        date={1978},
}

\bib{sokal-wagner}{article}{
      author={Choe, Young-Bin},
      author={Oxley, James~G.},
      author={Sokal, Alan~D.},
      author={Wagner, David~G.},
       title={Homogeneous multivariate polynomials with the half-plane
  property},
        date={2004},
        ISSN={0196-8858},
     journal={Adv. in Appl. Math.},
      volume={32},
      number={1-2},
       pages={88\ndash 187},
        note={Special issue on the Tutte polynomial},
      review={\MR{2037144 (2005d:05043)}},
}

\bib{MR1787335}{book}{
      editor={Colonius, Fritz},
      editor={Helmke, Uwe},
      editor={Prätzel-Wolters, Dieter},
      editor={Wirth, Fabian},
       title={Advances in mathematical systems theory},
      series={Systems \& Control: Foundations \& Applications},
   publisher={Birkh\"auser Boston Inc.},
     address={Boston, MA},
        date={2001},
        ISBN={0-8176-4162-9},
      review={\MR{1787335 (2001e:93003)}},
}

\bib{constan}{article}{
      author={Constantinescu, Florian},
       title={Relations entre les coefficients de deus polyn\^{o}mes dont les
  racines se s\'{e}parent},
        date={1964},
     journal={\v{C}asopis pro p\v{e}stov\'{a}ni matematiky},
       pages={1\ndash 3},
}

\bib{craven85}{article}{
      author={Craven, T.},
      author={Csordas, G.},
       title={On the betweenness condition of rolle's theorem},
        date={1985},
     journal={Rocky Mountain. J. of Math.},
      volume={15},
       pages={721\ndash 728},
}

\bib{cravenIII}{unpublished}{
      author={Craven, Thomas},
      author={Csordas, George},
       title={Composition theorems, multiplier sequences and complex zero
  decreasing sequences},
        note={to appear},
}

\bib{craven82}{article}{
      author={Craven, Thomas},
      author={Csordas, George},
       title={On the number of real roots of polynomials},
        date={1982},
     journal={Pacific Journal of Mathematics},
      volume={102},
      number={1},
       pages={15\ndash 28},
}

\bib{cravenII}{article}{
      author={Craven, Thomas},
      author={Csordas, George},
       title={Complex zero decreasing polynomials},
        date={1995},
     journal={Methods Appl. Anal.},
       pages={420\ndash 441},
}

\bib{csordas-waleffe}{article}{
      author={Csordas, George},
      author={Charalambides, Marios},
      author={Waleffe, Fabian},
       title={A new property of a class of jacobi polynomials},
        date={2005},
     journal={Proc. Amer. Math. Soc.},
      volume={133},
       pages={3551\ndash 3560},
}

\bib{csordas00}{article}{
      author={Csordas, George},
      author={Dimitrov, Dimitar~K.},
       title={Conjectures and theorems in the theory of entire functions},
        date={2000},
     journal={Numerical Algorithms},
       pages={1\ndash 14},
}

\bib{csordas75}{article}{
      author={Csordas, George},
      author={Williamson, Jack},
       title={The zeros of {J}ensen polynomials are simple},
        date={1975},
     journal={Proc. of the Amer. Math. Soc.},
      volume={49},
      number={1},
       pages={263\ndash 264},
}

\bib{dedieu94}{article}{
      author={Dedieu, Jean~Pierre},
      author={Gregorac, R.~J.},
       title={Corrigendum. {O}breshkoff's theorem revisited: {W}hat convex sets
  are contained in the set of hyperbolic polynomials?},
        date={1994},
     journal={Journal of Pure and Applied Algebra},
      volume={93},
       pages={111\ndash 112},
}

\bib{dedieu92}{article}{
      author={Dedieu, Jean~Pierre},
       title={Obreshkoff's theorem revisited: {W}hat convex sets are contained
  in the set of hyperbolic polynomials?},
        date={1992},
     journal={Journal of Pure and Applied Algebra},
      volume={81},
       pages={269\ndash 278},
}

\bib{dimitrov-pena}{misc}{
      author={Dimitrov, Dimitar~K.},
      author={Pena, Juan~Manuel},
       title={Almost strict total positivity, strong {T}uran inequalities and a
  class of {H}urwitz polynomials},
         how={Seminario mathematico, University of Zaragoza},
        date={2004},
}

\bib{dimitrov-ranga}{article}{
      author={Dimitrov, Dimitar~K.},
      author={Ranga, A.~Sri},
       title={Monotonicity of zeros of orthogonal {L}aurent polynomials},
        date={2002},
        ISSN={1073-2772},
     journal={Methods Appl. Anal.},
      volume={9},
      number={1},
       pages={1\ndash 11},
      review={\MR{1948463 (2003i:42034)}},
}

\bib{dimitrov98}{article}{
      author={Dimitrov, Dimitar~K.},
       title={Counterrexamples to a problem of {\Polya}\ and to a problem of
  {K}arlin},
        date={1998},
     journal={East. J. Approx.},
       pages={479\ndash 489},
}

\bib{dodonov}{article}{
      author={Dodonov, V.~V.},
      author={Man'ko, V.~I.},
       title={New relations for two-dimensional {H}ermite polynomials},
        date={1994},
     journal={J. Math. Phys.},
      volume={35},
      number={8},
       pages={4227\ndash 4288},
}

\bib{driver}{article}{
      author={Driver, Kathy},
      author={Duren, Peter},
       title={Trajectories of the zeros of hypergeometric polynomials
  f(-n,b;2b;z) for $b<-1/2$},
        date={2001},
     journal={Const. Approx.},
      volume={17},
       pages={169\ndash 179},
}

\bib{driver01}{article}{
      author={Driver, Kathy},
      author={Duren, Peter},
       title={Zeros of ultraspherical polynomials and the {H}ilbert-{K}lein
  formulas},
        date={2001},
     journal={J. Compu. Appl. Math.},
      volume={135},
       pages={293\ndash 301},
}

\bib{dubeau98}{article}{
      author={Dubeau, F.},
      author={Savoie, J.},
       title={More on the roots of {E}uler-{F}robenius polynomials},
        date={1998},
     journal={Ann. Sci. Math. Qu\'{e}bec},
      volume={22},
      number={1},
       pages={31\ndash 42},
}

\bib{dunkl-hu}{book}{
      author={Dunkl, Charles~F.},
      author={Hu, Yuan},
       title={Orthogonal polynomials of several variables},
      series={Encyclopedia of Mathematics and Its Applications},
   publisher={Cambridge University Press},
        date={2001},
      volume={81},
}

\bib{dunkl}{book}{
      author={Dunkl, Charles~F.},
      author={Xu, Yuan},
       title={Orthogonal polynomials of several variables},
   publisher={Cambridge University Press},
        date={2001},
}

\bib{ryavec-redmond}{article}{
      author={E{\u{g}}ecio{\u{g}}lu, {\"O}mer},
      author={Redmond, Timothy},
      author={Ryavec, Charles},
       title={From a polynomial {R}iemann hypothesis to alternating sign
  matrices},
        date={2001},
        ISSN={1077-8926},
     journal={Electron. J. Combin.},
      volume={8},
      number={1},
       pages={Research Paper 36, 51 pp. (electronic)},
      review={\MR{2002k:05251}},
}

\bib{ryavec}{inproceedings}{
      author={E{\u{g}}ecio{\u{g}}lu, {\"O}mer},
      author={Ryavec, Charles},
       title={Polynomial families satisfying a {R}iemann hypothesis},
        date={2001},
   booktitle={Proceedings of the thirty-second southeastern international
  conference on combinatorics, graph theory and computing (baton rouge, la,
  2001)},
      volume={149},
       pages={177\ndash 191},
      review={\MR{2002m:11019}},
}

\bib{faddeev}{book}{
      author={Faddeev, D.~K.},
      author={Sominskii, I.~S.},
       title={Problems in higher algebra},
      series={Translated by J. L. Brenner},
   publisher={W. H. Freeman and Co.},
     address={San Francisco},
        date={1965},
      review={\MR{31 \#1258}},
}

\bib{fernandez98}{article}{
      author={Fern{\'a}ndez, F.M.},
       title={Generating functions for {H}ermite polynomials of arbitrary
  order},
        date={1998},
     journal={Physics Letters A},
      volume={237},
       pages={189\ndash 191},
}

\bib{aspects}{misc}{
      author={Fisk, Steve},
       title={Aspects of stable polynomials},
        note={arXiv:0803.0286},
}

\bib{upper}{misc}{
      author={Fisk, Steve},
       title={An introduction to upper half plane polynomials},
        note={arXiv:0711.4043},
}

\bib{fisk}{book}{
      author={Fisk, Steve},
       title={Coloring theories},
      series={Contemporary Mathematics},
   publisher={American Mathematical Society},
        date={1989},
      number={103},
}

\bib{fomin}{unpublished}{
      author={Fomin, Sergey},
      author={Zelevinsky, Andrei},
       title={Total positivity: tests and parameterizations},
        note={math.L.A/9912128},
}

\bib{frenkel}{misc}{
      author={Frenkel, Peter~E.},
       title={Pfaffians, hafnians, and products of linear functionals},
        note={arXiv:math.CA/07040028},
}

\bib{fuhrmann}{book}{
      author={Fuhrmann, Paul~A.},
       title={A polynomial approach to linear algebra},
   publisher={Springer},
        date={1996},
}

\bib{fulton-horn}{article}{
      author={Fulton, William},
       title={Eigenvalues, invariant factors, highest weights, and {S}chubert
  calculus},
        date={2000},
        ISSN={0273-0979},
     journal={Bull. Amer. Math. Soc. (N.S.)},
      volume={37},
      number={3},
       pages={209\ndash 249 (electronic)},
      review={\MR{2001g:15023}},
}

\bib{gantm}{book}{
      author={Gantmacher, F.~R.},
       title={The theory of matrices},
   publisher={Chelsea},
        date={1984},
      volume={2},
}

\bib{garloff-wagner}{article}{
      author={Garloff, J{\"u}rgen},
      author={Wagner, David~G.},
       title={Hadamard products of stable polynomials are stable},
        date={1996},
     journal={Journal of mathematical analysis and applications},
      volume={202},
       pages={797\ndash 809},
}

\bib{gasper}{book}{
      author={Gasper, George},
      author={Rahman, Mizan},
       title={Basic hypergeometric series},
      series={Encyclopedia of mathematics and its applications},
   publisher={Cambridge University press},
        date={1990},
      number={31},
}

\bib{godsil}{book}{
      author={Godsil, C.~D.},
       title={Algebraic combinatorics},
   publisher={Chapman and Hall},
        date={1993},
}

\bib{gomez}{article}{
      author={Gomez-Calderon, Javier},
      author={Wells, David~M.},
       title={Why polynomials have roots},
        date={1996},
     journal={College Mathematics Journal},
      volume={27},
      number={2},
       pages={90\ndash 94},
}

\bib{goodman-sun}{article}{
      author={Goodman, T. N.~T.},
      author={Sun, Qiyu},
       title={Total positivity and refinable functions with general dilation},
        date={2004},
        ISSN={1063-5203},
     journal={Appl. Comput. Harmon. Anal.},
      volume={16},
      number={2},
       pages={69\ndash 89},
      review={\MR{MR2038266 (2004k:42056)}},
}

\bib{grace}{article}{
      author={Grace, J.~H.},
       title={The zeros of a polynomial},
        date={1902},
     journal={Proc. cambridge Philos. Soc.},
      volume={11},
       pages={352\ndash 357},
}

\bib{grads}{book}{
      author={Gradshteyn, I.~S.},
      author={Ryzhjk, J.~M.},
       title={Table of integrals, series, and products},
   publisher={Academic Press},
        date={1994},
}

\bib{gregor96}{article}{
      author={Gregor, Ji\v{r}\'{i}},
      author={Ti\v{s}er, Jaroslav},
       title={On convex combinations of {H}urwitz polynomials},
        date={1996},
     journal={Appl. Math. and Comp. Sci},
      volume={6},
      number={1},
       pages={33\ndash 47},
}

\bib{gregor98}{article}{
      author={Gregor, Ji\v{r}\'{i}},
      author={Ti\v{s}er, Jaroslav},
       title={On {H}adamard powers of polynomials},
        date={1998},
     journal={Math. Control Signals Systems},
       pages={372\ndash 378},
}

\bib{gurvits}{misc}{
      author={Gurvits, Leonid},
       title={{Combinatorial and algorithmic aspects of hyperbolic
  polynomials}},
        note={arXiv:math.CO/0404474},
}

\bib{haglund97}{incollection}{
      author={Haglund, James},
      author={Ono, Ken},
      author={Wagner, David~G.},
       title={Theorems and conjectures involving rook polynomials with only
  real zeros},
        date={1999},
   booktitle={Topics in number theory (university park, pa, 1997)},
      series={Math. Appl.},
      volume={467},
   publisher={Kluwer Acad. Publ.},
     address={Dordrecht},
       pages={207\ndash 221},
      review={\MR{1691320 (2000d:05097)}},
}

\bib{hahn}{article}{
      author={Hahn, W},
       title={\"{U}ber orthogonalpolynome, die q-diffenenzengleichungen},
        date={1949},
     journal={Math. Nachr.},
      volume={2},
       pages={4\ndash 34},
}

\bib{hlpolya}{book}{
      author={Hardy, G.~H.},
      author={Littlewood, J.~E.},
      author={P{\'o}lya, G.},
       title={Inequalities},
      series={Cambridge Mathematical Library},
   publisher={Cambridge University Press},
     address={Cambridge},
        date={1988},
        ISBN={0-521-35880-9},
        note={Reprint of the 1952 edition},
      review={\MR{944909 (89d:26016)}},
}

\bib{vinnikov-helton}{misc}{
      author={Helton, J.~William},
      author={Vinnikov, Victor},
       title={{Linear Matrix Inequality Representation of Sets}},
        note={arXiv:math.OC/0306180},
}

\bib{hoggatt}{article}{
      author={Hoggatt, V.~E., Jr.},
      author={Bicknell, Marjorie},
       title={Roots of {F}ibonacci polynomials},
        date={1973},
        ISSN={0015-0517},
     journal={Fibonacci Quart.},
      volume={11},
      number={3},
       pages={271\ndash 274},
      review={\MR{0323700 (48 \#2056)}},
}

\bib{holtz}{article}{
      author={Holtz, Olga},
       title={Hermite-{B}iehler, {R}outh-{H}urwitz, and total positivity},
        date={2003},
        ISSN={0024-3795},
     journal={Linear Algebra Appl.},
      volume={372},
       pages={105\ndash 110},
      review={\MR{1 999 142}},
}

\bib{holtz-m}{article}{
      author={Holtz, Olga},
       title={M-matrices satisfy {N}ewton's inequalities},
        date={2004August},
     journal={Proceeding of the AMS},
}

\bib{hormander}{book}{
      author={H\"{o}rmander, Lars},
       title={An introduction to complex analysis in several variables},
   publisher={Van Nostrand},
        date={1966},
}

\bib{horn-problem}{article}{
      author={Horn, Alfred},
       title={Eigenvalues of sums of {H}ermitian matrices},
        date={1962},
     journal={Pacific J. Math.},
      volume={12},
       pages={225\ndash 241},
      review={\MR{25 \#3941}},
}

\bib{horn-johnson-1}{book}{
      author={Horn, Roger~A.},
      author={Johnson, Charles~R.},
       title={Matrix analysis},
   publisher={Cambridge University Press},
        date={1985},
}

\bib{horn-johnson-2}{book}{
      author={Horn, Roger~A.},
      author={Johnson, Charles~R.},
       title={Topics in matrix analysis},
   publisher={Cambridge University Press},
        date={1991},
}

\bib{horwitz}{article}{
      author={Horwitz, Alan},
       title={On the ratio vectors of polynomials},
        date={1997},
        ISSN={0022-247X},
     journal={J. Math. Anal. Appl.},
      volume={205},
      number={2},
       pages={568\ndash 576},
      review={\MR{1428370 (97j:12001)}},
}

\bib{smith-hu}{article}{
      author={Hu, Shu{-}an},
      author={Smith, Ronald~L.},
       title={The {S}chur complement interlacing theorem},
        date={1995},
     journal={Siam J. Matrix Anal. Appl.},
      volume={16},
      number={3},
       pages={1013\ndash 1023},
}

\bib{iliev}{book}{
      author={Iliev, Ljubomir},
       title={Laguerre entire functions},
   publisher={Bulgarian {A}cademy of {S}ciences},
        date={1987},
}

\bib{iserles}{article}{
      author={Iserles, A.},
      author={Norsett, S.~P.},
       title={Zeros of transformed polynomials},
        date={1990},
     journal={Math. Anal.},
      number={2},
       pages={483\ndash 509},
}

\bib{iserles-saff}{article}{
      author={Iserles, A.},
      author={Saff, E.~B.},
       title={Zeros of expansions in orthogonal polynomials},
        date={1989},
     journal={Math. Proc. Camb. Phil. Soc.},
      volume={105},
       pages={559\ndash 573},
}

\bib{ismail-muldoon}{article}{
      author={Ismail, Mourad},
      author={Muldoon, Martin},
       title={A discrete approach to monotonicity of zeros of orthogonal
  polynomials},
        date={1991},
     journal={Transactions of the AMS},
      volume={323},
       pages={65\ndash 78},
}

\bib{jahnke}{book}{
      author={Jahnke-Emde-L\u{o}sch},
       title={Tables of higher functions},
     edition={sixth edition},
   publisher={McGraw-Hill},
        date={1960},
}

\bib{jodar98}{article}{
      author={J\'{o}dar, L.},
      author={Defez, E.},
       title={A connection between {L}aguerre's and {H}ermite's matrix
  polynomials},
        date={1998},
     journal={Appl. Math. Lett.},
      volume={11},
      number={1},
       pages={13\ndash 17},
}

\bib{johnson}{article}{
      author={Johnson, Charles~R.},
       title={Interlacing polynomials},
        date={1987},
     journal={Proceedings of the American Mathematical Society},
      volume={100},
       pages={401\ndash 404},
}

\bib{johnson89}{article}{
      author={Johnson, C.R.},
       title={A characteristic polynomial for matrix pairs},
        date={1989},
     journal={Linear and multilinear algebra},
      volume={25},
       pages={289\ndash 290},
}

\bib{wjohnson}{article}{
      author={Johnson, Warren~P.},
       title={The curious history of {F}a\`a di {B}runo's formula},
        date={2002},
        ISSN={0002-9890},
     journal={Amer. Math. Monthly},
      volume={109},
      number={3},
       pages={217\ndash 234},
      review={\MR{1903577 (2003d:01019)}},
}

\bib{kac}{book}{
      author={Kac, Victor},
      author={Cheung, Pokman},
       title={Quantum calculus},
      series={Universitext},
   publisher={Springer-Verlag},
     address={New York},
        date={2002},
        ISBN={0-387-95341-8},
      review={\MR{1865777 (2003i:39001)}},
}

\bib{karlin-szego}{article}{
      author={Karlin, Samuel},
      author={Szego, Gabor},
       title={On certain determinants whose entries are orthogonal
  polynomials},
        date={11960/61},
     journal={J. Analyse Math.},
      volume={8},
       pages={1\ndash 157},
}

\bib{karlin}{book}{
      author={Karlin, Samuel},
       title={Total positivity},
   publisher={Stanford University Press},
        date={1968},
}

\bib{kim-kwon}{article}{
      author={Kim, Deok},
      author={Kwon, Kil},
       title={{{On a conjecture of Karlin and Szego}}},
        date={1996},
     journal={Proceedings of the American Mathematical Society},
      volume={124},
       pages={227\ndash 231},
}

\bib{klyatcho}{article}{
      author={Klyachko, Alexander~A.},
       title={Random walks on symmetric spaces and inequalities for matrix
  spectra},
        date={2000},
        ISSN={0024-3795},
     journal={Linear Algebra Appl.},
      volume={319},
      number={1-3},
       pages={37\ndash 59},
      review={\MR{2002f:15023}},
}

\bib{koekoek}{techreport}{
      author={Koekoek, Roelof},
      author={Swarttouw, Ren\`{e}~F.},
       title={The {A}skey-scheme of hypergeometric orthogonal polynomials and
  its q-analog},
 institution={Technische Universiteit Delft},
        date={1994},
      number={94-05},
}

\bib{kostov-shapiro}{unpublished}{
      author={Kostov, V.},
      author={Shapiro, B.},
       title={On arrangements of roots for a real hyperbolic polynomial and its
  derivatives},
        note={to appear},
}

\bib{wolowski-diff}{article}{
      author={Kozitsky, Yu.},
      author={Oleszczuk, P.},
      author={Wolowski, L.},
       title={{Infinite Order Differential Operators in Spaces of Entire
  Functions}},
        date={2003},
     journal={J. Math. Anal. Appl.},
      volume={277},
       pages={423\ndash 437},
        note={arXiv:math.FA/0311213},
}

\bib{wolowski-convex}{misc}{
      author={Kozitsky, Yuri},
      author={Wolowski, Lech~B.},
       title={{Laguerre Entire Functions and Related Locally Convex Spaces}},
        note={arXiv:math.CV/9812111},
}

\bib{krasikov}{misc}{
      author={Krasikov, Ilia},
       title={On zeros of polynomials and allied functions satisfying second
  order differential equations},
        note={arXiv:math.CA/0211190},
}

\bib{kratt}{article}{
      author={Krattenthaler, C.},
       title={Advanced determinant calculus},
        date={1999},
     journal={Sém. Lothar. Combin},
      volume={42},
}

\bib{kuperberg}{misc}{
      author={Kuperberg, Greg},
       title={{Special moments}},
        note={arXiv:math.PR/0408360},
}

\bib{kurtz}{article}{
      author={Kurtz, David~C.},
       title={A sufficient condition for all the roots of a polynomial to be
  real},
        date={1992},
     journal={Amer. Math. Monthly},
       pages={259\ndash 263},
}

\bib{lancaster}{article}{
      author={Lancaster, Peter},
      author={Psarrakos, Panayiotis},
       title={The numerical range of self-adjoint quadratic matrix
  polynomials},
        date={2001/02},
        ISSN={0895-4798},
     journal={SIAM J. Matrix Anal. Appl.},
      volume={23},
      number={3},
       pages={615\ndash 631 (electronic)},
      review={\MR{1896809 (2003e:15032)}},
}

\bib{lax-conj}{misc}{
      author={Lewis, Adrian~S.},
      author={Parrilo, Pablo~A.},
      author={Ramana, Motakuri~V.},
       title={{The Lax conjecture is true}},
        note={arXiv:math.OC/0304104},
}

\bib{linear-preserver}{article}{
      author={Li, Chi{K}wong},
      author={Pierce, Stephen},
       title={Linear preserver problems},
        date={2001},
     journal={Amer. Math. Monthly},
      volume={108},
      number={7},
       pages={591\ndash 605},
}

\bib{lieb-seiringer}{article}{
      author={Lieb, Elliott~H.},
      author={Seiringer, Robert},
       title={Equivalent forms of the bessis-moussa-villani conjecture},
        date={2004},
     journal={Journal of Statistical Physics},
      volume={115},
      number={1-2},
       pages={185\ndash 190},
}

\bib{lieb-sokal}{article}{
      author={Lieb, Elliott~H.},
      author={Sokal, Alan~D.},
       title={A general {L}ee-{Y}ang theorem for one-component and
  multicomponent ferromagnets},
        date={1981},
        ISSN={0010-3616},
     journal={Comm. Math. Phys.},
      volume={80},
      number={2},
       pages={153\ndash 179},
      review={\MR{623156 (83c:82008)}},
}

\bib{little}{article}{
      author={Little, John~B.},
       title={On the zeros of two families of polynomials arising from certain
  rational integrals},
        date={2005},
        ISSN={0035-7596},
     journal={Rocky Mountain J. Math.},
      volume={35},
      number={4},
       pages={1205\ndash 1215},
      review={\MR{2178985}},
}

\bib{liu-wang}{misc}{
      author={Liu, Li},
      author={Wang, Yi},
       title={{A unified approach to polynomial sequences with only real
  zeros}},
        note={arXiv:math.CO/0509207},
}

\bib{liu-wang-log}{misc}{
      author={Liu, Li},
      author={Wang, Yi},
       title={{On the log-convexity of combinatorial sequences}},
        note={arXiv:math.CO/0602672},
}

\bib{locher}{article}{
      author={Locher, Franz},
       title={Stability tests for linear difference forms},
        date={1993},
     journal={International series of numerical mathematics},
      volume={112},
       pages={215\ndash 223},
}

\bib{macduffee}{book}{
      author={MacDuffee, C.C.},
       title={Theory of equations},
   publisher={Wiley},
        date={1954},
}

\bib{marden}{book}{
      author={Marden, Morris},
       title={The geometry of the zeros of a polynomial in a complex variable},
   publisher={American Mathematical Society},
        date={1949},
}

\bib{markus}{book}{
      author={Markus, A.S.},
       title={Introduction to the spectral theory of polynomial operator
  pencils},
      series={Translations of Mathematical Monographs},
   publisher={Americam Mathematical Society},
        date={1988},
      volume={71},
}

\bib{marquez}{article}{
      author={Marquez, Horacio~J.},
      author={Agathoklis, Panajotis},
       title={On the existence of robust strictly positive real rational
  functions},
        date={1998},
        ISSN={1057-7122},
     journal={IEEE Trans. Circuits Systems I Fund. Theory Appl.},
      volume={45},
      number={9},
       pages={962\ndash 967},
      review={\MR{1647619 (99h:93039)}},
}

\bib{masson}{article}{
      author={M{\'a}sson, Gisli},
      author={Shapiro, Boris},
       title={On polynomial eigenfunctions of a hypergeometric-type operator},
        date={2001},
        ISSN={1058-6458},
     journal={Experiment. Math.},
      volume={10},
      number={4},
       pages={609\ndash 618},
      review={\MR{2003e:33038}},
}

\bib{stanton-white}{misc}{
      author={de~M\'edicis, Anne},
      author={Stanton, Dennis~W.},
      author={White, Dennis~E.},
       title={{The combinatorics of q-Charlier polynomials}},
        note={arXiv:math.CA/9307208},
}

\bib{meijer}{article}{
      author={Meijer, H.~G.},
       title={{{Coherent pairs and zeros of Sobolev-type orthogonal
  polynomials}}},
        date={1993},
     journal={Indag. Mathem.},
      volume={2},
       pages={163\ndash 176},
}

\bib{mitev}{article}{
      author={Mitev, Todor~P.},
       title={New inequalities between elementary symmetric polynomials},
        date={2003},
        ISSN={1443-5756},
     journal={JIPAM. J. Inequal. Pure Appl. Math.},
      volume={4},
      number={2},
       pages={Article 48, 11 pp. (electronic)},
      review={\MR{1994261 (2004d:26015)}},
}

\bib{muldoon}{article}{
      author={Muldoon, M.~E.},
       title={Properties of zeros of orthogonal polynomials and related
  functions},
        date={1993},
     journal={J. Compu. Appl. Math.},
      volume={48},
       pages={167\ndash 186},
}

\bib{OEIS}{misc}{
      author={N.~J. A.~Sloane, editor~(2003)},
       title={The on-line encyclopedia of integer sequences},
        note={http://www.research.att.com/~njas/sequences/},
}

\bib{niculesscu}{article}{
      author={Niculescu, Constantin~P.},
       title={A new look at {N}ewton's inequalities},
        date={2000},
        ISSN={1443-5756},
     journal={JIPAM. J. Inequal. Pure Appl. Math.},
      volume={1},
      number={2},
       pages={Article 17, 14 pp. (electronic)},
      review={\MR{MR1786404 (2001h:26020)}},
}

\bib{nijenhuis}{article}{
      author={Nijenhuis, Albert},
       title={On permanants and the zeros of rook polynomials},
        date={1976},
     journal={J. Comb. Theory (A)},
      volume={21},
       pages={240\ndash 244},
}

\bib{nunemacher}{article}{
      author={Nunemacher, Jeffrey},
       title={Asymptotes, cubic curves, and the projective plane},
        date={1999},
     journal={Math. Magazine},
      volume={72},
      number={3},
       pages={183\ndash 192},
}

\bib{rao}{article}{
      author={P.~Rao, P.S.S.N.V},
       title={A characterization of merely subdefinite matrices and related
  results},
        date={1977},
     journal={Sankhya Ser. A},
      volume={39},
      number={4},
       pages={387\ndash 395},
}

\bib{parlett}{book}{
      author={Parlett, Beresford},
       title={The symmetric eigenvalue problem},
   publisher={Prentice-Hall},
        date={1980},
}

\bib{patrick}{article}{
      author={Patrick, Merrell~L.},
       title={Extensions of inequalities of the {L}aguerre and {T}ur\'an type},
        date={1973},
     journal={Pacific J. Math.},
      volume={44},
       pages={675\ndash 682},
      review={\MR{47 \#3725}},
}

\bib{peyser}{article}{
      author={Peyser, Gideon},
       title={On the roots of the derivative of a polynomial with real roots},
        date={1967},
        ISSN={0002-9890},
     journal={Amer. Math. Monthly},
      volume={74},
       pages={1102\ndash 1104},
      review={\MR{0227332 (37 \#2917)}},
}

\bib{polya-schur}{article}{
      author={P\'{o}lya, G.},
      author={Schur, J.},
       title={Uber zwei arten von factorenfolgen in der theorie der
  algebraischen gleichungen},
        date={1914},
     journal={Journal fur die reine und angewandte mathematik},
      volume={band 144},
       pages={89\ndash 113},
}

\bib{polya-szego2}{book}{
      author={P\'{o}lya, G.},
      author={Szego, G.},
       title={Problems and theorems in analysis},
   publisher={Springer-Verlag},
        date={1972},
      volume={2},
}

\bib{polya-szego1}{book}{
      author={P\'{o}lya, G.},
      author={Szego, G.},
       title={Problems and theorems in analysis},
   publisher={Springer-Verlag},
        date={1972},
      volume={1},
}

\bib{polya26}{article}{
      author={P\'{o}lya, G.},
       title={Bemerkung uber die integraldarstellung der {R}iemannssche
  $\zeta$-function},
        date={1926},
     journal={Acta Math.},
      volume={48},
       pages={305\ndash 317},
}

\bib{rahman}{book}{
      author={Rahman, Q.~I.},
      author={Schmeisser, G.},
       title={Analytic theory of polynomials},
      series={Mathematical Society Monographs},
   publisher={Oxford Science Publications},
        date={2002},
      number={26},
}

\bib{redmond}{unpublished}{
      author={Redmond, Timothy},
       title={A class of polynomial {R}iemann hypotheses},
        note={in prepartion},
}

\bib{reznick-book}{book}{
      author={Reznick, Bruce},
       title={Sums of even powers of real linear forms},
      series={Memoirs},
   publisher={American mathematical society},
        date={1992},
      number={463},
}

\bib{reznick}{article}{
      author={Reznick, Bruce},
       title={An inequality for products of polynomials},
        date={1993},
     journal={Proceedings of the American mathematical society},
      volume={117},
       pages={1063\ndash 1073},
}

\bib{rivlin}{book}{
      author={Rivlin, J.},
       title={The {C}hebychev polynomials},
   publisher={Wiley},
        date={1974},
}

\bib{roman}{book}{
      author={Roman, Steve},
       title={The umbral calculus},
   publisher={Academic Press},
        date={1984},
}

\bib{rosset}{article}{
      author={Rosset, Shmuel},
       title={Normalized symmetric functions, {N}ewton's inequalities and a new
  set of stronger inequalities},
        date={1989},
        ISSN={0002-9890},
     journal={Amer. Math. Monthly},
      volume={96},
      number={9},
       pages={815\ndash 819},
      review={\MR{1033349 (91e:05083)}},
}

\bib{rota}{book}{
      author={Rota, Gian-Carlo},
       title={Finite operator calculus},
   publisher={Academic Press},
        date={1975},
}

\bib{rubel}{book}{
      author={Rubel, Lee},
       title={Entire and meromorphic functions},
   publisher={Springer-Verlag},
        date={1996},
}

\bib{ruelle}{incollection}{
      author={Ruelle, David},
       title={Grace-like polynomials},
        date={2002},
   booktitle={Foundations of computational mathematics (hong kong, 2000)},
   publisher={World Sci. Publ., River Edge, NJ},
       pages={405\ndash 421},
      review={\MR{2021990 (2004k:30014)}},
}

\bib{shah}{article}{
      author={Shah, G.~M.},
       title={On the zeros of {V}an {V}leck polynomials},
        date={1968},
     journal={Proc. Amer. Math. Soc.},
      volume={19},
       pages={1421\ndash 1426},
      review={\MR{37 \#6511}},
}

\bib{bezout}{misc}{
      author={Shapiro, A.},
      author={Kaplan, S.},
      author={Teicher, M.},
       title={{Several Applications of Bezout Matrices}},
        note={arXiv:math.AG/0601047},
}

\bib{shapiro-01}{article}{
      author={Shapiro, B.},
      author={M\'asson, G.},
       title={On polynomial eigenfunctions of $\frac{d^k}{dx^k} q(x)$},
        date={2001},
     journal={IHES-preprint},
}

\bib{shapiro-welker}{article}{
      author={Shapiro, Boris},
      author={Welker, Volkmar},
       title={Combinatorics and topology of stratifications of the space of
  monic polynomials with real coefficients},
        date={1998},
     journal={Result. Math},
      volume={33},
       pages={338\ndash 355},
}

\bib{sheil-small}{book}{
      author={Sheil-Small, T.},
       title={Complex polynomials},
      series={Cambridge Studies in Advanced Mathematics},
   publisher={Cambridge University Press},
     address={Cambridge},
        date={2002},
      volume={75},
        ISBN={0-521-40068-6},
      review={\MR{1 962 935}},
}

\bib{shukla}{article}{
      author={Shukla, D.~P.},
       title={q-angelescu polynomials},
        date={1981},
     journal={Publications de l'{I}nstitut {M}ath/'{e}matique},
      volume={29},
      number={43},
       pages={205\ndash 213},
}

\bib{sidi}{article}{
      author={Sidi, Avram},
       title={Numerical quadrature and non-linear sequence transformations;
  unified rules for efficient computation of integrals with algebraic and
  logarithmic endpoint singularities},
        date={1989},
     journal={Math. Comp.},
      volume={35},
       pages={851\ndash 874},
}

\bib{simion}{article}{
      author={Simion, Rodica},
       title={A multiindexed {Strum} sequence of polynomials and unimodality of
  certain combinatorial sequences},
        date={1984},
     journal={J. Comb. Theory(A)},
      volume={36},
       pages={15\ndash 22},
}

\bib{singer}{article}{
      author={Singer, D.},
       title={A bijective proof of {B}orchardt's identity},
        date={2004},
     journal={Electon. J. Combin.},
      volume={11},
      number={R48},
}

\bib{smith}{article}{
      author={Smith, Ronald~L.},
       title={Some interlacing properties of the {S}chur complement of a
  {H}ermitian matrix},
        date={1992},
     journal={Linear algebra and its applications},
      volume={177},
       pages={137\ndash 144},
}

\bib{speyer}{unpublished}{
      author={Speyer, David~E.},
       title={Horn's problem, {V}innikov curves, and the hive cone},
        note={preprint},
}

\bib{stanley}{article}{
      author={Stanley, Richard},
       title={Log-concave and unimodal sequences in algebra, combinatorics, and
  geometry},
        date={1989},
     journal={Annals New York Academy of Sciences},
      volume={576},
       pages={500\ndash 535},
}

\bib{stanley97}{book}{
      author={Stanley, Richard},
       title={Enumerative combinatorics},
   publisher={Cambridge University Press},
        date={1997},
      volume={II},
}

\bib{steele}{book}{
      author={Steele, J.~Michael},
       title={The {C}auchy-{S}chwarz master class},
      series={MAA Problem Books Series},
   publisher={Mathematical Association of America},
     address={Washington, DC},
        date={2004},
        ISBN={0-521-83775-8; 0-521-54677-X},
        note={An introduction to the art of mathematical inequalities},
      review={\MR{2062704 (2005a:26035)}},
}

\bib{stoyanoff}{article}{
      author={Stoyanoff, A.},
       title={Sur un theorem de {M} {M}arcel {R}iesz},
        date={1926},
     journal={Nouvelles Annales de Mathematique},
      volume={1},
       pages={97\ndash 99},
}

\bib{suetin}{book}{
      author={Suetin, P.~K.},
       title={Orthogonal polynomials in two variables (analytical methods and
  special functions)},
   publisher={Taylor and Francis},
        date={1999},
}

\bib{szasz}{article}{
      author={Sz{\'a}sz, Otto},
       title={On sequences of polynomials and the distribution of their zeros},
        date={1943},
     journal={Bull. Amer. Math. Soc.},
      volume={49},
       pages={377\ndash 383},
      review={\MR{0008274 (4,273c)}},
}

\bib{szego}{book}{
      author={Szego, Gabor},
       title={Orthogonal polynomials},
   publisher={American Mathematical Society},
        date={1939},
}

\bib{szego-turan}{article}{
      author={Szego, Gabor},
       title={On an inequality of {P}. {T}uran concerning {L}egendre
  polynomials},
        date={1948},
     journal={Bull. Amer. Math. Soc},
       pages={401\ndash 5},
}

\bib{vanassche}{article}{
      author={Van~Assche, Walter},
      author={Coussement, Els},
       title={Some classical multiple orthogonal polynomials},
        date={2001},
        ISSN={0377-0427},
     journal={J. Comput. Appl. Math.},
      volume={127},
      number={1-2},
       pages={317\ndash 347},
        note={Numerical analysis 2000, Vol. V, Quadrature and orthogonal
  polynomials},
      review={\MR{2001i:33012}},
}

\bib{vassiliev89}{article}{
      author={Vassiliev, V.~A.},
       title={Topology of spaces of functions without compound singularities},
        date={1989},
     journal={Funct. Anal. Appl.},
      volume={23},
       pages={277\ndash 286},
}

\bib{vassiliev92}{book}{
      author={Vassiliev, V.~A.},
       title={Complements of discriminants of smooth maps: topology and
  applications},
      series={Translations of Math. Mono.},
   publisher={A.M.S.},
        date={1992},
      volume={98},
}

\bib{vegter}{article}{
      author={Vegter, Gert},
       title={The apolar bilinear form in geometric modeling},
        date={2000},
        ISSN={0025-5718},
     journal={Math. Comp.},
      volume={69},
      number={230},
       pages={691\ndash 720},
      review={\MR{1654006 (2000i:65018)}},
}

\bib{vinnikov}{article}{
      author={Vinnikov, Victor},
       title={Selfadjoint determinantal representations of real plane curves},
        date={1993},
        ISSN={0025-5831},
     journal={Math. Ann.},
      volume={296},
      number={3},
       pages={453\ndash 479},
      review={\MR{94e:14038}},
}

\bib{wagner-newton}{misc}{
      author={Wagner, David~G.},
       title={{Enumeration of spanning subgraphs with degree constraints}},
        note={arXiv:math.CO/0412059},
}

\bib{wagner-genus-polys}{unpublished}{
      author={Wagner, David~G.},
       title={Zeros of genus polynomials of graphs in some linear families},
        note={preprint},
}

\bib{wagner-zeros}{misc}{
      author={Wagner, David~G.},
       title={Zeros of reliability polynomials and $f$-vectors of matroids},
        note={arXiv:math.CO/9802047},
}

\bib{wagner}{article}{
      author={Wagner, David~G.},
       title={Total positivity of {Hadamard} products},
        date={1992},
     journal={J. Math. Anal. Appl.},
      volume={163},
       pages={459\ndash 483},
}

\bib{wagner96}{article}{
      author={Wagner, David~G.},
       title={Multipartition series},
        date={1996},
     journal={Siam J. Discrete Math.},
      volume={9},
      number={4},
       pages={529\ndash 544},
}

\bib{walker97}{article}{
      author={Walker, Peter},
       title={Separation of zeros of translates of polynomials and entire
  functions},
        date={1997},
        ISSN={0022-247X},
     journal={J. Math. Anal. Appl.},
      volume={206},
      number={1},
       pages={270\ndash 279},
      review={\MR{97m:30004}},
}

\bib{walker}{article}{
      author={Walker, Peter},
       title={Separation of zeros of translates of polynomials and entire
  functions},
        date={1997},
     journal={Journal of Mathematical Analysis and Applications},
      volume={206},
       pages={270\ndash 279},
}

\bib{wang-wang-yu}{misc}{
      author={Wang, Long},
      author={Wang, Zhizhen},
      author={Yu, Wensheng},
       title={Stability of polytopic polynomial matrices},
        date={2002},
        note={arXiv:math/0202238},
}

\bib{sampling-zeros}{article}{
      author={Weller, Steven},
      author={Moran, W.},
      author={Ninness, Brett},
      author={Pollington, A.D.},
       title={Sampling zeros and the {E}uler-{F}robenius polynomials},
        date={2001},
     journal={IEEE transactions on automatic control},
      volume={46},
      number={2},
       pages={340\ndash 343},
}

\bib{yamaguchi}{article}{
      author={Yamaguchi, Kohhei},
       title={Spaces of polynomials with real roots of bounded multiplicity},
        date={2002},
     journal={J. Math. Kyoto Univ.},
      volume={42},
      number={1},
       pages={105\ndash 115},
}

\bib{yu-wang-long-ack}{article}{
      author={Yu, Wensheng},
      author={Wang, Long},
      author={Ackermann, J{\"u}rgen},
       title={Robust strictly positive real synthesis for polynomial families
  of arbitrary order},
        date={2004},
        ISSN={1009-2757},
     journal={Sci. China Ser. F},
      volume={47},
      number={4},
       pages={475\ndash 489},
      review={\MR{2126240 (2006a:93044)}},
}

\bib{zabczyk}{book}{
      author={Zabczyk, Jerzy},
       title={Mathematical control theory: an introduction},
      series={Systems \& Control: Foundations \& Applications},
   publisher={Birkh\"auser Boston Inc.},
     address={Boston, MA},
        date={1992},
        ISBN={0-8176-3645-5},
      review={\MR{1193920 (93h:49001)}},
}

\bib{zemyan}{article}{
      author={Zemyan, Stephen~M.},
       title={On the zeroes of the {$N$}th partial sum of the exponential
  series},
        date={2005},
        ISSN={0002-9890},
     journal={Amer. Math. Monthly},
      volume={112},
      number={10},
       pages={891\ndash 909},
      review={\MR{2186832}},
}

\bib{zhedanov}{article}{
      author={Zhedanov, A.S.},
       title={Gauss sums and orthogonal polynomials},
        date={1997},
     journal={Intl. J. of Modern Physics A},
      volume={12},
      number={1},
       pages={289\ndash 294},
}

\end{biblist}
\end{bibdiv}

\printindex

\end{document}